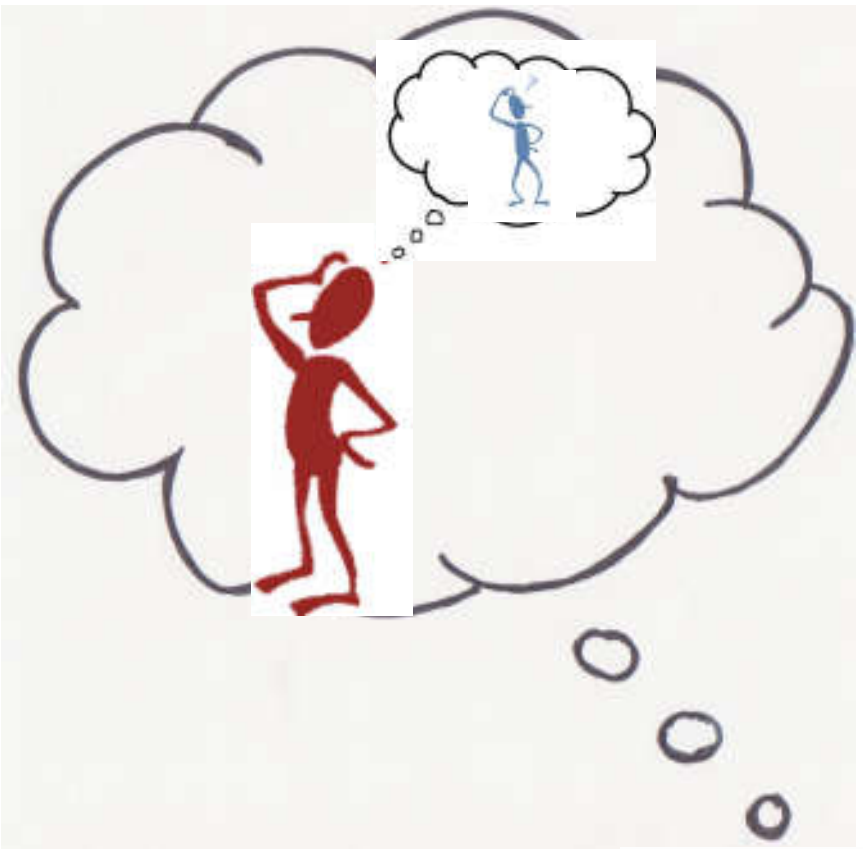

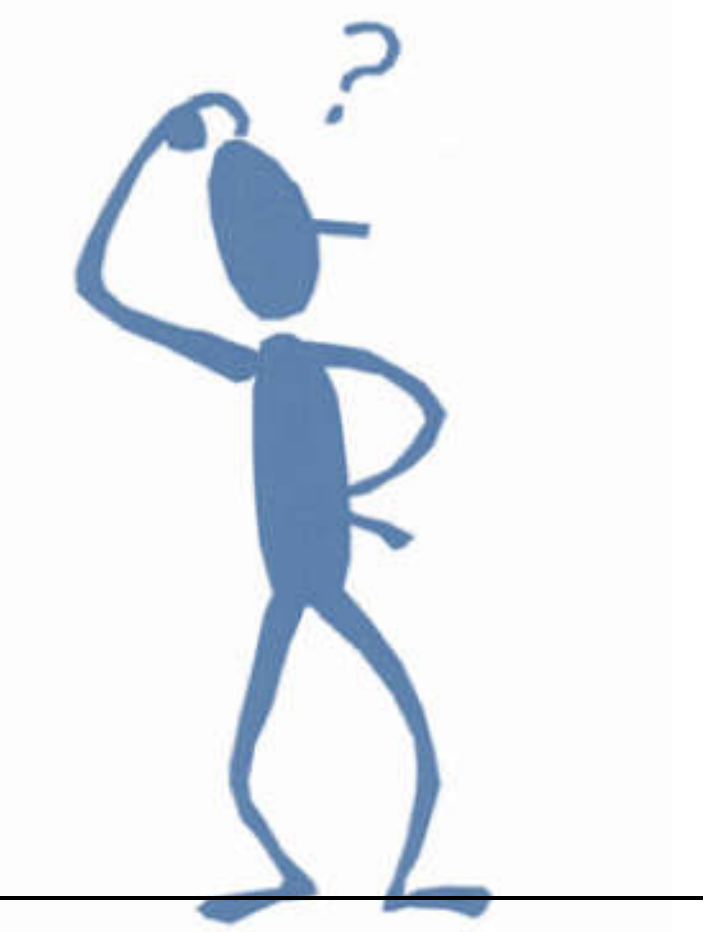

# GAME THEORY

## An open access textbook with 165 solved exercises

Giacomo Bonanno

University of California, Davis

http://www.econ.ucdavis.edu/faculty/bonanno/

# CONTENTS

## PART V: Advanced Topics III: Incomplete Information

# Preface and Acknowledgments

After teaching game theory (at both the undergraduate and graduate level) at the University of California, Davis for 25 years, I decided to organize all my teaching material in a textbook. There are many excellent textbooks in game theory and there is hardly any need for a new one. However, there are two distinguishing features of this textbook: (1) it is open access and thus free,[1] and (2) it contains an unusually large number of exercises (a total of 165) with complete and detailed answers.

I tried to write the book in such a way that it would be accessible to anybody with minimum knowledge of mathematics (high-school level algebra and some elementary notions of probability) and no prior knowledge of game theory. However, the book is intended to be rigorous and it includes several proofs. I believe it is appropriate for an advanced undergraduate class in game theory and also for a first-year graduate-level class.

I expect that there will be some typos and (hopefully minor) mistakes. If you come across any typos or mistakes, I would be grateful if you could inform me: I can be reached at gfbonanno@ucdavis.edu. I will maintain an updated version of the book on my web page at

http://www.econ.ucdavis.edu/faculty/bonanno/

I also intend to add, some time in the future, a further collection of exercises and exam questions with a solution manual to be given only to instructors. Details will appear on my web page.

I am very grateful to Elise Tidrick for meticulously going through each chapter of the book and for suggesting numerous improvements. Her insightful and constructive comments have considerably enhanced this book.

I would also like to thank Nicholas Bowden, Lester Lusher, Burkhard Schipper and Matthieu Stigler for pointing out typos.

Davis, October 27, 2015

[1] There may be several other free textbooks on game theory available. The only one I am aware of is the excellent book by Ariel Rubinstein and Martin Osborne, MIT Press, 1994, which can be downloaded for free from Ariel Rubinstein's web page: http://arielrubinstein.tau.ac.il. In my experience this book is too advanced for an undergraduate class in game theory.

# Chapter 0

# Introduction

The discipline of game theory was pioneered in the early 20th century by mathematicians Ernst Zermelo (1913) and John von Neumann (1928). The breakthrough came with John von Neumann and Oscar Morgenstern's book, *Theory of games and economic behavior*, published in 1944. This was followed by important work by John Nash (1950-51) and Lloyd Shapley (1953). Game theory had a major influence on the development of several branches of economics (industrial organization, international trade, labor economics, macroeconomics, etc.). Over time the impact of game theory extended to other branches of the social sciences (political science, international relations, philosophy, sociology, anthropology, etc.) as well as to fields outside the social sciences, such as biology, computer science, logic, etc. In 1994 the Nobel prize in economics was given to three game theorists, John Nash, John Harsanyi and Reinhardt Selten, for their theoretical work in game theory which was very influential in economics. At the same time, the US Federal Communications Commission was using game theory to help it design a $7-billion auction of the radio spectrum for personal communication services (naturally, the bidders used game theory too!). The Nobel prize in economics was awarded to game theorists three more times: in 2006 to Robert Aumann and Thomas Schelling, in 2007 to Leonid Hurwicz, Eric Maskin and Roger Myerson and in 2010 to Lloyd Shapley and Alvin Roth.

Game theory provides a formal language for the representation and analysis of *interactive situations*, that is, situations where several "entities", called *players*, take actions that affect each other. The nature of the players varies depending on the context in which the game theoretic language is invoked: in evolutionary biology (see, for example, John Maynard Smith, 1982) players are non-thinking living organisms;[2] in computer science (see, for example, Shoham-Leyton-Brown, 2008) players are artificial agents; in behavioral game theory (see, for example, Camerer, 2003) players are "ordinary" human beings, etc. Traditionally, however, game theory has focused on interaction among intelligent, sophisticated and rational individuals. For example, Aumann describes game theory as follows:

> "Briefly put, game and economic theory are concerned with the interactive behavior of *Homo rationalis* - rational man. *Homo rationalis* is the species that always acts both purposefully and logically, has well-defined goals, is motivated solely by the desire to approach these goals as closely as possible, and has the calculating ability required to do so." (Aumann, 1985, p. 35.)

This book is concerned with the traditional interpretation of game theory.

Game theory is divided into two main branches. The first is *cooperative* game theory, which assumes that the players can communicate, form coalitions and sign binding agreements. Cooperative game theory has been used, for example, to analyze voting behavior and other issues in political science and related fields. We will deal exclusively with the other main branch, namely *non-cooperative* game theory. Non-cooperative game theory models situations where the players are either unable to communicate or are able to communicate but cannot sign binding contracts. An example of the latter situation is the interaction among firms in an industry in an environment where antitrust laws make it illegal for firms to reach agreements concerning prices or production quotas or other forms of collusive behavior.

The book is divided into three parts.

**Part I** deals with games with *ordinal* payoffs, that is, with games where the players' preferences over the possible outcomes are only specified in terms of an ordinal ranking (outcome $o$ is better than outcome $o'$ or $o$ is just as good as $o'$). Chapter 1 covers strategic-form games, Chapter 2 deals with dynamic games with perfect information and Chapter 3 with dynamic games with (possibly) imperfect information.

[2] Evolutionary game theory has been applied not only to the analysis of animal and insect behavior but also to studying the " most successful strategies" for tumor and cancer cells (see, for example, Gerstung *et al.*, 2011).

**Part II** is devoted to games with *cardinal* payoffs, that is, with games where the players' preferences extend to uncertain prospects or lotteries: players are assumed to have a consistent ranking of the set of lotteries over basic outcomes. Chapter 4 reviews the theory of expected utility, Chapter 5 discusses the notion of mixed strategy in strategic-form games and of mixed-strategy Nash equilibrium , while Chapter 6 deals with mixed strategies in dynamic games.

Parts III, IV and V cover a number of advanced topics.

**Part III** deals with the notions of knowledge, common knowledge and belief. Chapter 7 explains how to model what an individual knows and what she is uncertain about and how to extend the analysis to the interactive knowledge of several individuals (e.g. what Individual 1 knows about what Individual 2 knows about some facts or about the state of knowledge of Individual 1). The chapter ends with the notion of common knowledge. Chapter 8 adds probabilistic beliefs to the knowledge structures of the previous chapter and discusses the notions of Bayesian updating, belief revision, like-mindedness and the possibility of "agreeing to disagree". Chapter 9 uses the interactive knowledge-belief structures of the previous two chapters to model the players' state of mind in a possible play of a given game and studies the implications of common knowledge of rationality in strategic-form games.

**Part IV** focuses on dynamic (or extensive-form) games and on the issue of how to refine the notion of subgame-perfect equilibrium (which was introduced in Chapters 3 and 6). Chapter 10 introduces a simple notion, called weak sequential equilibrium, which achieves some desirable goals (such as the elimination of strictly dominated choices) but fails to provide a refinement of subgame-perfect equilibrium. Chapter 11 explains the more complex notion of sequential equilibrium, which is extensively used in applications of game theory. That notion, however, leaves to be desired from a practical point of view (it is typically hard to show that an equilibrium is indeed a sequential equilibrium) and also from a conceptual point of view (it appeals to a topological condition, whose interpretation is not clear). Chapter 12 introduces an intermediate notion, called perfect Bayesian equilibrium, whose conceptual justification is anchored in the so called AGM theory of belief revision, extensively studied in philosophy and computer science, which was pioneered by Carlos Alchourrón (a legal scholar), Peter Gärdenfors (a philosopher) and David Makinson (a computer scientist) in 1985. In Chapter 12 we also provide an alternative characterization of sequential equilibrium based on the notion of perfect Bayesian equilibrium, which is free of topological conditions.

**Part V** deals with the so-called "theory of games of incomplete information", which was pioneered by John Harsanyi (1967-68). This theory is usually explained using the so-called "type-space" approach suggested by Harsanyi. However, we follow a different approach: the so-called "state-space" approach, which makes use of the interactive knowledge-belief structures developed in Part III. We find this approach both simpler and more elegant. For completeness, in Chapter 15 we explain the commonly used type-based structures and show how to convert a state-space structure into a type-space structure and *vice versa* (the two approaches are equivalent). Chapter 13 deals with situations of incomplete information that involve static (or strategic-form) games, while Chapter 14 deals with situations of incomplete information that involve dynamic (or strategic-form) games.

Appendix E ('E' for 'exercises') at the end of each chapter contains several exercises, which are divided according the that chapter's sections. Indeed, at the end of each chapter section the reader is invited to go to the appendix and try the exercises for that section. Complete and detailed answers for each exercise can be found in Appendix S ('S' for 'solutions') at the end of each chapter. For each chapter, the set of exercises culminates into a "challenging question", which is more difficult and more time consuming than the other exercises. In game theory, as in mathematics in general, it is essential to test one's understanding of the material by attempting to solve exercises and problems. The reader is encouraged to attempt solving exercises after the introduction of every new concept.

# PART I

# Games with Ordinal Payoffs

Chapter 1

# Ordinal games in strategic form

## 1.1. Game frames and games

Game theory deals with interactive situations where two or more individuals, called *players*, make decisions that jointly determine the final outcome. To see an example point your browser to the following video

https://www.youtube.com/watch?v=tBtr8-VMj0E

(if you search for 'Split or Steal' on youtube.com you will find several instances of this game). The video shows an excerpt from *Golden Balls*, a British daytime TV game show. In it each of two players, Sarah and Steve, has to pick one of two balls: inside one ball appears the word 'split' and inside the other the word 'steal' (each player is first asked to secretly check which of the two balls in front of him/her is the split ball and which is the steal ball). They make their decisions simultaneously. The possible outcomes are shown in Table 1.1, where each row is labeled with a possible choice for Sarah and each column with a possible choice for Steven. Each cell in the table thus corresponds to a possible pair of choices and the resulting outcome is written inside the cell.

| | | Steven: Split | | Steven: Steal | |
|---|---|---|---|---|---|
| Sarah | Split | Sarah gets $50,000 | Steven gets $50,000 | Sarah gets nothing | Steven gets $100,000 |
| Sarah | Steal | Sarah gets $100,000 | Steven gets nothing | Sarah gets nothing | Steven gets nothing |

Table 1.1

What should a rational player do in such a situation? It is tempting to reason as follows.

> Let us focus on Sarah's decision problem. She realizes that her decision alone is not sufficient to determine the outcome; she has no control over what Steven will choose to do. However, she can envision two scenarios: one where Steven chooses *Steal* and the other where he chooses *Split*.
>
> - If Steven decides to *Steal*, then it does not matter what Sarah does, because she ends up with nothing, no matter what she chooses.
> - If Steven picks *Split*, then Sarah will get either \$50,000 (if she also picks *Split*) or \$100,000 (if she picks *Steal*).
>
> Thus Sarah should choose *Steal*.

The trouble with the above argument is that it is not a valid argument because it makes an implicit assumption about how Sarah ranks the outcomes, namely that she is *selfish and greedy*, which may or may not be correct. Let us denote the outcomes as follows:

$o_1$ : Sarah gets \$50,000 and Steven gets \$50,000.
$o_2$ : Sarah gets nothing and Steven gets \$100,000.
$o_3$ : Sarah gets \$100,000 and Steven gets nothing.
$o_4$ : Sarah gets nothing and Steven gets nothing.

**Table 1.2**

If, indeed, Sarah is *selfish and greedy* – in the sense that, in evaluating the outcomes, she focuses exclusively on what she herself gets and prefers more money to less – then her ranking of the outcomes is as follows: $o_3 \succ o_1 \succ o_2 \sim o_4$ (which reads '$o_3$ is better than $o_1$, $o_1$ is better than $o_2$ and $o_2$ is just as good as $o_4$'). But there are other possibilities. For example, Sarah might be *fair-minded* and view the outcome where both get \$50,000 as better than all the other outcomes. For instance, her ranking could be $o_1 \succ o_3 \succ o_2 \succ o_4$; according to this ranking, besides valuing fairness, she also displays *benevolence* towards Steven, in the sense that – when comparing the two outcomes where she gets nothing, namely $o_2$ and $o_4$ – she prefers the one where at least Steven goes home with some money. If, in fact, Sarah is fair-minded and benevolent, then the logic underlying the above argument would yield the opposite conclusion, namely that she should choose *Split*.

Thus we cannot presume to know the answer to the question "What is the rational choice for Sarah?" if we don't know what her preferences are. It is a common mistake (unfortunately one that even game theorists sometimes make) to reason under the assumption that players are selfish and greedy. This is, typically, an unwarranted assumption. Research in experimental psychology, philosophy and economics has amply demonstrated that many people are strongly motivated by considerations of fairness. Indeed, fairness seems to motivate not only humans but also primates, as shown in the following video:[1]

http://www.ted.com/talks/frans_de_waal_do_animals_have_morals

The situation illustrated in Table 1.1 is not a game as we have no information about the preferences of the players; we use the expression *game-frame* to refer to it. In the case where there are only two players and each player has a small number of possible choices (also called strategies), a game-frame can be represented – as we did in Table 1.1 – by means of a table, with as many rows as the number of possible strategies of Player 1 and as many columns as the number of strategies of Player 2; each row is labeled with one strategy of Player 1 and each column with one strategy of Player 2; inside each cell of the table (which corresponds to a pair of strategies, one for Player 1 and one for Player 2) we write the corresponding outcome.

Before presenting the definition of game-frame, we remind the reader of what the Cartesian product of two or more sets is. Let $S_1$ and $S_2$ be two sets. Then the Cartesian product of $S_1$ and $S_2$, denoted by $S_1 \times S_2$, is the set of *ordered pairs* $(x_1, x_2)$ where $x_1$ is an element of $S_1$ $(x_1 \in S_1)$ and $x_2$ is an element of $S_2$ $(x_2 \in S_2)$. For example, if $S_1 = \{a, b, c\}$ and $S_2 = \{D, E\}$ then

$$S_1 \times S_2 = \{(a,D),(a,E),(b,D),(b,E),(c,D),(c,E)\}.$$

The definition extends to the general case of $n$ sets ($n \geq 2$): an element of $S_1 \times S_2 \times ... \times S_n$ is an ordered $n$-tuple $(x_1, x_2, ..., x_n)$ where, for each $i$ = 1,…,$n$, $x_i \in S_i$.

The definition of game-frame is as follows.

[1] Also available at https://www.youtube.com/watch?v=GcJxRqTs5nk

**Definition 1.1.** A *game-frame in strategic form* is a list of four items (a quadruple) $\langle I,(S_1,S_2,...,S_n),O,f\rangle$ where:

- $I = \{1,2,...,n\}$ is a set of *players* ($n \geq 2$).
- $(S_1,S_2,...,S_n)$ is a list of sets, one for each player. For every Player $i \in I$, $S_i$ is the set of *strategies* (or possible choices) of Player *i*. We denote by *S* the Cartesian product of these sets: $S = S_1 \times S_2 \times ... \times S_n$; thus an element of $S$ is a list $s = (s_1, s_2, ..., s_n)$ consisting of one strategy for each player. We call $S$ the set of *strategy profiles*.
- *O* is a set of *outcomes*.
- $f: S \to O$ is a function that associates with every strategy profile *s* an outcome $f(s) \in O$.

Using the notation of Definition 1.1, the situation illustrated in Table 1.1 is the following game-frame in strategic form:

- $I = \{1,2\}$ (letting Sarah be Player 1 and Steven Player 2),
- $(S_1,S_2) = (\{Split, Steal\}, \{Split, Steal\})$; thus $S_1 = S_2 = \{Split, Steal\}$, so that the set of strategy profiles is $S = \{(Split,Split),(Split,Steal),(Steal,Split),(Steal,Steal)\}$,
- *O* is the set of outcomes listed in Table 1.2,
- $f$ is the following function:

| $s:$ | $(Split,Split)$ | $(Split,Steal)$ | $(Steal,Split)$ | $(Steal,Steal)$ |
|---|---|---|---|---|
| $f(s):$ | $o_1$ | $o_2$ | $o_3$ | $o_4$ |

(that is, $f((Split,Split)) = o_1$, $f((Split,Steal)) = o_2$, etc.).

From a game-frame one obtains a game by adding, for each player, her preferences over (or ranking of) the possible outcomes. We use the following notation.

| Notation | Meaning |
|---|---|
| $o \succ_i o'$ | Player *i* considers outcome $o$ to be *better than* outcome $o'$ |
| $o \sim_i o'$ | Player *i* considers $o$ to be *just as good as* $o'$ (that is, Player *i* is indifferent between $o$ and $o'$) |
| $o \succsim_i o'$ | Player *i* considers $o$ to be *at least as good as* $o'$ (that is, either better than or just as good as) |

**Table 1.3**

For example, if $M$ denotes 'Mexican food' and $J$ 'Japanese food', $M \succ_{Alice} J$ means that Alice prefers Mexican food to Japanese food and $M \sim_{Bob} J$ means that Bob is indifferent between the two.

**Remark 1.1.** The "at least as good" relation $\succsim$ is sufficient to capture also strict preference $\succ$ and indifference $\sim$. In fact, starting from $\succsim$, one can define strict preference as follows: $o \succ_i o'$ if and only if $o \succsim_i o'$ and $o' \not\succsim_i o$ and one can define indifference as follows: $o \sim_i o'$ if and only if $o \succsim_i o'$ and $o' \succsim_i o$.

We will assume throughout this book that the "at least as good" relation $\succsim_i$ of Player $i$ – which embodies her preferences over (or ranking of) the outcomes – is *complete* (for every two outcomes $o_1$ and $o_2$, either $o_1 \succsim_i o_2$ or $o_2 \succsim_i o_1$) and *transitive* (if $o_1 \succsim_i o_2$ and $o_2 \succsim_i o_3$ then $o_1 \succsim_i o_3$).[2]

There are (at least) three ways of representing, or expressing, a complete and transitive ranking of a set of outcomes $O$. For example, suppose that $O = \{o_1, o_2, o_3, o_4, o_5\}$ and that we want to represent the following ranking: $o_3$ is the best outcome, it is better than $o_5$, which is just as good as $o_1$, $o_1$ is better than $o_4$, which, in turn, is better than $o_2$ (so that $o_2$ is the worst outcome). We can represent this ranking in one of the following ways.

1. By using the notation of Table 1.3: $o_3 \succ o_5 \sim o_1 \succ o_4 \succ o_2$.

2. By listing the outcomes in a column, starting with the best at the top and proceeding down to the worst, thus using the convention that if outcome $o$ is listed above outcome $o'$ then $o$ is preferred to $o'$, while if $o$ and $o'$ are written next to each other, then they are considered to be just as good:

| | |
|---|---|
| best | $o_3$ |
| | $o_1, o_5$ |
| | $o_4$ |
| worst | $o_2$ |

[2] Transitivity of the "at least as good" relation $\succsim$ implies transitivity of the indifference relation (if $o_1 \sim o_2$ and $o_2 \sim o_3$ then $o_1 \sim o_3$) as well as transitivity of the strict preference relation (not only in the sense that (1) if $o_1 \succ o_2$ and $o_2 \succ o_3$ then $o_1 \succ o_3$, but also that (2) if $o_1 \succ o_2$ and $o_2 \sim o_3$ then $o_1 \succ o_3$ and (3) if $o_1 \sim o_2$ and $o_2 \succ o_3$ then $o_1 \succ o_3$).

3. By assigning a number to each outcome, with the convention that *if the number assigned to $o$ is greater than the number assigned to $o'$ then $o$ is preferred to $o'$*, and if two outcomes are assigned the same number then they are considered to be just as good. For example, we could choose the following numbers: $\begin{matrix} o_1 & o_2 & o_3 & o_4 & o_5 \\ 6 & 1 & 8 & 2 & 6 \end{matrix}$. Such an assignment of numbers is called a *utility function*. A useful way of thinking of utility is as an "index of satisfaction": the higher the index the better the outcome; however, this suggestion is just to aid memory and should be taken with a grain of salt, because a utility function does not measure anything and, furthermore, as explained below, the actual numbers used as utility indices are completely arbitrary.[3]

**Definition 1.2.** Given a complete and transitive ranking $\succsim$ of a finite set of outcomes $O$, a function $U : O \to \mathbb{R}$ (where $\mathbb{R}$ denotes the set of real numbers) is said to be an *ordinal utility function that represents the ranking* $\succsim$ if, for every two outcomes $o$ and $o'$, $U(o) > U(o')$ if and only if $o \succ o'$ and $U(o) = U(o')$ if and only if $o \sim o'$. The number $U(o)$ is called the *utility of outcome* $o$.

**Remark 1.2.** Note that the statement "for Alice the utility of Mexican food is 10" is in itself a meaningless statement; on the other hand, what would be a meaningful statement is "for Alice the utility of Mexican food is 10 and the utility of Japanese food is 5", because such a statement conveys the information that she prefers Mexican food to Japanese food. However, the two numbers 10 and 5 have no other meaning besides the fact that 10 is greater than 5: for example we *cannot*, and should not, infer from these numbers that she considers Mexican food twice as good as Japanese food. The reason for this is that we could have expressed the same fact, namely that she prefers Mexican food to Japanese food, by assigning utility 100 to Mexican food and $-25$ to Japanese food, or with any other two numbers.

---

[3] Note that assigning a utility of 1 to an outcome *o* does *not* mean that *o* is the "first choice". Indeed, in this example a utility of 1 is assigned to the worst outcome: $o_2$ is the worst outcome because it has the lowest utility (which happens to be 1, in this example).

It follows from the above remark that there is an infinite number of utility functions that represent the same ranking. For instance, the following are equivalent ways of representing the ranking $o_3 \succ o_1 \succ o_2 \sim o_4$ (*f*, *g* and *h* are three out of the many possible utility functions):

| *outcome* → <br> *utility function* ↓ | $o_1$ | $o_2$ | $o_3$ | $o_4$ |
|---|---|---|---|---|
| $f$ | 5 | 2 | 10 | 2 |
| $g$ | 0.8 | 0.7 | 1 | 0.7 |
| $h$ | 27 | 1 | 100 | 1 |

Utility functions are a particularly convenient way of representing preferences. In fact, by using utility functions one can give a more condensed representation of games, as explained in the last paragraph of the following definition.

**Definition 1.3.** An *ordinal game in strategic form* is a quintuple $\langle I, (S_1, ..., S_n), O, f, (\succsim_1, ..., \succsim_n) \rangle$ where:

- $\langle I, (S_1, ..., S_n), O, f \rangle$ is a game-frame in strategic form (Definition 1.1) and
- for every Player $i \in I$, $\succsim_i$ is a complete and transitive ranking of the set of outcomes *O*.

If we replace each ranking $\succsim_i$ with a utility function $U_i$ that represents it, and we assign, to each strategy profile *s*, Player *i*'s utility of $f(s)$ (recall that $f(s)$ is the outcome associated with *s*) then we obtain a function $\pi_i : S \to \mathbb{R}$ called Player *i*'s *payoff function*. Thus $\pi_i(s) = U_i(f(s))$.[4] Having done so, we obtain a triple $\langle I, (S_1, ..., S_n), (\pi_1, ..., \pi_n) \rangle$ called a *reduced-form ordinal game in strategic form* ('reduced-form' because some information is lost, namely the specification of the possible outcomes).

For example, take the game-frame illustrated in Table 1.1, let Sarah be Player 1 and Steven Player 2 and name the possible outcomes as shown in Table 1.2. Let us add the information that both players are selfish and greedy (that is, Player 1's ranking is $o_3 \succ_1 o_1 \succ_1 o_2 \sim_1 o_4$ and Player 2's ranking is $o_2 \succ_2 o_1 \succ_2 o_3 \sim_2 o_4$) and let us represent their rankings with the following utility functions (note, again, that the

[4] Note that, in this book, the symbol $\pi$ is not used to denote the irrational number used to compute the circumference and area of a circle, but rather as the Greek letter for 'p' which stands for 'payoff'.

choice of numbers 2, 3 and 4 for utilities is arbitrary: any other three numbers would do):

| *outcome* → / *utility function* ↓ | $o_1$ | $o_2$ | $o_3$ | $o_4$ |
|---|---|---|---|---|
| $U_1$ (Player 1) | 3 | 2 | 4 | 2 |
| $U_2$ (Player 2) | 3 | 4 | 2 | 2 |

Then we obtain the following reduced-form game, where in each cell the first number is the payoff of Player 1 and the second number is the payoff of Player 2.

| | | Player 2 (Steven) | |
|---|---|---|---|
| | | Split | Steal |
| Player 1 (Sarah) | Split | 3 3 | 2 4 |
| | Steal | 4 2 | 2 2 |

**Table 1.4**

On the other hand, if we add to the game-frame of Table 1.1 the information that Player 1 is fair-minded and benevolent (that is, her ranking is $o_1 \succ_1 o_3 \succ_1 o_2 \succ_1 o_1$), while Player 2 is selfish and greedy and represent these rankings with the following utility functions,

| *outcome* → / *utility function* ↓ | $o_1$ | $o_2$ | $o_3$ | $o_4$ |
|---|---|---|---|---|
| $U_1$ (Player 1) | 4 | 2 | 3 | 1 |
| $U_2$ (Player 2) | 3 | 4 | 2 | 2 |

then we obtain the following reduced-form game:

| | | Player 2 (Steven) | |
|---|---|---|---|
| | | Split | Steal |
| Player 1 (Sarah) | Split | 4 3 | 2 4 |
| | Steal | 3 2 | 1 2 |

**Table 1.5**

In general, a player will act differently in different games, even if they are based on the same game-frame, because her incentives and objectives (as captured by her ranking of the outcomes) will be different. For example, one can argue that in the game of Table 1.4 a rational Player 1 would choose *Steal*, while in the game of Table 1.5 the rational choice for Player 1 is *Split*.

This is a good time to test your understanding of the concepts introduced in this section, by going through the exercises in Section 1.E.1 of Appendix 1.E at the end of this chapter.

## 1.2. Strict and weak dominance

In this section we define two relations on the set of strategies of a player. Before introducing the formal definition, we shall illustrate these notion with an example. The first relation is called "strict dominance". Let us focus our attention on one player, say Player 1, and select two of her strategies, say $a$ and $b$. We say that $a$ **strictly** dominates $b$ (for Player 1) if, for every possible strategy profile of the other players, strategy $a$ of Player 1, in conjunction with the strategies selected by the other players, yields a payoff for Player 1 which is **greater than** the payoff associated with strategy $b$ (in conjunction with the strategies selected by the other players). For example, consider the following two-player game, where only the payoffs of Player 1 are shown:

| | | Player 2 | | |
|---|---|---|---|---|
| | | *E* | *F* | *G* |
| Player 1 | *A* | 3 … | 2 … | 1 … |
| | *B* | 2 … | 1 … | 0 … |
| | *C* | 3 … | 2 … | 1 … |
| | *D* | 2 … | 0 … | 0 … |

**Table 1.6**

In this game for Player 1 strategy *A* strictly dominates strategy *B*:

- if Player 2 selects *E* then *A* in conjunction with *E* gives Player 1 a payoff of 3, while *B* in conjunction with *E* gives her only a payoff of 2,

- if Player 2 selects $F$ then $A$ in conjunction with $F$ gives Player 1 a payoff of 2, while $B$ in conjunction with $F$ gives her only a payoff of 1,
- if Player 2 selects $G$ then $A$ in conjunction with $G$ gives Player 1 a payoff of 1, while $B$ in conjunction with $G$ gives her only a payoff of 0.

In the game of Table 1.6 we also have that $A$ strictly dominates $D$ and $C$ strictly dominates $D$; however, it is not the case that $B$ strictly dominates $D$ because, in conjunction with strategy $E$ of Player 2, $B$ and $D$ yield the same payoff for Player 1.

The second relation is called "weak dominance". The definition is similar to that of strict dominance, but we replace 'greater than' with 'greater than or equal to' while insisting on at least one strict inequality: $a$ **weakly** dominates $b$ (for Player 1) if, for every possible strategy profile of the other players, strategy $a$ of Player 1, in conjunction with the strategies selected by the other players, yields a payoff for Player 1 which is **greater than or equal to** the payoff associated with strategy $b$ (in conjunction with the strategies selected by the other players) and, furthermore, there is at least one strategy profile of the other players against which strategy $a$ gives a larger payoff to Player 1 than strategy $b$. In the example of Table 1.6, we have that, while it is not true that $B$ *strictly* dominates $D$, it is true that $B$ *weakly* dominates $D$:

- if Player 2 selects $E$, then $B$ in conjunction with $E$ gives Player 1 the same payoff as $D$ in conjunction with $E$ (namely 2),
- if Player 2 selects $F$, then $B$ in conjunction with $F$ gives Player 1 a payoff of 2, while $B$ in conjunction with $F$ gives her only a payoff of 1,
- if Player 2 selects $G$ then $B$ in conjunction with $G$ gives Player 1 the same payoff as $D$ in conjunction with $G$ (namely 0).

In order to give the definitions in full generality we need to introduce some notation. Recall that $S$ denotes the set of strategy profiles, that is, an element $s$ of $S$ is an ordered list of strategies $s = (s_1, ..., s_n)$, one for each player. We will often want to focus on one player, say Player $i$, and view $s$ as a pair consisting of the strategy of Player $i$ and the remaining strategies of all the other players. For example, suppose that there are three players and the strategy sets are as follows: $S_1 = \{a,b,c\}$, $S_2 = \{d,e\}$ and $S_3 = \{f,g\}$. Then one possible strategy profile is $s = (b,d,g)$ (thus $s_1 = b$, $s_2 = d$ and $s_3 = g$). If we focus on, say, Player 2 then we will denote by $s_{-2}$ the sub-profile consisting of the strategies of the players *other than* 2: in this case $s_{-2} = (b,g)$. This gives us an alternative way of denoting $s$, namely as $(s_2, s_{-2})$. Continuing our example where $s = (b,d,g)$, letting $s_{-2} = (b,g)$, we can denote $s$ also by $(d, s_{-2})$ and we can write the result of replacing Player 2's strategy $d$ with her strategy $e$ in $s$ by $(e, s_{-2})$; thus $(d, s_{-2}) = (b,d,g)$ while $(e, s_{-2}) = (b,e,g)$. In general, given a Player $i$, we denote by $S_{-i}$ the set of strategy profiles of the players other than $i$ (that is, $S_{-i}$ is the Cartesian product of the strategy sets of the other players; in the above example we have that $S_{-2} = S_1 \times S_3 = \{a,b,c\} \times \{f,g\} = \{(a,f),(a,g),(b,f),(b,g),(c,f),(c,g)\}$. We denote an element of $S_{-i}$ by $s_{-i}$.

**Definition 1.4.** Given an ordinal game in strategic form, let $i$ be a Player and $a$ and $b$ two of her strategies ($a,b \in S_i$). We say that, for Player $i$,

- *a strictly dominates b* (or *b is strictly dominated by a*) if, in every situation (that is, no matter what the other players do), $a$ gives Player $i$ a payoff which is *greater* than the payoff that $b$ gives. Formally: for every $s_{-i} \in S_{-i}$, $\pi_i(a,s_{-i}) > \pi_i(b,s_{-i})$.[5]
- *a weakly dominates b* (or *b is weakly dominated by a*) if, in every situation, $a$ gives Player $i$ a payoff which is *greater than or equal to* the payoff that $b$ gives and, furthermore, there is at least one situation where $a$ gives a greater payoff than $b$. Formally: for every $s_{-i} \in S_{-i}$, $\pi_i(a,s_{-i}) \geq \pi_i(b,s_{-i})$ and there exists an $\bar{s}_{-i} \in S_{-i}$ such that $\pi_i(a,\bar{s}_{-i}) \geq \pi_i(b,\bar{s}_{-i})$.[6]
- *a is equivalent to b* if, in every situation, $a$ and $b$ give Player $i$ the same payoff. Formally: for every $s_{-i} \in S_{-i}$, $\pi_i(a,s_{-i}) = \pi_i(b,s_{-i})$.[7]

| Player 1 \ Player 2 | E | F | G |
|---|---|---|---|
| A | 3 ... | 2 ... | 1 ... |
| B | 2 ... | 1 ... | 0 ... |
| C | 3 ... | 2 ... | 1 ... |
| D | 2 ... | 0 ... | 0 ... |

For example, in the game of Table 1.6 reproduced above, we have that

- A strictly dominates B.
- A and C are equivalent.
- A strictly dominates D.
- B is strictly dominated by C.

[5] Or, stated in terms of rankings instead of payoffs, $f(a,s_{-i}) \succ_i f(b,s_{-i})$ for every $s_{-i} \in S_{-i}$.

[6] Or, stated in terms of rankings, $f(a,s_{-i}) \succsim_i f(b,s_{-i})$, for every $s_{-i} \in S_{-i}$, and there exists an $\bar{s}_{-i} \in S_{-i}$ such that $f(a,\bar{s}_{-i}) \succ_i f(b,\bar{s}_{-i})$.

[7] Or, stated in terms of rankings, $f(a,s_{-i}) \sim_i f(b,s_{-i})$, for every $s_{-i} \in S_{-i}$.

- B weakly (but not strictly) dominates D.
- C strictly dominates D.

**Remark 1.3.** Note that if strategy *a* strictly dominates strategy *b* then it also satisfies the conditions for weak dominance, that is, '*a* strictly dominates *b*' implies '*a* weakly dominates *b*'. Throughout the book the expression '*a weakly* dominates *b*' will be interpreted as '*a* dominates *b weakly but not strictly* '.

The expression '*a* dominates *b*' can be understood as '*a* is better than *b*'. The next term we define is 'dominant' which can be understood as ' best'. Thus one cannot meaningfully say "*a* dominates" because one needs to name another strategy that is dominated by *a*; for example, one would have to say "*a* dominates *b*". On the other hand, one *can* meaningfully say "*a* is dominant" because it is like saying "*a* is best", which means "*a* is better than every other strategy".

**Definition 1.5.** Given an ordinal game in strategic form, let $i$ be a Player and $a$ one of her strategies ($a \in S_i$). We say that, for Player *i*,

- *a is a strictly dominant* strategy if *a* strictly dominates every other strategy of Player *i*.
- *a is a weakly dominant* strategy if, for every other strategy *x* of Player *i*, one of the following is true: (1) either *a* weakly dominates *x* or (2) *a* is equivalent to *x*.

For example, in the game shown in Table 1.6, *A* and *C* are both weakly dominant strategies for Player 1. Note that if a player has two or more strategies that are weakly dominant, then any two of those strategies must be equivalent. On the other hand, there can be at most one strictly dominant strategy.

**Remark 1.4.** The reader should convince herself/himself that the definition of weakly dominant strategy given in Definition 1.5 is equivalent to the following:
$a \in S_i$ is a weakly dominant strategy for Player $i$ if and only if, for every $s_{-i} \in S_{-i}$,

$$\pi_i(a, s_{-i}) \geq \pi_i(s_i, s_{-i}) \text{ for every } s_i \in S_i.$$ [8]

In accordance with the convention established in Remark 1.3, the expression '*a* is a weakly dominant strategy' will have the default interpretation '*a* is a *weakly but not strictly* dominant strategy'.

Note: if you claim that, for some player, "strategy *x* is (weakly or strictly) dominated" then you ought to name another strategy of that player that dominates *x*. Saying "*x* is dominated" is akin to saying "*x* is worse": worse than what? On the

[8] Or, stated in terms of rankings, for every $s_{-i} \in S_{-i}$, $f(a, s_{-i}) \succsim_i f(s_i, s_{-i})$ for every $s_i \in S_i$.

other hand, claiming that strategy $y$ is weakly dominant is akin to claiming that it is best, that is, better than, or just as good as, any other strategy.

**Definition 1.6.** Given an ordinal game in strategic form, let $s = (s_1, ..., s_n)$ be a strategy profile. We say that

- $s$ *is a strict dominant-strategy equilibrium* if, for every Player $i$, $s_i$ is a strictly dominant strategy.
- $s$ *is a weak dominant-strategy equilibrium* if, for every Player $i$, $s_i$ is a weakly dominant strategy and, furthermore, for at least one Player $j$, $s_j$ is not a strictly dominant strategy.

If we refer to a strategy profile as a dominant-strategy equilibrium, without qualifying it as weak or strict, then the default interpretation will be 'weak'.

In the game of Table 1.4, reproduced below, *Steal* is a weakly dominant strategy for each player and thus (*Steal*,*Steal*) is a weak dominant-strategy equilibrium.

| | | Player 2 (Steven) | |
|---|---|---|---|
| | | Split | Steal |
| Player 1 (Sarah) | Split | 3 3 | 2 4 |
| | Steal | 4 2 | 2 2 |

In the game of Table 1.5, reproduced below, *Split* is a strictly dominant strategy for Player 1, while *Steal* is a weakly (but not strictly) dominant strategy for Player 2 and thus (*Split*,*Steal*) is a weak dominant-strategy equilibrium.

| | | Player 2 (Steven) | |
|---|---|---|---|
| | | Split | Steal |
| Player 1 (Sarah) | Split | 4 3 | 2 4 |
| | Steal | 3 2 | 1 2 |

The *Prisoner's Dilemma* is an example of a game with a strict dominant-strategy equilibrium. For a detailed account of the history of this game and an in-depth analysis of it see

http://plato.stanford.edu/entries/prisoner-dilemma or
http://en.wikipedia.org/wiki/Prisoner's_dilemma.

An instance of the Prisoner's Dilemma is the following. Doug and Ed work for the same company and the annual party is approaching. They know that they are the only two candidates for the best-worker-of-the-year prize and at the moment they are tied; however, only one person can be awarded the prize and thus, unless one of them manages to outperform the other, nobody will receive the prize. Each chooses between exerting *Normal* effort or *Extra* effort (that is, work overtime) before the party. The game-frame is as shown in Table 1.7 below.

| | | Player 2 (Ed): Normal effort | Player 2 (Ed): Extra effort |
|---|---|---|---|
| Player 1 (Doug) | Normal effort | $o_1$ | $o_2$ |
| | Extra effort | $o_3$ | $o_4$ |

$o_1$ : nobody gets the prize and nobody sacrifices family time

$o_2$ : Ed gets the prize and sacrifices family time, Doug does not

$o_3$ : Doug gets the prize and sacrifices family time, Ed does not

$o_4$ : nobody gets the prize and both sacrifice family time

**Table 1.7**

Suppose that both Doug and Ed are willing to sacrifice family time to get the prize, but otherwise value family time; furthermore, they are envious of each other, in the sense that they prefer nobody getting the prize to the other person's getting the prize (even at the personal cost of sacrificing family time). That is, their rankings are as follows: $o_3 \succ_{Doug} o_1 \succ_{Doug} o_4 \succ_{Doug} o_2$ and $o_2 \succ_{Ed} o_1 \succ_{Ed} o_4 \succ_{Ed} o_3$. Using utility function with values from the set {0,1,2,3} we can represent the game in reduced form as follows:

| | | Player 2 (Ed): Normal effort | Player 2 (Ed): Extra effort |
|---|---|---|---|
| Player 1 (Doug) | Normal effort | 2, 2 | 0, 3 |
| | Extra effort | 3, 0 | 1, 1 |

**Table 1.8**
The Prisoner's Dilemma game

In this game exerting extra effort is a strictly dominant strategy for every player; thus (*Extra effort, Extra effort*) is a strict dominant-strategy equilibrium.

**Definition 1.7.** Given an ordinal game in strategic form, let $o$ and $o'$ be two outcomes. We say that $o$ is *strictly Pareto superior* to $o'$ if every player prefers $o$ to $o'$ (that is, if $o \succ_i o'$ for every Player *i*). We say that $o$ is *weakly Pareto superior* to $o'$ if every player considers $o$ to be at least as good as $o'$ and at least one player prefers $o$ to $o'$ (that is, if $o \succsim_i o'$ for every Player *i* and there is a Player *j* such that $o \succ_j o'$). In reduced-form games, this definition can be extended to strategy profiles as follows. If $s$ and $s'$ are two strategy profiles, then $s$ is *strictly Pareto superior* to $s'$ if $\pi_i(s) > \pi_i(s')$ for every Player $i$ and $s$ is *weakly Pareto superior* to $s'$ if $\pi_i(s) \geq \pi_i(s')$ for every Player $i$ and, furthermore, there is a Player $j$ such that $\pi_j(s) > \pi_j(s')$.

For example, in the Prisoner's Dilemma game of Table 1.8, outcome $o_1$ is strictly Pareto superior to $o_4$ or, in terms of strategy profiles, (*Normal effort, Normal effort*) is strictly Pareto superior to (*Extra effort, Extra effort*).

When a player has a strictly dominant strategy, it would be irrational for that player to choose any other strategy, since she would be guaranteed a lower payoff in every possible situation (that is, no matter what the other players do). Thus in the Prisoner's Dilemma individual rationality leads to (*Extra effort, Extra effort*) despite the fact that both players would be better off if they both chose *Normal effort*. It is obvious that if the players could reach a *binding* agreement to exert normal effort then they would do so; however, the underlying assumption in non-cooperative game theory is that such agreements are not possible (e.g. because of lack of communication or because such agreements are illegal or cannot be enforced in a court of law, etc.). Any non-binding agreement to choose *Normal effort* would not be viable: if one player expects the other player to stick to the agreement, then he will gain by cheating and choosing *Extra effort*. On the other hand, if a player does not believe that the other player will honor the agreement then he will gain by deviating from the agreement herself. The Prisoner's Dilemma game is often used to illustrate a conflict between *individual* rationality and *collective* rationality: (*Extra effort, Extra effort*) is the individually rational outcome while (*Normal effort, Normal effort*) would be the collectively rational one.

This is a good time to test your understanding of the concepts introduced in this section, by going through the exercises in Section 1.E.2 of Appendix 1.E at the end of this chapter.

## 1.3. Second-price auction

The *second-price auction,* or Vickrey auction, is an example of a game that has a weak dominant-strategy equilibrium. It is a "sealed-bid" auction where bidders submit bids without knowing the bids of the other participants in the auction. The object which is auctioned is then assigned to the bidder who submits the highest bid (the winner), but the winner pays not her own bid but rather the second-highest bid, that is the highest bid among the bids that remain after removing the winner's own bid. Tie-breaking rules must be specified for selecting the winner when the highest bid is submitted by two or more bidders (in which case the winner ends up paying her own bid, because the second-highest bid is equal to the winner's bid). We first illustrate this auction with an example.

Two oil companies bid for the right to drill a field. The possible bids are \$10 million, \$20 million, ..., \$50 million. In case of ties the winner is Player 2 (this was decided earlier by tossing a coin). Let us take the point of view of Player 1. Suppose that Player 1 ordered a geological survey and, based on the report, concludes that the oil field would generate a profit of \$30 million. Suppose also that Player 1 is indifferent between any two outcomes where the oil field is given to Player 2 and prefers to get the oil field herself if and only if it has to pay not more than \$30 million for it; furthermore, getting the oil field for \$30 million is just as good as not getting it. Then we can take as utility function for Player 1 the net gain to Player 1 from the oil field (defined as profits from oil extraction minus the price paid for access to the oil field) if Player 1 wins, and zero otherwise. In Table 1.9 we have written inside each cell only the payoff of Player 1. For example, why is Player 1's payoff 20 when it bids \$30M and Player 2 bids \$10M? Since Player 1's bid is higher than Player 2's bid, Player 1 is the winner and thus the drilling rights are assigned to Player 1; hence Player 1 obtains something worth \$30M and pays, not its own bid of \$30M, but the bid of Player 2, namely \$10M; thus Player 1's net gain is \$(30 – 10)M = \$20M. It can be verified that for Player 1 submitting a bid equal to the value it assigns to the object (namely, a bid of \$30 million) is a weakly dominant strategy: it always gives Player 1 the largest of the payoffs that are possible, given the bid of the other player. This does not imply that it is the only weakly dominant strategy; indeed, in this example bidding \$40M is also a weakly dominant strategy for Player 1 (in fact, it is equivalent to bidding \$30M).

| | | Player 2 | | | | |
|---|---|---|---|---|---|---|
| | | \$10M | \$20M | \$30M | \$40M | \$50M |
| **Player 1 (value \$30M)** | \$10M | 0 | 0 | 0 | 0 | 0 |
| | \$20M | 20 | 0 | 0 | 0 | 0 |
| | \$30M | 20 | 10 | 0 | 0 | 0 |
| | \$40M | 20 | 10 | 0 | 0 | 0 |
| | \$50M | 20 | 10 | 0 | –10 | 0 |

**Table 1.9**

A second-price auction where, in case of ties, the winner is Player 2.

Now we can describe the second-price auction in more general terms. Let $n \geq 2$ be the number of bidders. We assume that all non-negative bids are allowed and that the tie-breaking rule favors the player with the lowest index among those who submit the highest bid: for example, if the highest bid is \$250 and it is submitted by Players 5, 8 and 10, then the winner is Player 5. We shall denote the possible outcomes as pairs $(i, p)$, where $i$ is the winner and $p$ is the price that the winner has to pay. Finally we denote by $b_i$ the bid of Player $i$. We start by describing the case where there are only two bidders and then generalize to the case of an arbitrary number of bidders. We denote the set of non-negative numbers by $[0,\infty)$.

**The case where $n$ = 2:** in this case we have that $I$ = {1,2}, $S_1 = S_2 = [0,\infty)$, $O = \{(i,p) : i \in \{1,2\}, p \in [0,\infty)\}$ and $f : S \to O$ is given by

$$f\left((b_1,b_2)\right) = \begin{cases} (1,b_2) & \text{if } b_1 \geq b_2 \\ (2,b_1) & \text{if } b_1 < b_2 \end{cases}.$$

**The case where $n \geq 2$:** in the general case the second-price auction is the following game-frame:

- $I$ = {1,…,$n$}
- $S_i = [0,\infty)$ for every $i$ = 1,…,$n$. We denote an element of $S_i$ by $b_i$.
- $O = \{(i,p) : i \in I, p \in [0,\infty)\}$
- $f : S \to O$ is defined as follows. Let $H(b_1,...,b_n)$ be the set of bidders who submit the highest bid: $H(b_1,...,b_n) = \{i \in I : b_i \geq b_j \text{ for all } j \in I\}$ and let $\hat{i}(b_1,...,b_n)$ be the smallest number in the set $H(b_1,...,b_n)$, that is, the winner of the auction. Finally, let $b^{\max}$ denote the maximum bid (that is, $b^{\max}(b_1,...,b_n) = Max\{b_1,...,b_n\}$), and let $b^{second}(b_1,...,b_n)$ be the second-highest bid (that is, $b^{second}(b_1,...,b_n) = Max\{\{b_1,...,b_n\} \setminus \{b^{\max}(b_1,...,b_n)\}\}$.[9]

Then $f(b_1,...,b_n) = (\hat{i}(b_1,...,b_n)\ ,\ b^{second}(b_1,...,b_n))$.

How much should a player bid in a second-price auction? Since what we have described is a game-frame and not a game, we cannot answer the question unless we specify the player's preferences over the set of outcomes $O$. Let us say that Player $i$ in a second-price auction is *selfish and greedy* if she only cares about whether or not she wins and, conditional on winning, prefers to pay less; furthermore, she prefers winning to not winning if and only if she has to pay less than the true value of the object for her, which we denote by $v_i$, and is indifferent between not winning and winning if she has to pay $v_i$. Thus the ranking of a selfish and greedy player is as follows (together with everything that follows from transitivity):

$(i,p) \succ_i (i,p')$ if and only if $p < p'$

$(i,p) \succ_i (j,p')$ for all $j \neq i$ and for all $p'$, if and only if $p < v_i$

$(i,v_i) \sim_i (j,p')$ for all $j \neq i$ and for all $p'$

$(j,p) \sim_i (k,p')$ for all $j \neq i, k \neq i$ and for all $p$ and $p'$

An ordinal utility function that represents those preferences is[10]

$$U_i(j,p) = \begin{cases} v_i - p & \text{if } i = j \\ 0 & \text{if } i \neq j \end{cases}.$$

Using this utility function we get the following payoff function for Player $i$:

$$\pi_i(b_1,...,b_n) = \begin{cases} v_i - b^{second}(b_1,...,b_n) & \text{if } i = \hat{i}(b_1,...,b_n) \\ 0 & \text{if } i \neq \hat{i}(b_1,...,b_n) \end{cases}.$$

We can now state the following theorem. The proof is given in Appendix 1.A at the end of this chapter.

[9] For example, if $n = 5$, $b_1 = 10$, $b_2 = 14$, $b_3 = 8$, $b_4 = 14$ and $b_5 = 14$ then $H(10,14,8,14,14) = \{2,4,5\}$, $\hat{i}(10,14,8,14,14) = 2$, $b^{\max}(10,14,8,14,14) = b^{second}(10,14,8,14,14) = 14$.

[10] Of course there are many more. For example, also the following utility function represents those preferences: $U_i(j,p) = \begin{cases} 2^{(v_i - p)} & \text{if } i = j \\ 1 & \text{if } i \neq j \end{cases}$.

**Theorem 1.1 [Vickrey, 1961].** In a second-price auction, if Player *i* is selfish and greedy then it is a weakly dominant strategy for Player *i* to bid her true value, that is, to choose $b_i = v_i$.

Note that, for a player who is not selfish and greedy, Theorem 1.1 is not true. For example, if a player has the same preferences as above for the case where she wins, but, conditional on not winning, prefers the other player to pay as much as possible (she is spiteful) or as little as possible (she is generous), then bidding her true value is no longer a dominant strategy.

This is a good time to test your understanding of the concepts introduced in this section, by going through the exercises in Section 1.E.3 of Appendix 1.E at the end of this chapter.

# 1.4. The pivotal mechanism

An article in the *Davis Enterprise* (the local newspaper in Davis, California) on January 12, 2001 started with the following paragraph:

> By consensus, the Davis City Council agreed Wednesday to order a communitywide public opinion poll to gauge how much Davis residents would be willing to pay for a park tax and a public safety tax.

Opinion polls of this type are worthwhile only if there are reasons to believe that the people who are interviewed will respond honestly. But will they? If I would like more parks and believe that the final tax I will have to pay is independent of the amount I state in the interview, I would have an incentive to overstate my willingness to pay, hoping to swing the decision towards building a new park. On the other hand, if I fear that the final tax will be affected by the amount I report, then I might have an incentive to understate my willingness to pay.

The *pivotal mechanism*, or Clarke mechanism, is a game designed to give the participants an incentive to report their true willingness to pay.

A public project, say to build a park, is under consideration. The cost of the project is $C. There are *n* individuals in the community. If the project is carried out, individual *i* (*i* = 1,..., *n*) will have to pay $\$c_i$ (with $c_1 + c_2 + \ldots + c_n = C$); these amounts are specified as part of the project. Note that we allow for the possibility that some individuals might have to contribute a larger share of the total cost *C* than others (e.g. because they live closer to the projected park and would therefore

benefit more from it). Individual $i$ has an initial wealth of $\$\overline{m}_i > 0$. If the project is carried out, individual $i$ receives benefits from it that she considers equivalent to receiving $\$v_i$. Note that for some individual $i$, $v_i$ could be negative, that is, the individual could be harmed by the project (e.g. because she likes peace and quiet and a park close to her home would bring extra traffic and noise). We assume that individual $i$ has the following utility-of-wealth function:

$$U_i(\$m) = \begin{cases} m & \text{if the project is not carried out} \\ m + v_i & \text{if the project is carried out} \end{cases}$$

The socially efficient decision is to carry out the project if and only if $\sum_{i=1}^{n} v_i > C$ (recall that $\sum$ is the summation sign: $\sum_{i=1}^{n} v_i$ is a short-hand for $v_1 + v_2 + \ldots + v_n$). For example, suppose that $n = 2$, $\overline{m}_1 = 50$, $\overline{m}_2 = 60$, $v_1 = 19$, $v_2 = -15$, $C = 6$, $c_1 = 6$, $c_2 = 0$. In this case $\sum_{i=1}^{n} v_i = 19 - 15 = 4 < C = 6$ hence the project should *not* be carried out. To see this consider the following table.

| | If the project is **not** carried out | If the project is carried out |
|---|---|---|
| Utility of individual 1 | 50 | 50 + 19 –6 = 63 |
| Utility of individual 2 | 60 | 60 –15 = 45 |

If the project is carried out, Individual 1 has a utility gain of 13, while Individual 2 has a utility loss of 15. Since the loss is greater than the gain, we have a *Pareto inefficient* situation. Individual 2 could propose the following alternative to Individual 1: let us not carry out the project and I will pay you \$14. Then Individual 1's wealth and utility would be 50 + 14 = 64 and Individual 2's wealth and utility would be $60 - 14 = 46$ and thus they would both be better off.

Thus Pareto efficiency requires that the project be carried out if and only if $\sum_{i=1}^{n} v_i > C$. This would be a simple decision for the government if it knew the $v_i$'s. But, typically, these values are private information to the individuals. Can the government find a way to induce the individuals to reveal their true valuations? It seems that in general the answer is No: those who gain from the project would have an incentive to overstate their potential gains, while those who suffer would have an incentive to overstate their potential losses. Influenced by Vickrey's work on second-price auctions, Clarke suggested the following procedure or game. Each individual $i$ is asked to submit a number $w_i$ which will be interpreted as the gross benefit (if positive) or harm (if negative) that individual $i$ associates with the project. Note that, in principle, individual $i$ can lie and report a value $w_i$ which is different from the true value $v_i$. Then the decision will be:

$$\text{Carry out the project?} \begin{cases} Yes & if \sum_{j=1}^{n} w_j > C \\ No & if \sum_{j=1}^{n} w_j \leq C \end{cases}$$

However, this is not the end of the story. Each individual will be classified as either not pivotal or pivotal.

$$\text{Individual } i \text{ is } \textbf{not} \text{ pivotal if } \begin{cases} either & \left( \sum_{j=1}^{n} w_j > C \;\; and \;\; \sum_{j \neq i} w_j > \sum_{j \neq i} c_j \right) \\ or & \left( \sum_{j=1}^{n} w_j \leq C \;\; and \;\; \sum_{j \neq i} w_j \leq \sum_{j \neq i} c_j \right) \end{cases}$$

and she is pivotal otherwise. In other words individual $i$ is pivotal if the decision about the project that would be made in the restricted society resulting from removing individual $i$ is different from the decision that is made when she is included. If an individual is not pivotal then she has to pay no taxes. If individual $i$ is pivotal then she has to pay a tax in the amount of $\left| \sum_{j \neq i} w_j - \sum_{j \neq i} c_j \right|$, the absolute value of $\sum_{j \neq i} w_j - \sum_{j \neq i} c_j$ (recall that the absolute value of $a$ is equal to $a$, if $a$ is positive, and to $-a$, if $a$ is negative; for instance, $|4| = 4$ and $|-4| = -(-4) = 4$).

For example, let $n$ = 3, $C$ = 10, $c_1$ = 3, $c_2$ = 2, $c_3$ = 5. Suppose that they state the following benefits/losses (which may or may not be the true ones): $w_1 = -1$, $w_2 = 8$, $w_3 = 3$. Then $\sum_{i=1}^{3} w_i = 10 = C$. Thus the project will **not** be carried out. Who is pivotal? The answer is provided in the following table.

| Individual | Σwj (including i) | Σcj (including i) | Decision | Σwj j ≠ i (without i) | Σcj j ≠ i (without i) | Decision | Pivotal? | Tax |
|---|---|---|---|---|---|---|---|---|
| 1 | 10 | 10 | No | 8 + 3 = 11 | 2 + 5 = 7 | Yes | Yes | 11 – 7 = 4 |
| 2 | 10 | 10 | No | –1 + 3 = 2 | 3 + 5 = 8 | No | No | 0 |
| 3 | 10 | 10 | No | –1 + 8 = 7 | 3 + 2 = 5 | Yes | Yes | 7 – 5 = 2 |

It may seem that, since it involves paying a tax, being pivotal is a bad thing and one should try to avoid it. It is certainly possible for individual $i$ to make sure that she is not pivotal: all she has to do is to report $w_i = c_i$; in fact, if $\sum_{j \neq i} w_j > \sum_{j \neq i} c_j$ then adding $c_i$ to both sides yields $\sum_{j=1}^{n} w_j > C$ and if $\sum_{j \neq i} w_j \leq \sum_{j \neq i} c_j$ then adding $c_i$ to both sides yields $\sum_{j=1}^{n} w_j \leq C$ . It is not true, however, that it is best to avoid being pivotal. The following example shows that one can gain by being truthful even if it involves being pivotal and thus having to pay a tax. Let $n$ = 4, $C$ = 15, $c_1$ = 5, $c_2$ = 0, $c_3$ = 5 and $c_4$ = 5. Suppose that $\bar{m}_1$ = 40 and $v_1$ = 25. Imagine that you are individual 1 and, for whatever reason, you expect the following reports by the other individuals: $w_2$ = –40, $w_3$ = 15 and $w_4$ = 20. If you report $w_1 = c_1 = 5$ then you ensure that you are not pivotal. In this case $\sum_{j=1}^{4} w_j = 5 - 40 + 15 + 20 = 0 < C = 15$ and thus the project is not carried out and your utility is equal to $\bar{m}_1$ = 40. If you report truthfully, that is, you report $w_1 = v_1 = 25$ then $\sum_{j=1}^{4} w_j = 25 - 40 + 15 + 20 = 20 > C = 15$ and the project is carried out; furthermore, you are pivotal and have to pay a tax $t_1$ equal to $\left| \sum_{j=2}^{4} w_j - \sum_{j=2}^{4} c_j \right|$ = $|(-40 + 15 + 20) - (0 + 5 + 5)| = |-15| = 15$ and your utility will be $\bar{m}_1 + v_1 - c_1 - t_1$ = 40 + 25 – 5 – 15 =45; hence you are better off. Indeed, the following theorem states that no individual can ever gain by lying.

The proof of Theorem 1.2 is given in Appendix 1.A at the end of this chapter.

**Theorem 1.2 [Clarke, 1971].** In the pivotal mechanism (under the assumed preferences) truthful revelation (that is, stating $w_i = v_i$) is a weakly dominant strategy for every Player *i*.

This is a good time to test your understanding of the concepts introduced in this section, by going through the exercises in Section 1.E.4 of Appendix 1.E at the end of this chapter.

## 1.5. Iterated deletion procedures

If in a game a player has a (weakly or strictly) dominant strategy then the player ought to choose that strategy: in the case of strict dominance choosing any other strategy guarantees that the player will do worse and in the case of weak dominance no other strategy can give a better outcome, no matter what the other players do. Unfortunately, games that have a dominant-strategy equilibrium are not very common. What should a player do when she does not have a dominant strategy? We shall consider two iterative deletion procedures that can help solve some games.

**1.5.1. IDSDS.** *The Iterated Deletion of Strictly Dominated Strategies* (IDSDS) is the following procedure or algorithm. Given a finite ordinal strategic-form game $G$, let $G^1$ be the game obtained by removing from $G$, for every Player *i*, those strategies of Player *i* (if any) that are strictly dominated in $G$ by some other strategy; let $G^2$ be the game obtained by removing from $G^1$, for every Player *i*, those strategies of Player *i* (if any) that are strictly dominated in $G^1$ by some other strategy, and so on. Let $G^\infty$ be the output of this procedure. Since the initial game $G$ is finite, $G^\infty$ will be obtained in a finite number of steps. Figure 1.10 illustrates this procedure. If $G^\infty$ contains a single strategy profile (this is not the case in the example of Figure 1.10) then we call that strategy profile the *iterated strict dominant-strategy equilibrium*. If $G^\infty$ contains two or more strategy profiles then we refer to those strategy profiles merely as the *output of the IDSDS procedure*. For example, in the game of Figure 1.10 the output of the IDSDS procedure is the set of strategy profiles $\{(A,e),(A,f),(B,e),(B,f)\}$.

What is the significance of the output of the IDSDS procedure? Consider game $G$ of Figure 1.10. Since, for Player 2, *h* is strictly dominated by *g*, if Player 2 is rational she will not play *h*. Thus, if Player 1 believes that Player 2 is rational then he believes that Player 2 will not play *h*, that is, he restricts attention to game $G^1$; since, in $G^1$, $D$ is strictly dominated by $C$ for Player 1, if Player 1 is rational he will not play $D$. It follows that if Player 2 believes that Player 1 is rational and that Player 1 believes that Player 2 is rational, then Player 2 restricts attention to game $G^2$ where rationality requires that Player 2 not play *g*, etc. It will be shown in a later chapter that if the players are rational and there is common belief of rationality[11] then only strategy profiles that survive the IDSDS procedure can be played; the converse is also true: any strategy profile that survives the IDSDS procedure is compatible with common belief of rationality.

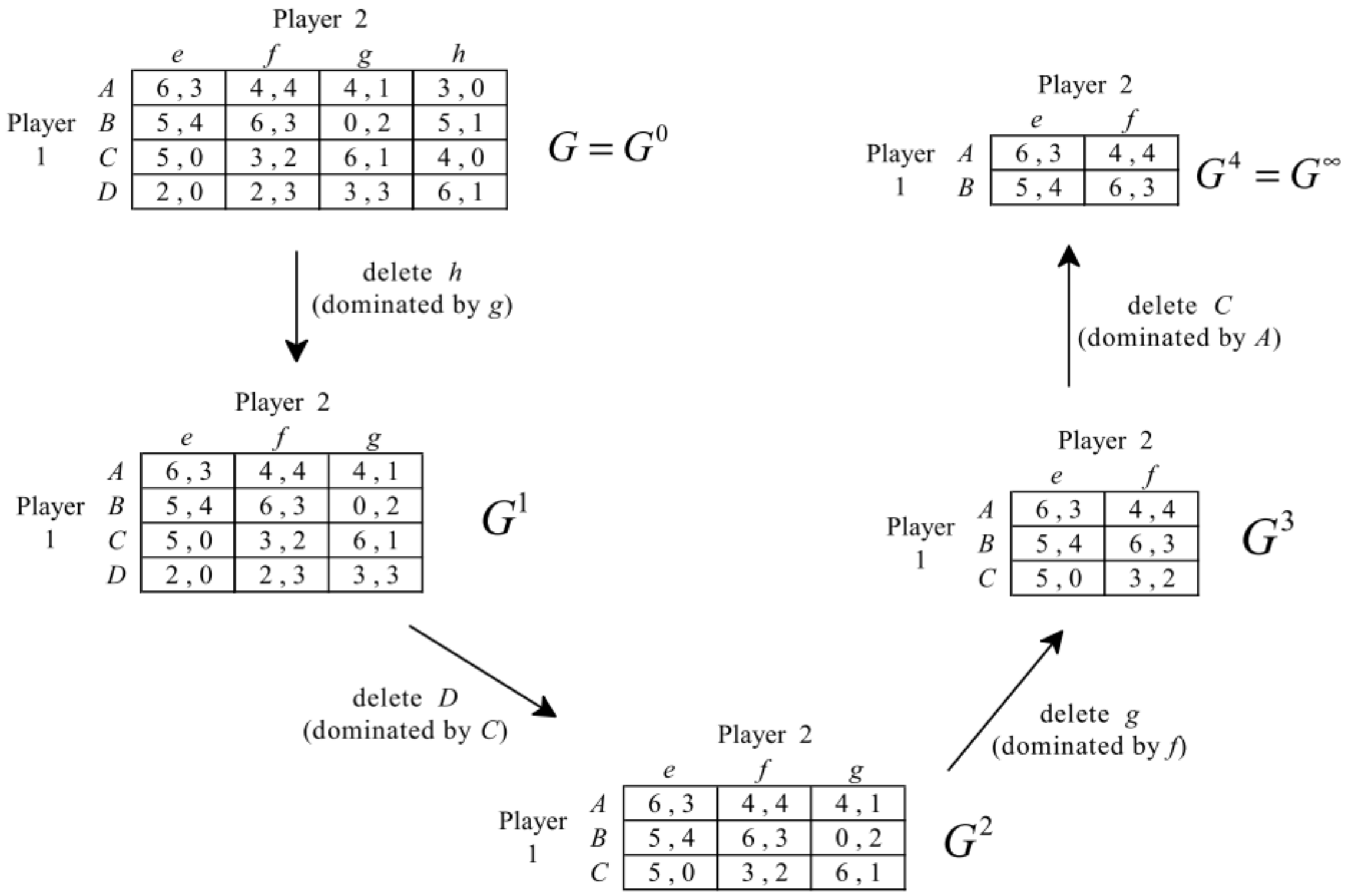

**Figure 1.10**

**Remark 1.5 .** In finite games, the order in which strictly dominated strategies are deleted is irrelevant, in the sense that any sequence of deletions of strictly dominated strategies leads to the same output.

[11] An event *E* is commonly believed if everybody believes *E* and everybody believes that everybody believes *E* and everybody believes that everybody believes that everybody believes *E*, and so on.

**1.5.1. IDWDS.** *The Iterated Deletion of Weakly Dominated Strategies* (IDWDS) is a weakening of IDSDS in that it allows the deletion also of *weakly* dominated strategies. However, this procedure has to be defined carefully, since in this case the order of deletion *can* matter. To see this, consider the game shown in Figure 1.11.

Player 2

| Player 1 | | L | R |
|---|---|---|---|
| | A | 4 , 0 | 0 , 0 |
| | T | 3 , 2 | 2 , 2 |
| | M | 1 , 1 | 0 , 0 |
| | B | 0 , 0 | 1 , 1 |

**Figure 1.11**

Since M is strictly dominated by T for Player 1, we can delete it and obtain

Player 2

| Player 1 | | L | R |
|---|---|---|---|
| | A | 4 , 0 | 0 , 0 |
| | T | 3 , 2 | 2 , 2 |
| | B | 0 , 0 | 1 , 1 |

Now L is weakly dominated by R for Player 2. Deleting L we are left with

Player 2

| Player 1 | | R |
|---|---|---|
| | A | 0 , 0 |
| | T | 2 , 2 |
| | B | 1 , 1 |

Now A and B are strictly dominated by T. Deleting them we are left with (T,R) with corresponding payoffs (2,2).

Alternatively, going back to the game of Figure 1.11, we could note that B is strictly dominated by T; deleting B we are left with

| Player 1 \ Player 2 | L | R |
|---|---|---|
| A | 4 , 0 | 0 , 0 |
| T | 3 , 2 | 2 , 2 |
| M | 1 , 1 | 0 , 0 |

Now R is weakly dominated by L for Player 2. Deleting R we are left with

| Player 1 \ Player 2 | L |
|---|---|
| A | 4 , 0 |
| T | 3 , 2 |
| M | 1 , 1 |

Now T and M are strictly dominated by A and deleting them leads to (A,L) with corresponding payoffs (4,0). Since one order of deletion leads to (T,R) with payoffs (2,2) and the other to (A,L) with payoffs (4,0), the procedure is not well defined: the output of a well-defined procedure should be unique.

**Definition 1.8 (IDWDS).** In order to avoid the problem illustrated above, the IDWDS procedure is defined as follows: *at every step identify, for every player, all the strategies that are weakly (or strictly) dominated and then delete all such strategies in that step*. If the output of the IDWDS procedure is a single strategy profile then we call that strategy profile the *iterated weak dominant-strategy equilibrium* (otherwise we just use the expression 'output of the IDWDS procedure').

For example, the IDWDS procedure when applied to the game of Figure 1.11 leads to the following output:

| Player 1 \ Player 2 | L | R |
|---|---|---|
| A | 4 , 0 | 0 , 0 |
| T | 3 , 2 | 2 , 2 |

Hence the game of Figure 1.11 does not have an iterated weak dominant-strategy equilibrium.

The interpretation of the output of the IDWDS procedure is not as simple as that of the IDSDS procedure: certainly common belief of rationality is not sufficient. In order to delete weakly dominated strategies one needs to appeal not only to rationality but also to some notion of caution: a player should not completely rule out any of her opponents' strategies. However, this notion of caution is in direct conflict with the process of deletion of strategies. In this book we shall not address the issue of how to justify the IDWDS procedure.

This is a good time to test your understanding of the concepts introduced in this section, by going through the exercises in Section 1.E.5 of Appendix 1.E at the end of this chapter.

## 1.6. Nash equilibrium

Games where either the IDSDS procedure or the IDWDS procedure leads to a unique strategy profile are not very common. How can one then "solve" games that are not solved by either procedure? The notion of Nash equilibrium offers a more general alternative. We first define Nash equilibrium for a two-player game.

**Definition 1.9.** Given an ordinal game in strategic form with two players, a strategy profile $s^* = (s_1^*, s_2^*) \in S_1 \times S_2$ is a Nash equilibrium if the following two conditions are satisfied:

(1) for every $s_1 \in S_1$, $\pi_1(s_1^*, s_2^*) \geq \pi_1(s_1, s_2^*)$ (or stated in terms of outcomes and preferences, $f(s_1^*, s_2^*) \succsim_1 f(s_1, s_2^*)$) and

(2) for every $s_2 \in S_2$, $\pi_2(s_1^*, s_2^*) \geq \pi_1(s_1^*, s_2)$ (or, $f(s_1^*, s_2^*) \succsim_2 f(s_1^*, s_2)$).

For example, in the game of Figure 1.12 there are two Nash equilibria: $(T, L)$ and $(B, C)$.

| | | Player 2 | | |
|---|---|---|---|---|
| | | $L$ | $C$ | $R$ |
| Player 1 | $T$ | 3 , 2 | 0 , 0 | 1 , 1 |
| | $M$ | 3 , 0 | 1 , 5 | 4 , 4 |
| | $B$ | 1 , 0 | 2 , 3 | 3 , 0 |

**Figure 1.12**

There are several possible interpretations of this definition.

**'No regret' interpretation**: $s^*$ is a Nash equilibrium if there is no player who, after observing the opponent's choice, regrets his own choice (in the sense that he could have done better with a different strategy of his, given the observed strategy of the opponent).

**'Self-enforcing agreement' interpretation**: imagine that the players are able to communicate before playing the game and reach a non-binding agreement expressed as a strategy profile $s^*$; then no player will have an incentive to deviate from the agreement, if she believes that the other player will follow the agreement, if and only if $s^*$ is a Nash equilibrium.

**'Viable recommendation' interpretation**: imagine that a third party makes a public recommendation to each player on what strategy to play; then no player will have an incentive to deviate from the recommendation, if she believes that the other players will follow the recommendation, if and only if the recommended strategy profile is a Nash equilibrium.

**'Transparency of reason' interpretation**: if players are all "equally rational" and Player 2 reaches the conclusion that she should play *y*, then Player 1 must be able to duplicate Player 2's reasoning process and come to the same conclusion; it follows that Player 1's choice of strategy is not rational unless it is a strategy *x* that is optimal against *y*. A similar argument applies to Player 2's choice of strategy (*y* must be optimal against *x*) and thus $(x,y)$ is a Nash equilibrium.

It is clear that all of the above interpretations are mere rewording of the formal definition of Nash equilibrium in terms of the inequalities of Definition 1.9.

The generalization of Definition 1.9 to games with more than two players is straightforward.

**Definition 1.10.** Given an ordinal game in strategic form with $n$ players, a strategy profile $s^* \in S$ is a Nash equilibrium if the following inequalities are satisfied:

$$\text{for every Player } i,\ \pi_i(s^*) \geq \pi_i(s_1^*,...,s_{i-1}^*,s_i,s_{i+1}^*,...,s_n^*) \text{ for all } s_i \in S_i.$$

The reader should convince himself/herself that a (weak or strict) dominant strategy equilibrium is a Nash equilibrium and the same is true of a (weak or strict) iterated dominant-strategy equilibrium.

**Definition 1.11.** Consider an ordinal game in strategic form, a Player $i$ and a strategy profile $\bar{s}_{-i} \in S_{-i}$ of the players other than $i$. A strategy $s_i \in S_i$ of Player $i$ is a *best reply (or best response) to* $\bar{s}_{-i}$ if $\pi_i(s_i,\bar{s}_{-i}) \geq \pi_i(s_i',\bar{s}_{-i})$ for every $s_i' \in S_i$.

For example, in the game of Figure 1.12, reproduced below, for Player 1 there are two best replies to $L$, namely $M$ and $T$, while the unique best reply to $C$ is $B$ and the unique best reply to $R$ is $M$; for Player 2 the best reply to $T$ is $L$, the best reply to $M$ is $C$ and the best reply to $B$ is $C$.

| | | Player 2 | | |
|---|---|---|---|---|
| | | $L$ | $C$ | $R$ |
| Player 1 | $T$ | 3 , 2 | 0 , 0 | 1 , 1 |
| | $M$ | 3 , 0 | 1 , 5 | 4 , 4 |
| | $B$ | 1 , 0 | 2 , 3 | 3 , 0 |

**Remark 1.6.** Using the notion of best reply, an alternative definition of Nash equilibrium is as follows: $\bar{s} \in S$ is a Nash equilibrium if and only if, for every Player $i$, $\bar{s}_i \in S_i$ is a best reply to $\bar{s}_{-i} \in S_{-i}$.

A quick way to find the Nash equilibria of a two-player game is as follows: in each column of the table underline the largest payoff of Player 1 in that column (if there are several instances, underline them all) and in each row underline the largest payoff of Player 2 in that row; if a cell has both payoffs underlined then the corresponding strategy profile is a Nash equilibrium. Underlining of the maximum payoff of Player 1 in a given column identifies the best reply of Player 1 to the strategy of Player 2 that labels that column and similarly for Player 2. This procedure is illustrated in Figure 1.13, where there is a unique Nash equilibrium, namely (*B*,*E*).

| | | Player 2 | | | |
|---|---|---|---|---|---|
| | | *E* | *F* | *G* | *H* |
| | *A* | 4 , 0 | 3 , 2 | 2 , 3 | 4 , 1 |
| Player | *B* | 4 , 2 | 2 , 1 | 1 , 2 | 0 , 2 |
| 1 | *C* | 3 , 6 | 5 , 5 | 3 , 1 | 5 , 0 |
| | *D* | 2 , 3 | 3 , 2 | 1 , 2 | 3 , 3 |

**Figure 1.13**

Exercise 1.2 in Appendix 1.E explains how to represent a three-player game by means of a set of tables. In a three-player game the procedure for finding the Nash equilibria is the same, with the necessary adaptation for Player 3: in each cell underline the payoff of Player 3 if and only if her payoff is the largest of all her payoffs in the same cell across different tables. This is illustrated in Figure 1.14, where there is a unique Nash equilibrium, namely (*B*,*R*,*W*).

| | | Player 2 | |
|---|---|---|---|
| | | *L* | *R* |
| Player | *T* | 0 , 0 , 0 | 2 , 8 , 6 |
| 1 | *B* | 5 , 3 , 2 | 3 , 4 , 2 |

Player 3 chooses *W*

| | Player 2 | |
|---|---|---|
| | *L* | *R* |
| *T* | 0 , 0 , 0 | 1 , 2 , 5 |
| *B* | 1 , 6 , 1 | 0 , 0 , 1 |

Player 3 chooses *E*

**Figure 1.14**

Unfortunately, when the game has too many players or too many strategies and it is thus impossible or impractical to represent it as a set of tables, there is no quick procedure for finding the Nash equilibria: one must simply apply the definition of Nash equilibrium. For example, consider the following game.

**Example 1.1.** There are 50 players. A benefactor asks them to simultaneously and secretly write on a piece of paper a request, which must be a multiple of \$10 up to a maximum of \$100 (thus the possible strategies of each player are \$10, \$20, ..., \$90, \$100). He will then proceed as follows: if not more than 10% of the players ask for \$100 then he will grant every player's request, otherwise every player will get nothing. Assume that every player is selfish and greedy (only cares about how much money she gets and prefers more money to less). What are the Nash equilibria of this game? There are several:

- every strategy profile where 7 or more players request \$100 is a Nash equilibrium (everybody gets nothing and no player can get a positive amount by unilaterally changing her request, since there will still be more than 10% requesting \$100; on the other hand, convince yourself that a strategy profile where exactly 6 players request \$100 is not a Nash equilibrium),
- every strategy profile where exactly 5 players request \$100 and the remaining players request \$90 is a Nash equilibrium.

Any other strategy profile is not a Nash equilibrium: (1) if fewer than 5 players request \$100, then a player who requested less than \$100 can increase her payoff by switching to a request of \$100, (2) if exactly 5 players request \$100 and among the remaining players there is one who is not requesting \$90, then that player can increase her payoff by increasing her request to \$90.

We conclude this section by noting that, since so far we have restricted attention to ordinal games, there is no guarantee that an arbitrary game will have at least one Nash equilibrium. An example of a game that has no Nash equilibria is the *Matching Penny* game. This is a simultaneous two-player game where each player has a coin and decides whether to show the Heads face or the Tails face. If both choose $H$ or both choose $T$ then Player 1 wins, otherwise Player 2 wins. Each player strictly prefers the outcome where she herself wins to the alternative outcome. The game is illustrated in Figure 1.15.

| | | Player 2 | |
|---|---|---|---|
| | | $H$ | $T$ |
| Player 1 | $H$ | 1 0 | 0 1 |
| | $T$ | 0 1 | 1 0 |

**Figure 1.15**

This is a good time to test your understanding of the concepts introduced in this section, by going through the exercises in Section 1.E.6 of Appendix 1.E at the end of this chapter.

# 1.7. Games with infinite strategy sets

Games where the strategy set of one or more players is infinite cannot be represented using a table or set of tables. However, all the concepts introduced in this chapter can still be applied. In this section we will focus on the notion of Nash equilibrium. We start with an example.

**Example 1.2.** There are two players. Each player has to write down a real number (not necessarily an integer) greater than or equal to 1; thus the strategy sets are $S_1 = S_2 = [1,\infty)$. Payoffs are as follows ($\pi_1$ is the payoff of Player 1, $\pi_2$ the payoff of Player 2, $x$ is the number written by Player 1 and $y$ the number written by Player 2):

$$\pi_1(x,y) = \begin{cases} x-1 & \text{if } x<y \\ 0 & \text{if } x \geq y \end{cases}, \qquad \pi_2(x,y) = \begin{cases} y-1 & \text{if } x>y \\ 0 & \text{if } x \leq y \end{cases}$$

What are the Nash equilibria of this game?

There is only one Nash equilibrium, namely (1,1) with payoffs (0,0). First of all, we must show that (1,1) is indeed a Nash equilibrium. If Player switched to some $x > 1$ then her payoff would remain 0 and the same is true for Player 2 if he unilaterally switched to some $y > 1$: $\pi_1(x,1)=0$, for all $x \in [1,\infty)$ and $\pi_2(1,y)=0$, for all $y \in [1,\infty)$. Now we show that no other pair $(x,y)$ is a Nash equilibrium. Consider first an arbitrary pair $(x,y)$ with $x = y > 1$. Then $\pi_1(x,y) = 0$, but if Player 1 switched to an $\hat{x}$ strictly between 1 and $x$ $(1<\hat{x}<x)$ her payoff would be $\pi_1(\hat{x},y)=\hat{x}-1>0$ (recall that, by hypothesis, $x = y$). Now consider an arbitrary $(x,y)$ with $x < y$. Then $\pi_1(x,y) = x - 1$, but if Player 1 switched to an $\hat{x}$ strictly between $x$ and $y$ $(x<\hat{x}<y)$ her payoff would be $\pi_1(\hat{x},y)=\hat{x}-1>x-1$. The argument for ruling out pairs $(x,y)$ with $y < x$ is similar.

Note the interesting fact that, for Player 1, $x = 1$ is a weakly dominated strategy: indeed it is weakly dominated by *any* other strategy: $x = 1$ guarantees a payoff of 0 for Player 1, while any $\hat{x}>1$ would yield a positive payoff to Player 1 in some cases (against any $y > \hat{x}$) and 0 in the remaining cases. The same is true for Player 2. Thus in this game there is *a unique Nash equilibrium where the strategy of each player is weakly dominated*!

[Note: the rest of this section makes use of calculus. The reader who is not familiar with calculus should skip this part.]

We conclude this section with an example based on the analysis of competition among firms proposed by Augustine Cournot in a book published in 1838. In fact, Cournot is the one who invented what we now call Nash equilibrium, although his analysis was restricted to a small class of games. Consider $n \geq 2$ firms which produce an identical product. Let $q_i$ be the quantity produced by Firm $i$ ($i = 1,...,n$). For Firm $i$ the cost of producing $q_i$ units of output is $c_i q_i$, where $c_i$ is a positive constant. For simplicity we will restrict attention to the case of two firms ($n = 2$) and identical cost functions: $c_1 = c_2 = c$. Let $Q$ be total industry output, that is, $Q = q_1 + q_2$. The price at which each firm can sell each unit of output is given by the inverse demand function $P = a - bQ$ where $a$ and $b$ are positive constants. Cournot assumed that each firm was only interested in its own profit and preferred higher profit to lower profit (that is, each firm is "selfish and greedy"). The profit function of Firm 1 is given by

$$\pi_1(q_1,q_2) = Pq_1 - cq_1 = \left[a - b(q_1 + q_2)\right]q_1 - cq_1 = (a-c)q_1 - b(q_1)^2 - bq_1q_2 .$$

Similarly, the profit function of Firm 2 is given by

$$\pi_2(q_1,q_2) = (a-c)q_2 - b(q_2)^2 - bq_1q_2$$

Cournot defined an equilibrium as a pair $\left(\bar{q}_1, \bar{q}_2\right)$ such that

$$\pi_1\left(\bar{q}_1, \bar{q}_2\right) \geq \pi_1\left(q_1, \bar{q}_2\right), \text{ for every } q_1 \geq 0 \qquad (♣)$$

and

$$\pi_2\left(\bar{q}_1, \bar{q}_2\right) \geq \pi_2\left(\bar{q}_1, q_2\right), \text{ for every } q_2 \geq 0 \qquad (♦)$$

Of course, this is the same as saying that $\left(\bar{q}_1, \bar{q}_2\right)$ is a Nash equilibrium of the game where the players are the two firms, the strategy sets are $S_1 = S_2 = [0,\infty)$ and the payoff functions are the profit functions. How do we find a Nash equilibrium? First of all, note that the profit functions are differentiable. Secondly note that (♣) says that, having fixed the value of $q_2$ at $\bar{q}_2$, the function $\pi_1\left(q_1, \bar{q}_2\right)$ viewed as a function of $q_1$ alone is maximized at the point $q_1 = \bar{q}_1$. A necessary condition for this (if $\bar{q}_1 > 0$) is that the derivative of this function be zero at the point $q_1 = \bar{q}_1$, that is, it must be that $\frac{\partial \pi_1}{\partial q_1}\left(\bar{q}_1, \bar{q}_2\right) = 0$. This condition is also sufficient since the second derivative of this function is always negative ($\frac{\partial^2 \pi_1}{\partial q_1^2}(q_1, q_2) = -2b$ for every $(q_1, q_2)$). Similarly, by (◆), it must be that $\frac{\partial \pi_2}{\partial q_2}(\bar{q}_1, \bar{q}_2) = 0$. Thus the Nash equilibrium is found by solving the system of two equations $\begin{cases} a - c - 2bq_1 - bq_2 = 0 \\ a - c - 2bq_2 - bq_1 = 0 \end{cases}$. The solution is $\bar{q}_1 = \bar{q}_2 = \frac{a-c}{3b}$. The corresponding price is $\bar{P} = a - b\left(2\frac{a-c}{3b}\right) = \frac{a+2c}{3}$ and the corresponding profits are $\pi_1(\frac{a-c}{3b}, \frac{a-c}{3b}) = \pi_2(\frac{a-c}{3b}, \frac{a-c}{3b}) = \frac{(a-c)^2}{9b}$. For example, if $a = 25,\ b = 2,\ c = 1$ then the Nash equilibrium is given by (4,4) with corresponding profits of 32 for each firm. The analysis can easily be extended to the case of more than two firms. The reader who is interested in further exploring the topic of competition among firms can consult any textbook on Industrial Organization.

This is a good time to test your understanding of the concepts introduced in this section, by going through the exercises in Section 1.E.7 of Appendix 1.E at the end of this chapter.

# Appendix 1.A: Proofs of theorems

**Theorem 1.1 [Vickrey, 1961].** In a second-price auction, if Player $i$ is selfish and greedy then it is a weakly dominant strategy for Player $i$ to bid her true value, that is, to choose $b_i = v_i$.

**Proof.** In order to make the notation simpler and the argument more transparent, we give the proof for the case where $n = 2$. We shall prove that bidding $v_1$ is a weakly dominant strategy for Player 1 (the proof for Player 2 is similar). Assume that Player 1 is selfish and greedy. Then we can take her payoff function to be as follows: $\pi_1(b_1,b_2) = \begin{cases} v_1 - b_2 & \text{if } b_1 \geq b_2 \\ 0 & \text{if } b_1 < b_2 \end{cases}$. We need to show that, whatever bid Player 2 submits, Player 1 cannot get a higher payoff by submitting a bid different from $v_1$. Two cases are possible (recall that $b_2$ denotes the actual bid of Player 2, which is unknown to Player 1).

**Case 1**: $b_2 \leq v_1$. In this case, bidding $v_1$ makes Player 1 the winner and his payoff is $v_1 - b_2 \geq 0$. Consider a different bid $b_1$. If $b_1 \geq b_2$ then Player 1 is still the winner and his payoff is still $v_1 - b_2 \geq 0$. Thus such a bid is as good as (hence not better than) $v_1$. If $b_1 < b_2$ then the winner is Player 2 and Player 1 gets a payoff of 0. Thus such a bid is also not better than $v_1$.

**Case 2**: $b_2 > v_1$. In this case, bidding $v_1$ makes Player 2 the winner and thus Player 1 gets a payoff of 0. Any other bid $b_1 < b_2$ gives the same outcome and payoff. On the other hand, any bid $b_1 \geq b_2$ makes Player 1 the winner, giving him a payoff of $v_1 - b_2 < 0$, thus making Player 1 worse off than with a bid of $v_1$. ■

**Theorem 1.2 [Clarke, 1971].** In the pivotal mechanism (under the assumed preferences) truthful revelation (that is, stating $w_i = v_i$) is a weakly dominant strategy for every Player $i$.

**Proof.** Consider an individual $i$ and possible statements $w_j$ for $j \neq i$. Several cases are possible.

**Case 1:** $\sum_{j \neq i} w_j > \sum_{j \neq i} c_j$ and $v_i + \sum_{j \neq i} w_j > c_i + \sum_{j \neq i} c_j = C$. Then

| | decision | *i*'s tax | *i*'s utility |
|---|---|---|---|
| if *i* states $v_i$ | Yes | 0 | $\bar{m}_i + v_i - c_i$ |
| if *i* states $w_i$ such that $w_i + \sum_{j \neq i} w_j > C$ | Yes | 0 | $\bar{m}_i + v_i - c_i$ |
| if *i* states $w_i$ such that $w_i + \sum_{j \neq i} w_j \leq C$ | No | $\sum_{j \neq i} w_j - \sum_{j \neq i} c_j$ | $\bar{m}_i - \left( \sum_{j \neq i} w_j - \sum_{j \neq i} c_j \right)$ |

Individual *i* cannot gain by lying if and only if $\bar{m}_i + v_i - c_i \geq \bar{m}_i - \left( \sum_{j \neq i} w_j - \sum_{j \neq i} c_j \right)$,
i.e. if $v_i + \sum_{j \neq i} w_j \geq C$ which is true by our hypothesis.

**Case 2:** $\sum_{j \neq i} w_j > \sum_{j \neq i} c_j$ and $v_i + \sum_{j \neq i} w_j \leq c_i + \sum_{j \neq i} c_j = C$. Then

| | decision | *i*'s tax | *i*'s utility |
|---|---|---|---|
| if states $v_i$ | No | $\sum_{j \neq i} w_j - \sum_{j \neq i} c_j$ | $\bar{m}_i - \left( \sum_{j \neq i} w_j - \sum_{j \neq i} c_j \right)$ |
| if states $w_i$ such that $w_i + \sum_{j \neq i} w_j \leq C$ | No | $\sum_{j \neq i} w_j - \sum_{j \neq i} c_j$ | $\bar{m}_i - \left( \sum_{j \neq i} w_j - \sum_{j \neq i} c_j \right)$ |
| if states $w_i$ such that $w_i + \sum_{j \neq i} w_j > C$ | Yes | 0 | $\bar{m}_i + v_i - c_i$ |

Individual *i* cannot gain by lying if and only if $\bar{m}_i - \left( \sum_{j \neq i} w_j - \sum_{j \neq i} c_j \right) \geq \bar{m}_i + v_i - c_i$ ,
i.e. if $v_i + \sum_{j \neq i} w_j \leq C$ which is true by our hypothesis.

**Case 3:** $\sum_{j\neq i} w_j \leq \sum_{j\neq i} c_j$ and $v_i + \sum_{j\neq i} w_j \leq c_i + \sum_{j\neq i} c_j = C$. Then

| | decision | *i*'s tax | *i*'s utility |
|---|---|---|---|
| if *i* states $v_i$ | No | 0 | $\bar{m}_i$ |
| if *i* states $w_i$ such that $w_i + \sum_{j\neq i} w_j \leq C$ | No | 0 | $\bar{m}_i$ |
| if *i* states $w_i$ such that $w_i + \sum_{j\neq i} w_j > C$ | Yes | $\left(\sum_{j\neq i} c_j - \sum_{j\neq i} w_j\right)$ (recall that $\sum_{j\neq i} w_j \leq \sum_{j\neq i} c_j$) | $\bar{m}_i + v_i - c_i - \left(\sum_{j\neq i} c_j - \sum_{j\neq i} w_j\right)$ |

Individual *i* cannot gain by lying if and only if $\bar{m}_i \geq \bar{m}_i + v_i - c_i - \left(\sum_{j\neq i} c_j - \sum_{j\neq i} w_j\right)$, i.e. if $v_i + \sum_{j\neq i} w_j \leq C$, which is our hypothesis.

**Case 4:** $\sum_{j\neq i} w_j \leq \sum_{j\neq i} c_j$ and $v_i + \sum_{j\neq i} w_j > c_i + \sum_{j\neq i} c_j = C$. Then

| | decision | *i*'s tax | *i*'s utility |
|---|---|---|---|
| if *i* states $v_i$ | Yes | $\left(\sum_{j\neq i} c_j - \sum_{j\neq i} w_j\right)$ (recall that $\sum_{j\neq i} w_j \leq \sum_{j\neq i} c_j$) | $\bar{m}_i + v_i - c_i - \left(\sum_{j\neq i} c_j - \sum_{j\neq i} w_j\right)$ |
| if *i* states $w_i$ such that $w_i + \sum_{j\neq i} w_j > C$ | Yes | $\left(\sum_{j\neq i} c_j - \sum_{j\neq i} w_j\right)$ (recall that $\sum_{j\neq i} w_j \leq \sum_{j\neq i} c_j$) | $\bar{m}_i + v_i - c_i - \left(\sum_{j\neq i} c_j - \sum_{j\neq i} w_j\right)$ |
| if *i* states $w_i$ such that $w_i + \sum_{j\neq i} w_j \leq C$ | No | 0 | $\bar{m}_i$ |

Individual *i* cannot gain by lying if and only if $\bar{m}_i + v_i - c_i - \left(\sum_{j\neq i} c_j - \sum_{j\neq i} w_j\right) \geq \bar{m}_i$, i.e. if $v_i + \sum_{j\neq i} w_j \geq C$, which is true by our hypothesis.

Since we have covered all the possible cases, the proof is complete. ■

# Appendix 1.E: Exercises

## 1.E.1. Exercises for Section 1.1: Game frames and games

The answers to the following exercises are in Appendix S at the end of this chapter.

**Exercise 1.1**. Antonia and Bob cannot decide where to go to dinner. Antonia proposes the following procedure: she will write on a piece of paper either the number 2 or the number 4 or the number 6, while Bob will write on his piece of paper either the number 1 or 3 or 5. They will write their numbers secretly and independently. They then will show each other what they wrote and choose a restaurant according to the following rule: if the sum of the two numbers is 5 or less, they will go to a Mexican restaurant, if the sum is 7 they will go to an Italian restaurant and if the number is 9 or more they will go to a Japanese restaurant.

**(a)** Let Antonia be Player 1 and Bob Player 2. Represent this situation as a game frame, first by writing out each element of the quadruple of Definition 1.1 and then by using a table (label the rows with Antonia's strategies and the columns with Bob's strategies, so that we can think of Antonia as choosing the row and Bob as choosing the column).

**(b)** Suppose that Antonia and Bob have the following preferences (where $M$ stands for 'Mexican', $I$ for 'Italian' and $J$ for 'Japanese). For Antonia: $M \succ_{Antonia} I \succ_{Antonia} J$; for Bob: $I \succ_{Bob} M \succ_{Bob} J$. Using utility function with values 1, 2 and 3 represent the corresponding reduced-form game as a table.

**Exercise 1.2**. Consider the following two-player game-frame where each player is given a set of cards and each card has a number on it. The players are Antonia (Player 1) and Bob (Player 2). Antonia's cards have the following numbers (one number on each card): 2, 4 and 6, whereas Bob's cards are marked 0, 1 and 2 (thus different numbers from the previous exercise). Antonia chooses one of her own cards and Bob chooses one of his own cards: this is done without knowing the other player's choice. The outcome depends on the sum of the points of the chosen cards, as follows. If the sum of the points on the two chosen cards is greater than or equal to 5, Antonia gets \$10 minus that sum; otherwise (that is, if the sum is less than 5) she gets nothing; furthermore, if the sum of points is an odd number, Bob gets as many dollars as that sum; if the sum of points turns out to be an even number and is less than or equal to 6, Bob gets \$2; otherwise he gets nothing.

**(a)** Represent the game-frame described above by means of a table. As in the previous exercise, assign the rows to Antonia and the columns to Bob.

**(b)** Using the game-frame of part (a) obtain a reduced-form game by adding the information that each player is selfish and greedy. This means that each player only cares about how much money he/she gets and prefers more money to less.

**Exercise 1.3**. Alice (Player 1), Bob (Player 2), and Charlie (Player 3) play the following simultaneous game. They are sitting in different rooms facing a keyboard with only one key and each has to decide whether or not to press the key. Alice wins if the number of people who press the key is odd (that is, all three of them or only Alice or only Bob or only Charlie) , Bob wins if exactly two people (he may be one of them) press the key and Charlie wins if nobody presses the key.

**(a)** Represent this situation as a game-frame. Note that we can represent a three-player game with a *set of tables*: Player 1 chooses the row, Player 2 chooses the column and Player 3 chooses the table (that is, we label the rows with Player 1's strategies, the columns with Player 2's strategies and the tables with Player 3's strategies).

**(b)** Using the game-frame of part (a) obtain a reduced-form game by adding the information that each player prefers winning to not winning and is indifferent between any two outcomes where he/she does not win. For each player use a utility function with values from the set {0,1}.

**(c)** Using the game-frame of part (a) obtain a reduced-form game by adding the information that (1) each player prefers winning to not winning, (2) Alice is indifferent between any two outcomes where she does not win, (3) conditional on not winning, Bob prefers if Charlie wins rather than Alice, (4) conditional on not winning, Charlie prefers if Bob wins rather than Alice. For each player use a utility function with values from the set {0,1,2}.

## 1.E.2. Exercises for Section 1.2: Strict/weak dominance

The answers to the following exercises are in Appendix S at the end of this chapter.

**Exercise 1.4**. There are two players. Each player is given an unmarked envelope and asked to put in it either nothing or \$300 of his own money or \$600 of his own money. A referee collects the envelopes, opens them, gathers all the money, then adds 50% of that amount (using his own money) and divides the total into two equal parts which he then distributes to the players.

**(a)** Represent this game frame with two alternative tables: the first table showing in each cell the amount of money distributed to Player 1 and the amount of money distributed to Player 2, the second table showing the change in wealth of each player (money received minus contribution).

**(b)** Suppose that Player 1 has some animosity towards the referee and ranks the outcomes in terms of how much money the referee loses (the more, the better), while Player 2 is selfish and greedy and ranks the outcomes in terms of her own net gain. Represent the corresponding game using a table.

**(c)** Is there a strict dominant-strategy equilibrium?

**Exercise 1.5**. **(a)** For the game of Part (b) of Exercise 1.1 determine, for each player, whether the player has *strictly* dominated strategies.
**(b)** For the game of Part (b) of Exercise 1.1 determine, for each player, whether the player has *weakly* dominated strategies.

**Exercise 1.6**. There are three players. Each player is given an unmarked envelope and asked to put in it either nothing or \$3 of his own money or \$6 of his own money. A referee collects the envelopes, opens them, gathers all the money and then doubles the amount (using his own money) and divides the total into three equal parts which he then distributes to the players. For example, if Players 1 and 2 put nothing and Player 3 puts \$6, then the referee adds another \$6 so that the total becomes \$12, divides this sum into three equal parts and gives \$4 to each player. Each player is selfish and greedy, in the sense that he ranks the outcomes exclusively in terms of his net change in wealth (what he gets from the referee minus what he contributed).

**(a)** Represent this game by means of a set of tables. (Do not treat the referee as a player.)

**(b)** For each player and each pair of strategies determine if one of the two dominates the other and specify if it is weak or strict dominance.

**(c)** Is there a strict dominant-strategy equilibrium?

## 1.E.3. Exercises for Section 1.3: Second price auction

The answers to the following exercises are in Appendix S at the end of this chapter.

**Exercise 1.7**. For the second-price auction partially illustrated in Table 1.9, complete the representation by adding the payoffs of Player 2, assuming that Player 2 assigns a value of \$50M to the field and, like Player 1, ranks the outcomes in terms of the net gain from the oil field (defined as profits minus the price paid, if Player 2 wins, and zero otherwise).

**Exercise 1.8**. Consider the following "third-price" auction. There are $n \geq 3$ bidders. A single object is auctioned and Player $i$ values the object \$$v_i$, with $v_i > 0$. The bids are simultaneous and secret. The utility of Player $i$ is: 0 if she does not win and $(v_i - p)$ if she wins and pays \$$p$. Every non-negative number is an admissible bid. Let $b_i$ denote the bid of Player $i$. The winner is the highest bidder. In case of ties the bidder with the lowest index among those who submitted the highest bid wins (e.g. if the highest bid is \$120 and it is submitted by players 6, 12 and 15, then the winner is Player 6). The losers don't get anything and don't pay anything. The winner gets the object and pays the **third** highest bid, which is defined as follows. Let $i$ be the winner and fix a Player $j$ such that $b_j = \max(\{b_1,...,b_n\} \setminus \{b_i\})$ [if $\max(\{b_1,...,b_n\} \setminus \{b_i\})$ contains more than one element, then we pick *one* of them].

Then the third price is defined as $\max(\{b_1,...,b_n\} \setminus \{b_i, b_j\})$. For example, if $n$ = 3 and the bids are $b_1 = 30$, $b_2 = 40$ and $b_3 = 40$ then the winner is Player 2 and she pays \$30; if $b_1 = b_2 = b_3 = 50$ then the winner is Player 1 and she pays \$50.

For simplicity, let us restrict attention to the case where $n$ = 3 and $v_1 > v_2 > v_3 > 0$. Does Player 1 have a weakly dominant strategy in this auction?

## 1.E.4. Exercises for Section 1.4: The pivotal mechanism

The answers to the following exercises are in Appendix S at the end of this chapter.

**Exercise 1.9**. The pivotal mechanism is used to decide whether a new park should be built. There are 5 individuals. According to the proposed project, the cost of the park would be allocated as follows:

| individual | 1 | 2 | 3 | 4 | 5 |
|---|---|---|---|---|---|
| share of cost | $c_1$ = \$30 | $c_2$= \$25 | $c_3$ = \$25 | $c_4$ = \$15 | $c_5$ = \$5 |

For every individual $i$ = 1, .., 5, let $v_i$ be the perceived gross benefit (if positive; perceived gross loss, if negative) from having the park built. The $v_i$'s are as follows:

| individual | 1 | 2 | 3 | 4 | 5 |
|---|---|---|---|---|---|
| gross benefit | $v_1$ = \$60 | $v_2$= \$15 | $v_3$ = \$55 | $v_4$ = \$–25 | $v_5$ = \$–20 |

(Thus the **net** benefit (loss) to individual $i$ is $v_i - c_i$). Individual $i$ has the following utility of wealth function (where $m_i$ denotes the wealth of individual $i$):

$$U_i = \begin{cases} m_i & \text{if the project is not carried out} \\ m_i + v_i & \text{if the project is carried out} \end{cases}$$

(Let $\bar{m}_i$ be the initial endowment of money of individual $i$ and assume that $\bar{m}_i$ is large enough that it exceeds $c_i$ plus any tax that the individual might have to pay.)

**(a)** What is the Pareto-efficient decision: to build the park or not?

Assume that the pivotal mechanism is used, so that each individual $i$ is asked to state a number $w_i$ which is going to be interpreted as the **gross** benefit to individual $i$ from carrying out the project. The are no restrictions on the number $w_i$: it can be positive, negative or zero. Suppose that the individuals make the following announcements:

| individual | 1 | 2 | 3 | 4 | 5 |
|---|---|---|---|---|---|
| stated benefit | $w_1$ = \$70 | $w_2$= \$10 | $w_3$ = \$65 | $w_4$ = \$−30 | $w_5$ = \$5 |

**(b)** Would the park be built based on the above announcements?

**(c)** Using the above announcements and the rules of the pivotal mechanism, fill in the following table:

| individual | 1 | 2 | 3 | 4 | 5 |
|---|---|---|---|---|---|
| Pivotal? | | | | | |
| Tax | | | | | |

**(d)** As you know, in the pivotal mechanism each individual has a dominant strategy. If all the individuals played their dominant strategies, would the park be built?

**(e)** Assuming that all the individuals play their dominant strategies, find out who is pivotal and what tax (if any) each individual has to pay?

**(f)** Show that if every other individual reports his/her true benefit, then it is best for Individual 1 to also report his/her true benefit.

## 1.E.5. Exercises for Section 1.5: Iterated deletion procedures

The answers to the following exercises are in Appendix S at the end of this chapter.

**Exercise 1.10. (a)** Apply the IDSDS procedure (Iterated Deletion of Strictly Dominated Strategies) to the game of Part (b) of Exercise 1.1.

**(b)** Apply the IDWDS procedure (Iterated Deletion of Weakly Dominated Strategies) to the game of Part (b) of Exercise 1.1.

**Exercise 1.11**. Apply the IDSDS procedure to the following game. Is there a strict iterated dominant-strategy equilibrium?

| | | Player 2 | | |
|---|---|---|---|---|
| | | $d$ | $e$ | $f$ |
| Player 1 | $a$ | 8 , 6 | 0 , 9 | 3 , 8 |
| | $b$ | 3 , 2 | 2 , 1 | 4 , 3 |
| | $c$ | 2 , 8 | 1 , 5 | 3 , 1 |

**Exercise 1.12**. Consider the following game. There is a club with three members: Ann, Bob and Carla. They have to choose which of the three is going to be president next year. Currently, Ann is the president. Each member is both a candidate and a voter. Voting is as follows: each member votes for one candidate (voting for oneself is allowed); if two or more people vote for the same candidate then that person is chosen as the next president; if there is complete disagreement, in the sense that there is exactly one vote for each candidate, then the person for whom Ann voted is selected as the next president.

**(a)** Represent this voting procedure as a game frame, indicating inside each cell of each table which candidate is elected.

**(b)** Assume that the players' preferences are as follows: $Ann \succ_{Ann} Carla \succ_{Ann} Bob$, $Carla \succ_{Bob} Bob \succ_{Bob} Ann$, $Bob \succ_{Carla} Ann \succ_{Carla} Carla$.
Using utility values 0, 1 and 2, convert the game frame into a game.

**(c)** Apply the IDWDS to the game of part (b). Is there an iterated weak dominant-strategy equilibrium?

**(d)** Does the extra power given to Ann (in the form of tie-breaking in case of complete disagreement) benefit Ann?

**Exercise 1.13**. Consider the following game:

| | | Player 2 | | | | | |
|---|---|---|---|---|---|---|---|
| | | $D$ | | $E$ | | $F$ | |
| Player 1 | $a$ | 2 | 3 | 2 | 2 | 3 | 1 |
| | $b$ | 2 | 0 | 3 | 1 | 1 | 0 |
| | $c$ | 1 | 4 | 2 | 0 | 0 | 4 |

**(a)** Apply the IDSDS procedure to it. Is there an iterated strict dominant-strategy equilibrium?

**(b)** Apply the IDWDS procedure to it. Is there an iterated weak dominant-strategy equilibrium?

## 1.E.6. Exercises for Section 1.6: Nash equilibrium.

The answers to the following exercises are in Appendix S at the end of this chapter.

**Exercise 1.14**. Find the Nash equilibria of the game of Exercise 1.2.

**Exercise 1.15**. Find the Nash equilibria of the games of Exercise 1.3 (b) and (c).

**Exercise 1.16**. Find the Nash equilibria of the games of Exercise 1.4 (b).

**Exercise 1.17**. Find the Nash equilibria of the game of Exercise 1.6.

**Exercise 1.18**. Find the Nash equilibria of the game of Exercise 1.7.

**Exercise 1.19**. Find one Nash equilibrium of the game of Exercise 1.8 for the case where
$n = 3$ and $v_1 > v_2 > v_3 > 0$.

**Exercise 1.20**. Find the Nash equilibria of the game of part (b) of Exercise 1.12.

**Exercise 1.21**. Find the Nash equilibria of the game of Exercise 1.13.

## 1.E.7. Exercises for Section 1.6: Games with infinite strategy sets.

**Exercise 1.22**. Consider a simultaneous $n$-player game where each player $i$ chooses an effort level $a_i \in [0,1]$. The payoff to player $i$ is given by $\pi_i(a_1,...,a_n) = 4\min\{a_1,...,a_n\} - 2a_i$ (efforts are complementary and each player's cost per unit of effort is 2).

**(a)** Find all the Nash equilibria and prove that they are indeed Nash equilibria.

**(b)** Are any of the Nash equilibria Pareto efficient?

**(c)** Find a Nash equilibrium where each player gets a payoff of 1.

◊◊◊◊◊◊◊◊◊◊◊◊

**Exercise 1.23: Challenging Question**. The Mondevil Corporation operates a chemical plant, which is located on the banks of the Sacramento river. Downstream from the chemical plant is a group of fisheries. The Mondevil plant emits by-products that pollute the river, causing harm to the fisheries. The profit Mondevil obtains from operating the chemical plant is $\$\Pi > 0$. The harm inflicted on the fisheries due to water pollution is equal to $\$L > 0$ of lost profit [without pollution the fisheries' profit is $\$A$, while with pollution it is $\$(A - L)$]. Suppose that the fisheries collectively sue the Mondevil Corporation. It is easily verified in court that Mondevil's plant pollutes the river. However, the values of $\Pi$ and $L$ cannot be verified by the court, although they are commonly known to the litigants. Suppose that the court requires the Mondevil attorney (Player 1) and the fisheries' attorney (Player 2) to play the following litigation game. Player 1 is asked to announce a number $x \geq 0$, which the court interprets as a claim about the plant's profits. Player 2 is asked to announce a number $y \geq 0$, which the court interprets as the fisheries' claim about their profit *loss*. The announcements are made simultaneously and independently. Then the court uses *Posner's nuisance rule* to make its decision (R. Posner, *Economic analysis of Law*, 9th edition, 1997). According to the rule, if $y > x$, then Mondevil must shut down its chemical plant. If $x \geq y$ , then the court allows Mondevil to operate the plant, but the court also requires Mondevil to pay the fisheries the amount $y$. Note that the court cannot force the attorneys to tell the truth: in fact, it would not be able to tell whether or not the lawyers were reporting truthfully. Assume that the attorneys want to maximize the payoff (profits) of their clients.

**(a)** Represent this situation as a normal-form game by describing the strategy set of each player and the payoff functions.

**(b)** Is it a dominant strategy for the Mondevil attorney to make a truthful announcement (i.e. to choose $x = \Pi$)? [Prove your claim.]

**(c)** Is it a dominant strategy for the fisheries' attorney to make a truthful announcement (i.e. to choose $y = L$)? [Prove your claim.]

**(d)** For the case where $\Pi > L$ (recall that $\Pi$ and $L$ denote the *true* amounts), find **all** the Nash equilibria of the litigation game. [Prove that what you claim to be Nash equilibria are indeed Nash equilibria and that there are no other Nash equilibria.]

**(e)** For the case where $\Pi < L$ (recall that $\Pi$ and $L$ denote the *true* amounts), find **all** the Nash equilibria of the litigation game. [Prove that what you claim to be Nash equilibria are indeed Nash equilibria and that there are no other Nash equilibria.]

**(f)** Does the court rule give rise to a Pareto efficient outcome? [Assume that the players end up playing a Nash equilibrium.]

# Appendix 1.S: Solutions to exercises

**Exercise 1.1.** **(a)** $I = \{1,2\}$, $S_1 = \{2,4,6\}$, $S_2 = \{1,3,5\}$, $O = \{M, I, J\}$ (where $M$ stands for 'Mexican', $I$ for 'Italian' and $J$ for 'Japanese); the set of strategy profiles is $S = \{(2,1),(2,3),(2,5),(4,1),(4,3),(4,5),(6,1),(6,3),(6,5)\}$; the outcome function is $f(2,1) = f(2,3) = f(4,1) = M$, $f(2,5) = f(4,3) = f(6,1) = I$ and $f(4,5) = f(6,3) = f(6,5) = J$. The representation as a table is as follows:

Player 2 (Bob)

| Player 1 (Antonia) | 1 | 3 | 5 |
|---|---|---|---|
| 2 | M | M | I |
| 4 | M | I | J |
| 6 | I | J | J |

**(b)** Using values 1, 2 and 3, the utility functions are as follows, where $U_1$ is the utility function of Player 1 (Antonia) and $U_2$ is the utility function of Player 2 (Bob): $\begin{pmatrix} & M & I & J \\ U_1: & 3 & 2 & 1 \\ U_2: & 2 & 3 & 1 \end{pmatrix}$. The reduced-form game is as follows:

Player 2 (Bob)

| Player 1 (Antonia) | 1 | | 3 | | 5 | |
|---|---|---|---|---|---|---|
| 2 | 3 | 2 | 3 | 2 | 2 | 3 |
| 4 | 3 | 2 | 2 | 3 | 1 | 1 |
| 6 | 2 | 3 | 1 | 1 | 1 | 1 |

**Exercise 1.2.**

**(a)** The game-frame is as follows:

| | | Player 2 (Bob): 0 | | 1 | | 2 | |
|---|---|---|---|---|---|---|---|
| Player 1 (Antonia) | 2 | Antonia gets nothing | Bob gets $2 | Antonia gets nothing | Bob gets $3 | Antonia gets nothing | Bob gets $2 |
| | 4 | Antonia gets nothing | Bob gets $2 | Antonia gets $5 | Bob gets $5 | Antonia gets $4 | Bob gets $2 |
| | 6 | Antonia gets $4 | Bob gets $2 | Antonia gets $3 | Bob gets $7 | Antonia gets $2 | Bob gets nothing |

**(b)** When the outcomes are sums of money and Player $i$ is selfish and greedy then we can take the following as $i$'s utility function: $U_i(\$x) = x$ (other utility functions would do too: the only requirement is that the utility of a larger sum of money is larger than the utility of a smaller sum of money). Thus the reduced-form game is as follows:

| | | Player 2 (Bob): 0 | | 1 | | 2 | |
|---|---|---|---|---|---|---|---|
| Player 1 (Antonia) | 2 | 0 | 2 | 0 | 3 | 0 | 2 |
| | 4 | 0 | 2 | 5 | 5 | 4 | 2 |
| | 6 | 4 | 2 | 3 | 7 | 2 | 0 |

**Exercise 1.3.**

**(a)** The game-frame is as follows:

P = Press

| | | Bob: P | not P |
|---|---|---|---|
| Alice | P | Alice wins | Bob wins |
| | not P | Bob wins | Alice wins |

Charlie: P

| | | Bob: P | not P |
|---|---|---|---|
| Alice | P | Bob wins | Alice wins |
| | not P | Alice wins | Charlie wins |

Charlie: not P

**(b)** The reduced-form game is as follows:

P = Press

| Alice \ Bob | P | | | not P | | |
|---|---|---|---|---|---|---|
| P | 1 | 0 | 0 | 0 | 1 | 0 |
| not P | 0 | 1 | 0 | 1 | 0 | 0 |

Charlie: P

| Alice \ Bob | P | | | not P | | |
|---|---|---|---|---|---|---|
| P | 0 | 1 | 0 | 1 | 0 | 0 |
| not P | 1 | 0 | 0 | 0 | 0 | 1 |

Charlie: not P

**(c)** The reduced-form game is as follows. For Alice we chose 1 and 0 as utilities, but one could also use 2 and 1 or 2 and 0.

P = Press

| Alice \ Bob | P | | | not P | | |
|---|---|---|---|---|---|---|
| P | 1 | 0 | 0 | 0 | 2 | 1 |
| not P | 0 | 2 | 1 | 1 | 0 | 0 |

Charlie: P

| Alice \ Bob | P | | | not P | | |
|---|---|---|---|---|---|---|
| P | 0 | 2 | 1 | 1 | 0 | 0 |
| not P | 1 | 0 | 0 | 0 | 1 | 2 |

Charlie: not P

**Exercise 1.4.** (**a**) The tables are as follows:

**Distributed money**

| | | Player 2 | | | | | |
|---|---|---|---|---|---|---|---|
| | | 0 | | 300 | | 600 | |
| Player 1 | 0 | 0 | 0 | 225 | 225 | 450 | 450 |
| | 300 | 225 | 225 | 450 | 450 | 675 | 675 |
| | 600 | 450 | 450 | 675 | 675 | 900 | 900 |

**Net amounts**

| | | Player 2 | | | | | |
|---|---|---|---|---|---|---|---|
| | | 0 | | 300 | | 600 | |
| Player 1 | 0 | 0 | 0 | 225 | –75 | 450 | –150 |
| | 300 | –75 | 225 | 150 | 150 | 375 | 75 |
| | 600 | –150 | 450 | 75 | 375 | 300 | 300 |

(**b**) For Player 1 we can take as his payoff the total money lost by the referee and for Player 2 her own net gain:

| | | Player 2 | | | | | |
|---|---|---|---|---|---|---|---|
| | | 0 | | 300 | | 600 | |
| Player 1 | 0 | 0 | 0 | 150 | –75 | 300 | –150 |
| | 300 | 150 | 225 | 300 | 150 | 450 | 75 |
| | 600 | 300 | 450 | 450 | 375 | 600 | 300 |

(**c**) For Player 1 contributing \$600 is a strictly dominant strategy and for Player 2 contribution \$0 is a strictly dominant strategy. Thus (\$600,\$0) is the strict dominant-strategy equilibrium.

**Exercise 1.5.** The game under consideration is the following:

| | | Player 2 (Bob) | | | | | |
|---|---|---|---|---|---|---|---|
| | | 1 | | 3 | | 5 | |
| Player 1 (Antonia) | 2 | 3 | 2 | 3 | 2 | 2 | 3 |
| | 4 | 3 | 2 | 2 | 3 | 1 | 1 |
| | 6 | 2 | 3 | 1 | 1 | 1 | 1 |

(**a**) For Player 1, 6 is strictly dominated by 2. There is no other strategy which is strictly dominated. Player 2 does not have any strictly dominated strategies.

**(b)** For Player 1, 6 is weakly dominated by 4 (and also by 2, since strict dominance implies weak dominance); 4 is weakly dominated by 2. Player 2 does not have any weakly dominated strategies.

**Exercise 1.6.** **(a)** The game is as follows:

Player 3: 0

| Player 1 \ Player 2 | 0 | 3 | 6 |
|---|---|---|---|
| 0 | 0 0 0 | 2 –1 2 | 4 –2 4 |
| 3 | –1 2 2 | 1 1 4 | 3 0 6 |
| 6 | –2 4 4 | 0 3 6 | 2 2 8 |

Player 3: 3

| Player 1 \ Player 2 | 0 | 3 | 6 |
|---|---|---|---|
| 0 | 2 2 –1 | 4 1 1 | 6 0 3 |
| 3 | 1 4 1 | 3 3 3 | 5 2 5 |
| 6 | 0 6 3 | 2 5 5 | 4 4 7 |

Player 3: 6

| Player 1 \ Player 2 | 0 | 3 | 6 |
|---|---|---|---|
| 0 | 4 4 –2 | 6 3 0 | 8 2 2 |
| 3 | 3 6 0 | 5 5 2 | 7 4 4 |
| 6 | 2 8 2 | 4 7 4 | 6 6 6 |

**(b)** For Player 1, 3 strictly dominates 6, 0 strictly dominates 6, 0 strictly dominates 3 (the same is true for every player).

**(c)** The strict dominant-strategy equilibrium is (0,0,0) (everybody contributes 0).

**Exercise 1.7.** The game is as follows:

| Player 1 (value \$30M) | \$10M | \$20M | \$30M | \$40M | \$50M |
|---|---|---|---|---|---|
| \$10M | 0 , 40 | 0 , 40 | 0 , 40 | 0 , 40 | 0 , 40 |
| \$20M | 20 , 0 | 0 , 30 | 0 , 30 | 0 , 30 | 0 , 30 |
| \$30M | 20 , 0 | 10 , 0 | 0 , 20 | 0 , 20 | 0 , 20 |
| \$40M | 20 , 0 | 10 , 0 | 0 , 0 | 0 , 10 | 0 , 10 |
| \$50M | 20 , 0 | 10 , 0 | 0 , 0 | –10 , 0 | 0 , 0 |

**Exercise 1.8.** No. Suppose, by contradiction, that $\hat{b}_1$ is a weakly dominant strategy for Player 1. It cannot be that $\hat{b}_1 > v_1$, because when $b_2 = b_3 = \hat{b}_1$ Player 1 wins and pays $\hat{b}_1$, thereby obtaining a payoff of $v_1 - \hat{b}_1 < 0$, whereas bidding 0 would give him a payoff of 0. It cannot be that $\hat{b}_1 = v_1$ because when $b_2 > \hat{b}_1$ and $b_3 < v_1$ the auction is won by Player 2 and Player 1 gets a payoff of 0, while a bid of Player 1 greater than $b_2$ would make him the winner with a payoff of $v_1 - b_3 > 0$. Similarly, it cannot be that $\hat{b}_1 < v_1$ because when $b_2 > \hat{b}_1$ and $b_3 < v_1$ then the auction is won by Player 2 and Player 1 gets a payoff of 0, while a bid greater than $b_2$ would make him the winner with a payoff of $v_1 - b_3 > 0$.

**Exercise 1.9.** (**a**) Since $\sum_{i=1}^{5} v_i = 85 < \sum_{i=1}^{5} c_i = 100$ the Pareto efficient decision is not to build the park.

(**b**) Since $\sum_{i=1}^{5} w_i = 120 > \sum_{i=1}^{5} c_i = 100$ the park would be built.

(**c**) Individuals 1 and 3 are pivotal and each of them has to pay a tax of $20. The other individuals are not pivotal and thus are not taxed.

(**d**) For each individual *i* it is a dominant strategy to report $v_i$ and thus, by part (a), the decision will be the Pareto efficient one, namely not to build the park.

(**e**) When every individual reports truthfully, Individuals 4 and 5 are pivotal and Individual 4 has to pay a tax of $25, while individual 5 has to pay a tax of $10. The others are not pivotal and do not have to pay a tax.

(**f**) Assume that all the other individuals report truthfully; then if Individual 1 reports truthfully, he is not pivotal, the project is not carried out and his utility is $\bar{m}_1$. Any other $w_1$ that leads to the same decision (not to build the park) gives him the same utility. If, on the other hand, he chooses a $w_1$ that leads to a decision to build the park, then Individual 1 will become pivotal and will have to pay a tax of $45 with a utility of $\bar{m}_1 + v_1 - 45 = \bar{m}_1 + 30 - 45 = \bar{m}_1 - 15$, so that he would be worse off relative to reporting truthfully.

**Exercise 1.10.** The game under consideration is the following:

| | | Player 2 (Bob) | | |
|---|---|---|---|---|
| | | 1 | 3 | 5 |
| Player 1 (Antonia) | 2 | 3, 2 | 3, 2 | 2, 3 |
| | 4 | 3, 2 | 2, 3 | 1, 1 |
| | 6 | 2, 3 | 1, 1 | 1, 1 |

(**a**) The first step of the procedure eliminates 6 for Player 1. After this step the procedure stops and thus the output is

| | | Player 2 (Bob) | | |
|---|---|---|---|---|
| | | 1 | 3 | 5 |
| Player 1 (Antonia) | 2 | 3  2 | 3  2 | 2  3 |
| | 4 | 3  2 | 2  3 | 1  1 |

**(b)** The first step of the procedure eliminates 4 and 6 for Player 1 and nothing for Player 2. The second step of the procedure eliminates 1 and 3 for Player 2. Thus the output is the strategy profile (2,5), which constitutes the iterated weak dominant-strategy equilibrium of this game.

**Exercise 1.11.** The game under consideration is

| | | Player 2 | | |
|---|---|---|---|---|
| | | *d* | *e* | *f* |
| Player 1 | *a* | 8 , 6 | 0 , 9 | 3 , 8 |
| | *b* | 3 , 2 | 2 , 1 | 4 , 3 |
| | *c* | 2 , 8 | 1 , 5 | 3 , 1 |

In this game *c* is strictly dominated by *b*; after deleting *c*, *d* becomes strictly dominated by *f*; after deleting *d*, *a* becomes strictly dominated by *b*; after deleting *a*, *e* becomes strictly dominated by *f*; deletion of *e* leads to only one strategy profile, namely (*b*,*f*). Thus (*b*,*f*) is the iterated strict dominant-strategy equilibrium.

**Exercise 1.12.** **(a)** The game-frame is as follows:

Carla votes for A

| ANN \ BOB | A | B | C |
|---|---|---|---|
| A | A | A | A |
| B | A | B | B |
| C | A | C | C |

Carla votes for B

| ANN \ BOB | A | B | C |
|---|---|---|---|
| A | A | B | A |
| B | B | B | B |
| C | C | B | C |

Carla votes for C

| ANN \ BOB | A | B | C |
|---|---|---|---|
| A | A | A | C |
| B | B | B | C |
| C | C | C | C |

**(b)** The game is as follows:

BOB

| ANN | A | B | C |
|---|---|---|---|
| A | 2,0,1 | 2,0,1 | 2,0,1 |
| B | 2,0,1 | 0,1,2 | 0,1,2 |
| C | 2,0,1 | 1,2,0 | 1,2,0 |

Carla votes for A

BOB

| ANN | A | B | C |
|---|---|---|---|
| A | 2,0,1 | 0,1,2 | 2,0,1 |
| B | 0,1,2 | 0,1,2 | 0,1,2 |
| C | 1,2,0 | 0,1,2 | 1,2,0 |

Carla votes for B

BOB

| ANN | A | B | C |
|---|---|---|---|
| A | 2,0,1 | 2,0,1 | 1,2,0 |
| B | 0,1,2 | 0,1,2 | 1,2,0 |
| C | 1,2,0 | 1,2,0 | 1,2,0 |

Carla votes for C

**(c)** For Ann, both *B* and *C* are weakly dominated by *A*, while the other two players do not have any dominated strategies. Thus in the first step of the IDWDS we delete *B* and *C* for Ann. Hence the game reduces to:

CARLA

| BOB | A | B | C |
|---|---|---|---|
| A | 0,1 | 0,1 | 0,1 |
| B | 0,1 | 1,2 | 0,1 |
| C | 0,1 | 0,1 | 2,0 |

Carla votes for A

In this game, for Carla, *B* weakly dominates both *A* and *C* and for Bob *A* is weakly dominated by *B* (and also by *C*). Thus in the second step of the IDWDS we delete *A* and *B* for Carla and *A* for Bob. In the third and final step we delete *C* for Bob. Thus we are left with a unique strategy profile, namely (*A*,*B*,*B*), that is, Ann votes for herself and Bob and Carla vote for Bob. This is the iterated weak dominant-strategy equilibrium.

**(d)** The elected candidate is Bob, who is Ann's least favorite; thus the extra power given to Ann (tie breaking in case of total disagreement) turns out to be detrimental for Ann!

**Exercise 1.13.** The game under consideration is

Player 2

| Player 1 | | *D* | | *E* | | *F* | |
|---|---|---|---|---|---|---|---|
| | *a* | 2 | 3 | 2 | 2 | 3 | 1 |
| | *b* | 2 | 0 | 3 | 1 | 1 | 0 |
| | *c* | 1 | 4 | 2 | 0 | 0 | 4 |

**(a)** The output of the IDSDS is as follows (first delete *c* and then *F*):

Player 2

| Player 1 | | *D* | | *E* | |
|---|---|---|---|---|---|
| | *a* | 2 | 3 | 2 | 2 |
| | *b* | 2 | 0 | 3 | 1 |

Thus there is no iterated strict dominant-strategy equilibrium.

**(b)** The output of the IDWDS is (*b*,*E*) (in the first step delete *c* and *F*, the latter because it is weakly dominated by *D*; in the second step delete *a* and in the third step delete *D*). Thus (*b*,*E*) is the iterated weak dominant-strategy equilibrium.

**Exercise 1.14.** The game under consideration is

Player 2 (Bob)

| Player 1 (Antonia) | | 0 | | 1 | | 2 | |
|---|---|---|---|---|---|---|---|
| | 2 | 0 | 2 | 0 | 3 | 0 | 2 |
| | 4 | 0 | 2 | 5 | 5 | 4 | 2 |
| | 6 | 4 | 2 | 3 | 7 | 2 | 0 |

There is only one Nash equilibrium, namely (4,1) with payoffs (5,5).

**Exercise 1.15.** The game of part (b) of Exercise 1.3 is as follows:

P = Press

Charlie: P

| | | Bob: P | | | Bob: not P | | |
|---|---|---|---|---|---|---|---|
| Alice | P | 1 | 0 | 0 | 0 | 1 | 0 |
| Alice | not P | 0 | 1 | 0 | 1 | 0 | 0 |

Charlie: not P

| | | Bob: P | | | Bob: not P | | |
|---|---|---|---|---|---|---|---|
| Alice | P | 0 | 1 | 0 | 1 | 0 | 0 |
| Alice | not P | 1 | 0 | 0 | 0 | 0 | 1 |

This game has only one Nash equilibrium, namely (not P, P, not P).

The game of part (c) of Exercise 1.3 is as follows:

P = Press

Charlie: P

| | | Bob: P | | | Bob: not P | | |
|---|---|---|---|---|---|---|---|
| Alice | P | 1 | 0 | 0 | 0 | 2 | 1 |
| Alice | not P | 0 | 2 | 1 | 1 | 0 | 0 |

Charlie: not P

| | | Bob: P | | | Bob: not P | | |
|---|---|---|---|---|---|---|---|
| Alice | P | 0 | 2 | 1 | 1 | 0 | 0 |
| Alice | not P | 1 | 0 | 0 | 0 | 1 | 2 |

This game does not have any Nash equilibria.

**Exercise 1.16.** The game under consideration is

| | | Player 2: 0 | | Player 2: 300 | | Player 2: 600 | |
|---|---|---|---|---|---|---|---|
| Player 1 | 0 | 0 | 0 | 150 | −75 | 300 | −150 |
| Player 1 | 300 | 150 | 225 | 300 | 150 | 450 | 75 |
| Player 1 | 600 | 300 | 450 | 450 | 375 | 600 | 300 |

This game has only one Nash equilibrium, namely (600,0).

**Exercise 1.17.** The game under consideration is:

Player 3: 0

| Player 1 \ Player 2 | 0 | 3 | 6 |
|---|---|---|---|
| 0 | 0 0 0 | 2 −1 2 | 4 −2 4 |
| 3 | −1 2 2 | 1 1 4 | 3 0 6 |
| 6 | −2 4 4 | 0 3 6 | 2 2 8 |

Player 3: 3

| Player 1 \ Player 2 | 0 | 3 | 6 |
|---|---|---|---|
| 0 | 2 2 −1 | 4 1 1 | 6 0 3 |
| 3 | 1 4 1 | 3 3 3 | 5 2 5 |
| 6 | 0 6 3 | 2 5 5 | 4 4 7 |

Player 3: 6

| Player 1 \ Player 2 | 0 | 3 | 6 |
|---|---|---|---|
| 0 | 4 4 −2 | 6 3 0 | 8 2 2 |
| 3 | 3 6 0 | 5 5 2 | 7 4 4 |
| 6 | 2 8 2 | 4 7 4 | 6 6 6 |

This game has only one Nash equilibrium, namely (0,0,0).

**Exercise 1.18.** The game under consideration is as follows:

| Player 1 (value \$30M) | \$10M | \$20M | \$30M | \$40M | \$50M |
|---|---|---|---|---|---|
| \$10M | 0 , 40 | 0 , 40 | 0 , 40 | 0 , 40 | 0 , 40 |
| \$20M | 20 , 0 | 0 , 30 | 0 , 30 | 0 , 30 | 0 , 30 |
| \$30M | 20 , 0 | 10 , 0 | 0 , 20 | 0 , 20 | 0 , 20 |
| \$40M | 20 , 0 | 10 , 0 | 0 , 0 | 0 , 10 | 0 , 10 |
| \$50M | 20 , 0 | 10 , 0 | 0 , 0 | −10 , 0 | 0 , 0 |

This game has 15 Nash equilibria: (10,30), (10,40), (10,50), (20,30), (20,40), (20,50), (30,30), (30,40), (30,50), (40,40), (40,50), (50,10), (50,20), (50,30), (50,50).

**Exercise 1.19.** A Nash equilibrium is $b_1 = b_2 = b_3 = v_1$ (with payoffs (0,0,0)). Convince yourself that this is indeed a Nash equilibrium. There are many more Nash equilibria: for example, any triple $(b_1, b_2, b_3)$ with $b_2 = b_3 = v_1$ and $b_1 > v_1$ is a Nash equilibrium (with payoffs (0,0,0)) and so is any triple $(b_1, b_2, b_3)$ with $b_2 = b_3 = v_2$ and $b_1 \geq v_2$ (with payoffs $(v_1 - v_2, 0, 0)$).

**Exercise 1.20.** The game under consideration is

| ANN \ BOB | A | B | C |
|---|---|---|---|
| A | 2,0,1 | 2,0,1 | 2,0,1 |
| B | 2,0,1 | 0,1,2 | 0,1,2 |
| C | 2,0,1 | 1,2,0 | 1,2,0 |

Carla votes for A

| ANN \ BOB | A | B | C |
|---|---|---|---|
| A | 2,0,1 | 0,1,2 | 2,0,1 |
| B | 0,1,2 | 0,1,2 | 0,1,2 |
| C | 1,2,0 | 0,1,2 | 1,2,0 |

Carla votes for B

| ANN \ BOB | A | B | C |
|---|---|---|---|
| A | 2,0,1 | 2,0,1 | 1,2,0 |
| B | 0,1,2 | 0,1,2 | 1,2,0 |
| C | 1,2,0 | 1,2,0 | 1,2,0 |

Carla votes for C

There are 5 Nash equilibria: (A,A,A), (B,B,B), (C,C,C), (A,C,A) and (A,B,B).

**Exercise 1.21.** The game under consideration is as follows:

| Player 1 \ Player 2 | *D* | | *E* | | *F* | |
|---|---|---|---|---|---|---|
| *a* | 2 | 3 | 2 | 2 | 3 | 1 |
| *b* | 2 | 0 | 3 | 1 | 1 | 0 |
| *c* | 1 | 4 | 2 | 0 | 0 | 4 |

There are 2 Nash equilibria: $(a,D)$ and $(b,E)$.

**Exercise 1.22.** **(a)** For every $e \in [0,1]$, $(e,e,...,e)$ is a Nash equilibrium. The payoff of Player $i$ is $\pi_i(e,e,...,e) = 2e$; if player $i$ increases her effort to $a > e$ (of course, this can only happen if $e < 1$), then her payoff decreases to $4e - 2a$ and if she decreases her effort to $a < e$ (of course, this can only happen if $e > 0$), then her payoff decreases to $2a$.

There is no Nash equilibrium where two players choose different levels of effort. Proof: suppose there is an equilibrium $(a_1, a_2, ..., a_n)$ where $a_i \neq a_j$ for two players $i$ and $j$. Let $a_{\min} = \min\{a_1, ..., a_n\}$ and let $k$ be a player such that $a_k > a_{\min}$ (such a player exists by our supposition). Then the payoff to player $k$ is $\pi_k = 4a_{\min} - 2a_k$ and if she reduced her effort to $a_{\min}$ her payoff would increase to $2a_{\min}$.

**(b)** Any symmetric equilibrium with $e < 1$ is Pareto inefficient, because all the players would be better off if they simultaneously switched to (1,1,...,1). On the other hand, the symmetric equilibrium (1,1,...,1) is Pareto efficient.

**(c)** The symmetric equilibrium $\left(\frac{1}{2}, \frac{1}{2}, ..., \frac{1}{2}\right)$.

**Exercise 1.23.** **(a)** The strategy sets are $S_1 = S_2 = [0,\infty)$. The payoff functions are as follows:

$$\pi_1(x,y) = \begin{cases} \Pi - y & \text{if } x \geq y \\ 0 & \text{if } y > x \end{cases} \quad \text{and} \quad \pi_2(x,y) = \begin{cases} A - L + y & \text{if } x \geq y \\ A & \text{if } y > x \end{cases}$$

**(b)** Yes, for player 1 choosing $x = \Pi$ is a weakly dominant strategy. **Proof.** Consider an arbitrary $y$. We must show that $x = \Pi$ gives at least a high a payoff against $y$ as any other $x$. Three cases are possible. **Case 1: $y < \Pi$.** In this case $x = \Pi$ or any other $x$ such that $x \geq y$ yields $\pi_1 = \Pi - y > 0$, while any $x < y$ yields $\pi_1 = 0$. **Case 2: $y = \Pi$.** In this case 1's payoff is zero no matter what $x$ he chooses. **Case 3: $y > \Pi$.** In this case $x = \Pi$ or any other $x$ such that $x < y$ yields $\pi_1 = 0$, while any $x \geq y$ yields $\pi_1 = \Pi - y < 0$.

**(c)** No, choosing $y = L$ is not a dominant strategy for Player 2. For example, if $x > L$ then choosing $y = L$ yields $\pi_2 = A$ while choosing a $y$ such that $L < y \leq x$ yields $\pi_2 = A - L + y > A$.

**(d)** **Suppose that $\Pi > L$.** If (x,y) is a Nash equilibrium **with $x \geq y$** then it must be that $y \leq \Pi$ (otherwise Player 1 could increase its payoff by reducing $x$ below $y$) and $y \geq L$ (otherwise Player 2 would be better off by increasing $y$ above $x$). Thus it must be $L \leq y \leq \Pi$, which is possible, given our assumption. However, it cannot be that $x > y$, because Player 2 would be getting a higher payoff by increasing $y$ to $x$. Thus it must be **$x \leq y$**, which (together with our hypothesis that $x \geq y$) implies that $x = y$. Thus the following are Nash equilibria:

all the pairs $(x,y)$ with $L \leq y \leq \Pi$ and $x = y$.

Now consider pairs (x,y) **with $x < y$**. Then it cannot be that $y < \Pi$, because Player 1 could increase its payoff by increasing $x$ to $y$. Thus it must be $y \geq \Pi$ (hence, by our supposition that $\Pi > L$, $y > L$). Furthermore, it must be that $x \leq L$ (otherwise Player 2 could increase its profits by reducing $y$ to (or below) $x$. Thus

(x,y) with $x < y$ is a Nash equilibrium if and only if $x \leq L$ and $y \geq \Pi$.

**(e)** **Suppose that $\Pi < L$.** For the same reasons given above, an equilibrium with $x \geq y$ requires $L \leq y \leq \Pi$. However, this is not possible given that $\Pi < L$. Hence, there is no Nash equilibrium $(x,y)$ with $x \geq y$.

Thus we must restrict attention to pairs (x,y) **with $x < y$**. As explained above, it must be that $y \geq \Pi$ and $x \leq L$. Thus,

$$(x,y) \text{ with } x < y \text{ is a Nash equilibrium if and only if } \Pi \leq y \text{ and } x \leq L.$$

**(f)** Pareto efficiency requires that the chemical plant be shut down if $\Pi < L$ and that it remain operational if $\Pi > L$. Now, when $\Pi < L$ all the equilibria have $x < y$ which leads to shut-down, hence a Pareto efficient outcome. When $\Pi > L$, there are two types of equilibria: one where $x = y$ and the plant remains operational (a Pareto efficient outcome) and the other where $x < y$ in which case the plant shuts down, yielding a Pareto inefficient outcome.

Chapter 2

# Dynamic games with perfect information

## 2.1 Trees, frames and games

Often interactions are not simultaneous but sequential. For example, in the game of Chess the two players, White and Black, take turns moving pieces on the board, having full knowledge of the opponent's (and their own) past moves. Games with sequential interaction are called *dynamic games* or *games in extensive form*. This chapter is devoted to the subclass of dynamic games characterized by *perfect information*, namely the property that, whenever it is her turn to move, a player knows all the preceding moves. Perfect-information games are represented by means of rooted directed trees.

A *rooted directed tree* consists of a set of nodes and directed edges joining them. The root of the tree has no directed edges leading to it (has indegree 0), while every other node has exactly one directed edge leading to it (has indegree 1). There is a unique path (that is, a unique sequence of directed edges) leading from the root to any other node. A node that has no directed edges out of it (has outdegree 0) is called a *terminal node*, while every other node is called a *decision node*. We shall denote the set of nodes by $X$, the set of decision nodes by $D$ and the set of terminal nodes by $Z$. Thus $X = D \cup Z$.

**Definition 2.1.** A *finite extensive form (or frame) with perfect information* consists of the following items.

- A finite rooted directed tree.
- A set of players $I = \{1,\ldots,n\}$ and a function that assigns one player to every decision node.
- A set of actions $A$ and a function that assigns one action to every directed edge, satisfying the restriction that no two edges out of the same node are assigned the same action.
- A set of outcomes $O$ and a function that assigns an outcome to every terminal node.

**Example 2.1.** Amy (Player 1) and Beth (Player 2) have decided to dissolve a business partnership whose assets have been valued at $100,000. The charter of the partnership prescribes that the senior partner, Amy, make an offer concerning the division of the assets to the junior partner, Beth. The junior partner can *Accept*, in which case the proposed division is implemented, or *Reject*, in which case the case goes to litigation. Litigating involves a cost of $20,000 in legal fees for each partner and the typical verdict assigns 60% of the assets to the senior partner and the remaining 40% to the junior partner. Suppose, for simplicity, that there is no uncertainty about the verdict (how to model uncertainty will be discussed in a later chapter). Suppose also that there are only two possible offers that Amy can make: a 50-50 split or a 70-30 split. This situation can be represented as a finite extensive form with perfect information as shown in Figure 2.1. Each outcome is represented as two sums of money: the top one is what Player 1 gets and the bottom one what Player 2 gets.

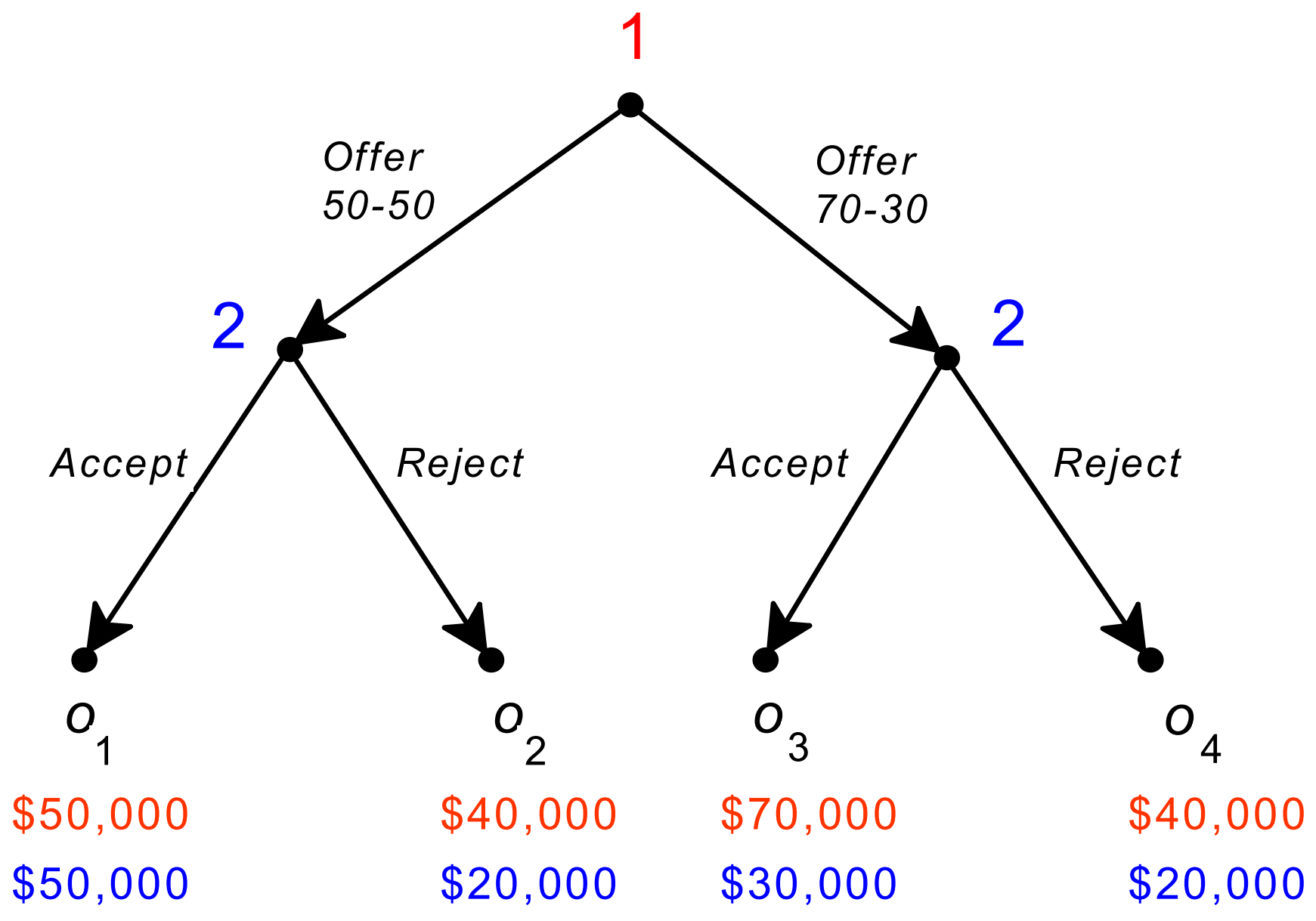

**Figure 2.1**
A perfect-information extensive form representing the situation described in Example 2.1

What should we expect the players to do in the above game? Consider the following reasoning, which is called *backward induction* reasoning, because it starts from the end of the game and proceeds backwards towards the root:

> If Player 2 (the junior partner) is offered a 50-50 split then, if she accepts, she will get \$50,000, while, if she rejects, she will get \$20,000 (the court-assigned 40% minus legal fees in the amount of \$20,000); thus, if rational, she will accept. Similarly, if Player 2 is offered a 70-30 split then, if she accepts, she will get \$30,000, while, if she rejects, she will get \$20,000 (the court-assigned 40% minus legal fees in the amount of \$20,000); thus, if rational, she will accept. Anticipating all of this, Player 1 realizes that, if she offers a 50-50 split then she will end up with \$50,000, while if she offers a 70-30 split then she will end up with \$70,000; thus, if Player 1 is rational and believes that Player 2 is rational, she will offer a 70-30 split and Player 2, being rational, will accept.

The above reasoning suffers from the same flaws as the reasoning described in Chapter 1 (Section 1.1): it is not a valid argument because it is based on an implicit assumption about how Player 2 ranks the outcomes, which may or may not be correct. For example, Player 2 may feel that she worked as hard as her senior partner and the only fair division is a 50-50 split; indeed she may feel so strongly about this that – if offered an unfair 70-30 split – she would be willing to sacrifice \$10,000 in order to "teach a lesson to Player 1"; in other words, she ranks outcome $o_4$ above outcome $o_3$.

Using the terminology introduced in Chapter 1, we say that the situation represented in Figure 2.1 is not a game but a *game-frame*. In order to convert that frame into a game we need to add a ranking of the outcomes for each player.

**Definition 2.2.** A *finite extensive game with perfect information* is a finite extensive form with perfect information together with a ranking $\succsim_i$ of the set of outcome $O$, for every player $i \in I$.

As usual, it is convenient to represent the ranking of Player $i$ by means of an ordinal utility function $U_i : O \to \mathbb{R}$. For example, take the extensive form of Figure 2.1 and assume that Player 1 is selfish and greedy, that is, her ranking is:

| | | |
|---|---|---|
| best | $o_3$ | |
| | $o_1$ | (or, in the alternative notation, $o_3 \succ_1 o_1 \succ_1 o_2 \sim_1 o_4$), |
| worst | $o_2, o_4$ | |

while Player 2 is concerned with fairness and her ranking is:

best $o_1$

$o_2, o_4$ (or, in the alternative notation, $o_1 \succ_2 o_2 \sim_2 o_4 \succ_2 o_3$).

worst $o_3$

Then we can represent the players' preferences using the following utility functions:

| outcome → <br> utility function ↓ | $o_1$ | $o_2$ | $o_3$ | $o_4$ |
|---|---|---|---|---|
| $U_1$ (Player 1) | 2 | 1 | 3 | 1 |
| $U_2$ (Player 2) | 3 | 2 | 1 | 2 |

and replace each outcome in Figure 2.1 with a pair of utilities or payoffs, as shown in Figure 2.2, thereby obtaining one of the many possible games based on the frame of Figure 2.1.

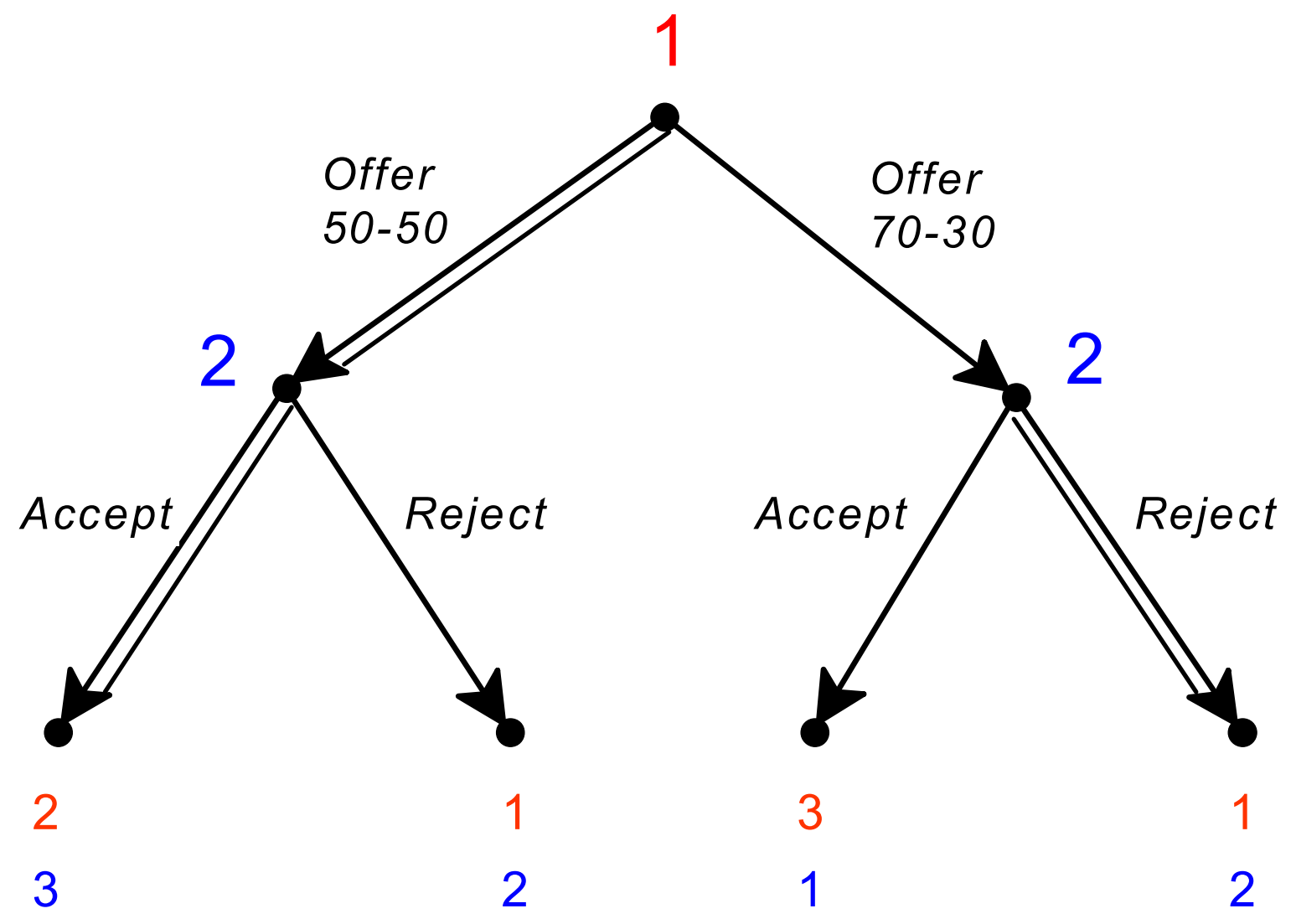

**Figure 2.2**
A perfect-information game based on the extensive form of Figure 2.1.

Now that we do have a game, rather than just a game-frame, we can indeed apply the backward-induction reasoning and conclude that Player 1 will offer a 50-50 split, anticipating that Player 2 would reject the offer of a 70-30 split, and Player 2 will accept Player 1's 50-50 offer. The choices selected by the backward-induction reasoning have been highlighted in Figure 2.2 by doubling the corresponding edges.

This is a good time to test your understanding of the concepts introduced in this section, by going through the exercises in Section 2.E.1 of Appendix 2.E at the end of this chapter.

## 2.2 Backward induction

The backward-induction reasoning mentioned above can be formalized as an algorithm for solving any perfect-information game, as follows. We say that a node is *marked* if a payoff or utility vector is associated with it. Initially all and only the terminal nodes are marked; the following procedure provides a way of marking all the nodes.

**Definition 2.3.** The *backward-induction algorithm* is the following procedure for solving a finite perfect-information game.

1. Select a decision node $x$ whose immediate successors are all marked. Let $i$ be the player who moves at $x$. Select a choice that leads to an immediate successor of $x$ with the highest payoff or utility for Player $i$ (highest among the utilities associated with the immediate successors of $x$). Mark $x$ with the payoff vector associated with the node following the selected choice.
2. Repeat the above step until all the nodes have been marked.

Note that, since the game is finite, the above procedure is well defined. In the initial steps one starts at those decision nodes that are followed only by terminal nodes, call them penultimate nodes. After all the penultimate nodes have been marked, there will be unmarked nodes whose immediate successors are all marked and thus the step can be repeated.

Note also that, in general, at a decision node there may be several choices that maximize the payoff of the player who moves at that node. If that is the case then the procedure requires that *one* such choice be selected. This arbitrary choice may lead to the existence of several backward-induction solutions. For example, consider the game of Figure 2.3.

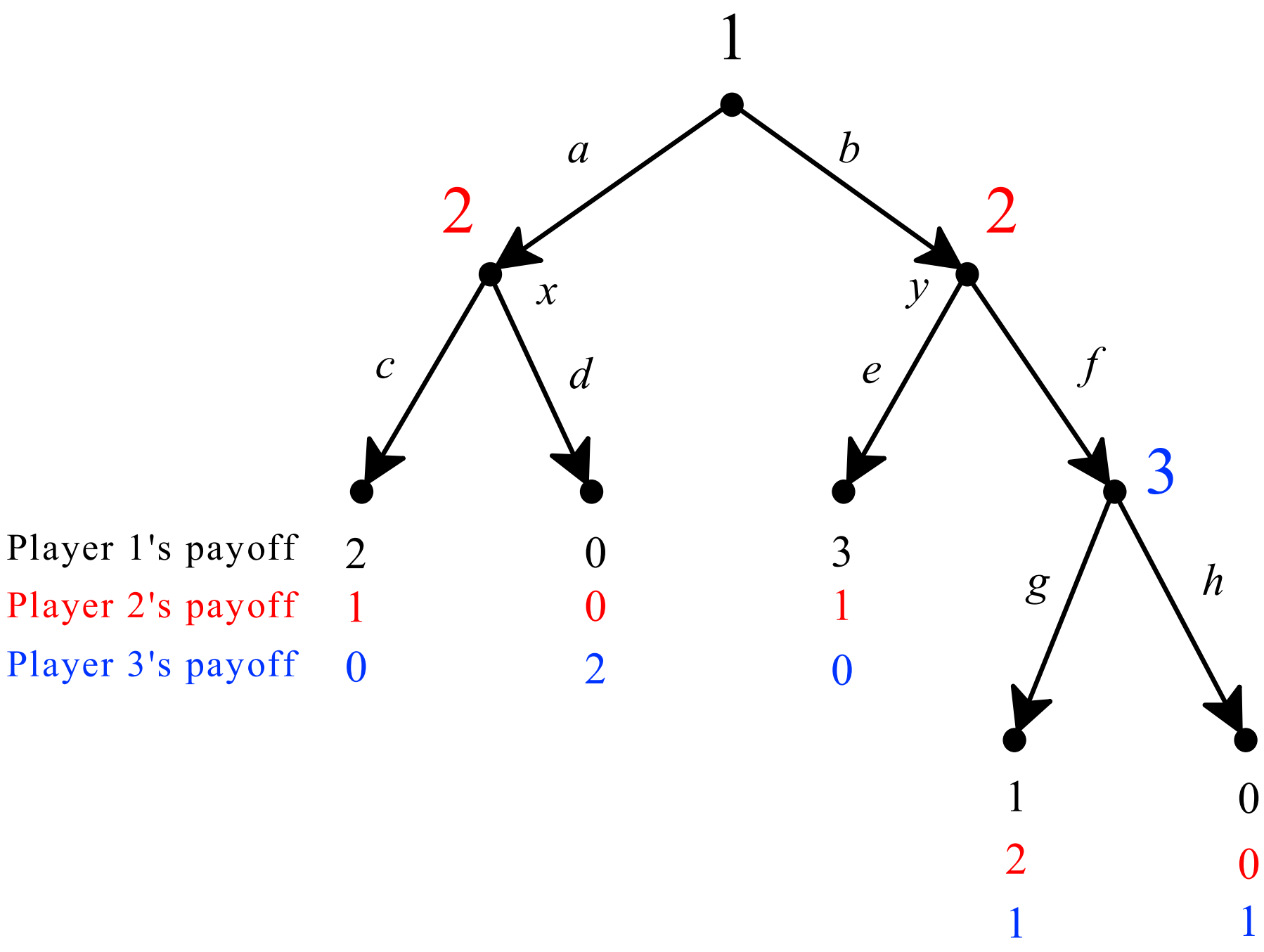

**Figure 2.3**

A perfect-information game with multiple backward-induction solutions.

In this game, starting at node *x* of Player 2 we select choice *c* (since it gives Player 2 a higher payoff than *d*). Then we move on to Player 3's node and we find that both choices there are payoff maximizing for Player 3; thus there are two ways to proceed. In Figure 2.4 we show the steps of the backward-induction algorithm with the selection of choice *g*, while Figure 2.5 shows the steps of the algorithm with the selection of choice *h*. As before, the selected choices are shown by double edges. In Figures 2.4 and 2.5 the marking of nodes is shown explicitly, but later on we will represent the backward-induction solution more succinctly by merely highlighting the selected choices.

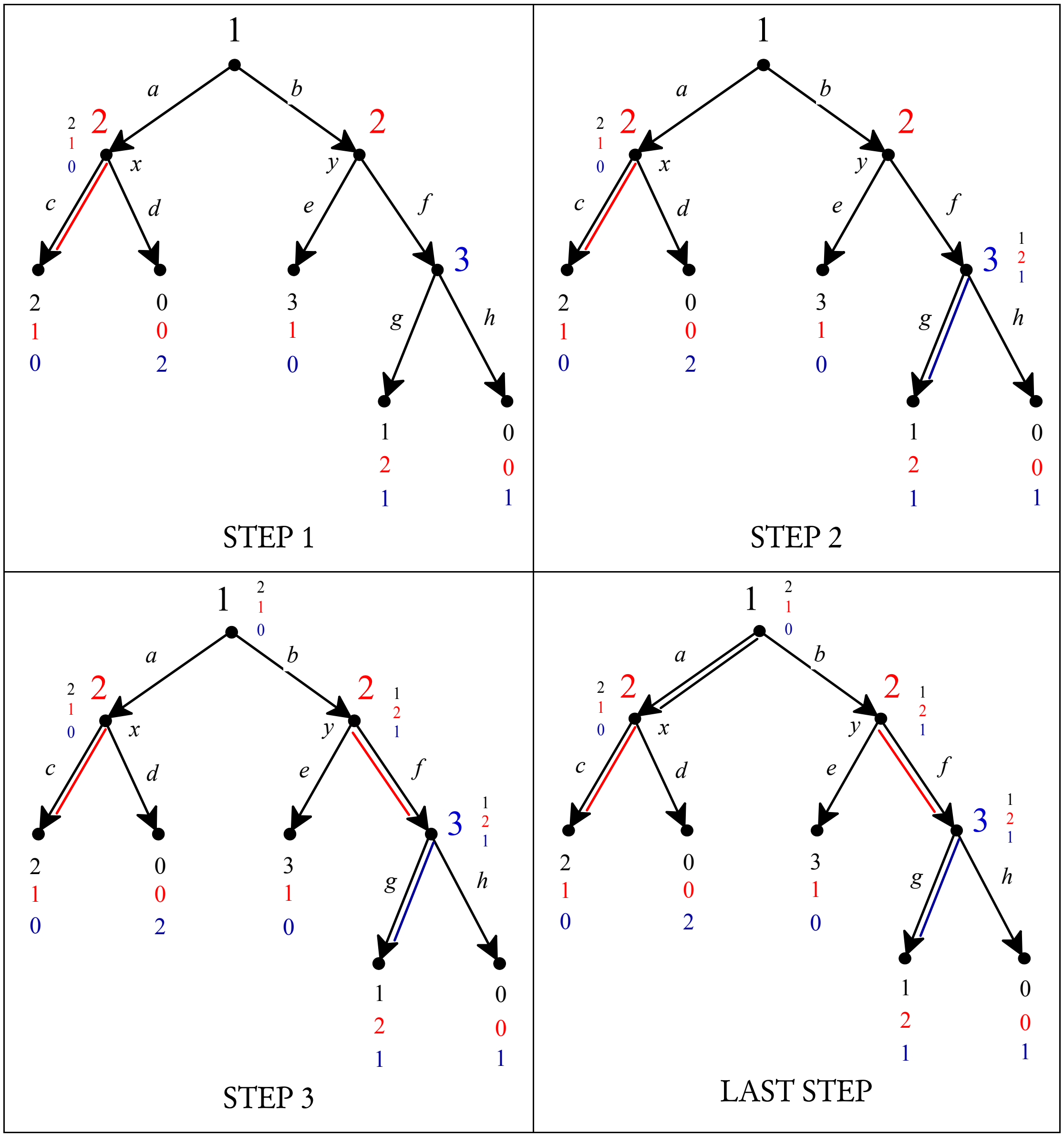

**Figure 2.4**
One possible output of the backward-induction algorithm to the game of Figure 2.3.

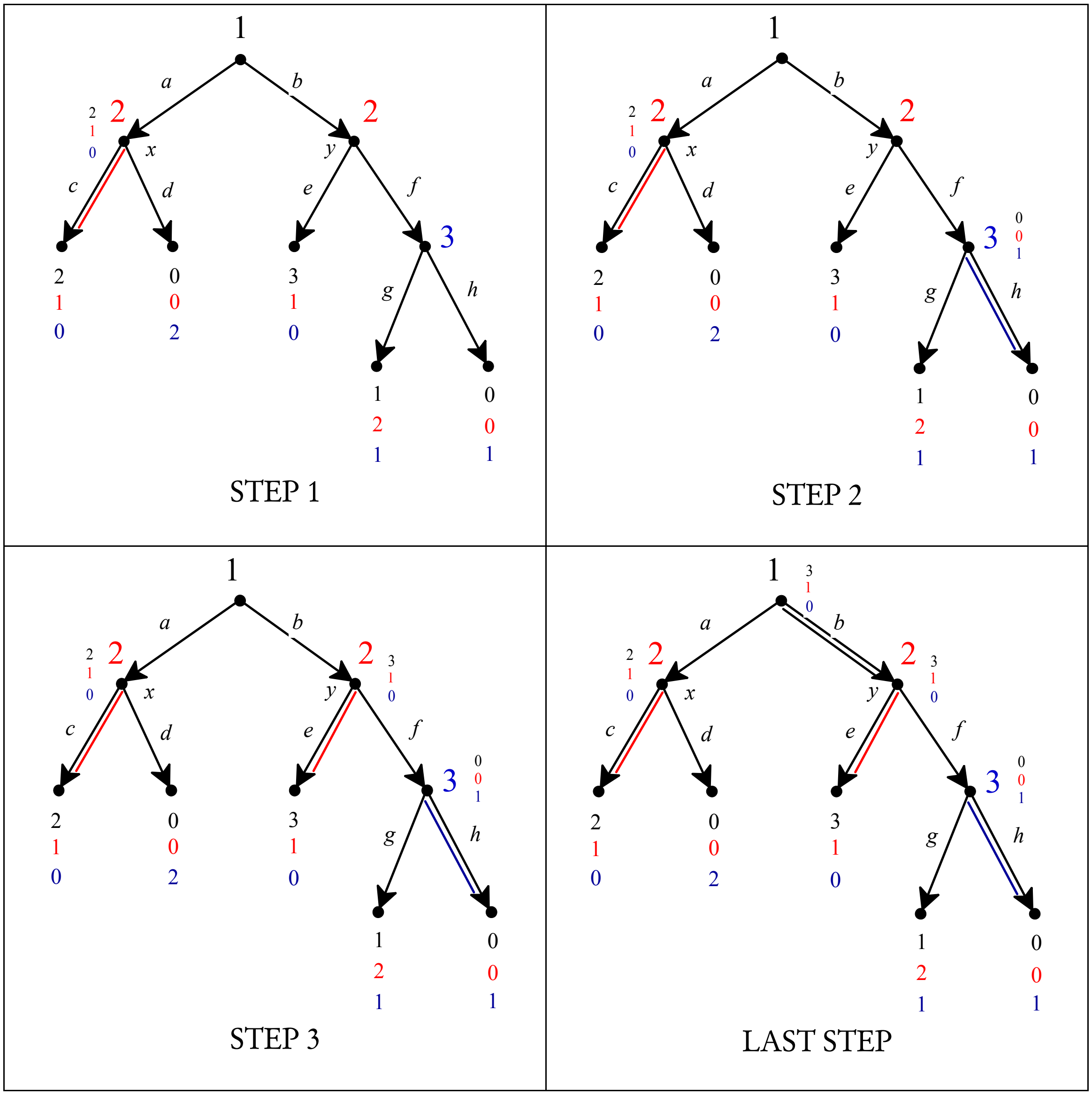

**Figure 2.5**
Another possible output of the backward-induction algorithm to the game of Figure 2.3

How should one define the output of the backward-induction algorithm and the notion of backward-induction solution? What kind of objects are they? Before we answer this question we need to introduce the notion of strategy in a perfect-information game.

This is a good time to test your understanding of the concepts introduced in this section, by going through the exercises in Section 2.E.2 of Appendix 2.E at the end of this chapter.

## 2.3 Strategies in perfect-information games

A strategy for a player in a perfect-information game is a complete, contingent plan on how to play the game. Consider, for example, the following game - which reproduces Figure 2.3 - and let us focus on Player 2.

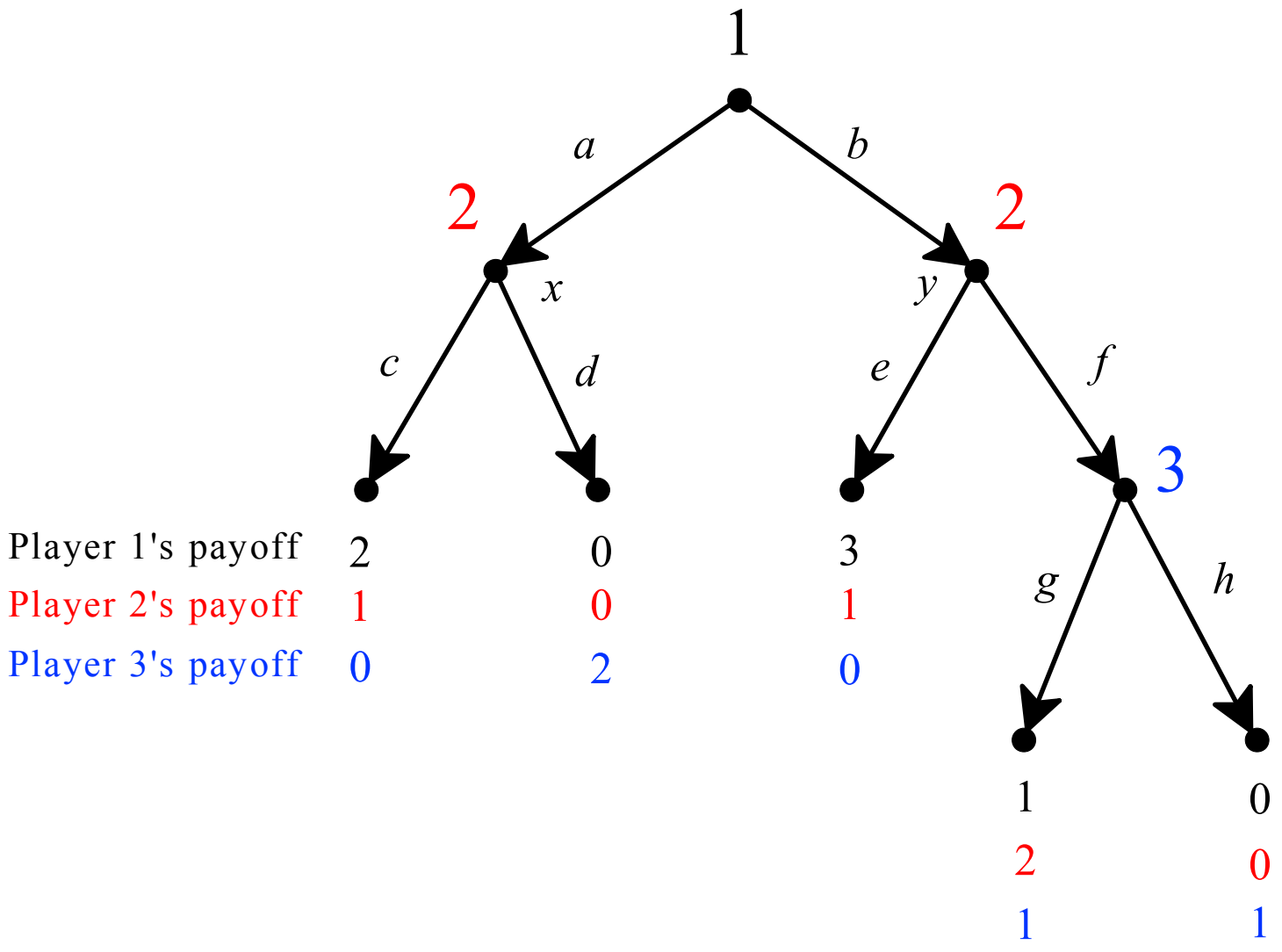

Before the game is played, Player 2 does not know what Player 1 will do and thus a complete plan needs to specify what she will do if Player 1 decides to play $a$ and what she will do if Player 1 decides to play $b$. A possible plan, or strategy, is "if Player 1 chooses $a$ then I will choose $c$ and if Player 1 chooses $b$ then I will choose $e$", which we can denote more succinctly as $(c,e)$. The other possible plans, or strategies, for Player 2 are $(c,f)$, $(d,e)$ and $(d,f)$. The formal definition of strategy is as follows.

**Definition 2.4.** A *strategy* for a player in a perfect-information game is a list of choices, one for each decision node of that player.

For example, suppose that Player 1 has three decision nodes in a given game: at one node she has three possible choices, $a_1, a_2$ and $a_3$, at another node she has two possible choices, $b_1$ and $b_2$, and at the third node she has four possible choices, $c_1, c_2, c_3$ and $c_4$. Then a strategy for Player 1 in that game can be thought of as a way of filling in three blanks: $\left( \underbrace{_____}_{\text{one of } a_1, a_2, a_3}, \underbrace{_____}_{\text{one of } b_1, b_2}, \underbrace{_____}_{\text{one of } c_1, c_2, c_3, c_4} \right)$. Since there are 3 choices for the first blank, 2 for the second and 4 for the third, the total number of possible strategies for Player 1 in this case would be $3 \times 2 \times 4 = 24$. One strategy is $(a_2, b_1, c_1)$, another strategy is $(a_1, b_2, c_4)$, etc.

It should be noted that the notion of strategy involves redundancies. To see this, consider the game of Figure 2.6 below. In this game a possible strategy for Player 1 is (*a*,*g*), which means that Player 1 is planning to choose *a* at the root of the tree and would choose *g* at her other node. But if Player 1 indeed chooses *a*, then her other node will *not* be reached and thus why should Player 1 make a plan on what to do there? One could justify this redundancy in the notion of strategy in a number of ways: (1) Player 1 is so cautious that she wants her plan to cover also the possibility that she might make mistakes in the implementation of parts of her plan (in this case, she allows for the possibility that – despite her intention to play *a* – she might end up playing *b*) or (2) we can think of a strategy as a set of instructions given to a third party on how to play the game on Player 1's behalf, in which case Player 1 might indeed be worried about the possibility of mistakes in the implementation and thus want to cover all contingencies, etc. An alternative justification relies on a different interpretation of the notion of strategy: not as a plan of Player 1 but as a belief in the mind of Player 2 concerning what Player 1 would do. For the moment we will set this issue aside and simply use the notion of strategy as given in Definition 2.4.

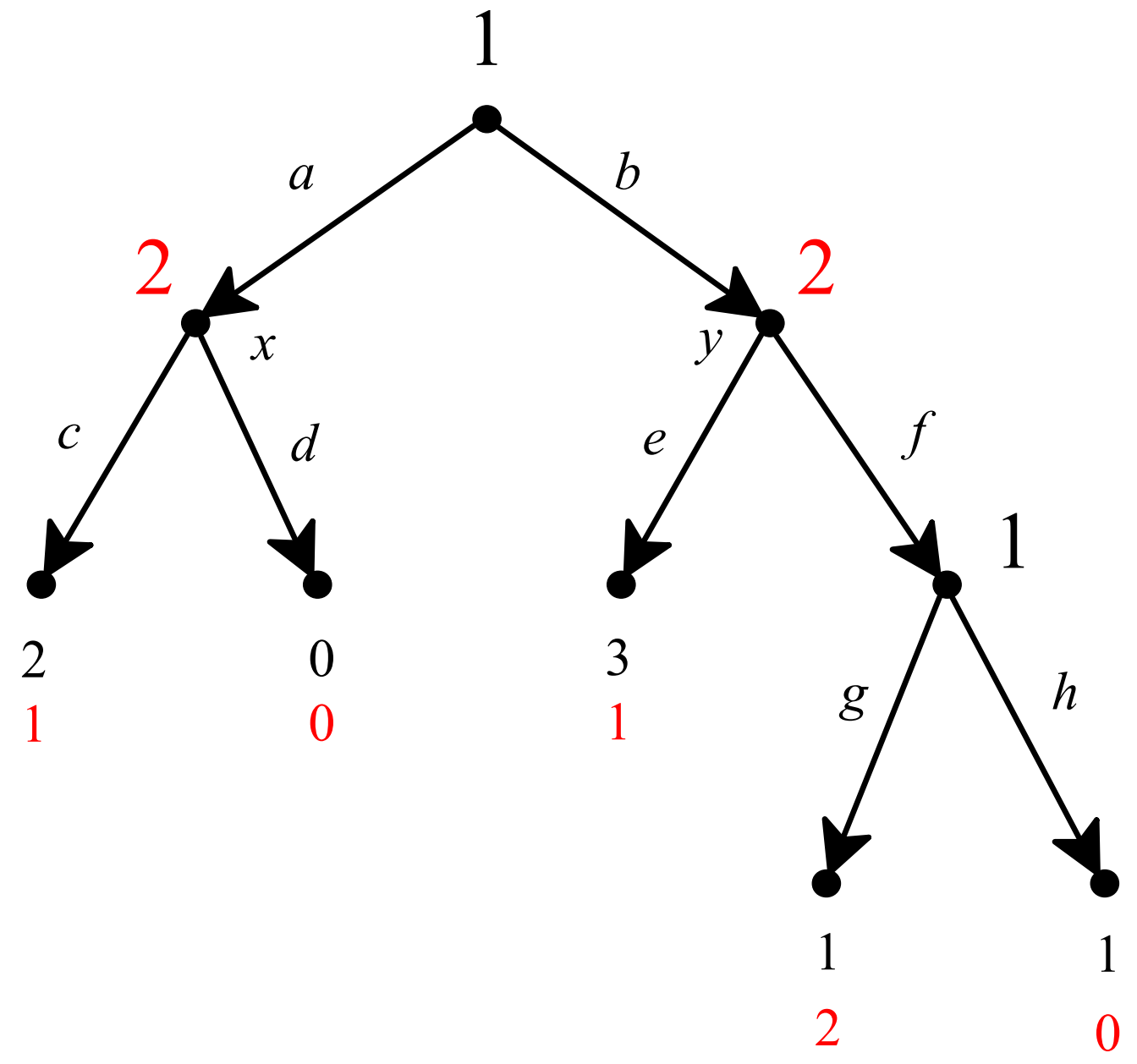

**Figure 2.6**
A perfect-information game.

Using Definition 2.4, one can associate with every perfect-information game a strategic-form (or normal-form) game: a strategy profile determines a unique terminal node that is reached if the players act according to that strategy profile and thus a unique vector of payoffs. Figure 2.7 shows the strategic-form associated with the perfect-information game of Figure 2.6, with the Nash equilibria highlighted.

Player 2

| Player 1 | *ce* | *cf* | *de* | *df* |
|---|---|---|---|---|
| *ag* | 2 1 | 2 1 | 0 0 | 0 0 |
| *ah* | 2 1 | 2 1 | 0 0 | 0 0 |
| *bg* | 3 1 | 1 2 | 3 1 | 1 2 |
| *bh* | 3 1 | 1 0 | 3 1 | 1 0 |

**Figure 2.7**
The strategic form of the perfect-information game of Figure 2.6
with the Nash equilibria highlighted.

Because of the redundancy discussed above, the strategic form also displays redundancies: in this case the top two rows are identical.

Armed with the notion of strategy, we can now revisit the notion of backward-induction solution. Figure 2.8 shows the two backward-induction solutions of the game of Figure 2.6.

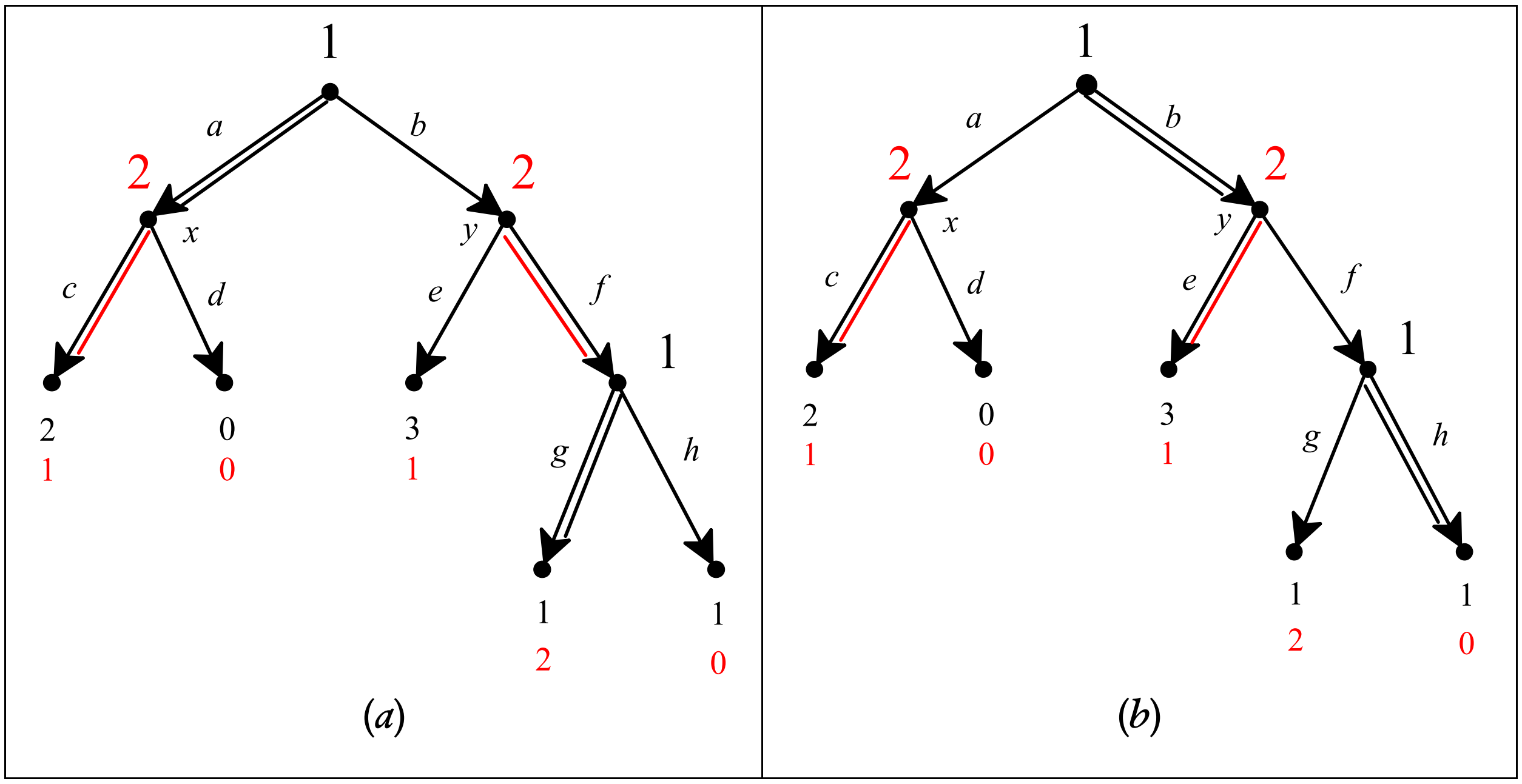

**Figure 2.8**
The backward-induction solutions of the game of Figure 2.6.

It is clear from the definition of the backward-induction algorithm (Definition 2.3) that the procedure selects a choice at every decision node and thus yields a strategy profile for the entire game: the backward-induction solution shown in Panel *a* of Figure 2.8 is the strategy profile $((a,g),(c,f))$, while the backward-induction solution shown in Panel *b* is the strategy profile $((b,h),(c,e))$. Both of them are Nash equilibria of the strategic form, but not all the Nash equilibria correspond to backward -induction solutions. The relationship between the two concepts will be explained in the following section.

**Remark 2.1.** A backward-induction solution is a strategy profile. Since strategies contain a description of what a player actually does and also of what the player would do in circumstances that do not arise, one often draws a distinction between the backward-induction solution and the *backward-induction outcome* which is defined as the sequence of actual moves. For example, the backward-induction outcome associated with the solution $((a,g),(c,f))$ is the play *ac* with corresponding payoff (2,1), while the backward-induction outcome associated with the solution $((b,h),(c,e))$ is the play *be* with corresponding payoff (3,1).

This is a good time to test your understanding of the concepts introduced in this section, by going through the exercises in Section 2.E.3 of Appendix 2.E at the end of this chapter.

## 2.4 Relationship between backward induction and other solutions

If you have gone through the exercises for the previous three sections, you will have seen that in all those games the backward-induction solutions are also Nash equilibria. This is always true, as stated in the following theorem.

**Theorem 2.1.** Every backward-induction solution of a perfect-information game is a Nash equilibrium of the associated strategic form.

In some games the set of backward-induction solutions coincides with the set of Nash equilibria (see, for example, Exercise 2.9), but typically the set of Nash equilibria is larger than (is a proper superset of) the set of backward-induction solutions (as shown for the game of Figure 2.6, which has two backward-induction solutions, shown in Figure 2.8, but five Nash equilibria, shown in Figure 2.7). Nash equilibria that are not backward-induction solutions often involve *incredible threats*. To see this, consider the following game. An industry is currently a monopoly and the incumbent monopolist is currently making a profit of \$5 million. A potential entrant is considering whether or not to enter this industry. If she does not enter, she will make \$1 million in an alternative investment; if she does enter, then the incumbent can either fight entry with a price war whose outcome is that both firms make zero profits, or it can accommodate entry, by sharing the market with the entrant, in which case both firms make a profit of \$2 million. This situation is illustrated in Figure 2.9 with the associated strategic form. Note that we are assuming that each player is selfish and greedy, that is, cares only about its own profit and prefers more money to less.

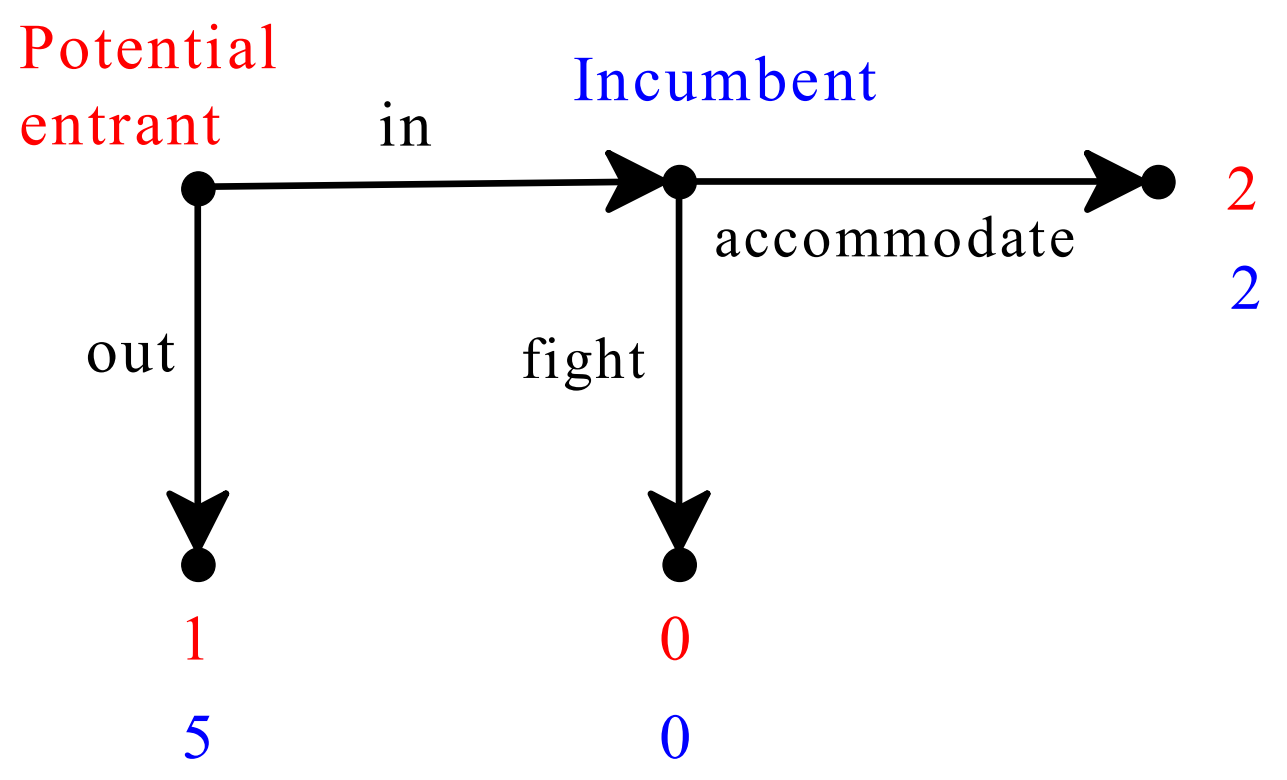

| | | Chain store: fight | Chain store: accommodate |
|---|---|---|---|
| Potential entrant | In | 0 , 0 | 2 , 2 |
| | Out | 1 , 5 | 1 , 5 |

**Figure 2.9**
The entry game.

The backward-induction solution is (in, accommodate) and it is also a Nash equilibrium. However, there is another Nash equilibrium, namely (out, fight). The latter should be discarded as a "rational solution" because it involves an incredible threat on the part of the incumbent, namely that it will fight entry if the potential entrant enters. It is true that, if the potential entrant believes the incumbent's threat, then she is better off staying out; however, she should ignore the incumbent's threat because she should realize that – when faced with the *fait accompli* of entry – the incumbent would not want to carry out the threat.

Reinhard Selten (who shared the 1994 Nobel prize in economics with two other game theorists, John Harsanyi and John Nash) discussed a repeated version of the above entry game, which has become known as Selten's Chain Store Game. The story is as follows. A chain store is a monopolist in an industry. It owns stores in $m$ different towns ($m \geq 2$). In each town the chain store makes \$5 million if left to enjoy its privileged position undisturbed. In each town there is a businesswoman who could enter the industry in that town, but earns \$1 million if she chooses not to enter; if she decides to enter, then the monopolist can either fight the entrant, leading to zero profits for both the chain store and the entrant in that town, or it can accommodate entry and share the market with the entrant, in which case both players make \$2 million in that town. Thus, in each town the interaction between the incumbent monopolist and the potential entrant is as illustrated in Figure 2.9 above. However, decisions are made sequentially, as follows. At date $t$ ($t$ = 1, ..., $m$) the businesswoman in town $t$ decides whether or not to enter and if she enters then the chain store decides whether or not to fight in that town. What happens in town $t$ at date $t$ becomes known to everybody. Thus, for example, the businesswoman in town 2 at date 2 knows what happened in town 1 at date 1 (either that there was no entry or that entry was met with a fight or that entry was accommodated). Intuition suggests that in this game the threat by the incumbent to fight early entrants might be credible, for the following reason. The incumbent could tell Businesswoman 1 the following:

> "It is true that, if you enter and I fight, I will make zero profits, while by accommodating your entry I would make \$2 million and thus it would seem that it cannot be in my interest to fight you. However, somebody else is watching us, namely Businesswoman 2. If she sees that I have fought your entry then she might fear that I would do the same with her and decide to stay out, in which case in town 2 I would make \$5 million, so that my total profits in towns 1 and 2 would be \$(0+5) = \$5 million. On the other hand, if I accommodate your entry, then she will be encouraged to entry herself and I will make \$2 million in each town, for a total profit of \$4 million. Hence – as you can see – it is indeed in my interest to fight you and thus you should stay out."

Does the notion of backward induction capture this intuition? To check this, let us consider the case where $m$ = 2, so that the extensive game is not too large to draw. It is shown in Figure 2.10, where at each terminal node the top number is the profit of the incumbent monopolist (it is the sum of the profits in the two towns), the middle number is the profit of Businesswoman 1 and the bottom number is the profit of Businesswoman 2. All profits are expressed in millions of dollars. We assume that all the players are selfish and greedy, so that we can take the profit of each player to be that player's payoff.

The backward-induction solution is unique and is shown by the thick directed edges in Figure 2.10. The corresponding outcome is that both businesswomen will enter and the incumbent monopolist accommodates entry in both towns.

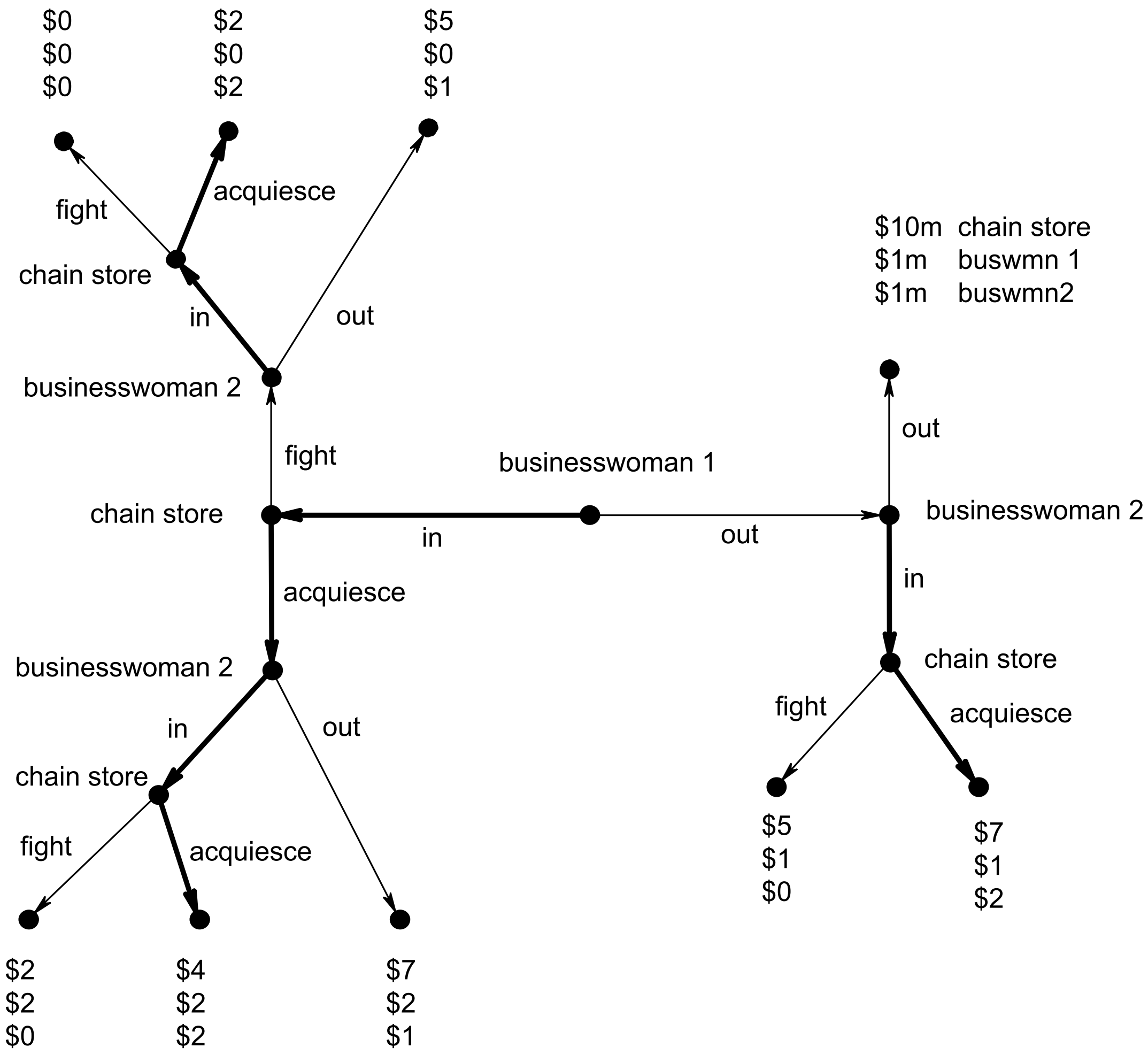

**Figure 2.10**
Selten's Chain-Store game.

Thus the backward-induction solution does not capture the "reputation" argument outlined above. However, the backward-induction solution does seem to capture the notion of rational behavior in this game. Indeed, Businesswoman 1 should reply to the incumbent with the following counter-argument:

> "Your reasoning is not valid. Whatever happens in town 1, it will be common knowledge between you and Businesswoman 2 that your interaction in town 2 will be the last; in particular, nobody else will be watching and thus there won't be an issue of establishing a reputation in the eyes of another player. Hence in town 2 it will be in your interest to accommodate entry, since – in essence – you will be playing the one-shot entry game of Figure 2.9. Hence a rational Businesswoman 2 will decide to enter in town 2 *whatever happened in town* 1: what you do against me will have no influence on her decision. Thus your "reputation" argument does not apply and it will in fact be in your interest not to fight my entry: your choice will be between a profit of \$(0+2) = \$2 million, if you fight me, and a profit of \$(2+2) = \$4 million, if you don't fight me. Hence I will enter and you will not fight me."

In order to capture the reputation argument described above we need to allow for some uncertainty in the mind of some of the players, as we will show in a later chapter. In a perfect-information game uncertainty is ruled out by definition.

By Theorem 2.1 the notion of backward induction can be seen as a refinement of the notion of Nash equilibrium. Another solution concept that is related to backward induction is the iterated elimination of weakly dominated strategies. Indeed the backward-induction algorithm could be viewed as a step-wise procedure that eliminates dominated choices at decision nodes, and thus strategies that contain those choices. What is the relationship between the two notions? In general this is all that can be said: applying the iterated deletion of weakly dominated strategies to the strategic form associated with a perfect-information game leads to a set of strategy profiles that contains at least one backward-induction solution; however, (1) it may also contain strategy profiles that are not backward-induction solutions and (2) it may fail to contain all the backward-induction solutions, as shown in Exercise 2.8.

## 2.5 Perfect-information games with two players

We conclude this chapter with a discussion of finite two-player extensive games with perfect information.

We will start with games that have only to outcomes, namely "Player 1 wins" (denoted by $W_1$) and "Player 2 wins" (denoted by $W_2$). We assume that Player 1 strictly prefers $W_1$ to $W_2$ and Player 2 strictly prefers $W_2$ to $W_1$. Thus we can use utility functions with values 0 and 1 and associate with each terminal node either the payoff vector (1,0) (if the outcome is $W_1$) or the payoff vector (0,1) (if the outcome is $W_2$). We call these games *win-lose games*. An example of such a game is the following.

**Example 2.2.** Two players take turns choosing a number from the set {1,2,...,10}, with Player 1 moving first. The first player who brings the sum of all the chosen numbers to 100 or more wins.

The following is one possible play of the game (the red numbers are the ones chosen by Player 1 and the blue numbers the ones chosen by Player 2):

10, 9, 9, 10, 8, 7, 10, 10, 1, 8, 1, 7, 6, 4

In this play Player 2 wins: at her last move the sum is 96 and with her choice of 4 she brings the total to 100. However, in this game Player 1 has a *winning strategy*, that is, a strategy that guarantees that he will win, no matter what numbers Player 2 chooses. To see this, we can use backward-induction reasoning. Drawing the tree is not a practical option, since the number of nodes is very large: one needs 10,000 nodes just to represent the first 4 moves! But we can imagine drawing the tree, placing ourselves towards the end of the tree and ask what partial sum represents a "losing position", in the sense that whoever is choosing in that position cannot win, while the other player can then win with his subsequent choice. With some thought one can see that 89 is the largest losing position: whoever moves there can take the sum to any number in the set {90, 91, ..., 99}, thus coming short of 100, while the other player can then take the sum to 100 with an appropriate choice. What is the largest losing position that precedes 89? The answer is 78: whoever moves at 78 must take the sum to a number in the set {79, 80, ..., 88} and then from there the other player can make sure to take the sum to 89 and then we know what happens from there! Repeating this reasoning we see that the losing positions are: 89, 78, 67, 56, 45, 34, 23, 12, 1. Since Player 1 moves first he can choose 1 and put Player 2 in the first losing position; then – whatever Player 2 chooses – Player 1 can put her in the next losing position, namely 12, etc. Recall that a strategy for Player 1 must specify what to do in every possible situation in which he might find himself.

In his game Player 1's winning strategy is as follows:

> Start with the number 1. Then, at every turn, choose the number $(11 - n)$, where $n$ is the number that was chosen by Player 2 in the immediately preceding turn.

Here is an example of a possible play of the game where Player 1 employs the winning strategy and does in fact win:

1, 9, 2, 6, 5, 7, 4, 10, 1, 8, 3, 3, 8, 9, 2, 5, 6, 1, 10.

We can now state a general result about this class of games.

**Theorem 2.2.** In every finite two-player, win-lose game with perfect information one of the two players has a winning strategy.

Although we will not give a detailed proof, the argument of the proof is rather simple. By applying the backward-induction algorithm we assign to every decision node either the payoff vector (1,0) or the payoff vector (0,1). Imagine applying the algorithm up to the point where the immediate successors of the root have been assigned a payoff vector. Two cases are possible.

**Case 1**: at least one of the immediate successors of the root has been assigned the payoff vector (1,0). In this case Player 1 is the one who has a winning strategy and his initial choice should be such that a node with payoff vector (1,0) is reached and then his future choices should also be such that only nodes with payoff vector (1,0) are reached.

**Case 2:** all the immediate successors of the root have been assigned the payoff vector (0,1). In this case it is Player 2 who has a winning strategy. An example of a game where it is Player 2 who has a winning strategy is given in Exercise 2.11.

We now turn to finite two-player games where there are three possible outcomes: "Player 1 wins" ($W_1$), "Player 2 wins" ($W_2$) and "Draw" ($D$). We assume that the rankings of the outcomes are as follows: $W_1 \succ_1 D \succ_1 W_2$ and $W_2 \succ_2 D \succ_2 W_1$. Examples of such games are Tic-Tac-Toe (http://en.wikipedia.org/wiki/Tic-tac-toe), Draughts or Checkers (http://en.wikipedia.org/wiki/Draughts) and Chess (although there does not seem to be agreement as to whether the rules of Chess guarantee that every possible play of the game is finite). What can we say about such games? The answer is provided by the following theorem.

**Theorem 2.3.** Every finite two-player, perfect-information game with three outcomes: Player 1 wins ($W_1$), Player 2 wins ($W_2$) and Draw ($D$) falls within one of the following three categories:

(1) Player 1 has a strategy that guarantees outcome $W_1$.

(2) Player 2 has a strategy that guarantees outcome $W_2$.

(3) Player 1 has a strategy that guarantees that the outcome will be $W_1$ or $D$ and Player 2 has a strategy that guarantees that the outcome will be $W_2$ or $D$, so that, if both players employ these strategies, the outcome will be $D$.

The logic of the proof is as follows. By applying the backward-induction algorithm we assign to every decision node either the payoff vector (2,0) (corresponding to outcome $W_1$) or the payoff vector (0,2) (corresponding to outcome $W_2$) or the payoff vector (1,1) (corresponding to outcome $D$). Imagine applying the algorithm up to the point where the immediate successors of the root have been assigned a payoff vector. Three cases are possible.

**Case 1**: at least one of the immediate successors of the root has been assigned the payoff vector (2,0); in this case Player 1 is the one who has a winning strategy.

**Case 2:** all the immediate successors of the root have been assigned the payoff vector (0,2); in this case it is Player 2 who has a winning strategy.

**Case 3:** there is at least one immediate successor of the root to which the payoff vector (1,1) has been assigned and all the other immediate successors of the root have been assigned either (1,1) or (0,2). In this case we fall within the third category of Theorem 2.10.

Both Tic-Tac-Toe and Draughts fall within the third category (http://en.wikipedia.org/wiki/Solved_game#Solved_games). As of 2015 it is not known to which category the game of Chess belongs.

This is a good time to test your understanding of the concepts introduced in this section, by going through the exercises in Section 2.E.5 of Appendix 2.E at the end of this chapter.

# Appendix 2.E: Exercises

## 2.E.1. Exercises for Section 2.1: trees, frames and games

The answers to the following exercises are in Appendix S at the end of this chapter.

**Exercise 2. 1.** How could they do that! They abducted Speedy, your favorite tortoise! They asked for $1,000 in unmarked bills and threatened to kill Speedy if you don't pay. Call the tortoise-napper Mr. T. Let the possible outcomes be as follows:

$$\begin{aligned} o_1 &: \text{ you don't pay and speedy is released} \\ o_2 &: \text{ you pay \$1,000 and speedy is released} \\ o_3 &: \text{ you don't pay and speedy is killed} \\ o_4 &: \text{ you pay \$1,000 and speedy is killed} \end{aligned}$$

You are attached to Speedy and would be willing to pay $1,000 to get it back. However, you also like your money and you prefer not to pay, conditional on the two separate events "Speedy is released" and "Speedy is killed". Thus your ranking of the outcomes is $o_1 \succ_{you} o_2 \succ_{you} o_3 \succ_{you} o_4$. On the other hand, you are not quite sure of what Mr. T's ranking is.

**(a)** Suppose first that Mr T has communicated that he wants you to go to Central Park tomorrow at 10:00am and leave the money in a garbage can; he also said that, two miles to the East and at the exact same time, he will decide whether or not to free Speedy in front of the police station and then go and collect his money in Central Park. What should you do?

**(b)** Suppose that Mr T is not as dumb as in part (a) and he instead gives you the following instructions: *first* you leave the money in a garbage can in Central Park and then he will go there to collect the money. He also told you that if you left the money there then he will free Speedy, otherwise he will kill it. Draw an extensive form or frame to represent this situation.

**(c)** Now we want to construct a game based on the extensive form of part (b). For this we need Mr T's preferences. There are two types of criminals in Mr T's line of work: the professionals and the one-timers. Professionals are in the business for the long term and thus worry about reputation; they want it to be known that (1) every time they were paid they honored their promise to free the hostage and (2) their threats are to be taken seriously, in the sense that every time they were not paid, the hostage was killed. The one-timers hit once and then they disappear; they don't try to establish a reputation and the only thing they worry about, besides money, is not to be caught: whether or not they get paid, they prefer to kill the hostage in order to eliminate any kind of evidence (DNA traces, fingerprints, etc.). Construct two games based on the extensive form of part (b) representing the two possible types of Mr T.

**Exercise 2.2**. A three-man board, composed of A, B, and C, has held hearings on a personnel case involving an officer of the company. This officer was scheduled for promotion but, prior to final action on his promotion, he made a decision that cost the company a good deal of money. The question is whether he should be (1) promoted anyway, (2) denied the promotion, or (3) fired. The board has discussed the matter at length and is unable to reach unanimous agreement. In the course of the discussion it has become clear to all three of them that their separate opinions are as follows:

- A considers the officer to have been a victim of bad luck, not bad judgment, and wants to go ahead and promote him but, failing that, would keep him rather than fire him.
- B considers the mistake serious enough to bar promotion altogether; he'd prefer to keep the officer, denying promotion, but would rather fire than promote him.
- C thinks the man ought to be fired but, in terms of personnel policy and morale, believes the man ought not to be kept unless he is promoted, i.e., that keeping an officer who has been declared unfit for promotion is even worse than promoting him.

To recapitulate, their preferences among the three outcomes are

| | PROMOTE | KEEP | FIRE |
|---|---|---|---|
| A: | best | middle | worst |
| B: | worst | best | middle |
| C: | middle | worst | best |

Assume that everyone's preferences among the three outcomes are fully evident as a result of discussion. The three must proceed to a vote.

Consider the following voting procedure. First A proposes an action (either promote or keep or fire). Then it is B's turn. If B accepts A's proposal, then this becomes the final decision. If B disagrees with A'a proposal, then C makes the final decision (which may be *any of the three*: promote, keep of fire). Represent this situation as an extensive game with perfect information.

## 2.E.2. Exercises for Section 2.2: backward induction.

The answers to the following exercises are in Appendix S at the end of this chapter.

**Exercise 2.3**. Apply the backward-induction algorithm to the two games of Exercise 2.1 Part *c*.

**Exercise 2.4**. Apply the backward-induction algorithm to the game of Exercise 2.2.

## 2.E.3. Exercises for Section 2.3: strategies in perfect-information games.

The answers to the following exercises are in Appendix S at the end of this chapter.

**Exercise 2.5**. Write the strategic form of the game of Figure 2.2, find all the Nash equilibria and verify that the backward-induction solution is a Nash equilibrium.

**Exercise 2.6**. Write the strategic form of the game of Figure 2.3, find all the Nash equilibria and verify that the backward-induction solutions are Nash equilibria.

**Exercise 2.7**. **(a)** Write down all the strategies of Player *B* in the game of Exercise 2.2.
**(b)** How many strategies does Player *C* have?

**Exercise 2.8**. Consider the following perfect-information game:

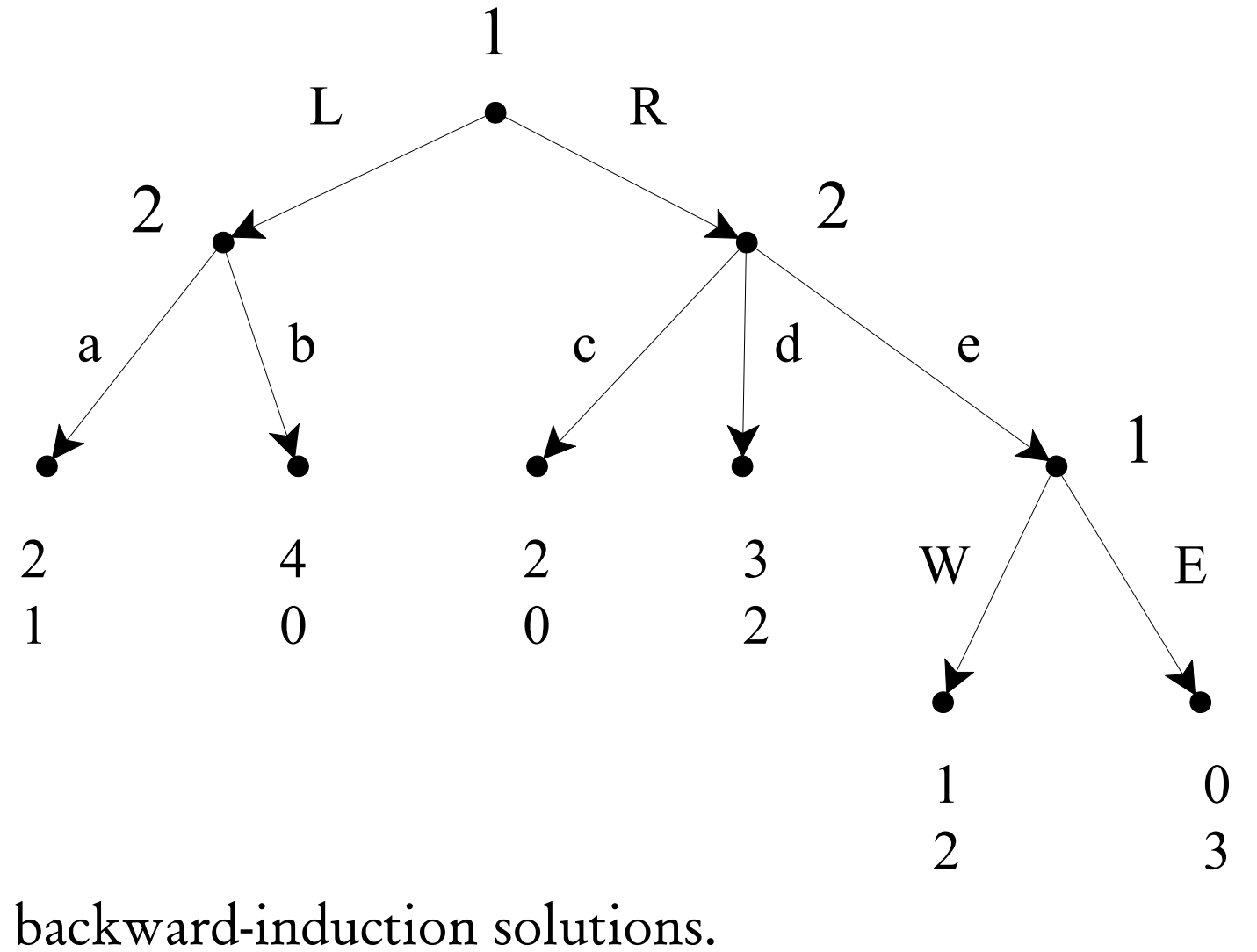

**(a)** Find the backward-induction solutions.
**(b)** Write down all the strategies of Player 1.
**(c)** Write down all the strategies of Player 2.
**(d)** Write the strategic form associated with this game.
**(e)** Does Player 1 have a dominant strategy?
**(f)** Does Player 2 have a dominant strategy?

**(g)** Is there a dominant-strategy equilibrium?
**(h)** Does Player 1 have any dominated strategies?
**(i)** Does Player 2 have any dominated strategies?
**(j)** What do you get when you apply the iterative elimination of weakly dominated strategies?
**(k)** What are the Nash equilibria?

**Exercise 2.9**. Consider an industry where there are two firms, a large firm, Firm 1, and a small firm, Firm 2. The two firms produce identical products. Let $x$ be the output of Firm 1 and $y$ the output of Firm 2. Industry output is $Q = x + y$. The price $P$ at which each unit of output can be sold is determined by the inverse demand function $P = 130 - 10Q$. For example, if Firm 1 produces 4 units and Firm 2 produces 2 units, then industry output is 6 and each unit is sold for $P$ = 130–60 = \$70. For each firm the cost of producing $q$ units of output is $C(q) = 10q + 62.5$. Each firm is only interested in its own profits. The profit of Firm 1 depends on both $x$ and $y$ and is given by $\Pi_1(x,y) = \underbrace{x[130 - 10(x+y)]}_{\text{revenue}} - \underbrace{(10x + 62.5)}_{\text{cost}}$ and similarly the profit function of Firm 2 is given by $\Pi_2(x,y) = \underbrace{y[130 - 10(x+y)]}_{\text{revenue}} - \underbrace{(10y + 62.5)}_{\text{cost}}$. The two firms play the following sequential game. First Firm 1 chooses its own output $x$ and commits to it; then Firm 2, after having observed Firm 1's output, chooses its own output $y$; then the price is determined according to the demand function and the two firms collect their own profits. In what follows assume, for simplicity, that $x$ can only be 6 or 6.5 units and $y$ can only be 2.5 or 3 units.

**(a)** Represent this situation as an extensive game with perfect information.

**(b)** Solve the game using backward induction.

**(c)** Write the strategic form associated with the perfect-information game.

**(d)** Find the Nash equilibria of this game and verify that the backward-induction solutions are Nash equilibria.

**Exercise 2.10**. Consider the following perfect-information game where $x$ is an integer.

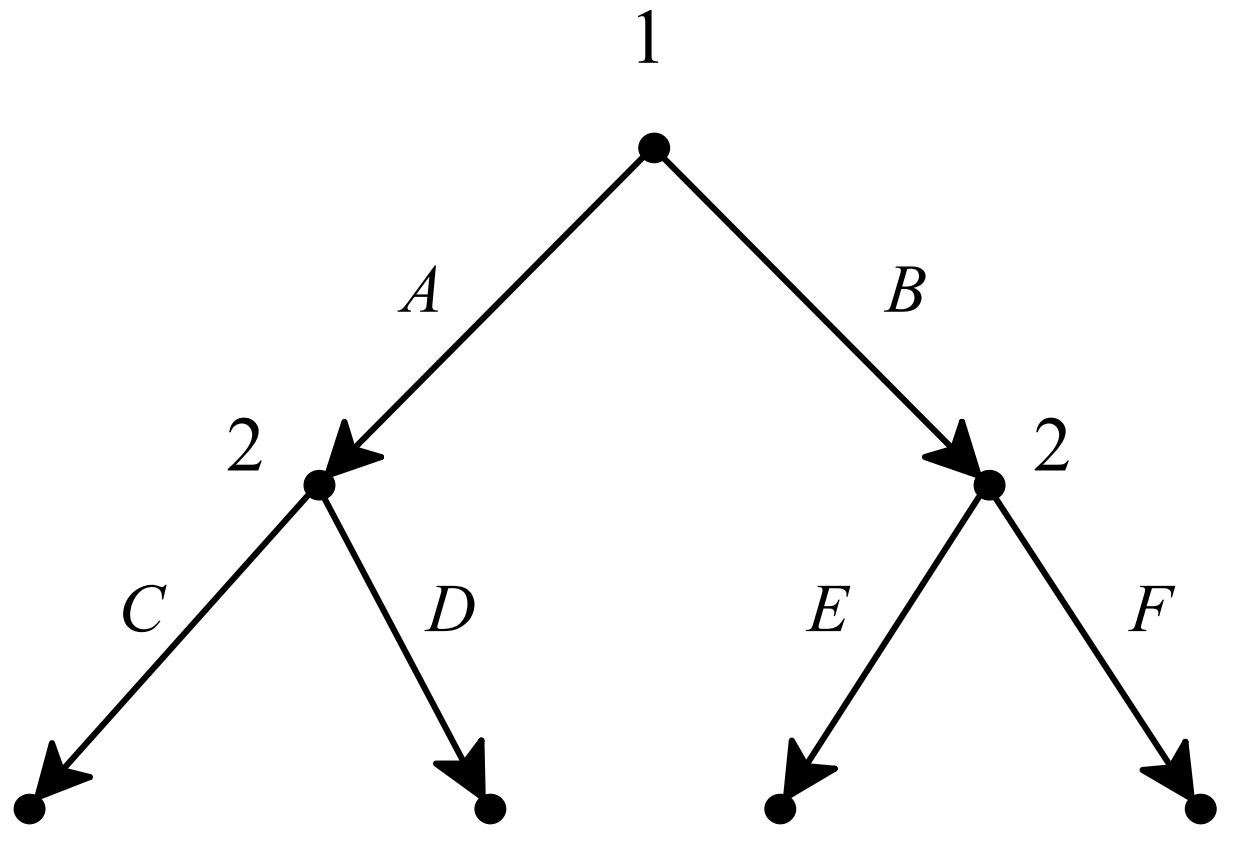

| | | | | |
|---|---|---|---|---|
| Player 1's payoff | 2 | 1 | 1 | $x$ |
| Player 2's payoff | 5 | 0 | 0 | 2 |

**(a)** For every value of $x$ find the backward induction solution(s).

**(b)** Write the corresponding strategic-form and find all the Nash equilibria.

## 2.E.5. Exercises for Section 2.5: two-player games.

The answers to the following exercises are in Appendix S at the end of this chapter.

**Exercise 2.11**. Consider the following perfect-information game. Player 1 starts by choosing a number from the set {1,2,3,4,5,6,7}, then Player 2 chooses a number from this set, then Player 1 again, followed by Player 2, etc. The first player who brings the cumulative sum of all the numbers chosen (up to and including the last one) to 48 or more wins. By Theorem 2.9 one of the two players has a winning strategy. Find out who that player is and fully describe the winning strategy.

**Exercise 2.12**. Consider the following two-player, perfect-information game. A coin is placed in the cell marked 'START' (cell A1). Player 1 moves first and can move the coin one cell up (to A2) or one cell to the left (to B1) or one cell diagonally in the left-up direction (to B2). Then Player 2 moves, according to the same rules (e.g. if the coin is in cell B2 then the admissible moves are shown by the directed edges). The players alternate moving the coin. Black cells are not accessible (so that, for example, from A3 the coin can only be moved to A4 or B3 and from F3 it can only be moved to G4, as shown by the directed edge). **The player who manages to place the coin in the cell marked 'END' wins.**

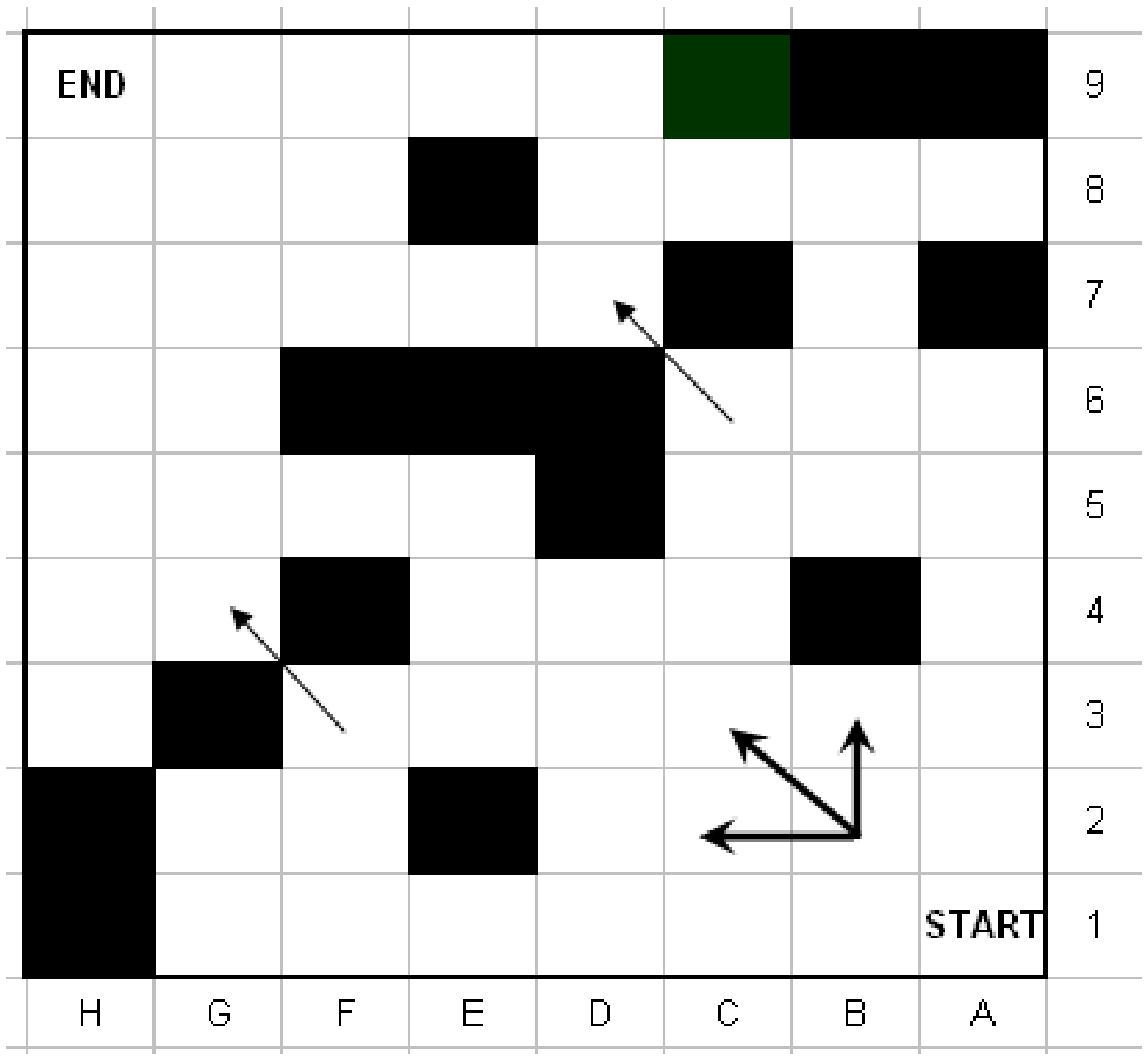

**(a)** Represent this game by means of an extensive form with perfect information by drawing the initial part of the tree that covers the first two moves (the first move of Player 1 and the first move of Player 2).

**(b)** Suppose that the coin is currently in cell G4 and it is Player 1's turn to move. Show that Player 1 has a strategy that allows her to win the game starting from cell G4. Describe the strategy in detail.

**(c)** Describe a play of the game (from cell A1) where Player 1 wins (describe it by means of the sequence of cells visited by the coin).

**(d)** Describe a play of the game (from cell A1) where Player 2 wins (describe it by means of the sequence of cells visited by the coin).

**(e)** Now go back to the beginning of the game. The coin is in cell A1 and player 1 has the first move. By Theorem 2.9 one of the two players has a winning strategy. Find out who that player is and fully describe the winning strategy.

◊◊◊◊◊◊◊◊◊◊◊◊

**Exercise 2.13: Challenging Question.** Two women, Anna and Bess, claim to be the legal owners of a diamond ring that – each claims – has great sentimental value. Neither of them can produce evidence of ownership and nobody else is staking a claim on the ring. Judge Sabio wants the ring to go to the legal owner, but he does not know which of the two women is in fact the legal owner. He decides to proceed as follows. First he announces a fine of $\$F > 0$ and then asks Anna and Bess to play the following game.

**Move 1**: Anna moves first. Either she gives up her claim to the ring (in which case Bess gets the ring, the game ends and nobody pays the fine) or she asserts her claim, in which case the game proceeds to Move 2.

**Move 2**: Bess either accepts Anna's claim (in which case Anna gets the ring, the game ends and nobody pays the fine) or challenges her claim. In the latter case, Bess must put in a bid, call it $B$, and Anna must pay the fine of $\$F$ to Sabio. The game goes on to Move 3.

**Move 3**: Anna now either matches Bess's bid (in which case Anna gets the ring, Anna pays $\$B$ to Sabio – in addition to the fine that she already paid – and Bess pays the fine of $\$F$ to Sabio) or chooses not to match (in which case Bess gets the ring and pays her bid of $\$B$ to Sabio and, furthermore, Sabio keeps the fine that Anna already paid).

Denote by $C_A$ the monetary equivalent of getting the ring for Anna (that is, getting the ring is as good, in Anna's mind, as getting $\$C_A$) and $C_B$ the monetary equivalent of getting the ring for Bess. Not getting the ring is considered by both as good as getting zero dollars.

**(a)** Draw an extensive game with perfect information to represent the above situation, assuming that there are only two possible bids: $B_1$ and $B_2$. Write the payoffs to Anna and Bess next to each terminal node.

**(b)** Find the backward-induction solution of the game you drew in part (a) for the case where $B_1 > C_A > C_B > B_2 > F > 0$.

Now consider the general case where the bid $B$ can be any non-negative number and assume that both Anna and Bess are very wealthy. Assume also that $C_A$, $C_B$ and F are positive numbers and that $C_A$ and $C_B$ are common knowledge between Anna and Bess . We want to show that, at the backward-induction solution of the game, the ring always goes to the legal owner. Since we (like Sabio) don't know who the legal owner is, we must consider two cases.

**Case 1:** the legal owner is Anna. Let us assume that this implies that $C_A > C_B$.

**Case 2:** the legal owner is Bess. Let us assume that this implies that $C_B > C_A$.

(c) Find the backward-induction solution for Case 1 and show that it implies that the ring goes to Anna.

(d) Find the backward-induction solution for Case 2 and show that it implies that the ring goes to Bess.

(e) How much money does Sabio make in equilibrium? How much money do Ann and Bess end up paying in equilibrium? (By 'equilibrium' we mean 'backward induction solution'.)

# Appendix 2.S: Solutions to exercises

**Exercise 2.1.** **(a)** For you it is a strictly dominant strategy to not pay and thus you should not pay.

**(b)** The extensive form is as follows:

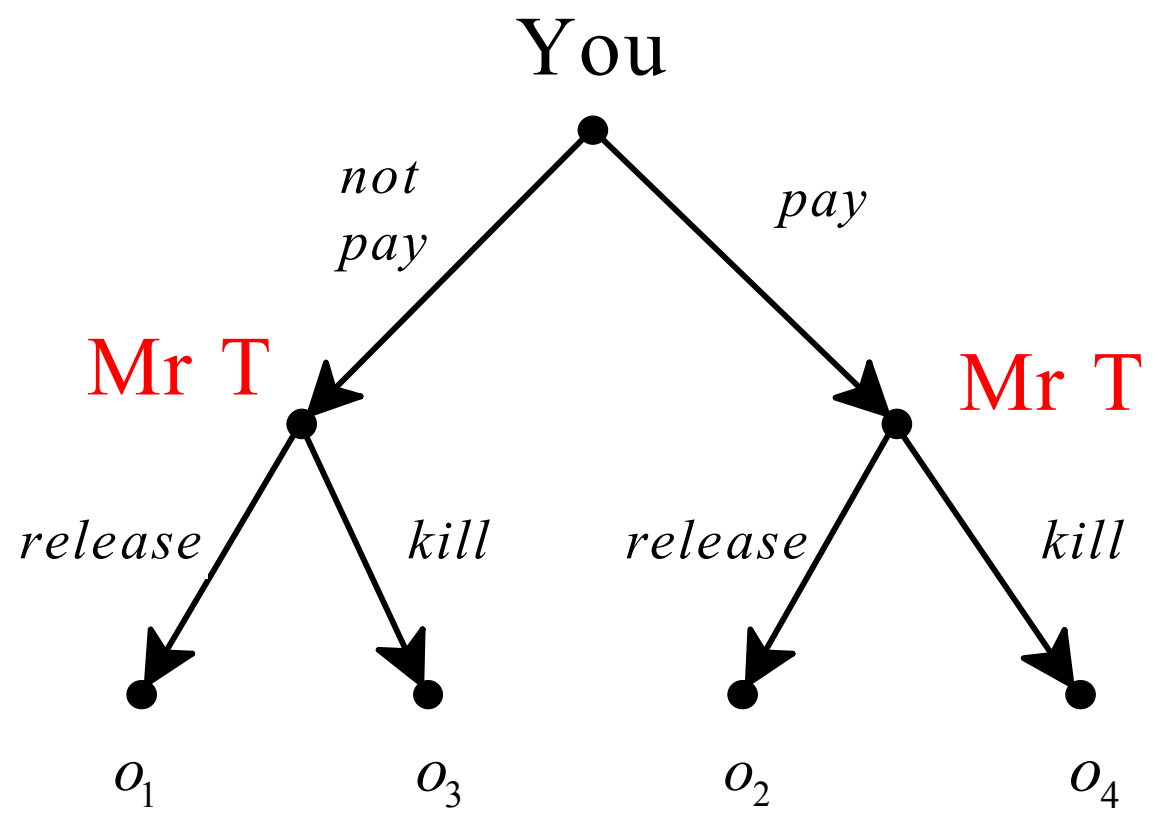

**(c)** For the professional, concern with reputation implies that $o_2 \succ_{MrT} o_4$ and $o_3 \succ_{MrT} o_1$. If we add the reasonable assumption that – after all – money is what they are after, then we can take the full ranking to be $o_2 \succ_{MrT} o_4 \succ_{MrT} o_3 \succ_{MrT} o_1$. Representing preferences with ordinal utility functions with values in the set {1,2,3,4}, we have

| outcome → <br> utility function ↓ | $o_1$ | $o_2$ | $o_3$ | $o_4$ |
|---|---|---|---|---|
| $U_{you}$ | 4 | 3 | 2 | 1 |
| $U_{MrT}$ | 1 | 4 | 2 | 3 |

The corresponding game is as follows:

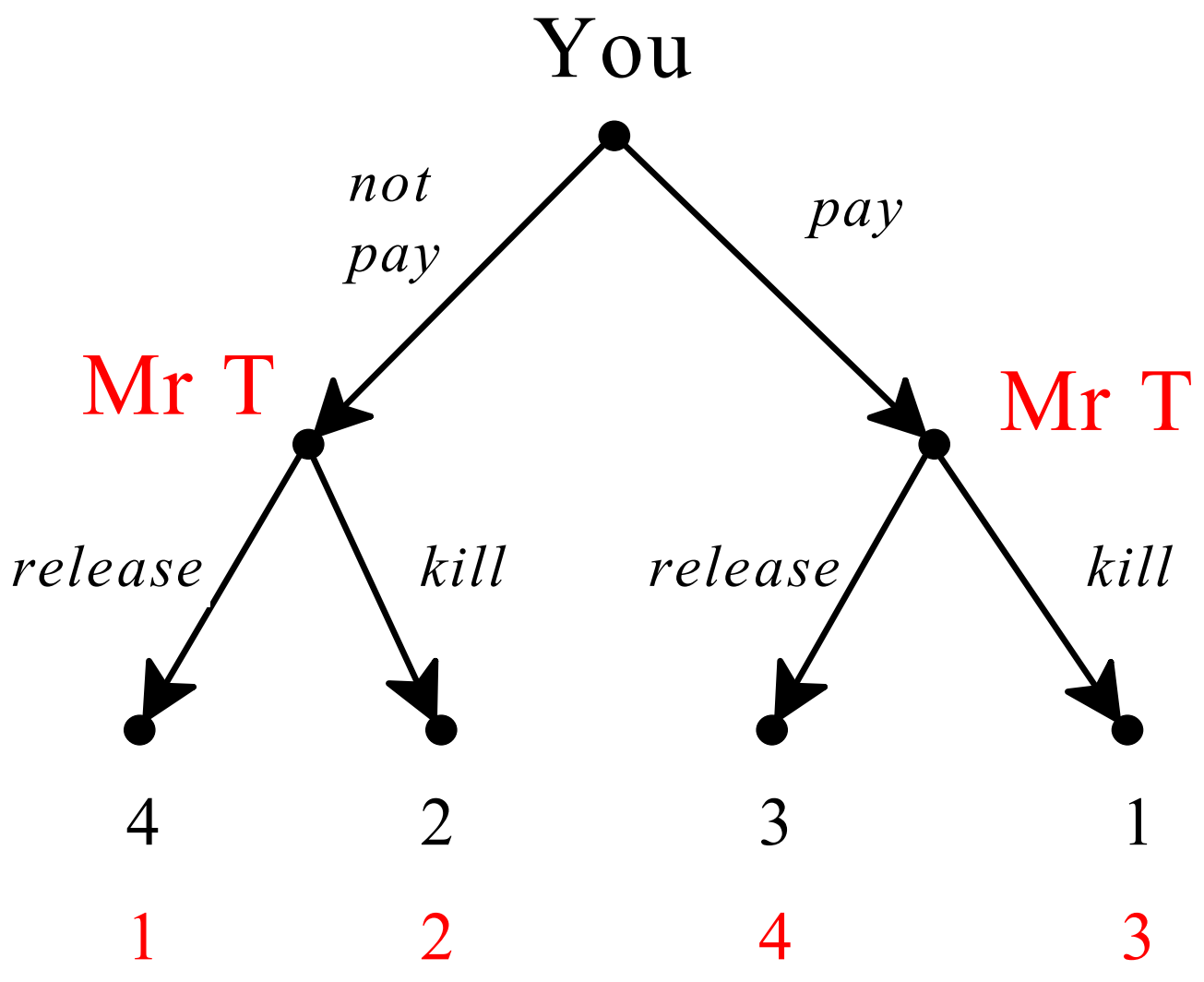

For the one-timer, the ranking can be taken to be (although this is not the only possibility) $o_4 \succ_{MrT} o_2 \succ_{MrT} o_3 \succ_{MrT} o_1$, with corresponding utility representation:

| outcome → <br> utility function ↓ | $o_1$ | $o_2$ | $o_3$ | $o_4$ |
|---|---|---|---|---|
| $U_{you}$ | 4 | 3 | 2 | 1 |
| $U_{MrT}$ | 1 | 3 | 2 | 4 |

and extensive game:

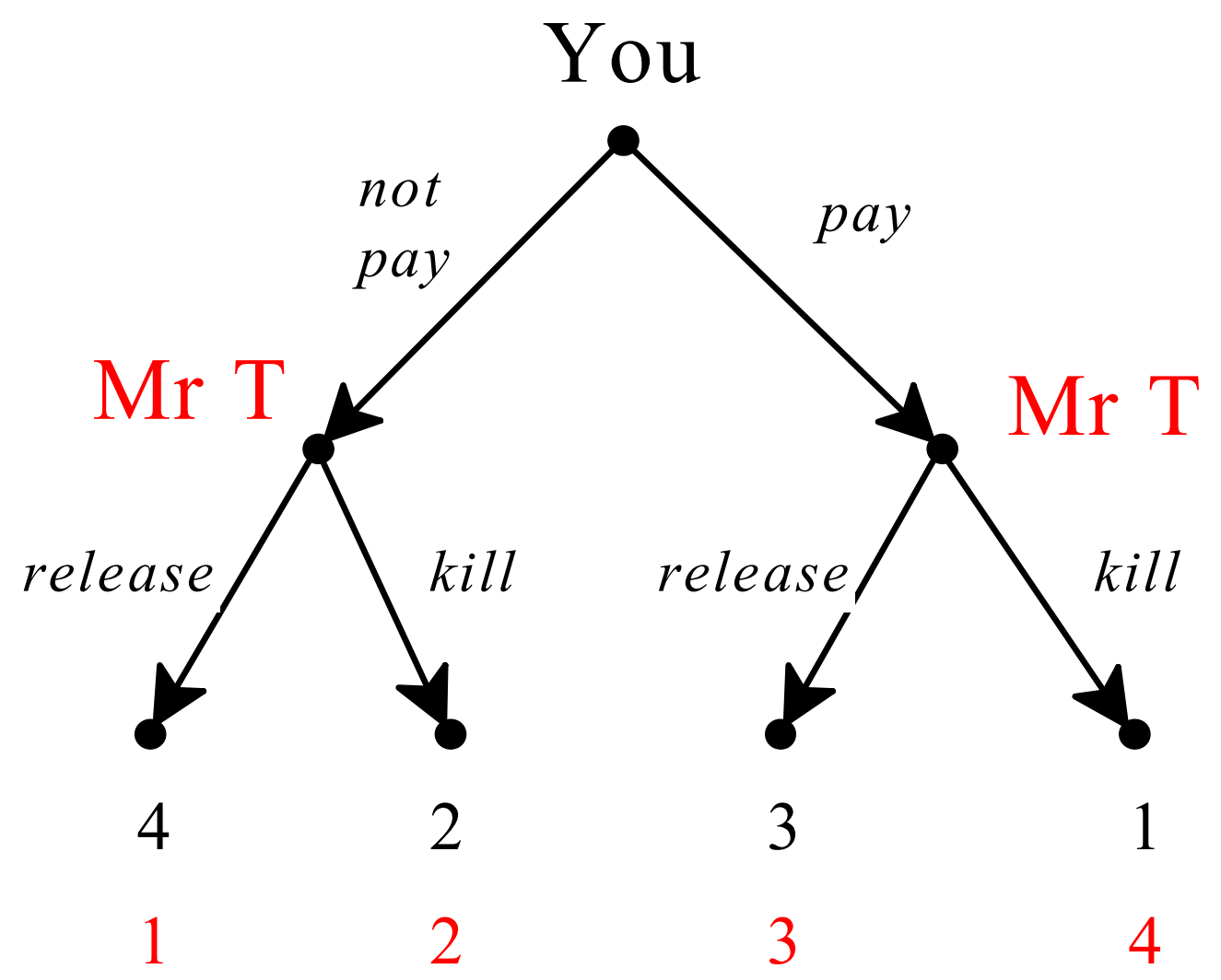

**Exercise 2.2.**
The game is as follows ("P" stands for promote, "K" for keep (without promoting), "F" for fire):

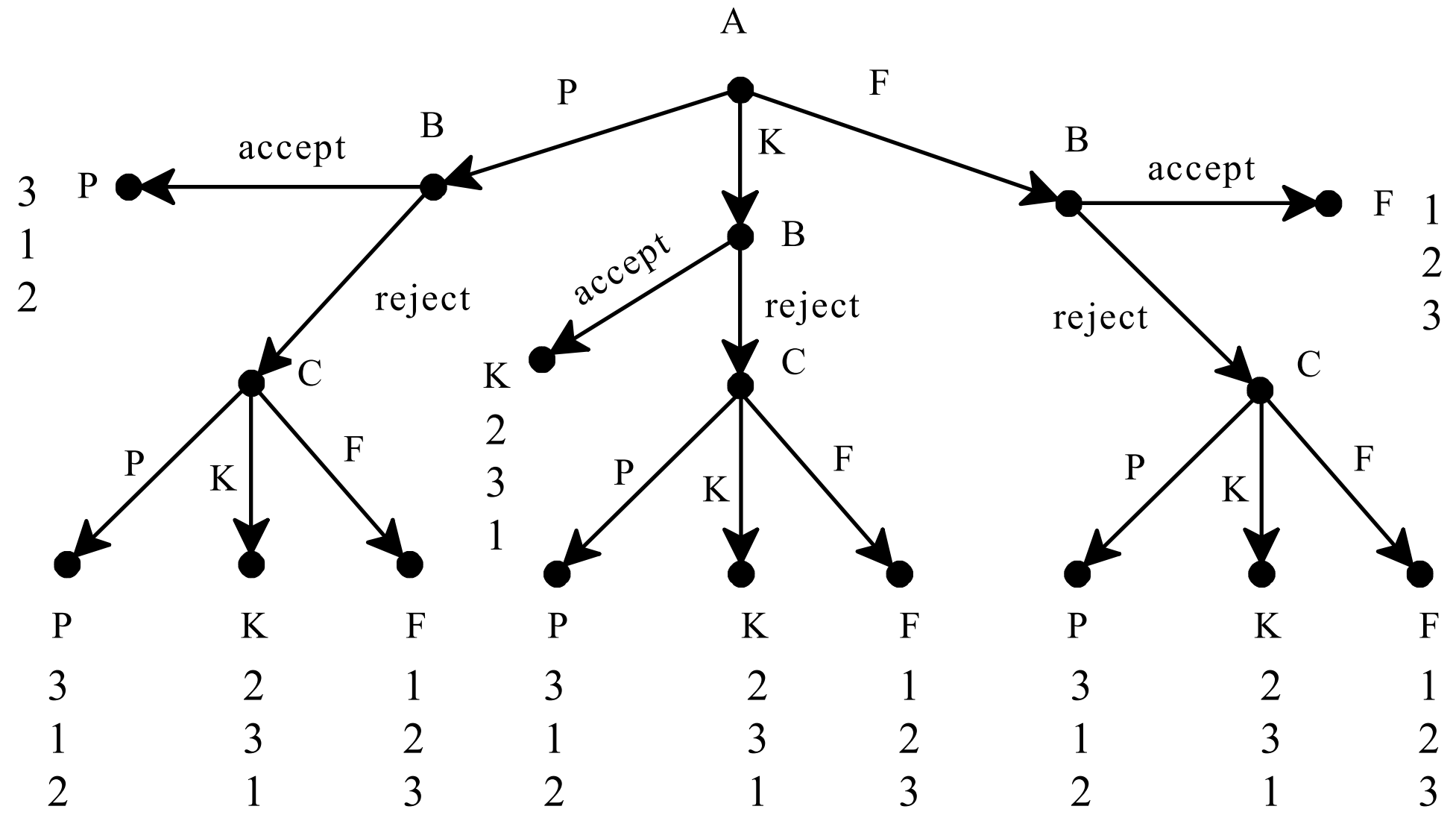

**Exercise 2.3**. The application of the backward-induction algorithm is shown by double edges in the following games, the first refers to the professional Mr. T and the second to the one-timer.

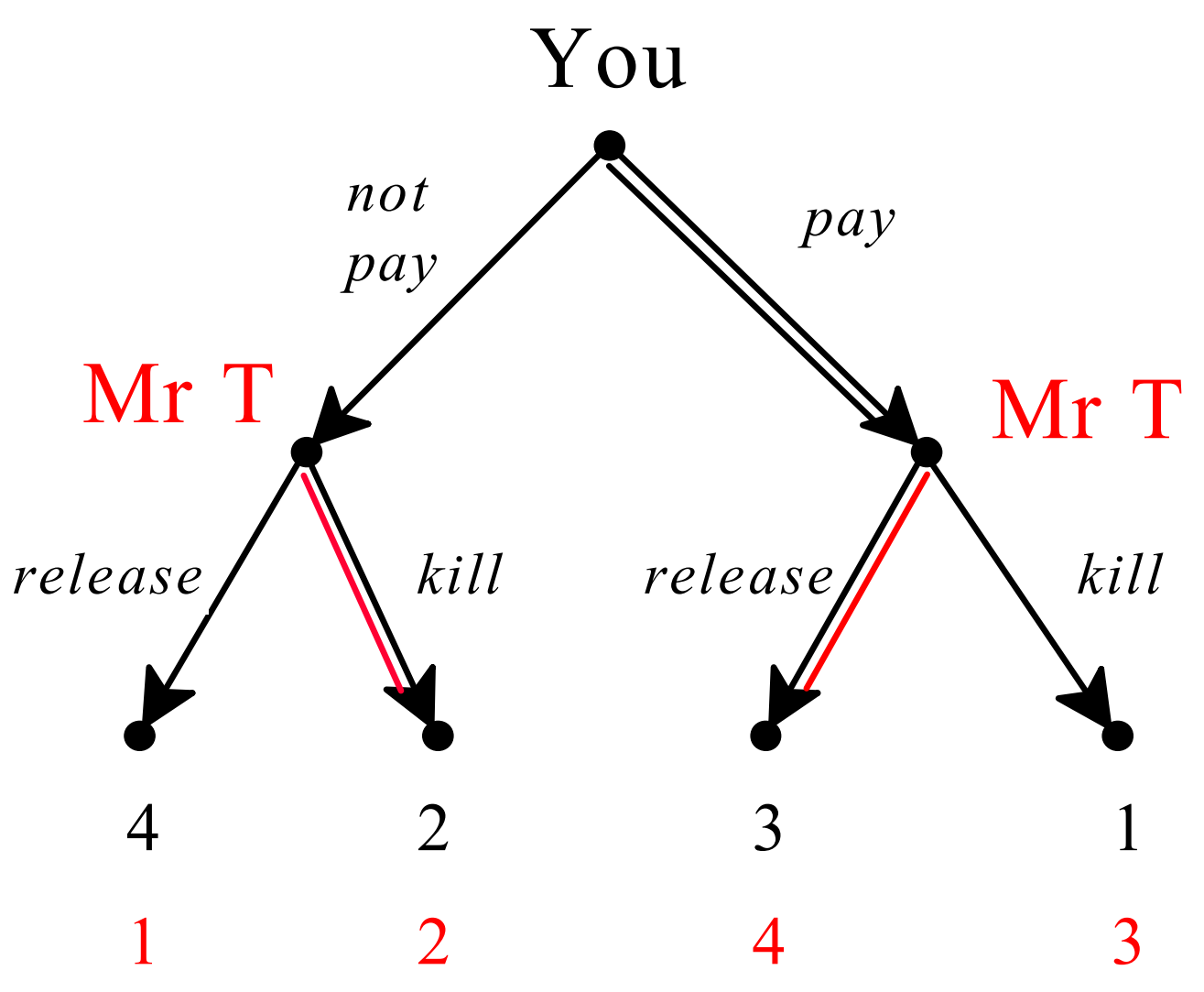

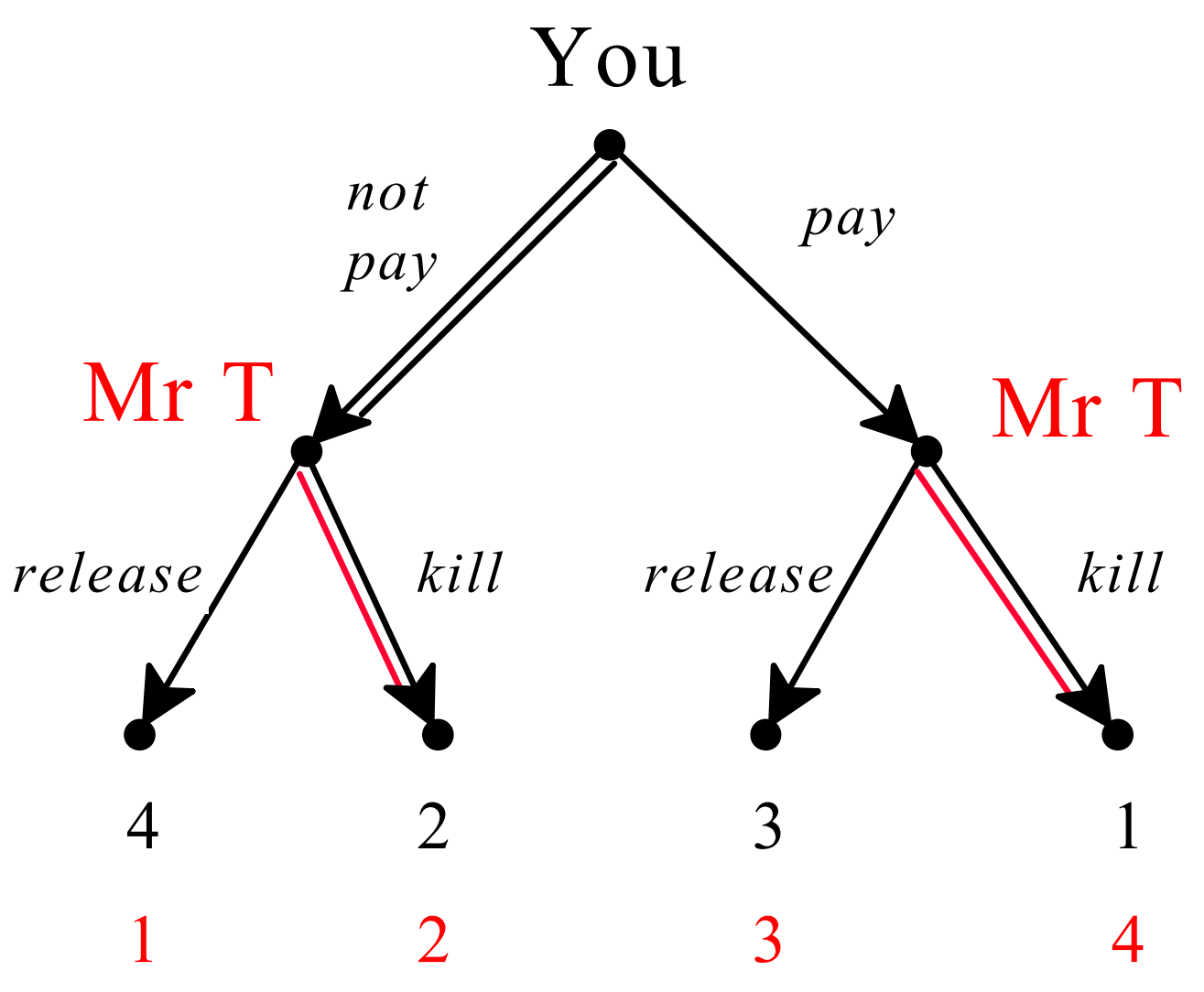

Thus, against a professional you would pay and against a one-timer you would not pay. With the professional you would get Speedy back, with the one-timer you will hold a memorial service for Speedy.

**Exercise 2.4**. The backward-induction algorithm yields two solutions, shown below. The difference between the two solutions lies in what Player *B* would do if Player *A* proposed *F*. In both solutions the officer is kept without promotion.

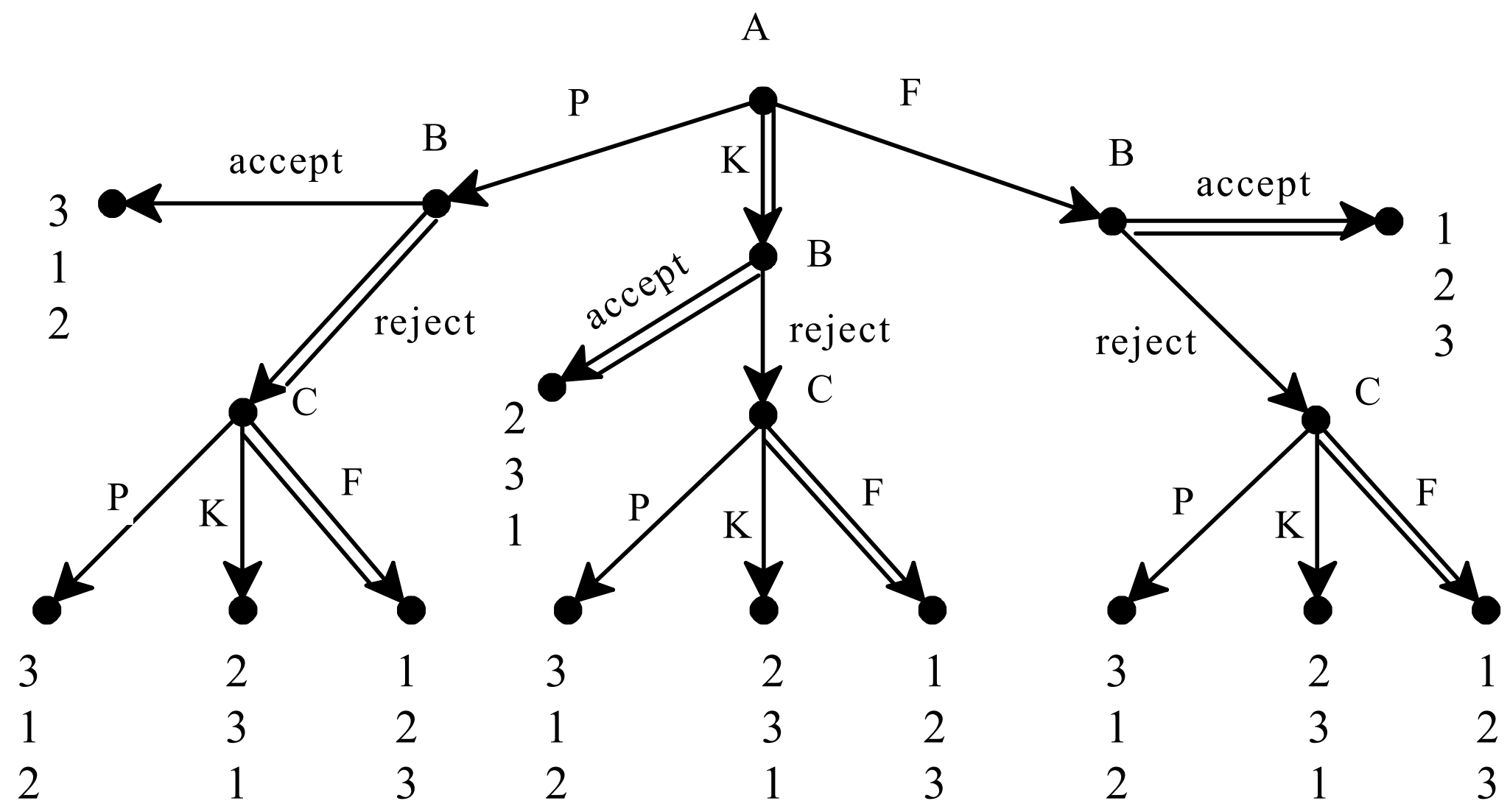

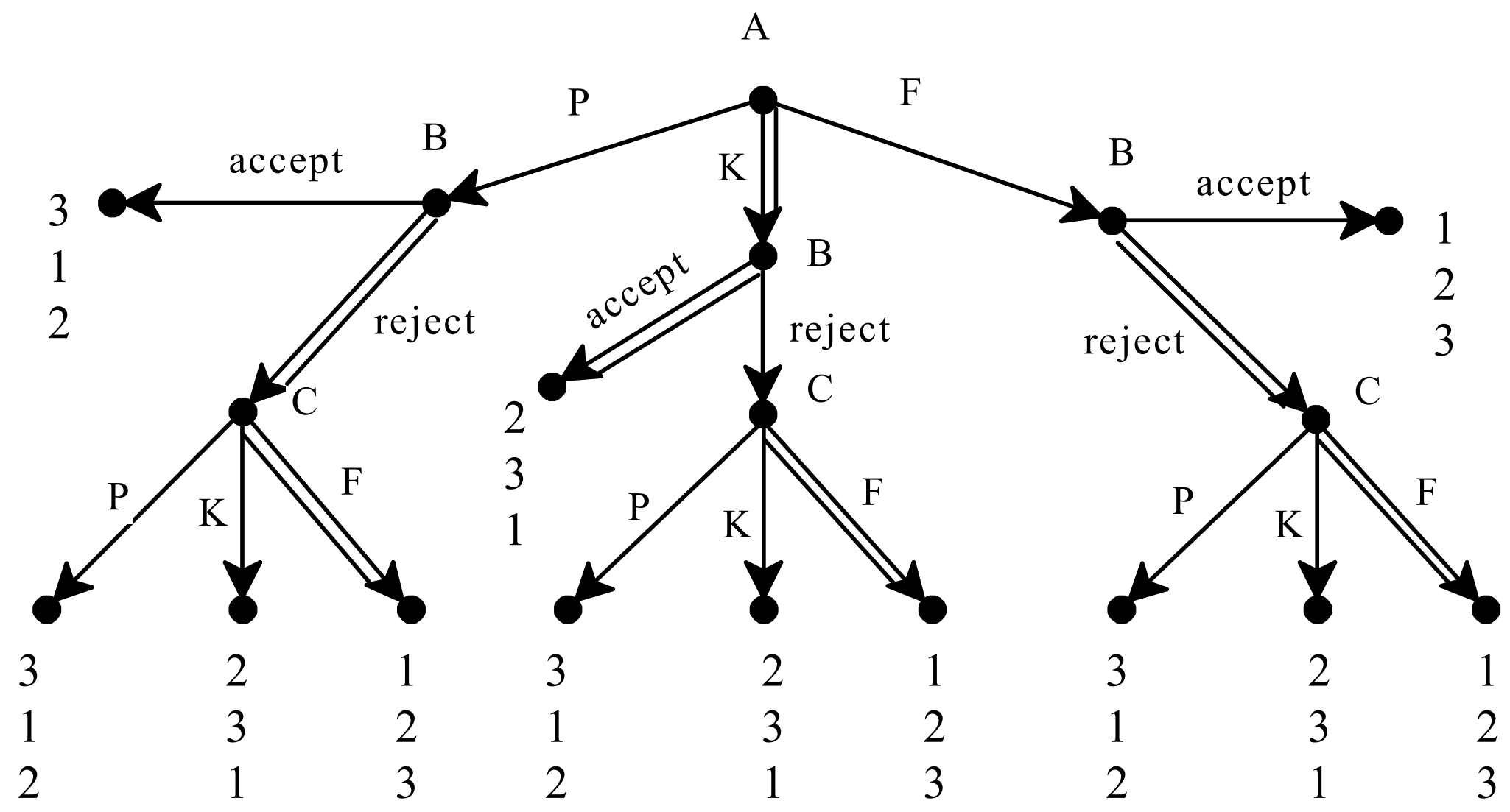

**Exercise 2.5.** The game of Figure 2.2 is reproduced below, with the unique backward-induction solution marked by double edges:

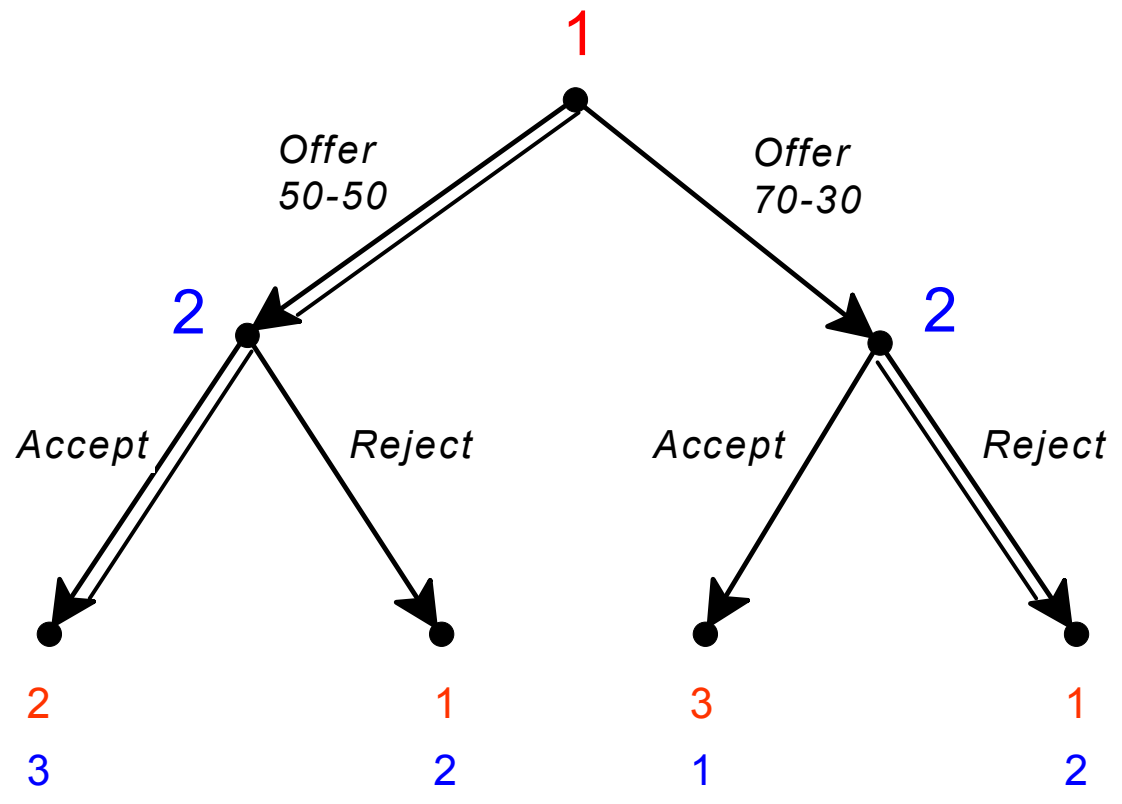

The corresponding strategic form is as follows (for each of Player 2's strategies, the first element in the pair is what Player 2 would do at her left node and the second element what she would do at her right node). The Nash equilibria are highlighted. One Nash equilibrium, namely (Offer 50-50,(Accept,Reject)), corresponds to the backward induction solution, while the other Nash equilibrium, namely (Offer 70-30,( Reject,Reject)) does not correspond to a backward-induction solution.

| | | Player 2 | | | | | | | |
|---|---|---|---|---|---|---|---|---|---|
| | | (Accept,Accept) | | (Accept,Reject) | | (Reject,Accept) | | (Reject,Reject) | |
| Player 1 | offer 50-50 | 2 | 3 | 2 | 3 | 1 | 2 | 1 | 2 |
| | offer 70-30 | 3 | 1 | 1 | 2 | 3 | 1 | 1 | 2 |

**Exercise 2.6.** The game of Figure 2.3 is reproduced below with the two backward-induction solutions marked by double edges.

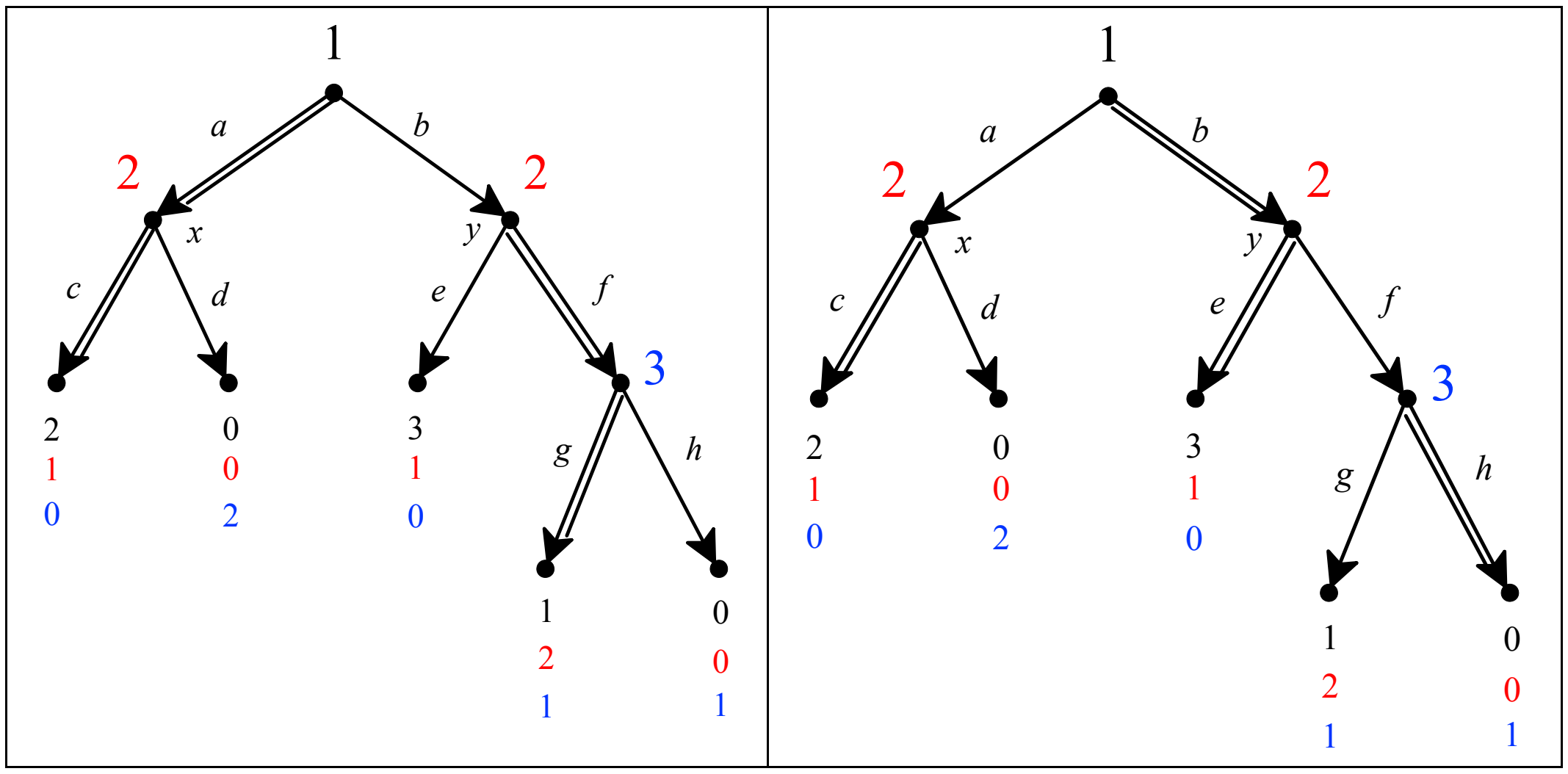

The corresponding strategic form is as follows. The Nash equilibria are highlighted. The backward-induction solutions are $(a,(c,f),g)$ and $(b,(c,e),h)$ and both of them are Nash equilibria. There are three more Nash equilibria which are not backward-induction solutions, namely $(a,(d,f),g)$, $(b,(c,f),h)$ and $(b,(d,e),h)$.

Player 3: *g*

| Player 1 \ Player 2 | *ce* | | | *cf* | | | *de* | | | *df* | | |
|---|---|---|---|---|---|---|---|---|---|---|---|---|
| *a* | 2 | 1 | 0 | 2 | 1 | 0 | 0 | 0 | 2 | 0 | 0 | 2 |
| *b* | 3 | 1 | 0 | 1 | 2 | 1 | 3 | 1 | 0 | 1 | 2 | 1 |

Player 3: *h*

| Player 1 \ Player 2 | *ce* | | | *cf* | | | *de* | | | *df* | | |
|---|---|---|---|---|---|---|---|---|---|---|---|---|
| *a* | 2 | 1 | 0 | 2 | 1 | 0 | 0 | 0 | 2 | 0 | 0 | 2 |
| *b* | 3 | 1 | 0 | 0 | 0 | 1 | 3 | 1 | 0 | 0 | 0 | 1 |

**Exercise 2.7.** The game of Exercise 2.2 is reproduced below:

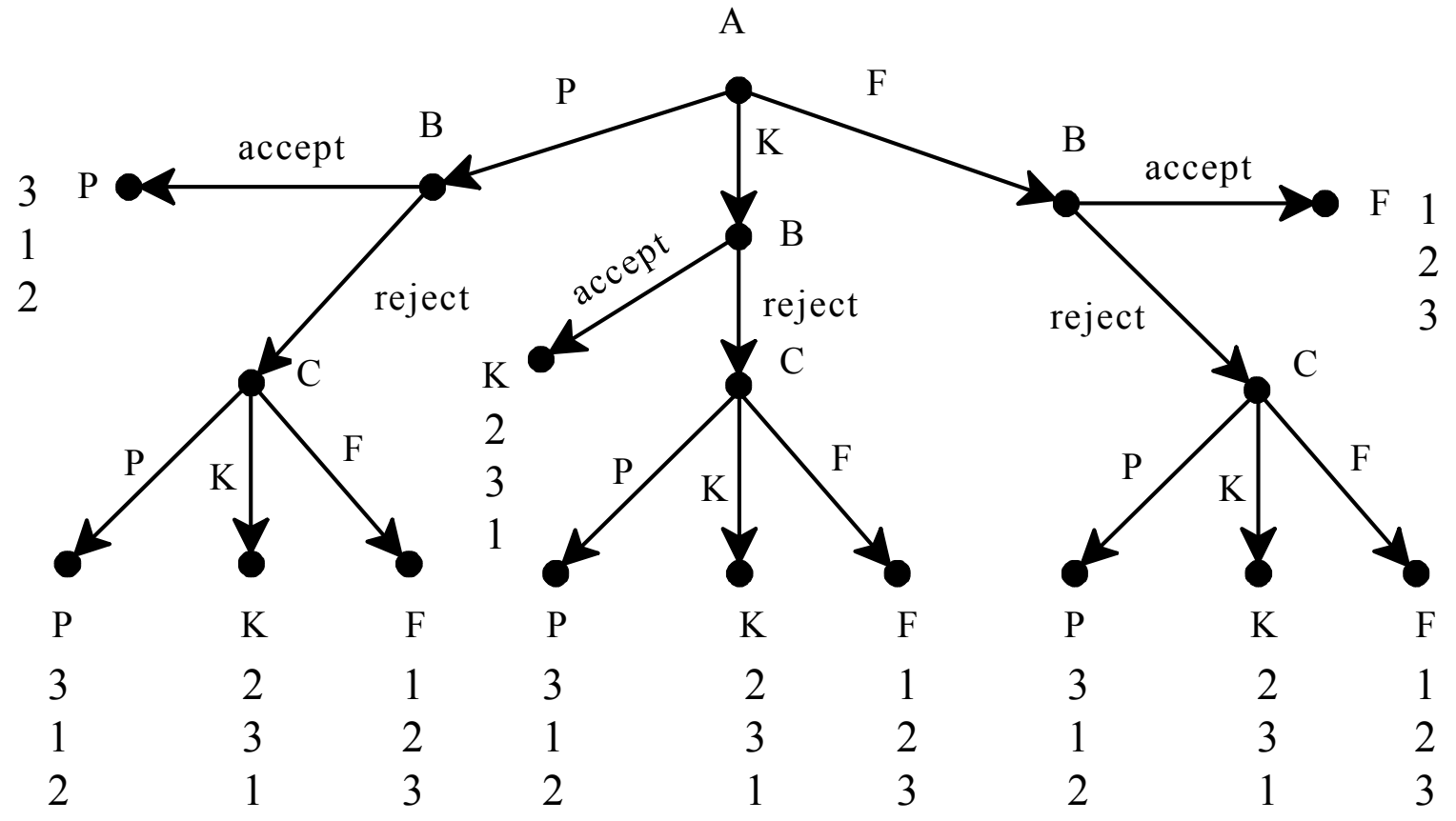

**(a)** All the possible strategies of Player *B* are shown in the following table:

| | If A chooses P | If A chooses K | If A chooses F |
|---|---|---|---|
| 1 | accept | accept | accept |
| 2 | accept | accept | reject |
| 3 | accept | reject | accept |
| 4 | accept | reject | reject |
| 5 | reject | accept | accept |
| 6 | reject | accept | reject |
| 7 | reject | reject | accept |
| 8 | reject | reject | reject |

**(b)** Player *C* has three decision nodes and three choices at each of her nodes. Thus she has $3 \times 3 \times 3 = 27$ strategies.

**Exercise 2.8.**

**(a)** One backward-induction solution is the strategy profile $((L,W),(a,e))$ shown in the following figure by doubling the corresponding edges. The corresponding backward-induction outcome is the play *La* with associated payoff vector (2,1).

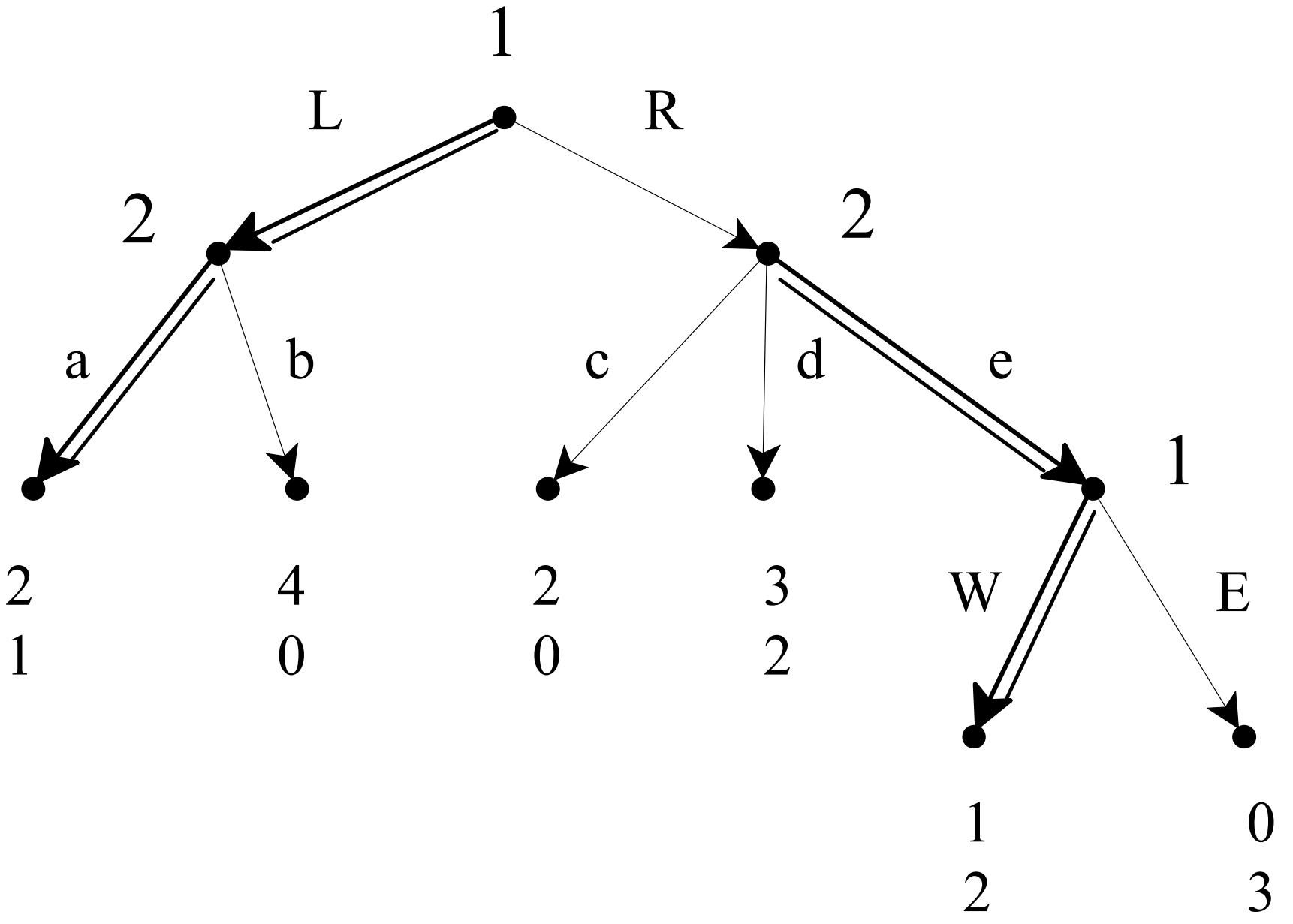

The other backward-induction solution is the strategy profile $((R,W),(a,d))$ shown in the following figure. The corresponding backward-induction outcome is the play *Rd* with associated payoff vector (3,2).

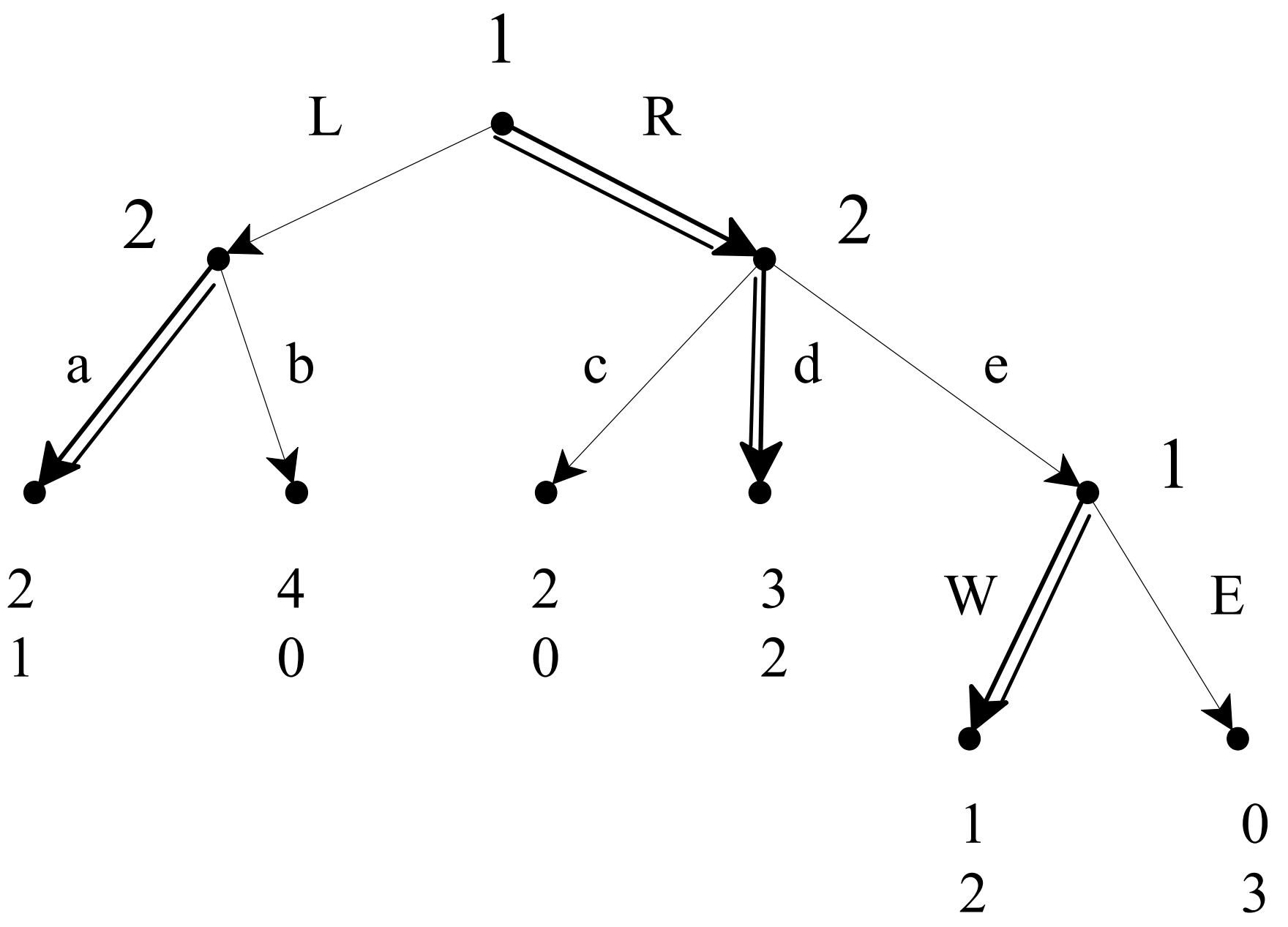

**(b)** Player 1 has four strategies are: *LW, LE, RW, RE*.

**(c)** Player 2 has six strategies: *ac, ad, ae, bc, bd, be*.

**(d)** The strategic form is as follows:

Player 2 (columns), Player 1 (rows)

| | *ac* | *ad* | *ae* | *bc* | *bd* | *be* |
|---|---|---|---|---|---|---|
| *LW* | 2 , 1 | 2 , 1 | 2 , 1 | 4 , 0 | 4 , 0 | 4 , 0 |
| *LE* | 2 , 1 | 2 , 1 | 2 , 1 | 4 , 0 | 4 , 0 | 4 , 0 |
| *RW* | 2 , 0 | 3 , 2 | 1 , 2 | 2 , 0 | 3 , 2 | 1 , 2 |
| *RE* | 2 , 0 | 3 , 2 | 0 , 3 | 2 , 0 | 3 , 2 | 0 , 3 |

**(e)** Player 1 does not have a dominant strategy.

**(f)** For Player 2 *ae* is a weakly dominant strategy.

**(g)** There is no dominant strategy equilibrium.

**(h)** For Player 1 *RE* is weakly dominated by *RW* (and *LW* and *LE* are equivalent).

**(i)** For Player 2 *ac* is weakly dominated by *ad* (or *ae*), *ad* is weakly dominated by *ae*, *bc* is (strictly or weakly) dominated by every other strategy, *bd* is weakly dominated by *be* (and by *ae* and *ad*), *be* is weakly dominated by *ae*. Thus the dominated strategies are: *ac*, *ad*, *bc*, *bd* and *be*.

**(j)** The iterative elimination of weakly dominated strategies yields the following reduced game (in Step 1 eliminate RE for Player 1 and *ac*, *ad*, *bc*, *bd* and *be* for Player 2; in Step 2 eliminate RW for Player 1:

| | *ae* |
|---|---|
| *LW* | 2 , 1 |
| *LE* | 2 , 1 |

Thus we are left with one of the two backward-induction solutions, namely $((L,W),(a,e))$ but also with $((L,E),(a,e))$ which is not a backward-induction solution.

**(k)** The Nash equilibria are highlighted in the following table. There are five Nash equilibria: (*LW*,*ac*), (*LE*,*ac*), (*RW*,*ad*), (*LW*,*ae*) and (*LE*,*ae*).

| 1 \ 2 | *ac* | *ad* | *ae* | *bc* | *bd* | *be* |
|---|---|---|---|---|---|---|
| *LW* | 2 , 1 | 2 , 1 | 2 , 1 | 4 , 0 | 4 , 0 | 4 , 0 |
| *LE* | 2 , 1 | 2 , 1 | 2 , 1 | 4 , 0 | 4 , 0 | 4 , 0 |
| *RW* | 2 , 0 | 3 , 2 | 1 , 2 | 2 , 0 | 3 , 2 | 1 , 2 |
| *RE* | 2 , 0 | 3 , 2 | 0 , 3 | 2 , 0 | 3 , 2 | 0 , 3 |

**Exercise 2.9.** **(a)** The extensive game is as follows:

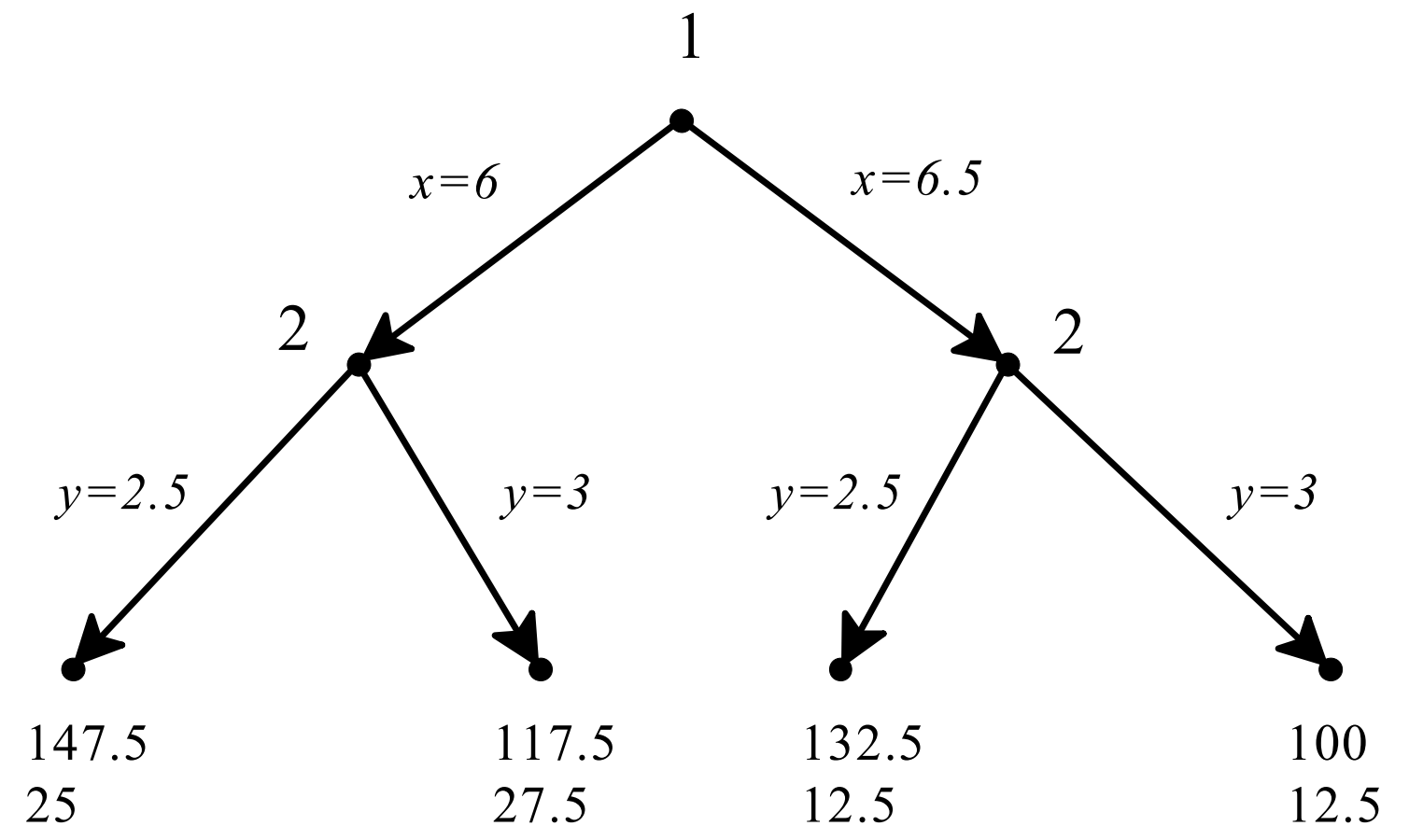

**(b)** There are two backward-induction solutions. The first is the strategy profile $(6,(3,3))$ shown in the following figure. The corresponding backward-induction outcome is given by Firm 1 producing 6 units and Firm 2 producing 3 units with profits 117.5 for Firm 1 and 27.5 for Firm 2.

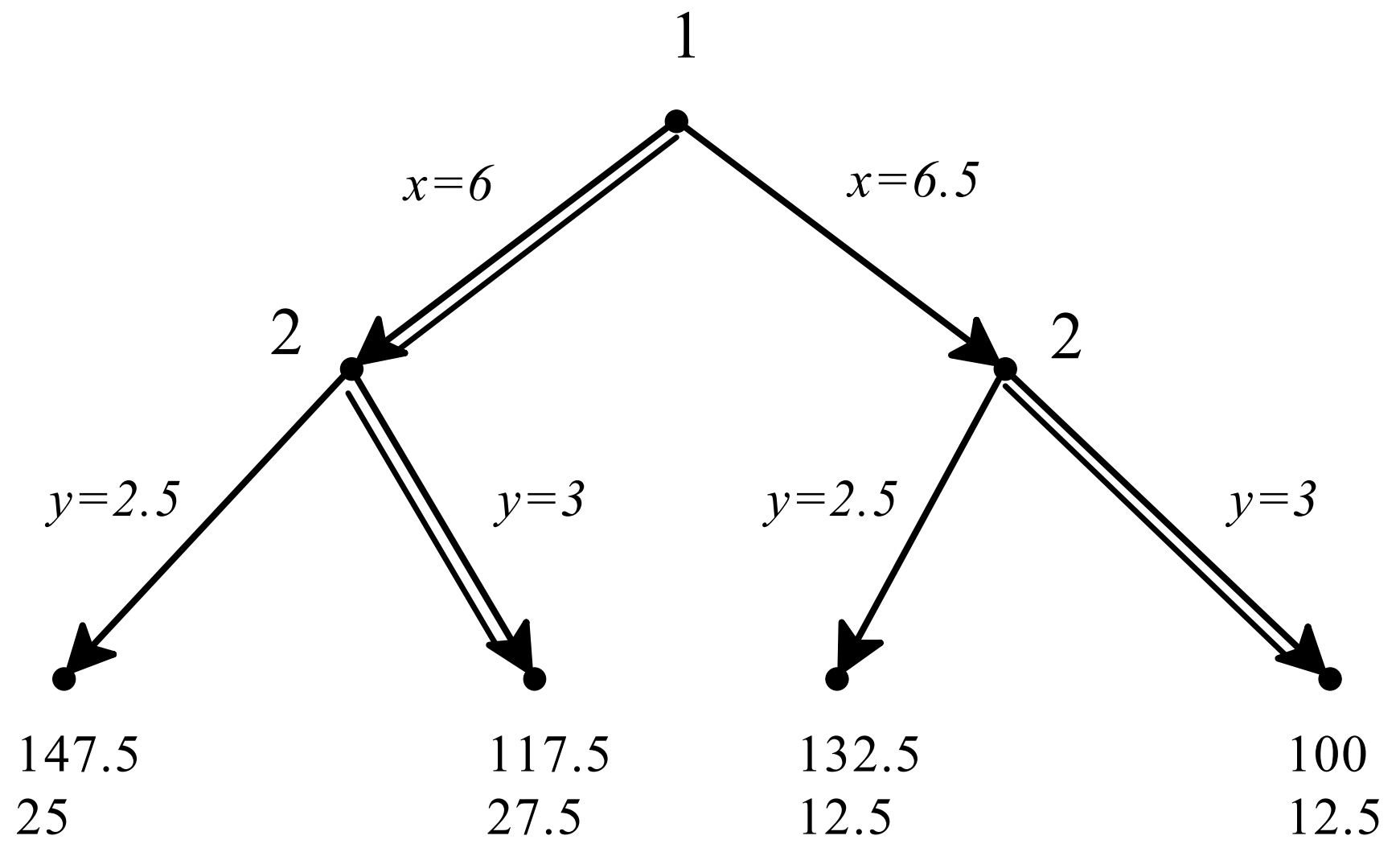

The other backward-induction solution is the strategy profile $(6.5,(3,2.5))$ shown in the following figure. The corresponding backward-induction outcome is given by Firm 1 producing 6.5 units and Firm 2 producing 2.5 units with profits 132.5 for Firm 1 and 12.5 for Firm 2.

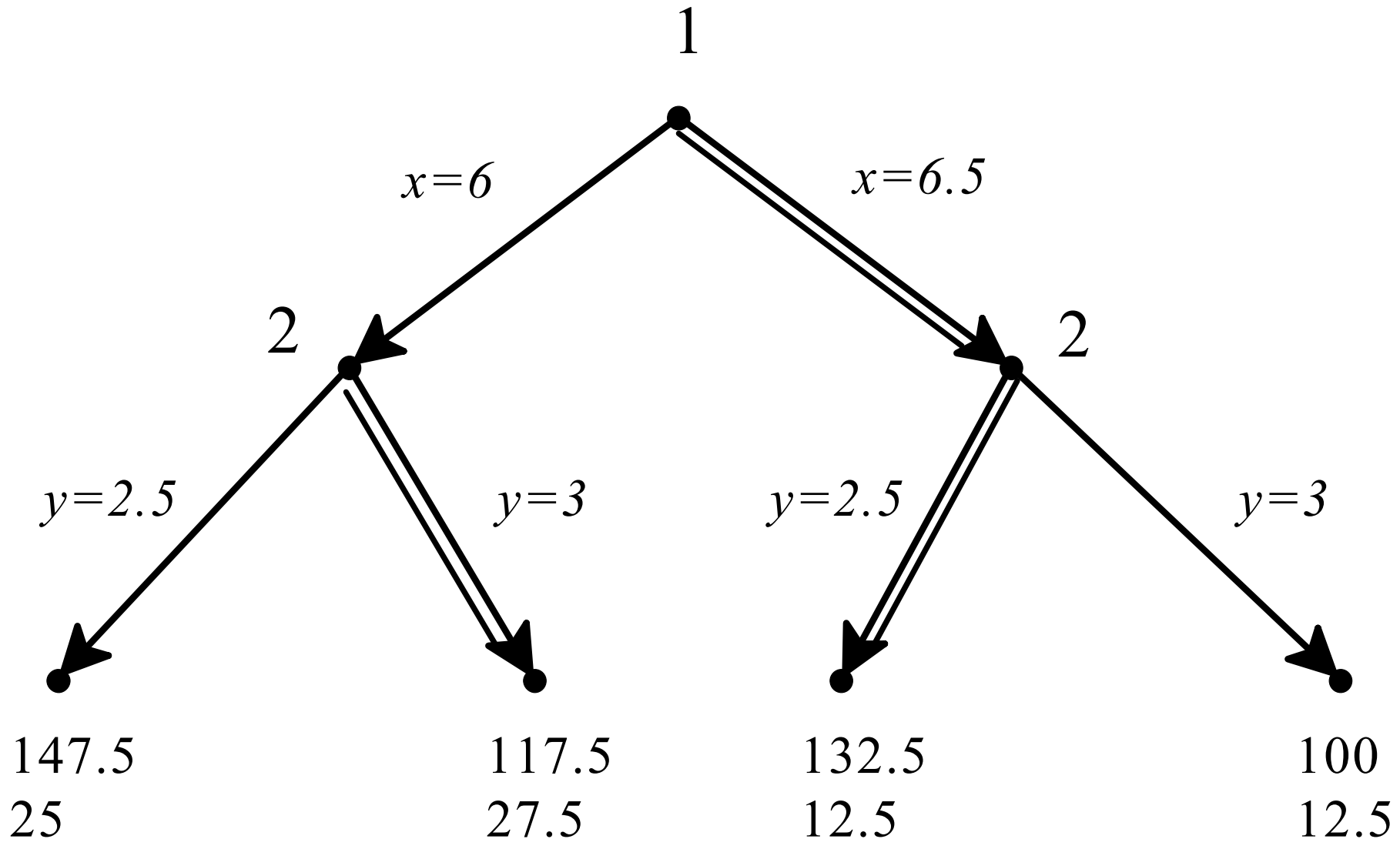

**(c)** The strategic form is as follows.

| | | Firm 2 | | | |
|---|---|---|---|---|---|
| | | (2.5,2.5) | (2.5,3) | (3,2.5) | (3,3) |
| Firm 1 | 6 | 147.5 25 | 147.5 25 | 117.5 27.5 | 117.5 27.5 |
| | 6.5 | 132.5 12.5 | 100 12.5 | 132.5 12.5 | 100 12.5 |

**(d)** The Nash equilibria are highlighted in the above table. In this game the set of Nash equilibria coincides with the set of backward-induction solutions.

**Exercise 2.10.** The game under consideration is the following, where $x$ is an integer:

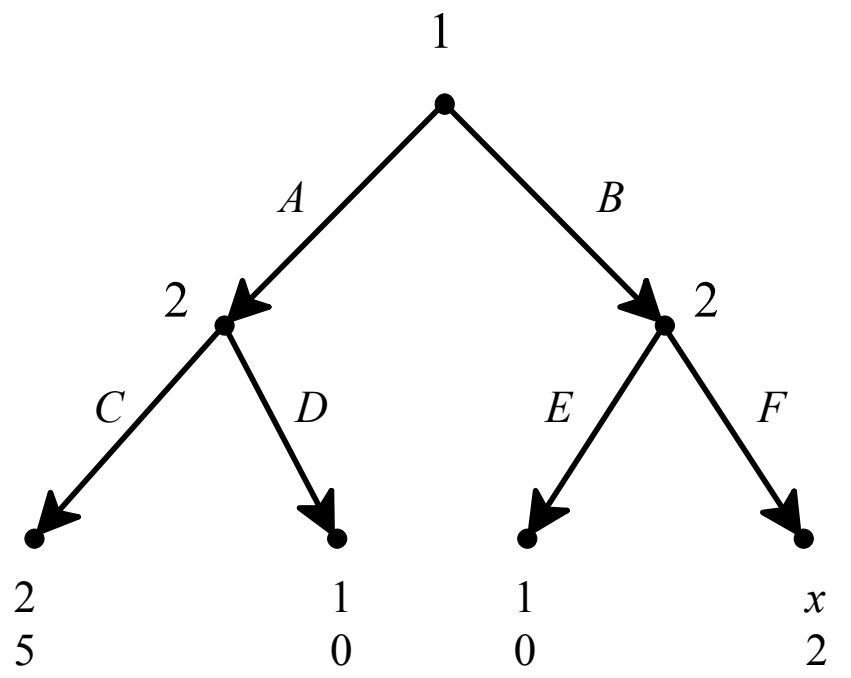

**(a)** The backward-induction strategy of Player 2 is the same, no matter what $x$ is, namely ($C$,$F$). Thus the backward induction solutions are as follows.

- If $x < 2$, there is only one: ($A$, ($C$,$F$)).
- If $x = 2$ there are two: ($A$, ($C$,$F$)) and ($B$, ($C$,$F$)).
- If $x > 2$, there is only one: ($B$, ($C$,$F$)).

**(b)** The strategic form is as follows:

| | | Player 2 | | | |
|---|---|---|---|---|---|
| | | CE | CF | DE | DF |
| Player 1 | A | 2 5 | 2 5 | 1 0 | 1 0 |
| | B | 1 0 | $x$ 2 | 1 0 | $x$ 2 |

First note that **($A$, ($C$,$E$)) is a Nash equilibrium for every value of $x$.** Now, depending on the value of $x$ the other Nash equilibria are as follows:

- If $x < 1$: ($A$, ($C$,$F$)).

- If $1 \leq x < 2$, $(A, (C,F))$ and $(B, (D,F))$.
- If $x = 2$, $(A, (C,F))$, $(B, (C,F))$ and $(B, (D,F))$.
- If $x > 2$, $(B, (C,F))$ and $(B, (D,F))$.

**Exercise 2.11**. Let us find the losing positions. If player $i$, with his choice, can bring the sum to **40** then he can win (the other player with her next choice will take the sum to a number between 41 and 47 and then player $i$ can win with his next choice). Working backwards, the previous losing position is **32** (from here the player who has to move will take the sum to a number between 33 and 39 and after this the opponent can take it to 40). Reasoning backwards, the earlier losing positions are **24**, **16**, **8** and **0**. Thus Player 1 starts from a losing position and therefore it is Player 2 who has a winning strategy. The winning strategy is: at every turn, if Player 1's last choice was $n$ then Player 2 should choose $(8 - n)$.

**Exercise 2.12**. **(a)** The initial part of the game is shown below.

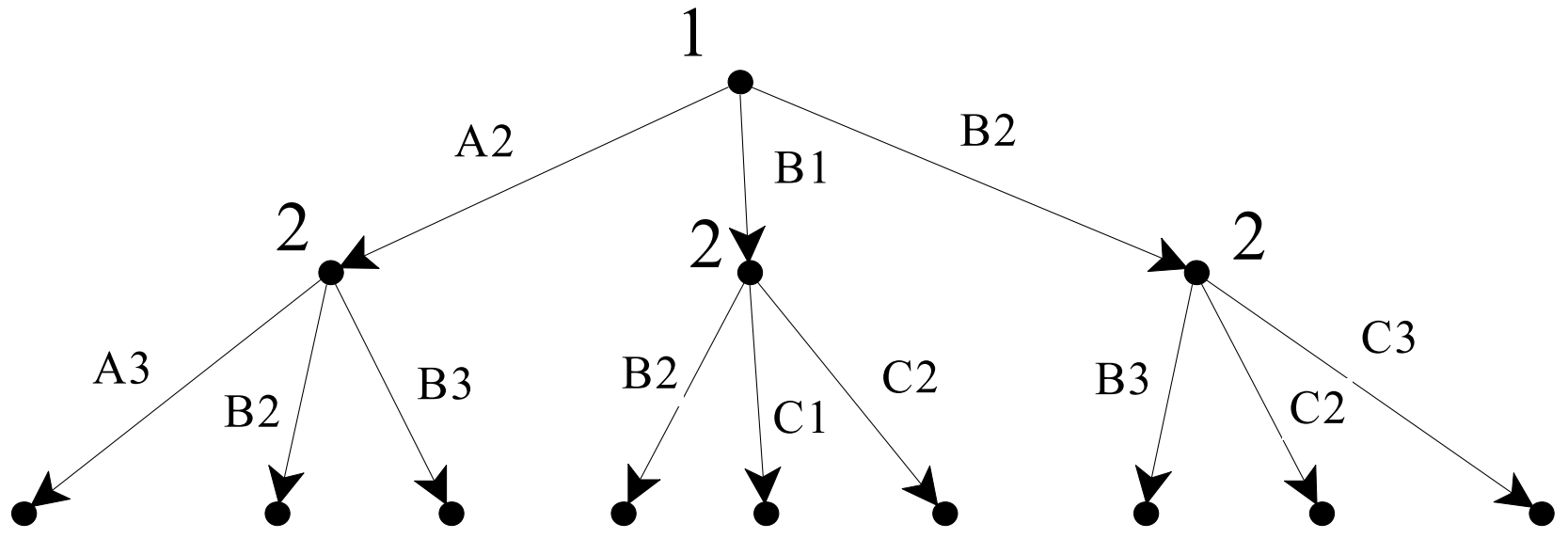

**(b)** From G4 Player 1 should move the coin to H5. From there Player 2 has to move it to H6 and Player 1 to H7 and Player 2 to H8 and from there Player 1 wins by moving it to H9.

**(c)** $A1 \xrightarrow{1} B2 \xrightarrow{2} C3 \xrightarrow{1} D4 \xrightarrow{2} E5 \xrightarrow{1} F5 \xrightarrow{2} G6 \xrightarrow{1} H7 \xrightarrow{2} H8 \xrightarrow{1} H9$

**(d)** $A1 \xrightarrow{1} B2 \xrightarrow{2} C3 \xrightarrow{1} D4 \xrightarrow{2} E5 \xrightarrow{1} F5 \xrightarrow{2} G6 \xrightarrow{1} G7 \xrightarrow{2} H7 \xrightarrow{1} H8 \xrightarrow{2} H9$

**(e)** Using backward induction we can label each cell with a W (meaning that the player who has to move when the coin is there has a winning continuation strategy) or with an L (meaning that the player who has to move when the coin is there can be made to lose). If all the cells that are accessible from a given cell are marked with a W then that cell must be marked with an L. If from a cell there is an accessible cell marked with an L then that cell should be marked with a W. See the following picture.

| | | | | | | | | |
|---|---|---|---|---|---|---|---|---|
| **END** | W | L | W | L | | | | 9 |
| W | W | W | | W | W | L | | 8 |
| L | W | L | W | L | | W | | 7 |
| W | W | | | | W | L | W | 6 |
| L | W | L | W | | L | W | W | 5 |
| W | W | | W | L | W | | L | 4 |
| L | | L | W | W | W | L | W | 3 |
| | W | W | | L | W | W | W | 2 |
| | L | W | L | W | W | L | **START** | 1 |
| H | G | F | E | D | C | B | A | |

From the picture it is clear that it is Player 1 who has a winning strategy. The winning strategy of Player 1 is: move the coin to cell B1 and from then on, after every move of Player 2, move the coin to a cell marked L.

**Exercise 2.13.** **(a)** The game is as follows:

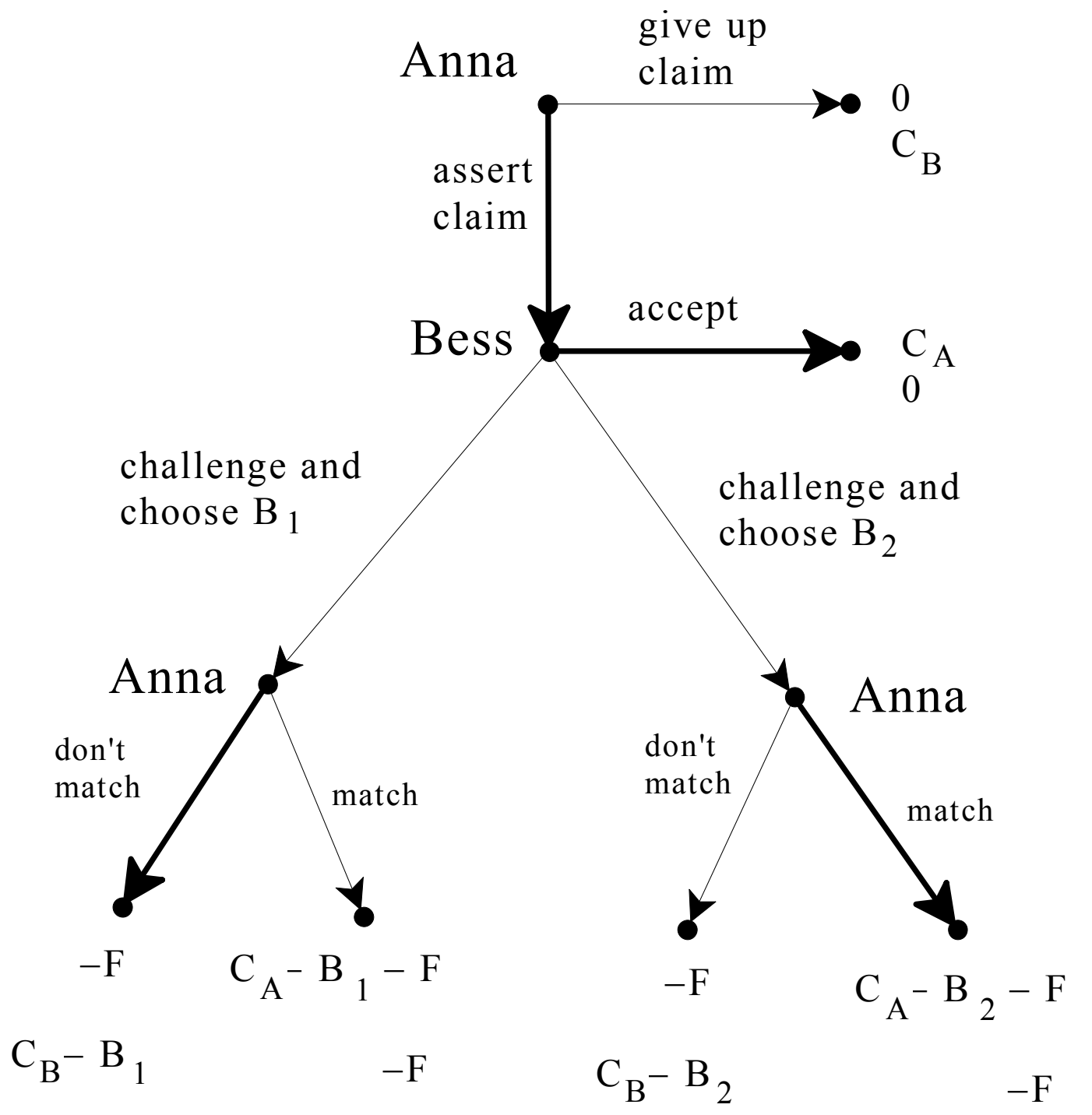

**(b)** The backward-induction solution is marked by thick arrows in the above figure.

**(c)** The structure of the game is as follows:

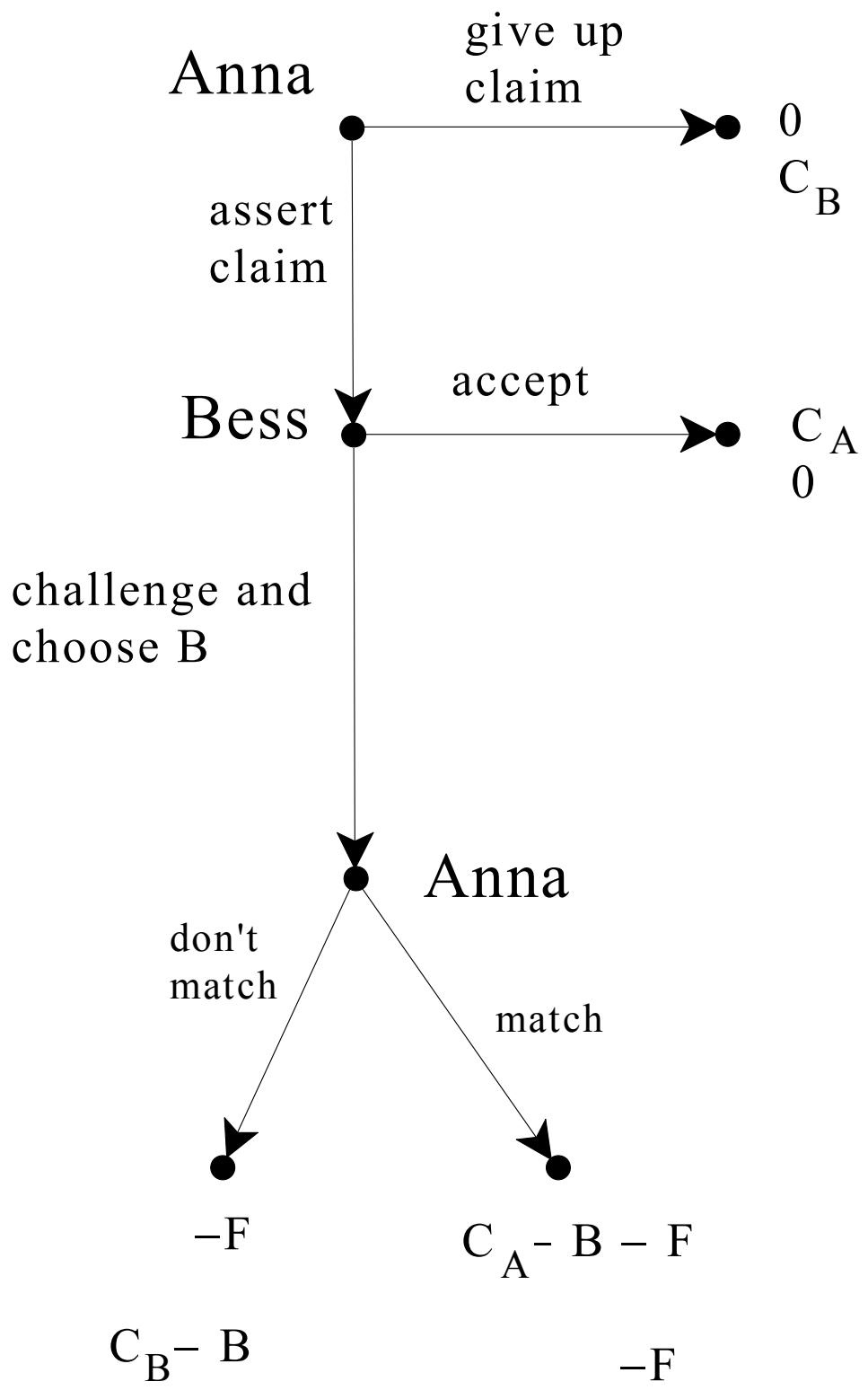

Suppose that Anna is the legal owner and values the ring more than Bess does: $C_A > C_B$. At the last node Anna will choose "match" if $C_A > B$ and "don't match" if $B \geq C_A$. In the first case Bess's payoff will be $-F$, while in the second case it will be $C_B - B$ which is negative since $B \geq C_A$ and $C_A > C_B$. Thus in either case Bess's payoff would be negative. Hence at her decision node Bess will choose "accept" (Bess can get the ring at this stage only if she bids more than the ring is worth to her). Anticipating this, Anna will assert her claim at the first decision node. Thus at the backward-induction solution the ring goes to Anna, the legal owner. The payoffs are $C_A$ for Anna and 0 for Bess. **Note that no money changes hands.**

**(d)** Suppose that Bess is the legal owner and values the ring more than Anna does: $C_B > C_A$. At the last node Anna will choose "match" if $C_A > B$ and "don't match" if $B \geq C_A$. In the first case Bess's payoff will be $-F$, while in the second case it will be $C_B - B$, which will be positive as long as $C_B > B$. Hence at her decision node Bess will choose to challenge and bid any amount B such that $C_B > B > C_A$. Anticipating this, at her first decision node Anna will give up (and get a payoff of 0), because if she asserted her claim then her final payoff would be $-F$. Thus at the backward-induction solution the ring goes to Bess, the legal owner. The payoffs are 0 for Anna and $C_B$ for Bess. **Note that no money changes hands**.

**(e)** As pointed out above, in both cases no money changes hands at the backward-induction solution. Thus Sabio collects no money at all and both Ann and Bess pay nothing.

Chapter 3

# General dynamic games

## 3.1 Imperfect information

There are many situations where players have to make decisions with only partial information about previous moves by other players. Here is an example from my professional experience: in order to discourage copying and cheating in exams, I prepare two versions of the exam, print one version on white paper and the other on pink paper and distribute the exams in such a way that if a student gets, say, the white version then the students on his left and right have the pink version. For simplicity let us assume that there is only one question in the exam. What matters for my purpose is not that the question is indeed different in the two versions, but rather that the students *believe* that they are different and thus refrain from copying from their neighbors. The students, however, are not naïve and realize that I might be bluffing; indeed, introducing differences between the two versions of the exam involves extra effort on my part. Consider a student who finds himself in the embarrassing situation of not having studied for the final exam and is tempted to copy from his neighbor, whom he knows to be a very good student. Let us assume that, if he does not copy, then he turns in a blank exam; in this case, because of his earlier grades in the quarter, he will get a C; on the other hand, if he copies he will get an A if the two versions are identical but will be caught cheating and get an F if the two versions are slightly different. How can we represent such a situation? Clearly this is a situation in which decisions are made sequentially: first the Professor decides whether to write identical versions (albeit printed on different-color paper) or different versions and then the Student chooses between copying and leaving the exam blank. We can easily represent this situation using a tree as we did with the case of perfect-information games, but the crucial element here is the fact that the Student *does not know* whether the two versions are identical or different. In order to represent this uncertainty (or lack of information) in the mind of the Student, we use the notion of *information set*. An information set for a player is a collection of decision nodes of that player and the interpretation is that the player does not know at which of these nodes he is making his decision. Graphically, we represent an information set by enclosing the corresponding nodes in a rounded rectangle. Figure 3.1 represents the situation described above.

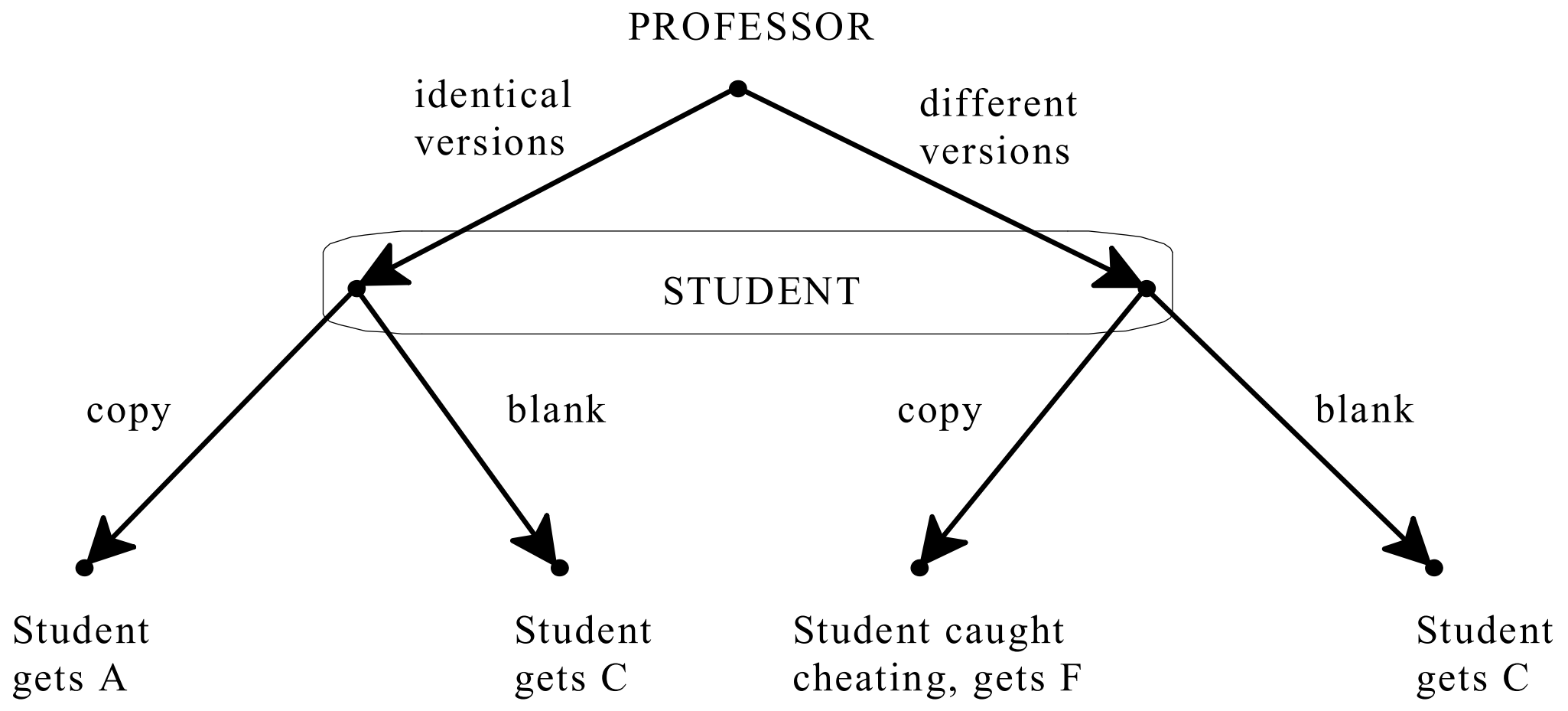

**Figure 3.1**
An extensive form, or frame, with imperfect information

As usual we need to distinguish between a game-frame and a game. Figure 3.1 depicts a game-frame: in order to obtain a game from it we need to add a ranking of the outcomes for each player. For the moment we shall ignore payoffs and focus on frames. A game-frame such as the one shown in Figure 3.1 is called an *extensive form (or frame) with imperfect information*: in this example it is the Student who has imperfect information, or uncertainty, about the earlier decision of the Professor.

We now give a general definition of extensive form that allows for perfect information (Definition 2.1, Chapter 2) as a special case. The first four items of Definition 3.1, marked by the bullet symbol •, coincide with Definition 2.1; what is new is the additional item marked by the symbol ♦.

First some additional terminology and notation. Given a directed tree and two nodes $x$ and $y$ we say that *$y$ is a successor of $x$* or *$x$ is a predecessor of $y$* if there is a sequence of directed edges from $x$ to $y$ (if the sequence consists of a single directed edge then we say that *$y$ is an immediate successor of $x$* or *$x$ is the immediate predecessor of $y$*). A *partition* of a set $H$ is a collection $\mathcal{H} = \{H_1, ..., H_m\}$ $(m \geq 1)$ of non-empty subsets of $H$ such that (1) any two elements of $\mathcal{H}$ are disjoint (if $H_j, H_k \in \mathcal{H}$ with $j \neq k$ then $H_j \cap H_k = \varnothing$) and (2) together the elements of $\mathcal{H}$ cover $H$: $H_1 \cup ... \cup H_m = H$.

**Definition 3.1.** A *finite extensive form (or frame) with perfect recall* consists of the following items.

- A finite rooted directed tree.
- A set of players $I = \{1, \ldots, n\}$ and a function that assigns one player to every decision node.
- A set of actions $A$ and a function that assigns one action to every directed edge, satisfying the restriction that no two edges out of the same node are assigned the same action.
- A set of outcomes $O$ and a function that assigns an outcome to every terminal node.

♦ For every player $i \in I$, a partition $\mathcal{D}_i$ of the set $D_i$ of decision nodes assigned to player $i$ (thus $\mathcal{D}_i$ is a collection of mutually disjoint subsets of $D_i$ whose union is equal to $D_i$). Each element of $\mathcal{D}_i$ is called an *information set of player i.* The elements of $\mathcal{D}_i$ satisfy the following restrictions:

(1) the actions available at any two nodes in the same information set must be the same (that is, for every $D \in \mathcal{D}_i$, if $x, y \in D$ then the outdegree of $x$ is equal to the outdegree of $y$ and the set of actions assigned to the directed edges out of $x$ is equal to the set of actions assigned to the directed edges out of $y$),

(2) if $x$ and $y$ are two nodes in the same information set then it is not the case that one node is a predecessor of the other,

(3) each player has *perfect recall* in the sense that if node $x \in D \in \mathcal{D}_i$ is a predecessor of node $y \in D' \in \mathcal{D}_i$ (thus, by (2), $D \neq D'$), and $a$ is the action assigned to the directed edge out of $x$ in the sequence of edges leading from $x$ to $y$, then for every node $z \in D'$ there is a predecessor $w \in D$ such that the action assigned to the directed edge out of $w$ in the sequence of edges leading from $w$ to $z$ is that same action $a$.

The perfect-recall restriction says that if a player takes action $a$ at an information set and later on has to move again, then at the later time she remembers that she took action $a$ at that earlier information set (because every node she is uncertain about at the later time comes after taking action $a$ at that information set). Perfect recall can be interpreted as requiring that *a player always remember what she knew in the past and what actions she herself took in the past*. Figure 3.2 shows two examples of violation of perfect recall. In the frame shown in Panel (*i*) Player 1 first chooses between $a$ and $b$ and then chooses between $c$ and $d$ having forgotten his previous choice: he does not remember *what* he chose previously. In the frame shown in Panel (*ii*) when Player 2 has to choose between $e$ and $f$ she is uncertain whether this is the first time she moves (left node) or the second time (right node): she is uncertain *whether* she moved in the past.

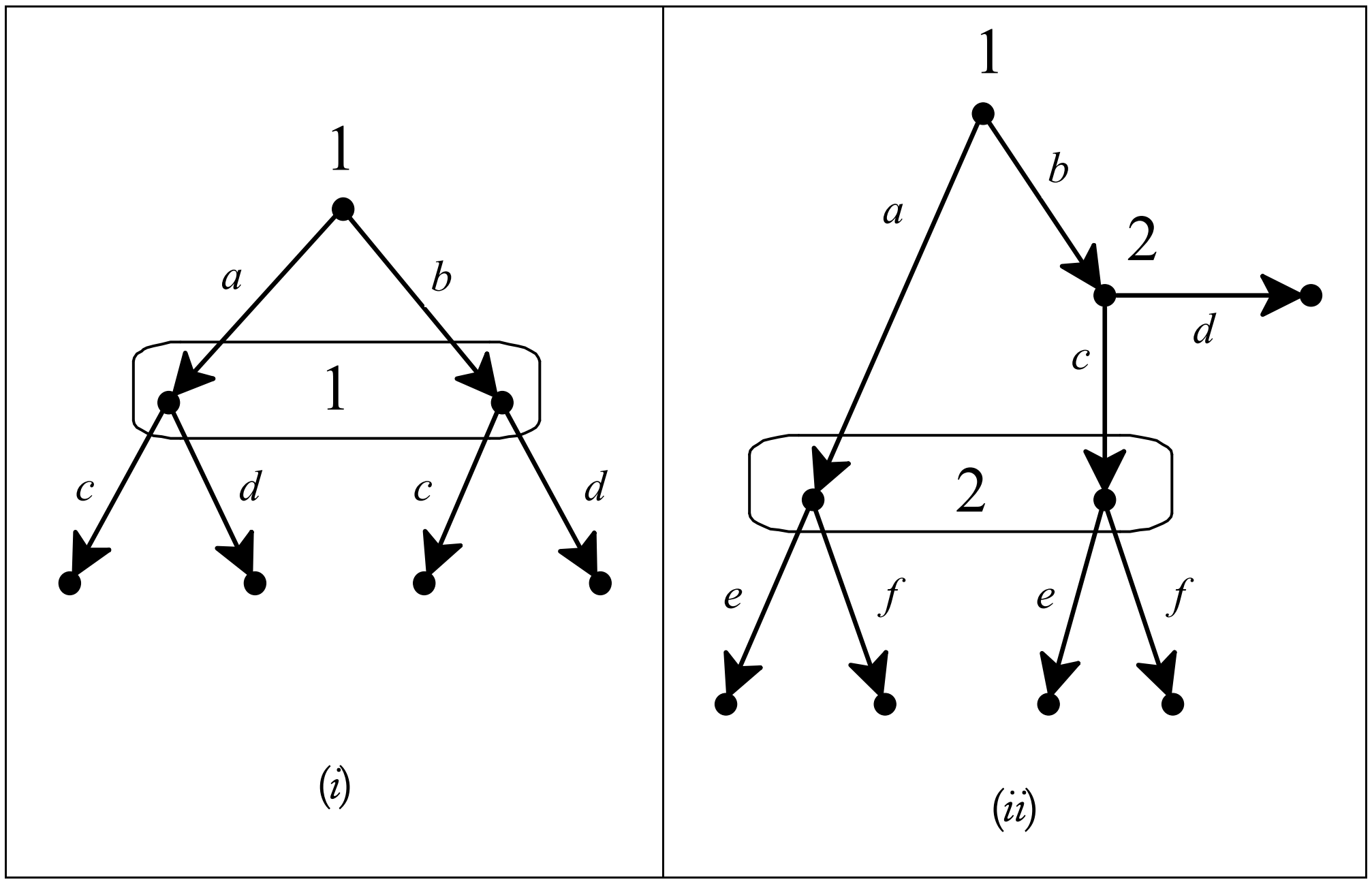

**Figure 3.2**
Examples of violations of perfect recall

If every information set of every player consists of a single node, then the frame is said to be a *perfect-information frame*: it is easy to verify that, in this case, the last item of Definition 3.1 (marked by the symbol ♦) is trivially satisfied and thus Definition 3.1 coincides with Definition 2.1 (Chapter 2). Otherwise (that is, if at least one player has at least one information set that consists of at least two nodes), the frame is said to have *imperfect information*. An example of an extensive frame with imperfect information is the one shown in Figure 3.1. We now give two more examples. In order to simplify the figures, when representing an extensive frame we enclose an information set in a rounded rectangle if and only if that information set contains at least two nodes.

**Example 3.1.** There are three players, Ann, Bob and Carla. Initially, Ann and Bob are in the same room and Carla is outside the room. Ann moves first, chooses either a red card or a black card from a full deck of cards, shows it to Bob and puts it, face down, on the table. Now Carla enters the room and Bob makes a statement to Carla: he either says "Ann chose a Red card" or he says "Ann chose a Black card"; thus Bob could be lying or could be telling the truth. After hearing Bob's statement Carla guesses the color of the card that was picked by Ann. The card is then turned and if Carla's guess was correct then Ann and Bob give $1 each to Carla, otherwise Carla gives $1 to each of Ann and Bob. When drawing an extensive frame to represent this situation, it is important to be careful about what Carla knows, when she makes her guess, and what she is uncertain about. The extensive frame is shown in Figure 3.3.

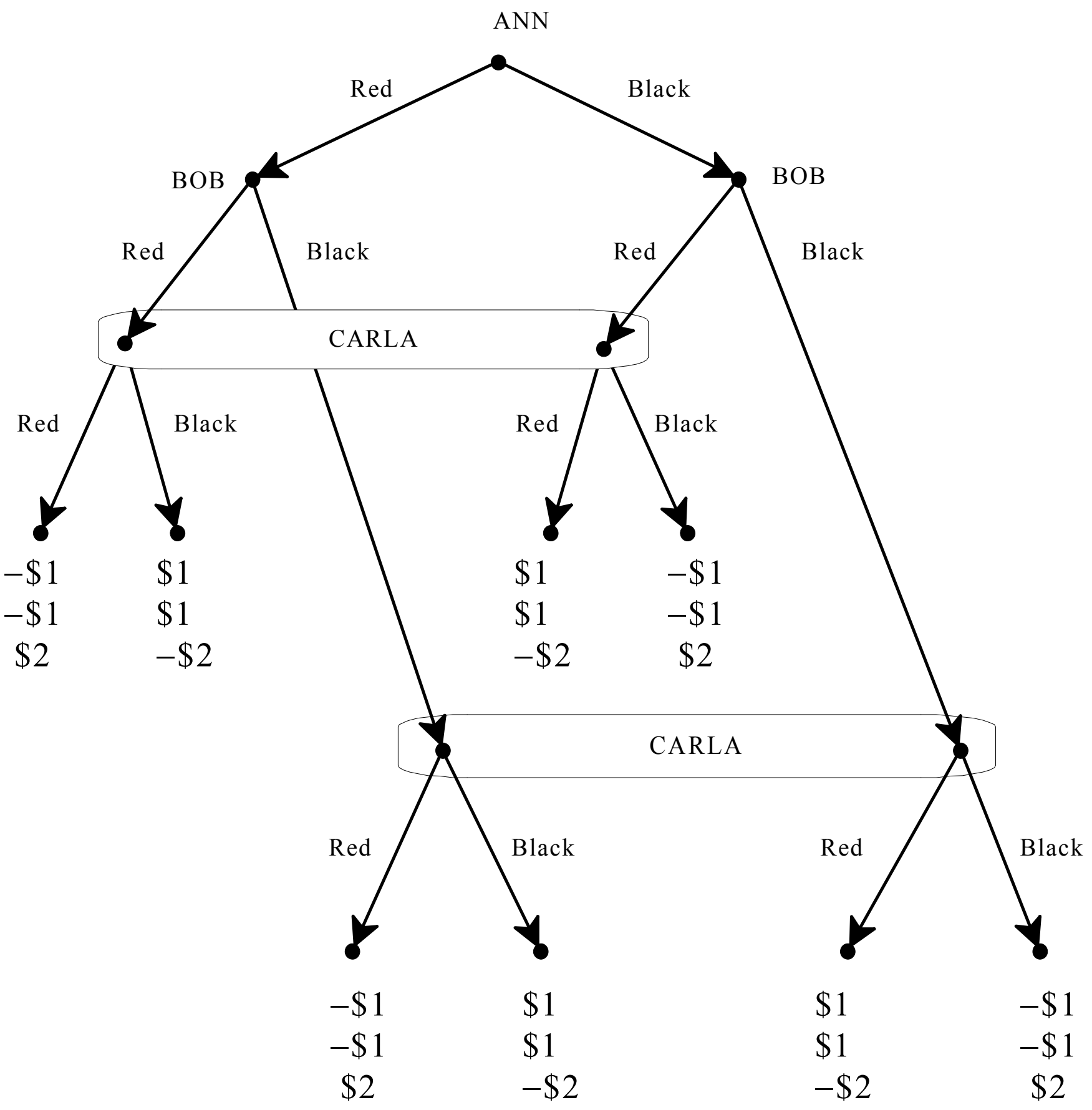

**Figure 3.3**

The extensive form, or frame, representing Example 3.1.

Carla's top information set captures the situation she is in after hearing Bob say "Ann chose a red card" and not knowing if he is telling the truth (left node) or he is lying (right node). Carla's bottom information set captures the alternative situation where she hears Bob say "Ann chose a black card" and does not know if he is lying (left node) or telling the truth (right node). In both situations Carla knows something, namely what Bob tells her, but lacks information about something else, namely what Ann chose. The fact that Bob knows the color of the card chosen by Ann is captured by giving Bob two information sets, each consisting of a single node: Bob's left node represents the situation he is in when he sees that Ann picked a red card, while his right node represents the situation he is in when he sees that Ann picked a black card.

**Example 3.2.** Yvonne and Fran were both interviewed for the same job, but only one person can be hired. The employer told each candidate: "don't call me, I will call you if I want to offer you the job". He also told them that he desperately needs to fill the position and thus, if turned down by one candidate, he will automatically make the offer to the other candidate, without revealing whether he is making a first offer or a "recycled" offer. This situation is represented in the extensive frame shown in Figure 3.4.

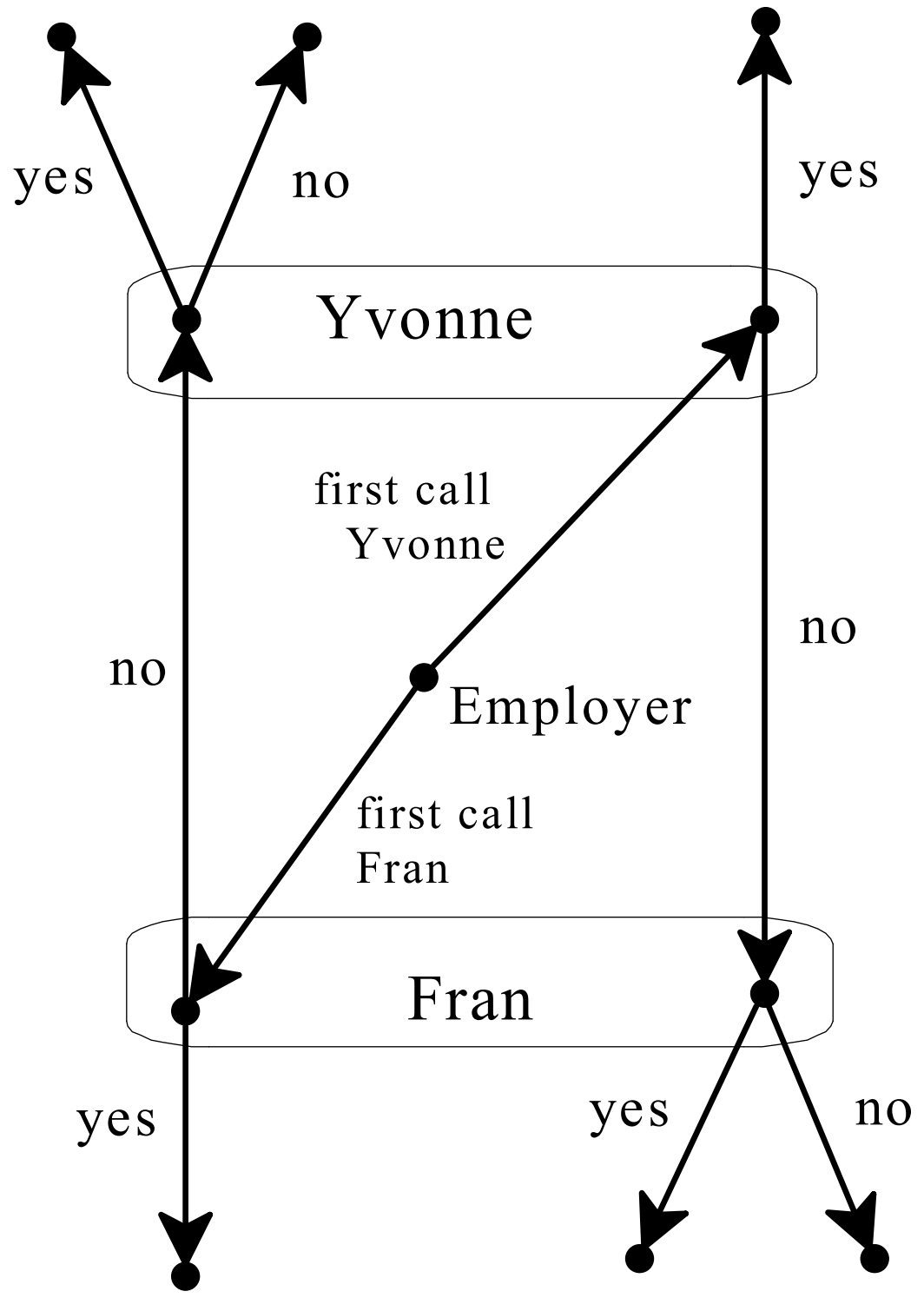

**Figure 3.4**
The extensive form, or frame, representing Example 3.2.

As before, in order to obtain a game from an extensive frame all we need to do is add a ranking of the outcomes for each player. As usual, the best way to represent such rankings is by means of an ordinal utility function for each player and thus represent an extensive-form game by associating a vector of utilities with each terminal node. For instance, expanding on Example 3.2, suppose that the employer only cares about whether the position is filled or not, prefers filling the position to not filling it, but is indifferent between filling it with Yvonne or with Fran; thus we can assign a utility of 1 for the employer to every outcome where one of the two candidates accepts the offer and a utility of 0 to every other outcome. Yvonne's favorite outcome is to be hired if she was the recipient of the first call by the employer; her second best outcome is not to be hired and her worst outcome is to accept a recycled offer (in the latter case Fran would have a blast telling Yvonne "You took that job?! It was offered to me but I turned it down. Who, in her right mind, would want that job? What's wrong with you?!"). Thus for Yvonne we can use utilities of 2 (if she accepts a first offer), 1 (if she is not hired) and 0 (if she accepts a recycled offer). Finally, suppose that Fran has preferences similar (but symmetric) to Yvonne's. Then the extensive frame of Figure 3.4 gives rise to the extensive game shown in Figure 3.5.

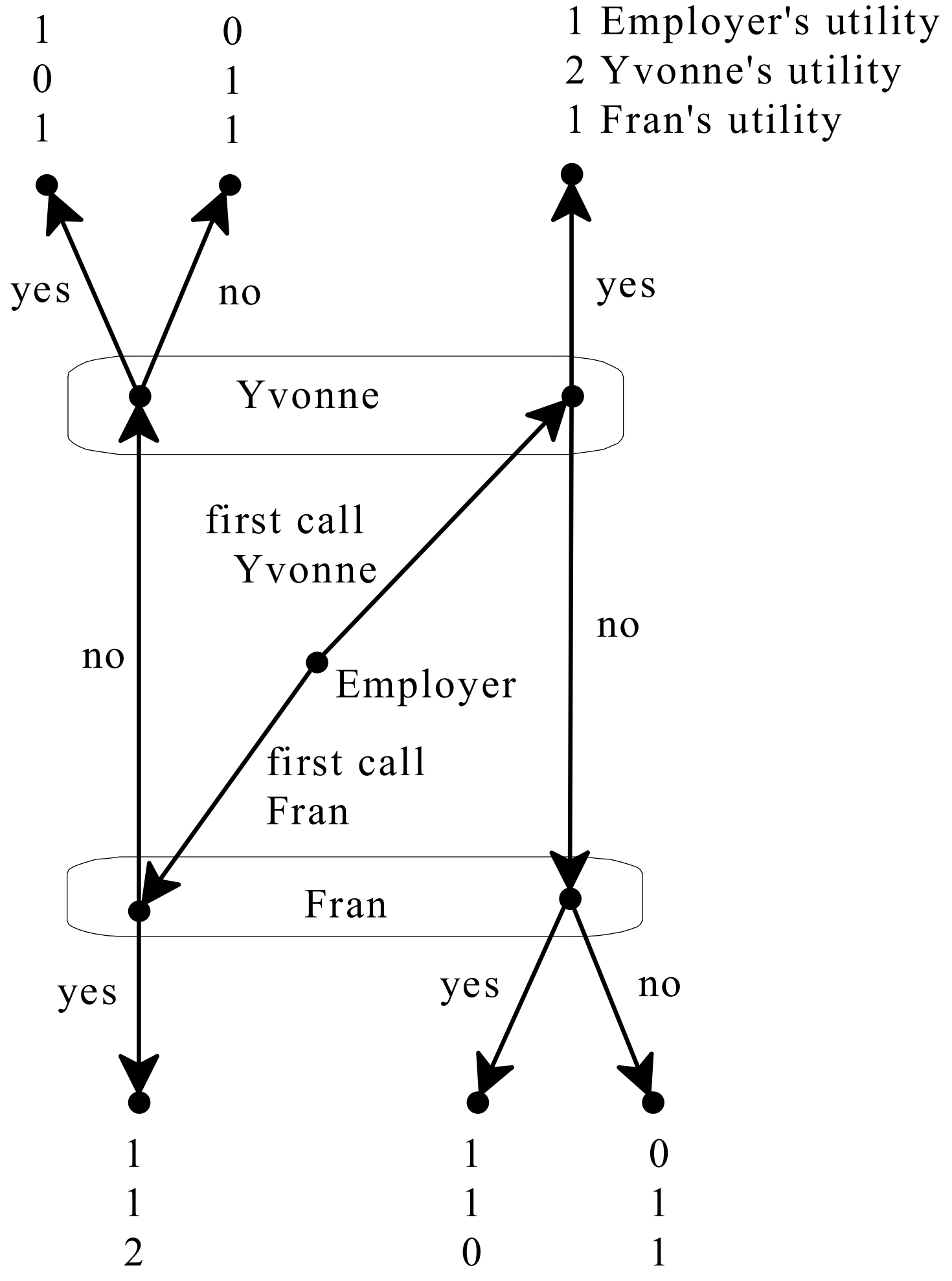

**Figure 3.5**
A game based on the extensive form of Figure 3.4.

This is a good time to test your understanding of the concepts introduced in this section, by going through the exercises in Section 3.E.1 of Appendix 3.E at the end of this chapter.

## 3.2. Strategies

The notion of strategy for general extensive games is the same as before: a strategy for Player $i$ is a complete, contingent plan that covers all the possible situations Player $i$ might find herself in. In the case of a perfect-information game a "possible situation" for a player is a decision node of that player; in the general case, where there may be imperfect information, a "possible situation" for a player is an *information set* of that player. The following definition reduces to Definition 2.4 (Chapter 2) if the game is a perfect-information game (where each information set consists of a single node).

**Definition 3.2.** A *strategy* for a player in an extensive-form game is a list of choices, one for every information set of that player.

For example, in the game of Figure 3.5, Yvonne has only one information set and thus a strategy for her is what to do at that information set, namely either say Yes or say No. Yvonne cannot make the plan "if the employer calls me first I will say Yes and if he calls me second I will say No", because when she receives the call she is not told if this is a first call or a recycled call and thus she cannot make her decision depend on information she does not have.

As in the case of perfect-information games, the notion of strategy allows us to associate with every extensive-form game a strategic-form game. For example, the strategic form associated with the game of Figure 3.5 is shown in Table 3.6 with the Nash equilibria highlighted.

Fran: Yes

| Employer \ Yvonne | Yes | | | No | | |
|---|---|---|---|---|---|---|
| first call Yvonne | 1 | 2 | 1 | 1 | 1 | 0 |
| first call Fran | 1 | 1 | 2 | 1 | 1 | 2 |

Fran: No

| Employer \ Yvonne | Yes | | | No | | |
|---|---|---|---|---|---|---|
| first call Yvonne | 1 | 2 | 1 | 0 | 1 | 1 |
| first call Fran | 1 | 0 | 1 | 0 | 1 | 1 |

**Table 3.6**
The strategic form of the game of Figure 3.5
with the Nash equilibria highlighted.

As another example, consider the extensive form of Figure 3.3 and view it as a game by assuming that each player is selfish and greedy (only cares about how much money he/she gets and prefers more money to less). Then the associated strategic forms is shown in Table 3.7, where Bob's strategy ($x$,$y$) means "I say $x$ if Ann chose a black card and I say $y$ if Ann chose a red card". Thus (R,B) means "if Ann chose a black card I say Red and if Ann chose a red card I say Black" (that is, Bob plans to lie in both cases). Similarly, Carla's strategy ($x$,$y$) means "I guess $x$ if Bob tells me Black and I guess $y$ if Bob tells me Red". Thus (B,R) means "if Bob tells me Black I guess Black and if Bob tells me Red I guess Red" (that is, Carla plans to repeat what Bob says).

| ANN \ BOB | B,B | R,R | B,R | R,B |
|---|---|---|---|---|
| B | -1,-1,2 | -1,-1,2 | -1,-1,2 | -1,-1,2 |
| R | 1,1,-2 | 1,1,-2 | 1,1,-2 | 1,1,-2 |

CARLA: B, B

| ANN \ BOB | B,B | R,R | B,R | R,B |
|---|---|---|---|---|
| B | 1,1,-2 | 1,1,-2 | 1,1,-2 | 1,1,-2 |
| R | -1,-1,2 | -1,-1,2 | -1,-1,2 | -1,-1,2 |

CARLA: R, R

| ANN \ BOB | B,B | R,R | B,R | R,B |
|---|---|---|---|---|
| B | -1,-1,2 | 1,1,-2 | -1,-1,2 | 1,1,-2 |
| R | 1,1,-2 | -1,-1,2 | -1,-1,2 | 1,1,-2 |

CARLA: B, R

| ANN \ BOB | B,B | R,R | B,R | R,B |
|---|---|---|---|---|
| B | 1,1,-2 | -1,-1,2 | 1,1,-2 | -1,-1,2 |
| R | -1,-1,2 | 1,1,-2 | 1,1,-2 | -1,-1,2 |

CARLA: R, B

**Table 3.7**

The strategic form of the game of Figure 3.3.

In order to "solve" an extensive-form game we could simply construct the associated strategic-form game and look for the Nash equilibria. However, we saw in Chapter 2 that in the case of perfect-information games not all Nash equilibria of the associated strategic form can be considered "rational solutions" and we introduced the notion of backward induction to select the "reasonable" Nash equilibria. What we now need is a generalization of the notion of backward induction that can be applied to general extensive-form games. This generalization is called *subgame-perfect equilibrium*. First we need to define the notion of subgame.

This is a good time to test your understanding of the concepts introduced in this section, by going through the exercises in Section 3.E.2 of Appendix 3.E at the end of this chapter.

# 3.3. Subgames

Roughly speaking, a subgame of an extensive-form game is a portion of the game that could be a game in itself. What we need to be precise about is the meaning of "portion of the game".

**Definition 3.3.** A *proper subgame* of an extensive-form game is obtained as follows:

(1) start from a decision node $x$, different from the root, whose information set consists of node $x$ only and enclose in an oval node $x$ itself and all its successors,

(2) if the oval does not "cut" any information sets (that is, there is no information set $S$ and two nodes $y, z \in S$ such that $y$ is a successor of $x$ while $z$ is not) then what is included in the oval is a proper subgame, otherwise it is not.

The reason why we use the qualifier 'proper' is that one could start from the root, in which case one would end up taking the entire game and consider this as a (trivial) subgame (just like any set is a subset of itself; a proper subgame is analogous to a proper subset).

Consider, for example, the extensive-form game of Figure 3.8 below. There are three possible starting points for identifying a proper subgame: nodes $x$, $y$ and $z$ (the other nodes fail to satisfy condition (1) of Definition 3.3).

- Starting from node $x$ and including all of its successors, we do indeed obtain a proper subgame, which is the portion included in the blue oval on the left.

- Starting from node $y$ and including all of its successors we obtain the portion of the game that is included in the red oval on the right; in this case, condition (2) of Definition 3.3 is violated, since we are cutting the top information set of Player 3; hence the portion of the game inside the red oval is *not* a proper subgame.

- Finally, starting from node $z$ and including all of its successors, we do obtain a proper subgame, which is the portion included in the purple oval at the bottom.

Thus the game of Figure 3.8 has two proper subgames.

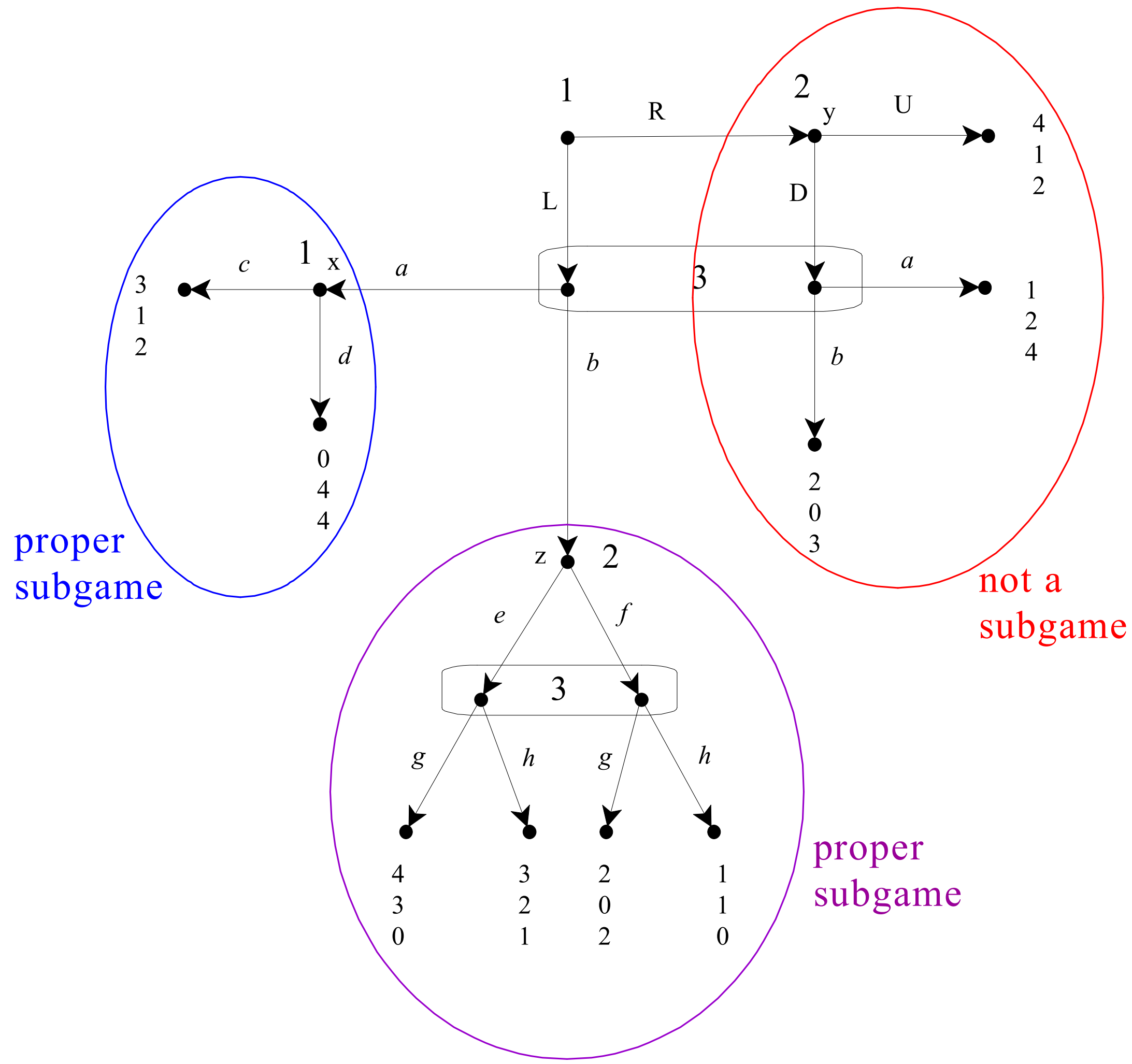

**Figure 3.8**
An extensive-form game and its proper subgames.

**Definition 3.4.** A proper subgame of an extensive-form game is called *minimal* if it does not strictly contain another proper subgame (that is, if there is no other proper subgame which is contained in it and does not coincide with it).

For example, the game shown in Figure 3.9 below has three proper subgames, one starting at node *x*, another at node *y* and the third at node *z*. The ones starting at nodes *x* and *z* are minimal subgames, while the one that starts at node *y* is *not* a minimal subgame, since it contains the one that starts at node *z*.

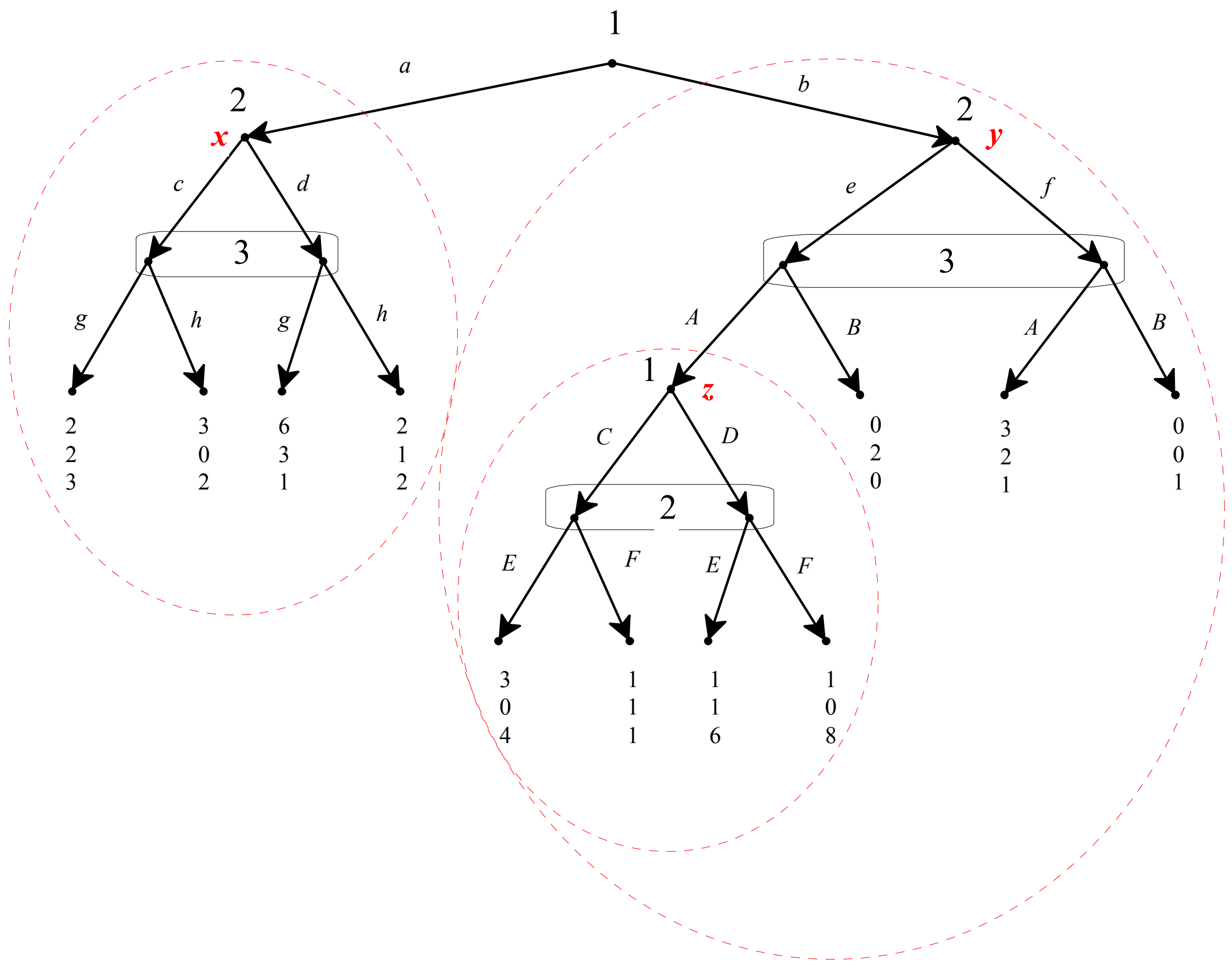

**Figure 3.9**

An extensive-form game with three proper subgames, two of which are minimal.

This is a good time to test your understanding of the concepts introduced in this section, by going through the exercises in Section 3.E.3 of Appendix 3.E at the end of this chapter.

# 3.4. Subgame-perfect equilibrium

A subgame-perfect equilibrium of an extensive-form game is a Nash equilibrium of the entire game which remains an equilibrium in every proper subgame. Consider an extensive-form game and let *s* be a strategy profile for that game. Let *G* be a proper subgame. Then the *restriction of s to G*, denoted by $s|_G$, is that part of *s* which prescribes choices at every information set of *G* and only at those information sets. For example, consider the extensive-form game of Figure 3.9 above and the strategy profile $\left( \underbrace{(a,C)}_{\text{1's strategy}}, \underbrace{(d,f,E)}_{\text{2's strategy}}, \underbrace{(h,B)}_{\text{3's strategy}} \right)$. Let *G* be the subgame that starts at node *y* of Player 2. Then $s|_G = \left( \underbrace{C}_{\substack{\text{1's strategy}\\ \text{in G}}}, \underbrace{(f,E)}_{\substack{\text{2's strategy}\\ \text{in G}}}, \underbrace{B}_{\substack{\text{3's strategy}\\ \text{in G}}} \right)$.

**Definition 3.5.** Given an extensive-form game, let *s* be a strategy profile for the entire game. Then *s* is a *subgame-perfect equilibrium* if

(1) *s* is a Nash equilibrium of the entire game and

(2) for every proper subgame *G*, $s|_G$ (the restriction of *s* to *G*) is a Nash equilibrium of *G*.

For example, consider again the extensive-form game of Figure 3.9 and the strategy profile $s = ((a,C),(d,f,E),(h,B))$. Then *s* is a Nash equilibrium of the entire game: Player 1's payoff is 2 and, if he were to switch to any strategy where he plays *b*, his payoff would be 0; Player 2's payoff is 1 and if she were to switch to any strategy where she plays *c* her payoff would be 0; Player 3's payoff is 2 and if he were to switch to any strategy where he plays *g* his payoff would be 0. However, *s* is not a subgame-perfect equilibrium, because the restriction of *s* to the proper subgame that starts at node *z* of Player 1, namely $(C,E)$, is not a Nash equilibrium of that subgame: in that subgame, for Player 2 the unique best reply to *C* is *F*.

One way of finding the subgame-perfect equilibria of a given game is to first find the Nash equilibria and then, for each of them, check if it satisfies condition (2) of Definition 3.5. However, this is not a practical way to proceed. A quicker and easier way is to apply the following algorithm, which generalizes the backward-induction algorithm for games with perfect information (Definition 2.3, Chapter 2).

**Definition 3.6.** Given an extensive-form game, the *subgame-perfect equilibrium algorithm* is the following procedure.

1. Start with a minimal proper subgame and select a Nash equilibrium of it.
2. Delete the selected proper subgame and replace it with the payoff vector associated with the selected Nash equilibrium, making a note of the strategies that constitute the Nash equilibrium. This yields a smaller extensive-form game.
3. Repeat Steps 1 and 2 in the smaller game so obtained.

For example, let us apply the algorithm to the game of Figure 3.9. Begin with the proper subgame that starts at node $x$ of Player 2, shown below with its associated strategic form, where the unique Nash equilibrium $(d,h)$ is highlighted. Note that this is a game only between Players 2 and 3 and thus in Figure 3.10 we only show the payoffs of these two players.

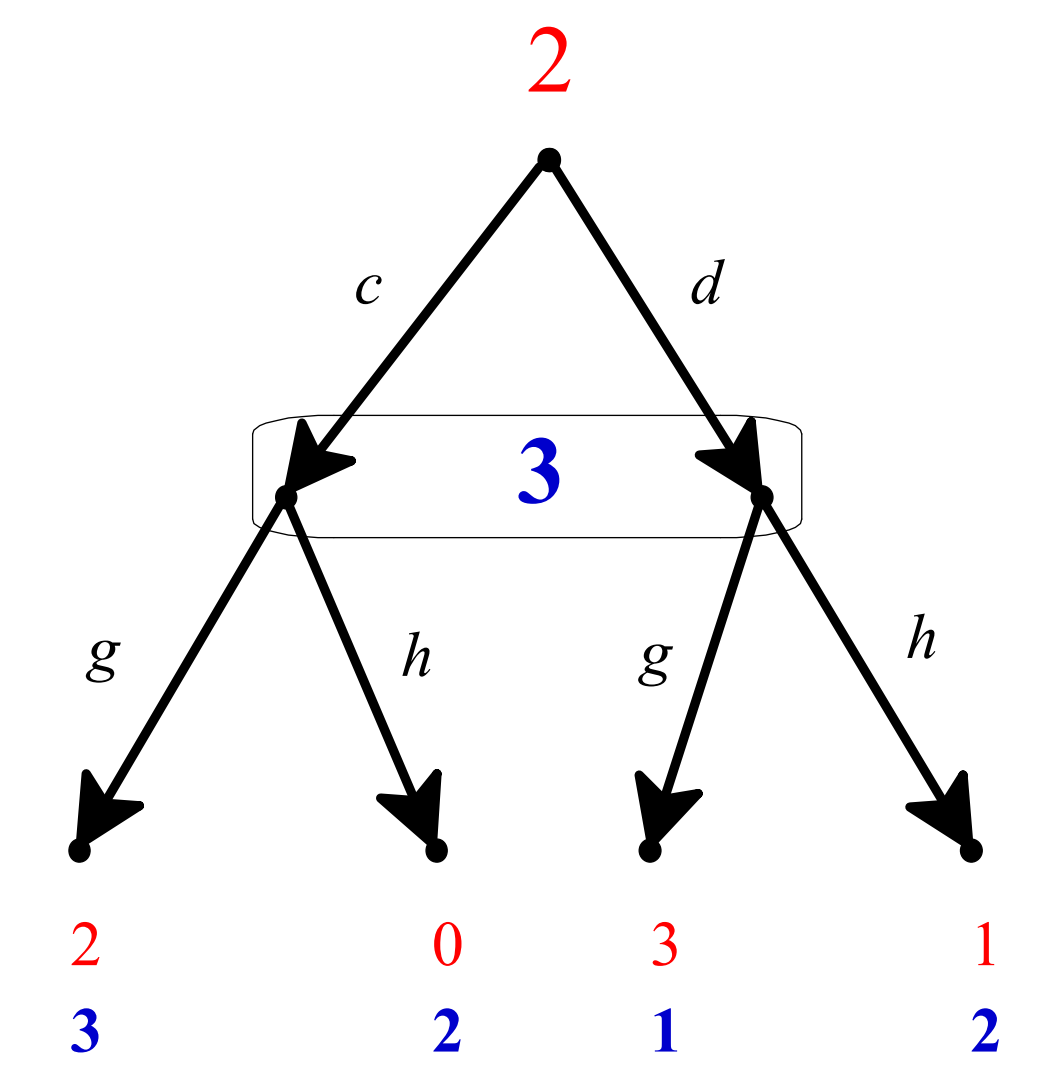

| | | Player 3: $g$ | Player 3: $h$ |
|---|---|---|---|
| Player 2 | $c$ | 2, 3 | 0, 2 |
| | $d$ | 3, 1 | 1, 2 |

**Figure 3.10**
A minimal proper subgame of the game of Figure 3.9 and its strategic form.

Now we delete the proper subgame, thereby turning node $x$ into a terminal node to which we attach the full payoff vector associated, in the original game, with the terminal node following history *adh*, namely (2,1,2). Hence we obtain the smaller game shown in Figure 3.11.

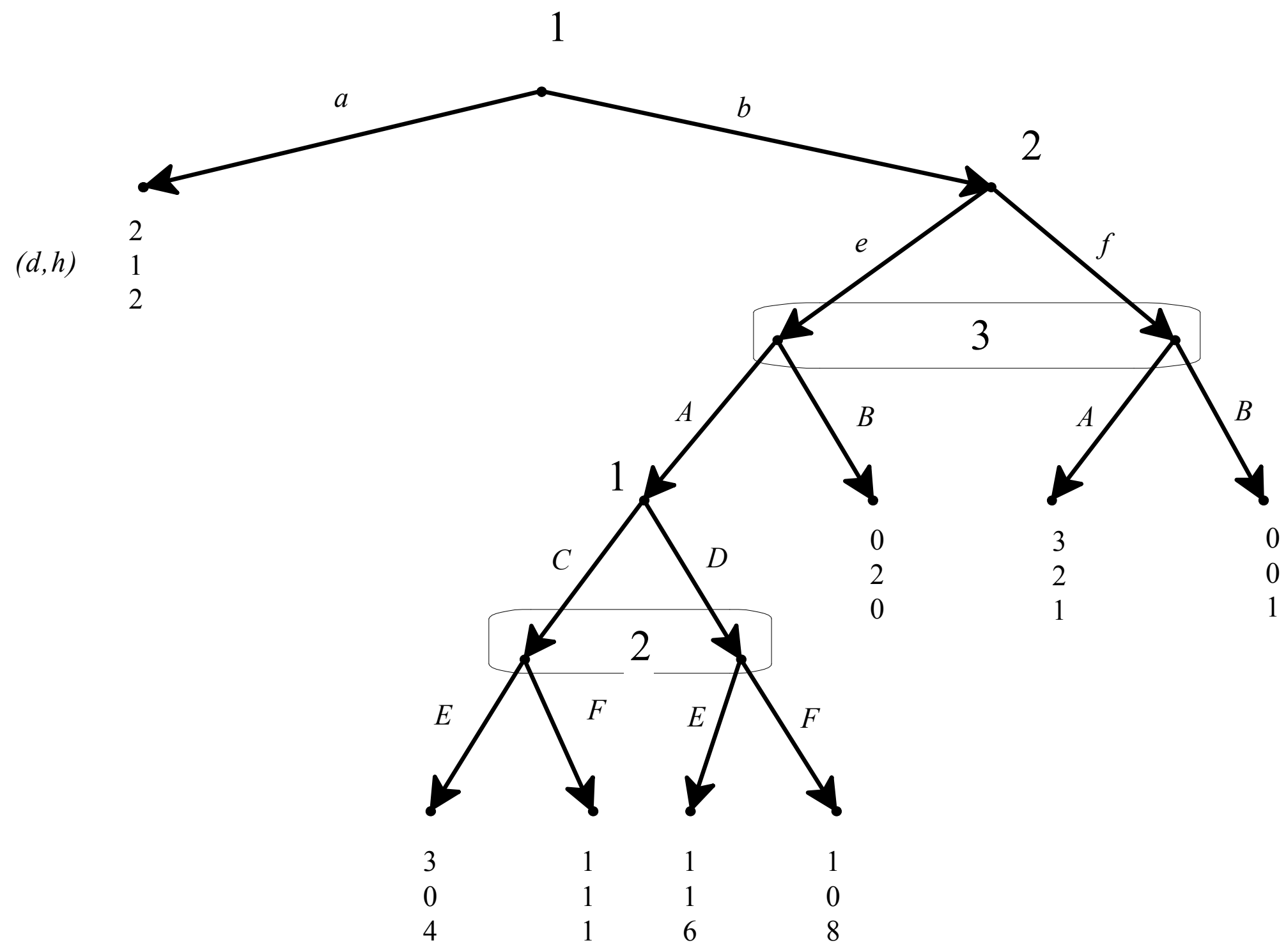

**Figure 3.11**
The reduced game after replacing a proper minimal subgame in the game of Figure 3.9.

Now, in the reduced game of Figure 3.11 we select the only minimal proper subgame, namely the one that starts at the bottom decision node of Player 1. This subgame is shown in Figure 3.12 together with its associated strategic form. The unique Nash equilibrium of that subgame is $(C,F)$.

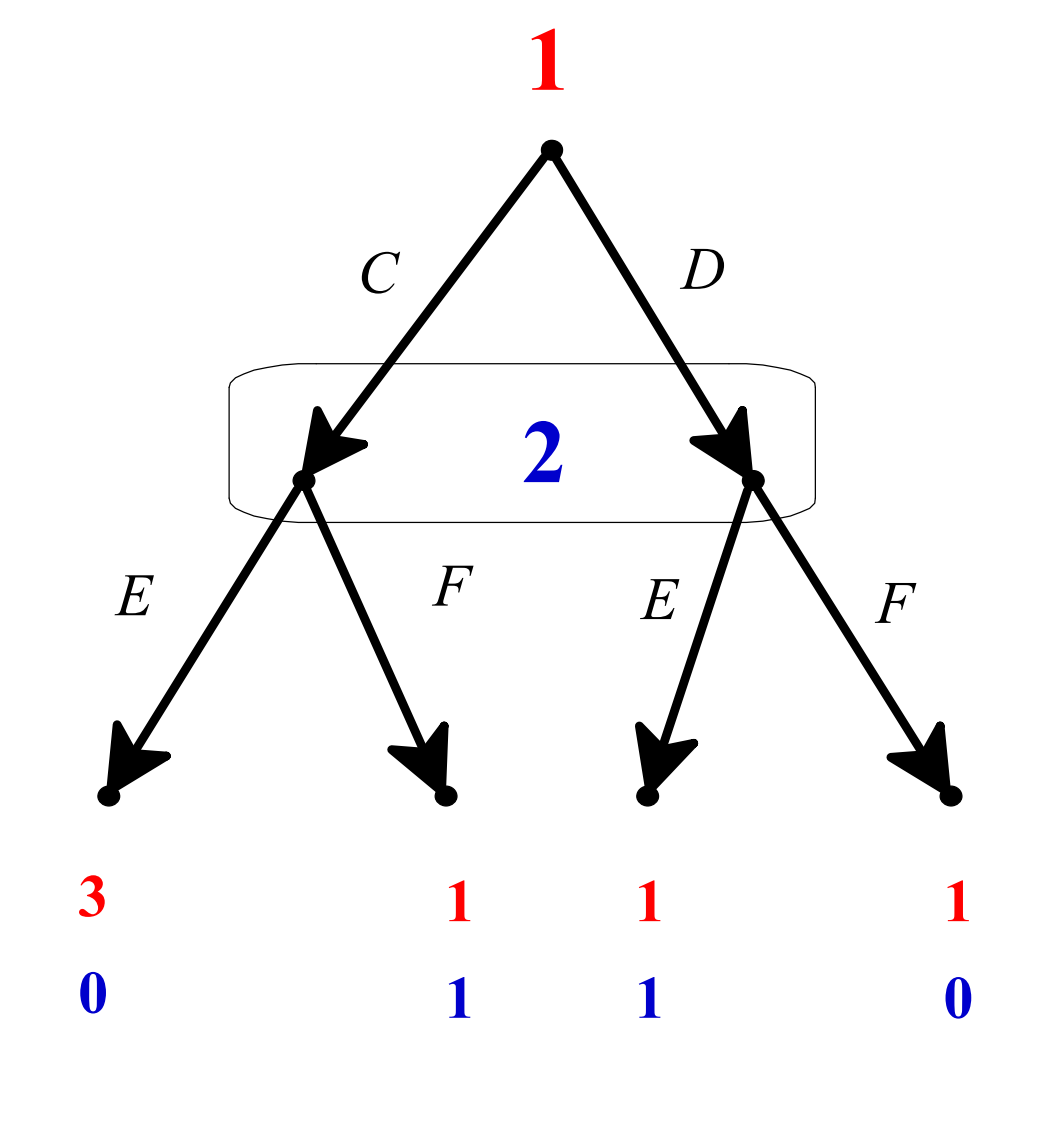

| | | Player 3: *E* | Player 3: *F* |
|---|---|---|---|
| **Player 2** | *C* | 3 , 0 | 1 , 1 |
| | *D* | 1 , 1 | 1 , 0 |

**Figure 3.12**
A minimal proper subgame of the game of Figure 3.11 and its strategic form.

Then, in the reduced game of Figure 3.11, we replace the selected proper subgame with the payoff vector associated with the history *beACF*, namely (1,1,1), thus obtaining the smaller game shown in Figure 3.13.

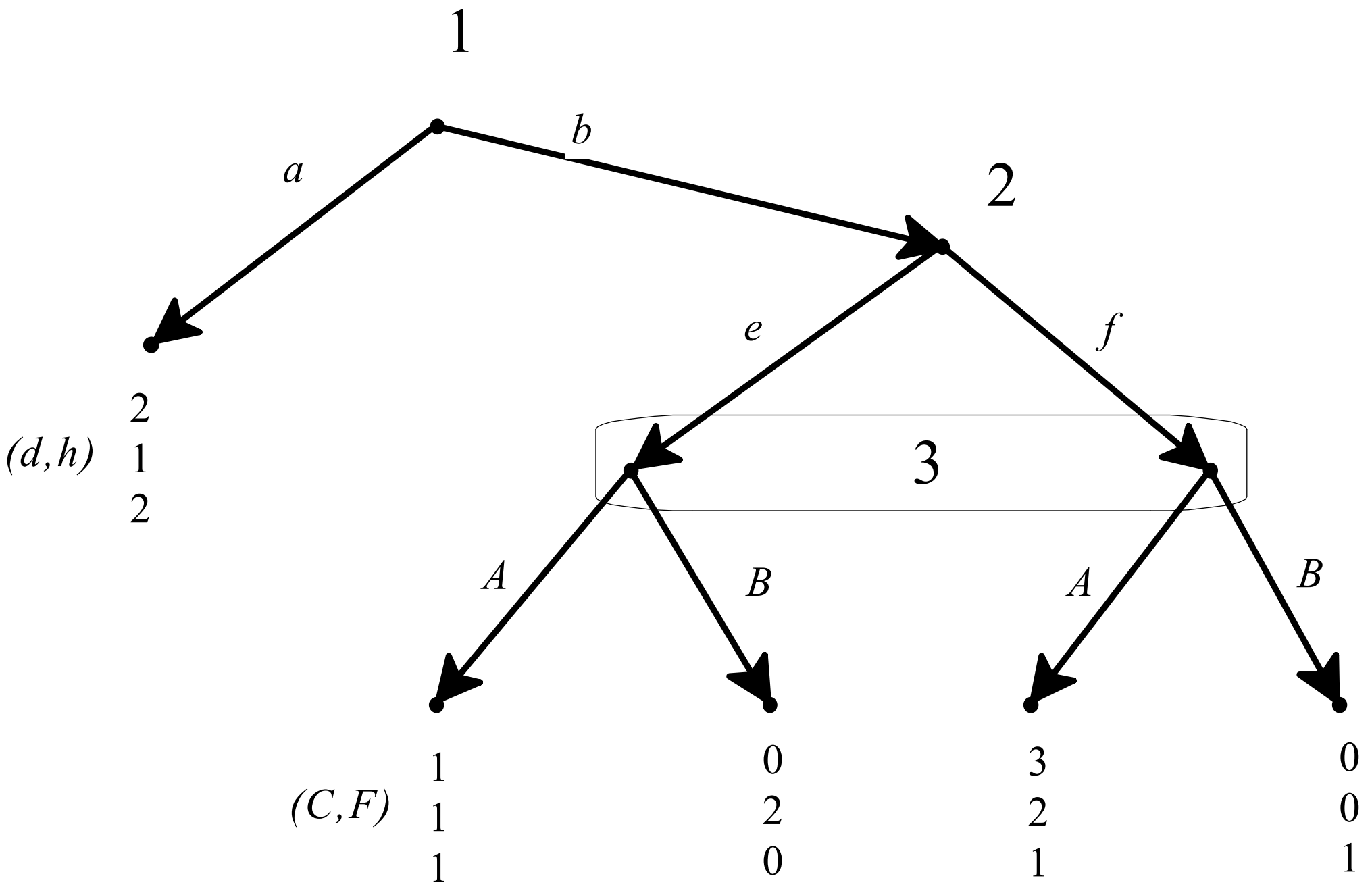

**Figure 3.13**
The reduced game after replacing a proper minimal subgame in the game of Figure 3.11.

The game of Figure 3.13 has a unique proper subgame, which has a unique Nash equilibrium, namely $(f, A)$. Replacing the subgame with the payoff vector associated with the history *bfA* we get the following smaller game.

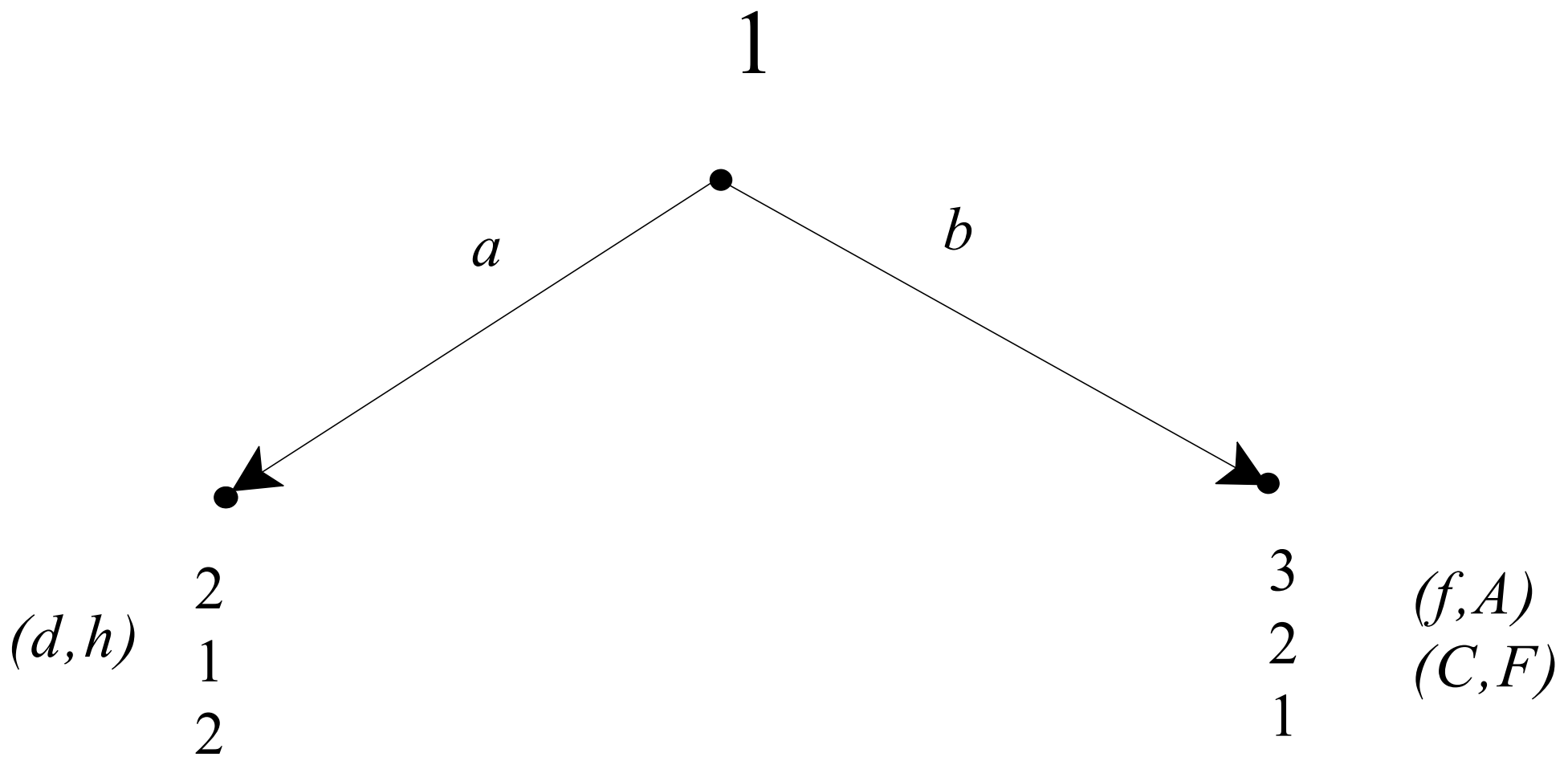

**Figure 3.14**
The reduced game after replacing the subgame in the game of Figure 3.13.

In the reduced game of Figure 3.14 the unique Nash equilibrium is *b*. Now patching together the choices selected during the application of the algorithm we get the following subgame-perfect equilibrium for the game of Figure 3.9: $((b,C),(d,f,F),(h,A))$.

As a second example, consider the game of Figure 3.15, which reproduces the game of Figure 3.8.

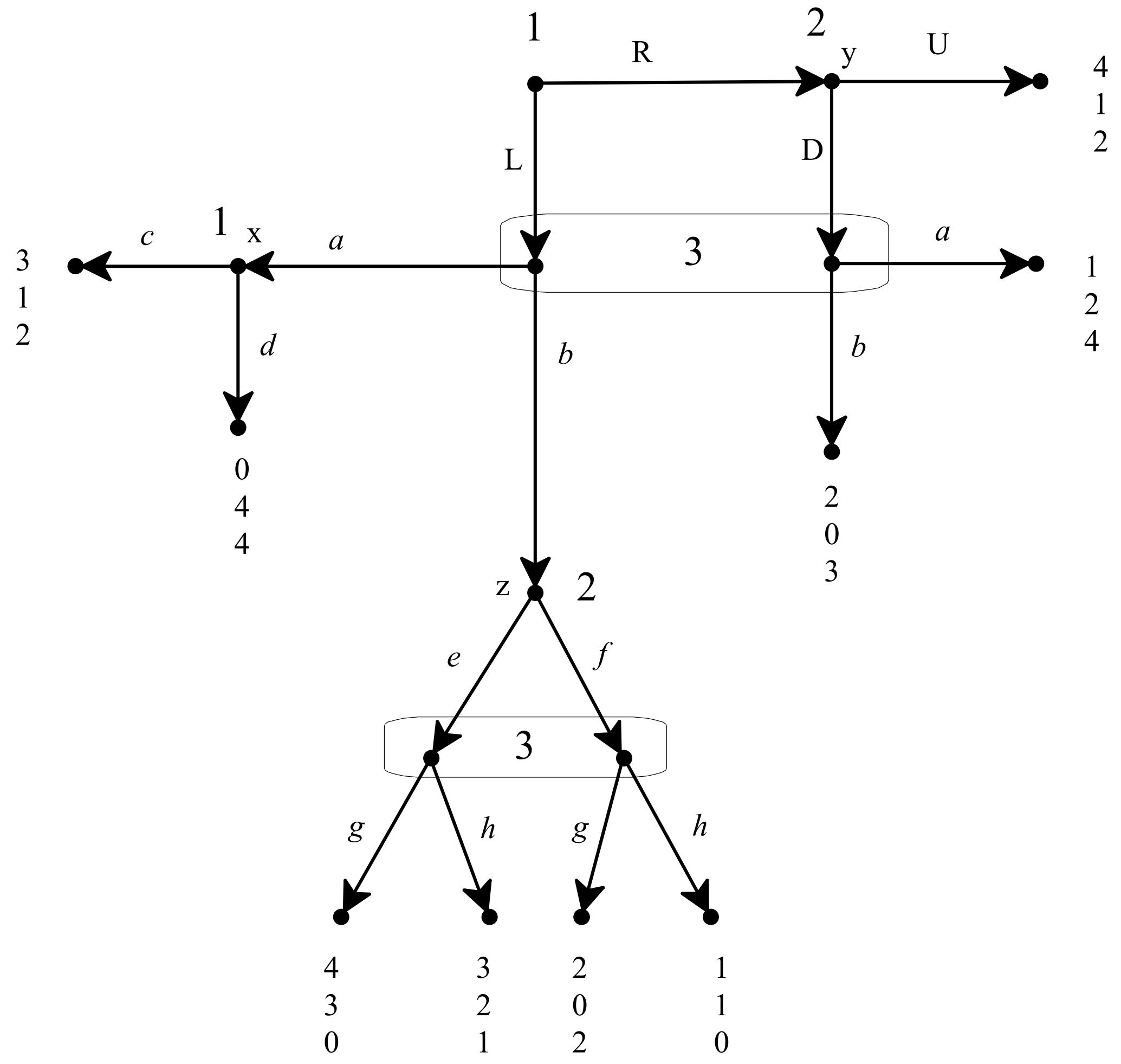

**Figure 3.15**
Reproduction of the game of Figure 3.8.

Begin with the subgame that starts at node *x* and replace it with the payoff vector (3,1,2). Next replace the subgame that starts at node *z* with the payoff vector (3,2,1) which corresponds to the Nash equilibrium (*e*,*h*) of that subgame, so that the game is reduced to the one shown in Figure 3.16, together with its strategic form.

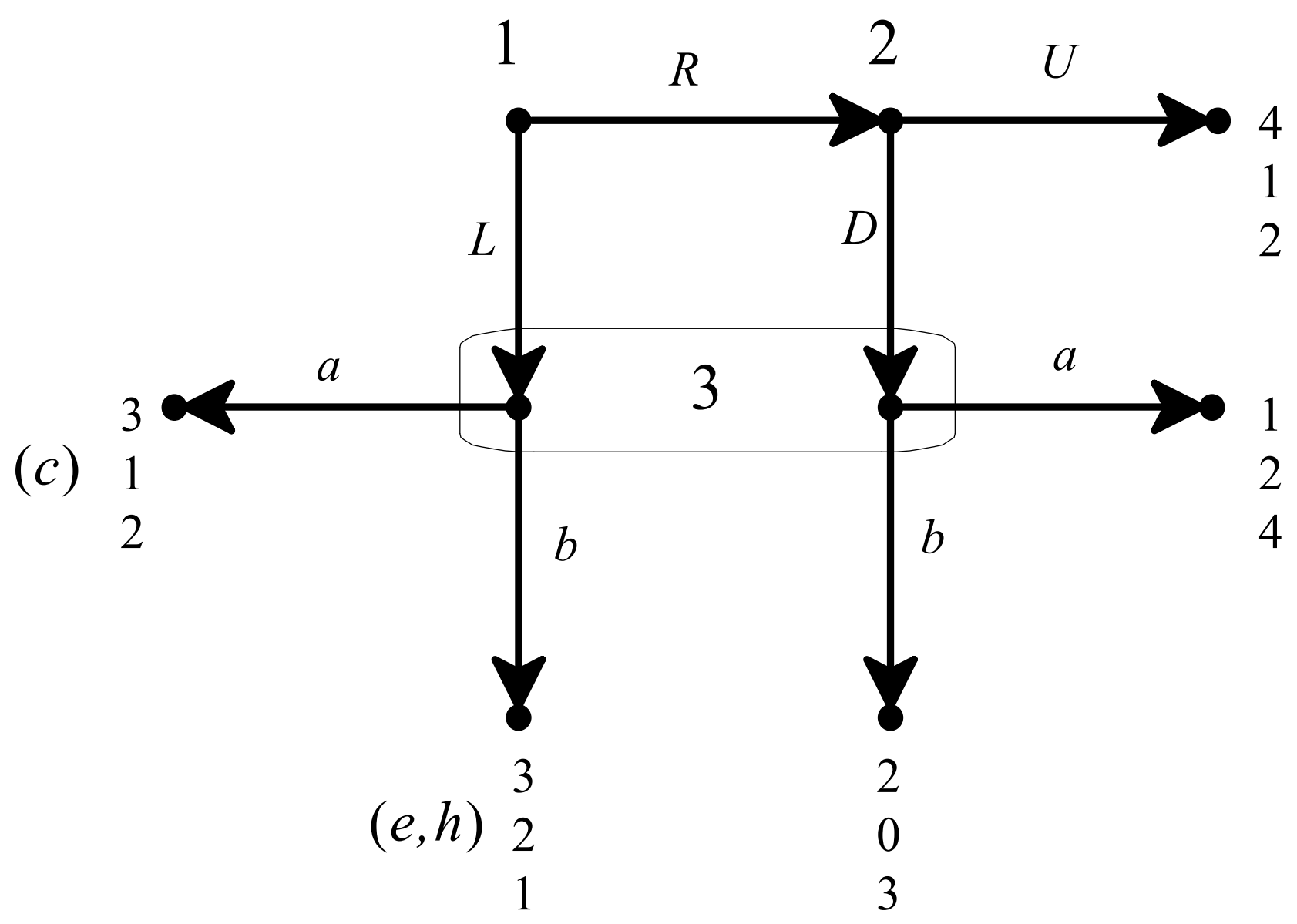

**Player 2**

| Player 1 | | U | | | D | | |
|---|---|---|---|---|---|---|---|
| | L | 3 | 1 | 2 | 3 | 1 | 2 |
| | R | 4 | 1 | 2 | 1 | 2 | 4 |

**Player 3:** ***a***

**Player 2**

| Player 1 | | U | | | D | | |
|---|---|---|---|---|---|---|---|
| | L | 3 | 2 | 1 | 3 | 2 | 1 |
| | R | 4 | 1 | 2 | 2 | 0 | 3 |

**Player 3: b**

**Figure 3.16**
The game of Figure 3.15 reduced after solving the proper subgames and the associated strategic form with the Nash equilibria highlighted.

The reduced game of Figure 3.16 has two Nash equilibria: (*L*,*D*,*a*) and (*R*,*U*,*b*). Thus the game of Figure 3.15 has two subgame-perfect equilibria:

$$\left( \underbrace{(L,c)}_{Player 1}, \underbrace{(D,e)}_{Player 2}, \underbrace{(a,h)}_{Player 3} \right) \text{ and } \left( \underbrace{(R,c)}_{Player 1}, \underbrace{(U,e)}_{Player 2}, \underbrace{(b,h)}_{Player 3} \right).$$

**Remark 3.1.** As shown in the last example, it is possible that – when applying the subgame-perfect equilibrium algorithm – one encounters a proper subgame or a reduced game that has several Nash equilibria. In this case one Nash equilibrium must be selected to continue the procedure and in the end one obtains one subgame-perfect equilibrium. One then has to repeat the procedure by selecting a different Nash equilibrium and thus obtain a different subgame-perfect equilibrium, and so on. This is similar to what happens with the backward-induction algorithm in perfect-information games.

**Remark 3.2.** It is also possible that – when applying the subgame-perfect equilibrium algorithm – one encounters a proper subgame or a reduced game that has no Nash equilibria.[12] In such a case the game under consideration does not have any subgame-perfect equilibria.

**Remark 3.3.** When applied to perfect-information games, the notion of subgame-perfect equilibrium coincides with the notion of backward-induction solution. Thus subgame-perfect equilibrium is a generalization of backward induction.

**Remark 3.4.** For extensive-form games that have no proper subgames (for example, the game of Figure 3.3) the set of Nash equilibria coincides with the set of subgame-perfect equilibria. In general, however, the notion of subgame-perfect equilibrium is a refinement of the notion of Nash equilibrium.

This is a good time to test your understanding of the concepts introduced in this section, by going through the exercises in Section 3.E.4 of Appendix 3.E at the end of this chapter.

[12] We will see in Part II that, when payoffs are cardinal and one allows for mixed strategies, then every finite game has at least one Nash equilibrium in mixed strategies.

# 3.5. Games with chance moves

So far we have only considered games where the outcomes do not involve any uncertainty. As a way of introducing the topic discussed in Part II, in this section we consider games where uncertain, probabilistic events are incorporated in the extensive form.

We begin with an example. There are three cards, one black and two red. They are shuffled well and put face down on the table. Adele picks the top card, looks at it without showing it to Ben and then tells Ben either "the top card is black" or "the top card is red": she could be telling the truth or she could be lying. Ben then has to guess the true color of the top card. If he guesses correctly he gets \$9 from Adele, otherwise he gives her \$9. How can we represent this situation? Whether the top card is black or red is not the outcome of a player's decision, but the outcome a random event, namely the shuffling of the cards. In order to capture this random event we introduce a fictitious player called *Nature* or *Chance*. We assign a probability distribution to Nature's "choices". In this case, since one card is black and the other two are red, the probability that the top card is black is $\frac{1}{3}$ and the probability that the top card is red is $\frac{2}{3}$. Note that we don't assign payoffs to Nature and thus the only real players are Adele and Ben. The situation can be represented as shown in Figure 3.17, where the numbers associated with the terminal nodes are dollar amounts.

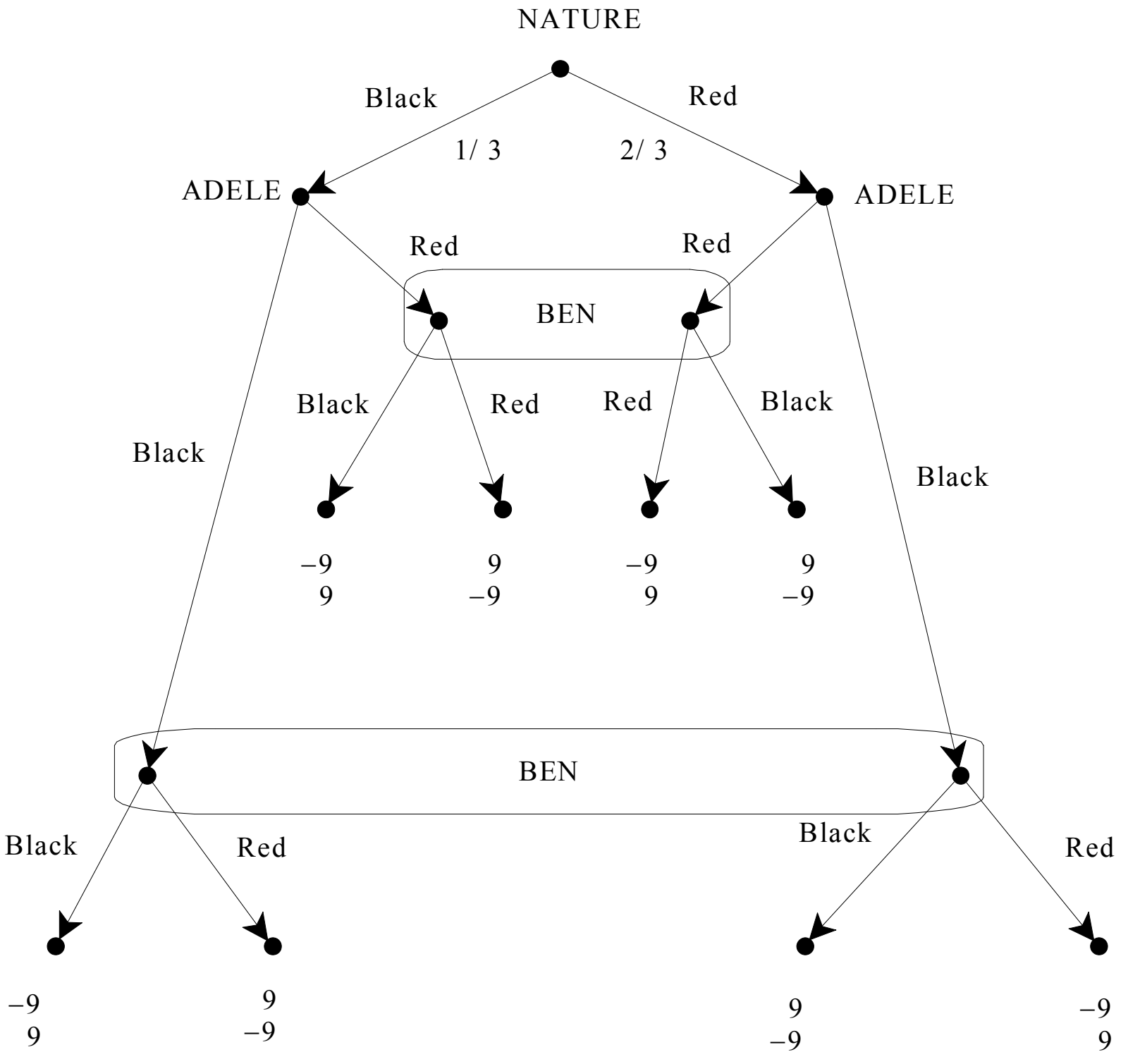

**Figure 3.17**
An extensive form with a chance move.

Clearly, the notion of strategy is not affected by the presence of chance moves. In the game of Figure 3.17 Adele has four strategies and so does Ben. However, we do encounter a difficulty when we try to write the associated strategic form. For example, consider the following strategy profile: $((B,R),(B,B))$ where Adele's strategy is to be truthful (say Black if she sees a black card and Red if she sees a red card) and Ben's strategy is to guess Black no matter what Adele says. What is the outcome in this case? It depends on what the true color of the top card is and thus the outcome is a probabilistic one:

$$\begin{pmatrix} \textit{outcome} & \text{Adele gives \$9 to Ben} & \text{Ben gives \$9 to Adele} \\ \textit{probability} & \frac{1}{3} & \frac{2}{3} \end{pmatrix}$$

We call such probabilistic outcomes *lotteries*. In order to convert the game-frame into a game we need to specify how the players rank probabilistic outcomes. Consider the case where Adele is selfish and greedy, in the sense that she only cares about her own wealth and she prefers more money to less. Then, from her point of view, the above probabilistic outcome reduces to the following monetary lottery $\begin{pmatrix} -\$9 & \$9 \\ \frac{1}{3} & \frac{2}{3} \end{pmatrix}$. If Ben is also selfish and greedy, then he views the same outcome as the lottery $\begin{pmatrix} \$9 & -\$9 \\ \frac{1}{3} & \frac{2}{3} \end{pmatrix}$.

How do we convert a lottery into a payoff or utility? The general answer to this question will be provided in Chapter 4. Here we consider one possibility, which is particularly simple.

**Definition 3.7.** Given a lottery whose outcomes are sums of money

$$\begin{pmatrix} \$x_1 & \$x_2 & \dots & \$x_n \\ p_1 & p_2 & \dots & p_n \end{pmatrix} \quad (p_i \geq 0, \text{ for all } i = 1,2,\dots,n \text{ and } \; p_1 + p_2 + \dots + p_n = 1)$$

the *expected value* of the lottery is the sum of money

$$\$\left(x_1 p_1 + x_2 p_2 + \dots + x_n p_n\right).$$

For example, the expected value of the lottery $\begin{pmatrix} \$5 & \$15 & \$25 \\ \frac{1}{5} & \frac{2}{5} & \frac{2}{5} \end{pmatrix}$ is $\$\left[5\left(\tfrac{1}{5}\right)+15\left(\tfrac{2}{5}\right)+25\left(\tfrac{2}{5}\right)\right] = \$\left(1+6+10\right) = \$17$ and the expected value of the lottery $\begin{pmatrix} -\$9 & \$9 \\ \frac{1}{3} & \frac{2}{3} \end{pmatrix}$ is \$3. We call lotteries whose outcomes are sums of money, *money lotteries*.

**Definition 3.8.** A player is defined to be *risk neutral* if she considers a money lottery to be just as good as its expected value. Hence a risk neutral person ranks money lotteries according to their expected value.

For example, consider the following money lotteries: $L_1 = \begin{pmatrix} \$5 & \$15 & \$25 \\ \frac{1}{5} & \frac{2}{5} & \frac{2}{5} \end{pmatrix}$, $L_2 = \begin{pmatrix} \$16 \\ 1 \end{pmatrix}$ and $L_3 = \begin{pmatrix} \$0 & \$32 & \$48 \\ \frac{5}{8} & \frac{1}{8} & \frac{1}{4} \end{pmatrix}$. The expected value of $L_1$ is \$17 and the expected value of both $L_2$ and $L_3$ is \$16. Thus a risk-neutral player would have the following ranking: $L_1 \succ L_2 \sim L_3$, that is, she would prefer $L_1$ to $L_2$ and be indifferent between $L_2$ and $L_3$.

For a selfish and greedy player who is risk neutral we can take the expected value of a money lottery as the utility of that lottery. For example, if we make the assumption that, in the extensive form of Figure 3.17, Adele and Ben are selfish, greedy *and risk neutral* then we can associate a strategic-form game to it as shown in Table 3.18. Note that inside each cell we have two numbers: the first is the utility (= expected value) of the underlying money lottery as perceived by Adele and the second number is the utility (= expected value) of the underlying money lottery as perceived by Ben.

| Adele \ Ben | *BB* | *BR* | *RB* | *RR* |
|---|---|---|---|---|
| *BB* | 3, –3 | –3, 3 | 3, –3 | –3, 3 |
| *BR* | 3, –3 | 9, –9 | –9, 9 | –3, 3 |
| *RB* | 3, –3 | –9, 9 | 9, –9 | –3, 3 |
| *RR* | 3, –3 | 3, –3 | –3, 3 | –3, 3 |

**Figure 3.18**

The strategic form of the game of Figure 3.17 when the two players are selfish, greedy and risk neutral. The first element of Adele's strategy is what she says if she sees a black card and the second element what she says if she sees a red card. The first element of Ben's strategy is what he guesses if Adele says "Red", the second element what he guesses if Adele says "Black".

We conclude this section with one more example.

**Example 3.3.** There are three unmarked, opaque envelopes. One contains \$100, one contains \$200 and the third contains \$300. They are shuffled well and then one envelope is given to Player 1 and another is given to Player 2 (the third one remains on the table). Player 1 opens her envelope and checks its content without showing it to Player 2. Then she either says "pass", in which case each player gets to keep his/her envelope, or she asks Player 2 to trade his envelope for hers. Player 2 is not allowed to see the content of his envelope and has to say either Yes or No. If he says No, then the two players get to keep their envelopes. If Player 2 says Yes, then they trade envelopes and the game ends. Each player is selfish, greedy and risk neutral.

This situation is represented by the extensive-form game shown in Figure 3.19, where (100,200) means that Player 1 gets the envelope with \$100 and Player 2 gets the envelope with \$200, etc.; P stands for "pass" and T for "suggest a trade"; Y for "Yes" and N for "No".

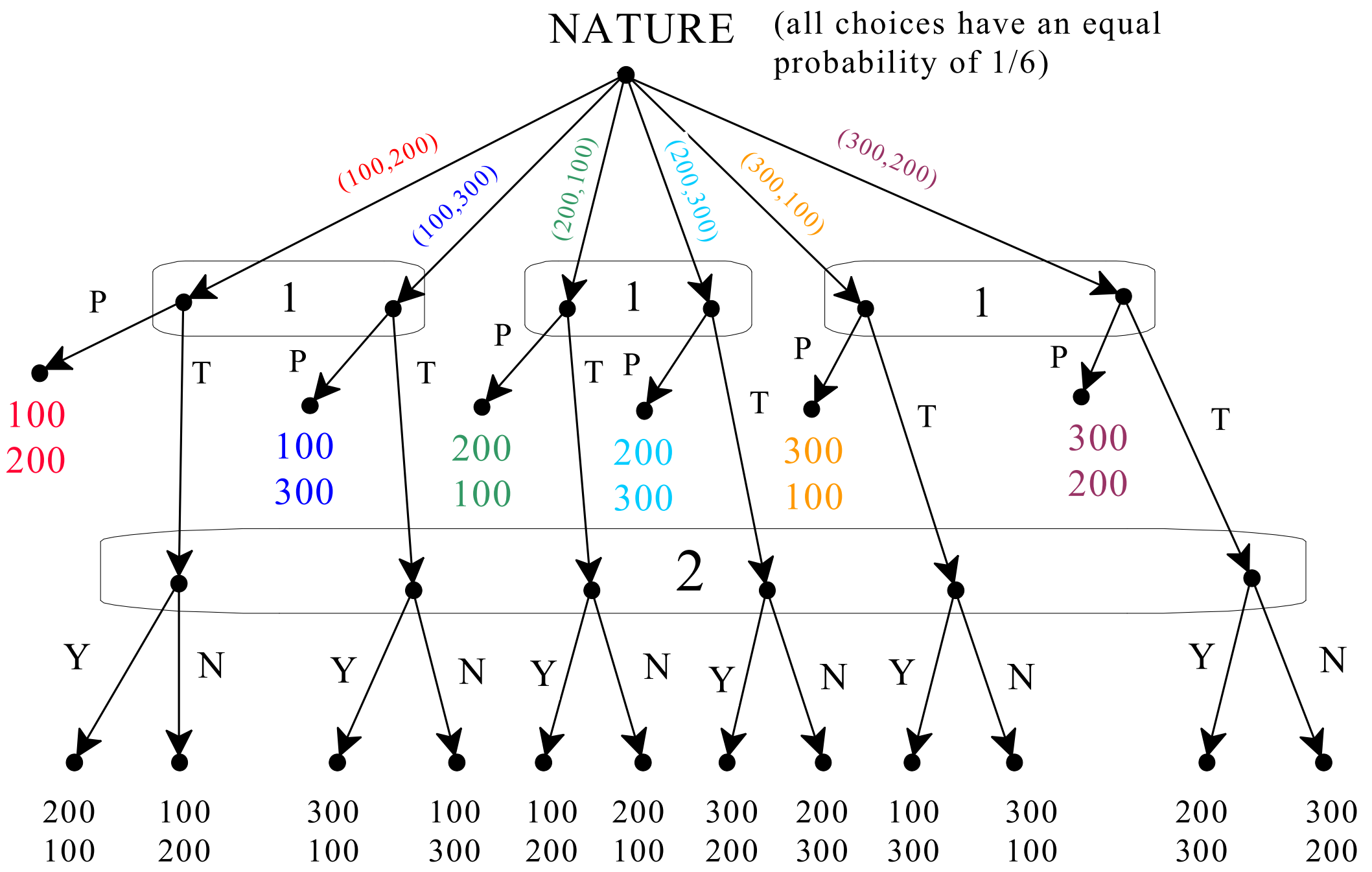

**Figure 3.19**
The extensive-form game of Example 3.3

In this game Player 1 has eight strategies. One possible strategy is: "if I get \$100 I will pass, if I get \$200 I will propose a trade, if I get \$300 I will pass": we will use the shorthand PTP for this strategy. Similarly for the other strategies. Player 2 has only two strategies: Yes and No. The strategic form associated with the game of Figure 3.19 is shown in Table 3.20, where the Nash equilibria are highlighted.

| | | Player 2 | |
|---|---|---|---|
| | | Y | N |
| P | PPP | 200 , 200 | 200 , 200 |
| l | PPT | 150 , 250 | 200 , 200 |
| a | PTP | 200 , 200 | 200 , 200 |
| y | PTT | 150 , 250 | 200 , 200 |
| e | TPP | 250 , 150 | 200 , 200 |
| r | TPT | 200 , 200 | 200 , 200 |
| | TTP | 250 , 150 | 200 , 200 |
| 1 | TTT | 200 , 200 | 200 , 200 |

**Table 3.20**
The strategic form of the game of Figure 3.19

How did we get those payoffs? Consider, for example, the first cell. Given the strategies PPP and Y, the outcomes are ($100,$200) with probability 1/6, ($100,$300) with probability 1/6, ($200,$100) with probability 1/6, ($200,$300) with probability 1/6, ($300,$100) with probability 1/6, ($300,$200) with probability 1/6. Being risk neutral, Player 1 views his corresponding money lottery as equivalent to getting its expected value $(100 + 100 + 200 + 200 + 300 + 300) (1/6) = $200. Similarly for Player 2 and for the other cells.

Since the game of Figure 3.19 has no proper subgames, all the Nash equilibria are also subgame-perfect equilibria. Are some of the Nash equilibria more plausible than others? For Player 1 all the strategies are weakly dominated, except for TPP and TTP. Elimination of the weakly dominated strategies leads to a game where Y is strictly dominated for Player 2. Thus one could argue that (TPP, N) and (TTP, N) are the most plausible equilibria; in both of them Player 2 refuses to trade.

This is a good time to test your understanding of the concepts introduced in this section, by going through the exercises in Section 3.E.5 of Appendix 3.E at the end of this chapter.

# Appendix 3.E: Exercises

## 3.E.1. Exercises for Section 3.1: Imperfect information

The answers to the following exercises are in Appendix S at the end of this chapter.

**Exercise 3.1**. Amy and Bill simultaneously write a bid on a piece of paper. The bid can only be either \$2 or \$3. A referee then looks at the bids, announces the amount of the lowest bid (without revealing who submitted it) and invites Amy to pass or double her initial bid. If Amy's final bid and Bill's bid are equal then Bill gets the object (and pays his bid), otherwise Amy gets the object (and pays her bid). Represent this situation by means of two alternative extensive frames. **Note**: (1) when there are simultaneous moves we have a choice as to which player we select as moving first: the important thing is that the second player does not know what the first player did; (2) when representing, by means of information sets, what a player is uncertain about, we typically assume that a player is smart enough to deduce relevant information, even if that information is not explicitly given to her.

**Exercise 3.2**. Consider the following situation. An incumbent monopolist decides at date 1 whether to build a small plant or a large plant. At date 2 a potential entrant observes the plant built by the incumbent and decides whether or not to enter. If she does not enter then her profits are 0 while the incumbent's profits are \$25 million with a small plant and \$20 million with a large plant. If the potential entrant decides to enter, she pays a cost of entry equal to \$$K$ million and at date 3 the two firms **simultaneously** decide whether to produce high output or low output. The profits of the firms are as follows (these figure do **not** include the cost of entry for the entrant; thus you need to subtract that cost for the entrant). In each cell, the first number is the profit of the entrant in millions of dollars and the second is the profit of the incumbent.

| | | incumbent | | | | incumbent | |
|---|---|---|---|---|---|---|---|
| | | low output | high output | | | low output | high output |
| Entrant | low output | 10 , 10 | 7 , 7 | Entrant | low output | 10 , 7 | 5 , 9 |
| | high output | 7 , 6 | 4 , 3 | | high output | 7 , 3 | 4 , 5 |
| | If incumbent has small plant | | | | If incumbent has large plant | | |

Draw an extensive form that represents this situation.

## 3.E.2. Exercises for Section 3.2: Strategies

The answers to the following exercises are in Appendix S at the end of this chapter.

**Exercise 3.3**. Write the strategic-form game-frame of the extensive form of Exercise 3.1 (that is, instead of writing payoffs in each cell you write the outcome). Verify that the strategic forms of the two possible versions of the extensive form are identical.

**Exercise 3.4**. Consider the extensive-form game of Exercise 3.2.

**(a)** Write down in words one of the strategies of the potential entrant.

**(b)** How many strategies does the potential entrant have?

**(c)** Write the strategic-form game associated with the extensive-form game.

**(d)** Find the Nash equilibria for the case where $K = 2$.

## 3.E.3. Exercises for Section 3.3: Subgames

The answers to the following exercises are in Appendix S at the end of this chapter.

**Exercise 3.5**. How many proper subgames does the extensive form of Figure 3.3 have?

**Exercise 3.6**. How many proper subgames does the extensive game of Figure 3.5 have?

**Exercise 3.7**. **(a)** How many proper subgames does the extensive game of Figure 2.3 (Chapter 2) have? **(b)** How many of those proper subgames are minimal?

## 3.E.4. Exercises for Section 3.4: Subgame-perfect equilibrium

**Exercise 3.8**. Find the Nash equilibria and the subgame-perfect equilibria of the following game.

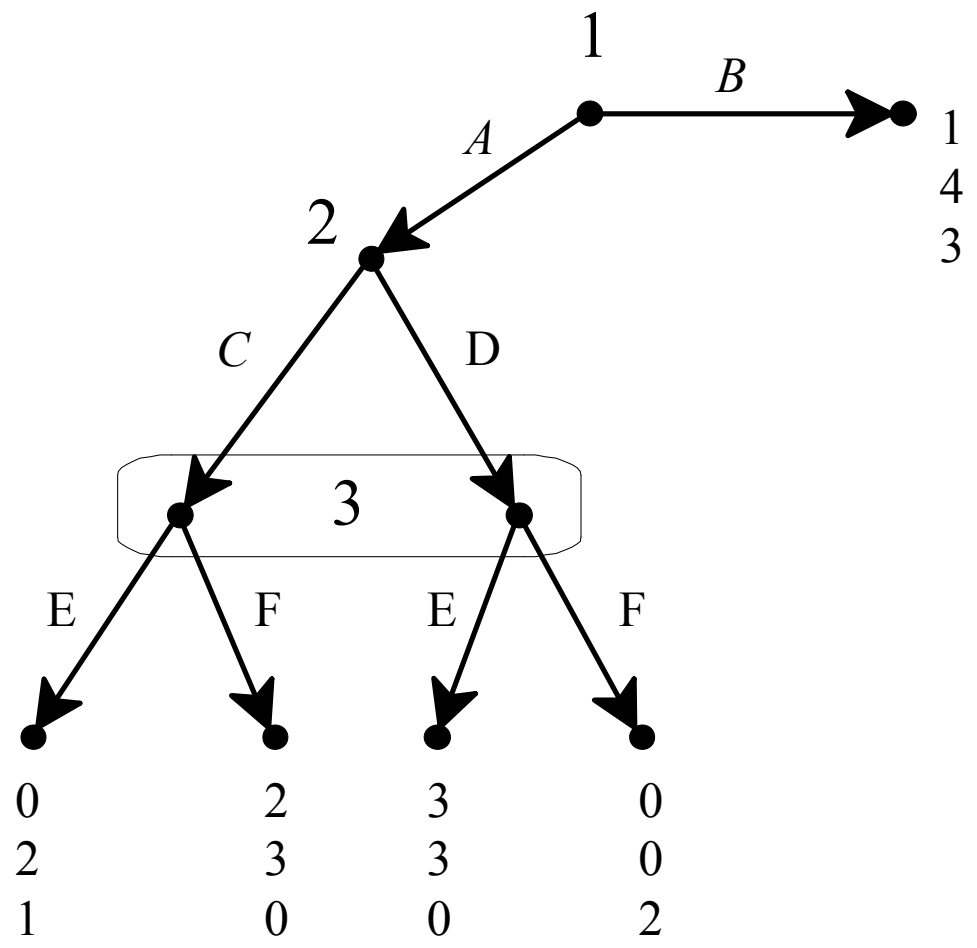


**Exercise 3.9**. Find the Nash equilibria and the subgame-perfect equilibria of the following game.

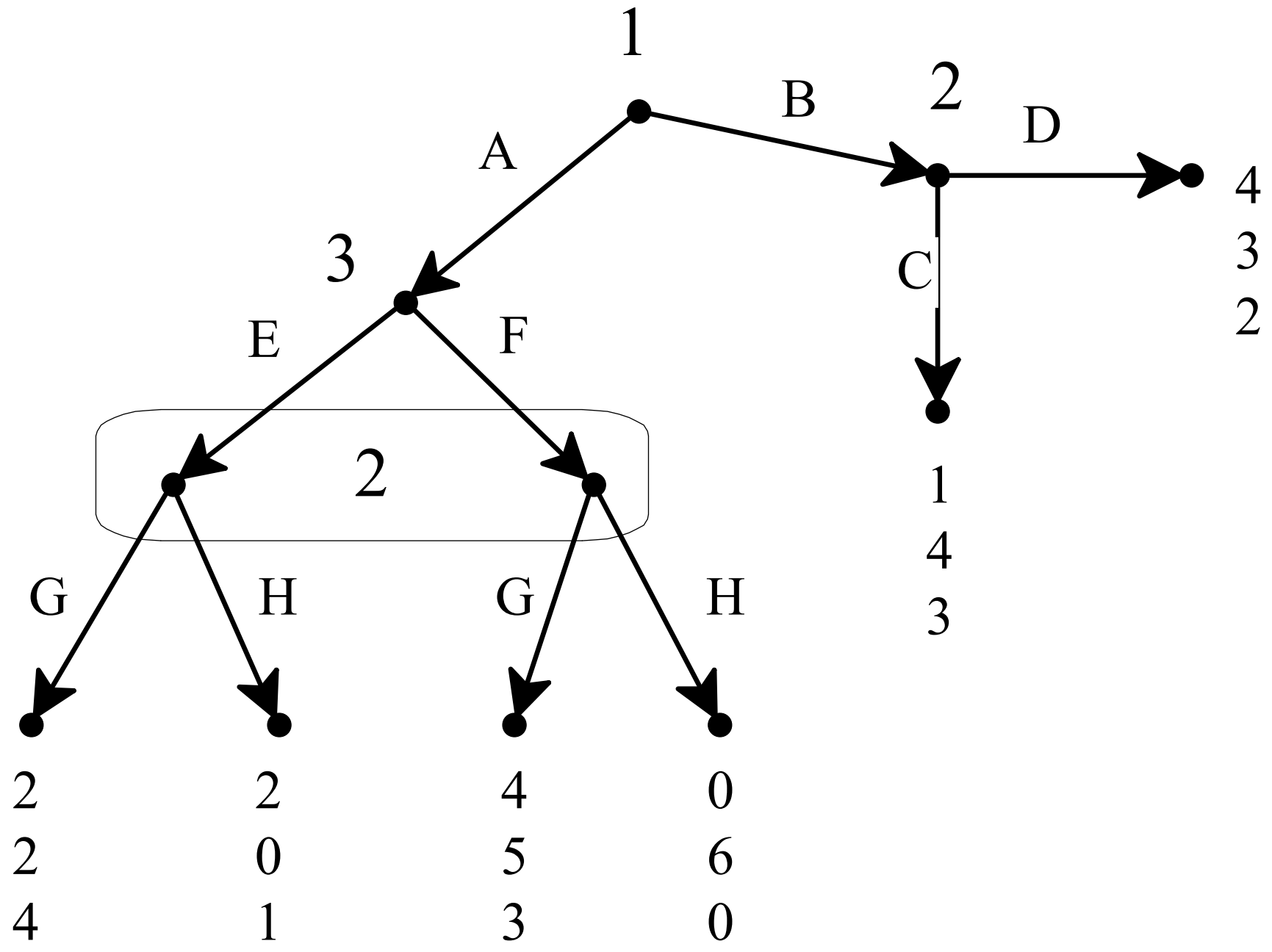

**Exercise 3.10**. Find the subgame-perfect equilibria of the following extensive-form game.

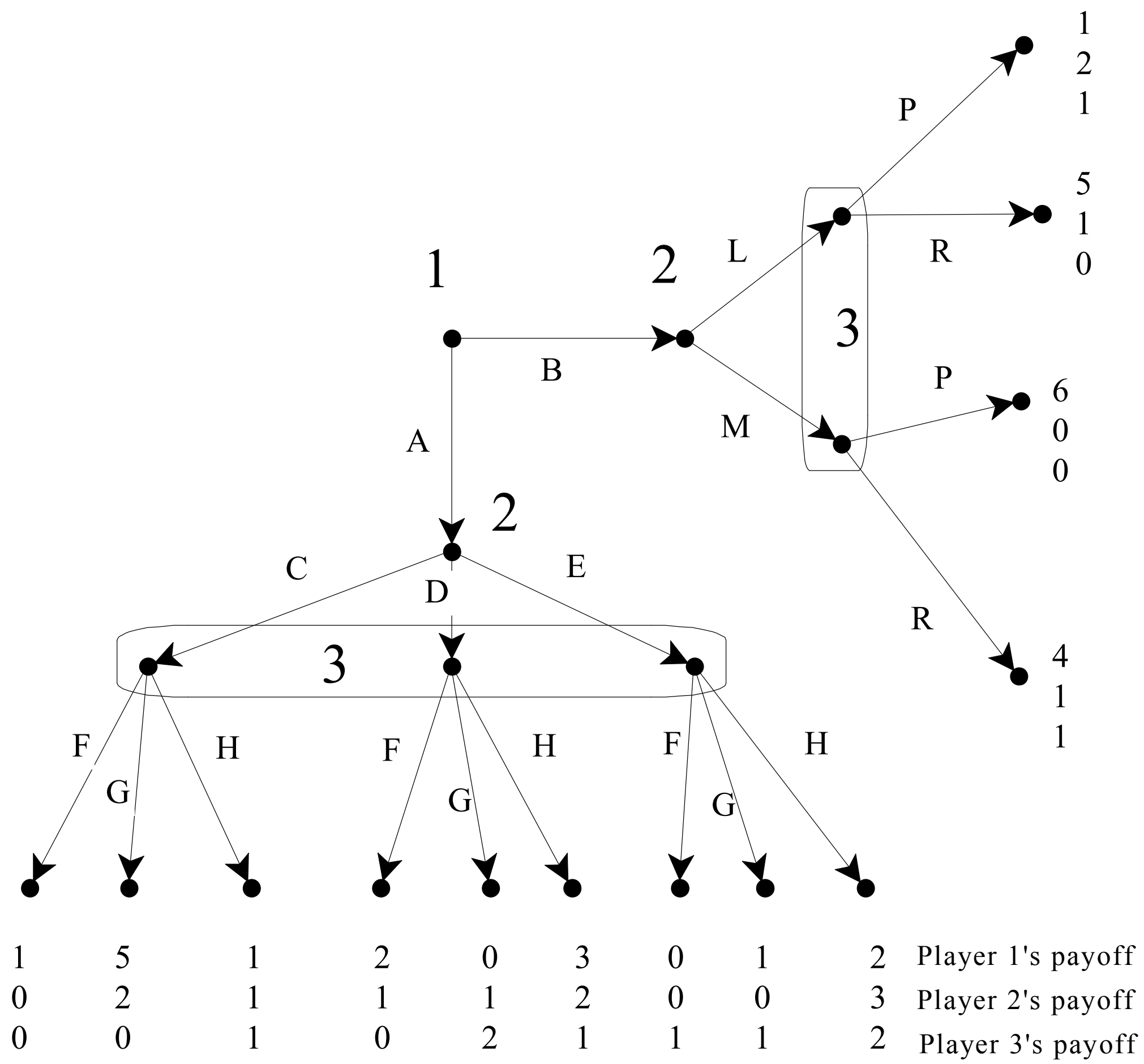


## 3.E.5. Exercises for Section 3.5: Games with chance moves

**Exercise 3.11**. Modify the game of Example 3.3 as follows: Player 2 is allowed to secretly check the content of his envelope before he decides whether or not to accept Player 1's proposal.

**(a)** Represent this situation as an extensive-form game.

**(b)** List all the strategies of Player 1 and all the strategies of Player 2.

**Exercise 3.12**. Three players, Avinash, Brian and John, play the following game. Two cards, one red and the other black, are shuffled well and put face down on the table. Brian picks the top card, looks at it without showing it to the other players (Avinash and John) and puts it back face down on the table. Then Brian whispers either "Black" or "Red" in Avinash's ear, making sure that John doesn't hear. Avinash then tells John either "Black" or "Red". Finally John announces either "Black" or "Red" and this exciting game ends. The payoffs are as follows: if John's final announcement matches the true color of the card Brian looked at, then Brian and Avinash give $2 each to John. In every other case John gives $2 each to Brian and Avinash.

**(a)** Represent this situation as an extensive-form game.
**(b)** Write the corresponding strategic form (or normal form) assuming that the players are selfish, greedy and risk neutral.

**Exercise 3.13**. Consider the following highly simplified version of Poker. There are three cards, marked A, B and C. A beats B and C, B beats C. There are two players, Yvonne and Zoe. Each player contributes $1 to the pot before the game starts. The cards are then shuffled and the top card is given, face down, to Yvonne and the second card (face down) to Zoe. Each player looks at, and only at, her own card: **she does not see the card of the other player nor the remaining card.** Yvonne, the first player, may *pass*, or *bet* $1. If she passes, the game ends, the cards are turned and the pot goes to the high-card holder (recall that A beats B and C, B beats C). If Yvonne bets, then Zoe can *fold*, in which case the game ends and the pot goes to Yvonne, or *see* by betting $1. If Yvonne sees the game ends, the cards are turned and the pot goes to the high-card holder. Both players are selfish, greedy and risk neutral.

**(a)** Draw the extensive-form game.
**(b)** How many strategies does Yvonne have?
**(c)** How many strategies does Zoe have?
**(d)** Consider the following strategies. For Yvonne: If A pass, if B pass, if C bet. For Zoe: if Yvonne bets and I get an A I will fold, if Yvonne bets and I get a B I will fold, if Yvonne bets and I get a C I will fold. Calculate the corresponding payoffs.
**(e)** Redo the same with the following strategies. For Yvonne: If A pass, if B pass, if C bet. For Zoe: see always (that is, no matter what card she gets).
**(f)** Now that you have understood how to calculate the payoffs, represent the entire game as a normal form game, assigning the rows to Yvonne and the columns to Zoe. [This might take you the entire night, but at your age sleep is not that important and also it will keep you out of trouble.]
**(g)** What strategies of Yvonne are weakly dominated? What strategies of Zoe are weakly dominated?
**(h)** What do you get when you apply the procedure of iterative elimination of weakly dominated strategies?

◊◊◊◊◊◊◊◊◊◊◊◊

**Exercise 3.14: Challenging Question.** In an attempt to reduce the deficit, the government of Italy has decided to sell a 14th century palace near Rome. The palace is in disrepair and is not generating any revenue for the government. From now on we'll call the government Player G. A Chinese millionaire has offered to purchase the palace for $*p*. Alternatively, Player G can organize an auction among *n* interested parties ($n \geq 2$). The participants to the auction (we'll call them players) have been randomly assigned labels 1,2,...,n. Player *i* is willing to pay up to $$p_i$ for the palace, where $p_i$ is a positive integer. For the auction assume the following: (1) it is a simultaneous, sealed-bid second-price auction, (2) bids must be non-negative integers, (3) each player only cares about his own wealth, (4) the tie-breaking rule for the auction is that the palace is given to that player who has the lowest index (e.g. if the highest bid was submitted by Players 3, 5 and 12 then the palace is given to Player 3). All of the above is commonly known among everybody involved, as is the fact that for every $i, j \in \{1,...,n\}$ with $i \neq j$, $p_i \neq p_j$. We shall consider four different scenarios. In all scenarios you can assume that the $p_i$'s are common knowledge.

**Scenario 1.** Player *G* first decides whether to sell the palace to the Chinese millionaire or make a public and irrevocable decision to auction it.

**(a)** Draw the extensive form of this game for the case where $n = 2$ and the only possible bids are $1 and $2. [List payoffs in the following order: first *G* then 1 then 2.]

**(b)** For the general case where $n \geq 2$ and every positive integer is a possible bid, find a pure-strategy subgame-perfect equilibrium of this game. What are the players' payoffs at the equilibrium?

**Scenario 2.** Here we assume that $n = 2$, and $p_1 > p_2 + 1 > 2$. First Player *G* decides whether to sell the palace to the Chinese or make a public and irrevocable decision to auction it. In the latter case he first asks Player 2 to publicly announce whether or not he is going to participate in the auction. If Player 2 says Yes then he has to pay $1 to Player *G* as a participation fee, which is non-refundable. If he says No then she is out of the game. After player 2 has made his announcement (and paid his fee if he decided to participate) Player 1 is asked to make the same decision (participate and pay a non-refundable fee of $1 to Player *G* or stay out); Player 1 knows player 2's decision when he makes his own decision. After both players have made their decisions, player *G* proceeds as follows: (1) if both 1 and 2 said Yes then he makes them play a simultaneous second-price auction, (2) if only one player said Yes then he is asked to put an amount $*x* of his choice in an envelope (where *x* is a positive integer) and give it to Player *G* in exchange for the palace, (3) if both 1 and 2 said No then *G* is no longer bound by his commitment to auction the palace and he sells it to the Chinese.

**(c)** Draw the extensive form of this game for the case where the only possible bids are $1 and $2 and also $x \in \{1,2\}$ [List payoffs in the following order: first *G* then 1 then 2.]

**(d)** For the general case where all possible bids are allowed (subject to being positive integers) and *x* can be any positive integer, find a pure-strategy subgame-perfect equilibrium of this game. What are the players' payoffs at the equilibrium?

**Scenario 3**. Same as Scenario 2; the only difference is that if both Players 1 and 2 decide to participate in the auction then Player *G* gives to the loser the fraction *a* (with $0 < a < 1$) of the amount paid by the winner *in the auction* (note that player *G* still keeps 100% of the participation fees). This is publicly announced at the beginning and is an irrevocable commitment.

**(e)** For the general case where all possible bids are allowed (subject to being positive integers) find a subgame-perfect equilibrium of this game. What are the players' payoff at the equilibrium?

**Scenario 4**. Player *G* tells the Chinese millionaire the following:

> "First you (= the Chinese) say Yes or No; if you say No I will sell you the palace for the price that you offered me, namely \$100 (that is, we now assume that $p = 100$); if you say Yes then we play the following perfect information game. I start by choosing a number from the set {1,2,3}, then you (= the Chinese) choose a number from this set, then I choose again, followed by you, etc. The first player who brings the cumulative sum of all the numbers chosen (up to and including the last one) to 40 **wins**. If you win I will sell you the palace for \$50, while if I win I will sell you the palace for \$200." Thus there is no auction in this scenario. Assume that the Chinese would actually be willing to pay up to \$300 for the palace.

**(f)** Find a pure-strategy subgame-perfect equilibrium of this game.

# Appendix 3.S: Solutions to exercises

**Exercise 3.1**. One possible extensive frame is as follows, where Amy moves first. Note that we have only one non-trivial information set for Amy, while each of the other three consists of a single node. The reason is as follows: if Amy initially bids \$3 and Bill bids \$2 then the referee announces "the lowest bid was \$2"; this announcement does not directly reveal to Amy that Bill's bid was \$2, but she can figure it out from her knowledge that her own bid was \$3; similarly, if the initial two bids are both \$3 then the referee announces "the lowest bid was \$3", in which case Amy is able to figure out that Bill's bid was also \$3. If we included those two nodes in the same information set for Amy, we would not show much faith in Amy's reasoning ability!

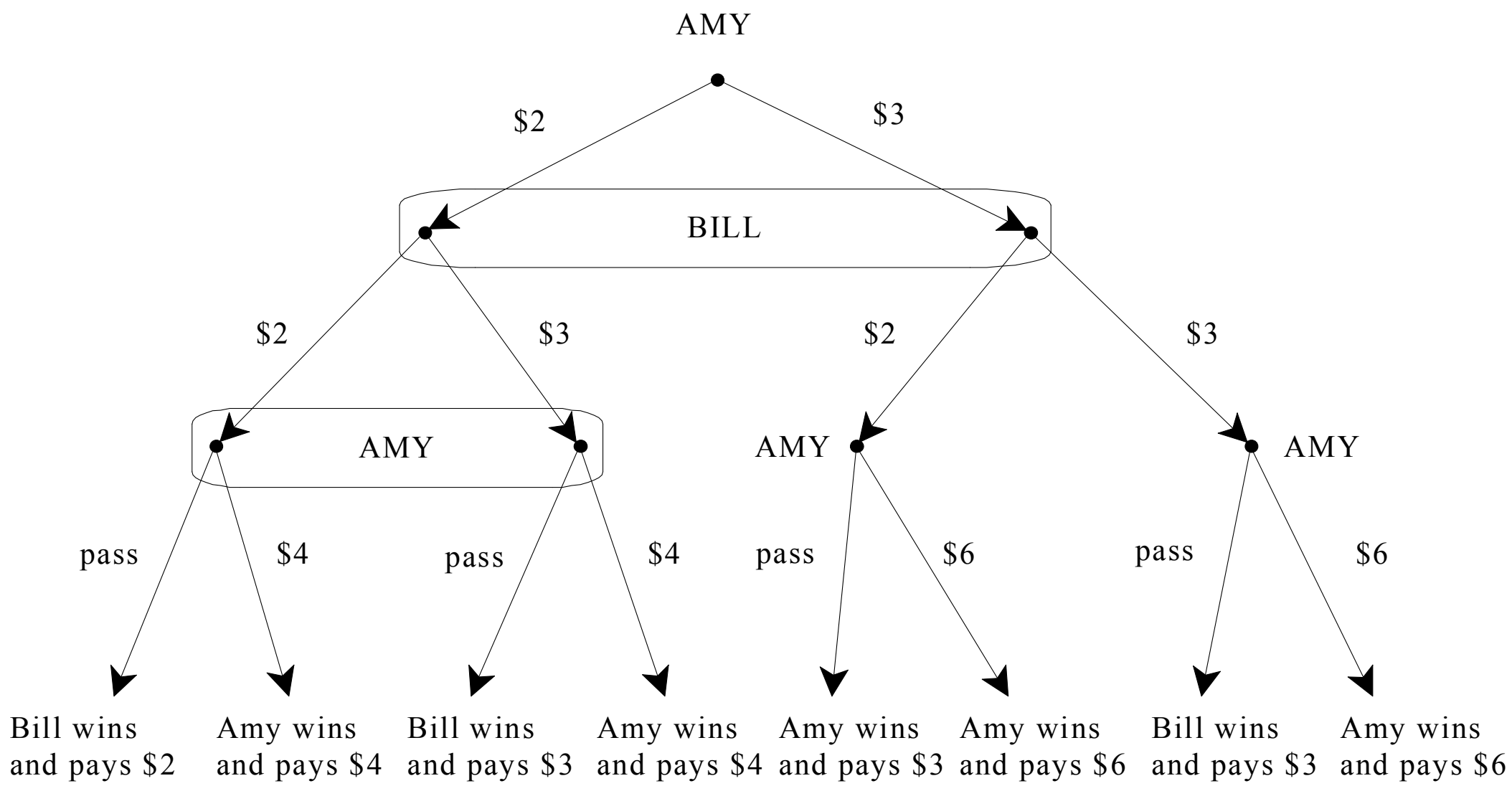

Another possible extensive frame is as follows, where Bill moves first.

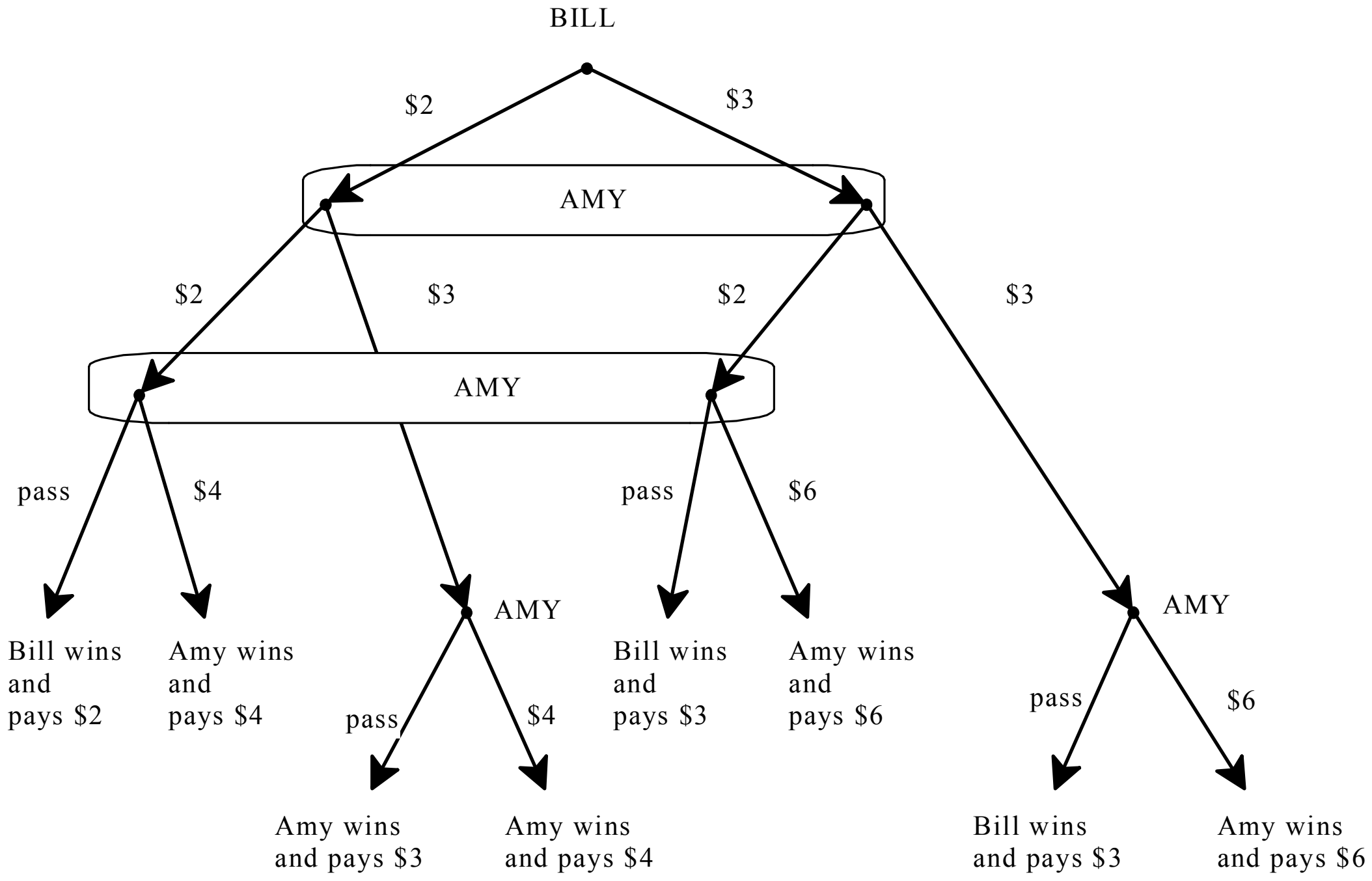

**Exercise 3.2.** The extensive form is as follows.

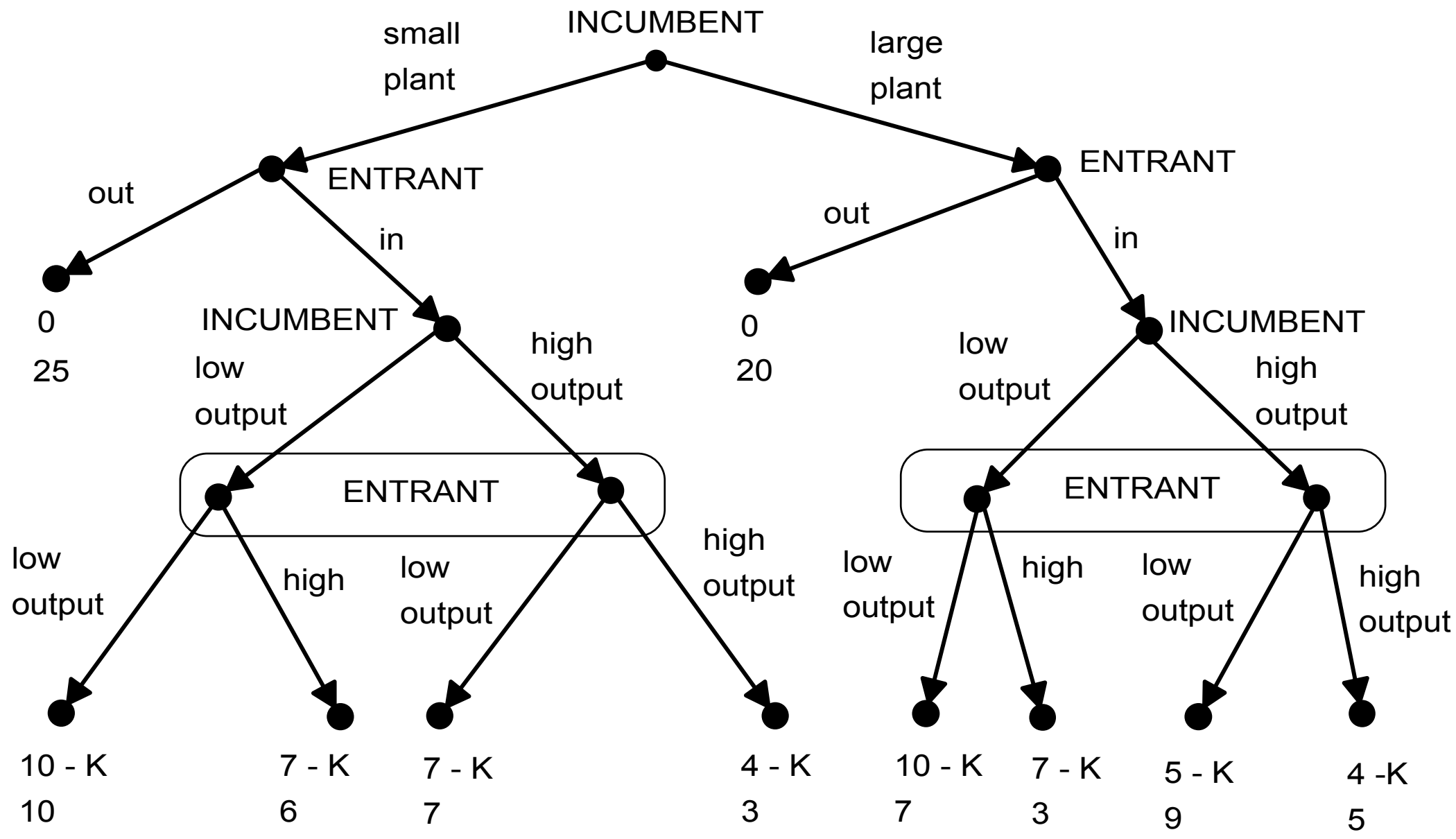

**Exercise 3.3**. Amy's strategy $(x,y,w,z)$ means: at the beginning I bid \$$x$, at the non-trivial information set on the left I choose $y$, at the singleton node in the middle I choose $w$ and at the singleton node on the right I choose $z$. The numbers are bid amounts and $P$ stands for "Pass".

| | | BILL | |
|---|---|---|---|
| | | bid \$2 | bid \$3 |
| | \$2, P, P, \$P | Bill wins pays 2 | Bill wins pays 3 |
| | 2 , P, P, 6 | Bill wins pays 2 | Bill wins pays 3 |
| | 2, P, 6, P | Bill wins pays 2 | Bill wins pays 3 |
| | 2, P, 6, 6 | Bill wins pays 2 | Bill wins pays 3 |
| A | 2, 4, P, P | Amy wins pays 4 | Amy wins pays 4 |
| M | 2, 4, P, 6 | Amy wins pays 4 | Amy wins pays 4 |
| Y | 2, 4, 6, P | Amy wins pays 4 | Amy wins pays 4 |
| | 2, 4, 6, 6 | Amy wins pays 4 | Amy wins pays 4 |
| | 3, P, P, P | Amy wins pays 3 | Bill wins pays 3 |
| | 3 , P, P, 6 | Amy wins pays 3 | Amy wins pays 6 |
| | 3, P, 6, P | Amy wins pays 6 | Bill wins pays 3 |
| | 3, P, 6, 6 | Amy wins pays 6 | Amy wins pays 6 |
| | 3, 4, P, P | Amy wins pays 3 | Bill wins pays 3 |
| | 3, 4, P, 6 | Amy wins pays 3 | Amy wins pays 6 |
| | 3, 4, 6, P | Amy wins pays 6 | Bill wins pays 3 |
| | 3, 4, 6, 6 | Amy wins pays 6 | Amy wins pays 6 |

**Exercise 3.4**. **(a)** The potential entrant has four information sets, hence a strategy has to say what she would do in each of the four situations. A possible strategy is: "if the incumbent chooses a small plant I stay out, if the incumbent chooses a large plant I enter, if small plant and I entered then I choose low output, if large plant and I entered then I choose high output".

**(b)** The potential entrant has $2^4 = 16$ strategies.

**(c)** The strategic form is as follows.

| | | INCUMBENT | | | | | | | |
|---|---|---|---|---|---|---|---|---|---|
| | | SLL | SLH | SHL | SHH | LLL | LLH | LHL | LHH |
| | OOLL | 0 , 25 | 0 , 25 | 0 , 25 | 0 , 25 | 0 , 20 | 0 , 20 | 0 , 20 | 0 , 20 |
| | OOLH | 0 , 25 | 0 , 25 | 0 , 25 | 0 , 25 | 0 , 20 | 0 , 20 | 0 , 20 | 0 , 20 |
| | OOHL | 0 , 25 | 0 , 25 | 0 , 25 | 0 , 25 | 0 , 20 | 0 , 20 | 0 , 20 | 0 , 20 |
| | OOLH | 0 , 25 | 0 , 25 | 0 , 25 | 0 , 25 | 0 , 20 | 0 , 20 | 0 , 20 | 0 , 20 |
| E | OILL | 0 , 25 | 0 , 25 | 0 , 25 | 0 , 25 | 10−K , 7 | 5−K, 9 | 10−K , 7 | 5−K, 9 |
| N | OILH | 0 , 25 | 0 , 25 | 0 , 25 | 0 , 25 | 7−K, 3 | 4−K, 5 | 7−K, 3 | 4−K, 5 |
| T | OIHL | 0 , 25 | 0 , 25 | 0 , 25 | 0 , 25 | 10−K , 7 | 5−K, 9 | 10−K , 7 | 5−K, 9 |
| R | OIHH | 0 , 25 | 0 , 25 | 0 , 25 | 0 , 25 | 7−K, 3 | 4−K, 5 | 7−K, 3 | 4−K, 5 |
| A | IOLL | 10−K, 10 | 10−K, 10 | 7−K, 7 | 7−K, 7 | 0 , 20 | 0 , 20 | 0 , 20 | 0 , 20 |
| N | IOLH | 10−K, 10 | 10−K, 10 | 7−K, 7 | 7−K, 7 | 0 , 20 | 0 , 20 | 0 , 20 | 0 , 20 |
| T | IOHL | 7−K, 6 | 7−K, 6 | 4−K,3 | 4−K, 3 | 0 , 20 | 0 , 20 | 0 , 20 | 0 , 20 |
| | IOHH | 7−K, 6 | 7−K, 6 | 4−K,3 | 4−K, 3 | 0 , 20 | 0 , 20 | 0 , 20 | 0 , 20 |
| | IILL | 10−K, 10 | 10−K, 10 | 7−K, 7 | 7−K, 7 | 10−K , 7 | 5−K, 9 | 10−K , 7 | 5−K, 9 |
| | IILH | 10−K, 10 | 10−K, 10 | 7−K, 7 | 7−K, 7 | 7−K, 3 | 4−K, 5 | 7−K, 3 | 4−K, 5 |
| | IIHL | 7−K, 6 | 7−K, 6 | 4−K,3 | 4−K, 3 | 10−K , 7 | 5−K, 9 | 10−K , 7 | 5−K, 9 |
| | IIHH | 7−K, 6 | 7−K, 6 | 4−K,3 | 4−K, 3 | 7−K, 3 | 4−K, 5 | 7−K, 3 | 4−K, 5 |

**(d)** For the case where $K = 2$ the Nash equilibria are highlighted in the above table.

**Exercise 3.5**. None.

**Exercise 3.6**. None.

**Exercise 3.7**. The game under consideration is the following:

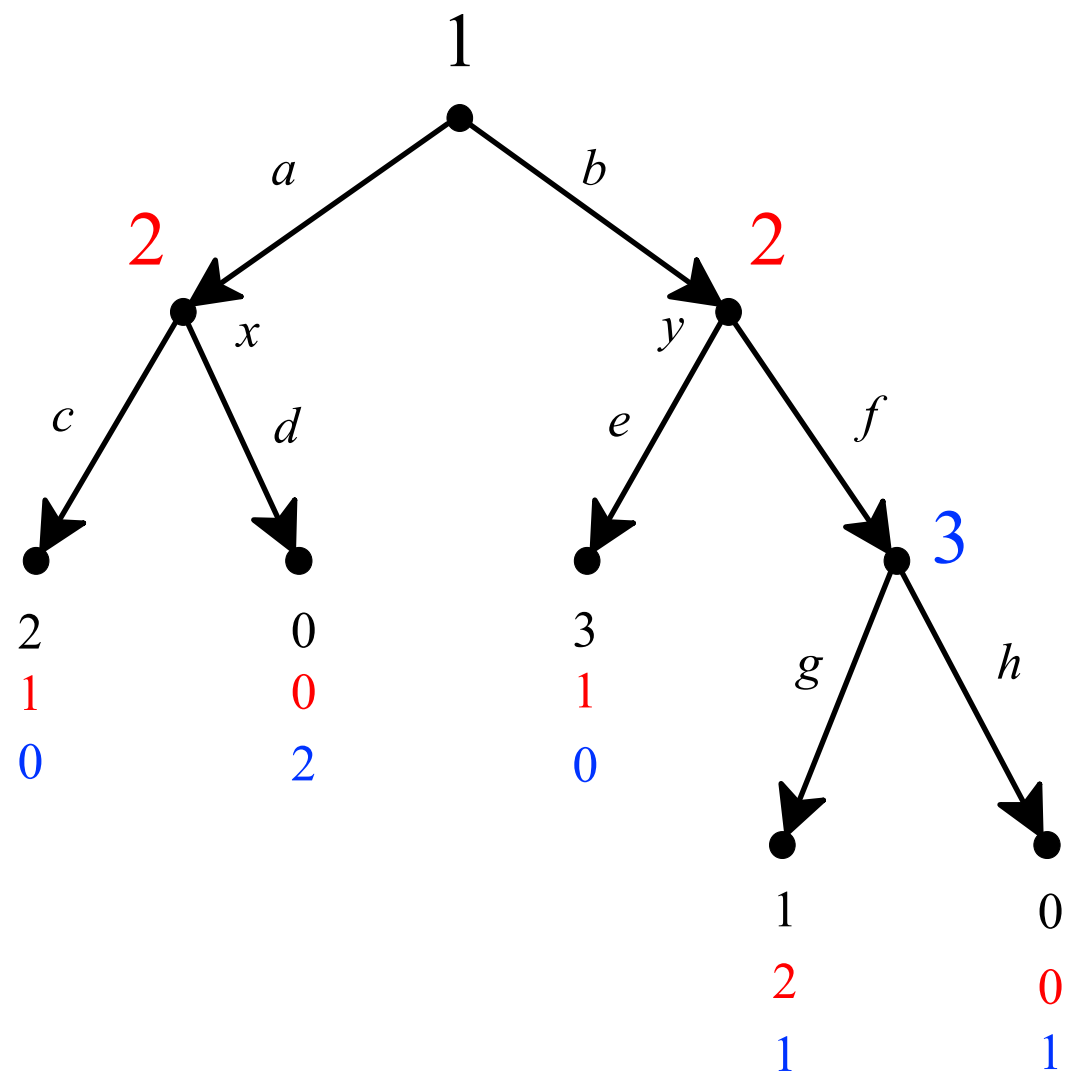


**(a)** There are three proper subgames: one starting from node *x*, one starting from node *y* and one starting at the node of Player 3.

**(b)** Two: the one starting from node *x* and the one starting at the decision node of Player 3. In a perfect-information game a minimal proper subgame is one that starts at a decision node followed only by terminal nodes.

**Exercise 3.8**. The Nash equilibria are (*B*,*C*,*E*) and (*B*,*D*,*F*). There are no subgame-perfect equilibria since the proper subgame has no Nash equilibria.

**Exercise 3.9**. The strategic form is as follows:

| Player 1 | Player 2: GC | GD | HC | HD |
|---|---|---|---|---|
| A | 2, 2, 4 | 2, 2, 4 | 2, 0, 1 | 2, 0, 1 |
| B | 1, 4, 3 | 4, 3, 2 | 1, 4, 3 | 4, 3, 2 |

Player 3 chooses E

| Player 1 | Player 2: GC | GD | HC | HD |
|---|---|---|---|---|
| A | 4, 5, 3 | 4, 5, 3 | 0, 6, 0 | 0, 6, 0 |
| B | 1, 4, 3 | 4, 3, 2 | 1, 4, 3 | 4, 3, 2 |

Player 3 chooses F

The Nash equilibria are: (*A*,(*G*,*C*),*E*) and (*B*,(*H*,*C*),*F*).

The extensive-form game has two proper subgames. The one on the left has a unique Nash equilibrium, (*G*,*E*), and the one on the right has a unique Nash equilibrium, *C*. Hence the game reduces to:

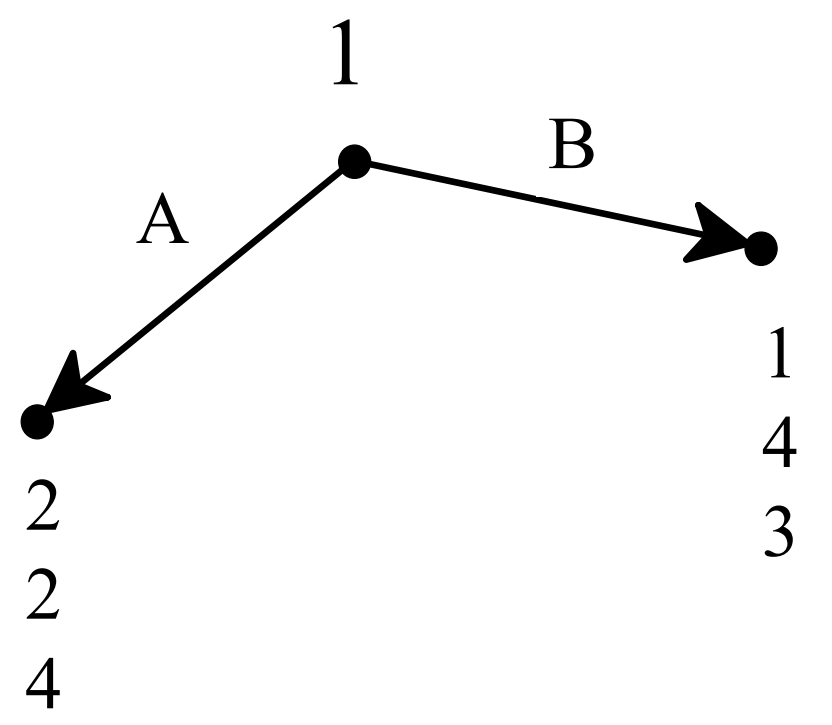


where $A$ is the unique optimal choice. Hence there is only one subgame-perfect equilibrium, namely $(A,(G,C),E)$.

**Exercise 3.10.** Consider first the subgame that starts at Player 2's decision node coming after choice $A$ of Player 1 (only the payoff of Players 2 and 3 are shown in the following table):

| | | Player 3: *F* | *G* | *H* |
|---|---|---|---|---|
| Player 2 | *C* | 0 , 0 | 2 , 0 | 1 , 1 |
| | *D* | 1 , 0 | 1 , 2 | 2 , 1 |
| | *E* | 0 , 1 | 0 , 1 | 3 , 2 |

The unique Nash equilibrium is (E, H).

Now consider the subgame that starts at Player 2's decision node coming after choice $B$ of Player 1

| | | Player 3: *P* | *R* |
|---|---|---|---|
| Player 2 | *L* | 2 , 1 | 1 , 0 |
| | *M* | 0 , 0 | 1 , 1 |

There are two Nash equilibria: $(L, P)$ and $(M, R)$.

Thus there are two subgame-perfect equilibria of the entire game:

(1) $(A,(E,L),((H,P))$. Player 1's strategy: $A$, Player 2's strategy: $E$ if $A$ and $L$ if $B$, Player 3's strategy: $H$ if $A$ and $P$ if $B$:

(2) $(B,(E,M),((H,R))$. Player 1's strategy: $B$, Player 2's strategy: $E$ if $A$ and $M$ if $B$, Player 3's strategy: $H$ if $A$ and $R$ if $B$.

**Exercise 3.11. (a)** The game is as follows.

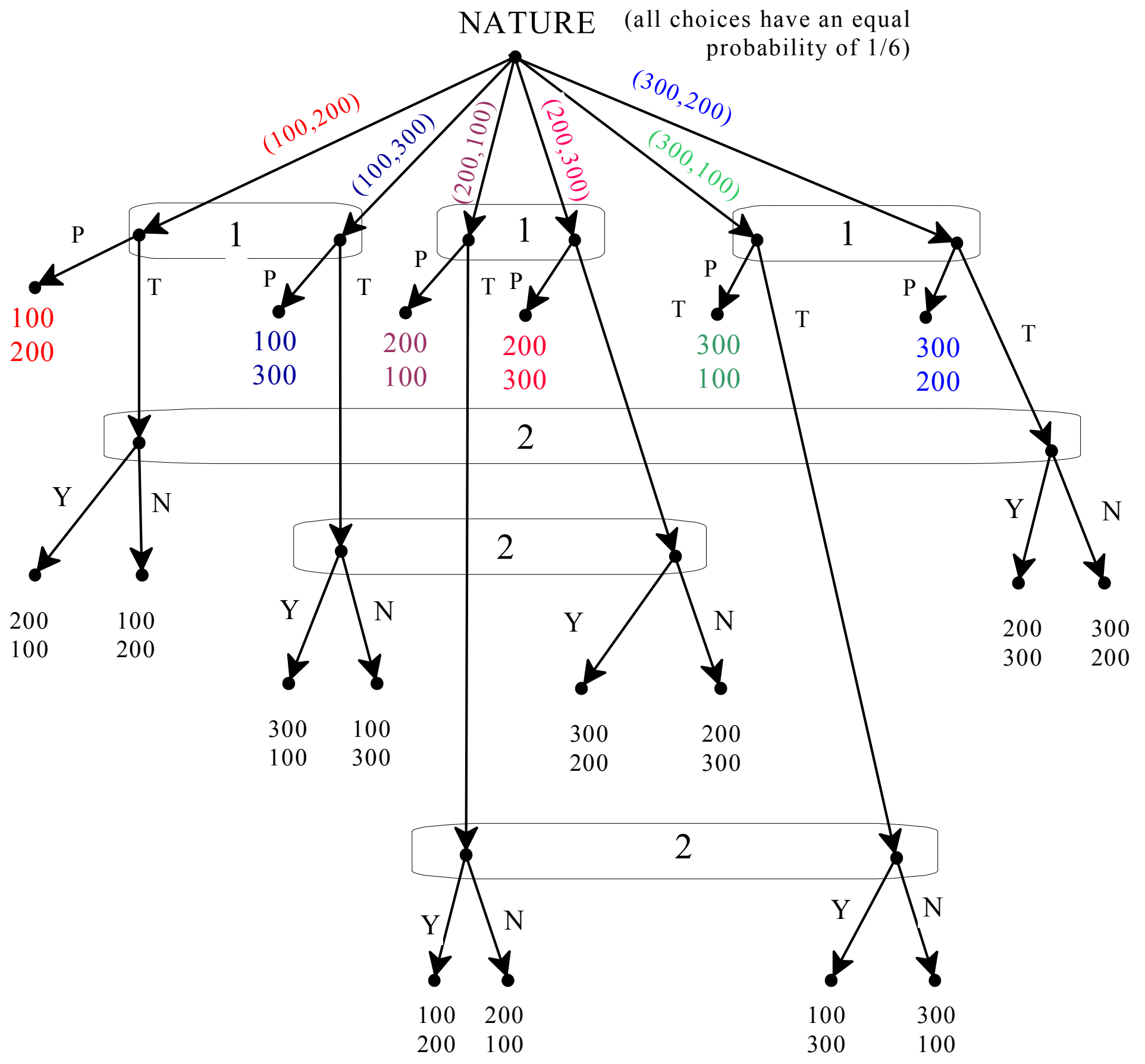


**(b)** Player 1's strategies are the same as in Example 3.3. Player 2 now has 8 strategies. Each strategy has to specify how to reply to Player 1's proposal depending on the sum he (Player 2) has. Thus one possible strategy is: if I have \$100 I say No, if I have \$200 I say Yes and if I have \$300 I say No.

**Exercise 3.12.** **(a)** The game is as follows.

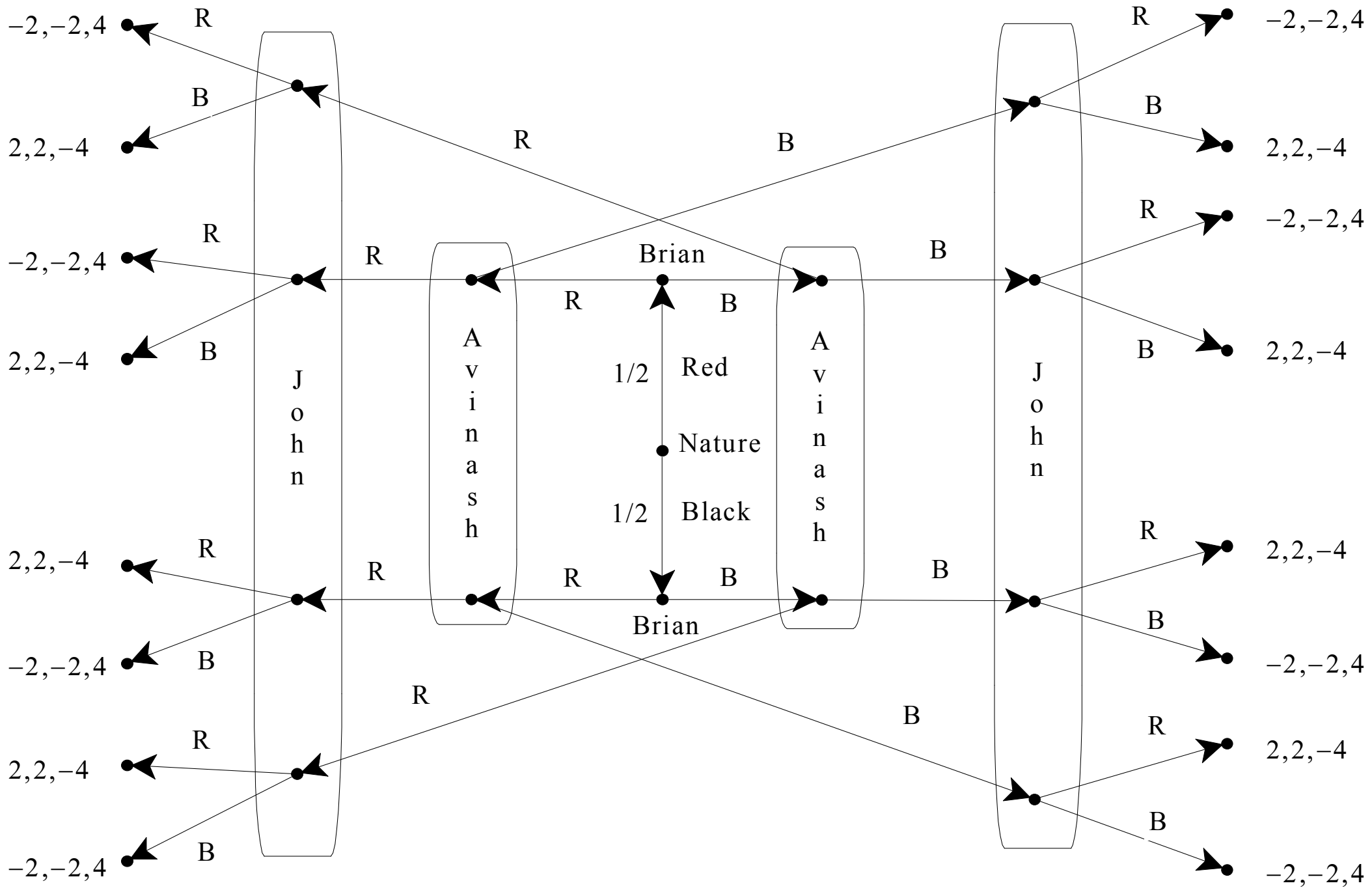


**(b)** Each player has two information sets, two choices at each information set, hence four strategies. The strategic form is given by the following set of tables.

**Avinash**

| | | if B, B, if R, R | if B, B, if R, B | if B, R, if R, R | if B, R, if R, B |
|---|---|---|---|---|---|
| **Brian** | if B, B, if R, R | −2,−2, 4 | 0, 0, 0 | 0, 0, 0 | 2, 2, −4 |
| | if B, B, if R, B | 0, 0, 0 | 0, 0, 0 | 0, 0, 0 | 0, 0, 0 |
| | if B, R, if R, R | 0, 0, 0 | 0, 0, 0 | 0, 0, 0 | 0, 0, 0 |
| | if B, R, if R, B | 2, 2, −4 | 0, 0, 0 | 0, 0, 0 | −2,−2, 4 |

**John** chooses: if B, B and if R, R

**Avinash**

| **Brian** | if B, B, if R, R | if B, B, if R, B | if B, R, if R, R | if B, R, if R, B |
|---|---|---|---|---|
| if B, B, if R, R | 0, 0, 0 | 0, 0, 0 | 0, 0, 0 | 0, 0, 0 |
| if B, B, if R, B | 0, 0, 0 | 0, 0, 0 | 0, 0, 0 | 0, 0, 0 |
| if B, R, if R, R | 0, 0, 0 | 0, 0, 0 | 0, 0, 0 | 0, 0, 0 |
| if B, R, if R, B | 0, 0, 0 | 0, 0, 0 | 0, 0, 0 | 0, 0, 0 |

**John** chooses: if B, B and if R, B

---

**Avinash**

| **Brian** | if B, B, if R, R | if B, B, if R, B | if B, R, if R, R | if B, R, if R, B |
|---|---|---|---|---|
| if B, B, if R, R | 0, 0, 0 | 0, 0, 0 | 0, 0, 0 | 0, 0, 0 |
| if B, B, if R, B | 0, 0, 0 | 0, 0, 0 | 0, 0, 0 | 0, 0, 0 |
| if B, R, if R, R | 0, 0, 0 | 0, 0, 0 | 0, 0, 0 | 0, 0, 0 |
| if B, R, if R, B | 0, 0, 0 | 0, 0, 0 | 0, 0, 0 | 0, 0, 0 |

**John** chooses: if B, R and if R, R

---

| | | **Avinash** | | | |
|---|---|---|---|---|---|
| | | if B, B, if R, R | if B, B, if R, B | if B, R, if R, R | if B, R, if R, B |
| | if B, B, if R, R | 2, 2, −4 | 0, 0, 0 | 0, 0, 0 | −2, −2, 4 |
| **Brian** | if B, B, if R, B | 0, 0, 0 | 0, 0, 0 | 0, 0, 0 | 0, 0, 0 |
| | if B, R, if R, R | 0, 0, 0 | 0, 0, 0 | 0, 0, 0 | 0, 0, 0 |
| | if B, R, if R, B | −2, −2, 4 | 0, 0, 0 | 0, 0, 0 | 2, 2, −4 |

**John** chooses: if B, R and if R, B

How can we fill in the payoffs without spending more than 24 hours on this problem? There is a quick way of doing it. First of all, when John's strategy is to guess Black, no matter what Avinash says, then he has a 50% chance of being right and a 50% chance of being wrong. Thus his expected payoff is $\frac{1}{2}(4)+\frac{1}{2}(-4)=0$ and, similarly, the expected payoff of each of the other two players is $\frac{1}{2}(2)+\frac{1}{2}(-2)=0$. This explains why the second table is filled with the same payoff vector, namely (0,0,0). The same reasoning applies to the case where when John's strategy is to guess Red, no matter what Avinash says (leading to the third table, filled with the same payoff vector (0,0,0)).

For the remaining strategies of John's, one can proceed as follows. Start with the two colors, B and R. Under B write T (for true) if Brian's strategy says "if B then B" and write F (for false) if Brian's strategy says "if B then R"; similarly, under R write T (for true) if Brian's strategy says "if R then R" and write F (for false) if Brian's strategy says "if R then B". In the next row, in the B column rewrite what is in the previous row if Avinash's strategy says "if B then B" and change a T into an F or an F into a T if Avinash's strategy says "if B then R". Similarly for the R column. Now repeat the same for John (in the B column a T remains a T and an F remains an F is John's strategy is "if B then B", while a T is changed into an F and an F is changed into a T if John's strategy is "if B then R"). Now in each column the payoffs are (−2, −2, 4) if the last row has a T and (2, 2, −4) if the last row has an F. The payoffs are then given by $\frac{1}{2}$ the payoff in the left column plus $\frac{1}{2}$ the payoff in the right column. For example,

for the cell in the second row, third column of the third table we have the following.

| | B | R | |
|---|---|---|---|
| Brian's strategy: if B, B and if R, B | T | F | |
| Avinash's strategy: if B, R and if R, R | F | T | |
| John's strategy: if B, R and if R, R | F | T | Payoffs: |
| Payoffs | (2, 2, −4) | (−2, −2, 4) | $\frac{1}{2}(2, 2, -4) + \frac{1}{2}(-2, -2, 4) = (0, 0, 0)$ |

**Exercise 3.13**. **(a)** The extensive-form representation of the simplified poker game is as follows (the top number is Yvonne's net take in dollars and the bottom number is Zoe's net take).

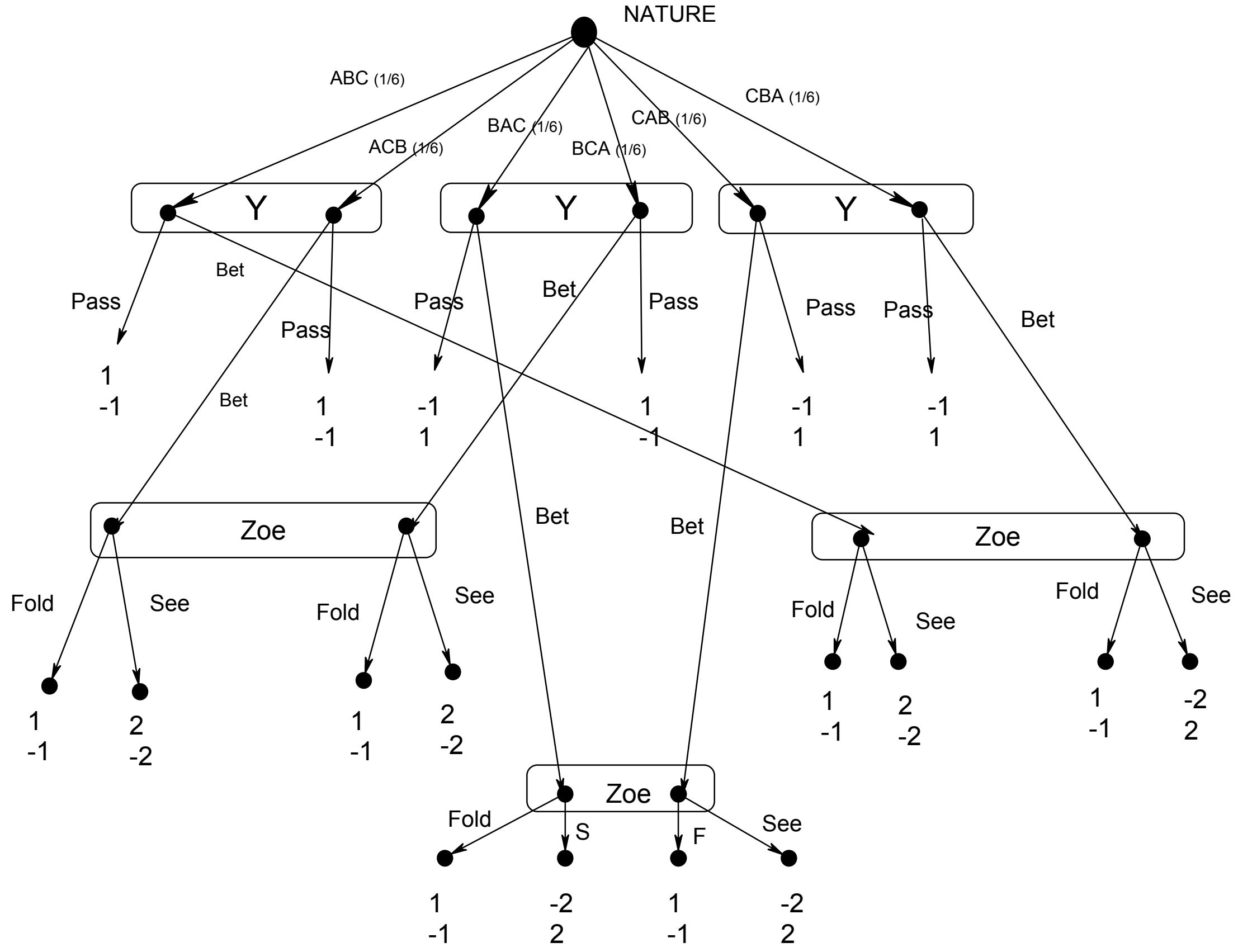

**(b) and (c)** Each player has eight strategies (three information sets, two choices at each information set, thus $2 \times 2 \times 2 = 8$ possible strategies).

**(d)** Yvonne uses the strategy "If A pass, if B pass, if C bet" and Zoe uses the strategy "If A fold, if B fold, if C fold"). We calculate the expected net payoff for Yvonne. Zoe's expected net payoff is the negative of that.

| Possible Cases | Top card is A<br>Second is B | A<br>C | B<br>A | B<br>C | C<br>A | C<br>B | Sum | Probability of each | Expected payoff |
|---|---|---|---|---|---|---|---|---|---|
| Y's payoff | 1 | 1 | –1 | 1 | 1 | 1 | 4 | $\frac{1}{6}$ | $\frac{4}{6}$ |
| Explain | pass | pass | pass | pass | bet + fold | bet + fold | | | |

**(e)** Yvonne uses the strategy "If A pass, if B pass, if C bet" and Zoe uses the strategy "see with any card". Once again, we calculate the expected net payoff for Yvonne. Zoe's expected net payoff is the negative of that.

| Possible Cases | Top card is A<br>Second is B | A<br>C | B<br>A | B<br>C | C<br>A | C<br>B | Sum | Probability of each | Expected payoff |
|---|---|---|---|---|---|---|---|---|---|
| Y's payoff | 1 | 1 | –1 | 1 | –2 | –2 | –2 | $\frac{1}{6}$ | $-\frac{2}{6}$ |
| Explain | pass | pass | pass | pass | bet + see | bet + see | | | |

**(f)** The strategic form is as follows.

| | **ZOE**<br>If A fold,<br>If B fold,<br>If C fold | If A see,<br>If B see,<br>If C see | If A see,<br>If B fold,<br>If C fold | If A fold,<br>If B see,<br>If C fold | If A fold,<br>If B fold,<br>If C see | If A see,<br>If B see,<br>If C fold | If A see,<br>If B fold,<br>If C see | If A fold,<br>If B see,<br>If C see |
|---|---|---|---|---|---|---|---|---|
| If A pass, if B pass if C pass | 0, 0 | 0, 0 | 0, 0 | 0, 0 | 0, 0 | 0, 0 | 0, 0 | 0, 0 |
| If A bet, if B bet if C bet | 1, −1 | 0, 0 | 0, 0 | 4/6, −4/6 | 8/6, −8/6 | −2/6, 2/6 | 2/6, −2/6 | 1, −1 |
| If A bet, if B pass if C pass | 0, 0 | 2/6, −2/6 | 0, 0 | 1/6, −1/6 | 1/6, −1/6 | 1/6, −1/6 | 1/6, −1/6 | 2/6, −2/6 |
| If A pass, if B bet if C pass | 2/6, −2/6 | 0, 0 | −1/6, 1/6 | 2/6, −2/6 | 3/6, −3/6 | −1/6, 1/6 | 0, 0 | 3/6, −3/6 |
| If A pass, if B pass if C bet | 4/6, −4/6 | −2/6, 2/6 | 1/6, −1/6 | 1/6, −1/6 | 4/6, −4/6 | −2/6, 2/6 | 1/6, −1/6 | 1/6, −1/6 |
| If A bet, if B bet if C pass | 2/6, −2/6 | 2/6, −2/6 | −1/6, 1/6 | 3/6, −3/6 | 4/6, −4/6 | 0, 0 | 1/6, −1/6 | 5/6, −5/6 |
| If A bet, if B pass if C bet | 4/6, −4/6 | 0, 0 | 1/6, −1/6 | 2/6, −2/6 | 5/6, −5/6 | −1/6, 1/6 | 2/6, −2/6 | 3/6, −3/6 |
| if A Pass, If B or C, Bet, | 1, −1 | −2/6, 2/6 | 0, 0 | 3/6, −3/6 | 7/6, −7/6 | −3/6, 3/6 | 1/6, −1/6 | 4/6, −4/6 |

**(g)** Let $\searrow$ denote weak dominance, that is, $a \searrow b$ means that $a$ weakly dominates $b$.

FOR YVETTE (row player): $3^{rd}$ row $\searrow$ $1^{st}$ row, $6^{th} \searrow 4^{th}$, $7^{th} \searrow 4^{th}$, $7^{th} \searrow 5^{th}$, $2^{nd} \searrow 8^{th}$.

FOR ZOE (column player): $3^{rd}$ col $\searrow$ $1^{st}$ col, $3^{rd} \searrow 4^{th}$, $3^{rd} \searrow 5^{th}$, $3^{rd} \searrow 7^{th}$, $3^{rd} \searrow 8^{th}$, $2^{nd} \searrow 8^{th}$, $4^{th} \searrow 5^{th}$, $4^{th} \searrow 8^{th}$, $6^{th} \searrow 2^{nd}$, $6^{th} \searrow 4^{th}$, $6^{th} \searrow 5^{th}$,
$6^{th} \searrow 7^{th}$, $6^{th} \searrow 8^{th}$, $7^{th} \searrow 4^{th}$, $7^{th} \searrow 5^{th}$, $7^{th} \searrow 8^{th}$.

**(h)** Eliminating rows 1, 4, 5 and 8 and all columns except 3 and 6 we are left with:

| | | **Zoe** | |
|---|---|---|---|
| | | See only if A | See with A or B |
| **Y** | Bet always | 0, 0 | **−2/6, 2/6** |
| | If A, Bet, otherwise pass | 0, 0 | **1/6, −1/6** |
| | If A or B, Bet, if C Pass | **−1/6, 1/6** | 0, 0 |
| | If A or C, Bet, if B Pass | **1/6, −1/6** | **−1/6, 1/6** |

In the above table, the second row dominates the first and the third. Eliminating them we have the following, which is a remarkable simplification of the original strategic form:

| | | **Zoe** | |
|---|---|---|---|
| | | See only if A | See with A or B |
| **Y** | If A, Bet, otherwise pass | 0, 0 | **1/6, −1/6** |
| | If A or C, Bet, if B Pass | **1/6, −1/6** | **−1/6, 1/6** |

**Exercise 3.14.** **(a)** The extensive form is as follows:

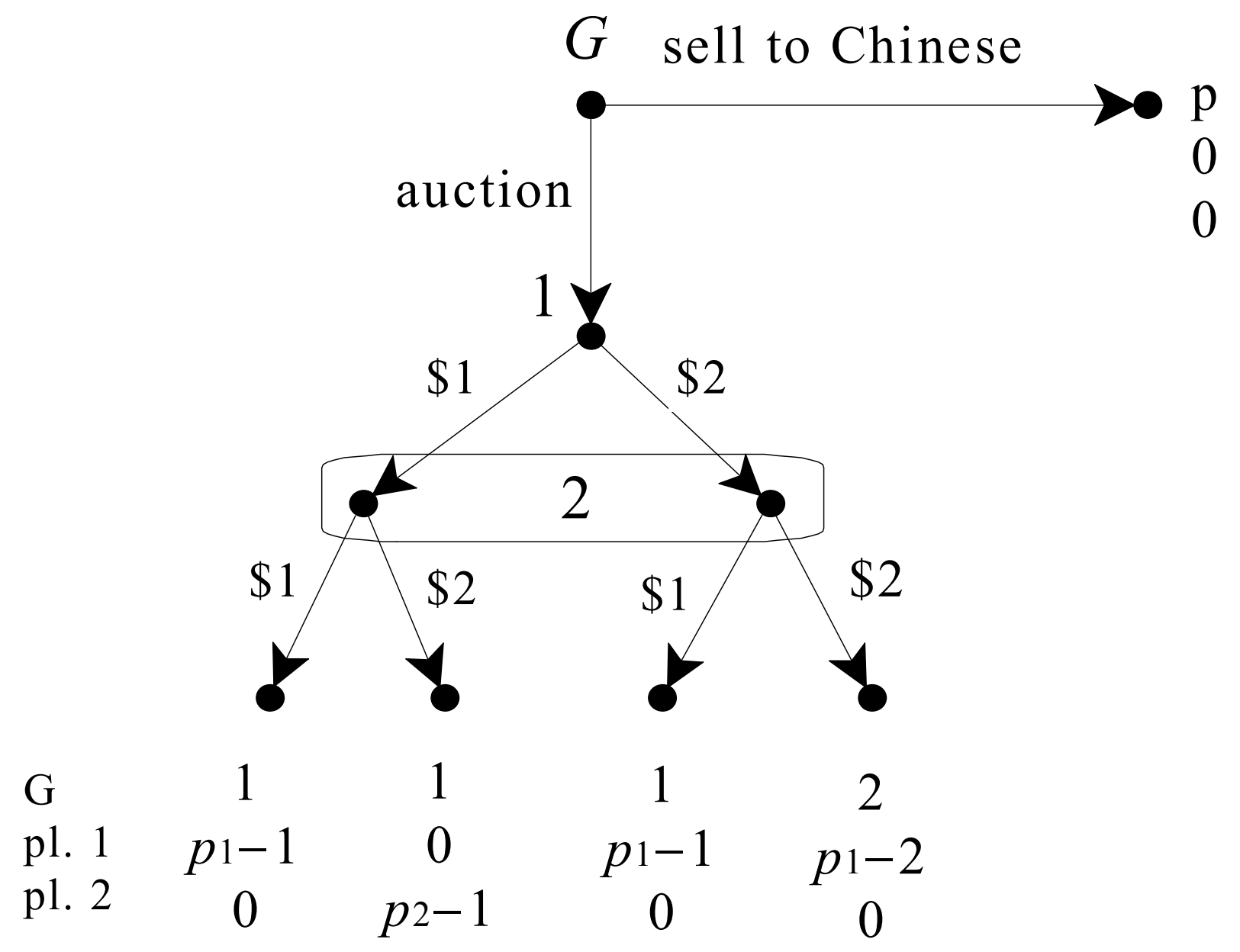


**(b)** In the auction subgame for every player it is a weakly dominant strategy to bid his own value. Let $p_j = \max_{i=1,\dots,n} \{p_i\}$ be the highest value and $p_k = \max_{\substack{i=1,\dots,n \\ i \neq j}} \{p_i\}$ be the second highest value. Then the auction, if it takes place, will be won by Player *j* and he will pay $p_k$. Hence there are three cases. **Case 1**: $p > p_k$. In this case Player G will sell to the Chinese (and the strategy of Player *i* in the subgame is to bid $p_i$), G's payoff is *p* and the payoff of Player *i* is 0. **Case 2:** $p < p_k$. In this case Player G announces the auction, the strategy of Player *i* in the subgame is to bid $p_i$, the winner is Player *j* and he pays $p_k$, so that the payoff of G is $p_k$, the payoff of player *j* is $p_j - p_k$ and the payoff of every other player is 0. **Case 3:** $p = p_k$. In this case there are two subgame-perfect equilibria: one as in Case 1 and the other as in Case 2 and G is indifferent between the two.

**(c)** The game is as follows:

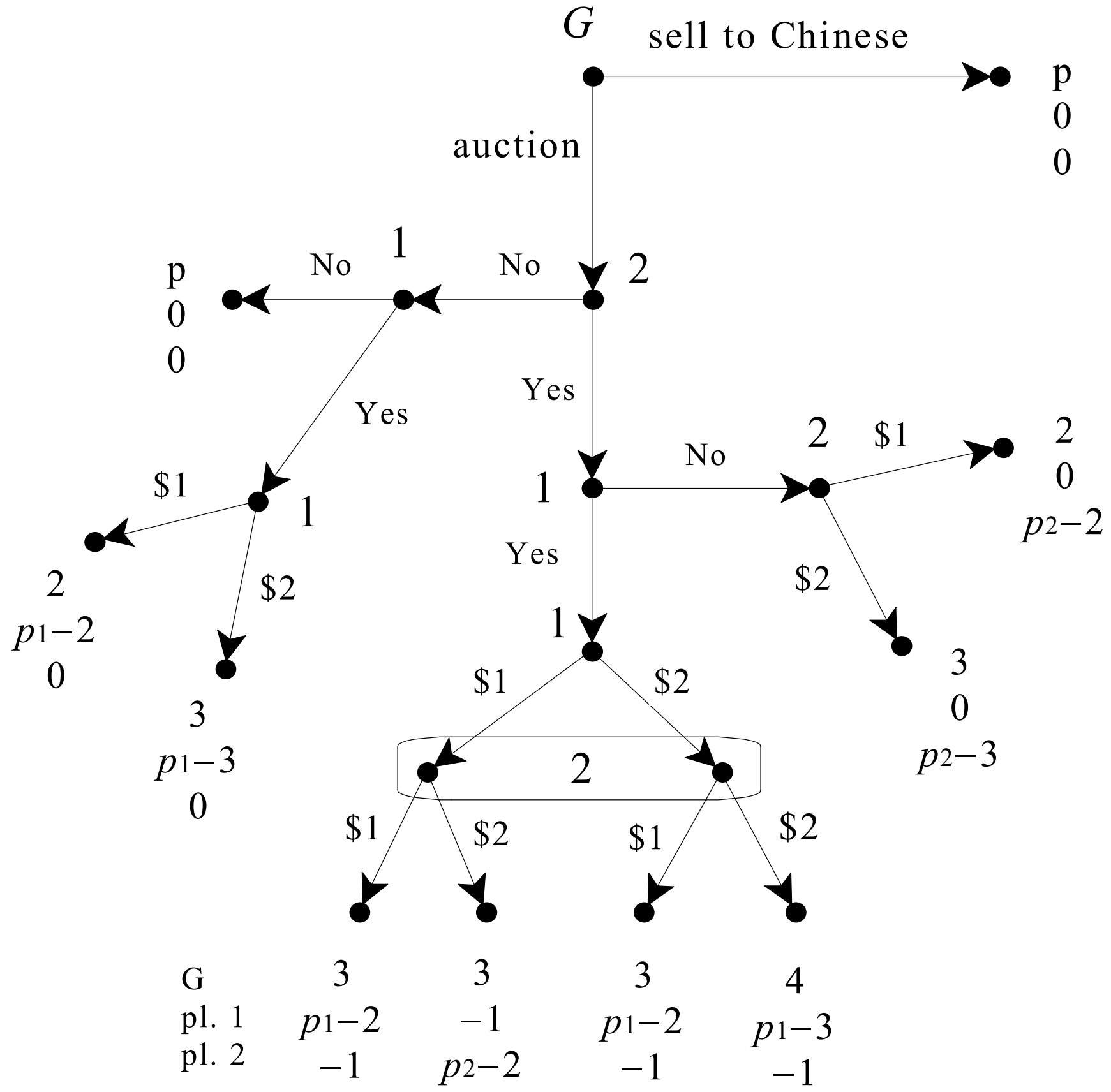


**(d)** In the simultaneous subgame after both players have said Yes, the participation fee paid is a sunk cost and for every player bidding the true value is a weakly dominant strategy. Thus the outcome there is as follows: Player 1 bids $p_1$, gets the palace by paying $p_2$, $G$'s payoff is $(p_2 + 2)$, 1's payoff is $(p_1 - p_2 - 1)$ and Player 2's payoff is $-1$. In the subgames where one player said No and the other said Yes the optimal choice is obviously $x = 1$, with payoffs of 2 for Player $G$, 0 for the player who said No and $p_i - 2$ for the player who said Yes. Thus the game reduces to the one shown in Figure A below:

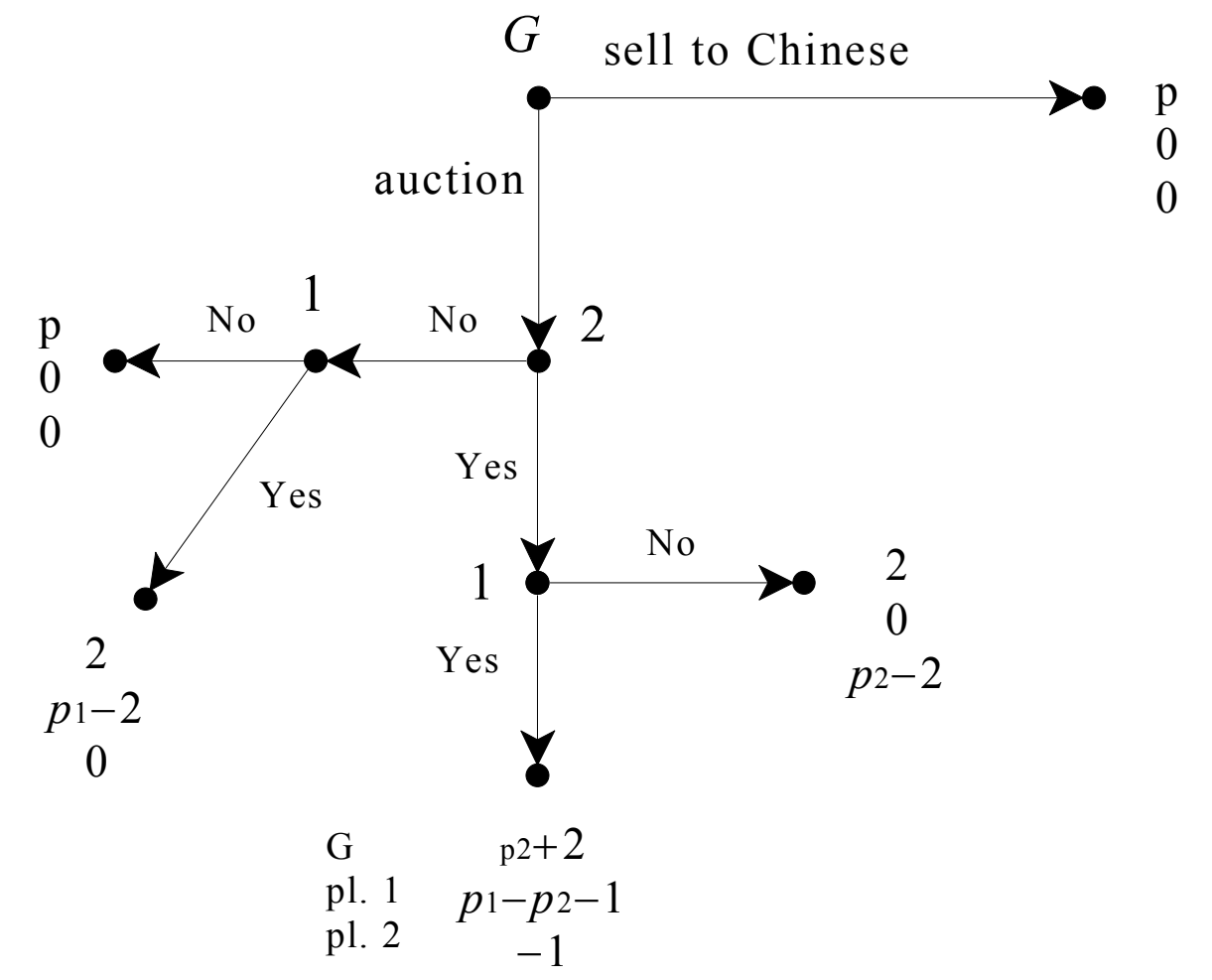


Figure A

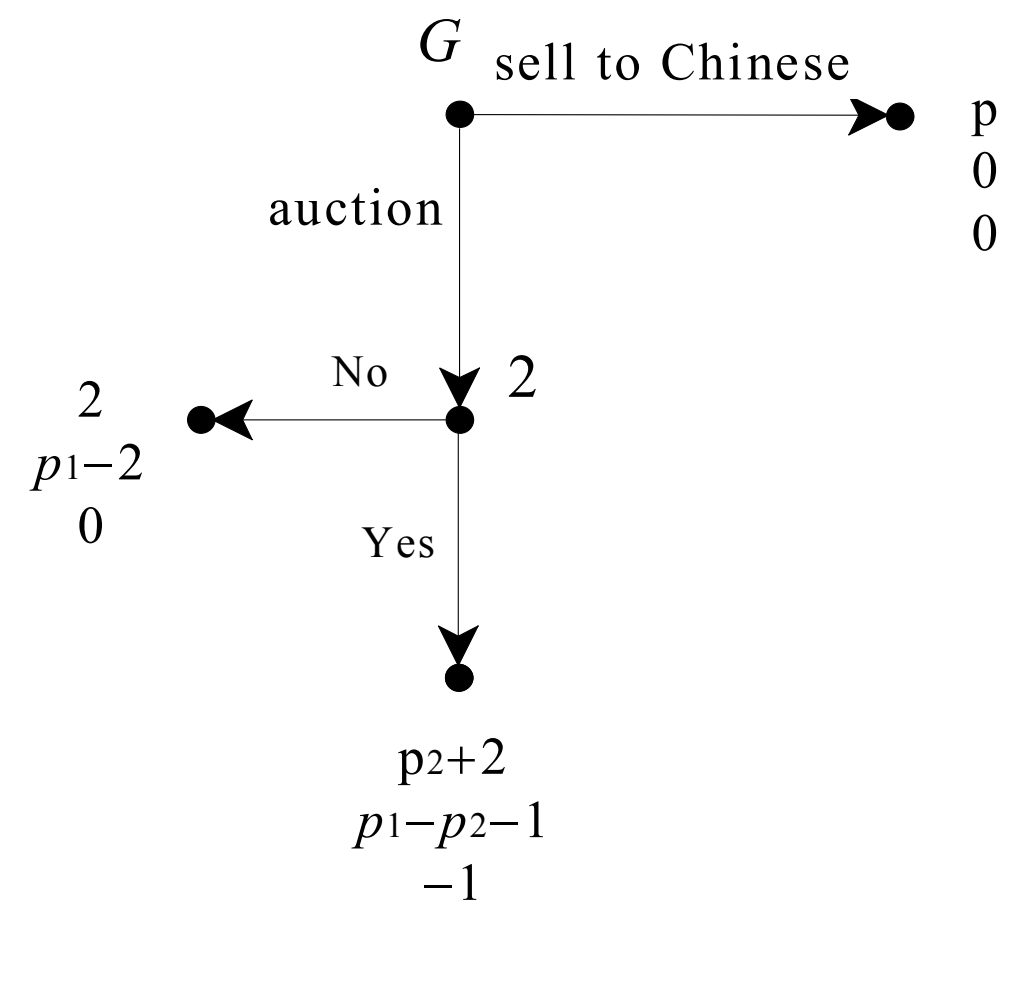


Figure B

By assumption, $p_1 > p_2 + 1 > 2$, so that $p_1 - p_2 - 1 > 0$ and $p_1 - 2 > 0$. Thus at the bottom node and at the left node player 1 prefers Yes to No. Thus the game reduces to the one shown in Figure B.
Hence Player 2 will say No. Thus the subgame-perfect equilibrium is as follows:
(1) if $p > 2$ then player $G$ will sell to the Chinese (and the choices off the equilibrium path are as explained above) and the payoffs are $(p, 0, 0)$;
(2) if $p < 2$ then $G$ chooses to auction, 2 says No, 1 says Yes and then offers \$1 and the payoffs are $(2, p_1 - 2, 0)$ (and the choices off the equilibrium path are as explained above);
(3) if $p = 2$ then there are two equilibria: one as in (1) and the other as in (2).

**(e)** When the loser is given the fraction $a$ of the amount paid by the winner (that is, the loser is given the fraction $a$ of his own bid), **it is no longer true that bidding one's true value is a dominant strategy**. In fact, $(p_1, p_2)$ is not even a Nash equilibrium any more. To see this, imagine that Player 1's true value is 10 and Player 2's true value is 6 and $a = 50\%$. Then if Player 1 bids 10 and 2 bids 6, Player 2 ends up losing the auction but being given \$3, while if he increased his bid to 8 then he would still lose the auction but receive \$4. This shows that there cannot be a Nash equilibrium where Player 2 bids less than Player 1. Now there are several Nash equilibria of the auction, for example, all pairs $(b_1, b_2)$ with $b_1 = b_2 = b$ and $p_2 \le b < p_1$ provided that $p_1 - b \ge a(b-1)$, that is, $b \le \frac{p_1 + a}{1 + a}$ (but there are more: for example all pairs $(b_1, b_2)$ with $b_1 = b_2 = b$ and $b < p_2$ provided that $p_1 - b \ge a(b-1)$ and $ab \ge p_2 - b$). Thus to find a subgame-perfect

equilibrium of the game one first has to select a Nash equilibrium of the auction game and then apply backward induction to see if the players would want to say Yes or No to the auction, etc.

**(f)** Let us start by considering the perfect-information game that is played if the Chinese says Yes. This is a game similar to the one discussed in Example 2.2 (Chapter 2, Section 2.5). We first determine the losing positions. Whoever has to move when the sum is 36 cannot win. Thus 36 is a losing position. Working backwards, the losing positions are 32, 28, 24, 20, 16, 12, 8, 4 and 0. Thus the first player (= player $G$) starts from a losing position: whatever his initial choice, he can be made to choose the second time when the sum is 4, and then 8, etc. Hence **the second player (= the Chinese) has a winning strategy**, which is as follows: if Player $G$ just chose $n$, then choose $(4 - n)$. If the Chinese says Yes and then follows this strategy he can guarantee that he will buy the palace for \$50. Thus the subgame-perfect equilibrium of this game is: the Chinese says Yes and uses the winning strategy in the ensuing game, while for Player $G$ we can pick any arbitrary choices.

# PART II

# Games with Cardinal Payoffs

Chapter 4

# Expected Utility

## 4.1 Money lotteries and attitudes to risk

At the end of Chapter 3 we showed how to incorporate random events in extensive-form games by means of chance moves. The introduction of chance moves gives rise to probabilistic outcomes, which we called lotteries. In Chapter 3 we restricted attention to lotteries whose outcomes are sums of money (money lotteries) and to one possible way of ranking such lotteries, based on the notion of risk neutrality. In this section we will continue to focus on money lotteries and define other possible attitudes to risk.[1]

Throughout this chapter we will restrict attention to finite lotteries. Recall that a money lottery is a probability distribution of the form

$$\begin{pmatrix} \$x_1 & \$x_2 & \dots & \$x_n \\ p_1 & p_2 & \dots & p_n \end{pmatrix} \quad (0 \le p_i \le 1, \text{ for all } i = 1,2,\dots,n \text{ and } p_1 + p_2 + \dots + p_n = 1)$$

and that (Definition 3.13, Chapter 3) its expected value is the number $(x_1 p_1 + x_2 p_2 + \dots + x_n p_n)$. If $L$ is a money lottery, we will denote by $\mathbb{E}[L]$ the expected value of $L$. Thus, for example, if $L = \begin{pmatrix} \$30 & \$45 & \$90 \\ \frac{1}{3} & \frac{5}{9} & \frac{1}{9} \end{pmatrix}$ then $\mathbb{E}[L] = \frac{1}{3}(30) + \frac{5}{9}(45) + \frac{1}{9}(90) = 45$. Recall also (Definition 3.14, Chapter 3) that a person is said to be risk neutral if she considers a money lottery to be just as good as its expected value. For example, a risk-neutral person would consider getting \$45 with certainty to be just as good a playing lottery

[1] In the next section we will consider more general lotteries, where the outcomes need not be sums of money, and also introduce the theory of expected utility which will then be used in the context of games in Chapters 5 and 6.

$L = \begin{pmatrix} \$30 & \$45 & \$90 \\ \frac{1}{3} & \frac{5}{9} & \frac{1}{9} \end{pmatrix}$. We can now consider different attitudes to risk, besides risk neutrality.

**Definition 4.1.** Let $L$ be a money lottery and consider the choice between $L$ and the degenerate lottery $\begin{pmatrix} \$\mathbb{E}[L] \\ 1 \end{pmatrix}$ (that is, the choice between facing the lottery $L$ or getting the expected value of $L$ with certainty). Then

- An individual who prefers $\$\mathbb{E}[L]$ for certain to $L$ is said to be *risk averse.*
- An individual who is indifferent between $\$\mathbb{E}[L]$ for certain and $L$ is said to be *risk neutral.*
- An individual who prefers $L$ to $\$\mathbb{E}[L]$ for certain is said to be *risk loving.*

Note the following: if an individual is risk neutral, has transitive preferences over money lotteries and prefers more money to less, then we can tell how that individual ranks any two money lotteries. For example, how would such an individual rank the two lotteries $L_1 = \begin{pmatrix} \$30 & \$45 & \$90 \\ \frac{1}{3} & \frac{5}{9} & \frac{1}{9} \end{pmatrix}$ and $L_2 = \begin{pmatrix} \$5 & \$100 \\ \frac{3}{5} & \frac{2}{5} \end{pmatrix}$? Since $\mathbb{E}[L_1] = 45$ and the individual is risk neutral, $L_1 \sim \$45$; since $\mathbb{E}[L_2] = 43$ and the individual is risk neutral, $\$43 \sim L_2$; since the individual prefers more money to less, $\$45 \succ \$43$; thus, by transitivity, $L_1 \succ L_2$. On the other hand, knowing that an individual is risk averse, has transitive preferences over money lotteries and prefers more money to less is not sufficient to predict how she will choose between two arbitrary money lotteries. For example, as we will see later (Exercise 4.11), it is possible that one such individual will prefer $L_3 = \begin{pmatrix} \$28 \\ 1 \end{pmatrix}$ (whose expected value is 28) to $L_4 = \begin{pmatrix} \$10 & \$50 \\ \frac{1}{2} & \frac{1}{2} \end{pmatrix}$ (whose expected value is 27), while another such individual will prefer $L_4$ to $L_3$. Similarly, knowing that an individual is risk loving, has transitive preferences over money lotteries and prefers more money to less is not sufficient to predict how she will choose between two arbitrary money lotteries.

**Remark 4.1.** Note that "rationality" does not and cannot dictate whether an individual should be risk neutral, risk averse or risk loving: an individual's attitude to risk is merely a reflection of that individual's preferences. It is a generally accepted principle that *de gustibus non est disputandum* (in matters of taste, there can be no disputes). According to this principle, there is no such thing as an irrational preference and thus there is no such thing as an irrational attitude to risk. As a matter of fact, most people reveal through their choices (e.g. the decision to buy insurance) that they are risk averse, at least when the stakes are high.

As noted above, with the exception of risk-neutral individuals, even if we restrict attention to money lotteries we are not able to say much about how an individual would choose among lotteries. What we need is a theory of "rational" preferences over lotteries that is general enough to cover lotteries whose outcomes are not just sums of money and is capable of accounting for different attitudes to risk in the case of money lotteries. One such theory is the *theory of expected utility*, to which we now turn.

Test your understanding of the concepts introduced in this section, by going through the exercises in Section 4.E.1 of Appendix 4.E at the end of this chapter.

## 4.2 Expected utility: theorems

The theory of expected utility was developed by the founders of game theory, namely John Von Neumann and Oskar Morgenstern, in their 1944 book *Theory of Games and Economic Behavior*. In a rather unconventional way, we shall first (in this section) state the main result of the theory (which we split into two theorems) and then (in the following section) explain the assumptions (or axioms) behind that result. The reader who is not interested in understanding the conceptual foundation of expected utility theory, but wants to understand what the theory says and how it can be used, can thus study this section and skip the next. First we review the notation.

Let $O$ be a set of *basic outcomes*. Note that a basic outcome need not be a sum of money: it could be the state of an individual's health, or whether the individual under consideration receives an award, or whether it will rain on the day of her planned outdoor party, etc. Let $\mathcal{L}(O)$ be the set of *simple lotteries* (or

probability distributions) over O. We will assume throughout that O is a finite set: $O = \{o_1, o_2, ..., o_m)$ . Thus an element of $\mathcal{L}(O)$ is of the form

$$\begin{pmatrix} o_1 & o_2 & \dots & o_m \\ p_1 & p_2 & \dots & p_m \end{pmatrix} \text{ with } 0 \le p_i \le 1, \text{ for all } i = 1, 2, ..., m \text{ and } p_1 + p_2 + ... + p_m = 1.$$

We will use the symbol $L$ (with or without subscript) to denote an element of $\mathcal{L}(O)$, that is, a simple lottery. Lotteries are used to represent situations of uncertainty. For example, if $m$ = 4 and the individual faces the lottery $L = \begin{pmatrix} o_1 & o_2 & o_3 & o_4 \\ \frac{2}{5} & 0 & \frac{1}{5} & \frac{2}{5} \end{pmatrix}$ then she knows that, eventually, the outcome will be one and only one of $o_1, o_2, o_3, o_4$, but does not know which one; furthermore, she is able to quantify her uncertainty by assigning probabilities to these outcomes. We take these probabilities either as objectively obtained from relevant past data or as subjective estimates by the individual. For example, an individual who is considering whether or not to insure her bicycle against theft for the following 12 months, knows that there are two relevant basic outcomes: either the bicycle will be stolen or it will not be stolen. Furthermore, she can look up data about past bicycle thefts in her area and use the proportion of bicycles that was stolen as an "objective" estimate of the probability that her bicycle will be stolen; alternatively, she can use a more subjective estimate: for example she might use a lower probability of theft than suggested by the data because she knows herself to be very conscientious and – unlike other people – to always lock her bicycle when left unattended.

The assignment of zero probability to a particular basic outcome is taken to be an expression of *belief* not of impossibility: the individual is *confident* that the outcome will not arise, but she cannot rule out that outcome on logical grounds or by appealing to the laws of nature.

Among the elements of $\mathcal{L}(O)$ there are the *degenerate lotteries* that assign probability 1 to one basic outcome: for example, if $m$ = 4 one degenerate lottery is $\begin{pmatrix} o_1 & o_2 & o_3 & o_4 \\ 0 & 0 & 1 & 0 \end{pmatrix}$. To simplify the notation we will often denote degenerate lotteries as basic outcomes, that is, instead of writing $\begin{pmatrix} o_1 & o_2 & o_3 & o_4 \\ 0 & 0 & 1 & 0 \end{pmatrix}$ we will simply write $o_3$. Thus, in general, the degenerate lottery $\begin{pmatrix} o_1 & \dots & o_{i-1} & o_i & o_{i+1} & \dots & o_m \\ 0 & 0 & 0 & 1 & 0 & 0 & 0 \end{pmatrix}$ will be denoted by $o_i$. As another simplification, we will often omit those outcomes that are assigned zero

probability. For example, if $m$ = 4, the lottery $\begin{pmatrix} o_1 & o_2 & o_3 & o_4 \\ \frac{1}{3} & 0 & \frac{2}{3} & 0 \end{pmatrix}$ will be written more simply as $\begin{pmatrix} o_1 & o_3 \\ \frac{1}{3} & \frac{2}{3} \end{pmatrix}$.

Throughout this chapter we shall call the individual under consideration the *Decision-Maker*, or DM for short. The theory of expected utility assumes that the DM has a complete and transitive ranking $\succsim$ of the elements of $\mathcal{L}(O)$ (indeed, this is one of the axioms listed in the next section). By completeness, given any two lotteries $L$ and $L'$, either $L \succ L'$ (the DM prefers $L$ to $L'$) or $L' \succ L$ (the DM prefers $L'$ to $L$) or $L \sim L'$ (the DM is indifferent between $L$ and $L'$). Furthermore, by transitivity, for any three lotteries $L_1, L_2$ and $L_3$, if $L_1 \succsim L_2$ (the DM considers $L_1$ to be at least as good as $L_2$) and $L_2 \succsim L_3$ then $L_1 \succsim L_3$. Besides completeness and transitivity, a number of other "rationality" constraints are postulated on the ranking $\succsim$ of the elements of $\mathcal{L}(O)$; these constraints are the so-called Expected Utility Axioms and are discussed in the next section.

**Definition 4.2.** A ranking $\succsim$ of the elements of $\mathcal{L}(O)$ that satisfies the Expected Utility Axioms is called a *von Neumann-Morgenstern ranking*.

The following two theorems are the key results in the theory of expected utility.

**Theorem 4.1 [von Neumann-Morgenstern, 1944].** Let $O = \{o_1, o_2, ..., o_m\}$ be a set of basic outcomes and $\mathcal{L}(O)$ the set of lotteries over O. If $\succsim$ is a von Neumann-Morgenstern ranking of the elements of $\mathcal{L}(O)$ then there exists a function $U : O \to \mathbb{R}$, called a *von Neumann-Morgenstern utility function*, that assigns a number (called utility) to every basic outcome and is such that, for any two lotteries $L = \begin{pmatrix} o_1 & o_2 & \dots & o_m \\ p_1 & p_2 & \dots & p_m \end{pmatrix}$ and $L' = \begin{pmatrix} o_1 & o_2 & \dots & o_m \\ q_1 & q_2 & \dots & q_m \end{pmatrix}$,

$$L \succ L' \text{ if and only if } \mathbb{E}[U(L)] > \mathbb{E}[U(L')] \text{ and}$$

$$L \sim L' \text{ if and only if } \mathbb{E}[U(L)] = \mathbb{E}[U(L')]$$

where

$$U(L) = \begin{pmatrix} U(o_1) & U(o_2) & \dots & U(o_m) \\ p_1 & p_2 & \dots & p_m \end{pmatrix}, \; U(L') = \begin{pmatrix} U(o_1) & U(o_2) & \dots & U(o_m) \\ q_1 & q_2 & \dots & q_m \end{pmatrix},$$

$\mathbb{E}[U(L)]$ is the expected value of the lottery $U(L)$ and $\mathbb{E}[U(L')]$ is the expected value of the lottery $U(L')$, that is,

$$\mathbb{E}[U(L)] = p_1 U(o_1) + p_2 U(o_2) + ... + p_m U(o_m)$$

$$\mathbb{E}[U(L')] = q_1 U(o_1) + q_2 U(o_2) + ... + q_m U(o_m).$$

$\mathbb{E}[U(L)]$ is called the *expected utility of lottery* $L$ (and $\mathbb{E}[U(L')]$ the expected utility of lottery $L'$).
We say that any function $U : O \to \mathbb{R}$ that satisfies the property that, for any two lotteries $L$ and $L'$, $L \succsim L'$ if and only if $\mathbb{E}[U(L)] \geq \mathbb{E}[U(L')]$ *represents the preferences (or ranking)* $\succsim$.

Before we comment on Theorem 4.1 we give an example of how one can use of it. Theorem 4.1 sometimes allows us to predict an individual's choice between two lotteries $C$ and $D$ if we know how that individual ranks two *different* lotteries $A$ and $B$. For example, suppose we observe that Susan is faced with the choice between lotteries $A$ and $B$ below and she says that she prefers $A$ to $B$

$$A = \begin{pmatrix} o_1 & o_2 & o_3 \\ 0 & 0.25 & 0.75 \end{pmatrix} \qquad B = \begin{pmatrix} o_1 & o_2 & o_3 \\ 0.2 & 0 & 0.8 \end{pmatrix}.$$

With this information we can predict which of the following two lotteries $C$ and $D$ she will choose, *if she has von Neumann-Morgenstern preferences*:

$$C = \begin{pmatrix} o_1 & o_2 & o_3 \\ 0.8 & 0 & 0.2 \end{pmatrix} \qquad D = \begin{pmatrix} o_1 & o_2 & o_3 \\ 0 & 1 & 0 \end{pmatrix} = o_2 .$$

Applying the theorem, let $U$ be a von Neumann-Morgenstern utility function whose existence is guaranteed by Theorem 4.1. Let $U(o_1) = a$, $U(o_2) = b$ and $U(o_3) = c$ (where $a$, $b$ and $c$ are numbers). Then, since Susan prefers $A$ to $B$, the expected utility of $A$ must be greater than the expected utility of $B$: $0.25b + 0.75c > 0.2a + 0.8c$. This inequality is equivalent to $0.25b > 0.2a + 0.05c$ or, dividing both sides by 0.25, $b > 0.8a + 0.2c$. It follows from this and Theorem 4.1 that Susan prefers $D$ to $C$, because the expected utility of $D$ is $b$ and the expected utility of $C$ is $0.8a + 0.2c$. Note that, in this example, we merely used the fact that a von Neumann-Morgenstern utility function *exists*, even though we don't know what the values of this function are.

Theorem 4.1 is an example of a *representation theorem* and is a generalization of a similar result for the case of the ranking of a finite set of basic outcomes $O$. It is not difficult to prove that if $\succsim$ is a complete and transitive ranking of $O$ then there exists a function $U : O \to \mathbb{R}$, called a *utility function* (see Chapter 1), such that, for any two basic outcomes $o, o' \in O$, $U(o) \geq U(o')$ if and only if $o \succsim o'$.

Now, it is quite possible that an individual has a complete and transitive ranking of *O*, is fully aware of her ranking and yet she is not able to answer the question "what is your utility function?", perhaps because she has never heard about utility functions. A utility function is a tool that we use to represent her ranking, nothing more than that. The same applies to von Neumann-Morgenstern rankings: Theorem 4.1 tells us that if an individual has a von Neumann-Morgenstern ranking of the set of lotteries $\mathcal{L}(O)$ then there exists a von Neumann-Morgenstern utility function that we can use to represent her preferences, but it would not make sense for us to ask the individual "what is your von Neumann-Morgenstern utility function?" (indeed this was a question that could not even be conceived before von Neumann and Morgenstern stated and proved Theorem 4.1 in 1944!). Theorem 4.1 tells us that such a function exists. The next theorem can be used to actually construct such a function by asking the individual to answer a few questions, formulated in a way that is fully comprehensible to her (that is, without using the word 'utility').

The next theorem says that, although there are many utility functions that represent a given von Neumann-Morgenstern ranking, once you know one function you know them all, in the sense that there is a simple operation that takes you from one function to the other.

**Theorem 4.2 [von Neumann-Morgenstern, 1944].** Let $\succsim$ be a von Neumann-Morgenstern ranking of the set of basic lotteries $\mathcal{L}(O)$, where $O = \{o_1, o_2, ..., o_m\}$. Then the following are true.

(A) If $U : O \to \mathbb{R}$ is a von Neumann-Morgenstern utility function that represents $\succsim$, then, for any two real numbers *a* and *b* with $a > 0$, the function $V : O \to \mathbb{R}$ defined by $V(o_i) = aU(o_i) + b$ ($i = 1,2,...,m$) is also a von Neumann-Morgenstern utility function that represents $\succsim$.

(B) If $U : O \to \mathbb{R}$ and $V : O \to \mathbb{R}$ are two von Neumann-Morgenstern utility functions that represent $\succsim$, then there exist two real numbers *a* and *b* with $a > 0$ such that $V(o_i) = aU(o_i) + b$ ($i = 1,2,...,m$).

**Proof**. The proof of (A) of Theorem 4.2 is very simple. Let *a* and *b* be two numbers, with $a > 0$. The hypothesis is that $U : O \to \mathbb{R}$ is a von Neumann-Morgenstern utility function that represents $\succsim$, that is, that, for any two lotteries $L = \begin{pmatrix} o_1 & \cdots & o_m \\ p_1 & \cdots & p_m \end{pmatrix}$ and $L' = \begin{pmatrix} o_1 & \cdots & o_m \\ q_1 & \cdots & q_m \end{pmatrix}$,

$$L \succsim L' \text{ if and only if}$$
$$p_1 U(o_1) + ... + p_m U(o_m) \geq q_1 U(o_1) + ... + q_m U(o_m) \quad (\clubsuit)$$

Multiplying both sides of the above inequality by $a > 0$ and adding $(p_1 + ... + p_m)b$ to both sides we obtain

$$p_1 U(o_1) + ... + p_m U(o_m) \geq q_1 U(o_1) + ... + q_m U(o_m)$$
$$\text{if and only if}$$
$$p_1 [aU(o_1) + b] + ... + p_m [aU(o_m) + b] \geq q_1 [aU(o_1) + b] + ... + q_m [aU(o_m) + b] \quad (\blacklozenge)$$

Defining $V(o_i) = aU(o_i) + b$, it follows from (♣) and (♦) that

$$L \succsim L' \text{ if and only if}$$
$$p_1 V(o_1) + ... + p_m V(o_m) \geq q_1 V(o_1) + ... + q_m V(o_m)$$

That is, the function $V$ is a von Neumann-Morgenstern utility function that represents the ranking $\succsim$.

The proof of Part (B) will be given later, after introducing some more notation and some observations. ■

Since among the lotteries there are the degenerate ones that assign probability 1 to a single basic outcome, this implies that the DM has a complete and transitive ranking of the basic outcomes. We shall write $o_{best}$ for a best basic outcome, that is, a basic outcome which is at least as good as any other outcome ($o_{best} \succsim o$, for every $o \in O$) and $o_{worst}$ for a worst basic outcome, that is, a basic outcome such that every other outcome is at least as good as it ($o \succsim o_{worst}$, for every $o \in O$). Note that there may be several best outcomes (then the DM would be indifferent among them) and several worst outcomes. We shall assume throughout that the DM is not indifferent among all the outcomes, that is, we shall assume that $o_{best} \succ o_{worst}$.

We now show that, in virtue of Theorem 4.2, among the von Neumann-Morgenstern utility functions that represent a given von Neumann-Morgenstern ranking $\succsim$ of $\mathcal{L}(O)$, there is one that assigns the value 1 to the best basic outcome(s) and the value 0 to the worst basic outcome(s). To see this, consider an arbitrary von Neumann-Morgenstern utility function $F : O \to \mathbb{R}$ that represents $\succsim$ and define $G : O \to \mathbb{R}$ as follows: for every $o \in O$, $G(o) = F(o) - F(o_{worst})$. Then, by Theorem 4.5, $G$ is also a utility function that represents $\succsim$ and, by construction, $G(o_{worst}) = F(o_{worst}) - F(o_{worst}) = 0$; note also

that, since $o_{best} \succ o_{worst}$, it follows that $G(o_{best}) > 0$. Finally, define $U : O \to \mathbb{R}$ as follows: for every $o \in O$, $U(o) = \dfrac{G(o)}{G(o_{best})}$. Then, by Theorem 4.2, $U$ is a utility function that represents $\succsim$ and, by construction, $U(o_{worst}) = 0$ and $U(o_{best}) = 1$. For example, if there are six basic outcomes and the ranking of the basic outcomes is $o_3 \sim o_6 \succ o_1 \succ o_4 \succ o_2 \sim o_5$, then one can take as $o_{best}$ either $o_3$ or $o_6$ and as $o_{worst}$ either $o_2$ or $o_5$; furthermore, if $F$ is given by

| outcome | $o_1$ | $o_2$ | $o_3$ | $o_4$ | $o_5$ | $o_6$ |
|---|---|---|---|---|---|---|
| $F$ | 2 | −2 | 8 | 0 | −2 | 8 |

then $G$ is the function

| $o_1$ | $o_2$ | $o_3$ | $o_4$ | $o_5$ | $o_6$ |
|---|---|---|---|---|---|
| 4 | 0 | 10 | 2 | 0 | 10 |

and $U$ is the function

| $o_1$ | $o_2$ | $o_3$ | $o_4$ | $o_5$ | $o_6$ |
|---|---|---|---|---|---|
| 0.4 | 0 | 1 | 0.2 | 0 | 1 |

.

**Definition 4.3.** Let $U : O \to \mathbb{R}$ be a utility function that represents a given von Neumann-Morgenstern ranking $\succsim$ of the set of lotteries $\mathcal{L}(O)$. We say that $U$ is the *normalized utility function* if $U(o_{worst}) = 0$ and $U(o_{best}) = 1$.

The transformations described above show how to normalize any given utility function.

Armed with the notion of normalized utility function we can now complete the proof of Theorem 4.2.

**Proof of Part (B) of Theorem 4.2.** Let $F : O \to \mathbb{R}$ and $G : O \to \mathbb{R}$ be two von Neumann-Morgenstern utility functions that represent a given von Neumann-Morgenstern ranking of $\mathcal{L}(O)$. Let $U : O \to \mathbb{R}$ be the normalization of $F$ and $V : O \to \mathbb{R}$ be the normalization of $G$. First we show that it must be that $U = V$, that is, $U(o) = V(o)$ for every $o \in O$. Suppose, by contradiction, that there is an $\hat{o} \in O$ such that $U(\hat{o}) \neq V(\hat{o})$. Without loss of generality we can assume that $U(\hat{o}) > V(\hat{o})$. Construct the following lottery: $L = \begin{pmatrix} o_{best} & o_{worst} \\ \hat{p} & 1-\hat{p} \end{pmatrix}$ with $\hat{p} = U(\hat{o})$. Then $\mathbb{E}[U(L)] = \mathbb{E}[V(L)] = U(\hat{o})$. Hence, according to $U$ it must be that $\hat{o} \sim L$ (this follows from Theorem 4.1), while according to $V$ it must be (again, by Theorem 4.1) that $L \succ \hat{o}$ (since $\mathbb{E}[V(L)] = U(\hat{o}) > V(\hat{o})$). Then $U$ and $V$ cannot be two representations of the same ranking.

Now let $a_1 = \frac{1}{F(o_{best}) - F(o_{worst})}$ and $b_1 = -\frac{F(o_{worst})}{F(o_{best}) - F(o_{worst})}$. Note that $a_1 > 0$. Then it is easy to verify that, for every $o \in O$, $U(o) = a_1 F(o) + b_1$. Similarly let $a_2 = \frac{1}{G(o_{best}) - G(o_{worst})}$ and $b_2 = -\frac{G(o_{worst})}{G(o_{best}) - G(o_{worst})}$; again, $a_2 > 0$ and, for every $o \in O$, $V(o) = a_2 G(o) + b_2$. We can invert the latter transformation and obtain that, for every $o \in O$, $G(o) = \frac{V(o)}{a_2} - \frac{b_2}{a_2}$. Thus we can transform *F* into *U*, which – as proved above – is the same as *V*, and then transform *V* into *G* as follows: $G(o) = aF(o) + b$ where $a = \frac{a_1}{a_2} > 0$ and $b = \frac{b_1 - b_2}{a_2}$. ■

**Remark 4.2.** Theorem 4.2 is often stated as follows: a utility function that represents a von Neumann-Morgenstern ranking ≿ of $\mathcal{L}(O)$ is *unique up to a positive affine transformation*.[2] Because of Theorem 4.2, a von Neumann-Morgenstern utility function is usually referred to as a *cardinal* utility function.

Theorem 4.1 guarantees the existence of a utility function that represents a given von Neumann-Morgenstern ranking ≿ of $\mathcal{L}(O)$ and Theorem 4.2 characterizes the set of such functions. Can one actually construct a utility function that represents a given ranking? The answer is affirmative: if there are *m* basic outcomes one can construct an individual's von Neumann-Morgenstern utility function by asking her at most $(m-1)$ questions. The first question is "what is your ranking of the basic outcomes?". Then we can construct the normalized utility function by first assigning the value 1 to the best outcome(s) and the value 0 to the worst outcome(s). This leaves us with at most $(m-2)$ values to determine. For this we appeal to one of the axioms discussed in the next section, namely the *Continuity Axiom*, which says that, for every basic outcome $o_i$ there is a probability $p_i \in [0,1]$ such that the DM is indifferent between $o_i$ for certain and the lottery that gives a best outcome with probability $p_i$ and a worst outcome with probability $(1-p_i)$: $o_i \sim \begin{pmatrix} o_{best} & o_{worst} \\ p_i & 1-p_i \end{pmatrix}$. Thus, for each basic outcome $o_i$ for which a utility has not

[2] An affine transformation is a function $f : \mathbb{R} \to \mathbb{R}$ of the form $f(x) = ax + b$ with $a, b \in \mathbb{R}$. The affine transformation is positive if $a > 0$.

been determined yet, we should ask the individual to tell us the value of $p_i$ such that $o_i \sim \begin{pmatrix} o_{best} & o_{worst} \\ p_i & 1-p_i \end{pmatrix}$; then we can then set $U_i(o_i) = p_i$, because the expected utility of the lottery $\begin{pmatrix} o_{best} & o_{worst} \\ p_i & 1-p_i \end{pmatrix}$ is $p_i U_i(o_{best}) + (1-p_i) U_i(o_{worst}) = p_i(1) + (1-p_i)0 = p_i$.

**Example 4.1.** Suppose that there are five basic outcomes, that is, $O = \{o_1, o_2, o_3, o_4, o_5\}$ and the DM, who has von Neumann-Morgenstern preferences, tells us that her ranking of the basic outcomes is as follows:

$$o_2 \succ o_1 \sim o_5 \succ o_3 \sim o_4 .$$

Then we can begin by assigning utility 1 to the best outcome $o_2$ and utility 0 to the worst outcomes $o_3$ and $o_4$: $\begin{pmatrix} \text{outcome} & o_1 & o_2 & o_3 & o_4 & o_5 \\ \text{utility} & ? & 1 & 0 & 0 & ? \end{pmatrix}$. There is only one value left to be determined, namely the utility of $o_1$ (which is also the utility of $o_5$, since $o_1 \sim o_5$). To find this value, we ask the DM to tell us what value of $p$ makes her indifferent between the lottery $L = \begin{pmatrix} o_2 & o_3 \\ p & 1-p \end{pmatrix}$ and outcome $o_1$ with certainty. Suppose that her answer is: 0.4. Then her normalized von Neumann-Morgenstern utility function is $\begin{pmatrix} \text{outcome} & o_1 & o_2 & o_3 & o_4 & o_5 \\ \text{utility} & 0.4 & 1 & 0 & 0 & 0.4 \end{pmatrix}$. Knowing this, we can predict her choice among any set of lotteries over the five basic outcomes.

Test your understanding of the concepts introduced in this section, by going through the exercises in Section 4.E.2 of Appendix 4.E at the end of this chapter.

# 4.3 Expected utility: the axioms

We can now turn to the list of "rationality" axioms proposed by von Neumann and Morgenstern. Let $O = \{o_1, o_2, ..., o_m\}$ be the set of basic outcomes and $\mathcal{L}(O)$ the set of simple lotteries, that is, the set of probability distributions over *O*. Let $\succsim$ be a binary relation on $\mathcal{L}(O)$. We say that $\succsim$ is a von Neumann-Morgenstern ranking of $\mathcal{L}(O)$ if it satisfies the following four axioms or properties.

**Axiom 1 [Completeness and transitivity].** $\succsim$ is complete and transitive.

As noted in the previous section, Axiom 1 implies that there is a complete and transitive ranking of the basic outcomes. Recall that $o_{best}$ denotes a best basic outcome and $o_{worst}$ denotes a worst basic outcome and that we are assuming that $o_{best} \succ o_{worst}$, that is, that the DM is not indifferent among all the basic outcomes.

**Axiom 2 [Monotonicity].** $\begin{pmatrix} o_{best} & o_{worst} \\ p & 1-p \end{pmatrix} \succsim \begin{pmatrix} o_{best} & o_{worst} \\ q & 1-q \end{pmatrix}$ if and only if $p \geq q$.

**Axiom 3 [Continuity].** For every basic axioms $o_i$ there is a $p_i \in [0,1]$ such that $o_i \sim \begin{pmatrix} o_{best} & o_{worst} \\ p_i & 1-p_i \end{pmatrix}$.

Before we introduce the last axiom we need one more definition.

**Definition 4.1.** A *compound lottery* is a lottery of the form $\begin{pmatrix} x_1 & x_2 & \cdots & x_r \\ p_1 & p_2 & \cdots & p_r \end{pmatrix}$ where each $x_i$ is either an element of *O* or an element of $\mathcal{L}(O)$. Given a compound lottery $C = \begin{pmatrix} x_1 & x_2 & \cdots & x_r \\ p_1 & p_2 & \cdots & p_r \end{pmatrix}$ the *corresponding simple lottery* $L(C) = \begin{pmatrix} o_1 & o_2 & \cdots & o_m \\ q_1 & q_2 & \cdots & q_m \end{pmatrix}$ is defined as follows. First of all, for *i* = 1,...,*m* and *j* = 1,...,*r*, define

$$o_i(x_j) = \begin{cases} 1 & \text{if } x_j = o_i \\ 0 & \text{if } x_j = o_k \text{ with } k \neq i \\ s_i & \text{if } x_j = \begin{pmatrix} o_1 & \dots & o_{i-1} & o_i & o_{i+1} & \dots & o_m \\ s_1 & \dots & s_{i-1} & s_i & s_{i+12} & \dots & s_m \end{pmatrix} \end{cases}.$$

Then $q_i = \sum_{j=1}^{r} p_j \, o_i(x_j)$.

For example, let $m = 4$. Then $L = \begin{pmatrix} o_1 & o_2 & o_3 & o_4 \\ \frac{2}{5} & 0 & \frac{1}{5} & \frac{2}{5} \end{pmatrix}$ is a simple lottery (an element of $\mathcal{L}(O)$), while $C = \begin{pmatrix} \begin{pmatrix} o_1 & o_2 & o_3 & o_4 \\ \frac{1}{3} & \frac{1}{6} & \frac{1}{3} & \frac{1}{6} \end{pmatrix} & o_1 & \begin{pmatrix} o_1 & o_2 & o_3 & o_4 \\ \frac{1}{5} & 0 & \frac{1}{5} & \frac{3}{5} \end{pmatrix} \\ \frac{1}{2} & \frac{1}{4} & \frac{1}{4} \end{pmatrix}$ is a compound lottery (with $r = 3$). The compound lottery $C$ can be viewed graphically as a tree, as shown in Figure 4.1.

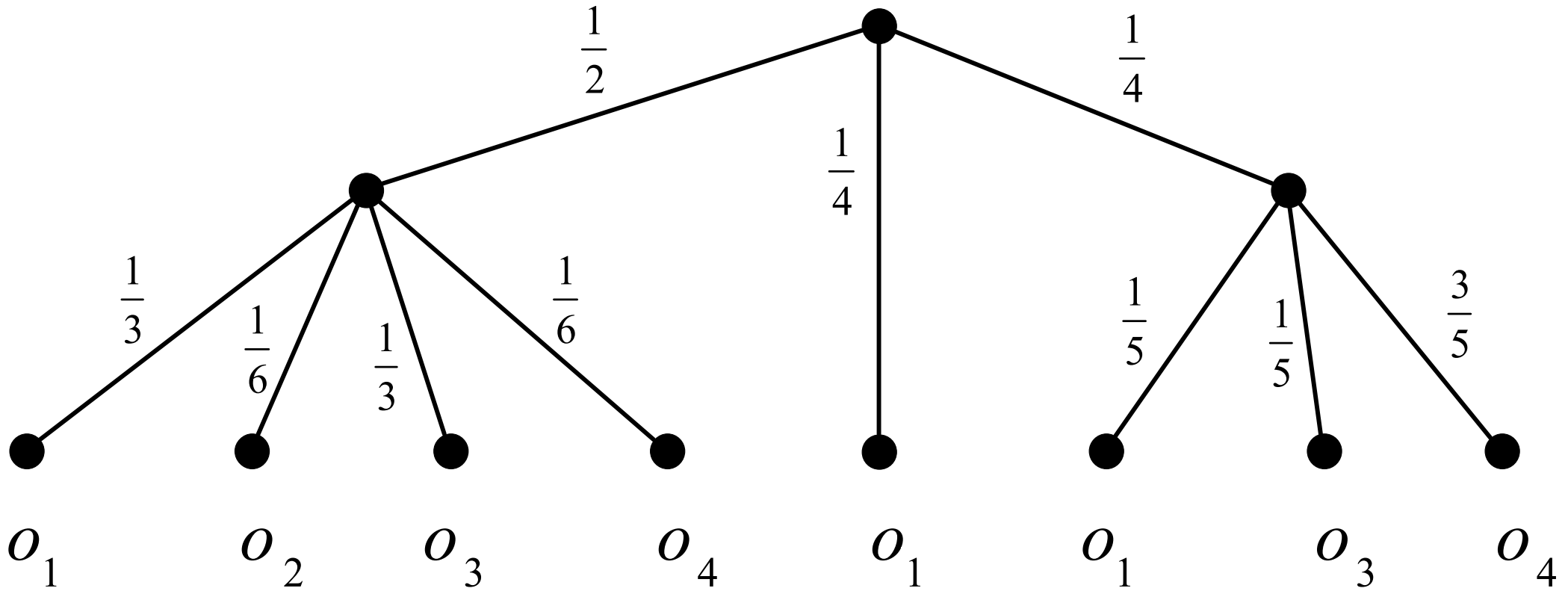


**Figure 4.1**
The compound lottery $C$.

Then $o_1(x_1) = \frac{1}{3}$, $o_1(x_2) = 1$ and $o_1(x_3) = \frac{1}{5}$ so that $q_1 = \frac{1}{2}\left(\frac{1}{3}\right) + \frac{1}{4}(1) + \frac{1}{4}\left(\frac{1}{5}\right) = \frac{28}{60}$. Similarly, $q_2 = \frac{1}{2}\left(\frac{1}{6}\right) + \frac{1}{4}(0) + \frac{1}{4}(0) = \frac{1}{12} = \frac{5}{60}$, $q_3 = \frac{1}{2}\left(\frac{1}{3}\right) + \frac{1}{4}(0) + \frac{1}{4}\left(\frac{1}{5}\right) = \frac{13}{60}$ and $q_4 = \frac{1}{2}\left(\frac{1}{6}\right) + \frac{1}{4}(0) + \frac{1}{4}\left(\frac{3}{5}\right) = \frac{14}{60}$. These numbers correspond to multiplying the probabilities along the edges of the tree of Figure 4.1 leading to an outcome, as shown in Figure 4.2 (a) and then adding up the probabilities of each outcome, as shown in Figure 4.2 (b). Thus the simple lottery $L(C)$ that corresponds to $C$ is $L(C) = \begin{pmatrix} o_1 & o_2 & o_3 & o_4 \\ \frac{28}{60} & \frac{5}{60} & \frac{13}{60} & \frac{14}{60} \end{pmatrix}$, namely the lottery shown in Figure 4.2 (b).

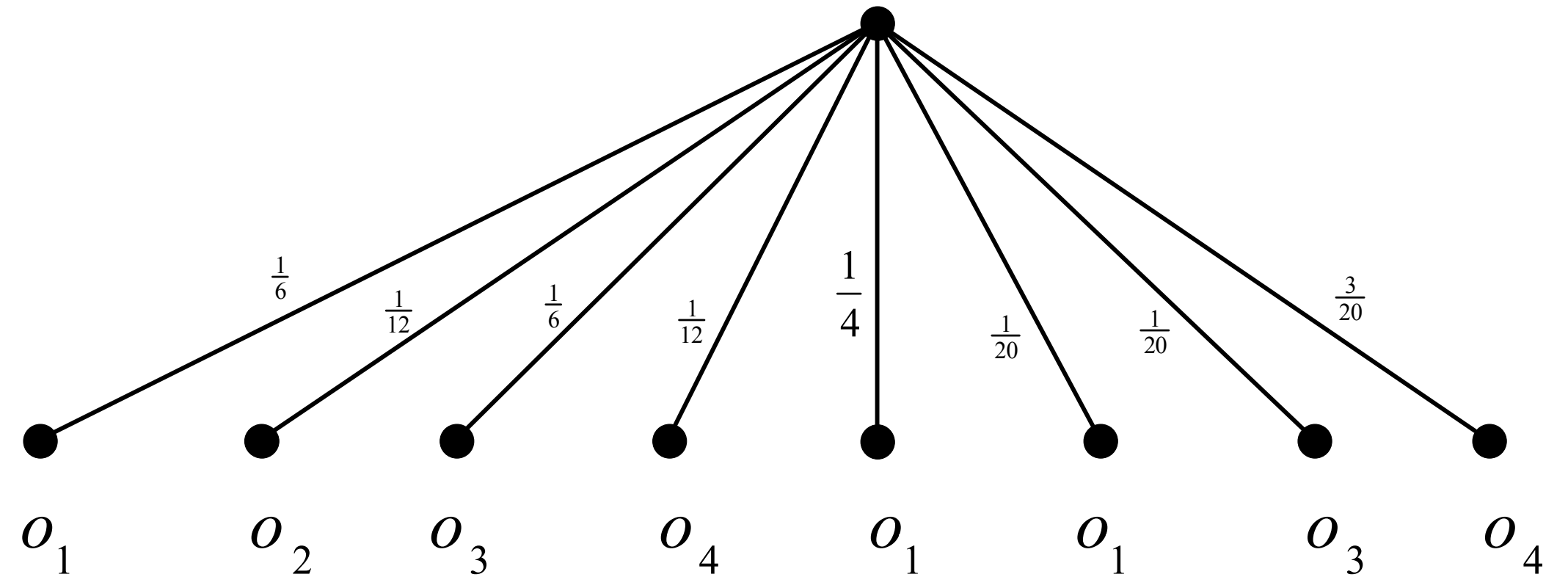


(a) Simplification of Figure 4.1 obtained by condensing paths into simple edges and associating with the simple edges the products of the probabilities along the path.

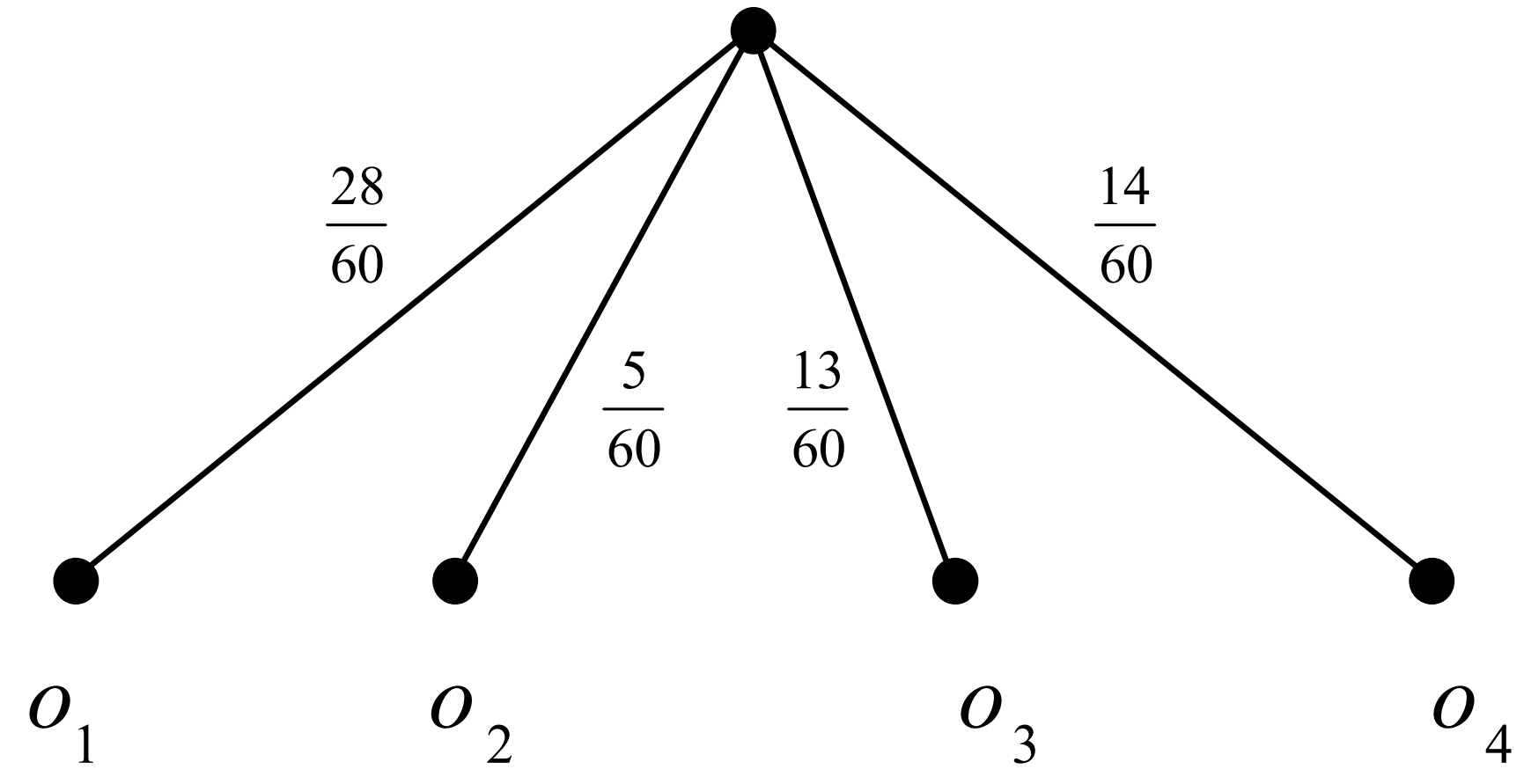


(b) Simplification of Part (a) obtained by adding, for each outcome, the probabilities of that outcome.

**Figure 4.2**

**Axiom 4 [Independence or substitutability].** Consider an arbitrary basic outcome $o_i$ and an arbitrary simple lottery $L = \begin{pmatrix} o_1 & \cdots & o_{i-1} & o_i & o_{i+1} & \cdots & o_m \\ p_1 & \cdots & p_{i-1} & p_i & p_{i+1} & \cdots & p_m \end{pmatrix}$. If $\hat{L}$ is a simple lottery such that $o_i \sim \hat{L}$ then $L \sim M$ where $M$ is the simple lottery corresponding to the compound lottery $C = \begin{pmatrix} o_1 & \cdots & o_{i-1} & \hat{L} & o_{i+1} & \cdots & o_m \\ p_1 & \cdots & p_{i-1} & p_i & p_{i+1} & \cdots & p_m \end{pmatrix}$ (Definition 4.9) obtained by replacing $o_i$ with $\hat{L}$ in $L$.

We can now prove the first theorem of the previous section.

**Proof of Theorem 4.1.** To simplify the notation, throughout this proof we will assume that we have renumbered the basic outcomes in such a way that

$$o_{best} = o_1 \text{ and } o_{worst} = o_m .$$

First of all, for every basic outcome $o_i$, let $u_i \in [0,1]$ be such that $o_i \sim \begin{pmatrix} o_1 & o_m \\ u_i & 1-u_i \end{pmatrix}$.

The existence of such values $u_i$'s is guaranteed by the Continuity Axiom (Axiom 3); clearly $u_1 = 1$ and $u_m = 0$. Now consider an arbitrary lottery $L_1 = \begin{pmatrix} o_1 & \ldots & o_m \\ p_1 & \ldots & p_m \end{pmatrix}$. First we show that

$$L_1 \sim \begin{pmatrix} o_1 & o_m \\ \sum_{i=1}^{m} p_i u_i & 1 - \sum_{i=1}^{m} p_i u_i \end{pmatrix} \qquad (\dagger)$$

This is done through a repeated application of the Independence Axiom (Axiom 4), as follows: consider the compound lottery $C_2 = \begin{pmatrix} o_1 & \begin{pmatrix} o_1 & o_m \\ u_2 & 1-u_2 \end{pmatrix} & o_3 & \ldots & o_m \\ p_1 & p_2 & p_3 & \ldots & p_m \end{pmatrix}$ obtained by replacing $o_2$ in lottery $L_1$ with the lottery $\begin{pmatrix} o_1 & o_m \\ u_2 & 1-u_2 \end{pmatrix}$ that the DM considers to be just as good as $o_2$. The simple lottery corresponding to $C_2$ is $L_2 = \begin{pmatrix} o_1 & o_3 & \ldots & o_{m-1} & o_m \\ p_1 + p_2 u_2 & p_3 & \ldots & p_{m-1} & p_m + p_2(1-u_2) \end{pmatrix}$. Note that $o_2$ is assigned probability 0 in $L_2$ and thus we have omitted it. By Axiom 4, $L_1 \sim L_2$. Now apply the same argument to $L_2$: let $C_3 = \begin{pmatrix} o_1 & \begin{pmatrix} o_1 & o_m \\ u_3 & 1-u_3 \end{pmatrix} & \ldots & o_m \\ p_1 + p_2 u_2 & p_3 & \ldots & p_m + p_2(1-u_2) \end{pmatrix}$ whose corresponding simple lottery is

$$L_3 = \begin{pmatrix} o_1 & \ldots & o_m \\ p_1 + p_2 u_2 + p_3 u_3 & \ldots & p_m + p_2(1-u_2) + p_3(1-u_3) \end{pmatrix}.$$

Note, again, that $o_3$ is assigned probability zero in $L_3$. By Axiom 4, $L_2 \sim L_3$; thus, by transitivity (since $L_1 \sim L_2$ and $L_2 \sim L_3$) we have that $L_1 \sim L_3$. Repeating this argument we get that $L_1 \sim L_{m-1}$, where

$$L_{m-1} = \begin{pmatrix} o_1 & o_m \\ p_1 + p_2 u_2 + \ldots + p_{m-1} u_{m-1} & p_m + p_2(1-u_2) + \ldots + p_{m-1}(1-u_{m-1}) \end{pmatrix}.$$

Since $u_1 = 1$ and $u_m = 0$,

$$p_1 + p_2 u_2 + \ldots + p_{m-1} u_{m-1} = p_1 u_1 + p_2 u_2 + \ldots + p_{m-1} u_{m-1} + p_m u_m = \sum_{i=1}^{m} p_i u_i$$

and

$$p_2(1-u_2) + \ldots + p_{m-1}(1-u_{m-1}) + p_m = \sum_{i=2}^{m} p_i - \sum_{i=2}^{m-1} p_i u_i = p_1 + \sum_{i=2}^{m} p_i - \sum_{i=2}^{m-1} p_i u_i - p_1$$

$$\underset{\substack{\text{since } u_1 = 1 \\ \text{and } u_m = 0}}{=} \sum_{i=1}^{m} p_i - \sum_{i=2}^{m-1} p_i u_i - p_1 u_1 - p_m u_m \underset{\text{since } \sum_{i=1}^{m} p_i = 1}{=} 1 - \sum_{i=1}^{m} p_i u_i$$

Thus $L_{m-1} = \begin{pmatrix} o_1 & o_m \\ \sum_{i=1}^{m} p_i u_i & 1 - \sum_{i=1}^{m} p_i u_i \end{pmatrix}$, which proves (†).

Now define the following utility function $U: \{o_1, \ldots, o_m\} \to [0,1]$: $U(o_i) = u_i$, where, as before, for every basic outcome $o_i$, $u_i \in [0,1]$ is such that $o_i \sim \begin{pmatrix} o_1 & o_m \\ u_i & 1-u_i \end{pmatrix}$. Consider two arbitrary lotteries $L = \begin{pmatrix} o_1 & \ldots & o_m \\ p_1 & \ldots & p_m \end{pmatrix}$ and $L' = \begin{pmatrix} o_1 & \ldots & o_m \\ q_1 & \ldots & q_m \end{pmatrix}$. We want to show that $L \succsim L'$ if and only if $\mathbb{E}[U(L)] \geq \mathbb{E}[U(L')]$, that is, if and only if $\sum_{i=1}^{m} p_i u_i \geq \sum_{i=1}^{m} q_i u_i$. By (†), $L \sim M$, where $M = \begin{pmatrix} o_1 & o_m \\ \sum_{i=1}^{m} p_i u_i & 1 - \sum_{i=1}^{m} p_i u_i \end{pmatrix}$ and also $L' \sim M'$, where $M' = \begin{pmatrix} o_1 & o_m \\ \sum_{i=1}^{m} q_i u_i & 1 - \sum_{i=1}^{m} q_i u_i \end{pmatrix}$. Thus, by transitivity of $\succsim$, $L \succsim L'$ if and only if

$M \succsim M'$; by the Monotonicity Axiom (Axiom 2), $M \succsim M'$ if and only if $\sum_{i=1}^{m} p_i u_i \geq \sum_{i=1}^{m} q_i u_i$. ■

The following example, known as the *Allais paradox*,[3] suggests that one should view expected utility theory as a *prescriptive* or *normative* theory (that is, as a theory about how rational people should choose) rather than as a *descriptive* theory (that is, as a theory about the actual behavior of individuals). In 1953 the French economist Maurice Allais published a paper regarding a survey he had conducted in 1952 concerning a hypothetical decision problem. Subjects "with good training in and knowledge of the theory of probability, so that they could be considered to behave rationally" were asked to rank the following pairs of lotteries:

$$A = \begin{pmatrix} \$5 \text{ million} & \$0 \\ \frac{89}{100} & \frac{11}{100} \end{pmatrix} \quad \textit{versus} \quad B = \begin{pmatrix} \$1 \text{ million} & \$0 \\ \frac{90}{100} & \frac{10}{100} \end{pmatrix}$$

and

$$C = \begin{pmatrix} \$5 \text{ million} & \$1 \text{ million} & \$0 \\ \frac{89}{100} & \frac{10}{100} & \frac{1}{100} \end{pmatrix} \quad \textit{versus} \quad D = \begin{pmatrix} \$1 \text{ million} \\ 1 \end{pmatrix}.$$

Most subjects reported the following ranking:

$$A \succ B \text{ and } D \succ C.$$

Such ranking violates the axioms of expected utility. To see this, let $O = \{o_1, o_2, o_3\}$ with $o_1 = \$5$ million, $o_2 = \$1$ million and $o_3 = \$0$. Let us assume that the individual in question prefers more money to less: $o_1 \succ o_2 \succ o_3$ and has a von Neumann-Morgenstern ranking of the lotteries over $O$. Let $u_2 \in (0,1)$ is such that $D \sim \begin{pmatrix} \$5 \text{ million} & \$0 \\ u_2 & 1-u_2 \end{pmatrix}$ (the existence of such $u_2$ is guaranteed by the Continuity Axiom). Then, since $D \succ C$, by transitivity

$$\begin{pmatrix} \$5 \text{ million} & \$0 \\ u_2 & 1-u_2 \end{pmatrix} \succ C \qquad (\clubsuit)$$

[3] See http://en.wikipedia.org/wiki/Allais_paradox .

Let $C'$ be the simple lottery corresponding to the compound lottery $\begin{pmatrix} \$5 \text{ million} & \begin{pmatrix} \$5 \text{ million} & \$0 \\ u_2 & 1-u_2 \end{pmatrix} & \$0 \\ \frac{89}{100} & \frac{10}{100} & \frac{1}{100} \end{pmatrix}$. Then $C' = \begin{pmatrix} \$5 \text{ million} & \$0 \\ \frac{89}{100}+\frac{10}{100}u_2 & 1-\left(\frac{89}{100}+\frac{10}{100}u_2\right) \end{pmatrix}$.

By the Independence Axiom, $C \sim C'$ and thus, by (♣) and transitivity,

$$\begin{pmatrix} \$5 \text{ million} & \$0 \\ u_2 & 1-u_2 \end{pmatrix} \succ \begin{pmatrix} \$5 \text{ million} & \$0 \\ \frac{89}{100}+\frac{10}{100}u_2 & 1-\left(\frac{89}{100}+\frac{10}{100}u_2\right) \end{pmatrix}$$

Hence, by the Monotonicity Axiom, $u_2 > \frac{89}{100}+\frac{10}{100}u_2$, that is,

$$u_2 > \frac{89}{90} \qquad (\blacklozenge)$$

Let $B'$ be the simple lottery corresponding to the compound lottery $\begin{pmatrix} \begin{pmatrix} \$5 \text{ million} & \$0 \\ u_2 & 1-u_2 \end{pmatrix} & \$0 \\ \frac{90}{100} & \frac{10}{100} \end{pmatrix}$. Then $B' = \begin{pmatrix} \$5 \text{ million} & \$0 \\ \frac{90}{100}u_2 & 1-\frac{90}{100}u_2 \end{pmatrix}$. By the Independence Axiom, $B \sim B'$; thus, since $A \succ B$, by transitivity $A \succ B'$ and therefore, by the Monotonicity Axiom, $\frac{89}{100} > \frac{90}{100}u_2$, that is, $u_2 < \frac{89}{90}$, contradicting (♦). Thus, if one finds the expected utility axioms compelling as axioms of rationality, then one cannot consistently express a preference for $A$ over $B$ and also a preference for $D$ over $C$.

Another well-known paradox is the *Ellsberg paradox*. Suppose you are told that an urn contains 30 red balls and 60 other balls that are either blue or yellow. You don't know how many blue or how many yellow balls there are, but the number of blue balls plus the number of yellow ball equals 60 (they could be all blue or all yellow or any combination of the two). The balls are well mixed so that each individual ball is as likely to be drawn as any other. You are given a choice between the bets $A$ and $B$, where

$A$ = you get \$100 if you pick a **red** ball and nothing otherwise,

$B$ = you get \$100 if you pick a **blue** ball and nothing otherwise.

Many subjects in experiments state a strict preference for $A$ over $B$: $A \succ B$.
Consider now the following bets:

C = you get \$100 if you pick a **red or yellow** ball and nothing otherwise,

D = you get \$100 if you pick a **blue or yellow** ball and nothing otherwise.

Do the axioms of expected utility constrain your ranking of *C* and *D*? Many subjects in experiments state the following ranking: $A \succ B$ and $D \succsim C$. All such people violate the axioms of expected utility. The fraction of red balls in the urn is $\frac{30}{90} = \frac{1}{3}$ Let $p_2$ be the fraction of blue balls and $p_3$ the fraction of yellow balls (either of these can be zero: all we know is that $p_2 + p_3 = \frac{60}{90} = \frac{2}{3}$ ). Then *A*, *B*, *C* and *D* can be viewed as the following lotteries:

$$A = \begin{pmatrix} \$100 & \$0 \\ \frac{1}{3} & p_2 + p_3 \end{pmatrix} \quad B = \begin{pmatrix} \$100 & \$0 \\ p_2 & \frac{1}{3} + p_3 \end{pmatrix} \quad C = \begin{pmatrix} \$100 & \$0 \\ \frac{1}{3} + p_3 & p_2 \end{pmatrix}$$

$$D = \begin{pmatrix} \$100 & \$0 \\ p_2 + p_3 = \frac{2}{3} & \frac{1}{3} \end{pmatrix}$$

Let $U$ be the normalized von Neumann-Morgenstern utility function that represents the individual's ranking; then $U(\$100) = 1$ and $U(0) = 0$. Thus $\mathbb{E}[U(A)] = \frac{1}{3}$, $\mathbb{E}[U(B)] = p_2$, $\mathbb{E}[U(C)] = \frac{1}{3} + p_3$ and $\mathbb{E}[U(D)] = p_2 + p_3 = \frac{2}{3}$. Hence, $A \succ B$ if and only if $\frac{1}{3} > p_2$, which implies that $p_3 > \frac{1}{3}$, so that $\mathbb{E}[U(C)] = \frac{1}{3} + p_3 > \mathbb{E}[U(D)] = \frac{2}{3}$ and thus $C \succ D$ (similarly, $B \succ A$ if and only if $\frac{1}{3} < p_2$, which implies that $\mathbb{E}[U(C)] < \mathbb{E}[U(D)]$ and thus $D \succ C$).

Test your understanding of the concepts introduced in this section, by going through the exercises in Section 4.E.3 of Appendix 4.E at the end of this chapter.

# Appendix 4.E: Exercises

## 4.E.1. Exercises for Section 4.1: Money lotteries and attitudes to risk.

The answers to the following exercises are in Appendix S at the end of this chapter.

**Exercise 4.1**. What is the expected value of the lottery: $\begin{pmatrix} 24 & 12 & 48 & 6 \\ \frac{1}{6} & \frac{2}{6} & \frac{1}{6} & \frac{2}{6} \end{pmatrix}$?

**Exercise 4.2**. Consider the following lottery: $\begin{pmatrix} o_1 & o_2 & o_3 \\ \frac{1}{4} & \frac{1}{2} & \frac{1}{4} \end{pmatrix}$, where

$o_1$ = you get an invitation to have dinner at the White House
$o_2$ = you get (for free) a puppy of your choice
$o_3$ = you get $600.

What is the expected value of this lottery?

**Exercise 4.3**. Consider the following money lottery

$$L = \begin{pmatrix} \$10 & \$15 & \$18 & \$20 & \$25 & \$30 & \$36 \\ \frac{3}{12} & \frac{1}{12} & 0 & \frac{3}{12} & \frac{2}{12} & 0 & \frac{3}{12} \end{pmatrix}.$$

**(a)** What is the expected value of the lottery?

**(b)** Ann says that, between getting $20 for certain and playing the above lottery, she would prefer the $20 for certain. What is her attitude to risk?

**(c)** Bob prefers more money to less and has transitive preferences. He says that, given the same choice as Ann, he would prefer playing the lottery. What is his attitude to risk?

**Exercise 4.4**. Sam has a debilitating illness and has been offered two mutually exclusive courses of action: (1) take some well-known drugs which have been tested for a long time and (2) take a new experimental drug. If he chooses (1) then for certain his pain will be reduced to a bearable level. If he chooses (2) then he has a 50% chance of being completely cured and a 50% chance of no benefits from the drug and possibly some harmful side effects. He chose (1). What is his attitude to risk?

## 4.E.2. Exercises for Section 4.2: Expected utility theory.

**Exercise 4.5**. Ben is offered a choice between the following two money lotteries: $A = \begin{pmatrix} \$4,000 & \$0 \\ 0.8 & 0.2 \end{pmatrix}$ and $B = \begin{pmatrix} \$3,000 \\ 1 \end{pmatrix}$. He says he strictly prefers $B$ to $A$. Which of the following two lotteries, $C$ and $D$, will Ben choose if he satisfies the axioms of expected utility and prefers more money to less?

$$C = \begin{pmatrix} \$4,000 & \$0 \\ 0.2 & 0.8 \end{pmatrix} \qquad D = \begin{pmatrix} \$3,000 & \$0 \\ 0.25 & 0.75 \end{pmatrix}$$

**Exercise 4.6**. There are three basic outcomes, $o_1, o_2$ and $o_3$. Ann satisfies the axioms of expected utility and her preferences over lotteries involving these three outcomes can be represented by the following von Neumann-Morgenstern utility function: $V(o_2) = a > V(o_1) = b > V(o_3) = c$. Normalize the utility function.

**Exercise 4.7**. Consider the following lotteries:

$$L_1 = \begin{pmatrix} \$3000 & \$2000 & \$1000 & \$500 \\ \frac{5}{6} & 0 & 0 & \frac{1}{6} \end{pmatrix}, \quad L_2 = \begin{pmatrix} \$3000 & \$2000 & \$1000 & \$500 \\ \frac{2}{3} & 0 & 0 & \frac{1}{3} \end{pmatrix}$$

$$L_3 = \begin{pmatrix} \$3000 & \$2000 & \$1000 & \$500 \\ \frac{1}{4} & \frac{1}{4} & \frac{1}{4} & \frac{1}{4} \end{pmatrix}, \quad L_4 = \begin{pmatrix} \$3000 & \$2000 & \$1000 & \$500 \\ 0 & \frac{1}{2} & \frac{1}{2} & 0 \end{pmatrix}$$

Jennifer says that she is indifferent between lottery $L_1$ and getting \$2,000 for certain. She is also indifferent between lottery $L_2$ and getting \$1,000 for certain. Finally, she says that between $L_3$ and $L_4$ she would choose $L_3$. Is she rational according to the theory of expected utility? [Assume that she prefers more money to less.]

**Exercise 4.8**. Consider the following basic outcomes: $o_1$ = a Summer internship at the White House, $o_2$ = a free 1-week vacation in Europe, $o_3$ = \$800, $o_4$ = a free ticket to a concert. Rachel says that her ranking of these outcomes is $o_1 \succ o_2 \succ o_3 \succ o_4$. She also says that (1) she is indifferent between $\begin{pmatrix} o_2 \\ 1 \end{pmatrix}$ and $\begin{pmatrix} o_1 & o_4 \\ \frac{4}{5} & \frac{1}{5} \end{pmatrix}$ and (2) she is indifferent between $\begin{pmatrix} o_3 \\ 1 \end{pmatrix}$ and $\begin{pmatrix} o_1 & o_4 \\ \frac{1}{2} & \frac{1}{2} \end{pmatrix}$. If she satisfies the axioms of expected utility theory, which of the two lotteries $L_1 = \begin{pmatrix} o_1 & o_2 & o_3 & o_4 \\ \frac{1}{8} & \frac{2}{8} & \frac{3}{8} & \frac{2}{8} \end{pmatrix}$ and $L_2 = \begin{pmatrix} o_1 & o_2 & o_3 \\ \frac{1}{5} & \frac{3}{5} & \frac{1}{5} \end{pmatrix}$ will she choose?

**Exercise 4.9**. Consider the following lotteries:

$$L_1 = \begin{pmatrix} \$30 & \$28 & \$24 & \$18 & \$8 \\ \frac{2}{10} & \frac{1}{10} & \frac{1}{10} & \frac{2}{10} & \frac{4}{10} \end{pmatrix} \text{ and } L_2 = \begin{pmatrix} \$30 & \$28 & \$8 \\ \frac{1}{10} & \frac{4}{10} & \frac{5}{10} \end{pmatrix}.$$

**(a)** Which lottery would a risk neutral person choose?

**(b)** Paul's von Neumann-Morgenstern utility-of-money function is $U(m) = \ln(m)$, where ln denotes the natural logarithm. Which lottery would Paul choose?

**Exercise 4.10**. There are five basic outcomes. Jane has a von Neumann-Morgenstern ranking of the set of lotteries over the set of basic outcomes that can be represented by either of the utility functions $U$ or $V$ given below:

$$\begin{pmatrix} & o_1 & o_2 & o_3 & o_4 & o_5 \\ U: & 44 & 170 & -10 & 26 & 98 \\ V: & 32 & 95 & 5 & 23 & 59 \end{pmatrix}.$$

**(a)** Show how to normalize each of $U$ and $V$ and verify that you get the same normalized utility function.

**(b)** Show how to transform $U$ into $V$ with a positive affine transformation of the form $x \mapsto ax + b$ with $a, b \in \mathbb{R}$ and $a > 0$.

**Exercise 4.11**. Consider the following lotteries: $L_3 = \begin{pmatrix} \$28 \\ 1 \end{pmatrix}$, $L_4 = \begin{pmatrix} \$10 & \$50 \\ \frac{1}{2} & \frac{1}{2} \end{pmatrix}$.

**(a)** Ann has the following von Neumann-Morgenstern utility function: $U_{Ann}(m) = \sqrt{m}$. How does she rank the two lotteries?

**(b)** Bob has the following von Neumann-Morgenstern utility function: $U_{Bob}(m) = 2m - \frac{m^4}{100^3}$. How does he rank the two lotteries?

**(c)** Verify that both Ann and Bob are risk averse, by determining what they would choose between lottery $L_4$ and its expected value for certain.

## 4.E.2. Exercises for Section 4.3: Expected utility axioms.

The answers to the following exercises are in Appendix S at the end of this chapter.

**Exercise 4.12**. Let $O = \{o_1, o_2, o_3, o_4\}$. Find the simple lottery corresponding to the following compound lottery

$$\begin{pmatrix} \begin{pmatrix} o_1 & o_2 & o_3 & o_4 \\ \frac{2}{5} & \frac{1}{10} & \frac{3}{10} & \frac{1}{5} \end{pmatrix} & o_2 & \begin{pmatrix} o_1 & o_3 & o_4 \\ \frac{1}{5} & \frac{1}{5} & \frac{3}{5} \end{pmatrix} & \begin{pmatrix} o_2 & o_3 \\ \frac{1}{3} & \frac{2}{3} \end{pmatrix} \\ \frac{1}{8} & \frac{1}{4} & \frac{1}{8} & \frac{1}{2} \end{pmatrix}$$

**Exercise 4.13**. Let $O = \{o_1, o_2, o_3, o_4\}$. Suppose that the DM has a von Neumann-Morgenstern ranking of $\mathcal{L}(O)$ and states the following indifference: $o_1 \sim \begin{pmatrix} o_2 & o_4 \\ \frac{1}{4} & \frac{3}{4} \end{pmatrix}$ and $o_2 \sim \begin{pmatrix} o_3 & o_4 \\ \frac{3}{5} & \frac{2}{5} \end{pmatrix}$. Find a lottery that the DM considers just as good as $L = \begin{pmatrix} o_1 & o_2 & o_3 & o_4 \\ \frac{1}{3} & \frac{2}{9} & \frac{1}{9} & \frac{1}{3} \end{pmatrix}$.

◊◊◊◊◊◊◊◊◊◊◊◊◊

**Exercise 4.14: Challenging Question.**

Would you be willing to pay more in order to reduce the probability of dying within the next hour from one sixth to zero or from four sixths to three sixths? Unfortunately, this is not a hypothetical question: you accidentally entered the office of an Italian mad scientist and have been overpowered and tied to a chair. The mad scientist has put six glasses in front of you, numbered 1 to 6, and tells you that one of them contains a deadly poison and the other five a harmless liquid. He says that he is going to roll a die and make you drink from the glass whose number matches the number that shows up from the rolling of the die. You beg to be exempted and he asks you "what is the largest amount of money that you would be willing to pay to replace the glass containing the poison with one containing a harmless liquid?". Interpret this question as "what sum of money $x$ makes you indifferent between (1) leaving the poison in whichever glass contains it and rolling the die, and (2) reducing your wealth by $\$x$ and rolling the die after the poison has been replaced by a harmless liquid". Your answer is: $\$X$. Then he asks you "suppose that instead of one glass with poison there had been four glasses with poison (and two with a harmless liquid); what is the largest amount of money that you would be willing to pay to replace one glass with poison with a glass containing a harmless liquid (and thus roll the die with 3 glasses with poison and 3 with a harmless liquid)?". Your answer is: $\$Y$. Show that if $X > Y$ then you do not satisfy the axioms of Expected Utility Theory. [Hint: think about what the basic outcomes are; assume that you do not care about how much money is left in your estate if you die and that, when alive, you prefer more money to less.]

# Appendix 4.S: Solutions to exercises

**Exercise 4.1.** The expected value of the lottery $\begin{pmatrix} 24 & 12 & 48 & 6 \\ \frac{1}{6} & \frac{2}{6} & \frac{1}{6} & \frac{2}{6} \end{pmatrix}$ is $\frac{1}{6}(24)+\frac{2}{6}(12)+\frac{1}{6}(48)+\frac{2}{6}(6)=18$

**Exercise 4.2.** This was a trick question! There is no expected value because the basic outcomes are not numbers.

**Exercise 4.3.** **(a)** The expected value is

$$\mathbb{E}[L]=\tfrac{3}{12}(10)+\tfrac{1}{12}(15)+(0)(18)+\tfrac{3}{12}(20)+\tfrac{2}{12}(25)+(0)(30)+\tfrac{3}{12}(36)=\frac{263}{12}=\$21.917$$

**(b)** Presumably Ann prefers more money to less, so that she prefers \$21.917 to \$20. She said that she prefers \$20 to the lottery *L*. Thus, assuming that her preferences are transitive, she prefers \$21.917 to *L*. Hence she is risk averse.

**(c)** The answer is: we cannot tell. First of all, since Bob prefers more money to less, he prefers \$21.917 to \$20. Bob could be risk neutral, because a risk neutral person would be indifferent between *L* and \$21.917; since Bob prefers \$21.917 to \$20 and has transitive preferences, if risk neutral he would prefer *L* to \$20 for certain. However, Bob could also be risk loving: a risk-loving person prefers *L* to \$21.917 and we know that he prefers \$21.917 to \$20; thus, by transitivity, if risk loving, he would prefer *L* to \$20 for certain. But Bob could also be risk averse: he could consistently prefer \$21.917 to *L* and *L* to \$20 (for example, he could consider *L* to be just as good as \$20.50).

**Exercise 4.4.** Just like Exercise 4.2, this was a trick question! Here the basic outcomes are not sums of money but states of health. Since the described choice is not one between money lotteries, the definitions of risk aversion/neutrality/love are not applicable.

**Exercise 4.5.** Since Ben prefers *B* to *A*, he must prefer *D* to *C*.

*Proof:* Let $U$ be a von Neumann-Morgenstern utility function that represents Ben's preferences. Let U(\$4,000) = *a*, U(\$3,000) = *b* and U(\$0) = *c*. Since Ben prefers more money to less, *a* > *b* > *c*. Now, $EU(A)=0.8U(\$4,000)+0.2U(\$0) \ =0.8a+0.2c$ and $EU(B)=U(\$3,000)=b$. Since Ben prefers *B* to *A*, it must be that $b>0.8a+0.2c$. Let us now compare

$C$ and $D$: $EU(C) = 0.2a + 0.8c$ and $EU(D) = 0.25b + 0.75c$. Since $b > 0.8a + 0.2c$, $0.25b > 0.25(0.8a + 0.2c) = 0.2a + 0.05c$ and thus, adding $0.75c$ to both sides, we get that $0.25b + 0.75c > 0.2a + 0.8c$, that is, $EU(D) > EU(C)$, so that $D \succ C$.

**Exercise 4.6.** Define the function $U$ as follows:

$$U(x) = \frac{1}{a-c}V(x) - \frac{c}{a-c} = \frac{V(x)-c}{a-c}.$$

Then $U$ represents the same preferences as $V$ (note that, by hypothesis, $a > c$ and thus $\frac{1}{a-c} > 0$). Then

$$U(o_2) = \frac{V(o_2)-c}{a-c} = \frac{a-c}{a-c} = 1, \quad U(o_1) = \frac{V(o_1)-c}{a-c} = \frac{b-c}{a-c}, \text{ and}$$

$$U(o_3) = \frac{V(o_3)-c}{a-c} = \frac{c-c}{a-c} = 0.$$

Note that, since $a > b > c$, $0 < \frac{b-c}{a-c} < 1$.

**Exercise 4.7.** Suppose that there is a von Neumann-Morgenstern utility function $U$ that represents Jennifer's preferences. We can normalize it so that $U(\$3000) = 1$ and $U(\$500) = 0$. Since Jennifer is indifferent between $L_1$ and \$2000, $U(\$2000) = \frac{5}{6}$. Since she is indifferent between $L_2$ and \$1000, $U(\$1000) = \frac{2}{3}$. Thus $EU(L_3) = \frac{1}{4}(1) + \frac{1}{4}\left(\frac{5}{6}\right) + \frac{1}{4}\left(\frac{2}{3}\right) + \frac{1}{4}(0) = \frac{5}{8}$ and $EU(L_4) = 0(1) + \frac{1}{2}\left(\frac{5}{6}\right) + \frac{1}{2}\left(\frac{2}{3}\right) + 0(0) = \frac{3}{4}$. Since $\frac{3}{4} > \frac{5}{8}$, Jennifer should prefer $L_4$ over $L_3$. Hence she is **not** rational according to the theory of expected utility.

**Exercise 4.8.** Normalize her utility function so that $U(z_1) = 1$ and $U(z_4) = 0$. Then, since Rachel is indifferent between $\begin{pmatrix} o_2 \\ 1 \end{pmatrix}$ and $\begin{pmatrix} o_1 & o_4 \\ \frac{4}{5} & \frac{1}{5} \end{pmatrix}$, we have that $U(o_2) = \frac{4}{5}$. Similarly, since she is indifferent between $\begin{pmatrix} o_3 \\ 1 \end{pmatrix}$ and $\begin{pmatrix} o_1 & o_4 \\ \frac{1}{2} & \frac{1}{2} \end{pmatrix}$,

$U(o_3) = \frac{1}{2}$. Then the expected utility of $L_1 = \begin{pmatrix} o_1 & o_2 & o_3 & o_4 \\ \frac{1}{8} & \frac{2}{8} & \frac{3}{8} & \frac{2}{8} \end{pmatrix}$ is $\frac{1}{8}(1)+\frac{2}{8}(\frac{4}{5})+\frac{3}{8}(\frac{1}{2})+\frac{2}{8}(0) = \frac{41}{80} =$ 0.5125, while the expected utility of $L_2 = \begin{pmatrix} o_1 & o_2 & o_3 \\ \frac{1}{5} & \frac{3}{5} & \frac{1}{5} \end{pmatrix}$ is $\frac{1}{5}(1)+\frac{3}{5}(\frac{4}{5})+\frac{1}{5}(\frac{1}{2}) = .\ \frac{39}{50} = 0.78$. Hence she prefers $L_2$ to $L_1$.

**Exercise 4.9.** **(a)** The expected value of $L_2$ is $\frac{1}{10}(30)+\frac{4}{10}(28)+\frac{5}{10}8 = 18.2$. The expected value of $L_1$ is $\frac{2}{10}(30)+\frac{1}{10}(28)+\frac{1}{10}(24)+\frac{2}{10}(18)+\frac{4}{10}(8) = 18$. Hence a risk-neutral person would prefer $L_2$ to $L_1$.

**(b)** The expected utility of $L_2$ is $\frac{1}{10}\cdot\ln(30) + \frac{2}{5}\cdot\ln(28) + \frac{1}{2}\cdot\ln(8) = 2.713$ while the expected utility of $L_1$ is $\frac{1}{5}\cdot\ln(30) + \frac{1}{10}\cdot\ln(28) + \frac{1}{10}\cdot\ln(24) + \frac{1}{5}\cdot\ln(18) + \frac{2}{5}\cdot\ln(8) = 2.741$. Thus Paul prefers $L_1$ to $L_2$.

**Exercise 4.10.** **(a)** To normalize $U$ first add 10 to each value and then divide by 180. Denote the normalization of $U$ by $\bar{U}$. Then

$$\begin{pmatrix} & o_1 & o_2 & o_3 & o_4 & o_5 \\ \bar{U}: & \frac{54}{180}=0.3 & \frac{180}{180}=1 & \frac{0}{180}=0 & \frac{36}{180}=0.2 & \frac{108}{180}=0.6 \end{pmatrix}.$$

To normalize $V$ first subtract 5 from each value and then divide by 100. Denote the normalization of $V$ by $\bar{V}$. Then

$$\begin{pmatrix} & o_1 & o_2 & o_3 & o_4 & o_5 \\ \bar{V}: & \frac{27}{90}=0.3 & \frac{90}{90}=1 & \frac{0}{90}=0 & \frac{18}{90}=0.2 & \frac{54}{90}=0.6 \end{pmatrix}.$$

**(b)** The transformation is of the form $V(o) = aU(o)+b$. To find the values of $a$ and $b$ plug in two sets of values and solve the system of equations $\begin{cases} 44a+b=32 \\ 170a+b=95 \end{cases}$. The solution is $a=\frac{1}{2},\ b=10$. Thus $V(o)=\frac{1}{2}U(o)+10$.

**Exercise 4.11.** **(a)** Ann prefers $L_3$ to $L_4$. In fact, $\mathbb{E}[U_{Ann}(L_3)] = \sqrt{28} = 5.2915$ while $\mathbb{E}[U_{Ann}(L_4)] = \frac{1}{2}\sqrt{10}+\frac{1}{2}\sqrt{50} = 5.1167$.

**(b)** Bob prefers $L_4$ to $L_3$. In fact, $\mathbb{E}[U_{Bob}(L_3)] = 2(28) - \frac{28^4}{100^3} = 55.38534$ while

$$\mathbb{E}[U_{Bob}(L_4)] = \frac{1}{2}\left[2(10)-\frac{10^4}{100^3}\right]+\frac{1}{2}\left[2(50)-\frac{50^4}{100^3}\right] = 56.87.$$

(c) The expected value of lottery $L_4$ is \$30; thus a risk-averse person would strictly prefer \$30 with certainty to the lottery $L_4$. We saw in part (a) that for Ann the expected utility of lottery $L_4$ is 5.1167; the utility of \$30 is $\sqrt{30} = 5.47723$. Thus Ann would indeed choose \$30 for certain over the lottery $L_4$. We saw in part (b) that for Bob the expected utility of lottery $L_4$ is 56.87; the utility of \$30 is $2(30) - \frac{30^4}{100^3} = 59.19$. Thus Bob would indeed choose \$30 for certain over the lottery $L_4$.

**Exercise 4.12.** The simple lottery is $\begin{pmatrix} o_1 & o_2 & o_3 & o_4 \\ \frac{18}{240} & \frac{103}{240} & \frac{95}{240} & \frac{24}{240} \end{pmatrix}$. For example, the probability of $o_2$ is computed as follows: $\frac{1}{8}\left(\frac{1}{10}\right) + \frac{1}{4}(1) + \frac{1}{8}(0) + \frac{1}{2}\left(\frac{1}{3}\right) = \frac{103}{240}$.

**Exercise 4.13.** Using the stated indifference, use lottery $L$ to construct the compound lottery $\begin{pmatrix} \begin{pmatrix} o_2 & o_4 \\ \frac{1}{4} & \frac{3}{4} \end{pmatrix} & \begin{pmatrix} o_3 & o_4 \\ \frac{3}{5} & \frac{2}{5} \end{pmatrix} & o_3 & o_4 \\ \frac{1}{3} & \frac{2}{9} & \frac{1}{9} & \frac{1}{3} \end{pmatrix}$, whose corresponding simple lottery is $L' = \begin{pmatrix} o_1 & o_2 & o_3 & o_4 \\ 0 & \frac{1}{12} & \frac{11}{45} & \frac{121}{180} \end{pmatrix}$. Then, by the Independence Axiom, $L \sim L'$.

**Exercise 4.14.** Let $W$ be your initial wealth. The basic outcomes are:

1. you do not pay any money, don't die and live to enjoy your wealth $W$ (denote this outcome by $A_0$),
2. you pay \$$Y$, don't die and live to enjoy your remaining wealth $W - Y$ (call this outcome $A_Y$),
3. you pay \$$X$, don't die and live to enjoy your remaining wealth $W - X$ (call this outcome $A_X$),
4. you die (call this outcome $D$); this could happen because (a) you do not pay any money, roll the die and drink the poison or (b) you pay \$$Y$, roll the die and drink the poison; we assume that you are indifferent between these two outcomes.

Since, by hypothesis, $X > Y$, your ranking of these outcomes must be $A_0 \succ A_Y \succ A_X \succ D$. If you satisfy the von Neumann-Morgenstern axioms, then your preferences can be represented by a von Neumann-Morgenstern utility function $U$ defined on the set of basic outcomes. We can normalize your utility function by setting $U(A_0) = 1$ and $U(D) = 0$. Furthermore, it must be that

$$U(A_Y) > U(A_X) \qquad (\heartsuit)$$

The maximum amount $p$ that you are willing to pay is that amount that makes you indifferent between (1) rolling the die with the initial number of poisoned glasses and (2) giving up $\$p$ and rolling the die with one less poisoned glass. Thus – based on your answers – you are indifferent between the two lotteries $\begin{pmatrix} D & A_0 \\ \frac{1}{6} & \frac{5}{6} \end{pmatrix}$ and $\begin{pmatrix} A_X \\ 1 \end{pmatrix}$ and you are indifferent between the two lotteries: $\begin{pmatrix} D & A_0 \\ \frac{4}{6} & \frac{2}{6} \end{pmatrix}$ and $\begin{pmatrix} D & A_Y \\ \frac{3}{6} & \frac{3}{6} \end{pmatrix}$. Thus $\underbrace{\tfrac{1}{6}U(D) + \tfrac{5}{6}U(A_0)}_{=\frac{1}{6}0+\frac{5}{6}1=\frac{5}{6}} = U(A_X)$ and $\underbrace{\tfrac{4}{6}U(D) + \tfrac{2}{6}U(A_0)}_{=\frac{4}{6}0+\frac{2}{6}1=\frac{2}{6}} = \underbrace{\tfrac{3}{6}U(D) + \tfrac{3}{6}U(A_Y)}_{=\frac{3}{6}0+\frac{3}{6}U(A_Y)}$. Hence $U(A_X) = \frac{5}{6}$ and $U(A_Y) = \frac{2}{3} = \frac{4}{6}$, so that $U(A_X) > U(A_Y)$, contradicting (♥).

Chapter 5

# Mixed Strategies in Strategic-Form Games

## 5.1 Strategic-form games with cardinal payoffs

At the end of Chapter 3 we discussed the possibility of incorporating random events in extensive-form games by means of chance moves. The introduction of chance moves gives rise to probabilistic outcomes and thus to the issue of how a player might rank such outcomes. Random events can also occur in strategic-form games, as the following example illustrates.

**Example 5.1.** Consider the following, very simple, first-price auction. Two players simultaneously submit a bid for a painting. Only two bids are possible: \$100 and \$200. If one player bids \$200 and the other \$100 then the high bidder wins the painting and has to pay her own bid. If the two players bid the same amount then a fair coin is tossed and if the outcome is Heads the winner is Player 1 (who then has to pay her own bid) while if the outcome is Tails the winner is Player 2 (who then has to pay her own bid). Define the following basic outcomes:

$o_1$ : Player 1 gets the painting and pays \$100
$o_2$ : Player 2 gets the painting and pays \$100
$o_3$ : Player 2 gets the painting and pays \$200
$o_4$ : Player 1 gets the painting and pays \$200

Then this situation can be represented by means of the game-frame in strategic form shown in Table 5.1.

| | | Player 2 | |
|---|---|---|---|
| | | bid \$100 | bid \$200 |
| Player 1 | bid \$100 | $\begin{pmatrix} o_1 & o_2 \\ \frac{1}{2} & \frac{1}{2} \end{pmatrix}$ | $o_3$ |
| | bid \$200 | $o_4$ | $\begin{pmatrix} o_3 & o_4 \\ \frac{1}{2} & \frac{1}{2} \end{pmatrix}$ |

where

$o_1$ : Player 1 gets the painting and pays \$100

$o_2$ : Player 2 gets the painting and pays \$100

$o_3$ : Player 2 gets the painting and pays \$200

$o_4$ : Player 1 gets the painting and pays \$200

**Table 5.1**

A game-frame in strategic form representing Example 5.1

Suppose that Player 1 ranks the basic outcomes as follows: $o_1 \succ o_4 \succ o_2 \sim o_3$, that is, she prefers winning to not winning; conditional on winning, she prefers to pay less and, conditional on not winning, she is indifferent as to how much Player 2 pays. Suppose also that Player 1 believes that Player 2 is going to submit a bid of \$100 (perhaps she has been informed of this by somebody spying on Player 2). What should we expect Player 1 to do? Knowing her ranking of the basic outcomes is of no help, because we need to know how she ranks the probabilistic outcome $\begin{pmatrix} o_1 & o_2 \\ \frac{1}{2} & \frac{1}{2} \end{pmatrix}$ relative to the basic outcome $o_4$.

The theory of expected utility introduced in Chapter 4 provides one possible answer to the question of how players rank probabilistic outcomes. With the aid of expected utility theory we can now generalize the definition strategic-form game. First we generalize the notion of game-frame in strategic form (Definition 1.1, Chapter 1) by allowing probabilistic outcomes, or lotteries, to be associated with strategy profiles. In the following definition, the bulleted items coincide with the first three items of Definition 1.1 (Chapter 1); the modified item is the last one, preceded by the symbol ♦.

**Definition 5.1.** A *game-frame in strategic form* is a quadruple $\langle I, (S_i)_{i \in I}, O, f \rangle$ where:

- $I = \{1,...,n\}$ is a set of *players* ($n \geq 2$).
- For every Player $i \in I$, $S_i$ is the set of *strategies* (or choices) of Player *i*. As before, we denote by $S = S_1 \times ... \times S_n$ the set of *strategy profiles*.
- *O* is a set of *basic outcomes*.
- ♦ $f : S \to \mathcal{L}(O)$ is a function that associates with every strategy profile *s* a lottery over the set of basic outcomes *O* (as in Chapter 4, we denote by $\mathcal{L}(O)$ the set of lotteries, or probability distributions, over *O*).

If, for every $s \in S$, $f(s)$ is a degenerate lottery (that is, a basic outcome) then we are back to Definition 1.1 (Chapter 1).

From a game-frame one obtains a game by adding, for every player $i \in I$, a von Neumann-Morgenstern ranking $\succsim_i$ of the elements of $\mathcal{L}(O)$. It is more convenient to represent such a ranking by means of a von Neumann-Morgenstern utility function $U_i : O \to \mathbb{R}$. We denote by $\mathbb{E}[U_i(f(s))]$ the expected utility of lottery $f(s) \in \mathcal{L}(O)$ for Player *i*. The following definition mirrors Definition 1.2 of Chapter 1.

**Definition 5.2.** A *game in strategic form with cardinal payoffs* is a quintuple $\langle I, (S_i)_{i \in I}, O, f, (\succsim_i)_{i \in I} \rangle$ where:

- $\langle I, (S_i)_{i \in I}, O, f \rangle$ is a game-frame in strategic form (Definition 5.1) and
- for every Player $i \in I$, $\succsim_i$ is a von Neumann-Morgenstern ranking of the set of lotteries $\mathcal{L}(O)$.

If we represent each ranking $\succsim_i$ by means of a von Neumann-Morgenstern utility function $U_i$ and define $\pi_i : S \to \mathbb{R}$ by $\pi_i(s) = \mathbb{E}[U_i(f(s))]$, then the triple $\langle I, (S_1,...,S_n), (\pi_1,...,\pi_n) \rangle$ is called a *reduced-form game in strategic form with cardinal payoffs* ('reduced-form' because some information is lost, namely the specification of the possible outcomes). The function $\pi_i : S \to \mathbb{R}$ is called the *von Neumann-Morgenstern payoff function* of Player *i*.

For example, consider the first-price auction of Example 5.1 whose game-frame in strategic form was shown in Table 5.1. Let $O = \{o_1, o_2, o_3, o_4\}$ and suppose that Player 1 has a von Neumann-Morgenstern ranking of $\mathcal{L}(O)$ that is represented by the following von Neumann-Morgenstern utility function $U_1$ (note that the implied ordinal ranking of the basic outcomes is indeed $o_1 \succ o_4 \succ o_2 \sim o_3$):

$$\begin{array}{lcccc} \text{outcome}: & o_1 & o_2 & o_3 & o_4 \\ U_1: & 4 & 1 & 1 & 2 \end{array}$$

Then, for Player 1, the expected utility of lottery $\begin{pmatrix} o_1 & o_2 \\ \frac{1}{2} & \frac{1}{2} \end{pmatrix}$ is 2.5 and the expected utility of lottery $\begin{pmatrix} o_3 & o_4 \\ \frac{1}{2} & \frac{1}{2} \end{pmatrix}$ is 1.5. Suppose also that Player 2 has (somewhat spiteful) preferences represented by the following von Neumann-Morgenstern utility function $U_2$:

$$\begin{array}{lcccc} \text{outcome}: & o_1 & o_2 & o_3 & o_4 \\ U_2: & 1 & 6 & 4 & 5 \end{array}$$

Thus, for Player 2, the expected utility of lottery $\begin{pmatrix} o_1 & o_2 \\ \frac{1}{2} & \frac{1}{2} \end{pmatrix}$ is 3.5 and the expected utility of lottery $\begin{pmatrix} o_3 & o_4 \\ \frac{1}{2} & \frac{1}{2} \end{pmatrix}$ is 4.5. Then we can represent the game in reduced form as shown in Table 5.2.

| | | Player 2: $100 | Player 2: $200 |
|---|---|---|---|
| Player 1 | $100 | 2.5, 3.5 | 1, 4 |
| | $200 | 2, 5 | 1.5, 4.5 |

**Table 5.2**

A cardinal game in reduced form based on the game-frame of Table 5.1

The game of Table 5.2 does not have any Nash equilibria. However, we will show in the next section that if we extend the notion of strategy, by allowing players to choose randomly, then the game of Table 5.2 does have a Nash equilibrium.

Test your understanding of the concepts introduced in this section, by going through the exercises in Section 5.E.1 of Appendix 5.E at the end of this chapter.

## 5.2 Mixed strategies

**Definition 5.3.** Consider a game in strategic form with cardinal payoffs and recall that $S_i$ denotes the set of strategies of Player *i*. From now on, we shall call $S_i$ the set of *pure strategies* of Player *i*. We assume that $S_i$ is a finite set (for every $i \in I$). A *mixed strategy* of Player *i* is a probability distribution over the set of pure strategies $S_i$. The set of mixed strategies of Player *i* is denoted by $\Sigma_i$.

**Remark 5.1.** Since among the mixed strategies of Player *i* there are the degenerate strategies that assign probability 1 to a pure strategy, the set of mixed strategies includes the set of pure strategies (viewed as degenerate probability distributions).

For example, one possible mixed strategy for Player 1 in the game of Table 5.2 is $\begin{pmatrix} \$100 & \$200 \\ \frac{1}{3} & \frac{2}{3} \end{pmatrix}$. The traditional interpretation of a mixed strategy is in terms of objective randomization: the player, instead of choosing a pure strategy herself, delegates the choice to a random device.[4] For example, Player 1 choosing the mixed strategy $\begin{pmatrix} \$100 & \$200 \\ \frac{1}{3} & \frac{2}{3} \end{pmatrix}$ is interpreted as a decision to let, say, a die determine whether she will bid \$100 or \$200: Player 1 will roll a die and if the outcome is 1 or 2 then she will bid \$100, while if the outcome is 3, 4, 5 or 6 then she will bid \$200. Suppose that Player 1 chooses this mixed strategy and Player 2 chooses the mixed strategy $\begin{pmatrix} \$100 & \$200 \\ \frac{3}{5} & \frac{2}{5} \end{pmatrix}$. Since the players rely on independent random devices, this pair of mixed strategies gives rise to the following probabilistic outcome:

$$\begin{pmatrix} \text{strategy profile} & (\$100,\$100) & (\$100,\$200) & (\$200,\$100) & (\$200,\$200) \\ \text{outcome} & \begin{pmatrix} o_1 & o_2 \\ \frac{1}{2} & \frac{1}{2} \end{pmatrix} & o_3 & o_4 & \begin{pmatrix} o_3 & o_4 \\ \frac{1}{2} & \frac{1}{2} \end{pmatrix} \\ \text{probability} & \frac{1}{3}\frac{3}{5}=\frac{3}{15} & \frac{1}{3}\frac{2}{5}=\frac{2}{15} & \frac{2}{3}\frac{3}{5}=\frac{6}{15} & \frac{2}{3}\frac{2}{5}=\frac{4}{15} \end{pmatrix}$$

[4] An alternative interpretation of mixed strategies in terms of beliefs will be discussed in a later chapter.

If the two players have von Neumann-Morgenstern preferences, then - by the Compound Lottery Axiom (Chapter 4) – they will view the above as the following lottery:

$$\begin{pmatrix} \text{outcome} & o_1 & o_2 & o_3 & o_4 \\ \text{probability} & \frac{3}{30} & \frac{3}{30} & \frac{8}{30} & \frac{16}{30} \end{pmatrix}.$$

Using the von Neumann-Morgenstern utility functions postulated above, namely

$$\begin{matrix} \text{outcome:} & o_1 & o_2 & o_3 & o_4 \\ U_1: & 4 & 1 & 1 & 2 \\ U_2: & 1 & 6 & 4 & 5 \end{matrix}$$

, the lottery $\begin{pmatrix} o_1 & o_2 & o_3 & o_4 \\ \frac{3}{30} & \frac{3}{30} & \frac{8}{30} & \frac{16}{30} \end{pmatrix}$ has an expected utility of $\frac{3}{30}4+\frac{3}{30}1+\frac{8}{30}1+\frac{16}{30}2=\frac{55}{30}$ for Player 1 and $\frac{3}{30}1+\frac{3}{30}6+\frac{8}{30}4+\frac{16}{30}5=\frac{133}{30}$ for Player 2. Thus we can define the payoffs of the two players from this mixed strategy profile by

$$\Pi_1\left(\begin{pmatrix} \$100 & \$200 \\ \frac{1}{3} & \frac{2}{3} \end{pmatrix}, \begin{pmatrix} \$100 & \$200 \\ \frac{3}{5} & \frac{2}{5} \end{pmatrix}\right)=\frac{55}{30} \text{ and } \Pi_2\left(\begin{pmatrix} \$100 & \$200 \\ \frac{1}{3} & \frac{2}{3} \end{pmatrix}, \begin{pmatrix} \$100 & \$200 \\ \frac{3}{5} & \frac{2}{5} \end{pmatrix}\right)=\frac{133}{30}.$$

Note that we can calculate these payoffs in a different – but equivalent – way by using the reduced-form game of Table 5.2, as follows.

$$\begin{pmatrix} \text{strategy profile} & (\$100,\$100) & (\$100,\$200) & (\$200,\$100) & (\$200,\$200) \\ \text{expected utilities} & (2.5,3.5) & (1,4) & (2,5) & (1.5,4.5) \\ \text{probability} & \frac{1}{3}\frac{3}{5}=\frac{3}{15} & \frac{1}{3}\frac{2}{5}=\frac{2}{15} & \frac{2}{3}\frac{3}{5}=\frac{6}{15} & \frac{2}{3}\frac{2}{5}=\frac{4}{15} \end{pmatrix}$$

so that the expected payoff of Player 1 is $\frac{3}{15}(2.5)+\frac{2}{15}(1)+\frac{6}{15}(2)+\frac{4}{15}(1.5)=\frac{55}{30}$ and the expected payoff of Player 2 is $\frac{3}{15}(3.5)+\frac{2}{15}(4)+\frac{6}{15}(5)+\frac{4}{15}(4.5)=\frac{133}{30}$.

The above example motivates the following definition. First some notation. Let $\sigma_i \in \Sigma_i$ be a mixed strategy of player *i*; then, for every pure strategy $s_i \in S_i$ of Player *i*, we denote by $\sigma_i(s_i)$ the probability that $\sigma_i$ assigns to $s_i$.[5] Let $\Sigma$ be the set of mixed-strategy profiles, that is, $\Sigma = \Sigma_1 \times ... \times \Sigma_n$. Consider a mixed-strategy profile $\sigma = (\sigma_1,...,\sigma_n) \in \Sigma$ and a pure-strategy profile $s = (s_1,...,s_n) \in S$;

[5] In the above example, if $\sigma_1 = \begin{pmatrix} \$100 & \$200 \\ \frac{1}{3} & \frac{2}{3} \end{pmatrix}$ then $\sigma_1(\$200)=\frac{2}{3}$.

then we denote by $\sigma(s)$ the product of the probabilities $\sigma_i(s_i)$, that is, $\sigma(s) = \prod_{i=1}^{n} \sigma_i(s_i) = \sigma_1(s_1) \times \ldots \times \sigma_n(s_n)$.[6]

**Definition 5.4.** Consider a reduced-form game in strategic form with cardinal payoffs $G = \langle I, (S_1,...,S_n), (\pi_1,...,\pi_n) \rangle$ (Definition 5.2), where, for every Player $i \in I$, the set of pure strategies $S_i$ is finite. Then the *mixed-strategy extension of* $G$ is the reduced-form game in strategic form $\langle I, (\Sigma_1,...,\Sigma_n), (\Pi_1,...,\Pi_n) \rangle$ where, for every Player $i \in I$,

- $\Sigma_i$ is the set of mixed strategies of Player $i$ in $G$ (that is, $\Sigma_i$ is the set of probability distributions over $S_i$).
- The payoff function $\Pi_i : \Sigma \to \mathbb{R}$ is defined by $\Pi_i(\sigma) = \sum_{s \in S} \sigma(s) \pi_i(s)$.[7]

**Definition 5.5.** Fix a reduced-form game in strategic form with cardinal payoffs $G = \langle I, (S_1,...,S_n), (\pi_1,...,\pi_n) \rangle$ (Definition 5.3), where, for every $i \in I$, the set of pure strategies $S_i$ is finite. A *Nash equilibrium in mixed-strategies of* $G$ is a Nash equilibrium of the mixed-strategy extension of $G$.

For example, consider the reduced-form game of Table 5.2, which is reproduced in Table 5.3 below with all the payoffs multiplied by 10 (this corresponds to representing the preferences of the players with different utility functions that are a obtained from the ones used above by multiplying them by 10). Is $\sigma = (\sigma_1, \sigma_2)$ with $\sigma_1 = \begin{pmatrix} \$100 & \$200 \\ \frac{1}{3} & \frac{2}{3} \end{pmatrix}$ and $\sigma_2 = \begin{pmatrix} \$100 & \$200 \\ \frac{3}{5} & \frac{2}{5} \end{pmatrix}$ a mixed-strategy Nash equilibrium of this game?

---

[6] In the above example, if $\sigma = (\sigma_1, \sigma_2)$ with $\sigma_1 = \begin{pmatrix} \$100 & \$200 \\ \frac{1}{3} & \frac{2}{3} \end{pmatrix}$ and $\sigma_2 = \begin{pmatrix} \$100 & \$200 \\ \frac{3}{5} & \frac{2}{5} \end{pmatrix}$ then $\sigma_1(\$200) = \frac{2}{3}$, $\sigma_2(\$100) = \frac{3}{5}$ and thus $\sigma((\$200, \$100)) = \frac{2}{3}\frac{3}{5} = \frac{6}{15}$.

[7] In the above example, if $\sigma_1 = \begin{pmatrix} \$100 & \$200 \\ \frac{1}{3} & \frac{2}{3} \end{pmatrix}$ and $\sigma_2 = \begin{pmatrix} \$100 & \$200 \\ \frac{3}{5} & \frac{2}{5} \end{pmatrix}$ then $\Pi_1(\sigma_1, \sigma_2) = \frac{3}{15}(2.5) + \frac{2}{15}(1) + \frac{6}{15}(2) + \frac{4}{15}(1.5) = \frac{55}{30}$.

| | | Player 2 | |
|---|---|---|---|
| | | $100 | $200 |
| Player 1 | $100 | 25 35 | 10 40 |
| | $200 | 20 50 | 15 45 |

**Table 5.3**
The game of Table 5.2 with the payoffs multiplied by 10

The payoff of Player 1 is $\Pi_1(\sigma_1,\sigma_2) = \frac{3}{15}(25) + \frac{2}{15}(10) + \frac{6}{15}(20) + \frac{4}{15}(15) = \frac{55}{3}$ $= 18.33$. If Player 1 switched from $\sigma_1 = \begin{pmatrix} \$100 & \$200 \\ \frac{1}{3} & \frac{2}{3} \end{pmatrix}$ to $\hat{\sigma}_1 = \begin{pmatrix} \$100 & \$200 \\ 1 & 0 \end{pmatrix}$, that is, to the pure strategy $100, then Player 1's payoff would be larger: $\Pi_1(\hat{\sigma}_1,\sigma_2) = \frac{3}{5}(25) + \frac{2}{5}(10) = 19$. Thus $(\sigma_1,\sigma_2)$ is not a Nash equilibrium.

John Nash, who shared the 1994 Nobel prize in economics with John Harsanyi and Reinhard Selten, proved the following theorem.

**Theorem 5.1 [Nash, 1951].** Every reduced-form game in strategic form with cardinal payoffs $\langle I, (S_1,...,S_n), (\pi_1,...,\pi_n) \rangle$ (Definition 5.2), where, for every Player $i \in I$, the set of pure strategies $S_i$ is finite, has at least one Nash equilibrium in mixed-strategies.

We will not give the proof of this theorem, since it is rather complex (it requires the use of fixed-point theorems).

Going back to the game of Table 5.3, let us verify that, on the other hand, $\sigma^* = (\sigma_1^*, \sigma_2^*)$ with $\sigma_1^* = \sigma_2^* = \begin{pmatrix} \$100 & \$200 \\ \frac{1}{2} & \frac{1}{2} \end{pmatrix}$ is a Nash equilibrium in mixed strategies. The payoff of Player 1 is

$$\Pi_1(\sigma_1^*,\sigma_2^*) = \tfrac{1}{4}(25) + \tfrac{1}{4}(10) + \tfrac{1}{4}(20) + \tfrac{1}{4}(15) = \tfrac{70}{4} = 17.5.$$

Could Player 1 obtain a larger payoff with some other mixed strategy $\sigma_1 = \begin{pmatrix} \$100 & \$200 \\ p & 1-p \end{pmatrix}$ for some $p \neq \frac{1}{2}$? Fix an arbitrary $p \in [0,1]$ and let us compute

Player 1's payoff if she uses the strategy $\sigma_1 = \begin{pmatrix} \$100 & \$200 \\ p & 1-p \end{pmatrix}$ against Player 2's mixed strategy $\sigma_2^* = \begin{pmatrix} \$100 & \$200 \\ \frac{1}{2} & \frac{1}{2} \end{pmatrix}$:

$$\Pi_1\left(\begin{pmatrix} \$100 & \$200 \\ p & 1-p \end{pmatrix}, \begin{pmatrix} \$100 & \$200 \\ \frac{1}{2} & \frac{1}{2} \end{pmatrix}\right) = \tfrac{1}{2}p(25) + \tfrac{1}{2}p(10) + \tfrac{1}{2}(1-p)(20) + \tfrac{1}{2}(1-p)(15) =$$
$$p\left(\tfrac{1}{2}25 + \tfrac{1}{2}10\right) + (1-p)\left(\tfrac{1}{2}20 + \tfrac{1}{2}15\right) = \tfrac{35}{2} = 17.5$$

Thus if Player 2 uses the mixed strategy $\sigma_2^* = \begin{pmatrix} \$100 & \$200 \\ \frac{1}{2} & \frac{1}{2} \end{pmatrix}$, then *Player 1 gets the same payoff no matter what mixed strategy she employs*. It follows that any mixed strategy of Player 1 is a best reply to $\sigma_2^* = \begin{pmatrix} \$100 & \$200 \\ \frac{1}{2} & \frac{1}{2} \end{pmatrix}$; in particular, $\sigma_1^* = \begin{pmatrix} \$100 & \$200 \\ \frac{1}{2} & \frac{1}{2} \end{pmatrix}$ is a best reply to $\sigma_2^* = \begin{pmatrix} \$100 & \$200 \\ \frac{1}{2} & \frac{1}{2} \end{pmatrix}$. It is easy to verify that the same applies to Player 2: any mixed strategy of Player 2 is a best reply to Player 1's mixed strategy $\sigma_1^* = \begin{pmatrix} \$100 & \$200 \\ \frac{1}{2} & \frac{1}{2} \end{pmatrix}$. Hence $\sigma^* = (\sigma_1^*, \sigma_2^*)$ is indeed a Nash equilibrium in mixed strategies.

We will see in the next section that this "indifference" phenomenon is true in general.

**Remark 5.2.** Since, among the mixed strategies of Player $i$ there are the degenerate strategies that assign probability 1 to a pure strategy (Remark 5.5), every Nash equilibrium in pure strategies is also a Nash equilibrium in mixed strategies. That is, the set of mixed-strategy Nash equilibria includes the set of pure-strategy Nash equilibria.

Test your understanding of the concepts introduced in this section, by going through the exercises in Section 5.E.2 of Appendix 5.E at the end of this chapter.

## 5.3 Computing the mixed-strategy Nash equilibria

How can we find the mixed-strategy equilibria of a given game? The first important observation is that if a pure strategy is strictly dominated by another pure strategy then it cannot be played with positive probability at a Nash equilibrium. Thus, for the purpose of finding Nash equilibria, one can delete all the strictly dominated strategies and focus on the resulting game. But then the same reasoning applies to the resulting game and one can delete all the strictly dominated strategies in that game, and so on. Thus we have the following observation.

**Remark 5.3.** In order to find the mixed-strategy Nash equilibria of a game one can first apply the iterated deletion of strictly dominated strategies (IDSDS: Section 1.5.1, Chapter 1) and then find the Nash equilibria of the resulting game (which can then be viewed as Nash equilibria of the original game where all the pure strategies that were deleted are assigned zero probability). Note, however, that – as we will see in Section 5.4 - one can perform more deletions than allowed by the IDSDS procedure.

For example, consider the game of Table 5.4.

| | | Player 2 | | |
|---|---|---|---|---|
| | | *E* | *F* | *G* |
| Player 1 | *A* | 2   4 | 3   3 | 6   0 |
| | *B* | 4   0 | 2   4 | 4   2 |
| | *C* | 3   3 | 4   2 | 3   1 |
| | *D* | 3   6 | 1   1 | 2   6 |

**Table 5.4**
A reduced-form game with cardinal payoffs

In this game there are no pure-strategy Nash equilibria; however, by Nash's theorem there will be at least one mixed-strategy equilibrium. To find it we can first note that, for Player 1, *D* is strictly dominated by *B*; deleting *D* we get a smaller game where, for Player 2, *G* is strictly dominated by *F*. Deleting *G* we are left with a smaller game where *A* is strictly dominated by *C*. Deleting *A* we are left with the game shown in Table 5.5.

| | | Player 2 | |
|---|---|---|---|
| | | $E$ | $F$ |
| Player 1 | $B$ | 4   0 | 2   4 |
| | $C$ | 3   3 | 4   2 |

**Table 5.5**
The result of applying the IDSDS procedure to the game of Table 5.4.

We will see below that the game of Table 5.5 has a unique Nash equilibrium in mixed strategies given by $\left(\begin{pmatrix} B & C \\ \frac{1}{5} & \frac{4}{5} \end{pmatrix}, \begin{pmatrix} E & F \\ \frac{2}{3} & \frac{1}{3} \end{pmatrix}\right)$. Thus the game of Table 5.4 has a unique Nash equilibrium in mixed strategies given by $\left(\begin{pmatrix} A & B & C & D \\ 0 & \frac{1}{5} & \frac{4}{5} & 0 \end{pmatrix}, \begin{pmatrix} E & F & G \\ \frac{2}{3} & \frac{1}{3} & 0 \end{pmatrix}\right)$.

Once we have simplified the game as suggested in Remark 5.3, in order to find the mixed-strategy Nash equilibria we can use the following result. First we recall some notation that was introduced in Chapter 1. Given a mixed-strategy profile $\sigma = (\sigma_1,...,\sigma_n)$ and a Player $i$, we denote by $\sigma_{-i} = (\sigma_1,...,\sigma_{i-1},\sigma_{i+1},...,\sigma_n)$ the profile of strategies of the players other than $i$ and use $(\sigma_i,\sigma_{-i})$ as an alternative notation for $\sigma$; furthermore, $(\tau_i,\sigma_{-i})$ denotes the result of replacing $\sigma_i$ with $\tau_i$ in $\sigma$, that is, $(\tau_i,\sigma_{-i}) = (\sigma_1,...,\sigma_{i-1},\tau_i,\sigma_{i+1},...,\sigma_n)$.

**Theorem 5.2.** Consider a reduced-form game in strategic form with cardinal payoffs. Suppose that $\sigma^* = (\sigma_1^*,...,\sigma_n^*)$ is a Nash equilibrium in mixed strategies. Consider an arbitrary Player $i$. Let $\Pi_i^* = \Pi_i(\sigma^*)$ be the payoff of Player $i$ at this Nash equilibrium and let $s_{ij}, s_{ik} \in S_i$ be two pure strategies of Player $i$ such that $\sigma_i^*(s_{ij}) > 0$ and $\sigma_i^*(s_{ik}) > 0$, that is, $s_{ij}$ and $s_{ik}$ are two pure strategies to which the mixed strategy $\sigma_i^*$ of Player $i$ assigns positive probability. Then $\Pi_i\left((s_{ij},\sigma_{-i}^*)\right) = \Pi_i\left((s_{ik},\sigma_{-i}^*)\right) = \Pi_i^*$. In other words, when the other players use the mixed-strategy profile $\sigma_{-i}^*$, Player $i$ gets the same payoff no matter whether she plays the mixed strategy $\sigma_i^*$ or the pure strategy $s_{ij}$ or the pure strategy $s_{ik}$.

The details of the proof of Theorem 5.2 will be omitted, but the idea is simple: if $s_{ij}$ and $s_{ik}$ are two pure strategies to which the mixed strategy $\sigma_i^*$ of Player $i$ assigns positive probability and $\Pi_i\left((s_{ij},\sigma_{-i}^*)\right) > \Pi_i\left((s_{ik},\sigma_{-i}^*)\right)$, then Player $i$ can increase her payoff from $\Pi_i^* = \Pi_i(\sigma^*)$ to a larger number by reducing the probability of $s_{ik}$ to zero and adding that probability to $\sigma_i^*(s_{ij})$, that is, by switching from $\sigma_i^*$ to the mixed strategy $\hat{\sigma}_i$ obtained as follows: $\hat{\sigma}_i(s_{ik}) = 0$, $\hat{\sigma}_i(s_{ij}) = \sigma_i^*(s_{ij}) + \sigma_i^*(s_{ik})$ and, for every other $s_i \in S_i$, $\hat{\sigma}_i(s_i) = \sigma_i^*(s_i)$. But this would contradict the hypothesis that $\sigma^* = (\sigma_1^*,...,\sigma_n^*)$ is a Nash equilibrium.

Let us now go back to the game of Table 5.5, reproduced below, and see how we can use Theorem 5.2 to find the Nash equilibrium in mixed strategies.

| | | Player 2 | |
|---|---|---|---|
| | | *E* | *F* |
| Player 1 | *B* | 4   0 | 2   4 |
| | *C* | 3   3 | 4   2 |

We want to find values of *p* and *q*, strictly between 0 and 1, such that $\left(\begin{pmatrix} B & C \\ p & 1-p \end{pmatrix}, \begin{pmatrix} E & F \\ q & 1-q \end{pmatrix}\right)$ is a Nash equilibrium. By Theorem 5.2, if Player 1 played the pure strategy *B* against $\begin{pmatrix} E & F \\ q & 1-q \end{pmatrix}$ she should get the same payoff as if she were to play the pure strategy *C*. The former would give her a payoff of $4q + 2(1-q)$ and the latter a payoff of $3q + 4(1-q)$. Thus we need *q* to be such that $4q + 2(1-q) = 3q + 4(1-q)$, that is, $q = \frac{2}{3}$. When $q = \frac{2}{3}$, both *B* and *C* give Player 1 a payoff of $\frac{10}{3}$ and thus any mixture of *B* and *C* would also give the same payoff of $\frac{10}{3}$. In other words, Player 1 is indifferent among all her mixed strategies and thus any mixed strategy is a best response to $\begin{pmatrix} E & F \\ \frac{2}{3} & \frac{1}{3} \end{pmatrix}$. Similar reasoning for Player 2 reveals that, by Theorem 5.2, we need *p* to be such that $0p + 3(1-p) = 4p + 2(1-p)$, that is, $p = \frac{1}{5}$. Against $\begin{pmatrix} B & C \\ \frac{1}{5} & \frac{4}{5} \end{pmatrix}$ any mixed strategy of Player 2 gives him the same payoff of $\frac{12}{5}$; thus any mixed strategy of Player 2 is a best reply to $\begin{pmatrix} B & C \\ \frac{1}{5} & \frac{4}{5} \end{pmatrix}$. It follows that $\left(\begin{pmatrix} B & C \\ \frac{1}{5} & \frac{4}{5} \end{pmatrix}, \begin{pmatrix} E & F \\ \frac{2}{3} & \frac{1}{3} \end{pmatrix}\right)$ is a Nash equilibrium.

**Remark 5.4.** It follows from Theorem 5.2, and was illustrated in the above example, that at a mixed strategy Nash equilibrium where Player *i* plays two or more pure strategies with positive probability, Player *i* does not have an incentive to use that mixed strategy: she would get the same payoff if, instead of randomizing, she played with probability 1 one of the pure strategies in the support of her mixed strategy (that is, if she increased the probability of any pure strategy from a positive number to 1).[8] *The only purpose of randomizing is to make the other player indifferent among two or more of his own pure strategies.*

Theorem 5.2 provides a necessary, *but not sufficient*, condition for a mixed-strategy profile to be a Nash equilibrium. To see that the condition is not sufficient, consider the game of Table 5.6 and the mixed-strategy profile $\left(\begin{pmatrix} A & B & C \\ \frac{1}{2} & \frac{1}{2} & 0 \end{pmatrix}, \begin{pmatrix} D & E \\ \frac{1}{2} & \frac{1}{2} \end{pmatrix}\right)$.

Player 2

| | | D | E |
|---|---|---|---|
| Player 1 | A | 3, 0 | 0, 2 |
| | B | 0, 2 | 3, 0 |
| | C | 2, 0 | 2, 1 |

**Table 5.6**
A reduced-form game with cardinal payoffs

Given that Player 2 plays the mixed strategy $\begin{pmatrix} D & E \\ \frac{1}{2} & \frac{1}{2} \end{pmatrix}$, Player 1 is indifferent between the two pure strategies that are in the support of her own mixed strategy, namely *A* and *B*: the payoff from playing *A* is 1.5 and so is the payoff from playing *B* (and 1.5 is also the payoff associated with the mixed strategy under consideration). However, the profile $\left(\begin{pmatrix} A & B & C \\ \frac{1}{2} & \frac{1}{2} & 0 \end{pmatrix}, \begin{pmatrix} D & E \\ \frac{1}{2} & \frac{1}{2} \end{pmatrix}\right)$ is not a Nash equilibrium, because Player 1 could get a payoff of 2 by switching to the pure strategy *C*. We know from Theorem 5.1 that this game does have a mixed-strategy Nash equilibrium. How can we find it?

[8] The support of a mixed strategy is the set of pure strategies that are assigned positive probability by that mixed strategy.

Let us calculate the best response of Player 1 to every possible mixed strategy $\begin{pmatrix} D & E \\ q & 1-q \end{pmatrix}$ of Player 2 (with $q \in [0,1]$). For Player 1 the payoff from playing $A$ against $\begin{pmatrix} D & E \\ q & 1-q \end{pmatrix}$ is $3q$, the payoff from playing $B$ is $3-3q$ and the payoff from playing $C$ is constant and equal to 2. These functions are shown in Figure 5.7.

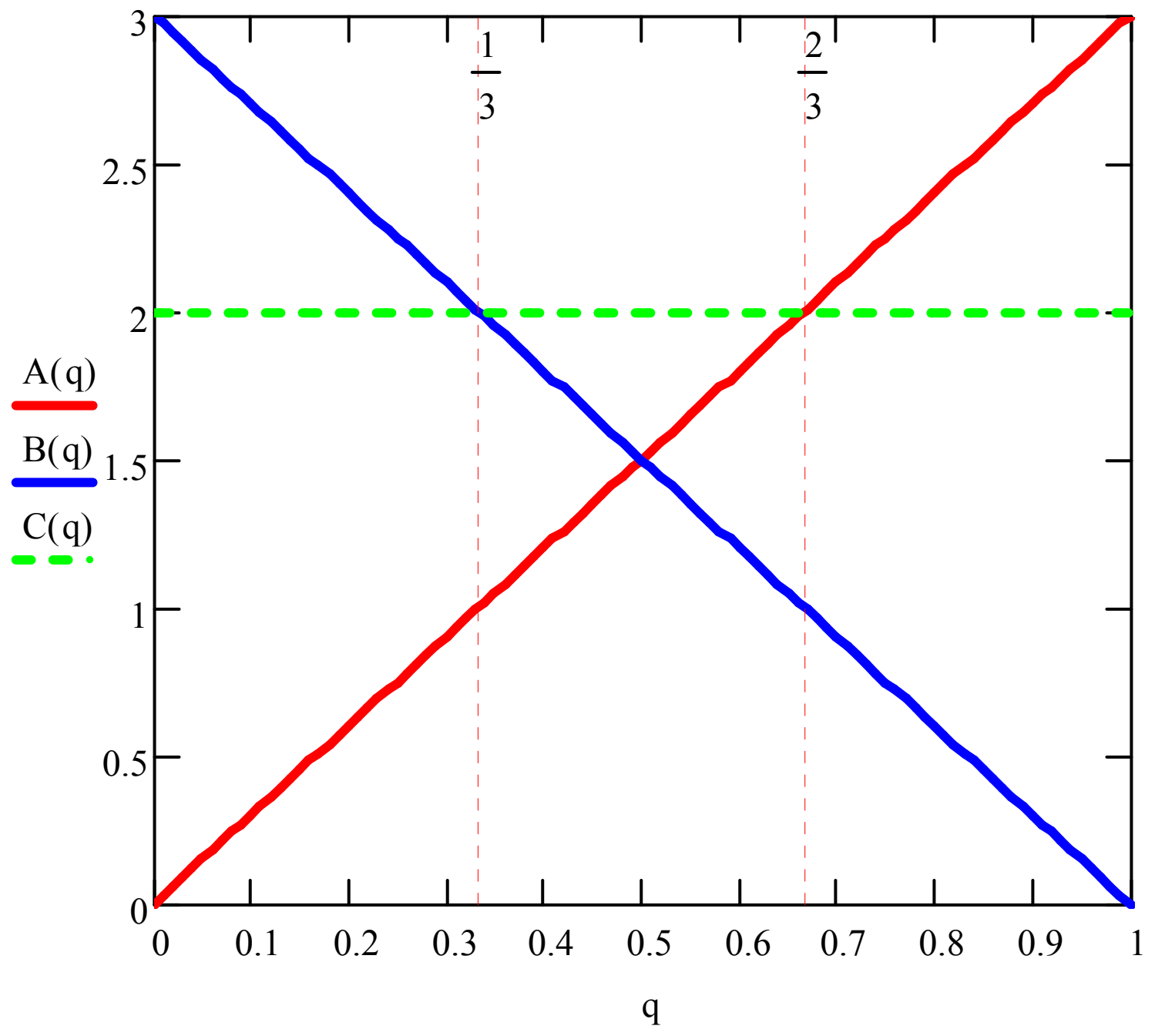

**Figure 5.7**
Player 1's payoff from each pure strategy against an arbitrary mixed strategy of Player 2.

The red upward-sloping line plots the function $A(q) = 3q$, the blue downward-sloping line plots the function $B(q) = 3-3q$ and the horizontal green dotted line the function $C(q) = 2$. The blue and green lines intersect when $q = \frac{1}{3}$ and the red and green lines intersect when $q = \frac{2}{3}$. The maximum payoff is given by the blue line up to $q = \frac{1}{3}$, then by the green line up to $q = \frac{2}{3}$ and then by the red line. Thus the best reply function of Player 1 is as follows:

$$\text{Player1's best reply} = \begin{cases} B & \text{if } 0 \le q < \frac{1}{3} \\ \begin{pmatrix} B & C \\ p & 1-p \end{pmatrix} \text{ for any } p \in [0,1] & \text{if } q = \frac{1}{3} \\ C & \text{if } \frac{1}{3} < q < \frac{2}{3} \\ \begin{pmatrix} A & C \\ p & 1-p \end{pmatrix} \text{ for any } p \in [0,1] & \text{if } q = \frac{2}{3} \\ A & \text{if } \frac{2}{3} < q \le 1 \end{cases}$$

Hence if there is a mixed-strategy Nash equilibrium it is either of the form $\left( \begin{pmatrix} A & B & C \\ 0 & p & 1-p \end{pmatrix}, \begin{pmatrix} D & E \\ \frac{1}{3} & \frac{2}{3} \end{pmatrix} \right)$ or of the form $\left( \begin{pmatrix} A & B & C \\ p & 0 & 1-p \end{pmatrix}, \begin{pmatrix} D & E \\ \frac{2}{3} & \frac{1}{3} \end{pmatrix} \right)$. The latter cannot be a Nash equilibrium for any *p*, because when Player 1 plays *B* with probability 0, *E* strictly dominates *D* for Player 2 and thus Player 2's mixed strategy $\begin{pmatrix} D & E \\ \frac{2}{3} & \frac{1}{3} \end{pmatrix}$ is not a best reply (*E* is the unique best reply). Thus the only candidate for a Nash equilibrium is of the form $\left( \begin{pmatrix} A & B & C \\ 0 & p & 1-p \end{pmatrix}, \begin{pmatrix} D & E \\ \frac{1}{3} & \frac{2}{3} \end{pmatrix} \right)$. In this case, by Theorem 5.11, we need *p* to be such that Player 2 is indifferent between *D* and *E*, that is, we need $2p = 1 - p$, that is, $p = \frac{1}{3}$. Hence the Nash equilibrium is $\left( \begin{pmatrix} A & B & C \\ 0 & \frac{1}{3} & \frac{2}{3} \end{pmatrix}, \begin{pmatrix} D & E \\ \frac{1}{3} & \frac{2}{3} \end{pmatrix} \right)$.

In games where the number of strategies or the number of players are larger than in the examples we have considered, finding the Nash equilibria involves lengthier calculations. However, computer programs have been developed that can be used to compute all the Nash equilibria of a finite game in a relatively short time.

Test your understanding of the concepts introduced in this section, by going through the exercises in Section 5.E.3 of Appendix 5.E at the end of this chapter.

## 5.4 Strict dominance and rationalizability

We remarked in the previous section that a pure strategy that is strictly dominated by another pure strategy cannot be played with positive probability at a Nash equilibrium. Thus, when looking for a Nash equilibrium, one can first simplify the game by applying the IDSDS procedure (Section 1.5.1, Chapter 1). When payoffs are cardinal (von Neumann-Morgenstern payoffs) it turns out that, *in a two-person game*, a pure strategy cannot be a best response to any mixed-strategy of the opponent not only when it is strictly dominated by another pure strategy but also when it is strictly dominated *by a mixed strategy*. To see this, consider the game of Table 5.8 below.

| | | Player 2: *D* | | Player 2: *E* | |
|---|---|---|---|---|---|
| Player 1 | *A* | 0 | 1 | 4 | 0 |
| | *B* | 1 | 2 | 1 | 4 |
| | *C* | 2 | 0 | 0 | 1 |

**Table 5.8**
A strategic-form game with cardinal payoffs.

The pure strategy *B* of Player 1 is not strictly dominated by another pure strategy and yet it cannot be a best reply to any mixed strategy of Player 2. To see this, consider an arbitrary mixed strategy $\begin{pmatrix} D & E \\ q & 1-q \end{pmatrix}$ of Player 2 with $q \in [0,1]$ . If Player 1 plays *B* against it, she gets a payoff of 1; if, instead, she plays the mixed strategy $\begin{pmatrix} A & B & C \\ \frac{1}{3} & 0 & \frac{2}{3} \end{pmatrix}$ then her payoff is $\frac{1}{3}4(1-q)+\frac{2}{3}2q = \frac{4}{3} > 1$.

**Theorem 5.3 [Pearce, 1984].** Consider a two-player reduced-form game in strategic form with cardinal payoffs, an arbitrary Player *i* and a pure strategy $s_i$ of Player *i*. Then there is no mixed-strategy of the opponent to which $s_i$ is a best response, if and only if $s_i$ is strictly dominated by a mixed strategy $\sigma_i$ of Player *i* (that is, there is a $\sigma_i \in \Sigma_i$ such that $\Pi_i(\sigma_i,\sigma_j) > \Pi_i(s_i,\sigma_j)$, for every $\sigma_j \in \Sigma_j$).

Note that, since the set of mixed strategies includes the set of pure strategies, strict dominance by a mixed strategy includes as a sub-case strict dominance by a pure strategy.

When the number of players is 3 or more, the generalization of Theorem 5.3 raises some subtle issues: see Exercise 5.14. However, we can appeal to the intuition behind Theorem 5.3 (see Remark 5.5 below) to refine the IDSDS procedure (Section 1.5.1, Chapter 1) for general $n$-player games with cardinal payoffs as follows.

**Definition 5.6 [Cardinal IDSDS].** The *Cardinal Iterated Deletion of Strictly Dominated Strategies* is the following algorithm. Given a finite $n$-player ($n \geq 2$) strategic-form game with cardinal payoffs $G$, let $G^1$ be the game obtained by removing from $G$, for every Player $i$, those *pure* strategies of Player $i$ (if any) that are strictly dominated in $G$ by some *mixed* strategy of Player $i$; let $G^2$ be the game obtained by removing from $G^1$, for every Player $i$, those pure strategies of Player $i$ (if any) that are strictly dominated in $G^1$ by some mixed strategy of Player $i$, and so on. Let $G^\infty$ be the output of this procedure. Since the initial game $G$ is finite, $G^\infty$ will be obtained in a finite number of steps. For every Player $i$, the pure strategies of Player $i$ in $G^\infty$ are called her *rationalizable* strategies.

Figure 5.9 illustrates this procedure as applied to the game in Panel (i).

| | | Player 2: $D$ | | $E$ | | $F$ | |
|---|---|---|---|---|---|---|---|
| Player 1 | $A$ | 3 | 4 | 2 | 1 | 1 | 2 |
| | $B$ | 0 | 0 | 1 | 3 | 4 | 1 |
| | $C$ | 1 | 4 | 1 | 4 | 2 | 6 |

(*i*) The game $G^0 = G$

| | | Player 2: $D$ | | $E$ | |
|---|---|---|---|---|---|
| Player 1 | $A$ | 3 | 4 | 2 | 1 |
| | $B$ | 0 | 0 | 1 | 3 |

(*iii*) The game $G^2$ after Step 2

| | | Player 2: $D$ | | $E$ | | $F$ | |
|---|---|---|---|---|---|---|---|
| Player 1 | $A$ | 3 | 4 | 2 | 1 | 1 | 2 |
| | $B$ | 0 | 0 | 1 | 3 | 4 | 1 |

(*ii*) The game $G^1$ after Step 1

| | | Player 2: $D$ | | $E$ | |
|---|---|---|---|---|---|
| Player 1 | $A$ | 3 | 4 | 2 | 1 |

| | | Player 2: $D$ | |
|---|---|---|---|
| Player 1 | $A$ | 3 | 4 |

(*iv*) The game $G^3$ at the top.
The game $G^4 = G^\infty$ at the bottom

**Figure 5.9**
Application of the cardinal IDSDS procedure.

In the first step, the pure strategy $C$ of Player 1 is deleted, because it is strictly dominated by the mixed strategy $\begin{pmatrix} A & B \\ \frac{1}{2} & \frac{1}{2} \end{pmatrix}$ thus yielding game $G^1$ shown in Panel (ii). In the second step, the pure strategy $F$ of Player 2 is deleted, because it is strictly dominated by the mixed strategy $\begin{pmatrix} D & E \\ \frac{1}{2} & \frac{1}{2} \end{pmatrix}$ thus yielding game $G^2$ shown in Panel (iii). In the third step, $B$ is deleted because it is strictly dominated by $A$ thus yielding game $G^3$ shown in the top part of Panel (iv). In the final step, $E$ is deleted because it is strictly dominated by $D$ so that the final output is the strategy profile $(A,D)$. Hence the only rationalizable strategies are $A$ for Player 1 and $D$ for Player 2.

Note that, in the game of Figure 5.9 (i), since the only rationalizable strategy profile is (*A*,*D*), it follows that (*A*,*D*) is also the unique Nash equilibrium.

As noted in Chapter 1 (Section 1.5.1) the significance of the output of the IDSDS procedure is as follows. Consider game *G* in Panel (i) of Figure 5.9. Since, for Player 1, *C* is strictly dominated, if Player 1 is rational she will not play *C*. Thus, if Player 2 believes that Player 1 is rational then he believes that Player 1 will not play *C*, that is, he restricts attention to game $G^1$; since, in $G^1$, *F* is strictly dominated for Player 2, if Player 2 is rational he will not play *F*. It follows that if Player 1 believes that Player 2 is rational and that Player 2 believes that Player 1 is rational, then Player 1 restricts attention to game $G^2$ where rationality requires that Player 1 not play *B*, etc.

**Remark 5.5.** Define a player to be rational if her chosen pure strategy is a best reply to her belief about what the opponent will do. In a two-player game a belief of Player 1 about what Player 2 will do can be expressed as a probability distribution over the set of pure strategies of Player 2; but this is the same object as a mixed strategy of Player 2. Thus, by Theorem 5.3, a rational Player 1 cannot choose a pure strategy that is strictly dominated by one of her own mixed strategies. The iterated reasoning outlined above can be captured by means of the notion of common belief of rationality. Indeed, it will be shown in a later chapter that if the players are rational and there is common belief of rationality then only rationalizable strategy profiles can be played. In a game with more than two players a belief of Player *i* about her opponents is no longer the same object as a mixed-strategy profile of the opponents, because a belief can allow for correlation in the behavior of the opponents, while the notion of mixed-strategy profile rules out such correlation (see Exercise 5.14).

**Remark 5.6.** The iterated reasoning outlined above *requires that the von Neumann-Morgenstern preferences of both players be common knowledge between them*. For example, if Player 2 believes that Player 1 is rational but only knows her ordinal ranking of the outcomes, then Player 2 will not be able to deduce that it is irrational for Player 1 to play *C* and thus it cannot be irrational for him to play *F*. Expecting a player to know the von Neumann-Morgenstern preferences of another player is often (almost always?) very unrealistic! Thus one should be aware of the implicit assumptions that one makes (and one should question the assumptions made by others in their analyses).

Test your understanding of the concepts introduced in this section, by going through the exercises in Section 5.E.4 of Appendix 5.E at the end of this chapter.

# Appendix 5.E: Exercises

## 5.E.1. Exercises for Section 5.1: Strategic-form games with cardinal payoffs

The answers to the following exercises are in Appendix S at the end of this chapter.

**Exercise 5.1.** Consider the following game-frame in strategic form, where $o_1, o_2, o_3$ and $o_4$ are basic outcomes:

Player 2

| | | $c$ | $d$ |
|---|---|---|---|
| Player 1 | $a$ | $o_1$ | $o_2$ |
| | $b$ | $o_3$ | $o_4$ |

Both players satisfy the axioms of expected utility. The best outcome for Player 1 is $o_3$; she is indifferent between outcomes $o_1$ and $o_4$ and ranks them both as worst; she considers $o_2$ to be worse than $o_3$ and better than $o_4$; she is indifferent between $o_2$ with certainty and the lottery $\begin{pmatrix} o_3 & o_1 \\ 0.25 & 0.75 \end{pmatrix}$. The best outcome for Player 2 is $o_4$, which he considers to be just as good as $o_1$; he considers $o_2$ to be worse than $o_1$ and better than $o_3$; he is indifferent between $o_2$ with certainty and the lottery $\begin{pmatrix} o_1 & o_3 \\ 0.4 & 0.6 \end{pmatrix}$.

Find the normalized von Neumann-Morgenstern utility functions for the two players and write the corresponding reduced-form game.

**Exercise 5.2.** Consider the following game-frame, where $o_1,\ldots,o_4$ are basic outcomes.

| | | Player 2 | |
|---|---|---|---|
| | | $C$ | $D$ |
| **Player 1** | A | $\begin{pmatrix} o_1 & o_4 \\ \frac{1}{4} & \frac{3}{4} \end{pmatrix}$ | $\begin{pmatrix} o_1 & o_2 \\ \frac{1}{2} & \frac{1}{2} \end{pmatrix}$ |
| | B | $o_3$ | $\begin{pmatrix} o_3 & o_4 \\ \frac{2}{5} & \frac{3}{5} \end{pmatrix}$ |

Both players have von Neumann-Morgenstern rankings of the basic outcomes. The ranking of Player 1 can be represented by the following von Neumann-Morgenstern utility function

| outcome: | $o_1$ | $o_2$ | $o_3$ | $o_4$ |
|---|---|---|---|---|
| $U_1$ : | 12 | 10 | 6 | 16 |

and the ranking of Player 2 can be represented by the following von Neumann-Morgenstern utility function

| outcome: | $o_1$ | $o_2$ | $o_3$ | $o_4$ |
|---|---|---|---|---|
| $U_2$ : | 6 | 14 | 8 | 10 |

.

Write the corresponding reduced-form game.

## 5.E.2. Exercises for Section 5.2: Mixed strategies.

The answers to the following exercises are in Appendix S at the end of this chapter.

**Exercise 5.3**. Consider the following reduced-form game with cardinal payoffs:

| | | Player 2: *D* | | Player 2: *E* | |
|---|---|---|---|---|---|
| Player 1 | *A* | 0 | 1 | 6 | 3 |
| | *B* | 4 | 4 | 2 | 0 |
| | *C* | 3 | 0 | 4 | 2 |

**(a)** Calculate the player's payoffs from the mixed strategy $\left(\begin{pmatrix} A & B & C \\ \frac{1}{4} & \frac{3}{4} & 0 \end{pmatrix}\begin{pmatrix} D & E \\ \frac{1}{2} & \frac{1}{2} \end{pmatrix}\right)$.

**(b)** Is $\left(\begin{pmatrix} A & B & C \\ \frac{1}{4} & \frac{3}{4} & 0 \end{pmatrix}\begin{pmatrix} D & E \\ \frac{1}{2} & \frac{1}{2} \end{pmatrix}\right)$ a Nash equilibrium?

**Exercise 5.4**. Consider the following reduced-form game with cardinal payoffs:

| | | Player 2: *D* | | Player 2: *E* | |
|---|---|---|---|---|---|
| Player 1 | *A* | 2 | 3 | 8 | 5 |
| | *B* | 6 | 6 | 4 | 2 |

Prove that $\left(\begin{pmatrix} A & B \\ \frac{2}{3} & \frac{1}{3} \end{pmatrix}\begin{pmatrix} D & E \\ \frac{1}{2} & \frac{1}{2} \end{pmatrix}\right)$ is a Nash equilibrium.

## 5.E.3. Exercises for Section 5.3: Computing the mixed-strategy Nash equilibria.

The answers to the following exercises are in Appendix S at the end of this chapter.

**Exercise 5.5**. **(a)** Find the mixed-strategy Nash equilibrium of the game of Exercise 5.1. **(b)** Calculate the payoffs of both players at the Nash equilibrium.

**Exercise 5.6**. Find the Nash equilibrium of the game of Exercise 5.2.

**Exercise 5.7**. Find all the mixed-strategy Nash equilibria of the game of Exercise 5.4. Calculate the payoffs of both players at every Nash equilibrium that you find.

**Exercise 5.8**. Find the mixed-strategy Nash equilibria of the following game:

| Player 1 \ Player 2 | L | R |
|---|---|---|
| T | 1 , 4 | 4 , 3 |
| C | 2 , 0 | 1 , 2 |
| B | 1 , 5 | 0 , 6 |

**Exercise 5.9**. Consider the following two-player game, where $z_1, z_2, \ldots, z_6$ are basic outcomes.

| Player 1 \ Player 2 | *d* | *e* |
|---|---|---|
| *a* | z1 | z2 |
| *b* | z3 | z4 |
| *c* | z5 | z6 |

Player 1 ranks the outcomes as indicated by $A$ below and Player 2 ranks the outcomes as indicated by $B$ below (if outcome $z$ is above outcome $z'$ then $z$ is strictly preferred to $z'$, and if $z$ and $z'$ are written next to each other then the player is indifferent between the two).

$$A = \begin{pmatrix} z_1 \\ z_6 \\ z_4, z_2 \\ z_5 \\ z_3 \end{pmatrix} \qquad B = \begin{pmatrix} z_3, z_4 \\ z_2 \\ z_1, z_5 \\ z_6 \end{pmatrix}$$

**(a)** One player has a strategy that is strictly dominated. Identify the player and the strategy.

[Hint: in order to answer the following questions, you can make your life a lot easier if you simplify the game on the basis of your answer to part (a).]

Player 1 satisfies the axioms of Expected Utility Theory and is indifferent between $z_6$ and the lottery $\begin{pmatrix} z_1 & z_5 \\ \frac{4}{5} & \frac{1}{5} \end{pmatrix}$ and is indifferent between $z_2$ and the lottery $\begin{pmatrix} z_6 & z_5 \\ \frac{1}{2} & \frac{1}{2} \end{pmatrix}$.

**(b)** Suppose that Player 1 believes that Player 2 is going to play $d$ with probability $\frac{1}{2}$ and $e$ with probability $\frac{1}{2}$. Which strategy should he play?

Player 2 satisfies the axioms of Expected Utility Theory and is indifferent between $z_5$ and the lottery $\begin{pmatrix} z_2 & z_6 \\ \frac{1}{4} & \frac{3}{4} \end{pmatrix}$.

**(c)** Suppose that Player 2 believes that Player 1 is going to play $a$ with probability $\frac{1}{4}$ and $c$ with probability $\frac{3}{4}$. Which strategy should she play?

**(d)** Find all the (pure- and mixed-strategy) Nash equilibria of this game.

**Exercise 5.10.** Consider the following game (where the payoffs are von Neumann-Morgenstern payoffs):

**Player 2**

| Player 1 | | $C$ | $D$ |
|---|---|---|---|
| | $A$ | $x$ $y$ | 3 0 |
| | $B$ | 6 2 | 0 4 |

**(a)** Suppose that $x = 2$ and $y = 2$. Find the mixed-strategy Nash equilibrium and calculate the payoffs of both players at the Nash equilibrium.

**(b)** For what values of $x$ and $y$ is $\left( \begin{pmatrix} A & B \\ \frac{1}{5} & \frac{4}{5} \end{pmatrix}, \begin{pmatrix} C & D \\ \frac{3}{4} & \frac{1}{4} \end{pmatrix} \right)$ a Nash equilibrium?

**Exercise 5.11.** Find the mixed-strategy Nash equilibria of the game of Exercise 5.3. Calculate the payoffs of both players at every Nash equilibrium that you find.

## 5.E.4. Exercises for Section 5.4: Strict dominance and rationalizability

The answers to the following exercises are in Appendix S at the end of this chapter.

**Exercise 5.12**. In the following game, for each player, find all the rationalizable pure strategies (that is, apply the cardinal IDSDS procedure).

| | | Player 2 L | M | R |
|---|---|---|---|---|
| Player 1 | A | 3 , 5 | 2 , 0 | 2 , 2 |
| | B | 5 , 2 | 1 , 2 | 2 , 1 |
| | C | 9 , 0 | 1 , 5 | 3 , 2 |

**Note:** The next three exercises are more difficult than the previous ones.

**Exercise 5.13**. Is the following statement true or false? (Either prove that it is true or give a counterexample.)

> "Consider a two-player strategic-form game with cardinal payoffs. Let *A* and *B* be two pure strategies of Player 1. Suppose that both *A* and *B* are rationalizable (that is, they survive the cardinal IDSDS procedure). Then any mixed strategy that attaches positive probability to both *A* and *B* and zero to every other strategy is a best reply to some mixed strategy of Player 2."

**Exercise 5.14**. Consider the following three-player game, where only the payoffs of Player 1 are shown.

Player 3: G

| | | Player 2 E | F |
|---|---|---|---|
| Player 1 | A | 3 | 0 |
| | B | 0 | 3 |
| | C | 0 | 0 |
| | D | 2 | 0 |

Player 3: H

| | | Player 2 E | F |
|---|---|---|---|
| Player 1 | A | 0 | 0 |
| | B | 3 | 0 |
| | C | 0 | 3 |
| | D | 0 | 2 |

**(a)** Show that if Player 1 assigns probability $\frac{1}{2}$ to the event "Player 2 will play *E* and Player 3 will play *G*" and probability $\frac{1}{2}$ to the event "Player 2 will play *F* and Player will play *H*", then playing *D* is a best reply.

Next we want to show that there is no mixed-strategy profile $\sigma_{-1} = \left( \begin{pmatrix} E & F \\ p & 1-p \end{pmatrix}, \begin{pmatrix} G & H \\ q & 1-q \end{pmatrix} \right)$ of Players 2 and 3 against which $D$ is a best reply for Player 1. We do this in steps. First define the following functions: $A(p,q) = \Pi_1(A, \sigma_{-1})$ (that is, $A(p,q)$ is Player 1's expected payoff if she plays the pure strategy $A$ against $\sigma_{-1}$), $B(p,q) = \Pi_1(B, \sigma_{-1})$, $C(p,q) = \Pi_1(C, \sigma_{-1})$ and $D(p,q) = \Pi_1(D, \sigma_{-1})$.

**(b)** In the $(p,q)$ plane (with $0 \le p \le 1$ and $0 \le q \le 1$) draw the curve corresponding to the equation $A(p,q) = D(p,q)$ and identify the region where $A(p,q) > D(p,q)$.

**(c)** In the $(p,q)$ plane draw the curve corresponding to the equation $C(p,q) = D(p,q)$ and identify the region where $C(p,q) > D(p,q)$.

**(d)** In the $(p,q)$ plane draw the two curves corresponding to the equation $B(p,q) = D(p,q)$ and identify the region where $B(p,q) > D(p,q)$.

**(e)** Infer from parts (b)-(c) that there is no mixed-strategy profile of Players 2 and 3 against which $D$ is a best reply for Player 1.

◊◊◊◊◊◊◊◊◊◊◊◊

**Exercise 5.15: Challenging Question**. A team of $n$ professional swimmers ($n \geq 2$) – from now on called players – are partying on the bank of the Sacramento river on a cold day in January. Suddenly a passerby shouts "Help! My dog fell into the water!". Each of the swimmers has to decide whether or not to jump into the icy cold water to rescue the dog. One rescuer is sufficient: the dog will be saved if at least one player jumps into the water; if nobody does, then the dog will die. Each player prefers somebody else to jump in, but each player prefers to jump in himself if nobody else does. Let us formulate this as a game. The strategy set of each player $i = 1,...,n$ is $S_i = \{J, \neg J\}$, where $J$ stands for 'jump in' and $\neg J$ for 'not jump in'. The possible basic outcomes can be expressed as subsets of the set $I = \{1,...,n\}$ of players: outcome $N \subseteq I$ is interpreted as 'the players in the set $N$ jump into the water'; if $N = \varnothing$ the dog dies, while if $N \neq \varnothing$ the dog is saved. Player $i$ has the following ordinal ranking of the outcomes: (1) $N \sim N'$, for every $N \neq \varnothing$, $N' \neq \varnothing$ with $i \notin N$ and $i \notin N'$, (2) $N \succ N'$ for every $N \neq \varnothing$, $N' \neq \varnothing$ with $i \notin N$ and $i \in N'$, (3) $\{i\} \succ \varnothing$.

**(a)** Find all the pure-strategy Nash equilibria.

**(b)** Suppose that each player $i$ has the following von Neumann-Morgenstern payoff function (which is consistent with the above ordinal ranking):

$$\pi_i(N) = \begin{cases} v & \text{if } N \neq \varnothing \text{ and } i \notin N \\ v - c & \text{if } N \neq \varnothing \text{ and } i \in N \\ 0 & \text{if } N = \varnothing \end{cases} \quad \text{with } 0 < c < v.$$

Find the symmetric mixed-strategy Nash equilibrium (symmetric means that all the players use the same mixed strategy).

**(c)** Assuming that the players behave according to the symmetric mixed-strategy Nash equilibrium of part (b), is it better for the dog if $n$ (the number of players) is large or if $n$ is small? Calculate the probability that the dog is saved at the mixed-strategy Nash equilibrium as a function of $n$, for all possible values of $c$ and $v$ (subject to $0 < c < v$), and plot it for the case where $c = 10$ and $v = 12$.

# Appendix 5.S: Solutions to exercises

**Exercise 5.1.** The normalized von Neumann-Morgenstern utility functions are:

Player 1:

| outcome | $U_1$ |
|---|---|
| $o_3$ | 1 |
| $o_2$ | 0.25 |
| $o_1, o_4$ | 0 |

Player 2:

| outcome | $U_2$ |
|---|---|
| $o_1, o_4$ | 1 |
| $o_2$ | 0.4 |
| $o_3$ | 0 |

Thus the reduced-form game is as follows:

| Player 1 \ Player 2 | c | d |
|---|---|---|
| a | 0, 1 | 0.25, 0.4 |
| b | 1, 0 | 0, 1 |

**Exercise 5.2.** The expected utility of the lottery $\begin{pmatrix} o_1 & o_4 \\ \frac{1}{4} & \frac{3}{4} \end{pmatrix}$ is $\frac{1}{4}(12)+\frac{3}{4}(16)=15$ for Player 1 and $\frac{1}{4}(6)+\frac{3}{4}(10)=9$ for Player 2. The expected utility of the lottery $\begin{pmatrix} o_1 & o_2 \\ \frac{1}{2} & \frac{1}{2} \end{pmatrix}$ is 11 for Player 1 and 10 for Player 2. The expected utility of the lottery $\begin{pmatrix} o_3 & o_4 \\ \frac{2}{5} & \frac{3}{5} \end{pmatrix}$ is 12 for Player 1 and 9.2 for Player 2. Thus the reduced-form game is as follows.

| Player 1 \ Player 2 | C | D |
|---|---|---|
| A | 15, 9 | 11, 10 |
| B | 6, 8 | 12, 9.2 |

**Exercise 5.3.** **(a)** $\Pi_1 = \frac{1}{8}0+\frac{1}{8}6+\frac{3}{8}4+\frac{3}{8}2=3$ and $\Pi_2 = \frac{1}{8}1+\frac{1}{8}3+\frac{3}{8}4+\frac{3}{8}0=2$.

**(b)** No, because if Player 1 switched to the pure strategy $C$ then her payoff would be $\frac{1}{2}3+\frac{1}{2}4=3.5 > 2$.

**Exercise 5.4.** Player 1's payoff is $\Pi_1 = \frac{2}{6}2 + \frac{2}{6}8 + \frac{1}{6}6 + \frac{1}{6}4 = 5$. If Player 1 switches to any other mixed strategy $\begin{pmatrix} A & B \\ p & 1-p \end{pmatrix}$, while Player 2's strategy is kept fixed at $\begin{pmatrix} C & D \\ \frac{1}{2} & \frac{1}{2} \end{pmatrix}$, then her payoff is $\Pi_1 = \frac{1}{2}p2 + \frac{1}{2}p8 + \frac{1}{2}(1-p)6 + \frac{1}{2}(1-p)4 = 5$. Thus any mixed strategy of Player 1 is a best response to $\begin{pmatrix} C & D \\ \frac{1}{2} & \frac{1}{2} \end{pmatrix}$. Similarly, Player 2's payoff is $\Pi_2 = \frac{2}{6}3 + \frac{2}{6}5 + \frac{1}{6}6 + \frac{1}{6}2 = 4$. If Player 2 switches to any other mixed strategy $\begin{pmatrix} C & D \\ q & 1-q \end{pmatrix}$, while Player 1's strategy is kept fixed at $\begin{pmatrix} A & B \\ \frac{2}{3} & \frac{1}{3} \end{pmatrix}$, then her payoff is $\Pi_2 = \frac{2}{3}q3 + \frac{2}{3}(1-q)5 + \frac{1}{3}q6 + \frac{1}{3}(1-q)2 = 4$. Thus any mixed strategy of Player 2 is a best response to $\begin{pmatrix} A & B \\ \frac{2}{3} & \frac{1}{3} \end{pmatrix}$. Hence $\begin{pmatrix} A & B \\ \frac{2}{3} & \frac{1}{3} \end{pmatrix}$ is a best reply to $\begin{pmatrix} C & D \\ \frac{1}{2} & \frac{1}{2} \end{pmatrix}$ and $\begin{pmatrix} C & D \\ \frac{1}{2} & \frac{1}{2} \end{pmatrix}$ is a best reply to $\begin{pmatrix} A & B \\ \frac{2}{3} & \frac{1}{3} \end{pmatrix}$, that is, $\left(\begin{pmatrix} A & B \\ \frac{2}{3} & \frac{1}{3} \end{pmatrix}, \begin{pmatrix} C & D \\ \frac{1}{2} & \frac{1}{2} \end{pmatrix}\right)$ is a Nash equilibrium.

**Exercise 5.5.** **(a)** We have to find the Nash equilibrium of the following game.

| | | Player 2 | |
|---|---|---|---|
| | | *c* | *d* |
| Player 1 | *a* | 0, 1 | 0.25, 0.4 |
| | *b* | 1, 0 | 0, 1 |

To make calculations easier, let us multiply all the payoffs by 100 (that is, we re-scale the von Neumann-Morgenstern utility functions by a factor of 100):

| | | Player 2 | |
|---|---|---|---|
| | | *c* | *d* |
| Player 1 | *a* | 0, 100 | 25, 40 |
| | *b* | 100, 0 | 0, 100 |

There are no pure-strategy Nash equilibria. To find the mixed-strategy Nash equilibrium, let $p$ be the probability with which Player 1 chooses $a$ and $q$ be the probability with which Player 2 chooses $c$. Then, for Player 1, the payoff from playing $a$ against $\begin{pmatrix} c & d \\ q & 1-q \end{pmatrix}$ must be equal to the payoff from playing $b$ against $\begin{pmatrix} c & d \\ q & 1-q \end{pmatrix}$. That is, it must be that $25(1-q)=100q$, which yields $q=\frac{1}{5}$. Similarly, for Player 2, the payoff from playing $c$ against $\begin{pmatrix} a & b \\ p & 1-p \end{pmatrix}$ must be equal to the payoff from playing $d$ against $\begin{pmatrix} a & b \\ p & 1-p \end{pmatrix}$. This requires $100p = 40p + 100(1-p)$, that is, $p=\frac{5}{8}$. Thus the Nash equilibrium is $\left(\begin{pmatrix} a & b \\ \frac{5}{8} & \frac{3}{8} \end{pmatrix}, \begin{pmatrix} c & d \\ \frac{1}{5} & \frac{4}{5} \end{pmatrix}\right)$. **(b)** At the Nash equilibrium the payoffs are 20 for Player 1 and 62.5 for Player 2. (If you worked with the original payoffs, then the Nash equilibrium payoffs would be 0.2 for Player 1 and 0.625 for Player 2.)

**Exercise 5.6.** We have to find the Nash equilibria of the following game.

| | | Player 2 | |
|---|---|---|---|
| | | *C* | *D* |
| Player 1 | *A* | 15 9 | 11 10 |
| | *B* | 6 8 | 12 9.2 |

For Player 2 *D* is a strictly dominant strategy, thus at a Nash equilibrium Player 2 must play *D* with probability 1. For Player 1, the unique best reply to *D* is *B*. Thus the pure-strategy profile (*B*,*D*) is the only Nash equilibrium.

**Exercise 5.7.** We have to find the Nash equilibria of the following game.

| | | Playe 2 | |
|---|---|---|---|
| | | *D* | *E* |
| Player 1 | *A* | 2 3 | 8 5 |
| | *B* | 6 6 | 4 2 |

(*B*,*D*) (with payoffs (6,6)) and (*A*,*E*) (with payoffs (8,5)) are both Nash equilibria. To see if there is also a mixed-strategy equilibrium we need to solve the following equations, where $p$ is the probability of *A* and $q$ is the probability of *D*: $2q+8(1-q)=6q+4(1-q)$ and $3p+6(1-p)=5p+2(1-p)$. The solution is $p=\frac{2}{3}$ and $q=\frac{1}{2}$ and indeed we verified in Exercise 5.4 that $\left(\begin{pmatrix} A & B \\ \frac{2}{3} & \frac{1}{3} \end{pmatrix}, \begin{pmatrix} C & D \\ \frac{1}{2} & \frac{1}{2} \end{pmatrix}\right)$ is a Nash equilibrium. The payoffs at this Nash equilibrium are 5 for Player 1 and 4 for Player 2.

**Exercise 5.8.** Since $B$ is strictly dominated (by $C$), it cannot be assigned positive probability at a Nash equilibrium. Let $p$ be the probability of $T$ and $q$ the probability of $L$. Then $p$ must be such that $4p+0(1-p) = 3p+2(1-p)$ and $q$ must be such that $q+4(1-q)=2q+(1-q)$. Thus $p=\frac{2}{3}$ and $q=\frac{3}{4}$. Hence there is only one mixed-strategy equilibrium, namely $\left(\begin{pmatrix} T & C & B \\ \frac{2}{3} & \frac{1}{3} & 0 \end{pmatrix}, \begin{pmatrix} L & R \\ \frac{3}{4} & \frac{1}{4} \end{pmatrix}\right)$.

**Exercise 5.9.** **(a)** Since Player 1 prefers $z_5$ to $z_3$ and prefers $z_6$ to $z_4$, strategy $b$ is strictly dominated by strategy c. Thus, at a Nash equilibrium, Player 1 will not play $b$ with positive probability and we can simplify the game to

Player 2

| | | $d$ | $e$ |
|---|---|---|---|
| Player 1 | $a$ | z1 | z2 |
| | $c$ | z5 | z6 |

**(b)** Of the remaining outcomes, for Player 1 $z_1$ is the best (we can assign utility 1 to it) and $z_5$ is the worst (we can assign utility 0 to it). Since he is indifferent between $z_6$ and the lottery $\begin{pmatrix} z_1 & z_5 \\ \frac{4}{5} & \frac{1}{5} \end{pmatrix}$, the utility of $z_6$ is $\frac{4}{5}$. Hence the expected utility of $\begin{pmatrix} z_6 & z_5 \\ \frac{1}{2} & \frac{1}{2} \end{pmatrix}$ is $\frac{1}{2}(\frac{4}{5})+\frac{1}{2}(0)=\frac{2}{5}$ and thus the utility of $z_2$ is also $\frac{2}{5}$. So, if Player 2 is going to play $d$ with probability $\frac{1}{2}$ and $e$ with probability $\frac{1}{2}$, then for Player 1 playing $a$ gives a payoff of $\frac{1}{2}(1)+\frac{1}{2}(\frac{2}{5})=\frac{7}{10}$, while playing $c$ gives a payoff of $\frac{1}{2}(0)+\frac{1}{2}(\frac{4}{5})=\frac{4}{10}$. Hence he should play $a$.

[If you did not take the hint to simplify the analysis as was done above, then you can still reach the same answer, although in a lengthier way. You would still set $U(z_1)=1$. Then the expected payoff from playing $a$ is

$$\Pi_1(a)=\tfrac{1}{2}U(z_1)+\tfrac{1}{2}U(z_2)=\tfrac{1}{2}+\tfrac{1}{2}U(z_2) \qquad (*)$$

Since $z_2$ is as good as $\begin{pmatrix} z_6 & z_5 \\ \frac{1}{2} & \frac{1}{2} \end{pmatrix}$,

$$U(z_2)=\tfrac{1}{2}U(z_6)+\tfrac{1}{2}U(z_5) \qquad (**)$$

Since $z_6$ is as good as $\begin{pmatrix} z_1 & z_5 \\ \frac{4}{5} & \frac{1}{5} \end{pmatrix}$,

$$U(z_6) = \tfrac{4}{5} + \tfrac{1}{5}U(z_5) \qquad (***)$$

Replacing (***) in (**) we get $U(z_2) = \tfrac{2}{5} + \tfrac{3}{5}U(z_5)$ and replacing this expression in (*) we get $\Pi_1(a) = \tfrac{7}{10} + \tfrac{3}{10}U(z_5)$. Similarly, $\Pi_1(c) = \tfrac{1}{2}U(z_5) + \tfrac{1}{2}U(z_6) = \tfrac{1}{2}U(z_5) + \tfrac{1}{2}\left(\tfrac{4}{5} + \tfrac{1}{5}U(z_5)\right) = \tfrac{4}{10} + \tfrac{6}{10}U(z_5)$. Now, $\Pi_1(a) > \Pi_1(c)$ if and only if $\tfrac{7}{10} + \tfrac{3}{10}U(z_5) > \tfrac{4}{10} + \tfrac{6}{10}U(z_5)$ if and only if $3 > 3U(z_5)$ if and only if $U(z_5) < 1$ which is the case because $z_5$ is worse than $z_1$ and $U(z_1) = 1$. Similar steps would be taken to answer (c) and (d).]

**(c)** In the reduced game, for Player 2 $z_2$ is the best outcome (we can assign utility 1 to it) and $z_6$ is the worst (we can assign utility 0 to it). Thus, since she is indifferent between $z_5$ and the lottery $\begin{pmatrix} z_2 & z_6 \\ \tfrac{1}{4} & \tfrac{3}{4} \end{pmatrix}$, the utility of $z_5$ is $\tfrac{1}{4}$ and so is the utility of $z_1$. Thus playing *d* gives an expected payoff of $\tfrac{1}{4}(\tfrac{1}{4}) + \tfrac{3}{4}(\tfrac{1}{4}) = \tfrac{1}{4}$, while playing *e* gives an expected utility of $\tfrac{1}{4}(1) + \tfrac{3}{4}(0) = \tfrac{1}{4}$. Thus she is indifferent between playing *d* and playing *e* (and any mixture of *d* and *e*).

**(d)** Using the calculations of parts (b) and (c) the game is as follows:

| | $d$ | $e$ |
|---|---|---|
| $a$ | $1, \tfrac{1}{4}$ | $\tfrac{2}{5}, 1$ |
| $c$ | $0, \tfrac{1}{4}$ | $\tfrac{4}{5}, 0$ |

There is no pure-strategy Nash equilibrium. At a mixed-strategy Nash equilibrium, each player must be indifferent between his/her two strategies. From part (c) we already know that Player 2 is indifferent if Player 1 plays *a* with probability $\tfrac{1}{4}$ and *c* with probability $\tfrac{3}{4}$. Now let *q* be the probability with which Player 2 plays *d*. Then we need $q + \tfrac{2}{5}(1-q) = \tfrac{4}{5}(1-q)$, hence $q = \tfrac{2}{7}$. Thus the Nash equilibrium is $\left(\begin{array}{ccc|cc} a & b & c & d & e \\ \tfrac{1}{4} & 0 & \tfrac{3}{4} & \tfrac{2}{7} & \tfrac{5}{7} \end{array}\right)$ which can be written more succinctly as $\left(\begin{array}{cc|cc} a & c & d & e \\ \tfrac{1}{4} & \tfrac{3}{4} & \tfrac{2}{7} & \tfrac{5}{7} \end{array}\right)$.

**Exercise 5.10.** **(a)** Let $p$ be the probability of $A$ and $q$ the probability of $B$. Then, Player 1 must be indifferent between playing $A$ and playing $B$: $2q+3(1-q)=6q$. This gives $q=\frac{3}{7}$. Similarly, Player 2 must be indifferent between playing $C$ and playing $D$: $2=4(1-p)$. This gives $p=\frac{1}{2}$. Thus the Nash equilibrium is given by $\left(\begin{pmatrix} A & B \\ \frac{1}{2} & \frac{1}{2} \end{pmatrix},\begin{pmatrix} C & D \\ \frac{3}{7} & \frac{4}{7} \end{pmatrix}\right)$. The equilibrium payoffs are $\frac{18}{7}$ = 2.57 for Player 1 and 2 for Player 2.

**(b)** Player 1 must be indifferent between playing $A$ and playing $B$: $\frac{3}{4}x+\frac{1}{4}3=\frac{3}{4}6$. Thus $x$ = 5. Similarly, Player 2 must be indifferent between playing $C$ and playing $D$: $\frac{1}{5}y+\frac{4}{5}2=\frac{4}{5}4$. Thus $y$ = 8.

**Exercise 5.11.** We have to find the Nash equilibria of the following game.

| | | Player 2: D | | E | |
|---|---|---|---|---|---|
| Player 1 | A | 0 | 1 | 6 | 3 |
| | B | 4 | 4 | 2 | 0 |
| | C | 3 | 0 | 4 | 2 |

There are two pure-strategy equilibria, namely ($B$,$D$) and ($A$,$E$). To see if there is a mixed-strategy equilibrium we calculate the best response of Player 1 to every possible mixed strategy $\begin{pmatrix} D & E \\ q & 1-q \end{pmatrix}$ of Player 2 (with $q\in[0,1]$). For Player 1 the payoff from playing $A$ against $\begin{pmatrix} D & E \\ q & 1-q \end{pmatrix}$ is $6-6q$, the payoff from playing $B$ is $4q+2(1-q)=2+2q$ and the payoff from playing $C$ is $3q+4(1-q)=4-q$. These functions are shown in the following diagram.

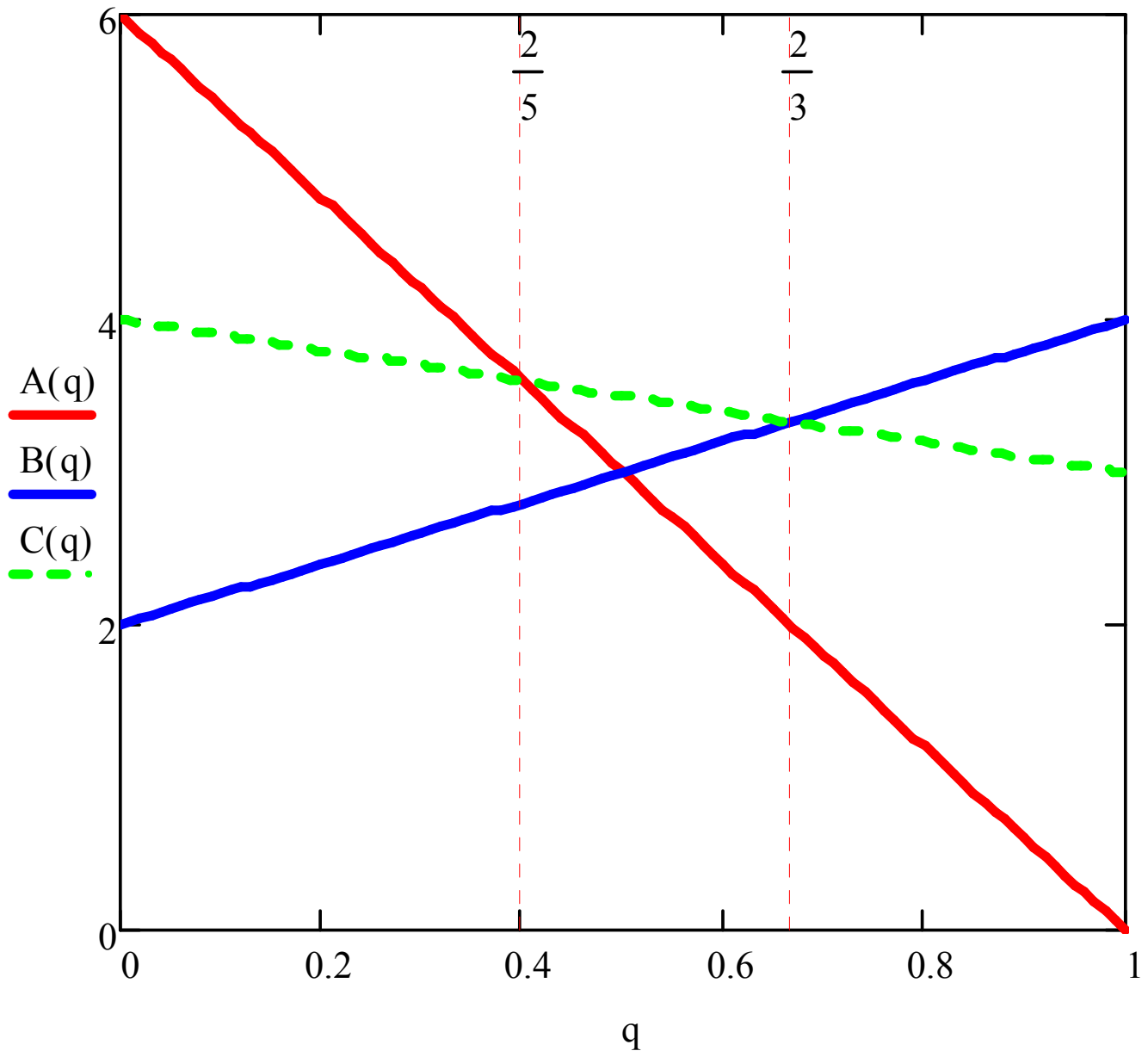

The red downward-sloping line plots the function where $A(q) = 6 - 6q$, the blue upward-sloping line plots the function $B(q) = 2 + 2q$ and the green dotted line the function $C(q) = 4 - q$. It can be seen from the diagram that

$$\text{Player1's best reply} = \begin{cases} A & \text{if } 0 \le q < \frac{2}{5} \\ \begin{pmatrix} A & C \\ p & 1-p \end{pmatrix} \text{ for any } p \in [0,1] & \text{if } q = \frac{2}{5} \\ C & \text{if } \frac{2}{5} < q < \frac{2}{3} \\ \begin{pmatrix} B & C \\ p & 1-p \end{pmatrix} \text{ for any } p \in [0,1] & \text{if } q = \frac{2}{3} \\ B & \text{if } \frac{2}{3} < q \le 1 \end{cases}$$

Thus if there is a mixed-strategy equilibrium it is either of the form $\left( \begin{pmatrix} A & C \\ p & 1-p \end{pmatrix}, \begin{pmatrix} D & E \\ \frac{2}{5} & \frac{3}{5} \end{pmatrix} \right)$ or of the form $\left( \begin{pmatrix} B & C \\ p & 1-p \end{pmatrix}, \begin{pmatrix} D & E \\ \frac{2}{3} & \frac{1}{3} \end{pmatrix} \right)$. In the first case, where Player 1 chooses *B* with probability zero, *E* strictly dominates *D* for Player 2 and thus $\begin{pmatrix} D & E \\ \frac{2}{5} & \frac{3}{5} \end{pmatrix}$ is not a best reply for Player 2, so that $\left( \begin{pmatrix} A & C \\ p & 1-p \end{pmatrix}, \begin{pmatrix} D & E \\ \frac{2}{5} & \frac{3}{5} \end{pmatrix} \right)$ is not a Nash equilibrium for any *p*. In the second case we need $D(p) = E(p)$, that is, $4p = 2(1-p)$, which yields $p = \frac{1}{3}$. Thus the mixed-strategy Nash equilibrium is $\left(\begin{pmatrix} B & C \\ \frac{1}{3} & \frac{2}{3} \end{pmatrix}, \begin{pmatrix} D & E \\ \frac{2}{3} & \frac{1}{3} \end{pmatrix}\right)$ with payoffs of $\frac{10}{3}$ for Player 1 and $\frac{4}{3}$ for Player 2.

**Exercise 5.12.** For Player 1, $B$ is strictly dominated by $\begin{pmatrix} A & C \\ \frac{1}{2} & \frac{1}{2} \end{pmatrix}$; for Player 2, $R$ is strictly dominated by $\begin{pmatrix} L & M \\ \frac{1}{2} & \frac{1}{2} \end{pmatrix}$. Eliminating $B$ and $R$ we are left with

| | | Player 2: L | M |
|---|---|---|---|
| Player 1 | A | 3 , 5 | 2 , 0 |
| | C | 9 , 0 | 1 , 5 |

In this game no player has a strictly dominated strategy. Thus for Player 1 both $A$ and $C$ are rationalizable and for Player 2 both $L$ and $M$ are rationalizable.

**Exercise 5.13.** The statement is false. Consider, for example, the following game:

| | | Player 2: L | R |
|---|---|---|---|
| Player 1 | A | 3 , 1 | 0 , 0 |
| | B | 0 , 0 | 3 , 1 |
| | C | 2 , 1 | 2 , 1 |

Here both $A$ and $B$ are rationalizable (indeed, they are both part of a Nash equilibrium; note that the cardinal IDSDS procedure leaves the game unchanged). However, the mixture $\begin{pmatrix} A & B \\ \frac{1}{2} & \frac{1}{2} \end{pmatrix}$ (which gives Player 1 a payoff of 1.5, no matter what Player 2 does) *cannot* be a best reply to any mixed strategy of Player 2, since it is strictly dominated by $C$.

**Exercise 5.14.** **(a)** if Player 1 assigns probability $\frac{1}{2}$ to the event "Player 2 will play $E$ and Player 3 will play $G$" and probability $\frac{1}{2}$ to the event "Player 2 will play $F$ and Player will play $H$", then $A$ gives Player 1 an expected payoff of 1.5, $B$ an expected payoff of 0, $C$ and expected payoff of 1.5 and $D$ an expected payoff of 2. Thus $D$ is a best reply to those beliefs.

The functions are as follows: $A(p,q) = 3pq$, $B(p,q) = 3(1-p)q + 3p(1-q)$, $C(p,q) = 3(1-p)(1-q)$, $D(p,q) = 2pq + 2(1-p)(1-q)$.

**(b)** $A(p,q) = D(p,q)$ at those points $(p,q)$ such that $q = \dfrac{2-2p}{2-p}$. The set of such points is the red curve in the diagram below. The region where $A(p,q) > D(p,q)$ is the region **above** the red curve.

**(c)** $C(p,q) = D(p,q)$ at those points $(p,q)$ such that $q = \dfrac{1-p}{1+p}$. The set of such points is the blue curve in the diagram below. The region where $C(p,q) > D(p,q)$ is the region **below** the blue curve.

**(c)** $B(p,q) = D(p,q)$ at those points $(p,q)$ such that $q = \dfrac{2-5p}{5-10p}$ (for $p \neq \frac{1}{2}$). The set of such points is given by the two green curves in the diagram below. The region where $B(p,q) > D(p,q)$ is the region **between** the two green curves.

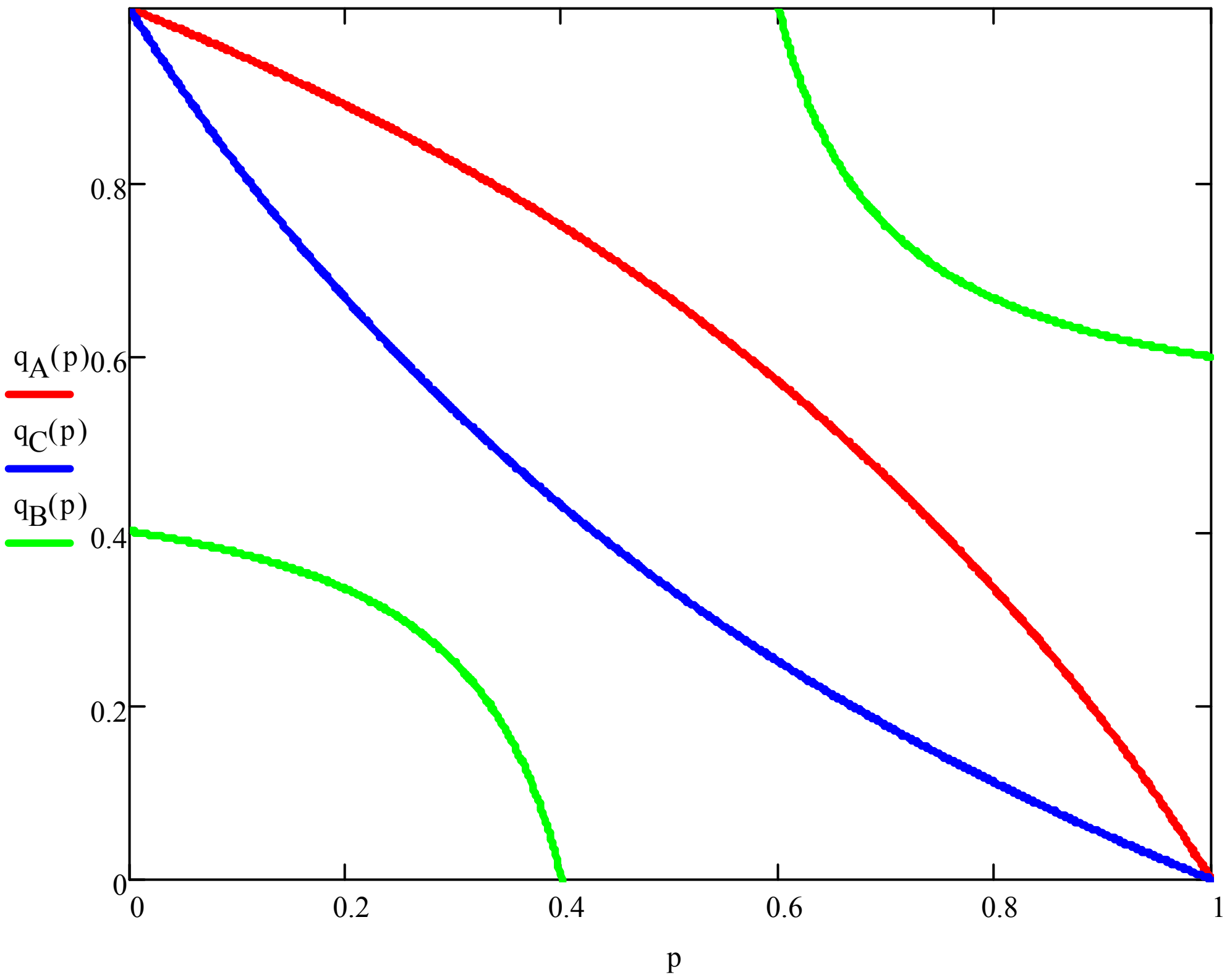

Thus in the region strictly above the red curve $A$ is better than $D$, in the region strictly below the blue curve $C$ is better than $D$ and in the region on and between the red and blue curves $B$ is better that $D$. Hence, at every point in the $(p,q)$ square there is a pure strategy of Player 1 which is strictly better than $D$. It follows that there is no mixed-strategy $\sigma_{-1}$ against which $D$ is a best reply.

**Exercise 5.15.** **(a)** There are $n$ pure-strategy Nash equilibria: at each equilibrium exactly one player jumps in.

**(b)** Let $p$ be the probability with which each player jumps into the water. Consider a Player $i$. The probability that none of the other players jump in is $(1-p)^{n-1}$ and thus the probability that somebody else jumps in is $\left[1-(1-p)^{n-1}\right]$. Player $i$'s payoff if he jumps in is $v-c$, while his expected payoff if he does not jump in is $v\left[1-(1-p)^{n-1}\right]+0(1-p)^{n-1}=v\left[1-(1-p)^{n-1}\right]$. Thus we need $v-c=v\left[1-(1-p)^{n-1}\right]$, that is, $\boxed{p=1-\left(\frac{c}{v}\right)^{\frac{1}{n-1}}}$, which is strictly between 0 and 1 because $c < v$.

**(c)** At the Nash equilibrium the probability that nobody jumps in is $(1-p)^n = \left(\frac{c}{v}\right)^{\frac{n}{n-1}}$; thus this is the probability that the dog dies. Hence the dog is rescued with the remaining probability $1-\left(\frac{c}{v}\right)^{\frac{n}{n-1}}$. This is a decreasing function of *n*. Thus the larger the number of swimmers who are present, the more likely it is that the dog dies. The plot of this function when $c$ = 10 and $v$ = 12 is shown below.

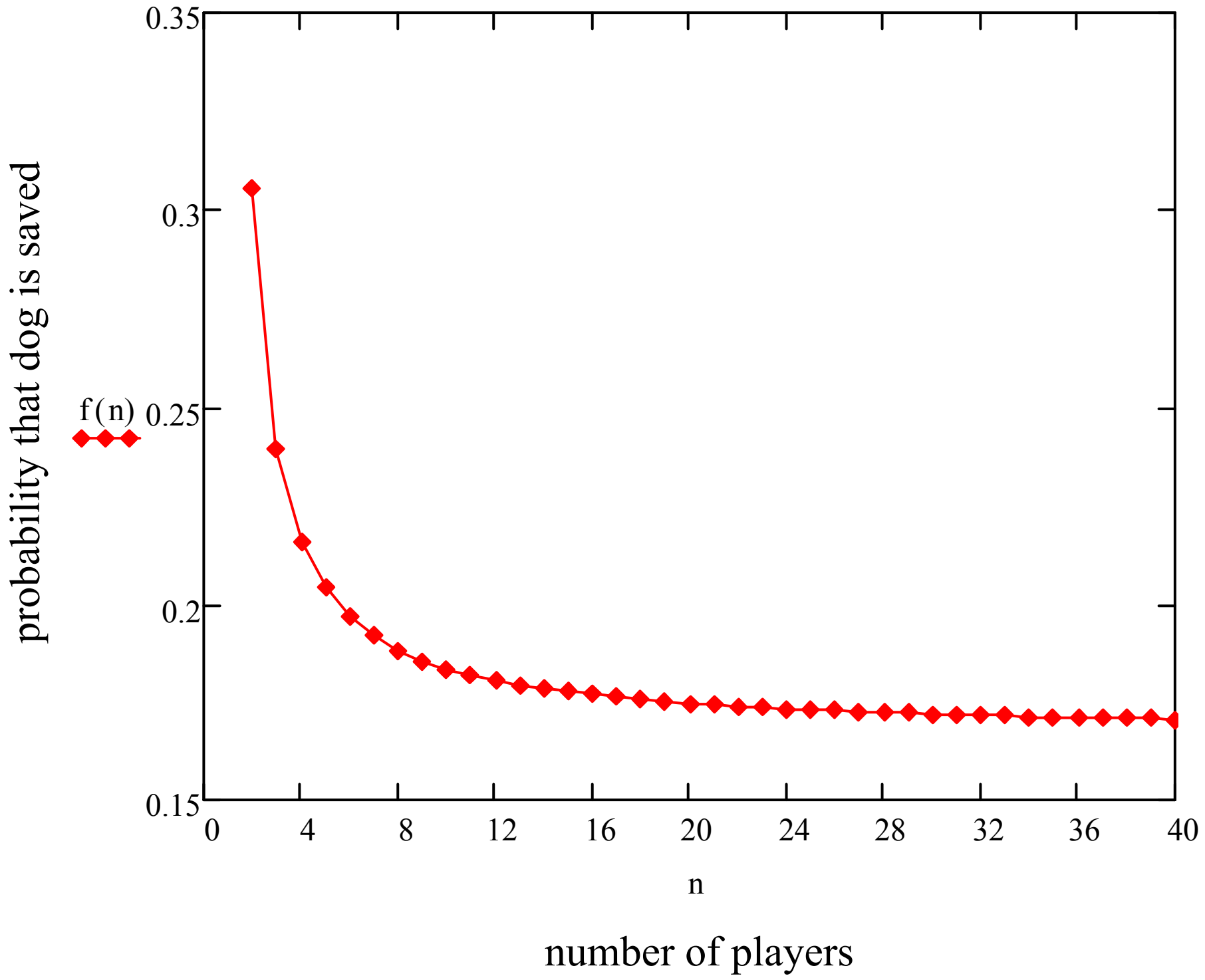

Chapter 6

# Dynamic Games with Cardinal Payoffs

## 6.1 Behavioral strategies in dynamic games

The definition of dynamic (or extensive-form) game with cardinal payoffs is just like the definition of extensive form game with ordinal payoffs (Definition 3.1, Chapter 3), the only difference being that we postulate von Neumann-Morgenstern preferences instead of merely ordinal preferences.

In Chapter 5 we generalized the notion of strategic-form frame by allowing for lotteries (rather than just simple outcomes) to be associated with strategy profiles. One could do the same for extensive-form frames, as shown in Figure 6.1.

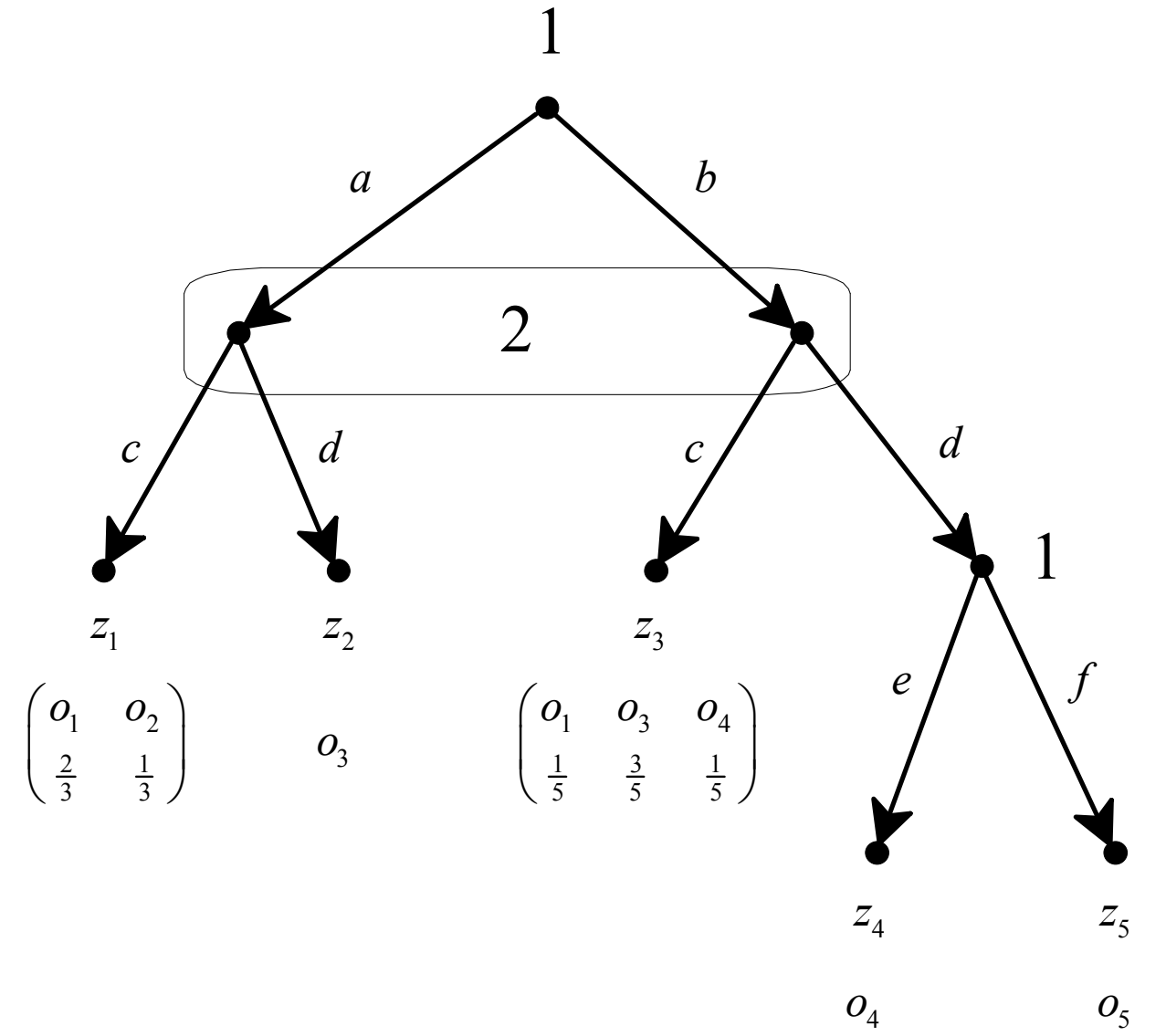

**Figure 6.1**
An extensive frame with probabilistic outcomes.
The $z_i$'s are terminal nodes and the $o_i$'s are basic outcomes.

In Figure 6.1 $\{z_1, z_2, ..., z_5\}$ is the set of terminal nodes and $\{o_1, o_2, ..., o_5\}$ is the set of basic outcomes. Associated with $z_1$ is the lottery $\begin{pmatrix} o_1 & o_2 \\ \frac{2}{3} & \frac{1}{3} \end{pmatrix}$, while the lottery associated with $z_3$ is $\begin{pmatrix} o_1 & o_3 & o_4 \\ \frac{1}{5} & \frac{3}{5} & \frac{1}{5} \end{pmatrix}$, etc.

However, as we saw at the end of Chapter 3, in extensive forms one can explicitly represent random events by means of chance moves (also called moves of Nature). Thus an alternative representation of the extensive-form frame of Figure 6.1 is the extensive form shown in Figure 6.2.

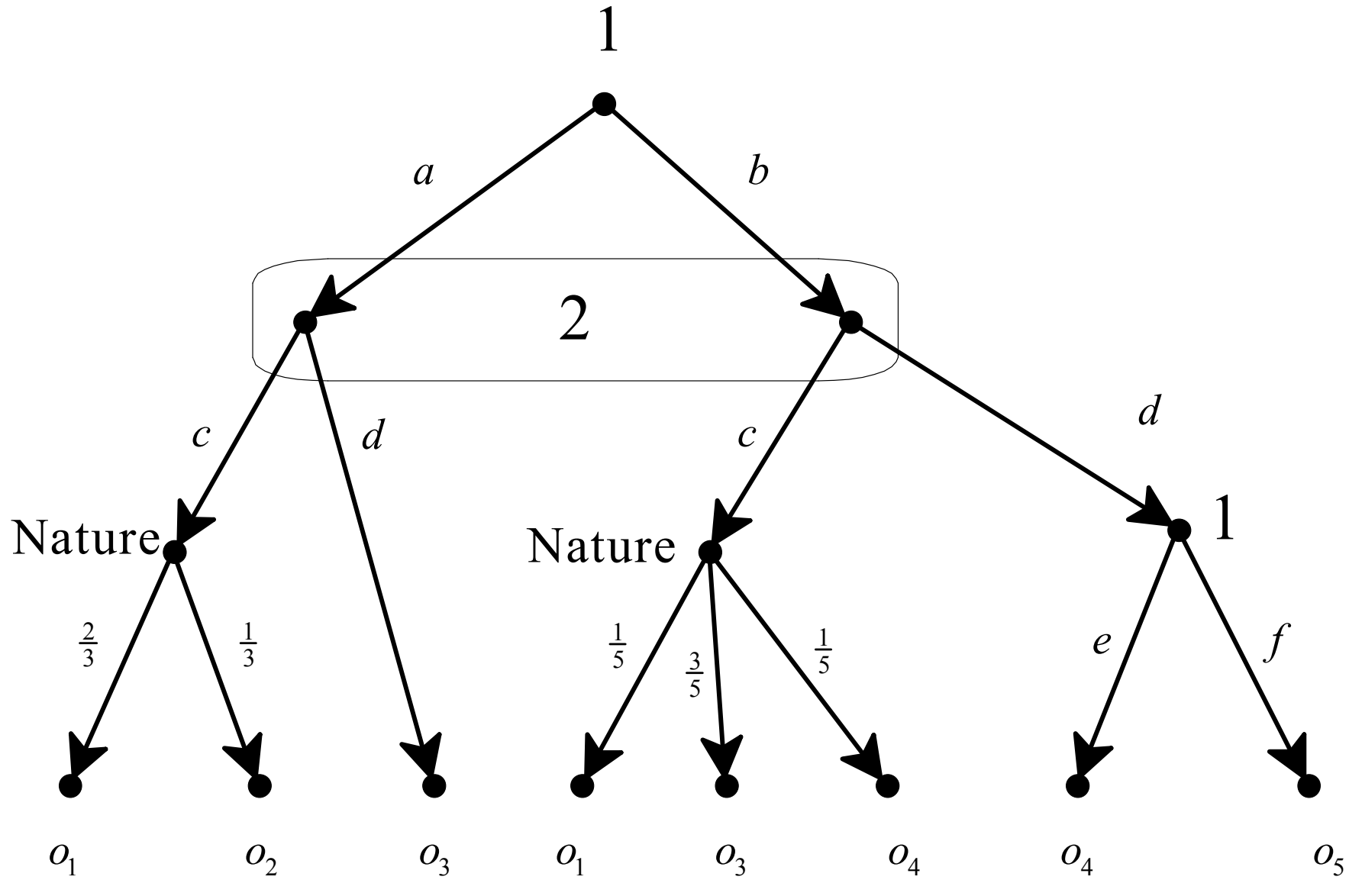

**Figure 6.2**
An alternative representation of the extensive frame of Figure 6.1.
The terminal nodes have not been labeled. The $o_i$'s are basic outcomes.

Hence we can continue to use the definition of extensive-form frame given in Chapter 3, but from now on we will allow for the possibility of chance moves.

The notion of strategy remains, of course, unchanged: a *strategy* for a player is a list of choices, one for every information set of that player (Definition 3.4, Chapter 3). For example, the set of strategies for Player 1 in the extensive frame of Figure 6.2 is $S_1 = \{(a,e),(a,f),(b,e),(b,f)\}$. Thus mixed strategies can easily be introduced also in extensive frames. For example, in the extensive frame of Figure 6.2, the set of mixed strategies for Player 1, $\Sigma_1$, is the set of probability distributions over $S_1$:

$$\Sigma_1 = \left\{ \begin{pmatrix} (a,e) & (a,f) & (b,e) & (b,f) \\ p & q & r & 1-p-q-r \end{pmatrix} : p,q,r \in [0,1] \text{ and } p+q+r \leq 1 \right\}.$$

However, it turns out that in extensive forms with perfect recall[9] one can use simpler objects than mixed strategies, namely behavioral strategies.

**Definition 6.1.** A *behavioral strategy* for a player in an extensive form is a list of probability distributions, one for every information set of that player; each probability distribution is over the set of choices at the corresponding information set.

For example, the set of behavioral strategies for Player 1 in the extensive frame of Figure 6.2 is $\left\{ \left( \begin{array}{cc|cc} a & b & e & f \\ p & 1-p & q & 1-q \end{array} \right) : p,q \in [0,1] \right\}$. A behavioral strategy is a simpler object than a mixed strategy: in this example, specifying a behavioral strategy for Player 1 requires specifying the values of two parameters (*p* and *q*), while specifying a mixed strategy requires specifying the values of three parameters (*p, r* and *q*). Can one then use behavioral strategies rather than mixed strategies? The answer is affirmative, as Theorem 6.2 below states. First we illustrate with an example based on the extensive form of Figure 6.3 (which is a simplified version of the extensive form of Figure 6.2, obtained by removing the moves of Nature; the $z_i$'s are terminal nodes and the outcomes have been omitted).

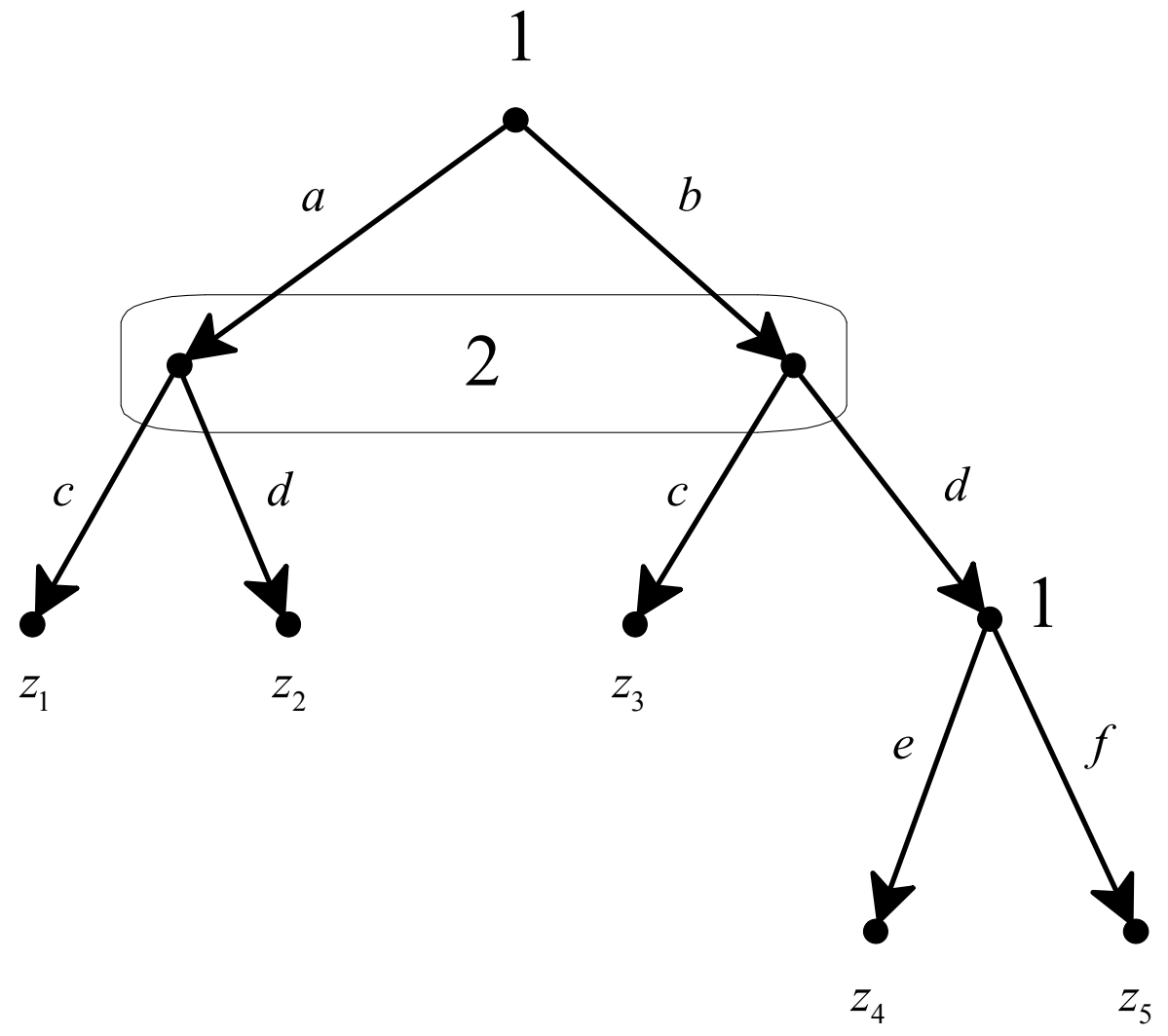

**Figure 6.3**

An extensive frame with the outcomes omitted. The $z_i$'s are terminal nodes.

[9] The definition of perfect recall was included as part of Definition 3.1 (Chapter 3).

Consider the mixed-strategy profile $\sigma = (\sigma_1, \sigma_2)$ with

$$\sigma_1 = \begin{pmatrix} (a,e) & (a,f) & (b,e) & (b,f) \\ \frac{1}{12} & \frac{4}{12} & \frac{2}{12} & \frac{5}{12} \end{pmatrix} \text{ and } \sigma_2 = \begin{pmatrix} c & d \\ \frac{1}{3} & \frac{2}{3} \end{pmatrix}$$

We can compute the probability of reaching every terminal node as follows:

$$P(z_1) = \sigma_1((a,e))\sigma_2(c) + \sigma_1((a,f))\sigma_2(c) = \tfrac{1}{12}\left(\tfrac{1}{3}\right) + \tfrac{4}{12}\left(\tfrac{1}{3}\right) = \tfrac{5}{36},$$
$$P(z_2) = \sigma_1((a,e))\sigma_2(d) + \sigma_1((a,f))\sigma_2(d) = \tfrac{1}{12}\left(\tfrac{2}{3}\right) + \tfrac{4}{12}\left(\tfrac{2}{3}\right) = \tfrac{10}{36}$$
$$P(z_3) = \sigma_1((b,e))\sigma_2(c) + \sigma_1((b,f))\sigma_2(c) = \tfrac{2}{12}\left(\tfrac{1}{3}\right) + \tfrac{5}{12}\left(\tfrac{1}{3}\right) = \tfrac{7}{36}$$
$$P(z_4) = \sigma_1((b,e))\sigma_2(d) = \tfrac{2}{12}\left(\tfrac{2}{3}\right) = \tfrac{4}{36}$$
$$P(z_5) = \sigma_1((b,f))\sigma_2(d) = \tfrac{5}{12}\left(\tfrac{2}{3}\right) = \tfrac{10}{36}$$

That is, the mixed-strategy profile $\sigma = (\sigma_1, \sigma_2)$ gives rise to the following probability distribution over terminal nodes: $\begin{pmatrix} z_1 & z_2 & z_3 & z_4 & z_5 \\ \frac{5}{36} & \frac{10}{36} & \frac{7}{36} & \frac{4}{36} & \frac{10}{36} \end{pmatrix}$. Now consider the following behavioral strategy of Player 1: $\left(\begin{array}{cc|cc} a & b & e & f \\ \frac{5}{12} & \frac{7}{12} & \frac{2}{7} & \frac{5}{7} \end{array}\right)$. What probability distribution over the set of terminal nodes would it determine in conjunction with Player 2's mixed strategy $\sigma_2 = \begin{pmatrix} c & d \\ \frac{1}{3} & \frac{2}{3} \end{pmatrix}$? The calculations are simple ($P(x)$ denotes the probability of choice $x$ for Player 1, according to the given behavioral strategy):

$$P(z_1) = P(a)\sigma_2(c) = \tfrac{5}{12}\left(\tfrac{1}{3}\right) = \tfrac{5}{36}, \quad P(z_2) = P(a)\sigma_2(d) = \tfrac{5}{12}\left(\tfrac{2}{3}\right) = \tfrac{10}{36}$$
$$P(z_3) = P(b)\sigma_2(c) = \tfrac{7}{12}\left(\tfrac{1}{3}\right) = \tfrac{7}{36}, \quad P(z_4) = P(b)\sigma_2(d)P(e) = \tfrac{7}{12}\left(\tfrac{2}{3}\right)\left(\tfrac{2}{7}\right) = \tfrac{4}{36}$$
$$P(z_5) = P(b)\sigma_2(d)P(f) = \tfrac{7}{12}\left(\tfrac{2}{3}\right)\left(\tfrac{5}{7}\right) = \tfrac{10}{36}.$$

Thus, against $\sigma_2 = \begin{pmatrix} c & d \\ \frac{1}{3} & \frac{2}{3} \end{pmatrix}$, Player 1's behavioral strategy $\left(\begin{array}{cc|cc} a & b & e & f \\ \frac{5}{12} & \frac{7}{12} & \frac{2}{7} & \frac{5}{7} \end{array}\right)$ and her mixed strategy $\begin{pmatrix} (a,e) & (a,f) & (b,e) & (b,f) \\ \frac{1}{12} & \frac{4}{12} & \frac{2}{12} & \frac{5}{12} \end{pmatrix}$ are equivalent, in the sense that they give rise to the same probability distribution $\begin{pmatrix} z_1 & z_2 & z_3 & z_4 & z_5 \\ \frac{5}{36} & \frac{10}{36} & \frac{7}{36} & \frac{4}{36} & \frac{10}{36} \end{pmatrix}$ over terminal nodes.

**Theorem 6.1 [Kuhn, 1953].** In extensive forms *with perfect recall,* behavioral strategies and mixed strategies are equivalent, in the sense that, for every mixed strategy there is a behavioral strategy that gives rise to the same probability distribution over terminal nodes.[10]

Without perfect recall, Theorem 6.1 does not hold. To see this, consider the one-player extensive form shown in Figure 6.4 and the mixed strategy $\begin{pmatrix} (a,c) & (a,d) & (b,c) & (b,d) \\ \frac{1}{2} & 0 & 0 & \frac{1}{2} \end{pmatrix}$ which induces the probability distribution $\begin{pmatrix} z_1 & z_2 & z_3 & z_4 \\ \frac{1}{2} & 0 & 0 & \frac{1}{2} \end{pmatrix}$ on the set of terminal nodes.

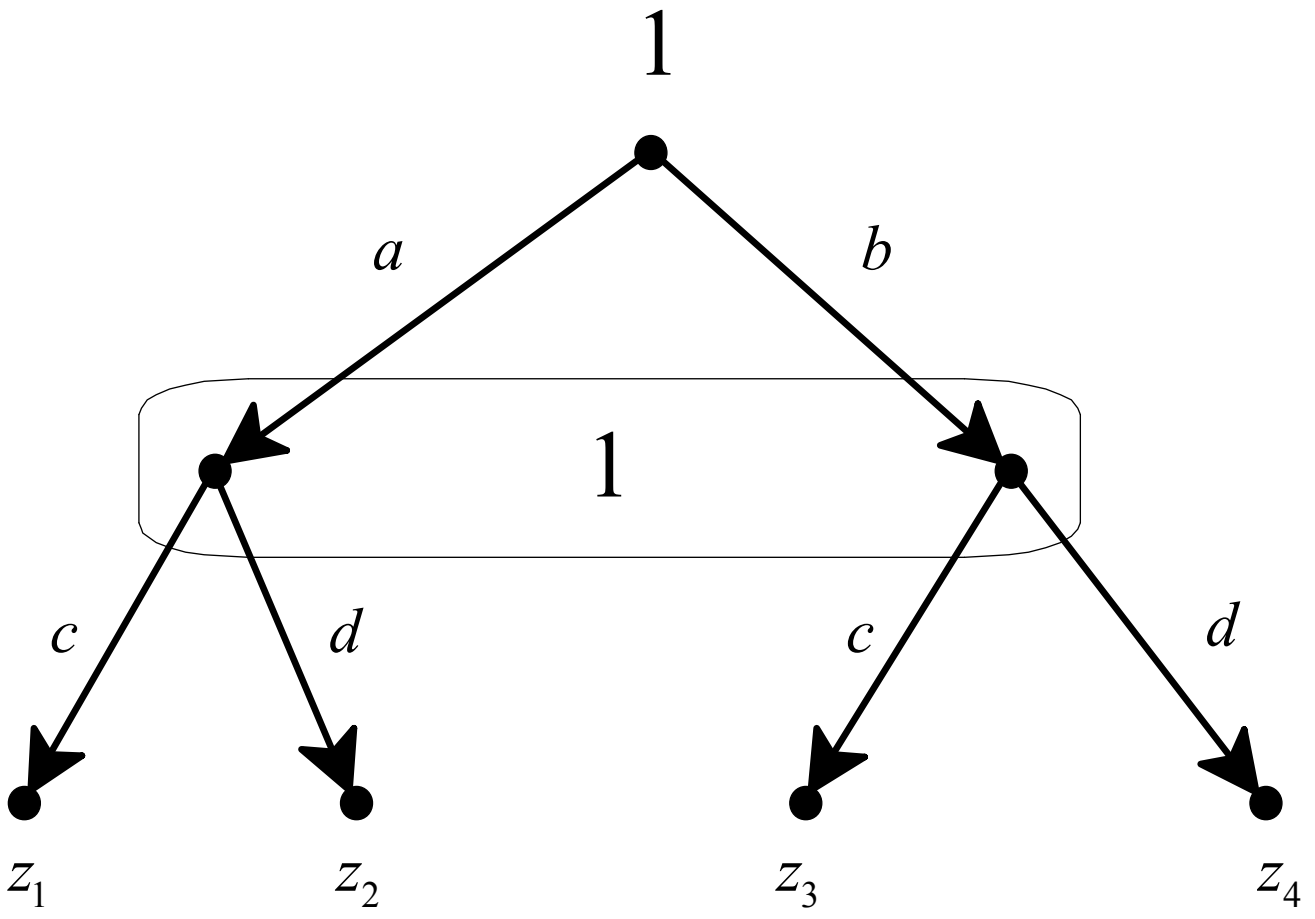

**Figure 6.4**
A one-player extensive frame without perfect recall.

Consider an arbitrary behavioral strategy $\left(\begin{array}{cc|cc} a & b & c & d \\ p & 1-p & q & 1-q \end{array}\right)$, whose corresponding probability distribution over the set of terminal nodes is $\begin{pmatrix} z_1 & z_2 & z_3 & z_4 \\ pq & p(1-q) & (1-p)q & (1-p)(1-q) \end{pmatrix}$. In order to have $P(z_2) = 0$ it must be

[10] A more precise statement is as follows. Consider an extensive form with perfect recall and a Player *i*. Let $x_{-i}$ be an arbitrary profile of strategies of the players other than *i*, where, for every $j \neq i$, $x_j$ is either a mixed or a behavioral strategy of Player *j*. Then, for every mixed strategy $\sigma_i$ of Player *i* there is a behavioral strategy $b_i$ of Player *i* such that $(\sigma_i, x_{-i})$ and $(b_i, x_{-i})$ give rise to the same probability distribution over the set of terminal nodes.

that either $p = 0$ or $q = 1$. If $p = 0$ then $P(z_1) = 0$ and if $q = 1$ then $P(z_4) = 0$. Thus the probability distribution $\begin{pmatrix} z_1 & z_2 & z_3 & z_4 \\ \frac{1}{2} & 0 & 0 & \frac{1}{2} \end{pmatrix}$ cannot be achieved with a behavioral strategy.

Since the focus of this book is on extensive-form games with perfect recall, by appealing to Theorem 6.3 from now on we can restrict attention to behavioral strategies.

As usual, one goes from a frame to a game by adding preferences over outcomes. Let $O$ be the set of basic outcomes (recall that with every terminal node is associated a basic outcome) and $\mathcal{L}(O)$ the set of lotteries (probability distributions) over $O$.

**Definition 6.2.** An *extensive-form game with cardinal payoffs* is an extensive frame (with, possibly, chance moves) together with a von Neumann-Morgenstern ranking $\succsim_i$ of the set of lotteries $\mathcal{L}(O)$, for every Player $i$.

As usual, it is convenient to represent a von Neumann-Morgenstern ranking by means of a von Neumann-Morgenstern utility function and replace the outcomes with a vector of utilities, one for every player.

For example, consider the extensive form of Figure 6.2 where the set of basic outcomes is $O = \{o_1, o_2, o_3, o_4, o_5\}$ and suppose that Player 1 has a von Neumann-Morgenstern ranking of $\mathcal{L}(O)$ that is represented by the following von Neumann-Morgenstern utility function: $\begin{cases} \text{outcome}: & o_1 & o_2 & o_3 & o_4 & o_5 \\ U_1: & 5 & 2 & 0 & 1 & 3 \end{cases}$. Suppose also that Player 2 has preferences represented by the von Neumann-Morgenstern utility function $\begin{cases} \text{outcome}: & o_1 & o_2 & o_3 & o_4 & o_5 \\ U_2: & 3 & 6 & 4 & 5 & 0 \end{cases}$.

Then from the extensive frame of Figure 6.2 we obtain the extensive-form game with cardinal payoffs shown in Figure 6.5. Since the expected utility of lottery $\begin{pmatrix} o_1 & o_2 \\ \frac{2}{3} & \frac{1}{3} \end{pmatrix}$ is 4 for both players, and the expected utility of lottery $\begin{pmatrix} o_1 & o_3 & o_4 \\ \frac{1}{5} & \frac{3}{5} & \frac{1}{5} \end{pmatrix}$ is 1.2 for Player 1 and 4 for Player 2, we can simplify the game by replacing the first move of Nature with the payoff vector (4,4) and the second move of Nature with the payoff vector (1.2, 4). The simplified game is shown in Figure 6.6.

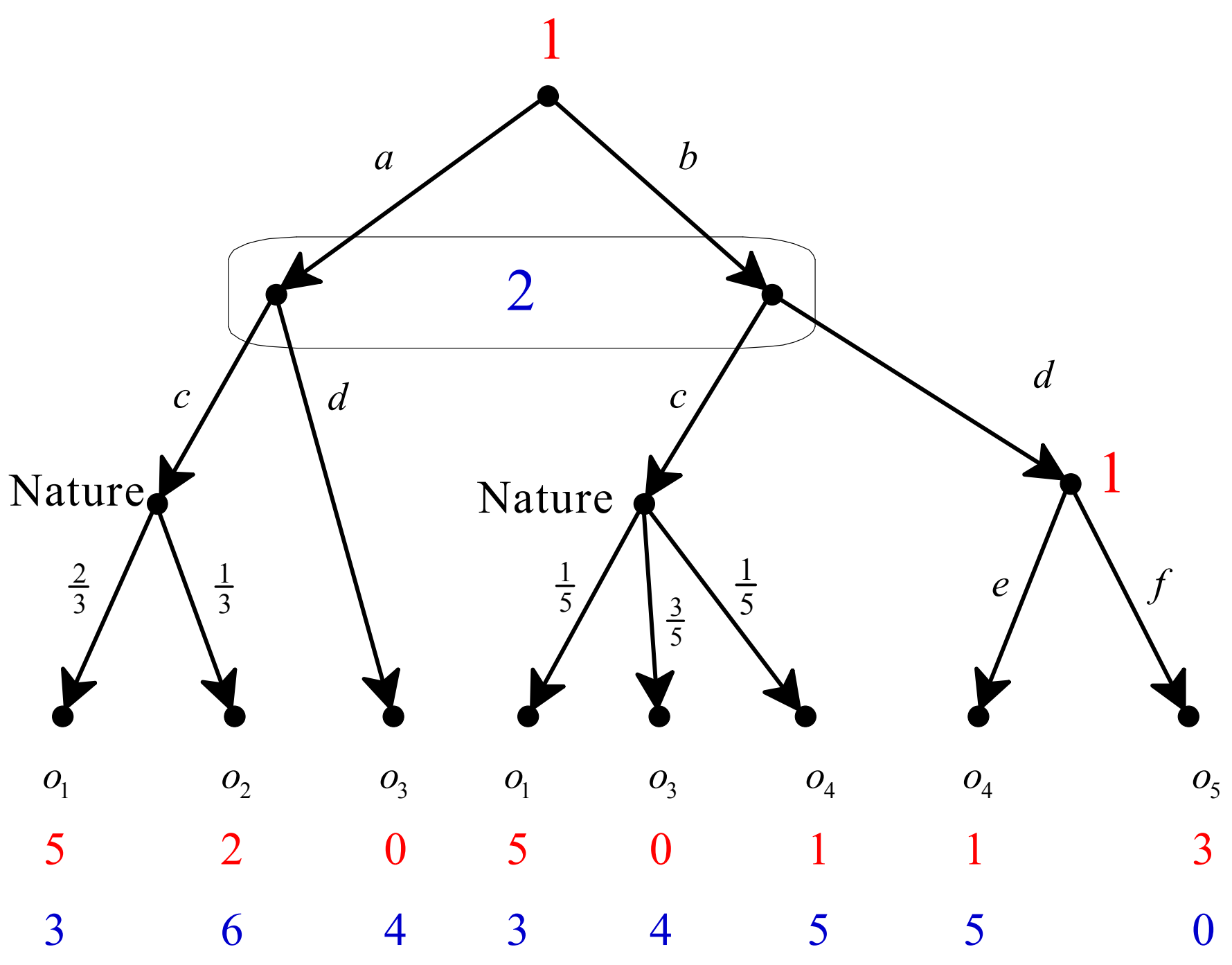

**Figure 6.5**
An extensive game based on the frame of Figure 6.2.
The terminal nodes have not been labeled. The $o_i$'s are basic outcomes.

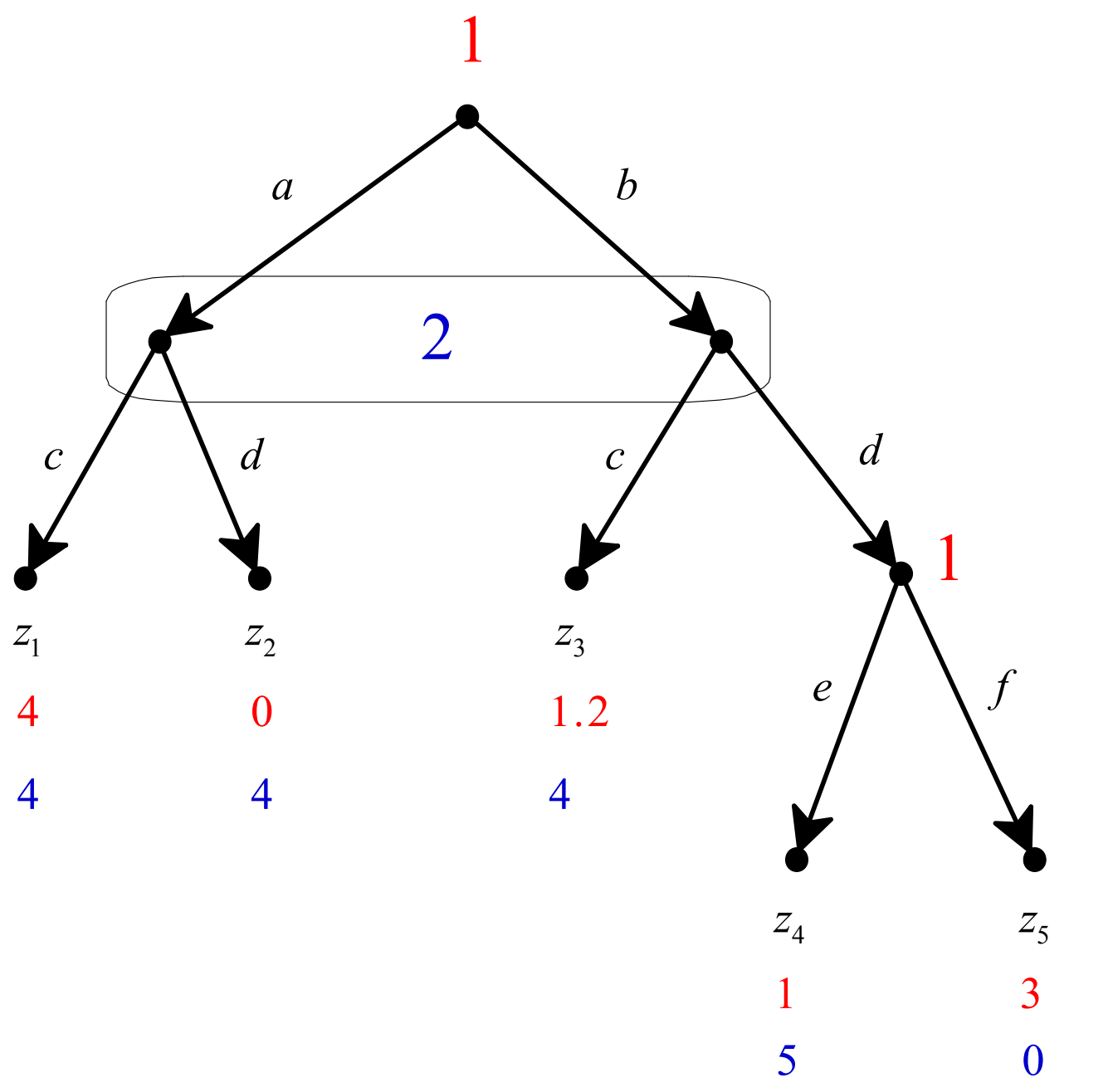

**Figure 6.6**
A simplified version of the game of Figure 6.5. The $z_i$'s are terminal nodes.
Note that this is a game based on the frame of Figure 6.3

Given an extensive game with cardinal payoffs, associated with every behavioral strategy profile is a lottery over basic outcomes and thus, using a von Neumann-Morgenstern utility function for each player, a payoff for each player. For example, the behavioral strategy profile $\left(\begin{pmatrix} a & b & e & f \\ \frac{5}{12} & \frac{7}{12} & \frac{2}{7} & \frac{5}{7} \end{pmatrix}, \begin{pmatrix} c & d \\ \frac{1}{3} & \frac{2}{3} \end{pmatrix}\right)$ for the extensive game of Figure 6.5 gives rise to the lottery $\begin{pmatrix} o_1 & o_2 & o_3 & o_4 & o_5 \\ \frac{71}{540} & \frac{25}{540} & \frac{213}{540} & \frac{81}{540} & \frac{150}{540} \end{pmatrix}$ (for instance, the probability of basic outcome $o_1$ is calculated as follows: $P(o_1) = P(a)P(c)\frac{2}{3} + P(b)P(c)\frac{1}{5} = \frac{5}{12}\frac{1}{3}\frac{2}{3} + \frac{7}{12}\frac{1}{3}\frac{1}{5} = \frac{71}{540}$). Using the utility function postulated above for Player 1, namely $\begin{matrix} \text{outcome}: & o_1 & o_2 & o_3 & o_4 & o_5 \\ U_1: & 5 & 2 & 0 & 1 & 3 \end{matrix}$, we get a corresponding payoff of $\frac{71}{540}5 + \frac{25}{540}2 + \frac{213}{540}0 + \frac{81}{540}1 + \frac{150}{540}3 = \frac{936}{540} = 1.733$. An alternative way of computing this payoff is by using the simplified game of Figure 6.6 where the behavioral strategy profile $\left(\begin{pmatrix} a & b & e & f \\ \frac{5}{12} & \frac{7}{12} & \frac{2}{7} & \frac{5}{7} \end{pmatrix}, \begin{pmatrix} c & d \\ \frac{1}{3} & \frac{2}{3} \end{pmatrix}\right)$ yields the probability distribution over terminal nodes $\begin{pmatrix} z_1 & z_2 & z_3 & z_4 & z_5 \\ \frac{5}{36} & \frac{10}{36} & \frac{7}{36} & \frac{4}{36} & \frac{10}{36} \end{pmatrix}$, which, in turn, yields the probability distribution $\begin{pmatrix} 4 & 0 & 1.2 & 1 & 3 \\ \frac{5}{36} & \frac{10}{36} & \frac{7}{36} & \frac{4}{36} & \frac{10}{36} \end{pmatrix}$ over utilities for Player 1. From the latter we get that the expected payoff for Player 1 is $\frac{5}{36}4 + \frac{10}{36}0 + \frac{7}{36}1.2 + \frac{4}{36}1 + \frac{10}{36}3 = \frac{936}{540} = 1.733$. The calculations for Player 2 are similar (see Exercise 6.3).

Test your understanding of the concepts introduced in this section, by going through the exercises in Section 6.E.1 of Appendix 6.E at the end of this chapter.

# 6.2 Subgame-perfect equilibrium revisited

The notion of subgame-perfect equilibrium was introduced in Chapter 3 (Definition 3.5) for extensive-form games with ordinal payoffs. When payoffs are ordinal, a subgame-perfect equilibrium might fail to exist because either the entire game or a proper subgame does not have any Nash equilibria. In the case of finite extensive-form games with cardinal payoffs, a subgame-perfect equilibrium always exists, because – by Nash's theorem (Theorem 5.1, Chapter 5) – every finite game has at least one Nash equilibrium in mixed strategies. Thus, in the case of cardinal payoffs, the subgame-perfect equilibrium algorithm (Definition 3.6, Chapter 3) never halts and the output of the algorithm is a subgame-perfect equilibrium. We shall illustrate this with the extensive-form game shown in Figure 6.7 below.

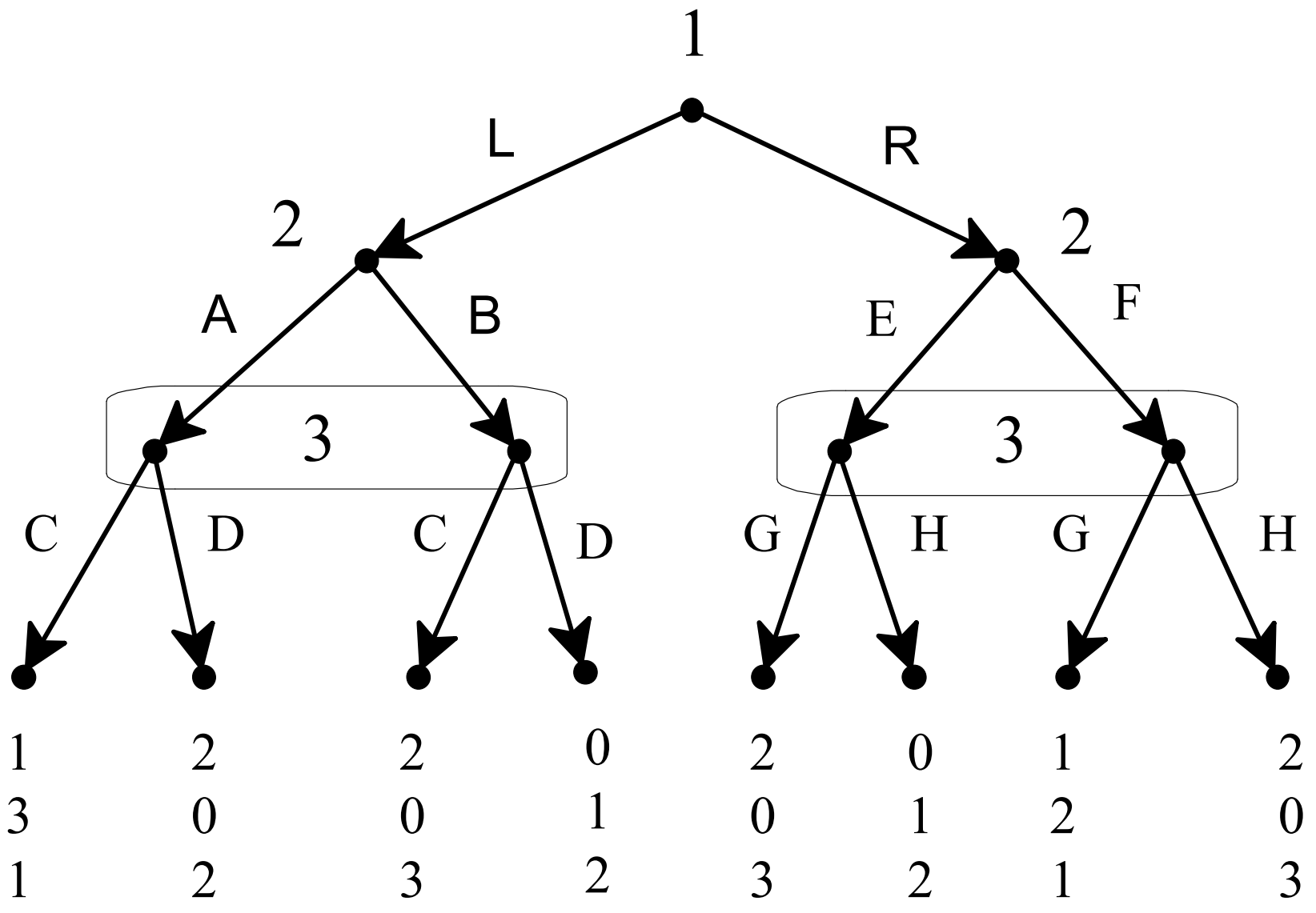

**Figure 6.7**
An extensive-form game with cardinal payoffs.

Let us apply the subgame-perfect equilibrium algorithm to this game. We start with the proper subgame that begins at Player 2's decision node on the left, whose strategic form is shown in Table 6.8. Note that this subgame has no pure-strategy Nash equilibria. Thus if payoffs were merely ordinal payoffs the algorithm would halt and we would conclude that the game of Figure 6.7 has no subgame-perfect equilibria. However, we will assume that payoffs are cardinal (that is, that they are von Neumann-Morgenstern utilities).

| | | Player 3 | |
|---|---|---|---|
| | | *C* | *D* |
| Player 2 | *A* | 3 , 1 | 0 , 2 |
| | *B* | 0 , 3 | 1 , 2 |

**Table 6.8**
The strategic form of the proper subgame on the left in the game of Figure 6.7.

To find the mixed-strategy Nash equilibrium of the game of Table 6.8, let $p$ be the probability of $A$ and $q$ the probability of $C$. Then we need $q$ to be such that $3q = 1 - q$, that is, $q = \frac{1}{4}$, and $p$ to be such that $p + 3(1 - p) = 2$, that is, $p = \frac{1}{2}$.

Thus the Nash equilibrium of this proper subgame is $\left(\begin{pmatrix} A & B \\ \frac{1}{2} & \frac{1}{2} \end{pmatrix}, \begin{pmatrix} C & D \\ \frac{1}{4} & \frac{3}{4} \end{pmatrix}\right)$, yielding the following payoffs: $\frac{1}{2}\frac{1}{4}1+\frac{1}{2}\frac{3}{4}2+\frac{1}{2}\frac{1}{4}2+\frac{1}{2}\frac{3}{4}0=1.125$ for Player 1, $\frac{1}{2}\frac{1}{4}3+\frac{1}{2}\frac{3}{4}0+\frac{1}{2}\frac{1}{4}0+\frac{1}{2}\frac{3}{4}1=0.75$ for Player 2 and $\frac{1}{2}\frac{1}{4}1+\frac{1}{2}\frac{3}{4}2+\frac{1}{2}\frac{1}{4}3+\frac{1}{2}\frac{3}{4}2=2$ for Player 3.[11] Thus we can simplify the game of Figure 6.7 as follows:

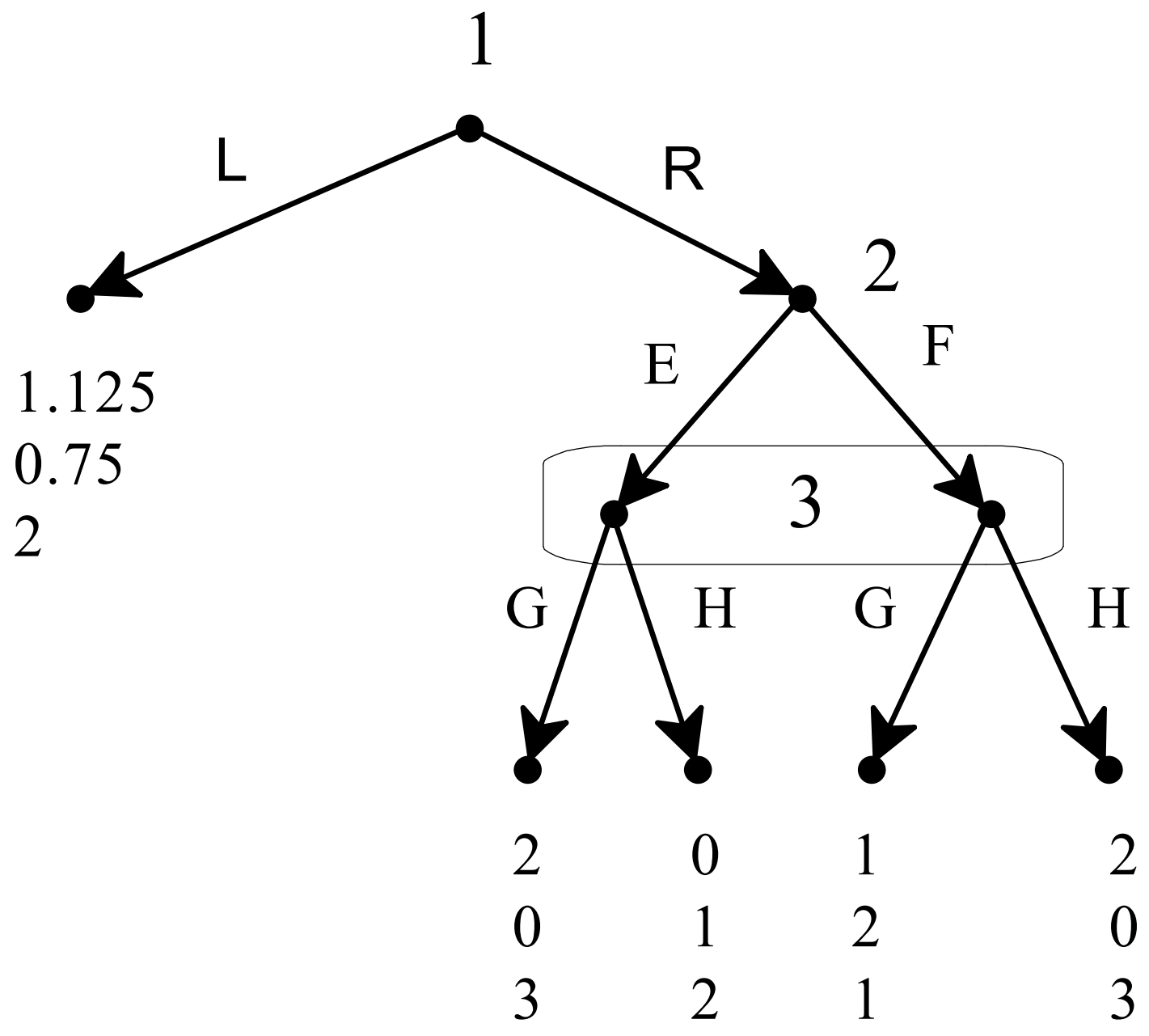

**Figure 6.9**
The game of Figure 6.7 after replacing the proper subgame on the left with the payoffs associated with its Nash equilibrium.

Now consider the proper subgame of the game of Figure 6.9 (the subgame that starts at Player 2's node). Its strategic form is shown in Table 6.10.

| | | Player 3 | |
|---|---|---|---|
| | | *G* | *H* |
| Player 2 | *E* | 0 , 3 | 1 , 2 |
| | *F* | 2 , 1 | 0 , 3 |

**Table 6.10**
The strategic form of the proper subgame of the game of Figure 6.9.

[11] Note that, by Theorem 5.2 (Chapter 5) the Nash equilibrium payoff of Player 2 can be computed much more quickly as the expected payoff from the pure strategy $A$, namely $\frac{3}{4}$, and the payoff of Player 3 as the expected payoff from the pure strategy $D$, namely 2.

Again, there is no pure-strategy Nash equilibrium. To find the mixed-strategy equilibrium let $p$ be the probability of $E$ and $q$ the probability of $G$. Then we need $q$ to be such that $1-q=2q$, that is, $q=\frac{1}{3}$, and $p$ to be such that $3p+1-p=2p+3(1-p)$, that is, $p=\frac{2}{3}$. Hence the Nash equilibrium is $\left(\begin{pmatrix} E & F \\ \frac{2}{3} & \frac{1}{3}\end{pmatrix},\begin{pmatrix} G & H \\ \frac{1}{3} & \frac{2}{3}\end{pmatrix}\right)$ yielding the following payoffs: $\frac{2}{3}\frac{1}{3}2+\frac{2}{3}\frac{2}{3}0+\frac{1}{3}\frac{1}{3}1+\frac{1}{3}\frac{2}{3}2=1$ for Player 1, $\frac{2}{3}\frac{1}{3}0+\frac{2}{3}\frac{2}{3}1+\frac{1}{3}\frac{1}{3}2+\frac{1}{3}\frac{2}{3}0=0.67$ for Player 2 and $\frac{2}{3}\frac{1}{3}3+\frac{2}{3}\frac{2}{3}2+\frac{1}{3}\frac{1}{3}1+\frac{1}{3}\frac{2}{3}3=2.33$ for Player 3.

Thus we can simplify the game of Figure 6.9 as shown in Figure 6.11 below, where the optimal choice for Player 1 is $L$.

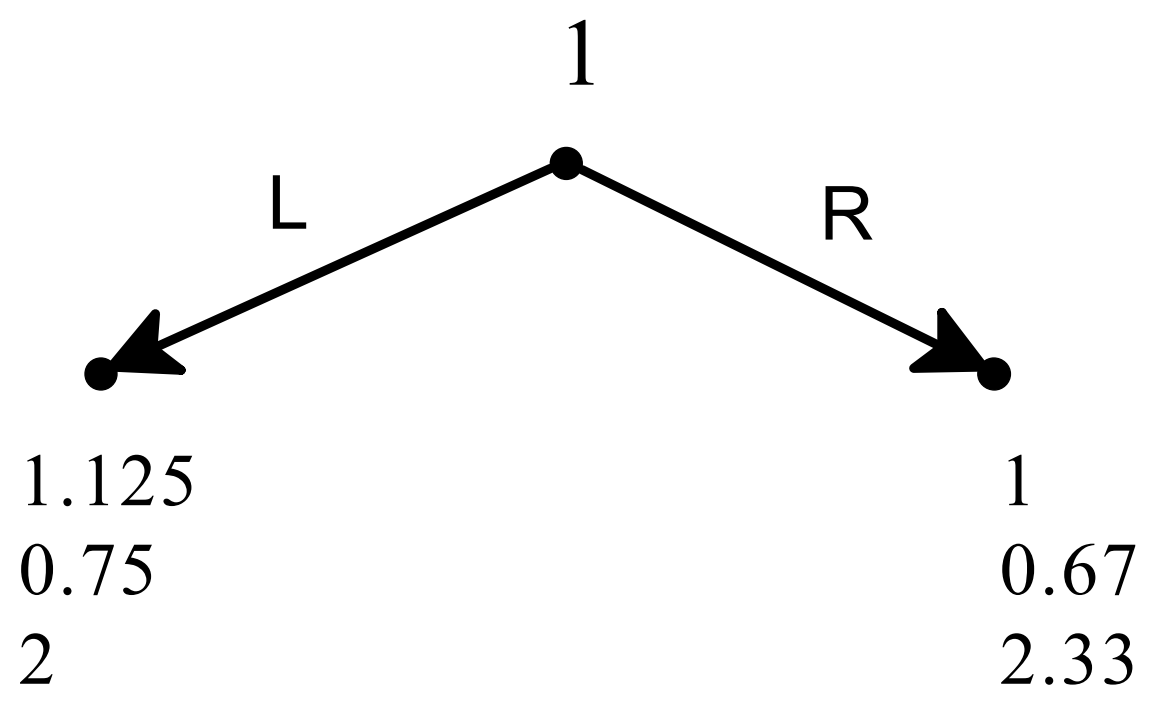

**Figure 6.11**
The game of Figure 6.9 after replacing the proper subgame
with the payoffs associated with the Nash equilibrium.

Hence the subgame-perfect equilibrium of the game of Figure 6.7 (expressed in terms of behavioral strategies) is:

$$\left(\begin{pmatrix} L & R \\ 1 & 0\end{pmatrix},\left(\begin{array}{cc|cc} A & B & E & F \\ \frac{1}{2} & \frac{1}{2} & \frac{2}{3} & \frac{1}{3}\end{array}\right),\left(\begin{array}{cc|cc} C & D & G & H \\ \frac{1}{4} & \frac{3}{4} & \frac{1}{3} & \frac{2}{3}\end{array}\right)\right)$$

We conclude this section with the following theorem, which is a corollary of Theorem 5.1 (Chapter 5).

**Theorem 6.2.** Every finite extensive-form game with cardinal payoffs has at least one subgame-perfect equilibrium

Test your understanding of the concepts introduced in this section, by going through the exercises in Section 6.E.2 of Appendix 6.E at the end of this chapter.

# 6.3 Problems with the notion of subgame-perfect equilibrium

The notion of subgame-perfect equilibrium is a refinement of Nash equilibrium. As explained in Chapter 2, in the context of perfect-information games, the notion of subgame-perfect equilibrium eliminates some "unreasonable" Nash equilibria that involve incredible threats. However, not every subgame-perfect equilibrium can be viewed as a "rational solution". To see this, consider the extensive-form game shown in Figure 6.12 below. This game has no proper subgames and thus the set of subgame-perfect equilibria coincides with the set of Nash equilibria. The pure-strategy Nash equilibria of this game are $(a, f, c)$, $(a, e, c)$, $(b, e, c)$ and $(b, f, d)$. It can be argued that neither $(a, f, c)$ nor $(b, f, d)$ can be considered "rational solutions".

Consider first the Nash equilibrium $(a, f, c)$. Player 2's plan to play *f* is rational only in the very limited sense that, given that Player 1 plays *a*, what Player 2 plans to do is irrelevant because it cannot affect anybody's payoff; thus *f* is as good as *e*. However, if we take Player 2's strategy as a "serious" plan specifying what Player 2 would actually do if she had to move, then –given that Player 3 plays *c* –*e* would give Player 2 a payoff of 2, while *f* would only give a payoff of 1. Thus *e* seems to be a better strategy than *f*, if Player 2 takes the contingency "seriously".

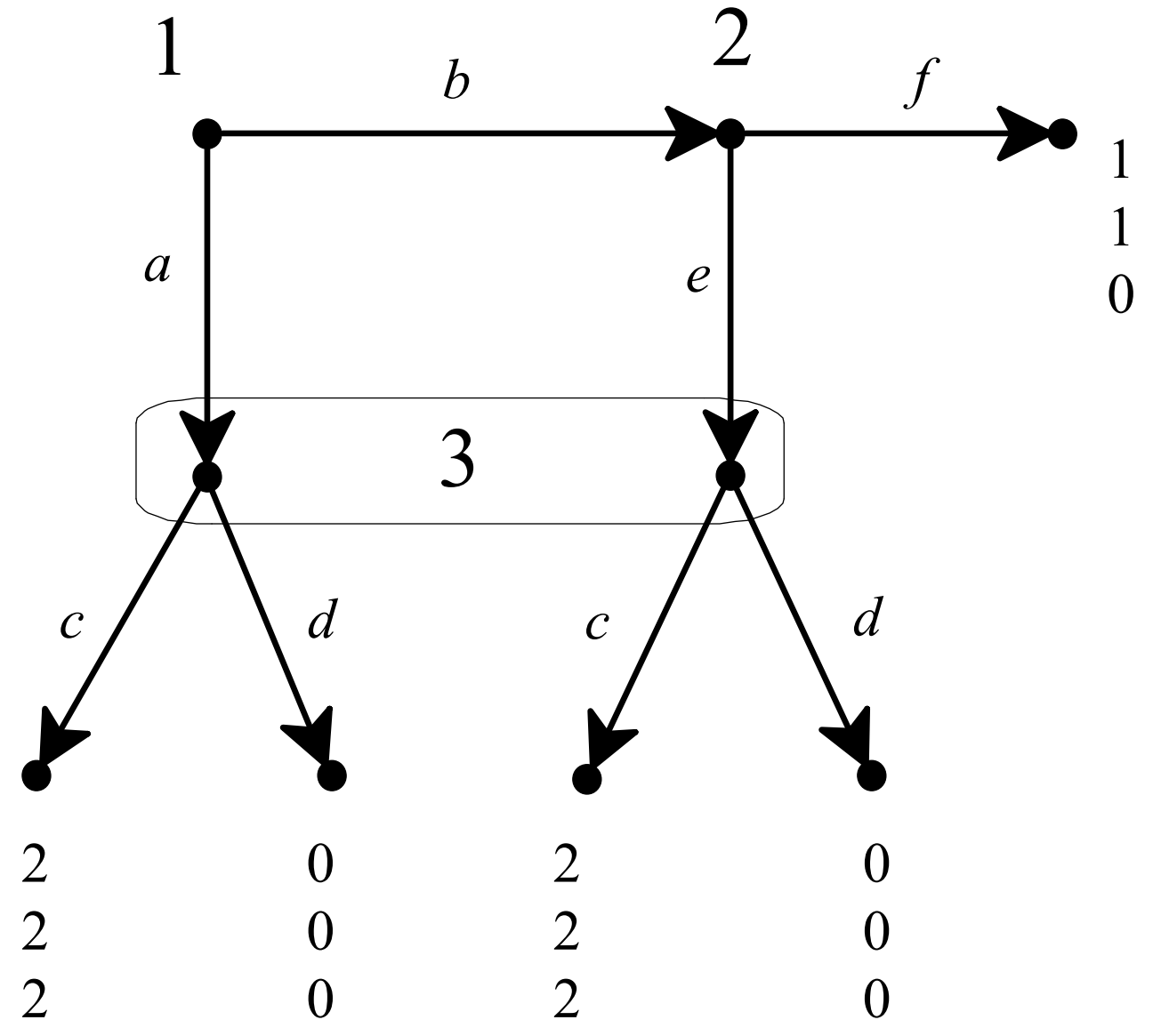

**Figure 6.12**
An extensive-form game showing the insufficiency of the notion of subgame-perfect equilibrium.

Consider now the Nash equilibrium $(b, f, d)$ and focus on Player 3. As before, Player 3's plan to play $d$ is rational only in the very limited sense that, given that Player 1 plays $a$ and Player 2 plays $f$, what Player 3 plans to do is irrelevant, so that $c$ is as good as $d$. However, if Player 3 did find himself having to play, it would not be rational for him to play $d$, since $d$ is a strictly dominated choice: no matter whether he is making his choice at the left node or at the right node of his information set, $c$ gives him a higher payoff than $d$. How can it be then that $d$ can be part of a Nash equilibrium? The answer is that $d$ is strictly dominated *conditional on Player 3's information set being reached* but not as a plan formulated before the play of the game starts. In other words, $d$ is strictly dominated *as a choice but not as a strategy*.

The notion of subgame-perfect equilibrium is not strong enough to eliminate "unreasonable" Nash equilibria such as $(a, f, c)$ and $(b, f, d)$ in the game of Figure 6.12. In order to do that we will need a stronger notion. This issue is postponed to a later chapter.

# Appendix 6.E: Exercises

## 6.E.1. Exercises for Section 6.1: Behavioral strategies in dynamic games

The answers to the following exercises are in Appendix S at the end of this chapter.

**Exercise 6.1**. What properties must an extensive-form frame satisfy in order for it to be the case that, for a given player, the set of mixed strategies coincides with the set of behavioral strategies? [Assume that there are at least two choices at every information set.]

**Exercise 6.2**. Suppose that, in a given extensive-form frame, Player 1 has four information sets: at one of them she has two choices and at each of the other three she has three choices.

**(a)** How many parameters are needed to specify a mixed strategy of Player 1?

**(b)** How many parameters are needed to specify a behavioral strategy of Player 1?

**Exercise 6.3**. From the behavioral strategy profile $\left(\begin{pmatrix} a & b & e & f \\ \frac{5}{12} & \frac{7}{12} & \frac{2}{7} & \frac{5}{7} \end{pmatrix}, \begin{pmatrix} c & d \\ \frac{1}{3} & \frac{2}{3} \end{pmatrix}\right)$ calculate the payoff of Player 2 in two ways: (1) using the game of Figure 6.5 and (2) using the simplified game of Figure 6.6.

**Exercise 6.4**. Consider the following extensive form, where $o_1,...,o_5$ are basic outcomes.

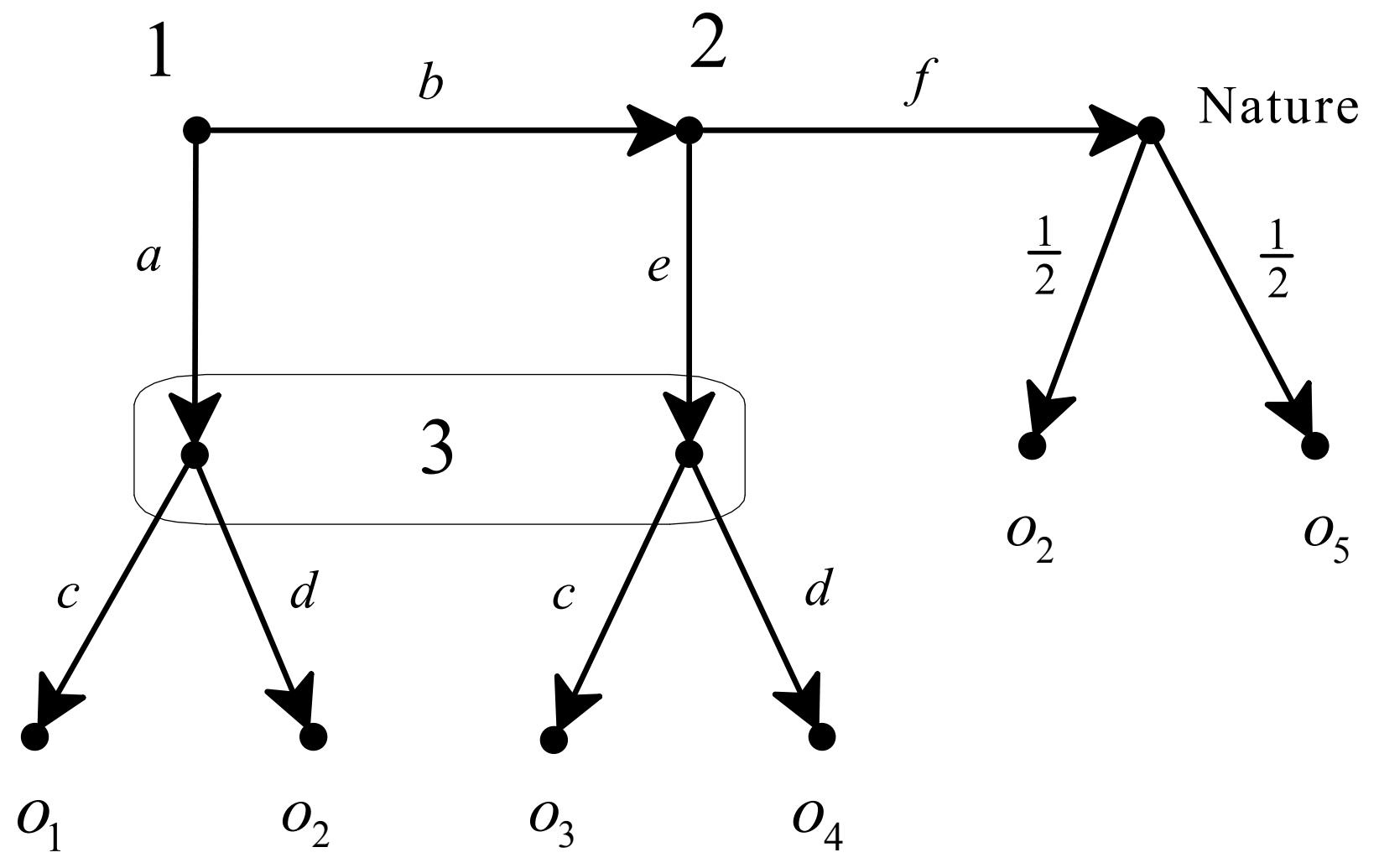

Player 1's ranking of $O$ is $o_1 \succ_1 o_5 \succ_1 o_4 \succ_1 o_2 \sim_1 o_3$; furthermore, she is indifferent between $o_5$ and the lottery $\begin{pmatrix} o_1 & o_2 & o_3 \\ \frac{6}{8} & \frac{1}{8} & \frac{1}{8} \end{pmatrix}$ and is also indifferent between $o_4$ and the lottery $\begin{pmatrix} o_2 & o_5 \\ \frac{2}{3} & \frac{1}{3} \end{pmatrix}$. Player 2's ranking of $O$ is $o_1 \sim_2 o_2 \sim_2 o_4 \succ_2 o_3 \succ_2 o_5$; furthermore, he is indifferent between $o_3$ and the lottery $\begin{pmatrix} o_1 & o_2 & o_5 \\ \frac{1}{10} & \frac{1}{10} & \frac{8}{10} \end{pmatrix}$. Finally, Player 3's ranking of $O$ is $o_2 \succ_3 o_4 \succ_3 o_3 \sim_3 o_5 \succ_3 o_1$; furthermore, she is indifferent between $o_4$ and the lottery $\begin{pmatrix} o_1 & o_2 & o_3 \\ \frac{1}{4} & \frac{1}{2} & \frac{1}{4} \end{pmatrix}$ and is also indifferent between $o_3$ and the lottery $\begin{pmatrix} o_1 & o_2 \\ \frac{3}{5} & \frac{2}{5} \end{pmatrix}$. Write the corresponding extensive-form game.

## 6.E.2. Exercises for Section 6.2: Subgame-perfect equilibrium revisited

The answers to the following exercises are in Appendix S at the end of this chapter.

**Exercise 6.5**. Consider the following extensive-form game with cardinal payoffs.

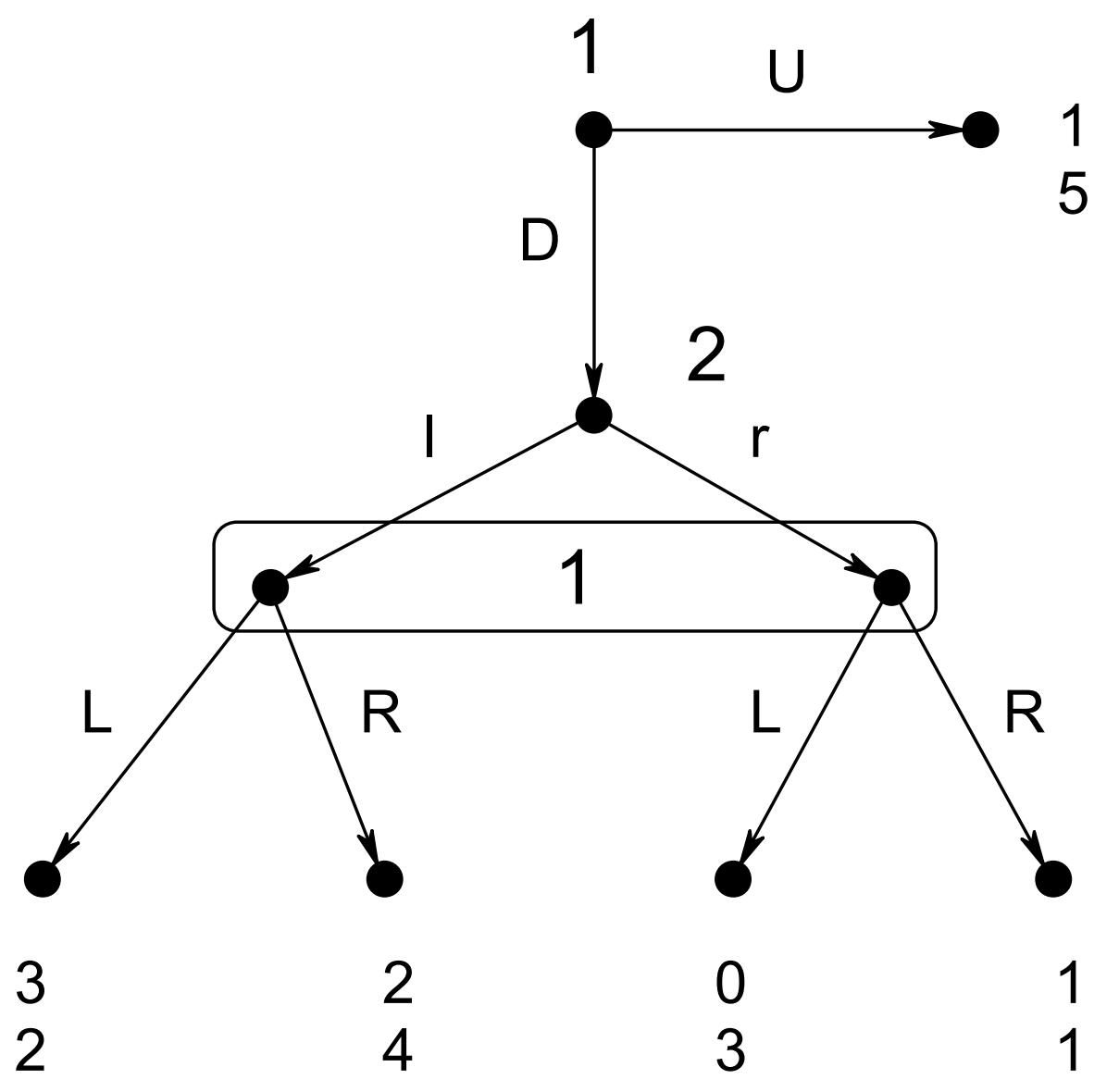

(a) Write the corresponding strategic-form game and find all the pure-strategy Nash equilibria.

(b) Find the subgame-perfect equilibrium.

**Exercise 6.6.** Consider the following extensive form (where the basic outcomes are denoted by *xj* instead of $o_j$, $j$ = 1,..., 10).

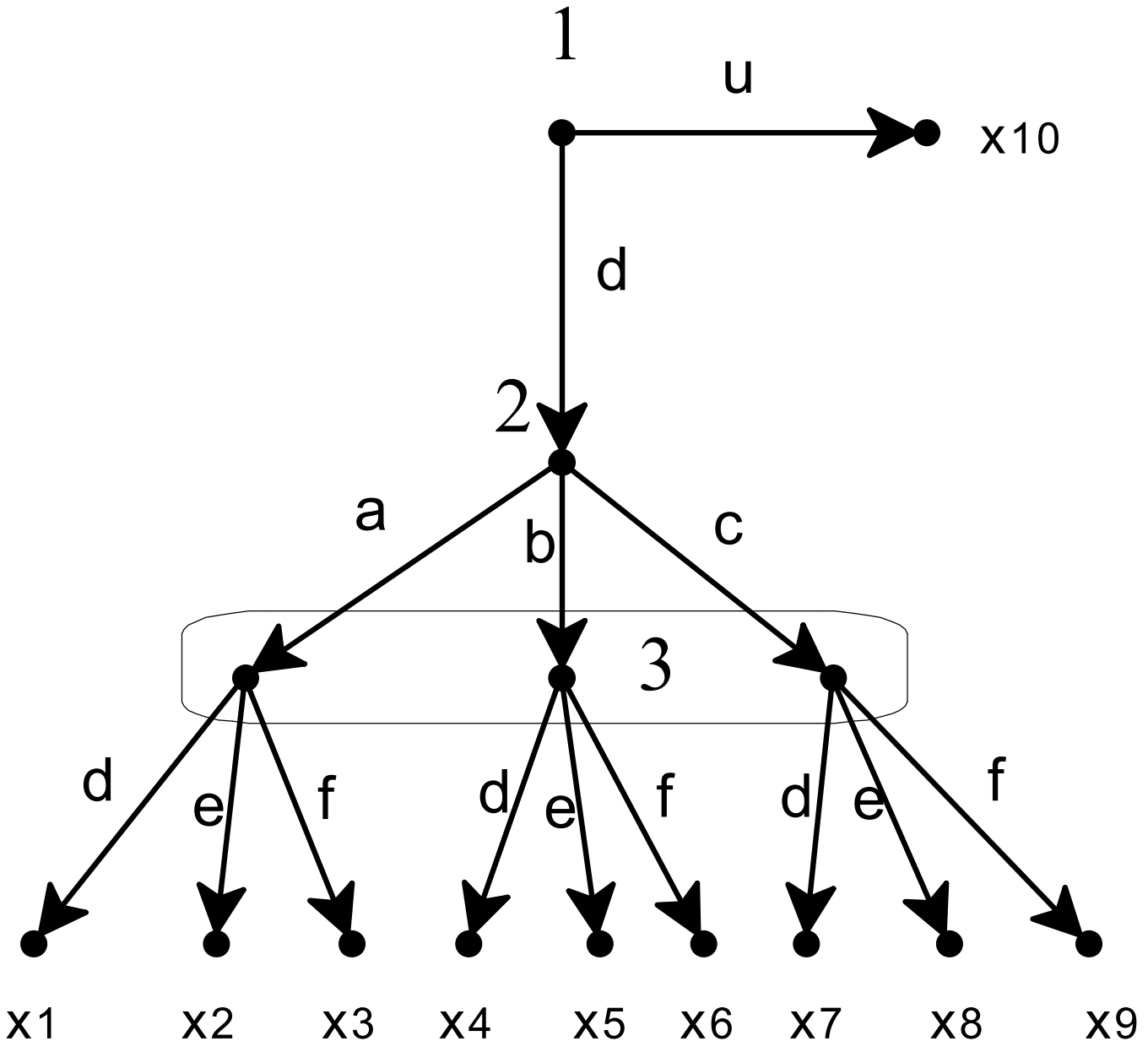

All the players satisfy the axioms of expected utility. They rank the outcomes as indicated below (if outcome $w$ is above outcome $y$ then $w$ is strictly preferred to $y$, and if $w$ and $y$ are written next to each other then the player is indifferent between the two):

$$Player\ 1: \begin{pmatrix} x_7, x_9 \\ x_1, x_2, x_4, x_5 \\ x_{10} \\ x_3, x_6, x_8 \end{pmatrix}, \quad Player\ 2: \begin{pmatrix} x_1, x_3 \\ x_4, x_5 \\ x_2, x_7, x_8 \\ x_6 \\ x_9 \end{pmatrix}, \quad Player\ 3: \begin{pmatrix} x_2, x_7 \\ x_8 \\ x_1, x_4, x_9 \\ x_3, x_5, x_6 \end{pmatrix}.$$

Furthermore, Player 2 is indifferent between $x_4$ and the lottery $\begin{pmatrix} x_1 & x_2 \\ \frac{1}{2} & \frac{1}{2} \end{pmatrix}$ and Player 3 is indifferent between $x_1$ and the lottery $\begin{pmatrix} x_2 & x_5 \\ \frac{1}{2} & \frac{1}{2} \end{pmatrix}$. Although the above information is not sufficient to determine the von Neumann-Morgenstern utility functions of the players, it is sufficient to compute the subgame-perfect equilibrium. Find the subgame-perfect equilibrium.

**Exercise 6.7**. Consider the following extensive-form game with cardinal payoffs.

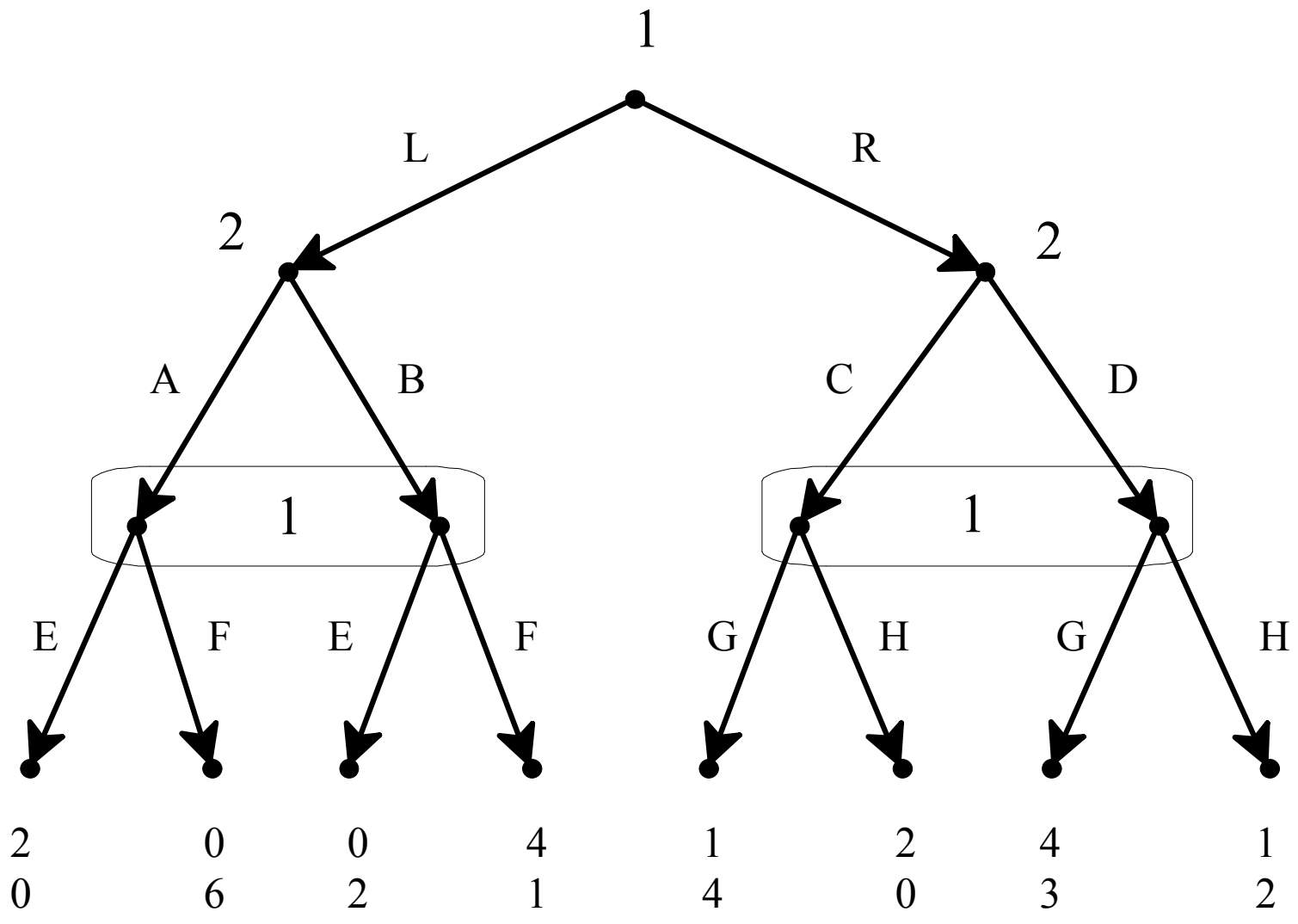

**(a)** Write the corresponding strategic-form game.
**(b)** Find all the pure-strategy Nash equilibria.
**(c)** Find the mixed-strategy subgame-perfect equilibrium.

**Exercise 6.8**. Consider the following extensive-form game with cardinal payoffs:

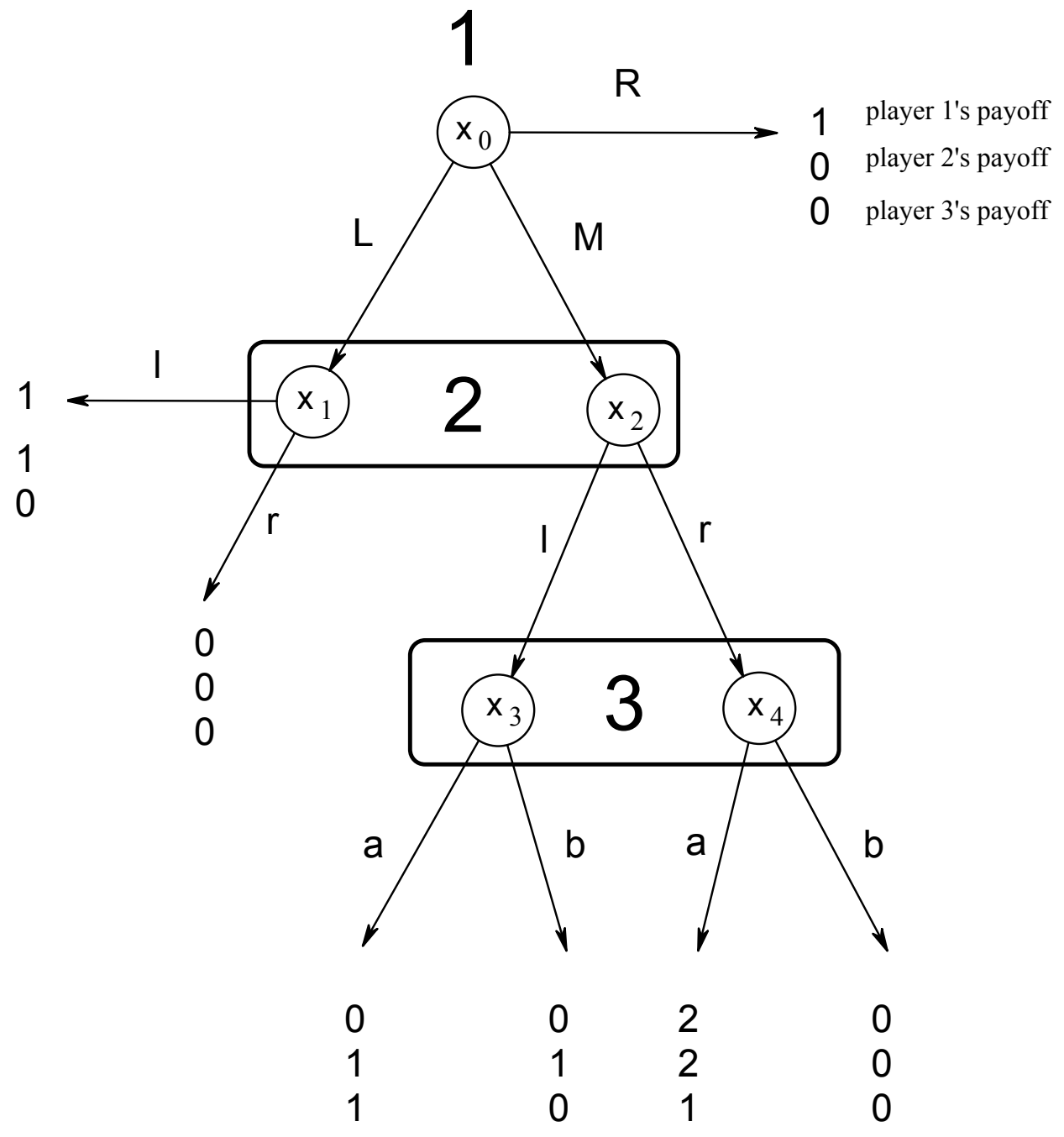

**(a)** Find all the pure-strategy Nash equilibria. Which ones are also subgame perfect?

**(b)** [This is a more challenging question] Prove that there is no mixed-strategy Nash equilibrium where Player 1 plays $M$ with probability strictly between 0 and 1.

◊◊◊◊◊◊◊◊◊◊◊◊

### Exercise 6.9: Challenging Question.

You have to go to a part of town where many people have been mugged recently. You consider whether you should leave your wallet at home or carry it with you. Of the four possible outcomes, your most preferred one is having your wallet with you and not being mugged. Being mugged is a very unpleasant experience, so your second favorite alternative is not carrying your wallet and not being mugged (although not having any money with you can be very inconvenient). If, sadly enough, your destiny is to be mugged, then you prefer to have your wallet with you (possibly with not too much money in it!) because you don't want to have to deal with a frustrated mugger. A typical potential mugger, on the other hand, does not care whether or not you are carrying a wallet, in case he decides not to mug you (that is, he is indifferent between the corresponding two outcomes). Of course his favorite outcome is the one where you have your wallet with you and he mugs you. His least preferred outcome is the one where he attempts to mug you and you don't have your wallet with you (he risks being caught for nothing). Denote the possible outcomes as follows:

| | | **Potential mugger** | |
|---|---|---|---|
| | | Not mug | Mug |
| **You** | Leave wallet at home | $z_1$ | $z_2$ |
| | Take wallet with you | $z_3$ | $z_4$ |

**(a)** What is the ordinal ranking of the outcomes for each player?

**(b)** Suppose now that both players have von Neumann-Morgenstern utility functions. **You** are indifferent between the following lotteries: $L_1 = \begin{pmatrix} z_1 & z_2 & z_3 & z_4 \\ \frac{3}{20} & \frac{14}{20} & \frac{3}{20} & 0 \end{pmatrix}$ and $L_2 = \begin{pmatrix} z_1 & z_2 & z_3 & z_4 \\ 0 & \frac{1}{2} & 0 & \frac{1}{2} \end{pmatrix}$; furthermore, you are indifferent between $L_3 = \begin{pmatrix} z_1 & z_2 & z_3 & z_4 \\ 0 & \frac{2}{3} & \frac{1}{3} & 0 \end{pmatrix}$ and $L_4 = \begin{pmatrix} z_1 & z_2 & z_3 & z_4 \\ \frac{1}{2} & \frac{1}{2} & 0 & 0 \end{pmatrix}$. The **potential mugger** is indifferent between the two lotteries $L_5 = \begin{pmatrix} z_1 & z_2 & z_3 & z_4 \\ \frac{1}{4} & \frac{1}{4} & \frac{1}{4} & \frac{1}{4} \end{pmatrix}$ and $L_6 = \begin{pmatrix} z_1 & z_2 & z_3 & z_4 \\ \frac{8}{128} & \frac{67}{128} & \frac{16}{128} & \frac{37}{128} \end{pmatrix}$. For each player find the normalized von Neumann-Morgenstern utility function.

You have to decide whether or not to leave your wallet at home. Suppose that, if you leave your wallet at home, with probability $p$ (with $0 < p < 1$) the potential mugger will notice that your pockets are empty and with probability $(1-p)$ he will not notice. If the potential mugger does not notice that your pockets are empty then he will nevertheless be aware of the fact that you might have left the wallet at home and he simply cannot tell.

**(c)** Represent this situation as an extensive game with imperfect information.

**(d)** Write the corresponding normal form.

**(e)** Find all the subgame-perfect equilibria (including the mixed-strategy ones, if any). [Hint: your answer should distinguish between different values of $p$].

# Appendix 6.S: Solutions to exercises

**Exercise 6.1.** It must be the case that the player under consideration has only one information set.

**Exercise 6.2.** **(a)** 53. The number of pure strategies is $2 \times 3 \times 3 \times 3 = 54$ and thus 53 probabilities are needed to specify a mixed strategy. **(b)** 7: one probability for the information set where she has two choices and two probabilities for each of the other three information sets.

**Exercise 6.3.** (1) The induced probability distribution on basic outcomes is $\begin{pmatrix} o_1 & o_2 & o_3 & o_4 & o_5 \\ \frac{71}{540} & \frac{25}{540} & \frac{213}{540} & \frac{81}{540} & \frac{150}{540} \end{pmatrix}$. Thus Player 2's expected utility is $\frac{71}{540}3 + \frac{25}{540}6 + \frac{213}{540}4 + \frac{81}{540}5 + \frac{150}{540}0 = \frac{1620}{540} = 3$. (2) The induced probability distribution on terminal nodes is $\begin{pmatrix} z_1 & z_2 & z_3 & z_4 & z_5 \\ \frac{5}{36} & \frac{10}{36} & \frac{7}{36} & \frac{4}{36} & \frac{10}{36} \end{pmatrix}$. Thus Player 2's expected payoff is $\frac{5}{36}4 + \frac{10}{36}4 + \frac{7}{36}14 + \frac{4}{36}5 + \frac{10}{36}0 = \frac{108}{36} = 3$. Not surprisingly, the same number as in part (1).

**Exercise 6.4.** The normalized von Neumann-Morgenstern utility functions are $\left(\begin{array}{c|ccccc} & o_1 & o_2 & o_3 & o_4 & o_5 \\ \hline U_1 & 1 & 0 & 0 & 0.25 & 0.75 \\ U_2 & 1 & 1 & 0.2 & 1 & 0 \\ U_3 & 0 & 1 & 0.4 & 0.6 & 0.4 \end{array}\right)$. Thus the extensive-form game is as follows:

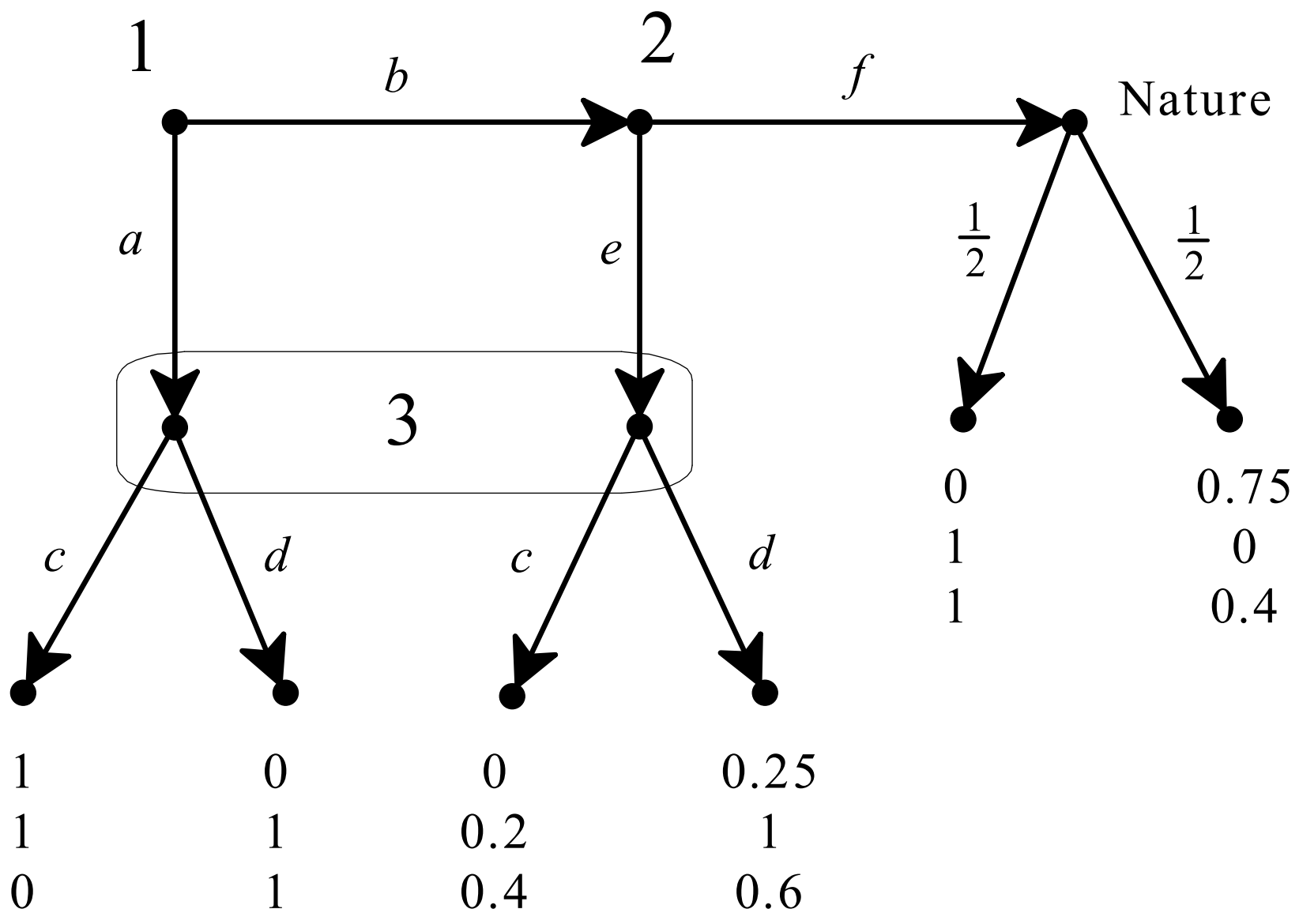

Or, in a simplified form obtained by removing the move of Nature, as follows:

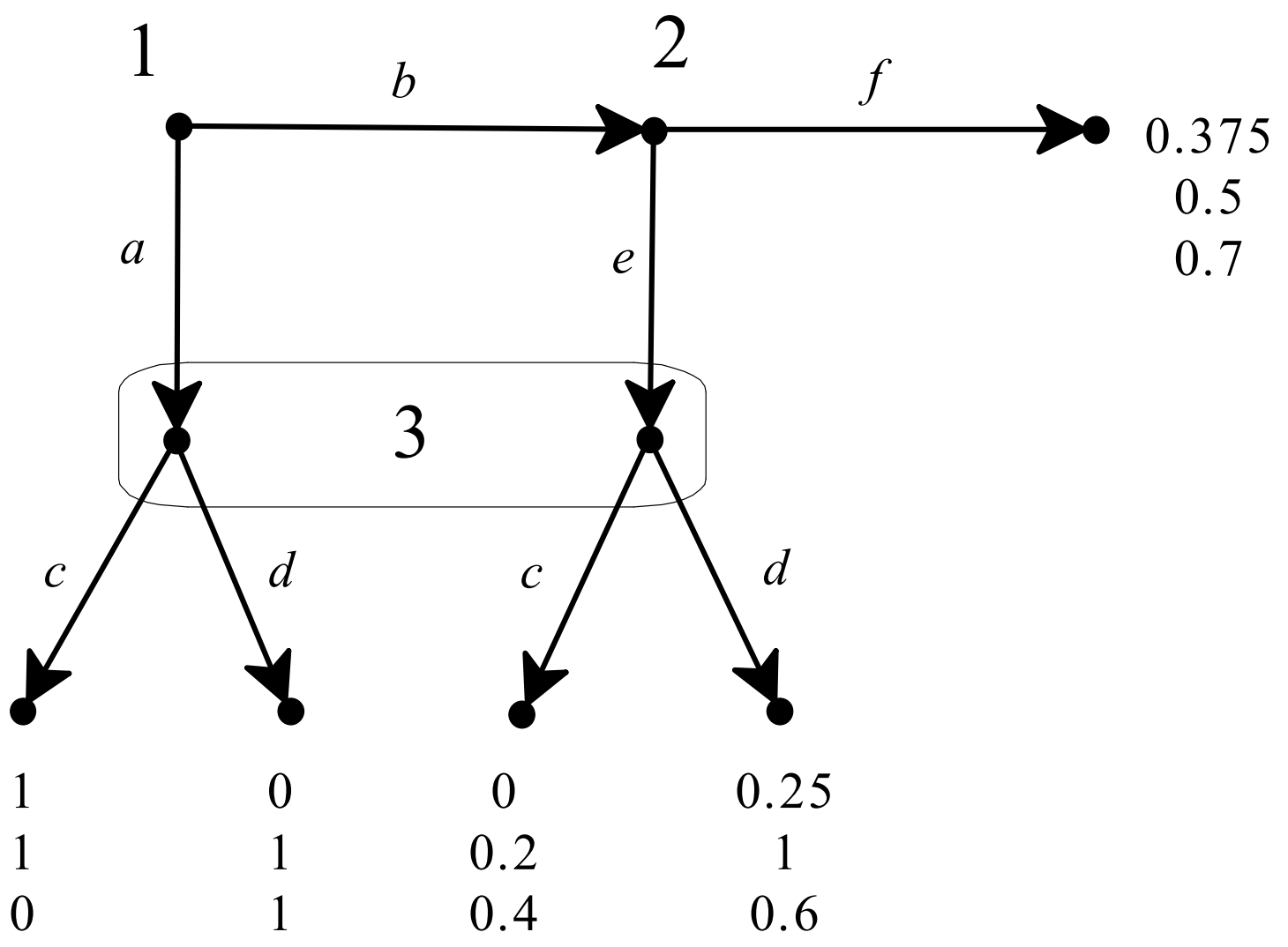

**Exercise 6.5**. **(a)** The strategic form is as follows.

| | | Player 2 $l$ | Player 2 $r$ |
|---|---|---|---|
| Player 1 | *UL* | 1 , 5 | 1 , 5 |
| | *UR* | 1 , 5 | 1 , 5 |
| | *DL* | 3 , 2 | 0 , 3 |
| | *DR* | 2 , 4 | 1 , 1 |

The pure-strategy Nash equilibria are (*UL*,*r*) and (*UR*,*r*).

**(b)** The strategic form of the proper subgame that starts at Player 2's node is:

| | | Player 2 $l$ | Player 2 $r$ |
|---|---|---|---|
| Player 1 | *L* | 3 , 2 | 0 , 3 |
| | *R* | 2 , 4 | 1 , 1 |

This game has a unique mixed-strategy Nash equilibrium given by $\left(\begin{pmatrix} L & R \\ \frac{3}{4} & \frac{1}{4} \end{pmatrix}, \begin{pmatrix} l & r \\ \frac{1}{2} & \frac{1}{2} \end{pmatrix}\right)$, yielding Player 1 an expected payoff of 1.5. Thus the unique subgame-perfect equilibrium, expressed as a behavioral-strategy profile, is $\left(\begin{pmatrix} U & D & L & R \\ 0 & 1 & \frac{3}{4} & \frac{1}{4} \end{pmatrix}, \begin{pmatrix} l & r \\ \frac{1}{2} & \frac{1}{2} \end{pmatrix}\right)$ or, expressed as a mixed-strategy profile, $\left(\begin{pmatrix} UL & UR & DL & DR \\ 0 & 0 & \frac{3}{4} & \frac{1}{4} \end{pmatrix}, \begin{pmatrix} l & r \\ \frac{1}{2} & \frac{1}{2} \end{pmatrix}\right)$.

**Exercise 6.6.** There is only one proper subgame starting from Player 2's node; its strategic-form frame is as follows:

| | | Player 3 | | |
|---|---|---|---|---|
| | | ***d*** | ***e*** | ***f*** |
| | ***a*** | **x1** | **x2** | **x3** |
| Player 2 | ***b*** | **x4** | **x5** | **x6** |
| | ***c*** | **x7** | **x8** | **x9** |

For Player 2 strategy *c* is strictly dominated by strategy *b* (she prefers $x_4$ to $x_7$, and $x_5$ to $x_8$ and $x_6$ to $x_9$) and for Player 3 strategy *f* is strictly dominated by strategy d (she prefers $x_1$ to $x_3$, and $x_4$ to $x_6$ and $x_7$ to $x_9$). Thus we can simplify the game to the following.

| | | **Player 3** | |
|---|---|---|---|
| | | *d* | *e* |
| **Player** | *a* | **x1** | **x2** |
| **2** | *b* | **x4** | **x5** |

Restricted to these outcomes the payers' rankings are: $Player\ 2: \begin{pmatrix} x_1 \\ x_4, x_5 \\ x_2 \end{pmatrix}$, $Player\ 3: \begin{pmatrix} x_2 \\ x_1, x_4 \\ x_5 \end{pmatrix}$. Let $U$ be Player 2's von Neumann-Morgenstern utility function. The we can set $U(x_1) = 1$ and $U(x_2) = 0$. Thus, since she is indifferent between $x_4$ and $x_5$ and also between $x_4$ and the lottery $\begin{pmatrix} x_1 & x_2 \\ \frac{1}{2} & \frac{1}{2} \end{pmatrix}$, $U(x_4) = U(x_5) = \frac{1}{2}$. Let $V$ be Player 3's von Neumann-Morgenstern utility function. The we can set $V(x_2) = 1$ and $V(x_5) = 0$. Thus, since she is indifferent between $x_1$ and $x_4$ and also between $x_1$ and the lottery $\begin{pmatrix} x_2 & x_5 \\ \frac{1}{2} & \frac{1}{2} \end{pmatrix}$, $V(x_1) = V(x_4) = \frac{1}{2}$. Hence the above game-frame becomes the following game:

Player 3

| Player 2 | | $d$ | $e$ |
|---|---|---|---|
| | $a$ | $1, \frac{1}{2}$ | $0, 1$ |
| | $b$ | $\frac{1}{2}, \frac{1}{2}$ | $\frac{1}{2}, 0$ |

There is no pure-strategy Nash equilibrium. Let $p$ be the probability of $a$ and $q$ the probability of $d$. Then for a Nash equilibrium we need $q = \frac{1}{2}$ and $p = \frac{1}{2}$. Hence in the subgame the outcome will be $\begin{pmatrix} x_1 & x_2 & x_4 & x_5 \\ \frac{1}{4} & \frac{1}{4} & \frac{1}{4} & \frac{1}{4} \end{pmatrix}$. Since all of these outcomes are better than $x_{10}$ for Player 1, Player 1 will play $d$. Thus the subgame-perfect equilibrium is $\left( \begin{pmatrix} d & u \\ 1 & 0 \end{pmatrix}, \begin{pmatrix} a & b & c \\ \frac{1}{2} & \frac{1}{2} & 0 \end{pmatrix}, \begin{pmatrix} d & e & f \\ \frac{1}{2} & \frac{1}{2} & 0 \end{pmatrix} \right)$.

**Exercise 6.7.** **(a)** The strategic form is as follows:

PLAYER 2

| PLAYER 1 | | $AC$ | $AD$ | $BC$ | $BD$ |
|---|---|---|---|---|---|
| | $LEG$ | 2 , 0 | 2 , 0 | 0 , 2 | 0 , 2 |
| | $LEH$ | 2 , 0 | 2 , 0 | 0 , 2 | 0 , 2 |
| | $LFG$ | 0 , 6 | 0 , 6 | 4 , 1 | 4 , 1 |
| | $LFH$ | 0 , 6 | 0 , 6 | 4 , 1 | 4 , 1 |
| | $REG$ | 1 , 4 | 4 , 3 | 1 , 4 | 4 , 3 |
| | $REH$ | 2 , 0 | 1 , 2 | 2 , 0 | 1 , 2 |
| | $RFG$ | 1 , 4 | 4 , 3 | 1 , 4 | 4 , 3 |
| | $RFH$ | 2 , 0 | 1 , 2 | 2 , 0 | 1 , 2 |

**(b)** There are no pure-strategy Nash equilibria.

**(c)** First let us solve the subgame on the left:

Player 2, Player 1

| | A | B |
|---|---|---|
| E | 2 , 0 | 0 , 2 |
| F | 0 , 6 | 4 , 1 |

There is no pure-strategy Nash equilibrium. Let us find the mixed-strategy equilibrium. Let $p$ be the probability assigned to E and $q$ the probability assigned to A. Then $p$ must be the solution to $6(1-p) = 2p + (1-p)$ and $q$ must be the solution to $2q = 4(1-q)$. Thus $p = \frac{5}{7}$ and $q = \frac{2}{3}$. The expected payoff of Player 1 is $\frac{4}{3} = 1.33$, while the expected payoff of player 2 is $\frac{12}{7} = 1.714$.

Next we solve the subgame on the right:

Player 2, Player 1

| | C | D |
|---|---|---|
| G | 1 , 4 | 4 , 3 |
| H | 2 , 0 | 1 , 2 |

There is no pure-strategy Nash equilibrium. Let us find the mixed-strategy equilibrium. Let $p$ be the probability assigned to G and $q$ the probability assigned to C. Then $p$ must be the solution to $4p = 3p + 2(1-p)$ and $q$ must be the solution to $q + 4(1-q) = 2q + (1-q)$. Thus $p = \frac{2}{3}$ and $q = \frac{3}{4}$. The expected payoff of Player 1 is $\frac{7}{4} = 1.75$. Thus the game reduces to:

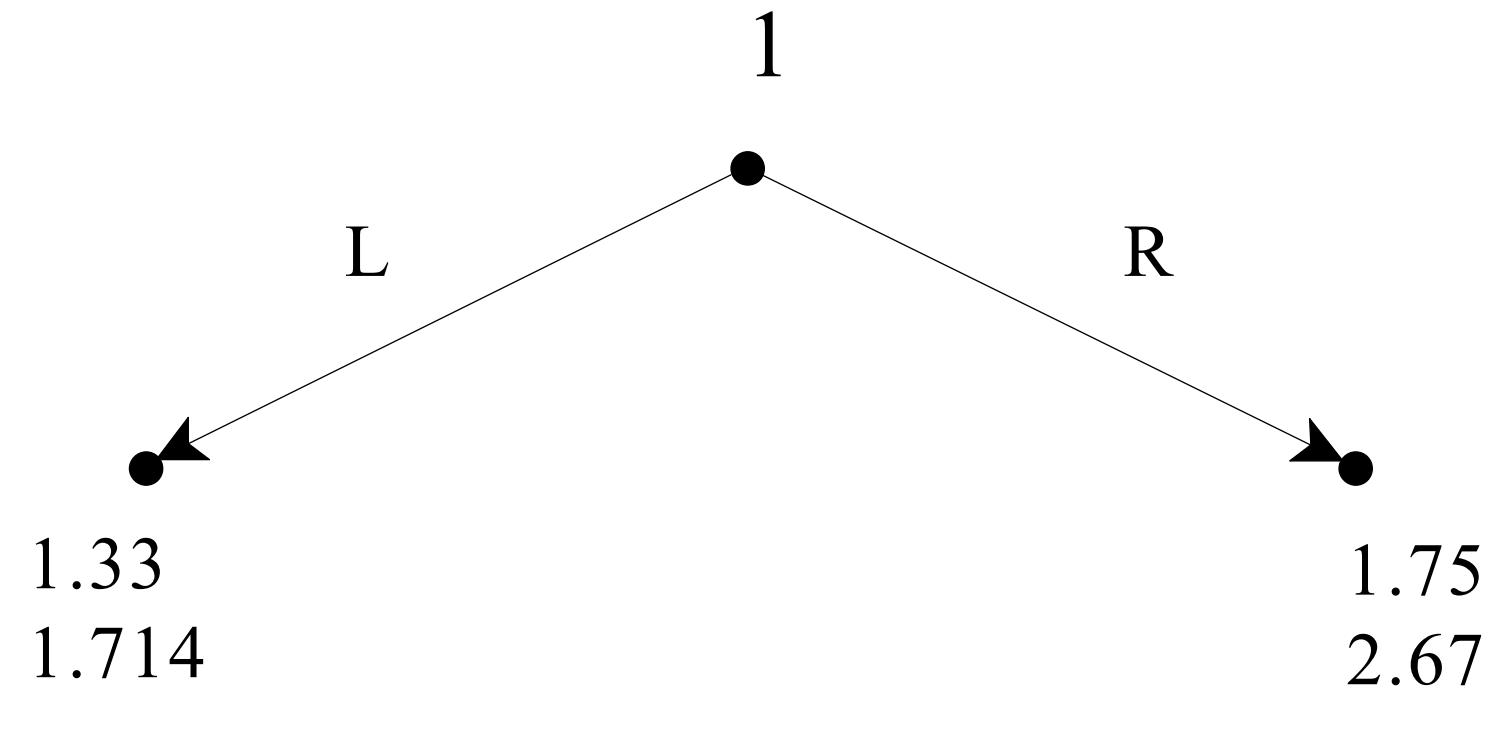

Hence the subgame-perfect equilibrium is:

$$\left( \begin{pmatrix} L & R & E & F & G & H \\ 0 & 1 & \frac{5}{7} & \frac{2}{7} & \frac{2}{3} & \frac{1}{3} \end{pmatrix}, \begin{pmatrix} A & B & C & D \\ \frac{2}{3} & \frac{1}{3} & \frac{3}{4} & \frac{1}{4} \end{pmatrix} \right)$$

**Exercise 6.8.** **(a)** The strategic form is as follows:

Player 3 chooses *a*

| Player 1 \ Player 2 | *l* | *r* |
|---|---|---|
| *R* | 1 , 0 , 0 | 1 , 0 , 0 |
| *M* | 0 , 1 , 1 | 2 , 2 , 1 |
| *L* | 1 , 1 , 0 | 0 , 0 , 0 |

Player 3 chooses *b*

| Player 1 \ Player 2 | *l* | *r* |
|---|---|---|
| *R* | 1 , 0 , 0 | 1 , 0 , 0 |
| *M* | 0 , 1 , 0 | 0 , 0 , 0 |
| *L* | 1 , 1 , 0 | 0 , 0 , 0 |

The pure-strategy Nash equilibria are highlighted: (*R*,*l*,*a*), (*M*,*r*,*a*), (*L*,*l*,*a*), (*R*,*l*,*b*), (*R*,*r*,*b*) and (*L*,*l*,*b*). They are all subgame perfect because there are no proper subgames.

**(b)** Since, for Player 3, *a* strictly dominates *b* conditional on his information set being reached, he will have to play *a* if his information set is reached with positive probability. Now, Player 3's information set is reached with positive probability if and only if player 1 plays *M* with positive probability. Thus when $\Pr(M) > 0$ the game essentially reduces to

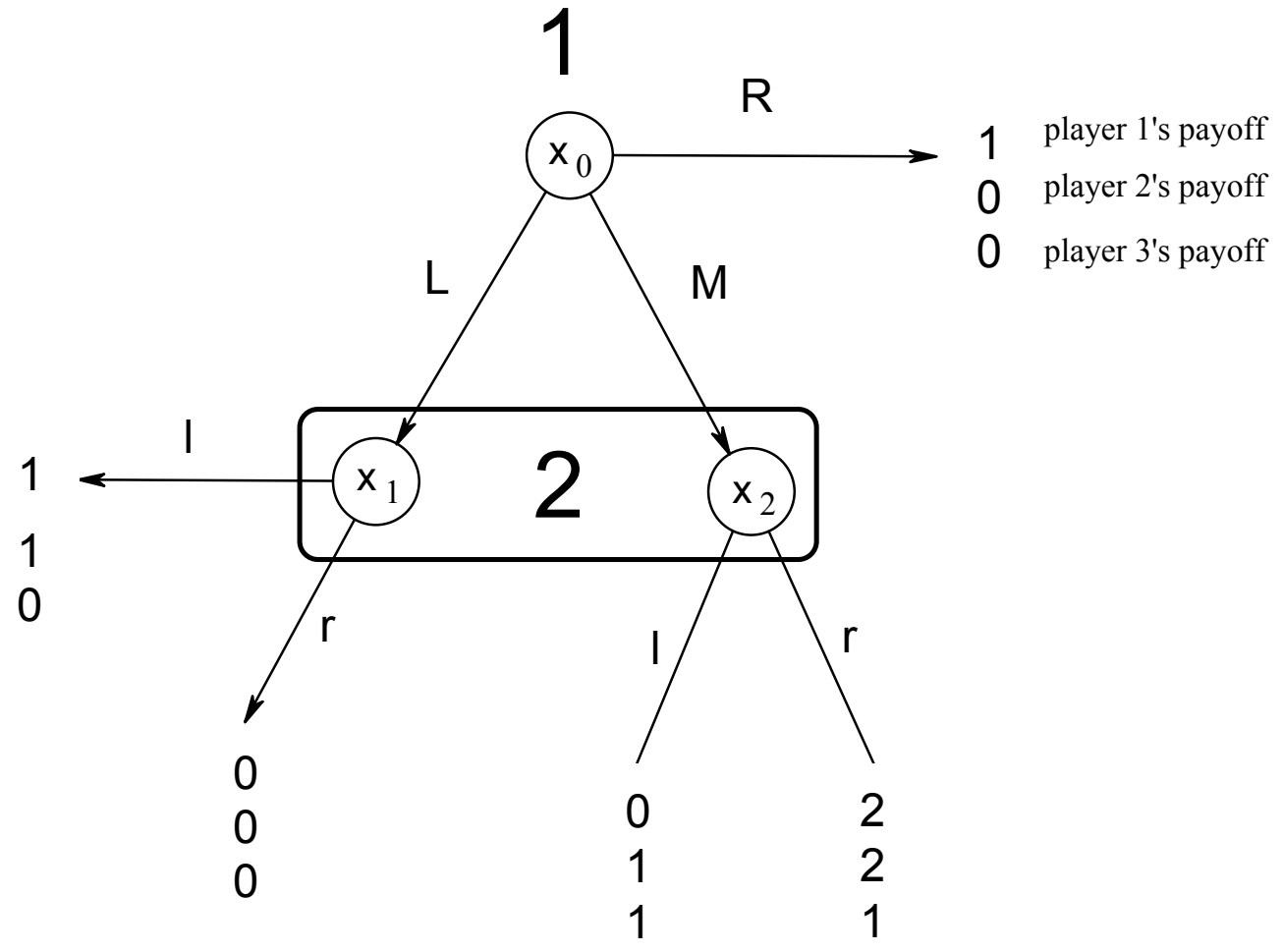

Now, in order for Player 1 to be willing to assign positive probability to *M* he must expect a payoff of at least 1 (otherwise *R* would be better) and the only way he can expect a payoff of at least 1 is if Player 2 plays *r* with probability at least $\frac{1}{2}$. Now if Player 2 plays *r* with probability greater than $\frac{1}{2}$, then *M* gives Player 1 a higher payoff than both *L* and *R* and thus he will choose *M* with probability 1, in which case Player 2 will choose *r* with probability 1 (and Player 3 will choose *a* with probability 1) and so we get the pure strategy equilibrium (*M*,*r*,*a*). If Player 2 plays *r* with probability exactly $\frac{1}{2}$ then Player 1 is indifferent between *M* and *R* (and can mix between the two), but finds *L* inferior and must give it probability 0. But then Player 2's best reply to a mixed strategy of Player 1 that assigns positive probability to *M* and *R* and zero probability to *L* is to play *r* with probability 1 (if his information set is reached it can only be reached at node $x_2$). Thus there cannot be a mixed-strategy equilibrium where Player 1 assigns to *M* probability $p$ with $0 < p < 1$ : it must be either Pr(*M*) = 0 or Pr(*M*) = 1.

**Exercise 6.9.** **(a)** Your ranking is (best at the top, worst at the bottom) $\begin{pmatrix} z_3 \\ z_1 \\ z_4 \\ z_2 \end{pmatrix}$

while the potential mugger's ranking is $\begin{pmatrix} z_4 \\ z_1, z_3 \\ z_2 \end{pmatrix}$.

**(b)** Let $U$ be your utility function. Let $U(z_3) = 1$, $U(z_1) = a$, $U(z_4) = b$ and $U(z_2) = 0$, with $0 < b < a < 1$. The expected utilities are as follows: $EU(L_1) = \frac{3}{20}a + \frac{3}{20}$, $EU(L_2) = \frac{1}{2}b$, $EU(L_3) = \frac{1}{3}$ and $EU(L_4) = \frac{1}{2}a$. From $EU(L_3) = E(L_4)$ we get that $a = \frac{2}{3}$. Substituting this into the equation $EU(L_1) = EU(L_2)$ gives $b = \frac{1}{2}$. Thus $U(z_3) = 1$, $U(z_1) = \frac{2}{3}$, $U(z_4) = \frac{1}{2}$ and $U(z_2) = 0$. Let $V$ be the mugger's utility function. Let $V(z_4) = 1$, $V(z_1) = V(z_3) = c$ and $V(z_2) = 0$ with $0 < c < 1$. The expected utilities are as follows: $EV(L_5) = \frac{1}{4}(2c+1)$ and $EV(L_6) = \frac{1}{128}(24c+37)$. Solving $EV(L_5) = EV(L_6)$ gives $c = \frac{1}{8}$. Thus, $V(z_4) = 1$, $V(z_1) = V(z_3) = \frac{1}{8}$ and $V(z_2) = 0$.

**(c)** The extensive game is as follows:

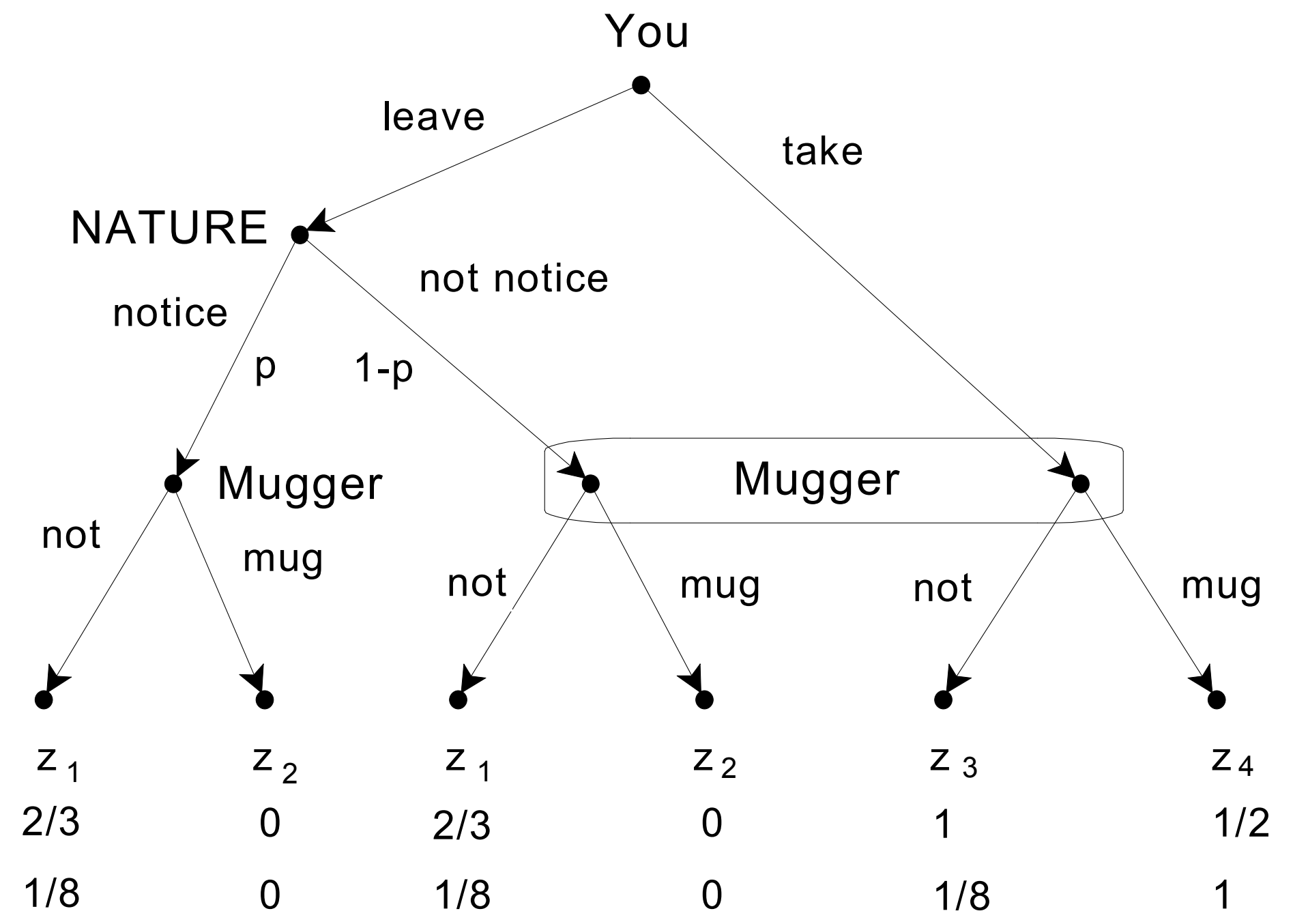

**(d)** The strategic form is as follows (for the mugger's strategy the first item refers to the left node, the second item to the information set on the right):

**Potential Mugger**

| | | NN | NM | MN | MM |
|---|---|---|---|---|---|
| **You** | L | 2/3 , 1/8 | (2/3)p , (1/8)p | 2/3(1-p) , 1/8 (1-p) | 0 , 0 |
| | T | 1 , 1/8 | 1/2 , 1 | 1 , 1/8 | 1/2 , 1 |

**(e)** At a subgame-perfect equilibrium the mugger will choose not to mug when he notices your empty pockets. Thus the normal form can be simplified as follows:

**Potential Mugger**

| | | NN | NM |
|---|---|---|---|
| **You** | L | 2/3 , 1/8 | (2/3)p , (1/8)p |
| | T | 1 , 1/8 | 1/2 , 1 |

- If $p < ¾$ then Take is a strictly dominant strategy for you and therefore there is a unique subgame-perfect equilibrium given by (Take, (Not mug,Mug)).

- If $p = ¾$ then there is a continuum of equilibria where the Mugger chooses (Not mug,Mug) with probability 1 and you choose L with probability $q$ and T with probability $(1–q)$ for any $q$ with $0 \leq q \leq \frac{28}{29}$ (obtained from $\frac{3}{32}q + 1 - q \geq \frac{1}{8}$)

- If p > ¾ then there is no pure-strategy subgame-perfect equilibrium. Let $q$ be the probability that you choose $L$ and $r$ the probability that the mugger chooses NN. Then the unique mixed strategy equilibrium is given by the solution to: $\frac{2}{3}r + \frac{2}{3}p(1-r) = r + \frac{1}{2}(1-r)$ and $\frac{1}{8} = \frac{1}{8}pq + (1-q)$ which is $q = \dfrac{7}{8-p}$ and $r = \dfrac{4p-3}{4p-1}$. Thus the unique subgame-perfect equilibrium is:

$$\left(\begin{array}{cc|cccc} L & T & NN & NM & MN & MM \\ \dfrac{7}{8-p} & \dfrac{1-p}{8-p} & \dfrac{4p-3}{4p-1} & \dfrac{2}{4p-1} & 0 & 0 \end{array}\right).$$

# PART III

## Advanced Topics I

### Knowledge, Common Knowledge, Belief

Chapter 7

# Knowledge and Common Knowledge

## 7.1 Individual knowledge

In extensive-form games with imperfect information we use information sets to represent what the players, when it is their turn to move, know about past choices. An information set of Player $i$ is a collection of nodes in the tree where it is Player $i$'s turn to move and the interpretation is that Player $i$ knows that she is making her choice at one of those nodes, but she does not know which of these nodes has actually been reached. In this chapter we extend the notion of information set to more general settings.

We start with an example. After listening to her patient's symptoms, a doctor reaches the conclusion that there are only five possible causes: (1) a bacterial infection, (2) a viral infection, (3) an allergic reaction to a drug, (4) an allergic reaction to food and (5) environmental factors. The doctor decides to do a lab test. If the lab test turns out to be positive then the doctor will be able to rule out causes (3)-(5), while a negative lab test will be an indication that causes (1) and (2) can be ruled out. To represent the doctor's possible states of information and knowledge we can use five states: *a*, *b*, *c*, *d*, and *e*. Each state represents a possible cause, as shown in Figure 7.1. We can partition the set $\{a,b,c,d,e\}$ into two sets: the set $\{a,b\}$, representing the state of knowledge of the doctor if she is informed that the lab test is positive, and the set $\{c,d,e\}$, representing the state of knowledge of the doctor if she is informed that the lab test is negative.

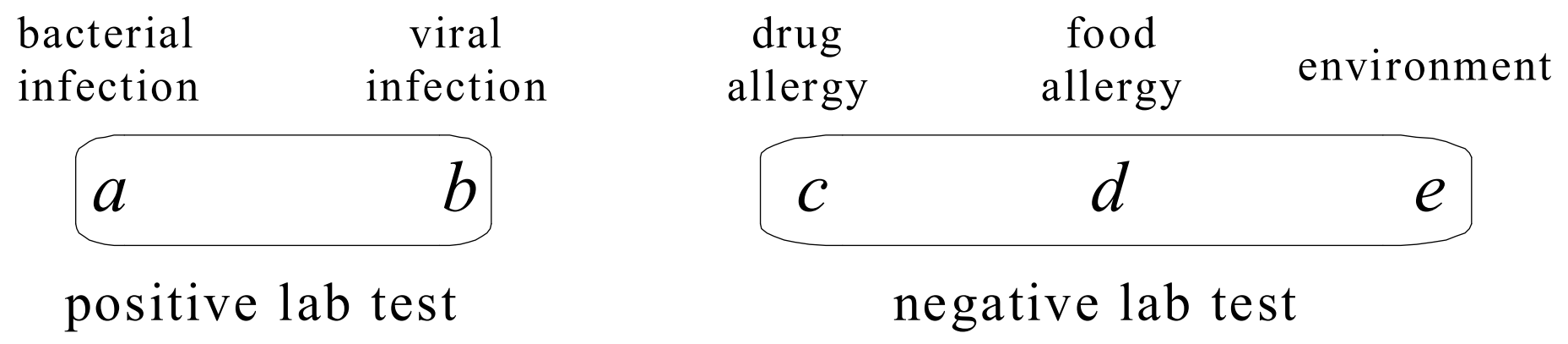

**Figure 7.1**

Consider the proposition "the cause of the patient's symptoms is either an infection or environmental factors". We can think of this proposition as the set of states $\{a,b,e\}$ where the proposition is in fact true; furthermore, we can ask the question "after receiving the result of the lab test, at which state would the doctor know the proposition represented by the set $\{a,b,e\}$?". If the true state is *a*, then – after viewing the result of the lab test – the doctor will think that it is possible that the state is either *a* or *b* and in either case it is indeed the case that the cause of the patient's symptoms is either an infection or environmental factors, so that the doctor will know this; the same is true if the true state is *b*. On the other hand, if the true state is *e* then the doctor will consider *c*, *d* and *e* as possibilities and thus will not be able to claim that she knows that the cause of the patient's symptoms is either an infection or environmental factors. Thus the answer to the question "after receiving the result of the lab test, at which state would the doctor know the proposition represented by the set $\{a,b,e\}$?" is "at states *a* and *b*". We can now turn to the general definitions.

**Definition 7.1.** Let *W* be a finite set of *states*, where each state is to be understood as a complete specification of the relevant *facts* about the world. An *information partition* is a partition $\mathcal{I}$ of *W* (that is, a collection of subsets of *W* that (1) are pairwise disjoint and (2) whose union covers the entire set *W*); the elements of the partition are called *information sets*. For every $w \in W$ we denote by $I(w)$ the information set that contains state *w*.

In the example of the doctor, $W = \{a,b,c,d,e\}$ and $\mathcal{I} = \{\{a,b\},\{c,d,e\}\}$; furthermore, $I(a) = I(b) = \{a,b\}$ and $I(c) = I(d) = I(e) = \{c,d,e\}$.

**Definition 7.2.** Let *W* be a set of states. We will call the subsets of *W* *events.* Let $\mathcal{I}$ be a partition of *W*, *E* an event (thus $E \subseteq W$) and $w \in W$ a state. We say that *at w the agent knows E* if and only if the information set to which *w* belongs is contained in *E*, that is, if and only if $I(w) \subseteq E$.

In the example of the doctor, where $W = \{a,b,c,d,e\}$ and $\mathcal{I} = \{\{a,b\},\{c,d,e\}\}$, let $E = \{a,b,d,e\}$; then at $a$ and $b$ the doctor knows $E$ because $I(a) = I(b) = \{a,b\} \subseteq E$, but at $d$ it is not true that the doctor knows $E$ because $I(d) = \{c,d,e\} \not\subseteq E$ (since $c \in I(d)$ but $c \notin E$) and, for the same reason, also at *c* and *e* it is not true that the doctor knows *E*.

Note that it is possible that there is no state where the agent knows a given event. In the doctor example, if we consider event $F = \{a,c\}$ then there is no state where the doctor knows *F*.

**Definition 7.3.** Using Definition 7.2, we can define a *knowledge operator K* on events that, given as input any event *E*, produces as output the event *KE* defined as the set of states at which the agent knows *E*. Let $2^W$ denote the set of events, that is the set of subsets of *W*.[1] Then the knowledge operator is the function $K: 2^W \to 2^W$ defined as follows: for every $E \subseteq W$, $KE = \{w \in W : I(w) \subseteq E\}$.

In the example of the doctor, where $W = \{a,b,c,d,e\}$ and $\mathcal{I} = \{\{a,b\},\{c,d,e\}\}$, let $E = \{a,b,d,e\}$ and $F = \{a,c\}$; then $KE = \{a,b\}$ and $KF = \varnothing$.

Given an event $G \subseteq W$ we denote by $\neg G$ the complement of *G*, that is, the set of states that are **not** in *G*. For example, if $W = \{a,b,c,d,e\}$ and $G = \{a,b,d\}$ then $\neg G = \{c,e\}$. Thus while *KG* is the event that the agent knows *G*, $\neg KG$ is the event that the agent does *not* know *G*. Note the important difference between event $\neg KG$ and event $K\neg G$: if $w \in \neg KG$ then at state *w* the agent does not know *G* but she might not know $\neg G$ either, that is, it may be that she considers *G* possible ($I(w) \cap G \neq \varnothing$) and she also considers $\neg G$ possible ($I(w) \cap \neg G \neq \varnothing$).[2] On the other hand, if $w \in K\neg G$ then at state *w* the agent knows that *G* is not true, because every state that she considers possible is in $\neg G$ ($I(w) \subseteq \neg G$). Thus $K\neg G \subseteq \neg KG$ but the converse inclusion does not hold. In the example of the doctor, where $W = \{a,b,c,d,e\}$ and $\mathcal{I} = \{\{a,b\},\{c,d,e\}\}$, again let $E = \{a,b,d,e\}$ (so that $\neg E = \{c\}$) and $F = \{a,c\}$ (so that $\neg F = \{b,d,e\}$); then $KE = \{a,b\}$, $\neg KE = \{c,d,e\}$, $K\neg E = \varnothing$, $KF = \varnothing$, $\neg KF = W$ and $K\neg F = \varnothing$.

[1] If *W* contains *n* elements, then there are $2^n$ subsets of *W*, hence the notation $2^W$. For example, if $W = \{a,b,c\}$ then $2^W = \{\varnothing,\{a\},\{b\},\{c\},\{a,b\},\{a,c\},\{b,c\},\{a,b,c\}\}$.

[2] In the example of the doctor, where $W = \{a,b,c,d,e\}$ and $\mathcal{I} = \{\{a,b\},\{c,d,e\}\}$, if $F = \{a,c\}$ (so that $\neg F = \{b,d,e\}$) then $KF = \varnothing$ and $K\neg F = \varnothing$; for instance, if the true state is *a* then the doctor considers *F* possible (because her information set is $I(a) = \{a,b\}$ and $I(a) \cap F = \{a\} \neq \varnothing$) but she also considers $\neg F$ possible (because $I(a) \cap \neg F = \{b\} \neq \varnothing$).

Note the interesting fact that, since we can apply the knowledge operator to any event, we can also compute the event that the agent knows that she knows $E$ (that is, the event $K(KE)$, which we will denote more simply as $KKE$) and the event that the agent knows that she does not know $E$ (that is, the event $K\neg KE$). Continuing the example of the doctor, where $W = \{a,b,c,d,e\}$, $\mathcal{I} = \{\{a,b\},\{c,d,e\}\}$ and $E = \{a,b,d,e\}$ $KKE$ is the set of states where the agent knows event $KE = \{a,b\}$; thus $KKE = \{a,b\}$. Furthermore, since $\neg KE = \{c,d,e\}$, $K\neg KE = \{c,d,e\}$. As noted in the following remark, this is not a coincidence.

**Remark 7.1.** The knowledge operator $K : 2^W \to 2^W$ satisfies the following properties (which you are asked to prove in Exercise 7.5): for every event $E \subseteq W$,

- $KE \subseteq E$
- $KKE = KE$
- $K\neg KE = \neg KE$

This is a good time to test your understanding of the concepts introduced in this section, by going through the exercises in Section 7.E.1 of Appendix 7.E at the end of this chapter.

# 7.2 Interactive knowledge

We can now extend our analysis to the case of several agents and talk about not only what an individual knows about relevant facts but also about what she knows about what other individuals know and what they know about what she knows, etc. There is an entertaining episode of the TV series *Friends* in which Phoebe and Rachel reason about whether Chandler and Monica know that they (=Phoebe and Rachel) know that Chandler and Monica are having an affair. Search for the string "Friends-They Don't Know That We Know They Know We Know" on www.youtube.com or point your browser to https://www.youtube.com/watch?v=LUN2YN0bOi8

Again we start with a set of states $W$, where each state represents a complete description of the relevant *facts*. Let there be $n$ individuals. To represent the possible states of mind of each individual we use an information partition: $I_i$ denotes the partition of individual $i \in \{1,...,n\}$. As before, we call the subsets of $W$ events. Using Definition 7.3 we can define a knowledge operator for every individual. Let $K_i$ be the knowledge operator of individual $i$; thus, for every event $E \subseteq W$, $K_iE = \{w \in W : I_i(w) \subseteq E\}$. Now consider an event $E$ and an individual, say Individual 1; since $K_1E$ is an event (it is the set of states where Individual 1 knows event $E$), we can compute the event $K_2K_1E$, which is the event that Individual 2 knows event $K_1E$, that is, the event that 2 knows that 1 knows $E$. But there is no need to stop there: we can also compute the event $K_3K_2K_1E$ (the event that 3 knows that 2 knows that 1 knows $E$) and the event $K_1K_3K_2K_1E$ (the event that 1 knows that 3 knows that 2 knows that 1 knows $E$), etc. A few examples will be useful.

We begin with an abstract example, without interpreting the states in terms of specific facts. Let the set of states be $W = \{a, b, c, d, e, f, g, h\}$. There are three individuals, Ann, Bob and Carol, with the information partitions shown in Figure 7.2.

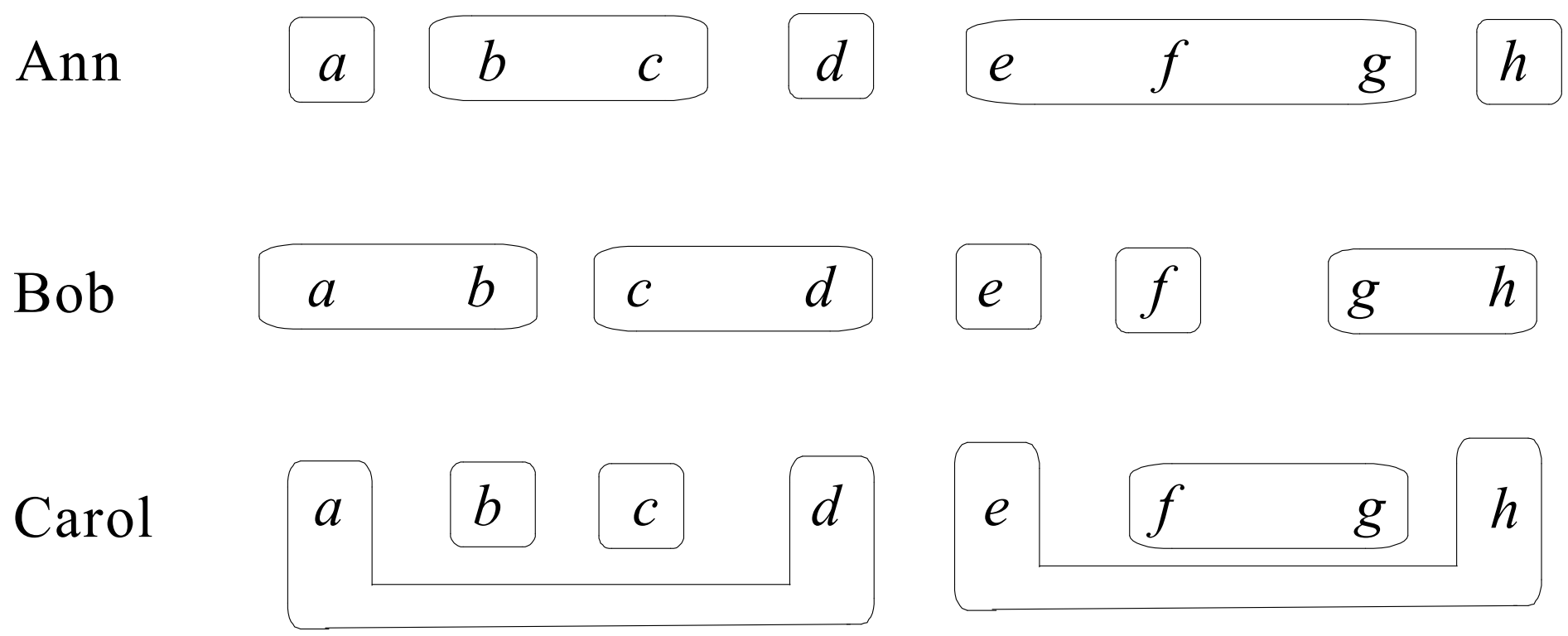

**Figure 7.2**

Consider the event $E = \{a, b, c, f, g\}$. Let us compute the following events: $K_{Ann}E$ (the event that Ann knows $E$), $K_{Bob}E$, $K_{Carol}E$, $K_{Carol}K_{Ann}E$ (the event that Carol knows that Ann knows $E$), $K_{Bob}K_{Carol}K_{Ann}E$ and $K_{Ann}\neg K_{Bob}K_{Carol}E$ (the event that Ann knows that Bob does not know that Carol knows $E$). All we need to do is apply Definition 7.3. First of all,

- $K_{Ann}E = \{a,b,c\}$

(for example, $b \in K_{Ann}E$ because $I_{Ann}(b) = \{b,c\}$ and $\{b,c\} \subseteq E$, while $f \notin K_{Ann}E$ because $I_{Ann}(f) = \{e,f,g\}$ and $\{e,f,g\}$ is not a subset of $E$). Similarly,

- $K_{Bob}E = \{a,b,f\}$, and
- $K_{Carol}E = \{b,c,f,g\}$.

To compute $K_{Carol}K_{Ann}E$ we need to find the set of states where Carol knows event $\{a,b,c\}$, since $K_{Ann}E = \{a,b,c\}$. Thus

- $K_{Carol}\underbrace{K_{Ann}E}_{=\{a,b,c\}} = \{b,c\}$.

Hence $K_{Bob}\underbrace{K_{Carol}K_{Ann}E}_{=\{b,c\}} = \varnothing$, that is, there is no state where Bob knows that Carol knows that Ann knows $E$. To compute $K_{Ann}\neg K_{Bob}K_{Carol}E$ first we start with $K_{Carol}E$, which we have already computed: $K_{Carol}E = \{b,c,f,g\}$; then we compute $K_{Bob} \underbrace{K_{Carol}E}_{=\{b,c,f,g\}}$, which is $\{f\}$: $K_{Bob}K_{Carol}E = \{f\}$; then we take the complement of this: $\neg K_{Bob}K_{Carol}E = \{a,b,c,d,e,g,h\}$ and finally we compute $K_{Ann}$ of this event:

- $K_{Ann} \underbrace{\neg K_{Bob}K_{Carol}E}_{=\{a,b,c,d,e,g,h\}} = \{a,b,c,d,h\}$.

Thus, for example, at state *a* it is true that Ann knows that Bob does not know that Carol knows *E*, while at state *e* this is not true.

Next we discuss a concrete example. The professor in a game theory class calls three students to the front of the room, shows them a large transparent box that contains many hats, some red and some white (there are no other colors) and tells them the following:

> "I will blindfold you and put a hat on each of you, then I will remove the box from the room and, after that, I will remove your blindfolds, so that each of you can see the color(s) of the hats worn by the other two students, but you will not be able to see the color of your own hat. Then I will ask you questions starting with Student 1 then Student 2 then Student 3 then back to Student 1 and so on."

After having placed the hats and removed the blindfolds, the professor asks Student 1 "do you know the color of your own hat?". She replies "No". Then he asks Student 2 the same question: "do you know the color of your own hat?". Student 2 says "No". Then he asks Student 3 and she says "No". Then he asks the same question again to Student 1 and again the answer is "No", and so on. After asking the same question over and over and always hearing the answer "No" he gets tired and tells them "I'll give you a piece of information: **I did not pick three white hats**". He then resumes the questioning: first he asks Student 1 "do you know the color of your own hat?". She replies "No". Then he asks Student 2 the same question: "do you know the color of your own hat?". Student 2 says "No". Then he asks Student 3 and she says "Yes I do!" What color hat does she have? What color hats do Students 1 and 2 have?

To answer these questions, we begin by defining the set of possible states. We can think of a state as a triple $(C_1, C_2, C_3)$, where $C_i \in \{R,\ W\}$ is the color of the hat of Student *i* (*R* means Red and *W* means White). Thus, for example, $(R, W, R)$ is the state where Students 1 and 3 have a red hat, while Student 2 has a white hat. The possible states of information of the three students *before the professor announces that he did not pick three white hats* are represented by the information partitions shown in Figure 7.3, where we have connected with a line states that are in the same information set. Whatever the state (that is, whatever hats the professor picks), each student is only uncertain about the color of her own hat: she can see, and thus knows, the colors of the hats of the other two students. Thus each information set contains only two elements.

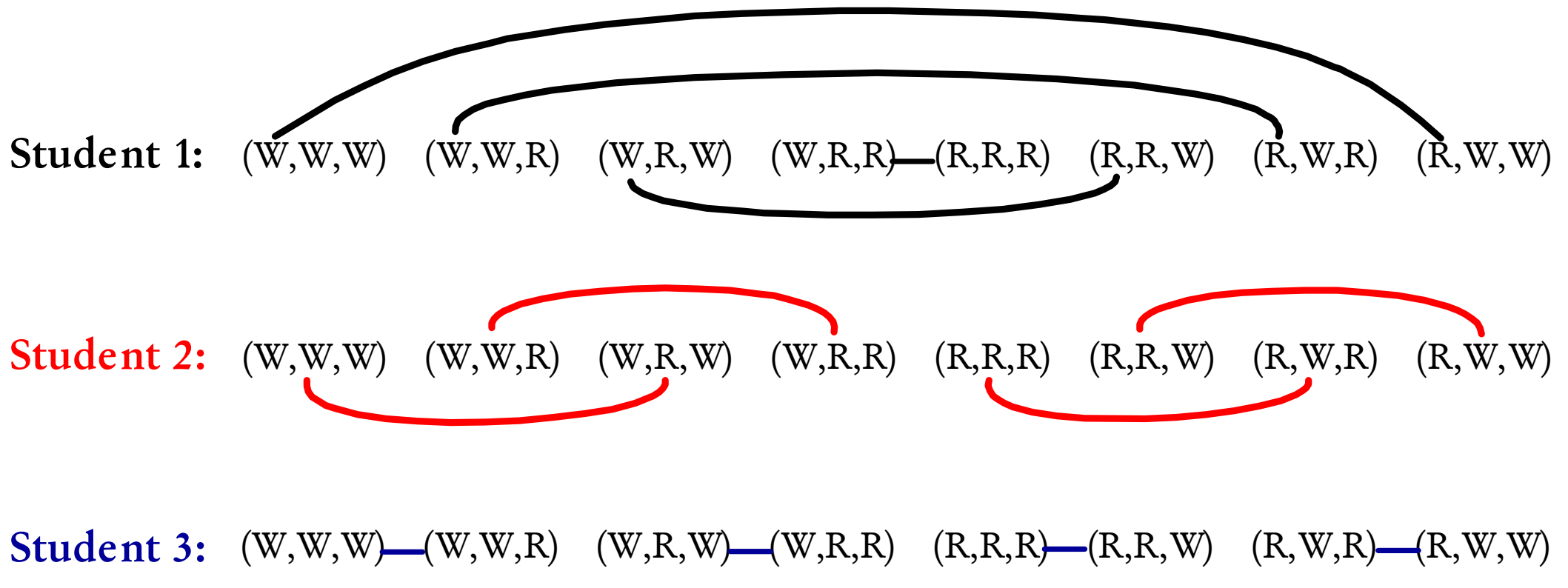

**Figure 7.3**

Consider a particular state, say the state where all three hats are red: $(R,R,R)$. At that state, obviously, each student knows that not all hats are white: he/she can actually see two red hats. Furthermore, each student knows that every other student knows that not all hats are white.[3] Take, for example, Student 1. She sees that the hats of the other two students are red and thus she reasons that Student 2 also sees that the hat of Student 3 is red and hence Student 2 knows that not all hats are white (similarly, she reasons that Student 3 knows that not all hats are white). But does Student 1 know that Student 2 knows that Student 3 knows that not all hats are white? The answer is No. Seeing two red hats, Student 1 must consider it possible that her own hat is white, in which case Student 2 would, like Student 1, see a red hat on Student 3 but (unlike Student 1) would also see a white hat on

[3] It is unfortunate that many people would use, incorrectly, the expression "all hats are not white" to mean that "it is not the case that all hats are white". The latter expression is equivalent to "at least one hat is red (possibly one, possibly two, possibly all three)", while the former is equivalent to "every hat is not white", that is, in this context, "every hat is red". [It is also unfortunate that when asked "how are you?" most people answer "I am good", instead of "I am well" or "I am fine"!]

Student 1; thus Student 2 would have to consider the possibility that her own hat is white in which case, putting herself in the shoes of Student 3, would reason that Student 3 would see two white hats and consider it possible that his own hat was also white, that is, consider it possible that all hats were white. We can see this more clearly by using the information partitions and the associated knowledge operators. To simplify matters, let us assign names to the states and rewrite the information partitions using these names, as shown in Figure 7.4.

| (W,W,W) | (W,W,R) | (W,R,W) | (W,R,R) | (R,R,R) | (R,R,W) | (R,W,R) | (R,W,W) |
|---|---|---|---|---|---|---|---|
| $a$ | $b$ | $c$ | $d$ | $e$ | $f$ | $g$ | $h$ |

Student 1's partition: $\{\{a,h\},\{b,g\},\{c,f\},\{d,e\}\}$

Student 2's partition: $\{\{a,c\},\{b,d\},\{e,g\},\{f,h\}\}$

Student 3's partition: $\{\{a,b\},\{c,d\},\{e,f\},\{g,h\}\}$

**Figure 7.4**

The proposition "not all hats are white" corresponds to event $E = \{b,c,d,e,f,g,h\}$ (the set of all states, excluding only state *a*). Using Definition 7.3 we get the following events:

$$\begin{aligned}
&K_1E = \{b,c,d,e,f,g\},\\
&K_2E = \{b,d,e,f,g,h\},\\
&K_3E = \{c,d,e,f,g,h\}\\
&K_1K_2E = K_2K_1E = \{b,d,e,g\},\\
&K_1K_3E = K_3K_1E = \{c,d,e,f\},\\
&K_2K_3E = K_3K_2E = \{e,f,g,h\}
\end{aligned}$$

Note that the intersection of all these events is the singleton set $\{e\}$. Thus at state *e* (where all the hats are red), and only at state *e*, everybody knows that not all hats are white and everybody knows that everybody knows that not all hats are white. Proceeding one step further, we have that

$$K_1K_2K_3E = K_1K_3K_2E = K_2K_1K_3E = K_2K_3K_1E = K_3K_1K_2E = K_3K_2K_1E = \varnothing .$$

Thus there is no state (not even *e*) where a student knows that another student knows that the third student knows that not all hats are white.

Let us now continue with the formal analysis of the story of the three hats. At some stage the professor makes the public announcement "I did not choose three white hats" (that is, he announces event $E = \{b, c, d, e, f, g, h\}$). This announcement makes it commonly known that state *a* is to be ruled out. Thus, after the announcement, the information partitions are reduced to the ones shown in Figure 7.5, obtained from Figure 7.3 by deleting the state $(W,W,W)$.

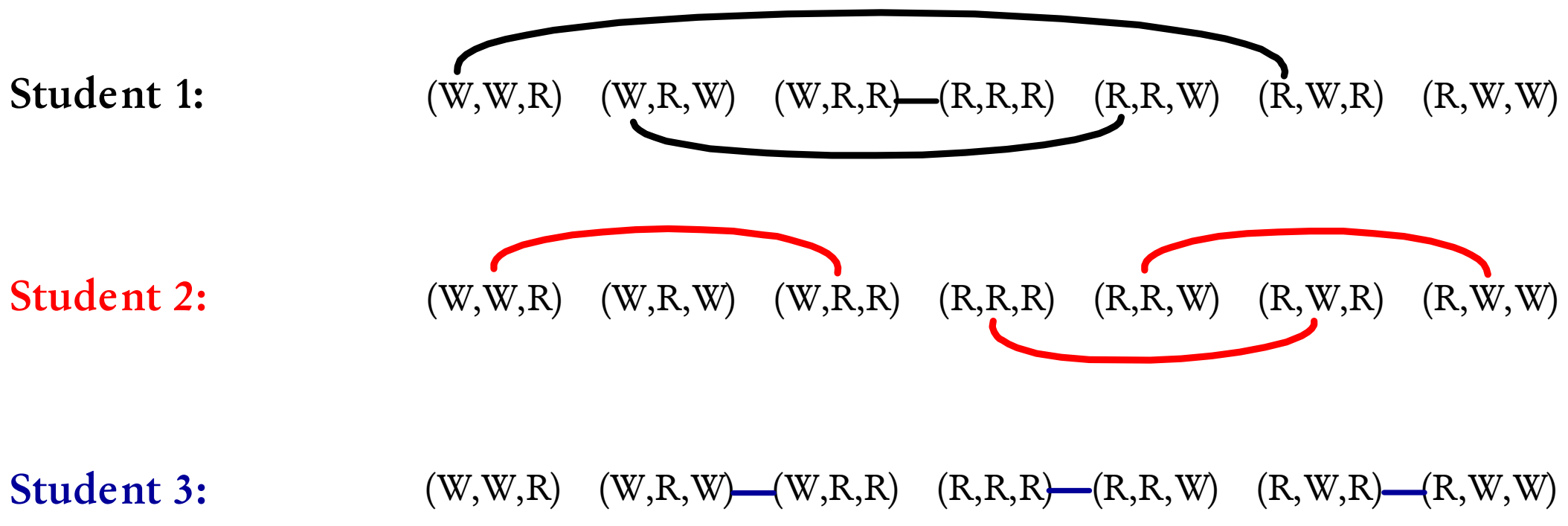

**Figure 7.5**

Note that, at this stage, if the true state is $(R,W,W)$ Student 1 knows that her hat is red (she sees two white hats and, having been informed that the professor did not choose three white hats, she can deduce that hers must be red). Similarly, if the true state is $(W,R,W)$, Student 2 knows that her own hat is red and, if the true state is $(W,W,R)$, Student 3 knows that her own hat is red. According to the story, after announcing that not all hats are white, the professor first asks Student 1 if she knows the color of her hat and she answers No. From this answer everybody can deduce that the state is not $(R,W,W)$ and thus this state can be deleted and the information partitions reduce to the ones shown in Figure 7.6.

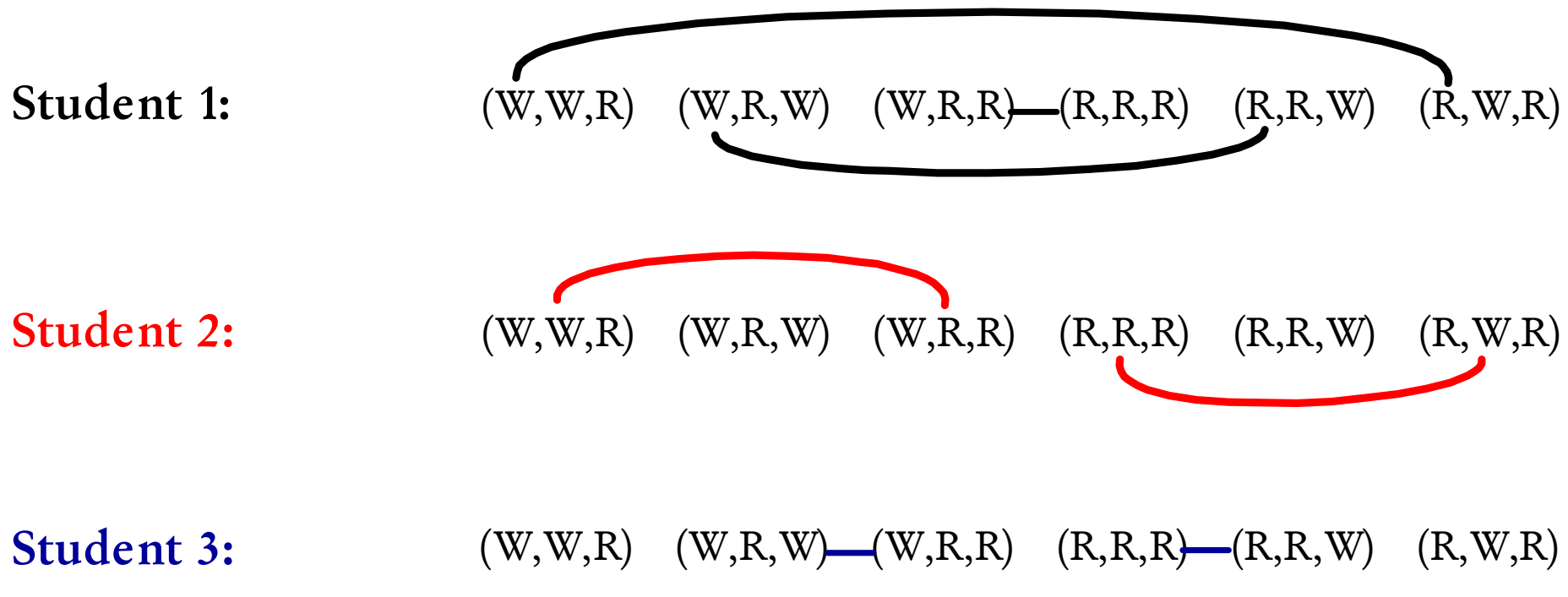

**Figure 7.6**

Now, according to the story, Student 2 is asked whether she knows the color of her hat. Let us see what the possibilities are.

1. Student 2 will answer **Yes** if state is either $(W,R,W)$ (in fact, in this case, she knew even before hearing Student 1's answer) or $(R,R,W)$ (in this case, before hearing Student 1's answer, she thought the state might be either $(R,R,W)$ or $(R,W,W)$ but then, after hearing Student 1's answer, she was able to eliminate $(R,W,W)$ as a possibility). In either of these two states Student 2's hat is red and thus she knows that her hat is red. Furthermore, upon hearing Student 2 say Yes, Student 3 learns that the state is either $(W,R,W)$ or $(R,R,W)$ and in both of these states his hat is white, thus he acquires the knowledge that his own hat is white.

2. In each of the remaining states, namely $(W,W,R)$, $(W,R,R)$, $(R,R,R)$ and $(R,W,R)$, Student 2 will answer **No**. Then everybody learns that the state is neither $(W,R,W)$ nor $(R,R,W)$ and thus the information partitions reduce to the ones shown in Figure 7.7.

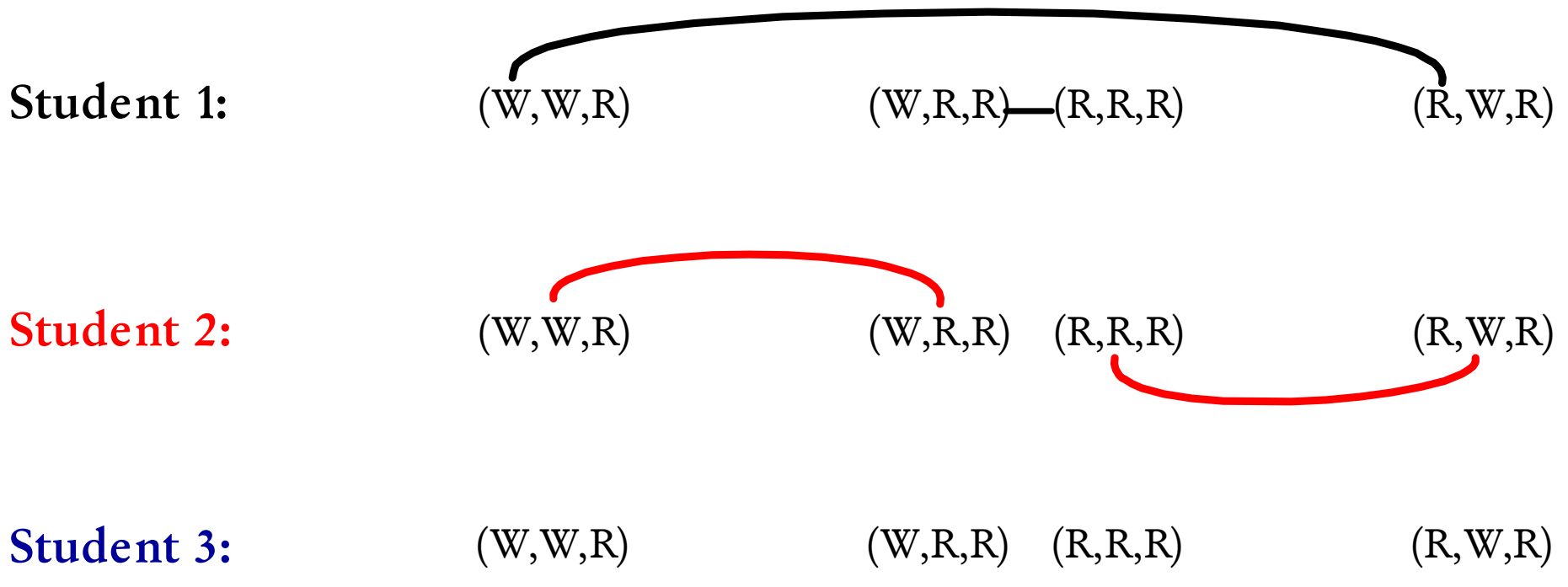

**Figure 7.7**

> Note that each of the remaining states is now a singleton information set for Student 3 and thus, upon hearing Student 2 say No, he learns what the state is: in particular he learns that his own hat is red (at each of these states Student 3's hat is red). In the original story, Students 1 and 2 answered No and Student 3 answered Yes and thus we fall in this second case.

Now consider a blindfolded witness, who cannot see the color of anybody's hat, but knows that they are either red or white and hears the professor's announcement and subsequent questions, as well as the answers by the three students. What would the blindfolded witness learn about the true state? The initial information set of the witness would consist of the set of all states; then, (1) after the announcement of the professor, he would be able to eliminate state $(W,W,W)$, (2) after the negative answer of Student 1 he would be able to eliminate state $(R,W,W)$ and (3) after the negative answer of Student 2 he would be able to eliminate states $(W,R,W)$ and $(R,R,W)$. The affirmative answer of Student 3 is not informative, because in each of these states Student 3 knows that the color of her own hat is red. Thus at the end of the exchange the witness's information set is

$$\{(W,W,R),(W,R,R),(R,R,R),(R,W,R)\}.$$

Hence all the witness knows is that the hat of Student 3 is red (on the other hand, Student 3 knows the color of all the hats, because she has figured out the color of her own hat and can see the hats of the other two students).

Let us now focus on state $(R,R,R)$ where all hats are red. Initially, before the professor makes the announcement, no student is ever able to figure out the color of her own hat, no matter how many times the students are asked. However, as we saw, once the professor announces that not all hats are white, then after Students 1 and 2 reply negatively to the question whether they know the color of their own hat, Student 3 is able to deduce that her hat is red. Thus the professor's announcement provides crucial information. This, however, seems puzzling, because the professor merely tells the students *what they already knew*: each student, seeing two red hats, knows that not all hats are white (furthermore, as we saw above, each student also knows that every other students knows this). So how can giving the students a piece of information that they already possess make any difference? The answer is that the professor's *public* announcement makes it a matter of *common knowledge* that not all hats are white. Indeed, we saw above that – at the beginning – if the true state is $(R,R,R)$, although everybody knows that not all hats are white and also everybody knows that everybody knows that not all hats are white, it is not the case that Student 1 knows that Student 2 knows that Student 3 knows that not all hats are white. Thus it is not common knowledge that not all hats are white. The notion of common knowledge is discussed in the next section.

This is a good time to test your understanding of the concepts introduced in this section, by going through the exercises in Section 7.E.2 of Appendix 7.E at the end of this chapter.

# 7.3 Common knowledge

Common knowledge is the strongest form of interactive knowledge: an event $E$ is common knowledge if everybody knows $E$ and everybody knows that everybody knows $E$ and everybody knows that everybody knows that everybody knows $E$, and so on. For example, in the case of two individuals, we say that at state $w$ event $E$ is common knowledge if

$$w \in K_1E \cap K_2E \cap K_1K_2E \cap K_2K_1E \cap K_1K_2K_1E \cap K_2K_1K_2E \cap \dots$$

We denote by $CKE$ the event that (that is, the set of states where) event $E$ is common knowledge. Thus, in the case of two individuals,

$$CKE = K_1E \cap K_2E \cap K_1K_2E \cap K_2K_1E \cap K_1K_2K_1E \cap K_2K_1K_2E \cap \dots$$

Given the definition of common knowledge, it may seem impossible to check if an event is common knowledge, because it requires checking an infinite number of conditions. We will see that, on the contrary, it is very easy to determine if an event is common knowledge at any state. We begin with an example.

**Example 7.1**. Abby proposes the following to Bruno and Caroline.

> "Tomorrow I will put you in two separate rooms, so that there will be no possibility of communication between you. I will then pick randomly an even number from the set {2,4,6}. Let $n$ be that number. Then I will write the number $n-1$ on a piece of paper and the number $n+1$ on another piece of paper, shuffle the two pieces of paper and hand one to Bruno and the other to Caroline. For example, if I happen to pick the number 6, then I will write 5 on a piece of paper and 7 on another piece of paper, shuffle and give one piece of paper to each of you. After seeing the number handed to you, each of you will then write a pair of numbers on your piece of paper and return it to me. If (1) you write the same pair of numbers and (2) at least one of the two numbers is equal to the number that was actually given to Bruno then I will give \$1,000 to each of you, otherwise each of you will give me \$1,000."

Should Bruno and Caroline accept to play this game? They can agree today on how they should act tomorrow under various contingencies, bearing in mind that they will be unable to communicate with each other tomorrow.

We will see that Bruno and Caroline should indeed accept, because they have a strategy that *guarantees* that they will each get \$1,000 from Abby. The first step is to represent the set of possible states and the information partitions. We will describe a state by a triple *abc,* where *a* is the number picked by Abby, *b* is the number given to Bruno and *c* is the number given to Caroline. Bruno only observes *b* and Caroline only observes *c*. Thus the information partitions are as shown in Figure 7.8.

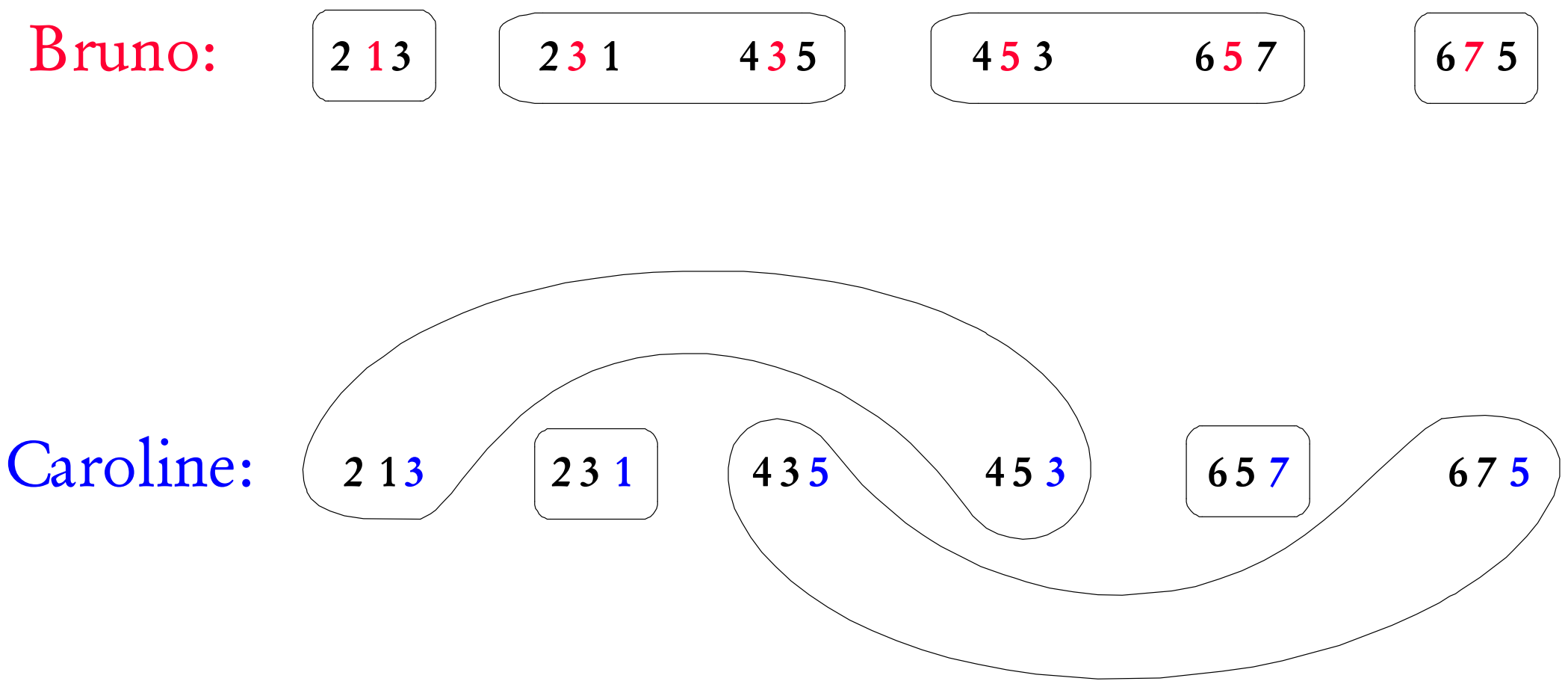

**Figure 7.8**

Let us see if the following would be a successful strategy for Bruno and Caroline: "if Bruno gets a 1 or a 3 we will write the pair (1,3) and if Bruno gets a 5 or a 7 we will write the pair (5,7)". Consider the event "Bruno gets a 1 or a 3"; call this event *E*. Then $E = \{213,\ 231,\ 435\}$. If this event occurs (e.g. because the actual state is 213), will it be common knowledge between Bruno and Caroline that it occurred? It is straightforward to check that (*B* stands for Bruno and *C* for Caroline) $K_B E = E$, $K_C E = \{231\}$ and thus $K_B K_C E = \varnothing$. Hence, while Bruno, of course, will know if he gets a 1 or a 3, Caroline might know (if the state is 231) or she might not know (that happens if the state is 213 or 435), but Bruno will never know that Caroline knows. Thus event *E* is far from being common knowledge, if it occurs. It is easy to check that the same is true of the event "Bruno gets a 5 or a 7". In order to be successful, a coordination strategy must be based on events that, when they occur, are commonly known.

Let us now consider an alternative strategy, namely

"if Bruno gets a 1 or a 5 let us write the pair (1,5) (♦)
and if Bruno gets a 3 or a 7 let us write the pair (3,7)"

Let *F* be the event "Bruno gets a 1 or a 5", that is, $F = \{213,\ 453,\ 657\}$. Then $K_B F = F$ and $K_C F = F$, so that $K_B \underbrace{K_C F}_{=F} = F$, $K_C \underbrace{K_B F}_{=F} = F$, $K_B \underbrace{K_C K_B F}_{=F} = F$, $K_C \underbrace{K_B K_C F}_{=F} = F$ and so on. Hence $CKF = F$, that is, if event $F$ occurs then it is common knowledge between Bruno and Caroline that it occurred. Similarly, letting $G = \{231, 435, 675\}$ be the event "Bruno gets a 3 or a 7" we have that $CKG = G$. Hence strategy (♦) will be a successful coordination strategy, since the conditioning events, when they occur, are common knowledge between Bruno and Caroline.

In the above example we were able to show directly that an event was common knowledge at a given state (for example, we showed that for every $w \in F$, $w \in CKF$). We now show a faster method for computing, for every event *E*, the event $CKE$. The crucial step is to derive from the individuals' information partitions a new partition of the set of states which we call the *common knowledge partition.*

**Definition 7.4.** Consider a set of states *W* and *n* partitions $\mathcal{I}_1, \mathcal{I}_2, ..., \mathcal{I}_n$ of *W*. As usual, if *w* is a state, we denote by $I_i(w)$ the element (information set) of the partition $\mathcal{I}_i$ that contains *w*. Given two states $w, w' \in W$, we say that

- $w'$ is *reachable from w in one step* if there exists an $i \in \{1,...,n\}$ such that $w' \in I_i(w)$.
- $w'$ is *reachable from w in two steps* if there exists a state $x \in W$ such that *x* is reachable from *w* in one step and $w'$ is reachable from *x* in one step.[4]
- In general, $w'$ is *reachable from w in m steps* ($m \geq 1$) if there is a sequence of states $\langle w_1, w_2, ..., w_m \rangle$ such that (1) $w_1 = w$, (2) $w_m = w'$ and (3) for every $k = 2,...,m$, $w_k$ is reachable from $w_{k-1}$ in one step.

Finally, we say that $w'$ is *reachable from w* if, for some $m \geq 1$, $w'$ is *reachable from w in m steps.*

[4] Thus $w'$ is *reachable from w in two steps* if there exist $x \in W$ and $i, j \in \{1,...,n\}$ such that $x \in I_i(w)$ and $w' \in I_j(x)$.

**Definition 7.5.** Consider a set of states $W$ and $n$ partitions $\mathcal{I}_1, \mathcal{I}_2, ..., \mathcal{I}_n$ of $W$. Given a state $w$, the *common knowledge information set that contains w*, denoted by $I_{CK}(w)$, is the set of states reachable from $w$. The *common knowledge information partition* is the collection of common knowledge information sets.

**Example 7.1 continued.** Let us go back to Example 7.1, where there are two individuals, Bruno and Caroline, with the following information partitions:

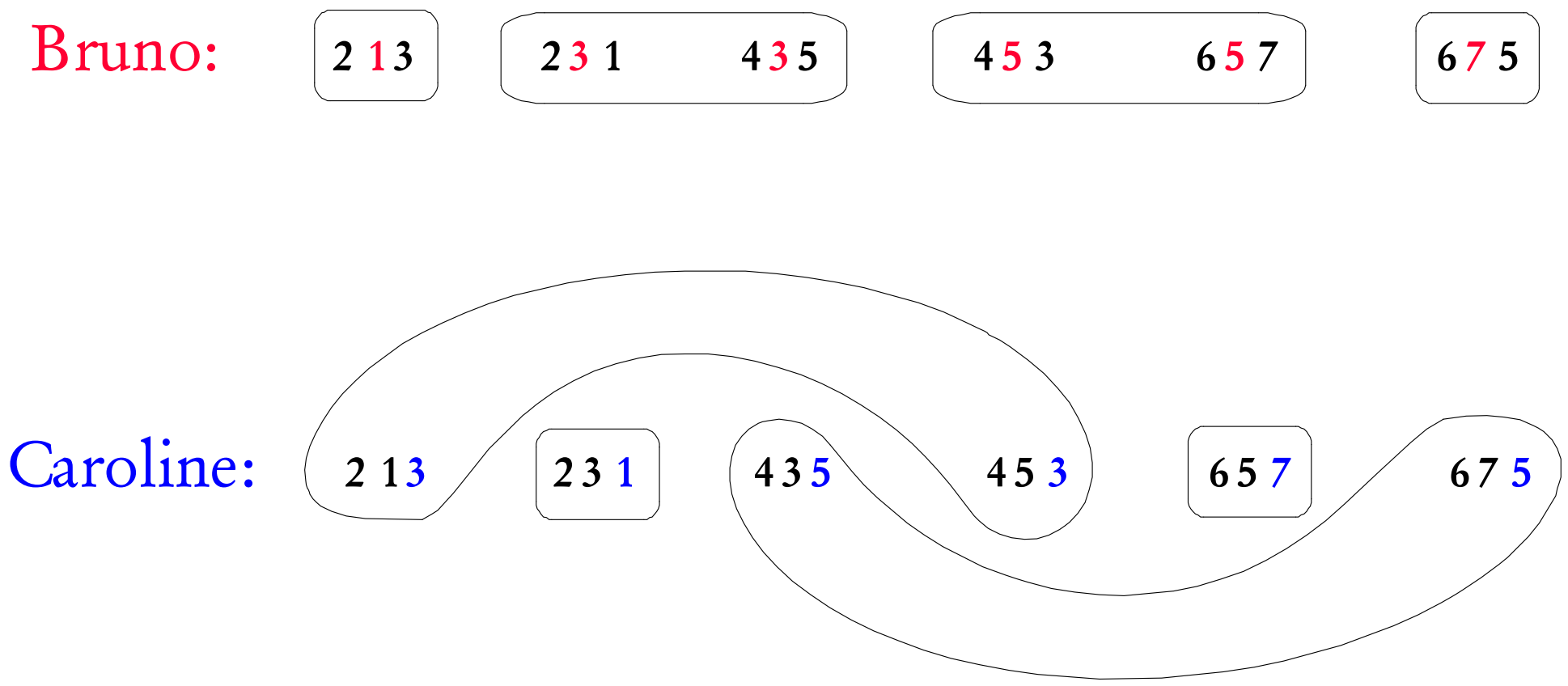

Applying Definition 7.5 we have that $I_{CK}(213) = I_{CK}(453) = I_{CK}(657) = \{213, 453, 657\}$[5] and $I_{CK}(231) = I_{CK}(435) = I_{CK}(675) = \{231, 435, 675\}$. Thus the common knowledge partition is as shown in Figure 7.9

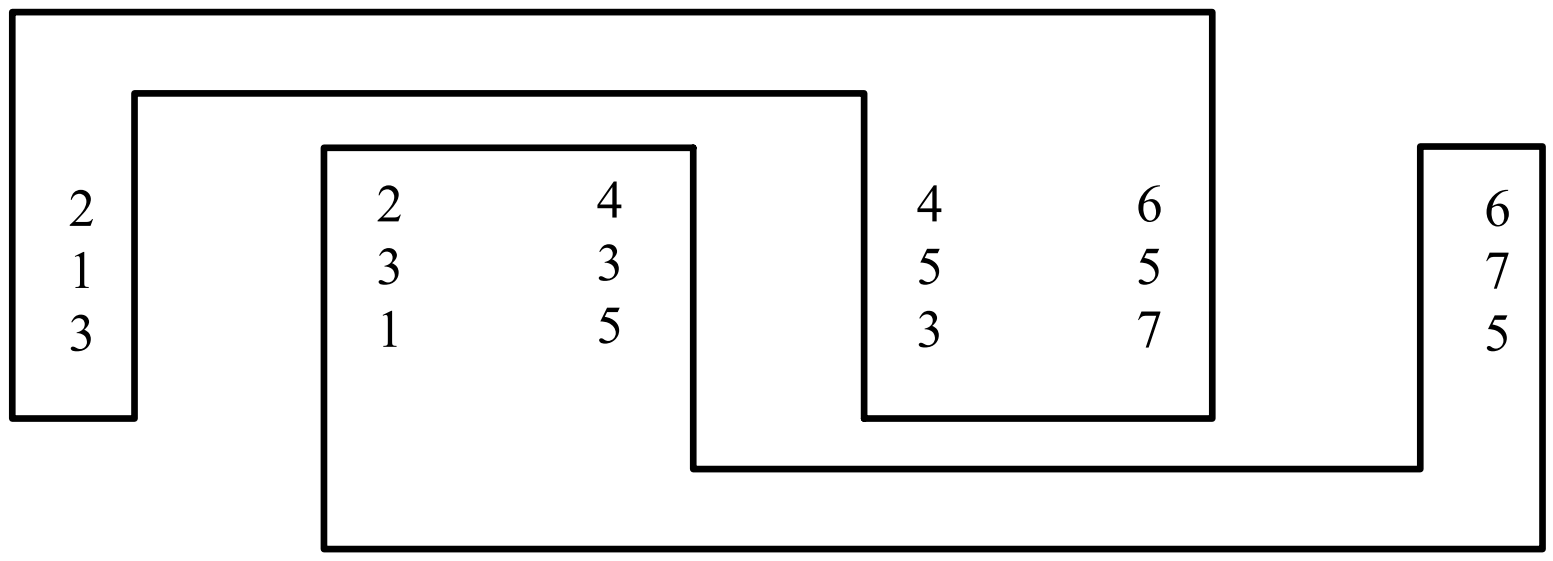

**Figure 7.9**

[5] In fact, 213 is reachable from itself in one step (through either Bruno or Caroline), 453 is reachable in one step from 213 (through Bruno) and 657 is reachable in two steps from 213 (the first step – to 453 – through Caroline and the second step – from 453 – through Bruno).

As a second example, consider the information partitions of Figure 7.2, reproduced below:

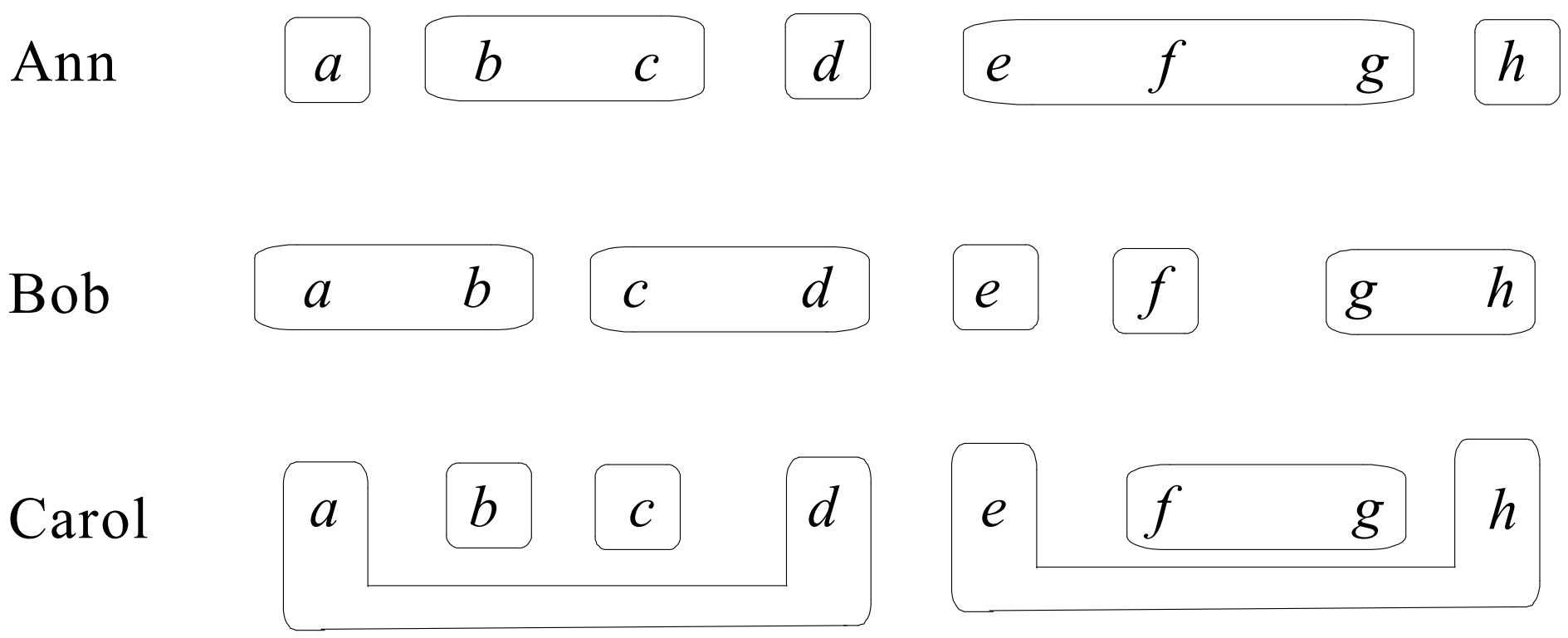

The corresponding common knowledge partition is shown in Figure 7.10.

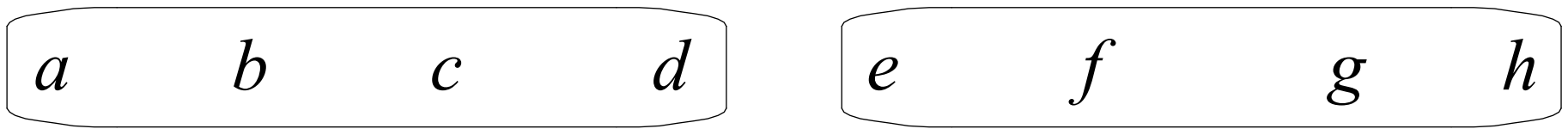

**Figure 7.10**

The following theorem states that, determining whether an event $E$ is common knowledge at a state $w$, is equivalent to determining whether an individual whose information partition coincided with the common knowledge partition would know $E$ at $w$.

**Theorem 7.1.** At state $w \in W$ event $E \subseteq W$ is common knowledge (that is, $w \in CKE$) if and only if $I_{CK}(w) \subseteq E$. Thus the we can define the common knowledge operator $CK : 2^W \to 2^W$ as follows: $CKE = \{w \in W : I_{CK}(w) \subseteq E\}$.

**Example 7.1 continued.** Let us go back to Example 7.1 about Bruno and Caroline. Let $F$ be the event that Bruno gets a 1 or a 5: $F = \{213,\ 453,\ 657\}$. Then, using Theorem 7.1, $CKF = F$ because $I_{CK}(213) = I_{CK}(453) = I_{CK}(657) = \{213, 453, 657\}$; thus – confirming what we found in Example 7.1 – at any state where Bruno gets a 1 or a 5 it is common knowledge between Bruno and Caroline that Bruno got a 1 or a 5.

Now let $H$ be the event "Bruno did *not* get a 5", that is, $H = \{213, 231, 435, 675\}$. Then, using Theorem 7.1 we have that $CKH = \{231, 435, 675\}$. Thus while at state 231 Bruno does not get a 5 and this is common knowledge between Bruno and Caroline, at state 213 Bruno does not get a 5 but this is *not* common knowledge between Bruno and Caroline (in fact $213 \notin K_{Caroline}H = \{231, 435, 675\}$).

As a last example, consider again the information partitions of Figure 7.2, whose common knowledge partition was shown in Figure 7.10 and is reproduced below:

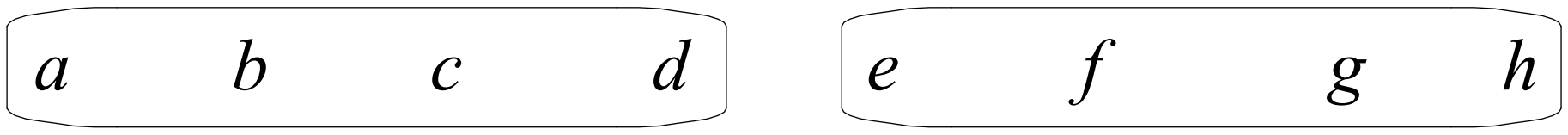

Let $E$ = $\{a,b,c,d,e,f\}$. Then $CKE$ = $\{a,b,c,d\}$. Let $F$ = $\{a,b,f,g,h\}\}$. Then $CKF$ = $\varnothing$. The smallest event that is common knowledge at state $a$ is $I_{CK}(a) = \{a, b, c, d\}$ and the smallest event that is common knowledge at state $g$ is $I_{CK}(g) = \{e, f, g, h\}$.

This is a good time to test your understanding of the concepts introduced in this section, by going through the exercises in Section 7.E.3 of Appendix 7.E at the end of this chapter.

# Appendix 7.E: Exercises

## 7.E.1. Exercises for Section 7.1: Individual knowledge

The answers to the following exercises are in Appendix S at the end of this chapter.

**Exercise 7.1**. You are in a completely dark room, where you cannot see anything. You open a drawer that you know to contain **individual** socks, all of the same size but of different colors. You know that 5 are blue and 7 are white.

**(a)** First use your intuition to answer the following question. What is the **smallest** number of socks you need to take out in order to know (that is, to be absolutely certain) that you have a matching pair (i.e. either a pair of blue socks or a pair of white socks)?

**(b)** Now represent this situation using states and information sets. Do this to (1) represent the situation after you have taken out one sock, (2) represent the situation after you have taken out two socks and (3) represent the situation after you have taken out three socks. Remember that a state should encode a complete description of the relevant aspects of the situation; in particular, the state should tell us how many socks you have taken out and the color of each sock that you have taken out (thus the set of states changes over time as you take out more socks).

**Exercise 7.2**. Consider the situation described in Exercise 7.1: the room is dark and **you have taken out three socks**. Consider the following alternative scenarios.

**(a)** Somebody tells you the color of the third sock (but you still don't know the color of the other two socks). Represent your state of knowledge by means of an information set.

**(b)** Somebody tells you the color of the matching pair (but you don't know what socks you picked). Represent your state of knowledge by means of an information set.

**Exercise 7.3**. Let the set of states be $W = \{a, b, c, d, e, f, g, h, k, m, n\}$ and the information partition of an individual be $\{\{a,b,c\}, \{d\}, \{e,f,g,h\}, \{k,m\}, \{n\}\}$. Consider the following event: $E = \{a,b,d,k,n\}$.

**(a)** Does the individual know $E$ at state $a$?

**(b)** Does the individual know $E$ at state $c$?

**(c)** Does the individual know $E$ at state $d$?

**(d)** Does the individual know $E$ at state $h$?

**(e)** Does the individual know $E$ at state $k$?

**(f)** Let $KE$ denote the event that the individual knows $E$ (that is, the set of states where the individual knows $E$). What is $KE$?

(g) [Recall that, given an event, say $F$, denote by $\neg F$ be the complement of $F$, that is, the set of states that are **not** in $F$.] Once again, let $E = \{a,b,d,k,n\}$. What is the event $\neg KE$, that is the event that the individual does not know $E$? What is the event $K\neg KE$, that is, the event that the individual knows that she does not know $E$?

**Exercise 7.4**. The famous pirate Sandokan has captured you and put you in front of three chests containing coins. Chest number 1 is labeled "gold," chest number 2 is labeled "bronze," and chest number 3 is labeled "gold or silver." One chest contains gold coins only, another contains silver coins only, and the third bronze coins only.

(a) Represent the set of possible states in the case where the labels might or might not be correct (a state must describe the label and content of each box).
(b) Let $E$ be the event "Box number 1 is mislabeled" (that is, what the label says is false). What states are in event $E$?
(c) Let $F$ be the event "Box number 2 is mislabeled". What states are in event $F$?
(d) What is the event "both boxes 1 and 2 are mislabeled"?
(e) Suppose now that Sandokan tells you that *all* the chests are falsely labeled, that is, what the label says is false (for example, if it says "gold" then you can be sure that the chest does *not* contain gold coins). If you correctly announce the content of **all** the chests you will be given a total of $1,000. If you make a mistake (e.g. state that a chest contains, say, gold coins while in fact it contains bronze coins) then you don't get any money at all. You can open any number of chests you like in order to inspect the content. However, the first time you open a chest, you have to pay $500, the second time $300, the third time $100. **(e.1)** What is the set of possible states (assuming that Sandokan told you the truth)? **(e.2)** What is the maximum amount of money you can be absolutely certain to make?

**Exercise 7.5**. Prove each of the following properties of the knowledge operator. The proofs are straightforward applications of Definition 7.2.

**Truth:** $KE \subseteq E$, that is, if at a state $w$ one knows $E$, then $E$ is indeed true at $w$.

**Consistency:** $KE \cap K\neg E = \varnothing$, that is, one never simultaneously knows E and also $\neg E$ ($\neg E$ denotes the complement of $E$).

**Positive introspection:** $KE \subseteq KKE$, that is, if one knows $E$ then one knows that one knows $E$ (one is fully aware of what one knows). [Note that it follows from this and Truth that $KE = KKE$, because from Truth we get that $KKE \subseteq KE$.]

**Negative Introspection:** $\neg KE \subseteq K\neg KE$, that is, if one does not know $E$, then one knows that one does not know $E$ (one is fully aware of what one does not know). [Note that it follows from this and Truth that $\neg KE = K\neg KE$, because from Truth we get that $K\neg KE \subseteq \neg KE$.]

**Monotonicity:** If $E \subseteq F$, then $KE \subseteq KF$, that is, if $E$ implies $F$ then if one knows $E$ then one knows $F$.

**Conjunction:** $KE \cap KF = K(E \cap F)$, that is, if one knows $E$ and one knows $F$, then one knows $E$ and $F$, and *vice versa*.

## 7.E.2. Exercises for Section 7.2: Interactive knowledge

The answers to the following exercises are in Appendix S at the end of this chapter.

**Exercise 7.6**. Let the set of states be $W = \{a, b, c, d, e, f, g, h\}$. There are three individuals with the following information partitions:

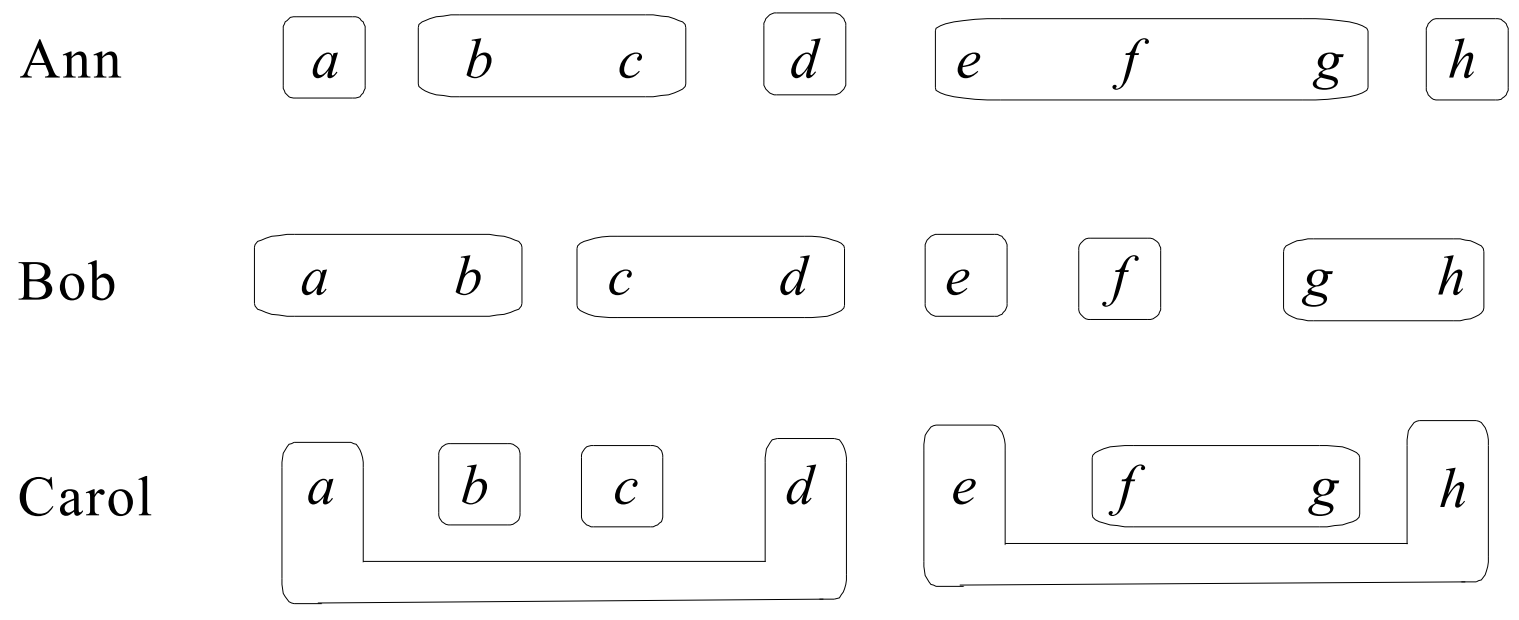

Consider the event $E = \{a,b,c,f,g\}$. Find the following events.

**(a)** $K_{Ann}E$ (the event that Ann knows $E$). **(b)** $K_{Bob}E$, **(c)** $K_{Carol}E$,

**(d)** $K_{Carol}K_{Ann}E$ (the event that Carol knows that Ann knows $E$),

**(e)** $K_{Bob}K_{Carol}K_{Ann}E$, **(f)** $K_{Ann}\neg K_{Bob}K_{Carol}E$ (the event that Ann knows that Bob does not know that Carol knows $E$).

**Exercise 7.7.** Dan is at the airport. He calls his office and leaves a voice message that says: "My flight was cancelled and I am waiting to see if they can re-route me through Los Angeles or San Francisco. I will call one of you at home tonight at 8:00 sharp to let that person know whether I am in San Francisco or Los Angeles." Dan's office staff, consisting of Ann, Barb and Carol, were out for lunch. When they come back they listen to the message together. They leave the office at 5:00 pm and each goes to her home.

**(a)** Using information partitions, represent the possible states of knowledge of Ann, Barb and Carol concerning Dan's whereabouts at 8:15 pm, after Dan's call (there has been no communication among Ann, Barb and Carol after they left the office).

**(b)** Let $E$ be the event that Dan calls either Barb or Carol. What states are in $E$?

**(c)** For the event $E$ of part (b), find $K_A E$ (Ann knows $E$), $K_B K_A E$ (Barb knows that Ann knows $E$) and $\neg K_C E$ (Carol does not know $E$ or it is not the case that Carol knows $E$).

**(d)** For the event $E$ of part (b), find a state $x$ where all of the following are true: (1) at $x$ Ann knows $E$, (2) at $x$ Barb knows that Ann knows $E$, (3) at $x$ it is not the case that Carol knows $E$.

**Exercise 7.8**. A set of lights is controlled by two switches, each of which can be in either the Up position or in the Down position. One switch is in room number 1, where Ann is; the other switch is in room number 2, where Bob is. The lights are in room number 3, where Carla is. There are two lights: one red and one green. The red light is on if the two switches are in different positions (one up and the other down: it doesn't matter which is up and which is down), while the green light is on if the two switches are in the same position (both up or both down). All this is common knowledge among Ann, Bob and Carla.

**(a)** Represent the possible states (you need to specify the position of each switch and which light is on).

**(b)** Represent the possible states of information of Ann, Bob and Carla by means of information partitions.

**(c)** Let $G$ be the event "the green light is on". Find the events $G$, $K_A G$ (Ann knows $G$), $K_B G$ (Bob knows $G$), $K_C G$ (Carla knows $G$).

**(d)** Let $L$ be the event "either the green light is on or the red light is on". Find the events $L$, $K_A L$ (Ann knows $L$), $K_B L$ (Bob knows $L$), $K_C L$ (Carla knows $L$).

## 7.E.3. Exercises for Section 7.3: Common knowledge

The answers to the following exercises are in Appendix S at the end of this chapter.

**Exercise 7.9**. In Exercise 7.6,

**(a)** find the common knowledge partition,

**(b)** find the event *CKE* (that is, the event that $E$ is common knowledge; recall that $E = \{a,b,c,f,g\}$)

**(c)** Letting $F = \{a,b,c,d,e,g\}$, find *CKF*, that is, the event that $F$ is common knowledge.

**Exercise 7.10**. In Exercise 7.7,

**(a)** find the common knowledge partition,

**(b)** find the event *CKE* (where $E$ is the event that Dan calls either Barb or Carol).

**Exercise 7.11**. In Exercise 7.8,

**(a)** find the common knowledge partition,

**(b)** find the event *CKG* (where $G$ is the event "the green light is on"),

**(c)** find the event *CKL* (where $L$ is the event "either the green light is on or the red light is on").

**Exercise 7.12**. The set of states is $W = \{a, b, c, d, e, f, g, h\}$. There are four individuals with the following information partitions:

Individual 1: $\{ \{a, b\}, \{c\}, \{d\}, \{e, f\}, \{g\}, \{h\} \}$
Individual 2: $\{ \{a\}, \{b, c\}, \{d, e\}, \{f\}, \{g\}, \{h\} \}$
Individual 3: $\{ \{a, c\}, \{b\}, \{d\}, \{e\}, \{g\}, \{f, h\} \}$
Individual 4: $\{ \{a\}, \{b, c\}, \{d, e\}, \{f, g\}, \{h\} \}$

**(a)** Let $E = \{a, c, d, e\}$. Find the following events: $K_1E$, $K_2E$, $K_3E$, $K_4E$ and $K_1 K_2 \neg K_3E$ (Recall that $\neg$ denotes the complement of a set, that is, $\neg F$ is the set of all states that are *not* in $F$).

**(b)** Find the common knowledge partition.

**(c)** At what states is event $E = \{a, c, d, e\}$ common knowledge?

**(d)** Let $F = \{a, b, c, d, g, h\}$. Find the event that $F$ is common knowledge, that is, find *CKF*.

**Exercise 7.13**. Amy, Beth and Carla are now in Room 1. They are asked to proceed, one at a time, to Room 3 through Room 2. In Room 2 there are two large containers, one with many red hats and the other with many white hats. They have to choose one hat, put it on their head and then go and sit in the chair in Room 3 that has their name on it.

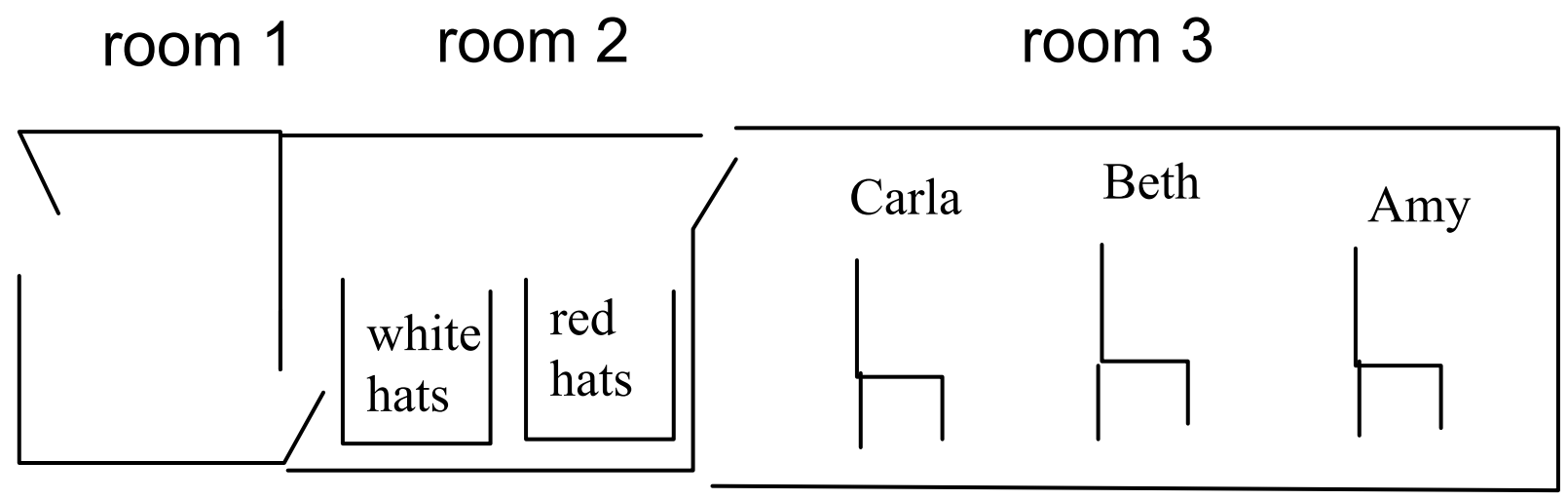

Amy goes first, then Beth then Carla. The chairs are turned with the back to the door. Thus a person entering Room 3 can see whomever is already seated there, but cannot be seen by them. Beth and Carla don't know how many hats there were in each box.

**(a)** Use information partitions to represent the possible states of knowledge of Amy, Beth and Carla after they are seated in Room 3.

**(b)** Suppose that Amy chose a white hat, Beth a red hat and Carla a white hat. Find the smallest event that is common knowledge among them. Give also a verbal description of this event.

**(c)** Repeat Parts (a) and (b) for the modified setting where there is a mirror in Room 3 that allows Amy to see the hat of Carla (but not that of Beth).

◊◊◊◊◊◊◊◊◊◊◊◊

**Exercise 7.14. Challenging Question.**

Pete and his three daughters Elise, Andrea and Theresa are in the same room. Pete gives a sealed envelope to Elise and tells her (in a loud voice, so that everybody can hear) "in this envelope I put a sum of money; I don't remember how much I put in it, but I know for sure that it was either \$4 or \$8 or \$12. Now, my dear Elise, go to your room by yourself and open it. Divide the money into two equal amounts, put the two sums in two different envelopes, seal them and give one to Andrea and one to Theresa". Unbeknownst to her sisters, Elise likes one of them more than the other and decides to disobey her father: after dividing the sum into two equal parts, she takes \$1 from one envelope and puts it in the other envelope. She then gives the envelope with more money to her favorite sister and the envelope with the smaller amount to the other sister. Andrea and Theresa go to their respective rooms and privately open their envelopes, to discover, to their surprise, an odd number of dollars. So they realize that Elise did not follow their father's instructions. Neither Andrea nor Theresa suspect that Elise kept some money for herself; in fact, it is common knowledge between them that Elise simply rearranged the money, without taking any for herself. Of course, *neither Andrea nor Theresa know in principle how much money Elise took from one envelope* (although in some cases they might be able to figure it out). Thus it is *not* common knowledge between Andrea and Theresa that Elise took only \$1 from one of the two envelopes. Your answers should reflect this.

(a) Use states and information partitions to represent the possible states of knowledge of Andrea and Theresa

(b) Let $E$ be the event that Andrea is Elise's favorite sister. Find the events $K_A E$ (the event that Andrea knows it), $K_T E$ (the event that Theresa knows it), $K_A K_T E$ and $K_T K_A E$.

(c) Find the common knowledge partition.

(d) Is there a state at which it is common knowledge between Andrea and Theresa that Elise's favorite sister is Andrea?

(e) The night before, Pete was looking through his digital crystal ball and saw what Elise was planning to do. However the software was not working properly (the screen kept freezing) and he could not tell whether Elise was trying to favor Andrea or Theresa. He knew he couldn't stop Elise and wondered how much money he should put in the envelope to make sure that the mistreated sister (whether it be Andrea or Theresa) would not know that she had been mistreated. How much money should he put in the envelope?

# Appendix 7.S: Solutions to exercises

**Exercise 7.1.** **(a)** The answer is that you need to take out three socks. **(b)** Describe a state by the number of socks you have taken out and the color of each sock that you have taken out.

**After you have taken out only one sock**:

| | |
|---|---|
| 1 | 1 |
| B | W |

(the information set captures the fact that you cannot see the color of the sock because it is dark)

**After you have taken out two socks:**

| | | | |
|---|---|---|---|
| 2 | 2 | 2 | 2 |
| B | B | W | W |
| B | W | B | W |

**After you have taken out three socks:**

| | | | | | | | |
|---|---|---|---|---|---|---|---|
| 3 | 3 | 3 | 3 | 3 | 3 | 3 | 3 |
| B | B | B | B | W | W | W | W |
| B | B | W | W | B | B | W | W |
| B | W | B | W | B | W | B | W |

Now at every state there are (at least) two socks of the same color, thus you know that you have a matching pair, even though you don't know what state you are in and hence you don't know the color of the matching pair.
An alternative (and equivalent) way to proceed would be to describe the state as a quadruple of numbers as follows:

| | |
|---|---|
| # of blue socks in drawer | $x_1$ |
| # of white socks in drawer | $x_2$ |
| # of blue socks in your hand | $x_3$ |
| # of white socks in your hand | $x_4$ |

Clearly, it must always be that $x_1 + x_3 = 5$ and $x_2 + x_4 = 7$. The initial state is $x_1 = 5$, $x_2 = 7$, $x_3 = x_4 = 0$. Taking one sock from the drawer will modify the state as follows:

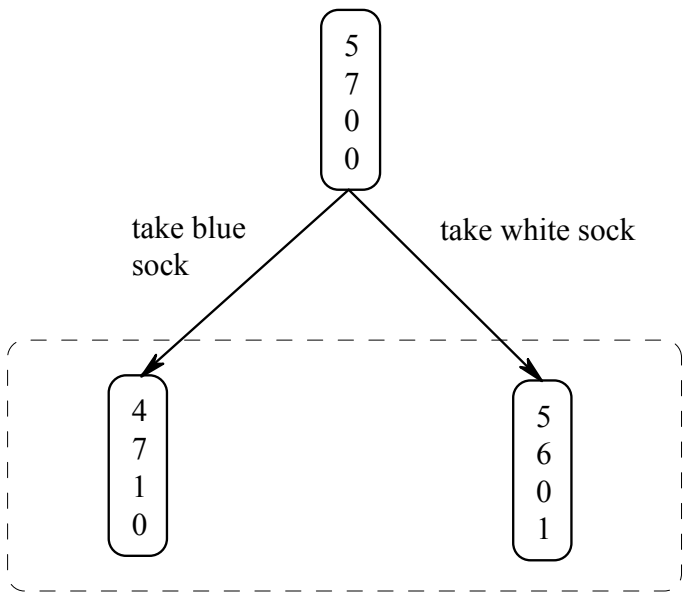

Since it is dark and you don't know what color sock you took, you cannot distinguish between the two states and we have represented this by enclosing the two states in a dashed rectangle, representing an information set.

Now taking a second sock will modify the state as follows.

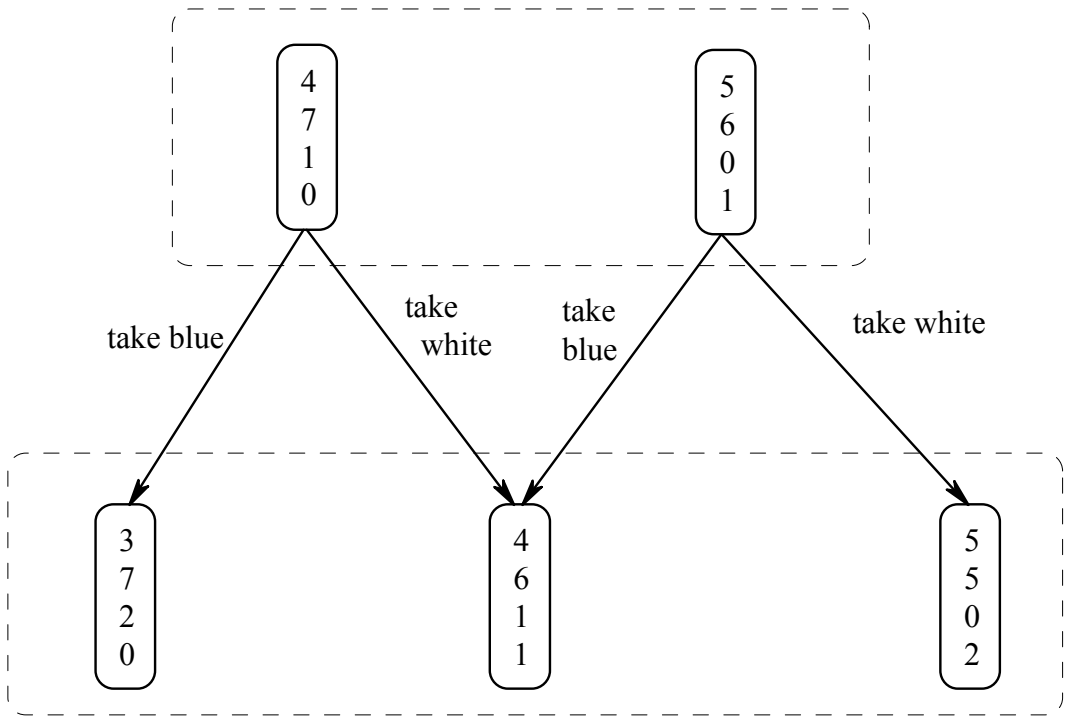

Since you could not distinguish between the two initial states, you cannot now distinguish among the three states that represent all the possibilities. In two of the states you have a matching pair, but you cannot be certain that you have a matching pair because you might be in the middle state where you have one blue sock and one white sock.

Now taking a third sock will lead to:

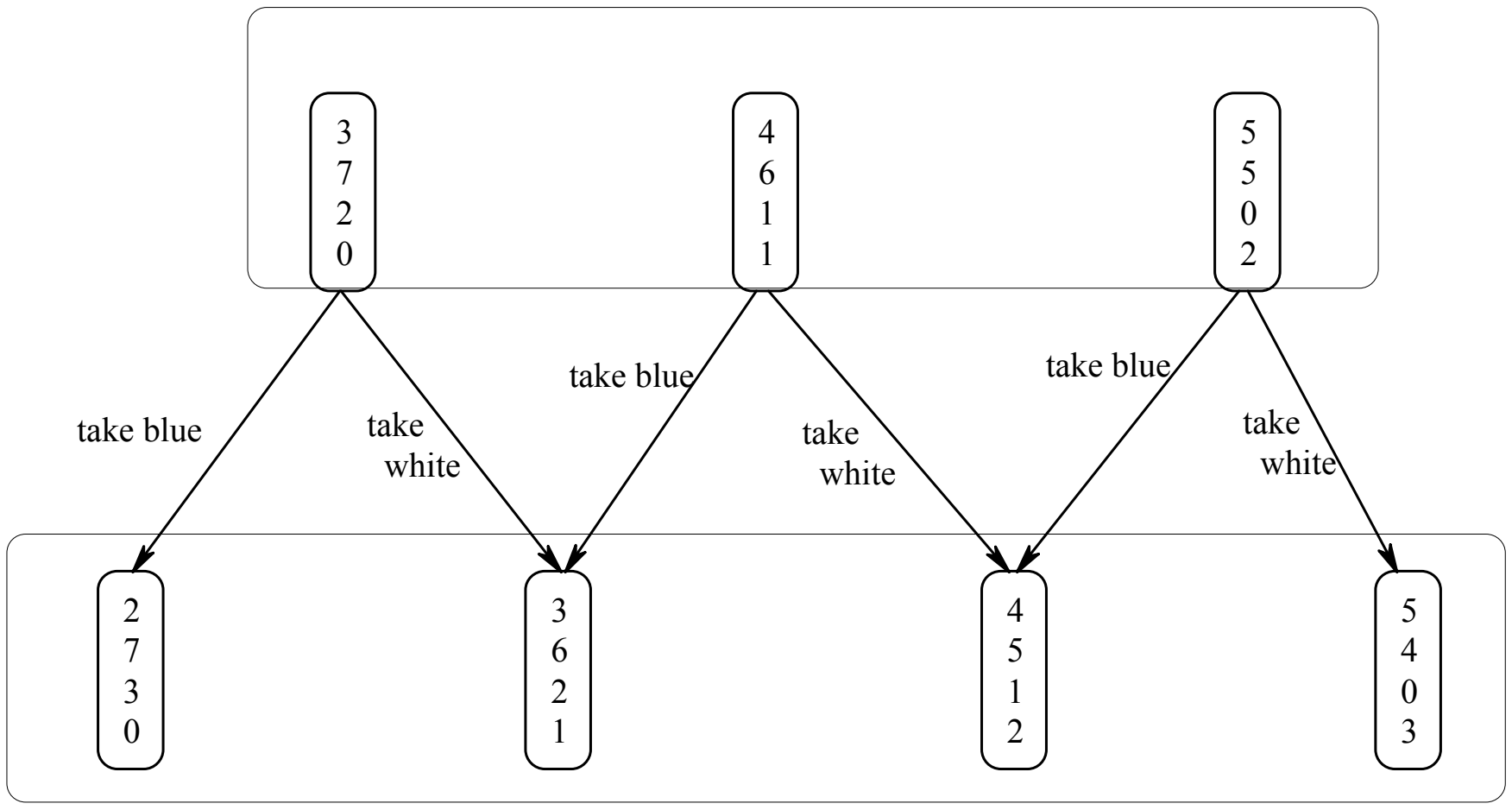

Now you know that you have a matching pair because in every possible state you have at least 2 socks of the same color. Thus the answer is indeed that you need to take out three socks.

**Exercise 7.2.**

**(a)**

| | | | |
|---|---|---|---|
| 3 | 3 | 3 | 3 |
| B | B | W | W |
| B | W | B | W |
| B | B | B | B |
| 3 | 3 | 3 | 3 |
| B | B | W | W |
| B | W | B | W |
| W | W | W | W |

**(b)**

| | | | |
|---|---|---|---|
| 3 | 3 | 3 | 3 |
| B | B | W | B |
| B | W | B | B |
| B | B | B | W |

you have a blue matching pair

| | | | |
|---|---|---|---|
| 3 | 3 | 3 | 3 |
| W | B | W | W |
| W | W | B | W |
| B | W | W | W |

you have a white matching pair

**Exercise 7.3.** **(a)** No, because the information set containing *a* is {*a,b,c*} which is not contained in *E* = {*a,b,d,k,n*}.

**(b)** No, because the information set containing *c* is {*a,b,c*} which is not contained in *E* = {*a,b,d,k,n*}.

**(c)** Yes, because the information set containing *d* is {*d*} which *is* contained in *E* = {*a,b,d,k,n*}.

**(d)** No, because the information set containing *h* is {*e,f,g,h*} which is not contained in *E* = {*a,b,d,k,n*}.

**(e)** No, because the information set containing *k* is {*k,m*} which is not contained in *E* = {*a,b,d,k,n*}.

**(f)** *KE* = {*d,n*}.

**(g)** $\neg KE$ = { *a,b,c,e,f,g,h,k,m* }. $K\neg KE$ = { *a,b,c,e,f,g,h,k,m* } = $\neg KE$.

**Exercise 7.4.**

**(a)** Describe a state by a triple $((x_1,y_1),(x_2,y_2),(x_3,y_3))$ where $x_i$ is the **label** on box number *i* (*i* = 1,2,3) and $y_i$ is the **content** of box number *i*. Then the set of possible states is

$$z_1 = ((G,B),(B,G),(G \text{ or } S,S)),\ z_2 = ((G,B),(B,S),(G \text{ or } S,G)),$$
$$z_3 = ((G,G),(B,B),(G \text{ or } S,S)),\ z_4 = ((G,G),(B,S),(G \text{ or } S,B)),$$
$$z_5 = ((G,S),(B,B),(G \text{ or } S,G)),\ z_6 = ((G,S),(B,G),(G \text{ or } S,B)).$$

Thus, for example, state $z_4$ is one where box number 1 (which is labeled 'gold') in fact contains gold coins, box number 2 (which is labeled 'bronze') as a matter of fact contains silver coins, and box number 3 (which is labeled 'gold or silver') as a matter of fact contains bronze coins. More simply, we could write a state as a triple $(y_1,y_2,y_3)$ denoting the content of each box (since we are told what the labels are). In this simpler notation the states are:

$$z_1 = (B,G,S),\ z_2 = (B,S,G),\ z_3 = (G,B,S),$$
$$z_4 = (G,S,B),\ z_5 = (S,B,G),\ z_6 = (S,G,B).$$

**(b)** $E = \{(B,G,S), (B,S,G), (S,B,G), (S,G,B)\}$ (recall that the label says "gold").

**(c)** $F = \{(B,G,S), (B,S,G), (G,S,B), (S,G,B)\}$ (recall that the label says "bronze").

**(d)** $E \cap F = \{(B,G,S), (B,S,G), (S,G,B)\}$.

**(e)** Of all the states listed above, only state $z_6$ is such that **all** the labels are false. Hence Sandokan's statement reduces the set of states to only one: $z_6 = (S,G,B)$. This is because the label "gold or silver" must be on the chest containing the bronze coins. Hence we are only left with gold and silver. Then silver must be in the chest labeled "gold". Hence gold must be in the chest labeled "bronze". Thus, by looking at the labels you can correctly guess the content without opening any chests: you know that the true state is (S,G,B). So your expected payoff is $1,000 (you are not guessing at random, you are deducing by reasoning).

**Exercise 7.5.**

**Truth:** $KE \subseteq E$. Consider an arbitrary $w \in KE$. We have to show that $w \in E$. Since
$w \in KE$, by Definition 7.2 $I(w) \subseteq E$. Since $w \in I(w)$, it follows that $w \in E$.

**Consistency:** $KE \cap K\neg E = \varnothing$. Suppose that $w \in KE \cap K\neg E$ for some $w$ and some $E$. Then, by Definition 7.2, $I(w) \subseteq E$ (because $w \in KE$) and $I(w) \subseteq \neg E$ (because $w \in \neg E$) and thus $I(w) \subseteq E \cap \neg E$. Since $E \cap \neg E = \varnothing$, this implies that $I(w) = \varnothing$, which is not true because $w \in I(w)$.

**Positive introspection:** $KE \subseteq KKE$. Consider an arbitrary $w \in KE$. We need to show that $w \in KKE$, that is, that $I(w) \subseteq KE$ which means that $w' \in KE$ for every $w' \in I(w)$. Since $w \in KE$, by Definition 7.2, $I(w) \subseteq E$. Consider an arbitrary $w' \in I(w)$. By definition of partition, $I(w') = I(w)$. Thus $I(w') \subseteq E$, that is, by Definition 7.2, $w' \in KE$.

**Negative Introspection:** $\neg KE \subseteq K\neg KE$. Consider an arbitrary $w \in \neg KE$. We need to show that $w \in K\neg KE$, that is, that $I(w) \subseteq \neg KE$. By Definition 7.2, since $I(w) \subseteq \neg KE$, $I(w) \cap \neg E \neq \varnothing$. Consider an arbitrary $w' \in I(w)$; then, since (by definition of partition) $I(w') = I(w)$, $I(w') \cap \neg E \neq \varnothing$ so that $w' \in \neg KE$. Thus we have shown that, for every $w' \in I(w)$, $w' \in \neg KE$, that is, $I(w) \subseteq \neg KE$ which, by Definition 7.2, yields $w \in K\neg KE$.

**Monotonicity:** if $E \subseteq F$, then $KE \subseteq KF$. Consider an arbitrary $w \in KE$. We need to show that $w \in KF$. Since $w \in KE$, by Definition 7.2, $I(w) \subseteq E$. Hence, since – by hypothesis – $E \subseteq F$, $I(w) \subseteq F$, that is, by Definition 7.2, $w \in KF$.

**Conjunction:** $KE \cap KF = K(E \cap F)$. Let $w \in KE \cap KF$. Then $w \in KE$ and $w \in KF$; by Definition 7.2, the former implies that $I(w) \subseteq E$ and the latter that $I(w) \subseteq F$, so that $I(w) \subseteq E \cap F$ and hence, by Definition 7.2, $w \in K(E \cap F)$. Conversely, suppose that $w \in K(E \cap F)$. Then, by Definition 7.2, $I(w) \subseteq E \cap F$ and thus $I(w) \subseteq E$ and $I(w) \subseteq F$ so that, by Definition 7.2, $w \in KE$ and $w \in KF$; hence $w \in KE \cap KF$.

**Exercise 7.6.** **(a)** $K_{Ann}E = \{a,b,c\}$, **(b)** $K_{Bob}E = \{a,b,f\}$, **(c)** $K_{Carol}E = \{b,c,f,g\}$,

**(d)** $K_{Carol}\, K_{Ann}E = K_{Carol}\{a,b,c\} = \{b,c\}$, **(e)** $K_{Bob}K_{Carol}\, K_{Ann}E = K_{Bob}\{b,c\} = \varnothing$

**(f)** $K_{Ann}\, \neg K_{Bob}K_{Carol}E = K_{Ann}\, \neg K_{Bob}\{b,c,f,g\} = K_{Ann}\, \neg\{f\} = K_{Ann}\{a,b,c,d,e,g,h\} = \{a,b,c,d,h\}$.

**Exercise 7.7.** **(a)** We can represent a state as a pair $(x,y)$ where $x$ is the city where Dan is (SF for San Francisco or LA for Los Angeles) and $y$ is the person that Dan called ($A$ for Ann, $B$ for Barb and $C$ for Carol). The information partitions are as follows:

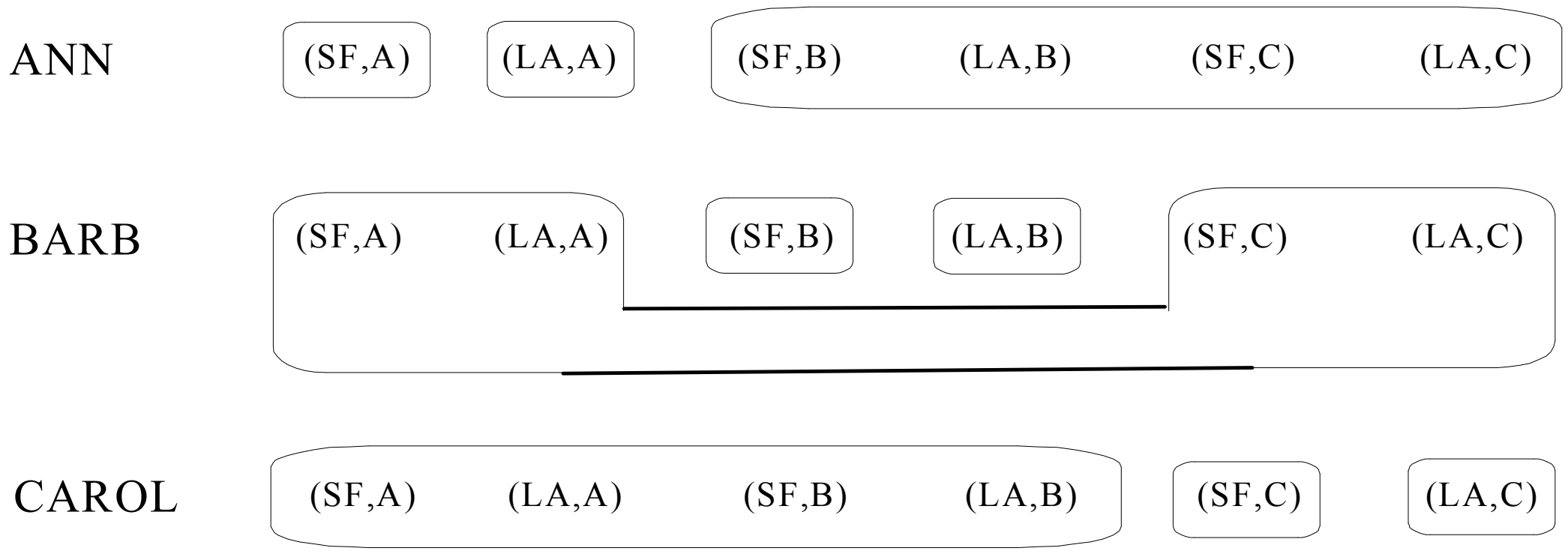

**(b)** $E$ = {(SF,B), (LA,B), (SF,C), (LA,C)}.

**(c)** $K_{Ann}E = E$, $K_{Barb}K_{Ann}E$ = {(SF,B), (LA,B)} and $K_{Carol}E$ = {(SF,C), (LA,C)} so that $\neg\, K_{Carol}E$ = {(SF,A), (LA,A), (SF,B), (LA,B)}.

**(d)** We want a state $x$ such that $x \in K_{Ann}E$, $x \in K_{Barb}K_{Ann}E$ and $x \in \neg K_{Carol}E$ (that is, $x \notin K_{Carol}E$). There are only two such states: (SF,B) and (LA,B). Thus either $x$ = (SF,B) or $x$ = (LA,B).

**Exercise 7.8.** **(a)** We can represent a state as a triple $\begin{matrix} x \\ y \\ z \end{matrix}$ where $x$ is the position of the switch in room 1 (Up or Down), $y$ is the position of the switch in room 2 and $z$ is the light which is on (Green or Red):

$$\begin{matrix} U & U & D & D \\ U & D & U & D \\ G & R & R & G \end{matrix}$$

**(b)**

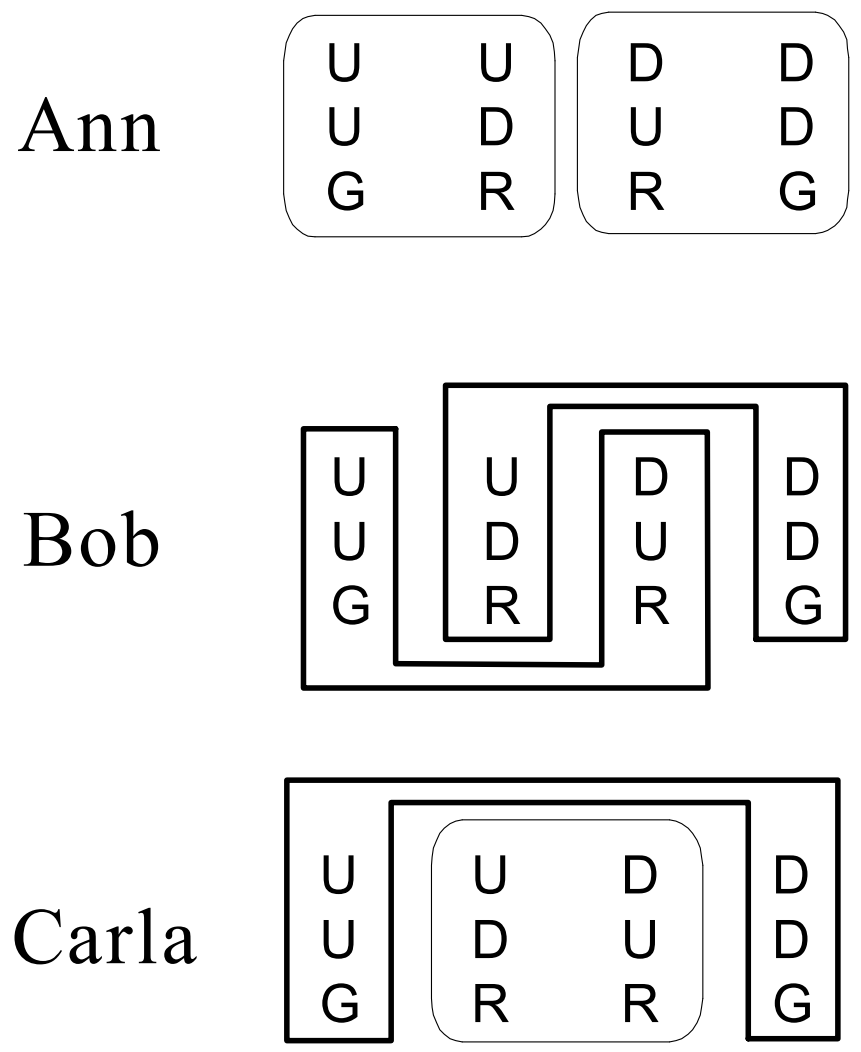

**(c)** $G = \left\{ \begin{matrix} U & D \\ U & D \\ G & G \end{matrix} \right\}$. $K_A G = \varnothing$, $K_B G = \varnothing$, $K_C G = \left\{ \begin{matrix} U & D \\ U & D \\ G & G \end{matrix} \right\}$.

**(d)** $L$ is the set of all states. Hence $K_A L = K_B L = K_C L = L$.

**Exercise 7.9.** **(a)** The common knowledge partition is:

$\{a \quad b \quad c \quad d\}$ $\{e \quad f \quad g \quad h\}$

**(b)** $CKE = \varnothing$ (where $E = \{a,b,c,f,g\}$).

**(c)** $CKF = \{a,b,c,d\}$ (where $F = \{a,b,c,d,e,g\}$).

**Exercise 7.10.** **(a)** The common knowledge partition is:

{(SF,A) (LA,A) (SF,B) (LA,B) (SF,C) (LA,C)}

**(b)** $CKE = \varnothing$.

**Exercise 7.11.** **(a)** The common knowledge partition is:

| | | | |
|---|---|---|---|
| U | U | D | D |
| U | D | U | D |
| G | R | R | G |

**(b)** $CKG = \varnothing$.
**(c)** $CKL = L$.

**Exercise 7.12.** **(a)** Let $E = \{a, c, d, e\}$. Then: $K_1E = \{c, d\}$, $K_2E = \{a, d, e\}$, $K_3E = \{a, c, d, e\}$, $K_4E = \{a, d, e\}$ and $K_1K_2\neg K_3E = \{g, h\}$ [in fact, $\neg K_3E = \{b, f, g, h\}$, so that $K_2\neg K_3E = \{f, g, h\}$ and $K_1K_2\neg K_3E = \{g, h\}$].

**(b)** The common knowledge partition is $\{ \{a, b, c\}, \{d, e, f, g, h\} \}$

**(c)** At no state is event $E = \{a, c, d, e\}$ common knowledge: $CKE = \varnothing$.

**(d)** $CKF = \{a,b,c\}$ (where $F = \{a, b, c, d, g, h\}$).

**Exercise 7.13.** **(a)** Represent the state as a triple of letters, where the top letter denotes the color of Amy's hat, the second letter the color of Beth's hat and the bottom letter the color of Carla's hat. Each of them knows which hat she chose; furthermore, Beth can see Amy's hat and Carla can see everything. Thus the partitions are as follows.

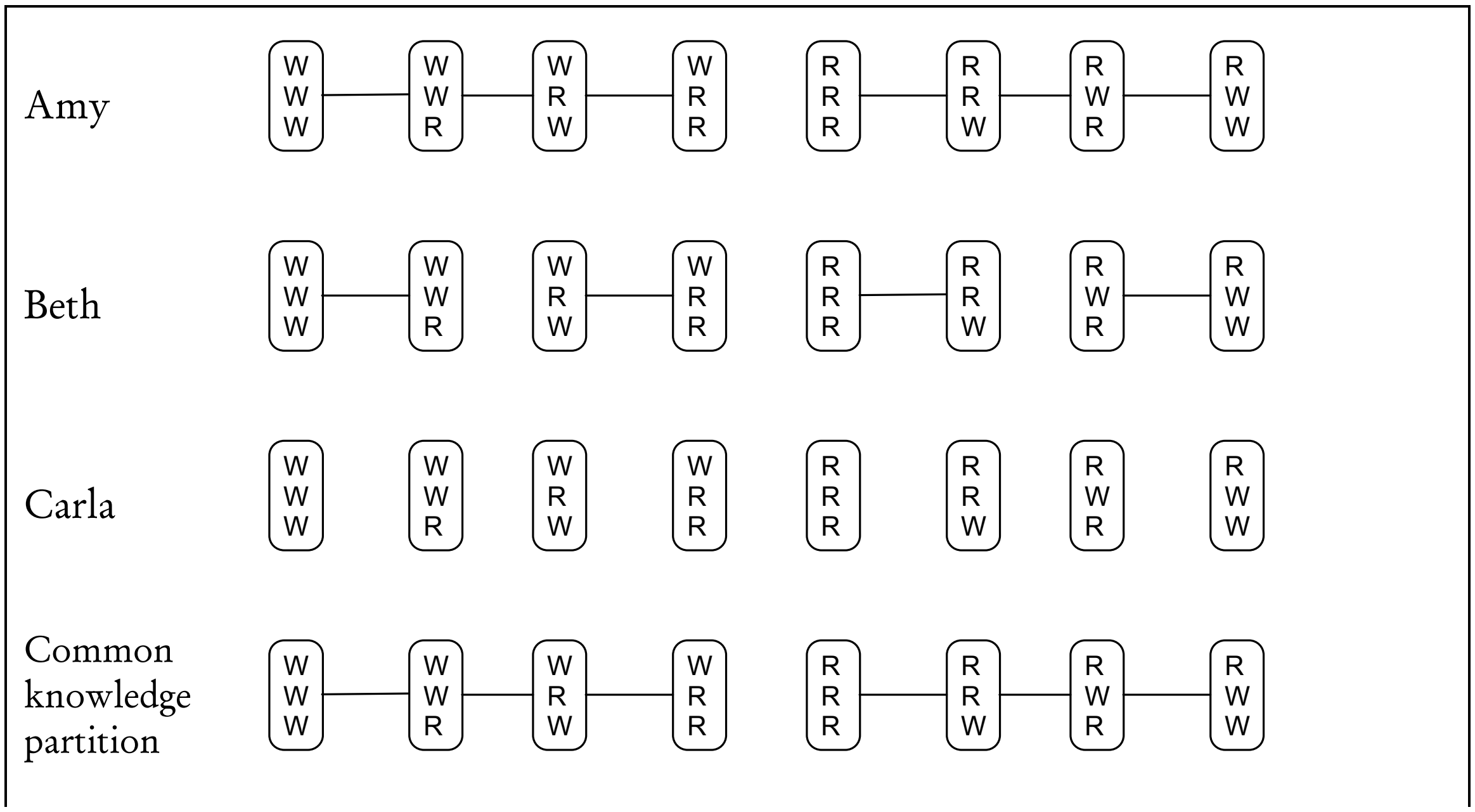

**(b)** When the true state is (*W, R, W*), the smallest event that is common knowledge among them is the first information set of the common knowledge partition, that is, the fact that Amy has a white hat.

**(c)** In the modified setting the information partitions are as follows:

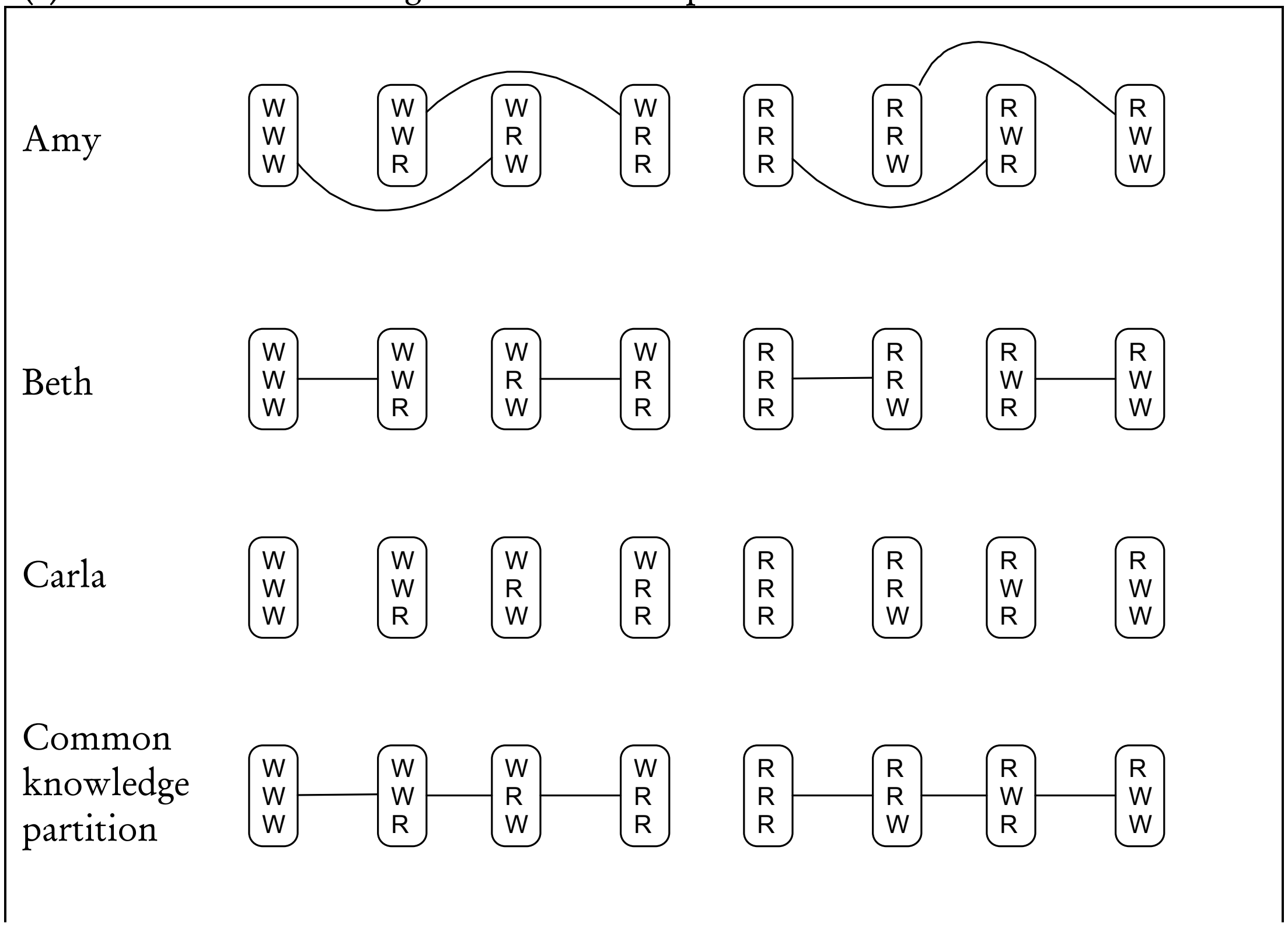

The common knowledge partition is the same as before, hence the smallest event that is common knowledge among all three of them is that Amy has a white hat.

**Exercise 7.14 (Challenging Question).** **(a)** Describe the state by three numbers where the top number is the amount of money that Dan put in the envelope, the middle umber is the amount given to Andrea and the bottom number is the amount given to Theresa. As a matter of fact, each sister can only find either \$1 or \$3 or \$5 or \$7 in the envelope. Thus the **objectively** possible states are: $\begin{pmatrix}4\\1\\3\end{pmatrix},\begin{pmatrix}4\\3\\1\end{pmatrix},\begin{pmatrix}8\\3\\5\end{pmatrix},\begin{pmatrix}8\\5\\3\end{pmatrix},\begin{pmatrix}12\\5\\7\end{pmatrix}$ and $\begin{pmatrix}12\\7\\5\end{pmatrix}$. However, besides these there are also **subjectively** possible states, namely $\begin{pmatrix}8\\1\\7\end{pmatrix},\begin{pmatrix}8\\7\\1\end{pmatrix},\begin{pmatrix}12\\1\\11\end{pmatrix},\begin{pmatrix}12\\11\\1\end{pmatrix},\begin{pmatrix}12\\3\\9\end{pmatrix}$ and $\begin{pmatrix}12\\9\\3\end{pmatrix}$ (because Andrea and Theresa don't know that Elise only transferred \$1 from one envelope to the other). The information partitions are as follows:

$$\text{ANDREA: } \left\{\left\{\begin{pmatrix}4\\1\\3\end{pmatrix},\begin{pmatrix}8\\1\\7\end{pmatrix},\begin{pmatrix}12\\1\\11\end{pmatrix}\right\},\left\{\begin{pmatrix}4\\3\\1\end{pmatrix},\begin{pmatrix}8\\3\\5\end{pmatrix},\begin{pmatrix}12\\3\\9\end{pmatrix}\right\},\left\{\begin{pmatrix}8\\5\\3\end{pmatrix},\begin{pmatrix}12\\5\\7\end{pmatrix}\right\},\left\{\begin{pmatrix}8\\7\\1\end{pmatrix},\begin{pmatrix}12\\7\\5\end{pmatrix}\right\},\left\{\begin{pmatrix}12\\9\\3\end{pmatrix}\right\},\left\{\begin{pmatrix}12\\11\\1\end{pmatrix}\right\}\right\}$$

$$\text{THERESA:} \left\{\left\{\begin{pmatrix}4\\1\\3\end{pmatrix},\begin{pmatrix}8\\5\\3\end{pmatrix},\begin{pmatrix}12\\9\\3\end{pmatrix}\right\},\left\{\begin{pmatrix}4\\3\\1\end{pmatrix},\begin{pmatrix}8\\7\\1\end{pmatrix},\begin{pmatrix}12\\11\\1\end{pmatrix}\right\},\left\{\begin{pmatrix}8\\3\\5\end{pmatrix},\begin{pmatrix}12\\7\\5\end{pmatrix}\right\},\left\{\begin{pmatrix}8\\1\\7\end{pmatrix},\begin{pmatrix}12\\5\\7\end{pmatrix}\right\},\left\{\begin{pmatrix}12\\3\\9\end{pmatrix}\right\},\left\{\begin{pmatrix}12\\1\\11\end{pmatrix}\right\}\right\}$$

**(b)** The event that Andrea is Elise's favorite sister is

$$E = \left\{\begin{pmatrix}4\\3\\1\end{pmatrix},\begin{pmatrix}8\\5\\3\end{pmatrix},\begin{pmatrix}8\\7\\1\end{pmatrix},\begin{pmatrix}12\\7\\5\end{pmatrix},\begin{pmatrix}12\\9\\3\end{pmatrix},\begin{pmatrix}12\\11\\1\end{pmatrix}\right\}.$$

$K_AE = \left\{\begin{pmatrix}8\\7\\1\end{pmatrix},\begin{pmatrix}12\\7\\5\end{pmatrix},\begin{pmatrix}12\\9\\3\end{pmatrix},\begin{pmatrix}12\\11\\1\end{pmatrix}\right\}$, $K_TE = \left\{\begin{pmatrix}4\\3\\1\end{pmatrix},\begin{pmatrix}8\\7\\1\end{pmatrix},\begin{pmatrix}12\\11\\1\end{pmatrix}\right\}$, $K_AK_TE = \left\{\begin{pmatrix}12\\11\\1\end{pmatrix}\right\}$, $K_TK_AE = \varnothing$.

**(c)** The common knowledge partition is:

$$\left\{\left\{\begin{pmatrix}4\\1\\3\end{pmatrix},\begin{pmatrix}8\\1\\7\end{pmatrix},\begin{pmatrix}12\\1\\11\end{pmatrix},\begin{pmatrix}8\\5\\3\end{pmatrix},\begin{pmatrix}12\\9\\3\end{pmatrix},\begin{pmatrix}12\\5\\7\end{pmatrix}\right\},\left\{\begin{pmatrix}4\\3\\1\end{pmatrix},\begin{pmatrix}8\\3\\5\end{pmatrix},\begin{pmatrix}12\\3\\9\end{pmatrix},\begin{pmatrix}8\\7\\1\end{pmatrix},\begin{pmatrix}12\\11\\1\end{pmatrix},\begin{pmatrix}12\\7\\5\end{pmatrix}\right\}\right\}$$

**(d)** At no state. In fact, $CKE = \varnothing$.

**(e)** Dan should put either \$8 or \$12 in the envelope. If he were to put \$4, then one sister would end up with \$1 and know that she was mistreated. In no other case does a mistreated sister know that she got less money than the other.

Chapter 8

# Adding Beliefs to Knowledge

## 8.1 Probabilistic beliefs

An information set contains all the states that an individual considers possible, that is, the states that the individual cannot rule out, given her information. However, of all the states that are possible, the individual might consider some to be more likely than others and might even dismiss some states as "extremely unlikely" or "implausible". For example, suppose that there are only three students in a class: Ann, Bob and Carla. The professor tells them that in the last exam one of them got 95 points (out of 100), another 78 and the third 54. We can think of a state as a triple (*a,b,c*), where *a* is Ann's score, *b* Bob's score and *c* Carla's score. Then, based on the information given by the professor, Ann must consider all of the following states as possible: (95,78,54), (95,54,78), (78,95,54), (78,54,95), (54,95,78) and (54,78,95). Suppose, however, that in almost all the past exams Ann and Bob always obtained a higher score than Carla and often Ann outperformed Bob. Then she might consider states (95,78,54) and (78,95,54) much more likely than (78,54,95) and (54,78,95). To represent such judgments of relative likelihood we add to an information set a probability distribution over the states in the information set.[6] The probability distribution expresses the individual's *beliefs*, while the information set represents what the individual *knows*. In this example, Ann's beliefs could be as shown in Table 8.1.

[6] Recall that a probability distribution over a finite set $\{x_1, x_2, ..., x_m\}$ is an assignment of a number $p_i$ to every element $x_i$ of the set satisfying two properties: (1) for every $i \in \{1,2,...,m\}$, $0 \le p_i \le 1$ and (2) $p_1 + p_2 + ... + p_m = 1$.

| *state* | (95,78,54) | (95,54,78) | (78,95,54) | (54,95,78) | (78,54,95) | (54,78,95) |
|---|---|---|---|---|---|---|
| *probability* | $\frac{9}{16}$ | $\frac{4}{16}$ | $\frac{2}{16}$ | $\frac{1}{16}$ | 0 | 0 |

**Table 8.1**

According to the beliefs expressed by the probability distribution shown in Table 8.1, Ann considers it very likely that she got the highest score, is willing to dismiss the possibility that Carla received the highest score (she assigns probability 0 to the two states where Carla's score is 95) and considers it much more likely that she, rather than Bob, received the highest score.

Let $W$ be a set of states (which we can think of as an information set of an individual). A *probability distribution on W* is a function $P: W \to [0,1]$ such that $\sum_{w \in W} P(w) = 1$. As usual, we call the subsets of *W* events. The *probability of event* $E \subseteq W$, denoted by $P(E)$, is defined as the sum of the probabilities of the elements of *E*: $P(E) = \sum_{w \in E} P(w)$. For instance, in the example of Table 8.1, the proposition "Ann received the highest score" corresponds to event $E$ = {(95,78,54),(95,54,78)} and – according to Ann's beliefs – the probability of that event is $P(E) = P\big((95,78,54)\big) + P\big((95,54,78)\big) = \frac{9}{16} + \frac{4}{16} = \frac{13}{16}$; on the other hand, event $F$ = {(95,78,54),(78,95,54),(54,95,78)} corresponds to the proposition "Bob's score was higher than Carla's score" and – according to Ann's beliefs – the probability of that event is $P(F) = P\big((95,78,54)\big) + P\big((78,95,54)\big) + P\big((54,95,78)\big) = \frac{9}{16} + \frac{2}{16} + \frac{1}{16} = \frac{12}{16} = \frac{3}{4} =$ 75%.

**Definition 8.1.** We say that an individual is *certain* of event *E* if she attaches probability 1 to *E*, that is, if $P(E) = 1$.

In the above example, Ann is certain of event $G$ = {(95,78,54), (95,54,78), (78,95,54), (54,95,78)}, corresponding to the proposition "Carla did not get the highest score".

**Remark 8.1.** Note the important difference between knowledge and certainty: if at a state $w$ the individual knows an event $E$ ($w \in KE$) then at that state $E$ is indeed true ($w \in E$),[7] that is, it is never the case that an individual knows something which is false; on the other hand, an individual can be certain of something which is false, that is, it is possible that the individual assigns probability 1 to an event $E$ even though the true state is not in $E$.[8]

Continuing the example of the exam scores, suppose that Ann has the beliefs shown in Table 8.1 and now the professor makes a new announcement: "I was surprised to see that, unlike past exams, Ann did not get the highest score." This new announcement by the professor informs the students that the true state is neither (95,78,54) nor (95,54,78). Thus we can view the effect of the new piece of information as shrinking Ann's information set from {(95,78,54), (95,54,78), (78,95,54), (54,95,78), (78,54,95), (54,78,95)} to {(78,95,54), (54,95,78), (78,54,95), (54,78,95)}. How should Ann revise her beliefs in response to the new information? The answer cannot be that we simply drop states (95,78,54) and (95,54,78) from Table 8.1, because the result would be

| *state* | (78, 95, 54) | (54, 95, 78) | (78, 54, 95) | (54, 78, 95) |
|---|---|---|---|---|
| *probability* | $\frac{2}{16}$ | $\frac{1}{16}$ | 0 | 0 |

which is not a probability distribution, since the probabilities do not add up to 1 (they add up to $\frac{3}{16}$). The topic of belief revision is addressed in next section.

This is a good time to test your understanding of the concepts introduced in this section, by going through the exercises in Section 8.E.1 of Appendix 8.E at the end of this chapter.

[7] Indeed it was proved in Exercise 7.5 (Chapter 7) that, for every event $E$, $KE \subseteq E$.

[8] In the above example, if the true state is (78,54,95) then it does not belong to event $G$, representing the proposition "Carla did not get the highest score" and yet Ann assigns probability 1 to $G$, that is, she is certain of $G$.

# 8.2 Conditional probability, belief updating, Bayes' rule

We start by reviewing the notion of conditional probability.

**Definition 8.2.** Let $W$ be a finite set of states and $P: W \to [0,1]$ a probability distribution on $W$. Let $w \in W$ be a state and $E \subseteq W$ an event such that $P(E) > 0$. Then the *probability of w conditional on E* (or *given E*), denoted by $P(w \mid E)$, is defined as

$$P(w \mid E) = \begin{cases} 0 & \text{if } w \notin E \\ \dfrac{P(w)}{P(E)} & \text{if } w \in E \end{cases}.$$

For example, let $W = \{a,b,c,d,e,f,g\}$ and let $P$ be as shown in Table 8.2.

| *state* | $a$ | $b$ | $c$ | $d$ | $e$ | $f$ | $g$ |
|---|---|---|---|---|---|---|---|
| *probability* | $\frac{3}{20}$ | 0 | $\frac{7}{20}$ | $\frac{1}{20}$ | 0 | $\frac{4}{20}$ | $\frac{5}{20}$ |

**Table 8.2**

Let $E = \{a,d,e,g\}$; thus $P(E) = P(a) + P(d) + P(e) + P(g) = \frac{3}{20} + \frac{1}{20} + 0 + \frac{5}{20} = \frac{9}{20}$. Then $P(b \mid E) = P(c \mid E) = P(f \mid E) = 0$ (since each of these states does not belong to $E$), $P(a \mid E) = \frac{\frac{3}{20}}{\frac{9}{20}} = \frac{3}{9}$, $P(d \mid E) = \frac{\frac{1}{20}}{\frac{9}{20}} = \frac{1}{9}$, $P(e \mid E) = \frac{0}{\frac{9}{20}} = 0$ and $P(g \mid E) = \frac{\frac{5}{20}}{\frac{9}{20}} = \frac{5}{9}$. Thus the operation of conditioning on event $E$ yields the *updated (on E)* probability distribution shown in Table 8.3.

| *state* | $a$ | $b$ | $c$ | $d$ | $e$ | $f$ | $g$ |
|---|---|---|---|---|---|---|---|
| *probability* | $\frac{3}{9}$ | 0 | 0 | $\frac{1}{9}$ | 0 | 0 | $\frac{5}{9}$ |

**Table 8.3**

Philosophers and logicians have argued that, when forming new beliefs in response to new information, the rational way to proceed is to modify the initial probabilities by using the conditional probability rule. Of course, this can only be done if the new information is *not surprising*, in the sense that it is represented by an event $E$ that, given the initial beliefs, has positive probability ($P(E) > 0$).[9]

**Definition 8.3.** We use the expression *'belief updating'* or *'Bayesian updating'* to refer to the modification of initial beliefs (expressed by an initial probability distribution $P$) obtained by applying the conditional probability rule; this assumes that the belief change is prompted by the arrival of new information, represented by an event $E$ such that $P(E) > 0$.

**Definition 8.4.** The conditional probability formula can be extended from states to events as follows. Let $A$ and $B$ be two events, with $P(B) > 0$. Then the *probability of A conditional on B* (or given $B$), denoted by $P(A \mid B)$ is defined by

$$P(A \mid B) = \frac{P(A \cap B)}{P(B)}.$$

It is straightforward to check (the reader should convince himself/herself of this) that an equivalent way of defining $P(A \mid B)$ is as follows: $P(A \mid B) = \sum_{w \in A \cap B} P(w \mid B)$.

For example, consider again the initial beliefs shown in Table 8.2, which is reproduced below.

| *state* | $a$ | $b$ | $c$ | $d$ | $e$ | $f$ | $g$ |
|---|---|---|---|---|---|---|---|
| *probability* | $\frac{3}{20}$ | 0 | $\frac{7}{20}$ | $\frac{1}{20}$ | 0 | $\frac{4}{20}$ | $\frac{5}{20}$ |

As before, let $E = \{a, d, e, g\}$ and consider the event $D = \{a, b, c, f, g\}$. What is the probability of $D$ given $E$? First compute the intersection of $D$ and $E$: $D \cap E = \{a, g\}$; thus $P(D \cap E) = P(a) + P(g) = \frac{3}{20} + \frac{5}{20} = \frac{8}{20}$. $P(E)$ was computed before: $P(E) = P(a) + P(d) + P(e) + P(g) = \frac{3}{20} + \frac{1}{20} + 0 + \frac{5}{20} = \frac{9}{20}$. Thus, using the formula given in Definition 8.4, we obtain

[9] When $P(E) = 0$, $P(w|E)$ is not well defined because division by zero has no meaning in arithmetics.

$P(D \mid E) = \frac{P(D \cap E)}{P(E)} = \frac{\frac{8}{20}}{\frac{9}{20}} = \frac{8}{9} \simeq 89\%$. Thus, before receiving information *E*, the individual was assigning probability $P(D) = P(a) + P(b) + P(c) + P(f) + P(g) = \frac{3}{20} + 0 + \frac{7}{20} + \frac{4}{20} + \frac{5}{20} = \frac{19}{20} = 95\%$ to event *D* and, after receiving information *E*, she updates that probability to 89%.

According to Definition 8.4, the probability of *A* conditional on *B* is given by

$$P(A \mid B) = \frac{P(A \cap B)}{P(B)}. \qquad (\clubsuit)$$

Similarly, the probability of *B* conditional on *A* is given by $P(B \mid A) = \frac{P(B \cap A)}{P(A)}$, which [since $A \cap B = B \cap A$ and thus $P(A \cap B) = P(B \cap A)$] can be written as

$$P(A \cap B) = P(B \mid A)\, P(A). \qquad (\blacklozenge)$$

From (♣) and (♦) we get that $P(A \mid B) = \frac{P(B \mid A)\, P(A)}{P(B)}$ which is known as *Bayes' rule*.

**Definition 8.5.** Let *E* and *F* be two events, with $P(F) > 0$. Then *Bayes' rule* says that the probability of *E* conditional on *F* is equal to the probability of *F* conditional on *E* times the probability of *E*, divided by the probability of *F*:

$$P(E \mid F) = \frac{P(F \mid E)\, P(E)}{P(F)}.$$

As an example of the application of Bayes' rule, suppose that you are a doctor examining a middle-aged man who complains of lower-back pain. You know that 25% of men in your patient's age group suffer from lower-back pain. There are various causes of lower back pain; one of them is chronic inflammation of the kidneys. This is not a very common disease: it affects only 4% of men in the age group you are considering. Among those who suffer from chronic inflammation of the kidneys, 85% complain of lower-back pain. What is the probability that your patient has chronic inflammation of the kidneys?

Let $I$ denote inflammation of the kidneys and $L$ denote lower-back pain. The information you have is that $P(I)=\frac{4}{100}=4\%$, $P(L)=\frac{25}{100}=25\%$ and $P(L\mid I)=\frac{85}{100}=85\%$. Thus, using Bayes' rule we get that $P(I\mid L)=\frac{P(L\mid I)\,P(I)}{P(L)}=\frac{\frac{85}{100}\frac{4}{100}}{\frac{25}{100}}=0.136=13.6\%$. This can also be seen as follows. Suppose there are 2,500 men in your patient's age group; 4% (= 100) have inflammation of the kidneys and the remaining 96% (= 2,400) don't:

| | Lower-back pain | No lower-back pain | Total | |
|---|---|---|---|---|
| Inflammation | | | 100 | |
| No inflammation | | | 2,400 | |
| | | | 2,500 | Total |

We also know that 25% (= 625) of the population has lower-back pain and the remaining 75% (=1,875) don't:

| | Lower-back pain | No lower-back pain | Total | |
|---|---|---|---|---|
| Inflammation | | | 100 | |
| No inflammation | | | 2,400 | |
| | 625 | 1,875 | **2,500** | Total |

Of the 100 people with inflammation of the kidneys, 85% (= 85) suffer from lower-back pain and the remaining 15 don't:

| | Lower-back pain | No lower-back pain | Total | |
|---|---|---|---|---|
| Inflammation | 85 | 15 | 100 | |
| No inflammation | | | 2,400 | |
| | 625 | 1,875 | **2,500** | Total |

Thus we can fill in the highlighted cells by taking the differences: 625 −85 = 540 and 1,875 − 15 = 1,860:

| | Lower-back pain | No lower-back pain | Total | |
|---|---|---|---|---|
| Inflammation | 85 | 15 | 100 | |
| No inflammation | 540 | 1,860 | 2,400 | |
| | 625 | 1,875 | **2,500** | Total |

Thus the fraction of people, *among those who suffer from lower-back pain,* who have inflammation of the kidneys is $\frac{85}{625}=0.136=13.6\%$.

So far we have focused on belief *updating*. The next section deals with belief *revision*.

This is a good time to test your understanding of the concepts introduced in this section, by going through the exercises in Section 8.E.2 of Appendix 8.E at the end of this chapter.

## 8.3 Belief revision

[Note: the material of this section will not be needed until Chapter 12. It is presented here for completeness on the topic of beliefs.]

How should a "rational" individual revise her beliefs when receiving information that is surprising, that is, when informed of an event $E$ to which her initial beliefs assigned zero probability ($P(E)=0$)? The best known theory of "rational" belief revision is the so-called *AGM theory*, which takes its name from its originators: Carlos Alchourrón (a legal scholar), Peter Gärdenfors (a philosopher) and David Makinson (a computer scientist); their pioneering contribution was published in 1985. Just like the theory of expected utility (Chapter 4), the AGM theory is an axiomatic theory: it provides a list of "rationality" axioms for belief revision and provides representation theorems.[10] Although the AGM theory was developed within the language of propositional logic, it can be restated in terms of a set of states and a collection of possible items of information represented as events. We first introduce the non-probabilistic version of the theory and then add graded beliefs, that is, probabilities.

Let $W$ be a finite set of states and $\mathcal{E} \subseteq 2^W$ a collection of events (subsets of $W$) representing possible items of information; we assume that $W \in \mathcal{E}$ and $\varnothing \notin \mathcal{E}$. To represent initial beliefs and revised beliefs we introduce a function $f : \mathcal{E} \to 2^W$, which we call a belief revision function.

[10] We will not list and discuss the axioms here. The interested reader can consult http://plato.stanford.edu/entries/formal-belief/ or, for a discussion which is closer to the approach followed in this section, see Bonanno (2009).

**Definition 8.6.** Let $W$ be a finite set of states and $\mathcal{E}$ a collection of events such that $W \in \mathcal{E}$ and $\varnothing \notin \mathcal{E}$. A *belief revision function* is a function $f:\mathcal{E} \to 2^W$ that satisfies the following properties: for every $E \in \mathcal{E}$, (1) $f(E) \subseteq E$ and (2) $f(E) \neq \varnothing$.

The interpretation of a belief revision function is as follows. First of all, $f(W)$ represents the initial beliefs, namely the set of states that the individual initially considers possible.[11] Secondly, for every $E \in \mathcal{E}$, $f(E)$ is the set of states that the individual would consider possible if informed that the true state belongs to $E$; thus $f(E)$ represents the individual's *revised beliefs* after receiving information $E$.[12]

One of the implications of the AGM axioms for belief revision is the following condition, which is known as *Arrow's Axiom* (proposed by the Nobel laureate Ken Arrow in the context of rational choice, rather than rational belief revision):

$$\text{if } E,F \in \mathcal{E},\ E \subseteq F \text{ and } E \cap f(F) \neq \varnothing \text{ then } f(E) = E \cap f(F).$$

Arrow's Axiom says that if information $E$ implies information $F$ ($E \subseteq F$) and there are states in $E$ that would be considered possible upon receiving information $F$ ($E \cap f(F) \neq \varnothing$), then the states that the individual would consider possible if informed that $E$ are precisely those that belong to both $E$ and $f(F)$ ($f(E) = E \cap f(F)$). Although necessary for a belief revision policy that satisfies the AGM axioms, Arrow's Axiom is not sufficient. Before stating the necessary and sufficient condition for rational belief revision, we remind the reader of the notion of a complete and transitive relation on a set $W$

[11] If the initial beliefs were to be expressed probabilistically, by means of a probability distribution $P$ over $W$, then $f(W)$ would be the support of $P$, that is, the set of states to which $P$ assigns positive probability. Thus $f(W)$ would be the smallest event of which the individual would initially be certain: she would initially be certain of any event $F$ such that $f(W) \subseteq F$.

[12] If the revised beliefs after receiving information $E$ were to be expressed probabilistically, by means of a probability distribution $P_E$ over $W$, then $f(E)$ would be the support of $P_E$, that is, the set of states to which $P_E$ assigns positive probability. Thus $f(E)$ would be the smallest event of which the individual would be certain after having been informed that $E$: according to her revised beliefs she would be certain of any event $F$ such that $f(E) \subseteq F$. [Note that, since - by assumption - $f(E) \subseteq E$, the individual is assumed to be certain of the information received (e.g. because she trusts the source of the information).]

(Remark 1.1, Chapter 1). In Chapter 1 the relation was denoted by $\succsim$ and was interpreted in terms of preference: $o \succsim o'$ was interpreted as "the individual considers outcome $o$ to be at least as good as outcome $o'$". In the present context the relation will be denoted by $\precsim$ and the interpretation is in terms of "plausibility": $w \precsim w'$ is to be read as "the individual considers state $w$ to be *at least as plausible* as state $w'$" ($w \prec w'$ means that $w$ is considered to be *more plausible* than $w'$ and $w \sim w'$ means that $w$ is considered to be *just as plausible* as $w'$).[13] An alternative reading of $w \precsim w'$ is "*$w$ weakly precedes $w'$ in terms of plausibility*".

**Definition 8.7.** A *plausibility order* on a set of states $W$ is a binary relation $\precsim$ on $W$ that is complete (for every two states $w_1$ and $w_2$, either $w_1 \precsim w_2$ or $w_2 \precsim w_1$ or both) and *transitive* (if $w_1 \precsim w_2$ and $w_2 \precsim w_3$ then $w_1 \precsim w_3$). We define $w \prec w'$ as "$w \precsim w'$ and $w' \not\precsim w$" and $w \sim w'$ as "$w \precsim w'$ and $w' \precsim w$").

**Theorem 8.1** (Grove, 1988; Bonanno, 2009). Let $W$ be a finite set of states, $\mathcal{E}$ a collection of events (representing possible items of information), with $W \in \mathcal{E}$ and $\varnothing \notin \mathcal{E}$, and $f : \mathcal{E} \to 2^W$ a belief revision function (Definition 8.6). Then the belief revision policy represented by the function $f$ is compatible with the AGM axioms of belief revision if and only if there exists a plausibility order $\precsim$ on $W$ that *rationalizes* $f$ in the sense that, for every $E \in \mathcal{E}$, $f(E)$ is the set of most plausible states in $E$: $f(E) = \{w \in E : w \precsim w' \text{ for every } w' \in E\}$.

**Definition 8.8.** A belief revision function $f : \mathcal{E} \to 2^W$ which is rationalized by a plausibility order is called an *AGM belief revision function*.

**Remark 8.2.** An AGM belief revision function satisfies Arrow's Axiom (the reader is asked to prove this in Exercise 8.6). The converse is not true: it is possible for a belief revision function $f : \mathcal{E} \to 2^W$ to satisfy Arrow's Axiom and yet fail to be rationalized by a plausibility order.

[13] One reason for the reversal of notation is that, in this context, when representing the order $\precsim$ numerically, it will be convenient to assign *lower* values to *more* plausible histories. Another reason is that it is the standard notation in the extensive literature that deals with AGM belief revision (for a recent survey of this literature the reader is referred to the special issue of the *Journal of Philosophical Logic*, Vol. 40, Issue 2, April 2011).

Within the context of probabilistic beliefs, let $P$ be the probability distribution on the set of states $W$ that represents the initial beliefs and $P_E$ be the probability distribution representing the updated beliefs after receiving not surprising information $E$ (not surprising in the sense that $P(E) > 0$). The *support* of a probability distribution $P$, denoted by $Supp(P)$, is the set of states to which $P$ assigns positive probability: $Supp(P) = \{w \in W : P(w) > 0\}$. The rule for *updating* beliefs upon receiving information $E$ (Definition 8.2) implies the following:

if $E \cap Supp(P) \neq \varnothing$ (that is, if $P(E) > 0$) then $Supp(P_E) = E \cap Supp(P)$.

We call this the *qualitative belief updating rule*. It is easy to check that the qualitative belief updating rule is implied by Arrow's Axiom.[14] Thus, by Remark 8.2, an AGM belief revision function has incorporated in it the qualitative belief updating rule. In other words, *belief updating is included in the notion of AGM belief revision*. A belief revision function, however, goes beyond belief updating because it also encodes new beliefs after receipt of surprising information (that is, after being informed of an event $E$ such that $P(E) = 0$).

What is the probabilistic version of AGM belief revision? It turns out that in order to obtain probabilistic beliefs we only need to make a simple addition to an AGM belief revision function $f : \mathcal{E} \to 2^W$. Let $P_0$ be any full-support probability distribution on $W$ (that is, $P_0(w) > 0$, for every $w \in W$). Then, for every $E \in \mathcal{E}$ let $P_E$ be the probability distribution obtained by conditioning $P_0$ on $f(E)$:

$$P_E(w) = P_0\left(w \mid f(E)\right) = \begin{cases} \dfrac{P_0(w)}{\sum\limits_{w' \in f(E)} P_0(w')} & \text{if } w \in f(E) \\ 0 & \text{if } w \notin f(E) \end{cases}$$

Then $P_W$ gives the initial probabilistic beliefs and, for every other $E \in \mathcal{E}$, $P_E$ gives the revised probabilistic beliefs after receiving information $E$. The collection $\{P_E\}_{E\in\mathcal{E}}$ of probability distributions on $W$ so obtained gives the individual's *probabilistic* belief revision policy (while the function $f:\mathcal{E}\to 2^W$ gives the *qualitative* probabilistic belief revision policy).

[14] For every event $E$ representing a possible item of information, let $P_E$ be the probability distribution on $E$ representing the *revised* beliefs of the individual after receiving information $E$. Let $P$ be the probability distribution on $W$ representing the individual's *initial* beliefs. Define the following belief revision function $f$ : $f(W) = Supp(P)$ and $f(E) = Supp(P_E)$. Suppose that $f$ satisfies Arrow's Axiom. Then, for every event $E$, if $E \cap f(W) \neq \varnothing$ [that is, if $E \cap Supp(P) \neq \varnothing$] then $f(E) = E \cap f(W)$ [that is, $Supp(P_E) = E \cap Supp(P)$].

**Definition 8.9.** Let $W$ be a finite set of states and $\mathcal{E}$ a collection of events such that $W\in\mathcal{E}$ and $\varnothing\notin\mathcal{E}$. A *probabilistic belief revision policy* is a collection $\{P_E\}_{E\in\mathcal{E}}$ of probability distributions on $W$ such that, for every $E\in\mathcal{E}$, $Supp(P_E)\subseteq E$. $P_W$ represents the initial beliefs and, for every other $E\in\mathcal{E}$, $P_E$ represents the revised beliefs after receiving information $E$.
The collection $\{P_E\}_{E\in\mathcal{E}}$ is called an *AGM probabilistic belief revision policy* if it satisfies the following properties:

(1) there exists a plausibility order $\precsim$ on $W$ such that, for every $E\in\mathcal{E}$, $Supp(P_E)$ is the set of most plausible states in $E$, that is, $Supp(P_E)=\{w\in E: w\precsim w' \text{ for every } w'\in E\}$,[15] and

(2) there exists a full-support probability distribution $P_0$ on $W$ such that, for every $E\in\mathcal{E}$, $P_E$ is equal to the probability distribution obtained by conditioning $P_0$ on $Supp(P_E)$.

As we will see in Part IV, belief revision is very important in dynamic (or extensive-form) games. In such games a player may find herself at an information set that, according to her initial beliefs, had zero probability of being reached and thus will have to form new beliefs reflecting the unexpected information.

This is a good time to test your understanding of the concepts introduced in this section, by going through the exercises in Section 8.E.3 of Appendix 8.E at the end of this chapter.

[15] This condition says that if one defines the function $f:\mathcal{E}\to 2^W$ by $f(E)=Supp(P_E)$ then this function $f$ is an AGM belief revision function.

## 8.4 Harsanyi consistency of beliefs or like-mindedness

[Note: the material of this section will not be needed until Chapter 13. It is presented here for completeness on the topic of beliefs.]

We can easily extend the analysis to the case of two or more individuals. We already know how to model interactive knowledge by means of information partitions; the addition of beliefs is a simple step: we merely add, for every individual and for every information set, a probability distribution over the elements of that information set. A two-person example is shown in Figure 8.4 below.

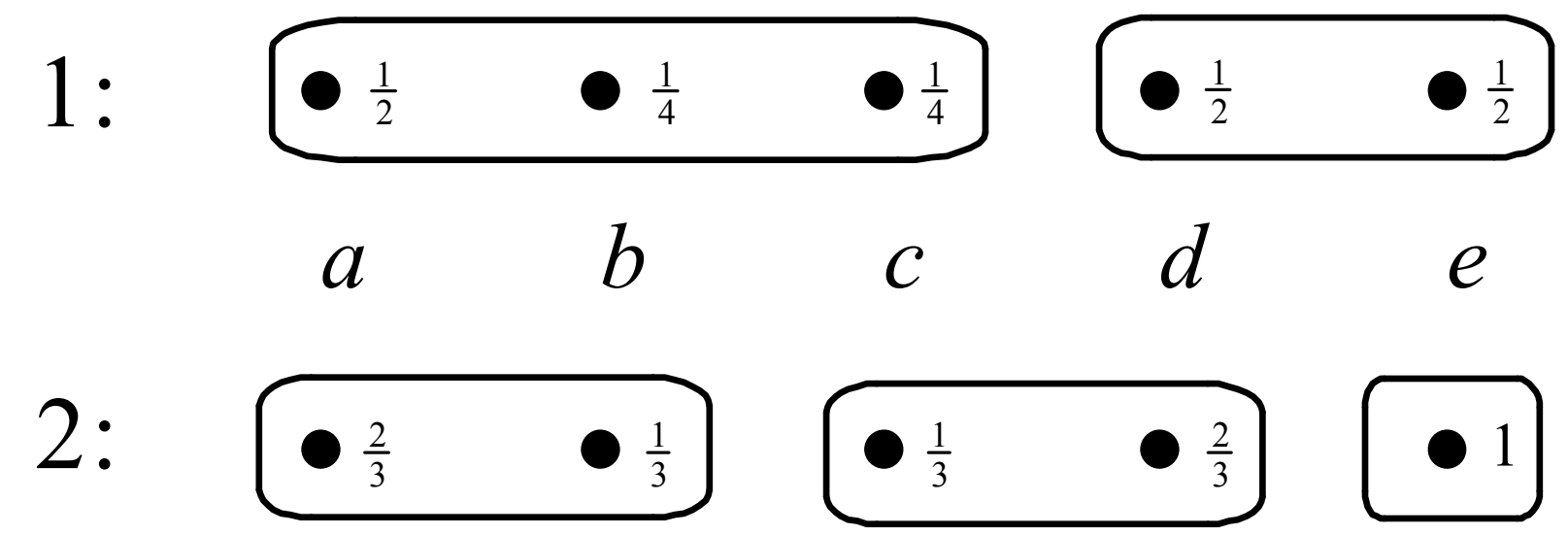

**Figure 8.4**

In this example, at every state, the two individuals hold different beliefs. For example, consider event $E = \{b,c\}$ and state *a*. Individual 1 attaches probability $\frac{1}{2}$ to *E* while Individual 2 attaches probability $\frac{1}{3}$ to *E*. Can two "equally rational" individuals hold different beliefs? The answer is: of course! In the above example it is not surprising that the two individuals assign different probabilities to the same event *E*, because they have different information. If the true state is *a*, then Individual 1's information is that the true state is either *a* or *b* or *c*, while Individual 2 considers only *a* and *b* possible (Individual 2 knows more than Individual 1). Is there a precise way of expressing the fact that two individuals assign different probabilities to an event *precisely because they have different information?* In the above example we could ask the hypothetical question: if Individual 1 had the same information as Individual 2, would he agree with Individual 2's assessment that the probability of event $E = \{b,c\}$ is $\frac{1}{3}$? This is, of course, a counterfactual question. The answer to this counterfactual question is affirmative: imagine giving Individual 1 the information that the true state is either *a* or *b*; then - according to Definition 8.3 - he would update his beliefs from $\begin{pmatrix} a & b & c \\ \frac{1}{2} & \frac{1}{4} & \frac{1}{4} \end{pmatrix}$ to $\begin{pmatrix} a & b \\ \frac{2}{3} & \frac{1}{3} \end{pmatrix}$ and thus have the same beliefs as Individual 2.

We say that two individuals are *like-minded* if it is the case that they would have the same beliefs if they had the same information. It is not straightforward how to turn this into a precise definition. Consider, again, the example of Figure 8.4 and event $E = \{b,c\}$. Above we asked the question "what would Individual 1 believe if he knew *as much as* Individual 2?". This is a simple question because we can imagine giving more information to Individual 1 and have him update his beliefs based on that information. However, we could also have asked the question "what would Individual 2 believe if he knew *as little as* Individual 1?". In this case we would have to imagine "taking away information" from Individual 2, by increasing his information set from $\{a,b\}$ to $\{a,b,c\}$. This is *not* the same as asking Individual 2 to update his beliefs based on information $\{a,b,c\}$, because updating on something you already know leaves your beliefs unchanged.

There is a sense in which the beliefs of the two individuals of Figure 8.4 are "in agreement": for example, they both consider state $a$ twice as likely as state $b$. One could try to use this condition to define like-mindedness: for every two states $x$ and $y$, whenever two individuals consider both $x$ and $y$ as possible (given their information) then they agree on the relative likelihood of $x$ versus $y$. Unfortunately, this condition is too weak. To see this, consider the three-individual example of Figure 8.5.

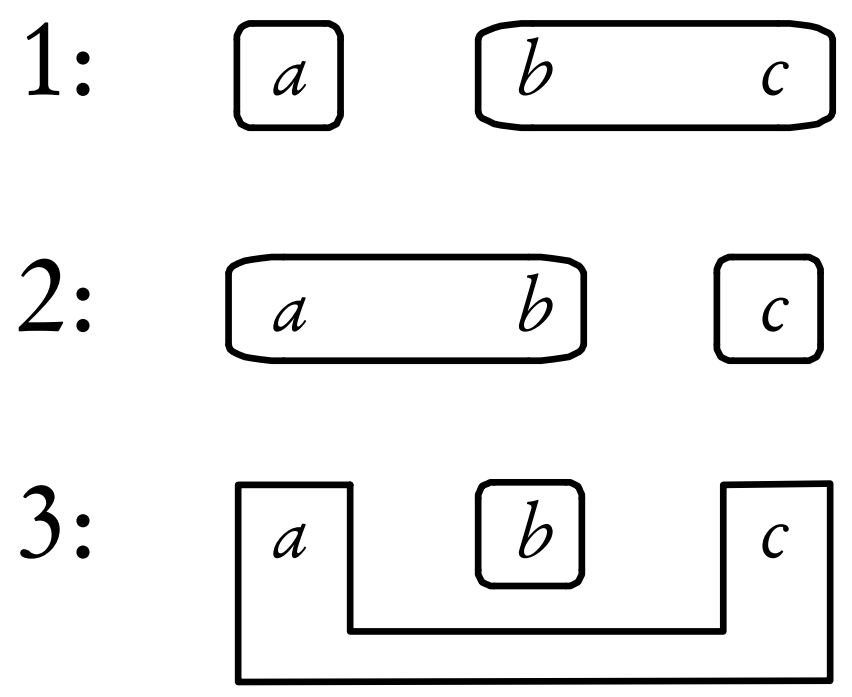

**Figure 8.5**

No matter what the true state is, we cannot find two individuals and two states $x$ and $y$ that both individuals consider possible. Thus any beliefs would make the three individuals like-minded. For example, Individual 1 at his information set $\{b,c\}$ could assign equal probability to $b$ and $c$, Individual 2 at her information set $\{a,b\}$ could assign probability $\frac{3}{4}$ to $a$ and $\frac{1}{4}$ to $b$ and Individual 3 at his information set $\{a,c\}$ could assign probability $\frac{1}{4}$ to $a$ and $\frac{3}{4}$ to $c$. But then, putting together the beliefs of Individuals 1 and 2 we would have that $a$ is judged to be three times more likely than $c$ (according to Individual 1, $a$ is three times more likely than $b$ and, according to Individual 2, $b$ is just as likely as $c$), while Individual 3 has the opposite judgment that $c$ is three times more likely than $a$.

In order to give a precise definition of like-mindedness we need to introduce some notation. Let there be $n$ individuals ($n \geq 2$). Let $W$ be a set of states and let $\mathcal{I}_i$ be the information partition of individual $i$ ($i \in \{1,...,n\}$). As usual, if $w$ is a state, we denote by $I_i(w)$ the element (information set) of the partition $\mathcal{I}_i$ that contains $w$. Let $P_{i,w}$ be the beliefs of individual $i$ at state $w$, that is, $P_{i,w}$ is a probability distribution over $I_i(w)$. Clearly, we can think of $P_{i,w}$ as a probability distribution over the entire set of states $W$ satisfying the property that if $w' \notin I_i(w)$ then $P_{i,w}(w') = 0$. Note that, for every individual $i$, there is just one probability distribution over $I_i(w)$ and thus if $w' \in I_i(w)$ then $P_{i,w'} = P_{i,w}$.

**Definition 8.10.** A probability distribution $P$ over $W$ is called a *common prior* if, for every individual $i$ and for every state $w$, (1) $P\left(I_i(w)\right) > 0$ and (2) updating $P$ on $I_i(w)$ (see Definition 8.4) yields precisely $P_{i,w}$, that is, $P\left(w' \mid I_i(w)\right) = P_{i,w}(w')$, for every $w' \in W$.

When a common prior exists we say that the individuals are *like-minded* or that the individuals' beliefs are *Harsanyi consistent*.[16]

[16] John Harsanyi, who in 1994 won the Nobel prize in Economics (together with Reinhardt Selten and John Nash), introduced the theory of games of incomplete information which will be the object of Part V. In that theory the notion of Harsanyi consistency plays a crucial role.

For instance, in the example of Figure 8.4, reproduced below, a common prior exists and thus the two individuals are like-minded.

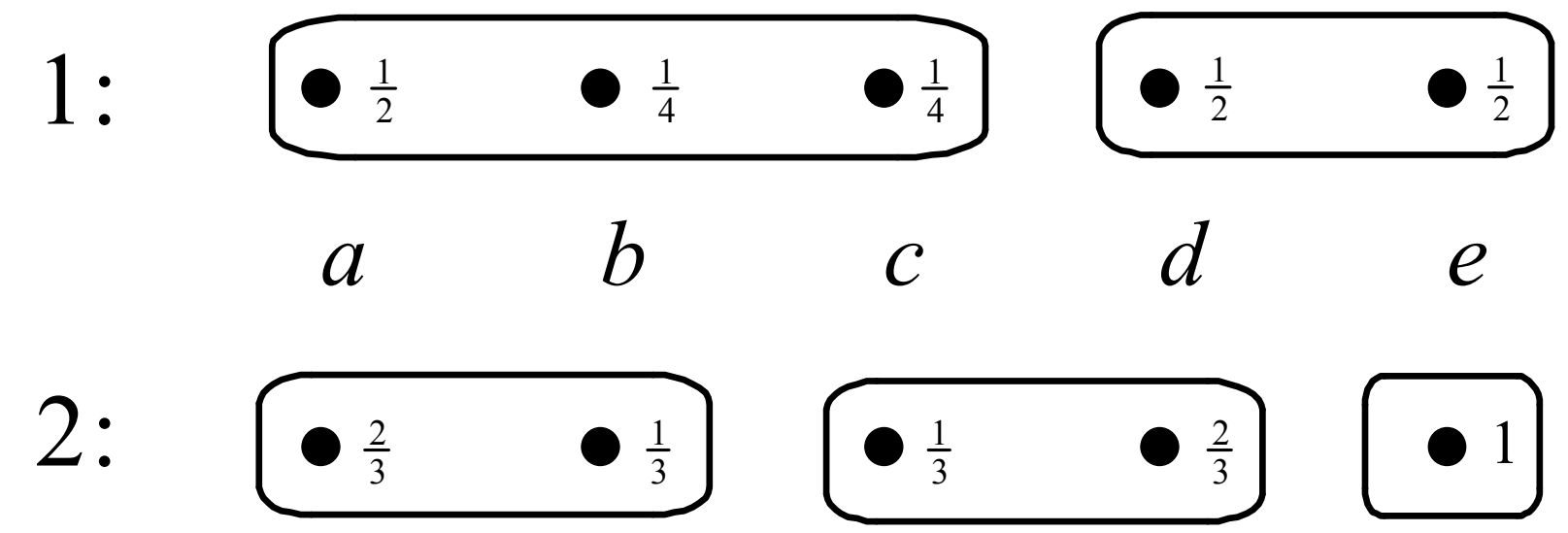

A common prior is $\begin{pmatrix} a & b & c & d & e \\ \frac{2}{8} & \frac{1}{8} & \frac{1}{8} & \frac{2}{8} & \frac{2}{8} \end{pmatrix}$ (the reader should convince himself/herself that, indeed, updating this probability distribution on each information set yields the probability distribution written inside that information set).

How can we determine if a common prior exists? The issue of existence of a common prior can be reduced to the issue of whether a system of equations has a solution. To see this, let us go back to the example of Figure 8.4, reproduced above. A common prior would be a probability distribution $\begin{pmatrix} a & b & c & d & e \\ p_a & p_b & p_c & p_d & p_e \end{pmatrix}$ that satisfies the following conditions.

1. Updating on information set $\{a,b,c\}$ of Individual 1 we need

$$\frac{p_b}{p_a + p_b + p_c} = \tfrac{1}{4} \text{ and } \frac{p_c}{p_a + p_b + p_c} = \tfrac{1}{4}$$

   Note that from these two the third condition follows, namely

$$\frac{p_a}{p_a + p_b + p_c} = \tfrac{1}{2}.$$

2. Updating on information set $\{d,e\}$ of Individual 1 we need $\frac{p_d}{p_d + p_e} = \tfrac{1}{2}$,

   from which it follows that $\frac{p_e}{p_d + p_e} = \tfrac{1}{2}$.

3. Updating on information set $\{a,b\}$ of Individual 2 we need $\frac{p_a}{p_a + p_b} = \frac{2}{3}$, from which it follows that $\frac{p_b}{p_a + p_b} = \frac{1}{3}$.

4. Updating on information set $\{c,d\}$ of Individual 2 we need $\frac{p_c}{p_c + p_d} = \frac{1}{3}$, from which it follows that $\frac{p_d}{p_c + p_d} = \frac{2}{3}$.

From the first we get $\boxed{p_b = p_c}$, from the second $\boxed{p_d = p_e}$, from the third $\boxed{p_a = 2p_b}$ and from the fourth $\boxed{p_d = 2p_c}$. Adding to these three equalities the requirement that $p_a + p_b + p_c + p_d + p_e = 1$, we have a system of five equations in five unknowns, which admits a unique solution, namely $\begin{pmatrix} a & b & c & d & e \\ \frac{2}{8} & \frac{1}{8} & \frac{1}{8} & \frac{2}{8} & \frac{2}{8} \end{pmatrix}$.

It is not always the case that a common prior exists. For instance, if we add to the example of Figure 8.5 the beliefs shown in Figure 8.6, then we get a situation where the individuals' beliefs are not Harsanyi consistent.

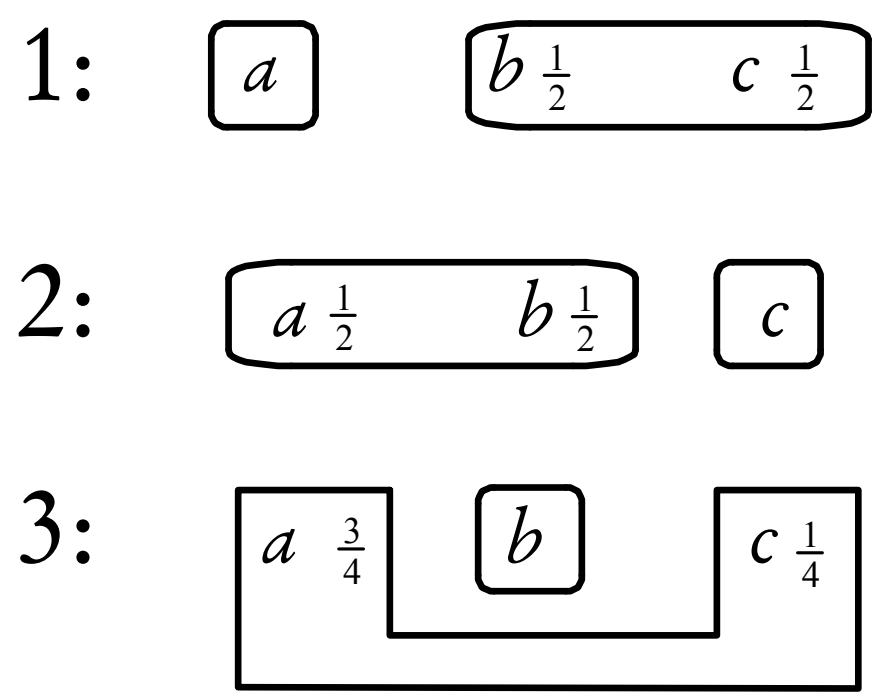

**Figure 8.6**

In this case, from the updating conditions of Individuals 1 and 2 one would get that $p_a = p_b$ and $p_b = p_c$, from which it follows that $p_a = p_c$; however, from the updating condition for Individual 3 we get that $p_a = 3p_c$, yielding a contradiction.

This is a good time to test your understanding of the concepts introduced in this section, by going through the exercises in Section 8.E.4 of Appendix 8.E at the end of this chapter.

## 8.5 Agreeing to disagree

Can two rational and like-minded individuals agree to disagree? This question was raised in 1976 by Robert Aumann (who received the Nobel Prize in Economics in 2005, together with Thomas Schelling).

As remarked above, it is certainly quite possible for two rational individuals to have different beliefs about a particular event, because they might have different information. Let us go back to the example of Figure 8.4, reproduced below in Figure 8.7, together with the common prior (showing that the individuals are like-minded). Suppose that the true state is $a$ and consider what the two individuals believe about event $E = \{b,c\}$. Individual 1's information set is $\{a,b,c\}$ and, given his beliefs at that information set, he attaches probability $\frac{1}{4}+\frac{1}{4}=\frac{1}{2}$ to $E$: let us denote this by "at state $a$, $P_1(E)=\frac{1}{2}$". Individual 2's information set is $\{a,b\}$ and, given her beliefs at that information set, she attaches probability $\frac{1}{3}$ to $E$: let us denote this by "at state $a$, $P_2(E)=\frac{1}{3}$". Thus the two individuals disagree about the probability of event $E$. Furthermore, *they know that they disagree*. To see this, let $\|P_1(E)=\frac{1}{2}\|$ be the event that Individual 1 attaches probability $\frac{1}{2}$ to $E$; then $\|P_1(E)=\frac{1}{2}\| = \{a, b, c\}$. Similarly, let $\|P_2(E)=\frac{1}{3}\|$ be the event that Individual 2 attaches probability $\frac{1}{3}$ to $E$; then $\|P_2(E)=\frac{1}{3}\| = \{a, b, c, d\}$. These are events and thus we can check at what states the two individuals know them. Using Definition 7.3 (Chapter 7), we have that $K_1\|P_2(E)=\frac{1}{3}\| = \{a,b,c\}$ and $K_2\|P_1(E)=\frac{1}{2}\| = \{a,b\}$. Hence

$$a \in \|P_1(E)=\tfrac{1}{2}\| \cap \|P_2(E)=\tfrac{1}{3}\| \cap K_1\|P_2(E)=\tfrac{1}{3}\| \cap K_2\|P_1(E)=\tfrac{1}{2}\|$$

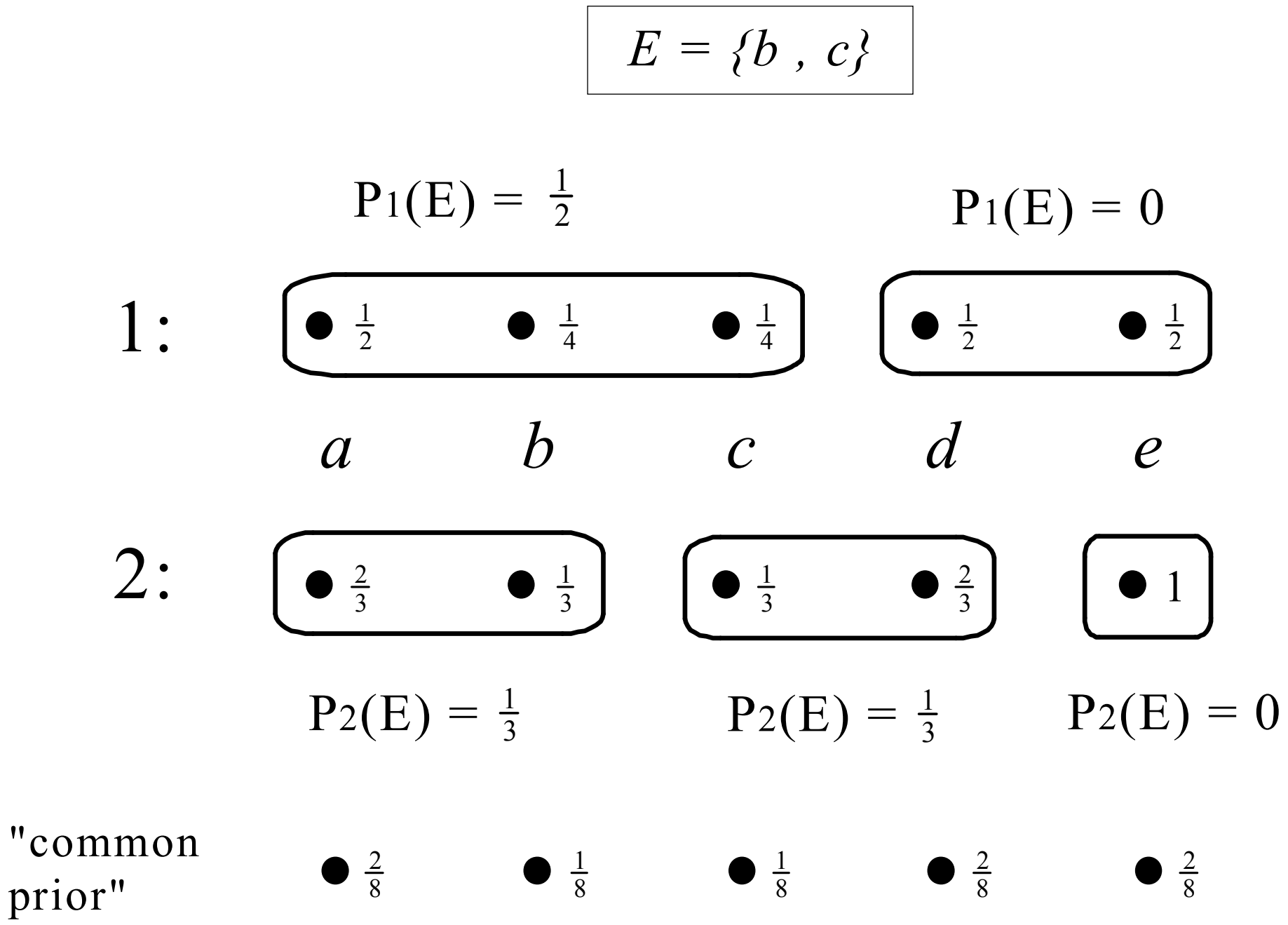

**Figure 8.7**

Thus at state $a$ not only do the individuals disagree, but they know that they disagree. However, their disagreement is not common knowledge. Indeed, $K_1K_2\|P_1(E)=\frac{1}{2}\|=\varnothing$ and thus $a \notin K_1K_2\|P_1(E)=\frac{1}{2}\|=\varnothing$, that is, at state $a$ it is not the case that Individual 1 knows that Individual 2 knows that Individual 1 assigns probability $\frac{1}{2}$ to event $E$. As the following theorem states, the opinions of two like-minded individuals about an event $E$ cannot be in disagreement and, at the same time, commonly known.

The following theorem is proved in Appendix 8.A at the end of this chapter.

**Theorem 8.1** (Agreement Theorem; Aumann, 1976). Let $W$ be a set of states and consider a knowledge-belief structure with two individuals, 1 and 2. Let $E$ be an event and let $p,q \in [0,1]$. Suppose that at some state $w$ it is common knowledge that Individual 1 assigns probability $p$ to $E$ and Individual 2 assigns probability $q$ to $E$. Then, if the individuals' beliefs are Harsanyi consistent (Definition 8.10), $p = q$. In other words, two like-minded individuals *cannot agree to disagree* about the probability of an event. Formally, if there exists a common prior and $w \in CK\left(\|P_1(E)=p\| \cap \|P_2(E)=q\|\right)$ then $p = q$.

Another way to think of this result is to imagine that the two individuals communicate their opinions to each other. Hearing that the other person has a different opinion is in itself a valuable piece of information, which ought to be incorporated (by updating) into one's beliefs. Thus sequential communication leads to changing beliefs. If, at the end of this process, the beliefs become common knowledge, then they must be identical. We shall illustrate this with an example.

Imagine two scientists who agree that the laws of Nature are such that the true state of the world must be one of seven, call them *a, b, c, d, e, f, g.* They also agree on the relative likelihood of these possibilities, which they take to be as shown in Figure 8.8.

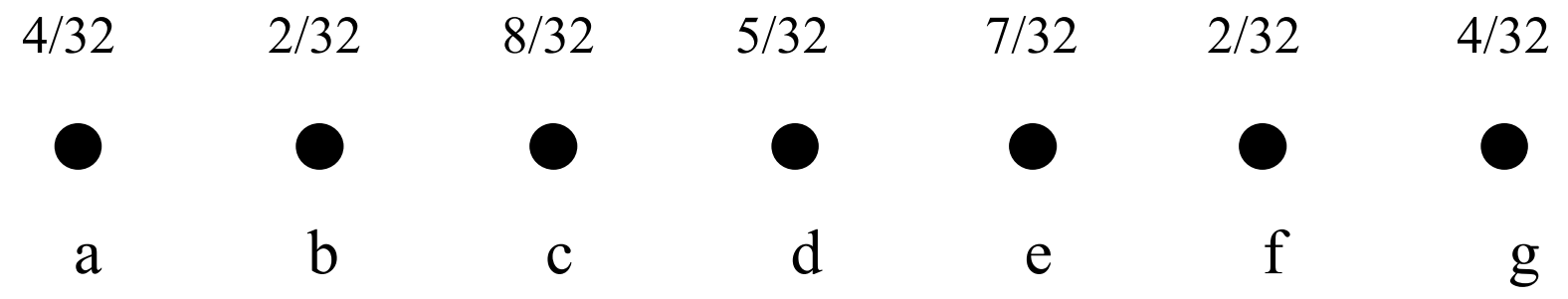

**Figure 8.8**

Experiments can be conducted to learn more. An experiment leads to a partition of the set of states. For example, if the true state is *a* and you perform an experiment then you might learn that the true state cannot be *d* or *e* or *f* or *g* but you still would not know which is the true state among the remaining ones. Suppose that the scientists agree that Scientist 1, from now on denoted by $S_1$, will perform experiment 1 and Scientist 2, denoted by $S_2$, will perform experiment 2. They also agree that each experiment would lead to a partition of the set of states as shown in Figure 8.9.

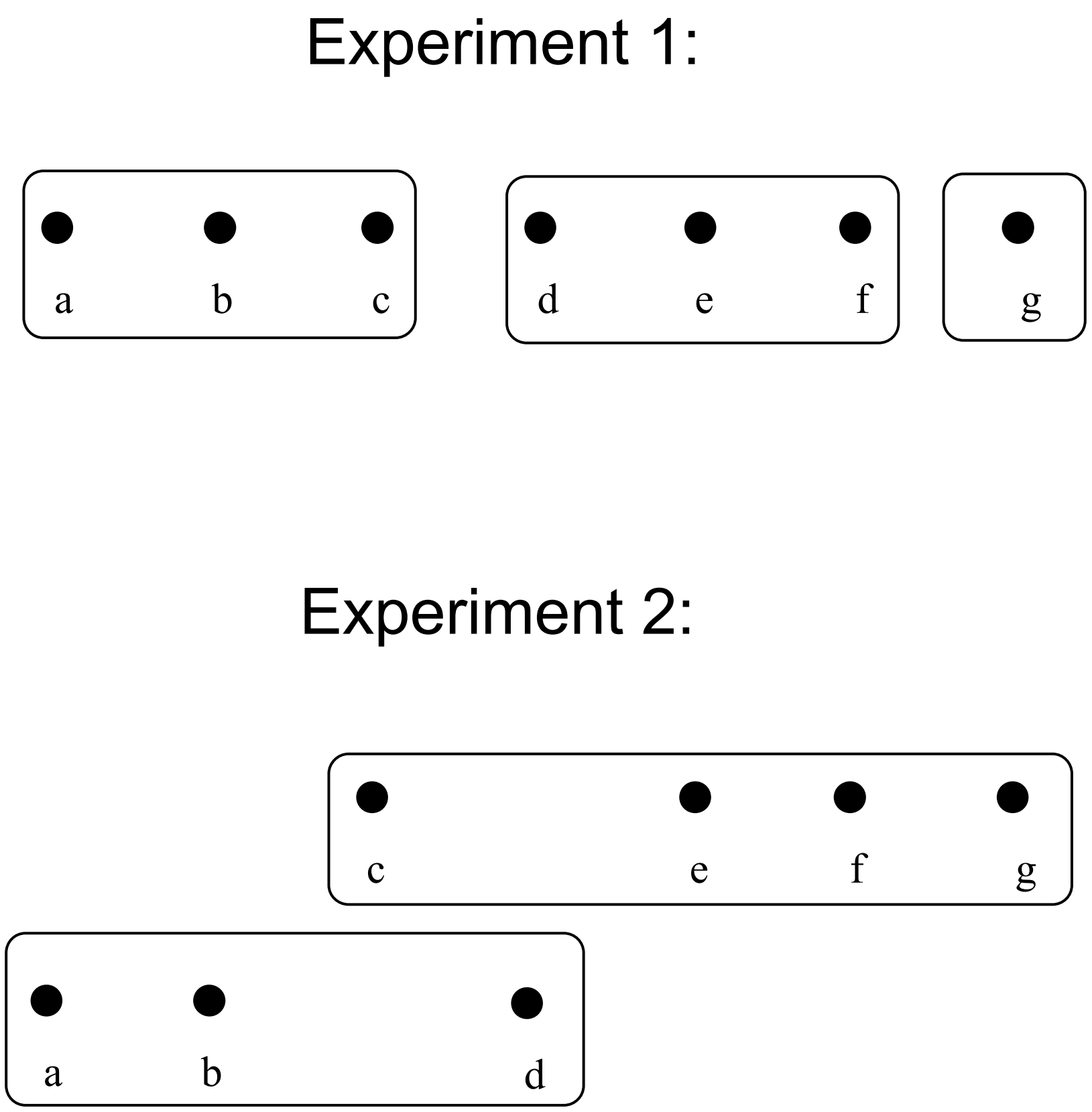

**Figure 8.9**

Suppose that the scientists are interested in establishing the truth of a proposition that is represented by the event $E = \{a, c, d, e\}$. Initially (given the shared probabilistic beliefs shown in Figure 8.8) they agree that the probability that $E$ is true is 75%:

$$P(E) = P(a) + P(c) + P(d) + P(e) = \tfrac{4}{32} + \tfrac{8}{32} + \tfrac{5}{32} + \tfrac{7}{32} = \tfrac{24}{32} = \tfrac{3}{4} = 75\%.$$

Before they perform the experiments they also realize that, depending on what the true state is, after the experiment they will have an updated probability of event $E$ conditional on what the experiment has revealed. For example, they agree that if one performs experiment 1 and the true state is $b$ then the experiment will yield the information $F = \{a,b,c\}$ and $P(E \mid F)$ is given by

$$P(E \mid F) = \frac{P(E \cap F)}{P(F)} = \frac{P(a) + P(c)}{P(a) + P(b) + P(c)} = \frac{\frac{4}{32} + \frac{8}{32}}{\frac{4}{32} + \frac{2}{32} + \frac{8}{32}} = \frac{12}{14} = \frac{6}{7} = 86\%.$$

[Note the interesting fact that sometimes experiments, although they are informative - that is, they reduce one's state of uncertainty - might actually induce one to become more confident of the truth of something that is false: in this case one would increase one's subjective probability that *E* is true from 75% to 86%, although *E* is actually false if the true state is *b*, as we hypothesized).]

We can associate with every cell of each experiment (that is, with every possible state of information yielded by the experiment) a new updated probability of event *E*, as shown in Figure 8.10.

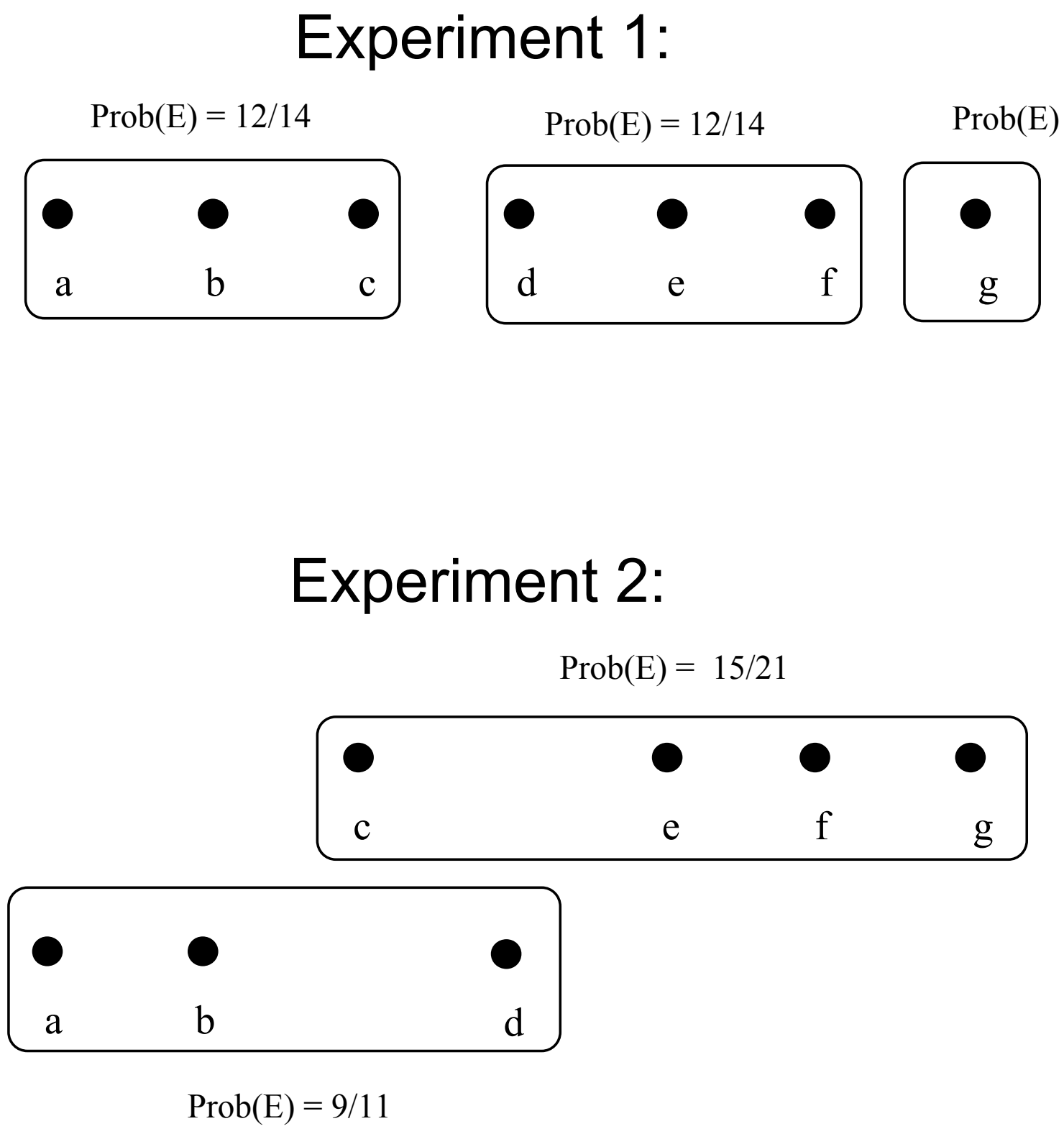

**Figure 8.10**

Suppose now that each scientist goes to her laboratory and performs the respective experiment (Scientist 1 experiment 1 and Scientist 2 experiment 2).

**Assume** also **that** **the true state is $f$** . Suppose that the two scientists send each other an e-mail communicating their new subjective estimates of the truth of $E$. Scientist 1 writes that he now attaches probability $\frac{12}{14}$ to $E$, while Scientist 2 says that she attaches probability $\frac{15}{21}$ to $E$. So their estimates disagree (not surprisingly, since they have performed different experiments and have collected different information). Should they be happy with these estimates? The answer is negative.

Consider first $S_1$ (Scientist 1). He hears that $S_2$ has a revised probability of $\frac{15}{21}$ (recall our assumption that the true state is $f$ ). What does he learn from this? He learns that the true state cannot be $d$ (if it had been, then he would have received an e-mail from $S_2$ saying "my new probability is $\frac{9}{11}$ "). $S_1$'s new state of knowledge and corresponding probabilities after receiving $S_2$'s e-mail are then as shown in Figure 8.11.

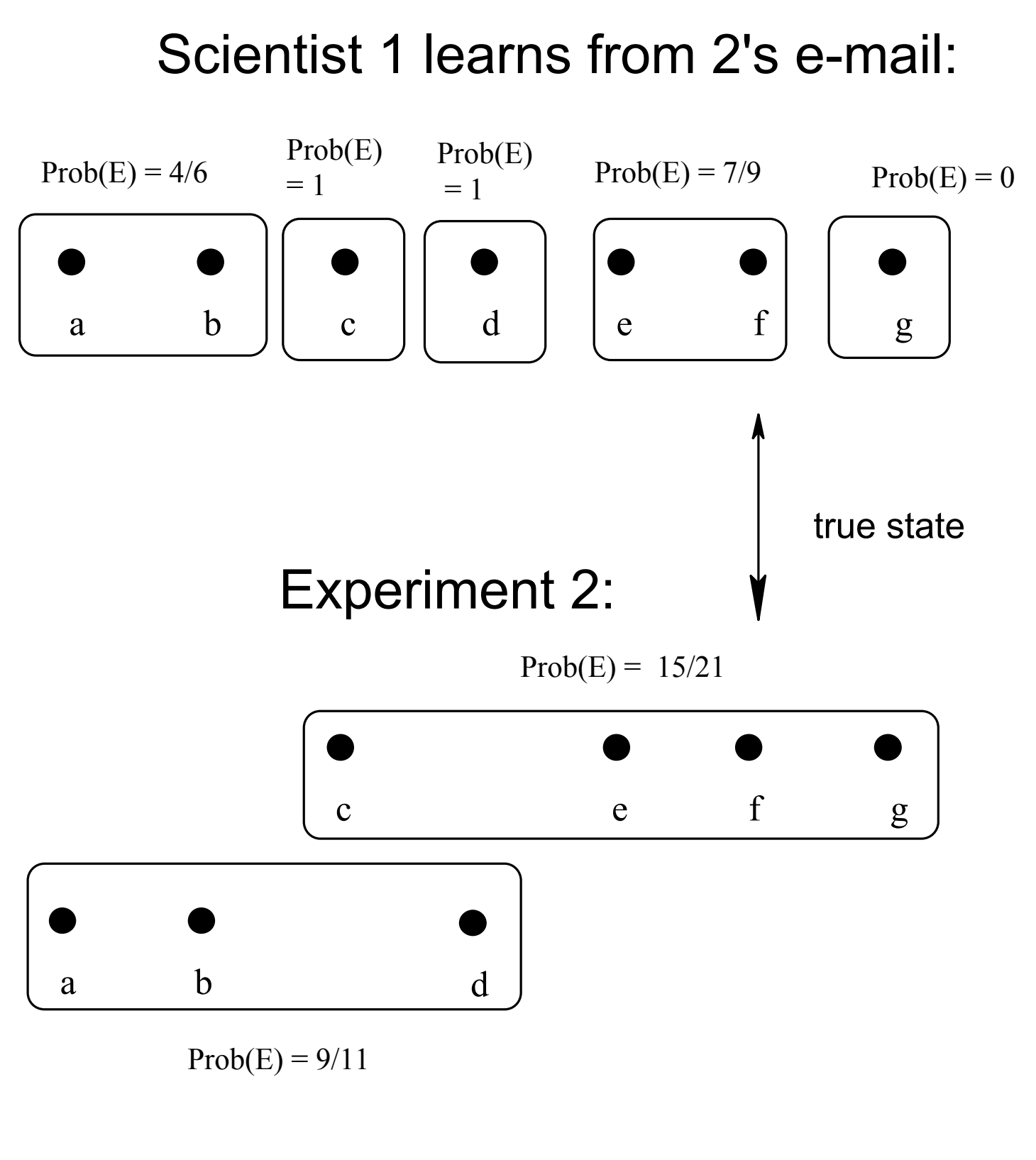

If $g$ is the true state, then experiment 1 reveals it.

If either $d$, $e$ or $f$ is the true state then experiment 1 reveals $\{d,e,f\}$. With $e$ or $f$ $S_2$'s e-mail is "15/21", with $d$ it is "9/11". Thus with $e$ or $f$ one can narrow down to $\{e,f\}$ while with $d$ one learns that the true state is $d$.

If the true state is $a$, $b$ or $c$ then experiment 1 reveals $\{a,b,c\}$. With $a$ or $b$ $S_2$ sends "9/11" while with $c$ she sends "15/21". Thus the e-mail splits $\{a,b,c\}$ into $\{a,b\}$ and $\{c\}$.

**Figure 8.11**

Consider now Scientist 2. From $S_1$'s e-mail she learns that $S_1$ has a new updated probability of $\frac{12}{14}$. What can she deduce from this? That the true state is *not* $g$ (if it had been $g$ then she would have received an e-mail from $S_1$ saying that her revised probability of $E$ was zero). Thus she can revise her knowledge partition by eliminating $g$ from her information set. A similar reasoning applies to the other states. $S_2$'s new state of knowledge and corresponding probabilities after receiving $S_1$'s e-mail are shown in Figure 8.12.

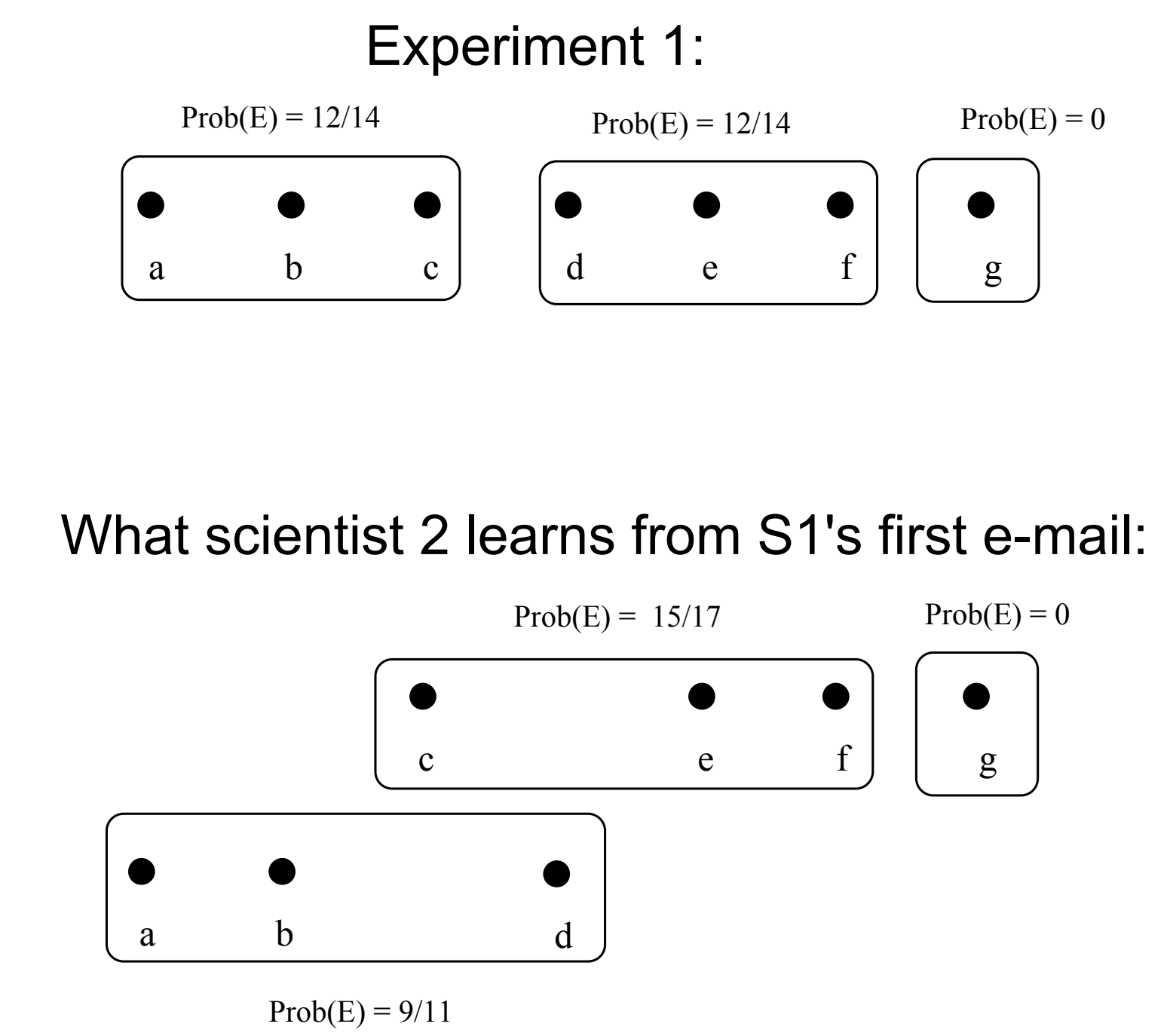

If either $c$, $e$, $f$ or $g$ is the true state then experiment 2 reveals $\{c,e,f,g\}$. With $c$, $e$ or $f$ $S_1$'s e-mail is "12/14", with $g$ it is "0". Thus with $c$, $e$ or $f$ one can narrow down to $\{c,e,f\}$ while with $g$ one learns that the true state is $g$.

If the true state is $a$, $b$ or $d$ then experiment 2 reveals $\{a,b,d\}$. In all three cases $S_1$'s e-mail is "12/14" and therefore nothing is learned.

**Figure 8.12**

Thus the new situation after the first exchange of e-mails is as shown in Figure 8.13 below.

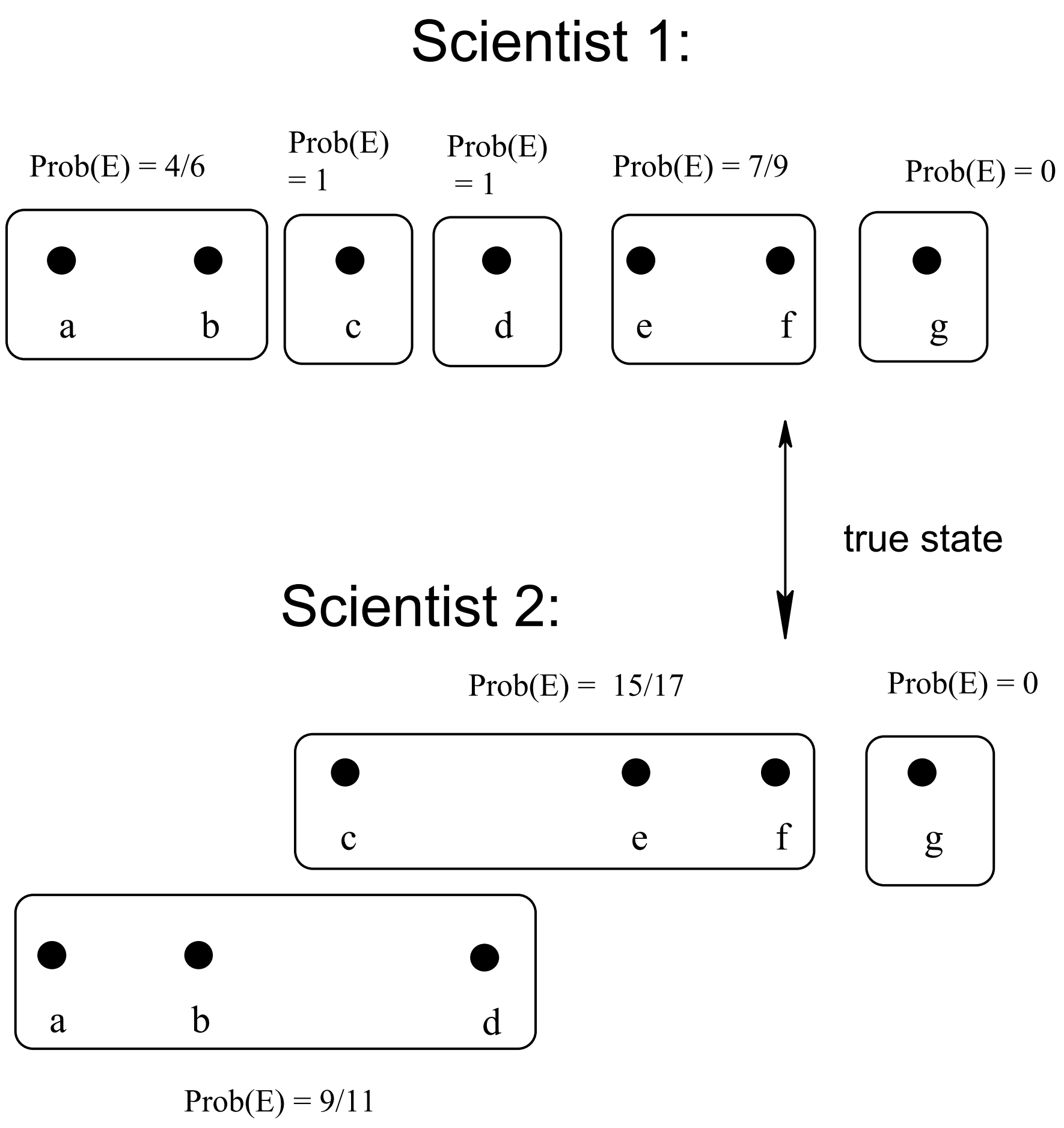

**Figure 8.13**

Now there is a second round of e-mails. $S_1$ communicates "in light of your e-mail, my new $P(E)$ is 7/9", while $S_2$ writes "after your e-mail I changed my $P(E)$ to 15/17". While $S_1$ learns nothing new from $S_2$'s second e-mail, $S_2$ learns that the true state cannot be $c$ (the second e-mail would have been "$P(E) = 1$" if $c$ had been the true state; in the hypothetical case where $S_2$'s revised information was $\{a,b,d\}$ then after $S_1$'s second e-mail it would split into $\{a,b\}$ and $\{d\}$). Thus the new situation is a shown in Figure 8.14.

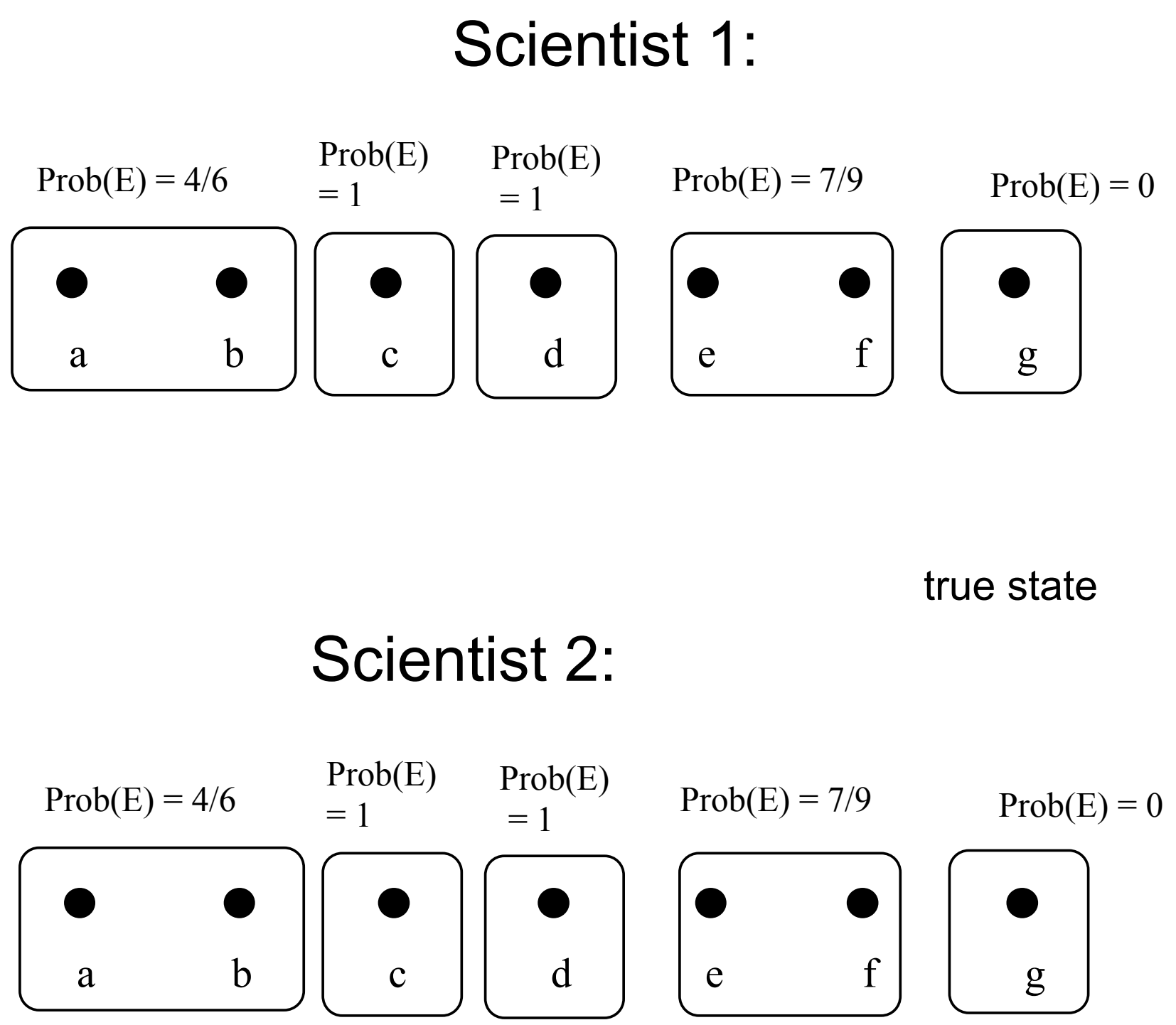

**Figure 8.14**

Now they have reached complete agreement: $P(E) = \frac{7}{9}$. Further exchanges do not convey any new information. Indeed it has become common knowledge (at state $f$) that both scientists estimate the probability of event $E$ to be $\frac{7}{9} = 78\%$ (before the experiments the probability of $E$ was judged to be $\frac{24}{32} = 75\%$; note, again, that with the experiments and the exchange of information they have gone farther from the truth than at the beginning!).

Notice that before the last step it was never common knowledge between the two what probability each scientist attached to $E$. When one scientist announced his subjective estimate, the other scientist found that announcement informative and revised her own estimate accordingly. After the exchange of two e-mails, the further announcement by one scientist of his/her estimate of the probability of $E$ would not make the other scientist change his/her own estimate: the announcement would reveal nothing new.

This is a good time to test your understanding of the concepts introduced in this section, by going through the exercises in Section 8.E.5 of Appendix 8.E at the end of this chapter.

# Appendix 8.A: Proof of the Agreement Theorem

First we prove the following.

**Lemma.** Let $W$ be a finite set of states, $P$ a probability distribution on $W$ and $E,F \subseteq W$ two events. Let $\{F_1, ..., F_m\}$ be a partition of $F$ (thus $F = F_1 \cup ... \cup F_m$ and any two $F_j$ and $F_k$ with $j \neq k$ are non-empty and disjoint). Suppose that $P(E|F_j) = q$ for all $j = 1, ..., m$. Then $P(E|F) = q$.

*Proof.* By definition of conditional probability, $P(E \mid F_j) = \dfrac{P(E \cap F_j)}{P(F_j)}$. Hence, since $P(E|F_j) = q$, we have that $P(E \cap F_j) = q\,P(F_j)$. Adding over $j$, the left-hand side becomes $P(E \cap F)$ $\big[$because $E \cap F = (E \cap F_1) \cup \ldots \cup (E \cap F_m)$ and for any $j$ and $k$ with $j \neq k$, $(E \cap F_j) \cap (E \cap F_k) = \varnothing$, so that $P(E \cap F) = P(E \cap F_1) + \ldots + P(E \cap F_m)\big]$ and the right-hand side becomes $qP(F)$. Hence $P(E \mid F) = \dfrac{P(E \cap F)}{P(F)} = \dfrac{q\,P(F)}{P(F)} = q$.

**Proof of Theorem 8.1**. Suppose that $CK\left(\|P_1(E) = p\| \cap \|P_2(E) = q\|\right) \neq \varnothing$. Let $P$ be a common prior. Select an arbitrary $w \in CK\left(\|P_1(E) = p\| \cap \|P_2(E) = q\|\right)$ and let $I_{CK}(w)$ be the cell of the common knowledge partition containing $w$. Consider Individual 1. $I_{CK}(w)$ is equal to the union of a collection of cells (information sets) of 1's information partition. On each such cell 1's conditional probability of $E$, using the common prior $P$, is $p$. Hence, by the above lemma, $P(E|I_{CK}(w)) = p$. A similar reasoning for Individual 2 leads to $P(E|I_{CK}(w)) = q$. Hence $p = q$. ■

# Appendix 8.E: Exercises

## 8.E.1. Exercises for Section 8.1: Probabilistic beliefs

The answers to the following exercises are in Appendix S at the end of this chapter.

**Exercise 8.1**. Let the set of states be $W = \{a, b, c, d, e, f, g, h, i\}$. Amy's initial beliefs are given by the following probability distribution:

| *state* | $a$ | $b$ | $c$ | $d$ | $e$ | $f$ | $g$ | $h$ | $i$ |
|---|---|---|---|---|---|---|---|---|---|
| *probability* | $\frac{11}{60}$ | 0 | $\frac{7}{60}$ | $\frac{9}{60}$ | $\frac{16}{60}$ | $\frac{5}{60}$ | $\frac{4}{60}$ | $\frac{8}{60}$ | 0 |

**(a)** Let $E = \{a, f, g, h, i\}$. What is the probability of *E*?

**(b)** Find all the events that Amy is certain of.

**Exercise 8.2**. Let *W* be a finite set of states and let *A* and *B* be two events. Recall that $A \cup B$ is the union of *A* and *B*, that is, the set of states that belong to either *A* or *B* (or both) and $A \cap B$ is the intersection of *A* and *B*, that is, the set of states that belong to both *A* and *B*. Explain why $P(A \cup B) = P(A) + P(B) - P(A \cap B)$.

## 8.E.2. Exercises for Section 8.2: Conditional probability, belief updating, Bayes' rule

The answers to the following exercises are in Appendix S at the end of this chapter.

**Exercise 8.3**. In a remote rural clinic with limited resources, a patient arrives complaining of low-abdomen pain. Based on all the information available, the doctor thinks that there are only four possible causes: a bacterial infection (*b*), a viral infection (*v*), cancer (*c*), internal bleeding (*i*). Of the four, only the bacterial infection and internal bleeding are treatable at the clinic. In the past the doctor has seen 600 similar cases and they eventually turned out to be as follows:

| Bacterial infection $b$ | Viral infection $v$ | Cancer $c$ | Internal bleeding $i$ |
|---|---|---|---|
| 140 | 110 | 90 | 260 |

Her probabilistic estimates are based on those past cases.

**(a)** What is the probability that the patient has a treatable disease?

There are two possible ways of gathering more information: a blood test and an ultrasound. A positive blood test will reveal that there is an infection, however it could be either a bacterial or a viral infection; a negative blood test rules out an infection (either bacterial or viral) and thus leaves cancer and internal bleeding as the only possibilities. The ultrasound, on the other hand, will reveal if there is internal bleeding.

**(b)** Suppose that the patient gets an ultrasound and it turns out that there is no internal bleeding. What is the probability that he does **not** have a treatable disease? What is the probability that he has cancer?

**(c)** If instead of getting the ultrasound he had taken the blood test and it had been positive, what would have been the probability that he had a treatable disease?

**(d)** Now let us go back to the hypothesis that the patient only gets the ultrasound and it turns out that there is no internal bleeding. He then asks the doctor: "if I were to take the blood test too, how likely is it that it would be positive?". What should the doctor's answer be?

**(e)** Finally, suppose that the patient gets both the ultrasound and the blood test and the ultrasound reveals that there is no internal bleeding, while the blood test is positive. How likely is it that he has a treatable disease?

**Exercise 8.4**. $A$ and $B$ are two events such that $P(A) = 0.2$, $P(B) = 0.5$ and $P(B \mid A) = 0.1$. Calculate $P(A \mid B)$ using Bayes' rule.

**Exercise 8.5**. A lab technician was asked to mark some specimens with two letters, one from the set {A, B, C} and the other from the set {E,F,G}. For example, a specimen could be labeled as AE or BG, etc. He had a **total of 220 specimens**. He has to file a report to his boss, consisting of the following table.

| LABEL | number |
|---|---|
| AE | |
| AF | |
| AG | |
| BE | |
| BF | |
| BG | |
| CE | |
| CF | |
| CG | |

Unfortunately, he does not remember all the figures. He had written some notes to himself, which are reproduced below. Help him fill in the above table with the help of his notes and the conditional probability formula. Here are the technician's notes:

**(a)** Of all the ones that he marked with an E, one fifth were also marked with an A and one fifth were marked with a B.

**(b)** He marked 36 specimens with the label CE.

**(c)** Of all the specimens that he marked with a C, $\frac{12}{33}$ were marked with a G.

**(d)** Of all the specimens, $\frac{23}{55}$ were marked with a C.

**(e)** The number of specimens marked BG was twice the number of specimens marked BE.

**(f)** Of all the specimens marked with an A, $\frac{3}{20}$ were marked with an E.

**(g)** Of all the specimens marked with an A, one tenth were marked with a G.

## 8.E.3. Exercises for Section 8.3: Belief revision

The answers to the following exercises are in Appendix S at the end of this chapter.

**Exercise 8.6.** Prove that an AGM belief revision function satisfies Arrow's Axiom (Remark 8.2).

**Exercise 8.7.** Let $W = \{a,b,c,d,e,g,h,k,m\}$ and let $\precsim$ be the following plausibility order on $W$ (we use the convention that if the row to which state $w$ belongs is above the row to which state $w'$ belongs then $w \prec w'$, and if $w$ and $w'$ belong to the same row then $w \sim w'$).

| | |
|---|---|
| most plausible | $b, g$ |
| | $c, k, m$ |
| | $d, h$ |
| | $e$ |
| least plausible | $a$ |

Let $\mathcal{E} = \{\{a,e\},\{d,e,k,m\},\{b,d,e,k\},W\}$. Find the belief revision function $f : \mathcal{E} \to 2^W$ that is rationalized by $\precsim$.

**Exercise 8.8.** As in Exercise 8.7, let $W = \{a,b,c,d,e,g,h,k,m\}$ and $\mathcal{E} = \left\{ \underbrace{\{a,e\}}_{E}, \underbrace{\{d,e,k,m\}}_{F}, \underbrace{\{b,d,e,k\}}_{G}, W \right\}$. Using the plausibility order of Exercise 8.7, find a collection of probability distributions $\{P_E, P_F, P_G, P_W\}$ that provides an AGM probabilistic belief revision policy (Definition 8.9). [There are many: find one.]

## 8.E.4. Exercises for Section 8.4: Harsanyi consistency of beliefs

The answers to the following exercises are in Appendix S at the end of this chapter.

**Exercise 8.9**. Consider the following knowledge-belief structure:

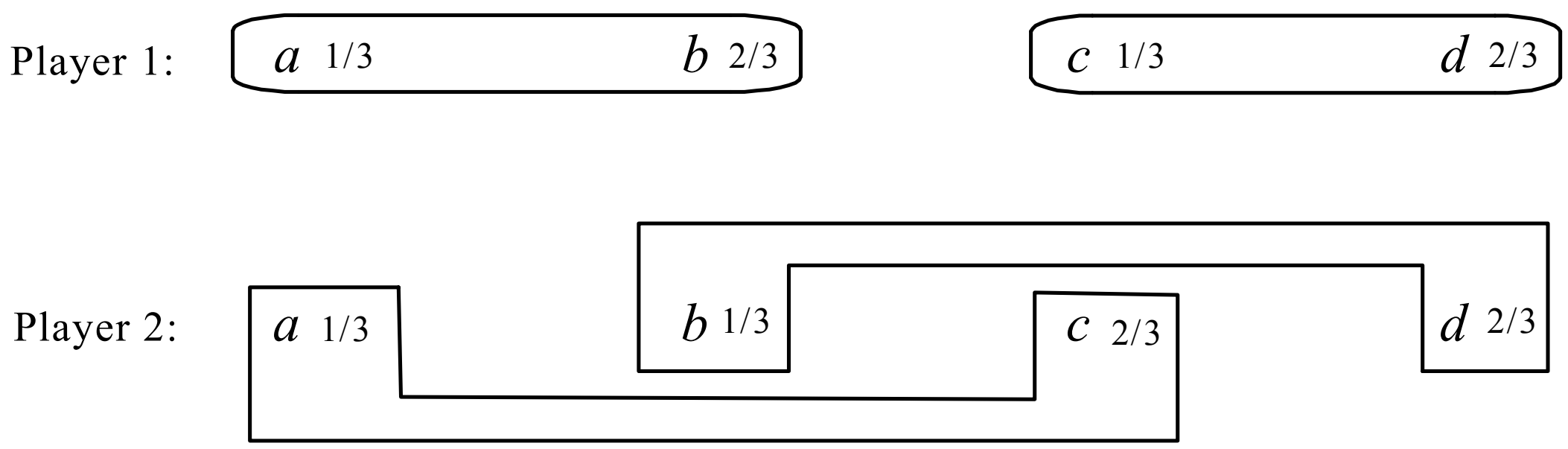

Are the individuals' beliefs Harsanyi consistent?

**Exercise 8.10**. Consider the following knowledge-belief structure:

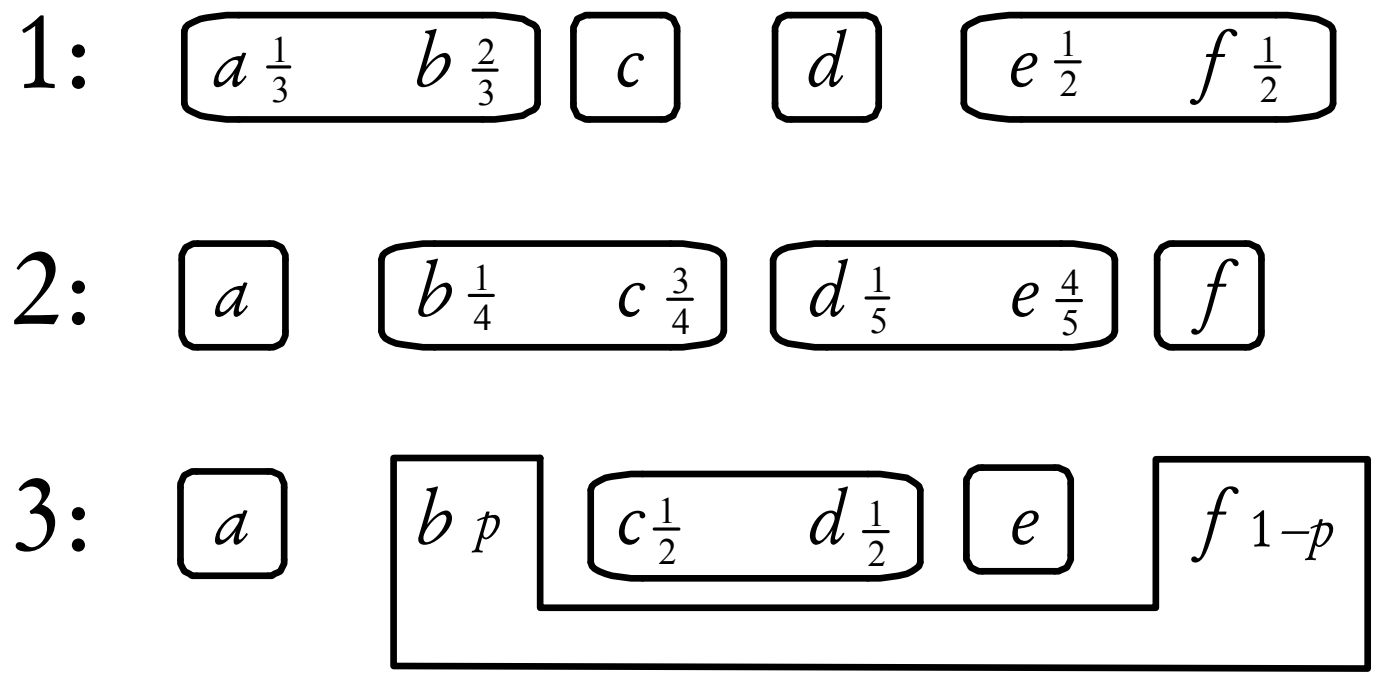

(**a**) Show that if $p = \frac{1}{13}$ then the beliefs are Harsanyi consistent.

(**b**) Show that if $p \neq \frac{1}{13}$ then the beliefs are *not* Harsanyi consistent.

## 8.E.5. Exercises for Section 8.5: Agreeing to disagree

The answers to the following exercises are in Appendix S at the end of this chapter.

**Exercise 8.11**. Consider the following knowledge-belief structure:

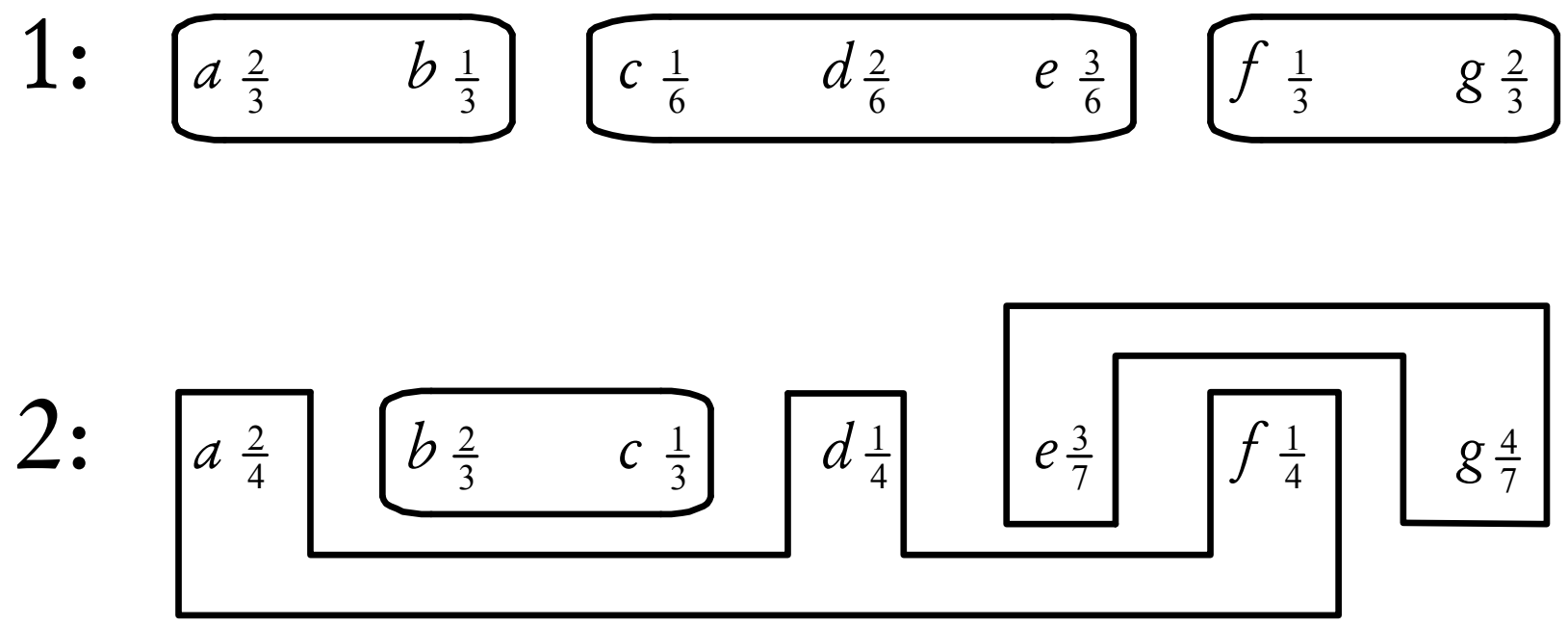

**(a)** Find the common knowledge partition.

**(b)** Find a common prior.

**(c)** Let $E = \{b,d,f\}$. **(c.1)** Is the probability that Individual 1 assigns to $E$ common knowledge? **(c.2)** Is the probability that Individual 2 assigns to $E$ common knowledge?

**(d)** Let $E = \{b,d,f\}$. **(d.1)** At state $b$ does Individual 1 know what probability Individual 2 assigns to $E$? **(d.2)** At state $c$ does Individual 1 know what probability Individual 2 assigns to $E$? **(d.3)** At state $f$ does Individual 1 know what probability Individual 2 assigns to $E$?

◊◊◊◊◊◊◊◊◊◊◊◊◊

**Exercise 8.12. Challenging Question.**

[This is known as the Monty Hall problem.]

You are a contestant in a show. You are shown three doors, numbered 1, 2 and 3. Behind one of them is a new car, which will be yours if you choose to open that door. The door behind which the car was placed was chosen randomly with equal probability (a die was thrown, if it came up 1 or 2 then the car was placed behind door 1, if it came up 3 or 4 then the car was placed behind door 2 and if it came up 5 or 6 then the car was placed behind door 3). You have to choose one door. Suppose that you have chosen door number 1. **Before the door is opened** the show host tells you that he knows where the car is and, to help you, he will open *one of the other two doors*, making sure that he opens a door behind which there is no car; if there are two such doors, then he will choose randomly with equal probability. Afterwards he will give you a chance to change your mind and switch to the other closed, **but you will have to pay $50 if you decide to switch**. Suppose that (initially you chose door 1 and) the host opens door 3 to show you that the car is not there. Should you switch from door 1 to door 2? [Assume that, if switching increases the probability of getting the car, then you find it worthwhile to pay $50 to switch.]

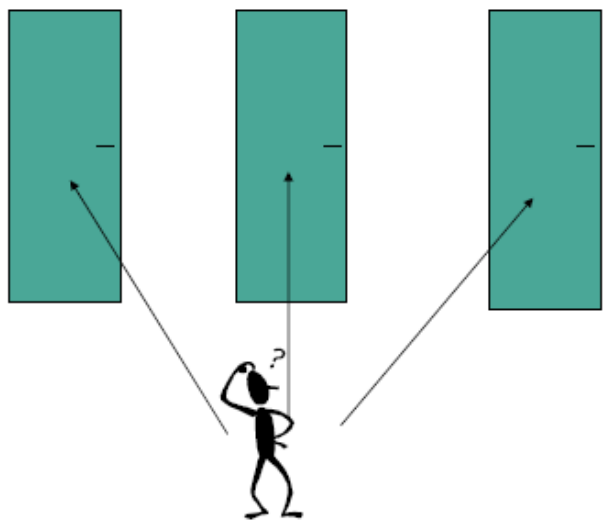

Let us approach this problem in two different ways.

**(a)** Method 1. Draw (part of) an extensive form with imperfect information where Nature moves first and chooses where the car is, then you choose one door and then the host chooses which door to open (of course, the host's choice is made according to the rules specified above). Reasoning about the information set you are in after you have pointed to door 1 and the host has opened door 3, determine if you should switch from door 1 to door 2.

**(b)** Method 2. Let $D_n$ denote the event that the car is behind door $n \in \{1,2,3\}$. Let $O_n$ denote the event that the host opens door $n$. The prior probabilities are $P(D_1) = P(D_2) = P(D_3) = \frac{1}{3}$. Compute $P(D_1 | O_3)$ using Bayes' rule (if $P(D_1 | O_3) = \frac{1}{2}$ then you should not switch, since there is a cost in switching).

# Appendix 8.S: Solutions to exercises

**Exercise 8.1**. The probability distribution is

| *state* | $a$ | $b$ | $c$ | $d$ | $e$ | $f$ | $g$ | $h$ | $i$ |
|---|---|---|---|---|---|---|---|---|---|
| *probability* | $\frac{11}{60}$ | 0 | $\frac{7}{60}$ | $\frac{9}{60}$ | $\frac{16}{60}$ | $\frac{5}{60}$ | $\frac{4}{60}$ | $\frac{8}{60}$ | 0 |

**(a)** Let $E=\{a,f,g,h,i\}$. Then $P(E)=P(a)+P(f)+P(g)+P(h)+P(i)=\frac{11}{60}+\frac{5}{60}+\frac{4}{60}+\frac{8}{60}+0=\frac{28}{60}=\frac{7}{15}$.

**(b)** Amy is certain of all events that have probability 1, namely $\{a,c,d,e,f,g,h\}$, $\{a,b,c,d,e,f,g,h\}$, $\{a,c,d,e,f,g,h,i\}$ and $\{a,b,c,d,e,f,g,h,i\}$.

**Exercise 8.2**. Since $P(A)=\sum_{w\in A}P(w)$ and $P(B)=\sum_{w\in B}P(w)$, when adding $P(A)$ to $P(B)$ the elements that belong to both $A$ and $B$ (that is, the elements of $A\cap B$) are added twice and thus we need to subtract $\sum_{w\in A\cap B}P(w)$ from $P(A)+P(B)$ in order to get $\sum_{w\in A\cup B}P(w)$, which is equal to $P(A\cup B)$.

**Exercise 8.3**. The probabilities are as follows:

| $b$ | $v$ | $c$ | $i$ |
|---|---|---|---|
| $\frac{140}{600}=\frac{7}{30}$ | $\frac{110}{600}=\frac{11}{60}$ | $\frac{90}{600}=\frac{3}{20}$ | $\frac{260}{600}=\frac{13}{30}$ |

**(a)** The event that the patient has a treatable disease is $\{b,i\}$. $P(\{b,i\})=P(b)+P(i)=\frac{7}{30}+\frac{13}{30}=\frac{2}{3}$.

**(b)** A negative result of the ultrasound is represented by the event $\{b,v,c\}$. A non-treatable disease is the event $\{v,c\}$. Thus

$$P(\{v,c\}\mid\{b,v,c\})=\frac{P(\{v,c\}\cap\{b,v,c\})}{P(\{b,v,c\})}=\frac{P(\{v,c\})}{P(\{b,v,c\})}=\frac{\frac{11}{60}+\frac{3}{20}}{\frac{7}{30}+\frac{11}{60}+\frac{3}{20}}=\frac{10}{17}.$$

$$P(c\mid\{b,v,c\})=\frac{P(c)}{P(\{b,v,c\})}=\frac{\frac{3}{20}}{\frac{7}{30}+\frac{11}{60}+\frac{3}{20}}=\frac{9}{34}.$$

**(c)** A positive blood test is represented by the event $\{b,v\}$. A treatable disease is the event $\{b,i\}$. Thus

$$P(\{b,i\} \mid \{b,v\}) = \frac{P(\{b,i\} \cap \{b,v\})}{P(\{b,v\})} = \frac{P(b)}{P(\{b,v\})} = \frac{\frac{7}{30}}{\frac{7}{30}+\frac{11}{60}} = \frac{14}{25}.$$

**(d)** Here we want $P(\{b,v\} \mid \{b,v,c\}) = \frac{P(\{b,v\})}{P(\{b,v,c\})} = \frac{\frac{7}{30}+\frac{11}{60}}{\frac{7}{30}+\frac{11}{60}+\frac{3}{20}} = \frac{25}{34}$.

**(e)** The information is $\{b,v\} \cap \{b,v,c\} = \{b,v\}$. Thus we want $P(\{b,i\} \mid \{b,v\})$ which was calculated in Part (c) as $\frac{14}{25}$.

**Exercise 8.4.** $P(A \mid B) = \frac{P(B \mid A)P(A)}{P(B)} = \frac{(0.1)(0.2)}{0.5} = 0.04$.

**Exercise 8.5.** Let P(xy) be the fraction of the total number of specimens that were marked xy and P(x|y) be the number of specimens marked with an x as a fraction of the specimens that were marked with a y (this is like a conditional probability). This is how we can re-write the information contained in the technician's notes:

**(a)** P(A|E) = P(B|E) = $\frac{1}{5}$. It follows that the remaining three fifths were marked with a *C*, that is, P(*C*|*E*) = $\frac{3}{5}$.

**(b)** The number of CEs is 36. Thus P(CE) = $\frac{36}{220}$. Thus $P(C \mid E) = \frac{P(CE)}{P(E)}$.

Using (a) we get $\frac{3}{5} = \frac{\frac{36}{220}}{P(E)}$, that is, P(E) = $\frac{36}{220}\frac{5}{3} = \frac{3}{11}$. Hence the number of specimens marked with an E is $\frac{3}{11}220 = 60$. Furthermore, $P(A \mid E) = \frac{P(AE)}{P(E)}$. Thus, using (a), $\frac{1}{5} = \frac{P(AE)}{\frac{3}{11}}$, that is, P(AE) = $\frac{3}{55}$. Thus the number of specimens marked AE is $\frac{3}{55}220 =$ 12. The calculation for P(B|E) is identical; thus the number of specimens marked BE is also 12. Thus, so far, we have:

| LABEL | number |
|---|---|
| AE | 12 |
| AF | |
| AG | |
| BE | 12 |
| BF | |
| BG | |
| CE | 36 |
| CF | |
| CG | |

**(c)** $P(G\,|\,C) = \frac{12}{23}$. Thus $P(G\,|\,C) = \dfrac{P(CG)}{P(C)}$, that is, $\frac{12}{23} = \dfrac{P(CG)}{P(C)}$.

**(d)** P(C) = $\frac{23}{55}$. Thus, using (c) we get $\frac{12}{23} = \dfrac{P(CG)}{\frac{23}{55}}$, that is, P(CG) = $\frac{12}{55}$. Hence the number of specimens marked CG is $\frac{12}{55}220 = 48$. Since P(C) = $\frac{23}{55}$, the total number of specimens marked with a C is $\frac{23}{55}220$ = 92. Since 36 were marked CE and 48 were marked CG, it follows that the number of specimens marked CF is 92 –48 –36 = 8. Thus, so far, we have:

| LABEL | number |
|---|---|
| AE | 12 |
| AF | |
| AG | |
| BE | 12 |
| BF | |
| BG | |
| CE | 36 |
| CF | 8 |
| CG | 48 |

**(e)** The number of BGs is twice the number of BEs. Since the latter is 12, the number of BGs is 24.

**(f)** P(E|A) = $\frac{3}{20}$. Since $P(E \mid A) = \frac{P(AE)}{P(A)}$ and, from (b), P(AE) = $\frac{3}{55}$, we have that $\frac{3}{20} = \frac{\frac{3}{55}}{P(A)}$. Hence $P(A) = \frac{3}{55}\,\frac{20}{3} = \frac{4}{11}$. Thus the number of specimens marked with an A is $\frac{4}{11}220 = 80$. Since P(A) = $\frac{4}{11}$ and (from (d)) P(C) = $\frac{23}{55}$, it follows that P(B) = 1 $-\frac{4}{11}$ $-\frac{23}{55}$ = $\frac{12}{55}$. Thus the number of specimens marked with a B is $\frac{12}{55}$220 = 48. Of these, 12 were marked BE and 24 were marked BG. Thus the number of specimens marked BF is 48 −12 −24 = 12. Thus, so far, we have:

| LABEL | number |
|---|---|
| AE | 12 |
| AF | |
| AG | |
| BE | 12 |
| BF | 12 |
| BG | 24 |
| CE | 36 |
| CF | 8 |
| CG | 48 |

**(g)** P(G|A) = $\frac{1}{10}$. $P(G \mid A) = \frac{P(AG)}{P(A)}$. From (f) we have that P(A) = $\frac{4}{11}$. Thus $\frac{1}{10} = \frac{P(AG)}{\frac{4}{11}}$. Thus P(AG) = $\frac{1}{10}\left(\frac{4}{11}\right) = \frac{4}{110}$. Thus the number of specimens marked AG is $\frac{4}{110}220 = 8$. Since the number marked with an A is $\frac{4}{11}220 = 80$ and those marked AE are 12 and those marked AG are 8, we get that the number of specimens marked AF is 80 −12 −8 = 60. Thus we have completed the table:

| LABEL | number |
|---|---|
| AE | 12 |
| AF | 60 |
| AG | 8 |
| BE | 12 |
| BF | 12 |
| BG | 24 |
| CE | 36 |
| CF | 8 |
| CG | 48 |

**Exercise 8.6**. Let $f : \mathcal{E} \to 2^W$ be an AGM belief revision function. Let $E, F \in \mathcal{E}$ be such that $E \subseteq F$ and $E \cap f(F) \neq \varnothing$. We need to show that $f(E) = E \cap f(F)$. By definition of AGM belief revision function,

$$f(F) = \{w \in F : w \precsim w' \text{ for every } w' \in F\} \tag{1}$$

and

$$f(E) = \{w \in E : w \precsim w' \text{ for every } w' \in E\} \tag{2}$$

Choose an arbitrary $w \in E \cap f(F)$. Then, by (1) and the fact that $E \subseteq F$, $w \precsim w'$ for every $w' \in E$ and thus, by (2), $w \in f(E)$. Hence $E \cap f(F) \subseteq f(E)$. Conversely, choose an arbitrary $w_1 \in f(E)$. Then, since $f(E) \subseteq E$, $w_1 \in E$. We want to show that $w_1 \in f(F)$ [so that $w_1 \in E \cap f(F)$ and, therefore, $f(E) \subseteq E \cap f(F)$]. Suppose not. Then, by (1), there exists a $w_2 \in F$ such that $w_2 \prec w_1$. Select a $w_3 \in E \cap f(F)$ (recall that, by hypothesis, $E \cap f(F) \neq \varnothing$). Then, by (1) (since $w_2, w_3 \in f(F)$), $w_3 \precsim w_2$, from which it follows (by transitivity of $\precsim$ and the fact that $w_2 \prec w_1$) that $w_3 \prec w_1$. But then, since $w_3 \in E$, it is not true that $w_1 \precsim w'$ for every $w' \in E$, contradicting - by (2) - the hypothesis that $w_1 \in f(E)$.

**Exercise 8.7**. We have that $\mathcal{E} = \{\{a,e\},\{d,e,k,m\},\{b,d,e,k\},W\}$ and the plausibility order is

| | |
|---|---|
| most plausible | $b,g$ |
| | $c,k,m$ |
| | $d,h$ |
| | $e$ |
| least plausible | $a$ |

The belief revision function rationalized by this plausibility order is given by $f(\{a,e\}) = \{e\}$, $f(\{d,e,k,m\}) = \{k,m\}$, $f(\{b,d,e,k\}) = \{b\}$ and $f(W) = \{b,g\}$.

**Exercise 8.8**. From Exercise 8.7 we get that $\{P_E, P_F, P_G, P_W\}$ must be such that $Supp(P_E) = \{e\}$, $Supp(P_F) = \{k,m\}$, $Supp(P_G) = \{b\}$ and $Supp(P_W) = \{b,g\}$. For every full support probability distribution $P_0$ on *W*, there is a corresponding collection $\{P_E, P_F, P_G, P_W\}$. For example, if $P_0$ is the uniform distribution on *W* (that assigns probability $\frac{1}{9}$ to every state) then the corresponding $\{P_E, P_F, P_G, P_W\}$ is given by:

| state | $a$ | $b$ | $c$ | $d$ | $e$ | $g$ | $h$ | $k$ | $m$ |
|---|---|---|---|---|---|---|---|---|---|
| $P_E$ | 1 | 0 | 0 | 0 | 0 | 0 | 0 | 0 | 0 |
| $P_F$ | 0 | 0 | 0 | 0 | 0 | 0 | 0 | $\frac{1}{2}$ | $\frac{1}{2}$ |
| $P_G$ | 0 | 1 | 0 | 0 | 0 | 0 | 0 | 0 | 0 |
| $P_W$ | 0 | $\frac{1}{2}$ | 0 | 0 | 0 | $\frac{1}{2}$ | 0 | 0 | 0 |

If $P_0$ is the following probability distribution

| state | $a$ | $b$ | $c$ | $d$ | $e$ | $g$ | $h$ | $k$ | $m$ |
|---|---|---|---|---|---|---|---|---|---|
| $P_0$ | $\frac{1}{50}$ | $\frac{3}{50}$ | $\frac{11}{50}$ | $\frac{4}{50}$ | $\frac{8}{50}$ | $\frac{9}{50}$ | $\frac{5}{50}$ | $\frac{2}{50}$ | $\frac{7}{50}$ |

then the corresponding $\{P_E, P_F, P_G, P_W\}$ is given by: $P_E$ and $P_G$ the same as above, and $P_F$ and $P_W$ as follows:

| state | $a$ | $b$ | $c$ | $d$ | $e$ | $g$ | $h$ | $k$ | $m$ |
|---|---|---|---|---|---|---|---|---|---|
| $P_F$ | 0 | 0 | 0 | 0 | 0 | 0 | 0 | $\frac{2}{9}$ | $\frac{7}{9}$ |
| $P_W$ | 0 | $\frac{1}{4}$ | 0 | 0 | 0 | $\frac{3}{4}$ | 0 | 0 | 0 |

**Exercise 8.9**. Yes. The following is a common prior: $\begin{pmatrix} a & b & c & d \\ \frac{1}{9} & \frac{2}{9} & \frac{2}{9} & \frac{4}{9} \end{pmatrix}$.

**Exercise 8.10**. **(a)** Assume that $p = \frac{1}{13}$ (so that $1-p = \frac{12}{13}$ ). From the updating conditions for Individual 1 we get $p_b = 2p_a$ and $p_e = p_f$; from the updating conditions for Individual 2 we get $p_c = 3p_b$ and $p_e = 4p_d$; from the updating conditions for Individual 3 we get $p_c = p_d$ and $p_f = 12p_b$. All these equations have a solution, which constitutes a common prior, namely $\begin{pmatrix} a & b & c & d & e & f \\ \frac{1}{63} & \frac{2}{63} & \frac{6}{63} & \frac{6}{63} & \frac{24}{63} & \frac{24}{63} \end{pmatrix}$.

**(b)** From the equations of part (a) we have that $p_f = p_e$, $p_e = 4p_d$, $p_d = p_c$ and $p_c = 3p_b$, from which it follows that $p_f = 12p_b$. Thus, updating on the information set {*b*,*f*} of Individual 3 we get that we must have $(1-p) = \dfrac{p_f}{p_b + p_f} = \dfrac{12p_b}{p_b + 12p_b} = \frac{12}{13}$ and thus $p = \frac{1}{13}$. Hence, if $p \neq \frac{1}{13}$ there is no common prior.

**Exercise 8.11**. **(a)** and **(b)** The common knowledge partition and the common prior are shown below.

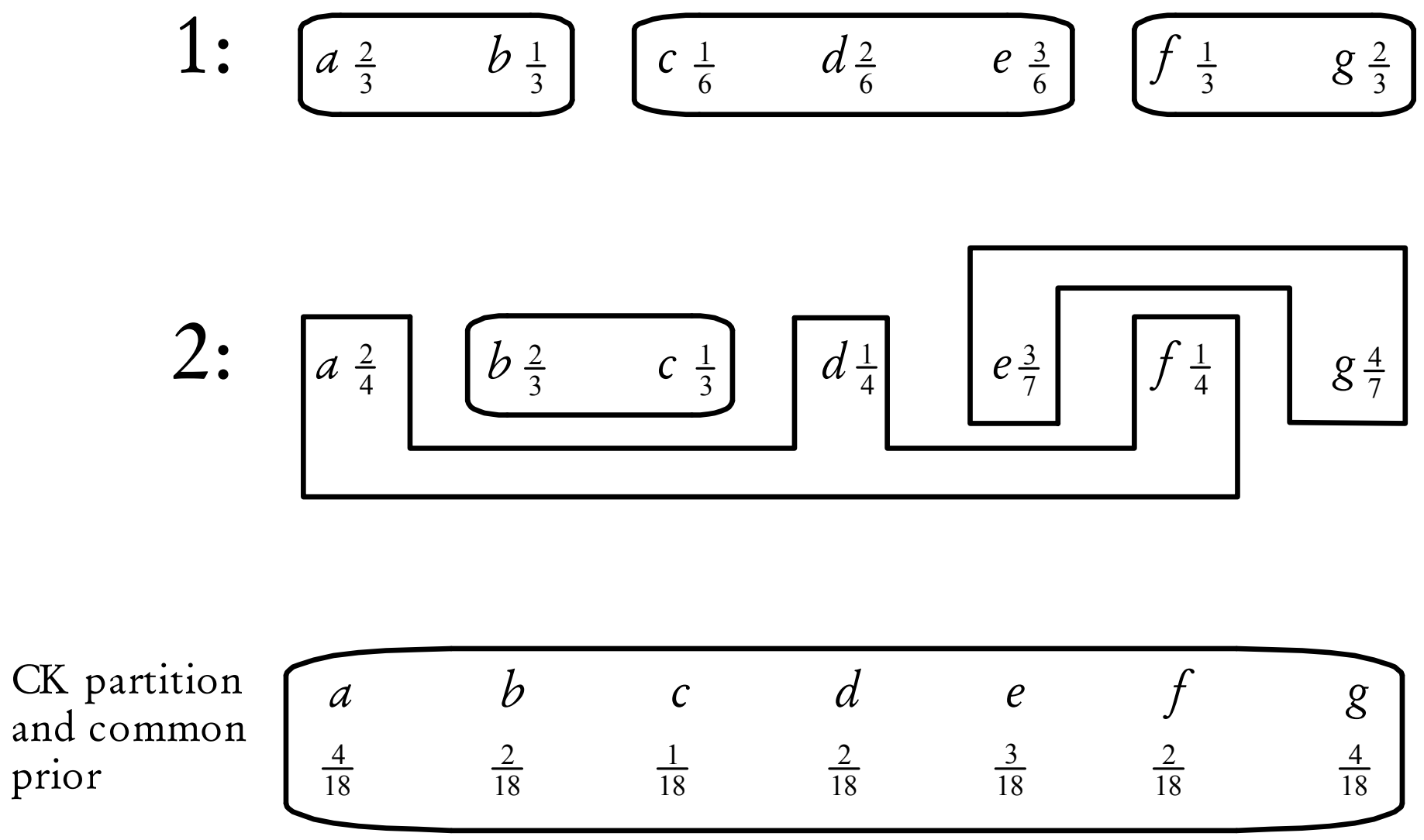

**(c.1)** Individual 1 assigns probability $\frac{1}{3}$ to $E$ at every state. Thus his assessment of the probability of $E$ is common knowledge at every state.

**(c.2)** At state $a$ Individual 2 assigns probability $\frac{1}{2}$ to $E$, while at state $b$ she assigns probability $\frac{2}{3}$ to $E$. Thus her assessment of the probability of $E$ is not common knowledge at any state.

**(d.1)** At state $a$ Individual 1's information set is $\{a,b\}$; since at $a$ Individual 2 assigns probability $\frac{1}{2}$ to $E$, while at $b$ she assigns probability $\frac{2}{3}$ to $E$, it follows that at $a$ it is not the case that Individual 1 knows Individual 2's assessment of the probability of $E$.

**(d.2)** At state $c$ Individual 1's information set is $\{c,d,e\}$; since at $c$ Individual 2 assigns probability $\frac{2}{3}$ to $E$, while at $d$ she assigns probability $\frac{1}{2}$ to $E$, it follows that at $c$ it is not the case that Individual 1 knows Individual 2's assessment of the probability of $E$.

**(d.3)** At state $f$ Individual 1's information set is $\{f,g\}$; since at $f$ Individual 2 assigns probability $\frac{1}{2}$ to $E$, while at $g$ she assigns probability 0 to $E$, it follows that at $f$ it is not the case that Individual 1 knows Individual 2's assessment of the probability of $E$.

**Exercise 8.12** (**Challenging question**).

**(a)** The extensive form is as follows ('cb$n$' means 'the car is behind door $n$', 'ch$n$' means 'you choose door $n$').

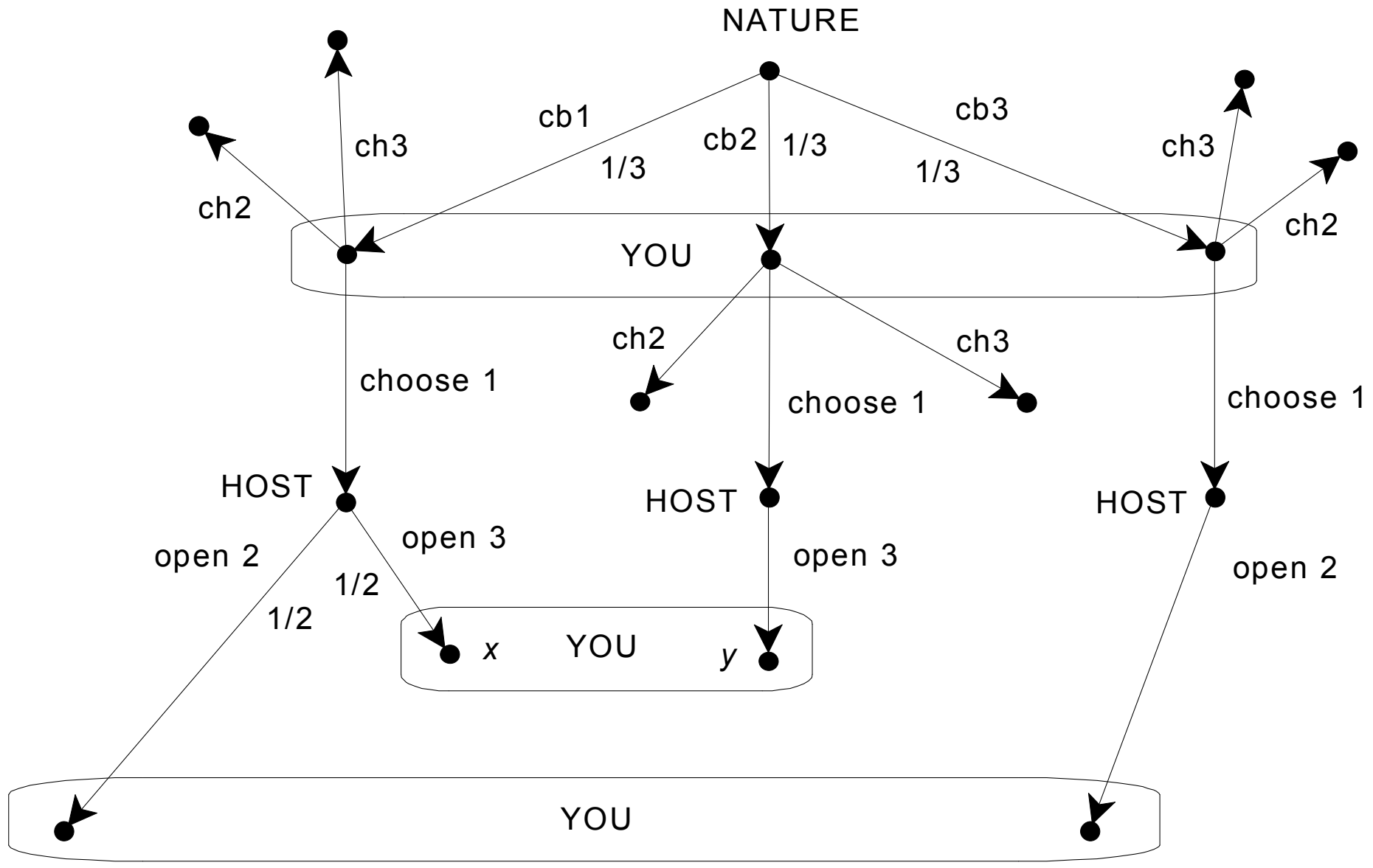

The hypothesized sequence of events leads to either node $x$ or node $y$ (you first choose door 1 and then the host opens door 3). The prior probability of getting to $x$ is $\frac{1}{3}\times\frac{1}{2}=\frac{1}{6}$ while the prior probability of getting to $y$ is $\frac{1}{3}\times 1=\frac{1}{3}$. The information you have is $\{x,y\}$ and, using the conditional probability rule, $P(\{x\}\,|\,\{x,y\})=\dfrac{P(\{x\}\cap\{x,y\})}{P(\{x,y\})}=\dfrac{P(\{x\})}{P(\{x,y\})}=\dfrac{\frac{1}{6}}{\frac{1}{6}+\frac{1}{3}}=\frac{1}{3}$. Thus you should switch.

**(b)** By Bayes' rule, $P(D_1\,|\,O_3)=\dfrac{P(O_3\,|\,D_1)\;P(D_1)}{P(O_3)}$. We know that $P(D_1)=\frac{1}{3}$ and $P(O_3\,|\,D_1)=\frac{1}{2}$ (when the car is behind door 1 then the host has a choice between opening door 2 and opening door 3 and he chooses with equal probability). Thus

$$P(D_1\,|\,O_3)=\frac{\frac{1}{2}\times\frac{1}{3}}{P(O_3)}=\frac{\frac{1}{6}}{P(O_3)}. \qquad (\blacklozenge)$$

Thus we only need to compute $P(O_3)$. By Bayes' rule,

$$P(O_3)=P(O_3\,|\,D_1)\,P(D_1)+P(O_3\,|\,D_2)\,P(D_2)+P(O_3\,|\,D_3)\,P(D_3)=$$

$$P(O_3\,|\,D_1)\tfrac{1}{3}+P(O_3\,|\,D_2)\tfrac{1}{3}+P(O_3\,|\,D_3)\tfrac{1}{3}=\tfrac{1}{2}\times\tfrac{1}{3}+1\times\tfrac{1}{3}+0\times\tfrac{1}{3}=\tfrac{1}{6}+\tfrac{1}{3}=\tfrac{1}{2}$$

because $P(O_3\,|\,D_1)=\frac{1}{2}$, $P(O_3\,|\,D_2)=1$ and $P(O_3\,|\,D_3)=0$. Substituting $\frac{1}{2}$ for $P(O_3)$ in ($\blacklozenge$) we get that $P(D_1\,|\,O_3)=\frac{1}{3}$. Hence the updated probability that the car is behind the other door (door 2) is $\frac{2}{3}$ and therefore you should switch.

Chapter 9

# Common Knowledge of Rationality

## 9.1 Models of strategic-form games

In Chapter 1 (Section 1.5) and Chapter 5 (Section 5.4) we discussed the procedure of iterative elimination of strictly dominated strategies and claimed that it captures the notion of common knowledge of rationality in strategic-form games. We now have the tools to state this precisely.

The *epistemic foundation program* in game theory aims to identify, for every game, the strategies that might be chosen by rational and intelligent players who know the structure of the game and the preferences of their opponents and who recognize each other's rationality. The two central questions are thus: (1) Under what circumstances is a player rational? and (2) What does 'mutual recognition of rationality' mean? A natural interpretation of the latter notion is in terms of 'common knowledge of rationality'.[17] Thus we only need to define what it means for a player to be rational. Intuitively, a player is rational if she chooses an action which is "best" given what she believes. In order to make this more precise we need to introduce the notion of a model of a game. We shall focus on strategic-form games with cardinal payoffs (Definition 5.2, Chapter 5).

[17] As pointed out in Chapter 7 (Remark 7.1), a defining characteristic of knowledge is truth: if, at a state, a player knows event $E$ then, at that state, $E$ must be true. A more general notion is that of *belief*, which allows for the possibility of mistakes: believing $E$ is compatible with $E$ being false. Thus a more appealing notion is that of *common belief* of rationality; however, in order to simplify the exposition, we shall restrict attention to the notions of knowledge and common knowledge developed in Chapter 7. For the analysis of common belief and an overview of the epistemic foundation program in game theory the reader is referred to Battigalli and Bonanno (1999).

The definition of a strategic-form game specifies the choices available to the players and what motivates those choices (their preferences over the possible outcomes); however, it leaves out an important factor in the determination of players' choices, namely what they believe about the other players. Adding a specification of the players' knowledge and beliefs determines the context in which a particular game is played; this can be done with the help of an interactive knowledge-belief structure.

Recall that an interactive knowledge-belief structure consists of a set of states *W*, *n* partitions $\mathcal{I}_1, \mathcal{I}_2, ..., \mathcal{I}_n$ of *W* and a collection of probability distributions on *W*, one for each information set of each partition, whose support is a subset of that information set. The probability distributions encode the beliefs of the players in each possible state of knowledge. Recall also that if $w$ is a state we denote by $I_i(w)$ the element of the partition $\mathcal{I}_i$ that contains *w*, that is, the information set of Player *i* at state $w$ and by $P_{i,w}$ the probability distribution on $I_i(w)$, representing the beliefs of Player *i* at state *w*.

**Definition 9.1.** Let *G* be a *n*-player strategic-form game with cardinal payoffs. A *model of G* is an interactive knowledge-belief structure together with a function $\sigma_i : W \to S_i$ (for every $i \in \{1,\ldots,n\}$) that associates with every state a pure strategy of Player *i* (recall that $S_i$ denotes the set of pure strategies of Player *i*).[18] The interpretation of $s_i = \sigma_i(w)$ is that, at state *w*, Player *i* plays (or has chosen) strategy $s_i$. We impose the restriction that a player always knows what strategy she has chosen, that is, the function $\sigma_i$ is constant on each information set of Player *i*: if $w' \in I_i(w)$ then $\sigma_i(w') = \sigma_i(w)$.

The addition of the functions $\sigma_i$ to an interactive knowledge-belief structure yields an interpretation of events in terms of propositions about what strategies the players choose, thereby giving content to players' knowledge and beliefs. Figure 9.1 shows a strategic-form game and a model of it.

[18] In Chapter 5 we used the symbol $\sigma_i$ to denote a mixed strategy of Player *i*. Mixed strategies will play no role in this chapter and thus there will be no possibility of confusion.

Player 2

| Player 1 | | $L$ | $C$ | $R$ |
|---|---|---|---|---|
| | $T$ | 4 , 6 | 3 , 2 | 8 , 0 |
| | $M$ | 0 , 9 | 0 , 0 | 4 , 12 |
| | $B$ | 8 , 3 | 2 , 4 | 0 , 0 |

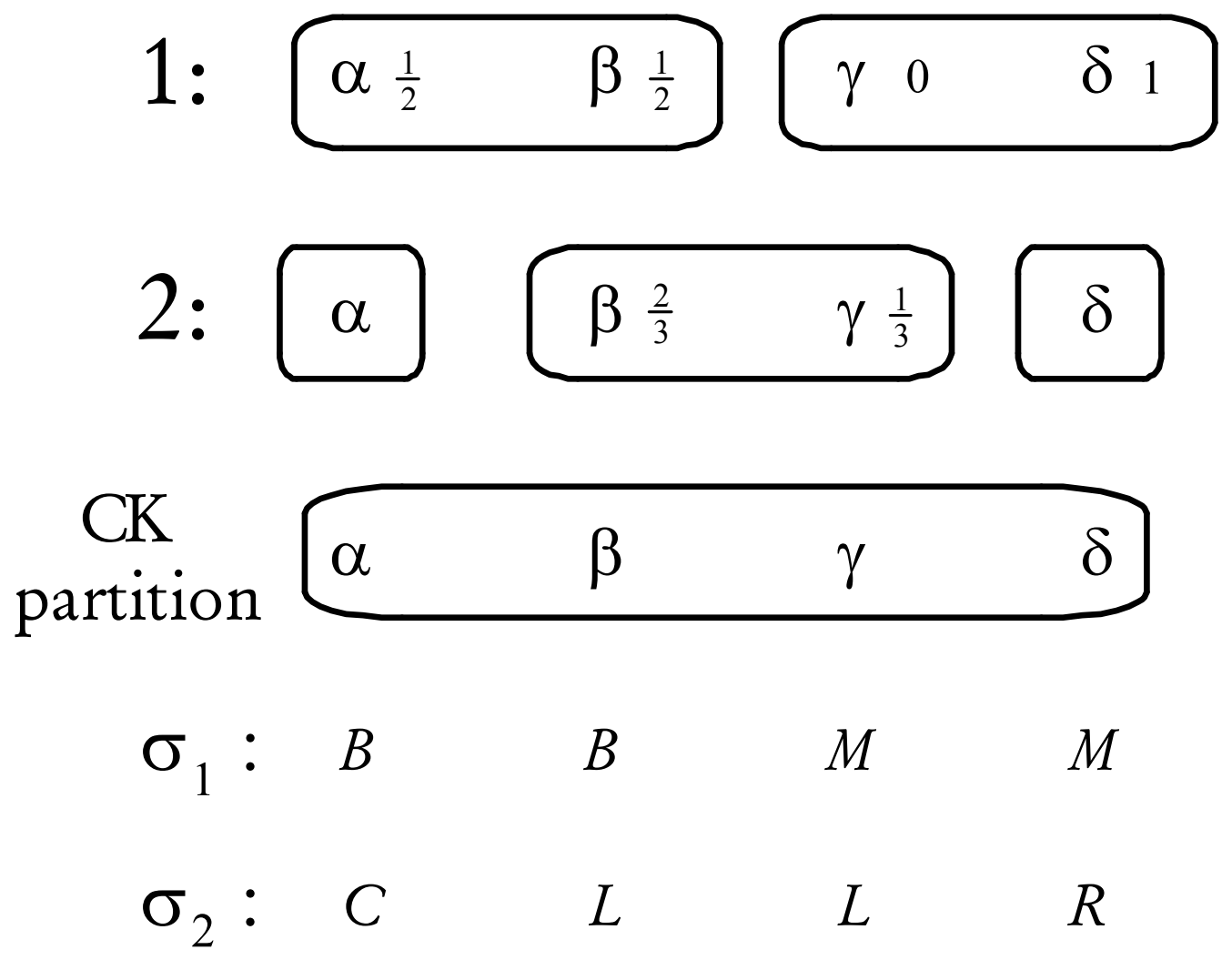

**Figure 9.1**

State $\beta$ in Figure 9.1 describes the following situation.

- Player 1 plays $B$ $(\sigma_1(\beta) = B)$ and Player 2 plays $L$ $(\sigma_2(\beta) = L)$.
- Player 1 (whose information set is $\{\alpha,\beta\}$) is uncertain as to whether Player 2 has chosen to play $C$ $(\sigma_2(\alpha) = C)$ or $L$ $(\sigma_2(\beta) = L)$; furthermore, he attaches probability $\frac{1}{2}$ to each of these two possibilities.
- Player 2 (whose information set is $\{\beta,\gamma\}$) is uncertain as to whether Player 1 has chosen to play $B$ $(\sigma_1(\beta) = B)$ or $M$ $(\sigma_1(\gamma) = M)$; furthermore, she attaches probability $\frac{2}{3}$ to Player 1 playing $B$ and $\frac{1}{3}$ to Player 1 playing $M$.

Are the two players' choices rational at state $\beta$?

Consider first Player 1. Given his beliefs and his choice of $B$, his expected payoff is $\frac{1}{2}\,\pi_1(B,C) + \frac{1}{2}\,\pi_1(B,L) = \frac{1}{2}(2) + \frac{1}{2}(8) = 5$. Given his beliefs about Player 2, could he obtain a higher expected payoff with a different choice? With $M$ his expected payoff would be $\frac{1}{2}\,\pi_1(M,C)+\frac{1}{2}\,\pi_1(M,L)=\frac{1}{2}(0)+\frac{1}{2}(0)=0$, while with $T$ his expected payoff would be $\frac{1}{2}\,\pi_1(T,C)+\frac{1}{2}\,\pi_1(T,L)=\frac{1}{2}(3)+\frac{1}{2}(4)=3.5$. Thus, given his beliefs, Player 1's choice of $B$ is optimal and we can say that Player 1 is rational at state $\beta$.

Consider now Player 2. Given her beliefs and her choice of $L$, her expected payoff is $\frac{2}{3}\,\pi_2(B,L)+\frac{1}{3}\,\pi_2(M,L)=\frac{2}{3}(3)+\frac{1}{3}(9)=5$. Given her beliefs about Player 1, could she obtain a higher expected payoff with a different choice? With $C$ her expected payoff would be $\frac{2}{3}\,\pi_2(B,C)+\frac{1}{3}\,\pi_2(M,C)=\frac{2}{3}(4)+\frac{1}{3}(0)=\frac{8}{3}$, while with $R$ her expected payoff would be $\frac{2}{3}\,\pi_2(B,R)+\frac{1}{3}\,\pi_2(M,R)=\frac{2}{3}(0)+\frac{1}{3}(12)=4$. Thus, given her beliefs, Player 2's choice of $L$ is optimal and we can say that also Player 2 is rational at state $\beta$.

One can view a strategic-form game as only a partial description of an interactive situation: it specifies who the players are, what actions they can take and how they rank the possible outcomes. A model of the game completes this description by also specifying what each player actually does and what she believes about what the other players are going to do. Once we know what a player does and what she believes then we are in a position to judge her choice to be either rational or irrational.

The following definition says that at a state (of a model) a player is rational if her choice at that state is optimal, given her beliefs, that is, if there is no other pure strategy that would give her a higher expected payoff, given what she believes about the choices of her opponents.

Given a state $w$ and a Player $i$, we denote by $\sigma_{-i}(w)$ the profile of strategies chosen at $w$ by the players other than $i$:

$$\sigma_{-i}(w)=\left(\sigma_1(w),\ldots,\sigma_{i-1}(w),\sigma_{i+1}(w),\ldots,\sigma_n(w)\right).$$

Furthermore, as in Chapter 8, we denote by $P_{i,w}$ the beliefs of individual $i$ at state $w$, that is, $P_{i,w}$ is a probability distribution over $I_i(w)$ (although we can think of $P_{i,w}$ as a probability distribution over the entire set of states $W$ satisfying the property that if $w'\notin I_i(w)$ then $P_{i,w}(w')=0$). Note that, for every Player $i$, there is just one probability distribution over $I_i(w)$, that is, if $w'\in I_i(w)$ then $P_{i,w'}=P_{i,w}$.

**Definition 9.2.** Player *i* is *rational* at state *w* if, given her beliefs $P_{i,w}$, there is no pure strategy $s_i' \in S_i$ of hers that yields a higher expected payoff than $\sigma_i(w)$ (recall that $\sigma_i(w)$ is the strategy chosen by Player *i* at state *w*):

$$\sum_{w' \in I_i(w)} P_{i,w}(w') \times \pi_i\left(\sigma_i(w), \sigma_{-i}(w')\right) \geq \sum_{w' \in I_i(w)} P_{i,w}(w') \times \pi_i\left(s_i', \sigma_{-i}(w')\right),$$
$$\text{for all } s_i' \in S_i.$$

Given a model of a game, using Definition 9.2 we can determine, for every player, the set of states where that player is rational. Let $\boldsymbol{R}_i$ be the event that (that is, the set of states at which) Player *i* is rational and let $\boldsymbol{R} = \boldsymbol{R}_1 \cap \ldots \cap \boldsymbol{R}_n$ be the event that all players are rational.

In the example of Figure 9.1 we have that $\boldsymbol{R}_1 = \{\alpha, \beta\}$. Indeed, it was shown above that $\beta \in \boldsymbol{R}_1$; since Player 1's choice and beliefs at $\alpha$ are the same as at $\beta$, it follows that also $\alpha \in \boldsymbol{R}_1$. On the other hand, Player 1 is not rational at state $\gamma$ because he is certain (that is, attaches probability 1 to the fact) that Player 2 plays *R* and his choice of *M* is not optimal against *R* (the unique best response to *R* is *T*). The same is true of state $\delta$. For Player 2 we have that $\boldsymbol{R}_2 = \{\alpha, \beta, \gamma, \delta\}$, that is, Player 2 is rational at every state (the reader should convince himself/herself of this). Thus $\boldsymbol{R} = \{\alpha, \beta\}$.

This is a good time to test your understanding of the concepts introduced in this section, by going through the exercises in Section 9.E.1 of Appendix 9.E at the end of this chapter.

## 9.2 Common knowledge of rationality in strategic-form games

Since we have defined the rationality of a player and the rationality of all the players as events, we can apply the knowledge operator and the common knowledge operator to these events; thus we can determine, for example, if there are any states where the rationality of all the players is common knowledge and see what strategy profiles are compatible with common knowledge of rationality.

**REMARK 1.** Note that, for every player i, $K_i\boldsymbol{R}_i = \boldsymbol{R}_i$, that is, every player is rational if and only if she knows it (this is a consequence of the assumption that a player always knows what strategy she has chosen and knows her own beliefs: if $w' \in I_i(w)$ then $\sigma_i(w') = \sigma_i(w)$ and $P_{i,w'} = P_{i,w}$). Note also that $CK\boldsymbol{R} \subseteq \boldsymbol{R}$, that is, if it is common knowledge that all the players are rational, then they are indeed rational.[19]

In Chapters 1 (for ordinal strategic-form games) and 5 (for cardinal strategic-form games) we claimed that the iterated deletion of strictly dominated strategies corresponds to the notion of common knowledge of rationality. We are now able to state this precisely. The following two theorems establish the correspondence between the two notions for the case of strategic-form games with cardinal payoffs. A corresponding characterization holds for strategic-form games with ordinal payoffs.[20] The proofs of Theorems 9.1 and 9.2 are given in Appendix 9.A.

**Theorem 9.1.** Given a finite strategic-form game with cardinal payoffs *G*, the following is true: for any model of *G* and any state *w* in that model, if $w \in CK\boldsymbol{R}$ (that is, at *w* there is common knowledge of rationality) then the pure-strategy profile associated with *w* must be one that survives the iterated deletion of strictly dominated strategies (Definition 5.6, Chapter 5).

**Theorem 9.2.** Given a finite strategic-form game with cardinal payoffs *G*, the following is true: if $s \in S$ is a pure-strategy profile that survives the iterated deletion of strictly dominated strategies, then there is a model of *G* and a state *w* in that model, such that $w \in CK\boldsymbol{R}$ and the strategy profile associated with *w* is *s*.

As an application of Theorem 9.1 consider again the game of Figure 9.1, reproduced below in Figure 9.2.

---

[19] Indeed, the common knowledge operator, like the individual knowledge operators, satisfies the Truth property (see Remark 7.1 and Exercise 7.5, Chapter 7). Note that, although $CK\boldsymbol{R} \subseteq \boldsymbol{R}$, the converse (that is, the inclusion $\boldsymbol{R} \subseteq CK\boldsymbol{R}$) is *not* true: it is possible for all the players to be rational at a state without this fact being common knowledge among them at that state. For example, this is the case at state $\alpha$ in the model of Figure 9.1.

[20] See Bonanno (2015).

| | | Player 2 | | |
|---|---|---|---|---|
| | | *L* | *C* | *R* |
| Player 1 | *T* | 4 , 6 | 3 , 2 | 8 , 0 |
| | *M* | 0 , 9 | 0 , 0 | 4 , 12 |
| | *B* | 8 , 3 | 2 , 4 | 0 , 0 |

**Figure 9.2**

In this game the iterated deletion of strictly dominated strategies leads to the following strategy profiles: (*T*, *L*), (*T*, *C*), (*B*, *L*) and (*B*, *C*).[21] Thus these are the only strategy profiles that are compatible with common knowledge of rationality. Hence, by Theorem 9.1, at a state in a model of this game where there is common knowledge of rationality the players can play only one of these strategy profiles. Furthermore, by Theorem 9.2, any of these four strategy profiles can in fact be played in a situation where there is common knowledge of rationality. It is worth stressing that common knowledge of rationality does *not* imply that the players play a Nash equilibrium: indeed, none of the above four strategy profiles is a Nash equilibrium. However, since a pure-strategy Nash equilibrium always survives the iterated deletion of strictly dominated strategies, the pure-strategy Nash equilibria (if any) are always compatible with common knowledge of rationality.

This is a good time to test your understanding of the concepts introduced in this section, by going through the exercises in Section 9.E.2 of Appendix 9.E at the end of this chapter.

[21] For Player 1, *M* is strictly dominated by any mixed strategy $\begin{pmatrix} T & B \\ p & 1-p \end{pmatrix}$ with $p > \frac{1}{2}$. After deletion of *M*, for Player 2 *R* becomes strictly dominated by either of the other two strategies.

# 9.3 Common knowledge of rationality in extensive-form games

So far we have only discussed the implications of common knowledge of rationality in strategic-form games, where the players make their choices simultaneously (or in ignorance of the other players' choices). We now turn to a brief discussion of the issues that arise when one attempts to determine the implications of common knowledge of rationality in dynamic (or extensive-form) games.

How should a model of a dynamic game be constructed? One approach in the literature has been to consider models of the corresponding strategic-form. However, there are several conceptual issues that arise in this context. In the models considered in Section 9.1 the interpretation of $s_i = \sigma_i(w)$ is that at state $w$ player $i$ "plays" or "chooses" strategy $s_i$. Now consider the perfect-information game shown in Figure 9.3 below and a model of the strategic-form game associated with it. Let $w$ be a state in that model where $\sigma_1(w) = (d_1, a_3)$. What does it mean to say that Player 1 "chooses" strategy $(d_1, a_3)$? The first part of the strategy, namely $d_1$, can be interpreted as a description of Player 1's actual behavior (what he actually does: he plays $d_1$), but the second part of the strategy, namely $a_3$, has no such interpretation: if Player 1 in fact plays $d_1$ then he knows that he will not have to make any further choices and thus it is not clear what it means for him to "choose" to play $a_3$ in a situation that is made impossible by his decision to play $d_1$.

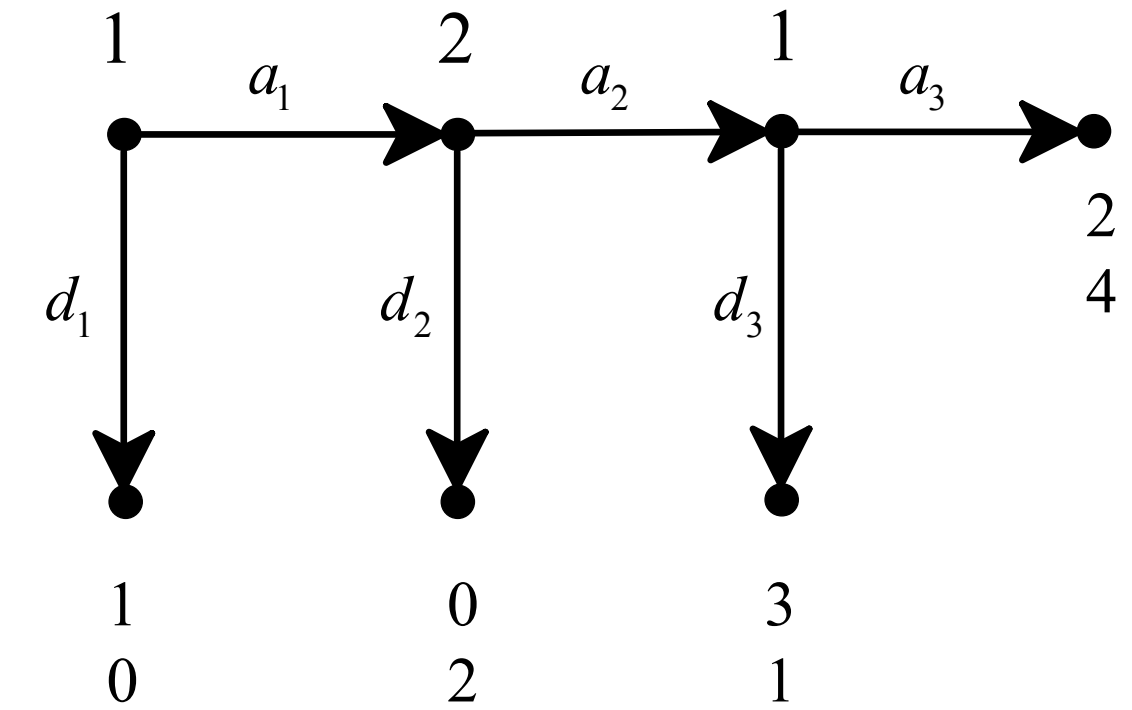

**Figure 9.3**

Thus it does not seem to make sense to interpret $\sigma_1(w) = (d_1, a_3)$ as 'at state *w* Player 1 chooses $(d_1, a_3)$'. Perhaps the correct interpretation is in terms of a more complex sentence such as 'Player 1 chooses to play $d_1$ and if - contrary to this plan - he were to play $a_1$ and Player 2 were to follow with $a_2$, then Player 1 would play $a_3$'. Thus, while in a simultaneous game the association of a strategy of Player *i* to a state can be interpreted as a description of Player *i*'s actual behavior at that state, in the case of dynamic games this interpretation is no longer valid, since one would end up describing not only the actual behavior of Player *i* at that state but also his counterfactual behavior. Methodologically, this is not satisfactory: if it is considered to be necessary to specify what a player would do in situations that do not occur at the state under consideration, then one should model the counterfactual explicitly. But why should it be necessary to specify at state *w* (where Player 1 is playing $d_1$) what he would do at the counterfactual node following actions $a_1$ and $a_2$? Perhaps what matters is not so much what Player 1 would actually do in that situation but what Player 2 *believes* that Player 1 would do: after all, Player 2 might not know that Player 1 has decided to play $d_1$ and she needs to consider what to do in the eventuality that Player 1 actually ends up playing $a_1$. So, perhaps, the strategy of Player 1 is to be interpreted as having two components: (1) a description of Player 1's behavior and (2) a conjecture in the mind of Player 2 about what Player 1 would do. If this is the correct interpretation, then one could object – again from a methodological point of view - that it would be preferable to disentangle the two components and model them explicitly.[22]

An alternative - although less common - approach in the literature dispenses with strategies and considers models of games where (1) states are described in terms of players' *actual behavior* and (2) players' conjectures concerning the actions of their opponents in various hypothetical situations are modeled by means of a generalization of the knowledge-belief structures considered in Section 9.1. The generalization is obtained by encoding not only the initial beliefs of the players (at each state) but also their *dispositions to revise those beliefs* under various hypotheses.[23]

[22] For an extensive discussion of these issues see Bonanno (forthcoming).

[23] The interested reader is referred to Perea (2012).

A third approach has been to move away from static models (like the models considered in Section 9.1) and consider dynamic models where time is introduced explicitly into the analysis. These are behavioral models where strategies play no role and the only beliefs that are specified are the actual beliefs of the players at the time of choice. Thus players' beliefs are modeled as temporal, rather than conditional, beliefs and rationality is defined in terms of actual choices, rather than hypothetical plans.[24]

A discussion of the implications of common knowledge (or common belief) of rationality in dynamic games would require the introduction of several new concepts and definitions. For the sake of brevity, we shall not pursue this topic in this book.

[24] The interested reader is referred to Bonanno (2014).

# Appendix 9.A: Proofs

In order to prove Theorem 9.1 we will need the following extension of Theorem 5.3 (Chapter 5) to the case of an arbitrary number of players. A proof of Theorem 9.3 can be found in Osborne and Rubinstein (1994, Lemma 60.1, p. 60).

**Theorem 9.3.** Consider a finite *n*-player game in strategic form with cardinal payoffs $\langle I,(S_i)_{i\in I},(\pi_i)_{i\in I}\rangle$. Select an arbitrary Player *i* and a pure strategy $s_i$ of Player *i*. For every Player *j* ≠ *i*, let $S'_j \subseteq S_j$ be a subset of *j*'s set of pure strategies and let $S'_{-i} = S'_1 \times ... \times S'_{i-1} \times S'_{i+1} \times ... \times S'_n$ be the Cartesian product of the sets $S'_j$ ($j \neq i$). Then the following are equivalent:

(1) There is a belief of Player *i* on $S'_{-i}$ (that is, a probability distribution $P: S'_{-i} \to [0,1]$) that makes $s_i$ a best response (that is, for every $x \in S_i$, $\sum_{s'_{-i}\in S'_{-i}} \pi_i(s_i, s'_{-i})\, P(s'_{-i}) \geq \sum_{s'_{-i}\in S'_{-i}} \pi_i(x, s'_{-i})\, P(s'_{-i})$),

(2) $s_i$ is not strictly dominated by a mixed strategy of Player *i* in the restriction of the game to the sets of pure strategies $S'_1,...,S'_{i-1},S_i,S'_{i+1},S'_n$.

Given a finite strategic-form game with cardinal payoffs $G = \langle I,(S_i)_{i\in I},(\pi_i)_{i\in I}\rangle$ we shall denote by $S_i^\infty$ the set of pure strategies of Player *i* that survive the iterated deletion of strictly dominated strategies (Definition 5.6, Chapter 5) and by $S^\infty$ the corresponding set of strategy profiles (that is, $S^\infty = S_1^\infty \times ... \times S_n^\infty$).

**Proof of Theorem 9.1**. Consider a finite strategic-form game *G* with cardinal payoffs and an epistemic model of *G*. Let *w* be a state in that model such that $w \in \boldsymbol{CKR}$. We want to show that $\sigma(w) = (\sigma_1(w),...,\sigma_n(w)) \in S^\infty$. We shall prove the stronger claim that, for every Player *i* and for every $w' \in I_{CK}(w)$ (recall that $I_{CK}(w)$ is the cell of the common knowledge partition that contains *w*) and for every $m \geq 0$, $\sigma_i(w') \notin D_i^m$, where $D_i^m \subseteq S_i$ is the set of pure strategies of Player *i* that are strictly dominated in game $G^m$ ($G^m$ is the subgame of *G* obtained at step *m* of the iterated deletion procedure; we define $G^0$ to be $G$ itself). Note that this is equivalent to stating that, for every $w' \in I_{CK}(w)$, $\sigma(w') \in S^\infty$. Since $w \in I_{CK}(w)$, it follows from this that $\sigma(w) \in S^\infty$. The proof is by induction.

1. BASE STEP ($m = 0$). Select an arbitrary $w' \in I_{CK}(w)$ and an arbitrary Player $i$. Since $w \in CK\boldsymbol{R}$, $I_{CK}(w) \subseteq \boldsymbol{R}$ and thus $w' \in \boldsymbol{R}$; furthermore, since $\boldsymbol{R} \subseteq \boldsymbol{R}_i$, $w' \in \boldsymbol{R}_i$. Thus $\sigma_i(w')$ is a best response to the beliefs of Player $i$ at state $w'$. Hence, by Theorem 9.3, $\sigma_i(w') \notin D_i^0$.

2. INDUCTIVE STEP: assuming that the statement is true for all $k \leq m$ (for some $m \geq 0$) we prove that it is true for $k+1$. The hypothesis is that, for every Player $i$ and for every $w' \in I_{CK}(w)$, $\sigma_i(w') \notin D_i^k$; that is, for every $w' \in I_{CK}(w)$, $\sigma_i(w') \in S^{k+1}$. Select an arbitrary $w' \in I_{CK}(w)$ and an arbitrary Player $i$. Since $w \in CK\boldsymbol{R}$, $w' \in \boldsymbol{R}_i$ and thus $\sigma_i(w')$ is a best reply to Player $i$'s beliefs at $w'$ which, by hypothesis, attach positive probability only to strategy profiles of the other players that belong to $S_{-i}^{k+1}$ (note that, since $w' \in I_{CK}(w)$, $I_{CK}(w') = I_{CK}(w)$ and, by definition of common knowledge partition, $I_i(w') \subseteq I_{CK}(w')$). By Theorem 9.3, it follows that $\sigma_i(w') \notin D_i^{k+1}$.

**Proof of Theorem 9.2**. Given a strategic-form game with cardinal payoffs $G = \left\langle I, \left(S_i\right)_{i \in I}, \left(\pi_i\right)_{i \in I}\right\rangle$, construct the following epistemic model of $G$: $W = S^\infty$ and, for every player $i$ and every $s \in S^\infty$, $I_i(s) = \left\{s' \in S^\infty : s'_i = s_i\right\}$ (that is, the information set of Player $i$ that contains $s$ consists of all the strategy profiles in $S^\infty$ that have the same strategy of Player $i$ as in $s$); finally, let $\sigma_i(s) = s_i$. Select an arbitrary $s \in S^\infty$ and an arbitrary Player $i$. By definition of $S^\infty$, it is not the case that $s_i$ is strictly dominated in the restriction of $G$ to the sets of pure strategies $S_1^\infty, ..., S_{i-1}^\infty, S_i, S_{i-1}^\infty, S_n^\infty$. Thus, by Theorem 9.3, there is a probability distribution over $S_{-i}^\infty$ that makes $s_i = \sigma_i(s)$ a best reply. Choose one such probability distribution and let that probability distribution give Player $i$'s beliefs at $s$. Then $s \in \boldsymbol{R}_i$. Since $i$ was chosen arbitrarily, $s \in \boldsymbol{R}$; hence, since $s \in S^\infty$ was chosen arbitrarily, $\boldsymbol{R} = S^\infty$. It follows that $s \in CK\boldsymbol{R}$ for every $s \in S^\infty$.

# Appendix 9.E: Exercises

## 9.E.1. Exercises for Section 9.1: Models of strategic-form games

The answers to the following exercises are in Appendix S at the end of this chapter.

**Exercise 9.1**. Consider the following game (where the payoffs are von Neumann-Morgenstern payoffs)

| Player 1 \ Player 2 | *L* | *M* | *R* |
|---|---|---|---|
| *A* | 3 , 5 | 2 , 0 | 2 , 2 |
| *B* | 5 , 2 | 1 , 2 | 2 , 1 |
| *C* | 9 , 0 | 1 , 5 | 3 , 2 |

and the following model of this game

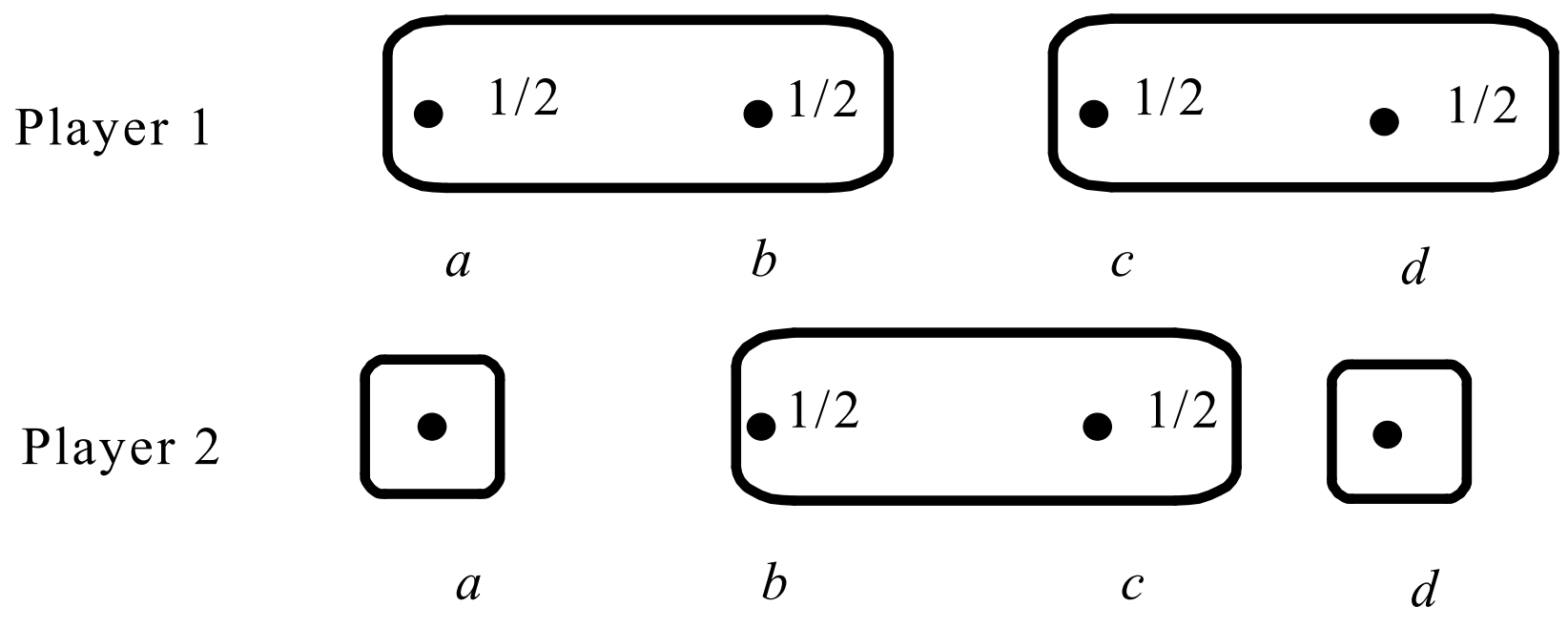

| | | | | |
|---|---|---|---|---|
| 1's strategy: | C | C | A | A |
| 2's strategy: | M | L | L | R |

**(a)** Find the event $R_1$ (that is, the set of states where Player 1 is rational).

**(b)** Find the event $R_2$ (that is, the set of states where Player 2 is rational).

**(c)** Find the event $R$ (that is, the set of states where both players are rational).

**Exercise 9.2**. Consider the following three-player game and model of it.

Player 3 chooses *E*

| Player 1 \ Player 2 | *C* | *D* |
|---|---|---|
| *A* | 2 , 3 , 2 | 1 , 0 , 4 |
| *B* | 3 , 6 , 4 | 0 , 8 , 0 |

Player 3 chooses *F*

| Player 1 \ Player 2 | *C* | *D* |
|---|---|---|
| *A* | 0 , 0 , 3 | 2 , 9 , 7 |
| *B* | 2 , 1 , 3 | 0 , 3 , 1 |

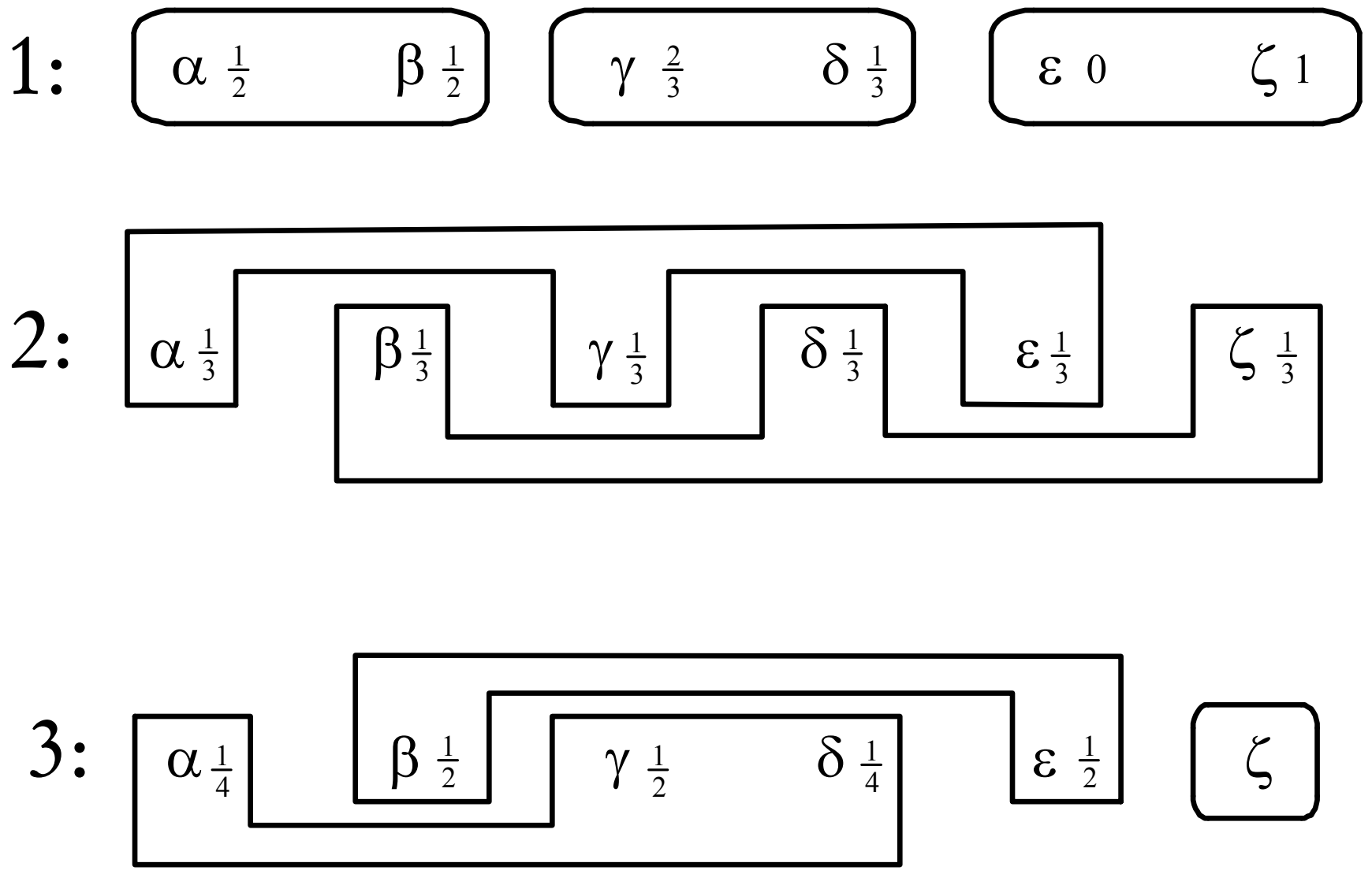

| | | | | | | |
|---|---|---|---|---|---|---|
| $\sigma_1$ : | *A* | *A* | *B* | *B* | *B* | *B* |
| $\sigma_2$ : | *C* | *D* | *C* | *D* | *C* | *D* |
| $\sigma_3$ : | *E* | *F* | *E* | *E* | *F* | *F* |

(a) Find the event $R_1$ (that is, the set of states where Player 1 is rational).

(b) Find the event $R_2$ (that is, the set of states where Player 2 is rational).

(c) Find the event $R_3$ (that is, the set of states where Player 3 is rational).

(d) Find the event $R$ (that is, the set of states where all players are rational).

## 9.E.2. Exercises for Section 9.2: Common knowledge of rationality

The answers to the following exercises are in Appendix S at the end of this chapter.

**Exercise 9.3.** (**a**) For the game and model of Exercise 9.1 find the following events: (*i*) $K_1\boldsymbol{R}_2$, (*ii*) $K_2\boldsymbol{R}_1$, (*iii*) $K_2K_1\boldsymbol{R}_2$, (*iv*) $K_1K_2\boldsymbol{R}_1$, (*v*) $CK\boldsymbol{R}_1$, (*vi*) $CK\boldsymbol{R}_2$, (*vii*) $CK\boldsymbol{R}$.

(**b**) Suppose that you found a model and a state *w* in that model such that $w \in CK\boldsymbol{R}$. What strategy profile could you find at *w*?

**Exercise 9.4.** (**a**) For the game and model of Exercise 9.2 find the following events: (*i*) $K_1\boldsymbol{R}$, (*ii*) $K_2\boldsymbol{R}$, (*iii*) $K_3\boldsymbol{R}$, (*iv*) $CK\boldsymbol{R}_3$, (*v*) $CK\boldsymbol{R}$.

(**b**) Suppose that you found a model and a state *w* in that model such that $w \in CK\boldsymbol{R}$. What strategy profile could you find at *w*?

**Exercise 9.5.** For the game of Exercise 9.1 construct a model where there is a state at which there is common knowledge of rationality and the strategy profile played there is (*C*,*L*). [Hints: (1) four states are sufficient, (2) it is easiest to postulate degenerated beliefs where a player assigns probability 1 to a particular state in his information set.]

**Exercise 9.6.** For the game of Exercise 9.2 construct a model where there is a state at which there is common knowledge of rationality. [Hint: think carefully, you don't need many states!]

◊◊◊◊◊◊◊◊◊◊◊◊◊

**Exercise 9.7 (Challenging Question).** Consider the following extensive-form game (where $o_1, o_2, o_3$ and $o_4$ are the possible outcomes)

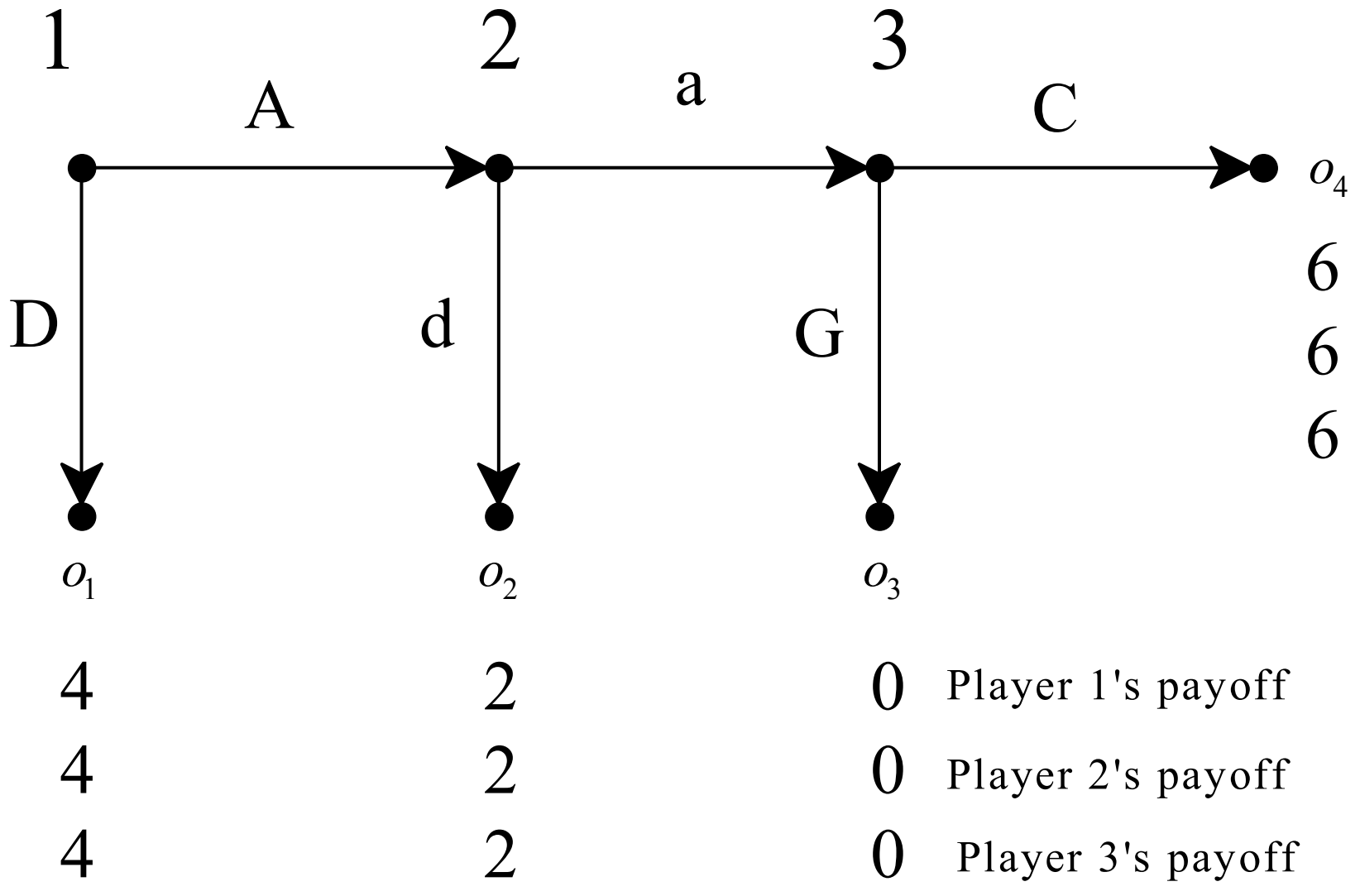

**(a)** What is the backward-induction outcome?

**(b)** Write the strategic-form game associated with the extensive-form game and find all the outcomes that are supported by a pure-strategy Nash equilibrium.

**(c)** Considering models of the strategic form, what strategy profiles are compatible with common knowledge of rationality?

**(d)** Choose a strategy profile which is not a Nash equilibrium and construct a model of the strategic form where at a state there is common knowledge of rationality and the associated strategy profile is the one you selected.

# Appendix 9.S: Solutions to exercises

**Exercise 9.1.** **(a)** $\boldsymbol{R}_1 = \{a,b\}$. At states *a* and *b* Player 1 believes that Player 2 is playing either *L* or *M* with equal probability. Thus his expected payoff from playing *A* is $\frac{1}{2}(3)+\frac{1}{2}(2)=2.5$, from playing *B* is $\frac{1}{2}(5)+\frac{1}{2}(1)=3$ and from playing *C* is $\frac{1}{2}(9)+\frac{1}{2}(1)=5$. Hence the best choice is *C* and this is indeed his choice at those two states. Thus $a,b \in \boldsymbol{R}_1$. At states *c* and *d* Player 1 plays *A* with an expected payoff of $\frac{1}{2}(3)+\frac{1}{2}(2)=2.5$ but he could get a higher expected payoff with *C* (namely $\frac{1}{2}(9)+\frac{1}{2}(3)=6$); thus $c,d \notin \boldsymbol{R}_1$. Hence $\boldsymbol{R}_1 = \{a,b\}$.
**(b)** $\boldsymbol{R}_2 = \{a,b,c\}$. **(c)** $\boldsymbol{R} = \{a,b\}$.

**Exercise 9.2.** **(a)** For each state let us calculate Player 1's expected payoff (denoted by $\mathbb{E}\pi_1$) from playing *A* and from playing *B*. **At states $\alpha$ and $\beta$,**

$$\mathbb{E}\pi_1(A) = \tfrac{1}{2}\pi_1(A,C,E) + \tfrac{1}{2}\pi_1(A,D,F) = \tfrac{1}{2}(2)+\tfrac{1}{2}(2) = 2,$$
$$\mathbb{E}\pi_1(B) = \tfrac{1}{2}\pi_1(B,C,E) + \tfrac{1}{2}\pi_1(B,D,F) = \tfrac{1}{2}(3)+\tfrac{1}{2}(0) = 1.5.$$

Thus *A* is optimal and hence, since $\sigma_1(\alpha)=\sigma_1(\beta)=A$, $\alpha,\beta \in \boldsymbol{R}_1$.

**At states $\gamma$ and $\delta$,** $\mathbb{E}\pi_1(A) = \frac{2}{3}\pi_1(A,C,E)+\frac{1}{3}\pi_1(A,D,E) = \frac{2}{3}(2)+\frac{1}{3}(1)=\frac{5}{3}$ and $\mathbb{E}\pi_1(B) = \frac{2}{3}\pi_1(B,C,E)+\frac{1}{3}\pi_1(B,D,E) = \frac{2}{3}(3)+\frac{1}{3}(0)=\frac{6}{3}$. Thus *B* is optimal and hence, since $\sigma_1(\gamma)=\sigma_1(\delta)=B$, $\gamma,\delta \in \boldsymbol{R}_1$.

**At states $\varepsilon$ and $\zeta$,** Player 1 assigns probability 1 to (*D*,*F*) against which *A* is the unique best reply and yet $\sigma_1(\varepsilon)=\sigma_1(\zeta)=B$; hence $\varepsilon,\zeta \notin \boldsymbol{R}_1$.

Thus $\boldsymbol{R}_1 = \{\alpha,\beta,\gamma,\delta\}$.

**(b)** For each state let us calculate Player 2's expected payoff (denoted by $\mathbb{E}\pi_2$) from playing *C* and from playing *D*. **At states $\alpha$, $\gamma$ and $\varepsilon$,**

$$\mathbb{E}\pi_2(C) = \tfrac{1}{3}\pi_2(A,C,E)+\tfrac{1}{3}\pi_2(B,C,E)+\tfrac{1}{3}\pi_2(B,C,F) = \tfrac{1}{3}(3)+\tfrac{1}{3}(6)+\tfrac{1}{3}(1)=\tfrac{10}{3},$$
$$\mathbb{E}\pi_2(D) = \tfrac{1}{3}\pi_2(A,D,E)+\tfrac{1}{3}\pi_2(B,D,E)+\tfrac{1}{3}\pi_2(B,D,F) = \tfrac{1}{3}(0)+\tfrac{1}{3}(8)+\tfrac{1}{3}(3)=\tfrac{11}{3}.$$

Thus *D* is optimal and hence, since $\sigma_2(\alpha)=\sigma_2(\gamma)=\sigma_2(\varepsilon)=C$, $\alpha,\gamma,\varepsilon \notin \boldsymbol{R}_2$.

**At states $\beta$, $\delta$ and $\zeta$,**

$$\mathbb{E}\pi_2(C) = \tfrac{1}{3}\pi_2(A,C,F)+\tfrac{1}{3}\pi_2(B,C,E)+\tfrac{1}{3}\pi_2(B,C,F) = \tfrac{1}{3}(0)+\tfrac{1}{3}(6)+\tfrac{1}{3}(1)=\tfrac{7}{3},$$
$$\mathbb{E}\pi_2(D) = \tfrac{1}{3}\pi_2(A,D,F)+\tfrac{1}{3}\pi_2(B,D,E)+\tfrac{1}{3}\pi_2(B,D,F) = \tfrac{1}{3}(9)+\tfrac{1}{3}(8)+\tfrac{1}{3}(3)=\tfrac{20}{3}.$$

Thus *D* is optimal and hence, since $\sigma_2(\beta)=\sigma_2(\delta)=\sigma_2(\zeta)=D$, $\beta,\delta,\zeta \in \boldsymbol{R}_2$.

Thus $\boldsymbol{R}_2 = \{\beta,\delta,\zeta\}$.

**(c)** For each state let us calculate Player 3's expected payoff (denoted by $\mathbb{E}\pi_3$) from playing *E* and from playing *F*. **At states $\alpha$, $\gamma$ and $\delta$,**

$\mathbb{E}\pi_3(E) = \frac{1}{4}\pi_3(A,C,E) + \frac{1}{2}\pi_3(B,C,E) + \frac{1}{4}\pi_3(B,D,E) = \frac{1}{4}(2) + \frac{1}{2}(4) + \frac{1}{4}(0) = 2.5$,

$\mathbb{E}\pi_3(F) = \frac{1}{4}\pi_3(A,C,F) + \frac{1}{2}\pi_3(B,C,F) + \frac{1}{4}\pi_3(B,D,F) =, = \frac{1}{4}(3) + \frac{1}{2}(3) + \frac{1}{4}(1) = 2.5$.

Thus both $E$ and $F$ are optimal and hence $\alpha, \gamma, \delta \in \boldsymbol{R}_3$.

**At states $\beta$ and $\varepsilon$,**

$$\mathbb{E}\pi_3(E) = \tfrac{1}{2}\pi_3(A,D,E) + \tfrac{1}{2}\pi_3(B,C,E) = \tfrac{1}{2}(4) + \tfrac{1}{2}(4) = 4,$$

$$\mathbb{E}\pi_3(F) = \tfrac{1}{2}\pi_3(A,D,F) + \tfrac{1}{2}\pi_3(B,C,F) =, = \tfrac{1}{2}(7) + \tfrac{1}{2}(3) = 5,$$

Thus $F$ is optimal and hence, since $\sigma_3(\beta) = \sigma_3(\varepsilon) = F$, $\beta, \varepsilon \in \boldsymbol{R}_3$.

**At state $\zeta$** Player 3 knows that Players 1 and 2 play (*B,D*) and she is best replying with *F*. Thus $\zeta \in \boldsymbol{R}_3$. Hence $\boldsymbol{R}_3 = \{\alpha, \beta, \gamma, \delta, \varepsilon, \zeta\}$

**(d)** $\boldsymbol{R} = \boldsymbol{R}_1 \cap \boldsymbol{R}_2 \cap \boldsymbol{R}_3 = \{\beta, \delta\}$.

**Exercise 9.3**. In Exercise 9.1 we determined that $\boldsymbol{R}_1 = \{a,b\}$, $\boldsymbol{R}_2 = \{a,b,c\}$ and $\boldsymbol{R} = \{a,b\}$. Thus

**(a)** (*i*) $K_1\boldsymbol{R}_2 = \{a,b\}$, (*ii*) $K_2\boldsymbol{R}_1 = \{a\}$, (*iii*) $K_2K_1\boldsymbol{R}_2 = \{a\}$, (*iv*) $K_1K_2\boldsymbol{R}_1 = \varnothing$. The common knowledge partition consists of a single information set containing all the states. Thus (*v*) $CK\boldsymbol{R}_1 = \varnothing$, (*vi*) $CK\boldsymbol{R}_2 = \varnothing$, (*vii*) $CK\boldsymbol{R} = \varnothing$.

**(b)** By Theorem 9.1, at a state at which there is common knowledge of rationality one can only find a strategy profile that survives the iterated deletion of strictly dominates strategies. In this game, for Player 1 strategy *B* is strictly dominated by the mixed strategy $\begin{pmatrix} A & C \\ \frac{1}{2} & \frac{1}{2} \end{pmatrix}$; after deleting *B*, for Player 2 strategy *R* is strictly dominated by the mixed strategy $\begin{pmatrix} L & M \\ \frac{1}{2} & \frac{1}{2} \end{pmatrix}$. Thus the iterated deletion of strictly dominated strategies yields the set of strategy profiles $\{(A,L),(A,M),(C,L),(C,M)\}$. Hence at a state where there is common knowledge of rationality one could only find one of these four strategy profiles.

**Exercise 9.4.** **(a)** **(*i*)** $K_1\boldsymbol{R} = \varnothing$, **(*ii*)** $K_2\boldsymbol{R} = \varnothing$, **(*iii*)** $K_3\boldsymbol{R} = \varnothing$. The common knowledge partition consists of a single information set containing all the states. Thus **(*iv*)** $CK\boldsymbol{R}_3 = \{\alpha,\beta,\gamma,\delta,\varepsilon,\zeta\}$, **(*v*)** $CK\boldsymbol{R} = \varnothing$.

**(b)** By Theorem 9.1, at a state at which there is common knowledge of rationality one can only find a strategy profile that survives the iterated deletion of strictly dominated strategies. Since in this game there are no strictly dominated strategies, at such a state one could find *any* strategy profile.

**Exercise 9.5.** In the following model there is common knowledge of rationality at every state. At state $\alpha$ the strategy profile played is $(C,L)$.

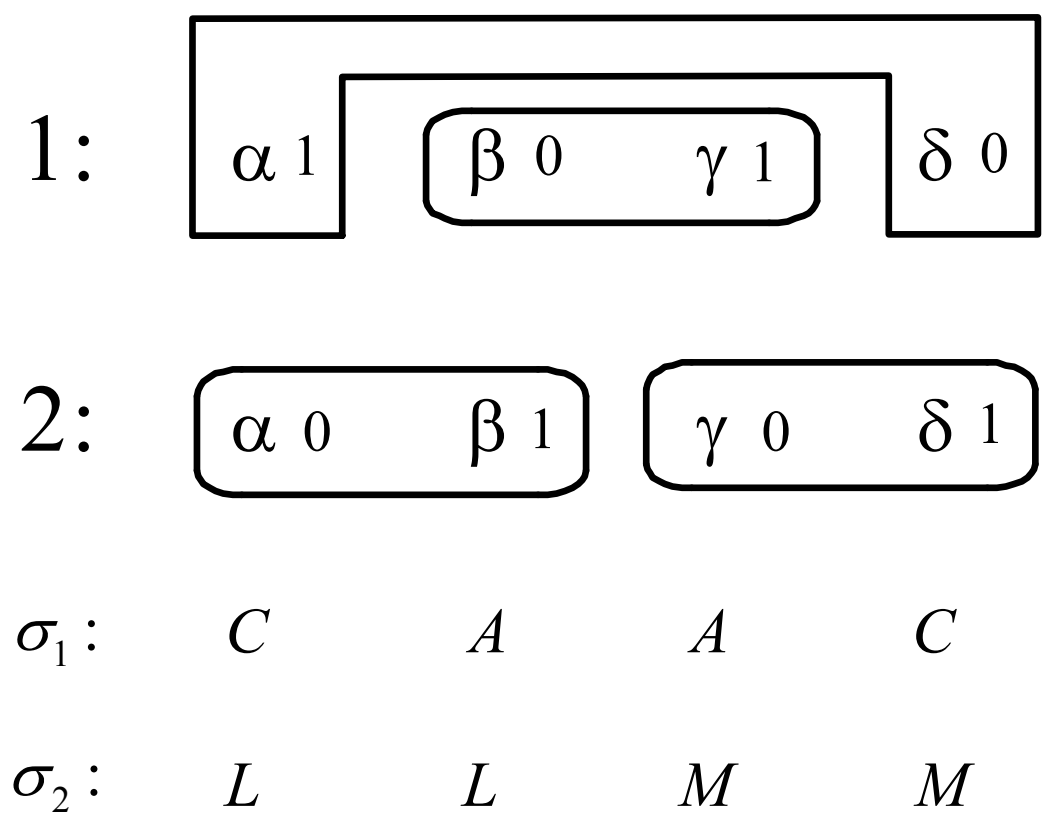

**Exercise 9.6.** Whenever a game has a Nash equilibrium in pure strategies, a one-state model where that Nash equilibrium is played is such that there is common knowledge of rationality at that state. Thus the following model provides a possible answer:

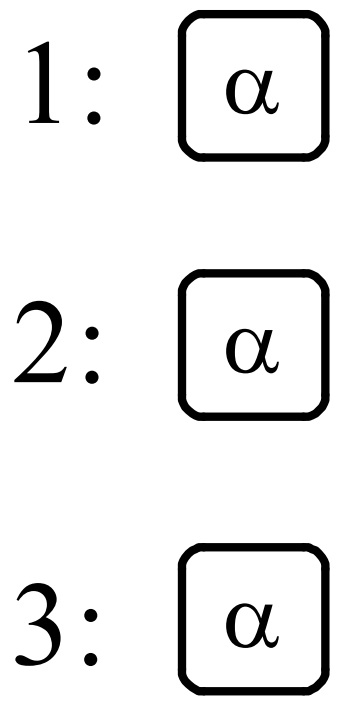

*A*
*D*
*F*

**Exercise 9.7** (Challenging Question).

**(a)** The backward induction outcome is $o_4$, with associated strategy profile (*A*,*a*,*C*) and payoffs (6,6,6).

**(b)** The strategic form is as follows (the Nash equilibria have been highlighted):

**Player 3** chooses *G*

| | | Player 2: *d* | Player 2: *a* |
|---|---|---|---|
| **Player 1** | *D* | 4 , 4 , 4 | 4 , 4 , 4 |
| | *A* | 2 , 2 , 2 | 0 , 0 , 0 |

**Player 3** chooses *C*

| | | Player 2: *d* | Player 2: *a* |
|---|---|---|---|
| **Player 1** | *D* | 4 , 4 , 4 | 4 , 4 , 4 |
| | *A* | 2 , 2 , 2 | 6 , 6 , 6 |

The Nash equilibria are: (*D*,*d*,*G*), (*D*,*a*,*G*), (*D*,*d*,*C*) and (*A*,*a*,*C*). Thus the only two outcomes sustained by a Nash equilibrium are $o_1$, with payoffs (4,4,4), and $o_4$, with payoffs (6,6,6).

(c) By Theorem 9.1, the outcomes that are compatible with common knowledge of rationality are those associated with strategy profiles that survive the iterated deletion of strictly dominated strategies. Since no player has any strictly dominated strategies, all the outcomes are compatible with common knowledge of rationality.

(d) In the following model at state $\alpha$ the players choose (*A*, *d*, *G*), which is not a Nash equilibrium; furthermore, there is no Nash equilibrium whose associated outcome is $o_2$ (with payoffs (2,2,2)). In this model $\alpha \in$ CK**R** (in particular, Player 1's choice of *A* is rational, given his belief that Player 2 plays *d* and *a* with equal probability).

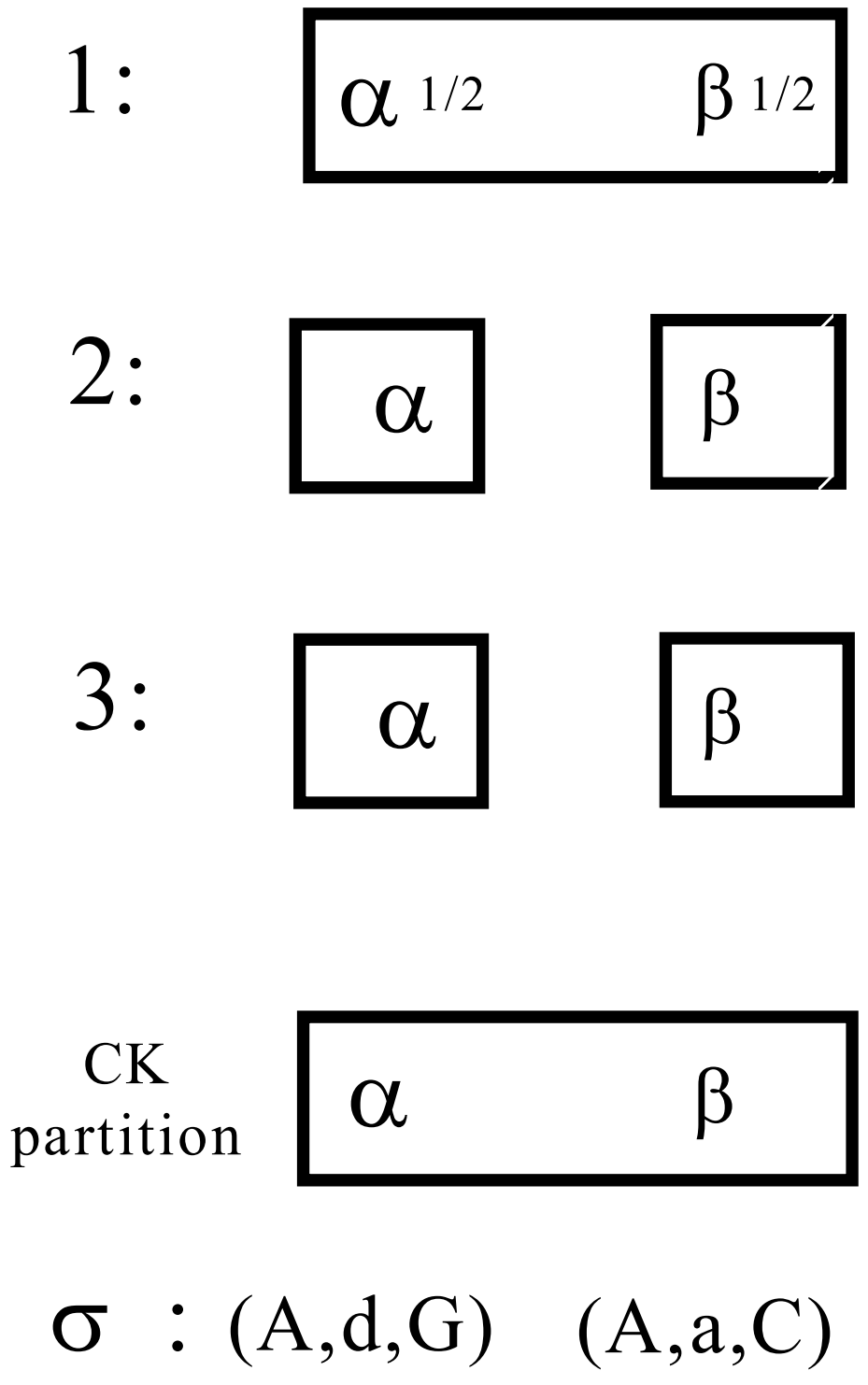

# PART IV

# Advanced Topics II:

## Refinements of Subgame-Perfect Equilibrium

Chapter

10

# A First Attempt: Weak Sequential Equilibrium

## 10.1 Assessments and sequential rationality

At the end of Chapter 6 (Section 6.3) we showed that, although the notion of subgame-perfect equilibrium is a refinement of Nash equilibrium, it is not strong enough to eliminate all "unreasonable" Nash equilibria. One reason for this is that a subgame-perfect equilibrium $\sigma$ allows a player's strategy to include a strictly dominated choice at an information set that is not reached by the play induced by $\sigma$. In order to eliminate such possibilities we need to define the notion of equilibrium for dynamic games in terms of a more complex object than merely a strategy profile.

**Definition 10.1.** Given an extensive-form game $G$, an *assessment* for $G$ is a pair $(\sigma, \mu)$, where $\sigma$ is a profile of behavioral strategies and $\mu$ is a list of probability distributions, one for every information set, over the nodes in that information set. We call $\mu$ a *system of beliefs*.

The system of beliefs $\mu$ specifies, for every information set, the beliefs – about past moves – that the relevant player would have if informed that her information set had been reached. Consider, for example, the extensive-form game of Figure 10.1.

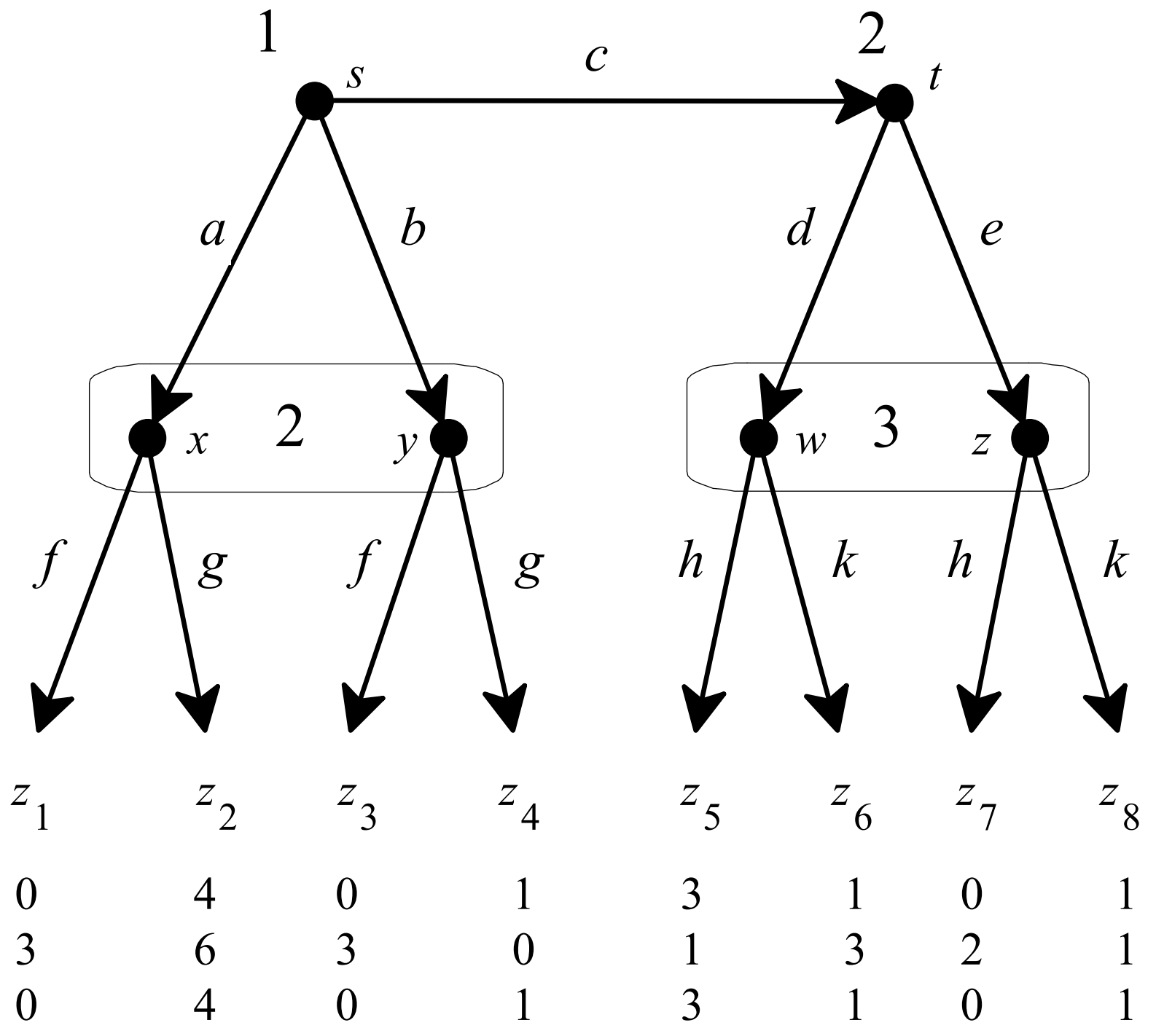

**Figure 10.1**

A possible assessment for this game is $(\sigma, \mu)$ with $\sigma = \left(\begin{array}{ccc|cc|cc|cc} a & b & c & f & g & d & e & h & k \\ \frac{1}{8} & \frac{3}{8} & \frac{4}{8} & 1 & 0 & \frac{3}{4} & \frac{1}{4} & \frac{1}{5} & \frac{4}{5} \end{array}\right)$ and $\mu = \left(\begin{array}{cc|cc} x & y & w & z \\ \frac{2}{3} & \frac{1}{3} & \frac{1}{2} & \frac{1}{2} \end{array}\right)$; note that, typically, we will not bother to include in $\mu$ the trivial probability distributions over singleton information sets.[1] The interpretation of this assessment is that Player 1 plans to play $a$ with probability $\frac{1}{8}$, $b$ with probability $\frac{3}{8}$ and $c$ with probability $\frac{4}{8}$; Player 2 plans to play $f$ if her information set $\{x, y\}$ is reached and to mix between $d$ and $e$ with probabilities $\frac{3}{4}$ and $\frac{1}{4}$, respectively, if her decision node $t$ is reached; Player 3 plans to mix between $h$ and $k$ with probabilities $\frac{1}{5}$ and $\frac{4}{5}$, respectively, if his information set $\{w, z\}$ is reached. Furthermore, Player 2 – if informed that her information set $\{x, y\}$ had been

[1] A complete specification of $\mu$ would be $\mu = \left(\begin{array}{c|c|cc|cc} s & t & x & y & w & z \\ 1 & 1 & \frac{2}{3} & \frac{1}{3} & \frac{1}{2} & \frac{1}{2} \end{array}\right)$ where $s$ is the decision node of Player 1 (the root) and $t$ is the decision node of Player 2 following choice $c$.

reached – would attach probability $\frac{2}{3}$ to Player 1 having played *a* (node *x*) and probability $\frac{1}{3}$ to Player 1 having played *b* (node *y*); Player 3 – if informed that his information set $\{w,z\}$ had been reached – would attach probability $\frac{1}{2}$ to node *w* (that is, to the sequence of moves *cd*) and probability $\frac{1}{2}$ to node *z* (that is, to the sequence of moves *ce*).

In order for an assessment $(\sigma,\mu)$ to be considered "reasonable" we will impose two requirements:

1. The choices specified by $\sigma$ should be optimal given the beliefs specified by $\mu$. We call this requirement *sequential rationality*.
2. The beliefs specified by $\mu$ should be consistent with the strategy profile $\sigma$. We call this requirement *Bayesian updating*.

Before we give a precise definition of these concepts, we shall illustrate them with reference to the above assessment $\sigma = \left(\begin{array}{ccc|cc|cc|cc} a & b & c & f & g & d & e & h & k \\ \frac{1}{8} & \frac{3}{8} & \frac{4}{8} & 1 & 0 & \frac{3}{4} & \frac{1}{4} & \frac{1}{5} & \frac{4}{5} \end{array}\right)$ and $\mu = \left(\begin{array}{cc|cc} x & y & w & z \\ \frac{2}{3} & \frac{1}{3} & \frac{1}{2} & \frac{1}{2} \end{array}\right)$ for the game of Figure 10.1. This assessment fails to satisfy sequential rationality, because, for example, at Player 3's information set {*w*,z} the planned mixed strategy $\begin{pmatrix} h & k \\ \frac{1}{5} & \frac{4}{5} \end{pmatrix}$ yields Player 3 – given her beliefs $\begin{pmatrix} w & z \\ \frac{1}{2} & \frac{1}{2} \end{pmatrix}$ – an expected payoff of $\frac{1}{2}\left[\frac{1}{5}(3) + \frac{4}{5}(1)\right] + \frac{1}{2}\left[\frac{1}{5}(0) + \frac{4}{5}(1)\right] = \frac{11}{10}$ while she could get a higher expected payoff, namely $\frac{1}{2}(3) + \frac{1}{2}(0) = \frac{15}{10}$, with the pure strategy *h*. This assessment also fails the rule for belief updating (Definition 8.3, Chapter 8); for example, given $\sigma$, the prior probability that node *x* is reached is $P(x) = \frac{1}{8}$ (it is the probability with which Player 1 plays *a*) and the prior probability that node *y* is reached is $P(y) = \frac{3}{8}$ (the probability with which Player 1 plays *b*), so that updating on information $\{x,y\}$ one gets $P(x \mid \{x,y\}) = \dfrac{P(x)}{P(\{x,y\})} = \dfrac{P(x)}{P(x)+P(y)} = \dfrac{\frac{1}{8}}{\frac{1}{8}+\frac{3}{8}} = \frac{1}{4}$ and $P(y \mid \{x,y\}) = \frac{3}{4}$. Thus, in order to be consistent with the rule for belief updating, Player 2's beliefs should be $\begin{pmatrix} x & y \\ \frac{1}{4} & \frac{3}{4} \end{pmatrix}$.

Now we can turn to the formal definitions. First we need to introduce some notation. If $\sigma$ is a profile of behavior strategies and $a$ is a choice of Player $i$ (at some information set of Player $i$), we denote by $\sigma(a)$ the probability that $\sigma_i$ (the strategy of Player $i$ that is part of $\sigma$) assigns to $a$. For example, for the game of Figure 10.1, if $\sigma = \left(\begin{array}{ccc|cc|cc|cc} a & b & c & f & g & d & e & h & k \\ \frac{1}{8} & \frac{3}{8} & \frac{4}{8} & 1 & 0 & \frac{3}{4} & \frac{1}{4} & \frac{1}{5} & \frac{4}{5} \end{array}\right)$ then $\sigma(b) = \frac{3}{8}$, $\sigma(g) = 0$, $\sigma(d) = \frac{3}{4}$, etc. Similarly, if $\mu$ is a system of beliefs and $x$ is a decision node, then we denote by $\mu(x)$ the probability that the relevant part of $\mu$ assigns to $x$. For example, if $\mu = \left(\begin{array}{c|c|cc|cc} s & t & x & y & w & z \\ 1 & 1 & \frac{2}{3} & \frac{1}{3} & \frac{1}{2} & \frac{1}{2} \end{array}\right)$ or, written more succinctly, $\mu = \left(\begin{array}{cc|cc} x & y & w & z \\ \frac{2}{3} & \frac{1}{3} & \frac{1}{2} & \frac{1}{2} \end{array}\right)$ then $\mu(s) = 1$, $\mu(y) = \frac{1}{3}$, $\mu(w) = \frac{1}{2}$, etc.

Recall that $Z$ denotes the set of terminal nodes and, for every Player $i$, $\pi_i : Z \to \mathbb{R}$ is the payoff function of Player $i$. Given a decision node $x$, let $Z(x) \subseteq Z$ be the set of terminal nodes that can be reached starting from $x$. For example, in the game of Figure 10.1, $Z(t) = \{z_5, z_6, z_7, z_8\}$. Given a behavior strategy profile $\sigma$ and a decision node $x$, let $\mathbb{P}_{x,\sigma}$ be the probability distribution over $Z(x)$ induced by $\sigma$, that is, if $z \in Z(x)$ and $\langle a_1, \ldots, a_m \rangle$ is the sequence of choices that leads from $x$ to $z$ then $\mathbb{P}_{x,\sigma}(z)$ is the product of the probabilities of those choices: $\mathbb{P}_{x,\sigma}(z) = \sigma(a_1) \times \sigma(a_2) \times \ldots \times \sigma(a_m)$. For example, in the game of Figure 10.1, if $\sigma$ is the strategy profile $\sigma = \left(\begin{array}{ccc|cc|cc|cc} a & b & c & f & g & d & e & h & k \\ \frac{1}{8} & \frac{3}{8} & \frac{4}{8} & 1 & 0 & \frac{3}{4} & \frac{1}{4} & \frac{1}{5} & \frac{4}{5} \end{array}\right)$ and $t$ is Player 2's decision node after choice $c$ of Player 1, then $\mathbb{P}_{t,\sigma}(z_5) = \sigma(d)\sigma(h) = \frac{3}{4}\left(\frac{1}{5}\right) = \frac{3}{20}$. If $H$ is an information set of Player $i$ we denoted by $\pi_i\left(H \mid \sigma, \mu\right)$ the expected payoff of Player $i$ starting from information set $H$, given the beliefs specified by $\mu$ at $H$ and given the choices prescribed by $\sigma$ at $H$ and at the information sets that come after $H$, that is,

$$\pi_i\left(H \mid \sigma, \mu\right) = \sum_{x \in H} \left[ \mu(x) \left( \sum_{z \in Z(x)} \mathbb{P}_{x,\sigma}(z)\, \pi_i(z) \right) \right].$$

For example, in the game of Figure 10.1, if

$$\sigma = \left(\begin{array}{ccc|cc|cc|cc} a & b & c & f & g & d & e & h & k \\ \frac{1}{8} & \frac{3}{8} & \frac{4}{8} & 1 & 0 & \frac{3}{4} & \frac{1}{4} & \frac{1}{5} & \frac{4}{5} \end{array}\right) \text{ and } \mu = \left(\begin{array}{cc|cc} x & y & w & z \\ \frac{2}{3} & \frac{1}{3} & \frac{1}{2} & \frac{1}{2} \end{array}\right),$$

then (as we computed above)

$$\pi_3\left(\{w,z\} \mid \sigma, \mu\right) =$$
$$\mu(w)\left(\mathbb{P}_{w,\sigma}(z_5)\pi_3(z_5) + \mathbb{P}_{w,\sigma}(z_6)\pi_3(z_6)\right) + \mu(z)\left(\mathbb{P}_{z,\sigma}(z_7)\pi_3(z_7) + \mathbb{P}_{z,\sigma}(z_8)\pi_3(z_8)\right) =$$
$$\tfrac{1}{2}\left[\tfrac{1}{5}(3) + \tfrac{4}{5}(1)\right] + \tfrac{1}{2}\left[\tfrac{1}{5}(0) + \tfrac{4}{5}(1)\right] = \tfrac{11}{10}.$$

Recall that if $\sigma$ is a strategy profile and $i$ is a player, then $\sigma_{-i}$ denotes the profile of strategies of the players other than $i$ and we often use $(\sigma_i, \sigma_{-i})$ as an alternative way of denoting $\sigma$; furthermore, if $\tau_i$ is a strategy of Player $i$, we denote by $(\tau_i, \sigma_{-i})$ the strategy profile obtained from $\sigma$ by replacing $\sigma_i$ with $\tau_i$ (and leaving everything else unchanged).

**Definition 10.2.** Fix an extensive-form game and an *assessment* $(\sigma, \mu)$. We say that Player *i*'s behavior strategy $\sigma_i$ is *sequentially rational* if, for every information set *H* of Player *i*,

$$\pi_i\left(H \mid (\sigma_i, \sigma_{-i}), \mu\right) \geq \pi_i\left(H \mid (\tau_i, \sigma_{-i}), \mu\right),$$

for every behavior strategy $\tau_i$ of Player *i*.

We say that $\sigma$ is *sequentially rational* if, for every Player *i*, $\sigma_i$ is sequentially rational.

Note that for Player *i*'s strategy $\sigma_i$ to be sequential rationality it is not sufficient (although it is necessary) that at every information set *H* of Player *i* the choice(s) at $H$ prescribed by $\sigma_i$ be optimal (given the choices of the other players prescribed by $\sigma_{-i}$): we need to check if Player *i* could improve her payoff by changing her choice(s) at *H and possibly also at information sets of hers that follow H.* To see this, consider the game of Figure 10.2 below.

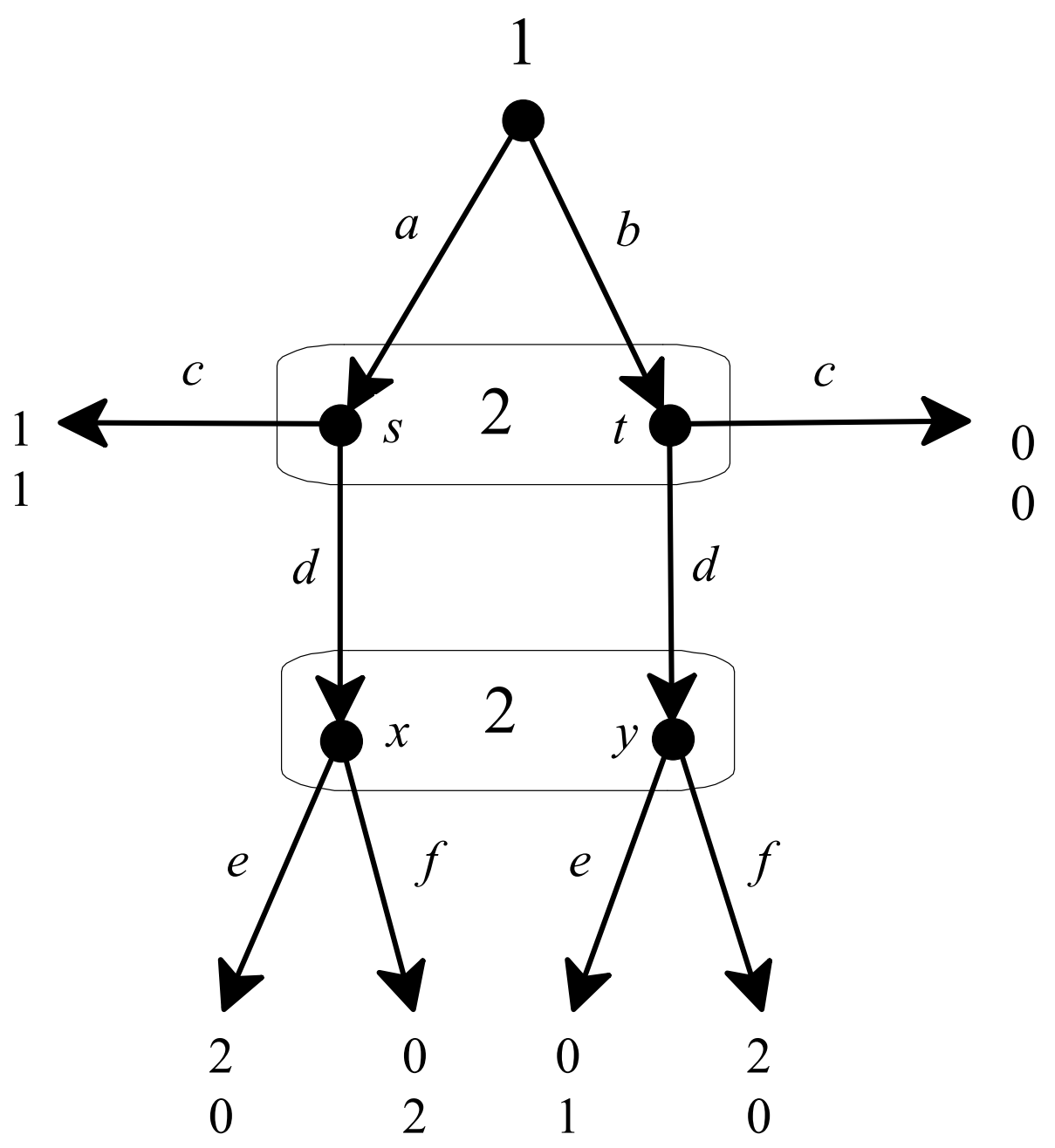

**Figure 10.2**

Let $(\sigma, \mu)$ be the assessment where $\sigma$ is the pure-strategy profile $(a,(c,e))$ and $\mu = \begin{pmatrix} s & t & | & x & y \\ 1 & 0 & | & 0 & 1 \end{pmatrix}$. For Player 2 *e* is rational at information set {*x*,*y*} because, given her belief that she is making her choice at node *y*, *e* gives her a payoff of 1 while *f* would give her a payoff of 0; furthermore, at information set {*s*,*t*} – given her belief that she is making her choice at node *s and given her future choice of e at* {*x*,*y*} – choice *c* is better than choice *d* because the former gives her a payoff of 1 while the latter a payoff of 0. However, the *strategy* (*c*,*e*) (while sequentially rational at {*x*,*y*}) is not sequentially rational at {*s*,*t*} because – given her belief that she is making her choice at node *s* – with (*c*,*e*) she gets a payoff of 1 but if she switched to (*d*,*f*) she would get a payoff of 2; in other words, Player 2 can improve her payoff by changing her choice at {*s*,*t*} from *c* to *d* and also changing her future planned choice at {*x*,*y*} from *e* to *f*.[2] Note that the pure-strategy profile $\left(a,(c,e)\right)$ is not a Nash equilibrium: for Player 2 the unique best reply to *a* is (*d,f*).

[2] It is possible to impose restrictions on the system of beliefs $\mu$ such that sequential rationality as defined in Definition 10.2 is equivalent to the weaker condition that, at every information set *H*, the corresponding player cannot increase her payoff by changing her choice(s) at *H* only. We will not discuss this condition here. The interested reader is referred to Hendon *et al* (1996) and Perea (2002).

An example of a sequentially rational assessment for the game of Figure 10.3 below is $\sigma = \begin{pmatrix} L & M & R & | & A & B & | & c & d \\ 0 & 0 & 1 & | & \frac{1}{2} & \frac{1}{2} & | & 1 & 0 \end{pmatrix}$ and $\mu = \begin{pmatrix} x & y & | & u & v & w \\ \frac{1}{4} & \frac{3}{4} & | & \frac{1}{5} & \frac{3}{5} & \frac{1}{5} \end{pmatrix}$.

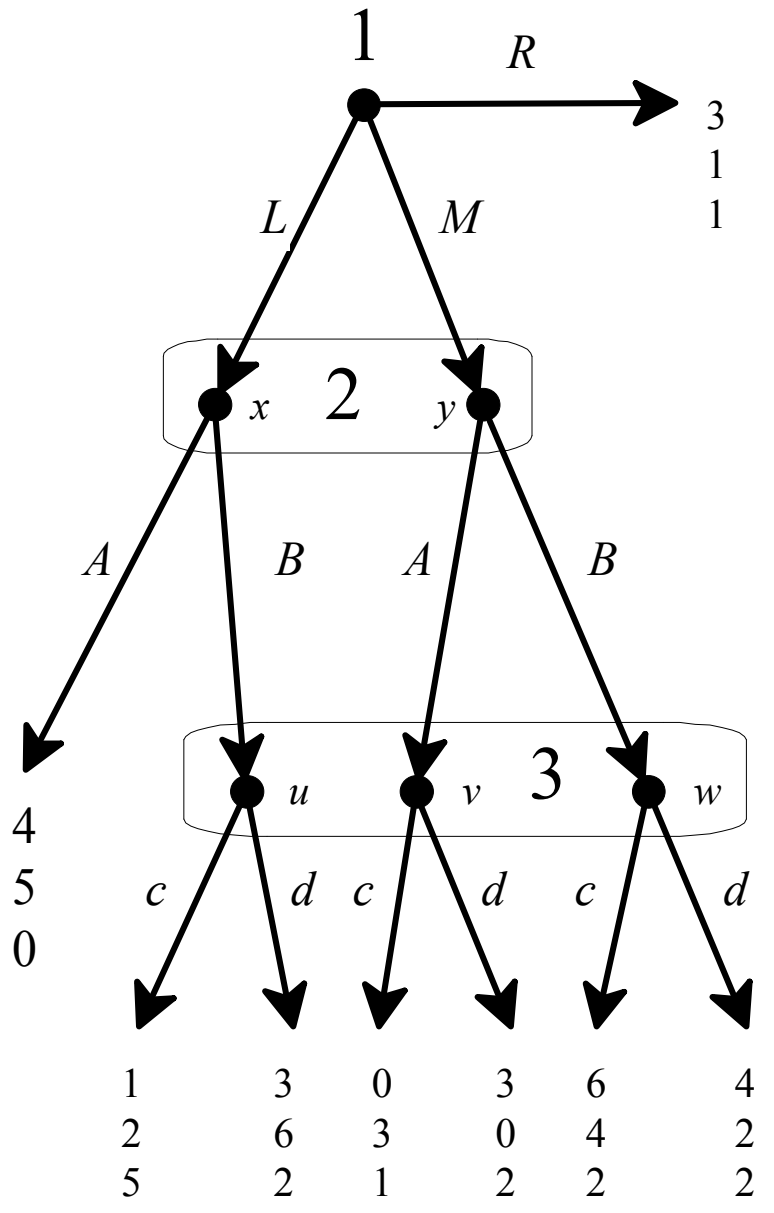

**Figure 10.3**

Let us first verify sequential rationality of Player 3's strategy. At her information set $\{u,v,w\}$, given her beliefs, *c* gives a payoff of $\frac{1}{5}(5)+\frac{3}{5}(1)+\frac{1}{5}(2)=2$, while *d* gives a payoff of $\frac{1}{5}(2)+\frac{3}{5}(2)+\frac{1}{5}(2)=2$. Thus *c* is optimal (as would be *d* and any randomization over *c* and *d*). Now consider Player 2: at his information set $\{x,y\}$, given his beliefs and given the strategy of Player 3, *A* gives a payoff of $\frac{1}{4}(5)+\frac{3}{4}(3)=3.5$ and *B* gives a payoff of $\frac{1}{4}(2)+\frac{3}{4}(4)=3.5$; thus any mixture of *A* and *B* is optimal, in particular, the mixture $\begin{pmatrix} A & B \\ \frac{1}{2} & \frac{1}{2} \end{pmatrix}$ is optimal. Finally, at the root, *R* gives Player 1 a payoff of 3, *L* a payoff of $\frac{1}{2}(4)+\frac{1}{2}(1)=2.5$ and *M* a payoff of $\frac{1}{2}(0)+\frac{1}{2}(6)=3$; thus *R* is optimal (as would be any mixture over *M* and *R*).

This is a good time to test your understanding of the concepts introduced in this section, by going through the exercises in Section 10.E.1 of Appendix 10.E at the end of this chapter.

## 10.2 Bayesian updating at reached information sets.

The second requirement for an assessment to be "reasonable" is that the beliefs encoded in $\mu$ should be consistent with the behavior postulated by $\sigma$, in the sense that the rule for belief updating (Definition 8.3, Chapter 8), should be used to form those beliefs, whenever it is applicable. We shall call this rule *Bayesian updating.*

Let $x$ be a node that belongs to information set $H$ and let $\mathbb{P}_{root,\sigma}(x)$ be the probability that node $x$ is reached (from the root of the tree) if $\sigma$ is implemented.[3] Let $\mathbb{P}_{root,\sigma}(H)=\sum_{y\in H}\mathbb{P}_{root,\sigma}(y)$ be the probability that information set $H$ is reached (that is, the probability that some node in $H$ is reached). Then Bayesian updating requires that the probability that is assigned to $x$ given the information that $H$ has been reached be given by the conditional probability $\mathbb{P}_{root,\sigma}(x\,|\,H)=\dfrac{\mathbb{P}_{root,\sigma}(x)}{\mathbb{P}_{root,\sigma}(H)}$; of course, this conditional probability is well defined if and only if $\mathbb{P}_{root,\sigma}(H)>0$.

**Definition 10.3.** Consider an extensive-form game and a behavioral strategy profile $\sigma$. We say that an information set $H$ is *reached by* $\sigma$ if $\mathbb{P}_{root,\sigma}(H)>0$.

[3] That is, if $\langle a_1,...,a_m\rangle$ is the sequence of choices that leads from the root to $x$ then $\mathbb{P}_{root,\sigma}(x)=\sigma(a_1)\times...\times\sigma(a_m)$.

For example, in the game of Figure 10.1 partially reproduced below (we omitted the payoffs), if $\sigma = \begin{pmatrix} a & b & c & f & g & d & e & h & k \\ \frac{1}{3} & \frac{2}{3} & 0 & 1 & 0 & \frac{3}{4} & \frac{1}{4} & 1 & 0 \end{pmatrix}$ then $\{x, y\}$ is reached by $\sigma$, while $\{w, z\}$ is not.[4]

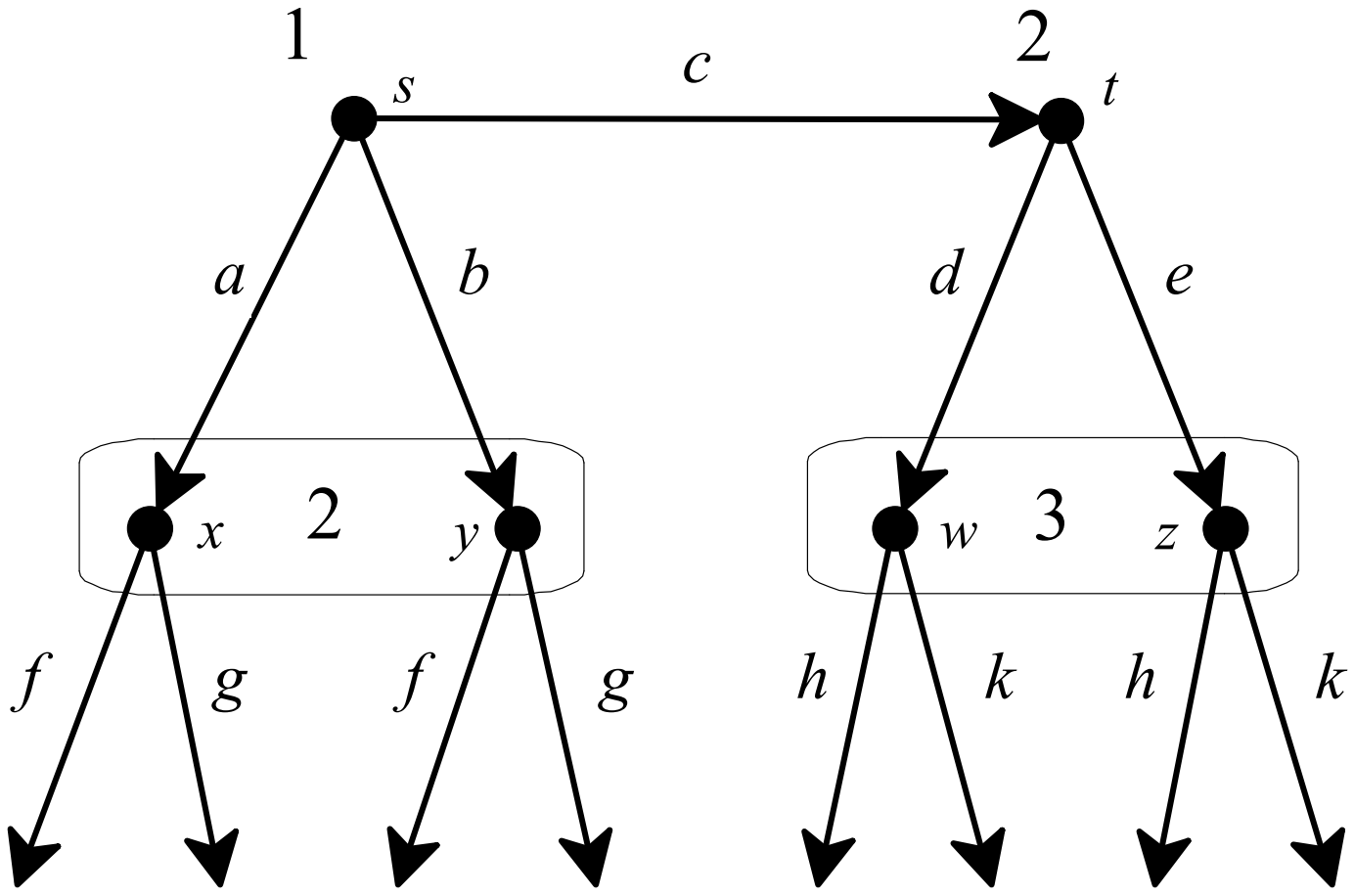

**Definition 10.4.** Given an extensive-form game and an assessment $(\sigma, \mu)$, we say that $(\sigma, \mu)$ satisfies *Bayesian updating at reached information sets* if, for every information set $H$ and every node $x \in H$, if $H$ is reached by $\sigma$ (that is, if $\mathbb{P}_{root,\sigma}(H) > 0$) then $\mu(x) = \dfrac{\mathbb{P}_{root,\sigma}(x)}{\mathbb{P}_{root,\sigma}(H)}$.

For example, in the game of Figure 10.1 partially reproduced above, the assessment $\sigma = \begin{pmatrix} a & b & c & f & g & d & e & h & k \\ \frac{1}{9} & \frac{5}{9} & \frac{3}{9} & 0 & 1 & \frac{3}{4} & \frac{1}{4} & 1 & 0 \end{pmatrix}$, $\mu = \begin{pmatrix} x & y & w & z \\ \frac{1}{6} & \frac{5}{6} & \frac{3}{4} & \frac{1}{4} \end{pmatrix}$ satisfies Bayesian updating at reached information sets. In fact, we have that $\mathbb{P}_{root,\sigma}(x) = \frac{1}{9}$, $\mathbb{P}_{root,\sigma}(y) = \frac{5}{9}$, $\mathbb{P}_{root,\sigma}(\{x, y\}) = \frac{1}{9} + \frac{5}{9} = \frac{6}{9}$, $\mathbb{P}_{root,\sigma}(w) = \frac{3}{9}\left(\frac{3}{4}\right) = \frac{9}{36}$, $\mathbb{P}_{root,\sigma}(z) = \frac{3}{9}\left(\frac{1}{4}\right) = \frac{3}{36}$ and $\mathbb{P}_{root,\sigma}(\{w, x\}) = \frac{9}{36} + \frac{3}{36} = \frac{12}{36}$. Thus

[4] $\mathbb{P}_{root,\sigma}(x) = \frac{1}{3}$, $\mathbb{P}_{root,\sigma}(y) = \frac{2}{3}$ and thus $\mathbb{P}_{root,\sigma}(\{x, y\}) = \frac{1}{3} + \frac{2}{3} = 1$; on the other hand, $\mathbb{P}_{root,\sigma}(w) = \sigma(c)\sigma(d) = 0\left(\frac{3}{4}\right) = 0$ and $\mathbb{P}_{root,\sigma}(z) = \sigma(c)\sigma(e) = 0\left(\frac{1}{4}\right) = 0$, so that $\mathbb{P}_{root,\sigma}(\{w, z\}) = 0$.

$$\frac{\mathbb{P}_{root,\sigma}(x)}{\mathbb{P}_{root,\sigma}(\{x,y\})} = \frac{\frac{1}{9}}{\frac{6}{9}} = \frac{1}{6} = \mu(x) \quad \text{and, similarly,} \quad \frac{\mathbb{P}_{root,\sigma}(y)}{\mathbb{P}_{root,\sigma}(\{x,y\})} = \frac{\frac{5}{9}}{\frac{6}{9}} = \frac{5}{6} = \mu(y),$$

$$\frac{\mathbb{P}_{root,\sigma}(w)}{\mathbb{P}_{root,\sigma}(\{w,z\})} = \frac{\frac{9}{36}}{\frac{12}{36}} = \frac{3}{4} = \mu(w) \text{ and } \frac{\mathbb{P}_{root,\sigma}(z)}{\mathbb{P}_{root,\sigma}(\{w,z\})} = \frac{\frac{3}{36}}{\frac{12}{36}} = \frac{1}{4} = \mu(z).$$

**Remark 10.1.** Note that the condition "Bayesian updating at reached information sets" is trivially satisfied at every information set $H$ that is *not* reached by $\sigma$, that is, at every information set $H$ such that $\mathbb{P}_{root,\sigma}(H) = 0$.[5]

This is a good time to test your understanding of the concepts introduced in this section, by going through the exercises in Section 10.E.2 of Appendix 10.E at the end of this chapter.

# 10.3 Weak sequential equilibrium

Our objective is to find a refinement of subgame-perfect equilibrium that rules out "unreasonable" equilibria. Moving from strategy profiles to assessments allowed us to judge the rationality of a choice at an information set independently of whether that information set is reached by the strategy profile under consideration: the requirement of sequential rationality rules out choices that are strictly dominated at an information set. As a first attempt we define a notion of equilibrium based on the two requirements of sequential rationality and Bayesian updating at reached information sets.

**Definition 10.5.** A *weak sequential equilibrium* is an assessment $(\sigma, \mu)$ which satisfies two requirements: (1) sequential rationality and (2) Bayesian updating at reached information sets.

[5] In logic a proposition of the form "if $A$ then $B$" ($A$ is called the *antecedent* and $B$ is called the *consequent*) is false only when $A$ is true and $B$ is false. Thus, in particular, the proposition "if $A$ then $B$" is true whenever $A$ is false (whatever the truth value of $B$). In our case the antecedent is "$\mathbb{P}_{root,\sigma}(H) > 0$".

Before we discuss the properties of weak sequential equilibria we illustrate the type of reasoning that one needs to go through in order to find weak sequential equilibria. Consider first the game of Figure 10.4.

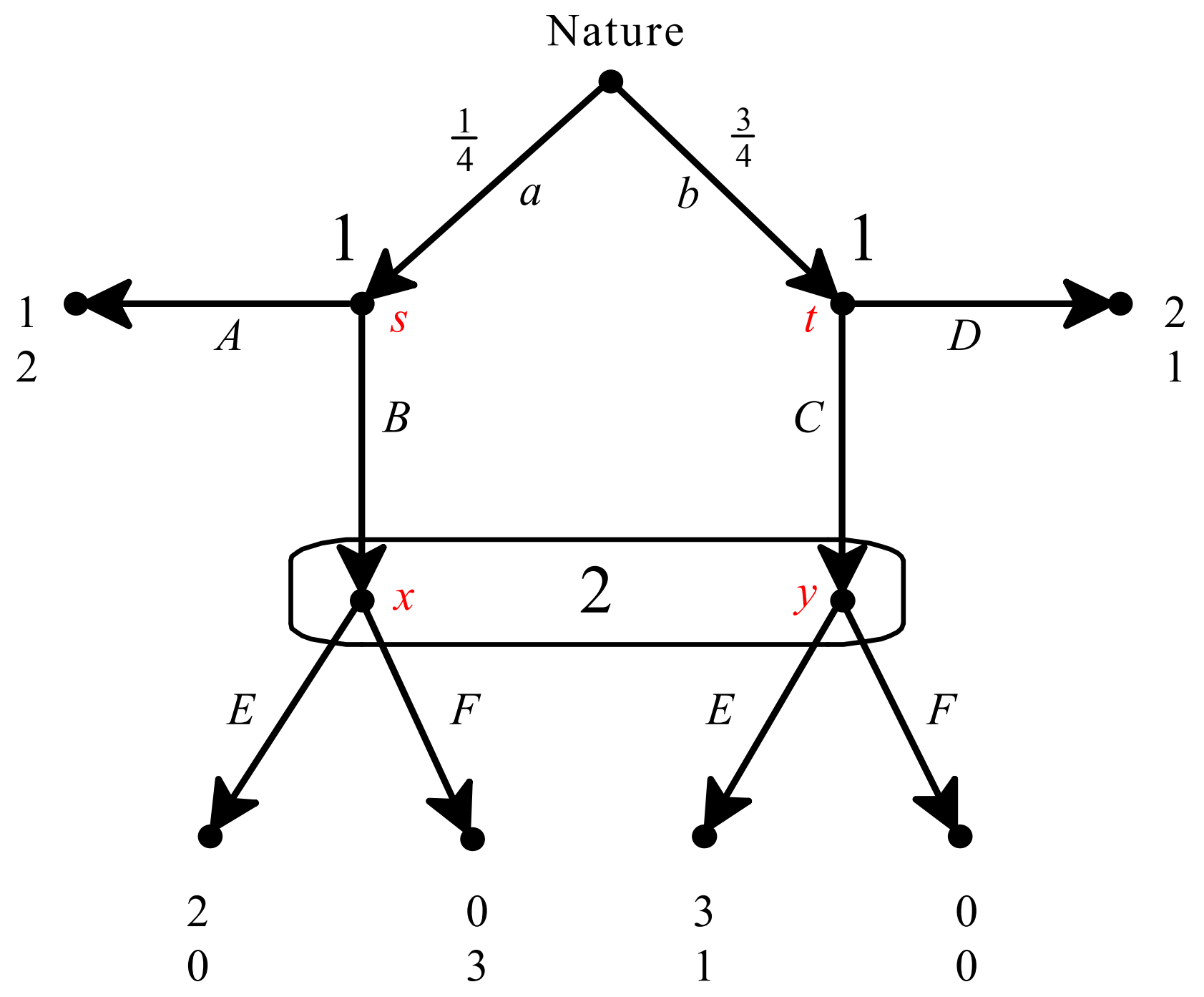

**Figure 10.4**

Let us see if there is a weak sequential equilibrium $(\sigma,\mu)$ where Player 1's strategy in $\sigma$ is a pure strategy. The set of pure strategies of Player 1 is

$$S_1 = \{(A,C),(B,D),(B,C),(A,D)\}$$

The strategy of Player 1 determines the beliefs of Player 2 at her information set $\{x,y\}$. Let us consider the four possibilities.

- If Player 1's strategy is $(A,C)$, then Player 2's information set $\{x,y\}$ is reached with positive probability and the only beliefs that are consistent with Bayesian updating are $\begin{pmatrix} x & y \\ 0 & 1 \end{pmatrix}$, so that – by sequential rationality – Player 2 must choose $E$. However, if Player 2's strategy is $E$ then at node $s$ it is not sequentially rational for Player 1 to choose $A$. Thus there is no weak sequential equilibrium where Player 1's strategy is $(A,C)$.

- If Player 1's strategy is (*B,D*), then Player 2's information set $\{x,y\}$ is reached with positive probability and the only beliefs that are consistent with Bayesian updating are $\begin{pmatrix} x & y \\ 1 & 0 \end{pmatrix}$, so that – by sequential rationality – Player 2 must choose *F*. However, if Player 2's strategy is *F* then at node *s* it is not sequentially rational for Player 1 to choose *B*. Thus there is no weak sequential equilibrium where Player 1's strategy is (*B,D*).

- If Player 1's strategy is (*B,C*), then Player 2's information set $\{x,y\}$ is reached (with probability 1) and the only beliefs that are consistent with Bayesian updating are $\begin{pmatrix} x & y \\ \frac{1}{4} & \frac{3}{4} \end{pmatrix}$. Given these beliefs, Player 2's payoff from playing *E* is $\frac{1}{4}(0)+\frac{3}{4}(1)=\frac{3}{4}$ and her payoff from playing *F* is $\frac{1}{4}(3)+\frac{3}{4}(0)=\frac{3}{4}$. Thus any mixed strategy is sequentially rational for Player 2, that is, for any $p\in[0,1]$, $\begin{pmatrix} E & F \\ p & 1-p \end{pmatrix}$ is sequentially rational. At node *s* choice *B* is sequentially rational for Player 1 if and only if the expected payoff from playing *B* is at least 1 (which is the payoff from playing *A*): $2p+0(1-p)\geq 1$, that is, $p\geq\frac{1}{2}$. At node *t* choice *C* is sequentially rational for Player 1 if and only if the expected payoff from playing *C* is at least 2 (which is the payoff from playing *D*): $3p+0(1-p)\geq 2$, that is, $p\geq\frac{2}{3}$. Hence, if $p\geq\frac{2}{3}$ then both *B* and *C* are sequentially rational. Thus we have an infinite number of weak sequential equilibria at which Player 1's strategy is (*B,C*): for every $p\in\left[\frac{2}{3},1\right]$, $(\sigma,\mu)$ is a weak sequential equilibrium, where

$$\sigma=\left(\begin{array}{cc|cc|cc} A & B & C & D & E & F \\ 0 & 1 & 1 & 0 & p & 1-p \end{array}\right),\quad \mu=\begin{pmatrix} x & y \\ \frac{1}{4} & \frac{3}{4} \end{pmatrix}.$$

- If Player 1's strategy is (*A*,*D*), then Player 2's information set $\{x, y\}$ is not reached and thus, according to the notion of weak sequential equilibrium, any beliefs are allowed there. Let $\begin{pmatrix} E & F \\ p & 1-p \end{pmatrix}$ be Player 2's strategy. From previous calculations we have that both $A$ and $D$ are sequentially rational if and only if $p \leq \frac{1}{2}$. One possibility is to set $p$ = 0; this means that Player 2 chooses the pure strategy *F* and this is sequentially rational if and only if her beliefs are $\begin{pmatrix} x & y \\ q & 1-q \end{pmatrix}$ with $\underbrace{3q + 0(1-q)}_{\text{payoff from } F} \geq \underbrace{0q + 1(1-q)}_{\text{payoff from } E}$ , that is, $q \geq \frac{1}{4}$. Thus we have an infinite number of weak sequential equilibria: $(\sigma, \mu)$ is a weak sequential equilibrium, where

  $$\sigma = \left(\begin{array}{cc|cc|cc} A & B & C & D & E & F \\ 0 & 1 & 1 & 0 & 0 & 1 \end{array}\right), \quad \mu = \begin{pmatrix} x & y \\ q & 1-q \end{pmatrix}, \quad \text{for every } q \in \left[\tfrac{1}{4}, 1\right].$$

  Next consider the case where $0 < p \leq \frac{1}{2}$, so that Player 2 uses a mixed strategy. This can be sequentially rational if and only if the payoff from *E* is equal to the payoff from *F*, that is (using the previous calculations) if and only if her beliefs are $\begin{pmatrix} x & y \\ \frac{1}{4} & \frac{3}{4} \end{pmatrix}$. Thus we have an infinite number of weak sequential equilibria: $(\sigma, \mu)$ is a weak sequential equilibrium, where

  $$\sigma = \left(\begin{array}{cc|cc|cc} A & B & C & D & E & F \\ 0 & 1 & 1 & 0 & p & 1-p \end{array}\right), \quad \mu = \begin{pmatrix} x & y \\ \frac{1}{4} & \frac{3}{4} \end{pmatrix}, \quad \text{for every } p \in (0, \tfrac{1}{2}]$$

Next we consider the more complex game shown in Figure 10.5.

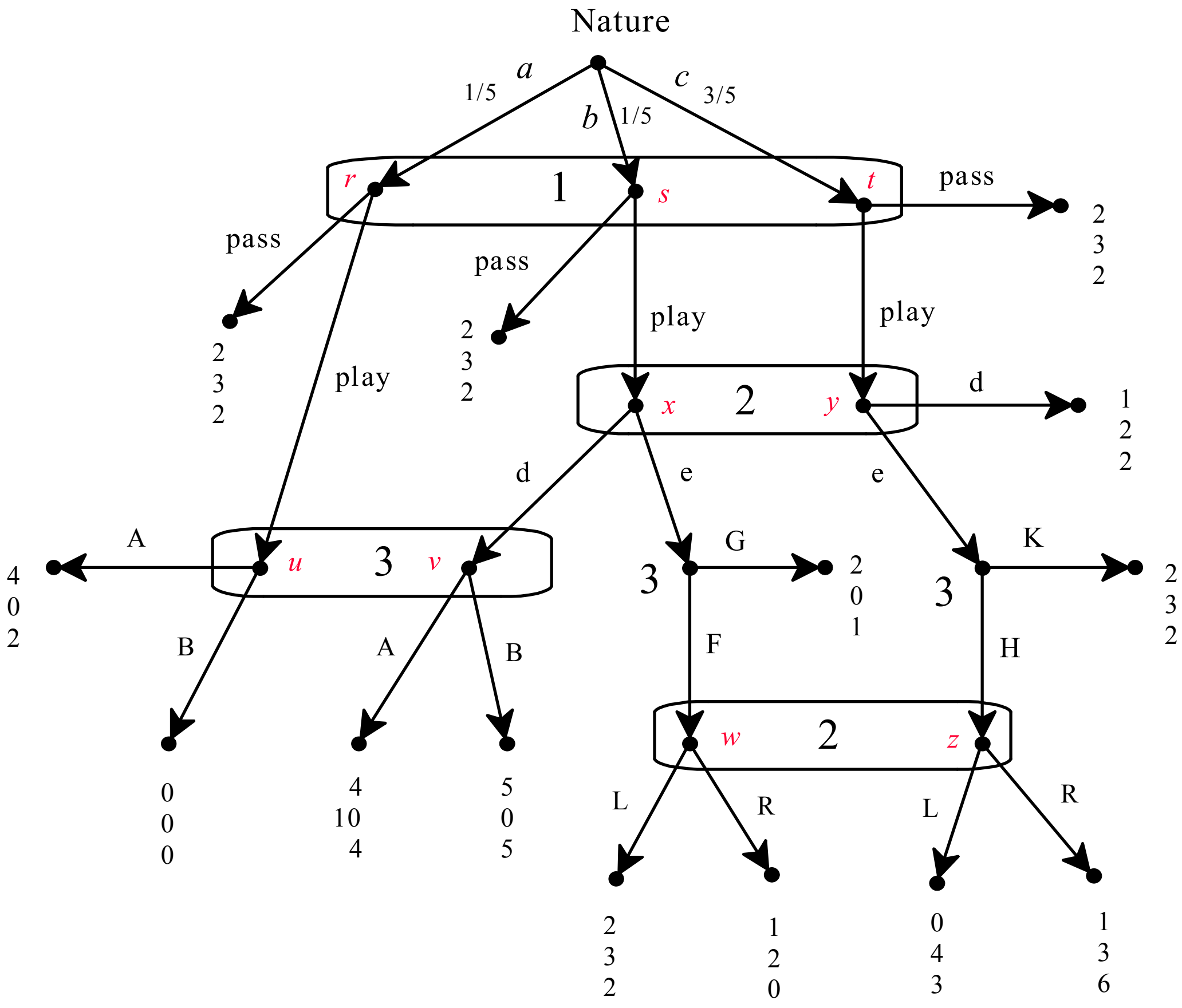

**Figure 10.5**

For simplicity, let us limit the search to pure-strategy weak sequential equilibria. How should we proceed? The first thing to do is to see if the game itself can be simplified; in particular, one should check if there are any information sets where there is a strictly dominant choice: if there is such a choice then, no matter what beliefs the relevant player has at that information set, sequential rationality requires that choice to be selected. In the game of Figure 10.5 there is indeed such an information set, namely information set $\{w,z\}$ of Player 2: here $L$ is strictly better than $R$ for Player 2 at both nodes, that is, $R$ is strictly dominated by $L$ as a choice. Thus we can simplify the game by (1) removing this information set, (2) converting nodes $w$ and $z$ into terminal nodes and (3) assigning to these newly created terminal nodes the payoffs associated with choice $L$. The simplified game is shown in Figure 10.6 below.

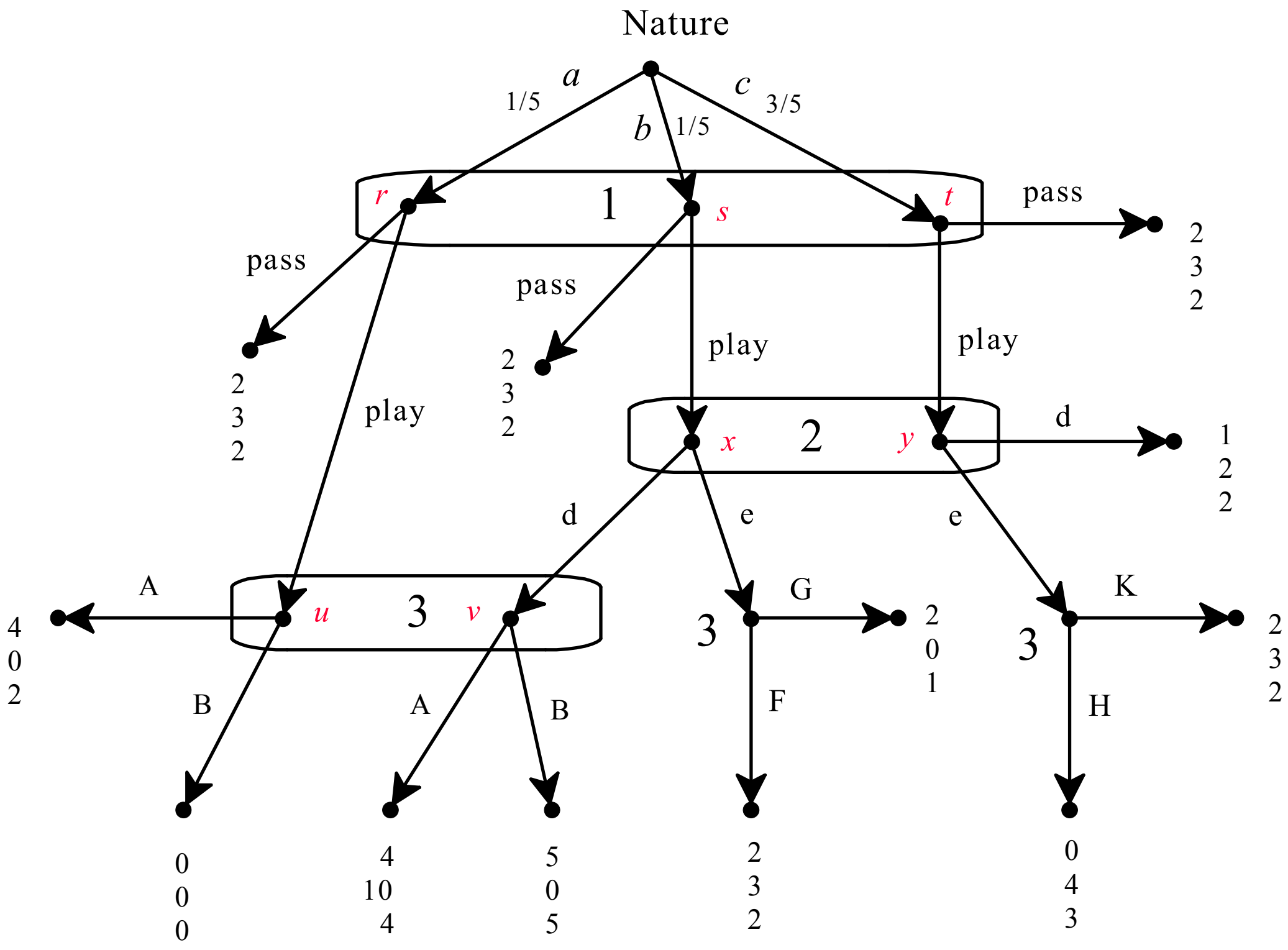

**Figure 10.6**

Applying the same reasoning to the simplified game of Figure 10.6, we can delete the two singleton decision nodes of Player 3 and replace one with the payoff associated with choice *F* and the other with the payoff associated with choice *H*, thus obtaining the further simplification shown in Figure 10.7 below. In this game there are no more information sets with strictly dominant choices. Thus we have to proceed by trial and error. Note first that, at any weak sequential equilibrium, by Bayesian updating Player 1's beliefs must be $\begin{pmatrix} r & s & t \\ \frac{1}{5} & \frac{1}{5} & \frac{3}{5} \end{pmatrix}$.

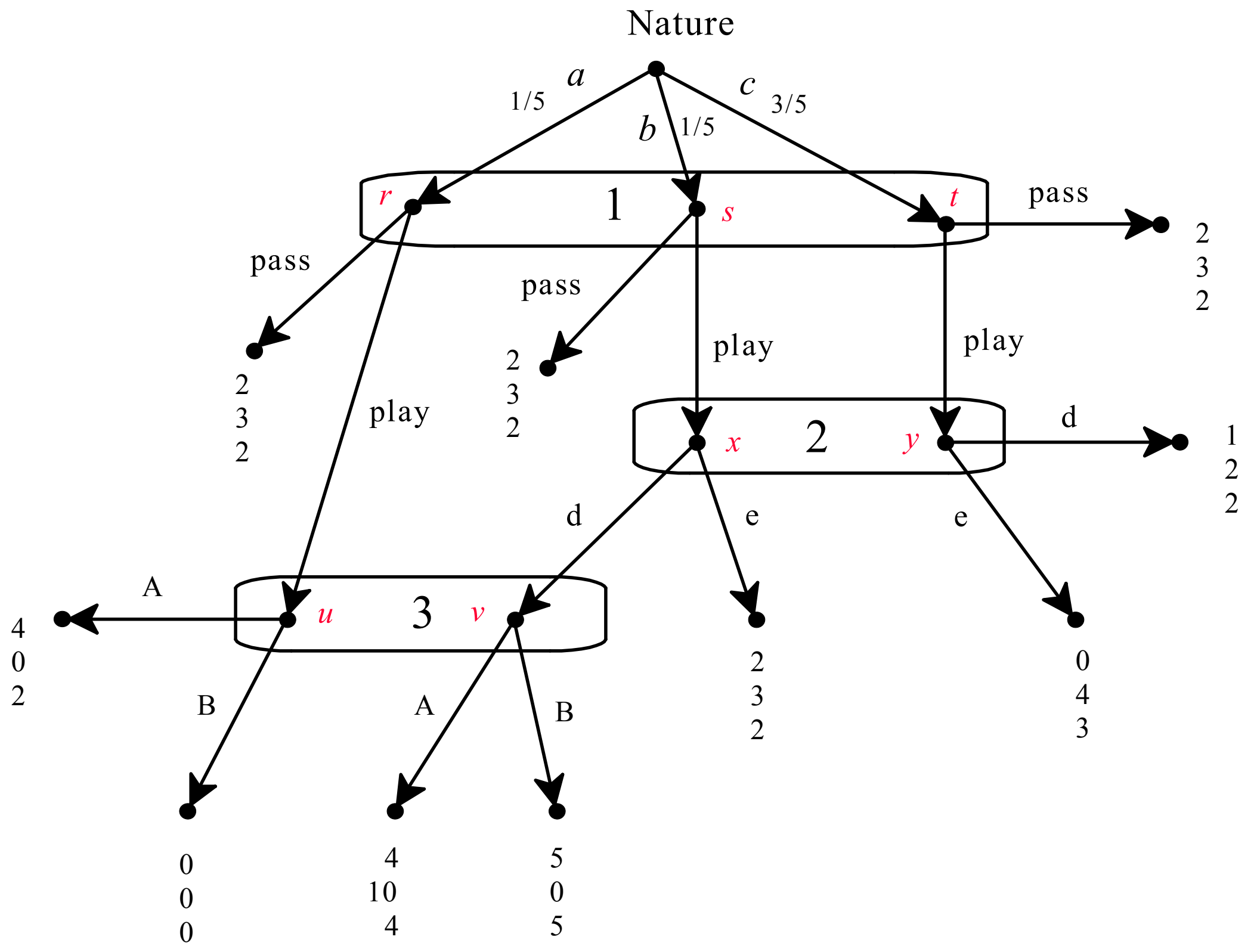

**Figure 10.7**

Let us see if in the simplified game of Figure 10.7 there is a pure-strategy weak sequential equilibrium where Player 1's strategy is "play". When Player 1 chooses "play" then, by Bayesian updating, Player 2's beliefs must be $\begin{pmatrix} x & y \\ \frac{1}{4} & \frac{3}{4} \end{pmatrix}$.

- Let us first try the hypothesis that there is a weak sequential equilibrium where (Player 1's strategy is "play" and) Player 2's strategy is *e*. Then, by Bayesian updating, Player 3's beliefs must be $\begin{pmatrix} u & v \\ 1 & 0 \end{pmatrix}$, making *A* the only sequentially rational choice at information set $\{u,v\}$. However, if Player 3's strategy is *A* then Player 2, at her information set $\{x,y\}$, gets a payoff of $\frac{1}{4}(10)+\frac{3}{4}(2)=\frac{16}{4}$ from playing *d* and a payoff of $\frac{1}{4}(3)+\frac{3}{4}(4)=\frac{15}{4}$ from playing *e*. Thus *e* is not sequentially rational and we have reached a contradiction.

- Let us now try the hypothesis that there is a weak sequential equilibrium where (Player 1's strategy is "play" and) Player 2's strategy is *d*. Then, by Bayesian updating, Player 3's beliefs must be $\begin{pmatrix} u & v \\ \frac{1}{2} & \frac{1}{2} \end{pmatrix}$; hence Player 3's expected payoff from playing *A* is $\frac{1}{2}(2)+\frac{1}{2}(4)=3$ and his expected payoff from playing *B* is $\frac{1}{2}(0)+\frac{1}{2}(5)=2.5$. Thus *A* is the only sequentially rational choice at information set {*u*,*v*}. Hence, by the previous calculations for Player 2, *d* is indeed sequentially rational for Player 2. Thus it only remains to check if Player 1 indeed wants to choose "play". Given the strategies *d* and *A* of Players 2 and 3, respectively, Player 1 gets a payoff of 2 from "pass" and a payoff of $\frac{1}{5}(4)+\frac{1}{5}(4)+\frac{3}{5}(1)=\frac{11}{5}>2$ from "play". Thus "play" is indeed sequentially rational. Thus we have found a pure-strategy weak sequential equilibrium of the game of Figure 10.7, namely

$$\sigma = \left(\begin{array}{cc|cc|cc} pass & play & d & e & A & B \\ 0 & 1 & 1 & 0 & 1 & 0 \end{array}\right), \quad \mu = \left(\begin{array}{ccc|cc|cc} r & s & t & x & y & u & v \\ \frac{1}{5} & \frac{1}{5} & \frac{3}{5} & \frac{1}{4} & \frac{3}{4} & \frac{1}{2} & \frac{1}{2} \end{array}\right).$$

This equilibrium can be extended to a weak sequential equilibrium of the original game of Figure 10.5 by adding the choices that led to the simplified game of Figure 10.7 and arbitrary beliefs at information set {*w*,*z*}:

$$\sigma = \left(\begin{array}{cc|cc|cc|cc|cc|cc} pass & play & d & e & A & B & F & G & H & K & L & R \\ 0 & 1 & 1 & 0 & 1 & 0 & 1 & 0 & 1 & 0 & 1 & 0 \end{array}\right),$$

$$\mu = \left(\begin{array}{ccc|cc|cc|cc} r & s & t & x & y & u & v & w & z \\ \frac{1}{5} & \frac{1}{5} & \frac{3}{5} & \frac{1}{4} & \frac{3}{4} & \frac{1}{2} & \frac{1}{2} & p & 1-p \end{array}\right), \text{ for any } p \in [0,1]$$

Are there any other pure-strategy weak sequential equilibria? This question is addressed in Exercise 10.7.

We now explore the relationship between the notion of weak sequential equilibrium and other equilibrium concepts.

**Theorem 10.1.** Given an extensive-form game with cardinal payoffs *G*, if $(\sigma,\mu)$ is a weak sequential equilibrium of *G* then $\sigma$ is a Nash equilibrium of *G*.

In general, not every Nash equilibrium can be part of a weak sequential equilibrium. To see this, consider the extensive-form game of Figure 10.8 that reproduces the game of Figure 6.12 (Chapter 6). The pure-strategy profile (*b*,*f*,*d*) is a Nash equilibrium, but it cannot be part of any weak sequential equilibrium, because – no matter what beliefs Player 3 has at her information set – choice *d* is not sequentially rational (it is strictly dominated by *c* at both nodes in that information set).

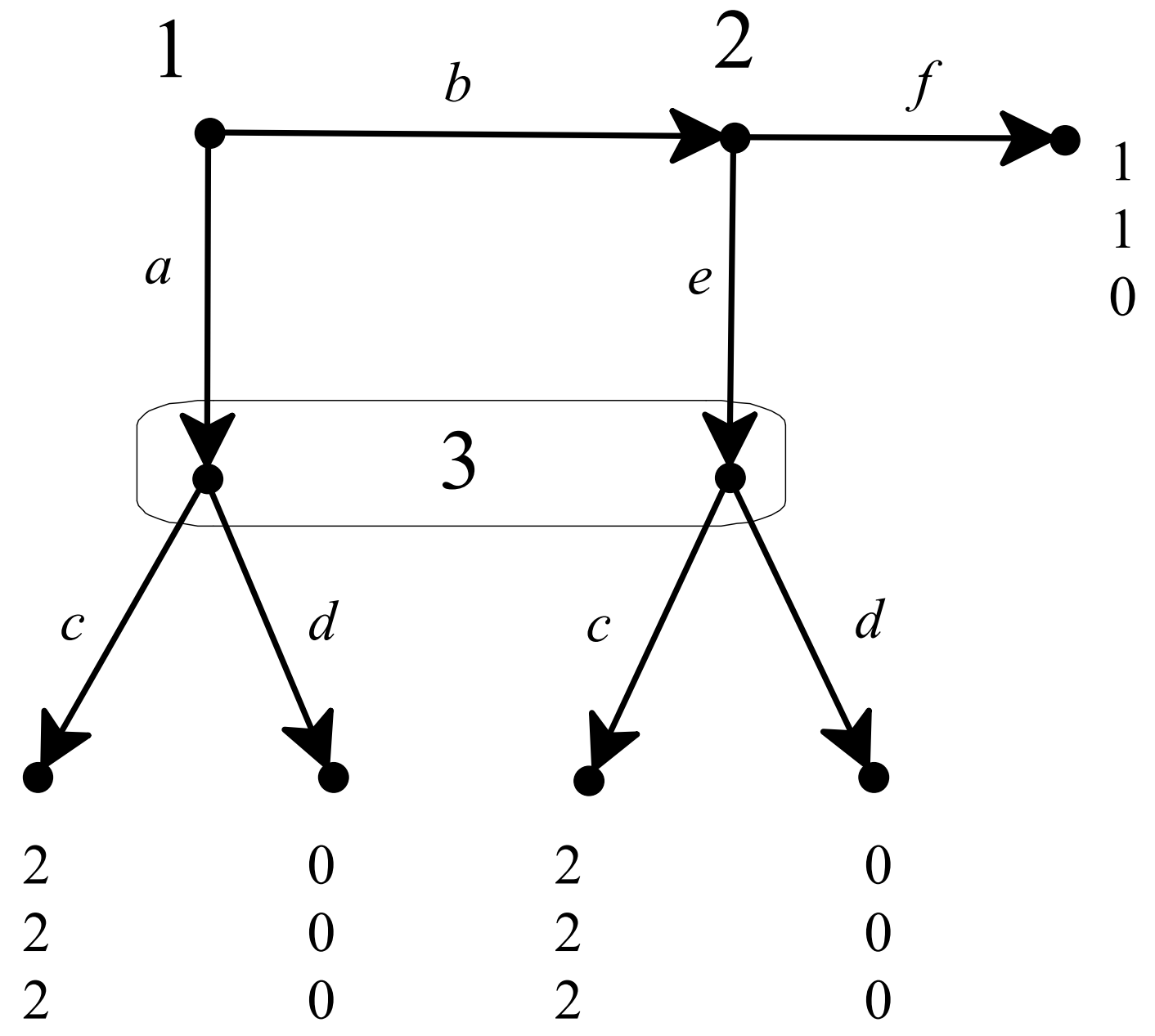

**Figure 10.8**

Thus weak sequential equilibrium is a strict refinement of Nash equilibrium. Does it also refine the notion of subgame-perfect equilibrium? Unfortunately, the answer is negative: it is possible for $(\sigma,\mu)$ to be a weak sequential equilibrium without $\sigma$ being a subgame-perfect equilibrium. To see this, consider the game of Figure 10.9, which reproduces Figure 10.1.

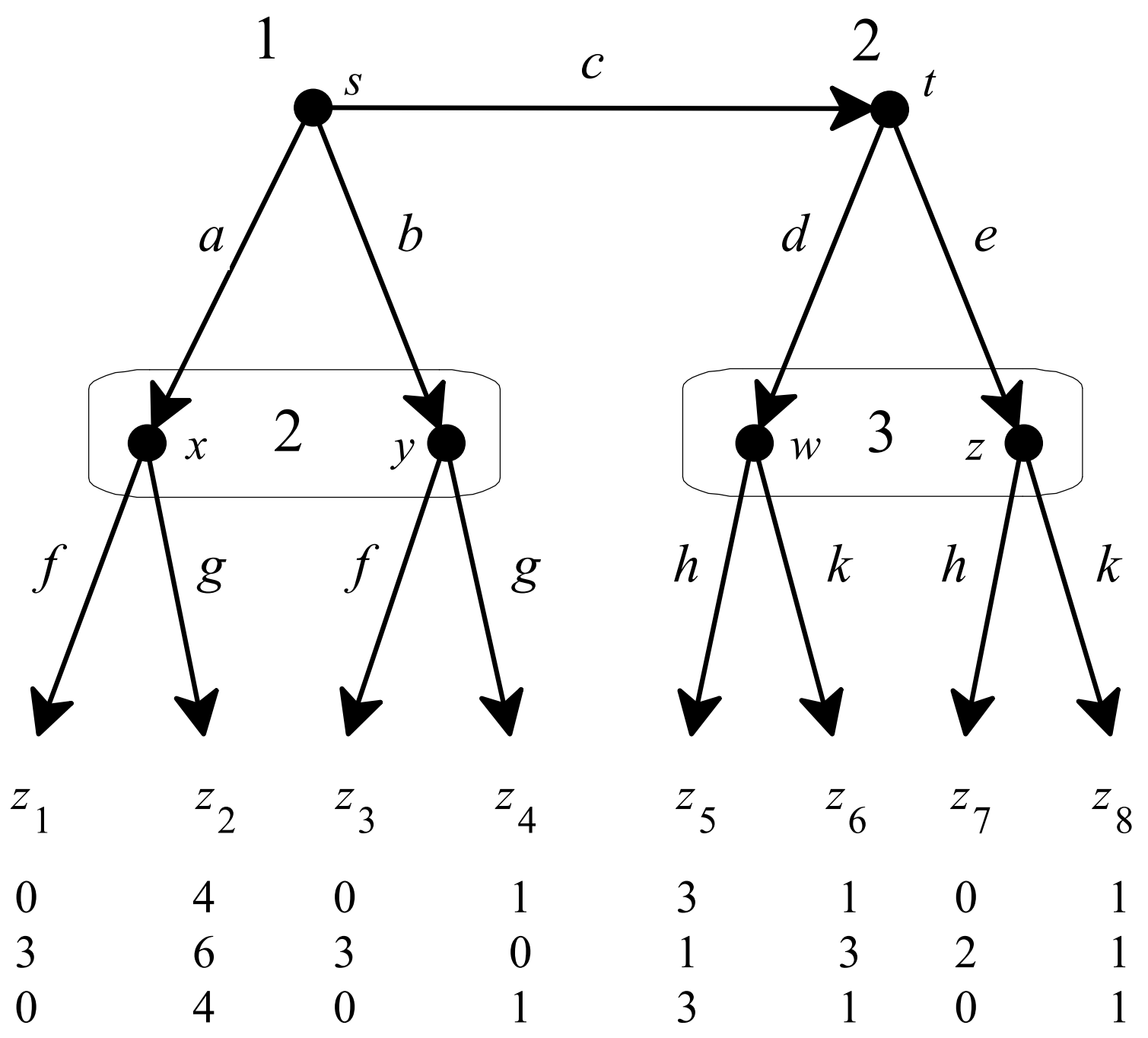

**Figure 10.9**

The assessment $\sigma = (b,(f,e),h)$ with $\mu = \begin{pmatrix} x & y & | & w & z \\ 0 & 1 & | & 1 & 0 \end{pmatrix}$ is a weak sequential equilibrium (the reader should verify this; in particular, note that information set {*w,z*} is not reached and thus Bayesian updating allows for arbitrary beliefs at that information set). However, $(b,(f,e),h)$ is not a subgame-perfect equilibrium, because the restriction of this strategy profile to the proper subgame that starts at note $t$ of Player 2, namely $(e,h)$, is not a Nash equilibrium of that subgame: *h* is not a best reply to *e*.

The relationship between the three notions of Nash equilibrium, subgame-perfect equilibrium and weak sequential equilibrium is illustrated in the Venn diagram of Figure 10.10.

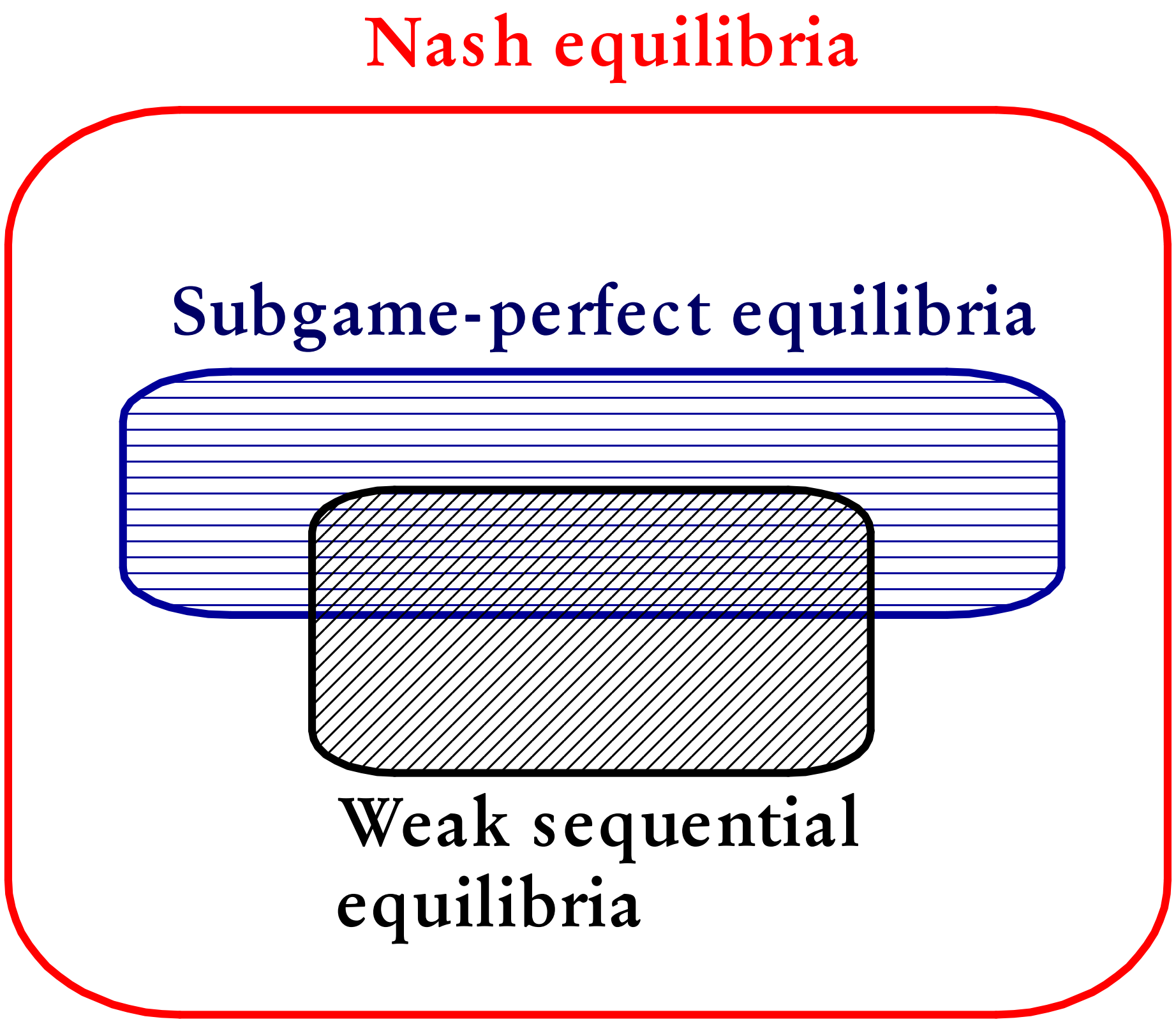

**Figure 10.10**

- Every subgame-perfect equilibrium is a Nash equilibrium. However, there are games in which there is a subgame-perfect equilibrium that is not part of any weak sequential equilibrium. For example, in the game of Figure 10.8, $(b, f, d)$ is a Nash equilibrium which is also subgame-perfect (because there are no proper subgames); however, choice $d$ is strictly dominated and thus is not sequentially rational for Player 3, no matter what beliefs he has at his information set. Thus $(b, f, d)$ cannot be part of a weak sequential equilibrium.

- Every weak sequential equilibrium is a Nash equilibrium. However, as shown in the example of Figure 10.9, there are games in which there is a weak sequential equilibrium whose strategy is not a subgame-perfect equilibrium.

As we will see in Chapter 11, every sequential equilibrium (Definition 11.2) is a weak sequential equilibrium (Remark 11.1) and every finite extensive-form game with cardinal payoffs has at least one sequential equilibrium (Theorem 11.2); hence the following theorem is a corollary of Remark 11.1 and Theorem 11.2.

**Theorem 10.2.** Every finite extensive-form game with cardinal payoffs has at least one weak sequential equilibrium (possibly in mixed strategies).

This is a good time to test your understanding of the concepts introduced in this section, by going through the exercises in Section 10.E.3 of Appendix 10.E at the end of this chapter.

# Appendix 10.E: Exercises

## 10.E.1. Exercises for Section 10.1: Assessments and sequential rationality

The answers to the following exercises are in Appendix S at the end of this chapter.

**Exercise 10.1**. For the game of Figure 10.1, reproduced below, check whether the assessment $\sigma = \begin{pmatrix} a & b & c & h & k & d & e & f & g \\ \frac{1}{8} & \frac{3}{8} & \frac{4}{8} & 1 & 0 & 0 & 1 & 1 & 0 \end{pmatrix}$, $\mu = \begin{pmatrix} x & y & w & z \\ \frac{1}{3} & \frac{2}{3} & \frac{1}{2} & \frac{1}{2} \end{pmatrix}$ is sequentially rational.

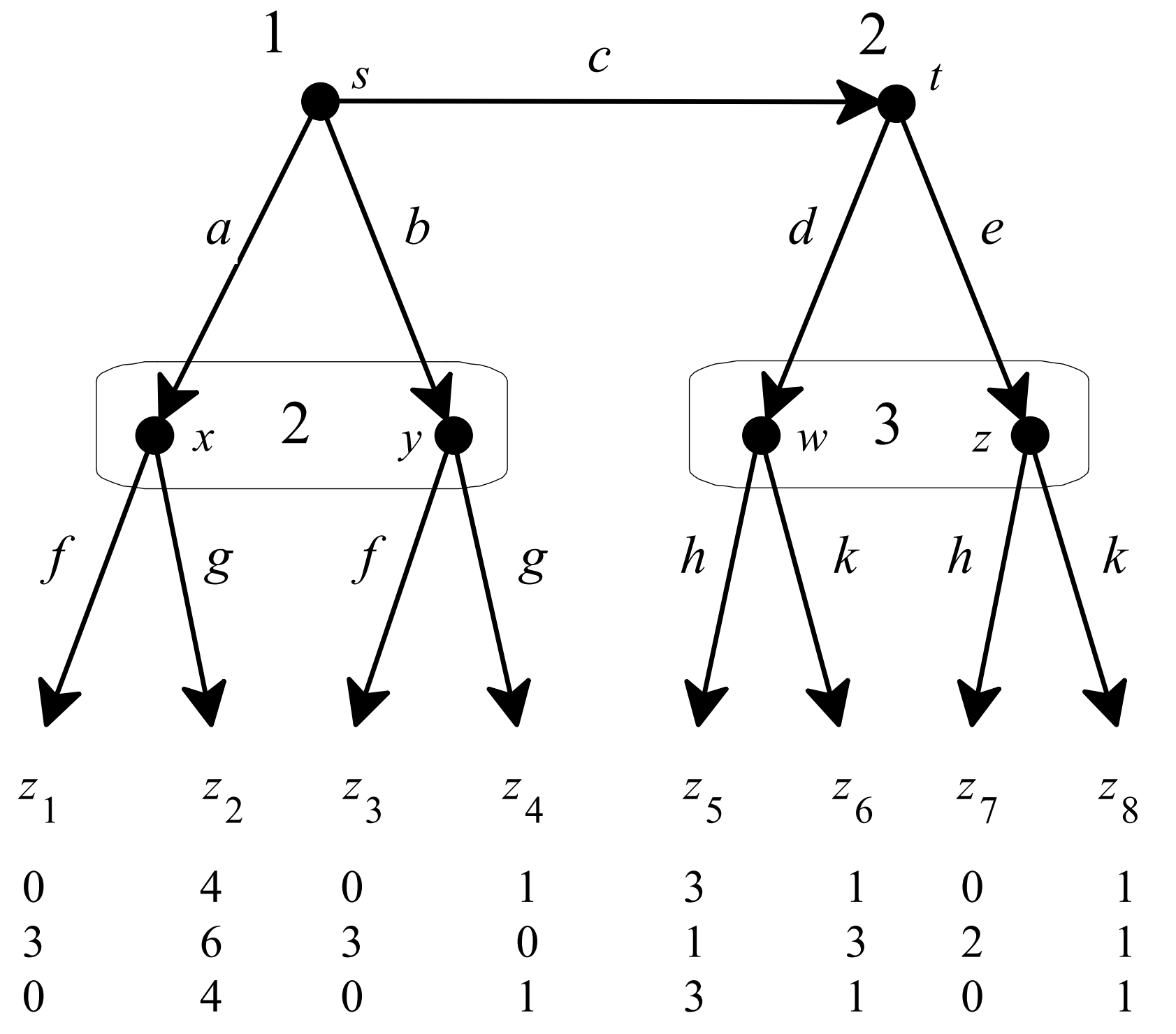

**Exercise 10.2.** Consider the following modification of the game of Figure 10.1 obtained by replacing Player 3 with Player 1 at information set $\{w, z\}$ (note that in the game of Figure 10.1 the payoffs of Players 1 and 3 are identical):

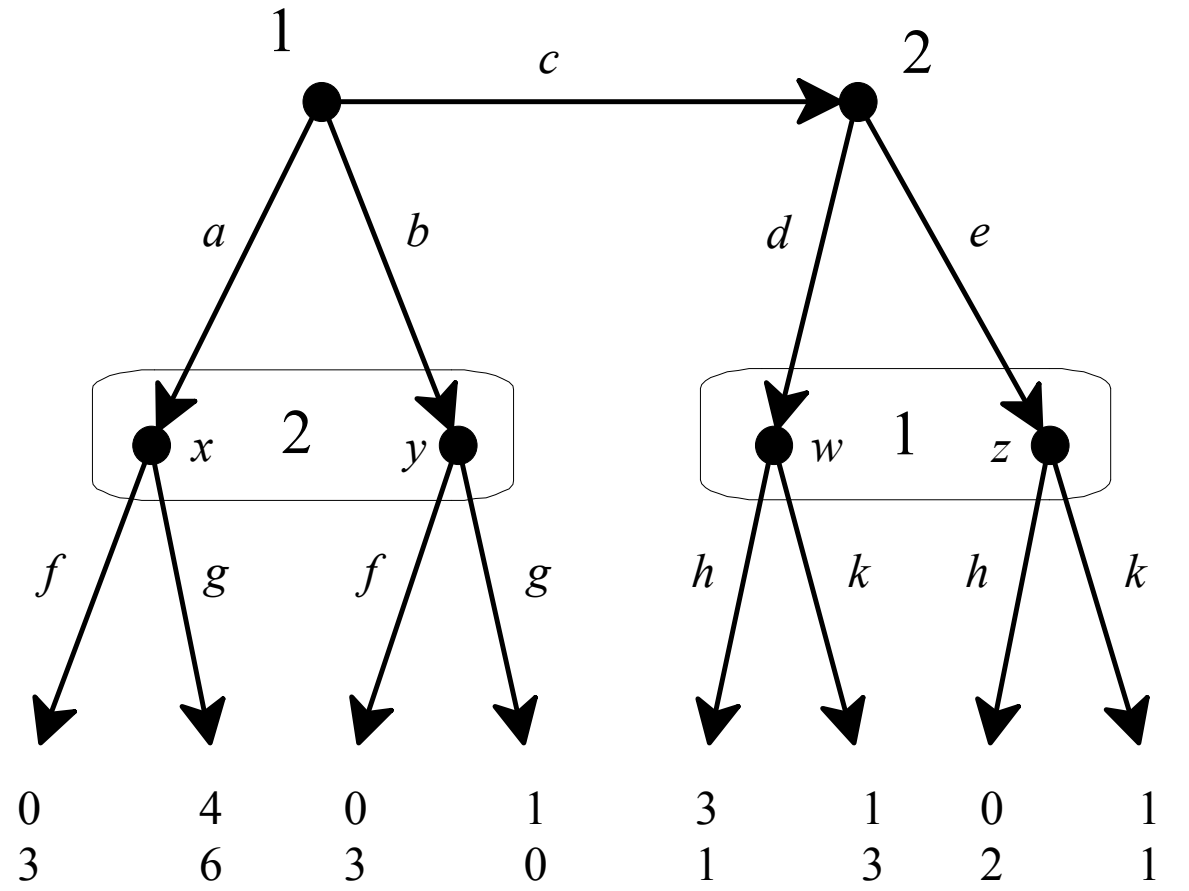

Is the assessment $\sigma = \begin{pmatrix} a & b & c & h & k & d & e & f & g \\ \frac{1}{8} & \frac{3}{8} & \frac{4}{8} & 1 & 0 & 0 & 1 & 1 & 0 \end{pmatrix}$, $\mu = \begin{pmatrix} x & y & w & z \\ \frac{1}{3} & \frac{2}{3} & \frac{1}{2} & \frac{1}{2} \end{pmatrix}$ sequentially rational?

**Exercise 10.3.** Consider the following game:

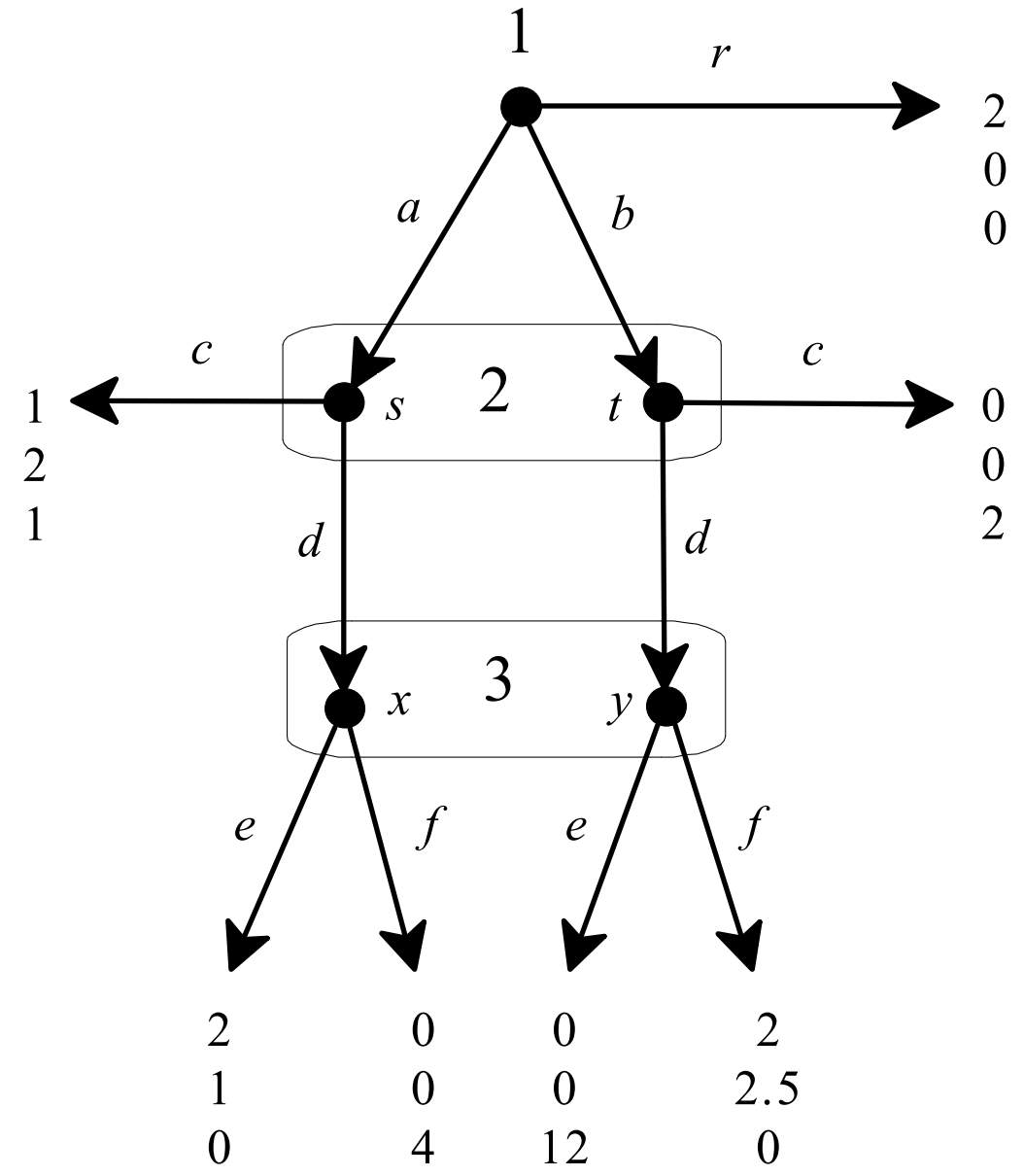

Is the assessment $\sigma = \begin{pmatrix} a & b & r & c & d & e & f \\ 0 & 0 & 1 & \frac{2}{5} & \frac{3}{5} & \frac{1}{3} & \frac{2}{3} \end{pmatrix}$, $\mu = \begin{pmatrix} s & t & x & y \\ \frac{1}{2} & \frac{1}{2} & \frac{3}{4} & \frac{1}{4} \end{pmatrix}$ sequentially rational?

## 10.E.2. Exercises for Section 10.2: Bayesian updating at reached information sets.

The answers to the following exercises are in Appendix S at the end of this chapter.

**Exercise 10.4.** For the game of Figure 10.1 find a system of beliefs $\mu$ such that $(\sigma,\mu)$ satisfies Bayesian updating at reached information sets (Definition 10.4), where $\sigma = \left(\begin{array}{ccc|cc|cc|cc} a & b & c & f & g & d & e & h & k \\ \frac{1}{8} & \frac{3}{8} & \frac{4}{8} & 1 & 0 & \frac{3}{4} & \frac{1}{4} & \frac{1}{5} & \frac{4}{5} \end{array}\right)$.

**Exercise 10.5.** In the following game, let $\sigma = \left(\begin{array}{ccc|cc|cc} a & b & r & c & d & e & f \\ \frac{2}{10} & \frac{1}{10} & \frac{7}{10} & 1 & 0 & \frac{1}{3} & \frac{2}{3} \end{array}\right)$. Find all the systems of beliefs which, combined with $\sigma$, yield assessments that satisfy Bayesian updating at reached information sets.

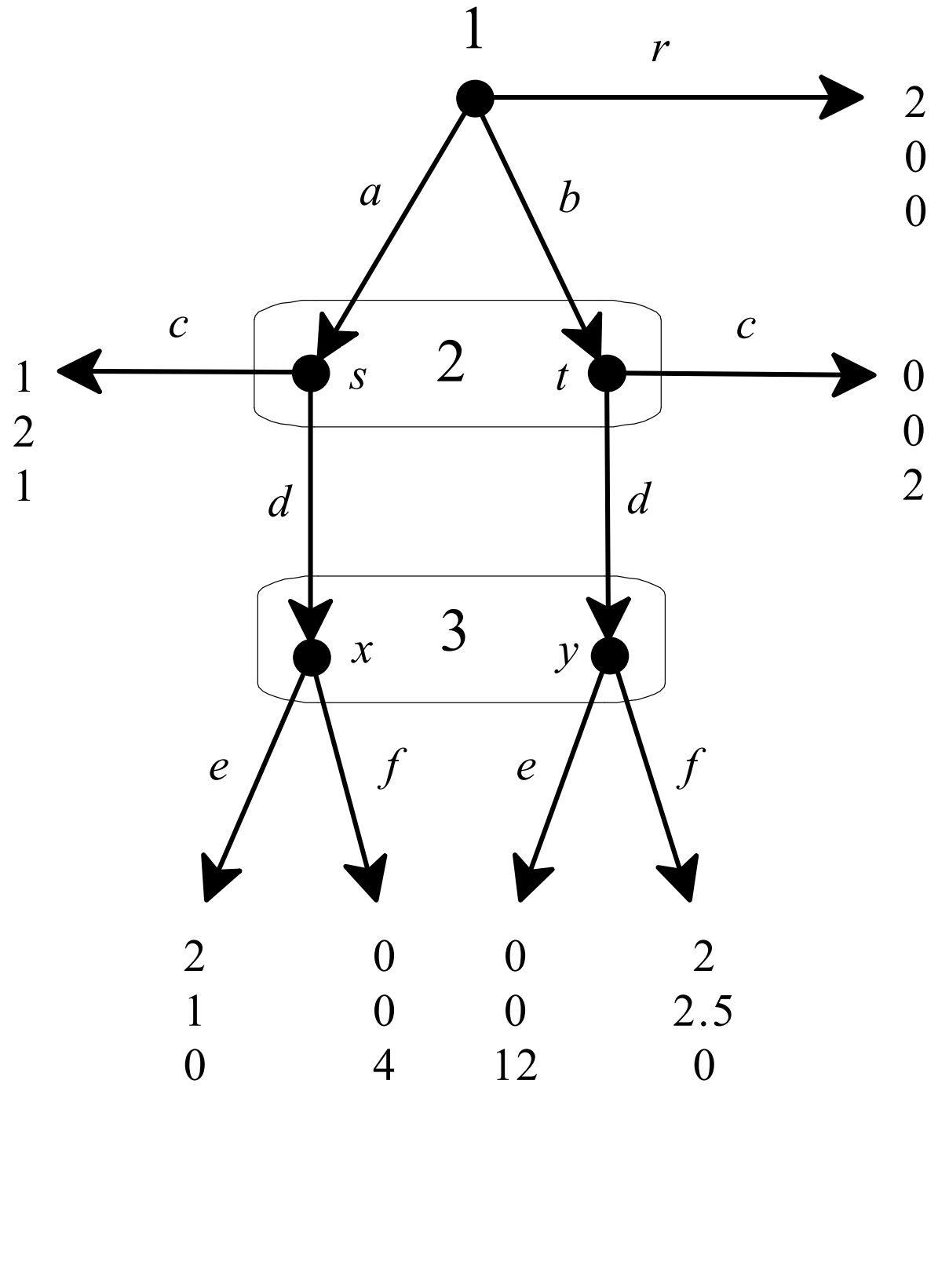

## 10.E.3. Exercises for Section 10.3: Weak Sequential Equilibrium.

The answers to the following exercises are in Appendix S at the end of this chapter.

**Exercise 10.6**. Find all the pure-strategy weak sequential equilibria of the following game:

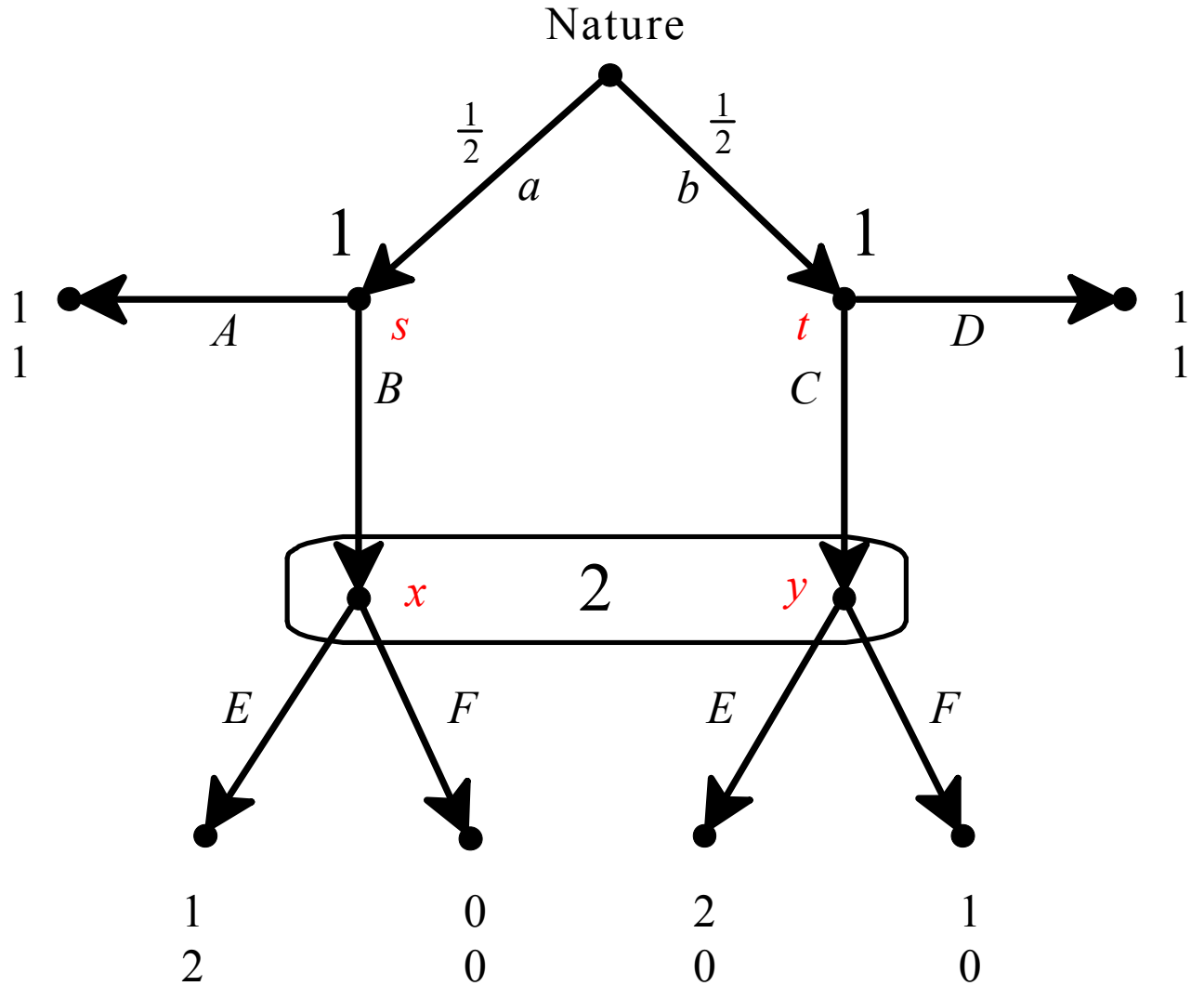

**Exercise 10.7**. In the game of Figure 10.7, reproduced below, is there a pure-strategy weak sequential equilibrium where Player 1 chooses "pass"?

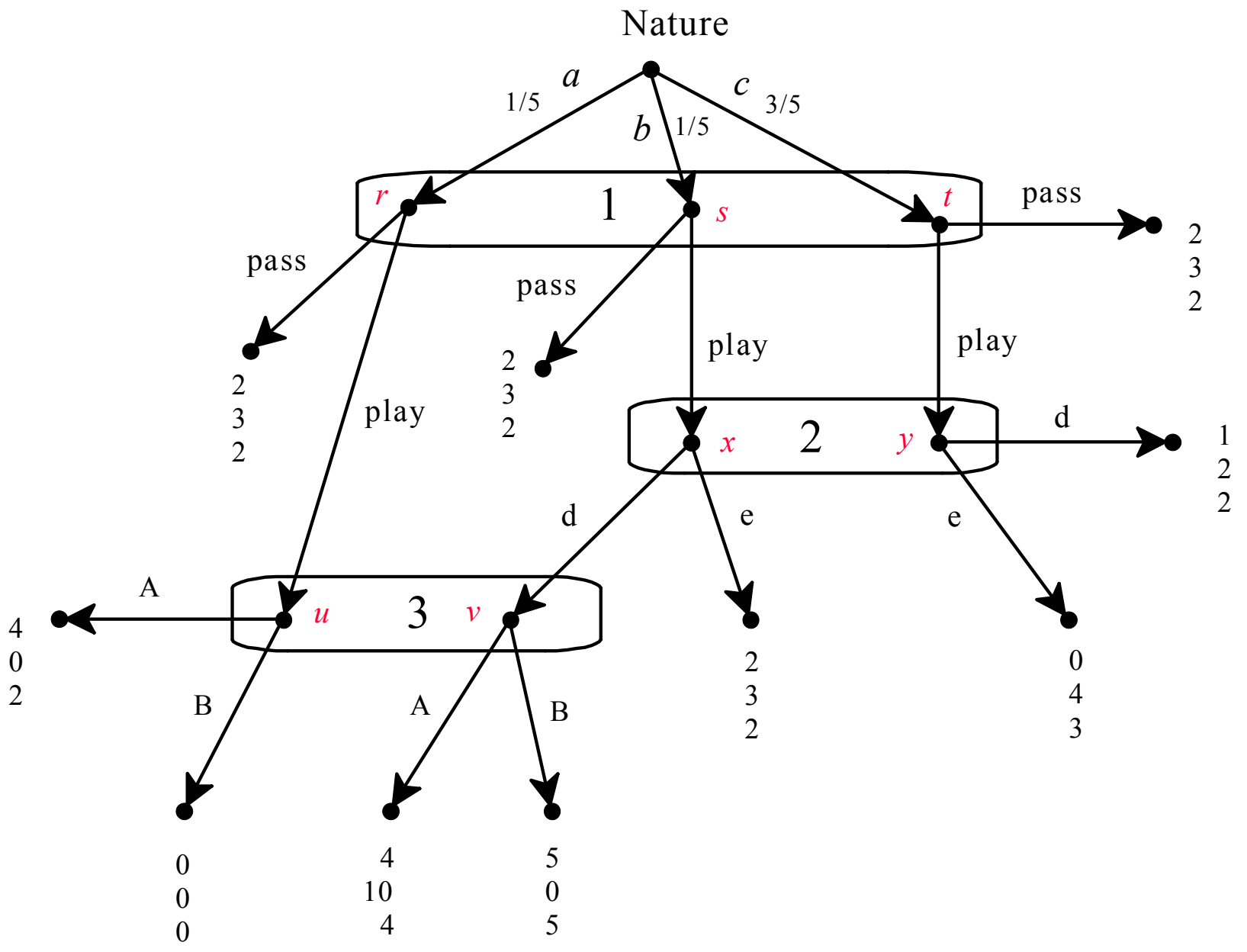

◊◊◊◊◊◊◊◊◊◊◊◊

**Exercise 10.8 (Challenging Question).**

Player 1 can take action $C$ or $L$ and Player 2 can take action $c$ or $f$. The von Neumann-Morgenstern payoffs are as follows:

| | | Player 2: $c$ | Player 2: $f$ |
|---|---|---|---|
| Player 1 | $C$ | 4 , 4 | 6 , 3 |
| | $L$ | 3 , 1 | 5 , 2 |

The game, however, is more complex than the above strategic form. Player 1 moves first and chooses between $C$ and $L$. He then sends an e-mail to Player 2 telling her truthfully what choice he made. However, it is commonly known between the two that a hacker likes to intercept e-mails and change the text. The hacker is a computer program that, with probability $(1-\varepsilon)$, leaves the text unchanged and, with probability $\varepsilon$, changes the sentence "I chose $C$" into the sentence "I chose $L$" and the sentence "I chose $L$" into the sentence "I chose $C$". This is commonly known. The value of $\varepsilon$ is also commonly known. Assume that $\varepsilon \in \left(0, \frac{1}{4}\right)$.

**(a)** Draw the extensive-form game.

**(b)** Find all the pure-strategy weak sequential equilibria.

**(c)** Are there any weak sequential equilibria (pure or mixed) in which Player 2, when he receives a message from Player 1 saying "I chose L" plays $f$ with probability 1?

# Appendix 10.S: Solutions to exercises

**Exercise 10.1**. The assessment $\sigma = \left(\begin{array}{ccc|cc|cc|cc} a & b & c & h & k & d & e & f & g \\ \frac{1}{8} & \frac{3}{8} & \frac{4}{8} & 1 & 0 & 0 & 1 & 1 & 0 \end{array}\right)$,

$\mu = \left(\begin{array}{cc|cc} x & y & w & z \\ \frac{1}{3} & \frac{2}{3} & \frac{1}{2} & \frac{1}{2} \end{array}\right)$ is sequentially rational. In fact,

- at the root, *a* gives Player 1 a payoff of 0 and so do *b* and *c* (given the strategies of Players 2 and 3). Thus any mixture of *a*, *b* and *c* is sequentially rational; in particular, $\begin{pmatrix} a & b & c \\ \frac{1}{8} & \frac{3}{8} & \frac{4}{8} \end{pmatrix}$;
- at Player 2's node *t*, given Player 3's strategy, *d* gives Player 2 a payoff of 1 while *e* gives a payoff of 2; thus *e* is sequentially rational;
- given $\mu$, at information set {*x*,*y*}, *f* gives Player 2 a payoff of $\frac{1}{3}(3) + \frac{2}{3}(3) = 3$ while *g* gives $\frac{1}{3}(6) + \frac{2}{3}(0) = 2$; thus *f* is sequentially rational;
- given $\mu$, at information set {*w*,*z*}, *h* gives Player 3 a payoff of $\frac{1}{2}(3) + \frac{1}{2}(0) = 1.5$, while *k* gives 1; thus *h* is sequentially rational.

**Exercise 10.2**. It might seem that the answer is the same as in the previous exercise, because the calculations at the various information sets remain the same; however, in this game checking sequential rationality at the root involves checking whether, given the strategy of Player 2 (namely $(f, e)$), Player 1 can increase his payoff by changing his *entire strategy*, that is by changing his choices at both the root and at information set {*w*,*z*}. Indeed, if Player 1 changes his strategy from $\left(\begin{array}{ccc|cc} a & b & c & h & k \\ \frac{1}{8} & \frac{3}{8} & \frac{4}{8} & 1 & 0 \end{array}\right)$ to $\left(\begin{array}{ccc|cc} a & b & c & h & k \\ 0 & 0 & 1 & 0 & 1 \end{array}\right)$ his payoff increases from 0 to 1. Hence $\left(\begin{array}{ccc|cc} a & b & c & h & k \\ \frac{1}{8} & \frac{3}{8} & \frac{4}{8} & 1 & 0 \end{array}\right)$ is not sequentially rational at the root.

**Exercise 10.3.** The game under consideration is the following:

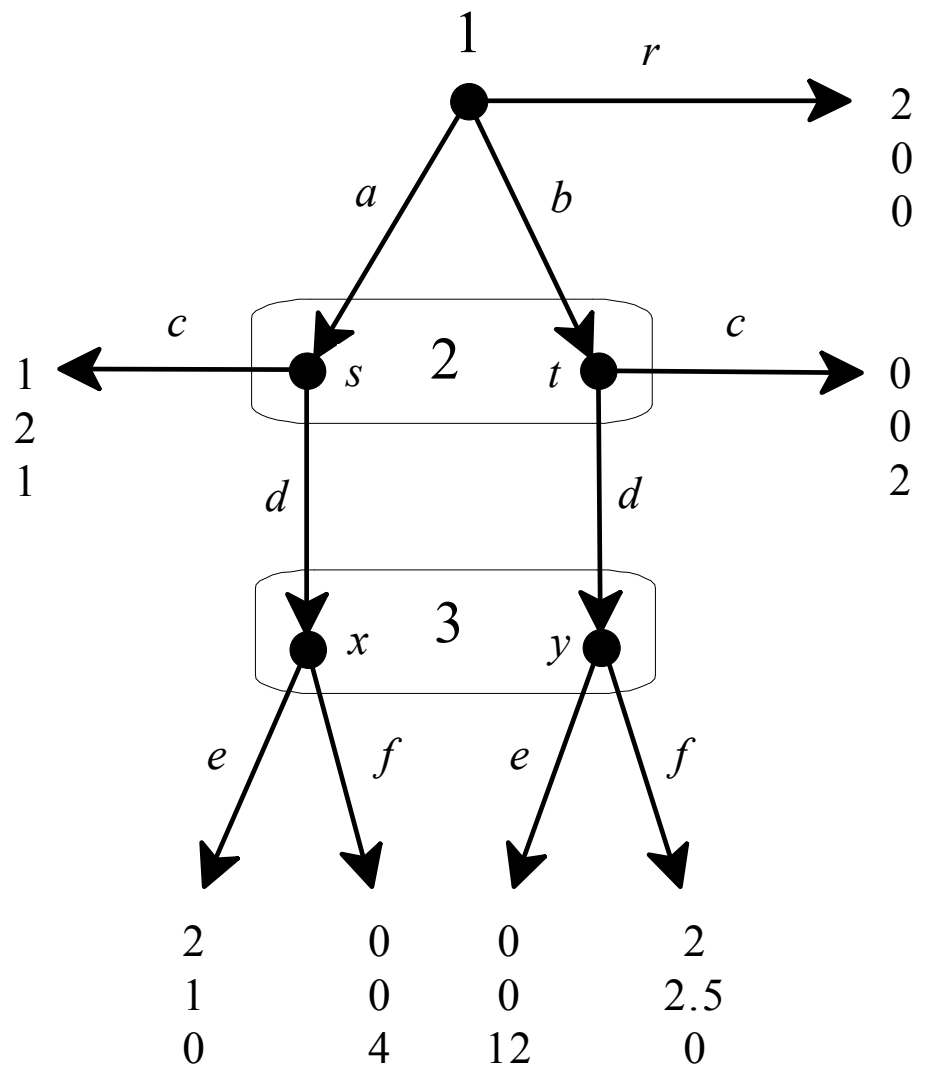

The assessment $\sigma = \begin{pmatrix} a & b & r & c & d & e & f \\ 0 & 0 & 1 & \frac{2}{5} & \frac{3}{5} & \frac{1}{3} & \frac{2}{3} \end{pmatrix}$, $\mu = \begin{pmatrix} s & t & x & y \\ \frac{1}{2} & \frac{1}{2} & \frac{3}{4} & \frac{1}{4} \end{pmatrix}$ is sequentially rational, as the following calculations show.

- At Player 3's information set $\{x,y\}$, $e$ gives Player 3 a payoff of $\frac{3}{4}(0)+\frac{1}{4}(12)=3$ and $f$ a payoff of $\frac{3}{4}(4)+\frac{1}{4}(0)=3$; thus both $e$ and $f$ are optimal and so is any mixture of $e$ and $f$; in particular, the mixture $\begin{pmatrix} e & f \\ \frac{1}{3} & \frac{2}{3} \end{pmatrix}$.

- At Player 2's information set $\{s,t\}$, $c$ gives Player 2 a payoff of $\frac{1}{2}(2)+\frac{1}{2}(0)=1$ and (given the strategy of Player 3) $d$ a payoff of $\frac{1}{2}\left[\frac{1}{3}(1)+\frac{2}{3}(0)\right]+\frac{1}{2}\left[\frac{1}{3}(0)+\frac{2}{3}\left(\frac{5}{2}\right)\right]=1$; thus both $c$ and $d$ are optimal and so is any mixture of $c$ and $d$; in particular the mixture $\begin{pmatrix} c & d \\ \frac{2}{5} & \frac{3}{5} \end{pmatrix}$.

- At the root, $r$ gives Player 1 a payoff of 2 and (given the strategies of Players 2 and 3) $a$ gives a payoff of $\frac{2}{5}(1)+\frac{3}{5}\left[\frac{1}{3}(2)+\frac{2}{3}(0)\right]=\frac{4}{5}$ and $b$ a payoff of $\frac{2}{5}(0)+\frac{3}{5}\left[\frac{1}{3}(0)+\frac{2}{3}(2)\right]=\frac{4}{5}$. Thus $r$ is sequentially rational.

**Exercise 10.4**. The game under consideration is

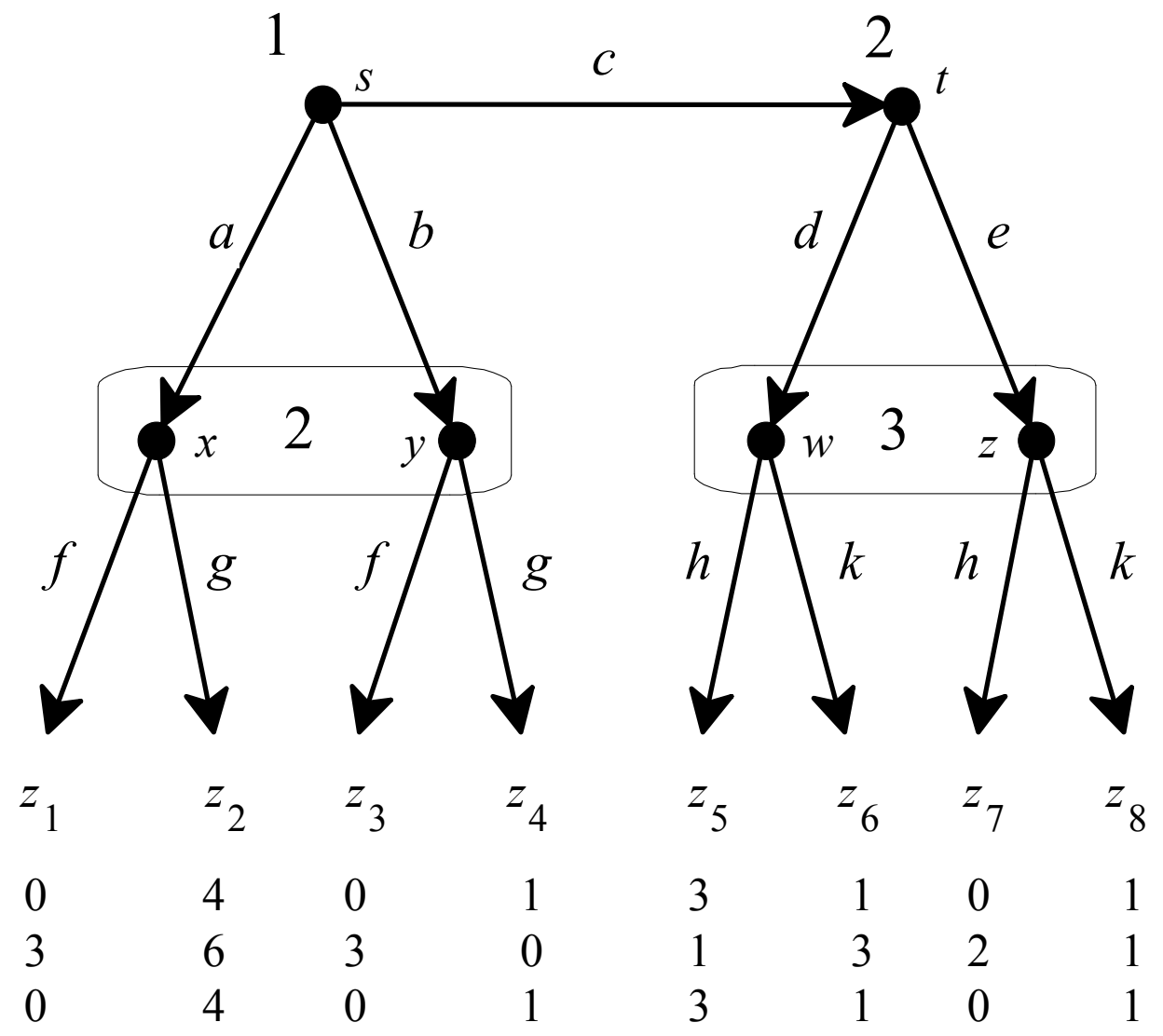

The system of beliefs is $\mu = \left(\begin{array}{cc|cc} x & y & w & z \\ \frac{1}{4} & \frac{3}{4} & \frac{3}{4} & \frac{1}{4} \end{array}\right)$. In fact, we have that $\mathbb{P}_{root,\sigma}(x) = \frac{1}{8}$, $\mathbb{P}_{root,\sigma}(y) = \frac{3}{8}$, $\mathbb{P}_{root,\sigma}(\{x,y\}) = \frac{1}{8} + \frac{3}{8} = \frac{4}{8}$, $\mathbb{P}_{root,\sigma}(w) = \frac{4}{8}\left(\frac{3}{4}\right) = \frac{3}{8}$, $\mathbb{P}_{root,\sigma}(z) = \frac{4}{8}\left(\frac{1}{4}\right) = \frac{1}{8}$ and $\mathbb{P}_{root,\sigma}(\{w,x\}) = \frac{3}{8} + \frac{1}{8} = \frac{4}{8}$. Thus

$$\mu(x) = \frac{\mathbb{P}_{root,\sigma}(x)}{\mathbb{P}_{root,\sigma}(\{x,y\})} = \frac{\frac{1}{8}}{\frac{4}{8}} = \frac{1}{4}$$

$$\mu(y) = \frac{\mathbb{P}_{root,\sigma}(y)}{\mathbb{P}_{root,\sigma}(\{x,y\})} = \frac{\frac{3}{8}}{\frac{4}{8}} = \frac{3}{4},$$

$$\mu(w) = \frac{\mathbb{P}_{root,\sigma}(w)}{\mathbb{P}_{root,\sigma}(\{w,z\})} = \frac{\frac{3}{8}}{\frac{4}{8}} = \frac{3}{4},$$

$$\mu(z) = \frac{\mathbb{P}_{root,\sigma}(z)}{\mathbb{P}_{root,\sigma}(\{w,z\})} = \frac{\frac{1}{8}}{\frac{4}{8}} = \frac{1}{4}.$$

**Exercise 10.5.** The game under consideration is

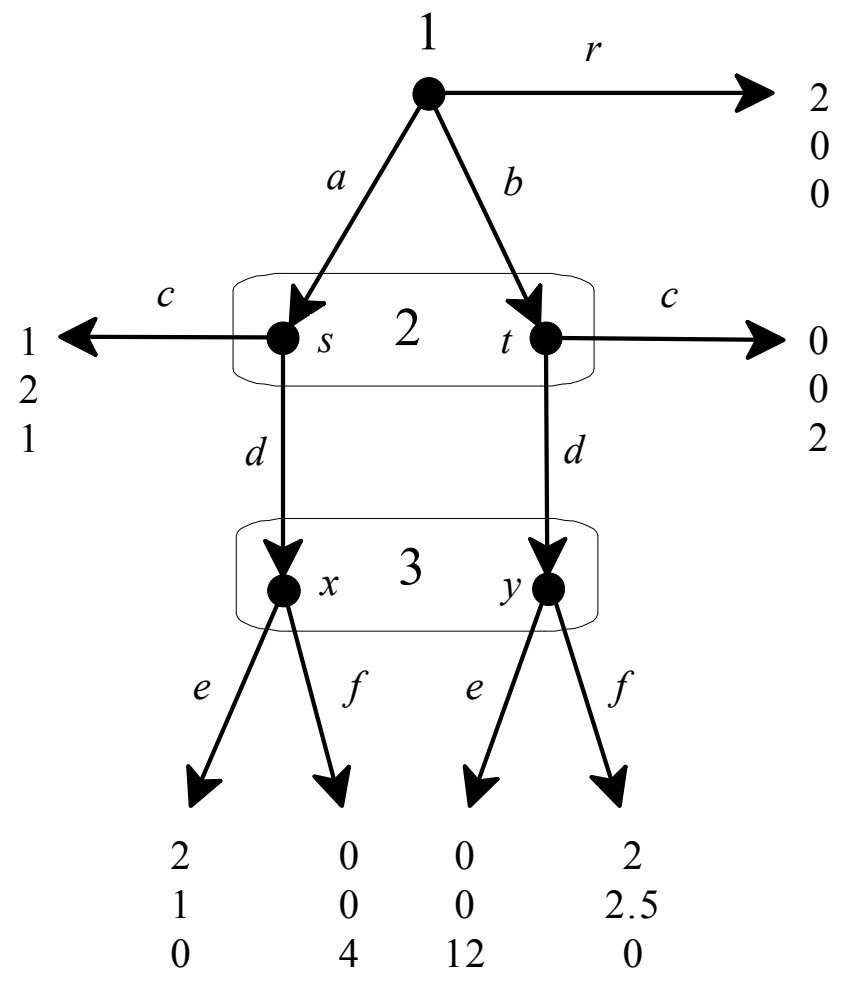

Let $\sigma = \begin{pmatrix} a & b & r & c & d & e & f \\ \frac{2}{10} & \frac{1}{10} & \frac{7}{10} & 1 & 0 & \frac{1}{3} & \frac{2}{3} \end{pmatrix}$. Since only information set $\{s,t\}$ is reached by $\sigma$, no restrictions are imposed on the beliefs at information set $\{x,y\}$. Thus, for every $p$ such that $0 \leq p \leq 1$, the system of beliefs $\mu = \begin{pmatrix} s & t & x & y \\ \frac{2}{3} & \frac{1}{3} & p & 1-p \end{pmatrix}$ combined with $\sigma$ yields an assessment that satisfies Bayesian updating at reached information sets.

**Exercise 10.6.** The game under consideration is

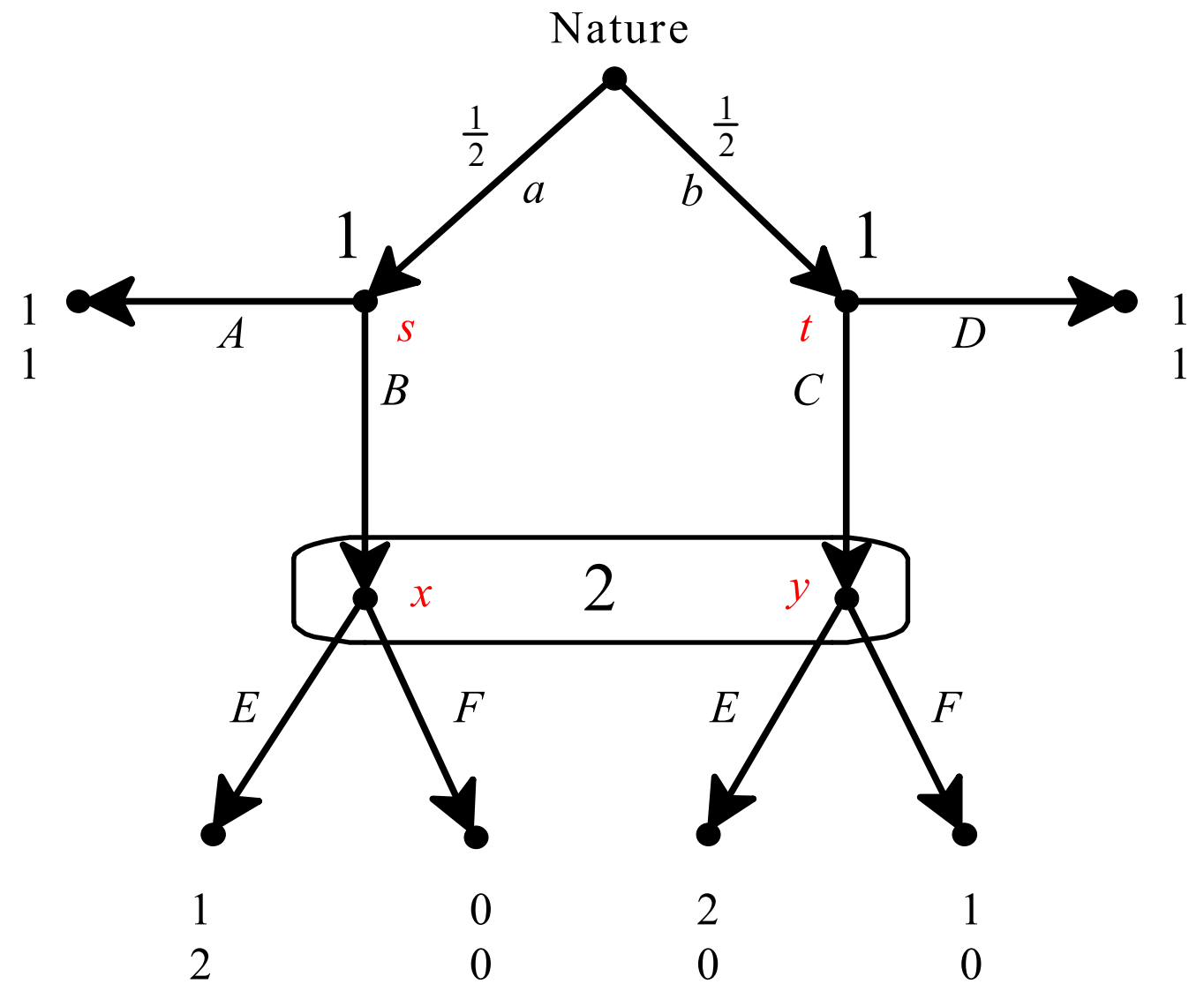

We restrict attention to pure strategies. First of all, note that – for Player 1 – $A$ is sequentially rational no matter what strategy Player 2 chooses and, similarly, $C$ is sequentially rational no matter what strategy Player 2 chooses. The strategy of Player 1 determines the beliefs of Player 2 at her information set $\{x, y\}$. Let us consider the four possibilities (recall that $S_1 = \{(A,C),(B,D),(B,C),(A,D)\}$).

- If Player 1's strategy is $(A,C)$, then Player 2's information set $\{x, y\}$ is reached with positive probability and the only beliefs that are consistent with Bayesian updating are $\begin{pmatrix} x & y \\ 0 & 1 \end{pmatrix}$, so that both $E$ and $F$ are sequentially rational for Player 2. By our preliminary observation it follows that $((A,C),E)$ with beliefs $\begin{pmatrix} x & y \\ 0 & 1 \end{pmatrix}$ is a weak sequential equilibrium and so is $((A,C),F)$ with beliefs $\begin{pmatrix} x & y \\ 0 & 1 \end{pmatrix}$.
- If Player 1's strategy is $(B,D)$, then Player 2's information set $\{x, y\}$ is reached with positive probability and the only beliefs that are consistent with Bayesian updating are $\begin{pmatrix} x & y \\ 1 & 0 \end{pmatrix}$, so that – by sequential rationality – Player 2 must choose $E$. However, if Player 2's strategy is $E$ then at node $t$ it is not sequentially rational for Player 1 to choose $D$. Thus there is no pure-strategy weak sequential equilibrium where Player 1's strategy is $(B,D)$.
- If Player 1's strategy is $(B,C)$, then Player 2's information set $\{x, y\}$ is reached (with probability 1) and the only beliefs that are consistent with Bayesian updating are $\begin{pmatrix} x & y \\ \frac{1}{2} & \frac{1}{2} \end{pmatrix}$. Given these beliefs, $E$ is the only sequentially rational choice for Player 2 (her payoff from playing $E$ is $\frac{1}{2}(2) + \frac{1}{2}(0) = 1$, while her payoff from playing $F$ is 0). Thus $((B,C),E)$ with beliefs $\begin{pmatrix} x & y \\ \frac{1}{2} & \frac{1}{2} \end{pmatrix}$ is a weak sequential equilibrium.
- If Player 1's strategy is $(A,D)$, then Player 2's information set $\{x, y\}$ is not reached and thus, according to the notion of weak sequential equilibrium, any beliefs are allowed there. In order for $D$ to be sequentially rational for Player 1 it must be that Player 2's pure strategy is $F$. In order for $F$ to be sequentially rational for Player 2, her beliefs must be $\begin{pmatrix} x & y \\ 0 & 1 \end{pmatrix}$. Thus $((A,D),F)$ with beliefs $\begin{pmatrix} x & y \\ 0 & 1 \end{pmatrix}$ is a weak sequential equilibrium.

Summarizing, there are four pure-strategy weak sequential equilibria:

(1) $((A,C),E)$ with beliefs $\begin{pmatrix} x & y \\ 0 & 1 \end{pmatrix}$, (2) $((A,C),F)$ with beliefs $\begin{pmatrix} x & y \\ 0 & 1 \end{pmatrix}$, (3) $((B,C),E)$ with beliefs $\begin{pmatrix} x & y \\ \frac{1}{2} & \frac{1}{2} \end{pmatrix}$ and (4) $((A,D),F)$ with beliefs $\begin{pmatrix} x & y \\ 0 & 1 \end{pmatrix}$.

**Exercise 10.7.** The game under consideration is

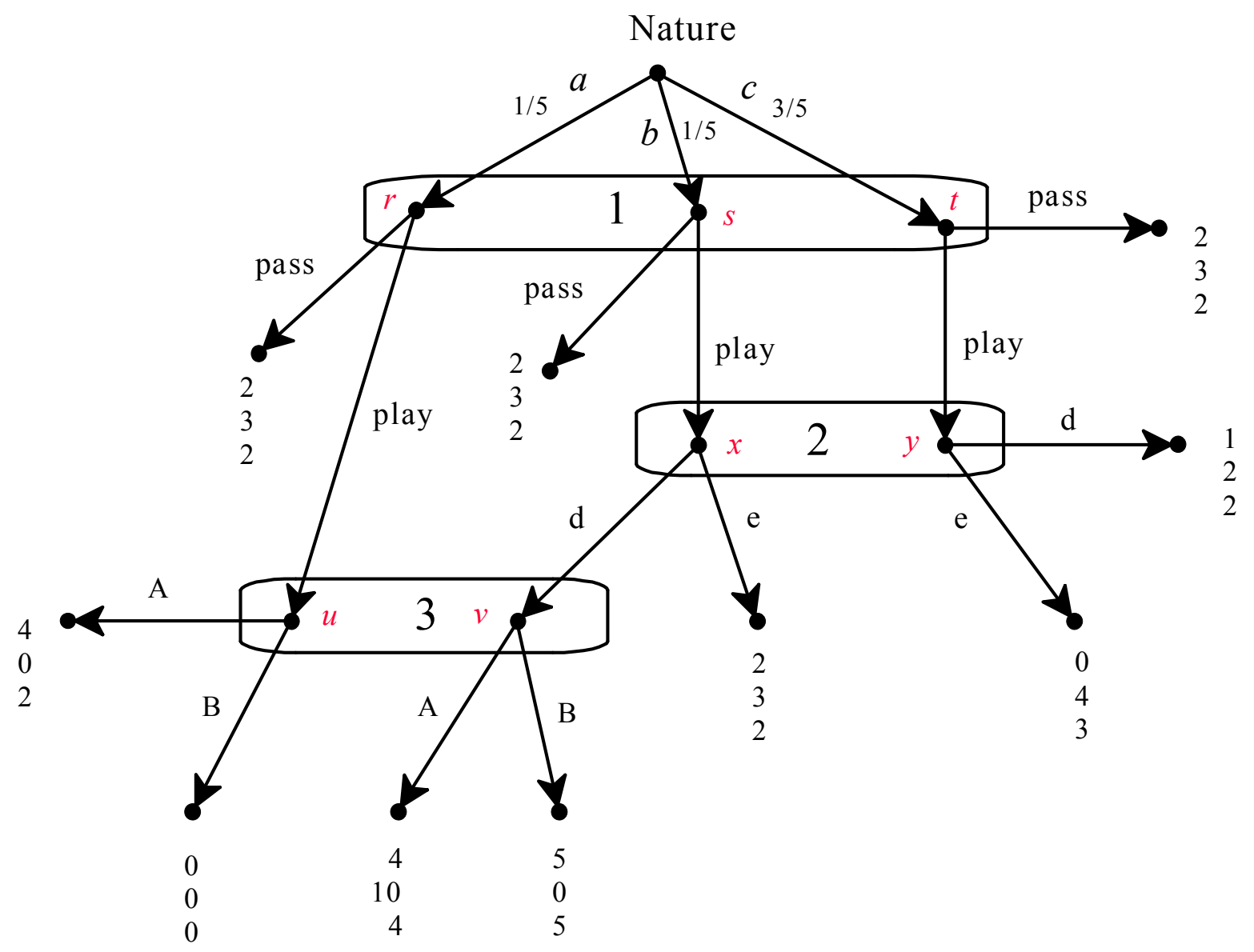

When Player 1's strategy is "pass" then it is much easier to construct a weak sequential equilibrium, because there are no restrictions on the beliefs at the information sets of Players 2 and 3. For example, we can choose beliefs $\begin{pmatrix} x & y \\ 0 & 1 \end{pmatrix}$ for Player 2, which make $e$ the only sequentially rational choice, and beliefs $\begin{pmatrix} u & v \\ 0 & 1 \end{pmatrix}$ for Player 3, which make $B$ the only sequentially rational choice. It only remains to check sequential rationality of "pass": if Player 1 chooses "pass" he gets a payoff of 2, while if he chooses "play" he gets a payoff of $\frac{1}{5}(0)+\frac{1}{5}(2)+\frac{3}{5}(0)=\frac{2}{5}<2$, so that "pass" is indeed the better choice. Thus we have found the following weak sequential equilibrium (note, however, that there are several weak sequential equilibria): $\sigma=\begin{pmatrix} pass & play & d & e & A & B \\ 1 & 0 & 0 & 1 & 0 & 1 \end{pmatrix}$,

$\mu=\begin{pmatrix} r & s & t & x & y & u & v \\ \frac{1}{5} & \frac{1}{5} & \frac{3}{5} & 0 & 1 & 0 & 1 \end{pmatrix}$.

**Exercise 10.8.** [Challenging question].

**(a)** The extensive form is as follows:

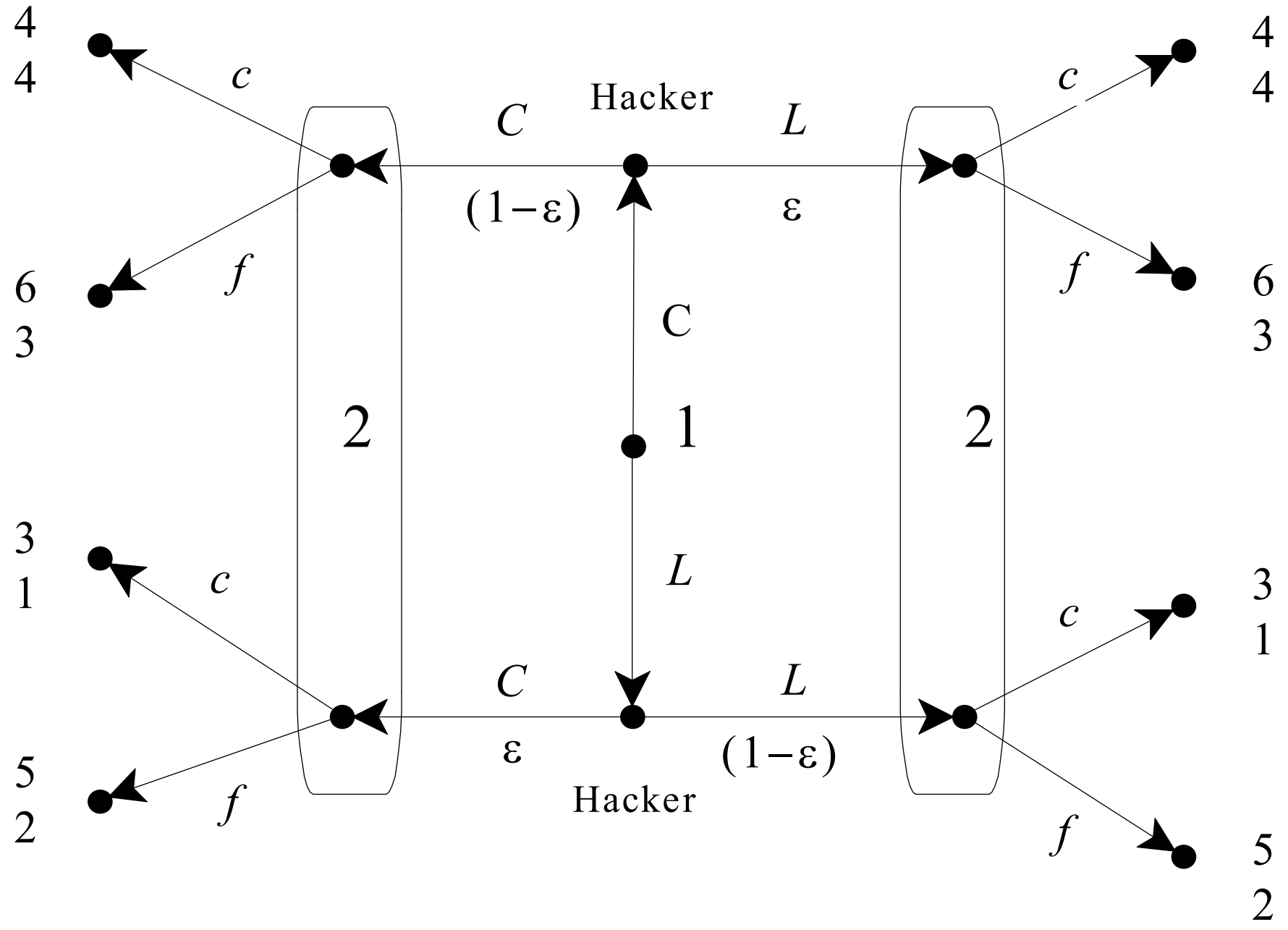

**(b)** Let us consider the pure strategies of Player 1. If Player 1 plays $L$ then Player 2 assigns probability 1 to the bottom node of each information set and responds with $f$ with probability 1, but this makes Player 1 want to deviate to $C$ to get a payoff of 6.[6] Thus there is no pure-strategy sequential equilibrium where Player 1 plays $L$. On the other hand, if Player 1 plays $C$, then Player 2 assigns probability 1 to the top node of each information set and thus sequential rationality requires her to respond with $c$ at each information set, which makes playing $C$ optimal for Player 1. Thus $(C,(c,c))$, with beliefs that assign probability one to the top node of each information set, is the only pure-strategy weak sequential equilibrium.

[6] In other words, if Player 1 plays $L$, then Bayesian updating requires Player 2 to assign probability 1 to the bottom node of each information set and then sequential rationality requires Player 2 to play $f$ at each information set, so that $(L,(f,f))$ is the only candidate for a weak sequential equilibrium where Player 1 plays $L$. But $L$ is not a best reply to $(f,f)$ and thus $(L,(f,f))$ is not a Nash equilibrium and hence cannot be part of a weak sequential equilibrium.

**(c)** Suppose that Player 2 plays $f$ after reading the message "I chose $L$", that is, at her information set on the right. We know from the argument in part (b) that there are no equilibria of this kind in which Player 1 chooses a pure strategy, so Player 1 must be playing a mixed strategy. For him to be willing to do so, he must receive the same payoff from playing $C$ or $L$. If we let $p$ be the probability with which Player 2 plays $c$ if she receives the message "I chose $C$", then 1 is indifferent when (the symbol $\Leftrightarrow$ stands for 'if and only if')

$$\underbrace{(1-\varepsilon)\left[4p+6(1-p)\right]+\varepsilon 6}_{=\pi_1(C)} = \underbrace{\varepsilon\left[3p+5(1-p)\right]+(1-\varepsilon)5}_{=\pi_1(L)} \Leftrightarrow p=\frac{1}{2-4\varepsilon}$$

Since $\varepsilon \in \left(0,\frac{1}{4}\right)$, $p \in \left(\frac{1}{2},1\right)$, so Player 2 randomizes after reading "I chose $C$". For Player 2 to be willing to do this, she must be indifferent between $c$ and $f$ in this event. Let $q \in (0,1)$ be the probability with which Player 1 plays $C$ (so that $1-q$ is the probability of $L$); then Bayesian updating requires that Player 2 assign probability $\frac{(1-\varepsilon)q}{(1-\varepsilon)q+\varepsilon(1-q)}$ to the top node of her information set on the left and probability $\frac{\varepsilon(1-q)}{(1-\varepsilon)q+\varepsilon(1-q)}$ to the bottom node. Then, for Player 2, the expected payoff from playing $c$ at the information set on the left is

$\pi_2(c)=\frac{(1-\varepsilon)q}{(1-\varepsilon)q+\varepsilon(1-q)}4+\frac{\varepsilon(1-q)}{(1-\varepsilon)q+\varepsilon(1-q)}1$ and the expected payoff from playing $f$ is $\pi_2(f)=\frac{(1-\varepsilon)q}{(1-\varepsilon)q+\varepsilon(1-q)}3+\frac{\varepsilon(1-q)}{(1-\varepsilon)q+\varepsilon(1-q)}2$. Player 2 is indifferent if these two are equal, that is if $4q(1-\varepsilon)+(1-q)\varepsilon=3q(1-\varepsilon)+2(1-q)\varepsilon$, which is true if and only if $q=\varepsilon$.

We have now specified behavior at all information sets. To ensure that the specified behavior is an equilibrium, we need to check that *f* is optimal for Player 2 if she receives the message "I chose *L*". This will be true if

$$\pi_2(c\,|\,\text{"I chose L"}) \le \pi_2(f\,|\,\text{"I chose L"}) \underbrace{\Leftrightarrow}_{\text{since } \varepsilon = q}$$
$$4q\varepsilon + 1(1-q)(1-\varepsilon) \le 3q\varepsilon + 2(1-q)(1-\varepsilon) \Leftrightarrow$$
$$4\varepsilon^2 + (1-\varepsilon)^2 \le 3\varepsilon^2 + 2(1-\varepsilon)^2 \;\Leftrightarrow\; \varepsilon \le \tfrac{1}{2}.$$

Since we have assumed that $\varepsilon < \frac{1}{4}$, Player 2 strictly prefers to play *f* after receiving the message "I chose *L*". Thus we have constructed the following weak sequential equilibrium:

- Behavior strategy of Player 1: $\begin{pmatrix} C & L \\ \varepsilon & 1-\varepsilon \end{pmatrix}$.
- Behavior strategy of Player 2: at the information set on the left (where she receives the message "I chose *C*"): $\begin{pmatrix} c & f \\ \frac{1}{2-4\varepsilon} & \frac{1-4\varepsilon}{2-4\varepsilon} \end{pmatrix}$, at the information set on the right (where she receives the message "I chose *L*"): $\begin{pmatrix} c & f \\ 0 & 1 \end{pmatrix}$.
- Player 2's beliefs at the information set on the left assign probability $\frac{\varepsilon(1-\varepsilon)}{\varepsilon(1-\varepsilon)+\varepsilon(1-\varepsilon)} = \frac{1}{2}$ to the top node and probability $\frac{1}{2}$ to the bottom node and her beliefs at the information set on the right assign probability $\frac{\varepsilon^2}{\varepsilon^2+(1-\varepsilon)^2}$ to the top node and probability $\frac{(1-\varepsilon)^2}{\varepsilon^2+(1-\varepsilon)^2}$ to the bottom node.

Chapter

11

# Sequential Equilibrium

## 11.1 Consistent assessments

The stated goal at the beginning of Chapter 10 was to seek a refinement of the notion of subgame-perfect equilibrium that would rule out strictly dominated choices at unreached information sets. The notion of weak sequential equilibrium achieved the goal of ruling out strictly dominated choices, but turned out not to be a refinement of subgame-perfect equilibrium. As shown at the end of Section 10.3 (Chapter 10), it is possible for the strategy profile in a weak sequential equilibrium not to be a subgame-perfect equilibrium. The reason for this is that the only restriction on beliefs that is incorporated in the notion of weak sequential equilibrium is Bayesian updating at reached information sets. At an information set that is not reached by the strategy profile under consideration any beliefs whatsoever are allowed, even if those beliefs are at odds with the strategy profile.

To see this, consider the game of Figure 11.1 and the assessment consisting of the pure-strategy profile $\sigma = (c, d, f)$ (highlighted as double edges) and the system of beliefs that assigns probability 1 to node $u$: $\begin{pmatrix} s & t & u \\ 0 & 0 & 1 \end{pmatrix}$.

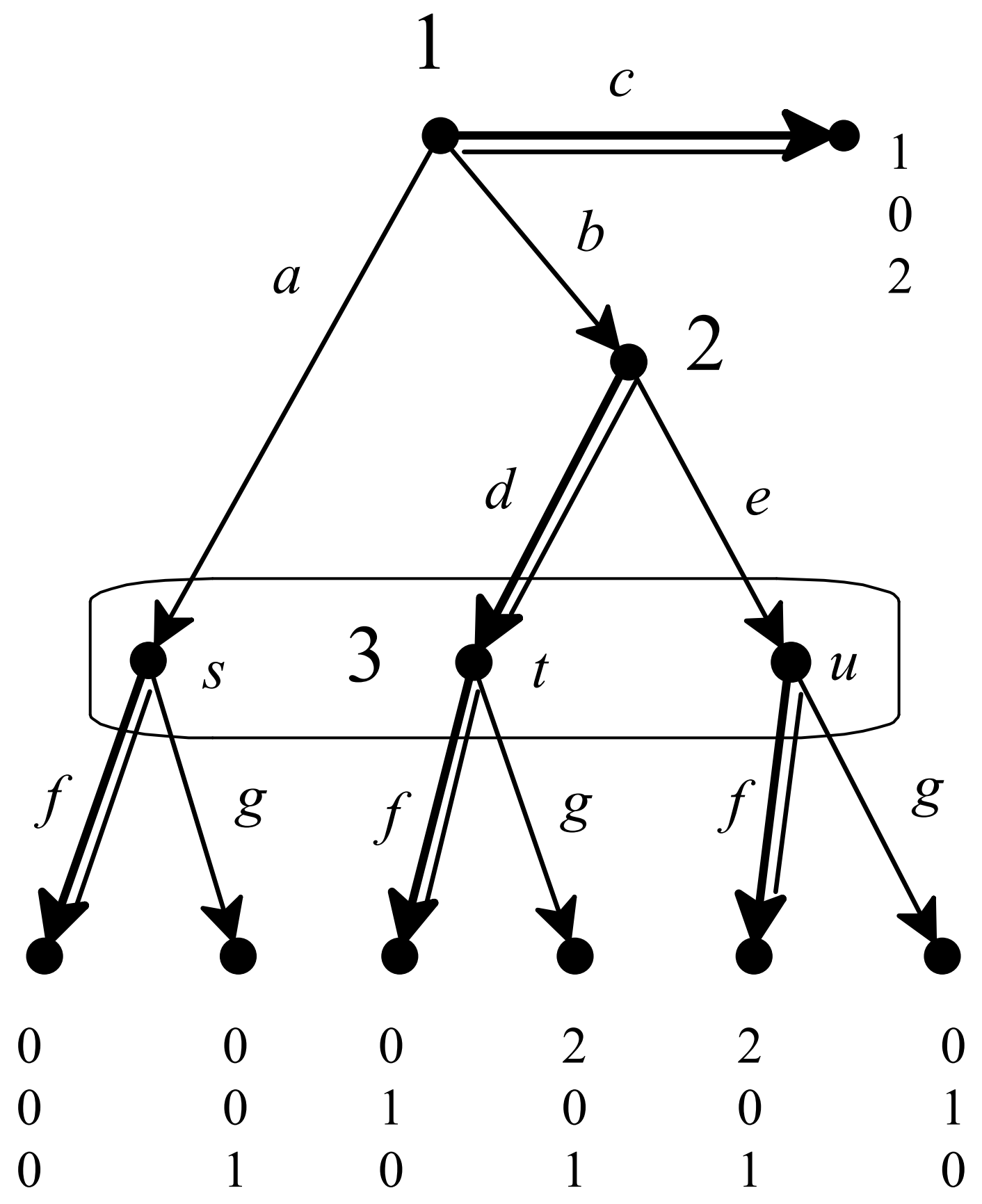

**Figure 11.1**

Since Player 3's information set $\{s,t,u\}$ is not reached by $\sigma$, Bayesian updating imposes no restrictions on beliefs at that information set. However, attaching probability 1 to node $u$ is at odds with $\sigma$, because in order for node $u$ to be reached the play must have gone through Player 2's node and there, according to $\sigma$, Player 2 should have played $d$ with probability 1, making it impossible for node $u$ to be reached.

Thus we need to impose some restrictions on beliefs to ensure that they are consistent with the strategy profile with which they are paired (in the assessment under consideration). At reached information sets this is achieved by requiring Bayesian updating, but so far we have imposed no restriction on beliefs at unreached information sets. We want these restrictions to be "just like Bayesian updating". Kreps and Wilson (1982) proposed a restriction on beliefs that they called *consistency*. To motivate it, note that if $\sigma$ is a *completely mixed* strategy profile (in the sense that $\sigma(a) > 0$, for every choice $a$) then the issue disappears, because every information set is reached with positive probability and Bayesian updating yields unique beliefs at every information set. For example, in the game of Figure 11.1, if Player 1 uses the completely mixed strategy $\begin{pmatrix} a & b & c \\ p_a & p_b & 1-p_a-p_b \end{pmatrix}$ with $p_a, p_b \in (0,1)$ and $p_a + p_b < 1$ and Player 2 uses the completely mixed strategy $\begin{pmatrix} d & e \\ p_d & 1-p_d \end{pmatrix}$ with $p_d \in (0,1)$ then, by Bayesian updating, Player 3's beliefs must be $\mu(s) = \dfrac{p_a}{p_a + p_b p_d + p_b(1-p_d)}$, $\mu(t) = \dfrac{p_b p_d}{p_a + p_b p_d + p_b(1-p_d)}$ and $\mu(u) = \dfrac{p_b(1-p_d)}{p_a + p_b p_d + p_b(1-p_d)}$.

In the case of a completely mixed strategy profile $\sigma$, it is clear what it means for a system of beliefs $\mu$ to be consistent with the strategy profile $\sigma$: $\mu$ must be the unique system of beliefs obtained from $\sigma$ by applying Bayesian updating. What about assessments $(\sigma,\mu)$ where $\sigma$ is such that some information sets are not reached? How can we decide, in such cases, whether $\mu$ is consistent with $\sigma$? Kreps and Wilson (1982) proposed the following criterion: the assessment $(\sigma,\mu)$ is consistent if there is a completely mixed strategy profile $\sigma'$ which is arbitrarily close to $\sigma$ and whose associated unique system of beliefs $\mu'$ (obtained by applying Bayesian updating) is arbitrarily close to $\mu$. In mathematics "arbitrary closeness" is captured by the notion of limit.

**Definition 11.1.** Given an extensive game, an assessment $(\sigma,\mu)$ is *consistent* if there is a sequence of completely mixed strategy profiles $\langle \sigma_1, \sigma_2, ..., \sigma_n, ... \rangle$ such that:
(1) the sequence converges to $\sigma$ as $n$ tends to infinity, that is, $\lim_{n\to\infty} \sigma_n = \sigma$, and
(2) letting $\mu_n$ be the unique system of beliefs obtained from $\sigma_n$ by using Bayesian updating, the sequence $\langle \mu_1, \mu_2, ..., \mu_n, ... \rangle$ converges to $\mu$ as $n$ tends to infinity, that is, $\lim_{n\to\infty} \mu_n = \mu$.

For example, consider the extensive form of Figure 11.2 below and the following assessment: $\sigma = (c,d,f)$, $\mu = \begin{pmatrix} s & t \\ \frac{3}{8} & \frac{5}{8} \end{pmatrix}$.

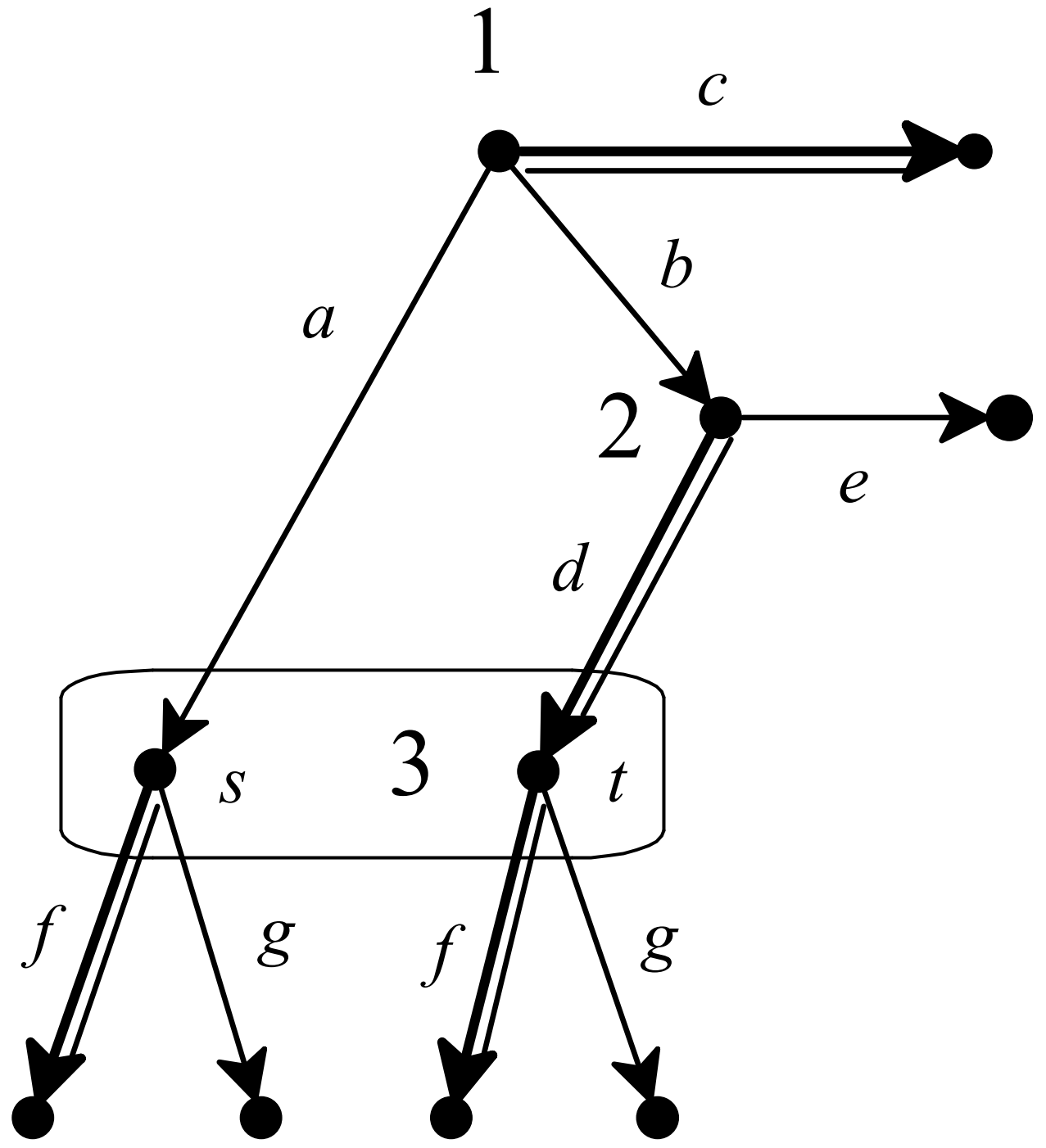

**Figure 11.2**

The assessment $\sigma = (c, d, f),\ \mu = \begin{pmatrix} s & t \\ \frac{3}{8} & \frac{5}{8} \end{pmatrix}$ is consistent. To see this, let $\sigma_n = \left(\begin{array}{ccc|cc|cc} a & b & c & d & e & f & g \\ \frac{3}{n} & \frac{5}{n} & 1-\frac{8}{n} & 1-\frac{1}{n} & \frac{1}{n} & 1-\frac{1}{n} & \frac{1}{n} \end{array}\right)$. Then, as $n$ tends to infinity, all of $\frac{3}{n}, \frac{5}{n}, \frac{1}{n}$ tend to 0 and both $1-\frac{8}{n}$ and $1-\frac{1}{n}$ tend to 1 so that $\lim_{n\to\infty} \sigma_n =$ $\left(\begin{array}{ccc|cc|cc} a & b & c & d & e & f & g \\ 0 & 0 & 1 & 1 & 0 & 1 & 0 \end{array}\right) = \sigma$. Furthermore, $\mu_n(s) = \dfrac{\frac{3}{n}}{\frac{3}{n} + \frac{5}{n}\left(1-\frac{1}{n}\right)} = \dfrac{3}{8-\frac{5}{n}}$, which tends to $\frac{3}{8}$ as $n$ tends to infinity and $\mu_n(t) = \dfrac{\frac{5}{n}\left(1-\frac{1}{n}\right)}{\frac{3}{n} + \frac{5}{n}\left(1-\frac{1}{n}\right)} = \dfrac{5-\frac{5}{n}}{8-\frac{5}{n}}$, which tends to $\frac{5}{8}$ as $n$ tends to infinity, so that $\lim_{n\to\infty} \mu_n = \begin{pmatrix} s & t \\ \frac{3}{8} & \frac{5}{8} \end{pmatrix} = \mu$.

The notion of consistent assessment $(\sigma, \mu)$ was meant to capture an extension of the requirement of Bayesian updating that would apply also to information sets that have zero probability of being reached (when the play of the game is according to the strategy profile $\sigma$). Definition 11.1 is rather technical and not easy to apply. Showing that an assessment is consistent requires displaying an appropriate sequence and showing that it converges to the given assessment. This is relatively easy as compared to the considerably more difficult task of proving that an assessment is *not* consistent: this requires showing that every possible sequence that one could construct would not converge. To see the kind of reasoning involved, consider the extensive form of Figure 11.3.

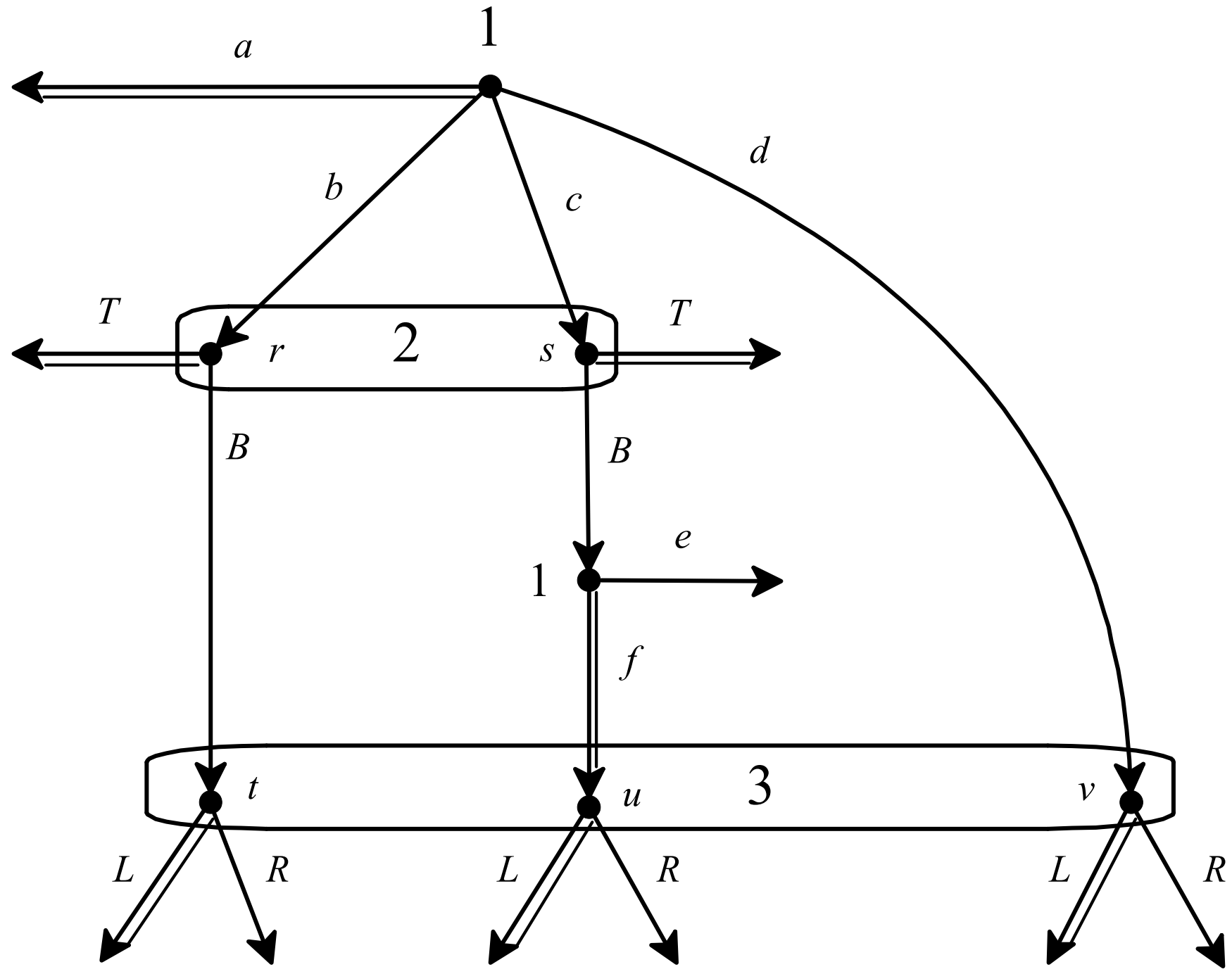

**Figure 11.3**

We want to show that the following assessment is *not* consistent:

$$\sigma = \begin{pmatrix} a & b & c & d & T & B & e & f & L & R \\ 1 & 0 & 0 & 0 & 1 & 0 & 0 & 1 & 1 & 0 \end{pmatrix}, \quad \mu = \begin{pmatrix} r & s & t & u & v \\ \frac{1}{2} & \frac{1}{2} & \frac{1}{5} & \frac{2}{5} & \frac{2}{5} \end{pmatrix}.$$

Suppose that $(\sigma, \mu)$ is consistent. Then there must be a sequence $\langle \sigma_n \rangle_{n=1,2,...}$ of completely mixed strategies that converges to $\sigma$, whose corresponding sequence of systems of beliefs $\langle \mu_n \rangle_{n=1,2,...}$ (obtained by applying Bayesian updating) converges to $\mu$. Let the $n^{\text{th}}$ term of this sequence of completely mixed strategies be:

$$\sigma_n = \left(\begin{array}{cccc|cc|cc|cc} a & b & c & d & T & B & e & f & L & R \\ p_n^a & p_n^b & p_n^c & p_n^d & p_n^T & p_n^B & p_n^e & p_n^f & p_n^L & p_n^R \end{array}\right).$$

Thus $p_n^a + p_n^b + p_n^c + p_n^d = p_n^T + p_n^B = p_n^e + p_n^f = p_n^L + p_n^R = 1$, $p_n^x \in (0,1)$ for all $x \in \{a,b,c,d,B,T,e,f,L,R\}$, $\lim_{n\to\infty} p_n^x = 0$ for $x \in \{b,c,d,B,e,R\}$ and $\lim_{n\to\infty} p_n^x = 1$ for $x \in \{a,T,f,L\}$. Let us compute the corresponding system of beliefs $\mu_n$ obtained by Bayesian updating: $\begin{pmatrix} r & s \\ \dfrac{p_n^b}{p_n^b + p_n^c} & \dfrac{p_n^b}{p_n^b + p_n^c} \end{pmatrix}$ and

$$\begin{pmatrix} t & u & v \\ \dfrac{p_n^b p_n^B}{p_n^b p_n^B + p_n^c p_n^B p_n^f + p_n^d} & \dfrac{p_n^c p_n^B p_n^f}{p_n^b p_n^B + p_n^c p_n^B p_n^f + p_n^d} & \dfrac{p_n^d}{p_n^b p_n^B + p_n^c p_n^B p_n^f + p_n^d} \end{pmatrix}$$

Note that $\dfrac{\mu_n(s)}{\mu_n(r)} = \dfrac{p_n^c}{p_n^b}$. By hypothesis, $\lim_{n\to\infty} \mu_n(s) = \mu(s) = \frac{1}{2}$ and $\lim_{n\to\infty} \mu_n(r) = \mu(r) = \frac{1}{2}$ and thus $\lim_{n\to\infty} \dfrac{p_n^c}{p_n^b} = \lim_{n\to\infty} \dfrac{\mu_n(s)}{\mu_n(r)} = \dfrac{\lim_{n\to\infty}\mu_n(s)}{\lim_{n\to\infty}\mu_n(r)} = \dfrac{\mu(s)}{\mu(r)} = \dfrac{\frac{1}{2}}{\frac{1}{2}} = 1$.

On the other hand, $\dfrac{\mu_n(u)}{\mu_n(t)} = \dfrac{p_n^c}{p_n^b} p_n^f$, so that $\lim_{n\to\infty} \dfrac{\mu_n(u)}{\mu_n(t)} = \lim_{n\to\infty}\left(\dfrac{p_n^c}{p_n^b} p_n^f\right)$.

By hypothesis, $\lim_{n\to\infty} \mu_n(u) = \mu(u) = \frac{2}{5}$ and $\lim_{n\to\infty} \mu_n(t) = \mu(t) = \frac{1}{5}$,

so that $\lim_{n\to\infty} \dfrac{\mu_n(u)}{\mu_n(t)} = \dfrac{\lim_{n\to\infty}\mu_n(u)}{\lim_{n\to\infty}\mu_n(t)} = \dfrac{\mu(u)}{\mu(t)} = 2$.

However, $\lim_{n\to\infty}\left(\dfrac{p_n^c}{p_n^b} p_n^f\right) = \underbrace{\left(\lim_{n\to\infty}\dfrac{p_n^c}{p_n^b}\right)}_{=1}\underbrace{\left(\lim_{n\to\infty} p_n^f\right)}_{=0} = (1)(0) = 0$, yielding a contradiction.

This is a good time to test your understanding of the concepts introduced in this section, by going through the exercises in Section 11.E.1 of Appendix 11.E at the end of this chapter.

# 11.2 Sequential equilibrium

The notion of sequential equilibrium was introduced by Kreps and Wilson (1982).

**Definition 11.2.** Given an extensive game, an assessment $(\sigma, \mu)$ is *a sequential equilibrium* if it is consistent (Definition 11.1) and sequentially rational (Definition 10.2).

For an example of a sequential equilibrium consider the extensive game of Figure 11.4 below.

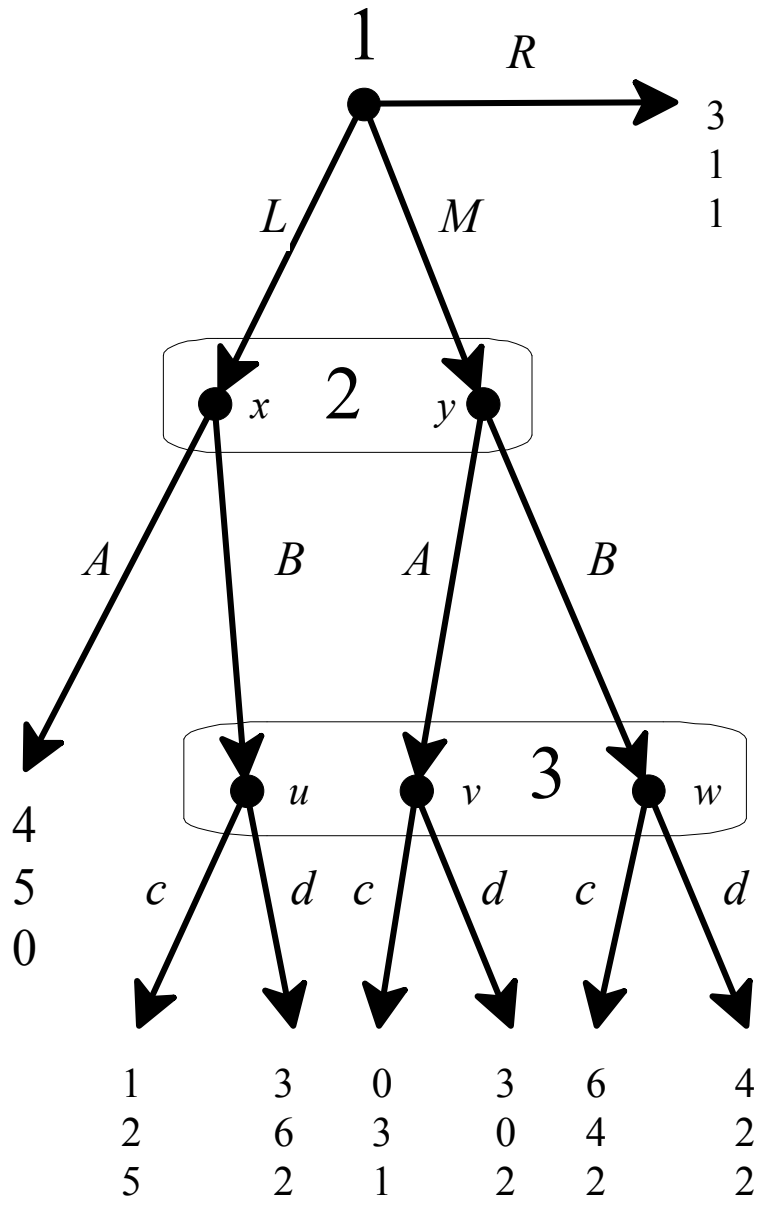

**Figure 11.4**

Let us show that the following assessment is a sequential equilibrium: $\sigma = \begin{pmatrix} L & M & R & | & A & B & | & c & d \\ 0 & 0 & 1 & | & \frac{1}{2} & \frac{1}{2} & | & 1 & 0 \end{pmatrix}$ and $\mu = \begin{pmatrix} x & y & | & u & v & w \\ \frac{1}{4} & \frac{3}{4} & | & \frac{1}{7} & \frac{3}{7} & \frac{3}{7} \end{pmatrix}$. Let us first verify sequential rationality of Player 3's strategy: at her information set $\{u, v, w\}$, given her beliefs, $c$ gives a payoff of $\frac{1}{7}(5) + \frac{3}{7}(1) + \frac{3}{7}(2) = 2$ and $d$ gives a payoff of $\frac{1}{7}(2) + \frac{3}{7}(2) + \frac{3}{7}(2) = 2$. Thus $c$ is optimal (as would be $d$ and any randomization over $c$ and $d$). For Player 2, at his information set $\{x, y\}$, given his beliefs and given the strategy of Player 3, $A$ gives a payoff of $\frac{1}{4}(5) + \frac{3}{4}(3) = 3.5$ and $B$ gives a payoff of $\frac{1}{4}(2) + \frac{3}{4}(4) = 3.5$; thus any mixture of

$A$ and $B$ is optimal, in particular, the mixture $\begin{pmatrix} A & B \\ \frac{1}{2} & \frac{1}{2} \end{pmatrix}$ is optimal. Finally, at the root, $R$ gives Player 1 a payoff of 3, $L$ a payoff of $\frac{1}{2}(4)+\frac{1}{2}(1)=2.5$ and $M$ a payoff of $\frac{1}{2}(0)+\frac{1}{2}(6)=3$; thus $R$ is optimal (as would be any mixture of $M$ and $R$).

Next we show consistency. Consider the sequence of completely mixed strategies $\langle \sigma_n \rangle_{n=1,2,...}$ where $\sigma_n = \left(\begin{array}{ccc|cc|cc} L & M & R & A & B & c & d \\ \frac{1}{n} & \frac{3}{n} & 1-\frac{4}{n} & \frac{1}{2} & \frac{1}{2} & 1-\frac{1}{n} & \frac{1}{n} \end{array}\right)$. Clearly $\lim_{n\to\infty} \sigma_n = \sigma$. The corresponding sequence of systems of beliefs $\langle \mu_n \rangle_{n=1,2,...}$ is given by the following constant sequence, which obviously converges to $\mu$:

$$\mu_n = \left(\begin{array}{cc|ccc} x & y & u & v & w \\ \frac{\frac{1}{n}}{\frac{1}{n}+\frac{3}{n}}=\frac{1}{4} & \frac{\frac{3}{n}}{\frac{1}{n}+\frac{3}{n}}=\frac{3}{4} & \frac{\frac{1}{n}\left(\frac{1}{2}\right)}{\frac{1}{n}\left(\frac{1}{2}\right)+\frac{3}{n}\left(\frac{1}{2}\right)+\frac{3}{n}\left(\frac{1}{2}\right)}=\frac{1}{7} & \frac{\frac{3}{n}\left(\frac{1}{2}\right)}{\frac{1}{n}\left(\frac{1}{2}\right)+\frac{3}{n}\left(\frac{1}{2}\right)+\frac{3}{n}\left(\frac{1}{2}\right)}=\frac{3}{7} & \frac{3}{7} \end{array}\right)$$

We now turn to the properties of sequential equilibria.

**Remark 11.1.** Since consistency implies Bayesian updating at reached information sets, every sequential equilibrium is a weak sequential equilibrium.

**Theorem 11.1** (Kreps and Wilson, 1982). Given an extensive game, if $(\sigma, \mu)$ is a sequential equilibrium then $\sigma$ is a subgame-perfect equilibrium.

**Theorem 11.2** (Kreps and Wilson, 1982). Every finite extensive-form game with cardinal payoffs has at least one sequential equilibrium .

By Theorem 11.1, the notion of sequential equilibrium achieves the objective of refining the notion of subgame-perfect equilibrium. The relationship between the various solution concepts considered so far is shown in the Venn diagram of Figure 11.5.

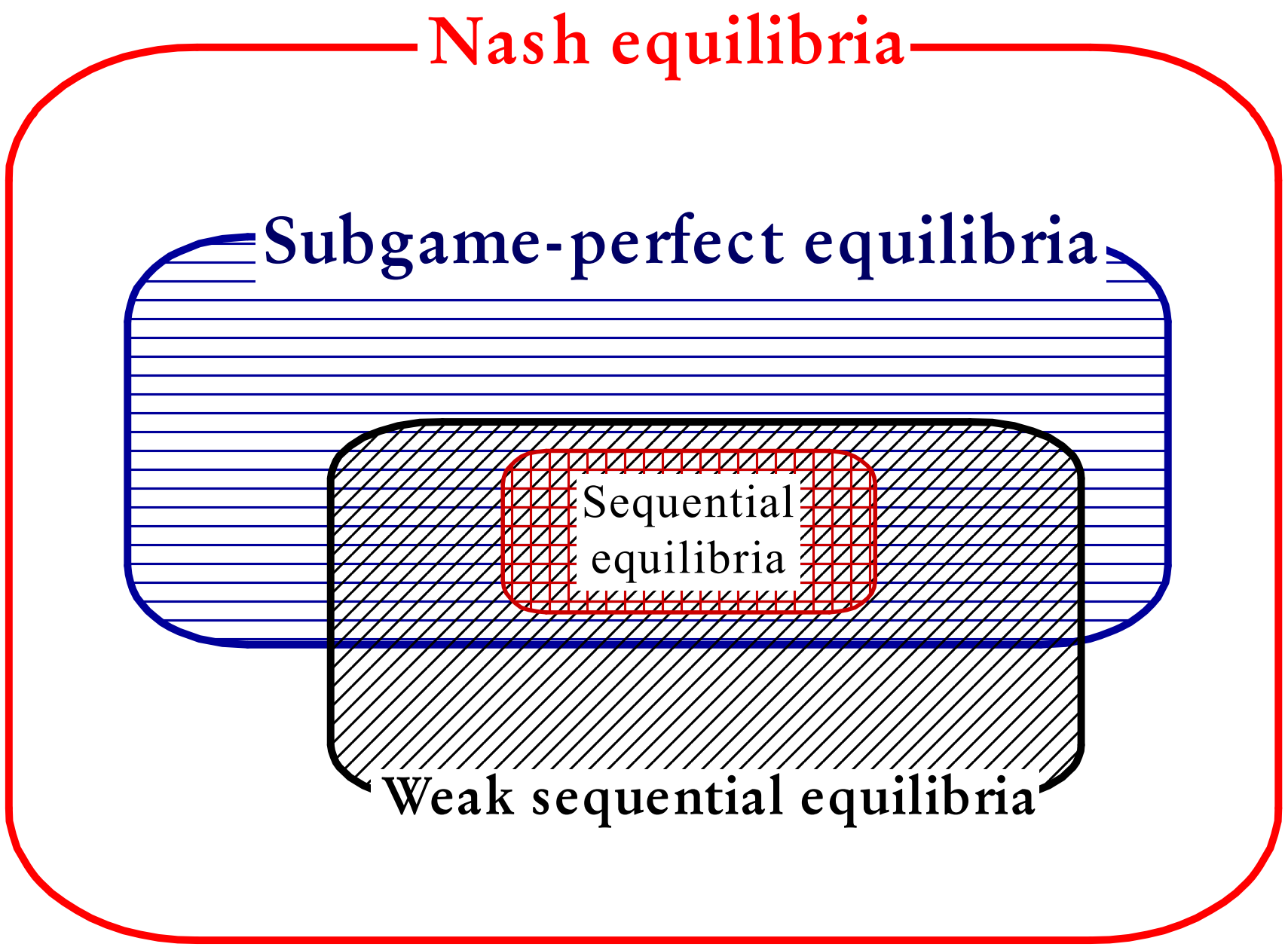

**Figure 11.5**

This is a good time to test your understanding of the concepts introduced in this section, by going through the exercises in Section 11.E.2 of Appendix 11.E at the end of this chapter.

## 11.3 Is 'consistency' a good notion?

The notion of consistency (Definition 11.1) is unsatisfactory in two respects.

From a practical point of view, consistency is computationally hard to prove, since one has to construct a sequence of completely mixed strategies, calculate the corresponding Bayesian beliefs and take the limit of the two sequences. The larger and more complex the game, the harder it is to establish consistency.

From a conceptual point of view, it is not clear how one should interpret, or justify, the requirement of taking the limit of sequences of strategies and beliefs. Indeed, Kreps and Wilson themselves express dissatisfaction with their definition of sequential equilibrium:

> "We shall proceed here to develop the properties of sequential equilibrium as defined above; however, we do so with some doubts of our own concerning what 'ought' to be the definition of a consistent assessment that, with sequential rationality, will give the 'proper' definition of a sequential equilibrium." (Kreps and Wilson, 1982, p. 876.)

In a similar vein, Osborne and Rubinstein (1994, p. 225) write

> "We do not find the consistency requirement to be natural, since it is stated in terms of limits; it appears to be a rather opaque technical assumption."

In the next chapter we will introduce a simpler refinement of subgame-perfect equilibrium which has a clear interpretation in terms of the AGM theory of belief revision: we call it Perfect Bayesian Equilibrium. We will also show that one can use this notion to provide a characterization of sequential equilibrium that does not require the use of limits of sequences of completely mixed strategies.

We conclude this chapter by observing that while the notion of sequential equilibrium eliminates strictly dominated choices at information sets (even if they are reached with zero probability), it is not strong enough to eliminate *weakly* dominated choices. To see this, consider the "simultaneous" game shown in Figure 11.6.

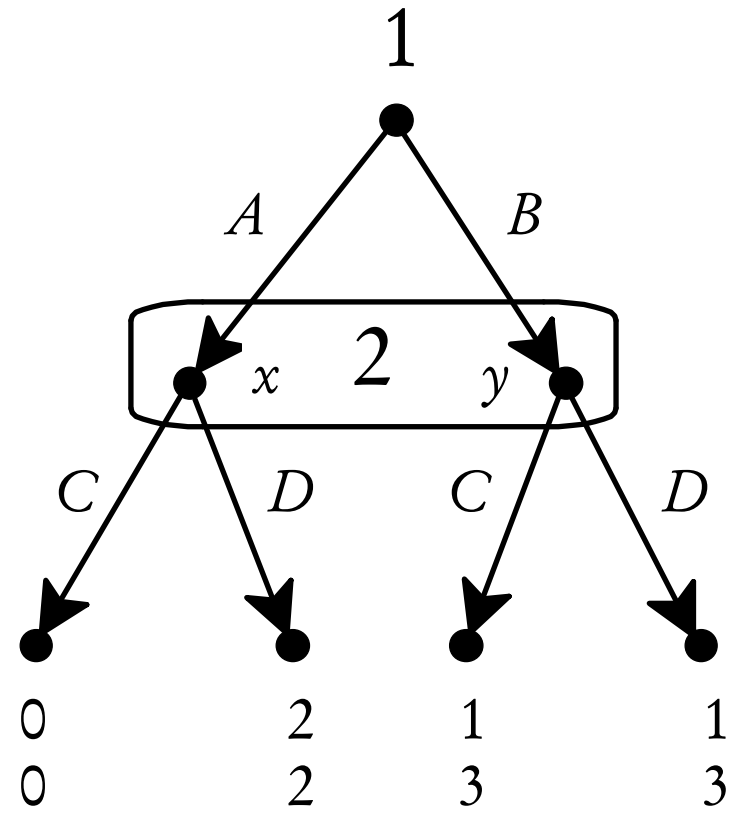

**Figure 11.6**

In this game there are two Nash (and subgame-perfect) equilibria: (*A*, *D*) and (*B*, *C*). However, *C* is a weakly dominated strategy for Player 2. The only beliefs of Player 2 that rationalize choosing *C* is that Player 1 chose *B* with probability 1 (if Player 2 attaches any positive probability, no matter how small, to Player 1 choosing *A*, then *D* is the only sequentially rational choice). Nevertheless, both Nash equilibria are sequential equilibria. For example, it is straightforward to check that the "unreasonable" Nash equilibrium σ = (*B*, *C*), when paired with beliefs $\mu\begin{pmatrix} x & y \\ 0 & 1 \end{pmatrix}$, constitutes a sequential equilibrium. Consistency of this assessment is easily verified by considering the sequence $\sigma_n = \left(\begin{array}{cc|cc} A & B & C & D \\ \frac{1}{n} & 1-\frac{1}{n} & 1-\frac{1}{n} & \frac{1}{n} \end{array}\right)$ whose associated beliefs are $\mu_n = \begin{pmatrix} x & y \\ \frac{1}{n} & 1-\frac{1}{n} \end{pmatrix}$ and sequential rationality is clearly satisfied.

Many game theorists feel that it is "irrational" to choose a weakly dominated strategy; thus further refinements beyond sequential equilibrium have been proposed. A stronger notion of equilibrium, which is a strict refinement of sequential equilibrium, is the notion of *trembling-hand perfect equilibrium*. This notion, due to Reinhardt Selten (who also introduced the notion of subgame-perfect equilibrium) precedes chronologically the notion of sequential equilibrium (Selten, 1975). Trembling-hand perfect equilibrium does in fact eliminate weakly dominated strategies. This topic is outside the scope of this book.[7]

[7] The interested reader is referred to van Damme (2002).

# Appendix 11.E: Exercises

## 11.E.1. Exercises for Section 11.1: Consistent assessments

The answers to the following exercises are in Appendix S at the end of this chapter.

**Exercise 11.1.** Consider the following extensive form.

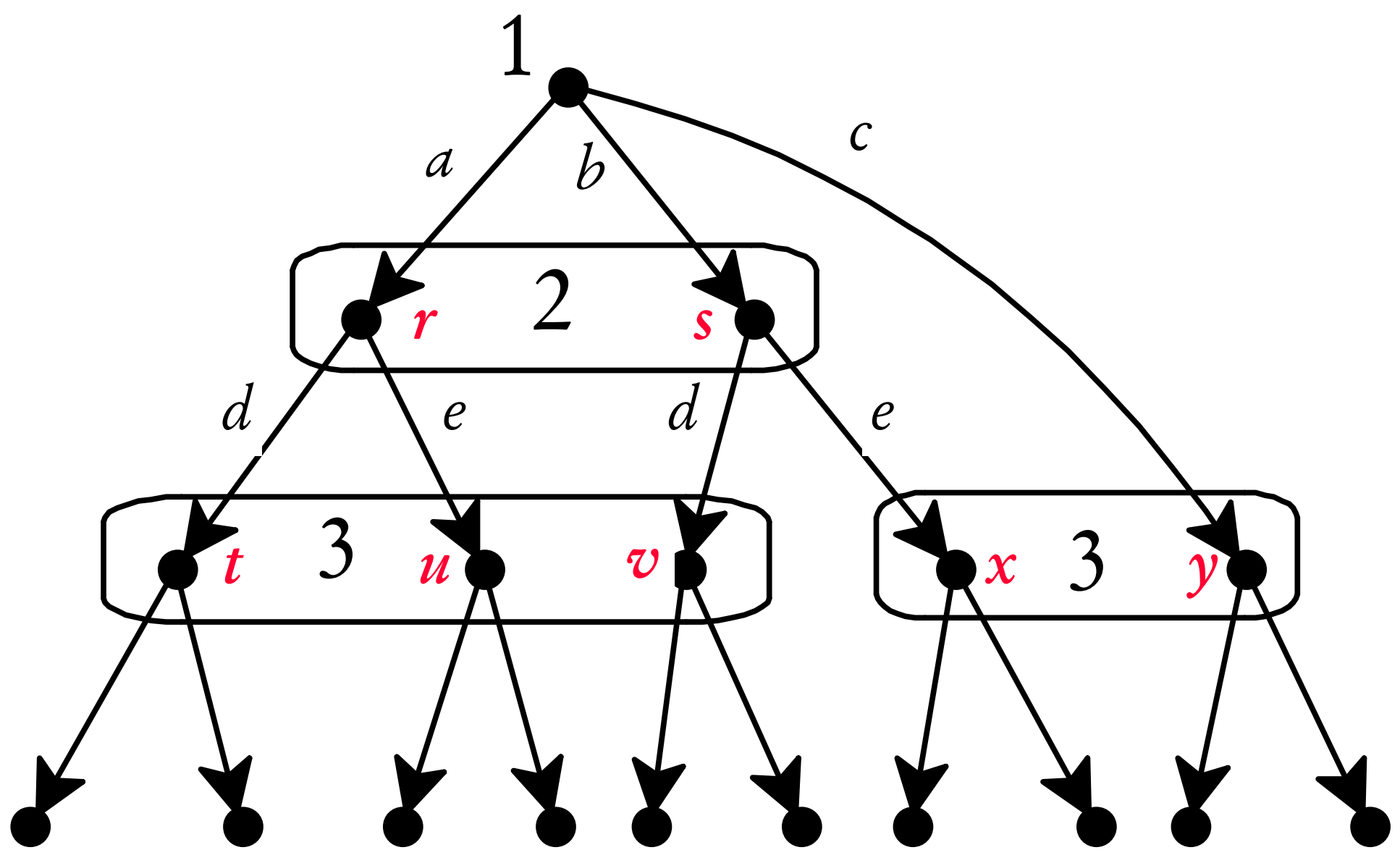

Consider the following (partial) behavior strategy profile

$$\sigma = \left( \begin{array}{ccc|cc} a & b & c & d & e \\ \frac{1}{5} & \frac{3}{5} & \frac{1}{5} & \frac{1}{4} & \frac{3}{4} \end{array} \right).$$

Find the corresponding system of beliefs obtained by Bayesian updating.

**Exercise 11.2**. Consider the following game.

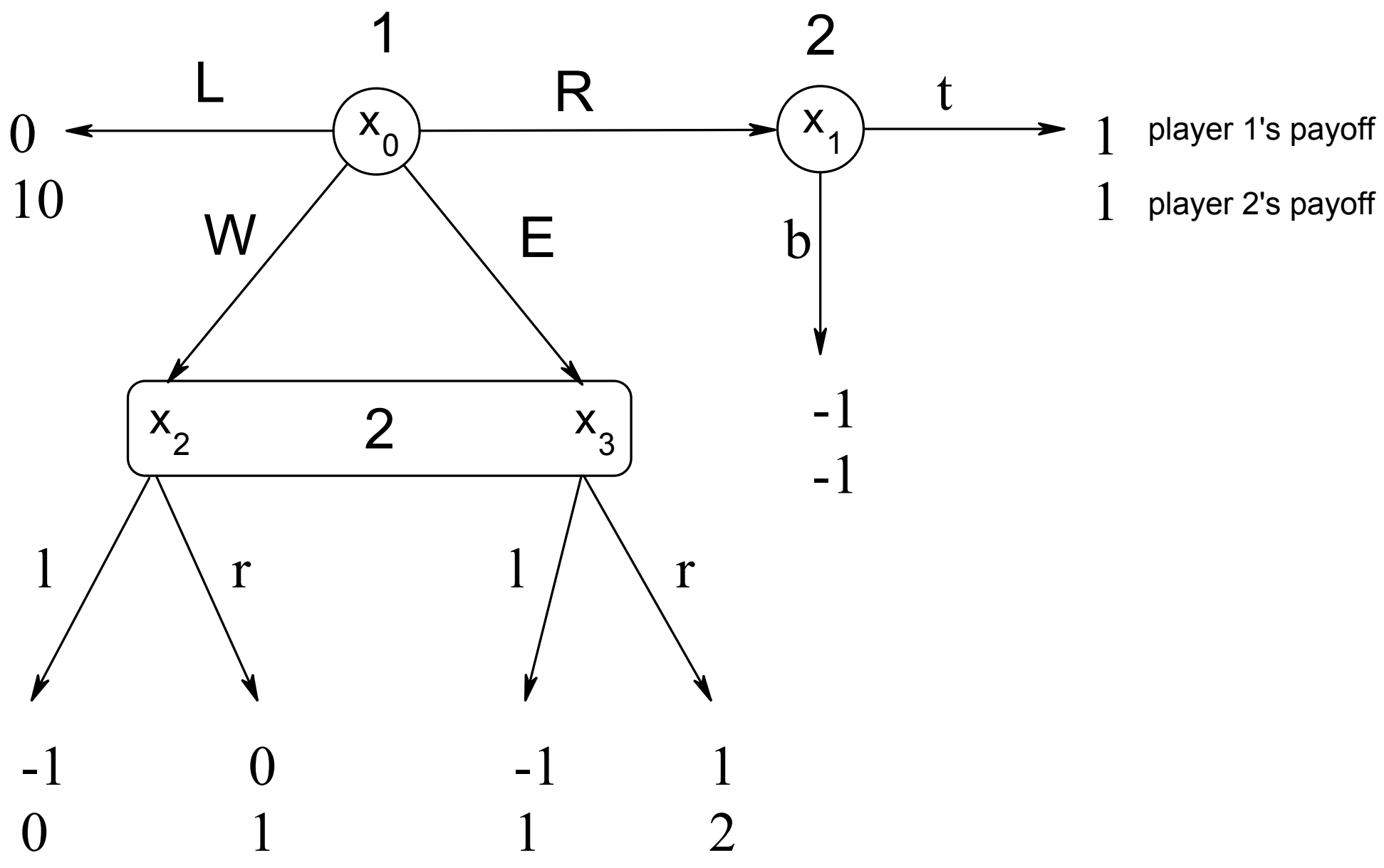

**(a)** Write the corresponding strategic form.

**(b)** Find all the pure-strategy Nash equilibria.

**(c)** Find all the pure-strategy subgame-perfect equilibria.

**(d)** Which of the pure-strategy subgame-perfect equilibria can be part of a consistent assessment? Give a proof for each of your claims.

## 11.E.2. Exercises for Section 11.2: Sequential equilibrium

**Exercise 11.3**. Consider the following extensive-form game. For each pure-strategy Nash equilibrium determine whether it is part of an assessment which is a sequential equilibrium.

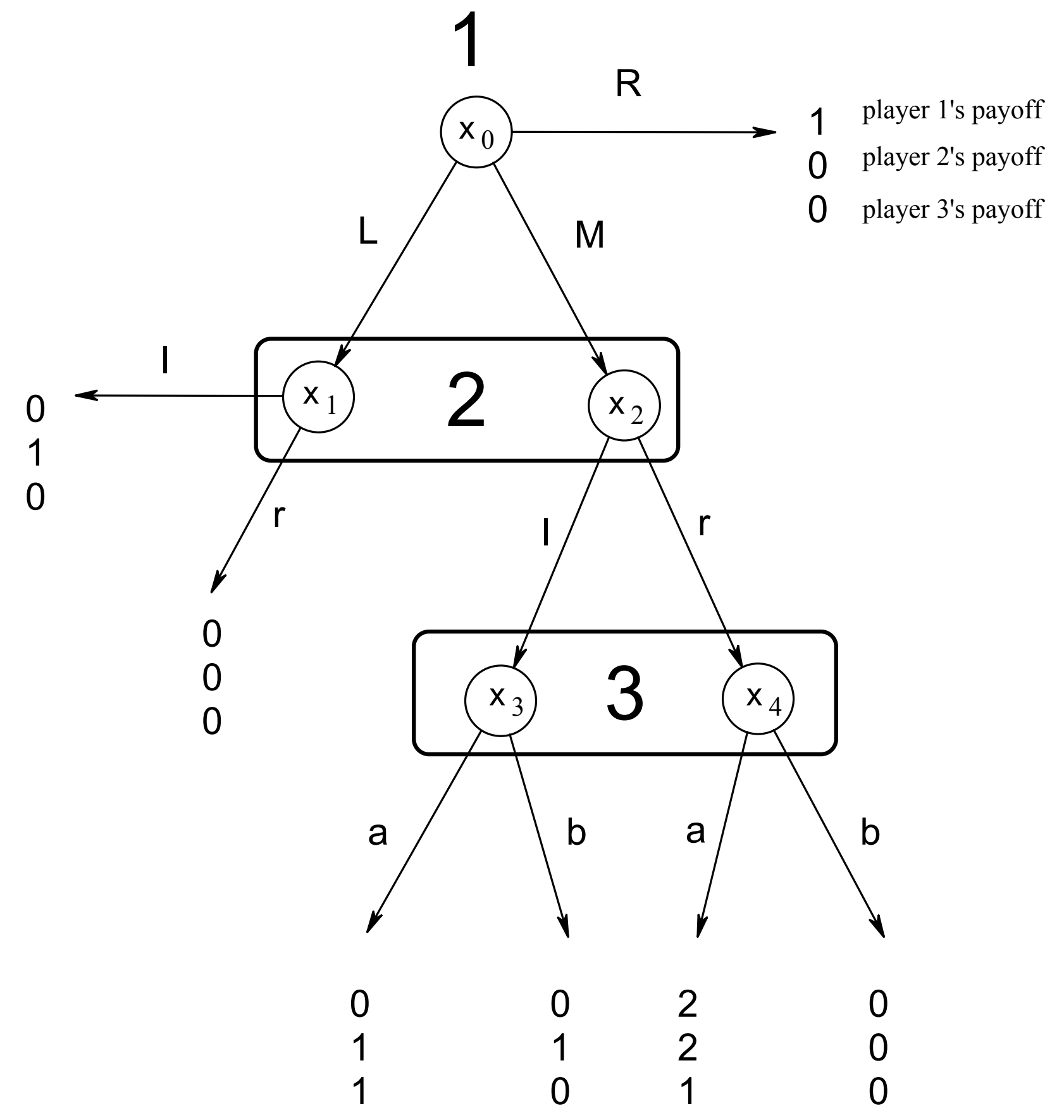

**Exercise 11.4.** Consider the following game:

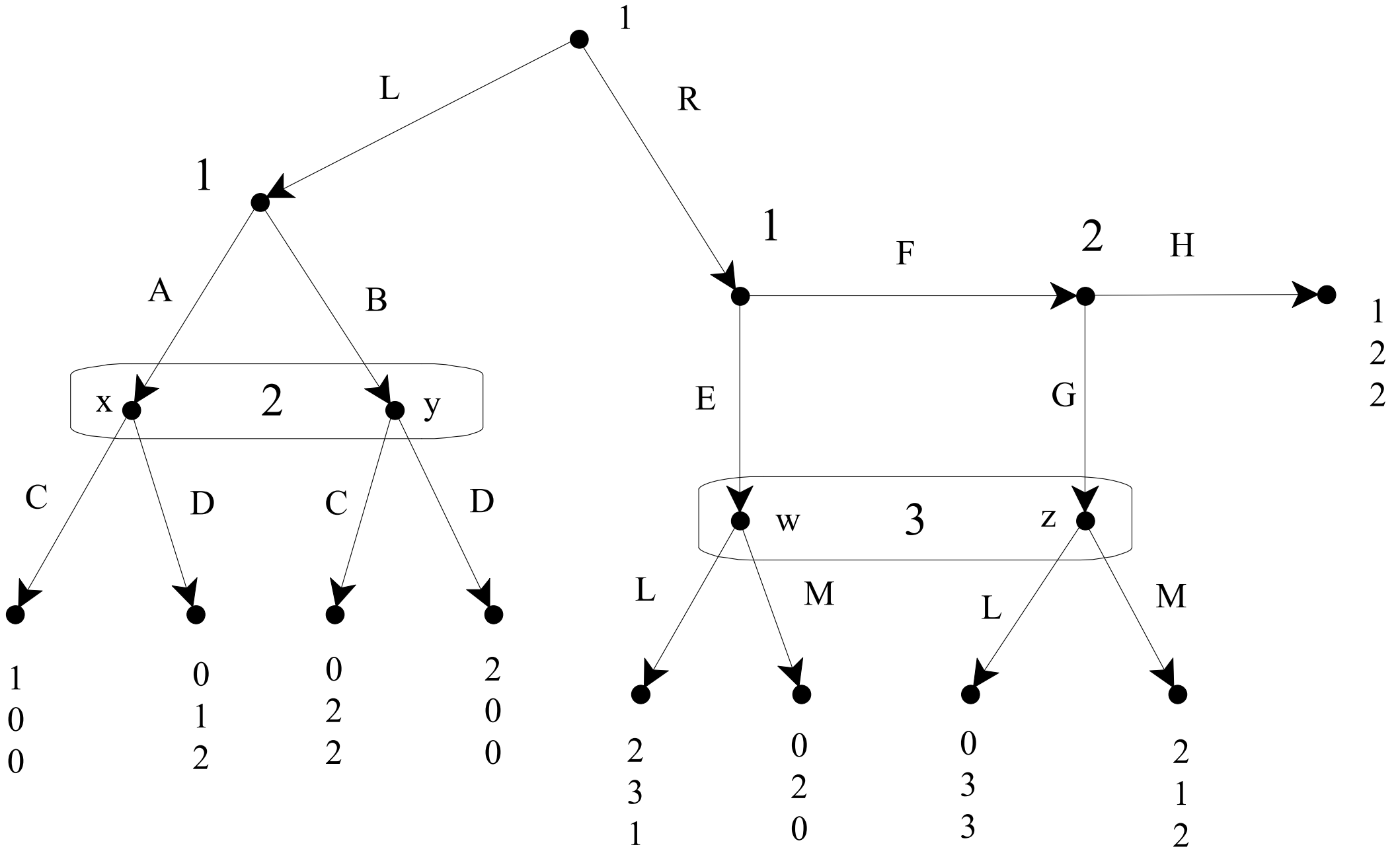

**(a)** Find three subgame-perfect equilibria. [Use pure strategies wherever possible.]

**(b)** For each of the equilibria you found in part (a), explain if it can be written as part of a weak sequential equilibrium.

**(c)** Find a sequential equilibrium. [Use pure strategies wherever possible.]

**Exercise 11.5**. An electric circuit connects two switches and a light. One switch is in Room 1, the second switch is in Room 2 and the light is in Room 3. Player 1 is in Room 1, Player 2 in Room 2 and Player 3 in Room 3. The two switches are now in the Off position. The light in Room 3 comes on if and only if **both** switches are in the On position. Players 1 and 2 act simultaneously and independently: each is allowed only two choices, namely whether to leave her switch in the Off position or turn it to the On position; these decisions have to be executed within one minute. If the light comes on in Room 3 then the game ends and Players 1 and 2 get \$100 each while Player 3 gets \$300. If the light in Room 3 stays off after one minute, then Player 3 (not knowing what the other players did) has to make a guess as to what Players 1 and 2 did (thus, for example, one possible guess is "both players left their respective switches in the Off position"). The payoffs are as follows: **(i)** if Player 3's guess turns out to be correct then each player gets \$200, **(ii)** if Player 3 makes one correct guess but the other wrong (e.g. he guesses that both Player 1 and Player 2 chose "Off" and, as a matter of fact, Player 1 chose "Off" while Player 2 chose "On"), then Player 3 gets \$50, the player whose action was guessed correctly gets \$100 and the remaining player gets nothing (in the previous example, Player 1 gets \$100, Player 2 gets nothing and Player 3 gets \$50) and **(iii)** if Player 3's guess is entirely wrong then all the players get nothing. All the players are selfish and greedy (that is, each player only cares about how much money he/she gets and prefers more money to less) and risk neutral.

**(a)** Represent this situation as an extensive-form game.

**(b)** Write the corresponding strategic form, assigning the rows to Player 1, the columns to Player 2, etc.

**(c)** Find all the pure-strategy Nash equilibria.

**(d)** For at least one pure-strategy Nash equilibrium prove that it cannot be part of a sequential equilibrium.

◊◊◊◊◊◊◊◊◊◊◊◊

**Exercise 11.6. (Challenging Question).**

A buyer and a seller are bargaining over an object owned by the seller. The value of the object to the buyer is known to her but not to the seller. The buyer is drawn randomly from a population with the following characteristics: the fraction $\lambda$ value the object at \$$H$ while the remaining fraction value the object at \$$L$, with $H > L > 0$. The bargaining takes place over two periods. In the first period the seller makes a take-it-or-leave-it offer (i.e. names the price) and the buyer accepts or rejects. If the buyer accepts, the transaction takes place and the game ends. If the buyer rejects, then the seller makes a new take-it-or-leave-it offer and the buyer accepts or rejects. In either case the game ends. Payoffs are as follows: (1) if the seller's offer is accepted (whether it was made in the first period or in the second period), the seller's payoff is equal to the price agreed upon and the buyer's payoff is equal to the difference between the value of the object to the buyer and the price paid; (2) if the second offer is rejected both players get a payoff of 0. Assume that both players discount period 2 payoffs with a discount factor $\delta \in (0,1)$, that is, from the point of view of period 1, getting \$$x$ in period 2 is considered to be the same as getting \$$\delta x$ in period 1. For example, if the seller offers price $p$ in the first period and the offer is accepted, then the seller's payoff is $p$, whereas if the same price $p$ is offered and accepted in the second period, then the seller's payoff, viewed from the standpoint of period 1, is $\delta p$. Assume that these payoffs are von Neumann-Morgenstern payoffs.

**(a)** Draw the extensive form of this game for the case where, in both periods, the seller can only offer one of two prices: \$10 or \$12; assume further that $H = 20$, $L = 10$, $\delta = \frac{3}{4}$ and $\lambda = \frac{2}{3}$. Nature moves first and selects the value for the buyer; the buyer is informed, while the seller is not. It is common knowledge between buyer and seller that Nature will pick $H$ with probability $\lambda$ and $L$ with probability $(1-\lambda)$.

**(b)** For the game of part (a) find a pure-strategy sequential equilibrium. Prove that what you suggest is indeed a sequential equilibrium.

# Appendix 11.S: Solutions to exercises

**Exercise 11.1.** The extensive form under consideration is

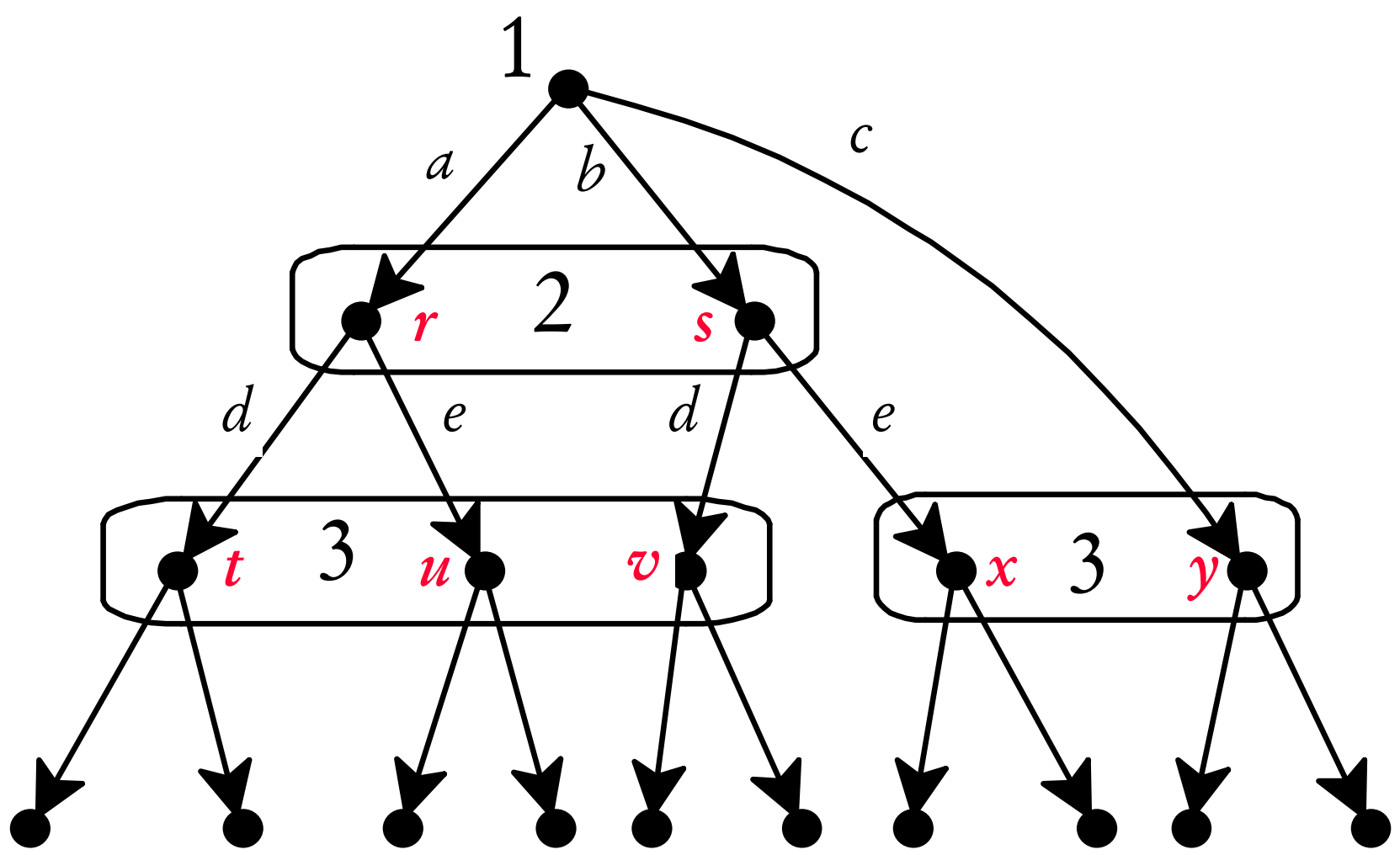

The system of beliefs obtained, by Bayesian updating, from the (partial) behavior strategy profile $\sigma = \begin{pmatrix} a & b & c & d & e \\ \frac{1}{5} & \frac{3}{5} & \frac{1}{5} & \frac{1}{4} & \frac{3}{4} \end{pmatrix}$ is as follows:

$$\mu(r) = \frac{\frac{1}{5}}{\frac{1}{5}+\frac{3}{5}} = \frac{1}{4},\ \mu(s) = \frac{\frac{3}{5}}{\frac{1}{5}+\frac{3}{5}} = \frac{3}{4},\ \mu(x) = \frac{\frac{3}{5}\left(\frac{3}{4}\right)}{\frac{3}{5}\left(\frac{3}{4}\right)+\frac{1}{5}} = \frac{9}{13},\quad \mu(y) = \frac{\frac{1}{5}}{\frac{3}{5}\left(\frac{3}{4}\right)+\frac{1}{5}} = \frac{4}{13},$$

$$\mu(t) = \frac{\frac{1}{5}\left(\frac{1}{4}\right)}{\frac{1}{5}\left(\frac{1}{4}\right)+\frac{1}{5}\left(\frac{3}{4}\right)+\frac{3}{5}\left(\frac{1}{4}\right)} = \frac{1}{7},\quad \mu(u) = \frac{\frac{1}{5}\left(\frac{3}{4}\right)}{\frac{1}{5}\left(\frac{1}{4}\right)+\frac{1}{5}\left(\frac{3}{4}\right)+\frac{3}{5}\left(\frac{1}{4}\right)} = \frac{3}{7}$$

$$\mu(v) = \frac{\frac{3}{5}\left(\frac{1}{4}\right)}{\frac{1}{5}\left(\frac{1}{4}\right)+\frac{1}{5}\left(\frac{3}{4}\right)+\frac{3}{5}\left(\frac{1}{4}\right)} = \frac{3}{7}.$$

**Exercise 11.2.** The game under consideration is

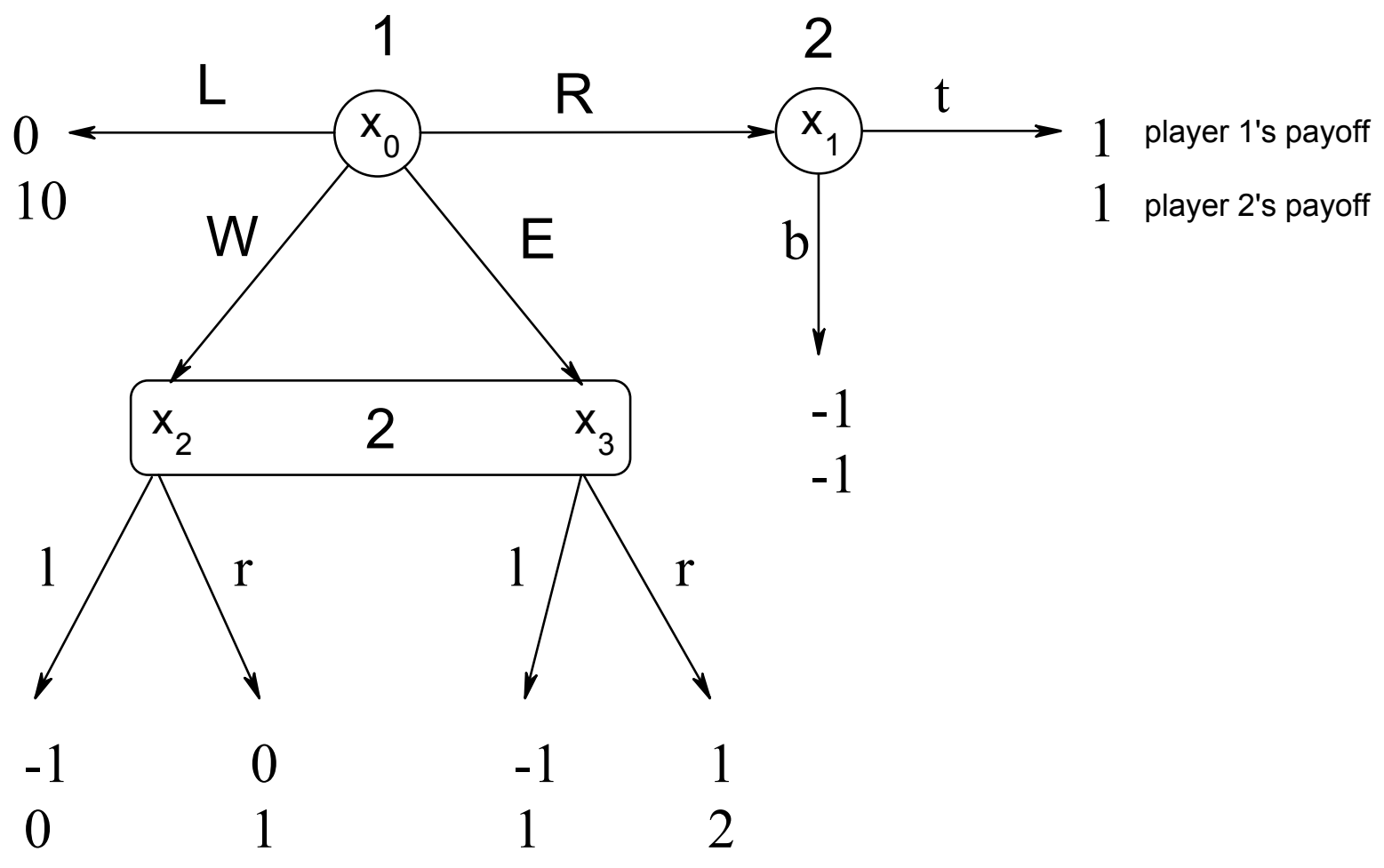

**(a)** The corresponding strategic form is as follows:

Player 2

| Player 1 | | *tl* | *tr* | *bl* | *br* |
|---|---|---|---|---|---|
| | *L* | 0 , 10 | 0 , 10 | 0 , 10 | 0 , 10 |
| | *R* | 1 , 1 | 1 , 1 | −1, −1 | −1, −1 |
| | *W* | −1, 0 | 0 , 1 | −1, 0 | 0 , 1 |
| | *E* | −1, 1 | 1 , 2 | −1, 1 | 1 , 2 |

**(b)** The pure-strategy Nash equilibria are: (*L*, *bl*), (*R*, *tl*), (*R*, *tr*), (*E*, *tr*), (*E*, *br*).

**(c)** There is only one proper subgame, namely the one that starts at node $x_1$ ; in that subgame the unique Nash equilibrium is *t*. Thus only the following are subgame-perfect equilibria: (*R*, *tl*), (*R*, *tr*), (*E*, *tr*).

**(c)** Each of the above subgame-perfect equilibria can be part of a consistent assessment. Let $\sigma = (R, tl)$, $\mu = \begin{pmatrix} x_2 & x_3 \\ \frac{1}{2} & \frac{1}{2} \end{pmatrix}$; this assessment is consistent as the following sequences of completely mixed strategies and corresponding system of beliefs show ($n = 1, 2, \ldots$):

$$\sigma_n = \begin{pmatrix} L & W & E & R & | & b & t & | & l & r \\ \frac{1}{n} & \frac{1}{n} & \frac{1}{n} & 1-\frac{3}{n} & | & \frac{1}{n} & 1-\frac{1}{n} & | & 1-\frac{1}{n} & \frac{1}{n} \end{pmatrix} \quad , \quad \mu_n = \begin{pmatrix} x_2 & x_3 \\ \frac{\frac{1}{n}}{\frac{1}{n}+\frac{1}{n}} = \frac{1}{2} & \frac{\frac{1}{n}}{\frac{1}{n}+\frac{1}{n}} = \frac{1}{2} \end{pmatrix}.$$

Clearly $\lim_{n\to\infty} \sigma_n = \sigma$ and $\lim_{n\to\infty} \mu_n = \mu$.

Let $\sigma = (R, tr)$, $\mu = \begin{pmatrix} x_2 & x_3 \\ \frac{1}{2} & \frac{1}{2} \end{pmatrix}$; this assessment is consistent as the following sequences of completely mixed strategies and corresponding system of beliefs show: $\sigma_n = \begin{pmatrix} L & W & E & R & | & b & t & | & l & r \\ \frac{1}{n} & \frac{1}{n} & \frac{1}{n} & 1-\frac{3}{n} & | & \frac{1}{n} & 1-\frac{1}{n} & | & \frac{1}{n} & 1-\frac{1}{n} \end{pmatrix}$,

$\mu_n = \begin{pmatrix} x_2 & x_3 \\ \frac{\frac{1}{n}}{\frac{1}{n}+\frac{1}{n}} = \frac{1}{2} & \frac{\frac{1}{n}}{\frac{1}{n}+\frac{1}{n}} = \frac{1}{2} \end{pmatrix}$. Clearly $\lim_{n\to\infty} \sigma_n = \sigma$ and $\lim_{n\to\infty} \mu_n = \mu$.

Let $\sigma = (E, tl)$, $\mu = \begin{pmatrix} x_2 & x_3 \\ 0 & 1 \end{pmatrix}$; this assessment is consistent as the following sequences of completely mixed strategies and corresponding beliefs show:

$$\sigma_n = \begin{pmatrix} L & W & E & R & | & b & t & | & l & r \\ \frac{1}{n} & \frac{1}{n} & 1-\frac{3}{n} & \frac{1}{n} & | & \frac{1}{n} & 1-\frac{1}{n} & | & 1-\frac{1}{n} & \frac{1}{n} \end{pmatrix}, \qquad \mu_n = \begin{pmatrix} x_2 & x_3 \\ \frac{\frac{1}{n}}{\frac{1}{n}+1-\frac{3}{n}} & \frac{1-\frac{3}{n}}{\frac{1}{n}+1-\frac{3}{n}} \end{pmatrix}.$$

Clearly $\lim_{n\to\infty} \sigma_n = \sigma$ and $\lim_{n\to\infty} \mu_n = \mu$.

**Exercise 11.3.** The game under consideration is the following:

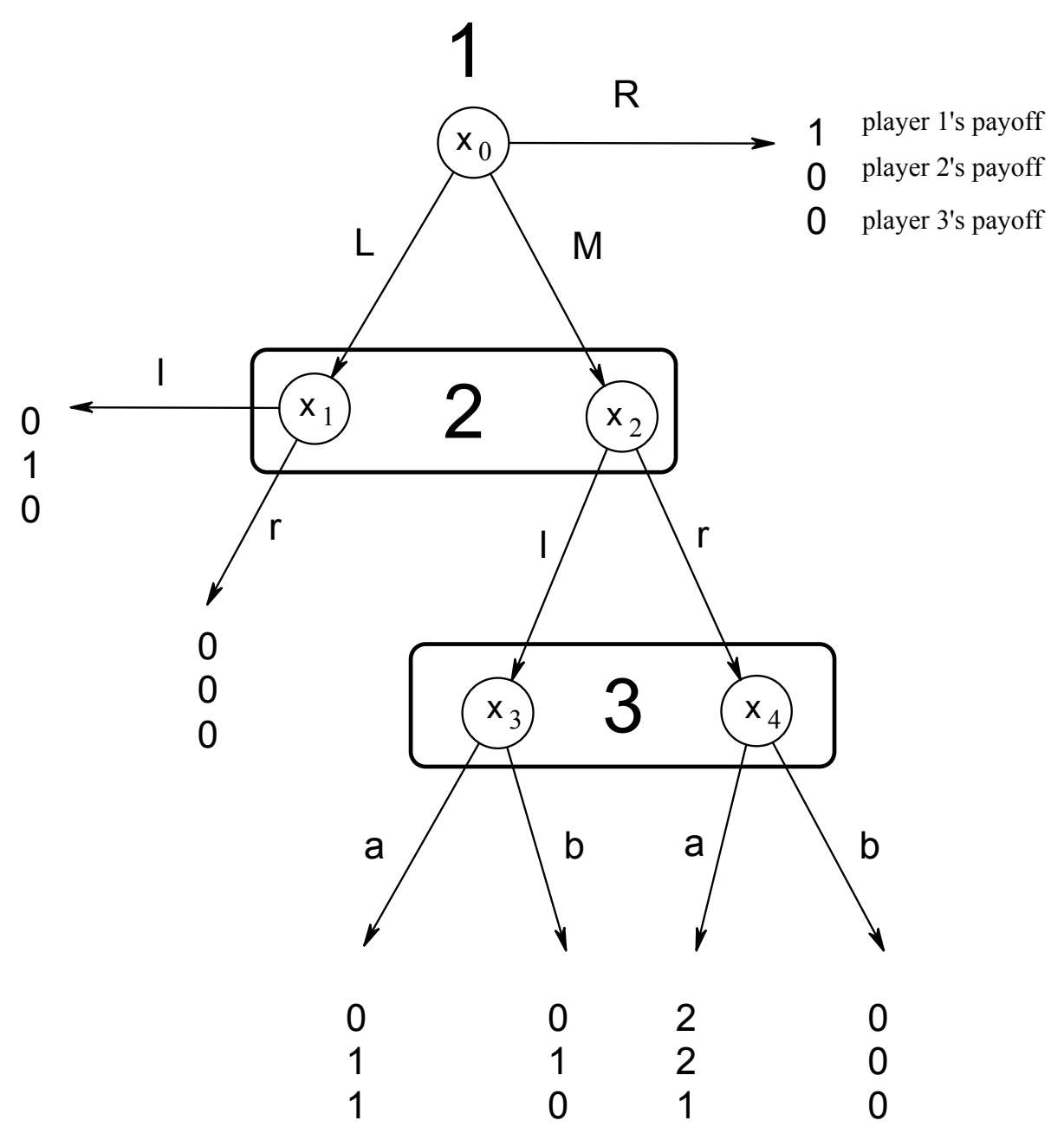

**(a)** The corresponding strategic form is:

Player 3 chooses *a*

| | | Player 2: *l* | Player 2: *r* |
|---|---|---|---|
| Player 1 | *R* | **1 , 0 , 0** | 1 , 0 , 0 |
| | *M* | 0 , 1 , 1 | **2 , 2 , 1** |
| | *L* | 0 , 1 , 0 | 0 , 0 , 0 |

Player 3 chooses *b*

| | | Player 2: *l* | Player 2: *r* |
|---|---|---|---|
| Player 1 | *R* | **1 , 0 , 0** | **1 , 0 , 0** |
| | *M* | 0 , 1 , 0 | 0 , 0 , 0 |
| | *L* | 0 , 1 , 0 | 0 , 0 , 0 |

**(b)** The pure-strategy Nash equilibria are: (*R*,*l*,*a*), (*M*,*r*,*a*), (*R*,*l*,*b*) and (*R*,*r*,*b*). They are all subgame perfect because there are no proper subgames. (*R*,*l*,*b*) and (*R*,*r*,*b*) cannot be part of a sequential equilibrium because *b* is a strictly dominated choice at Player 3's information set and, therefore, (*R*,*l*,*b*) and (*R*,*r*,*b*) would violate sequential rationality (with any system of beliefs). On the other hand, both (*R*,*l*,*a*) and (*M*,*r*,*a*) can be part of an assessment which is a sequential equilibrium.

First we show that $\sigma$ = $(R,l,a)$ together with the system of beliefs $\mu = \begin{pmatrix} x_1 & x_2 & x_3 & x_4 \\ \frac{2}{3} & \frac{1}{3} & 1 & 0 \end{pmatrix}$ is a sequential equilibrium. Consider the sequence of completely mixed strategy profiles whose $n^{\text{th}}$ term is $\sigma_n = \begin{pmatrix} L & M & R & l & r & a & b \\ \frac{2}{n} & \frac{1}{n} & 1-\frac{3}{n} & 1-\frac{1}{n} & \frac{1}{n} & 1-\frac{1}{n} & \frac{1}{n} \end{pmatrix}$. Clearly $\lim_{n\to\infty}\sigma_n = \sigma$. The corresponding Bayesian system of beliefs has $n^{\text{th}}$ term

$$\mu_n = \begin{pmatrix} x_1 & x_2 & x_3 & x_4 \\ \frac{\frac{2}{n}}{\frac{2}{n}+\frac{1}{n}} = \frac{2}{3} & \frac{\frac{1}{n}}{\frac{2}{n}+\frac{1}{n}} = \frac{1}{3} & \frac{\frac{1}{n}\left(1-\frac{1}{n}\right)}{\frac{1}{n}\left(1-\frac{1}{n}\right)+\frac{1}{n}\left(\frac{1}{n}\right)} = 1-\frac{1}{n} & \frac{\frac{1}{n}\left(\frac{1}{n}\right)}{\frac{1}{n}\left(1-\frac{1}{n}\right)+\frac{1}{n}\left(\frac{1}{n}\right)} = \frac{1}{n} \end{pmatrix}$$

Clearly $\lim_{n\to\infty}\mu_n = \mu$. Thus the assessment is consistent (Definition 11.1). Sequential rationality is easily checked: given the strategy profile and the system of beliefs,
(1) for Player 3, $a$ yields 1, while $b$ yields 0, (2) for Player 2, $l$ yields 1, while $r$ yields $\frac{2}{3}(0)+\frac{1}{3}(2)=\frac{2}{3}$ and (3) for Player 1, $R$ yields 1, while $L$ and M yield 0.

Next we show that $\sigma$ = $(M,r,a)$ together with $\mu = \begin{pmatrix} x_1 & x_2 & x_3 & x_4 \\ 0 & 1 & 0 & 1 \end{pmatrix}$ is a sequential equilibrium. Consider the sequence of completely mixed strategy profiles whose $n^{\text{th}}$ term is $\sigma_n = \begin{pmatrix} L & M & R & l & r & a & b \\ \frac{1}{2n} & 1-\frac{1}{n} & \frac{1}{2n} & \frac{1}{n} & 1-\frac{1}{n} & 1-\frac{1}{n} & \frac{1}{n} \end{pmatrix}$. Clearly $\lim_{n\to\infty}\sigma_n = \sigma$. The corresponding Bayesian system of beliefs has $n^{\text{th}}$ term

$$\mu_n = \begin{pmatrix} x_1 & x_2 & x_3 & x_4 \\ \frac{\frac{1}{2n}}{\frac{1}{2n}+\left(1-\frac{1}{n}\right)} = \frac{1}{2n-1} & 1-\frac{1}{2n-1} & \frac{\frac{1}{n}\left(1-\frac{1}{n}\right)}{\frac{1}{n}\left(1-\frac{1}{n}\right)+\left(1-\frac{1}{n}\right)^2} = \frac{1}{n} & 1-\frac{1}{n} \end{pmatrix}$$

Clearly $\lim_{n\to\infty}\mu_n = \mu$. Thus the assessment is consistent (Definition 11.1). Sequential rationality is easily checked: given the strategy profile and the system of beliefs,
(1) for Player 3, $a$ yields 1, while $b$ yields 0, (2) for Player 2, $l$ yields 1, while $r$ yields 2 and (3) for Player 1, $M$ yields 2, while $R$ yields 1 and $L$ yields 0.

**Exercise 11.4.** The game under consideration is the following:

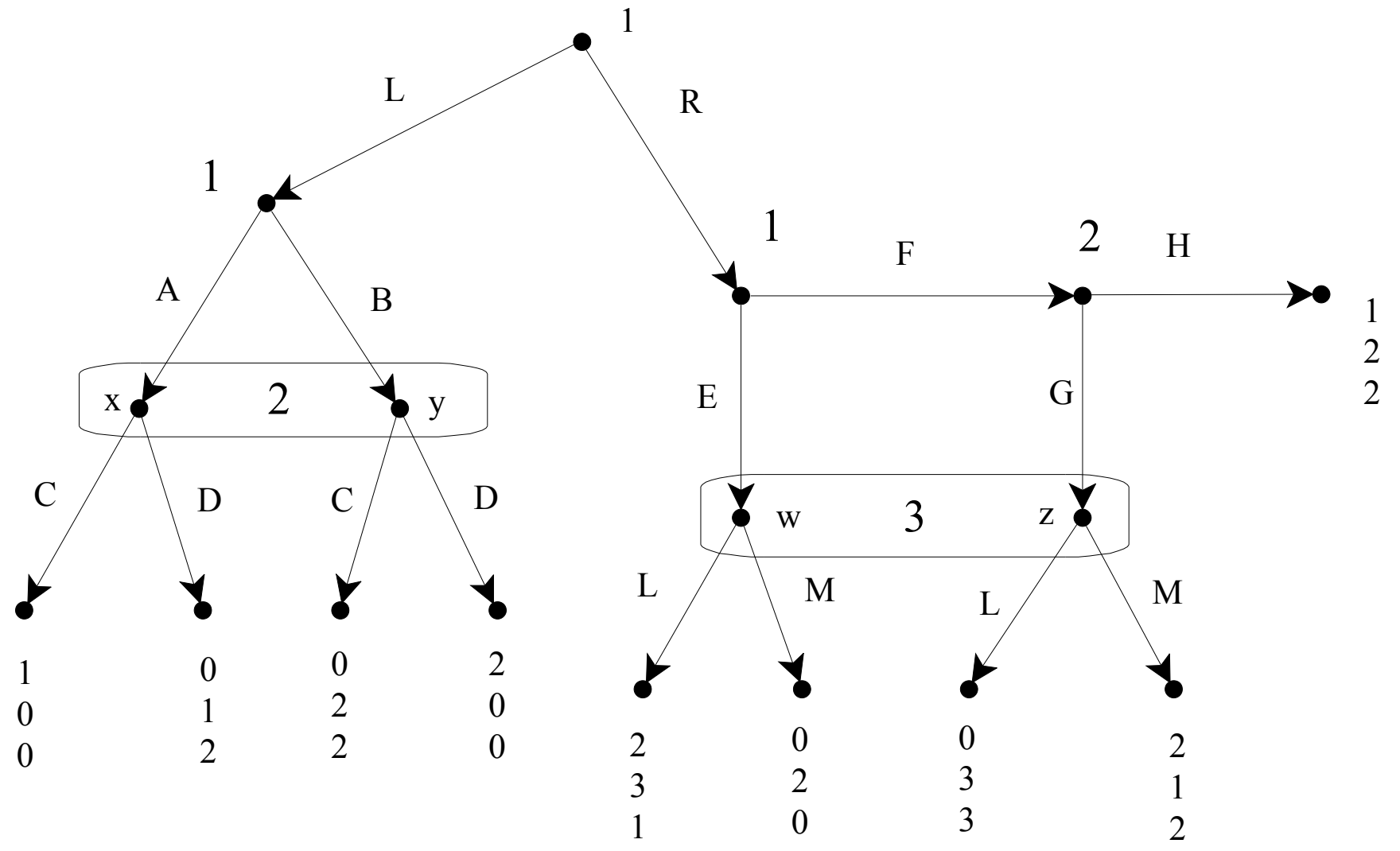

**(a)** First we solve the subgame on the left, whose strategic form is:

| | | Player 2 | |
|---|---|---|---|
| | | *C* | *D* |
| Player 1 | *A* | 1 , 0 | 0 , 1 |
| | *B* | 0 , 2 | 2 , 0 |

There is no pure-strategy Nash equilibrium. Let $p$ be the probability of $A$ and $q$ the probability of $C$. Then at a Nash equilibrium it must be that $q = 2(1-q)$ and $2(1-p) = p$. Thus there is a unique Nash equilibrium given by $\begin{pmatrix} A & B & | & C & D \\ \frac{2}{3} & \frac{1}{3} & | & \frac{2}{3} & \frac{1}{3} \end{pmatrix}$ with an expected payoff of $\frac{2}{3}$ for both players.

Next consider the subgame on the right. In this subgame the following are pure-strategy Nash equilibria: $(F,H,M)$ (where Player 1's payoff is 1), $(E,L,H)$ (where Player 1's payoff is 2), and $(E,L,G)$ (where Player 1's payoff is 2).

Thus the following are subgame-perfect equilibria:

$$\begin{pmatrix} L & R & | & A & B & | & C & D & | & E & F & | & G & H & | & L & M \\ 0 & 1 & | & \frac{2}{3} & \frac{1}{3} & | & \frac{2}{3} & \frac{1}{3} & | & 0 & 1 & | & 0 & 1 & | & 0 & 1 \end{pmatrix}$$

$$\begin{pmatrix} L & R & A & B & C & D & E & F & G & H & L & M \\ 0 & 1 & \frac{2}{3} & \frac{1}{3} & \frac{2}{3} & \frac{1}{3} & 1 & 0 & 0 & 1 & 1 & 0 \end{pmatrix}$$

$$\begin{pmatrix} L & R & A & B & C & D & E & F & G & H & L & M \\ 0 & 1 & \frac{2}{3} & \frac{1}{3} & \frac{2}{3} & \frac{1}{3} & 1 & 0 & 1 & 0 & 1 & 0 \end{pmatrix}$$

(b) $\begin{pmatrix} L & R & A & B & C & D & E & F & G & H & L & M \\ 0 & 1 & \frac{2}{3} & \frac{1}{3} & \frac{2}{3} & \frac{1}{3} & 0 & 1 & 0 & 1 & 0 & 1 \end{pmatrix}$ cannot be part of a weak sequential equilibrium, because choice *M* is strictly dominated by *L* and thus there are no beliefs at Player 3's information set that justify choosing *M*.

$\begin{pmatrix} L & R & A & B & C & D & E & F & G & H & L & M \\ 0 & 1 & \frac{2}{3} & \frac{1}{3} & \frac{2}{3} & \frac{1}{3} & 1 & 0 & 0 & 1 & 1 & 0 \end{pmatrix}$ cannot be part of a weak sequential equilibrium, because – given that Player 3 chooses *L* – *H* is not a sequentially rational choice for Player 2 at his singleton node.

$\begin{pmatrix} L & R & A & B & C & D & E & F & G & H & L & M \\ 0 & 1 & \frac{2}{3} & \frac{1}{3} & \frac{2}{3} & \frac{1}{3} & 1 & 0 & 1 & 0 & 1 & 0 \end{pmatrix}$ is a weak sequential equilibrium with the following system of beliefs: $\begin{pmatrix} x & y & w & z \\ \frac{2}{3} & \frac{1}{3} & 1 & 0 \end{pmatrix}$.

(c) $\begin{pmatrix} L & R & A & B & C & D & E & F & G & H & L & M \\ 0 & 1 & \frac{2}{3} & \frac{1}{3} & \frac{2}{3} & \frac{1}{3} & 1 & 0 & 1 & 0 & 1 & 0 \end{pmatrix}$ together with beliefs $\mu = \begin{pmatrix} x & y & w & z \\ \frac{2}{3} & \frac{2}{3} & 1 & 0 \end{pmatrix}$ is a sequential equilibrium. Consistency can be verified with the sequence of completely mixed strategies

$$\sigma_n = \begin{pmatrix} L & R & A & B & C & D & E & F & G & H & L & M \\ \frac{1}{n} & 1-\frac{1}{n} & \frac{2}{3} & \frac{1}{3} & \frac{2}{3} & \frac{1}{3} & 1-\frac{1}{n} & \frac{1}{n} & 1-\frac{1}{n} & \frac{1}{n} & 1-\frac{1}{n} & \frac{1}{n} \end{pmatrix}$$

For example, $\mu_n(w) = \dfrac{\left(1-\frac{1}{n}\right)\left(1-\frac{1}{n}\right)}{\left(1-\frac{1}{n}\right)\left(1-\frac{1}{n}\right)+\left(1-\frac{1}{n}\right)\frac{1}{n}\left(1-\frac{1}{n}\right)} = \dfrac{1}{1+\frac{1}{n}} \to 1$ as $n \to \infty$.

Sequential rationality is easily verified.

**Exercise 11.5.** **(a)** The extensive form is as follows, where FF means 1Off-2Off, FN means 1Off-2On, etc.

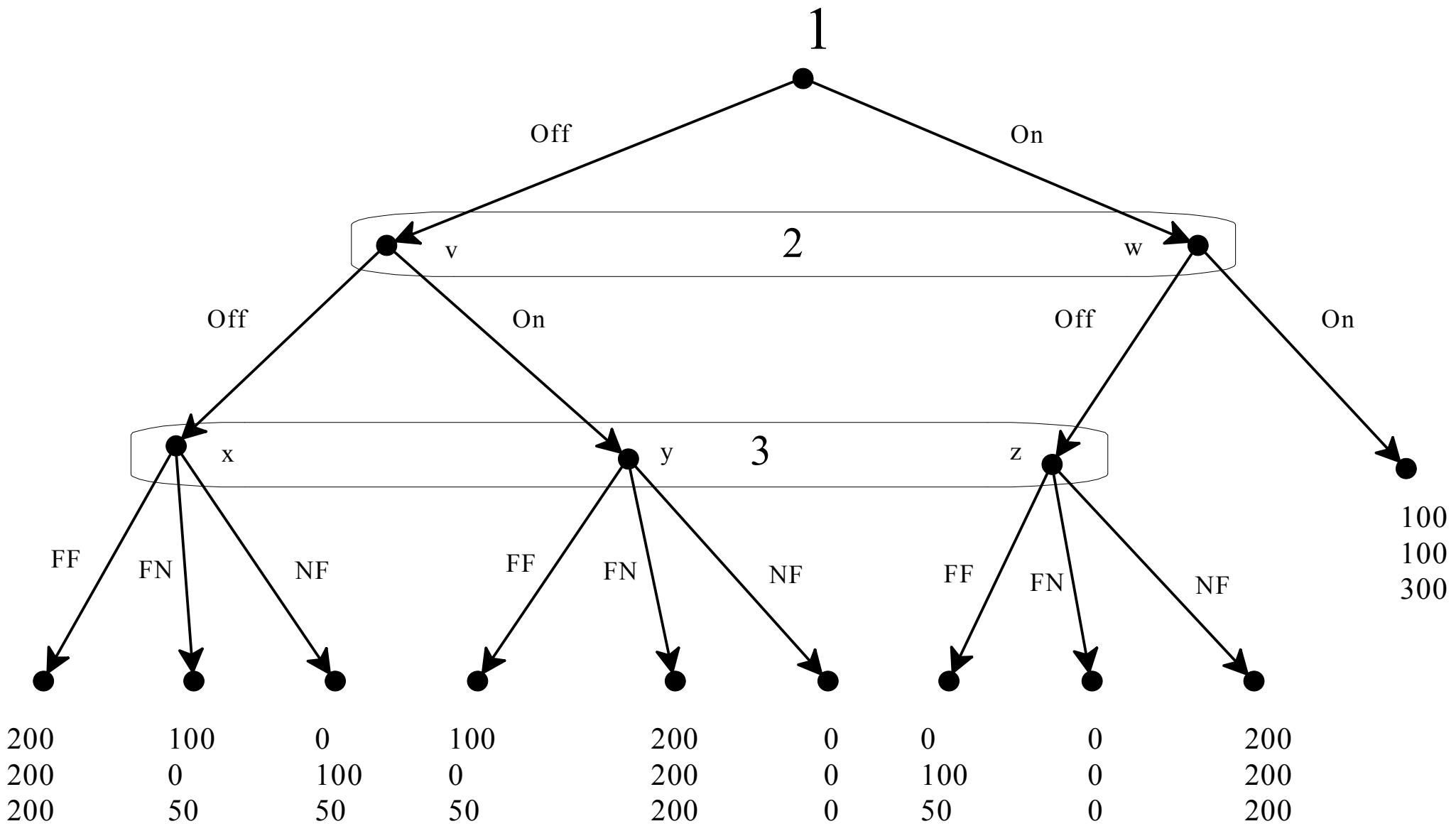

**(b)** The corresponding strategic form is as follows:

**Player 3: 1On-2Off**

| Player 1 \ Player 2 | *On* | *Off* |
|---|---|---|
| *On* | 100 , 100 , 300 | 200 , 200 , 200 |
| *Off* | 0 , 0 , 0 | 0 , 100 , 50 |

**Player 3: 1Off-2On**

| Player 1 \ Player 2 | *On* | *Off* |
|---|---|---|
| *On* | 100 , 100 , 300 | 0 , 0 , 0 |
| *Off* | 200 , 200 , 200 | 100 , 0 , 50 |

**Player 3: both Off**

| Player 1 \ Player 2 | *On* | *Off* |
|---|---|---|
| *On* | 100 , 100 , 300 | 0 , 100 , 100 |
| *Off* | 100 , 0 , 50 | 200 , 200 , 200 |

**(c)** The Nash equilibria are highlighted in the strategic form: (*On, Off,* 1*On*-2*Off*), (*Off, On,* 1*Off*-2*On*), (*On, On, both-Off*) and (*Off, Off, both-Off*).

**(d)** (*On, On, both-Off*) cannot be part of a sequential equilibrium. First of all, for Player 3 '*both-Off*' is a sequentially rational choice only if Player 3 attaches (sufficiently high) positive probability to node *x*. However, consistency does not allow beliefs with $\mu(x) > 0$. To see this, consider a sequence of completely mixed strategies $\{p_n, q_n\}$ for Players 1 and 2, where $p_n$ is the probability with which Player 1 chooses *On* and $q_n$ is the probability with which Player 2 chooses *On* and $\lim_{n\to\infty} p_n = \lim_{n\to\infty} q_n = 0$. Then, by Bayesian updating,

$$P_n\left(x \mid \{x,y,z\}\right) = \frac{p_n q_n}{p_n q_n + p_n(1-q_n) + q_n(1-p_n)}. \qquad (\blacklozenge)$$

If $q_n$ goes to 0 as fast as, or faster than, $p_n$ (that is, if $\lim_{n\to\infty} \frac{q_n}{p_n}$ is finite, then divide numerator and denominator of ($\blacklozenge$) by $p_n$ to get $P_n\left(x \mid \{x,y,z\}\right) = \dfrac{q_n}{q_n + (1-q_n) + \frac{q_n}{p_n}(1-p_n)}$. Taking the limit as $n \to \infty$ we get $\dfrac{0}{0+1+\left(\lim_{n\to\infty}\frac{q_n}{p_n}\right)(1)} = 0$. [If $p_n$ goes to 0 as fast as or faster than $q_n$ then repeat the above argument by dividing by $q_n$.] Thus a consistent assessment must assign zero probability to node *x*.

**Exercise 11.6** (Challenging Question). **(a)** The game is as follows, with $p$ = \$12 and $q$ = \$10.:

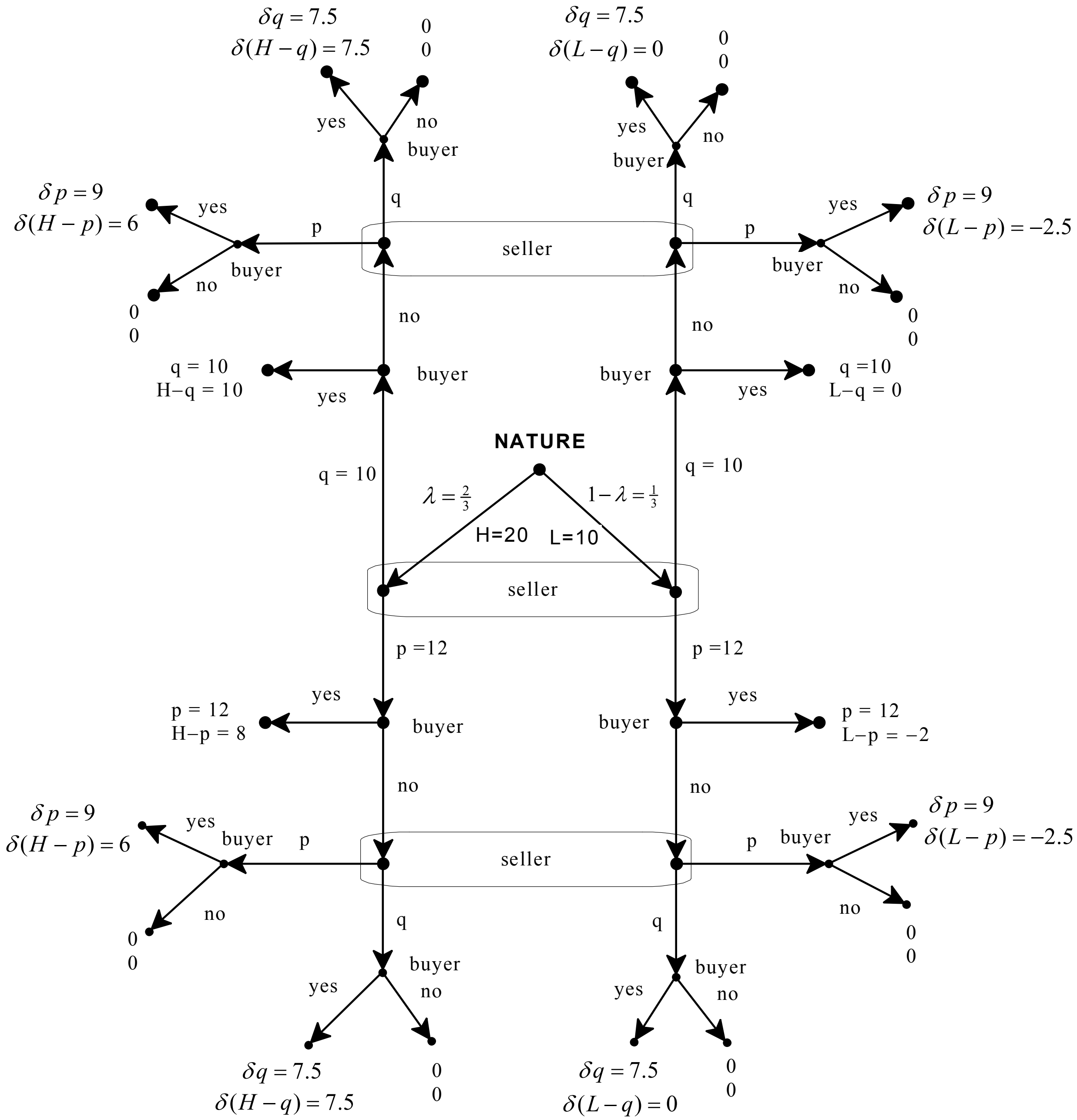

**(b)** The following assessment $(\sigma,\mu)$ is a sequential equilibrium. The pure-strategy profile $\sigma$ is as follows:

**The seller** offers $p$ = 12 in period 1 and, if his offer is rejected, he offers $q$ = 10 in period 2; furthermore, if the first-period offer had been $q$ = 10 and it had been rejected then he would have offered $q$ = 10 again in the second period.

**The *H* buyer** (that is, the buyer at information sets that follow Nature's choice of *H*) always says yes to any offer of the seller.

**The *L* buyer** (that is, the buyer at information sets that follow Nature's choice of *L*) always says no to an offer of $p$ and always says yes to an offer of $q$.

The system of beliefs is as follows (where *TL* means the left node of the top information set of the seller, *TR* the right node of that information set, *ML* means the left node of the middle information set of the seller, *MR* the right node of that information set, *BL* means the left node of the bottom information set of the seller, *BR* the right node of that information set):

$$\mu = \left(\begin{array}{cc|cc|cc} TL & TR & ML & MR & BL & BR \\ 0 & 1 & \frac{2}{3} & \frac{1}{3} & 0 & 1 \end{array}\right).$$

Let us first check sequential rationality. The seller's payoff is $\left(\frac{2}{3}\right)12+\left(\frac{1}{3}\right)\left(\frac{3}{4}\right)10=\frac{63}{6}$ which is greater than the payoff he would get if he offered $q$ = 10, namely a payoff of 10. The *H*-type's payoff is 20 −12 = 8, while if she said no to $p$ = 12 and then yes to $q$ = 10 in period 2 her payoff would be $\frac{3}{4}(20-10)=\frac{15}{2}=7.5$. Furthermore, for the *H* type, at every node of hers, saying yes is always strictly better than saying no. The *L*-type's payoff is 0, while if she said yes to $p$ = 12 then her payoff would be −2. At every node of the *L* type after having been offered $p$ = 12 saying no is strictly better than saying yes and at every node after having been offered $q$ = 10 saying no gives the same payoff as saying yes, namely 0.

To check consistency, construct the following completely mixed strategy profile $\langle\sigma_n\rangle_{n=1,2,\ldots}$ : (1) for the seller and for the *L*-buyer, any choice that has zero probability in $\sigma$ is assigned probability $\frac{1}{n}$ in $\sigma_n$ and any choice that has probability 1 in $\sigma$ is assigned probability $1-\frac{1}{n}$ in $\sigma_n$, (2) for the *H*- buyer, any choice that has zero probability in $\sigma$ is assigned probability $\frac{1}{n^2}$ in $\sigma_n$ and any choice that has probability 1 in $\sigma$ is assigned probability $1-\frac{1}{n^2}$ in $\sigma_n$. Let us compute the corresponding beliefs $\mu_n$ at the top and at the bottom information sets of the seller: $\mu_n(TL) = \dfrac{\frac{1}{n}\left(\frac{1}{n^2}\right)}{\frac{1}{n}\left(\frac{1}{n^2}\right)+\frac{1}{n}\left(\frac{1}{n}\right)} = \dfrac{1}{1+n}$ and $\mu_n(TR) = \dfrac{\frac{1}{n}\left(\frac{1}{n}\right)}{\frac{1}{n}\left(\frac{1}{n^2}\right)+\frac{1}{n}\left(\frac{1}{n}\right)} = \dfrac{1}{\frac{1}{n}+1}$; thus $\lim_{n\to\infty}\mu_n(TL) = 0 = \mu(TL)$ and $\lim_{n\to\infty}\mu_n(TR) = 1 = \mu(TR)$; $\mu_n(BL) = \dfrac{\left(1-\frac{1}{n}\right)\left(\frac{1}{n}\right)}{\left(1-\frac{1}{n}\right)\left(\frac{1}{n}\right)+\left(1-\frac{1}{n}\right)\left(1-\frac{1}{n}\right)} = \frac{1}{n}$ and $\mu_n(BR) == \dfrac{\left(1-\frac{1}{n}\right)\left(1-\frac{1}{n}\right)}{\left(1-\frac{1}{n}\right)\left(\frac{1}{n}\right)+\left(1-\frac{1}{n}\right)\left(1-\frac{1}{n}\right)} = 1-\frac{1}{n}$; thus $\lim_{n\to\infty}\mu_n(BL) = 0 = \mu(BL)$ and $\lim_{n\to\infty}\mu_n(BR) = 1 = \mu(BR)$.

Chapter

12

# Perfect Bayesian Equilibrium

## 12.1 Belief revision and *AGM* consistency

Any attempt to refine the notion of subgame-perfect equilibrium in extensive-form (or dynamic) games must deal with the issue of *belief revision*: how should a player revise her beliefs when informed that she has to make a choice at an information set to which she initially assigned zero probability? As we saw in Chapter 11, Kreps and Wilson (1982) suggested the notion of a consistent assessment (Definition 11.1) to deal with this issue. From now on, we shall refer to the notion of consistency proposed by Kreps and Wilson, as *KW-consistency* (KW stands for 'Kreps-Wilson'), in order to distinguish it from a different notion of consistency, called *AGM-consistency*, that will be introduced in this section. We shall make use of concepts introduced in Section 8.3 (Chapter 8): the reader might want to review that material before continuing.

In this chapter it will be more convenient to use the so called "history-based" definition of extensive-form game, which is spelled out in Appendix 12.A. Essentially it boils down to identifying a node in the tree with the sequence of actions leading from the root to it. We call a sequence of actions starting from the root of the tree a *history*. For example, consider the extensive form of Figure 12.1.

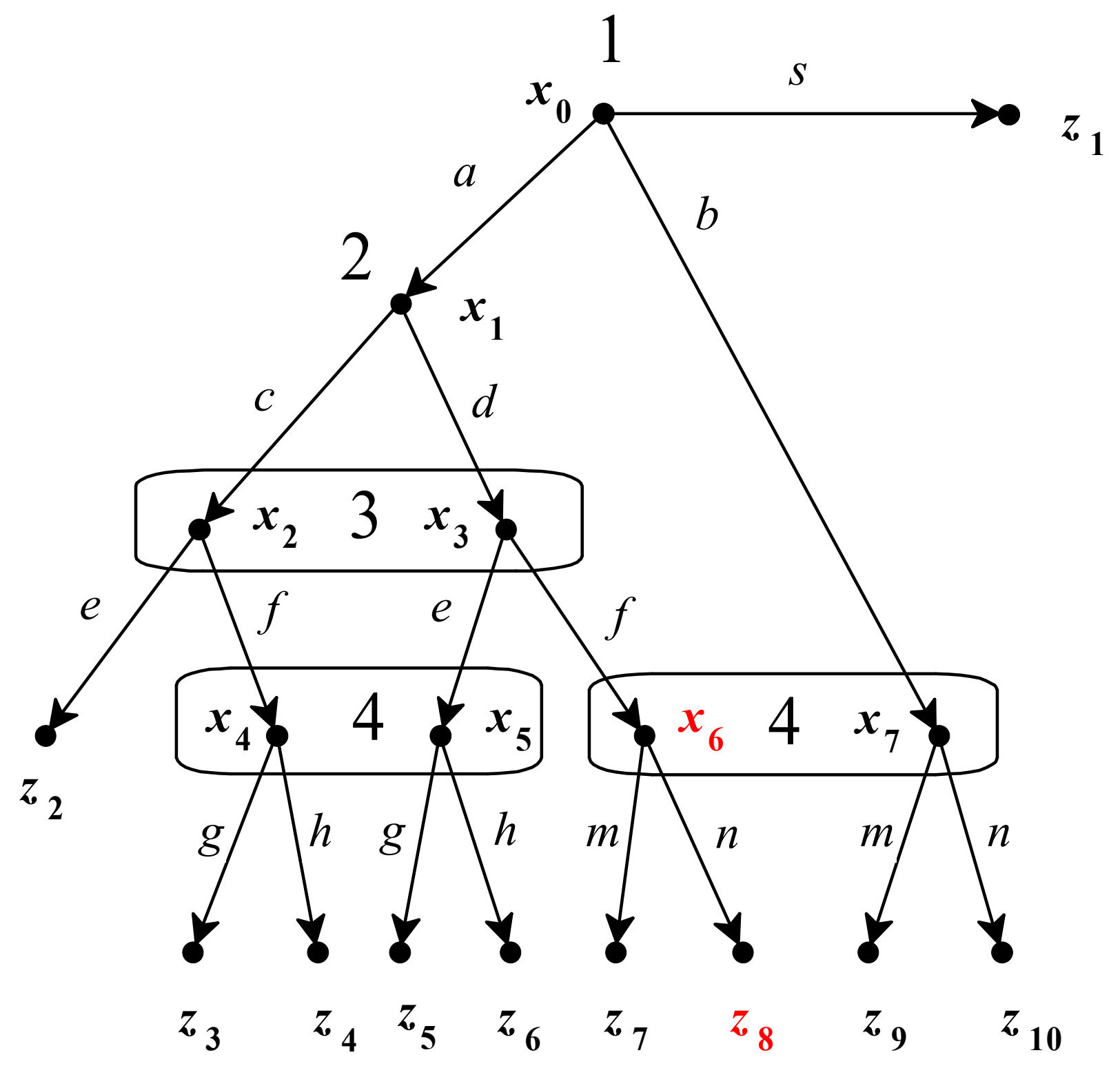

**Figure 12.1**

Node $x_0$ (the root of the tree) is identified with the null or empty history $\emptyset$, decision node $x_6$ with the history (or sequence of actions) *adf*, terminal node $z_8$ with history *adfn*, etc. Thus it will no longer be necessary to label the nodes of the tree, because we can refer to them by naming the corresponding histories.

If *h* is a decision history we denote by $A(h)$ the set of actions (or choices) available at *h*. For example, in the game of Figure 12.1, $A(a) = \{c, d\}$, $A(ac) = \{e, f\}$, $A(adf) = \{m, n\}$, etc. If *h* is a history and *a* is an action available at *h* (that is $a \in A(h)$), then we denote by *ha* the history obtained by appending *a* to *h*.

We saw in Chapter 8 (Section 8.3) that the AGM theory of belief revision (introduced by Alchourrón *et al.*, 1985) is intimately linked to the notion of a plausibility order. This is what motivates the following definition.[8]

[8] For an explicit analysis of the connection between the AGM theory of belief revision and the notion of plausibility order introduced in Definition 12.1 see Bonanno (2011).

First recall that a *total pre-order* on a set $H$ is a binary relation $\precsim$ which is complete (for all $h,h' \in H$, either $h \precsim h'$ or $h' \precsim h$ or both) and transitive (for all $h,h',h'' \in H$,, if $h \precsim h'$ and $h' \precsim h''$ then $h \precsim h''$). We write $h \sim h'$ as a short-hand for "$h \precsim h'$ and $h' \precsim h$" and $h \prec h'$ as a short-hand for "$h \precsim h'$ and $h' \not\precsim h$".

**Definition 12.1.** Given a finite extensive form, a *plausibility order* is a total pre-order $\precsim$ on the finite set of histories $H$ that satisfies the following properties (recall that $D$ is the set of decision histories, $A(h)$ is the set of actions available at decision history $h$ and $I(h)$ is the information set that contains decision history $h$): for all $h \in D$,

PL1. For all $a \in A(h)$, $h \precsim ha$.

PL2. (*i*) There exists an $a \in A(h)$ such that $h \sim ha$, and
(*ii*) for all $a \in A(h)$, if $h \sim ha$ then $h' \sim h'a$ for all $h' \in I(h)$,

PL3. if history $h$ is assigned to chance, then $h \sim ha$, for all $a \in A(h)$.

The interpretation of $h \precsim h'$ is that history $h$ is *at least as plausible* as history $h'$; thus $h \prec h'$ means that $h$ is *more plausible* than $h'$ and $h \sim h'$ means that $h$ is *just as plausible* as $h'$. An alternative reading of $h \precsim h'$ is "*$h$ weakly precedes $h'$ in terms of plausibility*".

Property PL1 of Definition 12.1 says that adding an action to a decision history $h$ cannot yield a more plausible history than $h$ itself. Property PL2 says that at every decision history $h$ there is at least one action $a$ which is "plausibility preserving" in the sense that adding $a$ to $h$ yields a history which is as plausible as $h$; furthermore, any such action $a$ performs the same role with any other history that belongs to the same information set as $h$. Property PL3 says that all the actions at a history assigned to chance are plausibility preserving.

**Definition 12.2.** Given an extensive-form, an assessment $(\sigma,\mu)$ is *AGM-consistent* if there exists a plausibility order $\precsim$ on the set of histories $H$ such that:

(*i*) the actions that are assigned positive probability by $\sigma$ are precisely the plausibility-preserving actions: for all $h \in D$ and for all $a \in A(h)$,

$$\sigma(a) > 0 \text{ if and only if } h \sim ha \tag{P1}$$

(*ii*) the histories that are assigned positive probability by $\mu$ are precisely those that are most plausible within the corresponding information set: for all $h \in D$,

$$\mu(h) > 0 \text{ if and only if } h \precsim h', \text{ for all } h' \in I(h). \tag{P2}$$

If $\precsim$ satisfies properties P1 and P2 with respect to $(\sigma,\mu)$, we say that $\precsim$ *rationalizes* $(\sigma,\mu)$.

In conjunction with sequential rationality, the notion of AGM-consistency is sufficient to rule out some subgame-perfect equilibria. Consider, for example, the extensive game of Figure 12.2 and the pure-strategy profile $\sigma = (c,d,f)$ (highlighted by double edges), which constitutes a Nash equilibrium of the game (and also a subgame-perfect equilibrium since there are no proper subgames).

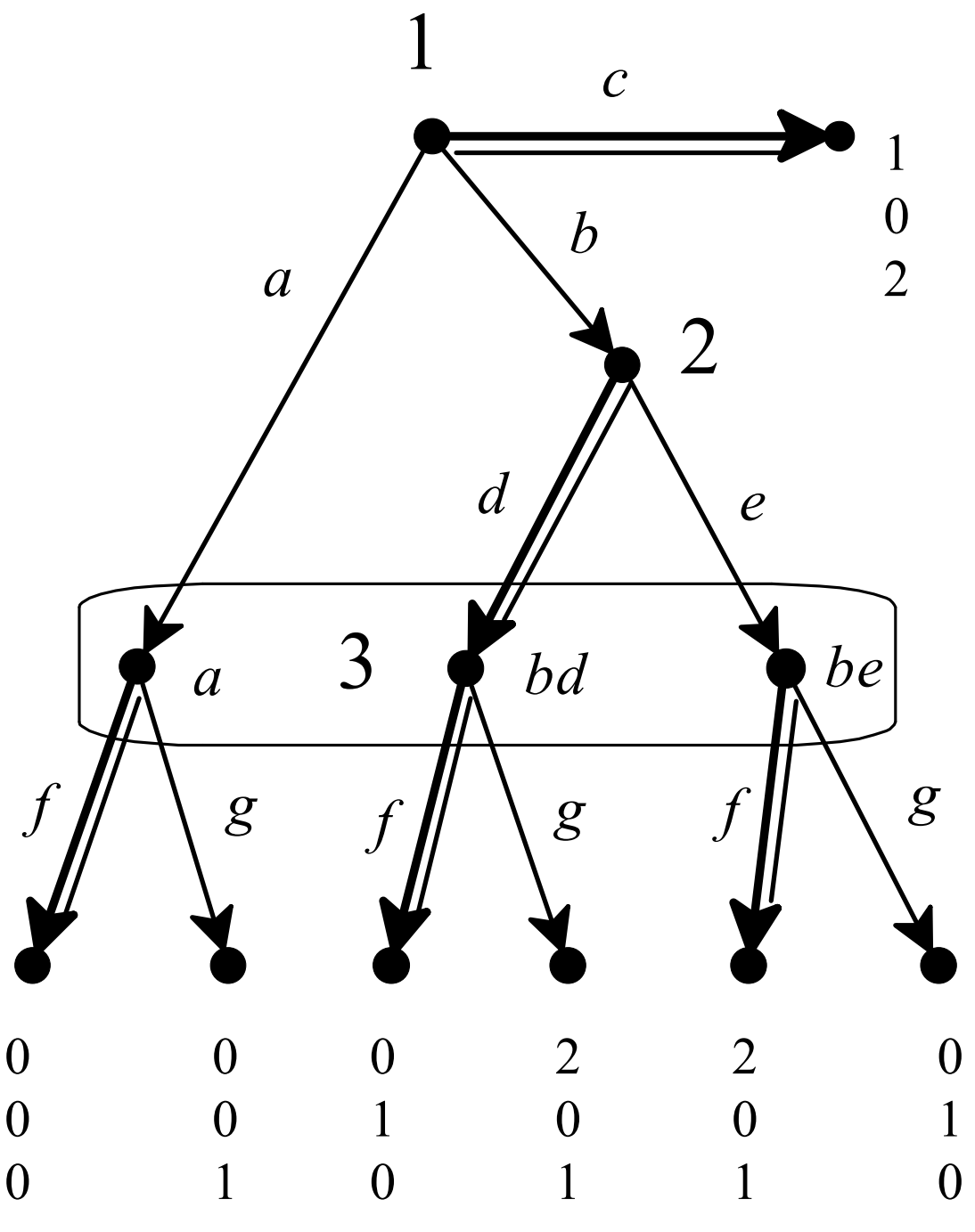

**Figure 12.2**

Can $\sigma$ be part of a sequentially rational AGM-consistent assessment $(\sigma,\mu)$? Since, for Player 3, choice *f* can be chosen rationally only if the player assigns (sufficiently high) positive probability to history *be*, sequential rationality requires that $\mu(be) > 0$; however, any such assessment is *not* AGM-consistent. In fact, if there were a plausibility order $\precsim$ that satisfied Definition 12.2, then, by P1, $b \sim bd$ (since $\sigma(d) = 1 > 0$ ) and $b \prec be$ (since $\sigma(e) = 0$)[9] and, by P2, $be \precsim bd$ (since – by hypothesis – $\mu$ assigns positive probability to *be*). By transitivity of $\precsim$, from $b \prec be$ and $b \sim bd$ it follows that $bd \prec be$, yielding a contradiction.

On the other hand, the Nash equilibrium $\sigma'(c) = 1$, $\sigma'(d) = \sigma'(e) = \frac{1}{2}$, $\sigma'(f) = \sigma'(g) = \frac{1}{2}$ together with beliefs $\mu'(bd) = \mu'(be) = \frac{1}{2}$ forms a sequentially rational, AGM-consistent assessment: it can be rationalized by the plausibility order shown in Figure 12.3, where each row represents an equivalence class ($\emptyset$ denotes the null history, that is, the root of the tree):

$$\begin{pmatrix} \emptyset, c & \text{most plausible} \\ b, bd, be, bdf, bdg, bef, beg & \\ a, af, ag & \text{least plausible} \end{pmatrix}$$

**Figure 12.3**

We use the following convention in representing a plausibility order: if the row to which history *h* belongs is above the row to which $h'$ belongs, then $h \prec h'$ (*h* is more plausible than $h'$) and if *h* and $h'$ belong to the same row then $h \sim h'$ (*h* is as plausible as $h'$).

This is a good time to test your understanding of the concepts introduced in this section, by going through the exercises in Section 12.E.1 of Appendix 12.E at the end of this chapter.

[9] By PL1 of Definition 12.1, $b \precsim be$ and, by P1 of Definition 12.2, it is not the case that $b \sim be$ because $e$ is not assigned positive probability by $\sigma$. Thus $b \prec be$.

# 12.2 Bayesian consistency

The definition of AGM-consistency deals with the supports of a given assessment, that is with the actions that are assigned positive probability by the strategy profile $\sigma$ and the histories that are assigned positive probability by the system of beliefs $\mu$. In this sense it is a *qualitative* (as opposed to quantitative) property: how the probabilities are distributed on those supports is irrelevant for AGM-consistency. However, we also need to impose quantitative restrictions concerning the actual probabilities. The reason for this is that we want the given assessment to satisfy "Bayesian updating as long as possible". By this we mean the following:

(1) when information causes no surprises, because the play of the game is consistent with the most plausible play(s) (that is, when information sets are reached that have positive prior probability), then beliefs should be formed using Bayesian updating (Definition 8.3, Chapter 8), and

(2) when information is surprising (that is, when an information set is reached that had zero prior probability) then new beliefs can be formed in an arbitrary way, but from then on Bayesian updating should be used to update those new beliefs, whenever further information is received that is consistent with those beliefs.

The next definition captures the above requirements.

**Definition 12.3.** Given a finite extensive form, let $\precsim$ be a plausibility order that rationalizes the assessment $(\sigma, \mu)$. We say that $(\sigma, \mu)$ is *Bayesian* (or *Bayes consistent*) *relative to* $\precsim$ if for every equivalence class $E$ of $\precsim$ that contains some decision history $h$ with $\mu(h) > 0$ (that is, $E \cap D_\mu^+ \neq \varnothing$, where $D_\mu^+ = \{h \in D : \mu(h) > 0\}$) there exists a probability distribution $\nu_E : H \to [0,1]$ (recall that $H$ is a finite set) such that:

B1. $\nu_E(h) > 0$ if and only if $h \in E \cap D_\mu^+$.

B2. If $h, h' \in E \cap D_\mu^+$ and $h' = ha_1 \ldots a_m$ (that is, $h$ is a prefix of $h'$) then $\nu_E(h') = \nu_E(h) \times \sigma(a_1) \times \ldots \times \sigma(a_m)$.

B3. If $h \in E \cap D_\mu^+$ then, for every $h' \in I(h)$, $\mu(h') = \nu_E\left(h' \mid I(h)\right) \overset{def}{=} \dfrac{\nu_E(h')}{\sum\limits_{h'' \in I(h)} \nu_E(h'')}$.

Property B1 requires that $\nu_E(h) > 0$ if and only if $h \in E$ and $\mu(h) > 0$. Property B2 requires $\nu_E$ to be consistent with the strategy profile $\sigma$ in the sense that if $h, h' \in E$, $\mu(h) > 0$, $\mu(h') > 0$ and $h' = ha_1 ... a_m$ then the probability that $\nu_E$ assigns to $h'$ is equal to the probability that $\nu_E$ assigns to $h$ multiplied by the probabilities (according to $\sigma$) of the actions that lead from $h$ to $h'$.[10] Property B3 requires the system of beliefs $\mu$ to satisfy the conditional probability rule in the sense that if $h \in E$ and $\mu(h) > 0$ (so that $E$ is the equivalence class of the most plausible elements of $I(h)$) then, for every history $h' \in I(h)$, $\mu(h')$ (the probability assigned to $h'$ by $\mu$) coincides with the probability of $h'$ conditional on $I(h)$ using the probability distribution $\nu_E$.

As an example of an AGM-consistent and Bayesian consistent assessment, consider the extensive form of Figure 12.4 and the assessment $\sigma = \left(\begin{array}{ccc|cc|cc} a & b & c & d & e & f & g \\ 0 & 0 & 1 & 1 & 0 & \frac{1}{3} & \frac{2}{3} \end{array}\right)$, $\mu = \left(\begin{array}{ccc|ccc} ad & ae & b & a & bf & bg \\ \frac{1}{4} & 0 & \frac{3}{4} & \frac{1}{4} & \frac{1}{4} & \frac{2}{4} \end{array}\right)$.

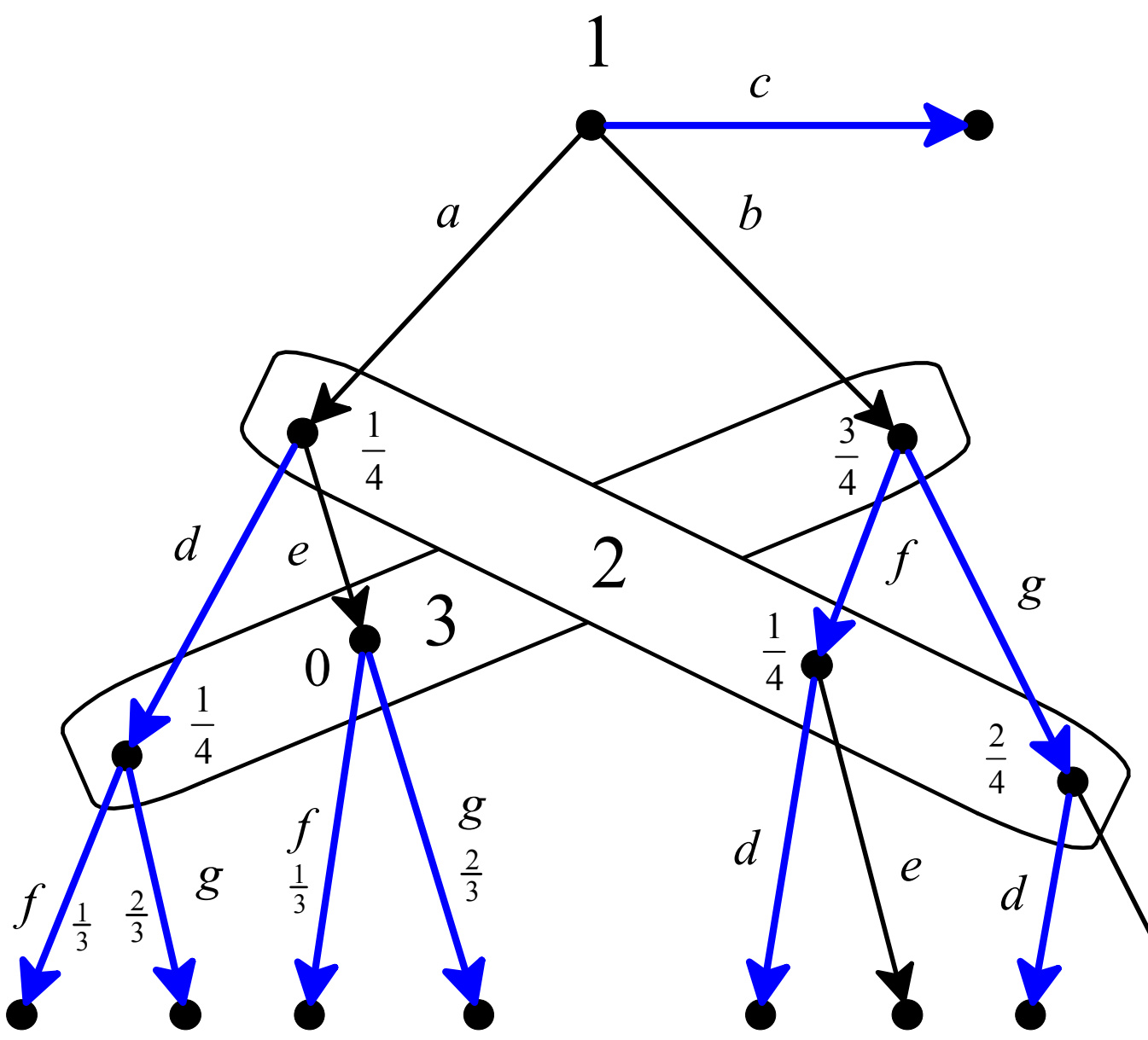

**Figure 12.4**

[10] Note that if $h, h' \in E$ and $h' = ha_1 ... a_m$ then $\sigma(a_j) > 0$ for all $j$ = 1,...,$m$. In fact, since $h \sim h'$, every action $a_j$ is plausibility preserving and thus, by Property P1 of Definition 12.2, $\sigma(a_j) > 0$ .

The assessment $\sigma = \begin{pmatrix} a & b & c & | & d & e & | & f & g \\ 0 & 0 & 1 & | & 1 & 0 & | & \frac{1}{3} & \frac{2}{3} \end{pmatrix}$, $\mu = \begin{pmatrix} ad & ae & b & | & a & bf & bg \\ \frac{1}{4} & 0 & \frac{3}{4} & | & \frac{1}{4} & \frac{1}{4} & \frac{2}{4} \end{pmatrix}$ is AGM consistent, since it is rationalized by the plausibility order shown in Figure 12.5

$$\begin{pmatrix} \emptyset, c & \text{most plausible} \\ a,b,ad,bf,bg,adf,adg,bfd,bgd & \\ ae,aef,aeg,bfe,bge & \text{least plausible} \end{pmatrix}$$

**Figure 12.5**

Furthermore, the assessment is Bayesian relative to this plausibility order. First of all, note that $D_\mu^+ = \{\emptyset, a, ad, b, bf, bg\}$.[11] Let $E$ be the top equivalence class ($E = \{\emptyset, c\}$), $F$ the middle one ($F = \{a,b,ad,bf,bg,adf,adg,bfd,bgd\}$) and $G$ the bottom one ($G = \{ae,aef,aeg,bfe,bge\}$). Then only $E$ and $F$ have a non-empty intersection with $D_\mu^+$ and thus, by Definition 12.3, we only need to specify two probability distributions: $\nu_E$ and $\nu_F$. The first one is trivial, since $D_\mu^+ \cap E = \{\emptyset\}$ and thus it must be $\nu_E(\emptyset) = 1$ (and $\nu_E(h) = 0$ for every other history *h*). Since $D_\mu^+ \cap F = \{a, ad, b, bf, bg\}$, by B1 of Definition 12.3 it must be that $\nu_F(h) > 0$ if and only if $h \in \{a, ad, b, bf, bg\}$. Consider the following probability distribution:

$$\nu_F = \begin{pmatrix} a & ad & b & bf & bg \\ \frac{1}{8} & \frac{1}{8} & \frac{3}{8} & \frac{1}{8} & \frac{2}{8} \end{pmatrix} \quad \text{(and } \nu_F(h) = 0 \text{ for every other history } h\text{).}$$

Property B2 of Definition 12.3 is satisfied, because $\nu_F(ad) = \underbrace{\tfrac{1}{8}}_{=\nu_F(a)} \times \underbrace{1}_{\sigma(d)}$, $\nu_F(bf) = \underbrace{\tfrac{3}{8}}_{=\nu_F(b)} \times \underbrace{\tfrac{1}{3}}_{\sigma(f)}$ and $\nu_F(bg) = \underbrace{\tfrac{3}{8}}_{=\nu_F(b)} \times \underbrace{\tfrac{2}{3}}_{\sigma(g)}$. To check that Property B3 of Definition 12.3 is satisfied, let $I_2 = \{a, bf, bg\}$ be the information set of Player 2 and $I_3 = \{b, ad, ae\}$ be the information set of Player 3. Then $\nu_F(I_2) = \nu_F(a) + \nu_F(bf) + \nu_F(bg) = \frac{4}{8}$ and $\nu_F(I_3) = \nu_F(b) + \nu_F(ad) + \nu_F(ae) = \frac{3}{8} + \frac{1}{8} + 0 = \frac{4}{8}$. Thus

[11] Recall that $D_\mu^+ \overset{def}{=} \{h \in D : \mu(h) > 0\}$.

$$\frac{\nu_F(a)}{\nu_F(I_2)} = \frac{\frac{1}{8}}{\frac{4}{8}} = \frac{1}{4} = \mu(a), \quad \frac{\nu_F(bf)}{\nu_F(I_2)} = \frac{\frac{1}{8}}{\frac{4}{8}} = \frac{1}{4} = \mu(bf), \quad \frac{\nu_F(bg)}{\nu_F(I_2)} = \frac{\frac{2}{8}}{\frac{4}{8}} = \frac{2}{4} = \mu(bg)$$

$$\frac{\nu_F(b)}{\nu_F(I_3)} = \frac{\frac{3}{8}}{\frac{4}{8}} = \frac{3}{4} = \mu(b), \quad \frac{\nu_F(ad)}{\nu_F(I_3)} = \frac{\frac{1}{8}}{\frac{4}{8}} = \frac{1}{4} = \mu(ad), \quad \frac{\nu_F(ae)}{\nu_F(I_3)} = \frac{0}{\frac{4}{8}} = 0 = \mu(ae).$$

This is a good time to test your understanding of the concepts introduced in this section, by going through the exercises in Section 12.E.2 of Appendix 12.E at the end of this chapter.

# 12.3 Perfect Bayesian equilibrium

By adding sequential rationality (Definition 10.2) to AGM-consistency (Definition 12.2) and Bayesian consistency (Definition 12.3) we obtain a new notion of equilibrium for extensive-form games.

**Definition 12.4.** Given a finite extensive-form game, an assessment $(\sigma, \mu)$ is a *perfect Bayesian equilibrium* if it is sequentially rational, it is rationalized by a plausibility order on the set of histories and is Bayesian relative to that plausibility order.

For an example of a perfect Bayesian equilibrium, consider the game of Figure 12.6 (which is based on the frame of Figure 12.4) and the assessment considered in the previous section, namely

$$\sigma = \left(\begin{array}{ccc|cc|cc} a & b & c & d & e & f & g \\ 0 & 0 & 1 & 1 & 0 & \frac{1}{3} & \frac{2}{3} \end{array}\right), \quad \mu = \left(\begin{array}{ccc|ccc} ad & ae & b & a & bf & bg \\ \frac{1}{4} & 0 & \frac{3}{4} & \frac{1}{4} & \frac{1}{4} & \frac{2}{4} \end{array}\right)$$

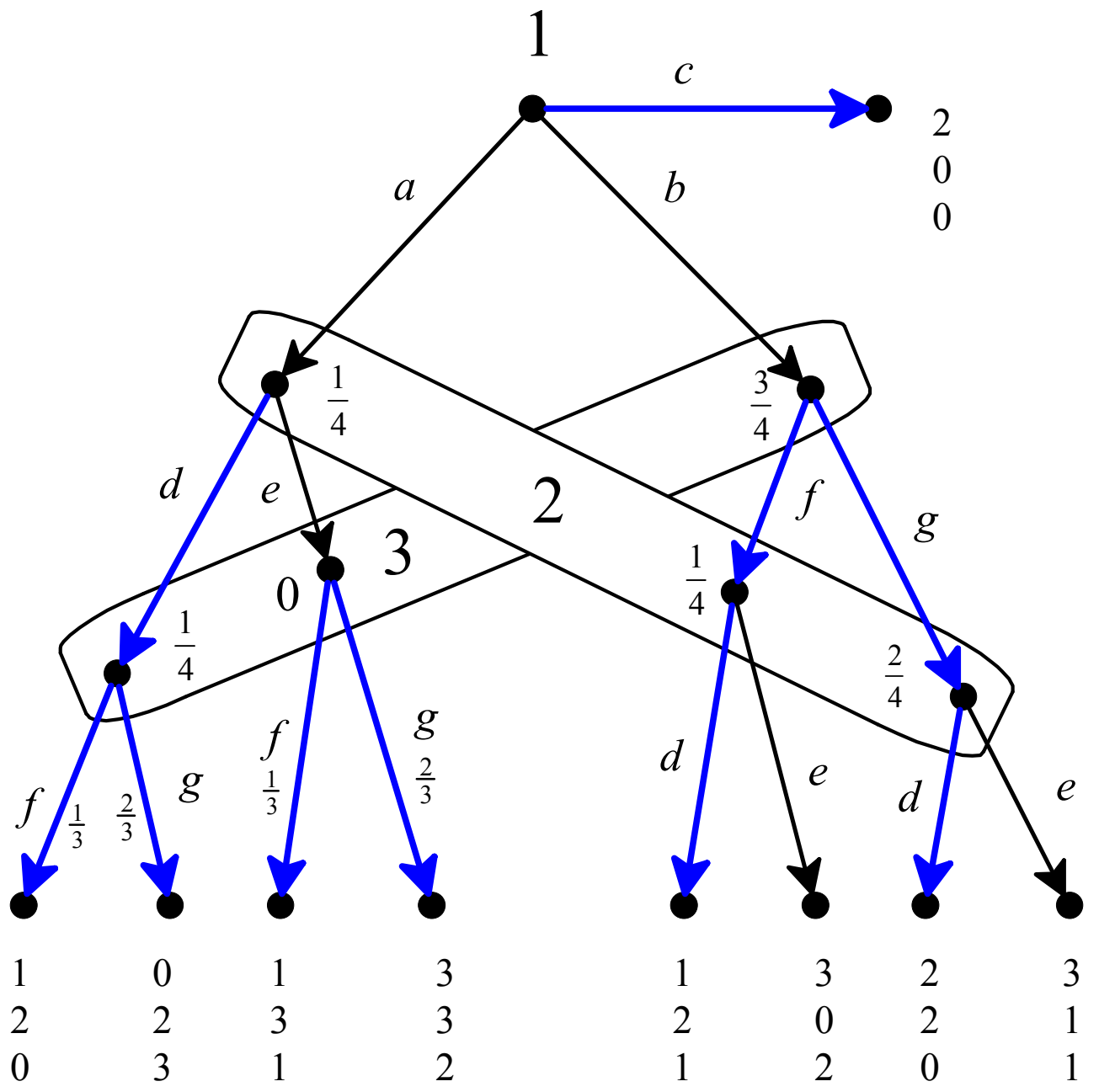

**Figure 12.6**

We showed in the previous section that the assessment under consideration is AGM and Bayesian consistent. Thus we only need to verify sequential rationality. For Player 3 the expected payoff from playing $f$ is $\left(\frac{1}{4}\right)0+(0)1+\left(\frac{3}{4}\right)1=\frac{3}{4}$ and the expected payoff from playing $g$ is $\left(\frac{1}{4}\right)3+(0)2+\left(\frac{3}{4}\right)0=\frac{3}{4}$; thus any mixture of $f$ and $g$ is sequentially rational, in particular the mixture $\begin{pmatrix} f & g \\ \frac{1}{3} & \frac{2}{3} \end{pmatrix}$. For Player 2 the expected payoff from playing $d$ is $\frac{1}{4}\left[\left(\frac{1}{3}\right)2+\left(\frac{2}{3}\right)2\right]+\left(\frac{1}{4}\right)2+\left(\frac{2}{4}\right)2=2$, while the expected payoff from playing $e$ is $\frac{1}{4}\left[\left(\frac{1}{3}\right)3+\left(\frac{2}{3}\right)3\right]+\left(\frac{1}{4}\right)0+\left(\frac{2}{4}\right)1=\frac{5}{4}$; thus $d$ is sequentially rational. For Player 1, $c$ gives a payoff of 2 while $a$ gives a payoff of $\frac{1}{3}$ and $b$ a gives a payoff of $\frac{5}{3}$; thus $c$ is sequentially rational.

We now turn to the properties of Perfect Bayesian equilibria.

**Theorem 12.1** (Bonanno, 2013). Consider a finite extensive-form game and an assessment $(\sigma,\mu)$. If $(\sigma,\mu)$ is a perfect Bayesian equilibrium then (1) $\sigma$ is a subgame-perfect equilibrium and (2) $(\sigma,\mu)$ is a weak sequential equilibrium.

The example of Figure 12.2 (Section 12.1) showed that not every subgame-perfect equilibrium can be part of a perfect Bayesian equilibrium. Thus, by Theorem 12.1, the notion of subgame-perfect equilibrium is a strict refinement of the notion of subgame-perfect equilibrium.

**Theorem 12.2** (Bonanno, 2013). Consider a finite extensive-form game and an assessment $(\sigma, \mu)$. If $(\sigma, \mu)$ is sequential equilibrium then $\sigma$ is a perfect Bayesian equilibrium.

Not every perfect Bayesian equilibrium is a sequential equilibrium. To see this, consider the game of Figure 12.7.

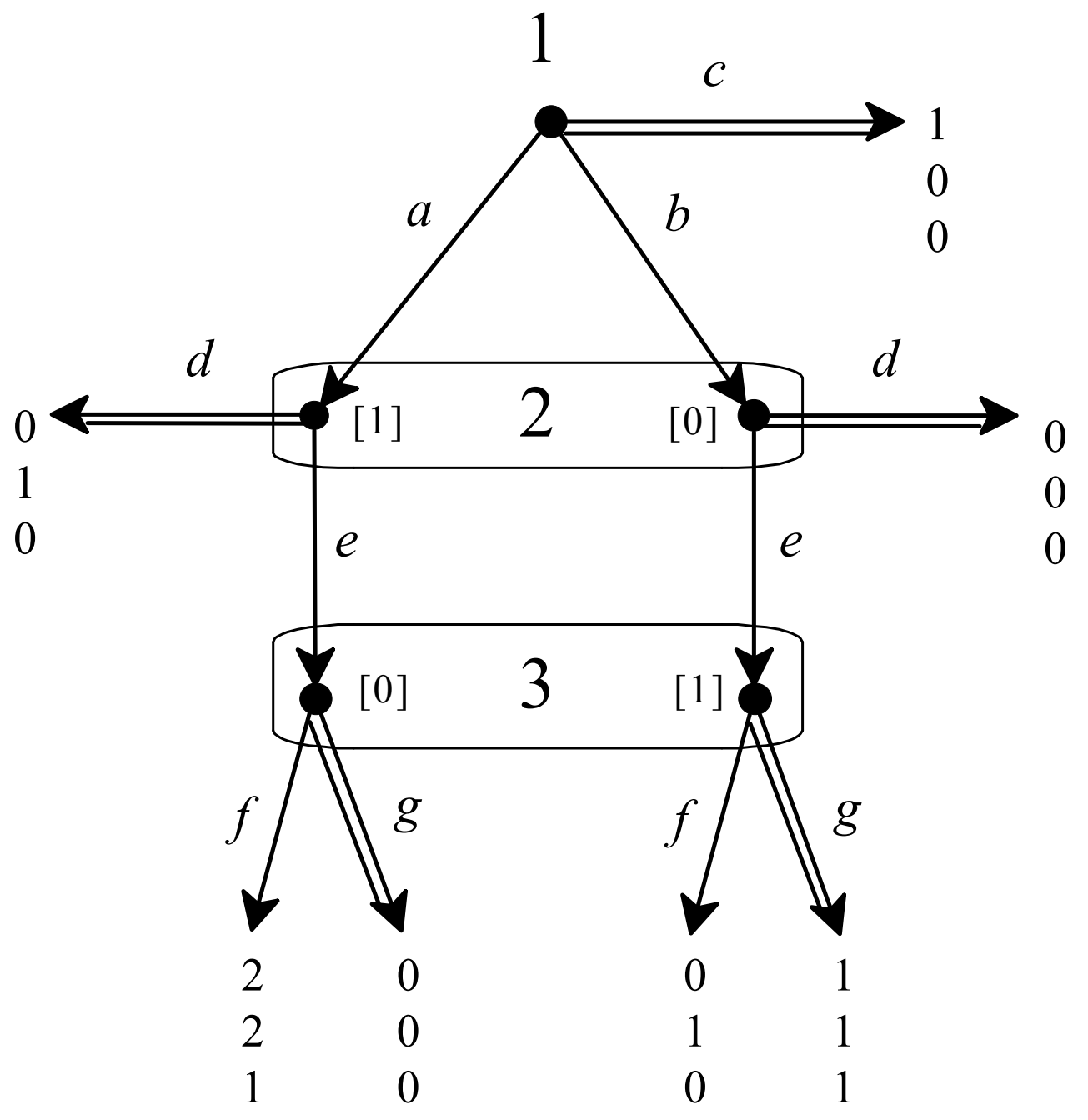

**Figure 12.7**

A perfect Bayesian equilibrium is given by the pure-strategy profile $\sigma = (c, d, g)$ (highlighted by double edges), together with the degenerate beliefs $\mu(a) = \mu(be) = 1$. In fact, $(\sigma, \mu)$ is sequentially rational and, furthermore, it is rationalized by the plausibility order shown in Figure 12.8 and is Bayesian relative to it (the probability distributions on the equivalence classes that contain histories $h$ with $\mu(h) > 0$ are written next to the order):

$$\begin{pmatrix} \emptyset, c & \text{most plausible} \\ a, ad & \\ b, bd & \\ be, beg & \\ ae, aeg & \\ bef & \\ aef & \text{least plausible} \end{pmatrix}, \quad \begin{pmatrix} \nu_{\{\emptyset,c\}}(\emptyset) = 1 \\ \nu_{\{a,ad\}}(a) = 1 \\ - \\ \nu_{\{be,beg\}}(be) = 1 \\ - \\ - \\ - \end{pmatrix}$$

**Figure 12.8**

The belief revision policy encoded in a perfect Bayesian equilibrium can be interpreted either as the point of view of an external observer or as a belief revision policy which is shared by all the players. For example, the perfect Bayesian equilibrium under consideration (for the game of Figure 12.7), namely $\sigma = (c, d, g)$ and $\mu(a) = \mu(be) = 1$, reflects the following belief revision policy: the initial beliefs are that Player 1 will play $c$; conditional on learning that Player 1 did not play *c*, the observer would become convinced that Player 1 played *a* (that is, she would judge *a* to be strictly more plausible than *b*) and would expect Player 2 to play $d$; upon learning that Player 1 did not play *c* and Player 2 did not play *d*, the observer would become convinced that Player 1 played *b* and Player 2 played *e*, hence judging *be* to be strictly more plausible than *ae*, thereby reversing her earlier belief that *a* was strictly more plausible than *b*. Such a belief revision policy is consistent with the AGM rationality axioms (Alchourrón *et al.*, 1985) but is incompatible with the notion of sequential equilibrium. In fact, $(\sigma, \mu)$ is not KW-consistent (Definition 11.1, Chapter 11). To see this, consider an arbitrary sequence $\langle \sigma_n \rangle_{n=1,2,\ldots}$ that converges to $\sigma$: $\sigma_n = \left(\begin{array}{ccc|cc|cc} a & b & c & d & e & f & g \\ p_n & q_n & 1-p_n-q_n & 1-r_n & r_n & t_n & 1-t_n \end{array}\right)$ with $\lim_{n\to\infty} p_n =$ $\lim_{n\to\infty} q_n = \lim_{n\to\infty} r_n = \lim_{n\to\infty} t_n = 0$. Then the corresponding sequence $\langle \mu_n \rangle_{n=1,2,\ldots}$ of beliefs obtained by Bayesian updating is given by

$$\mu_n = \left(\begin{array}{cc|cc} a & b & ae & be \\ \dfrac{p_n}{p_n+q_n} & \dfrac{q_n}{p_n+q_n} & \dfrac{p_n r_n}{p_n r_n + q_n r_n} = \dfrac{p_n}{p_n+q_n} & \dfrac{q_n r_n}{p_n r_n + q_n r_n} = \dfrac{q_n}{p_n+q_n} \end{array}\right).$$

Thus, $\mu_n(a) = \mu_n(ae)$; therefore if $\lim_{n\to\infty} \mu_n(a) = \mu(a) = 1$ then $\lim_{n\to\infty} \mu_n(ae) = 1 \neq \mu(ae) = 0$.

By Theorem 11.2 (Chapter 11), every finite extensive-form game with cardinal payoffs has at least one sequential equilibrium and, by Theorem 12.2, every sequential equilibrium is a perfect Bayesian equilibrium. Thus the following theorem follows as a corollary of these two results.

**Theorem 12.3.** Every finite extensive-form game with cardinal payoffs has at least one perfect Bayesian equilibrium.

The relationship among the different notions of equilibrium introduced so far (Nash equilibrium, subgame-perfect equilibrium, weak sequential equilibrium, perfect Bayesian equilibrium and sequential equilibrium) is shown in the Venn diagram of Figure 12.9.

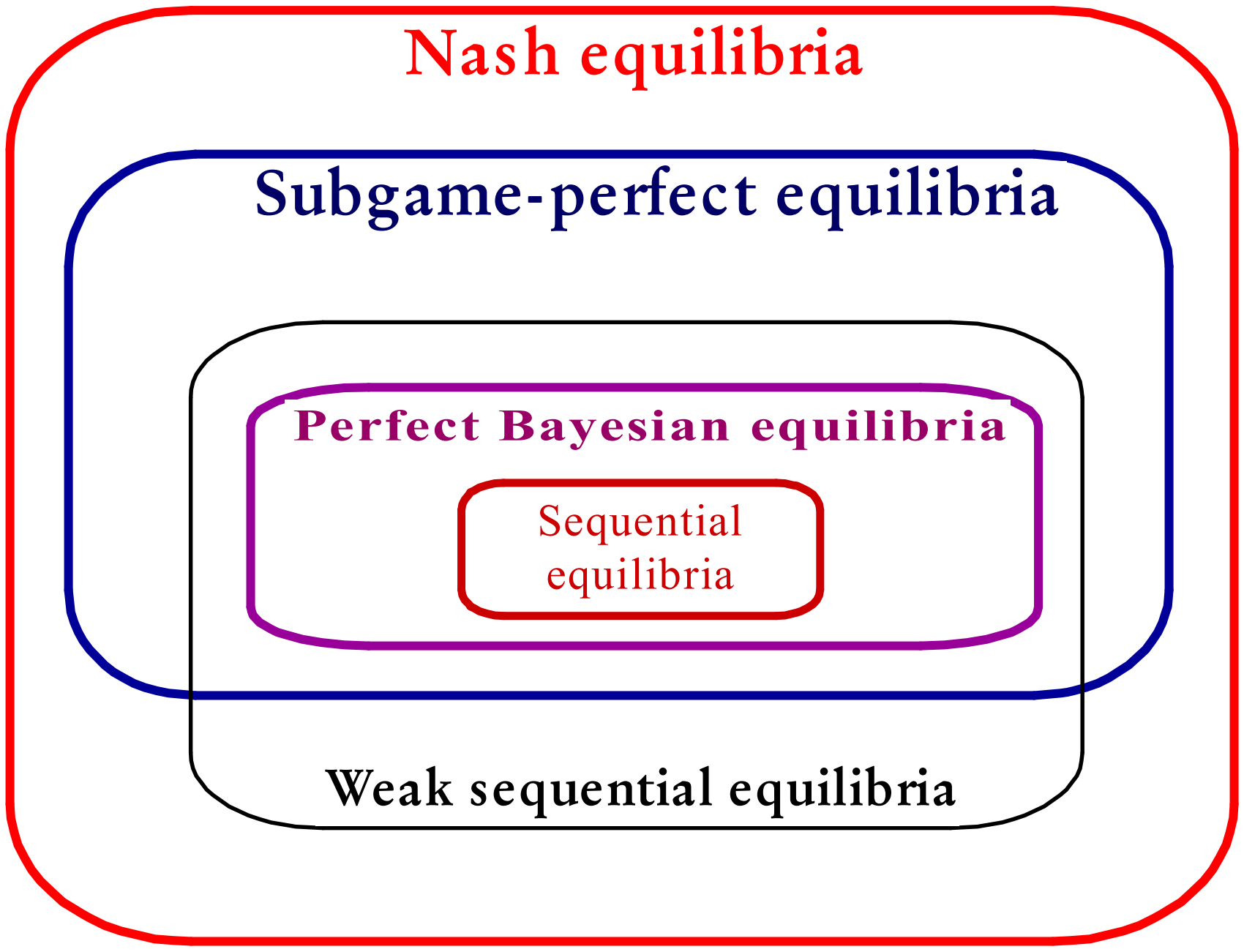

**Figure 12.9**

This is a good time to test your understanding of the concepts introduced in this section, by going through the exercises in Section 12.E.3 of Appendix 12.E at the end of this chapter.

# 12.4 Adding independence

The notion of perfect Bayesian equilibrium imposes relatively mild restrictions on beliefs at information sets that are reached with zero probability. Those restrictions require consistency between the strategy profile $\sigma$ and the system of belief $\mu$, as well as the requirement of Bayesian updating "as long as possible" (that is, also after beliefs have been revised at unreached information sets). The example of Figure 12.7 showed that perfect Bayesian equilibrium is compatible with a belief revision policy that allows a reversal of judgment concerning the behavior of one player after observing an unexpected move by a different player. One might want to rule out such forms of belief revision. In this section we introduce and discuss further restrictions on belief revision that incorporate some form of independence.

Situations like the one illustrated in Figure 12.7 are ruled out by the following restriction on the plausibility order: if $h$ and $h'$ belong to the same information set $(h' \in I(h))$ and $a$ is an action available at $h$ $(a \in A(h))$, then

$$h \precsim h' \text{ if and only if } ha \precsim h'a . \qquad (IND_1)$$

$IND_1$ says that if $h$ is deemed to be at least as plausible as $h'$, then the addition of any action $a$ must preserve this judgment, in the sense that $ha$ must be deemed to be at least as plausible as $h'a$, and *vice versa*. This restriction can be viewed as an independence condition, in the sense that observation of a new action cannot lead to a change in the relative plausibility of prior histories. Any plausibility order that rationalizes the assessment discussed for the game shown in Figure 12.7 violates $IND_1$, since it must be such that $a \prec b$ and also $be \prec ae$.

Another independence condition is the following, which says that if, conditional on history $h$, action $a$ is implicitly judged to be at least as plausible as action $b$ then the same judgment must be made conditional on any other history that belongs to the same information set as $h$: if $h' \in I(h)$ and $a, b \in A(h)$, then

$$ha \precsim hb \text{ if and only if } h'a \precsim h'b . \qquad (IND_2)$$

**Remark 12.1**. The two properties $IND_1$ and $IND_2$ are independent of each other, in the sense that there are plausibility orders that satisfy one of the two properties but not the other (see Exercises 12.7 and 12.8).

Property $IND_1$ is a qualitative property, to which we can add a corresponding quantitative condition on the probabilities: we call a system of beliefs $\mu$ *weakly independent* if it satisfies the following condition (recall that $D_\mu^+ = \{h \in D : \mu(h) > 0\}$):

$$\text{if } h' \in I(h),\ a \in A(h),\ h'a \in I(ha) \text{ and } \{h, h', ha, h'a\} \subseteq D_\mu^+ \quad \text{then } \frac{\mu(h)}{\mu(h')} = \frac{\mu(ha)}{\mu(h'a)} \qquad (IND_3)$$

The game illustrated in Figure 12.10 below shows that it is possible for an assessment $(\sigma, \mu)$ to satisfy Property $IND_1$ even though $\mu$ fails to satisfy Property $IND_3$.

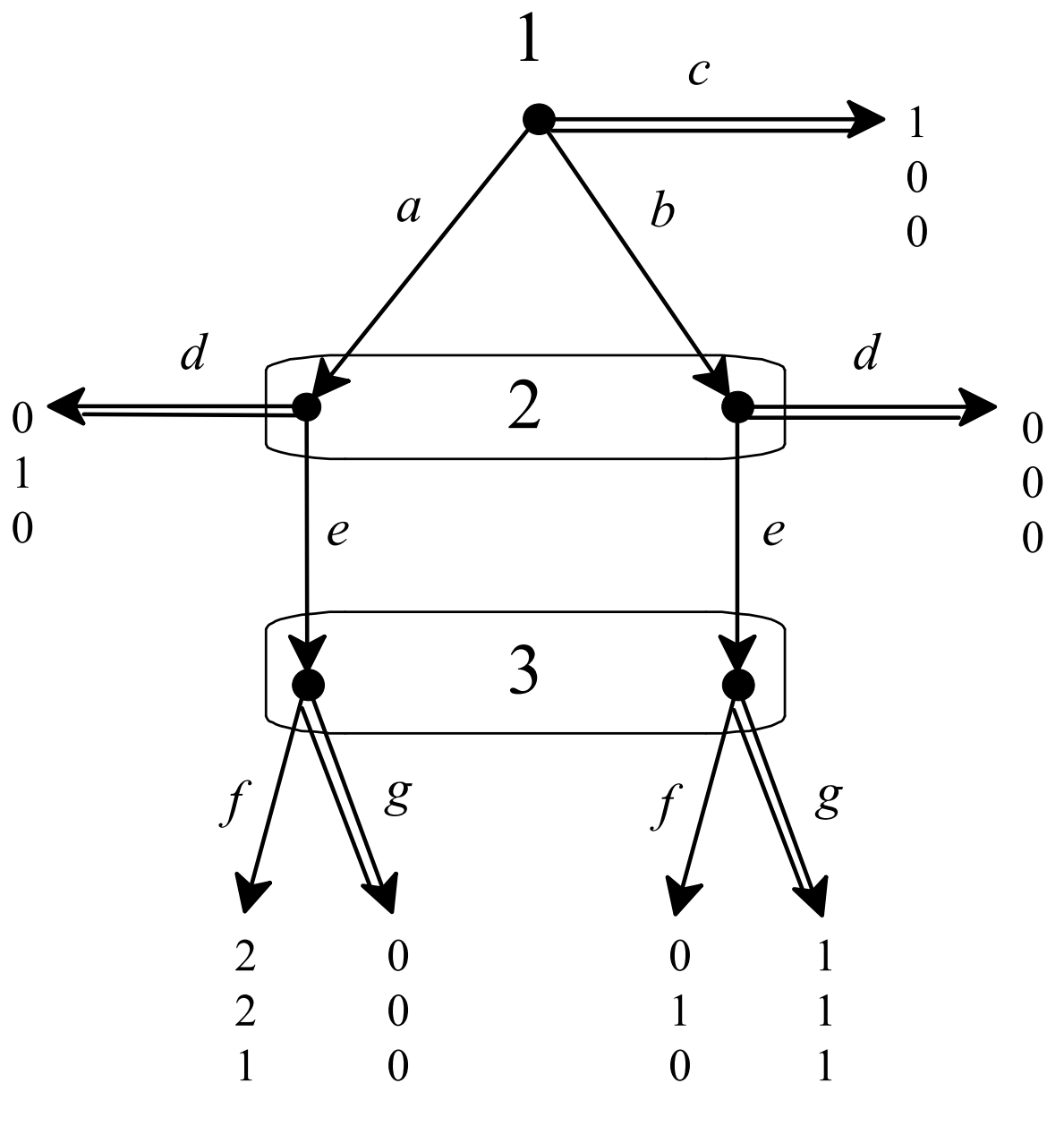

**Figure 12.10**

Consider the assessment consisting of the pure-strategy profile $\sigma = (c, d, g)$ (highlighted by double edges) and the system of beliefs $\mu = \left(\begin{array}{cc|cc} a & b & ae & be \\ \frac{1}{2} & \frac{1}{2} & \frac{1}{4} & \frac{3}{4} \end{array}\right)$.

This assessment is rationalized by the plausibility order

$$\begin{pmatrix} \emptyset, c & \text{most plausible} \\ a, b, ad, bd & \\ ae, be, aeg, beg & \\ aef, bef & \text{least plausible} \end{pmatrix}$$

which satisfies Property $IND_1$, since $a \sim b$ and $ad \sim bd$ and $ae \sim be$. However, $\mu$ fails to satisfy Property $IND_3$ since $\{a, b, ae, be\} \subseteq D_\mu^+$ and $\frac{\mu(a)}{\mu(b)} = 1 \neq \frac{\mu(ae)}{\mu(be)} = \frac{1}{3}$.

**Definition 12.5.** Given a finite extensive-form game, an assessment $(\sigma, \mu)$ is an *independent perfect Bayesian equilibrium* if it is a perfect Bayesian equilibrium that is rationalized by a plausibility order that satisfies properties $IND_1$ and $IND_2$ and, furthermore, $\mu$ satisfies property $IND_3$.

The following theorem is proved in Appendix 12.B at the end of this chapter, using a result given in the next section.

**Theorem 12.4.** Consider a finite extensive-form game and an assessment $(\sigma, \mu)$. If $(\sigma, \mu)$ is a sequential equilibrium then it is an *independent perfect Bayesian equilibrium.*

Exercise 12.10 shows that not every independent perfect Bayesian equilibrium is a sequential equilibrium. Thus the notion of sequential equilibrium is a strict refinement of the notion of independent perfect Bayesian equilibrium, which – in turn – is a strict refinement of perfect Bayesian equilibrium, as illustrated in the Venn diagram of Figure 12.11 below.

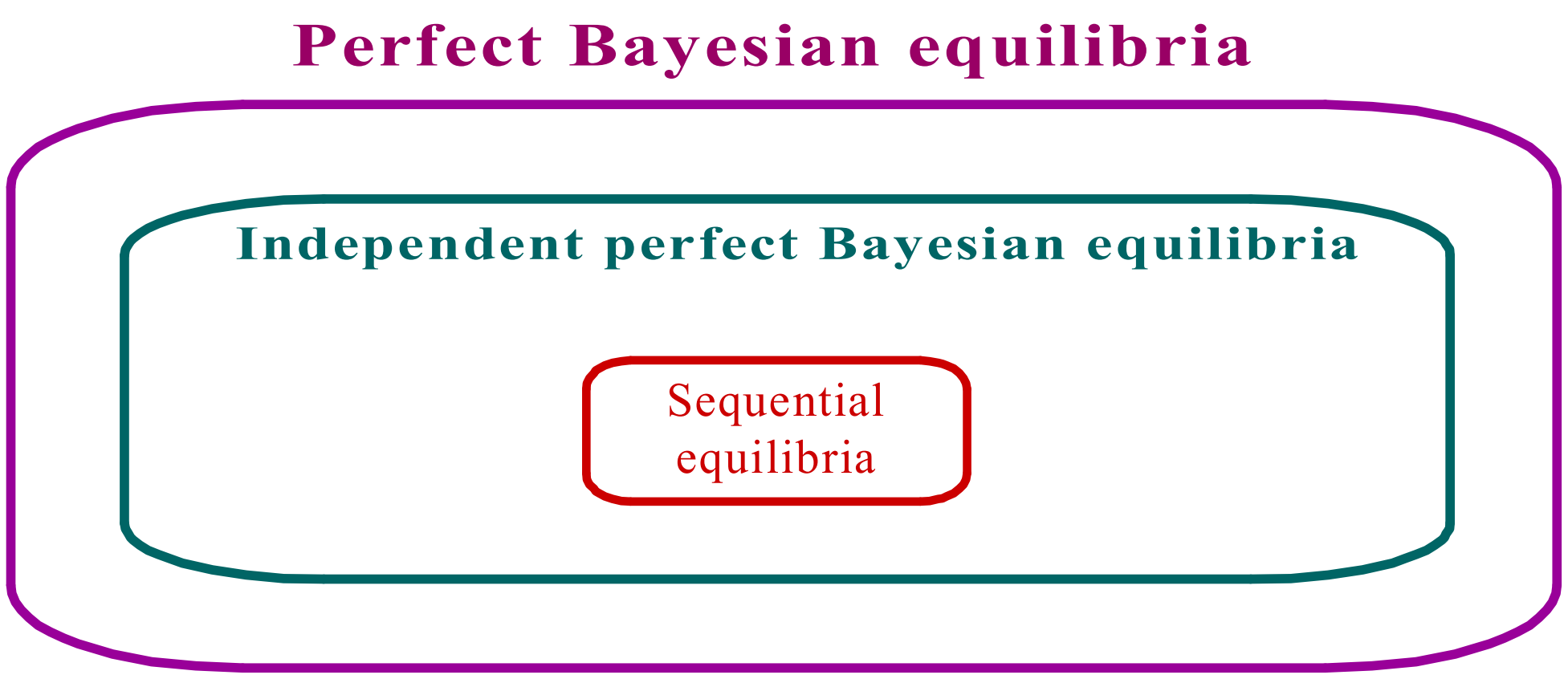

**Figure 12.11**

This is a good time to test your understanding of the concepts introduced in this section, by going through the exercises in Section 12.E.4 of Appendix 12.E at the end of this chapter.

## 12.5 Filling the gap between perfect Bayesian equilibrium and sequential equilibrium

Besides sequential rationality, the notion of perfect Bayesian equilibrium (Definition 12.4) is based on two elements: (1) the qualitative notions of plausibility order and AGM-consistency and (2) the notion of Bayesian consistency relative to the plausibility order. In this section we will show that by strengthening these two conditions one obtains a characterization of sequential equilibrium. The strengthening of the first condition is that the plausibility order that rationalizes the given assessment has a "cardinal" numerical representation that can be interpreted as measuring the plausibility distance between histories in a way that is preserved by the addition of a common action. The strengthening of the second condition imposes uniform consistency on the conditional probability distributions on the equivalence classes of the plausibility order, by requiring that the relative probabilities of two decision histories in the same information set be preserved when a common action is added.

**Definition 12.6.** Given a finite extensive form and a plausibility order $\precsim$ on the set of histories *H*, a function $F: H \to \mathbb{N}$ (where $\mathbb{N}$ denotes the set of non-negative integers) is an *integer-valued representation* of $\precsim$ if, for all $h, h' \in H$,

$$F(h) \le F(h') \text{ if and only if } h \precsim h'.$$

An integer-valued representation *F* of a plausibility order is analogous to an *ordinal* utility function for a preference relation: it is just a numerical representation of the order. Note that, in the case of a plausibility order, we find it more convenient to assign *lower* values to *more* plausible histories (while a utility function assigns higher values to more preferred alternatives or outcomes).

**Remark 12.2.** If $F: H \to \mathbb{N}$ is an integer-valued representation of a plausibility order $\precsim$ on the set of histories *H*, without loss of generality we can assume that $F(\emptyset) = 0$ (recall that $\emptyset$ denotes the null history, which is always one of the most plausible histories).[12]

**Remark 12.3.** Since $H$ is a finite set, an integer-valued representation of a plausibility order $\precsim$ on $H$ always exists. A natural integer-valued representation is the following. Define

$$H_0 = \{h \in H : h \precsim x, \text{ for all } x \in H\}$$

(thus $H_0$ is the set of most plausible histories in *H*), and

$$H_1 = \{h \in H \setminus H_0 \ : h \precsim x, \text{ for all } x \in H \setminus H_0\}$$

(thus $H_1$ is the set of most plausible histories among the ones that remain after removing the set $H_0$ from $H$) and, in general, for every integer $k \ge 1$ define $H_k = \{h \in H \setminus (H_0 \cup \ldots \cup H_{k-1}) : h \precsim x, \text{ for all } x \in H \setminus (H_0 \cup \ldots \cup H_{k-1})\}$. Since $H$ is finite, there is an $m \in \mathbb{N}$ such that $\{H_0, \ldots, H_m\}$ is a partition of $H$ and, for every $j, k \in \mathbb{N}$, with $j < k \le m$, and for every $h, h' \in H$, if $h \in H_j$ and $h' \in H_k$ then $h \prec h'$. Define $\hat{F}: H \to \mathbb{N}$ as follows: $\hat{F}(h) = k$ if and only if $h \in H_k$; then the function $\hat{F}$ so defined is an integer-valued representation of $\precsim$ and

[12] Let $\bar{F}: H \to \mathbb{N}$ be an integer-valued representation of a plausibility order $\precsim$ and define $F: H \to \mathbb{N}$ as follows: $F(h) = \bar{F}(h) - \bar{F}(\emptyset)$. Then *F* is also an integer-valued representation of $\precsim$ and $F(\emptyset) = 0$.

$\hat{F}(\emptyset) = 0$ .

**Definition 12.7.** Given a finite extensive form, a plausibility order $\precsim$ on the set of histories $H$ is *choice measurable* if it has at least one integer-valued representation $F$ that satisfies the following property: for all $h \in D$, for all $h' \in I(h)$ and for all $a \in A(h)$,[13]

$$F(h) - F(h') = F(ha) - F(h'a). \tag{CM}$$

If we think of $F$ as measuring the "plausibility distance" between histories, then we can interpret *choice measurability* (CM) as a distance-preserving condition: the plausibility distance between two histories in the same information set is preserved by the addition of the same action.

For example, consider the extensive form of Figure 12.12.

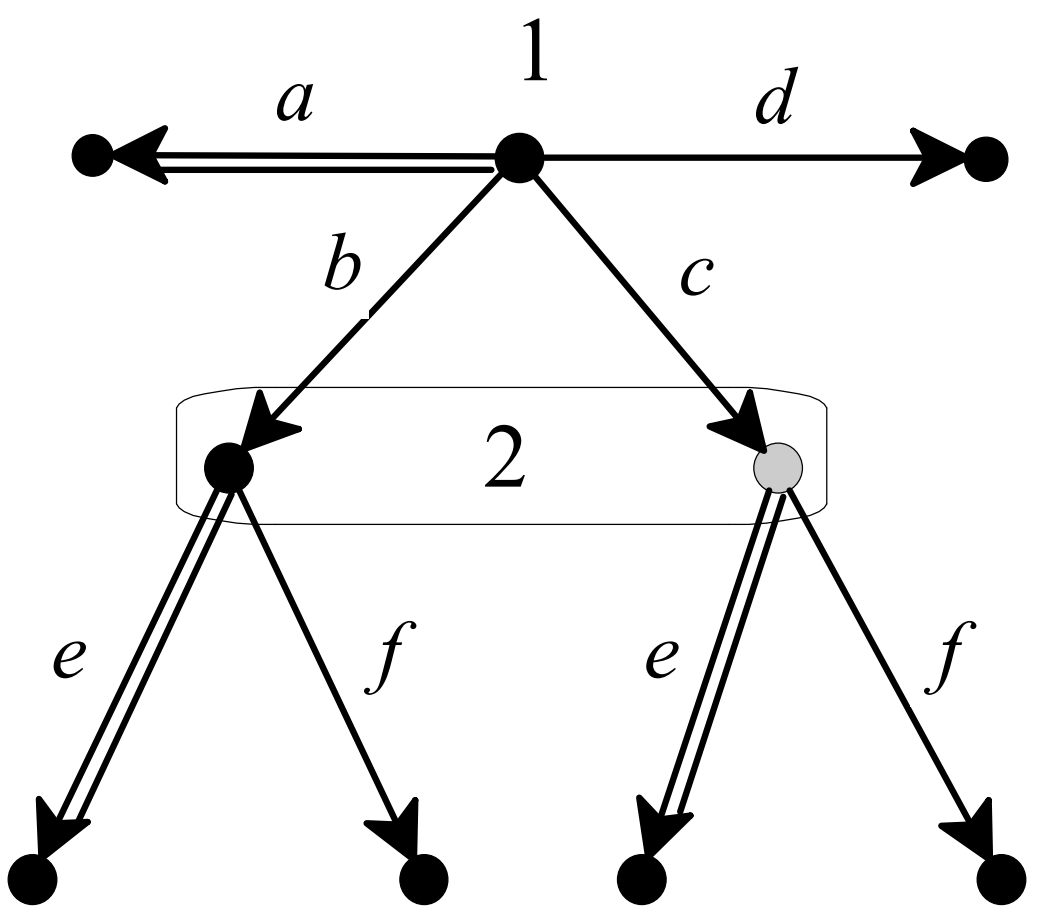

**Figure 12.12**

The assessment consisting of the pure-strategy profile $\sigma = (a, e)$ (highlighted by double edges) and the system of beliefs $\mu = \begin{pmatrix} b & c \\ 1 & 0 \end{pmatrix}$ (the grey node in Player 2's information set represents the history which is assigned zero probability) is rationalized by the plausibility order shown in Figure 12.13 below together with two integer-valued representations ($\hat{F}$ is the one described in Remark

[13] Recall that $D$ is the set of decision histories, $I(h)$ the information set that contains decision history $h$ and $A(h)$ the set of actions available at decision history $h$.

12.3). While $\hat{F}$ does not satisfy Property CM (since $\hat{F}(c) - \hat{F}(b) = 3 - 1 = 2$ $\neq \hat{F}(cf) - \hat{F}(bf) = 5 - 2 = 3$), $F$ does satisfy Property CM. Since there is at least one integer-valued representation that satisfies Property CM, by Definition 12.7 the plausibility order of Figure 12.13 is choice measurable.

| $\precsim$: | $\hat{F}$: | $F$: |
|---|---|---|
| $\emptyset, a$ | 0 | 0 |
| $b, be$ | 1 | 1 |
| $bf$ | 2 | 3 |
| $c, ce$ | 3 | 4 |
| $d$ | 4 | 5 |
| $cf$ | 5 | 6 |

**Figure 12.13**

On the other hand, the plausibility order shown in Figure 12.14 (which reproduces Figure 12.8 and refers to the game of Figure 12.7) is not choice measurable, because any integer-valued representation $F$ must be such that $F(b) - F(a) > 0$ while $F(be) - F(ae) < 0$.

$$\begin{pmatrix} \emptyset, c \\ a, ad \\ b, bd \\ be, beg \\ ae, aeg \\ bef \\ aef \end{pmatrix}$$

**Figure 12.14**

Rationalizability by a choice measurable plausibility order is one of the two properties that are needed to fill the gap between perfect Bayesian equilibrium and sequential equilibrium. The second property is a strengthening of the notion of Bayesian consistency (Definition 12.3).

To set the stage for the next definition, let $(\sigma, \mu)$ be an assessment which is rationalized by a plausibility order $\precsim$. As before, let $D^+_\mu$ be the set of decision histories to which $\mu$ assigns positive probability: $D^+_\mu = \{h \in D : \mu(h) > 0\}$. Let $\mathcal{E}^+_\mu$ be the set of equivalence classes of $\precsim$ that have a non-empty intersection with $D^+_\mu$. Then $\mathcal{E}^+_\mu$ is a non-empty, finite set. Suppose that $(\sigma, \mu)$ is Bayesian relative to $\precsim$ (Definition 12.3) and let $\{\nu_E\}_{E \in \mathcal{E}^+_\mu}$ be a collection of probability distributions that satisfy the properties of Definition 12.3. We call a probability distribution $\nu : D \to (0,1]$ (recall that $D$ is the set of all decision histories) a *full-support common prior* of $\{\nu_E\}_{E \in \mathcal{E}^+_\mu}$ if, for every $E \in \mathcal{E}^+_\mu$, the probability distribution $\nu_E$ coincides with the distribution obtained by conditioning $\nu$ on the set $E \cap D^+_\mu$, that is, for all $h \in E \cap D^+_\mu$, $\nu_E(h) = \nu(h \mid E \cap D^+_\mu) \overset{def}{=} \dfrac{\nu(h)}{\sum_{h' \in E \cap D^+_\mu} \nu(h')}$.

Note that a full support common prior $\nu$ assigns positive probability to *all* decision histories, not only to those in $D^+_\mu$.

**Remark 12.4.** Since any two elements of $\mathcal{E}^+_\mu$ are mutually disjoint, a full support common prior always exists; indeed, there is an infinite number of them. Furthermore, it can be shown (see Bonanno, 2015) that there always exists a full-support common prior $\nu$ that satisfies the following properties: if $a \in A(h)$ and $ha \in D$, then (1) $\nu(ha) \le \nu(h)$ and (2) if $\sigma(a) > 0$ then $\nu(ha) = \nu(h) \times \sigma(a)$.

The following definition requires that, among the many full-support common priors, there be one that, besides the properties of Remark 12.4, satisfies the additional property that the relative likelihood of any two histories in the same information set be preserved by the addition of the same action.

**Definition 12.8.** Consider an extensive form. Let $(\sigma,\mu)$ be an assessment that is rationalized by a plausibility order $\precsim$ on the set of histories $H$ and is Bayesian relative to it and let $\{\nu_E\}_{E\in\mathcal{E}_\mu^+}$ be a collection of probability distributions that satisfy the properties of Definition 12.3. We say that $(\sigma,\mu)$ is *uniformly Bayesian relative to* $\precsim$ if there exists a full-support common prior $\nu : D \to (0,1]$ of $\{\nu_E\}_{E\in\mathcal{E}_\mu^+}$ that satisfies the following properties:

UB1. If $a\in A(h)$ and $ha\in D$ then (1) $\nu(ha)\le \nu(h)$ and (2) if $\sigma(a)>0$ then $\nu(ha)=\nu(h)\times\sigma(a)$.

UB2. If $a\in A(h)$, $h$ and $h'$ belong to the same information set and $ha, h'a\in D$ then $\dfrac{\nu(h)}{\nu(h')}=\dfrac{\nu(ha)}{\nu(h'a)}$.

UB1 is the property of Remark 12.4, which can always be satisfied by an appropriate choice of a full-support common prior. UB2 requires that the relative probability, according to the common prior $\nu$, of any two histories that belong to the same information set remain unchanged by the addition of the same action.

Choice measurability and uniform Bayesian consistency are independent properties. For example, the perfect Bayesian equilibrium $\sigma=(c,d,g)$ and $\mu(a)=\mu(be)=1$ of the game of Figure 12.7 is such that any plausibility order that rationalizes it cannot be choice measurable[14] and yet $(\sigma,\mu)$ is uniformly Bayesian relative to the plausibility order shown in Figure 12.8 that rationalizes it.[15] On the other hand, the perfect Bayesian equilibrium $\sigma=(c,d,g)$, $\mu=\left(\begin{array}{cc|cc} a & b & ae & be \\ \frac{3}{4} & \frac{1}{4} & \frac{1}{4} & \frac{3}{4}\end{array}\right)$ of the game of Figure 12.10 is rationalized by the choice measurable plausibility order $\begin{pmatrix} \emptyset, c & \text{most plausible} \\ a, b, ad, bd & \\ ae, be, aeg, beg & \\ aef, bef & \text{least plausible} \end{pmatrix}$ but it cannot be uniformly Bayesian relative to any plausibility order that rationalizes it.[16]

[14] Because, by P2 of Definition 12.2, any such plausibility order $\precsim$ would have to satisfy $a\prec b$ and $be\prec ae$, so that any integer-valued representation $F$ of it would be such that $F(b)-F(a)>0$ and $F(be)-F(ae)<0$

[15] As can be seen by taking $\nu$ to be the uniform distribution over the set $D=\{\varnothing,a,b,ae,be\}$ (that is, $\nu(h)=\frac{1}{5}$ for every $h\in D$): UB1 is clearly satisfied and UB2 is also satisfied, since $\dfrac{\nu(a)}{\nu(b)}=\dfrac{\frac{1}{5}}{\frac{1}{5}}=\dfrac{\nu(ae)}{\nu(be)}$.

We can now state the main result of this section, namely that choice measurability and uniform Bayesian consistency are necessary and sufficient for a perfect Bayesian equilibrium to be a sequential equilibrium.

**Theorem 12.5** (Bonanno, 2015)**.** Consider a finite extensive-form game and an assessment $(\sigma, \mu)$. The following are equivalent:

(A) $(\sigma, \mu)$ is a sequential equilibrium,

(B) $(\sigma, \mu)$ is a perfect Bayesian equilibrium that is rationalized by a choice-measurable plausibility order and is uniformly Bayesian relative to it.

Thus Theorem 12.5 provides a characterization (or understanding) of sequential equilibrium which is free of the questionable requirement of taking a limit of sequences of completely mixed strategies and associated systems of beliefs (see the remarks in Chapter 11, Section 11.3).

Using Theorem 12.5 one can prove Theorem 12.4, namely that a sequential equilibrium is an independent perfect Bayesian equilibrium. The proof is given in Appendix B at the end of this chapter.

As an application of Theorem 12.5, consider the extensive-form game of Figure 12.15 below.

---

[16] Because, by P2 of Definition 12.2, any such plausibility order would have to have the following equivalence classes: $E = \{a, b, ad, bd\}$ and $F = \{ae, be, aeg, beg\}$ (thus $E \cap D_\mu^+ = \{a, b\}$ and $F \cap D_\mu^+ = \{ae, be\}$). Hence, if $\nu$ is any common prior then $\nu_E(a) = \frac{\nu(a)}{\nu(a)+\nu(b)}$ and $\nu_E(b) = \frac{\nu(b)}{\nu(a)+\nu(b)}$. By B3 of Definition 12.3, $\mu(a) = \frac{\nu_E(a)}{\nu_E(a)+\nu_E(b)}$ and $\mu(b) = \frac{\nu_E(b)}{\nu_E(a)+\nu_E(b)}$. Thus $\frac{\nu(a)}{\nu(b)} = \frac{\nu_E(a)}{\nu_E(b)} = \frac{\mu(a)}{\mu(b)} = 3$ and, similarly, $\frac{\nu(ae)}{\nu(be)} = \frac{\nu_F(ae)}{\nu_F(be)} = \frac{\mu(ae)}{\mu(be)} = \frac{1}{3}$, yielding a violation of UB2 of Definition 12.8.

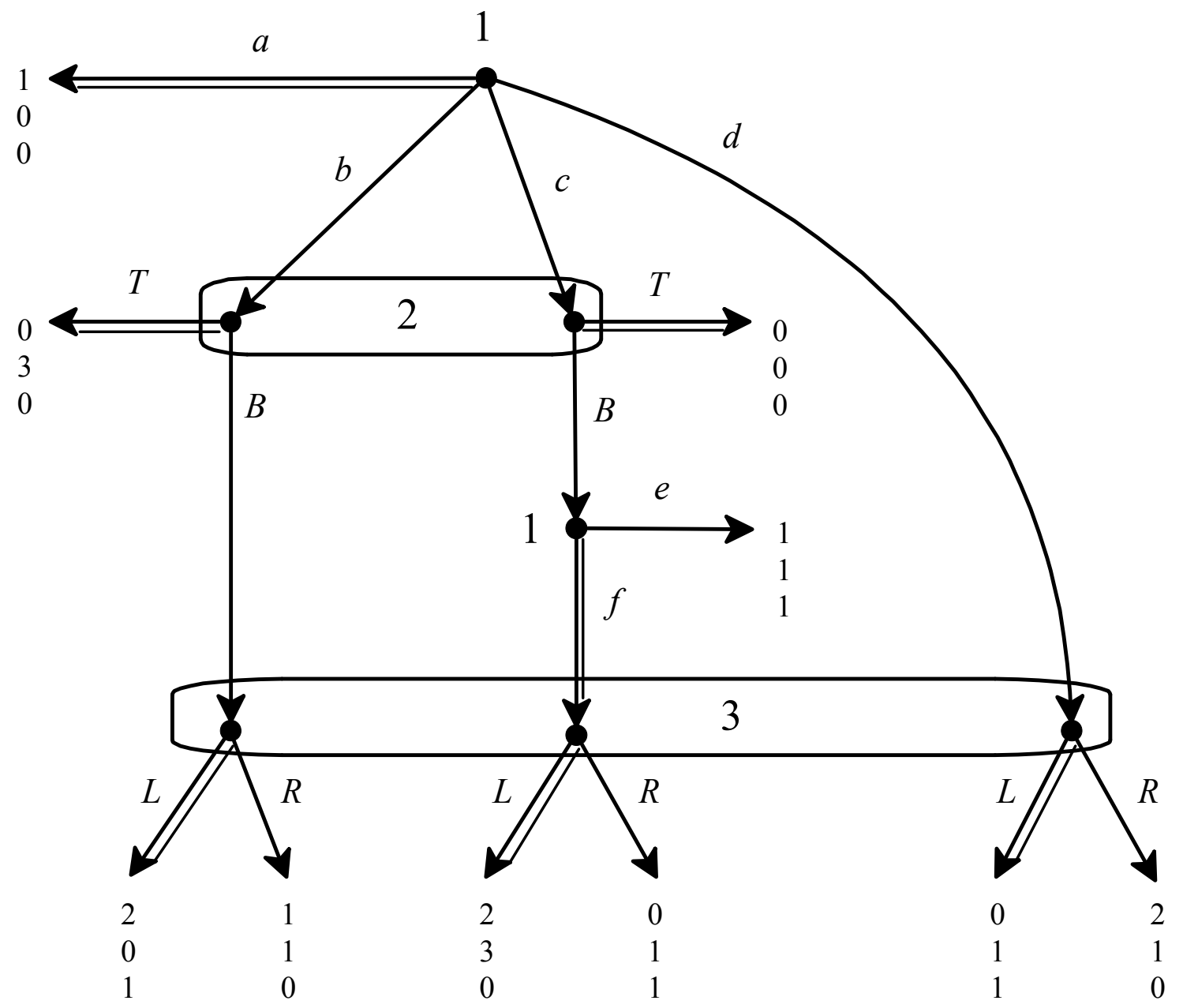

**Figure 12.15**

Using Theorem 12.5 let us show that the assessment $(\sigma,\mu)$ with $\sigma=(a,T,f,L)$ and $\mu=\left(\begin{array}{cc|ccc} b & c & bB & cBf & d \\ \frac{7}{10} & \frac{3}{10} & \frac{7}{18} & \frac{3}{18} & \frac{8}{18}\end{array}\right)$ is a sequential equilibrium.

First we check sequential rationality. For Player 1 at the root the possible payoffs are: 1 with *a*, 0 with *b*, 0 with *c* (more precisely, with either *ce* or *cf*) and 0 with *d*; thus *a* is sequentially rational. For Player 1 at history *cB* the possible payoffs are: 1 with *e* and 2 with *f*; thus *f* is sequentially rational. For Player 2 the possible payoffs are: $\frac{7}{10}(3)+\frac{3}{10}(0)=\frac{21}{10}$ with *T* and $\frac{7}{10}(0)+\frac{3}{10}(3)=\frac{9}{10}$ with *B*; thus *T* is sequentially rational. For Player 3 the possible payoffs are: $\frac{7}{18}(1)+\frac{3}{18}(0)+\frac{8}{18}(1)=\frac{15}{18}$ with *L* and $\frac{7}{18}(0)+\frac{3}{18}(1)+\frac{8}{18}(0)=\frac{3}{18}$ with *R*; thus *L* is sequentially rational. By Theorem 12.5 it only remains to show that $(\sigma,\mu)$ is rationalized by a choice-measurable plausibility order and is uniformly Bayesian relative to it. Indeed, $(\sigma,\mu)$ is rationalized by the plausibility order shown in Figure 12.16 together with a choice-measurable integer-valued representation *F*.

| $\precsim$ : | $F$ : |
|---|---|
| $\emptyset, a$ | 0 |
| $b, c, bT, cT$ | 1 |
| $d, bB, cB, cBf, dL, bBL, cBfL$ | 2 |
| $dR, bBR, cBfR, cBe$ | 3 |

**Figure 12.16**

To see that $(\sigma,\mu)$ is uniformly Bayesian relative to this plausibility order, let $E_1, E_2$ and $E_3$ be the top three equivalence classes of the order and consider the following probability distributions, which satisfy the properties of Definition 12.3 (so that $(\sigma,\mu)$ is Bayesian relative to this plausibility order): $\nu_{E_1}(\emptyset)=1$, $\nu_{E_2}=\begin{pmatrix} b & c \\ \frac{7}{10} & \frac{3}{10} \end{pmatrix}$ and $\nu_{E_3}=\begin{pmatrix} bB & cB & cBf & d \\ \frac{7}{21} & \frac{3}{21} & \frac{3}{21} & \frac{8}{21} \end{pmatrix}$; then a full support common prior that satisfies the properties of Definition 12.8 is the one shown in Figure 12.17 below:

$$\nu=\begin{pmatrix} \emptyset & b & bB & c & cB & cBf & d \\ \frac{9}{40} & \frac{7}{40} & \frac{7}{40} & \frac{3}{40} & \frac{3}{40} & \frac{3}{40} & \frac{8}{40} \end{pmatrix}$$

**Figure 12.17**

This is a good time to test your understanding of the concepts introduced in this section, by going through the exercises in Section 12.E.5 of Appendix 12.E at the end of this chapter.

# Appendix 12.A:
## History-based definition of extensive-form game

If $A$ is a set, we denote by $A^*$ the set of finite sequences in $A$. If $h = \langle a_1,...,a_k \rangle \in A^*$ and $1 \le j \le k$, the sequence $h' = \langle a_1,...,a_j \rangle$ is called a *prefix* of $h$ (a *proper prefix* of $h$ if $j < k$). If $h = \langle a_1,...,a_k \rangle \in A^*$ and $a \in A$, we denote the sequence $\langle a_1,...,a_k,a \rangle \in A^*$ by $ha$. A *finite extensive form* is a tuple $\langle A, H, N, \iota, \{\approx_i\}_{i \in N} \rangle$ whose elements are:

- A finite set of actions $A$.
- A finite set of histories $H \subseteq A^*$ which is closed under prefixes (that is, if $h \in H$ and $h' \in A^*$ is a prefix of $h$, then $h' \in H$). The null (or empty) history $\langle \rangle$, denoted by $\emptyset$, is an element of $H$ and is a prefix of every history (the null history $\emptyset$ represents the root of the tree). A history $h \in H$ such that, for every $a \in A$, $ha \notin H$, is called a *terminal history*. The set of terminal histories is denoted by $Z$. $D = H \setminus Z$ denotes the set of non-terminal or *decision histories*. For every decision history $h \in D$, we denote by $A(h)$ the set of actions available at $h$, that is, $A(h) = \{a \in A : ha \in H\}$.
- A finite set $N = \{1,...,n\}$ of players. In some cases there is also an additional, fictitious, player called *chance*.
- A function $\iota : D \to N \cup \{chance\}$ that assigns a player to each decision history. Thus $\iota(h)$ is the player who moves at history $h$. A game is said to be *without chance moves* if $\iota(h) \in N$ for every $h \in D$. For every $i \in N \cup \{chance\}$, let $D_i = \iota^{-1}(i)$ be the set of histories assigned to Player $i$. Thus $\{D_{chance}, D_1,...,D_n\}$ is a partition of $D$. If history $h$ is assigned to chance ($\iota(h) = chance$), then a probability distribution over $A(h)$ is given that assigns positive probability to every $a \in A(h)$.
- For every player $i \in N$, $\approx_i$ is an equivalence relation on $D_i$. The interpretation of $h \approx_i h'$ is that, when choosing an action at history $h \in D_i$, Player $i$ does not know whether she is moving at $h$ or at $h'$. The equivalence class of $h \in D_i$ is denoted by $I_i(h)$ and is called an *information set* of Player $i$; thus $I_i(h) = \{h' \in D_i : h \approx_i h'\}$. The following restriction applies: if $h' \in I_i(h)$ then $A(h') = A(h)$, that is, the set of actions available to a player is the same at any two histories that belong to the same information set of that player.

- The following property, known as *perfect recall*, is assumed: for every Player $i \in N$, if $h_1, h_2 \in D_i$, $a \in A(h_1)$ and $h_1 a$ is a prefix of $h_2$ then, for every $h' \in I_i(h_2)$, there exists an $h \in I_i(h_1)$ such that $ha$ is a prefix of $h'$. Intuitively, perfect recall requires a player to remember what she knew in the past and what actions she took previously.
- Given an extensive form, one obtains an *extensive game* by adding, for every Player $i \in N$, a *utility* (or *payoff*) *function* $U_i : Z \to \mathbb{R}$ (where $\mathbb{R}$ denotes the set of real numbers; recall that *Z* is the set of terminal histories).

Figure 12.A.1 below shows an extensive form without chance moves where $A = \{a,b,s,c,d,e,f,g,h,m,n\}$, $D = \{\emptyset, a, b, ac, ad, acf, ade, adf\}$ (to simplify the notation we write *a* instead of $\langle \emptyset, a \rangle$, *ac* instead of $\langle \emptyset, a, c \rangle$, etc.),
$Z = \{s, ace, acfg, acfh, adeg, adeh, adfm, adfn, bm, bn\}$, $H = D \cup Z$,
$A(\emptyset) = \{a,b,s\}$, $A(a) = \{c,d\}$, $A(ac) = A(ad) = \{e, f\}$,
$A(acf) = A(ade) = \{g,h\}$, $A(adf) = A(b) = \{m,n\}$,
$\iota(\emptyset) = 1$, $\iota(a) = 2$, $\iota(ac) = \iota(ad) = 3$, $\iota(acf) = \iota(ade) = \iota(adf) = \iota(b) = 4$,
$\approx_1 = \{(\emptyset,\emptyset)\}$, $\approx_2 = \{(a,a)\}$, $\approx_3 = \{(ac,ac),(ac,ad),(ad,ac),(ad,ad)\}$ and
$\approx_4 = \{(acf,acf),(acf,ade),(ade,acf),(ade,ade),(adf,adf),(adf,b),(b,adf),(b,b)\}$.
The information sets containing more than one history (for example, $I_4(b) = \{adf, b\}$) are shown as rounded rectangles. The root of the tree represents the null history $\emptyset$.

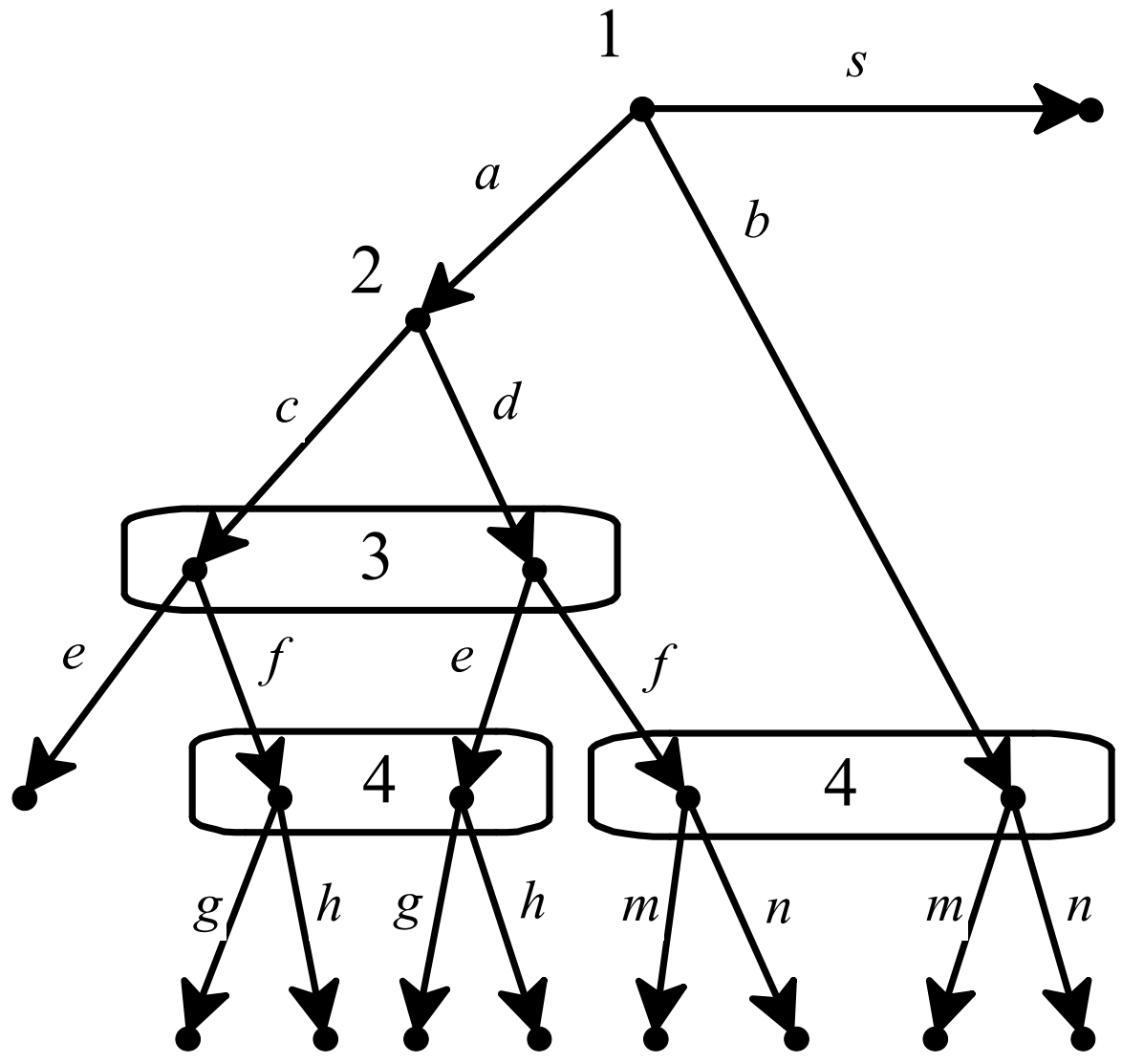

**Figure 12.A.1**

If $h$ and $h'$ are decision histories not assigned to chance, we write $h' \in I(h)$ as a short-hand for $h' \in I_{\iota(h)}(h)$. Thus $h' \in I(h)$ means that $h$ and $h'$ belong to the same information set (of the player who moves at $h$). If $h$ is a history assigned to chance, we use the convention that $I(h) = \{h\}$.

Given an extensive form, a *pure strategy* of Player $i \in N$ is a function that associates with every information set of Player $i$ an action at that information set, that is, a function $s_i : D_i \to A$ such that (1) $s_i(h) \in A(h)$ and (2) if $h' \in I_i(h)$ then $s_i(h') = s_i(h)$. For example, one of the pure strategies of Player 4 in the extensive form illustrated in Figure 12.A.1 is $s_4(acf) = s_4(ade) = g$ and $s_4(adf) = s_4(b) = m$. A *behavior strategy* of player $i$ is a collection of probability distributions, one for each information set, over the actions available at that information set; that is, a function $\sigma_i : D_i \to \Delta(A)$ (where $\Delta(A)$ denotes the set of probability distributions over $A$) such that (1) $\sigma_i(h)$ is a probability distribution over $A(h)$ and (2) if $h' \in I_i(h)$ then $\sigma_i(h') = \sigma_i(h)$. If the game does not have chance moves, we define a behavior strategy *profile* as an $n$-tuple $\sigma = (\sigma_1, ..., \sigma_n)$ where, for every $i \in N$, $\sigma_i$ is a behavior strategy of Player $i$. If the game has chance moves then we use the convention that a behavior strategy profile is an $(n+1)$-tuple $\sigma = (\sigma_1, ..., \sigma_n, \sigma_{chance})$ where, if $h$ is a history assigned to chance and $a \in A(h)$, then $\sigma_{chance}(h)(a)$ is the probability associated with $a$. When there is no risk of ambiguity (e.g. because no action is assigned to more than one information set) we shall denote by $\sigma(a)$ the probability assigned to action $a$ by the relevant component of the strategy profile $\sigma$.[17] Note that a pure strategy is a special case of a behavior strategy where each probability distribution is degenerate. A behavior strategy is *completely mixed at history* $h \in D$ if, for every $a \in A(h)$, $\sigma(a) > 0$. For example, in the extensive form of Figure 12.A.1 a possible behavior strategy for Player 1 is $\begin{pmatrix} a & b & s \\ 0 & 0 & 1 \end{pmatrix}$, which can be more simply denoted by $s$ (and which coincides with a pure strategy of Player 1) and a possible behavior strategy of Player 2 is $\begin{pmatrix} c & d \\ \frac{1}{3} & \frac{2}{3} \end{pmatrix}$ (which is a completely mixed strategy).

---

[17] If $h \in D_i$ and $\sigma_i$ is the $i^{th}$ component of $\sigma$, then $\sigma_i(h)$ is a probability distribution over $A(h)$; thus if $a \in A(h)$ then $\sigma_i(h)(a)$ is the probability assigned to action $a$ by $\sigma_i(h)$: we denote $\sigma_i(h)(a)$ more succinctly by $\sigma(a)$ .

# Appendix 12.B:
## Proof of Theorem 12.4 using Theorem 12.5

We prove Theorem 12.4, namely that a sequential equilibrium is an independent perfect Bayesian equilibrium, using Theorem 12.5. To simplify the notation, we will assume that no action is available at more than one information set, that is, that if $a \in A(h) \cap A(h')$ then $h' \in I(h)$ (this is without loss of generality, because we can always rename some of the actions).

Let $(\sigma, \mu)$ be a sequential equilibrium. First we show that Property $IND_1$ is satisfied. By Theorem 12.5 there is a choice measurable plausibility order $\precsim$ that rationalizes $(\sigma, \mu)$. Let $F$ be an integer-valued representation of $\precsim$ that satisfies Property CM of Definition 12.6 and let $h' \in I(h)$ and $a \in A(h)$. Suppose that $h \precsim h'$; then $F(h) \le F(h')$, that is, $F(h) - F(h') \le 0$ and thus, by CM, $F(ha) - F(h'a) \le 0$, that is, $F(ha) \le F(h'a)$, which implies that $ha \precsim h'a$; conversely, if $ha \precsim h'a$ then $F(ha) - F(h'a) \le 0$ and, by Property CM, $F(h) - F(h') \le 0$, which implies that $h \precsim h'$. Thus Property $IND_1$ is satisfied.

Next we shows that Property CM is equivalent to the following property:

$$\text{if } h' \in I(h) \text{ and } a, b \in A(h) \text{ then } F(hb) - F(ha) = F(h'b) - F(h'a) \qquad (\blacklozenge)$$

Let $\precsim$ be a plausibility order on the set of histories $H$ and let $F : H \to \mathbb{N}$ be an integer-valued representation of $\precsim$ that satisfies Property CM of Definition 12.6. Without loss of generality (see Remark 12.2), we can assume that $F(\emptyset) = 0$. For every decision history $h$ and action $a \in A(h)$, define $\lambda(a) = F(ha) - F(h)$. The function $\lambda : A \to \mathbb{N}$ is well defined, since, by assumption, no action is available at more than one information set and, by CM, if $h' \in I(h)$ then $F(h'a) - F(h') = F(ha) - F(h)$. Then, for every history $h = \langle a_1, a_2, ..., a_m \rangle$, $F(h) = \sum_{i=1}^{m} \lambda(a_i)$. In fact,

$$\lambda(a_1) + \lambda(a_2) + ... + \lambda(a_m) = \left(F(a_1) - F(\emptyset)\right) + \left(F(a_1 a_2) - F(a_1)\right) + ... \\ + \left(F(a_1 a_2 ... a_m) - F(a_1 a_2 ... a_{m-1})\right) = F(a_1 a_2 ... a_m) = F(h)$$

(recall that $F(\emptyset) = 0$). Thus, for every $h \in D$ and $a \in A(h)$, $F(ha) = F(h) + \lambda(a)$. Hence, $F(hb) - F(ha) = F(h) + \lambda(b) - \left(F(h) + \lambda(a)\right) = \lambda(b) - \lambda(a)$ and $F(h'b) - F(h'a) = F(h') + \lambda(b) - \left(F(h') + \lambda(a)\right) = \lambda(b) - \lambda(a)$ so that $F(hb) - F(ha) = F(h'b) - F(h'a)$. Thus we have shown that CM implies (♦). Now we show the converse, namely that (♦) implies (CM). Let $\precsim$ be a plausibility order on the set of histories $H$ and let $F: H \to \mathbb{N}$ be an integer-valued representation of $\precsim$ that satisfies (♦). Select arbitrary $h' \in I(h)$ and $a \in A(h)$. Let $b \in A(h)$ be a plausibility-preserving action at $h$ (there must be at least one such action: see Definition 12.1); then, $h \sim hb$ and $h' \sim h'b$. Hence, since $F$ is a representation of $\precsim$, $F(hb) = F(h)$ and $F(h'b) = F(h')$ so that $F(h') - F(h) = F(h'b) - F(hb)$. By (♦), $F(h'b) - F(hb) = F(h'a) - F(ha)$. From the last two equalities it follows that $F(h') - F(h) = F(h'a) - F(ha)$, that is, CM holds.

Since Property CM is equivalent to Property (♦) and the latter implies Property $IND_2$, we have also proved that a sequential equilibrium satisfies $IND_2$.

It only remains to show that a sequential equilibrium satisfies Property $IND_3$. Let $(\sigma, \mu)$ be a sequential equilibrium. Then by Theorem 12.5 there is a choice measurable plausibility order $\precsim$ that rationalizes $(\sigma, \mu)$, relative to which $(\sigma, \mu)$ is uniformly Bayesian. Let $h' \in I(h)$, $a \in A(h)$, $h'a \in I(ha)$ and $\{h, h', ha, h'a\} \subseteq D_\mu^+$. Since $h' \in I(h)$, $\mu(h) > 0$ and $\mu(h') > 0$, by P2 of Definition 12.2, $h$ and $h'$ are most plausible histories in $I(h)$ and thus $h \sim h'$. Similarly $ha$ and $h'a$ are most plausible histories in $I(ha)$ and thus $ha \sim h'a$. Let $E$ be the equivalence class to which $h$ and $h'$ belong, $F$ the equivalence class to which $ha$ and $h'a$ belong and let $\nu_E$ and $\nu_F$ be the probability distributions postulated by Definition 12.3. Let $\nu$ be a uniform common prior (Definition 12.8). Then, by UB2 of Definition 1.8, $\frac{\nu(h)}{\nu(h')} = \frac{\nu(ha)}{\nu(h'a)}$ and by definition of common prior $\nu_E(h) = \frac{\nu(h)}{\sum_{x \in E \cap D_\mu^+} \nu(x)}$, $\nu_E(h') = \frac{\nu(h')}{\sum_{x \in E \cap D_\mu^+} \nu(x)}$, $\nu_F(ha) = \frac{\nu(ha)}{\sum_{y \in F \cap D_\mu^+} \nu(y)}$ and $\nu_F(h'a) = \frac{\nu(h'a)}{\sum_{y \in F \cap D_\mu^+} \nu(y)}$. Thus $\frac{\nu_E(h)}{\nu_E(h')} = \frac{\nu(h)}{\nu(h')}$ and $\frac{\nu_F(ha)}{\nu_F(h'a)} = \frac{\nu(ha)}{\nu(h'a)}$, so that (since $\frac{\nu(h)}{\nu(h')} = \frac{\nu(ha)}{\nu(h'a)}$) $\frac{\nu_E(h)}{\nu_E(h')} = \frac{\nu_F(ha)}{\nu_F(h'a)}$. By Property B3 of Definition 12.3, $\frac{\nu_E(h)}{\nu_E(h')} = \frac{\mu(h)}{\mu(h')}$ and $\frac{\nu_F(ha)}{\nu_F(h'a)} = \frac{\mu(ha)}{\mu(h'a)}$; thus $\frac{\mu(h)}{\mu(h')} = \frac{\mu(ha)}{\mu(h'a)}$, that is, Property $IND_3$ is satisfied. ■

# Appendix 12.E: Exercises

## 12.E.1. Exercises for Section 12.1: Belief revision and *AGM* consistency

The answers to the following exercises are in Appendix S at the end of this chapter.

**Exercise 12.1**. Consider the following game:

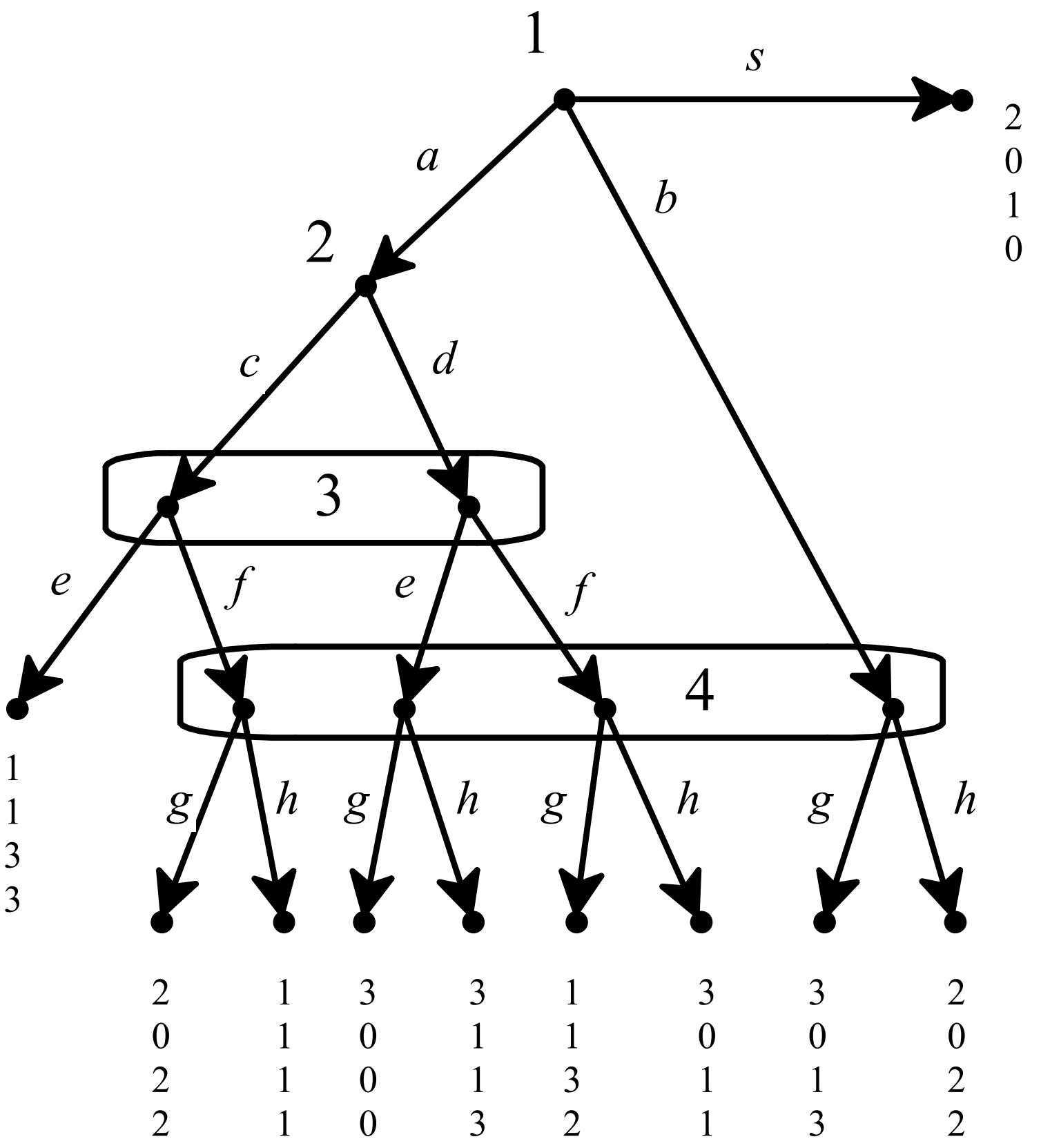

Determine if there is a plausibility order that rationalizes the following assessment (Definition 12.2):

$$\sigma = \left(\begin{array}{ccc|cc|cc|cc} a & b & s & c & d & e & f & g & h \\ \frac{1}{3} & 0 & \frac{2}{3} & 0 & 1 & \frac{1}{2} & \frac{1}{2} & 1 & 0 \end{array}\right), \quad \mu = \left(\begin{array}{cc|cccc} ac & ad & acf & ade & adf & b \\ 0 & 1 & \frac{1}{5} & \frac{1}{5} & \frac{3}{5} & 0 \end{array}\right)$$

**Exercise 12.2**. Consider again the game of Exercise 12.1. Find all the assessments that are rationalized by the following plausibility order:

$$\begin{pmatrix} \emptyset, s & \text{most plausible} \\ a, ac, ad, ace, ade, adeg, b, bg & \\ acf, adf, acfg, adfg & \\ adeh, bh & \\ acfh, adfh & \text{least plausible} \end{pmatrix}$$

## 12.E.2. Exercises for Section 12.2: Bayesian consistency

The answers to the following exercises are in Appendix S at the end of this chapter.

**Exercise 12.3**. Consider the game of Figure 12.4, reproduced below, and the assessment $\sigma = \left(\begin{array}{ccc|cc|cc} a & b & c & d & e & f & g \\ 0 & 0 & 1 & 1 & 0 & \frac{1}{3} & \frac{2}{3} \end{array}\right)$, $\mu = \left(\begin{array}{ccc|ccc} ad & ae & b & a & bf & bg \\ \frac{1}{4} & 0 & \frac{3}{4} & \frac{1}{4} & \frac{1}{4} & \frac{2}{4} \end{array}\right)$ which is rationalized by the plausibility order

$$\begin{pmatrix} \emptyset, c & \text{most plausible} \\ a, ad, adf, adg, b, bf, bg, bfd, bgd & \\ ae, aef, aeg, bfe, bge & \text{least plausible} \end{pmatrix}.$$

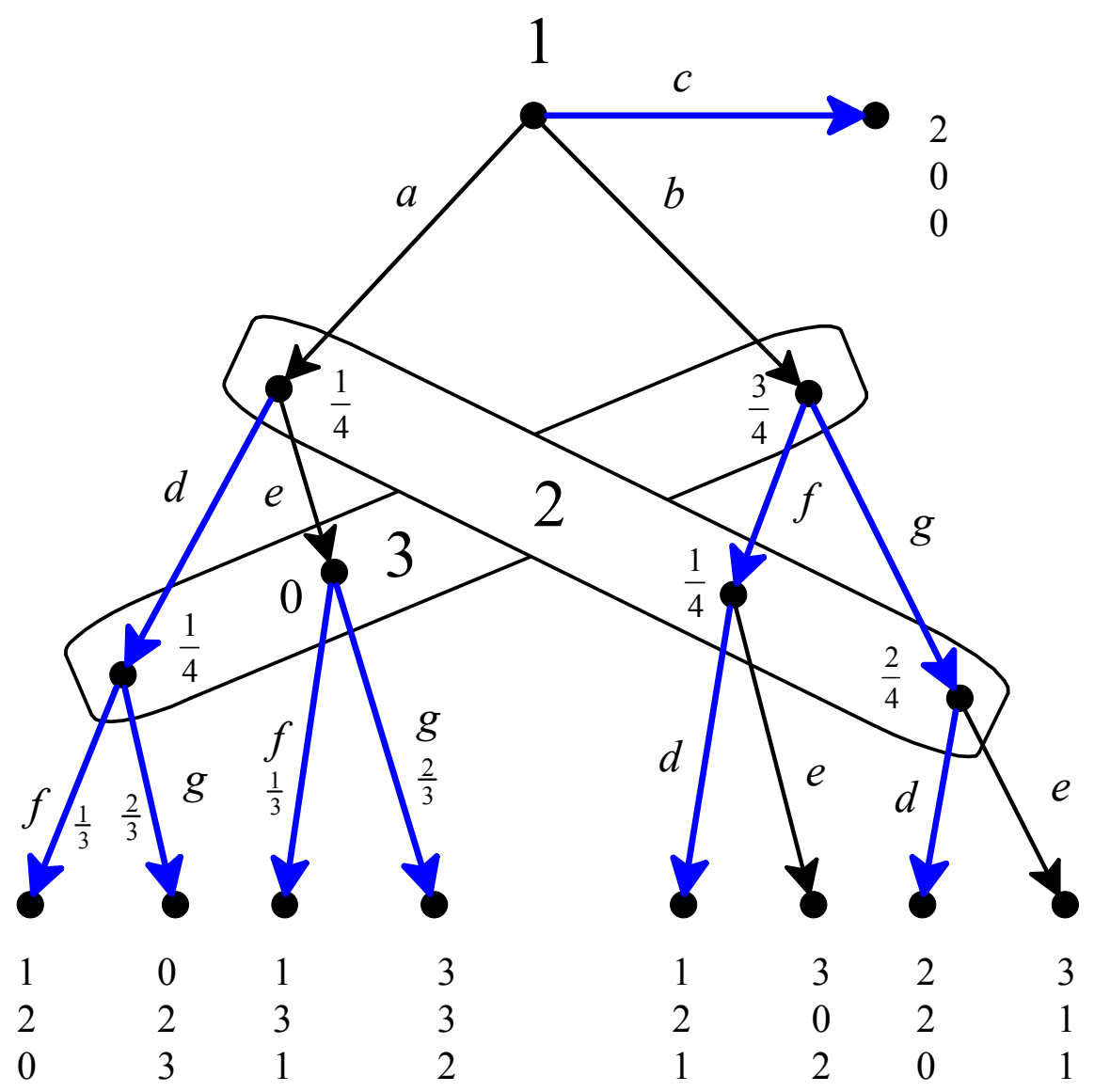

Let $F = \{a, ad, adf, adg, b, bf, bg, bfd, bgd\}$ be the middle equivalence class of the plausibility order. It was shown in Section 12.2 that the probability distribution $\nu_F = \begin{pmatrix} a & ad & b & bf & bg \\ \frac{1}{8} & \frac{1}{8} & \frac{3}{8} & \frac{1}{8} & \frac{2}{8} \end{pmatrix}$ satisfies the properties of Definition 12.3. Show that there is no other probability distribution on $F$ that satisfies those properties.

**Exercise 12.4**. Consider the following extensive-form and the assessment $\sigma = \left( f, A, \begin{pmatrix} L & R & | & \ell & r \\ \frac{1}{2} & \frac{1}{2} & | & \frac{1}{5} & \frac{4}{5} \end{pmatrix}, D \right)$, $\mu = \begin{pmatrix} a & b & c & | & d & e \\ 0 & \frac{1}{3} & \frac{2}{3} & | & \frac{3}{4} & \frac{1}{4} \end{pmatrix}$.

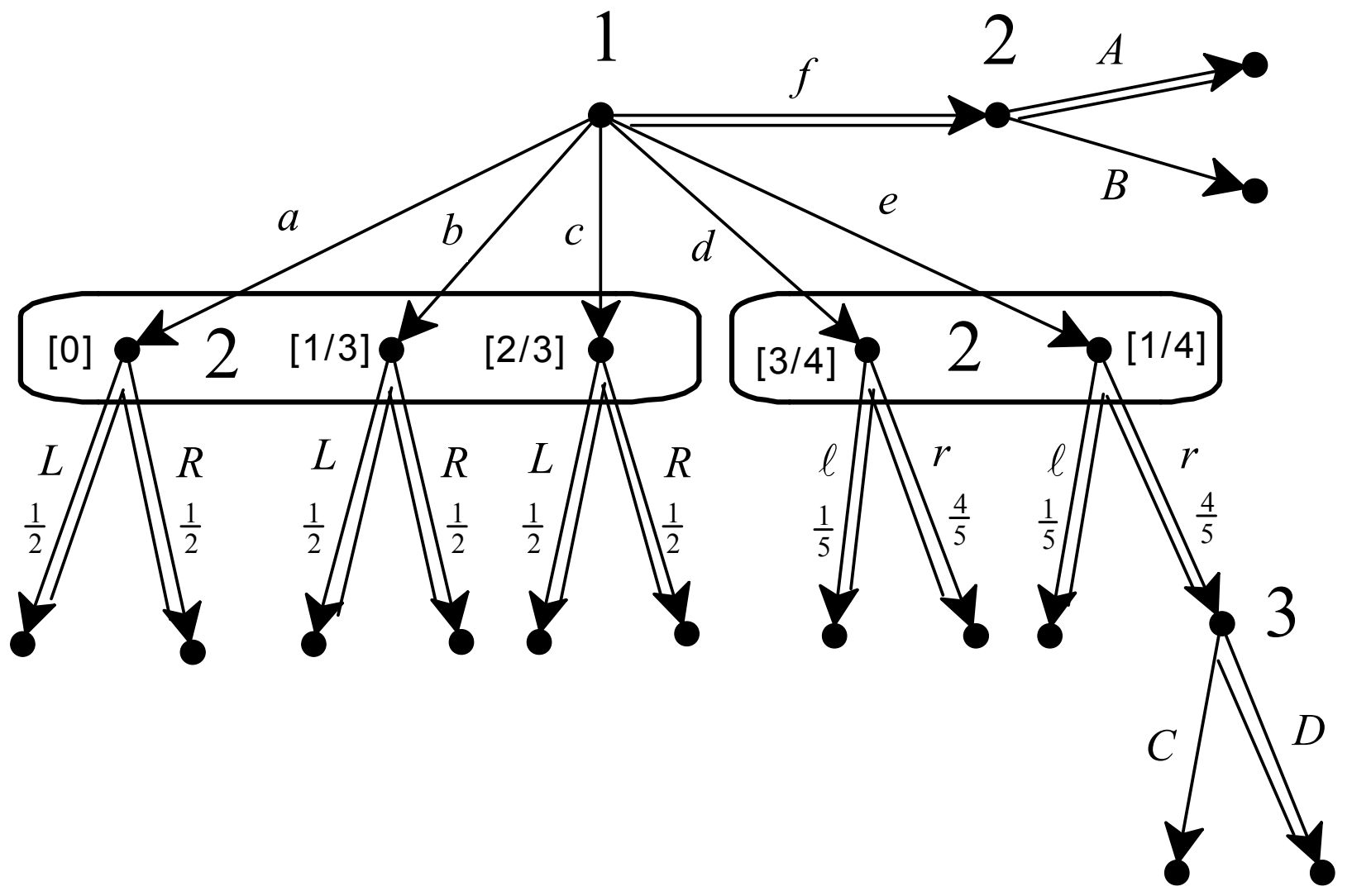

**(a)** Verify that the given assessment is rationalized by the plausibility order

$$\begin{pmatrix} \emptyset, f, fA \\ b, c, d, e, bL, bR, cL, cR, d\ell, dr, e\ell, er, erD \\ a, aL, aR, fB, erC \end{pmatrix}.$$

**(b)** Let $E$ be the top equivalence class of the plausibility order: $E = \{\emptyset, f, fA\}$. Show that there is a unique probability distribution $\nu_E$ on $E$ that satisfies the properties of Definition 12.3.

**(c)** Let $F$ be the middle equivalence class of the plausibility order: $F = \{b, c, d, e, bL, bR, cL, cR, d\ell, dr, e\ell, er, erD\}$. Show that both of the following probability distributions satisfy the properties of Definition 12.3:

$$\nu_F = \begin{pmatrix} b & c & d & e & er \\ \frac{20}{132} & \frac{40}{132} & \frac{45}{132} & \frac{15}{132} & \frac{12}{132} \end{pmatrix}, \quad \hat{\nu}_F = \begin{pmatrix} b & c & d & e & er \\ \frac{5}{87} & \frac{10}{87} & \frac{45}{87} & \frac{15}{87} & \frac{12}{87} \end{pmatrix}$$

## 12.E.3. Exercises for Section 12.3: Perfect Bayesian equilibrium

The answers to the following exercises are in Appendix S at the end of this chapter.

**Exercise 12.5**. Consider the following extensive-form game:

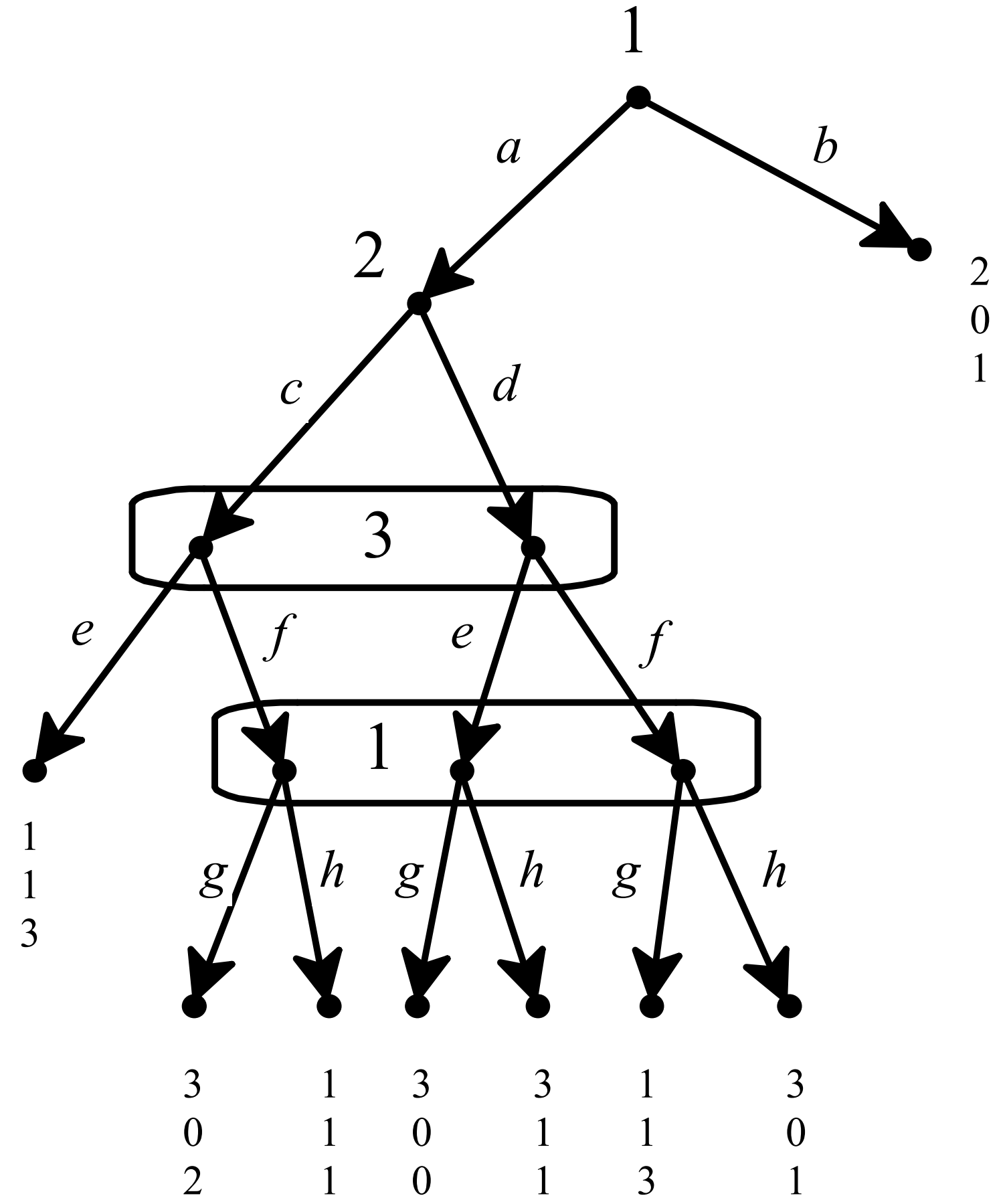

Find two perfect Bayesian equilibria $(\sigma, \mu)$ and $(\sigma', \mu')$ where both $\sigma$ and $\sigma'$ are pure-strategy profiles and $\sigma \neq \sigma'$.

**Exercise 12.6**. Consider the following game:

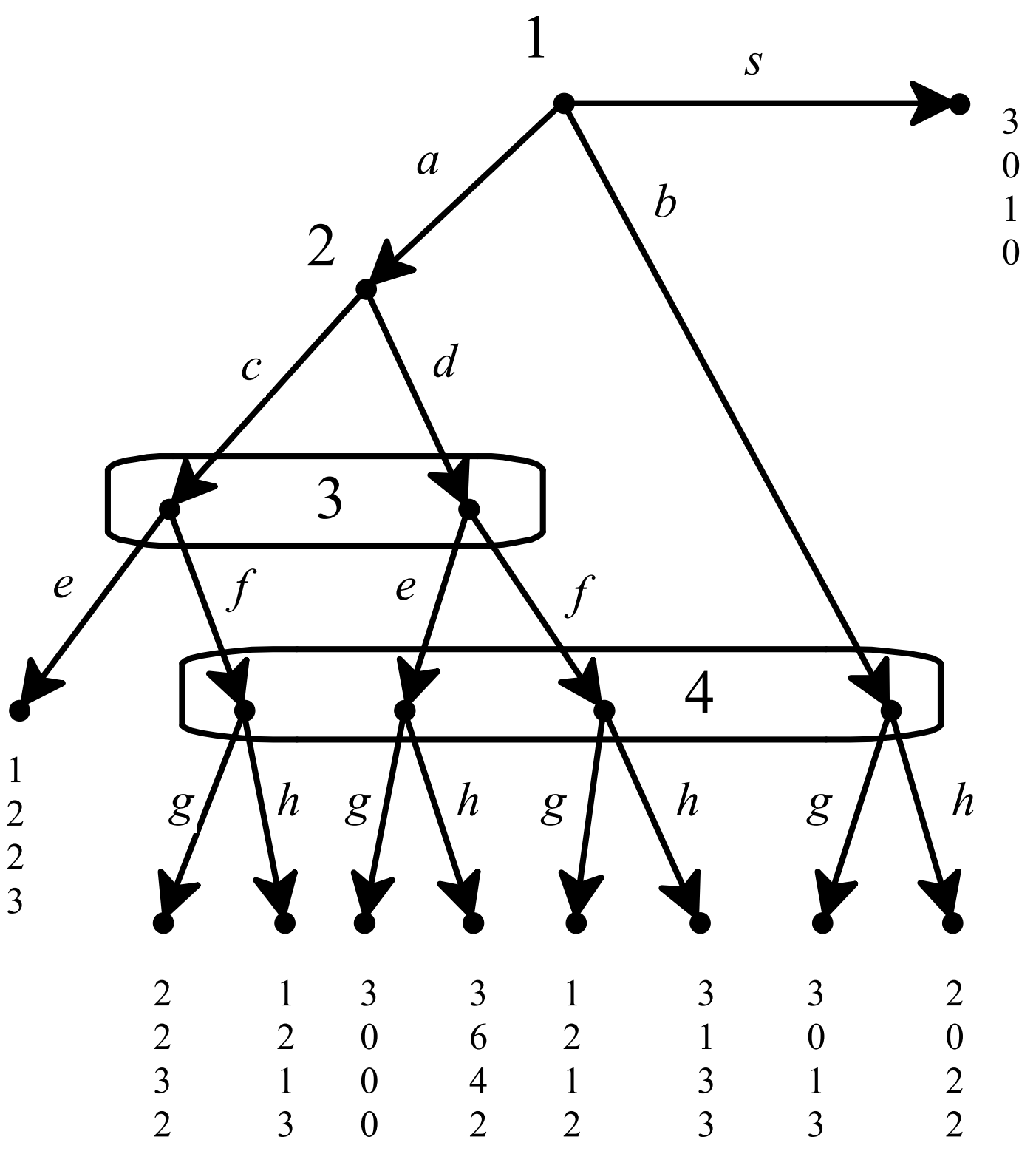

Prove that the following assessment is a perfect Bayesian equilibrium:

$$\sigma = \left( s, \left| \begin{matrix} c & d \\ \frac{1}{2} & \frac{1}{2} \end{matrix} \right| \begin{matrix} e & f \\ \frac{1}{3} & \frac{2}{3} \end{matrix} \left| \begin{matrix} g & h \\ \frac{1}{2} & \frac{1}{2} \end{matrix} \right. \right), \quad \mu = \left( \begin{matrix} ac & ad \\ \frac{1}{2} & \frac{1}{2} \end{matrix} \left| \begin{matrix} acf & ade & adf & b \\ \frac{2}{11} & \frac{1}{11} & \frac{2}{11} & \frac{6}{11} \end{matrix} \right. \right).$$

## 12.E.3. Exercises for Section 12.4: Adding independence

The answers to the following exercises are in Appendix S at the end of this chapter.

**Exercise 12.7**. Draw an extensive form where Player 1 moves first and Player 2 moves second without being informed of Player 1's choice. Player 1 chooses between *a* and *b*, while Player 2's choices are *c*, *d* and *e*. Find an assessment $(\sigma, \mu)$ which is rationalized by a plausibility order that satisfies Property $IND_1$ but fails Property $IND_2$.

**Exercise 12.8**. Find an extensive form and an assessment $(\sigma,\mu)$ which is rationalized by a plausibility order that violates Properties $IND_1$ but satisfies Property $IND_2$ and, furthermore, $\mu$ satisfies Property $IND_3$.

## 12.E.1. Exercises for Section 12.5: Filling the gap between perfect Bayesian equilibrium and sequential equilibrium

The answers to the following exercises are in Appendix S at the end of this chapter.

**Exercise 12.9**. Consider the following (partial) extensive form:

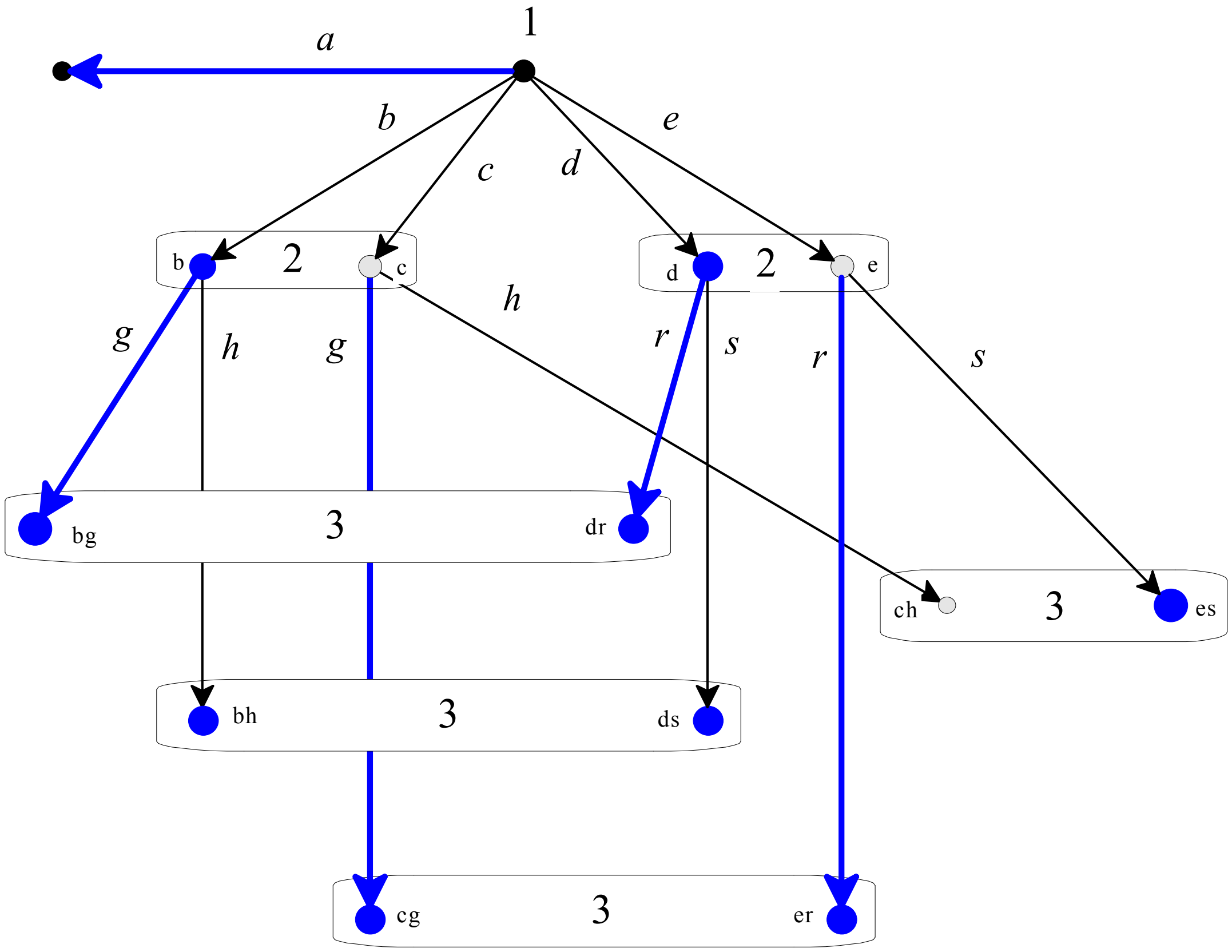

Using Theorem 12.5, prove that there is no sequential equilibrium $(\sigma,\mu)$ where $\sigma=(a,g,r,...)$ (that is, $\sigma$ assigns probability 1 to *a*, *g* and *r*), $\mu(c)=\mu(e)=\mu(ch)=0$ and $\mu(h)>0$ for every other decision history *h*. [Hint: consider all the possible plausibility orders that rationalize $(\sigma,\mu)$.]

**Exercise 12.10**. Consider the game of Figure 12.4, reproduced below:

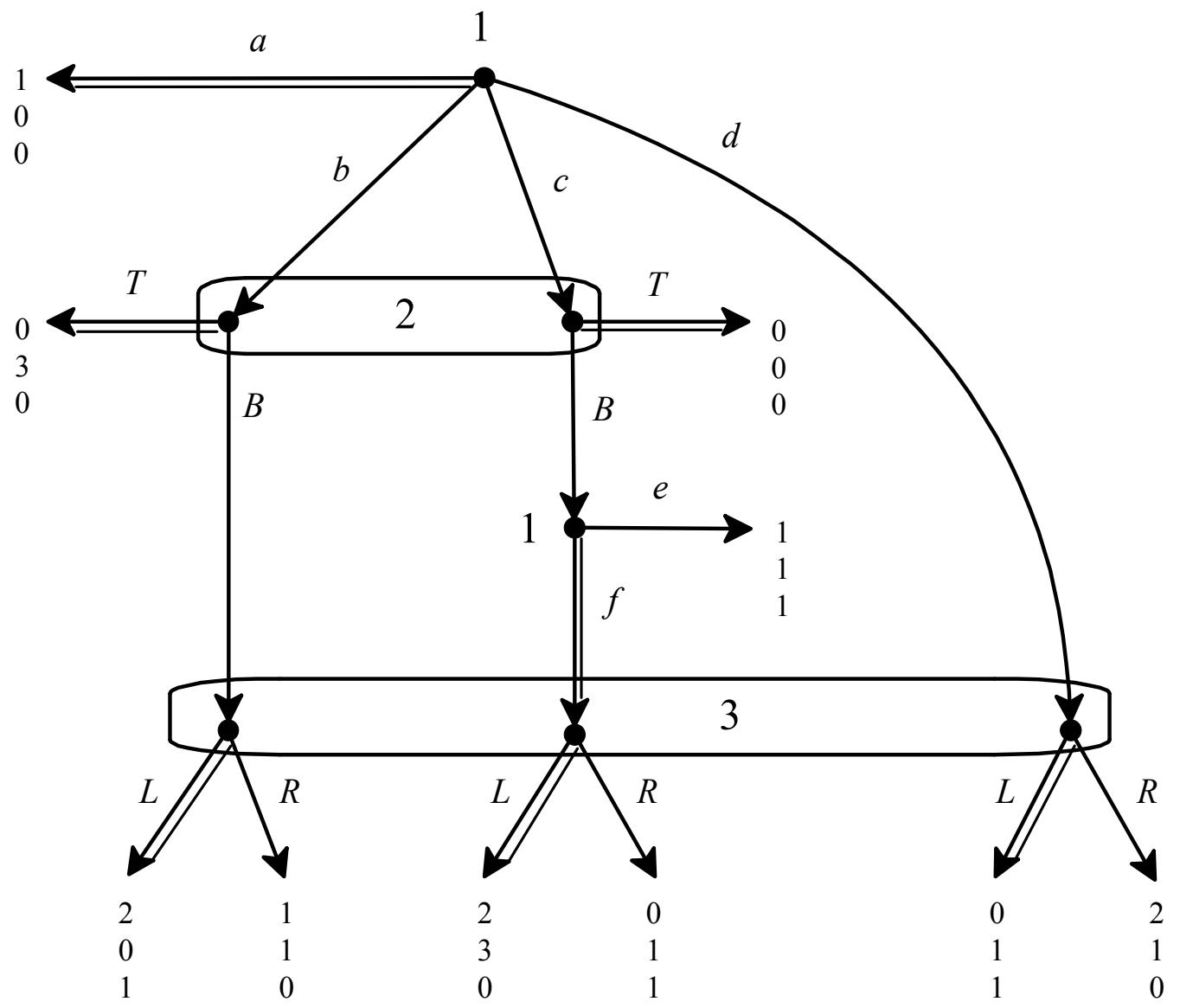

Let $(\sigma, \mu)$ be an assessment with $\sigma = (a, T, f, L)$ (highlighted by double edges; note that $\sigma$ is a subgame-perfect equilibrium), $\mu(b) > 0$ and $\mu(c) > 0$.

(a) Prove that $(\sigma, \mu)$ can be rationalized by a choice-measurable plausibility order only if $\mu$ satisfies the following condition:

$$\mu(bB) > 0 \text{ if and only if } \mu(cBf) > 0$$

(b) Prove that if, besides from being rationalized by a choice-measurable plausibility order $\precsim$, $(\sigma, \mu)$ is also uniformly Bayesian relative to $\precsim$ (Definition 12.8), then $\mu$ satisfies the following condition:

$$\text{if } \mu(bB) > 0 \text{ then } \frac{\mu(cBf)}{\mu(bB)} = \frac{\mu(c)}{\mu(b)}.$$

◊◊◊◊◊◊◊◊◊◊◊◊◊

**Exercise 12.11**. **(Challenging Question).** Consider the following game:

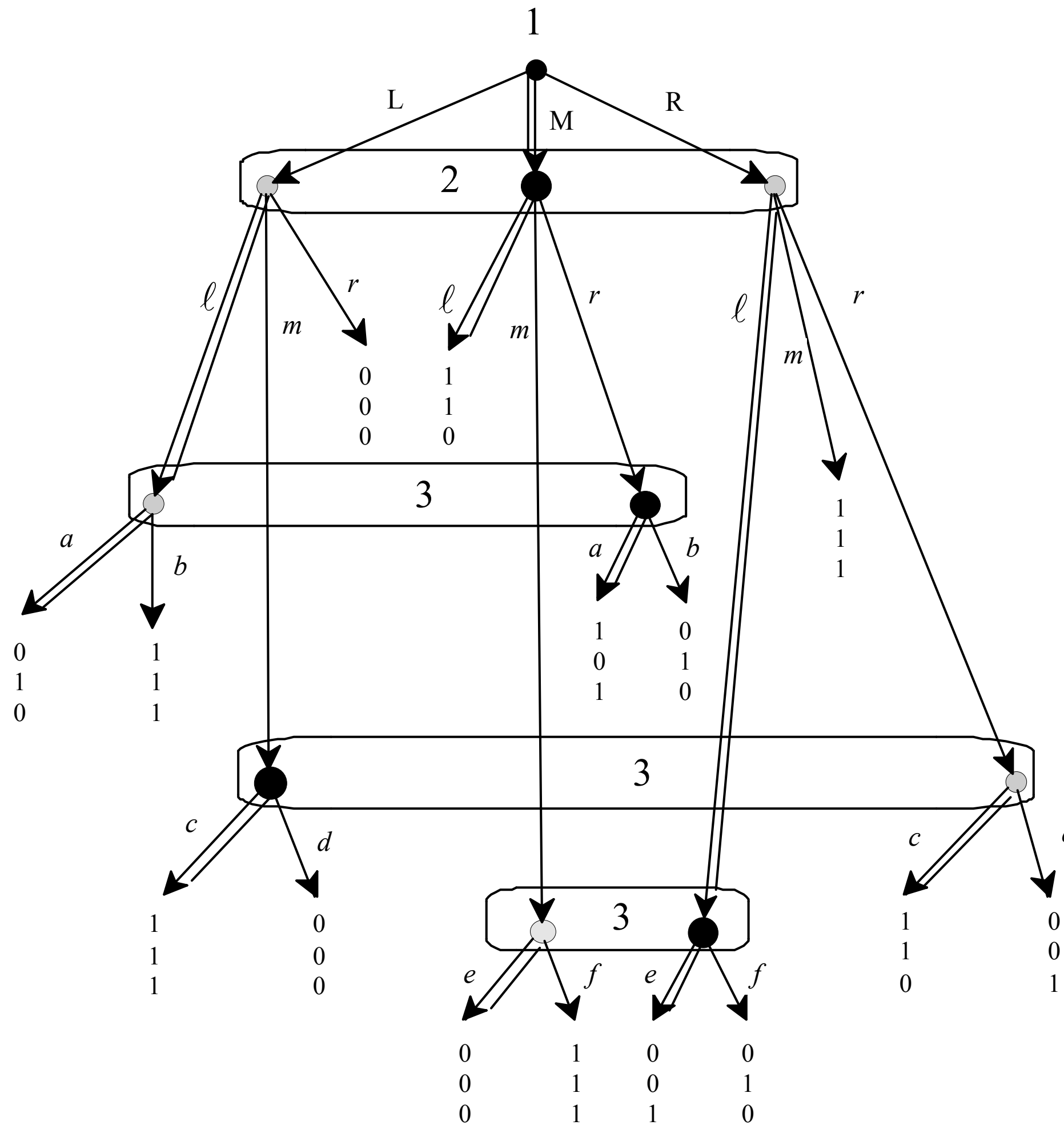

Let $(\sigma,\mu)$ be the following assessment: $\sigma = (M,\ell,a,c,e)$ (highlighted by double edges) and $\mu = \begin{pmatrix} L & M & R & L\ell & Mr & Lm & Rr & Mm & R\ell \\ 0 & 1 & 0 & 0 & 1 & 1 & 0 & 0 & 1 \end{pmatrix}$ (the decision histories $h$ such that $\mu(h)=1$ are denoted by large black dots and the decision histories $h$ such that $\mu(h)=0$ are denoted by small grey dots).

**(a)** Prove that $(\sigma,\mu)$ is an independent perfect Bayesian equilibrium (Definition 12.5).

**(b)** Using Theorem 12.5 prove that $(\sigma,\mu)$ is not a sequential equilibrium. [Hint: prove that $(\sigma,\mu)$ cannot be rationalized by a choice-measurable plausibility order (Definition 12.7).]

# Appendix 12.S: Solutions to exercises

**Exercise 12.1**. From $\sigma$ we get that the plausibility preserving actions are *a, s, d, e, f* and *g*. Thus *ade* is as plausible as *adf* (by transitivity, since each of them is as plausible as *ad*) and each of them is more plausible than both *acf* [since (1) *ad* is more plausible than *ac*, (2) *ade* and *adf* are as plausible as *ad* and (3) *acf* is as plausible as *ac*] and *b* [since (1) *a* is more plausible than *b* and (2) *ade* and *adf* are as plausible as *ad*, which in turn is as plausible as *a*]. Hence the most plausible histories in Player 4's information set are *ade* and *adf*. It follows that it must be that $\mu(ade) > 0$, $\mu(adf) > 0$, $\mu(acf) = 0$ and $\mu(b) = 0$. Hence the given assessment cannot be rationalized by a plausibility order.

**Exercise 12.2**. The plausibility-preserving actions are *s, c, d, e* and *g*. Thus the strategy profile must be of the form $\sigma = \left(\begin{array}{ccc|cc|cc|cc} a & b & s & c & d & e & f & g & h \\ 0 & 0 & 1 & p & 1-p & 1 & 0 & 1 & 0 \end{array}\right)$ with $0 < p < 1$. At Player 3's information set *ac* and *ad* are equally plausible and at Player 4's information set *ade* and *b* are equally plausible and each is more plausible than *acf* and *adf*. Thus the system of beliefs must be of the form $\mu = \left(\begin{array}{cc|cccc} ac & ad & acf & ade & adf & b \\ q & 1-q & 0 & r & 0 & 1-r \end{array}\right)$ with $0 < q < 1$ and $0 < r < 1$.

**Exercise 12.3**. Since $D_\mu^+ \cap F = \{a, ad, b, bf, bg\}$, Property B1 of Definition 12.3 requires that $\nu_F(h) > 0$ if and only if $h \in \{a, ad, b, bf, bg\}$. Let $\nu_F(d) = p \in (0,1)$ and $\nu_F(b) = q \in (0,1)$. The, by Property B2 of Definition 12.3, $\nu_F(ad) = \nu_F(a) \times \sigma(d) = p \times 1 = p$, $\nu_F(bf) = \nu_F(b) \times \sigma(f) = q \times \frac{1}{3}$ and $\nu_F(bg) =$. $\nu_F(b) \times \sigma(g) = q \times \frac{2}{3}$ Thus $\nu_F = \begin{pmatrix} a & ad & b & bf & bg \\ p & p & q & \frac{1}{3}q & \frac{2}{3}q \end{pmatrix}$ and the sum of these probabilities must be 1: $2p + 2q = 1$, that is $\boxed{p + q = \frac{1}{2}}$. let $I_2 = \{a, bf, bg\}$ be the information set of Player 2 and $I_3 = \{b, ad, ae\}$ the information set of Player 3. Then $\nu_F(I_2) = \nu_F(a) + \nu_F(bf) + \nu_F(bg) = p + \frac{1}{3}q + \frac{2}{3}q = p + q = \frac{1}{2}$ and $\nu_F(I_3) = \nu_F(b) + \nu_F(ad) + \nu_F(ae) = q + p + 0 = \frac{1}{2}$. Thus $\frac{\nu_F(a)}{\nu_F(I_2)} = \frac{p}{\frac{1}{2}} = 2p$ and, by Property B3, we need this to be equal to $\mu(a) = \frac{1}{4}$; solving $2p = \frac{1}{4}$ we get $p = \frac{1}{8}$. Similarly, $\frac{\nu_F(b)}{\nu_F(I_3)} = \frac{q}{\frac{1}{2}} = 2q$ and, by Property B3, we need this to be equal to $\mu(b)=\frac{3}{4}$; solving $2q=\frac{3}{4}$ we get $q=\frac{3}{8}$ (alternatively, we could have derived $q=\frac{3}{8}$ from $p=\frac{1}{8}$ and $p+q=\frac{1}{2}$).

**Exercise 12.4**. **(a)** The plausibility-preserving actions are $f, A, L, R, \ell, r$ and $D$ and these are precisely the actions that are assigned positive probability by $\sigma$. Furthermore, the most plausible actions in information set $\{a,b,c\}$ are $b$ and $c$ and the two actions in information set $\{d,e\}$ are equally plausible. This is consistent with the fact that $D_\mu^+ = \{\emptyset, f, b, c, d, e, er\}$. Thus the properties of Definition 12.2 are satisfied.

**(b)** By Property B1 of Definition 12.3 the support of $\nu_E$ must be $D_\mu^+ \cap E = \{\emptyset, f\}$ and by Property B2 it must be that $\nu_E(f) = \nu_E(\emptyset)\times\sigma(f) = \nu_E(\emptyset)\times 1 = \nu_E(\emptyset)$. Thus the only solution is $\nu_E(\emptyset) = \nu_E(f) = \frac{1}{2}$.

**(c)** Since $D_\mu^+ \cap F = \{b,c,d,e,er\}$, by Property B1 of Definition 12.3 the support of the probability distribution must coincide with the set $\{b,c,d,e,er\}$, which is indeed true for both $\nu_F$ and $\hat{\nu}_F$. By Property B2, the probability of *er* must be equal to the probability of *e* times $\sigma(r)=\frac{4}{5}$ and this is indeed true for both $\nu_F$ and $\hat{\nu}_F$. By Property B3, the conditional probability $\dfrac{\nu_F(b)}{\nu_F(a)+\nu_F(b)+\nu_F(c)} = \dfrac{\frac{20}{132}}{0+\frac{20}{132}+\frac{40}{132}} = \frac{1}{3}$ must be equal to $\mu(b)$ and this is indeed true; similarly, the conditional probability $\dfrac{\nu_F(d)}{\nu_F(d)+\nu_F(e)} = \dfrac{\frac{45}{132}}{\frac{45}{132}+\frac{15}{132}} = \frac{3}{4}$ must be equal to $\mu(d)$ and this is also true. Similar computations show that $\dfrac{\hat{\nu}_F(b)}{\hat{\nu}_F(a)+\hat{\nu}_F(b)+\hat{\nu}_F(c)} = \mu(b)$ and $\dfrac{\hat{\nu}_F(d)}{\hat{\nu}_F(d)+\hat{\nu}_F(e)} = \mu(d)$.

**Exercise 12.5**. The game under consideration is:

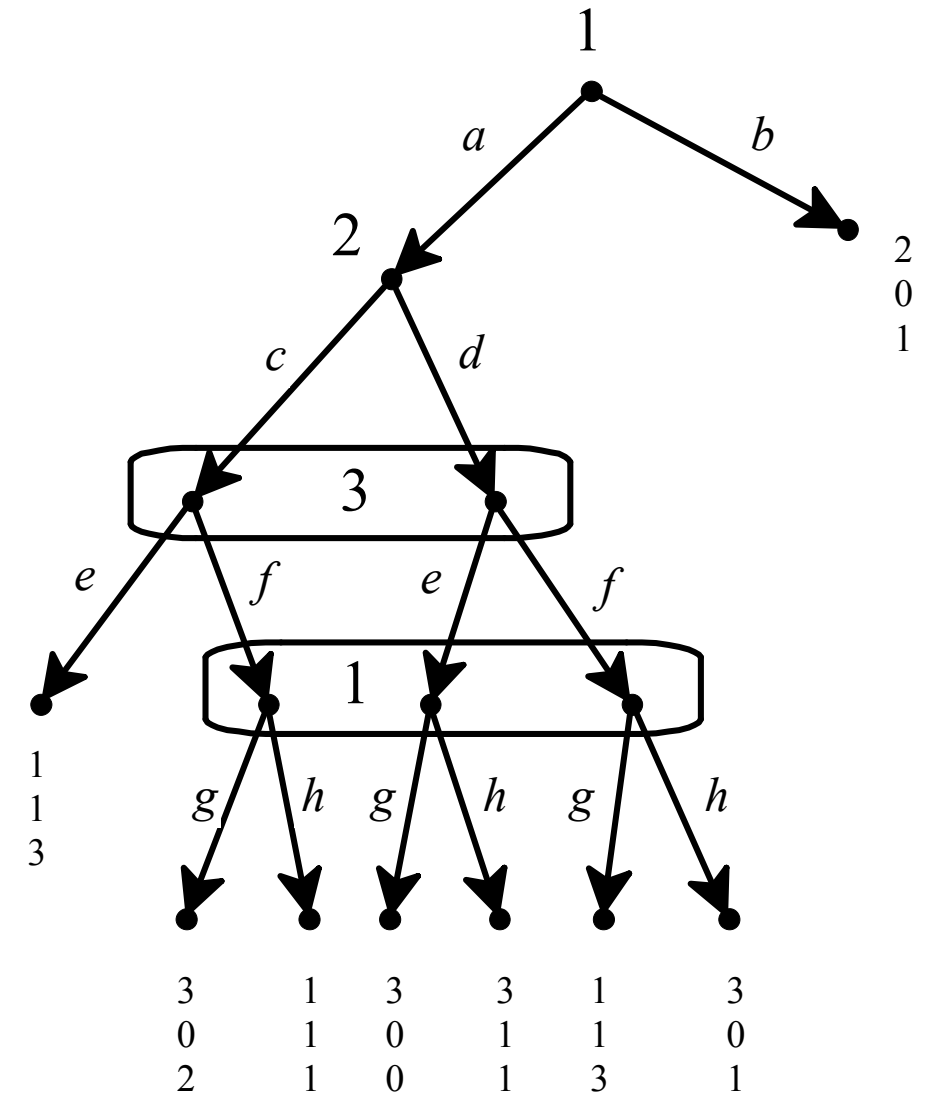

One perfect Bayesian equilibrium is

$$\sigma = (b, c, e, g), \qquad \mu = \begin{pmatrix} ac & ad & acf & ade & adf \\ 1 & 0 & \frac{1}{2} & \frac{1}{2} & 0 \end{pmatrix}.$$

Let us first verify sequential rationality. For Player 1, at the root, $b$ gives a payoff of 2, while $a$ gives a payoff of 1; thus $b$ is sequentially rational. For Player 2, $c$ gives a payoff of 1 and $d$ a payoff of 0; thus $c$ is sequentially rational. For Player 3, $e$ gives a payoff of 3, while $f$ gives a payoff of 2; thus $e$ is sequentially rational. For Player 1 at his bottom information set, $g$ gives a payoff of $\left(\frac{1}{2}\right)3 + \left(\frac{1}{2}\right)3 = 3$ and $h$ gives a payoff of $\left(\frac{1}{2}\right)1 + \left(\frac{1}{2}\right)3 = 2$; thus $g$ is sequentially rational.

The following plausibility order rationalizes the above assessment:

$$\begin{pmatrix} \emptyset, b & \text{most plausible} \\ a, ac, ace & \\ acf, acfg, ad, ade, adeg & \\ adf, adfg, acfh, adeh & \\ adfh & \text{least plausible} \end{pmatrix}$$

Name the equivalence classes of this order $E_1, E_2, \ldots, E_5$ (with $E_1$ being the top one and $E_5$ the bottom one). Then the following probability distributions on the sets $E_i \cap D_\mu^+$ ($i$ = 1,2,3: note $E_i \cap D_\mu^+ \neq \varnothing$ if and only if $i \in \{1, 2, 3\}$, since $D_\mu^+ = \{\emptyset, a, ac, acf, ade\}$) satisfy the properties of Definition 12.3: $\nu_{E_1}(\emptyset) = 1$, $\nu_{E_2} = \begin{pmatrix} a & ac \\ \frac{1}{2} & \frac{1}{2} \end{pmatrix}$, $\nu_{E_3} = \begin{pmatrix} acf & ad & ade \\ \frac{1}{3} & \frac{1}{3} & \frac{1}{3} \end{pmatrix}$.

Another perfect Bayesian equilibrium is:

$$\sigma = (a, d, e, h), \qquad \mu = \left(\begin{array}{cc|ccc} ac & ad & acf & ade & adf \\ 0 & 1 & 0 & 1 & 0 \end{array}\right).$$

Sequential rationality is easily verified. The assessment is rationalized by:

$$\begin{pmatrix} \emptyset, a, ad, ade, adeh & \text{most plausible} \\ b, ac, ace, adf, adfh & \\ acf, acfh, adeg & \\ acfg, adfg & \text{least plausible} \end{pmatrix}$$

Only the top equivalence class $E_1 = \{\emptyset, a, ad, ade, adeh\}$ has a non-empty intersection with $D_\mu^+ = \{\emptyset, a, ad, ade\}$. The following probability distribution satisfies the properties of Definition 12.3: $\nu_{E_1} = \begin{pmatrix} \emptyset & a & ad & ade \\ \frac{1}{4} & \frac{1}{4} & \frac{1}{4} & \frac{1}{4} \end{pmatrix}$.

**Exercise 12.6.** The game under consideration is:

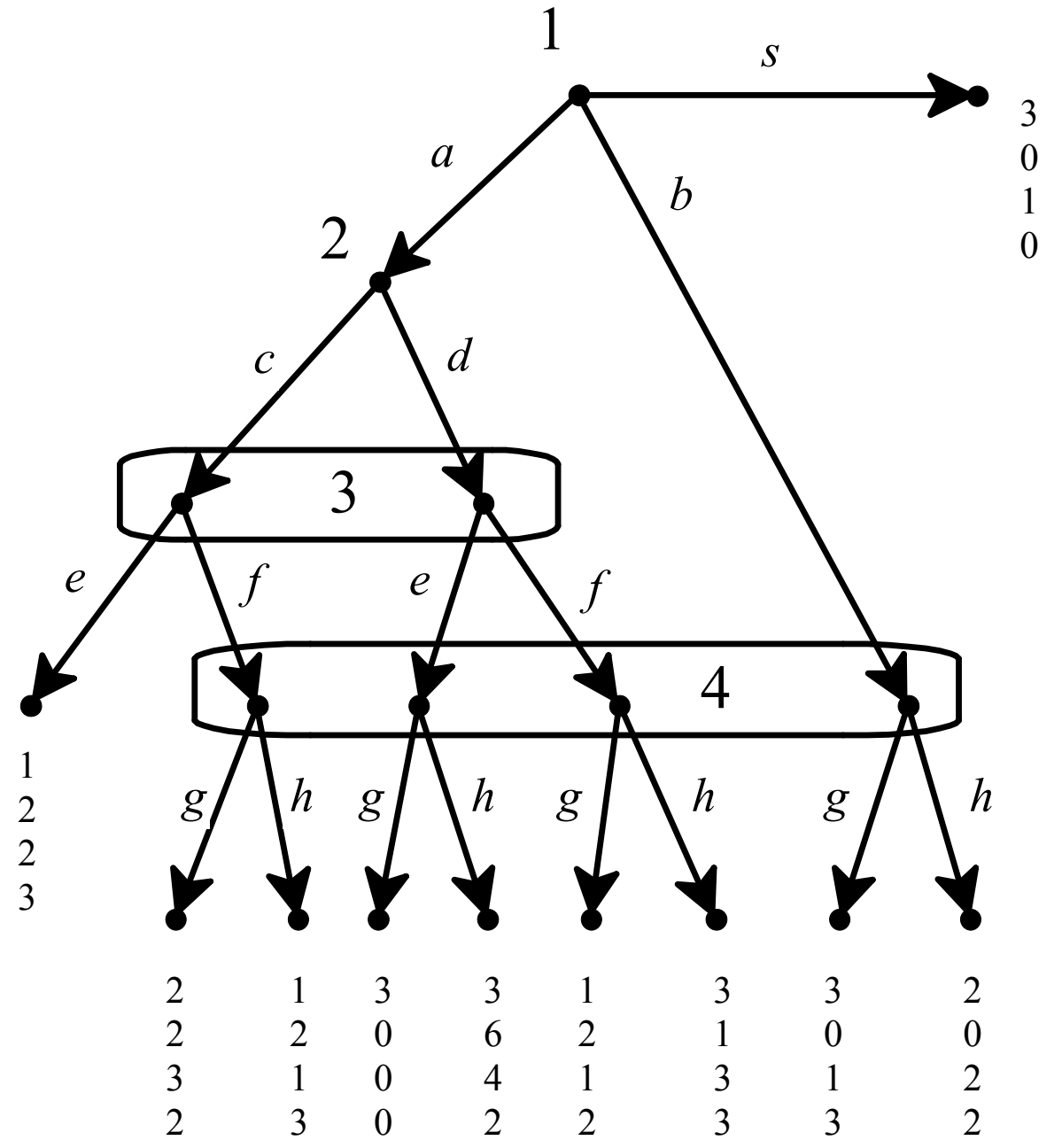

We need to show that the following is a perfect Bayesian equilibrium:

$$\sigma = \left( s, \left| \begin{matrix} c & d \\ \frac{1}{2} & \frac{1}{2} \end{matrix} \right| \begin{matrix} e & f \\ \frac{1}{3} & \frac{2}{3} \end{matrix} \left| \begin{matrix} g & h \\ \frac{1}{2} & \frac{1}{2} \end{matrix} \right. \right), \quad \mu = \left( \begin{matrix} ac & ad \\ \frac{1}{2} & \frac{1}{2} \end{matrix} \left| \begin{matrix} acf & ade & adf & b \\ \frac{2}{11} & \frac{1}{11} & \frac{2}{11} & \frac{6}{11} \end{matrix} \right. \right).$$

First we verify sequential rationality. For Player 1 the possible payoffs are:

from $s$: 3

from $b$: $\frac{1}{2}(3)+\frac{1}{2}(2)=2.5$ .

from $a$: $\frac{1}{2}\left[\frac{1}{3}(1)+\frac{2}{3}\left(\frac{1}{2}(2)+\frac{1}{2}(1)\right)\right]+\frac{1}{2}\left[\frac{1}{3}(3)+\frac{2}{3}\left(\frac{1}{2}(1)+\frac{1}{2}(3)\right)\right]=\frac{11}{6}$

Thus $s$ is sequentially rational. For Player 2 the possible payoffs are:

from $c$: $\frac{1}{3}(2)+\frac{2}{3}(2)=2$

from $d$: $\frac{1}{3}\left[\frac{1}{2}(0)+\frac{1}{2}(6)\right]+\frac{2}{3}\left[\frac{1}{2}(2)+\frac{1}{2}(1)\right]=2$ . Thus both $c$ and $d$ are sequentially rational and so is any mixture of the two, in particular the mixture $\begin{pmatrix} c & d \\ \frac{1}{2} & \frac{1}{2} \end{pmatrix}$.

For Player 3 the possible payoffs are:

from $e$: $\frac{1}{2}(2)+\frac{1}{2}\left[\frac{1}{2}(0)+\frac{1}{2}(4)\right]=2$

from $f$: $\frac{1}{2}\left[\frac{1}{2}(3)+\frac{1}{2}(1)\right]+\frac{1}{2}\left[\frac{1}{2}(1)+\frac{1}{2}(3)\right]=2$ . Thus both $e$ and $f$ are sequentially rational and so is any mixture of the two, in particular the mixture $\begin{pmatrix} e & f \\ \frac{1}{3} & \frac{2}{3} \end{pmatrix}$.

For Player 4 the possible payoffs are as follows: From *g*: $\frac{2}{11}(2)+\frac{1}{11}(0)+\frac{2}{11}(2)+\frac{6}{11}(3)=\frac{26}{11}$ and from *h*: $\frac{2}{11}(3)+\frac{1}{11}(2)+\frac{2}{11}(3)+\frac{6}{11}(2)=\frac{26}{11}$. Thus both *g* and *h* are rational and so is any mixture of the two, in particular the mixture $\begin{pmatrix} g & h \\ \frac{1}{2} & \frac{1}{2} \end{pmatrix}$.

For AGM consistency, note that all of the actions, except *a* and *b*, are plausibility preserving and, furthermore, all decision histories are assigned positive probability by $\mu$. Thus there is only one plausibility order that rationalizes the given assessment, namely the one that has only two equivalence classes: the top one being $\{\emptyset, s\}$ and the other one consisting of all the remaining histories. For Bayesian consistency, the probability distribution for the top equivalence class is the trivial one that assigns probability 1 to the null history $\emptyset$. Let $E = H \setminus \{\emptyset, s\}$ be the other equivalence class and note that $E \cap D_m^+ = \{a, ac, ad, acf, ade, adf, b\}$. In order for a probability distribution $\nu_E$ to satisfy the properties of Definition 12.3, the support of $\nu_E$ must be the set $\{a, ac, ad, acf, ade, adf, b\}$ (Property B1). Let $\nu_E(a) = p$; then by Property B2 it must be that $\nu_E(ac) = \nu_E(ad) = \frac{p}{2}$ (since $\sigma(c) = \sigma(d) = \frac{1}{2}$), $\nu_E(acf) = \nu_E(adf) = \frac{p}{3}$ (since $\sigma(c) \times \sigma(f) = \sigma(d) \times \sigma(f) = \frac{1}{2} \times \frac{2}{3} = \frac{1}{3}$) and $\nu_E(ade) = \frac{p}{6}$ (since $\sigma(d) \times \sigma(e) = \frac{1}{2} \times \frac{1}{3} = \frac{1}{6}$). Thus $\nu_E$ must be of the form $\begin{pmatrix} a & ac & ad & acf & ade & adf & b \\ p & \frac{p}{2} & \frac{p}{2} & \frac{p}{3} & \frac{p}{6} & \frac{p}{3} & q \end{pmatrix}$ with the sum equal to 1, that is, $\frac{17}{6}p + q = 1$. Furthermore, by Property B3, it must be that

$$\frac{\nu_E(b)}{\nu_E(acf) + \nu_E(ade) + \nu_E(adf) + \nu_E(b)} = \frac{q}{\frac{p}{3} + \frac{p}{6} + \frac{p}{3} + q} = \underbrace{\frac{6}{11}}_{=\mu(b)}.$$

The solution to the pair of equations $\frac{17}{6}p + q = 1$ and $\frac{q}{\frac{p}{3} + \frac{p}{6} + \frac{p}{3} + q} = \frac{6}{11}$ is $p = q = \frac{6}{23}$, yielding

$$\nu_E = \begin{pmatrix} a & ac & ad & acf & ade & adf & b \\ \frac{6}{23} & \frac{3}{23} & \frac{3}{23} & \frac{2}{23} & \frac{1}{23} & \frac{2}{23} & \frac{6}{23} \end{pmatrix}.$$

Thus we have shown that the given assessment is AGM- and Bayes-consistent and sequentially rational, hence a perfect Bayesian equilibrium .

**Exercise 12.7.** The extensive form is as follows:

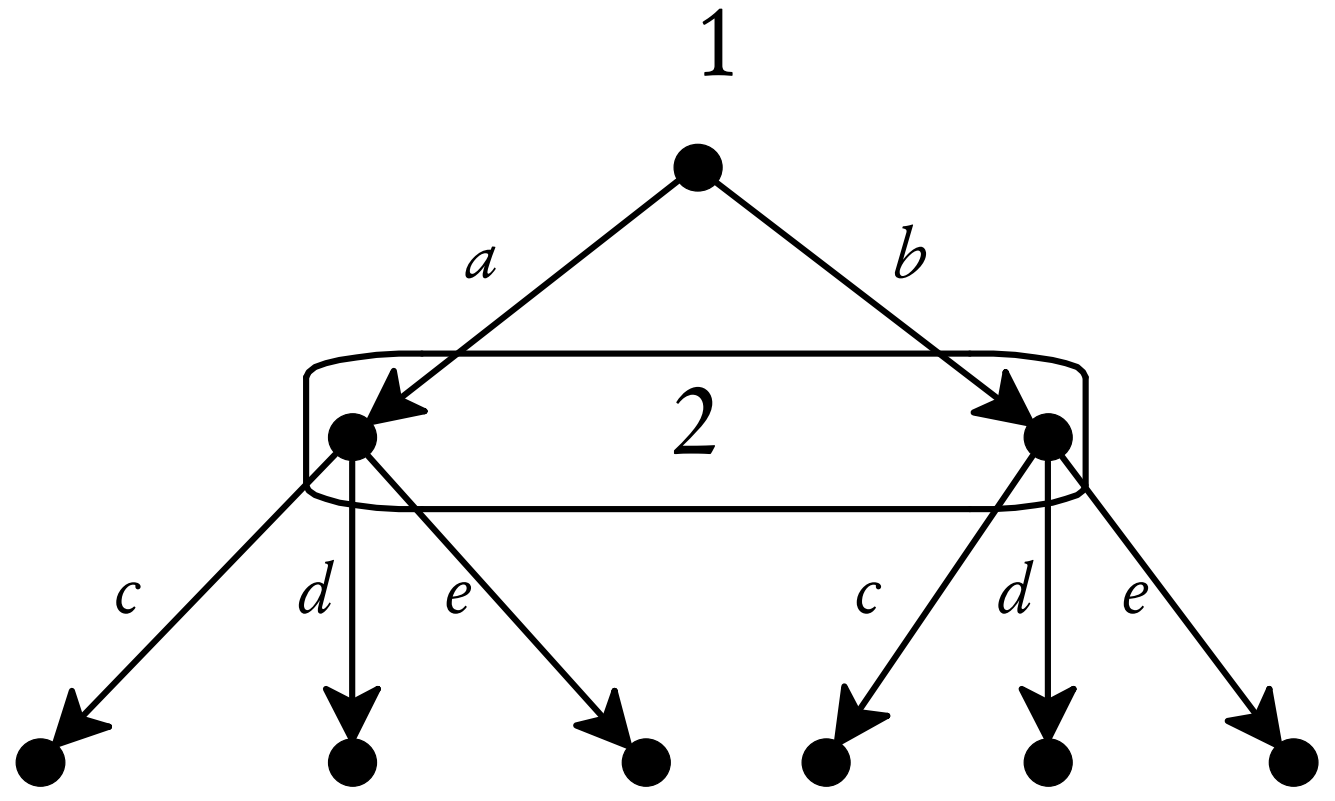

Let $\sigma = \left(\begin{array}{cc|ccc} a & b & c & d & e \\ 1 & 0 & 1 & 0 & 0 \end{array}\right)$, $\mu = \begin{pmatrix} a & b \\ 1 & 0 \end{pmatrix}$. This assessment is rationalized by the following plausibility order:

$$\begin{pmatrix} \emptyset, a, ac & \text{most plausible} \\ b, bc & \\ ad & \\ ae & \\ be & \\ bd & \text{least plausible} \end{pmatrix}$$

which satisfies Property $IND_1$ (since $a \prec b$ and $ac \prec bc$, $ad \prec bd$, $ae \prec be$) but fails Property $IND_2$ since $b \in I(a)$ and $ad \prec ae$ (implying that, conditional on *a*, *d* is more plausible than *e*) but $be \prec bd$ (implying that, conditional on *b*, *e* is more plausible than *d*).

**Exercise 12.8**. The game of Figure 12.7, reproduced below, provides such an example:

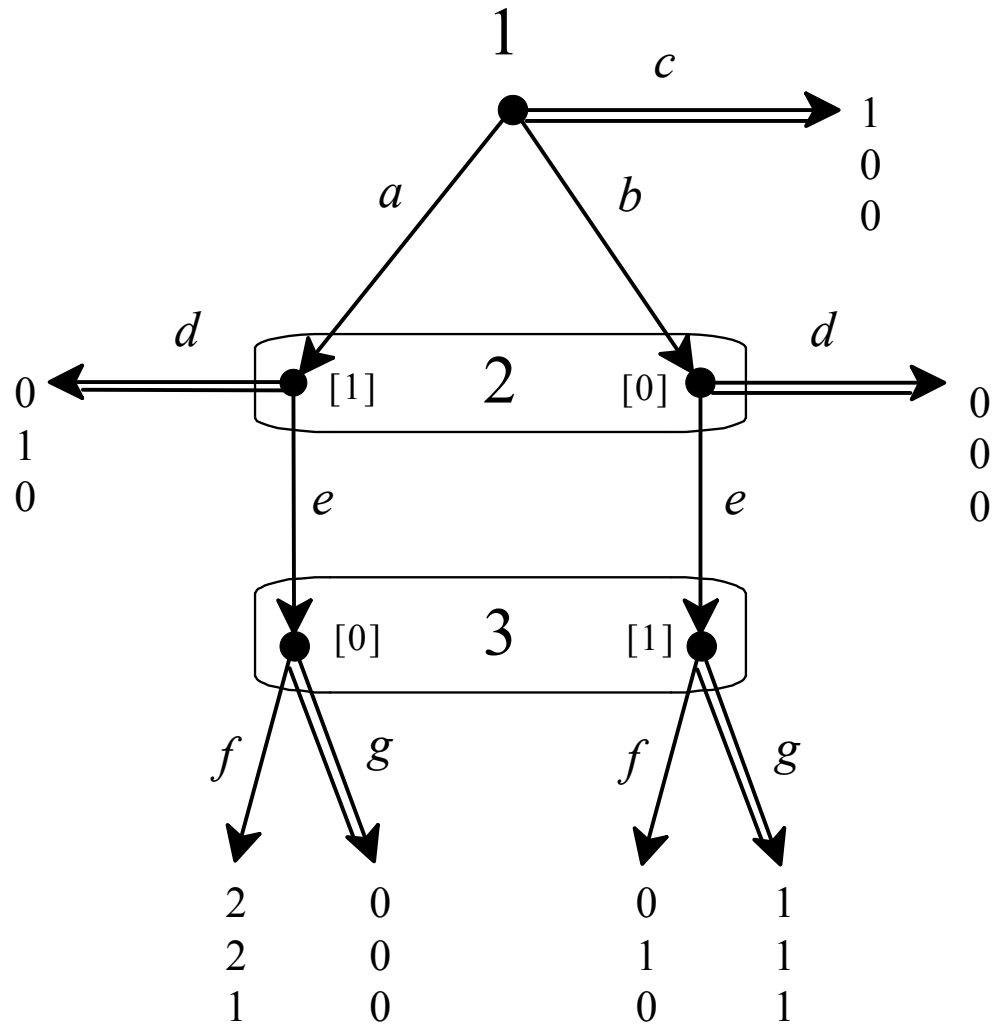

Consider the assessment $\sigma = (c, d, g)$ (highlighted by double edges), together with the system of beliefs $\mu = \begin{pmatrix} a & b & ae & be \\ 1 & 0 & 0 & 1 \end{pmatrix}$. This assessment is rationalized by the plausibility order $\begin{pmatrix} \emptyset, c & \text{most plausible} \\ a, ad & \\ b, bd & \\ be, beg & \\ ae, aeg & \\ bef & \\ aef & \text{least plausible} \end{pmatrix}$ which violates Property $IND_1$, since $a \prec b$ and $be \prec ae$. On the other hand, the above plausibility order satisfies Property $IND_2$ since (1) $ad \prec ae$ and $bd \prec be$ and (2) $aeg \prec aef$ and $beg \prec bef$. Furthermore, $\mu$ satisfies Property $IND_3$: trivially, since $a, be \in D_\mu^+$ but $ae, b \notin D_\mu^+$.

**Exercise 12.9.** The game under consideration is the following:

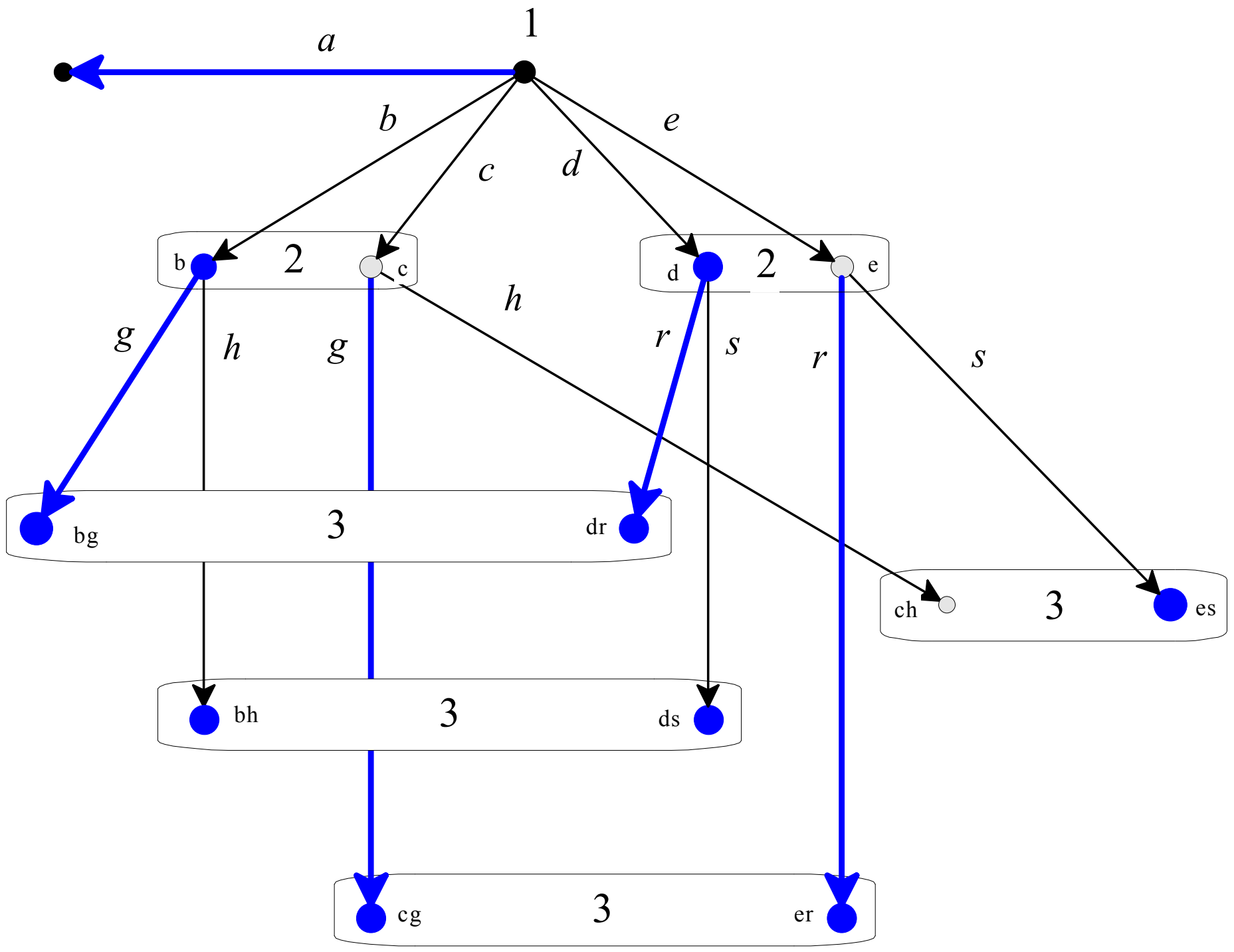

Consider an arbitrary assessment of the form $(\sigma, \mu)$ where $\sigma = (a, g, r, ...)$, $\mu(c) = \mu(e) = \mu(ch) = 0$ and $\mu(h) > 0$ for every other decision history *h*. There are plausibility orders that rationalize $(\sigma, \mu)$, for example the following:

$$\begin{pmatrix} \emptyset, a & \text{most plausible} \\ b, bg, d, dr & \\ bh, ds & \\ c, cg, e, er & \\ es & \\ ch & \text{least plausible} \end{pmatrix} \quad \text{or} \quad \begin{pmatrix} \emptyset, a & \text{most plausible} \\ b, bg, d, dr & \\ c, cg, e, er & \\ bh, ds & \\ es & \\ ch & \text{least plausible} \end{pmatrix}.$$

First we show that any plausibility order $\precsim$ that rationalizes $(\sigma, \mu)$ must satisfy the following properties:

- $c \sim e$ [because, by P2 of Definition 12.2, $cg \sim er$ (since $cg \in I(er)$, $\mu(cg) > 0$ and $\mu(er) > 0$) and, by P1 of Definition 12.2, $c \sim cg$ and $e \sim er$ (since both *g* and *r* are plausibility preserving) and thus, by transitivity, $c \sim e$],

- $es \prec ch$ [this follows from P2 of Definition 12.2, since *es* and *ch* belong to the same information set and $\mu(es) > 0$, while $\mu(ch) = 0$],

- $b \sim d$ [because, by P2 of Definition 12.2, $bg \sim dr$ (since $bg \in I(dr)$, $\mu(bg) > 0$ and $\mu(dr) > 0$) and, by P1 of Definition 12.2, $b \sim bg$ and $d \sim dr$ (since both *g* and *r* are plausibility preserving) and thus, by transitivity, $b \sim d$],

- $bh \sim ds$ [by P2 of Definition 12.2, because $bh \in I(ds)$, $\mu(bh) > 0$ and $\mu(ds) > 0$.

Next we show that no plausibility order that rationalizes $(\sigma, \mu)$ is choice measurable. Select an arbitrary plausibility order that rationalizes $(\sigma, \mu)$ and let $F$ be an integer valued representation of it. Then the following must be true:

(1) $F(e) - F(es) > F(c) - F(ch)$ (because $c \sim e$, implying that $F(c) = F(e)$, and $es \prec ch$, implying that $F(es) < F(ch)$),

(2) $F(b) - F(bh) = F(d) - F(ds)$ (because $bh \sim ds$, implying that $F(bh) = F(ds)$, and $b \sim d$, implying that $F(b) = F(d)$),

Thus if, as required by choice measurability, $F(c) - F(ch) = F(b) - F(bh)$ then, by (1) and (2), $F(e) - F(es) > F(d) - F(ds)$, which violates choice measurability. It follows from Theorem 12.5 that, since any plausibility ordering that rationalizes $(\sigma, \mu)$ is not choice measurable, $(\sigma, \mu)$ cannot be a sequential equilibrium.

**Exercise 12.10.** The game under consideration is the following:

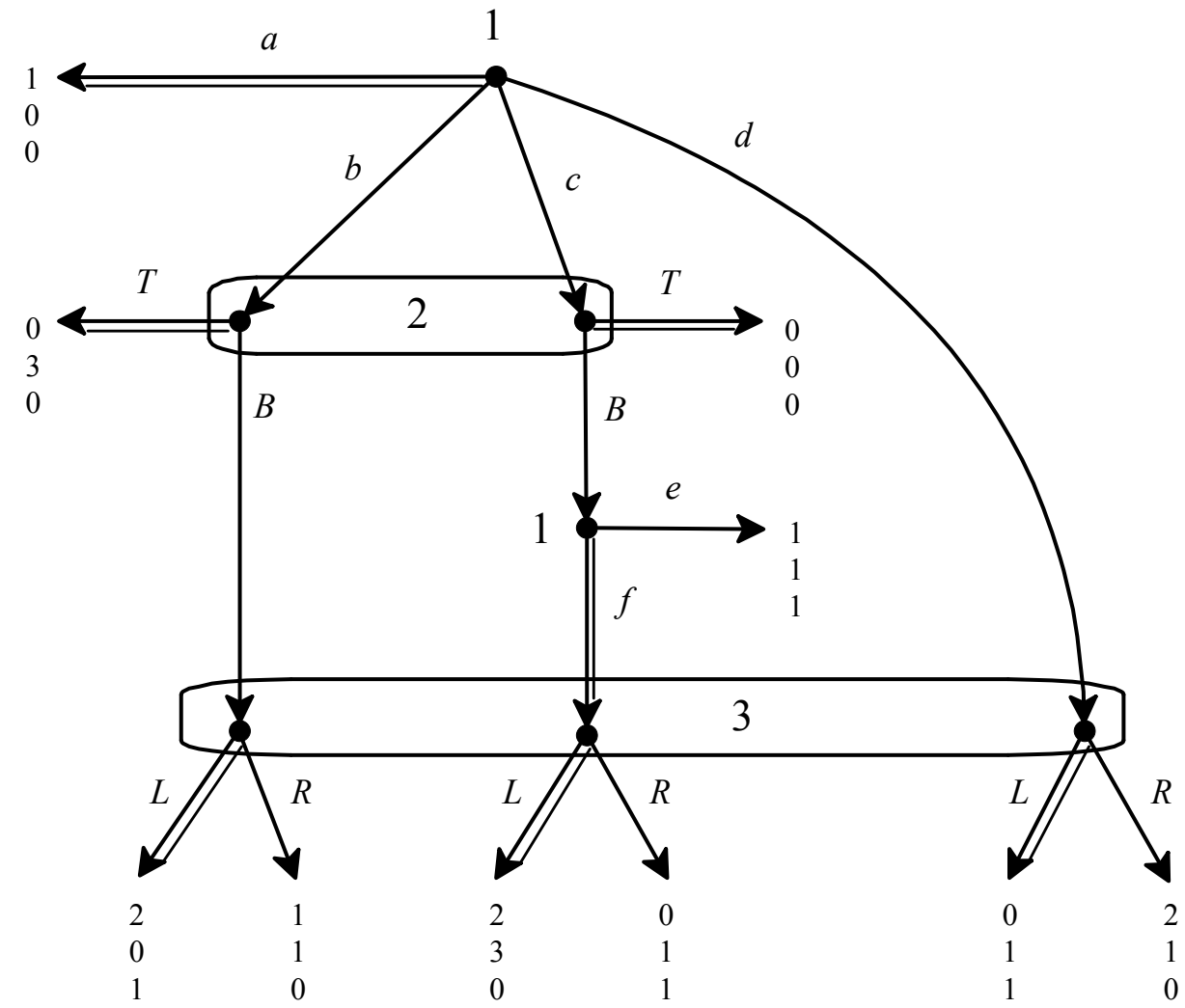

Let $(\sigma, \mu)$ be an assessment with $\sigma = (a, T, f, L)$, $\mu(b) > 0$ and $\mu(c) > 0$.

**(a)** We have to prove that $(\sigma, \mu)$ can be rationalized by a choice-measurable plausibility order only if $\mu$ satisfies the following condition:

$$\mu(bB) > 0 \text{ if and only if } \mu(cBf) > 0 .$$

Let $\precsim$ be a choice measurable plausibility order that rationalizes $(\sigma, \mu)$ and let $F$ be an integer-valued representation of $\precsim$ that satisfies choice measurability. Since $\mu(b) > 0$ and $\mu(c) > 0$, by P2 of Definition 12.2, $b \sim c$ and thus $F(b) = F(c)$; by choice measurability, $F(b) - F(c) = F(bB) - F(cB)$ and thus $F(bB) = F(cB)$, so that $bB \sim cB$. Since $\sigma(f) > 0$, by P1 of Definition 12.2, $cB \sim cBf$ and therefore, by transitivity of $\precsim$, $bB \sim cBf$. Hence if $\mu(bB) > 0$ then, by P2 of Definition 12.2, $bB \in Min_{\precsim}\{bB, cBf, d\}$ (for any set $S \subseteq H$, $Min_{\precsim} S = \{h \in S : h \precsim h', \text{ for all } h' \in S\}$) and thus $cBf \in Min_{\precsim}\{bB, cBf, d\}$ so that, by P2 of Definition 12.2, $\mu(cBf) > 0$. The proof that if $\mu(cBf) > 0$ then $\mu(bB) > 0$ is analogous.

**(b)** We have to prove that if, besides from being rationalized by a choice-measurable plausibility order $\precsim$, $(\sigma,\mu)$ is also uniformly Bayesian relative to $\precsim$ (Definition 12.8), then $\mu$ satisfies the following condition:

$$\text{if } \mu(bB) > 0 \text{ then } \frac{\mu(cBf)}{\mu(bB)} = \frac{\mu(c)}{\mu(b)}.$$

Suppose that $\mu(b) > 0$, $\mu(c) > 0$ (so that $b \sim c$) and $\mu(bB) > 0$. Let $\nu$ be a full-support common prior that satisfies the properties of Definition12.8. Then, by UB2, $\frac{\nu(c)}{\nu(b)} = \frac{\nu(cB)}{\nu(bB)}$ and, by UB1, since $\sigma(f) = 1$, $\nu(cBf) = \nu(cB) \times \sigma(f) = \nu(cB)$. Let $E$ be the equivalence class that contains *b*. Then $E \cap D_\mu^+ = \{b,c\}$. Since $\nu_E(\cdot) = \nu(\cdot \,|\, E \cap D_\mu^+)$, by B3 of Definition 12.3, $\mu(b) = \frac{\nu(b)}{\nu(b)+\nu(c)}$ and $\mu(c) = \frac{\nu(c)}{\nu(b)+\nu(c)}$, so that $\frac{\mu(c)}{\mu(b)} = \frac{\nu(c)}{\nu(b)}$. Let *G* be the equivalence class that contains *bB*. Then, since – by hypothesis – $\mu(bB) > 0$, it follows from the condition proved in part (a) that either $G \cap D_\mu^+ = \{bB, cBf\}$ or $G \cap D_\mu^+ = \{bB, cBf, d\}$. Since $\nu_G(\cdot) = \nu(\cdot \,|\, G \cap D_\mu^+)$, by B3 of Definition 12.3, in the former case $\mu(bB) = \frac{\nu(bB)}{\nu(bB)+\nu(cBf)}$ and $\mu(cBf) = \frac{\nu(cBf)}{\nu(bB)+\nu(cBf)}$ and in the latter case $\mu(bB) = \frac{\nu(bB)}{\nu(bB)+\nu(cBf)+\nu(d)}$ and $\mu(cBf) = \frac{\nu(cBf)}{\nu(bB)+\nu(cBf)+\nu(d)}$ ; thus in both cases $\frac{\mu(cBf)}{\mu(bB)} = \frac{\nu(cBf)}{\nu(bB)}$. Hence, since $\nu(cBf) = \nu(cB)$, $\frac{\mu(cBf)}{\mu(bB)} = \frac{\nu(cB)}{\nu(bB)}$ and, therefore, since (as shown above) $\frac{\nu(cB)}{\nu(bB)} = \frac{\nu(c)}{\nu(b)}$ and $\frac{\nu(c)}{\nu(b)} = \frac{\mu(c)}{\mu(b)}$, we have that $\frac{\mu(cBf)}{\mu(bB)} = \frac{\mu(c)}{\mu(b)}$.

**Exercise 12.11** (Challenging Question). The game under consideration is:

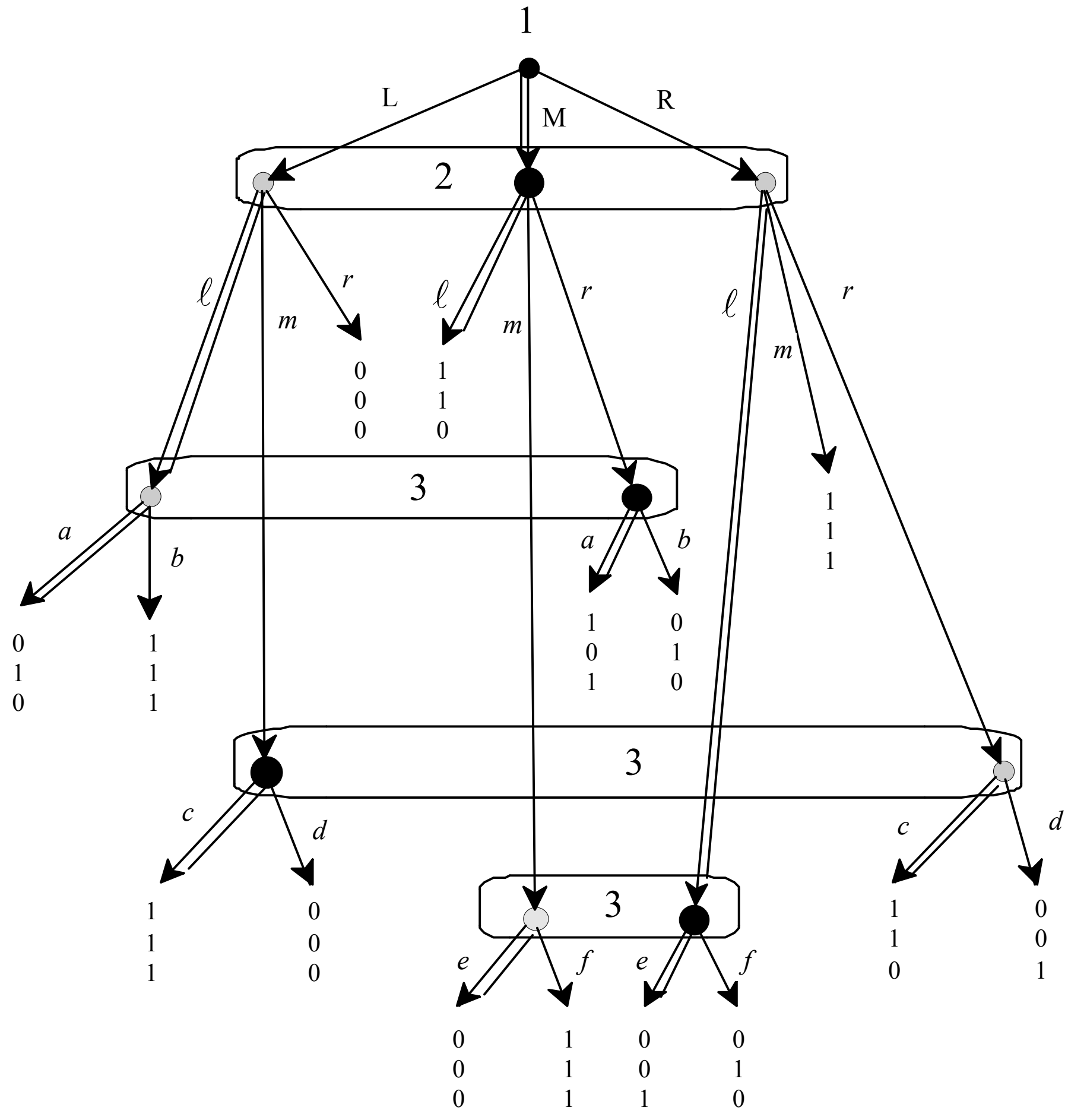

Let $\sigma = (M, \ell, a, c, e)$ and $\mu = \begin{pmatrix} L & M & R & L\ell & Mr & Lm & Rr & Mm & R\ell \\ 0 & 1 & 0 & 0 & 1 & 1 & 0 & 0 & 1 \end{pmatrix}$.

**(a)** We have to show that $(\sigma, \mu)$ is an independent perfect Bayesian equilibrium. Sequential rationality is straightforward to verify. It is also straightforward to verify that the following plausibility order rationalizes $(\sigma, \mu)$:

$$\begin{pmatrix} \emptyset, M, M\ell & \text{most plausible} \\ R, R\ell, R\ell e & \\ Mm, Mme & \\ Mr, Mra & \\ L, L\ell, L\ell a & \\ Rm & \\ Lm, Lmc & \\ Rr, Rrc & \\ Lr & \\ R\ell f & \\ Mmf & \\ Lmd & \\ Rrd & \\ Mrb & \\ L\ell b & \text{least plausible} \end{pmatrix}$$

Let us check that the above plausibility order satisfies properties $IND_1$ and $IND_2$. For $IND_1$ first note that there are no two decision histories $h$ and $h'$ that belong to the same information set and are such that $h \sim h'$; thus we only need to check that if $h$ and $h'$ belong to the same information set and $h \prec h'$ then $ha \prec h'a$ for every $a \in A(h)$. This is indeed true:

(1) $M \prec R$ and $M\ell \prec R\ell$, $Mm \prec Rm$ and $Mr \prec Rr$,
(2) $M \prec L$ and $M\ell \prec L\ell$, $Mm \prec Lm$ and $Mr \prec Lr$,
(3) $Mr \prec L\ell$ and $Mra \prec L\ell a$ and $Mrb \prec L\ell b$,
(4) $Lm \prec Rr$ and $Lmc \prec Rrc$ and $Lmd \prec Rrd$,
(5) $R\ell \prec Mm$ and $R\ell e \prec Mme$ and $R\ell f \prec Mmf$.

For $IND_2$ first note that there is no decision history $h$ that $ha \sim hb$ for $a, b \in A(h)$ with $a \neq b$; thus we only need to show that if $ha \prec hb$ and $h' \in I(h)$ then $h'a \prec h'b$. This is indeed true:

(1) $M\ell \prec Mm$ and $L\ell \prec Lm$ and $R\ell \prec Rm$,
(2) $Mm \prec Mr$ and $Lm \prec Lr$ and $Rm \prec Rr$,

and the rest is trivial, since at the other information sets there are only two actions, one of which is plausibility preserving and the other is not.

Finally, $\mu$ satisfies $IND_3$ trivially since at each information set only one history is assigned positive probability by $\mu$. Thus we have shown that $(\sigma, \mu)$ is an independent perfect Bayesian equilibrium.

**(b)** To prove that $(\sigma,\mu)$ is not a sequential equilibrium it would not be sufficient to prove that the plausibility order given above is not choice measurable (although it is indeed true that it is not choice measurable), because in principle there could be another plausibility order which is choice measurable and rationalizes $(\sigma,\mu)$. Thus we need to show that *any* plausibility order that rationalizes $(\sigma,\mu)$ is not choice measurable. Let $\precsim$ be a plausibility order that rationalizes $(\sigma,\mu)$; then it must satisfy the following properties:

- $Lm \prec Rr$ (because they belong to the same information set and $\mu(Lm) > 0$ while $\mu(Rr) = 0$). Thus if *F* is any integer-valued representation of $\precsim$ it must be that
$$F(Lm) < F(Rr). \tag{1}$$
- $Mr \prec L\ell \sim L$ (because $Mr$ and $L\ell$ belong to the same information set and $\mu(Mr) > 0$ while $\mu(L\ell) = 0$; furthermore, $\ell$ is a plausibility-preserving action since $\sigma(\ell) > 0$). Thus if *F* is any integer-valued representation of $\precsim$ it must be that
$$F(Mr) < F(L). \tag{2}$$
- $R \sim R\ell \prec Mm$ (because $\ell$ is a plausibility-preserving action, $R\ell$ and $Mm$ belong to the same information set and $\mu(R\ell) > 0$ while $\mu(Mm) = 0$). Thus if *F* is any integer-valued representation of $\precsim$ it must be that
$$F(R) < F(Mm). \tag{3}$$

Now suppose that $\precsim$ is choice measurable and let *F* be an integer-valued representation of it that satisfies choice measurability (Property CM of Definition 12.7). From (1) and (2) we get that

$$F(Lm) - F(L) < F(Rr) - F(Mr) \tag{4}$$

and by choice measurability

$$F(Rr) - F(Mr) = F(R) - F(M) \tag{5}$$

It follows from (4) and (5) that

$$F(Lm) - F(L) < F(R) - F(M) \tag{6}$$

Subtracting $F(M)$ from both sides of (3) we obtain

$$F(R) - F(M) < F(Mm) - F(M). \tag{7}$$

It follows from (6) and (7) that $F(Lm) - F(L) < F(Mm) - F(M)$, which can be written as $F(M) - F(L) < F(Mm) - F(Lm)$, yielding a contradiction, because choice measurability requires that $F(M) - F(L) = F(Mm) - F(Lm)$.

# PART V

## Advanced Topics III:

## Incomplete Information

Chapter

13

# Incomplete Information: Static Games

## 13.1 Interactive situations with incomplete information

An implicit assumption in game theory is that the game being played is common knowledge among the players. The expression "incomplete information" refers to situations where some of the elements of the game (e.g. the actions available to the players, the possible outcomes, the players' preferences, etc.) are *not* common knowledge. In such situations the knowledge and beliefs of the players about the game need to be made an integral part of the description of the situation. Pioneering work in this direction was done by John Harsanyi (1967, 1968), who (together with John Nash and Reinhard Selten) was the recipient of the 1994 Nobel prize in economics. Harsanyi suggested a method for converting a situation of incomplete information into an extensive game with imperfect information (this is the so-called *Harsanyi transformation*). The theory of games of incomplete information has been developed for the case of von Neumann-Morgenstern payoffs and the solution concept proposed by Harsanyi is *Bayesian Nash equilibrium*, which is merely Nash equilibrium of the imperfect-information game that is obtained by applying the Harsanyi transformation.

Although the approach put forward by Harsanyi was in terms of "types" of players and of probability distributions over types, we shall develop the theory using the so-called "state-space" approach, which makes use of the interactive knowledge-belief structures developed in Chapter 8. In Chapter 15 we will explain the "type-space" approach and show how to convert one type of structure into the other.

The distinction between complete and incomplete information is not at all the same as that between perfect and imperfect information. To say that a game (in extensive form) has *imperfect* information is to say that there is at least one player who may have to make a choice in a situation where she does not know what choices were made previously by other players. To say that there is *incomplete* information is to say that there is at least one player who does not quite know what game she is playing. Consider, for example, the players' preferences. In some situations it is not overly restrictive to assume that the players' preferences are common knowledge: for instance, in the game of chess, it seems quite plausible to assume that it is common knowledge that both players prefer winning to either drawing or losing (and prefer drawing to losing). But think of a contractual dispute. Somebody claims that you owe him a sum of money and threatens to sue you if you don't pay. We can view this situation as a two-player game: you have two strategies, "pay" and "not pay", and he has two strategies, "sue (if no payment)" and "not sue (if no payment)". In order to determine your best choice, you need to know how he will respond to your refusal to pay. A lawsuit is costly for him as well as for you; if he is the "tough" type he will sue you; if he is the "soft" type he will drop the dispute. If you don't know what type he is, you are in a situation of incomplete information.

As explained above, we have a situation of incomplete information whenever at least one of the players does not know some aspects of the game: it could be the preferences of her opponents, the choices available to her opponents or the set of possible outcomes. We shall focus almost exclusively on the case where the uncertainty concerns the preferences of the players, while everything else will be assumed to be common knowledge (for an exception see Exercise 13.2). Harsanyi argued that every situation can be reduced to this case. For example, he argued, if Player 1 does not know whether Player 2 has available only choices *a* and *b* or also choice *c*, we can model this as a situation where there are two possible "states of the world" in both of which Player 2 has three choices available, *a*, *b* and *c*, but in one of the two states choice *c* gives an extremely low payoff to Player 2, so that she would definitely not choose *c*.

The interactive knowledge-belief structures developed in Chapter 8 are sufficiently rich to model situations of incomplete information. The states in these structures can be used to express any type of uncertainty. In Chapter 9 we interpreted the states in terms of the actual choices of the players, thus expressing uncertainty in the mind of a player about the behavior of another player. In that case the game was assumed to be common knowledge among the players (that is, it was the same game at every state) and what varied from one state to another was the choice of at least one player. If we want to represent a situation where one player is not sure what game she is playing, all we have to do is interpret the states in terms of games, that is, assign different games to different states.

We shall begin with games in strategic form with cardinal payoffs where only one player is uncertain about what game is being played. We shall use the expression "one-sided incomplete information" to refer to these situations.

## 13.2 One-sided incomplete information

Let us model the following situation: the game being played is the one shown in Figure 13.1; call this the "true" game.

| | | 2 | |
|---|---|---|---|
| | | *L* | *R* |
| 1 | *T* | 6 , 3 | 0 , 9 |
| | *B* | 3 , 3 | 3 , 0 |

**Figure 13.1**

Player 1 knows that she is playing this game, while Player 2 is uncertain as to whether she is playing this game or a different game, shown in Figure 13.2, where the payoffs of Player 1 are different:

| | | 2 | |
|---|---|---|---|
| | | *L* | *R* |
| 1 | *T* | 0 , 3 | 3 , 9 |
| | *B* | 3 , 3 | 0 , 0 |

**Figure 13.2**

For convenience, let us refer to the "true game" of Figure 13.1 as the game where Player 1 is of type *b* and the game of Figure 13.2 as the game where Player 1 is of type $a$. Suppose that Player 2 assigns probability $\frac{2}{3}$ to Player 1 being of type $a$ and probability $\frac{1}{3}$ to Player 1 being of type $b$. The description is not complete yet, because we need to specify whether Player 1 has any uncertainty concerning the beliefs of Player 2; let us assume that the beliefs of Player 2 are common knowledge. We also need to specify whether Player 2 is uncertain about the state of mind of Player 1; let us assume that it is common knowledge that Player 2 knows that Player 1 knows what game is being played. This is a long verbal description! A picture is worth a thousand words and indeed the simple knowledge-belief structure shown in Figure 13.3 captures all of the above. There are two states, $\alpha$ and $\beta$: $\alpha$ is interpreted as a (counterfactual) state where the game of Figure 13.2 is played, while $\beta$ is interpreted as the state where the "true" game of Figure 13.1 is played. We capture the fact that latter is the "true" game, that is, the game which is actually played, by designating state $\beta$ as "the true state".[1]

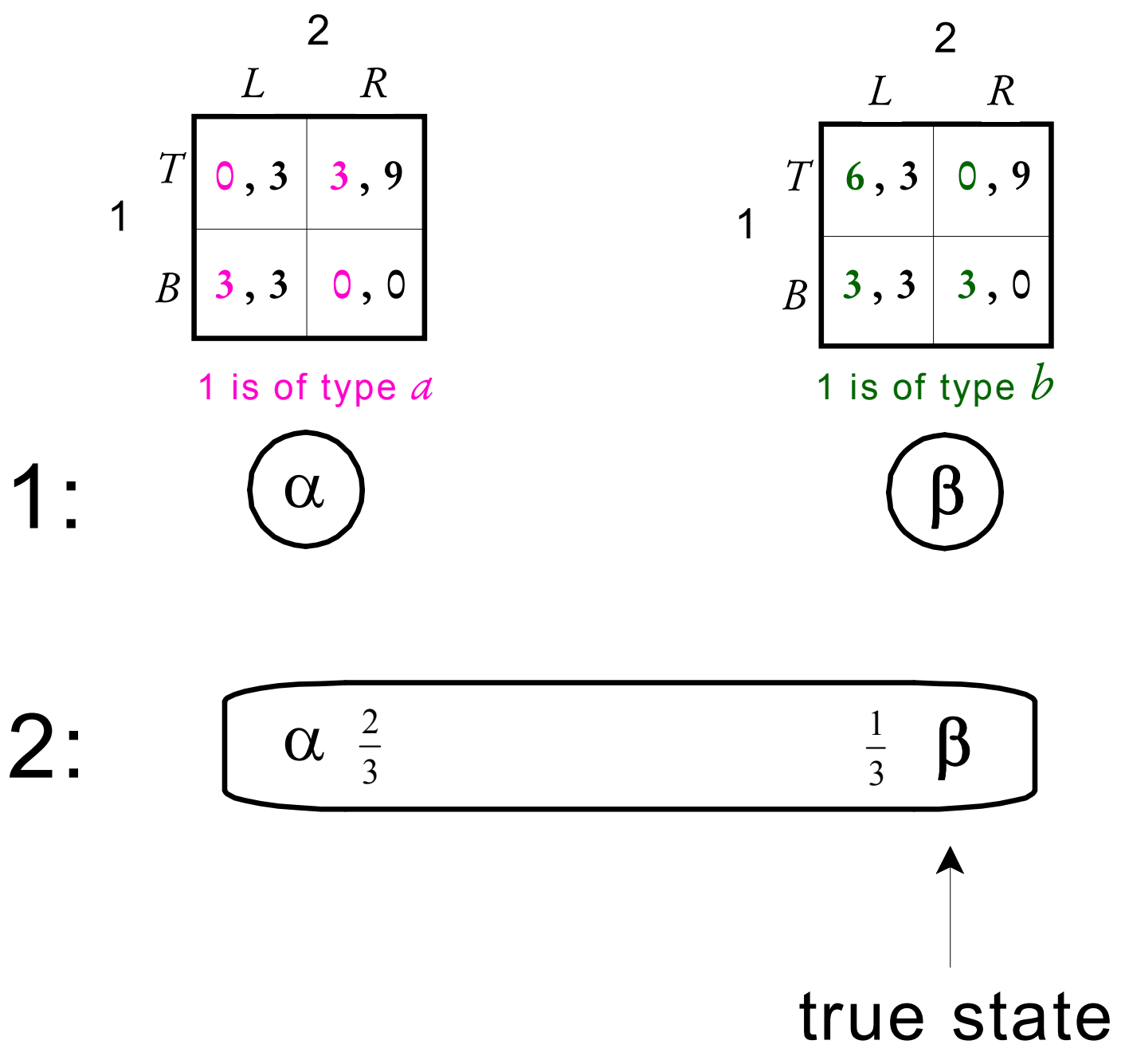

**Figure 13.3**

It is easy to check that, in the knowledge-belief structure of Figure 13.3, at state $\beta$ all of the elements of the verbal description given above are true. Anything that is constant across states is common knowledge: thus the payoffs of Player 2 are common knowledge, the beliefs of Player 2, namely $\begin{pmatrix} \alpha & \beta \\ \frac{2}{3} & \frac{1}{3} \end{pmatrix}$, are common knowledge and the fact that if the game being played is the one of Figure 13.2 (associated with state $\alpha$) then Player 1 knows that they are playing that game and that if the game being played is the one of Figure 13.1 (associated with state $\beta$) then Player 1 knows that they are playing that game, that is, Player 1 knows what game is being played. Note, however, that, at the true state $\beta$, Player 1 knows more than Player 2, namely which of the two games is actually being played.

[1] This is something that is almost never done in the literature, but it is an important element of a description of a situation of incomplete information: what is the actual state of affairs?

**Remark 13.1.** In the situation illustrated in Figure 13.3, at every state each player knows his/her own payoffs. It may seem natural to make this a requirement of rationality: shouldn't a rational player know her own preferences? The answer is Yes for preferences and No for payoffs: a rational player should know how she ranks the possible outcomes (her preferences over outcomes); however, it is perfectly rational to be uncertain about what outcome will follow a particular action and thus about one's own payoff. For example, you know that you prefer accepting a wrapped box containing a thoughtful gift (payoff of 1) to accepting a wrapped box containing an insulting gift (payoff of $-1$): you know your preferences over these two outcomes. However, you may be uncertain about the intentions of the gift-giver and thus you are uncertain about what outcome will follow if you accept the gift: you don't know if your payoff will be 1 or $-1$. Thus, although you know your payoff *function,* you may be uncertain about what your *payoff* will be if you accept the gift. This example is developed in Exercise 13.2.

Figure 13.3 illustrates a situation that has to do with games, but is not a game. Harsanyi's insight was that we can transform that situation into an extensive-form game with imperfect information, as follows. We start with a chance move, where Nature chooses the state; then Player 1 is informed of Nature's choice and makes her decision between $T$ and $B$; Player 2 then makes his choice between $L$ and $R$, without being informed of Nature's choice (to capture her uncertainty about the game) and without being informed of Player 1's choice (to capture the simultaneity of the game). The game is shown in Figure 13.4 below.

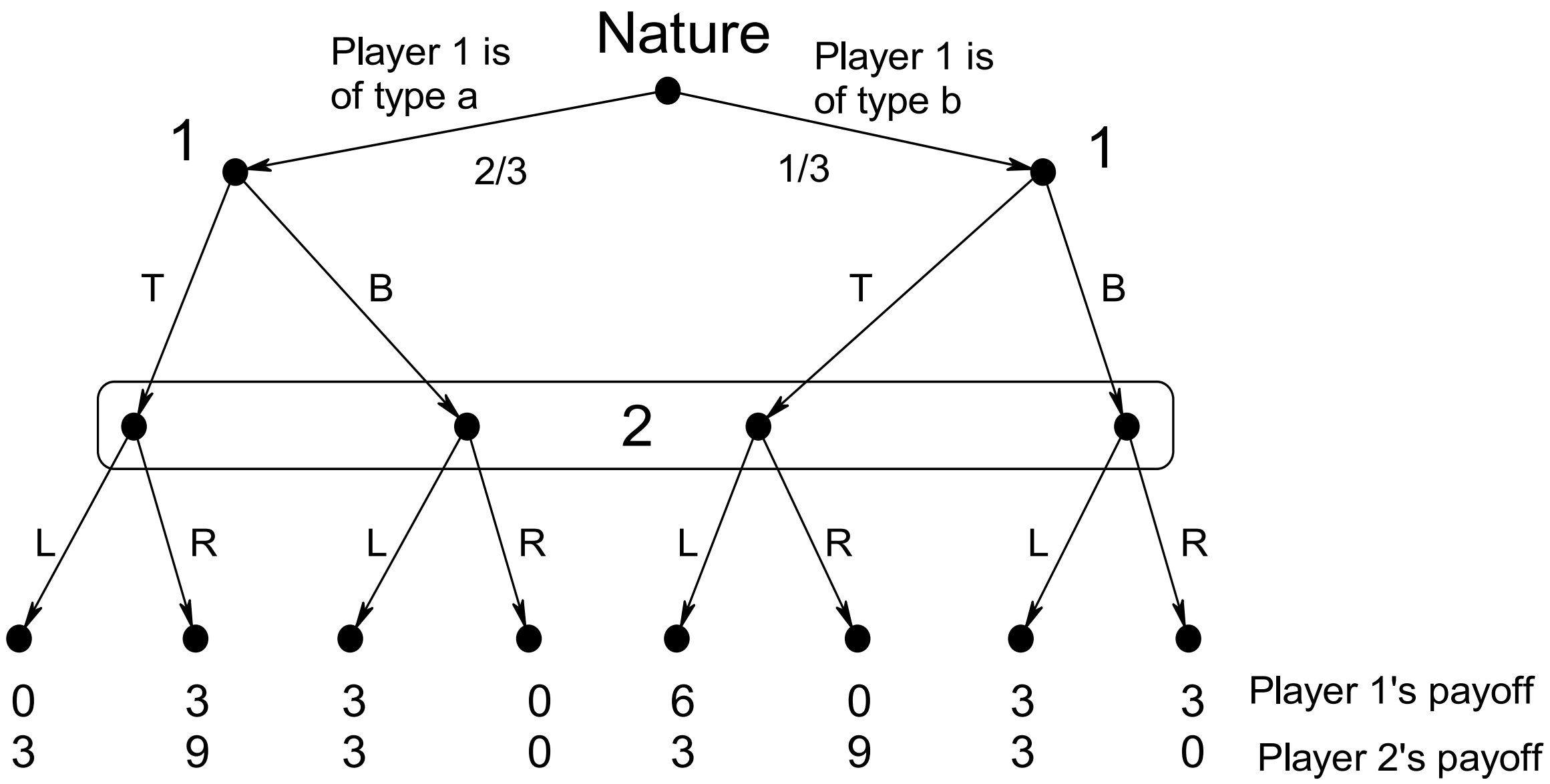

**Figure 13.4**

The reduction of the situation illustrated in Figure 13.3 to the extensive-form game of Figure 13.4 is called the *Harsanyi transformation.* Once a game has been obtained by means of the Harsanyi transformation, it can be solved using, for example, the notion of Nash equilibrium. The Nash equilibria of a game corresponding to a situation of incomplete information are called *Bayesian Nash equilibria.* Note that these are nothing more than Nash equilibria: the extra term 'Bayesian' is merely a hint that the game being solved is meant to represent a situation of incomplete information.

Note that the Harsanyi transformation involves some loss of information: in particular, in the resulting game one can no longer tell what the true state is, that is, what the actual game being played is.

To find the Bayesian Nash equilibria of the game of Figure 13.4 we can construct the corresponding strategic-form game, which is shown in Figure 13.5 below.[2]

[2] Recall that, when constructing the strategic-form game associated with an extensive-form game with chance moves, payoffs are *expected* payoffs: see Chapter 6.

| | | Player 2 | |
|---|---|---|---|
| | | *L* | *R* |
| **Player 1** | *T* if type *a*, *T* if type *b* | 2 , 3 | 2 , 9 |
| | *T* if type *a*, *B* if type *b* | 1 , 3 | **3 , 6** |
| | *B* if type *a*, *T* if type *b* | **4 , 3** | 0 , 3 |
| | *B* if type *a*, *B* if type *b* | 3 , 3 | 1 , 0 |

**Figure 13.5**

There are two pure-strategy Bayesian Nash equilibria: $((T,B),R)$ and $((B,T),L)$. How should we interpret them?

Let us consider, for example, the Bayesian Nash equilibrium $((T,B),R)$. This is a Nash equilibrium of the game of Figure 13.4; however, we assumed that the true state was $\beta$, where Player 1 is of type *b* and thus, in the actual game (the one associated with state $\beta$), the actual play is (*B*,*R*), which is not a Nash equilibrium of that game (because, while *B* is a best reply to *R*, *R* is not a best reply to *B*). This is not surprising, since Player 1 knows that she is playing that game, while Player 2 attaches only probability $\frac{1}{3}$ to that game being the actual game. Thus the first observation is that a Bayesian Nash equilibrium of a "game of incomplete information" does not imply that the players play a Nash equilibrium in the actual (or "true") game. The second observation has to do with Player 1's strategy. By definition, Player 1's strategy in the game of Figure 13.4 consists of a pair of choices, one for the case where she is informed that her type is *a* and the other for the case where she is informed that her type is *b*: in the Bayesian Nash equilibrium under consideration, Player 1's strategy (*T*,*B*) means "Player 1 plays *T* if of type *a* and plays *B* if of type *b*". However, Player 1 *knows* that she is of type *b*: she knows what game she is playing. Why, then, should she formulate a plan on how she would play in a counterfactual world where her type was *a* (that is, in a counterfactual game)? The answer is that Player 1's strategy (*T,B*) is best understood not as a contingent plan of action formulated by Player 1, but as a complex object that incorporates (1) the actual choice (namely to play *B*) made by Player 1 in the game that she knows she is playing and (2) a belief in the mind of Player 2 about what Player 1 would do in the two games that, as far as Player 2 knows, are actual possibilities.

An alternative (and easier) way to find the Bayesian Nash equilibria of the game of Figure 13.4 is to use the notion of weak sequential equilibrium for that game. For example, to verify that $((B,T),L)$ is a Nash equilibrium we can reason as follows: If the strategy of Player 1 is (*B,T*), then Bayesian updating requires Player 2 to have the following beliefs: probability $\frac{2}{3}$ on the second node from the left and probability $\frac{1}{3}$ on the third node from the left. Given these beliefs, playing *L* yields him an expected payoff of $\frac{2}{3}(3)+\frac{1}{3}(3)=3$ while playing *R* yields him an expected payoff of $\frac{2}{3}(0)+\frac{1}{3}(9)=3$; thus any strategy (pure or mixed) is sequentially rational: in particular *L* is sequentially rational. If Player 2 plays *L* then at her left node Player 1 gets 0 with *T* and 3 with *B*, so that *B* is sequentially rational; at her right node Player 1 gets 6 with *T* and 3 with *B*, so that *T* is sequentially rational. Hence $((B,T),L)$ with the stated beliefs is a weak sequential equilibrium, implying that $((B,T),L)$ is a Nash equilibrium (Theorem 10.1, Chapter 10).

**Definition 13.1.** In a "game of one-sided incomplete information" a pure-strategy Bayesian Nash equilibrium where the "informed" player makes the same choice at every singleton node is called a *pooling equilibrium*, while a pure-strategy Bayesian Nash equilibrium where the "informed" player makes different choices at different nodes is called a *separating equilibrium.*

Thus in the game of Figure 13.4 there are no pooling equilibria: all the pure-strategy Bayesian Nash equilibria are separating equilibria.

**Remark 13.2.** Although we have only looked at the case of two players, situations of one-sided incomplete information can involve any number of players, as long as only one player is uncertain about the game being played (while all the other players are not). An example of a situation of one-sided incomplete information with three players will be given in the next chapter with reference to Selten's chain-store game (which we studied in Chapter 2, Section 2.4).

We conclude this section with an example of a situation of one-sided incomplete information involving a two-player game where each player has an infinite number of strategies. This is an incomplete-information version of Cournot's game of competition between two firms (Chapter 1, Section 1.7). This example involves the use of calculus and some readers might want to skip it.

Consider a Cournot duopoly (that is, an industry consisting of two firms, which compete in output levels) with inverse demand given by

$$P(Q) = 34 - Q,$$

where $Q = q_1 + q_2$ is total industry output ($q_1$ is the output of Firm 1 and $q_2$ is the output of Firm 2). It is common knowledge between the two firms that Firm 1's cost function is given by:

$$C_1(q_1) = 6q_1.$$

Firm 2' s cost function is

$$C_2(q_2) = 9q_2.$$

Firm 2 knows this, while Firm 1 believes that Firm 2's cost function is $C_2(q_2) = 9q_2$ with probability $\frac{1}{3}$ and $C_2(q_2) = 3q_2$ with probability $\frac{2}{3}$ (Firm 2 could be a new entrant to the industry using an old technology, or could have just invented a new technology). Firm 1 knows that Firm 2 knows its own cost function. Everything that Firm 1 knows is common knowledge between the two firms. Thus we have the situation of one-sided incomplete information illustrated in Figure 13.6.

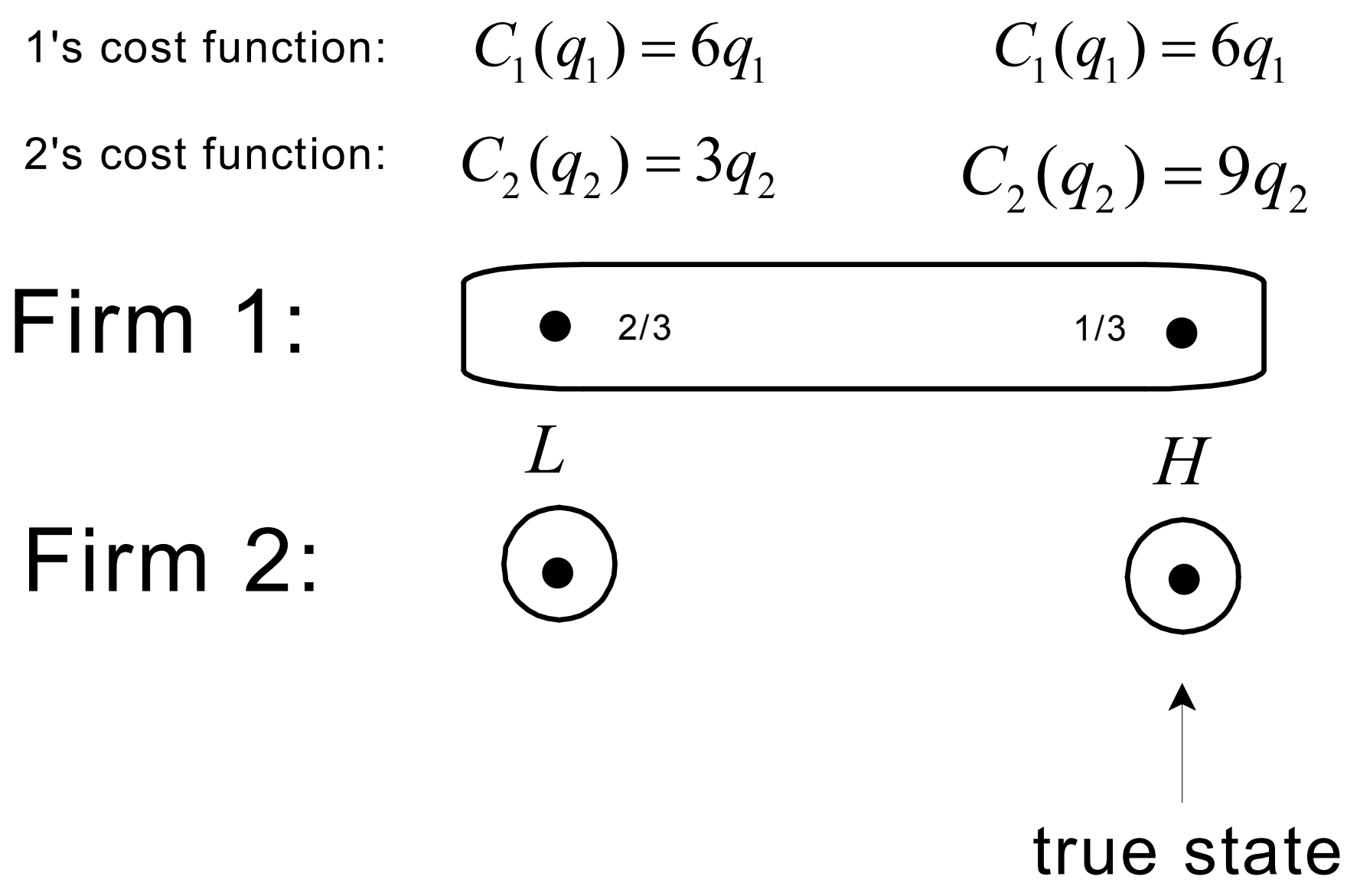

**Figure 13.6**

Let us find a Bayesian Nash equilibrium of the game obtained by applying the Harsanyi transformation to this situation of incomplete information and compare it to the Nash equilibrium of the complete information game where Firm 2's cost function is common knowledge.

In the complete information case where it is common knowledge that Firm 2's cost function is $C_2(q_2) = 9q_2$, the profit (= payoff) functions of the two firms are given by:

$$\pi_1(q_1,q_2) = \left[34 - (q_1 + q_2)\right]q_1 - 6q_1$$
$$\pi_2(q_1,q_2) = \left[34 - (q_1 + q_2)\right]q_2 - 9q_2$$

The Nash equilibrium is given by the solution to the following pair of equations:

$$\frac{\partial \pi_1(q_1,q_2)}{\partial q_1} = 34 - 2q_1 - q_2 - 6 = 0$$
$$\frac{\partial \pi_2(q_1,q_2)}{\partial q_2} = 34 - q_1 - 2q_2 - 9 = 0$$

The solution is $q_1 = \frac{31}{3} = 10.33$ and $\mathrm{q}_2 = \frac{22}{3} = 7.33$.

Next we go back to the incomplete-information situation illustrated in Figure 13.6. Although it is not possible to draw the extensive-form game that results from applying the Harsanyi transformation (because of the infinite number of possible output levels) we can nevertheless sketch the game as shown in Figure 13.7 (where *H* means "high cost" and *L* "low cost"):

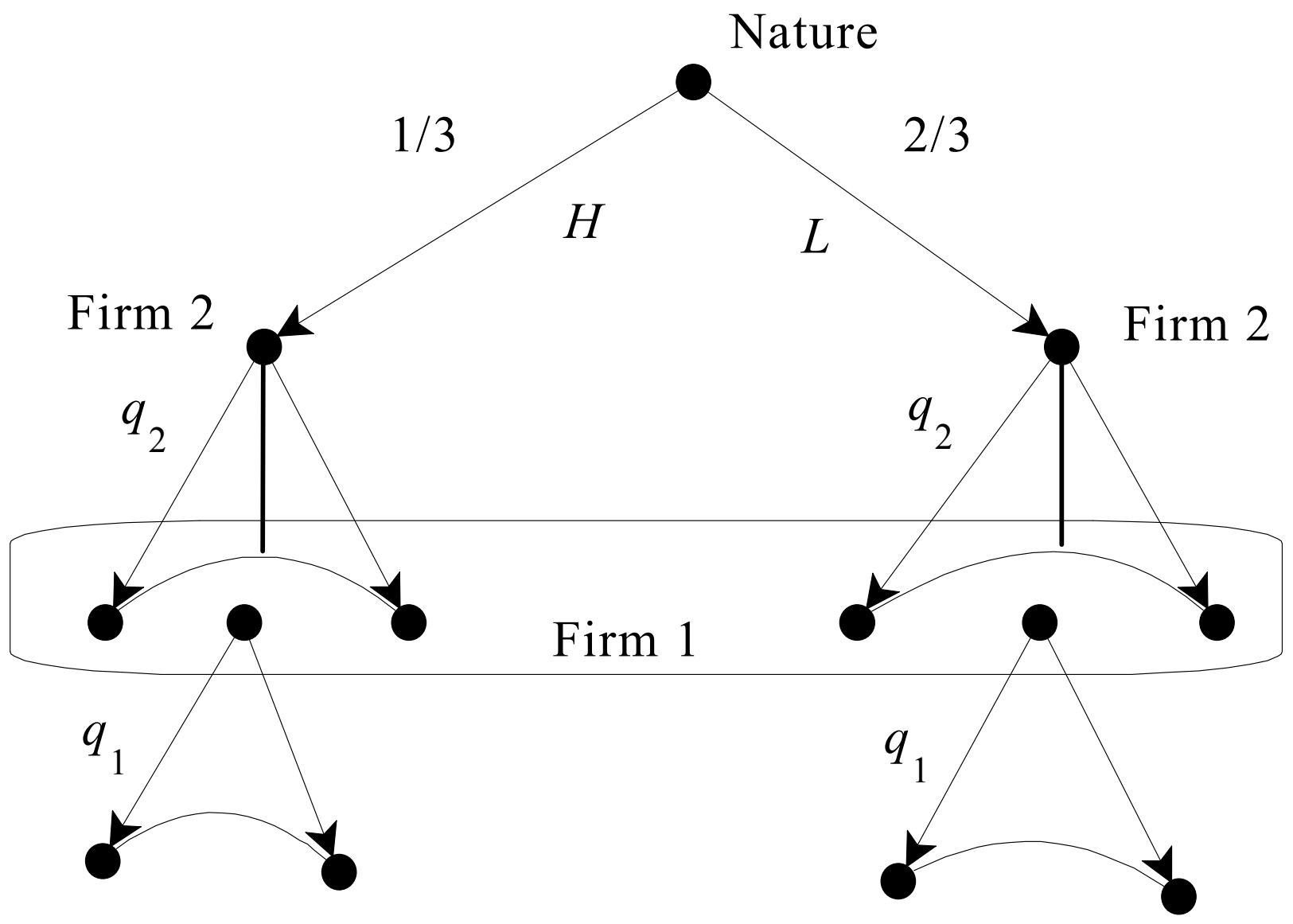

**Figure 13.7**

To find a Bayesian Nash equilibrium it is easiest to think in terms of weak sequential equilibrium: if the strategy profile $\left(\left(\hat{q}_2^H,\hat{q}_2^L\right),\hat{q}_1\right)$ is part of a weak sequential equilibrium,[3] then – by sequential rationality – $\hat{q}_2^H$ must maximize the expression $[(34-\hat{q}_1-q_2^H)q_2^H-9q_2^H]$ and $\hat{q}_2^L$ must maximize the expression $[(34-\hat{q}_1-q_2^L)q_2^L-3q_2^L]$ ; furthermore, by Bayesian updating, Firm 1 must assign probability $\frac{1}{3}$ to the node following choice $\hat{q}_2^H$ and probability $\frac{2}{3}$ to the node following choice $\hat{q}_2^L$ so that $\hat{q}_1$ must maximize the expression $\frac{1}{3}[(34-q_1-\hat{q}_2^H)q_1-6q_1]+\frac{2}{3}[(34-q_1-\hat{q}_2^L)q_1-6q_1]$. Thus $\left(\left(\hat{q}_2^H,\hat{q}_2^L\right),\hat{q}_1\right)$ must be a solution to the following system of three equations:

[3] $\hat{q}_2^H$ is the output chosen by Firm 2 at the left node, where its cost is high, and $\hat{q}_2^L$ is the output chosen by Firm 2 at the right node, where its cost is low.

$$\frac{\partial}{\partial q_2^H}[(34 - q_1 - q_2^H)q_2^H - 9q_2^H] = 0,$$

$$\frac{\partial}{\partial q_2^L}[(34 - q_1 - q_2^L)q_2^L - 3q_2^L] = 0,$$

$$\frac{\partial}{\partial q_1}\left(\tfrac{1}{3}[(34 - q_1 - \hat{q}_2^H)q_1 - 6q_1] + \tfrac{2}{3}[(34 - q_1 - \hat{q}_2^L)q_1 - 6q_1]\right) = 0.$$

There is a unique solution given by

$$\hat{q}_1 = 9,\ \hat{q}_2^H = 8,\ \hat{q}_2^L = 11.$$

Since the true state is where Firm 2's cost is high, the actual output levels are $\hat{q}_1 = 9$ for Firm 1 and $\hat{q}_2^H = 8$ for Firm 2. Thus in the incomplete information case Firm 2's output is higher than in the complete information case and Firm 1's output is lower than in the complete information case.

This is a good time to test your understanding of the concepts introduced in this section, by going through the exercises in Section 13.E.2 of Appendix 13.E at the end of this chapter.

# 13.3 Two-sided incomplete information

Let us now consider a situation involving two players, both of whom face some uncertainty. Such a situation is called a *two-sided incomplete-information situation.* Note that for such a situation to arise, it is not necessary that both players are uncertain about some "objective" aspect of the game (such as the preferences of the opponent): one of the two players might simply be uncertain about the beliefs of the other player. A situation of this sort is illustrated in Figure 13.8. Let $G$ be the game associated with states $\alpha$ and $\beta$ (it is the same game) and $G'$ the game associated with state $\gamma$. Suppose that the true state is $\alpha$. Then all of the following are true at state $\alpha$:

1) both Player 1 and Player 2 know that they are playing game $G$ (that is, neither player has any uncertainty about the objective aspects of the game),
2) Player 1 knows that Player 2 knows that they are playing game $G$,

3) Player 2 is uncertain as to whether Player 1 knows that they are playing $G$ (which is the case if the actual state is $\alpha$) or whether Player 1 is uncertain (if the actual state is $\beta$) between the possibility that they are playing game $G$ and the possibility that they are playing game $G'$ and considers the two possibilities equally likely; furthermore, Player 2 attaches probability $\frac{2}{3}$ to the first case (where Player 1 knows that they are playing game $G$) and probability $\frac{1}{3}$ to the second case (where Player 1 is uncertain between game $G$ and game $G'$),
4) Player 1 knows the state of uncertainty of Player 2 (concerning Player 1, as described in point 3 above),
5) The payoffs of Player 1 are common knowledge; furthermore, it is common knowledge that Player 2 knows his own payoffs.

| 1 \ 2 | C | D |
|---|---|---|
| A | 1,3 | 0,1 |
| B | 0,0 | 3,1 |

| 1 \ 2 | C | D |
|---|---|---|
| A | 1,3 | 0,1 |
| B | 0,0 | 3,1 |

| 1 \ 2 | C | D |
|---|---|---|
| A | 1,0 | 0,2 |
| B | 0,0 | 3,1 |

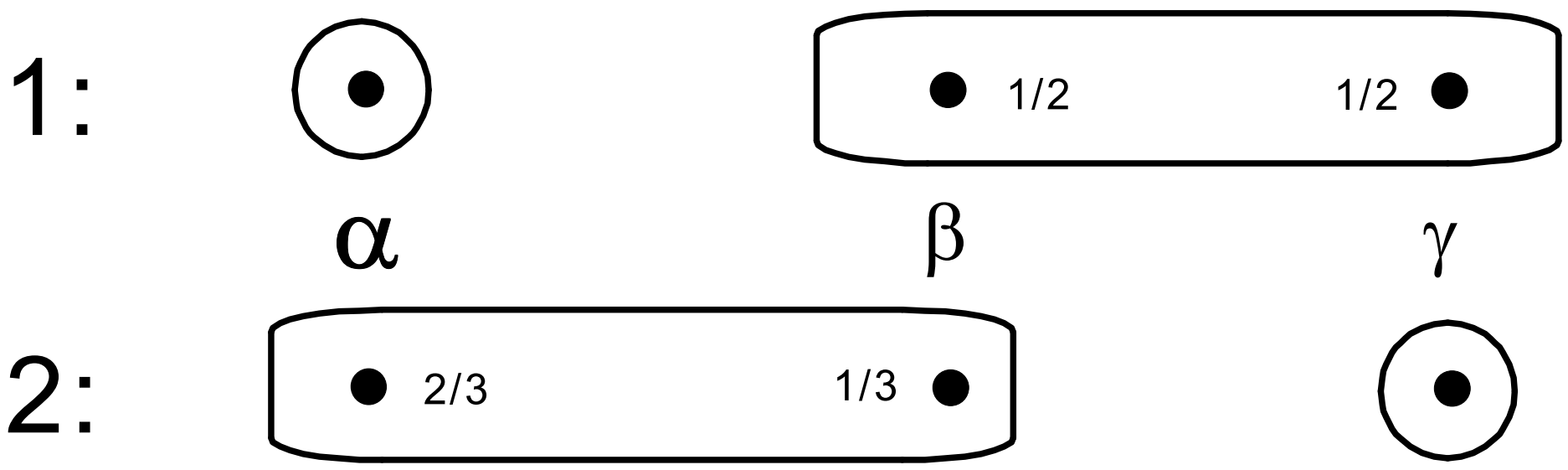

**Figure 13.8**

In principle, the Harsanyi transformation can be applied also to situations of two-sided incomplete information. However, in such cases there is an issue concerning the probabilities with which Nature chooses the states. In the case of one-sided incomplete information we can take those probabilities to be the beliefs of the uninformed player, but in the case of two-sided incomplete information we have two uninformed players and thus two sets of beliefs. For example, if we look at state $\beta$ in Figure 13.8, we have two different probabilities assigned to that state: $\frac{1}{2}$

by Player 1 and $\frac{1}{3}$ by Player 2. Which of the two should we take as Nature's probability for state $\beta$? The answer is: neither of them. What we should take as Nature's probabilities is a probability distribution over the set of states $\{\alpha,\beta,\gamma\}$ that reflects the beliefs of *both* players: we have encountered such a notion before, in Chapter 8: it is the notion of a *common prior* (Definition 8.10). A common prior is a probability distribution over the set of states that yields the players' beliefs upon conditioning on the information represented by a cell of an information partition. For example, in the situation illustrated in Figure 13.8 we are seeking a probability distribution $\nu:\{\alpha,\beta,\gamma\}\to[0,1]$ such that $\nu(\beta\,|\,\{\beta,\gamma\}) = \dfrac{\nu(\beta)}{\nu(\beta)+\nu(\gamma)} = \frac{1}{2}$, $\nu(\alpha\,|\,\{\alpha,\beta\}) = \dfrac{\nu(\alpha)}{\nu(\alpha)+\nu(\beta)} = \frac{2}{3}$ (and, of course, $\nu(\alpha)+\nu(\beta)+\nu(\gamma)=1$). In this case a common prior exists and is given by $\nu(\alpha)=\frac{2}{4}$ and $\nu(\beta)=\nu(\gamma)=\frac{1}{4}$. Using this common prior to assign probabilities to Nature's choices we obtain the imperfect-information game shown in Figure 13.9.

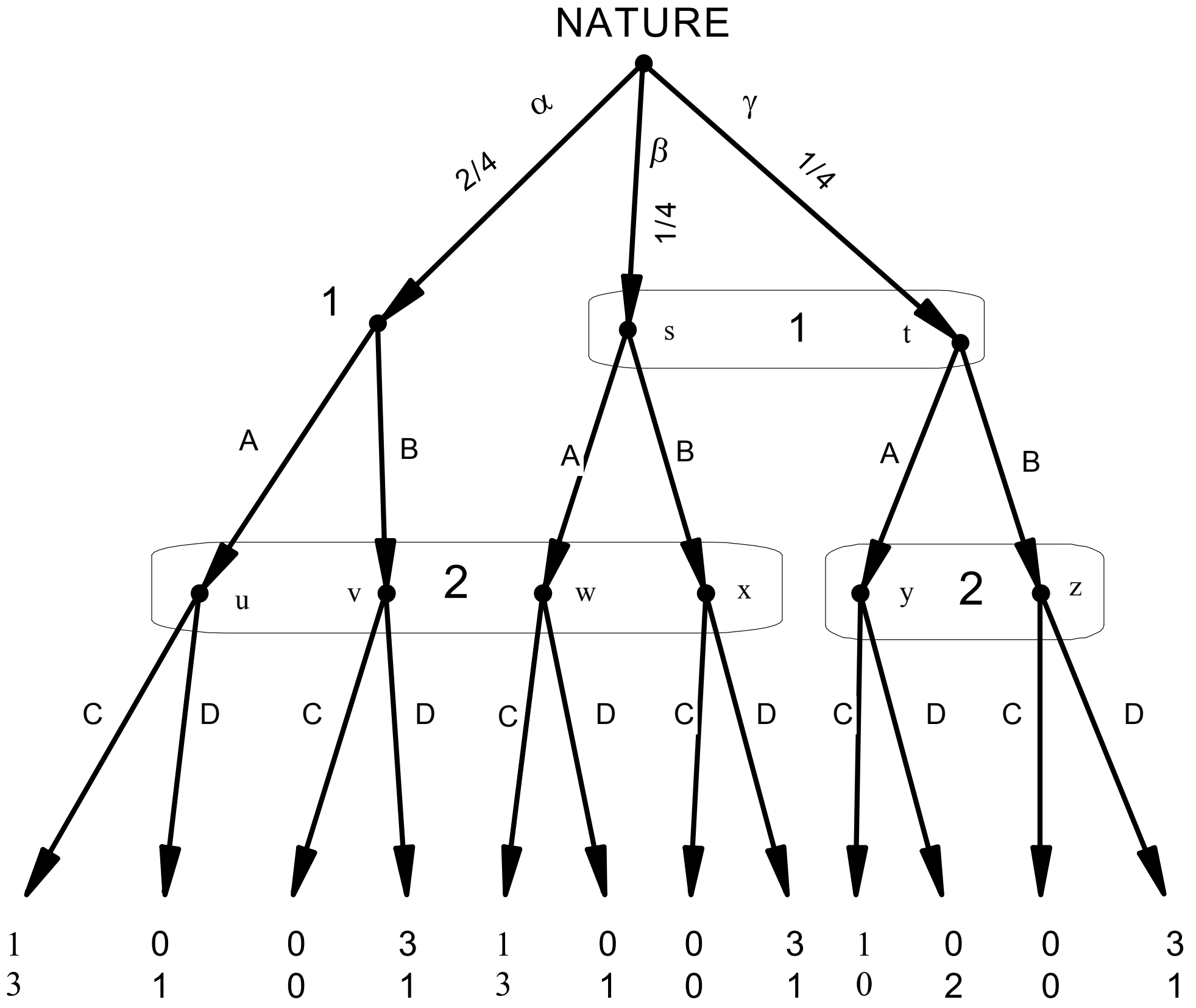

**Figure 13.9**

The following pure-strategy profile is Bayesian Nash equilibrium of the game of Figure 13.9: Player 1's strategy is *AB* (that is, he plays *A* if informed that the state is $\alpha$ and plays *B* if informed that the state is either $\beta$ or $\gamma$) and Player 2's strategy is *CD* (that is, she plays *C* at her information set on the left and *D* at her information set on the right). To verify that this is a Bayesian Nash equilibrium it is easier to verify that it is a weak sequential equilibrium together with the following system of belief, obtained by using Bayesian updating (note that every information set is reached by the strategy profile):

$$\mu = \begin{pmatrix} s & t & u & v & w & x & y & z \\ \frac{1}{2} & \frac{1}{2} & \frac{2}{3} & 0 & 0 & \frac{1}{3} & 0 & 1 \end{pmatrix}$$

Let us check sequential rationality. We begin with Player 1: (1) at the singleton node on the left, *A* gives Player 1 a payoff of 1 (given Player 2's choice of *C*) while *B* gives him a payoff of 0, hence *A* is sequentially rational; (2) at the information set on the right, given his beliefs and given Player 2's strategy *CD*, choosing *A* gives him an expected payoff of $\frac{1}{2}1+\frac{1}{2}0=\frac{1}{2}$ while *B* gives him an expected payoff of $\frac{1}{2}0+\frac{1}{2}3=\frac{3}{2}$, hence *B* is sequentially rational. Now consider Player 2: (1) at her information set on the left, given her beliefs, *C* gives her an expected payoff of $\frac{2}{3}3+\frac{1}{3}0=2$ while *D* gives her an expected payoff of $\frac{2}{3}1+\frac{1}{3}1=1$, hence *C* is sequentially rational; (2) at her information set on the right, given her beliefs, *C* gives her a payoff of 0 while *D* gives her a payoff of 1, hence *D* is sequentially rational.

The existence of a common prior is essential in order to be able to apply the Harsanyi transformation to a situation of two-sided incomplete information. In some cases a common prior does *not* exist (see Exercise 13.6) and thus the Harsanyi transformation cannot be carried out.

**Remark 13.3.** Besides the conceptual issues that arise in general with respect to the notion of Nash equilibrium, the notion of Bayesian Nash equilibrium for games with incomplete information raises the further issue of how one should understand or justify the notion of a common prior. This issue is not a trivial one and has been the object of debate in the literature.[4]

This is a good time to test your understanding of the concepts introduced in this section, by going through the exercises in Section 13.E.3 of Appendix 13.E at the end of this chapter.

[4] See, for example, Bonanno and Nehring (1999), Gul (1998) and Morris (1995).

# 13.4 Multi-sided incomplete information

So far we have only considered strategic-form games with two players. However, the analysis extends easily to games involving more than two players. If there are $n \geq 3$ players and only one player has uncertainty about some aspects of the game, while the others have no uncertainty whatsoever, then we have a situation of one-sided incomplete information; if two or more players have uncertainty (not necessarily all about the game but some possibly about the beliefs of other players) then we have a *multi-sided incomplete-information situation* (the two-sided case being a special case). We will consider an example of this below. Let $G_1$ and $G_2$ be the three-player games in strategic form (with cardinal payoffs) shown in Figure 13.10.

**GAME $G_1$**

Player 3 chooses E

| Player 1 \ Player 2 | C | | | D | | |
|---|---|---|---|---|---|---|
| A | 4 | 1 | 0 | 2 | 4 | 2 |
| B | 4 | 1 | 0 | 0 | 0 | 0 |

Player 3 chooses F

| Player 1 \ Player 2 | C | | | D | | |
|---|---|---|---|---|---|---|
| A | 2 | 0 | 2 | 0 | 1 | 0 |
| B | 0 | 0 | 0 | 0 | 2 | 0 |

**GAME $G_2$**

Player 3 chooses E

| Player 1 \ Player 2 | C | | | D | | |
|---|---|---|---|---|---|---|
| A | 4 | 4 | 0 | 2 | 1 | 2 |
| B | 4 | 4 | 0 | 0 | 1 | 0 |

Player 3 chooses F

| Player 1 \ Player 2 | C | | | D | | |
|---|---|---|---|---|---|---|
| A | 2 | 1 | 2 | 0 | 0 | 0 |
| B | 0 | 2 | 0 | 0 | 0 | 0 |

**Figure 13.10**

Consider the multi-sided situation of incomplete information illustrated in Figure 13.11, where with each state is associated one of the two games of Figure 13.10.

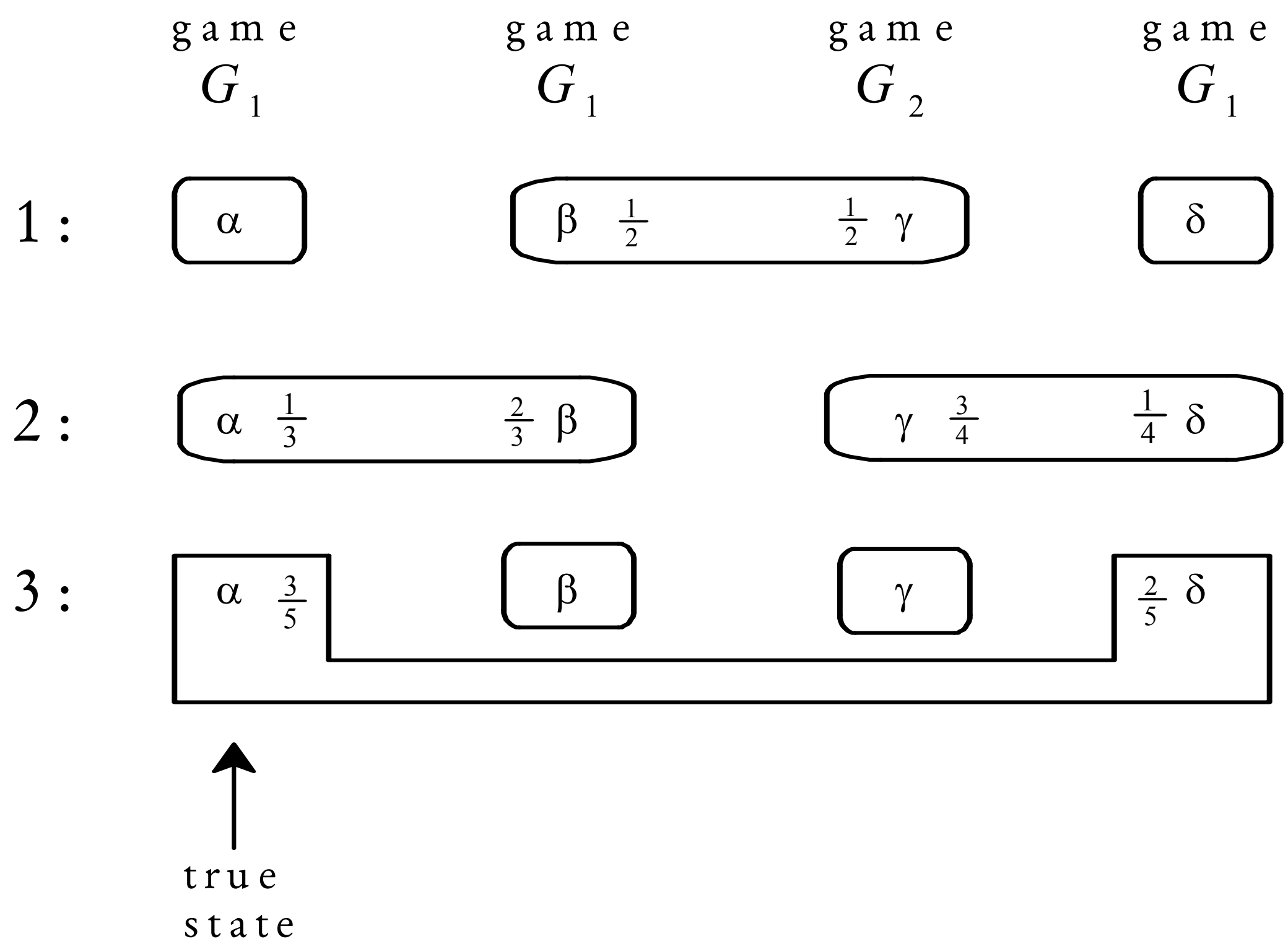

**Figure 13.11**

At the true state $\alpha$, all of the following are true:

- all three players know that they are playing game $G_1$,
- Player 1 knows that Players 2 and 3 know that the actual game is $G_1$;
- Player 2 knows that Player 3 knows that the actual game is $G_1$, but is uncertain as to whether Player 1 knows or is uncertain; furthermore, Player 2 assigns probability $\frac{2}{3}$ to Player 1 being uncertain;
- Player 3 knows that Player 1 knows that the actual game is $G_1$, but is uncertain as to whether Player 2 knows or is uncertain; furthermore, Player 3 assigns probability $\frac{2}{5}$ to Player 2 being uncertain;
- The payoffs of Players 1 and 3 are common knowledge.

The beliefs of the players in the situation illustrated in Figure 13.11 are compatible with each other, in the sense that there exists a common prior. In fact the following is a common prior: $\nu = \begin{pmatrix} \alpha & \beta & \gamma & \delta \\ \frac{3}{17} & \frac{6}{17} & \frac{6}{17} & \frac{2}{17} \end{pmatrix}$. Thus we can apply the Harsanyi transformation and obtain the game shown in Figure 13.12.

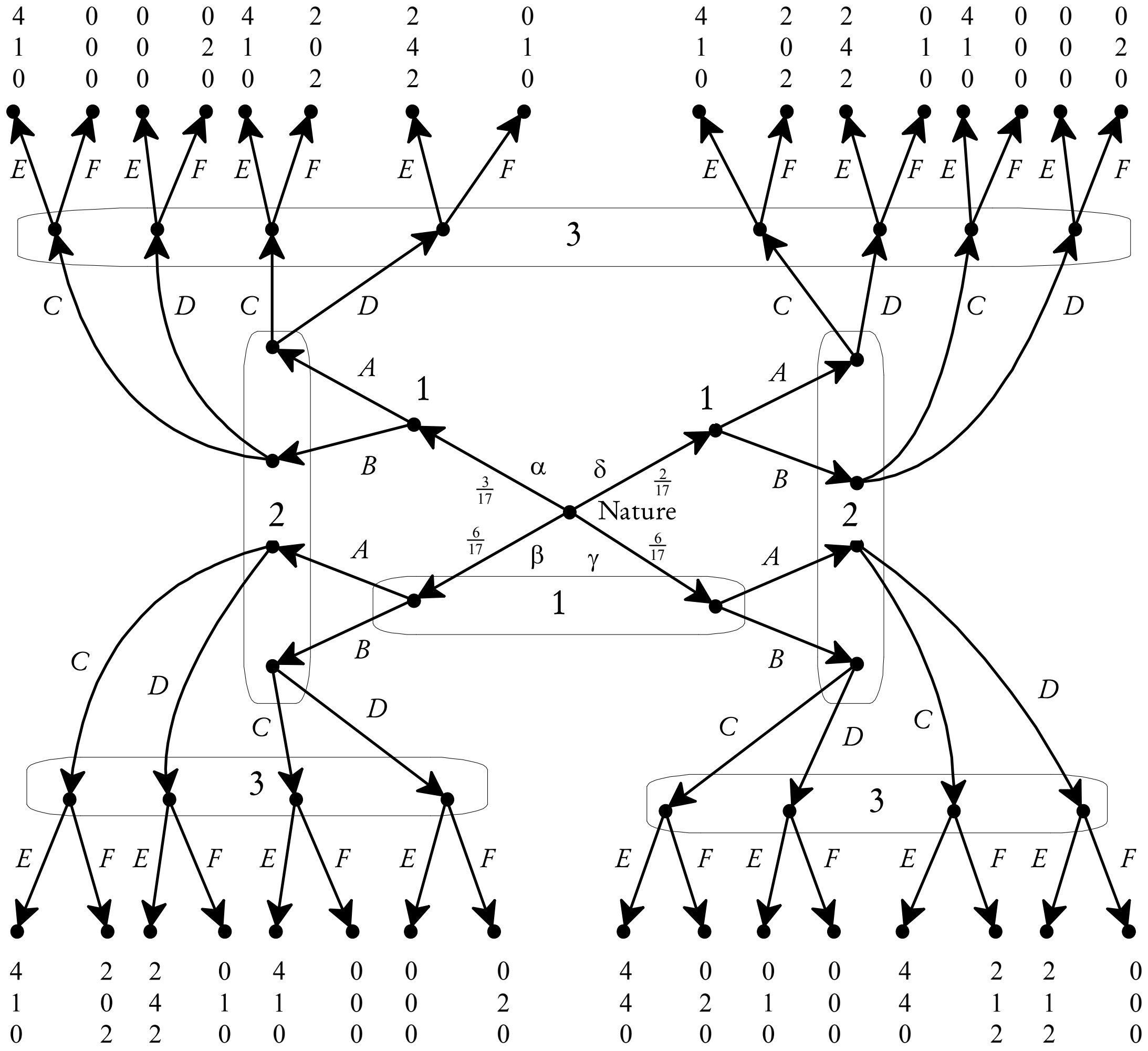

**Figure 13.12**

The task of finding Bayesian Nash equilibria for the game of Figure 13.12 is left as an exercise (see Exercise 13.7).

This is a good time to test your understanding of the concepts introduced in this section, by going through the exercises in Section 13.E.4 of Appendix 13.E at the end of this chapter.

# Appendix 13.E: Exercises

## 13.E.2. Exercises for Section 13.2: One-sided incomplete information.

The answers to the following exercises are in Appendix S at the end of this chapter.

**Exercise 13.1**. Albert and Bill play a Game of Chicken. Each player can choose "swerve" or "don't swerve." Choices are made simultaneously. If a player swerves then he is a chicken. If he does not swerve, then he is a rooster if the other player swerves, but he gets injured if the other player does not. A player can be a "normal" type or a "reckless" type. A "normal" type gets a payoff of 100 from being a rooster, 85 from being a chicken, and zero from being injured. A "reckless" type thinks being a rooster is worth 100, being injured is worth 50, but being a chicken is worth zero! As a matter of fact, both Albert and Bill are normal; it is common knowledge between them that Bill is normal, but Bill thinks that there is a 20% chance that Albert is reckless (and an 80% chance that Albert is normal). Bill knows that Albert knows whether he (= Albert) is reckless or normal. Everything that Bill knows is common knowledge between Albert and Bill.

**(a)** Construct an interactive knowledge-belief structure that represents the situation of incomplete information described above.

**(b)** Apply the Harsanyi transformation to obtain an extensive-form game.

**(c)** Construct the strategic-form associated with the extensive-form game of part (b).

**(d)** Find the pure-strategy Bayesian Nash equilibria and classify them as either pooling or separating.

**Exercise 13.2**. Bill used to be Ann's boyfriend. Today it is Ann's birthday. Bill can either not give a gift to Ann or give her a nicely wrapped gift. If offered a gift, Ann can either accept or reject. This would be a pretty simple situation, if it weren't for the fact that Ann does not know if Bill is still a friend or has become an enemy. If he is a friend, she expects a nice present from him. If Bill is an enemy, she expects to find a humiliating thing in the box (he is known to have given dead frogs to his "not so nice friends" when he was in third grade!).

The preferences are as follows: Bill's favorite outcome (payoff = 1) occurs when he offers a gift and it is accepted (in either case: if he is a friend, he enjoys seeing her unwrap a nice present, and if he is an enemy, he revels in the cruelty of insulting Ann with a humiliating "gift"). Whether he is a friend or an enemy, Bill prefers having not extended a gift (payoff = 0) to enduring the humiliation of a rejected gift (payoff = $-1$). Ann prefers accepting a gift coming from a friend (payoff = 1) to refusing a gift (payoff = 0); the worst outcome for her (payoff = $-1$) is one where she accepts a gift from an enemy. Ann attaches probability $p$ (with $0 < p < 1$) to the event that Bill is a friend (and $1-p$ to Bill being an enemy); however, Ann knows that Bill knows whether he is a friend or an enemy. Everything that Ann knows is common knowledge between Ann and Bill. As a matter of fact, Bill is a friend.

(a) Construct an interactive knowledge-belief structure that represents the situation of incomplete information described above.

(b) Apply the Harsanyi transformation to obtain an extensive-form game.

(c) Construct the strategic-form associated with the extensive-form game.

(d) Find all the pure-strategy Bayesian Nash equilibria and classify them as either pooling or separating.

(e) Suppose that $p = \frac{1}{4}$. Is the outcome associated with a pure-strategy Bayesian Nash equilibrium Pareto efficient?

## 13.E.3. Exercises for Section 13.2: Two-sided incomplete information.

The answers to the following exercises are in Appendix S at the end of this chapter.

**Exercise 13.3**. Consider the situation of two-sided incomplete information illustrated below (where the true state is $\alpha$).

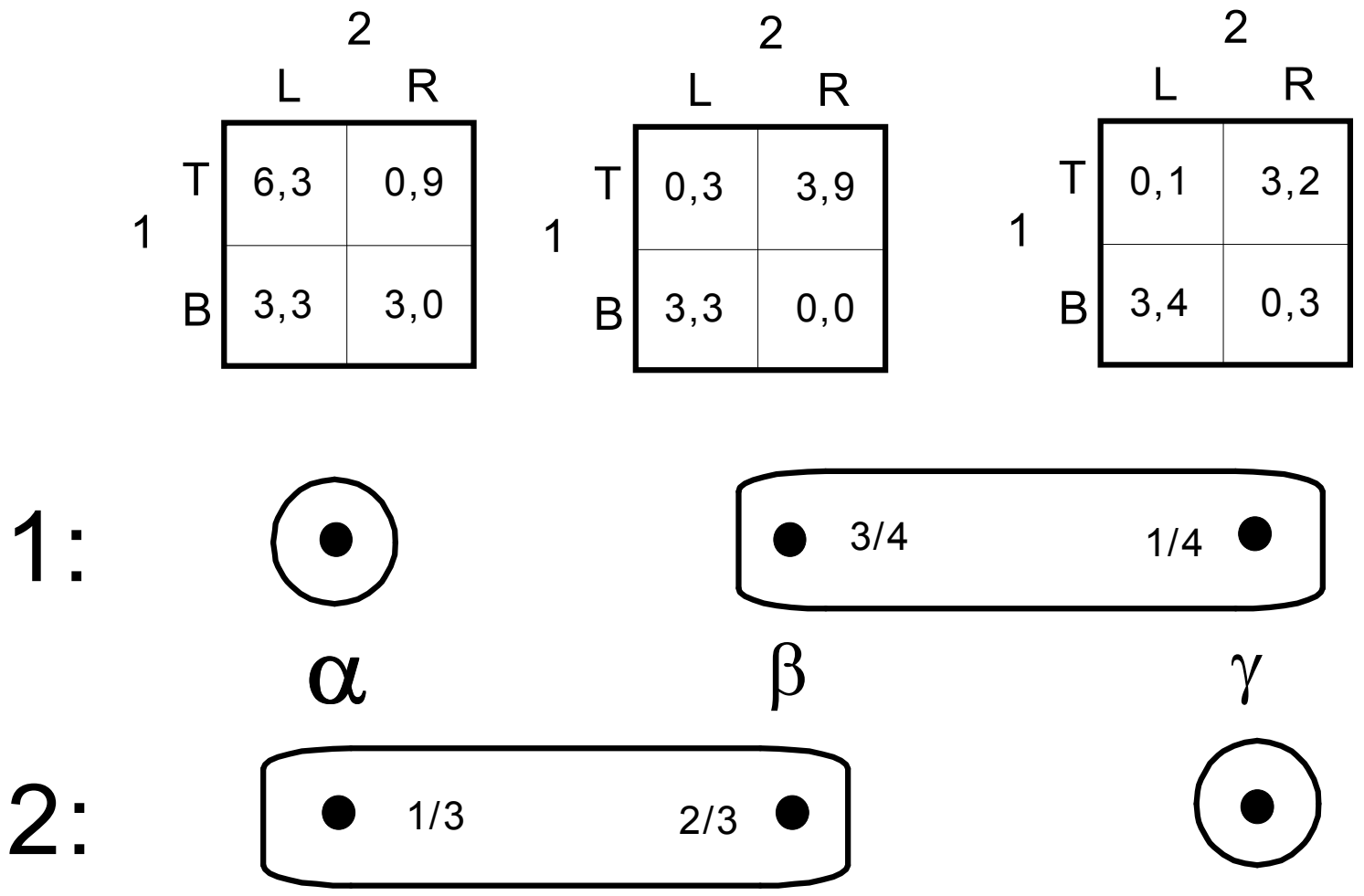

Use the Harsanyi transformation to represent this incomplete-information situation as an extensive-form game. Be explicit about how you calculated the probabilities for Nature.

**Exercise 13.4**. Consider the following congestion situation (a variant of the so called El Farol Bar Problem).[5] Two students at a college are simultaneously deciding between going to a bar and going home. The bar is extremely small and it gets congested when more than one person is there. In principle, there are two types of students. One type, call it the *b* type, prefers going to the bar if he is the only customer there (in which case he gets a utility of 20) but dislikes congestion and gets a utility of $-20$ if he goes to the bar and is not the only customer there; furthermore, for the *b* type, the utility of being at home is 0. The other type of student, call him the *not-b* type, prefers being at home (in which case his utility is 20) to being at the bar; being at the bar alone gives him a utility of 0 but being at the bar with other customers is very stressful and gives him a utility of $-40$. Let $G(b_1, b_2)$ be the game where both players are *b* types, $G(b_1, \text{not-}b_2)$ the game where Player 1 is a *b* type and Player 2 a *not-b* type, etc. Assume that all payoffs are von Neumann-Morgenstern payoffs.

**(a)** Write the four possible strategic-form games.

[5] See https://en.wikipedia.org/wiki/El_Farol_Bar_problem

**(b)** **(b.1)** For each of the games of part (a) find the pure-strategy Nash equilibria.
**(b.2)** For game $G(b_1, b_2)$ find also a mixed-strategy equilibrium where each choice is made with positive probability.

**(c)** Assume now that, as a matter of fact, both players are *b* types. However, it is a situation of incomplete information where it is common knowledge that each player knows his own type but is uncertain about the type of the other player and assigns probability $\frac{1}{5}$ to the other player being the same type as he is and probability $\frac{4}{5}$ to the other player being of the opposite type. Draw an interactive knowledge-belief structure that represents this situation of incomplete information.

**(d)** Use the Harsanyi transformation to represent the above situation of incomplete information as an extensive-form game.

**(e)** For the game of part (d), pick one strategy of Player 1 and explain in words what it means.

**(f)** For the game of part (d), write the corresponding strategic-form game.

**(g)** Find a pure-strategy Bayesian Nash equilibrium of the game of part (d).

**(h)** For the Bayesian Nash equilibrium of part (g), find **(h.1)** where the players are actually going and **(h.2)** the actual payoffs of the players in the game they are actually playing (that is, at the true state). **(h.3)** Do their actual choices yield a Nash equilibrium of the game that they are actually playing?

**(i)** If you didn't know what the true state was but you knew the game of part (d), what probability would you assign to the event that the players would end-up making actual choices that constitute a Nash equilibrium of the true game that they are playing?

**Exercise 13.5**. Consider the situation of two-sided incomplete information illustrated below (where the true state is $\alpha$).

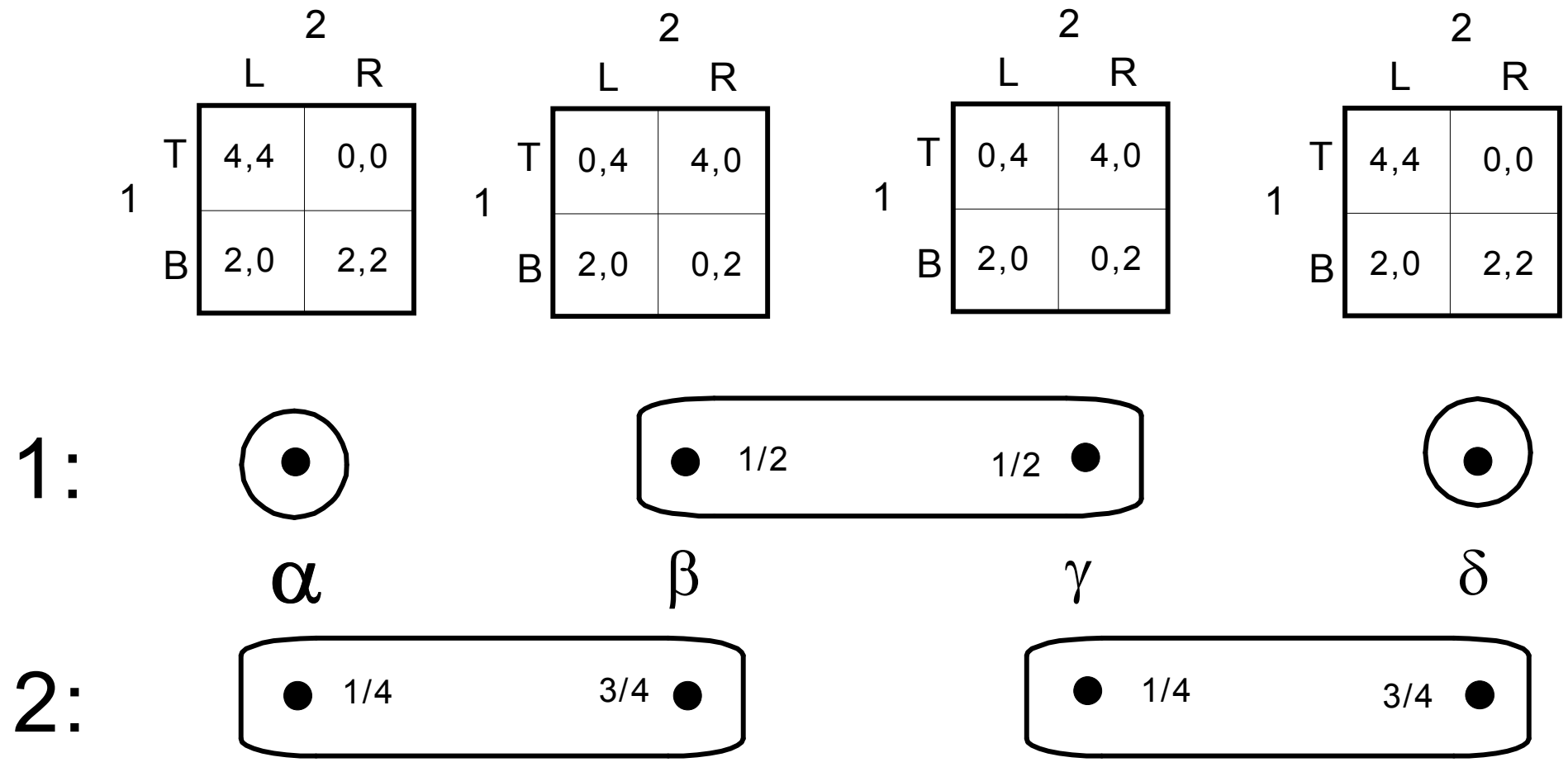

**(a)** Use the Harsanyi transformation to represent this situation as an extensive-form game. Be explicit about how you calculated the probabilities for Nature.

**(b)** Write down all the pure strategies of Player 1 and all the pure strategies of Player 2.

**(c)** Consider the following pure-strategy profile: Player 1 plays $T$ always and Player 2 plays $L$ always. What belief system, paired with this strategy profile, would satisfy Bayesian updating?

## 13.E.4. Exercises for Section 13.4: Multi-sided incomplete information.

The answers to the following exercises are in Appendix S at the end of this chapter.

**Exercise 13.6**. Consider the multi-sided situation of incomplete information shown below (what games $G_1$ and $G_2$ are is irrelevant to the following question). For what values of $p$ can the Harsanyi transformation be carried out?

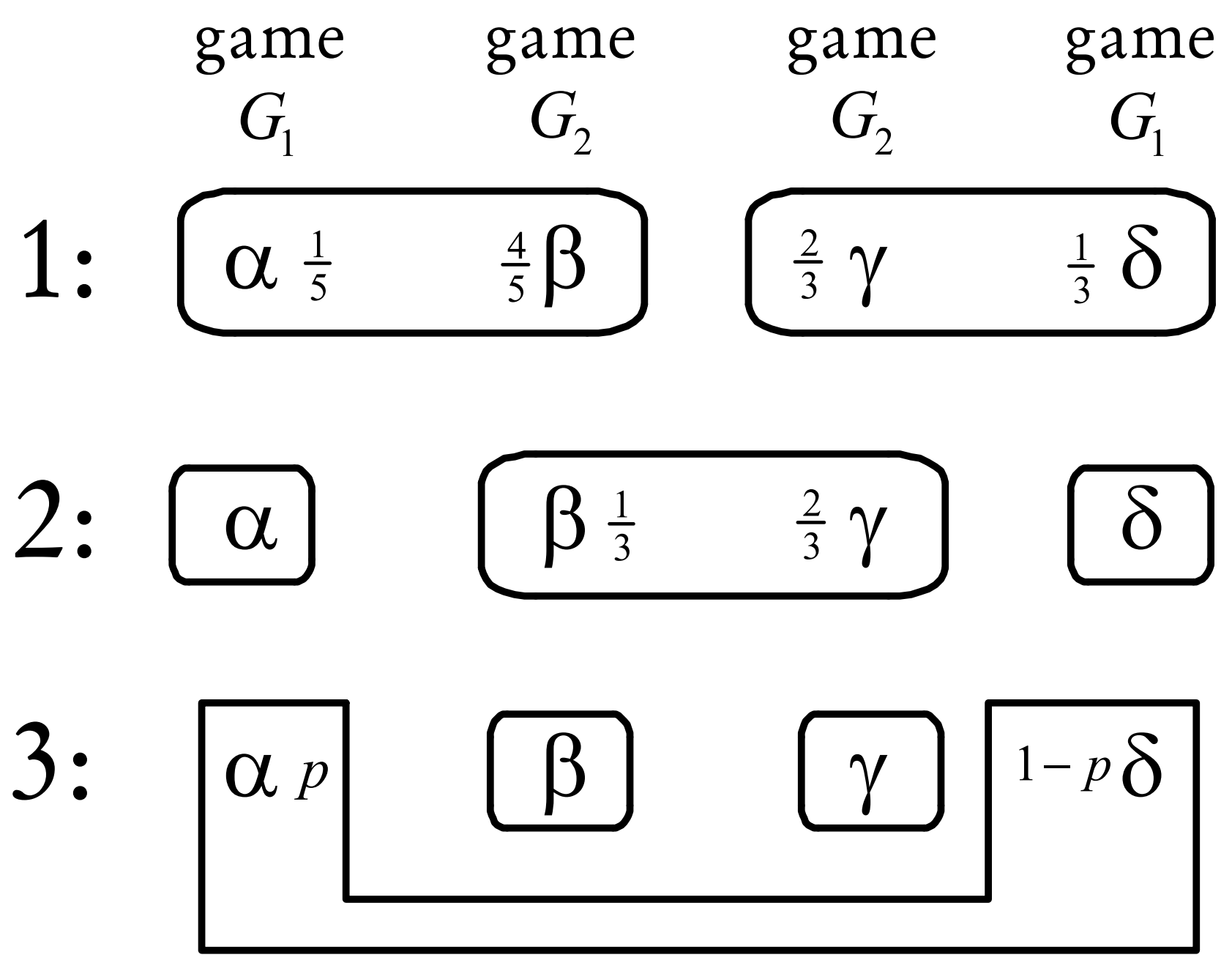

**Exercise 13.7**. Find a pure-strategy Bayesian Nash equilibrium of the game of Figure 13.12.

◊◊◊◊◊◊◊◊◊◊◊◊

**Exercise 13.8: Challenging Question.**

Consider the following situation. It is common knowledge between Players 1 and 2 that tomorrow one of three states will occur: *a*, *b* or *c*. It is also common knowledge between them that (1) if state *a* materializes, then Player 1 will only know that either *a* or *b* occurred and Player 2 will only know that either *a* or *c* occurred, (2) if state *b* materializes, then Player 1 will only know that either *a* or *b* occurred and Player 2 will know that *b* occurred, (3) if state *c* materializes, then Player 1 will know that *c* occurred while Player 2 will only know that either *a* or *c* occurred. Tomorrow they will play the following *simultaneous* game: each will report, confidentially, one of her two possible states of information to a third party (for example, Player 1 can only report $\{a,b\}$ or $\{c\}$). Note that lying is a possibility: for example, if the state is *a* Player 1 can choose to report $\{c\}$. Let $R_1$ be the report of Player 1 and $R_2$ the report of Player 2. The third party, who always knows the true state, then proceeds as follows:

(1) if the reports are compatible, in the sense that $R_1 \cap R_2 \neq \varnothing$, then he gives the players the following sums of money:

If the true state is *a*

| $R_1 \cap R_2 =$ | $a$ | $b$ | $c$ |
|---|---|---|---|
| 1 gets | 5 | 4 | 6 |
| 2 gets | 5 | 6 | 4 |

If the true state is *b*

| $R_1 \cap R_2 =$ | $a$ | $b$ | $c$ |
|---|---|---|---|
| 1 gets | 5 | 4 | 4 |
| 2 gets | 0 | 1 | 1 |

If the true state is *c*

| $R_1 \cap R_2 =$ | $a$ | $b$ | $c$ |
|---|---|---|---|
| 1 gets | 0 | 1 | 1 |
| 2 gets | 5 | 4 | 4 |

(2) if the reports are incompatible ($R_1 \cap R_2 = \varnothing$) then he gives the players the following sums of money:

| The true state is | $a$ | $b$ | $c$ |
|---|---|---|---|
| 1 gets | 5 | 4 | 1 |
| 2 gets | 5 | 1 | 4 |

**(a)** Represent this situation of incomplete information by means of an interactive knowledge structure (for the moment do not worry about beliefs).

**(b)** Apply the Harsanyi transformation to the situation represented in part (a) to obtain an extensive-form frame (again, at this stage, do not worry about probabilities).

**(c)** Suppose first that both players have no idea what the probabilities of the states are and are not willing to form subjective probabilities. It is common knowledge that each player is selfish (i.e. only cares about how much money she herself gets) and greedy (i.e. prefers more money to less) and ranks sets of outcomes according to the worst outcome, in the sense that she is indifferent between sets $X$ and $Y$ if and only if the worst outcome in $X$ is equal to the worst outcome in $Y$ and prefers $X$ to $Y$ if and only if the worst outcome in $X$ is better than the worst outcome in $Y$.
**(c.1)** Write the normal-form (or strategic-form) of the game of Part (a).
**(c.2)** Find all the pure-strategy Nash equilibria of this game.
**(c.3)** Among the Nash equilibria, is there one where each player tells the truth?

**(d)** Suppose now that it is common knowledge between the players that there are objective probabilities for the states as follows: $\begin{pmatrix} a & b & c \\ \frac{2}{5} & \frac{1}{5} & \frac{2}{5} \end{pmatrix}$. This time assume that it is common knowledge that both players are selfish, greedy and risk-neutral. (Thus ignore now the preferences described in part c.)

**(d.1)** Suppose that Player 2 expects Player 1 to report truthfully. Is it rational for Player 2 to also report truthfully?

**(d.2)** Is "always lying" for each player part of a pure-strategy weak sequential equilibrium? Prove your claim.

# Appendix 13.S: Solutions to exercises

**Exercise 13.1.** (a) The structure is as follows (*s* means "swerve", *ds* means "don't swerve"):

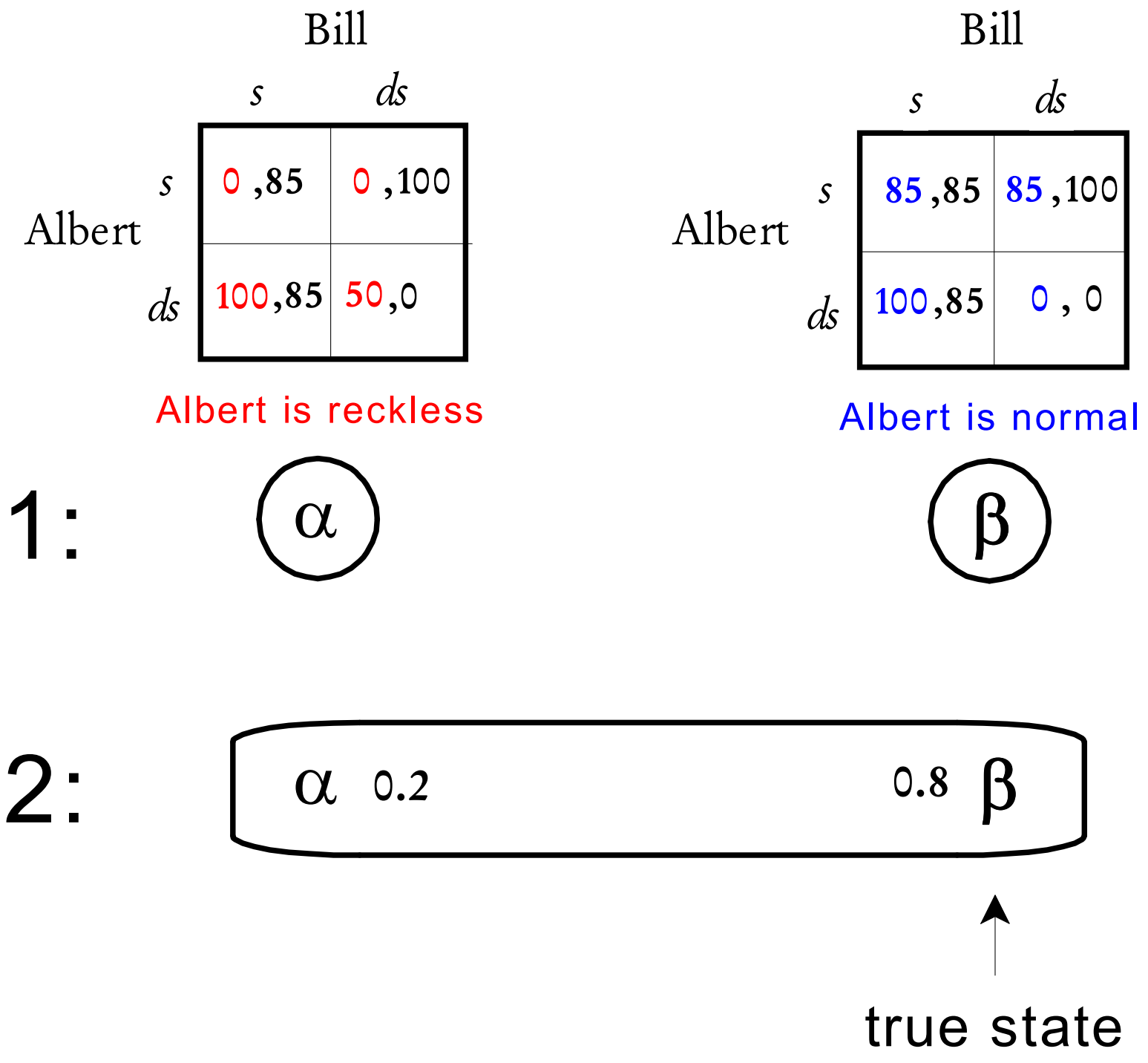

(b) The Harsanyi transformation yields the following game:

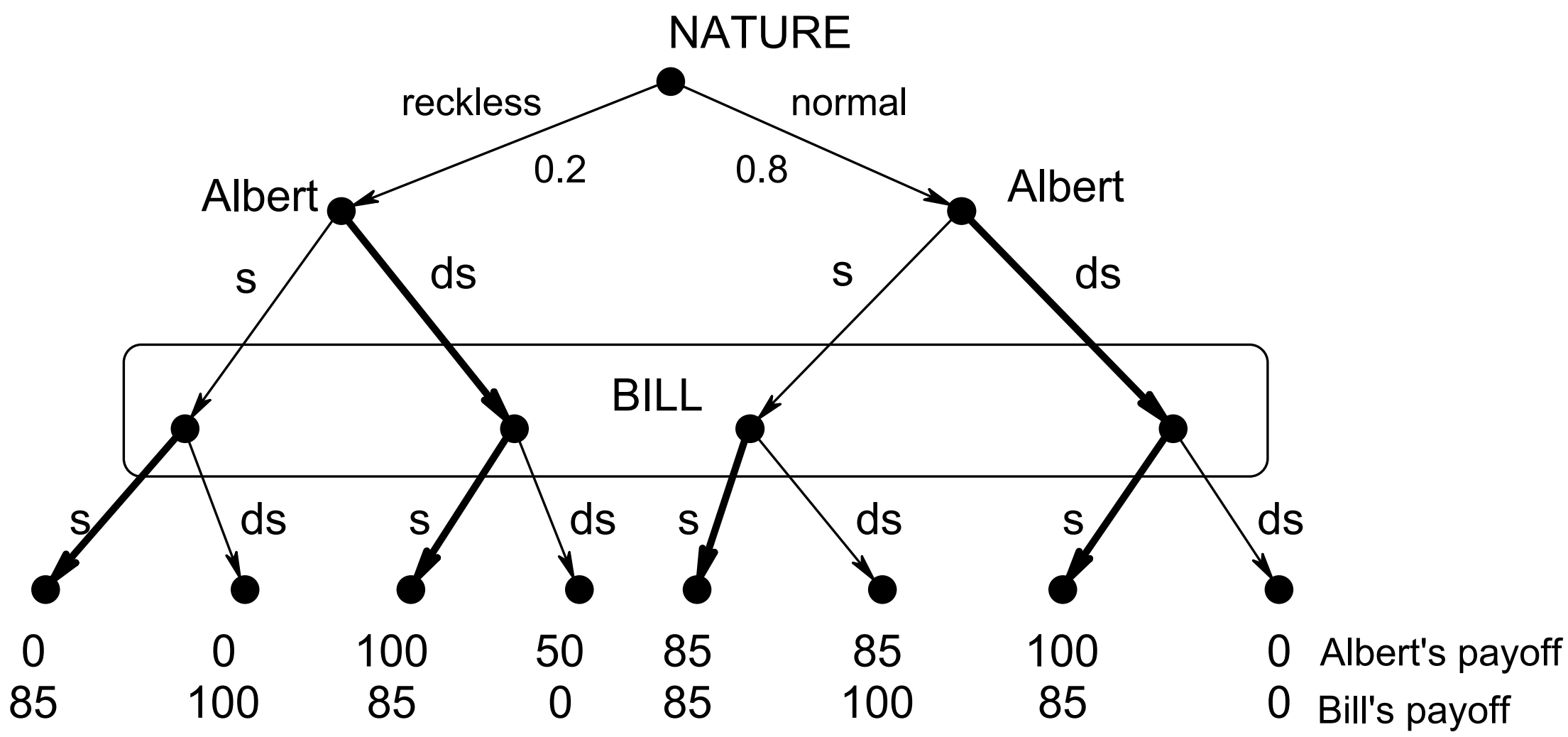

**(c)** The corresponding strategic-form game is as follows:

| Albert \ Bill | s | ds |
|---|---|---|
| s, s | (0.8)85 = **68 , 85** | (0.8)85 = **68 , 100** |
| s, ds | (0.8)100 = **80 , 85** | **0 ,** (0.2)100 = **20** |
| ds, s | (0.8)85 +(0.2)100 = **88 , 85** | (0.8)85 + (0.2)50 = **78 ,** (0.8)100 = **80** |
| ds, ds | **100 , 85** | (0.2)50 = **10 , 0** |

**(d)** There is only one pure-strategy Bayesian Nash equilibrium, namely ((*ds*,*ds*),*s*) and it is a pooling equilibrium.

**Exercise 13.2.** **(a)** The structure is as follows (*g* means "gift", *ng* means "no gift", *a* means "accept" and *r* means "reject"):

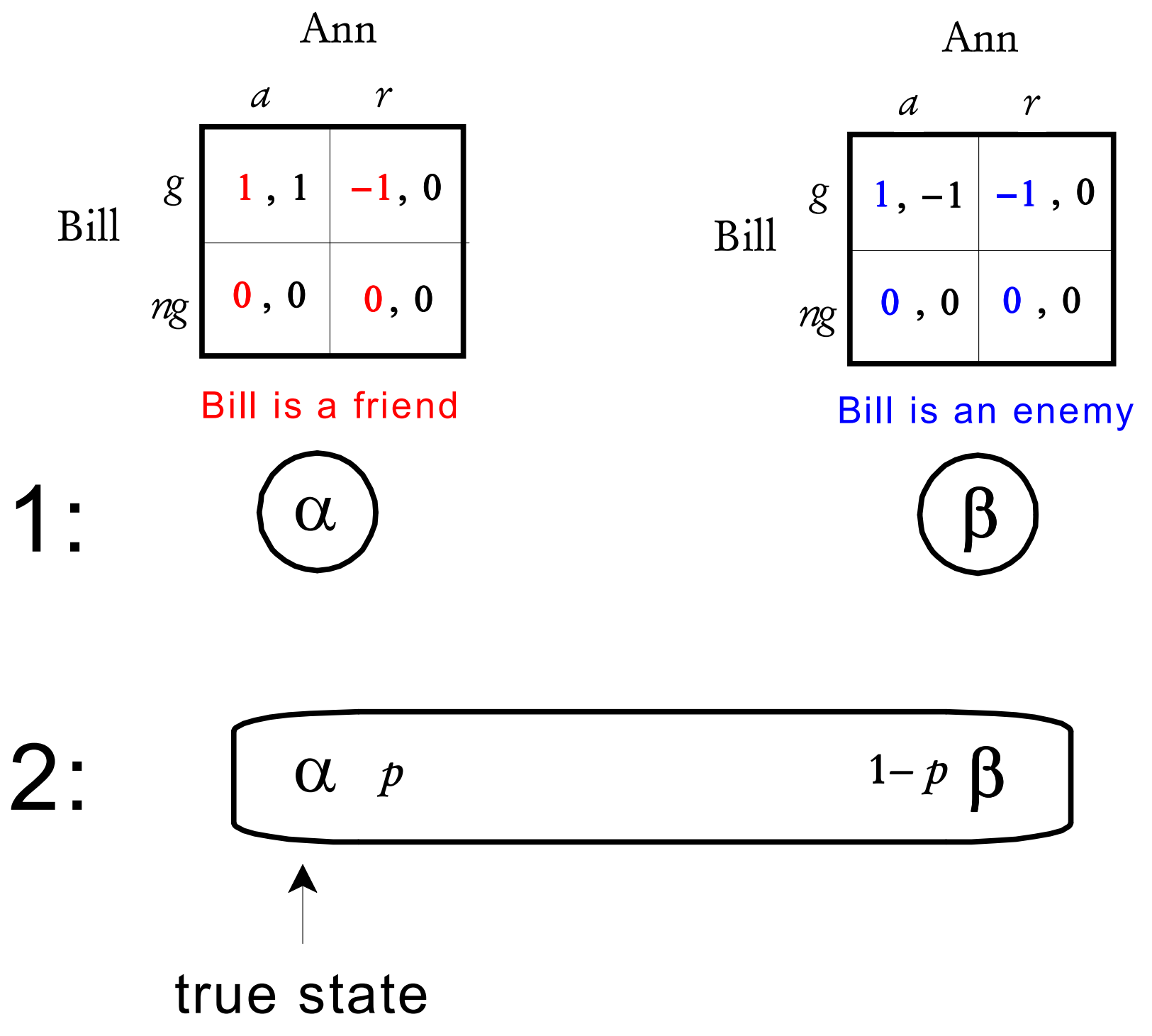

**(b)** The Harsanyi transformation yields the following game:

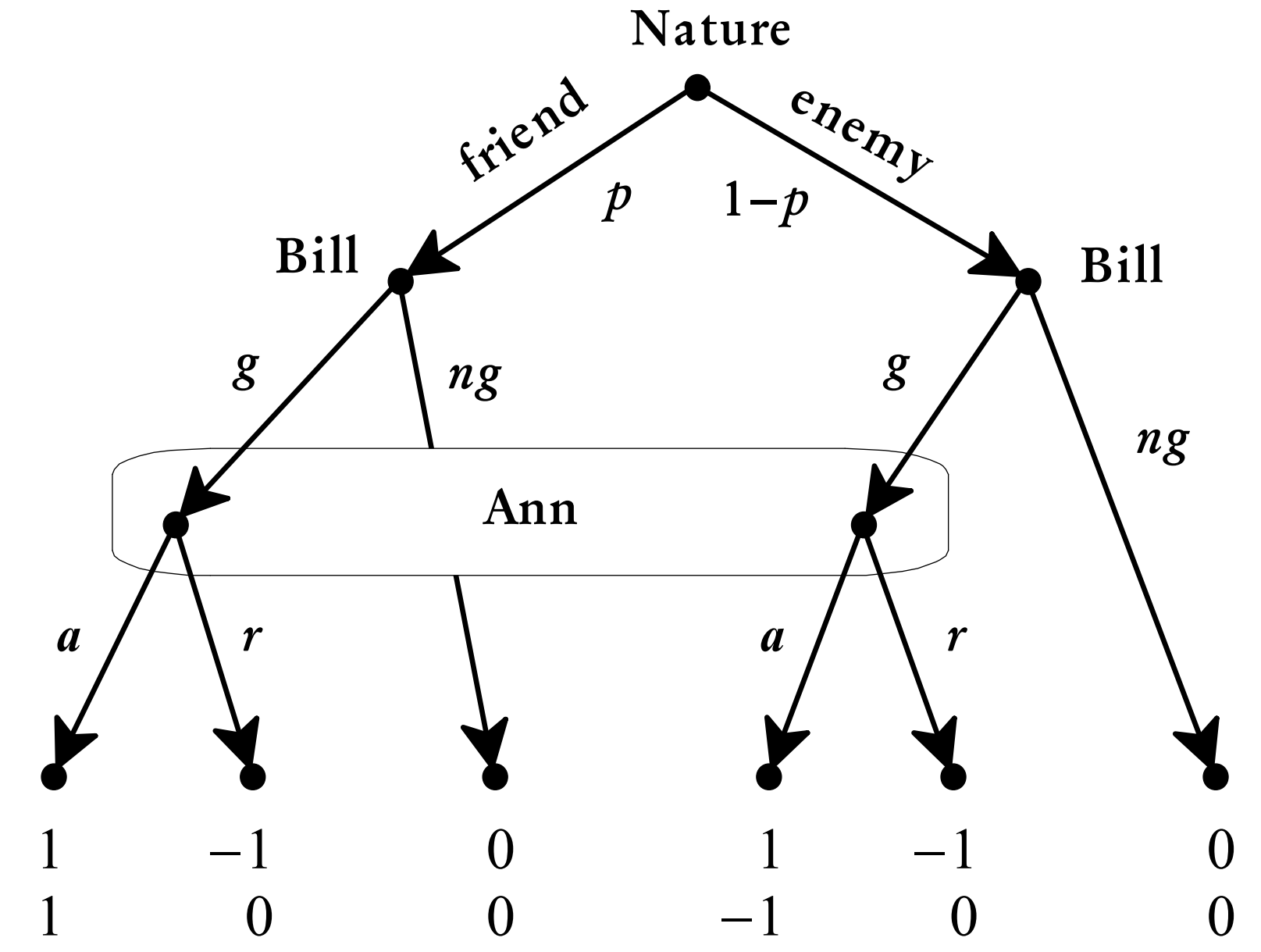

**(c)** The corresponding strategic-form game is as follows:

**Ann**

| Bill | *a* | *r* |
|---|---|---|
| *g g* | 1 , 2*p*−1 | −1 , 0 |
| *g ng* | *p* , *p* | −*p* , 0 |
| *ng g* | 1−*p* , *p*−1 | *p*−1 , 0 |
| *ng ng* | 0 , 0 | 0 , 0 |

**(d)** If $p \geq \frac{1}{2}$ then there are two pure-strategy Bayesian Nash equilibria: ((*g*,*g*), *a*) and ((*ng*,*ng*), *r*). Both of them are pooling equilibria. If p $<$ $\frac{1}{2}$ then there is only one Bayesian -Nash equilibrium: ((*ng*,*ng*), *r*) (a pooling equilibrium).

**(e)** If $p = \frac{1}{4}$ the only Bayesian-Nash equilibrium is ((*ng*,*ng*), *r*) and the outcome is that Bill does not offer a gift to Ann. This outcome is Pareto inefficient in the true game because, given that the true state of affairs is one where Bill is a friend, a Pareto superior outcome would be one where Bill offers a gift and Ann accepts.

**Exercise 13.3. (a) and (b.1):**

$G(b_1, b_2)$:

| | | Player 2 | |
|---|---|---|---|
| | | B | H |
| Player 1 | B | −20, −20 | 20, 0 |
| | H | 0, 20 | 0, 0 |

Nash equilibria: (H,B) and (B,H)

$G(b_1, \text{not-}b_2)$:

| | | Player 2 | |
|---|---|---|---|
| | | B | H |
| Player 1 | B | −20, −40 | 20, 20 |
| | H | 0, 0 | 0, 20 |

Nash equilibrium: (B,H)

$G(\text{not-}b_1, b_2)$:

| | | Player 2 | |
|---|---|---|---|
| | | B | H |
| Player 1 | B | −40, −20 | 0, 0 |
| | H | 20, 20 | 20, 0 |

Nash equilibrium: (H,B)

$G(\text{not-}b_1, \text{not-}b_2)$:

| | | Player 2 | |
|---|---|---|---|
| | | B | H |
| Player 1 | B | −40, −40 | 0, 20 |
| | H | 20, 0 | 20, 20 |

Nash equilibrium: (H,H).

**(b.2)** Let $p$ be the probability that Player 1 plays $B$ and $q$ the probability that Player 2 plays $B$. Then the mixed-strategy equilibrium is given by the solution to $-20q+20(1-q)=0$ and $-20p+20(1-p)=0$, which is $p=q=\frac{1}{2}$.

**(c)** The structure is as follows:

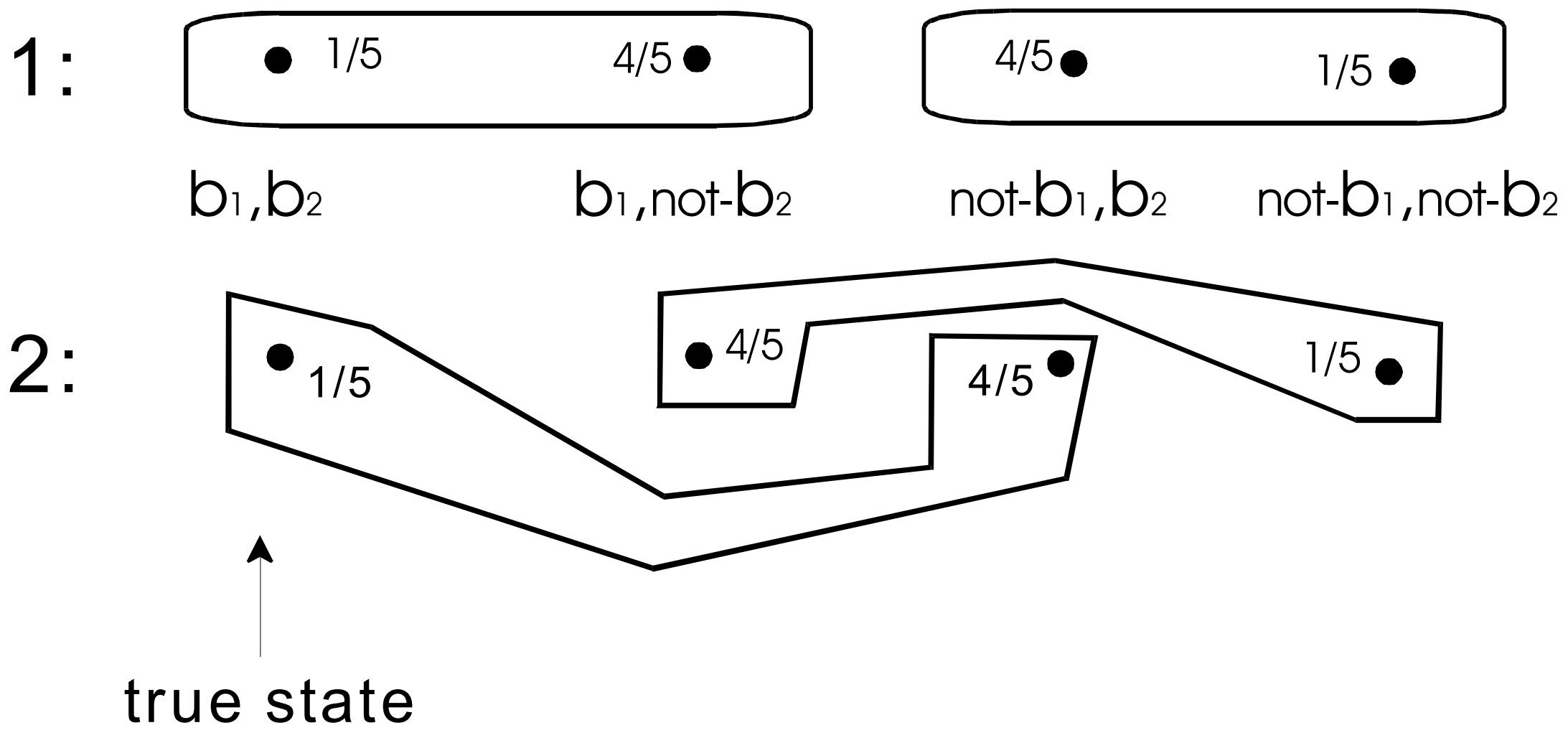

**(d)** First of all, note that the common prior is

$$\begin{pmatrix} b_1, b_2 & b_1, \text{not-}b_2 & \text{not-}b_1, b_2 & \text{not-}b_1, \text{not-}b_2 \\ \frac{1}{10} & \frac{4}{10} & \frac{4}{10} & \frac{1}{10} \end{pmatrix}.$$

Thus the extensive form is as follows:

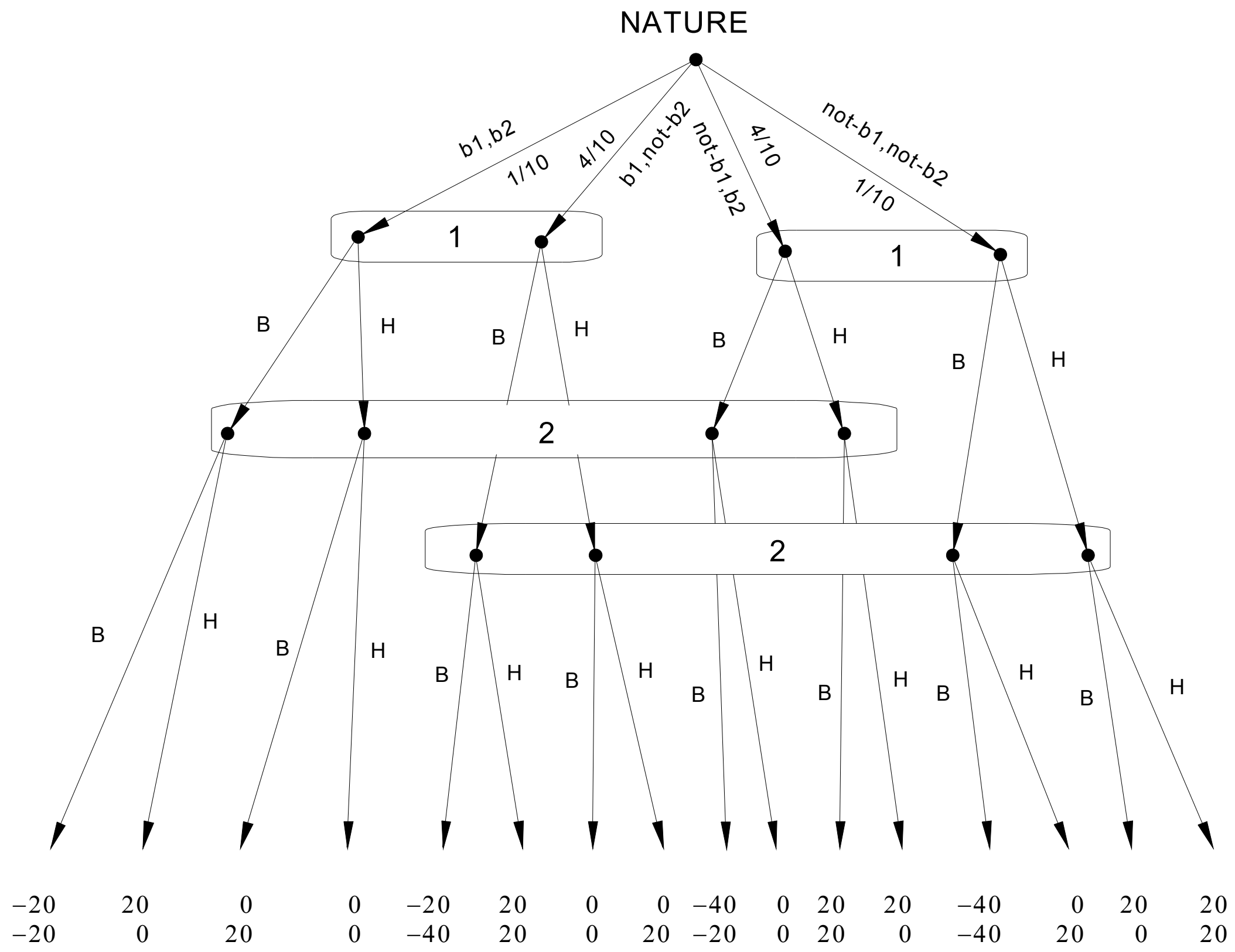

**(e)** One possible strategy for Player 1 is (*B*,*H*) which means "if I am the *b* type then I go to the bar and if I am the *not-b* type then I go home".

**(f)** The strategic form is as follows:

| | | **Player 2** | | | |
|---|---|---|---|---|---|
| | | *B*,*B* | *B*,*H* | *H*,*B* | *H*,*H* |
| **Player 1** | *B*,*B* | −30, −30 | −10 , 0 | −10 , −20 | 10, 10 |
| | *B*,*H* | 0, −10 | 16, 16 | 4 , −16 | 20, 10 |
| | *H*,*B* | −20 , −10 | −16 , 4 | −4 , −4 | 0, 10 |
| | *H*,*H* | 10 , 10 | 10 , 20 | 10 , 0 | 10 , 10 |

**(g)** There is a unique pure-strategy Nash equilibrium, namely $((B,H),(B,H))$ where each player goes to the bar if he is a $b$ type and goes home if he is a *not-b* type. This can be found either the long way, by filling in all the payoffs in the above matrix, or by reasoning as follows. For each player, going home is strictly better than going to the bar if the player is of type *not-b*, no matter what she anticipates the other player doing (in other words, $H$ strictly dominates $B$ at the information set where not-$b_i$ holds). Thus the question is what to do if you are of type $b$. You know that the other player is going home if he is of type *not-b*, thus you only need to consider his choice if he is of type $b$; if his plan is to go home also in that case, then $B$ gives you 20 for sure and $H$ gives you 0 for sure, hence $B$ is better; if his plan is to go to the bar, then $H$ gives you 0 for sure while $B$ gives you the lottery $\begin{pmatrix} -20 & 20 \\ \frac{1}{5} & \frac{4}{5} \end{pmatrix}$, that is, an expected payoff of 12; hence $B$ is better in that case too.

**(h)** At the true state, both players prefer going to the bar, thus **(h.1)** they both end up going to the bar and **(h.2)** get a payoff of −20. **(h.3)** $(B,B)$ is not a Nash equilibrium of the game that they are actually playing (game $G(b_1,b_2)$).

**(i)** If the true game is $G(b_1,b_2)$ they end up playing $(B,B)$ which is not a Nash equilibrium of that game, if the true game is $G(b_1,\text{not-}b_2)$ they end up playing $(B,H)$ which is a Nash equilibrium of that game, if the true game is $G(\text{not-}b_1,b_2)$ they end up playing $(H,B)$ which is a Nash equilibrium of that game and if the true game is $G(\text{not-}b_1,\text{not-}b_2)$ they end up playing $(H,H)$ which is a Nash equilibrium of that game. Since the probability of $G(b_1,b_2)$ is $\frac{1}{10}$, the probability that they end-up playing a Nash equilibrium of the actual game is $\frac{9}{10}$.

**Exercise 13.4.** The game is as follows:

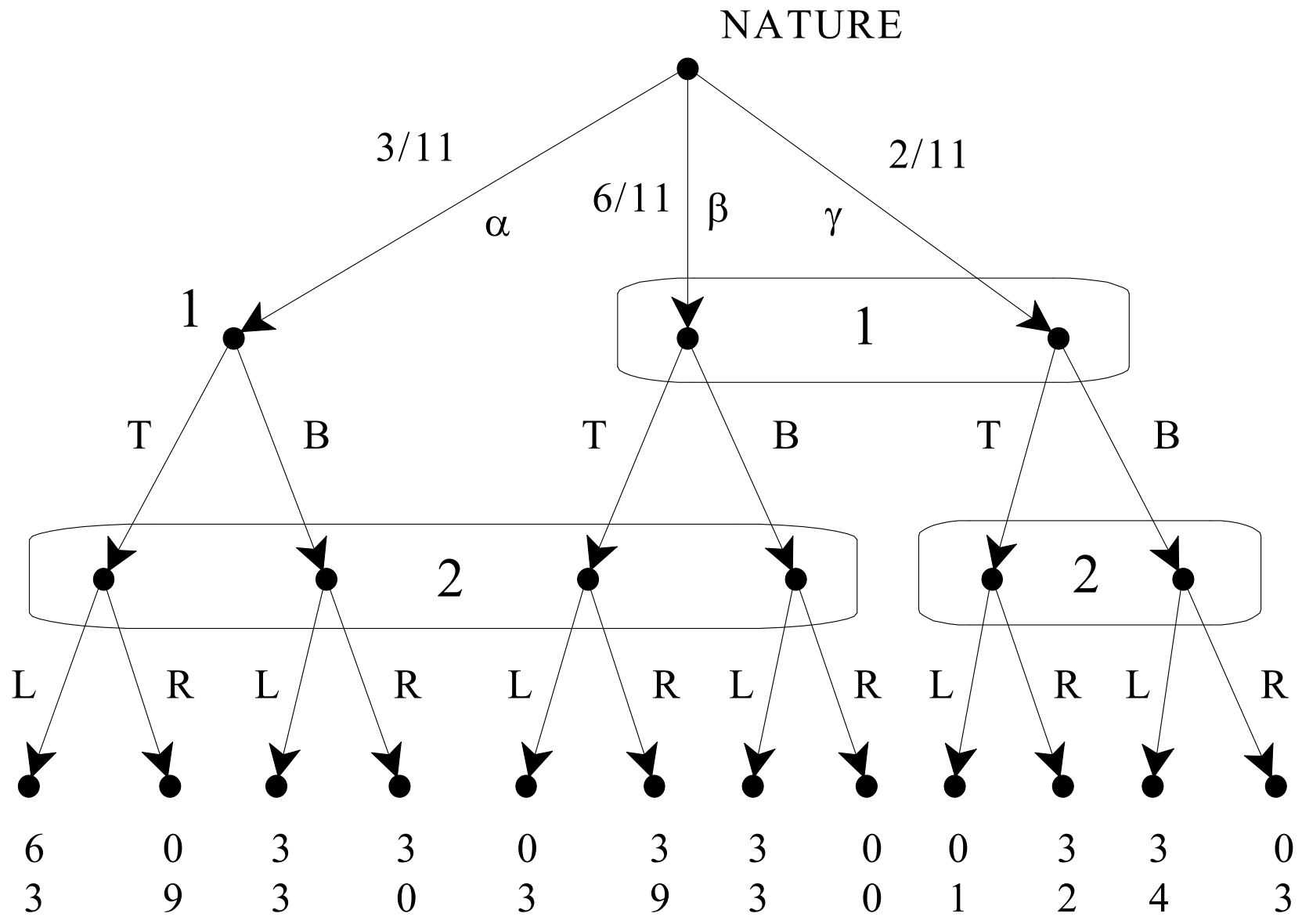

Nature's probabilities are obtained by solving $\dfrac{p_\alpha}{p_\alpha + p_\beta} = \frac{1}{3}$, $\dfrac{p_\beta}{p_\beta + p_\gamma} = \frac{3}{4}$ and $p_\alpha + p_\beta + p_\gamma = 1$.

**Exercise 13.5.** **(a)** The game is as follows:

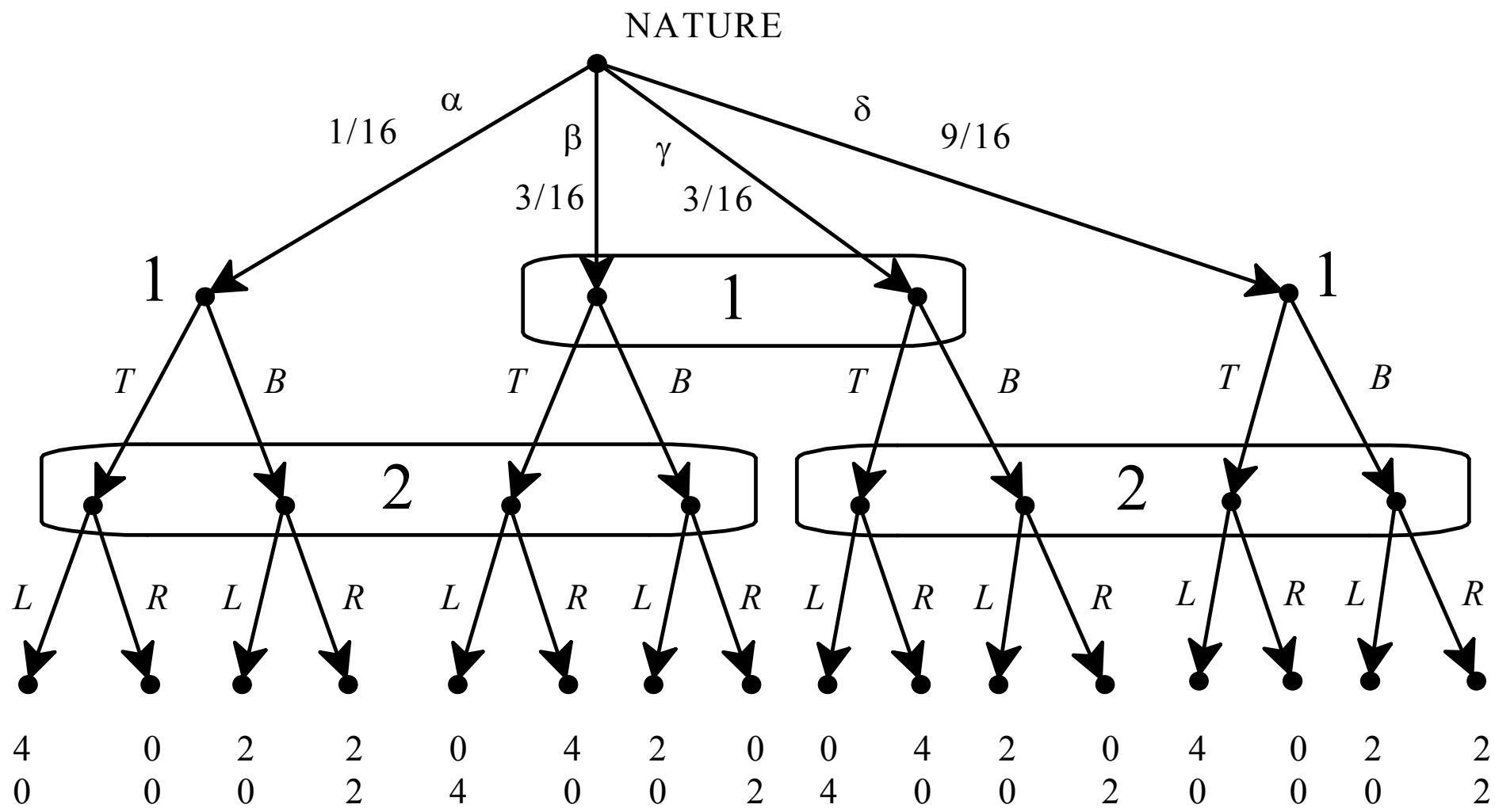

Nature's probabilities are obtained solving the following equations:

$$\frac{p_\alpha}{p_\alpha + p_\beta} = \tfrac{1}{4},\ \frac{p_\beta}{p_\beta + p_\gamma} = \tfrac{1}{2},\ \frac{p_\gamma}{p_\gamma + p_\delta} = \tfrac{3}{4},\ p_\alpha + p_\beta + p_\gamma + p_\delta = 1.$$

**(b)** Player 1 has 8 strategies: *TTT, TTB, TBT, TBB, BTT, BTB, BBT, BBB.* Player 2 has four strategies: *LL, LR, RL, RR.*

**(c)** Player 1's beliefs must be $\frac{1}{2}$ at the left node and $\frac{1}{2}$ at the right node of his middle information set. Player 2's beliefs at her information set on the left must be: $\frac{1}{4}$ at the left-most node and $\frac{3}{4}$ at the third node from the left. The same is true for the other information set of Player 2.

**Exercise 13.6**. The Harsanyi transformation requires that there be a common prior. Thus we need a probability distribution $\nu : \{\alpha, \beta, \gamma, \delta\} \to (0,1)$ such that:

$$(1)\ \frac{\nu(\alpha)}{\nu(\alpha)+\nu(\beta)} = \tfrac{1}{5},\quad (2)\ \frac{\nu(\gamma)}{\nu(\gamma)+\nu(\delta)} = \tfrac{2}{3},\quad (3)\ \frac{\nu(\beta)}{\nu(\beta)+\nu(\gamma)} = \tfrac{1}{3}$$

$$\text{and}\quad (4)\ \frac{\nu(\alpha)}{\nu(\alpha)+\nu(\delta)} = p.$$

From (1) we get $\nu(\beta) = 4\nu(\alpha)$, from (2) we get $\nu(\gamma) = 2\nu(\delta)$ and from (3) we get $\nu(\gamma) = 2\nu(\delta)$; these three equalities, together with $\alpha + \beta + \gamma + \delta = 1$ yield a unique solution, namely $\nu = \begin{pmatrix} \alpha & \beta & \gamma & \delta \\ \frac{1}{17} & \frac{4}{17} & \frac{8}{17} & \frac{4}{17} \end{pmatrix}$. This is a common prior if and only if $p = \dfrac{\nu(\alpha)}{\nu(\alpha)+\nu(\delta)} = \dfrac{\frac{1}{17}}{\frac{1}{17}+\frac{4}{17}} = \frac{1}{5}$. Thus $p = \frac{1}{5}$ is the only value that makes it possible to apply the Harsanyi transformation.

**Exercise 13.7**. The game under consideration is reproduced below, with the equilibrium described below highlighted in red. Given the complexity of the game, it is definitely not a good idea to construct the corresponding strategic form. It is best to think in terms of weak sequential equilibrium. Consider the assessment $(\sigma, \mu)$, highlighted in the figure above, where $\sigma = (BBB, CC, EEE)$ (that is, $\sigma$ is the pure-strategy profile where Player 1 chooses $B$ at each of his three information sets, Player 2 chooses $C$ at each of her two information sets and Player 3 chooses $E$ at each of her three information sets) and $\mu$ is the following system of beliefs: Player 1 attaches probability $\frac{1}{2}$ to each node in his bottom information set; Player 2, at her information set on the left, attaches probability $\frac{1}{3}$ to the second node from the top and probability $\frac{2}{3}$ to the bottom node and, at her information set on the right, attaches probability $\frac{1}{3}$ to the second node from the top and probability $\frac{2}{3}$ to the bottom node; Player 3, at her top information set, attaches probability $\frac{3}{5}$ to the left-most node and probability $\frac{2}{5}$ to the second node from the right and, at her bottom-left information set, attaches probability 1 to the third node from the left and, at her bottom-right information set, attaches probability 1 to the left-most node.

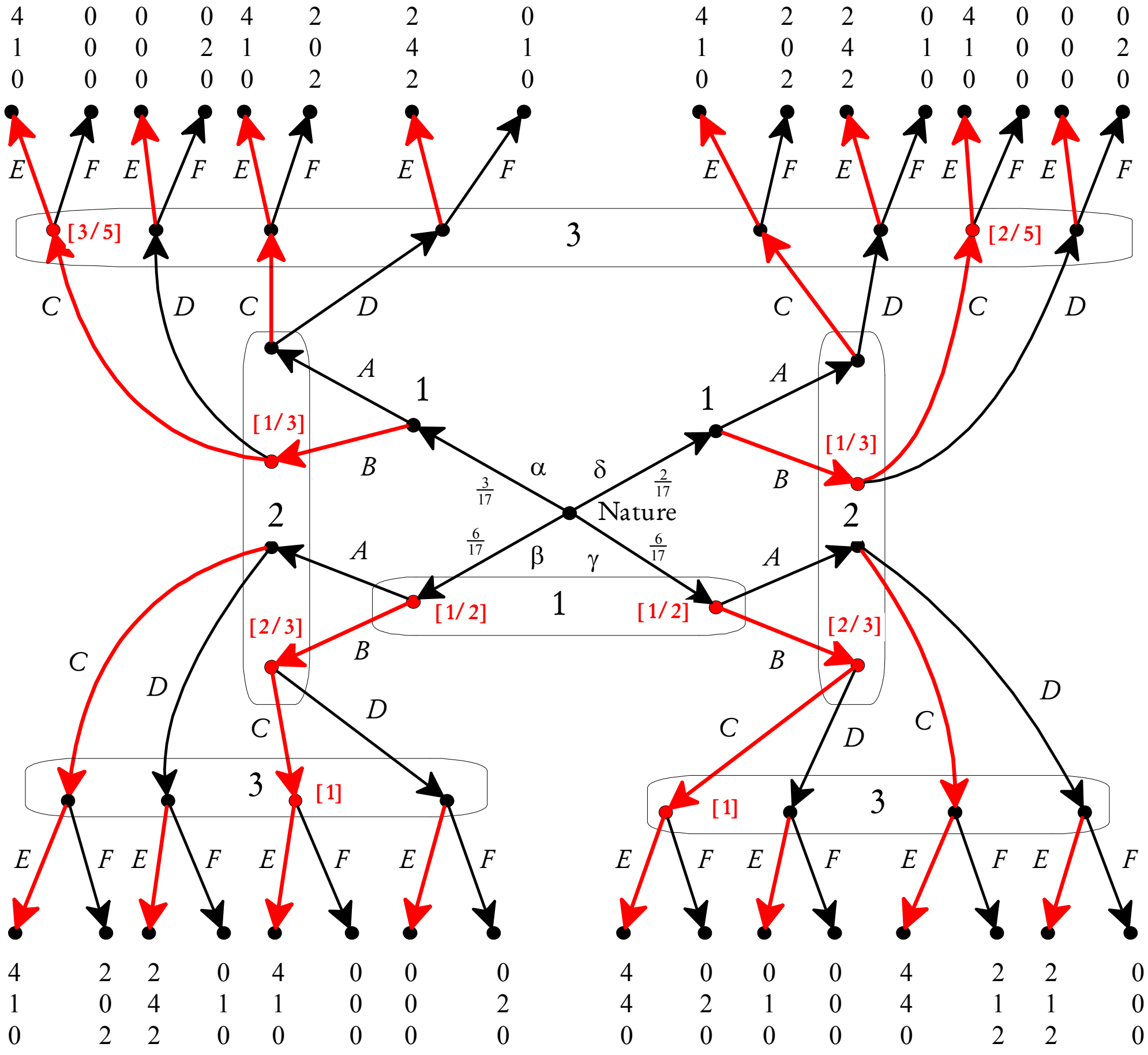

Let us verify that $(\sigma, \mu)$ is a weak sequential equilibrium. The beliefs described above are obtained from $\sigma$ using Bayesian updating. Thus we only need to check sequential rationality.

For Player 1: (1) at the top-left node, both $A$ and $B$ give the same payoff (namely, 4), thus $B$ is sequentially rational; (2) the same is true at the top-right node; (3) at the bottom information set both $A$ and $B$ give an expected payoff of $\frac{1}{2}(4)+\frac{1}{2}(4)=4$, thus $B$ is sequentially rational.

For Player 2: (1) at the left information set $C$ gives an expected payoff of $\frac{1}{3}(1)+\frac{2}{3}(1)=1$ while $D$ gives $\frac{1}{3}(0)+\frac{2}{3}(0)=0$, thus $C$ is sequentially rational; (2) at the information on the right $C$ gives an expected payoff of $\frac{1}{3}(1)+\frac{2}{3}(4)=3$ while $D$ gives an expected payoff of $\frac{1}{3}(0)+\frac{2}{3}(1)=\frac{2}{3}$, thus $C$ is sequentially rational.

For Player 3: (1) at the top information set both $E$ and $F$ give an expected payoff of $\frac{3}{5}(0)+\frac{2}{5}(0)=0$, thus $E$ is sequentially rational; (2) at the bottom-left information set both $E$ and $F$ give a payoff of 0, thus $E$ is sequentially rational; (3) the same is true at the bottom-right information set.

**Exercise 13.8** (Challenging Question). **(a)** The structure is as follows:

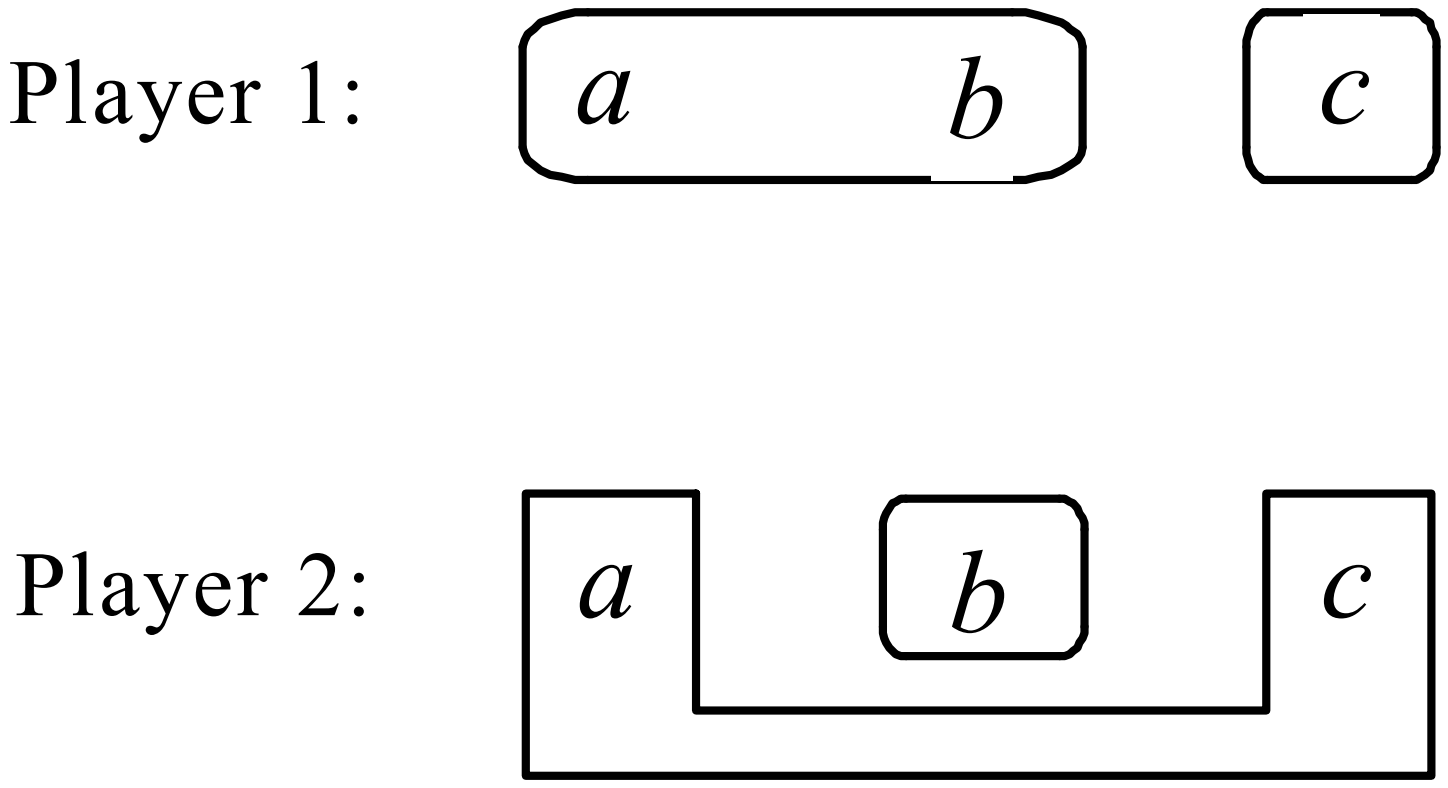

**(b)** The game is as follows:

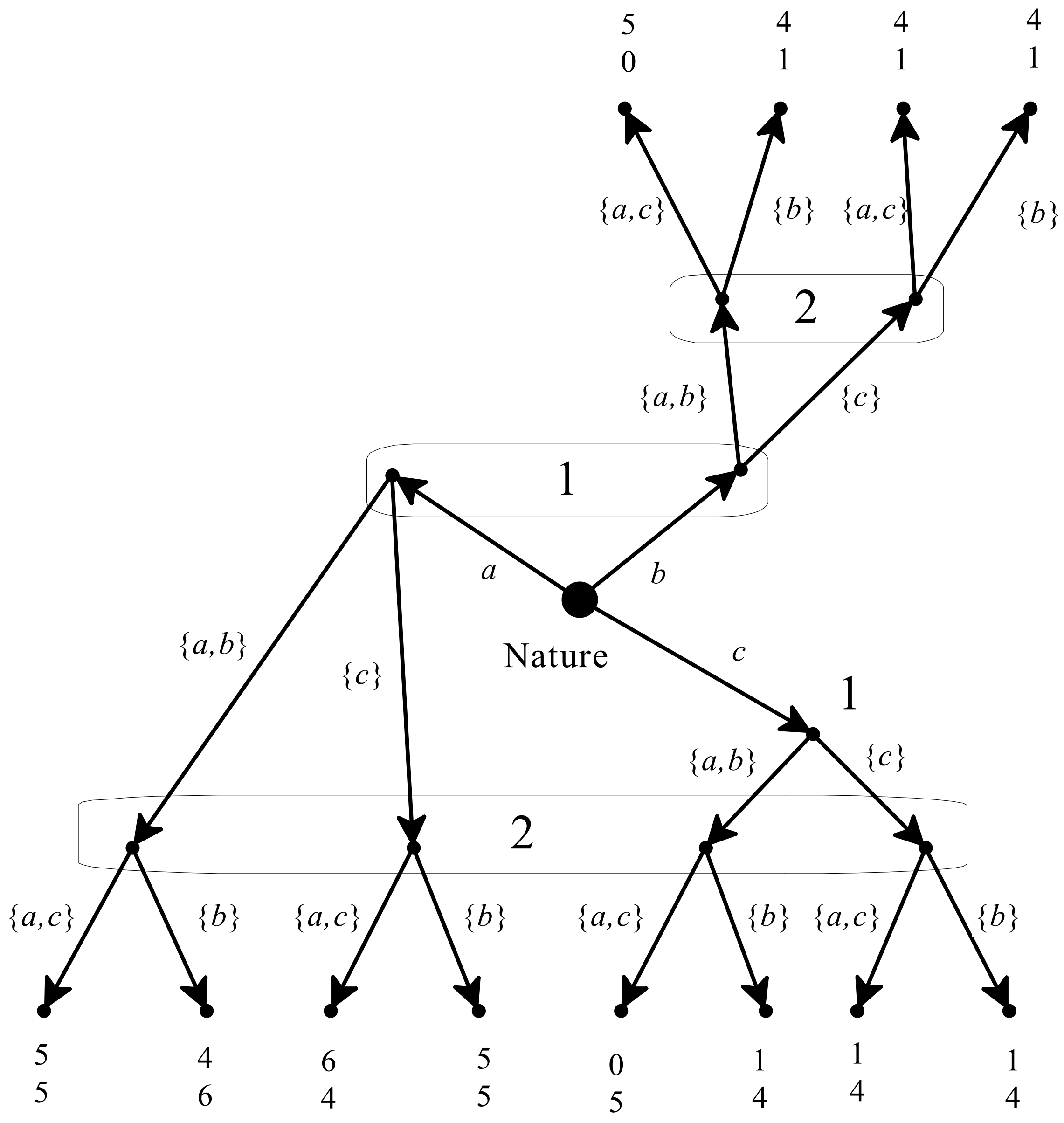

**(c)** **(c.1)** The strategic form is as follows (where the strategy $(x, y)$ for Player 1 means $x$ if $\{a,b\}$ and $y$ if $\{c\}$ and the strategy $(z, w)$ for Player 2 means $z$ if $\{a,c\}$ and $w$ if $\{b\}$). Inside each cell the sets of outcomes are given:

| | | | Player | 2 | |
|---|---|---|---|---|---|
| | | $(\{a,c\},\{a,c\})$ | $(\{a,c\},\{b\})$ | $(\{b\},\{a,c\})$ | $(\{b\},\{b\})$ |
| | $(\{a,b\},\{a,b\})$ | $\{(5,5),(5,0),(0,5)\}$ | $\{(5,5),(4,1),(0,5)\}$ | $\{(4,6),(5,0),(1,4)\}$ | $\{(4,6),(4,1),(1,4)\}$ |
| Pl | $(\{a,b\},\{c\})$ | $\{(5,5),(5,0),(1,4)\}$ | $\{(5,5),(4,1),(1,4)\}$ | $\{(4,6),(5,0),(1,4)\}$ | $\{(4,6),(4,1),(1,4)\}$ |
| 1 | $(\{c\},\{a,b\})$ | $\{(6,4),(4,1),(0,5)\}$ | $\{(6,4),(4,1),(0,5)\}$ | $\{(5,5),(4,1),(1,4)\}$ | $\{(5,5),(4,1),(1,4)\}$ |
| | $(\{c\},\{c\})$ | $\{(6,4),(4,1),(1,4)\}$ | $\{(6,4),(4,1),(1,4)\}$ | $\{(5,5),(4,1),(1,4)\}$ | $\{(5,5),(4,1),(1,4)\}$ |

Taking as payoffs the smallest sum of money in each cell (for the corresponding player) the game can be written as follows:

| | | | Player | 2 | |
|---|---|---|---|---|---|
| | | $(\{a,c\},\{a,c\})$ | $(\{a,c\},\{b\})$ | $(\{b\},\{a,c\})$ | $(\{b\},\{b\})$ |
| | $(\{a,b\},\{a,b\})$ | 0 , 0 | 0 , 1 | 1 , 0 | 1 , 1 |
| Pl | $(\{a,b\},\{c\})$ | 1 , 0 | 1 , 1 | 1 , 0 | 1 , 1 |
| 1 | $(\{c\},\{a,b\})$ | 0 , 1 | 0 , 1 | 1 , 1 | 1 , 1 |
| | $(\{c\},\{c\})$ | 1 , 1 | 1 , 1 | 1 , 1 | 1 , 1 |

**(c.2)** There are 9 Nash equilibria which are highlighted in red.

**(c.3)** Truth telling is represented by the strategy profile $((\{a,b\},\{c\}),(\{a,c\},\{b\}))$ and it is one of the Nash equilibria.

**(d)** **(d.1)** No. If the state is $b$ then it is a good idea for Player 2 to report truthfully because $\{a,c\}$ yields her 0 while $\{b\}$ yields her 1. But if the state is either $a$ or $c$ then, by Bayesian updating, Player 2 must assign probability $\frac{1}{2}$ to the left-most node and probability $\frac{1}{2}$ to the right-most node of her bottom information set; thus her expected payoff from reporting $\{a,c\}$ is $\frac{1}{2}(5)+\frac{1}{2}(4)=4.5$ while the expected payoff from reporting $\{b\}$ is $\frac{1}{2}(6)+\frac{1}{2}(4)=5$.

**(d.2)** Yes. "Always lie" corresponds to the strategy profile $((\{c\},\{a,b\}),(\{b\},\{a,c\}))$. By Bayesian updating the corresponding beliefs must be: for Player 1 $\left(\frac{2}{3},\frac{1}{3}\right)$ and for Player 2 $\left(0,\frac{1}{2},\frac{1}{2},0\right)$ at the bottom information set and (0,1) at the top information set. Sequential rationality is then satisfied at every information set: **for Player 1** at the top information set $\{c\}$ gives an expected payoff of $\frac{2}{3}(5)+\frac{1}{3}(4)=\frac{14}{3}$ while $\{a,b\}$ gives $\frac{2}{3}(4)+\frac{1}{3}(5)=\frac{13}{3}$ and at the singleton node on the right $\{a,b\}$ gives 1 and so does $\{c\}$; **for Player 2** at the bottom information set $\{b\}$ gives an expected payoff of $\frac{1}{2}(5)+\frac{1}{2}(4)=4.5$ and $\{a,c\}$ gives $\frac{1}{2}(4)+\frac{1}{2}(5)=4.5$ and at the top information set both $\{a,c\}$ and $\{b\}$ give 1.

Chapter

14

# Incomplete Information: Dynamic Games

## 14.1 One-sided incomplete information

At the conceptual level, situations of incomplete information involving dynamic (or extensive-form) games are essentially the same as situations involving static games: the only difference in the representation is that one would associate with every state a dynamic game instead of a static game.[6] As in the case of static games, we will distinguish between one-sided and multi-sided incomplete information. In this section we will go through two examples of the former, while the latter will be discussed in Section 14.2.

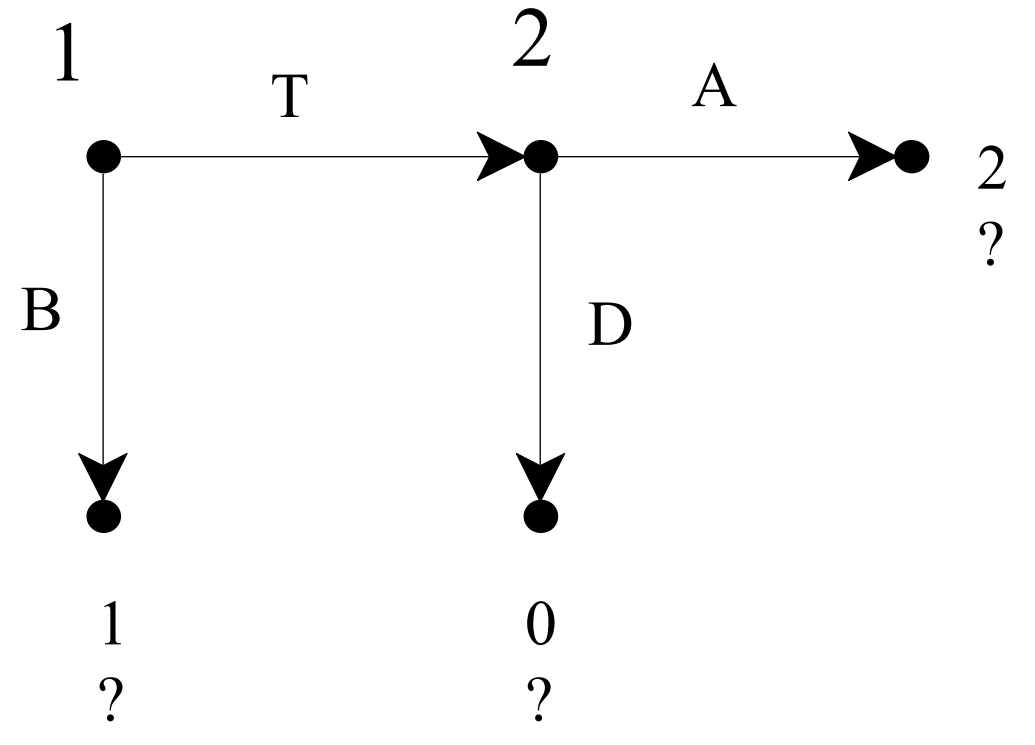

**Figure 14.1**

[6] The reader might have noticed that in Exercise 13.2 we "sneaked in" a dynamic game (where Bill moved first and decided whether or not to offer a gift and Ann – if offered a gift – decided whether or not to accept it). However, the game could also have been written as a simultaneous game, where Ann decided whether or not to accept before knowing whether Bill would offer a gift (without knowing Ann's decision).

Consider the perfect-information frame shown in Figure 14.1, where the question marks in place of Player 2's payoffs indicate that her payoffs are not common knowledge between the two players. Suppose that, on the other hand, Player 1's payoffs *are* common knowledge. Furthermore, Player 1 is uncertain between the two possibilities shown in Figure 14.2 and assigns probability $\frac{1}{3}$ to the one on the left and probability $\frac{2}{3}$ to the one on the right. The thick edges represent the backward-induction solutions of the two games.

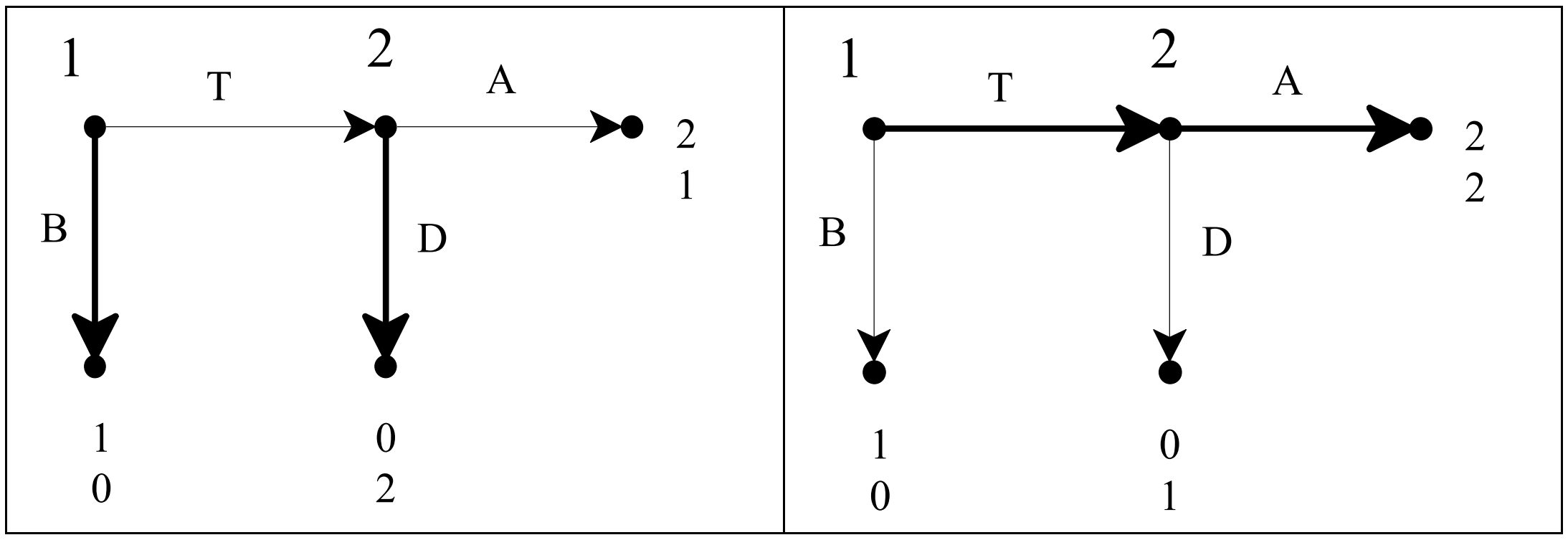

**Figure 14.2**

Thus if Player 1 knew that he was playing the game on the left, then he would choose *B*, while if he knew that he was playing the game on the right then he would choose *T*. If we assume that Player 1's beliefs are common knowledge between the players, then we have the situation of one-sided incomplete information shown in Figure 14.3. Assume that state $\alpha$ is the true state.

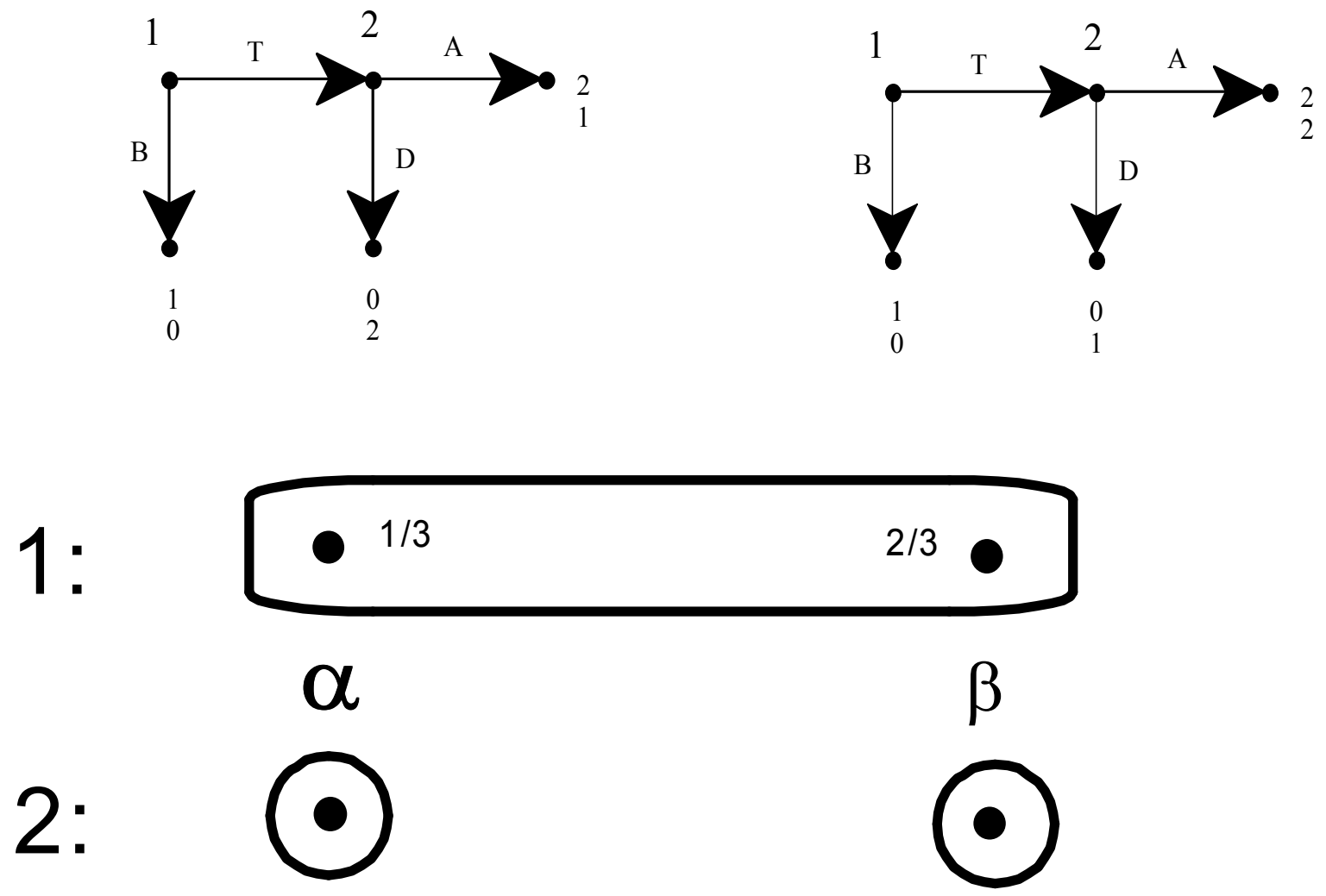

**Figure 14.3**

Using the Harsanyi transformation we can convert the situation illustrated in Figure 14.3 into the extensive-form game shown in Figure 14.4.

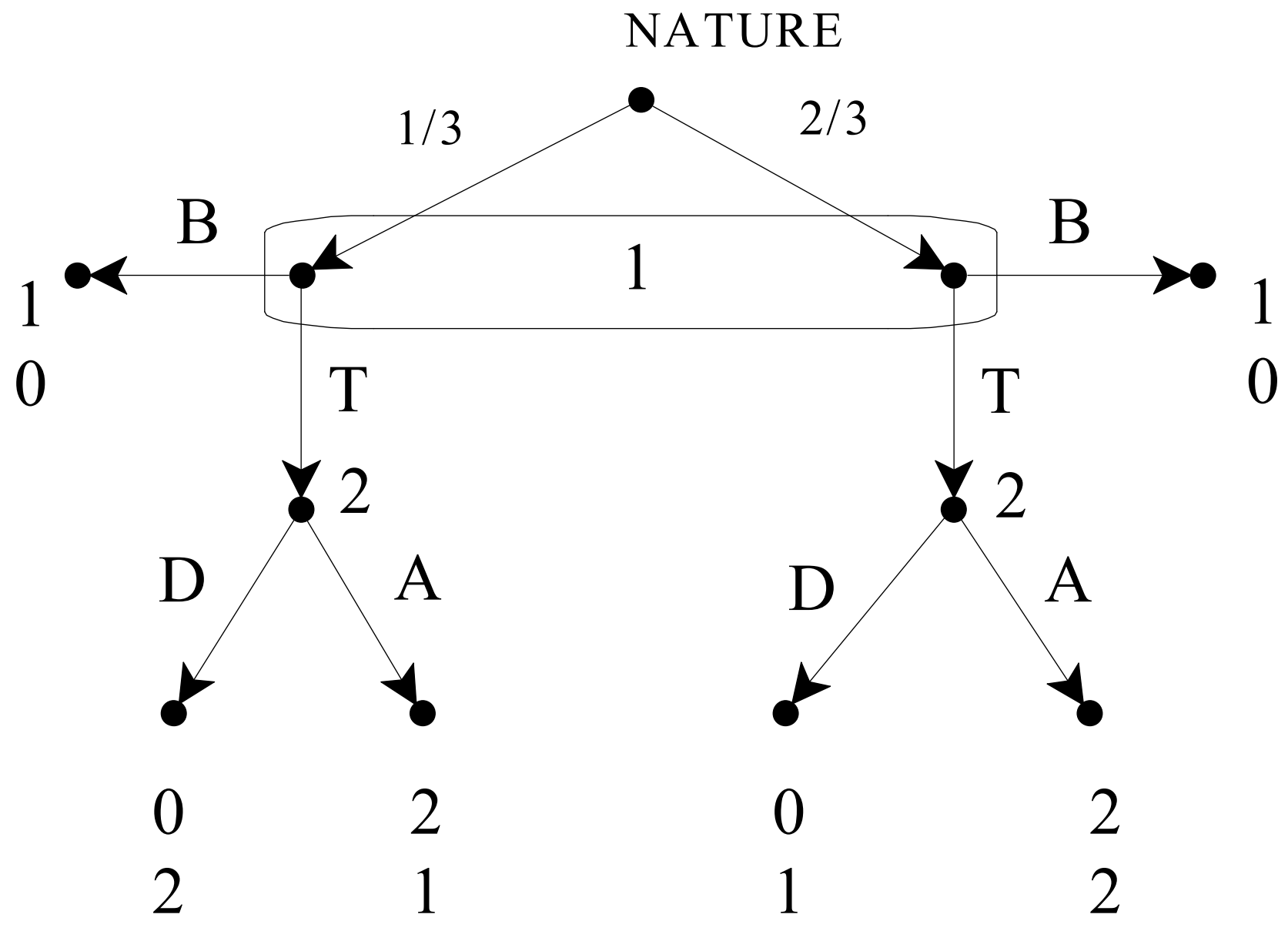

**Figure 14.4**

The solution concept used for the case of static games was Bayesian Nash equilibrium (that is, Nash equilibrium). In the case of dynamic games this is no longer an appropriate solution concept, since – as we know from Chapter 3 – it allows a player to "choose" a strictly dominated action at an unreached information set. To see this, consider the strategic-form game associated with the game of Figure 14.4, shown in Figure 14.5.

| | | Player 2 | | | |
|---|---|---|---|---|---|
| | | *DD* | *DA* | *AD* | *AA* |
| Player 1 | *B* | 1 **0** | 1 **0** | 1 **0** | 1 **0** |
| | *T* | 0 **4/3** | 4/3 **2** | 2/3 **1** | 2 **5/3** |

**Figure 14.5**

The Nash equilibria of the strategic-form game shown in Figure 14.5 (and thus the Bayesian Nash equilibria of the game shown in Figure 14.4) are: (*B* , *DD*), (*B, AD*) and (*T, DA*). Of these only (*T,DA*) is a subgame-perfect equilibrium.

Thus from now on we shall use either the notion of subgame-perfect equilibrium or the notion of weak sequential equilibrium (in the game of Figure 14.4 the two notions coincide). Thus we will take (*T,DA*) to be the solution of the game of Figure 14.4 Since we postulated that the true game was the one associated with state $\alpha$, as a matter of fact the players end up playing (*T,D*) which is not the backward induction solution of the true game. This is not surprising in light of what was already pointed out in Chapter 13.

Next we consider a more complex example, built on Selten's Chain-Store Game analyzed in Chapter 2 (Section 2.4). Recall that the story is as follows (we will consider the special case where $m$, the number of towns and thus of potential entrants, is 2). A chain store is a monopolist in an industry. It owns stores in two different towns. In each town the chain store makes \$5 million if left to enjoy its privileged position undisturbed. In each town there is a businesswoman who could enter the industry in that town, but earns \$1 million if she chooses not to enter; if she decides to enter, then the monopolist can either fight the entrant, leading to zero profits for both the chain store and the entrant in that town, or it can accommodate entry and share the market with the entrant, in which case both players make \$1.5 million in that town. Thus, in each town the interaction between the incumbent monopolist and the potential entrant is as illustrated in Figure 14.6 below.

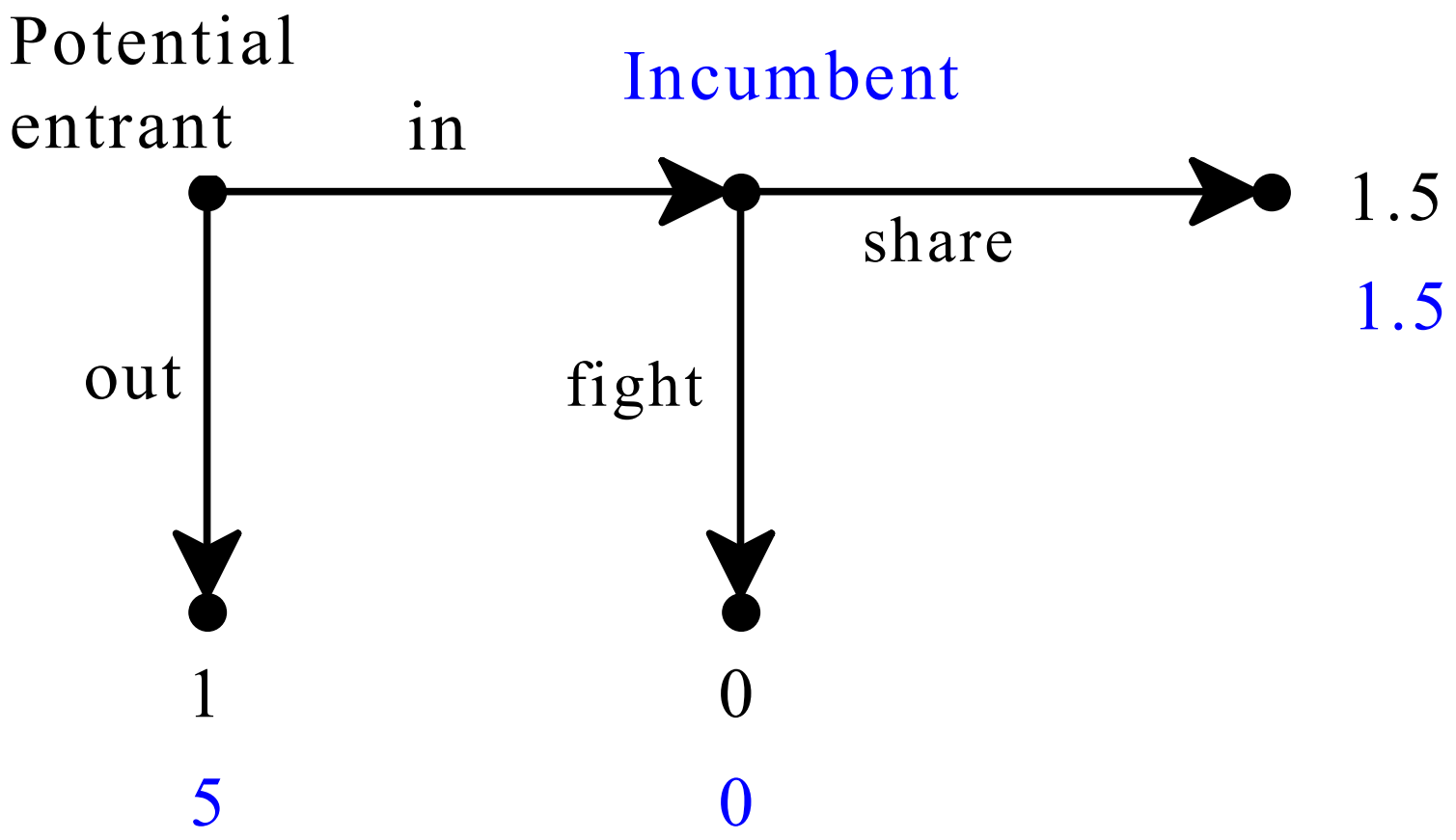

**Figure 14.6**

Decisions are made sequentially, as follows. At date $t$ ($t$ = 1,2) the businesswoman in town $t$ decides whether or not to enter and if she enters then the chain store decides whether or not to fight in that town. What happens in town $t$ (at date $t$) becomes known to everybody. Thus the businesswoman in town 2 at date 2 knows what happened in town 1 at date 1 (either that there was no entry or that entry was met with a fight or that entry was accommodated). Intuition suggests that in this game the threat by the Incumbent to fight early entrants might be credible, for the following reason. The Incumbent could tell Businesswoman 1 the following:

> "It is true that, if you enter and I fight, I will make zero profits, while by accommodating your entry I would make \$1.5 million and thus it would seem that it cannot be in my interest to fight you. However, somebody else is watching us, namely Businesswoman 2. If she sees that I have fought your entry then she might fear that I would do the same with her and decide to stay out, in which case in town 2 I would make \$5 million, so that my total profits in towns 1 and 2 would be $\$(0+5)$ = \$5 million. On the other hand, if I accommodate your entry, then she will be encouraged to entry herself and I will make \$1.5 million in each town, for a total profit of \$3 million. Hence, as you can see, it is indeed in my interest to fight you and thus you should stay out."

We showed in Chapter 2 that the notion of backward induction does not capture this intuition. In the game depicting the entire interaction (Figure 2.10, Chapter 2) there was a unique backward-induction solution whose corresponding outcome was that both businesswomen entered and the Incumbent accommodated entry in both towns. The reason why the backward-induction solution did not capture the "reputation" argument outlined above was explained in Chapter 2. We remarked there that, in order to capture the reputation effect, one would need to allow for some uncertainty in the mind of some of the players. This is what we will show below.

Suppose that, in principle, there are two types of incumbent monopolists: the rational type and the hotheaded type. The payoffs of a rational type are as shown in Figure 14.6, while a hotheaded type enjoys fighting; Figure 14.8 summarizes the payoffs for the two types of Incumbent (while there is only one type of potential entrant, whose payoffs are as shown earlier).

*If the Incumbent is rational*:

| | | potential entrant: in | potential entrant: out |
|---|---|---|---|
| Incumbent | fight | 0 , 0 | 5 , 1 |
| | share | 1.5 , 1.5 | 5 , 1 |

*If the Incumbent is hotheaded*:

| | | potential entrant: in | potential entrant: out |
|---|---|---|---|
| Incumbent | fight | 2 , 0 | 5 , 1 |
| | share | 1.5 , 1.5 | 5 , 1 |

**Figure 14.7**

Thus a hotheaded Incumbent enjoys fighting entry (he considers a fight to be as good as getting $2 million!).

Consider the following situation of one-sided incomplete information. As a matter of fact, the Incumbent is rational and this fact is common knowledge between the Incumbent and Potential Entrant 1 (from now on denoted by PE-1). However, Potential Entrant 2 (PE-2) is uncertain whether the Incumbent is rational or hotheaded and attaches probability $p$ to the latter case. The beliefs of PE-2 are common knowledge as are the payoffs of PE-1 and PE-2. This situation can be illustrated with the knowledge-belief structure of Figure 14.8 below.

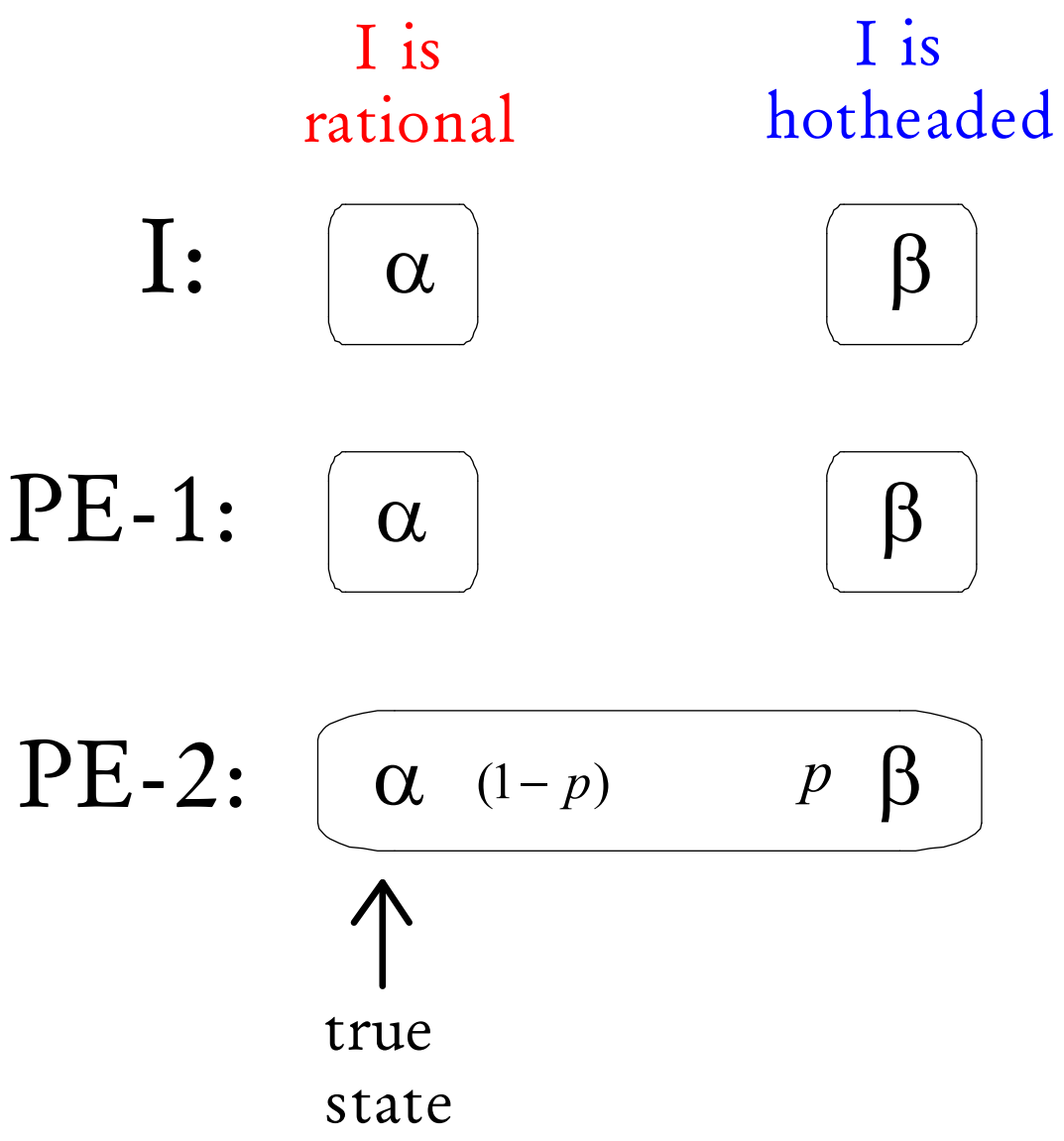

**Figure 14.8**

Applying the Harsanyi transformation to the situation depicted in Figure 14.8 yields the extensive-form game shown in Figure 14.9 below. We want to see if this situation of one-sided incomplete information can indeed capture the reputation effect discussed above; in particular, we want to check if there is a weak sequential equilibrium of the game of Figure 14.9 where the Incumbent's strategy would be to fight PE-1's entry in order to scare off PE-2 and, as a consequence, PE-1 decides to stay out and so does PE-2.

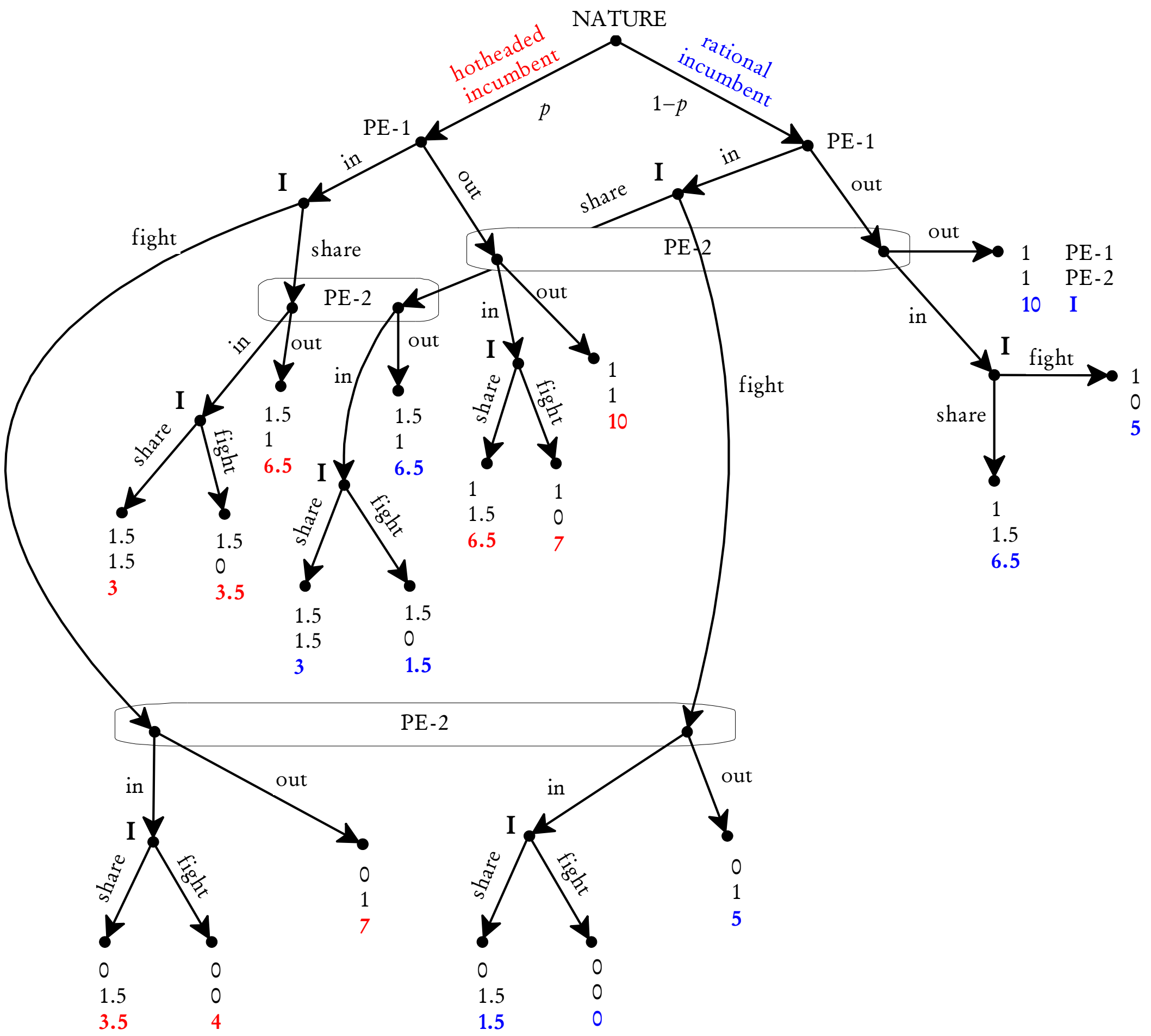

**Figure 14.9**

Using the notion of weak sequential equilibrium allows us to simplify the game, by selecting the strictly dominant choice for the Incumbent at each of his singleton nodes followed only by terminal nodes. It is straightforward to check that at such nodes a hotheaded Incumbent would choose "fight" while a rational Incumbent would choose "share". Thus we can delete those decision nodes and replace them with the payoff vectors associated with the optimal choice. The simplified game is shown in Figure 14.10 below.

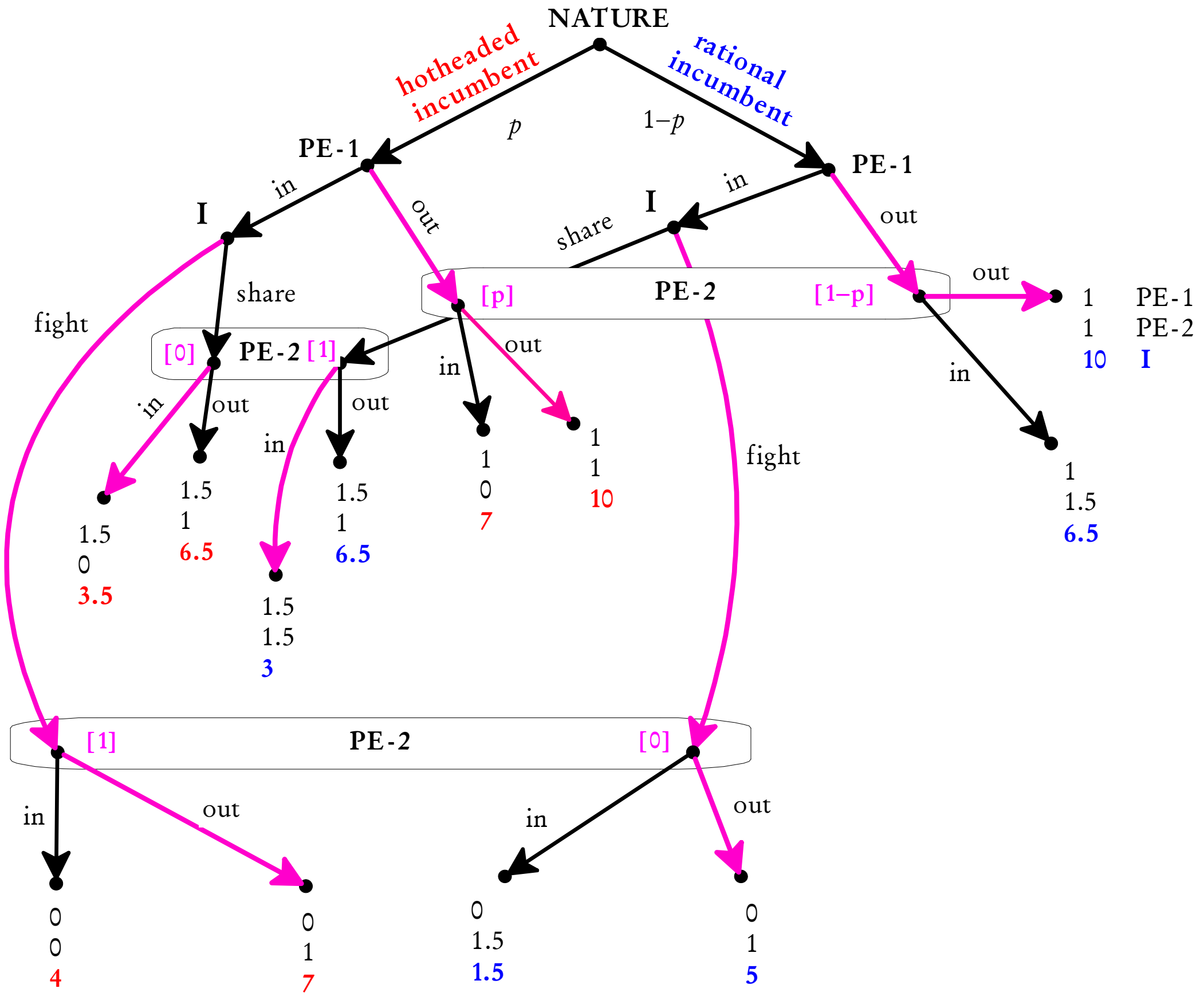

**Figure 14.10**

Consider the following pure-strategy profile, call it $\sigma$, for the simplified game of Figure 14.10 (it is highlighted in dark pink in Figure 14.10):

(1) PE-1's strategy is "out" at both nodes,

(2) PE-2's strategy is "out" at the top information set (after having observed that PE-1 stayed out), "in" at the middle information set (after having observed that PE-1's entry was followed by the Incumbent sharing the market with PE-1) and "out" at the bottom information set (after having observed that PE-1's entry was followed by the Incumbent fighting against PE-1),

(3) the Incumbent's strategy is to fight entry of PE-1 in any case (that is, whether the Incumbent himself is hotheaded or rational).

We want to show that $\sigma$, together with the system of beliefs $\mu$ shown in Figure 14.10 (where at her top information PE-2 assigns probability $p$ to the Incumbent being hotheaded, at her middle information set she assigns probability 1 to the Incumbent being rational and at her bottom information set she assigns probability 1 to the Incumbent being hotheaded) constitutes a weak sequential equilibrium for any value of $p \geq \frac{1}{3}$. Bayesian updating is satisfied at PE-2's top information set (the only information set reached by $\sigma$), and at the other two information sets of PE-2 any beliefs are allowed by the notion of weak sequential equilibrium. Sequential rationality also holds:

- For PE-1, at the left node "in" yields 0 and "out" yields 1, so that "out" is sequentially rational; the same is true at the right node.

- For the Incumbent, at the left node "fight" yields 7 and "share" yields 3.5, so that "fight" is sequentially rational; at the right node "fight" yields 5 and "share" yields 3, so that "fight" is sequentially rational.

- For PE-2, at the top information set "in" yields an expected payoff of $p(0)+(1-p)(1.5)$ and "out" yields 1, so that "out" is sequentially rational as long as $p \geq \frac{1}{3}$; at the middle information set, given her belief that $I$ is rational, "in" yields 1.5 and "out" yields 1, so that "in" is sequentially rational; at the bottom information set, given her belief that $I$ is hotheaded, "in" yields 0 and "out" yields 1, so that "out" is sequentially rational.

Thus the equilibrium described above captures the intuition suggested in Chapter 2, namely that – even though it is common knowledge between the Incumbent and PE-1 that the Incumbent is rational and thus would suffer a loss of 1.5 by fighting PE-1's entry – it is still credible for the Incumbent to threaten to fight PE-1's entry because it would influence the beliefs of PE-2 and induce her to stay out; understanding the credibility of this threat, it is optimal for PE-1 to stay out.

The above argument exploits the fact that the notion of weak sequential equilibrium allows for any beliefs whatsoever at unreached information sets. However, the beliefs postulated for PE-2 at her middle and bottom information sets seem highly plausible. Indeed, as shown in Exercise 14.2, the assessment described above is also a sequential equilibrium.

We now turn to an example that deals with the issue of labor-management negotiations and the inefficiency of strikes. It is not uncommon to observe a firm and a union engaging in unsuccessful negotiations leading to a strike by the workers, followed by renewed negotiations and a final agreement. Strikes are costly for the workers, in terms of lost wages, and for the firm, in terms of lost production. Why, then, don't the parties reach the agreement at the very beginning, thus avoiding the inefficiency of a strike? The answer in many cases has to do with the fact that there is incomplete information on the side of the labor union and enduring a strike is the only credible way for the firm to convince the union to reduce its demands. We shall illustrate this in a simple example of one-sided incomplete information.

Consider the following game between a new firm and a labor union. The union requests a wage (either high, $w_H$, or low, $w_L$) and the firm can accept or reject. If the union's request is accepted, then a contract is signed and production starts immediately. If the request is rejected, then the union goes on strike for one period and at the end of that period makes a second, and last, request to the firm, which the firm can accept or reject. If the firm rejects, then it cannot enter the industry. When no agreement is signed, both parties get a payoff of 0. Both firm and union have a discount factor of $\delta$, with $0 < \delta < 1$, which means that \$1 accrued one period into the future is considered to be equivalent to \$$\delta$ at the present time (this captures the cost of waiting and thus the desirability of avoiding a strike). The extensive-form game is shown in Figure 14.11, where $\pi$ denotes the firm's profit (gross of labor costs).

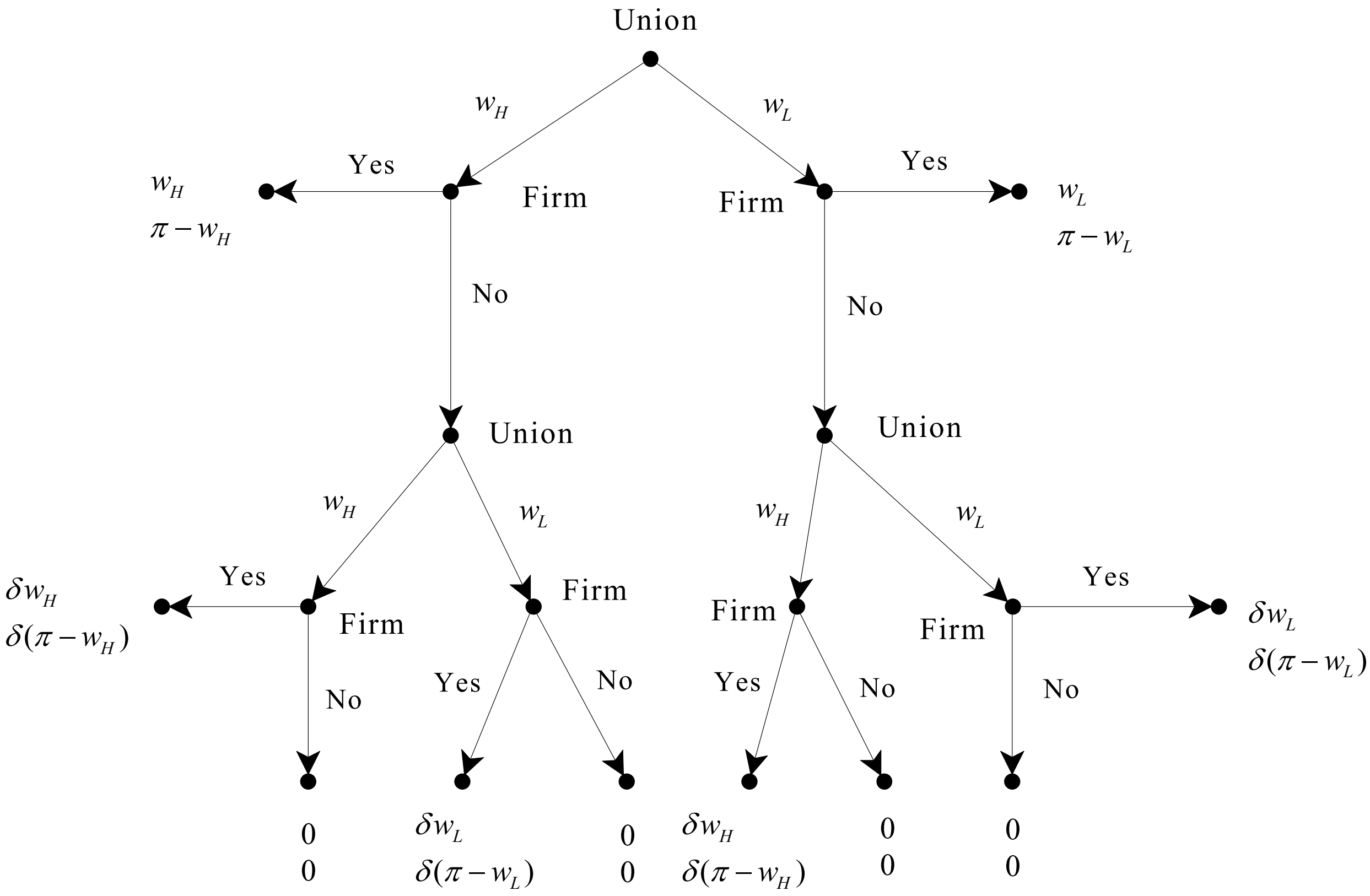


**Figure 14.11**

Suppose, however, that the union does not know the firm's expected profit $\pi$ (gross of labor costs). It is not unlikely that the management will know what profits to expect, while the union will not (because, for example, the union does not have enough information on the intensity of competition, the firm's non-labor costs, etc.). Suppose we have a one-sided incomplete information situation where the union believes that $\pi$ can have two values: $\pi_H$ (high) and $\pi_L$ (low) and assigns probability $\alpha$ to $\pi_H$ and $(1-\alpha)$ to $\pi_L$. This situation of one-sided incomplete information is illustrated in Figure 14.12.

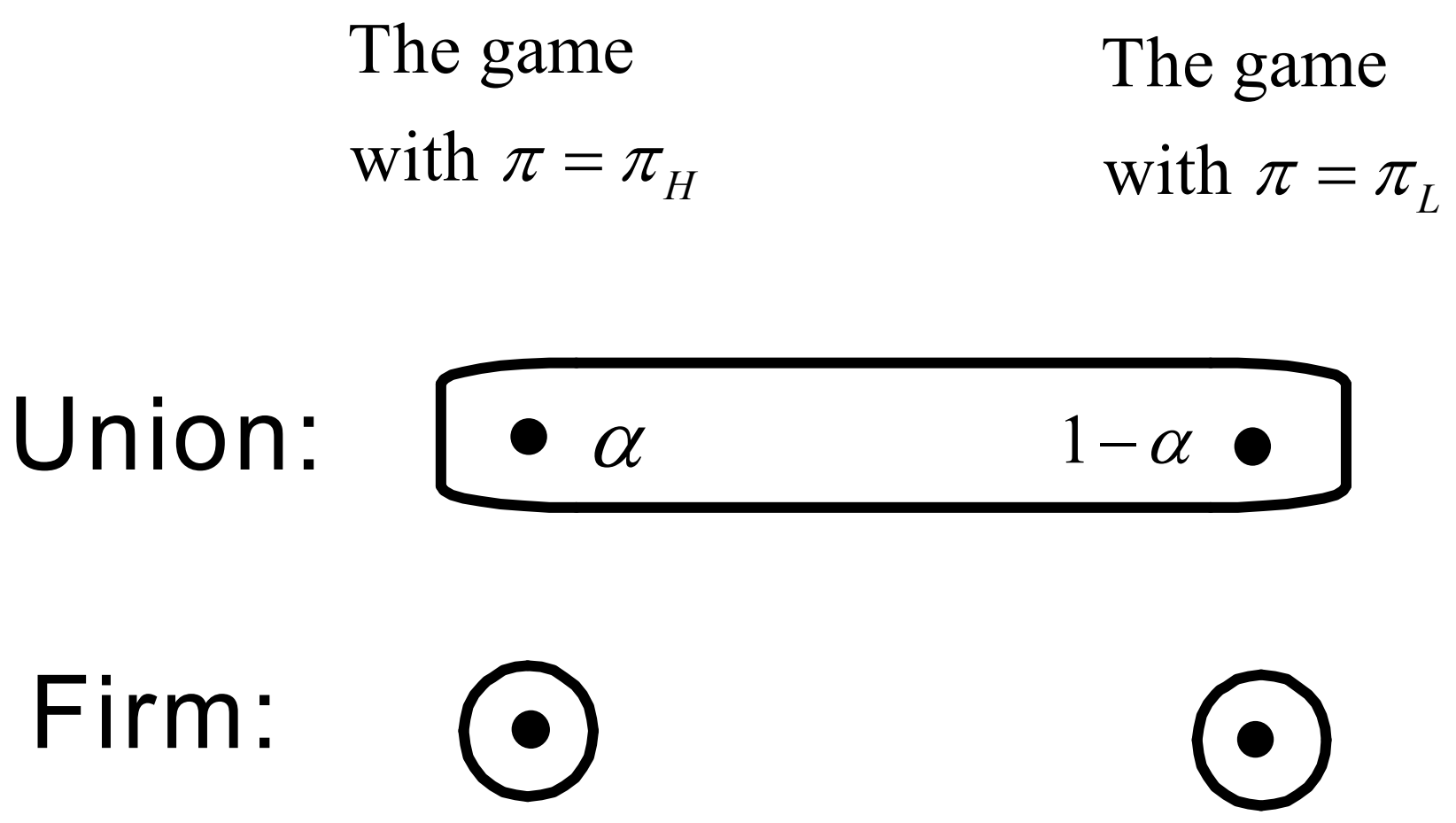


**Figure 14.12**

Let $\pi_H > \pi_L > 0$ and $w_H > w_L > 0$ (so that $H$ indeed means "high" and $L$ indeed means "low"). We also assume that $\pi_H - w_H > 0$ (so that the high-profit firm could in fact afford to pay a high wage), $\pi_L - w_L > 0$ (that is, the low wage is sufficiently low for the low-profit firm to be able to pay it) and $\pi_L - w_H < 0$ (that is, the low-profit firm cannot afford to pay a high wage). Finally we assume that the true state is the one where $\pi = \pi_L$, that is, the firm's potential profits are in fact low. These assumptions imply that it is in the interest of both the firm and the union to sign a contract (even if the firm's profits are low).

Using the Harsanyi transformation we can convert the situation shown in Figure 14.12 into the extensive-form game shown in Figure 14.13.

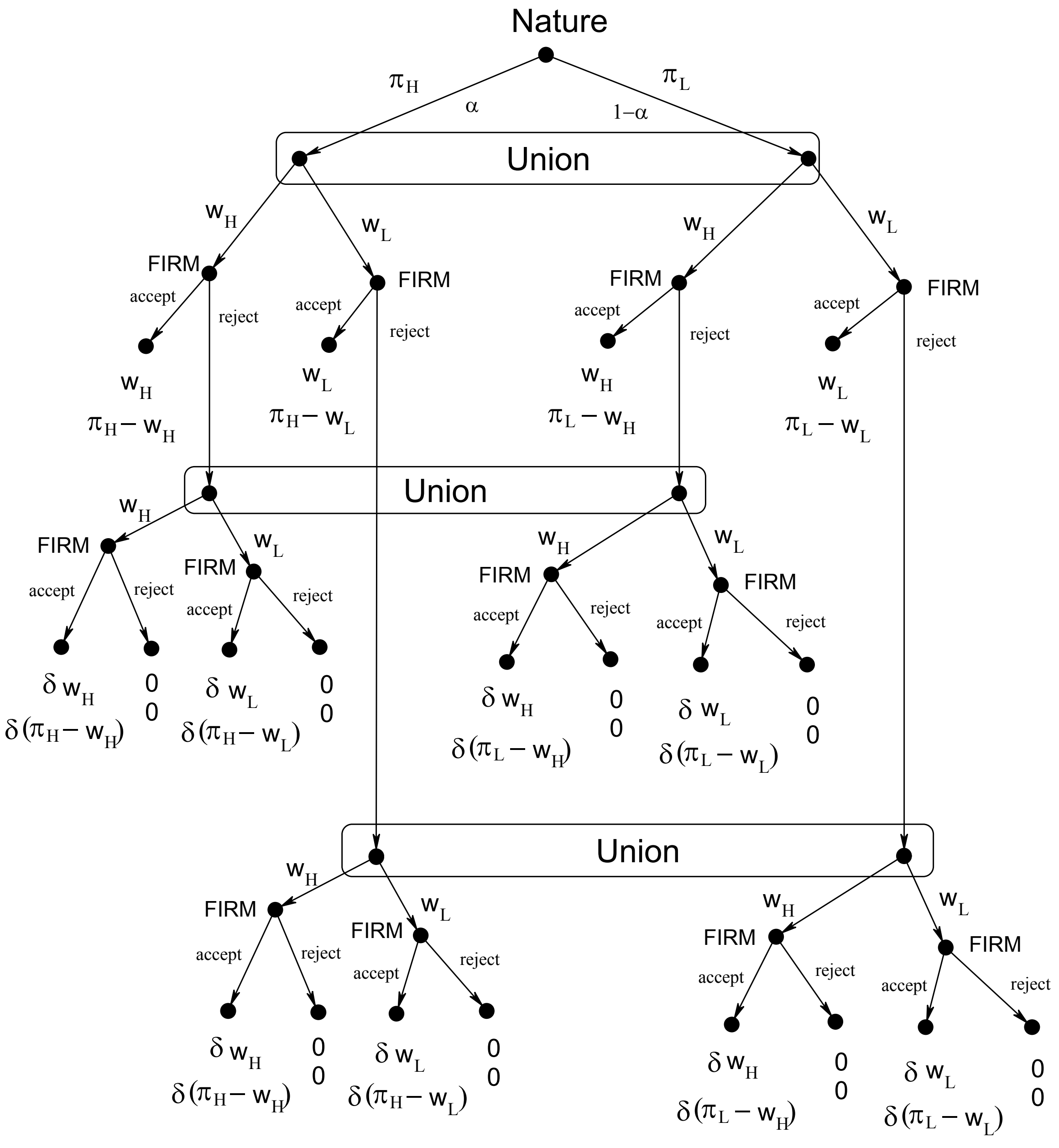
Nature
π_H
α
π_L
1−α
Union
w_H
w_L
w_H
w_L
FIRM
FIRM
FIRM
FIRM
accept
reject
accept
reject
accept
reject
accept
reject
w_H
π_H − w_H
w_L
π_H − w_L
w_H
π_L − w_H
w_L
π_L − w_L
Union
w_H
w_L
w_H
w_L
FIRM
FIRM
FIRM
FIRM
accept
reject
δ w_H
δ(π_H − w_H)
0
0
δ w_L
δ(π_H − w_L)
0
0
δ w_H
δ(π_L − w_H)
0
0
δ w_L
δ(π_L − w_L)
0
0
Union
w_H
w_L
w_H
w_L
FIRM
FIRM
FIRM
FIRM
accept
reject
δ w_H
δ(π_H − w_H)
0
0
δ w_L
δ(π_H − w_L)
0
0
δ w_H
δ(π_L − w_H)
0
0
δ w_L
δ(π_L − w_L)
0
0

**Figure 14.13**

Let us find conditions under which the following strategy profile is part of a separating weak sequential equilibrium:

(1) The union requests a high wage in the first period.

(2) The high-profit firm accepts while the low-profit firm rejects.

(3) After rejection of the first-period high-wage offer the union requests a low wage and the firm accepts.

That is, we want an equilibrium where the low-profit firm endures a strike to signal to the union that its expected profits are low and cannot afford to pay a high wage. The union reacts to the signal by lowering its demand. Of course, we need to worry about the fact that the high-profit firm might want to masquerade as a low-profit firm by rejecting the first-period offer (that is, by sending the same signal as the low-profit firm).

First of all, we can simplify the game by eliminating the second-period choice of the firm (under our assumptions any second-period offer will be accepted, except for $w_H$ by the low firm). The simplified game is shown in Figure 14.14.

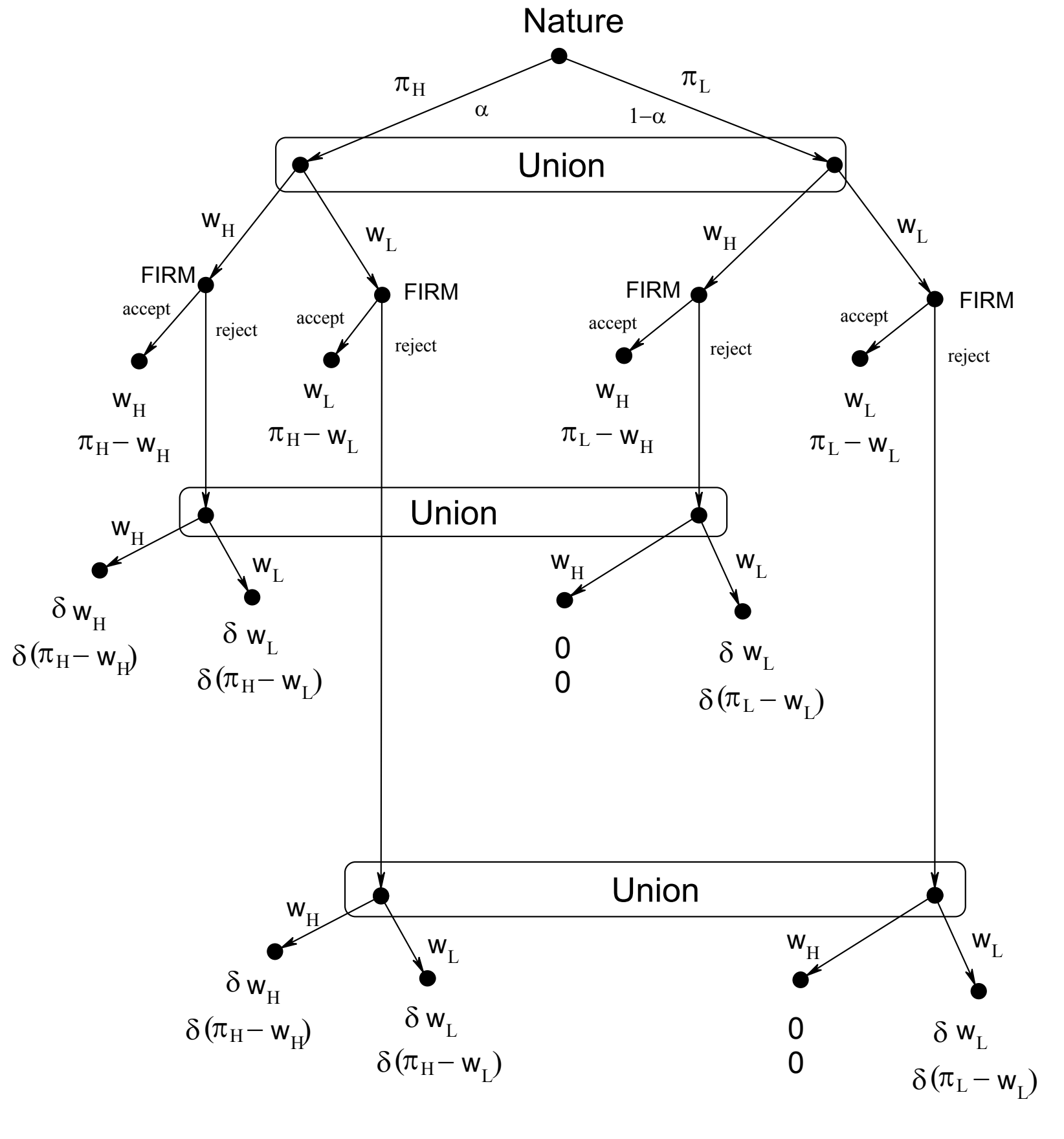


**Figure 14.14**

In order for a low-wage offer to be optimal for the union in Period 2 it is necessary for the union to assign sufficiently high probability to the firm being a low-profit one: if $p$ is this probability, then we need $\delta w_L \geq (1-p)\delta w_H$, that is, $p \geq 1 - \frac{w_L}{w_H}$ (for example, $p = 1$ would be fine). For the high-profit firm not to have an incentive to send the same signal as the low-profit firm it is necessary that $\pi_H - w_H \geq \delta(\pi_H - w_L)$, that is $\delta \leq \frac{\pi_H - w_H}{\pi_H - w_L}$. If this condition is satisfied, then the high-profit firm will accept $w_H$ in period 1, the low-profit firm will reject and Bayesian updating requires the union to assign probability 1 to the right node of its middle information set (that is, to the firm being a low-profit type), in which case adjusting its demand to $w_L$ is sequentially rational. On the other hand, requesting a low wage in period 1 will lead to both types of the firm accepting immediately. So in order for a high-wage offer to be optimal for the union in period 1 it is necessary that $w_L \leq \alpha w_H + (1-\alpha)\delta w_L$, which will be true if $\alpha$ is sufficiently large, that is, if $\alpha \geq \frac{(1-\delta)w_L}{w_H - \delta w_L}$. Finally, given the above equilibrium choices, the notion of weak sequential equilibrium imposes no restrictions on the beliefs (and hence choice) of the union at the bottom information set.

For example, all the above inequalities are satisfied if: $\pi_H = 100$, $\pi_L = 55$, $w_H = 60$, $w_L = 50$, $\alpha = 0.7$, $\delta = 0.6$. Then $\pi_H - w_H = 40 > \delta(\pi_H - w_L) = 0.6\,(50) = 30$, $w_L = 50 < \alpha w_H + (1-\alpha)\,\delta w_L = (0.7)\,60 + 0.3\,(0.6)\,50 = 51$.

We conclude this section with one more example that has to do with the effectiveness of truth-in-advertising laws.

Consider the case of a seller and a buyer. The seller knows the quality $x$ of his product, while the buyer does not, although she knows that he knows. The buyer can thus ask the seller to reveal the information he has. Suppose that there is a truth-in-advertising law which imposes harsh penalties for false claims in advertising. This is not quite enough, because the seller can tell the truth without necessarily revealing all the information. For example, if $x$ is the octane content of gasoline and $x = 89$, then the following are all true statements:

"the octane content of this gasoline is at least 70"

"the octane content of this gasoline is at least 85"

"the octane content of this gasoline is at most 89"

"the octane content of this gasoline is exactly 89"

etc.

Other examples are: (1) the fat content of food (the label on a package of ground meat might read "not more than 30% fat"), (2) fuel consumption for cars (a car might be advertised as yielding "at least 35 miles per gallon"), (3) the label on a mixed-nut package might read "not more than 40% peanuts", etc.

An interesting question is: Would the seller reveal all the information he has or would he try to be as vague as possible (while being truthful)? In the latter case, will the buyer be in a worse position than she would be in the case of complete information?

Milgrom and Roberts (1986) consider the following game (resulting from applying the Harsanyi transformation to the situation of one-sided incomplete information described above). First Nature selects the value of $x$, representing the seller's information, from a finite set $X$; as usual, the probabilities with which Nature chooses are a reflection of the buyer's beliefs. The seller "observes" Nature's choice and makes an assertion $A$ to the buyer; $A$ is a subset of $X$. The seller is restricted to make true assertions, that is, we require that $x \in A$. The buyer observes the seller's claim $A$ and then selects a quantity $q \geq 0$ to buy. Milgrom and Roberts show that (under reasonable hypotheses), if the buyer adopts a skeptical view concerning the seller's claim, that is, she always interprets the seller's claim in a way which is least favorable to the seller (for example "not more than 30% fat" is interpreted as "exactly 30% fat"), then there is an equilibrium where the outcome is the same as it would be in the case of complete information. We shall illustrate this phenomenon by means of a simple example.

Suppose that there are three possible quality levels for the good under consideration: low ($l$), medium ($m$) and high ($h$); thus $X = \{l,m,h\}$. The buyer believes that the three quality levels are equally likely. The buyer has to choose whether to buy one unit or two units. The seller can only make truthful claims. For example, if the quality is $l$, the seller's possible claims are $\{l\}$ (full revelation) or vague claims such as $\{l,m\}$, $\{l,h\}$, $\{l,m,h\}$. The extensive-form game is shown in Figure 14.15.

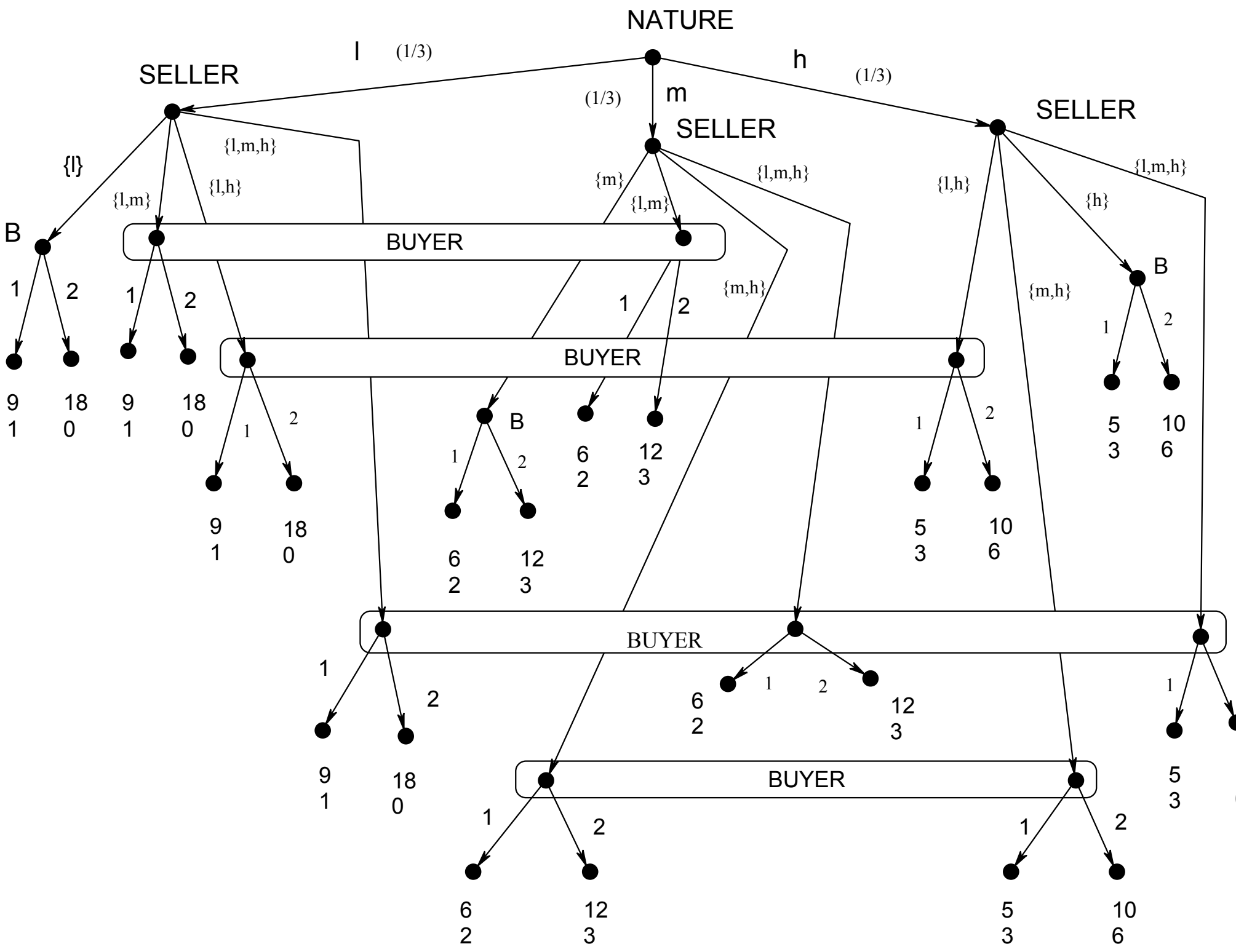


**Figure 14.15**

It is straightforward to check that the following is a weak sequential equilibrium:

- The seller claims $\{l\}$ if Nature chooses $l$, $\{m\}$ if Nature chooses $m$ and $\{h\}$ if Nature chooses $h$ (that is, the seller reveals the whole truth by choosing not to make vague claims).
- The buyer buys one unit if told $\{l\}$, two units if told $\{m\}$ or if told $\{h\}$; furthermore, the buyer adopts beliefs which are least favorable to the seller: if told $\{l,m\}$ or $\{l,h\}$ or $\{l,m,h\}$ she will believe $l$ with probability one and buy one unit, if told $\{m,h\}$ she will believe $m$ with probability one and buy two units. All these beliefs are admissible because they concern information sets that are not reached in equilibrium and therefore Bayesian updating does not apply.

On the other hand, as shown in Exercise 14.4, if the buyer is naïve then there is an equilibrium where the seller chooses to make vague claims.

This is a good time to test your understanding of the concepts introduced in this section, by going through the exercises in Section 14.E.1 of Appendix 14.E at the end of this chapter.

# 14.2 Multi-sided incomplete information

The case of situations of multi-sided incomplete information involving dynamic games is conceptually the same as the case of multi-sided situations of incomplete information involving static games. In this section we shall go through one example.

Consider the following situation of two-sided incomplete information. A seller (player $S$) owns an item that a buyer (player $B$) would like to purchase. The seller's reservation price is $s$ (that is, she is willing to sell if and only if the price paid by the buyer is at least $s$) and the buyer's reservation price is $b$ (that is, he is willing to buy if and only if the price is less than or equal to $b$). It is common knowledge between the two that (1) both $b$ and $s$ belong to the set $\{1,2,...,n\}$, (2) the buyer knows the value of $b$ and the seller knows the value of $s$, (3) both the buyer and the seller attach equal probability to all the possibilities among which they are uncertain.

Buyer and Seller play the following game. First the buyer makes an offer of a price $p \in \{1,...,n\}$ to the seller. If $p = n$ the game ends and the object is exchanged for \$$p$. If $p < n$ then the seller either accepts (in which case the game ends and the object is exchanged for \$$p$) or makes a counter-offer of $p' > p$, in which case either the buyer accepts (and the game ends and the object is exchanged for $\$\,p'$) or the buyer rejects, in which case the game ends without an exchange. Payoffs are as follows:

$$\pi_{seller} = \begin{cases} 0 & \text{if there is no exchange} \\ x - s & \text{if exchange takes place at price } \$x \end{cases}$$

$$\pi_{buyer} = \begin{cases} 0 & \text{if there is no exchange} \\ b - p & \text{if exchange takes place at price } \$p \text{ (the initial offer)} \\ b - p' - \varepsilon & \text{if exchange takes place at price } \$p' \text{ (the counter-offer)} \end{cases}$$

where $\varepsilon > 0$ is a measure of the buyer's "hurt feelings" for seeing his initial offer rejected. These are von Neumann-Morgenstern payoffs.

Let us start by focusing on the case $n = 2$. First we represent the situation described above by means of an interactive knowledge-belief structure. A possible state can be written as a pair $(b, s)$, where $b$ is the reservation price of the buyer and $s$ that of the seller. Thus when $n = 2$ the possible states are (1,1), (1,2), (2,1) and (2,2). Figure 14.16 below represents this two-sided situation of incomplete information.

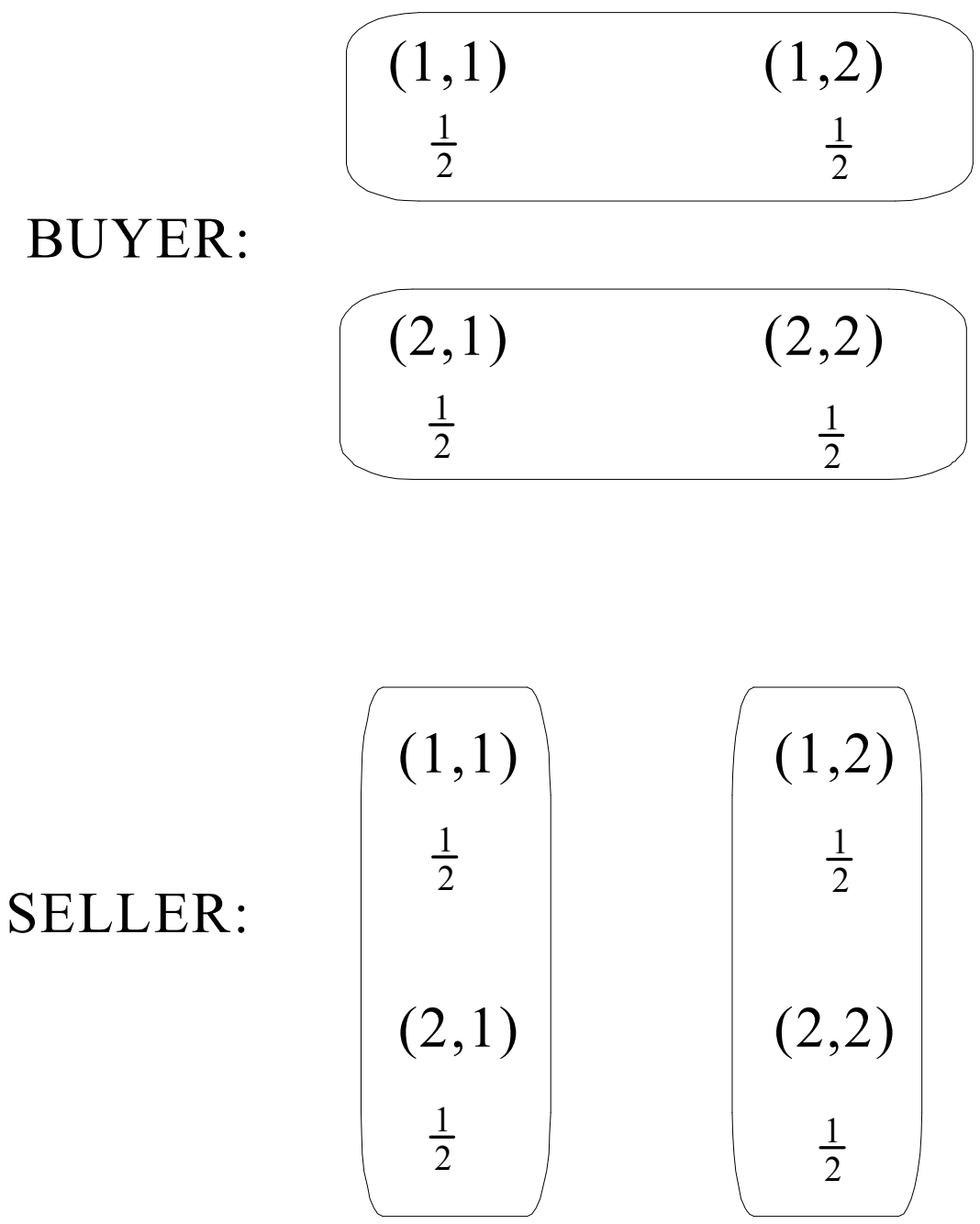


**Figure 14.16**

Next we apply the Harsanyi transformation to the situation illustrated in Figure 14.16. For each state there is a corresponding extensive-form game. The four possible games are shown in Figure 14.17 below (the top number is the Buyer's payoff and the bottom number the Seller's payoff).

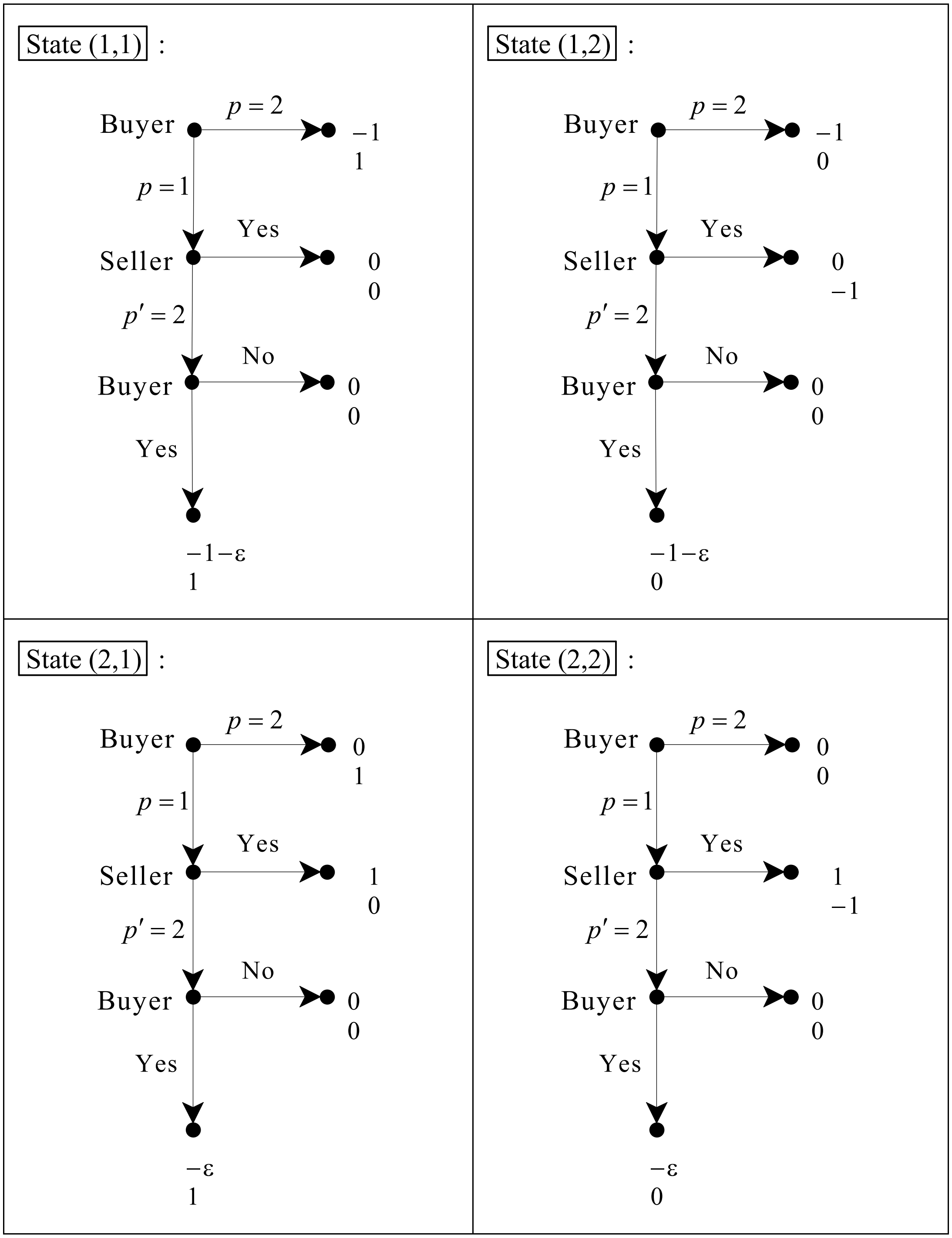
State (1,1) :
Buyer
p = 2
−1
1
p = 1
Yes
Seller
0
0
p′ = 2
No
Buyer
0
0
Yes
−1−ε
1
State (1,2) :
Buyer
p = 2
−1
0
p = 1
Yes
Seller
0
−1
p′ = 2
No
Buyer
0
0
Yes
−1−ε
0
State (2,1) :
Buyer
p = 2
0
1
p = 1
Yes
Seller
1
0
p′ = 2
No
Buyer
0
0
Yes
−ε
1
State (2,2) :
Buyer
p = 2
0
0
p = 1
Yes
Seller
1
−1
p′ = 2
No
Buyer
0
0
Yes
−ε
0

**Figure 14.17**

The extensive-form game that results from applying the Harsanyi transformation to the situation illustrated in Figure 14.16 is shown in Figure 14.18 below.

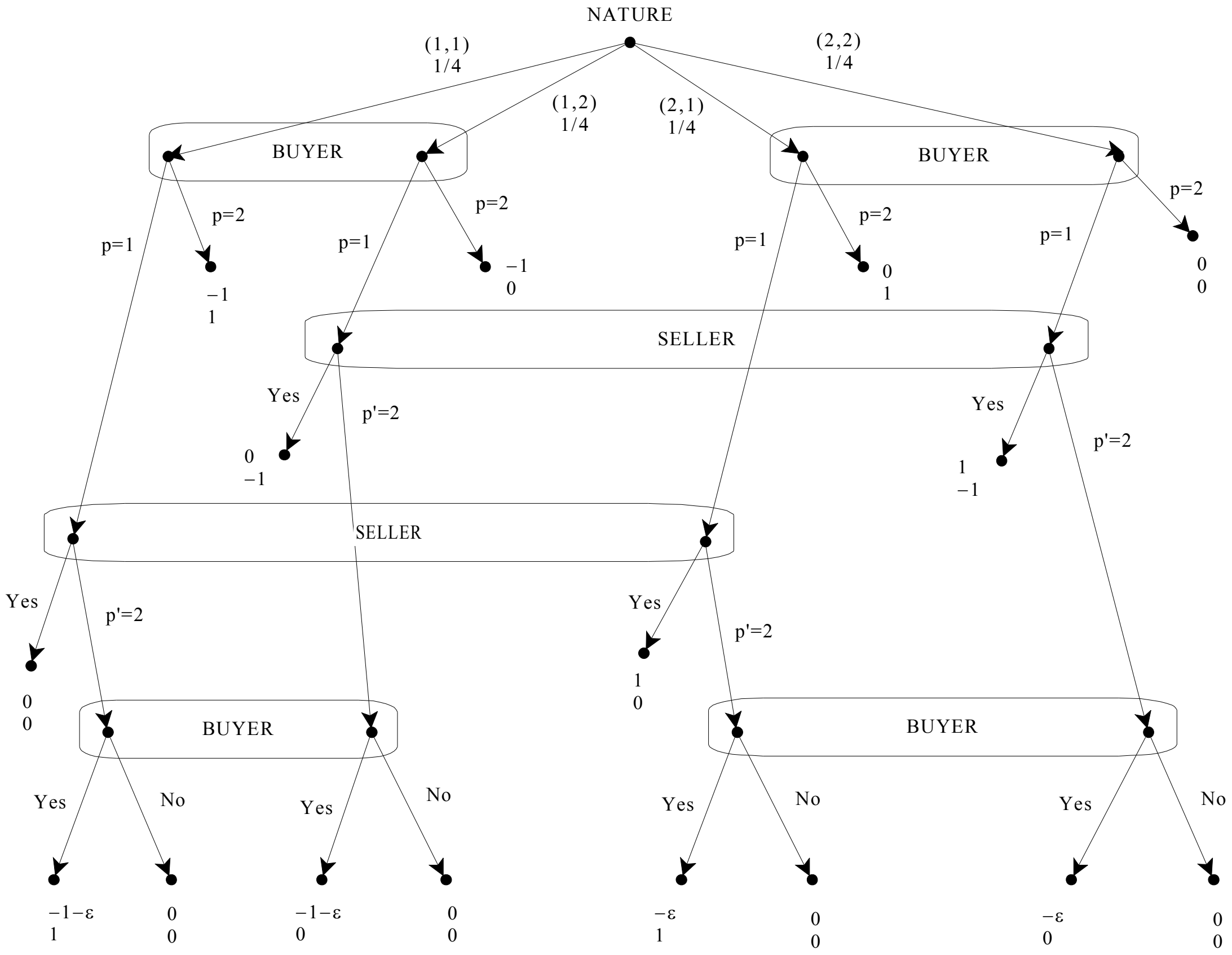


**Figure 14.18**

Let us find all the pure-strategy weak sequential equilibria of the game of Figure 14.18. First of all, note that at the bottom information sets of the Buyer, "Yes" is strictly dominated by "No" and thus a weak sequential equilibrium must select "No". Hence the game simplifies to the one shown in Figure 14.19 below.

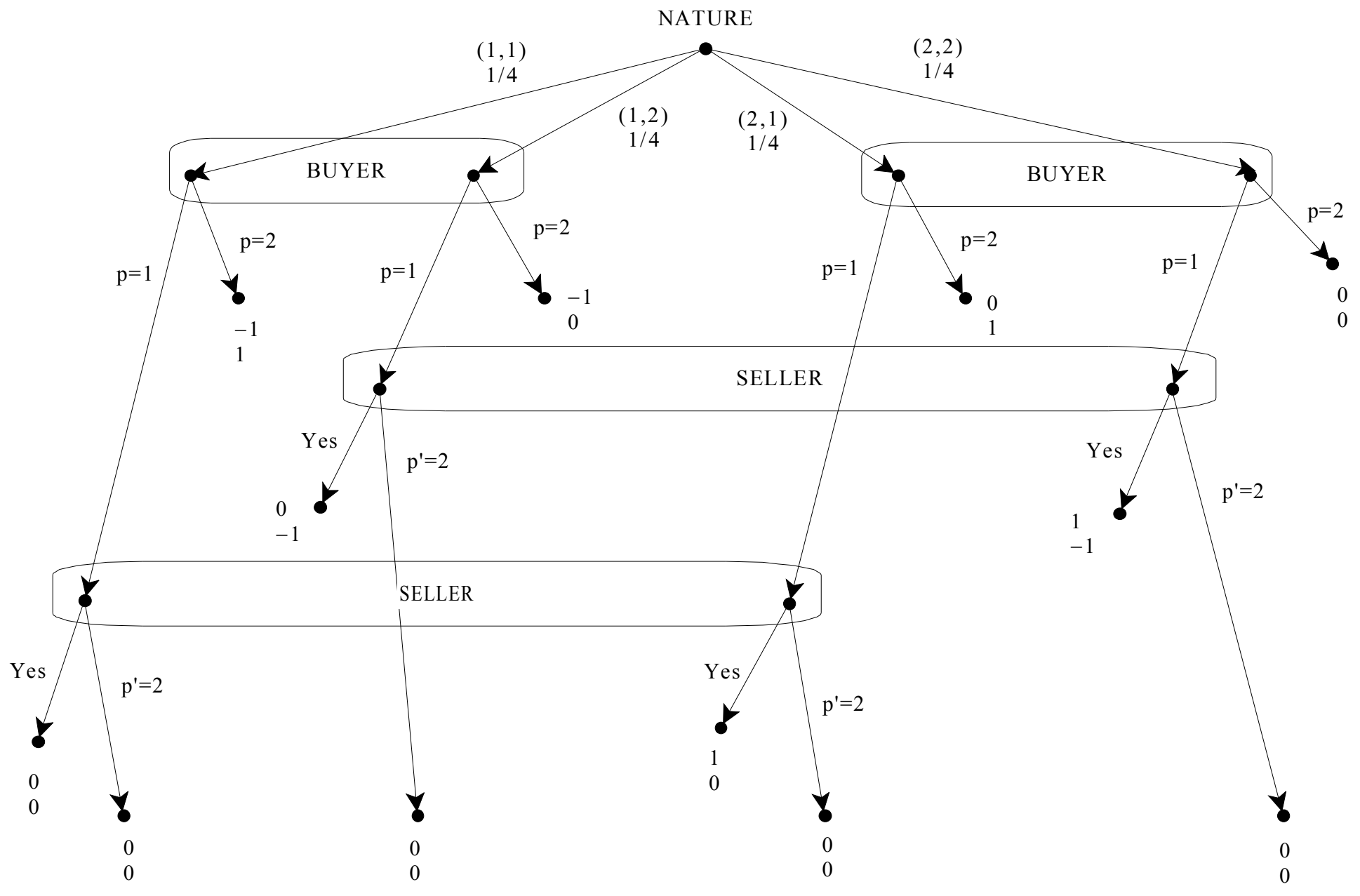


**Figure 14.19**

In the game of Figure 14.19, at the middle information of the Seller, making a counteroffer of $p' = 2$ strictly dominates "Yes". Thus the game can be further simplified as shown in Figure 14.20.

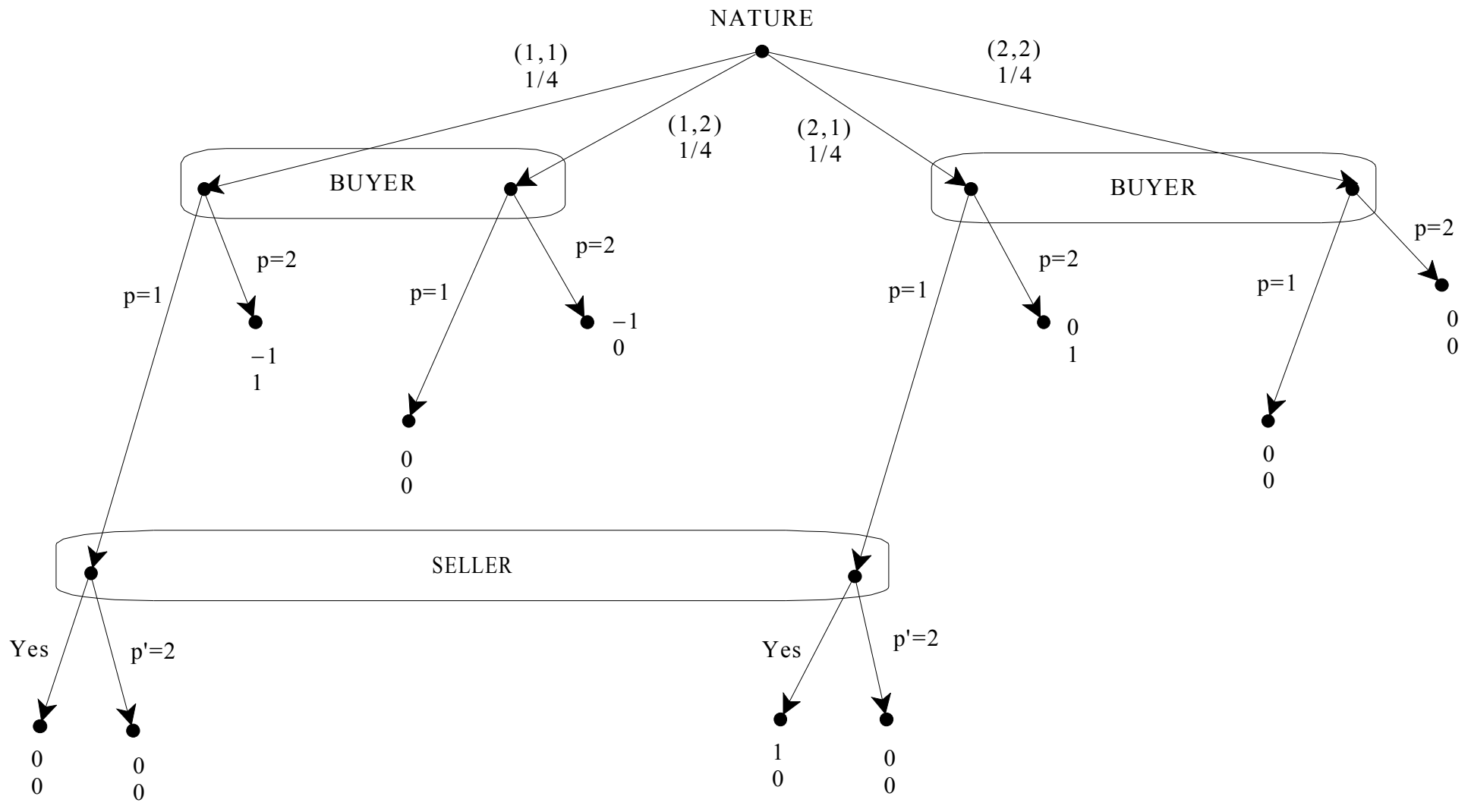


**Figure 14.20**

In the game of Figure 14.20 at the left information set of the Buyer $p = 1$ strictly dominates $p = 2$. At the right information set the Buyer's beliefs must

be $\frac{1}{2}$ on each node so that if the Seller's strategy is to say Yes, then $p = 1$ is the only sequentially rational choice there, otherwise both $p = 1$ and $p = 2$ are sequentially rational. Thus the pure-strategy weak sequential equilibria of the reduced game shown in Figure 14.20 are as follows:

(1) $\left((p=1, p=1), Yes\right)$ with beliefs given by probability ½ on each node at every information set.

(2) $\left((p=1, p=1),\ p'=2\right)$ with beliefs given by probability ½ on each node at every information set.

(3) $\left((p=1, p=2),\ p'=2\right)$ with beliefs given by (*i*) probability ½ on each node at both information sets of the Buyer and (*ii*) probability 1 on the left node at the information set of the Seller.

These equilibria can be extended to the original game of Figure 14.18 as follows:

(1) $\left((p=1, p=1, No, No), (p'=2, Yes)\right)$ with beliefs given by (*i*) probability ½ on each node at both information sets of the Buyer at the top and at both information sets of the Seller and (*ii*) probability 1 on right node at both information sets of Buyer at the bottom. The corresponding payoffs are: ¼ for the Buyer and 0 for the Seller.

(2) $\left((p=1, p=1, No, No), (p'=2, p'=2)\right)$ with beliefs given by probability ½ on each node at every information set. The corresponding payoffs are 0 for both Buyer and Seller.

(3) $\left((p=1, p=2, No, No), (p'=2, p'=2)\right)$ with beliefs given by (*i*) probability ½ on each node at both information sets of the Buyer at the top and at the lower left information set of the Buyer, (*ii*) any beliefs at the lower right information set of the Buyer and (*iii*) probability 1 on the left node at each information set of the Seller. The corresponding payoffs are: 0 for the Buyer and ¼ for the Seller.

Now let us consider the case $n$ = 100. Drawing the interactive knowledge-belief structure and the corresponding extensive-form game (obtained by applying the Harsanyi transformation) is clearly not practical. However, one can still reason about what could be a pure-strategy Bayesian Nash equilibrium of that game. As a matter of fact, there are many Bayesian Nash equilibria. One of them is the following: (*i*) each type of the Buyer offers a price equal to his reservation price, (*ii*) the Seller accepts if, and only if, that price is greater than or equal to her reservation price, (iii) at information sets of the Seller that are not reached the Seller rejects and counteroffers $100, (*iv*) at information sets of the Buyer that are not reached the Buyer says "No". The reader should convince himself/herself that this is indeed a Nash equilibrium.

This is a good time to test your understanding of the concepts introduced in this section, by going through the exercises in Section 14.E.2 of Appendix 14.E at the end of this chapter.

# Appendix 14.E: Exercises

## 14.E.1. Exercises for Section 14.1: One-sided incomplete information

The answers to the following exercises are in Appendix S at the end of this chapter.

**Exercise 14.1**. Consider the following situation of one-sided incomplete information. Players 1 and 2 are playing the extensive-form game shown below (where $z_1, z_2$, etc. are outcomes and the question mark stands for either outcome $z_5$ or outcome $z_6$). The outcome that is behind the question mark is actually outcome $z_5$ and Player 1 knows this, but Player 2 does not know. Player 2 thinks that the outcome behind the question mark is either $z_5$ or $z_6$ and assigns probability 25% to it being $z_5$ and probability 75% to it being $z_6$. Player 2 also thinks that whatever the outcome is, Player 1 knows (that is, if it is $z_5$, then Player 1 knows that it is $z_5$, and if it is $z_6$ then Player 1 knows that it is $z_6$). The beliefs of Player 2 are common knowledge between the two players.

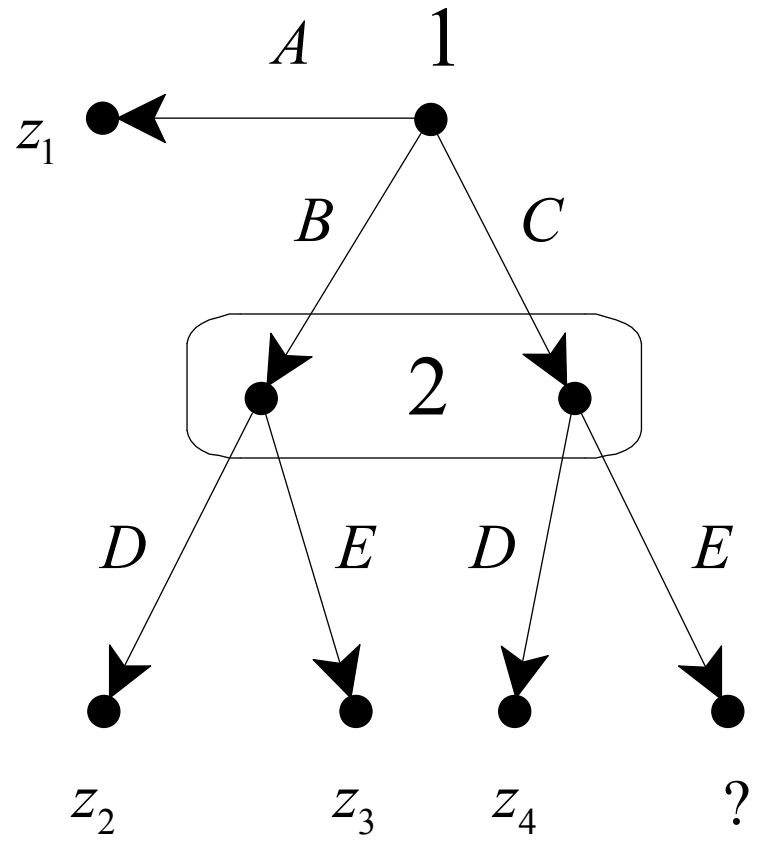


**(a)** Represent this situation of incomplete information using an interactive knowledge-belief structure.

**(b)** Apply the Harsanyi transformation to transform the situation represented in part (a) into an extensive-form frame. [Don't worry about payoffs for the moment.]

From now on assume that the following is common knowledge: (1) both players have von Neumann-Morgenstern preferences, (2) the ranking of Player

1 is $\begin{pmatrix} best & second & worst \\ z_4, z_6 & z_1 & z_2, z_3, z_5 \end{pmatrix}$ and he is indifferent between $z_1$ for sure and the lottery $\begin{pmatrix} z_6 & z_5 \\ 0.5 & 0.5 \end{pmatrix}$, (3) the ranking of player 2 is $\begin{pmatrix} best & second & third & worst \\ z_6 & z_4 & z_2, z_5 & z_1, z_3 \end{pmatrix}$ and she is indifferent between $z_4$ for sure and the lottery $\begin{pmatrix} z_6 & z_3 \\ 0.5 & 0.5 \end{pmatrix}$ and she is also indifferent between $z_2$ for sure and the lottery $\begin{pmatrix} z_6 & z_3 \\ 0.25 & 0.75 \end{pmatrix}$.

**(c)** Calculate the von Neumann-Morgenstern normalized utility functions for the two players.

**(d)** Is there a weak sequential equilibrium of the game of part (b) where Player 1 always plays $A$ (thus a pooling equilibrium)? Explain your answer.

**(e)** Is there a weak sequential equilibrium of the game of part (b) where player 1 always plays $C$ (thus a pooling equilibrium)? Explain your answer.

**(f)** Is there a pure-strategy weak sequential equilibrium of the game of part (b) where Player 1 does not always choose the same action (thus a separating equilibrium)?

**Exercise 14.2.** Show that the assessment highlighted in the following figure, which reproduces Figure 14.10, is a sequential equilibrium as long as $p \geq \frac{1}{3}$.

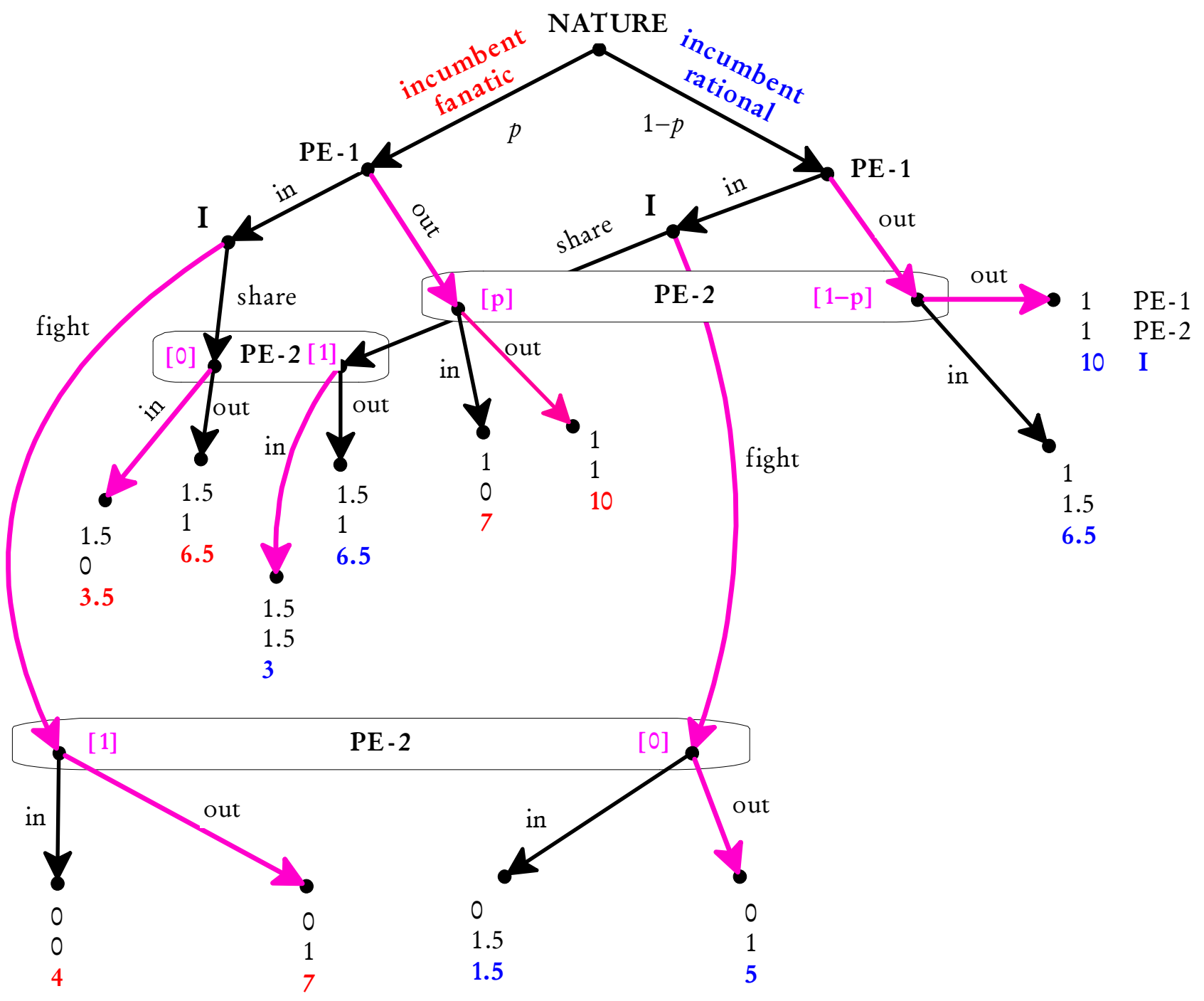

**Exercise 14.3**. **(a)** For the case where $p \leq \frac{1}{3}$, find a pure-strategy weak sequential equilibrium of the game shown in Exercise 14.1 (which reproduces Figure 14.10) where PE-1 stays out but PE-2 enters; furthermore, the Incumbent would fight PE-1 if she entered. Prove that what you propose is a weak sequential equilibrium.

**(b)** Either prove that the weak sequential equilibrium of part (a) is a sequential equilibrium or prove that it is not a sequential equilibrium.

**Exercise 14.4**. Show that in the "truth-in-advertising" game of Figure 14.15 there is a weak sequential equilibrium where the seller makes vague claims (that is, does not reveal the whole truth).

**Exercise 14.5**. Consider a simpler version of the "truth-in-advertising" game, where there are only two quality levels: *L* and *H*. The payoffs are as follows:

- If the quality is *L* and the buyer buys one unit, the seller's payoff is 9 and the buyer's payoff is 1,
- If the quality is *L* and the buyer buys two units, the seller's payoff is 18 and the buyer's payoff is 0,
- If the quality is *H* and the buyer buys one unit, the seller's payoff is 6 and the buyer's payoff is 2,
- If the quality is *H* and the buyer buys two units, the seller's payoff is 12 and the buyer's payoff is 3,

Let the buyer's initial beliefs be as follows: the good is of quality *L* with probability $p \in (0,1)$ (and of quality *H* with probability $1 - p$).

**(a)** Draw the extensive-form game that results from applying the Harsanyi transformation to this one-sided situation of incomplete information.

**(b)** Find the pure-strategy subgame-perfect equilibria of the game of part (a) for every possible value of $p \in (0,1)$.

## 14.E.2. Exercises for Section 14.2: Multi-sided incomplete information.

The answers to the following exercises are in Appendix S at the end of this chapter.

**Exercise 14.6**. Consider the following situation of two-sided incomplete information. Players 1 and 2 are having dinner with friends and during dinner Player 2 insinuates that Player 1 is guilty of unethical behavior. Player 1 can either demand an apology (*D*) or ignore (*I*) Player 2's remark; if Player 1 demands an apology, then Player 2 can either apologize (*A*) or refuse to apologize (*not-A*); if Player 2 refuses to apologize, Player 1 can either concede (*C*) or start a fight (*F*). Thus the extensive form is as follows (where $z_1,...,z_4$ are the possible outcomes):

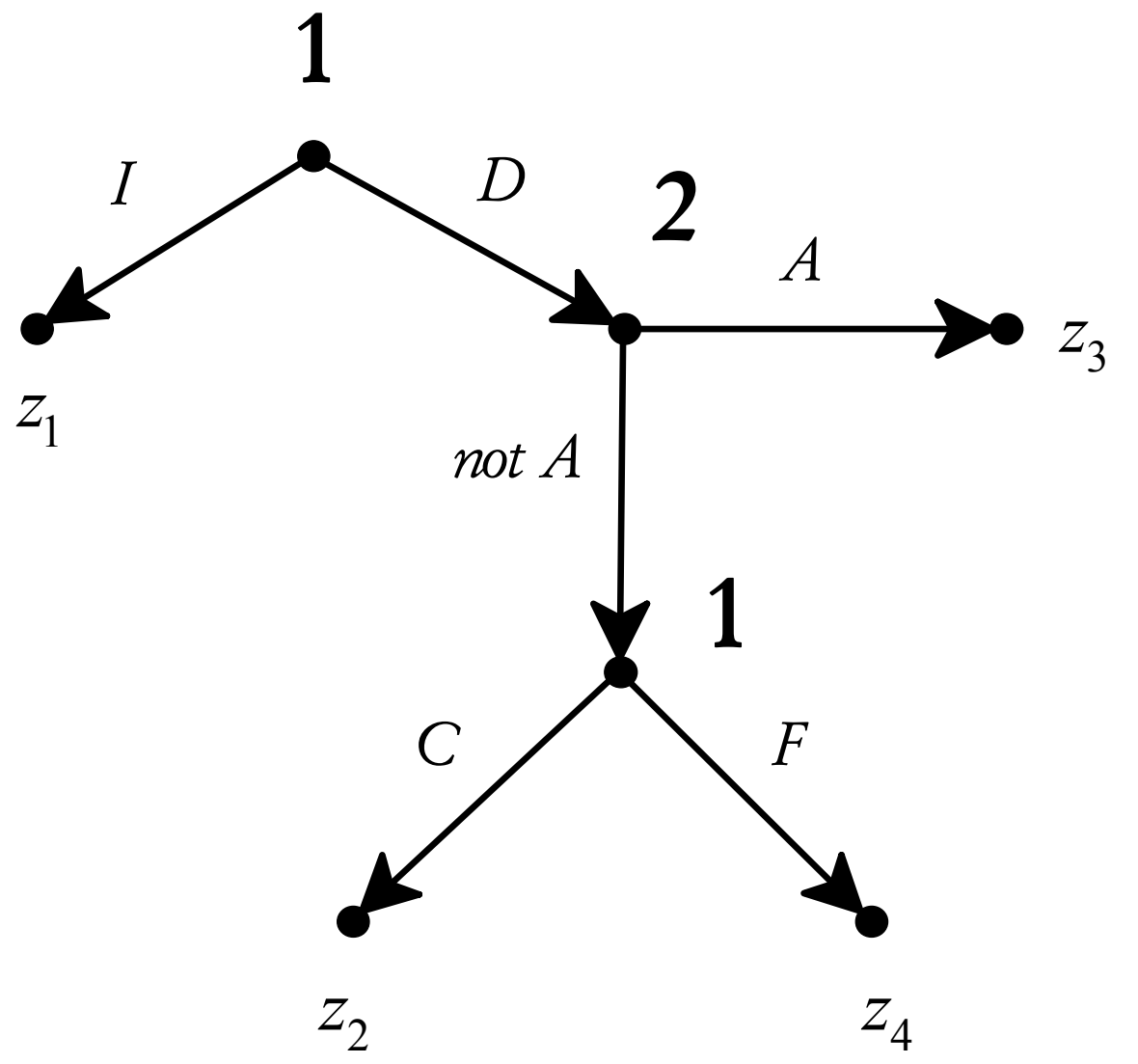


I= ignore, D = demand apology
A = apologize, not A = not apologize
C = concede, F = fight

Let $U_i$ be the von Neumann-Morgenstern utility function of Player *i* ($i$ = 1,2). The following is common knowledge between Players 1 and 2:

$$U_1(z_1) = 4,\ U_1(z_2) = 2, U_1(z_3) = 6$$

$$U_2(z_1) = 6,\ U_2(z_2) = 8, U_2(z_3) = 2,\ U_2(z_4) = 0\,.$$

The only uncertainty concerns the value of $U_1(z_4)$. As a matter of fact, $U_1(z_4) = 0$; Player 1 knows this, but Player 2 does not. It is common knowledge between Players 1 and 2 that Player 2 thinks that either $U_1(z_4) = 0$ or $U_1(z_4) = 3$. Furthermore, as a matter of fact, Player 2 assigns probability $\frac{1}{4}$ to $U_1(z_4) = 0$ and probability $\frac{3}{4}$ to $U_1(z_4) = 3$; Player 1 is uncertain as to whether Player 2's beliefs are $\begin{pmatrix} U_1(z_4)=0 & U_1(z_4)=3 \\ \frac{1}{4} & \frac{3}{4} \end{pmatrix}$ or $\begin{pmatrix} U_1(z_4)=0 \\ 1 \end{pmatrix}$ and he attaches probability $\frac{1}{2}$ to each.

**(a)** Construct a three-state interactive knowledge-belief structure that captures all of the above.

**(b)** Draw the extensive-form game that results from applying the Harsanyi transformation to the situation of part (a).

**(c)** Find the pure-strategy weak sequential equilibria of the game of part (b).

◊◊◊◊◊◊◊◊◊◊◊◊◊

**Exercise 14.7: Challenging Question.**

There are two parties to a potential lawsuit: the owner of a chemical plant and a supplier of safety equipment. The chemical plant owner, from now on called the *plaintiff*, alleges that the supplier, from now on called the *defendant,* was negligent in providing the safety equipment. The defendant knows whether or not he was negligent, while the plaintiff does not know; the plaintiff believes that there was negligence with probability $q$. These beliefs are common knowledge between the parties. The plaintiff has to decide whether or not to sue. If she does not sue then nothing happens and both parties get a payoff of 0. If the plaintiff sues then the defendant can either offer an out-of-court settlement of $\$S$ or resist. If the defendant offers a settlement, the plaintiff can either accept (in which case her payoff is $S$ and the defendant's payoff is $-S$) or go to trial. If the defendant resists then the plaintiff can either drop the case (in which case both parties get a payoff of 0) or go to trial. If the case goes to trial then legal costs are created in the amount of $\$P$ for the plaintiff and $\$D$ for the defendant. Furthermore (if the case goes to trial), the judge is able to determine if there was negligence and, if there was, requires the defendant to pay $\$W$ to the plaintiff (and each party has to pay its own legal costs), while if there was

no negligence the judge will drop the case without imposing any payments to either party (but each party still has to pay its own legal costs). It is common knowledge that each party is "selfish and greedy" (that is, only cares about its own wealth and prefers more money to less) and is risk neutral.

Assume the following about the parameters:

$$0 < q < 1, \quad 0 < D < S, \quad 0 < P < S < W - P.$$

**(a)** Represent this situation of incomplete information by means of an interactive knowledge-belief structure (the only two players are the plaintiff and the defendant).

**(b)** Apply the Harsanyi transformation to represent the situation in part (a) as an extensive-form game. [Don't forget to subtract the legal expenses from each party's payoff if the case goes to trial.]

**(c)** Write down all the strategies of the plaintiff.

**(d)** Prove that there is no pure-strategy weak sequential equilibrium which (1) is a separating equilibrium and (2) involves suing.

**(e)** For what values of the parameters ($q$, $S$, $P$, $W$, $D$) are there pure-strategy weak sequential equilibria which (1) are pooling equilibria and (2) involve suing? Consider all types of pooling equilibria and prove your claim.

**(f)** For the case where $q = \frac{1}{12}$, $P = 70$, $S = 80$, $W = 100$ find all the pure-strategy weak sequential equilibria which (1) are pooling equilibria and (2) involve suing. [Note that here we have dropped the assumption that $S < W - P$.]

# Appendix 14.S: Solutions to exercises

**Exercise 14.1.** **(a)** Let $G_1$ be the game with outcome $z_5$ and $G_2$ the game with outcome $z_6$. Then the structure is as follows, with $q = \frac{1}{4}$. The true state is *a*.

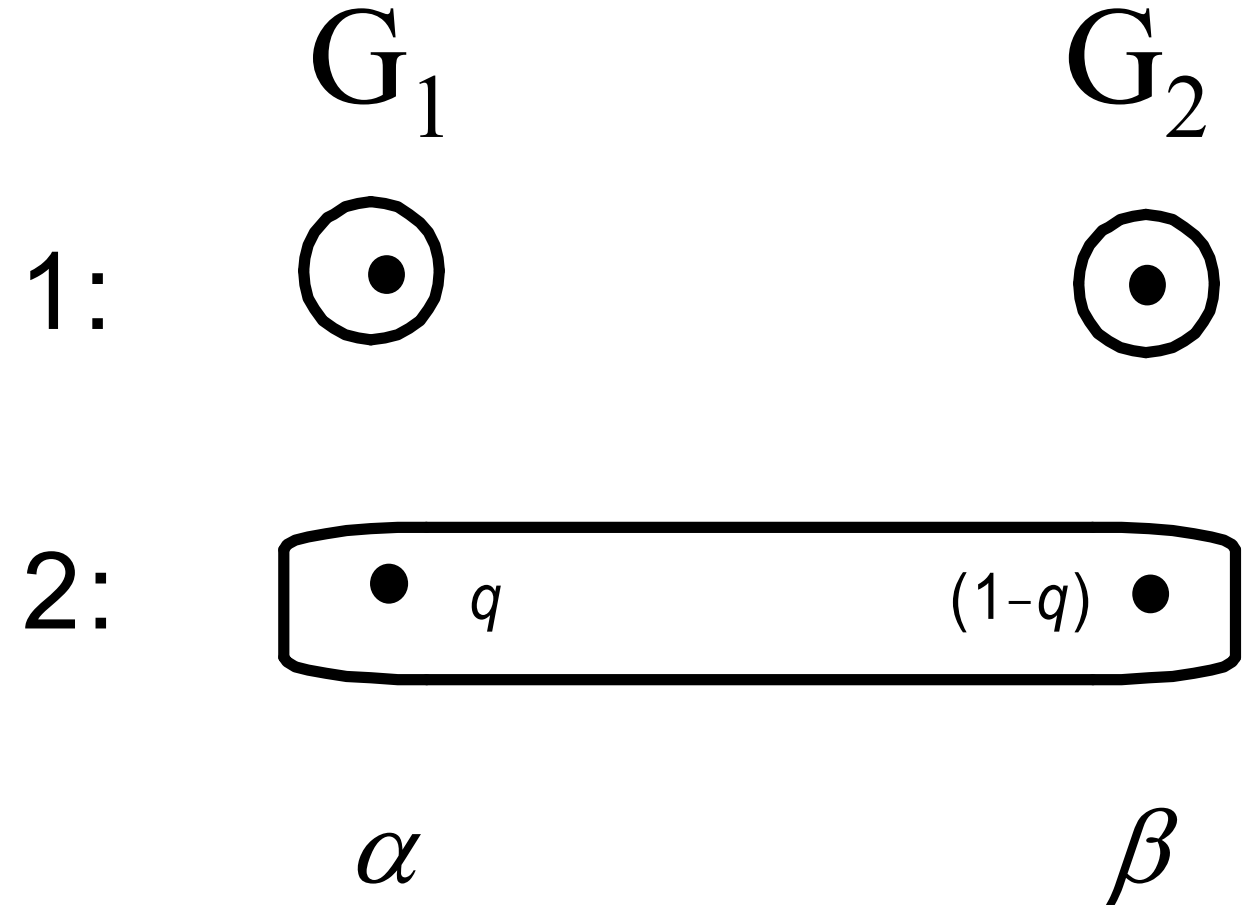


**(b)** The extensive form is as follows:

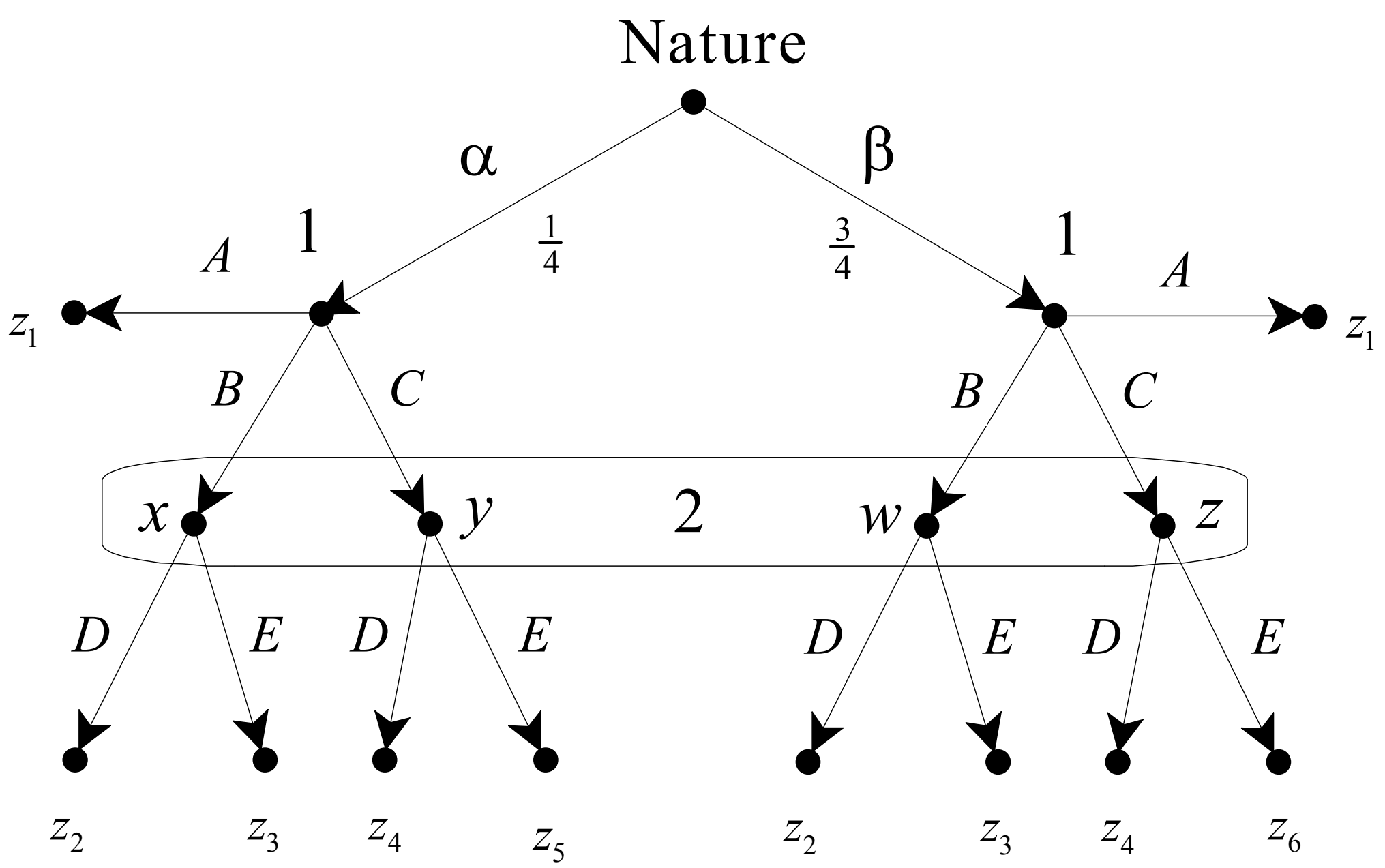

**(c)** $\begin{pmatrix} & z_1 & z_2 & z_3 & z_4 & z_5 & z_6 \\ U_1: & 0.5 & 0 & 0 & 1 & 0 & 1 \\ U_2: & 0 & 0.25 & 0 & 0.5 & 0.25 & 1 \end{pmatrix}$.

Adding these payoffs to the extensive form we obtain the following game:

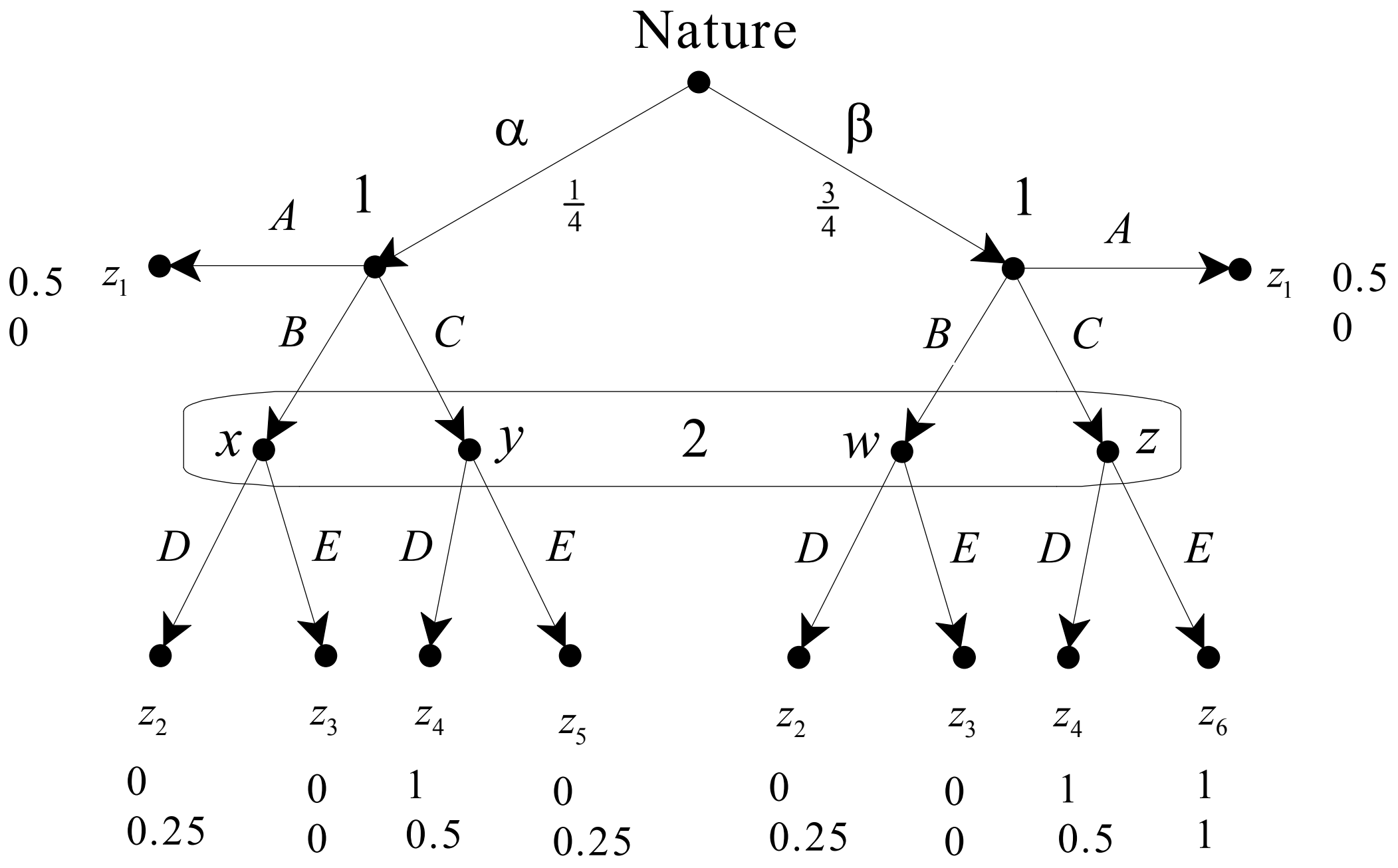


**(d)** No, because at the right node of Player 1, *C* gives a payoff of 1 no matter what Player 2 does and thus Player 1 would choose *C* rather than *A* (which only gives him a payoff of 0.5).

**(e)** If Player 1 always plays *C* (that is, his strategy is *CC*), then – by Bayesian updating – Player 2 should assign probability $\frac{1}{4}$ to node *y* and probability $\frac{3}{4}$ to node *z*, in which case *D* gives her a payoff of 0.5 while *E* gives her a payoff of $\frac{1}{4}(0.25)+\frac{3}{4}(1)=0.8125$; hence she must choose *E*. But then at his left node Player 1 with *C* gets 0 while with *A* he gets 0.5. Hence choosing *C* at the left node is not sequentially rational.

**(e)** Yes, the following is a weak sequential equilibrium: $\sigma=(AC,E)$ with beliefs $\mu=\begin{pmatrix} x & y & w & z \\ 0 & 0 & 0 & 1 \end{pmatrix}$. It is straightforward to verify that sequential rationality and Bayesian updating are satisfied.

**Exercise 14.2.** Sequential rationality was verified in Section 14.1. Thus we only need to show that the highlighted pure-strategy profile together with the following system of beliefs $\mu = \left( \begin{array}{cc|cc|cc} s & t & u & v & x & w \\ p & 1-p & 0 & 1 & 1 & 0 \end{array} \right)$ (for the names of the nodes refer to the following figure) constitutes a consistent assessment (Definition 11.1, Chapter 11). Let $\langle \sigma_n \rangle_{n=1,2,...}$ be the sequence of completely mixed strategies shown in the following figure.

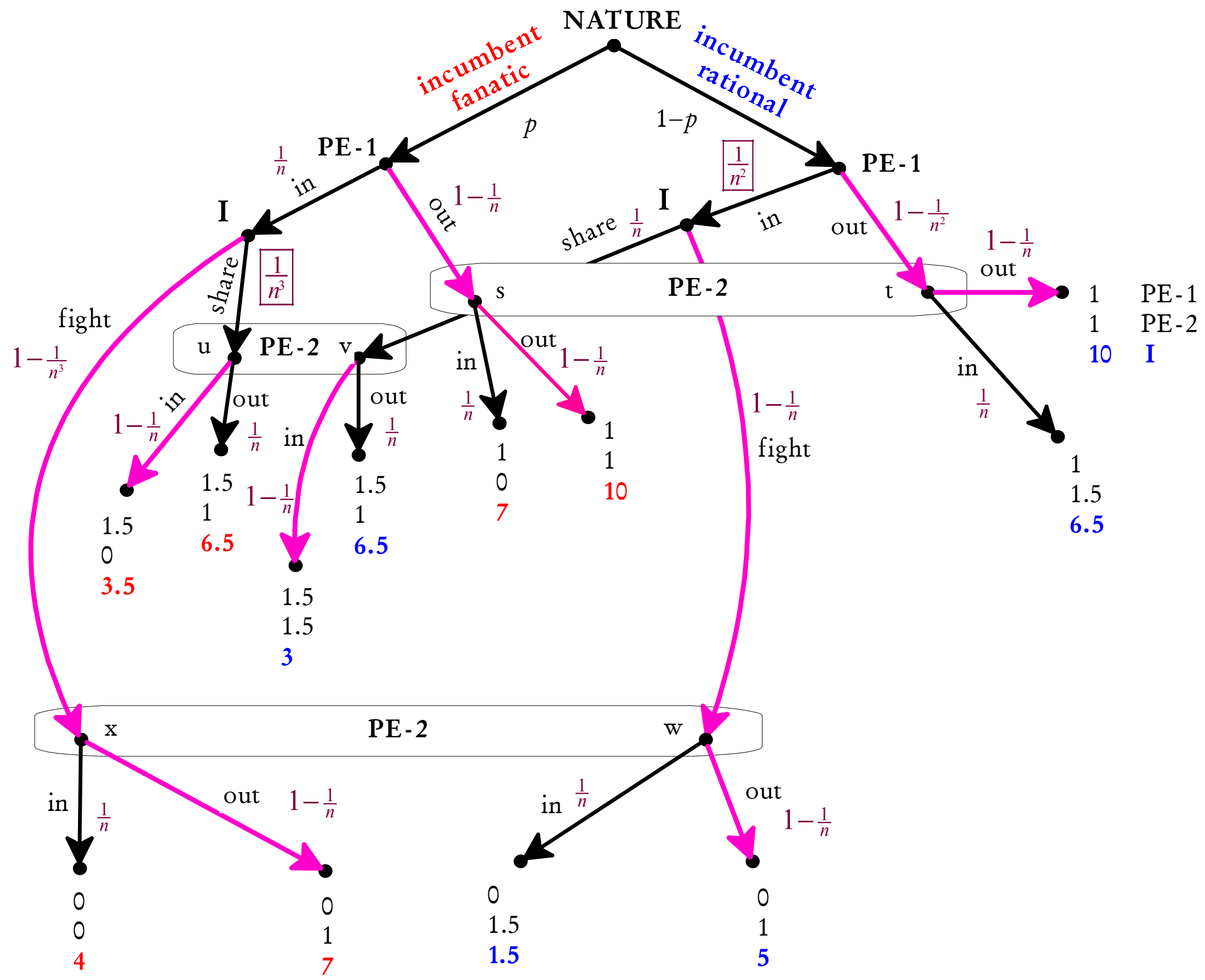


It is clear that $\lim_{n \to \infty} \sigma_n = \sigma$. Let us calculate the system of beliefs $\mu_n$ obtained from $\sigma_n$ by using Bayesian updating: $\mu_n(s) = \dfrac{p\left(1-\frac{1}{n}\right)}{p\left(1-\frac{1}{n}\right)+(1-p)\left(1-\frac{1}{n^2}\right)}$ and

$\lim_{n\to\infty} \mu_n(s) = \frac{p}{p+(1-p)} = p$;  $\mu_n(u) = \frac{p\left(\frac{1}{n}\right)\left(\frac{1}{n^3}\right)}{p\left(\frac{1}{n}\right)\left(\frac{1}{n^3}\right)+(1-p)\left(\frac{1}{n^2}\right)\left(\frac{1}{n}\right)} = \frac{p}{p+n(1-p)}$

and  $\lim_{n\to\infty} \mu_n(u) = 0$;  $\mu_n(x) = \frac{p\left(\frac{1}{n}\right)\left(1-\frac{1}{n^3}\right)}{p\left(\frac{1}{n}\right)\left(1-\frac{1}{n^3}\right)+(1-p)\left(\frac{1}{n^2}\right)\left(1-\frac{1}{n}\right)} =$

$\frac{p\left(1-\frac{1}{n^3}\right)}{p\left(1-\frac{1}{n^3}\right)+(1-p)\left(\frac{1}{n}\right)\left(1-\frac{1}{n}\right)}$ and $\lim_{n\to\infty} \mu_n(x) = 1$ (the numerator tends to $p(1) = p$ and the denominator tends to $p(1)+(1-p)(0)(1) = p$). Thus $(\sigma,\mu)$ is consistent and sequentially rational and therefore it is a sequential equilibrium.

**Exercise 14.3**. **(a)** Consider the following pure-strategy profile, call it $\sigma$, which is highlighted in dark pink in the figure below:

(1) PE-1's strategy is "out" at both nodes,

(2) PE-2's strategy is "in" at the top information set (after having observed that PE-1 stayed out), "in" at the middle information set (after having observed that PE-1's entry was followed by the Incumbent sharing the market with PE-1) and "out" at the bottom information set (after having observed that PE-1's entry was followed by the Incumbent fighting against PE-1),

(3) the Incumbent's strategy is to fight entry of PE-1 in any case (that is, whether the Incumbent himself is hotheaded or rational).

We want to show that $\sigma$, together with the system of beliefs $\mu$ shown in the Figure below (where at her top information PE-2 assigns probability $p$ to the Incumbent being hotheaded, at her middle information set she assigns probability 1 to the Incumbent being rational and at her bottom information set she assigns probability 1 to the Incumbent being hotheaded) constitutes a weak sequential equilibrium for any value of $p \leq \frac{1}{3}$. Bayesian updating is satisfied at PE-2's top information set (the only non-singleton information set reached by $\sigma$), and at the other two information sets of PE-2 any beliefs are allowed by the notion of weak sequential equilibrium. Sequential rationality also holds:

- For PE-1, at the left node "in" yields 0 and "out" yields 1, so that "out" is sequentially rational; the same is true at the right node.

- For the Incumbent, at the left node "fight" yields 7 and "share" yields 3.5, so that "fight" is sequentially rational; at the right node "fight" yields 5 and "share" yields 3, so that "fight" is sequentially rational.

- For PE-2, at the top information set "in" yields an expected payoff of $p(0)+(1-p)(1.5)=(1-p)(1.5)$ and "out" yields 1, so that "in" is sequentially rational as long as $p \le \frac{1}{3}$; at the middle information set "in" yields 1.5 and "out" yields 1, so that "in" is sequentially rational; at the bottom information set "in" yields 0 and "out" yields 1, so that "out" is sequentially rational.

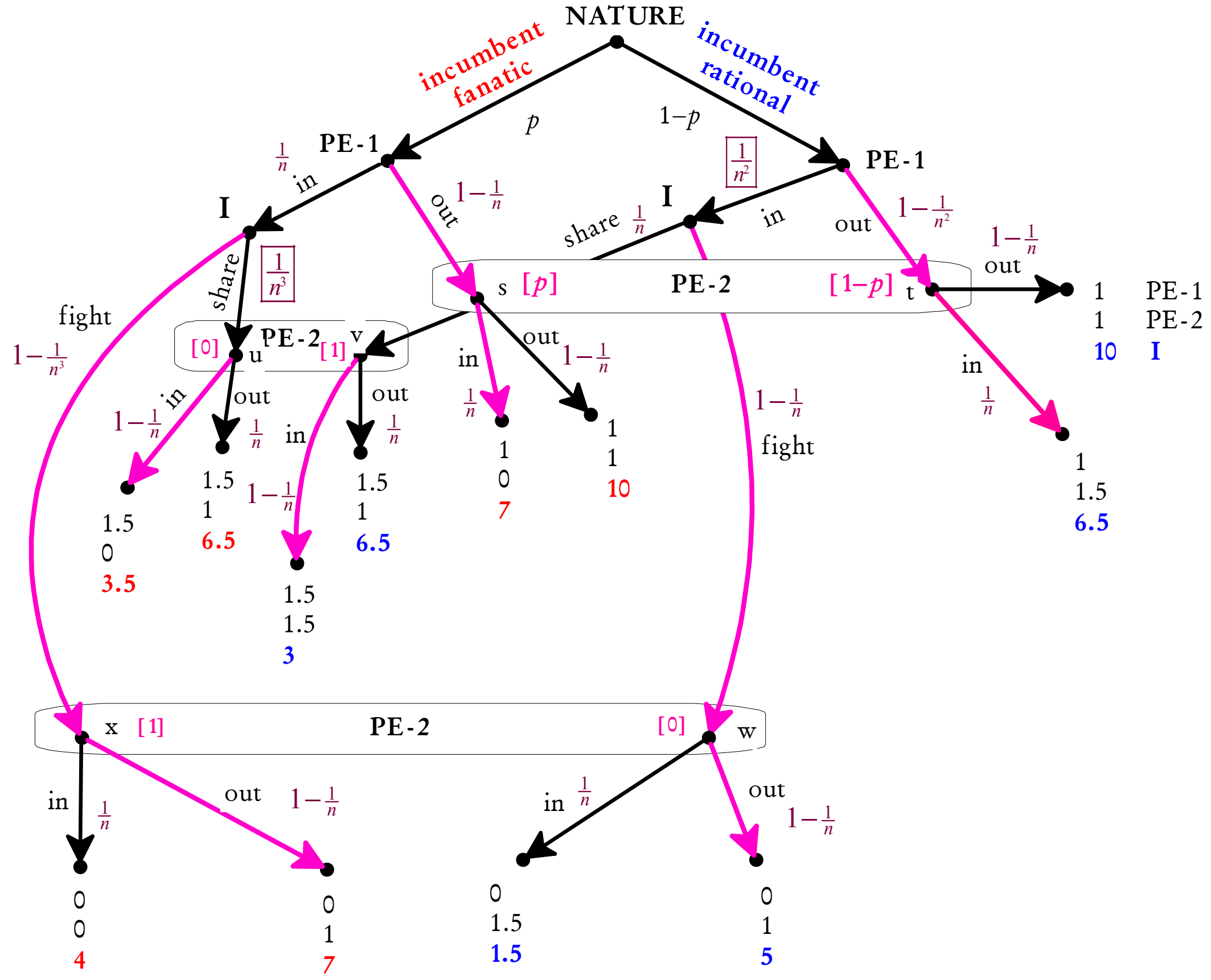


**(b)** The assessment described in part (a) is in fact a sequential equilibrium. This can be shown using the sequence of completely mixed strategies marked in the figure above, which coincides with the sequences of mixed strategies considered in Exercise 14.2; thus the calculations to show consistency are identical to the ones carried out in Exercise 14.2.

**Exercise 14.4.** The game under consideration is the following (which reproduces Figure 14.15):

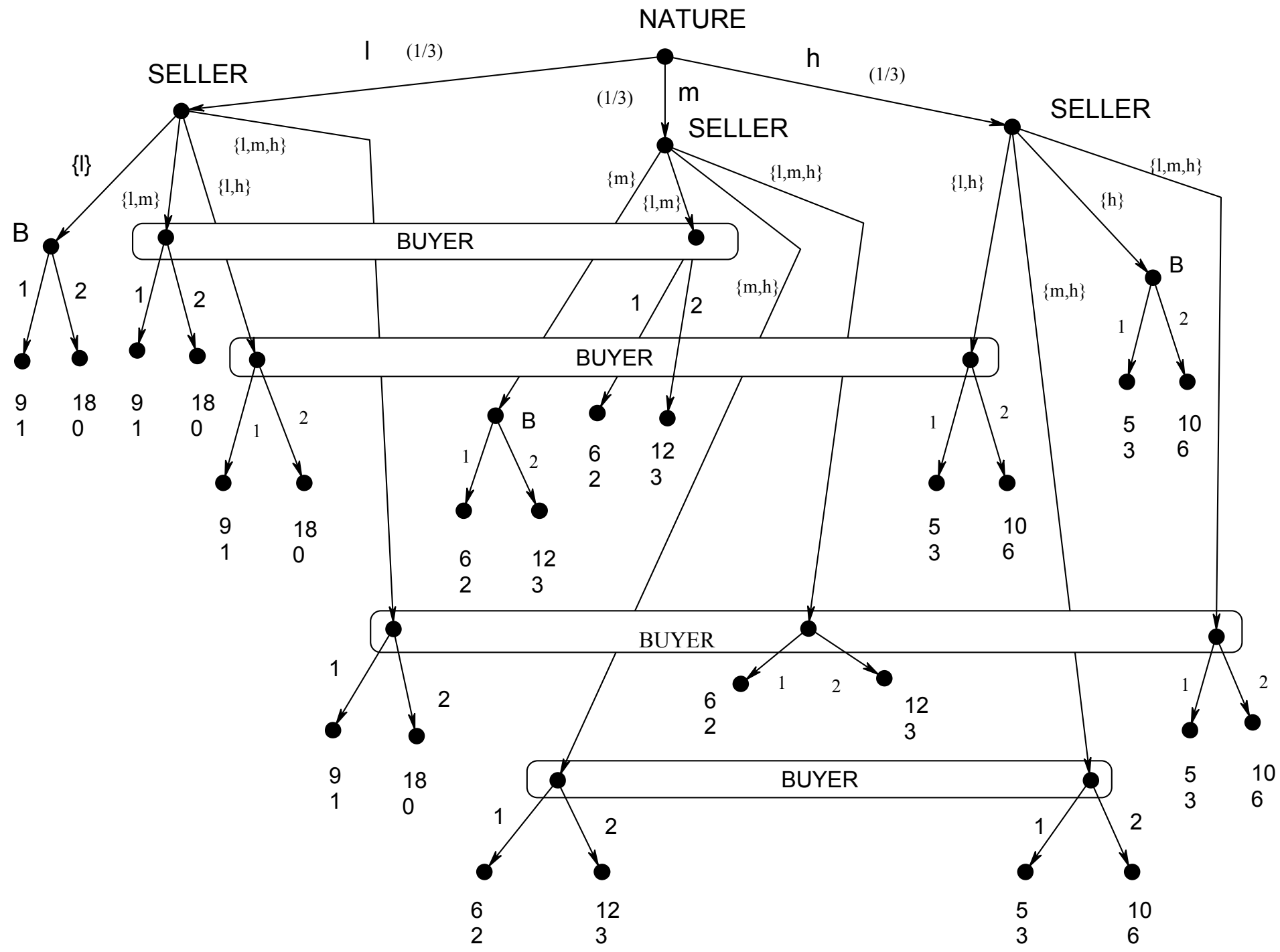


Let the buyer be naïve, in the sense that at every unreached information set she assigns equal probability to each node; furthermore, let the buyer's strategy be as follows: if told $\{l\}$ buy one unit, in every other case buy two units. Let the seller's strategy be as follows: if Nature chooses $l$ or $m$ claim $\{l,m\}$ and if Nature chooses $h$ then claim $\{h\}$. Let us verify that this assessment constitutes a weak sequential equilibrium. The only non-singleton information set that is reached by the strategy profile is the top information set (where the buyer hears the claim $\{l,m\}$) and Bayesian updating requires the buyer to assign equal probability to each node. Every other non-singleton information set is not reached and thus the notion of weak sequential equilibrium allows any beliefs: in particular beliefs that assign probability $\frac{1}{2}$ to each node. Now we check sequential rationality. Let us start with the buyer. At the singleton node following claim $\{l\}$, buying one unit is optimal and at the singleton nodes following claims $\{m\}$ and $\{h\}$ buying two units is optimal. At the top information set (after the claim $\{l,m\}$) choosing one unit gives an expected payoff of $\frac{1}{2}(1)+\frac{1}{2}(2)=1.5$ and choosing two units yields $\frac{1}{2}(0)+\frac{1}{2}(3)=1.5$, thus buying two units is sequentially rational. At the information set following claim $\{l,h\}$ choosing one unit gives an expected payoff of $\frac{1}{2}(1)+\frac{1}{2}(3)=2$ while

choosing two units yields $\frac{1}{2}(0)+\frac{1}{2}(6)=3$, thus buying two units is sequentially rational. At the information set following claim {*l,m,h*} choosing one unit gives an expected payoff of $\frac{1}{3}(1)+\frac{1}{3}(2)+\frac{1}{3}(3)=2$ while choosing two units yields $\frac{1}{3}(0)+\frac{1}{3}(3)+\frac{1}{3}(6)=3$, thus buying two units is sequentially rational. Finally, at the information set following claim {*m,h*} choosing one unit gives an expected payoff of $\frac{1}{2}(2)+\frac{1}{2}(3)=2.5$ while choosing two units yields $\frac{1}{2}(3)+\frac{1}{2}(6)=4.5$, thus buying two units is sequentially rational. Now let us check sequential rationality for the seller's strategy. At the left node (after Nature chooses *l*) claiming {*l*} yields a payoff of 9, while every other claim yields a payoff of 18. Thus {*l,m*} is sequentially rational. At the other two nodes, the seller is indifferent among all the claims, because they yield the same payoff (12 at the node after Nature chooses *m* and 10 at the node after Nature chooses *h*) thus the postulated choices are sequentially rational.

**Exercise 14.5.** **(a)** The extensive form is as follows:

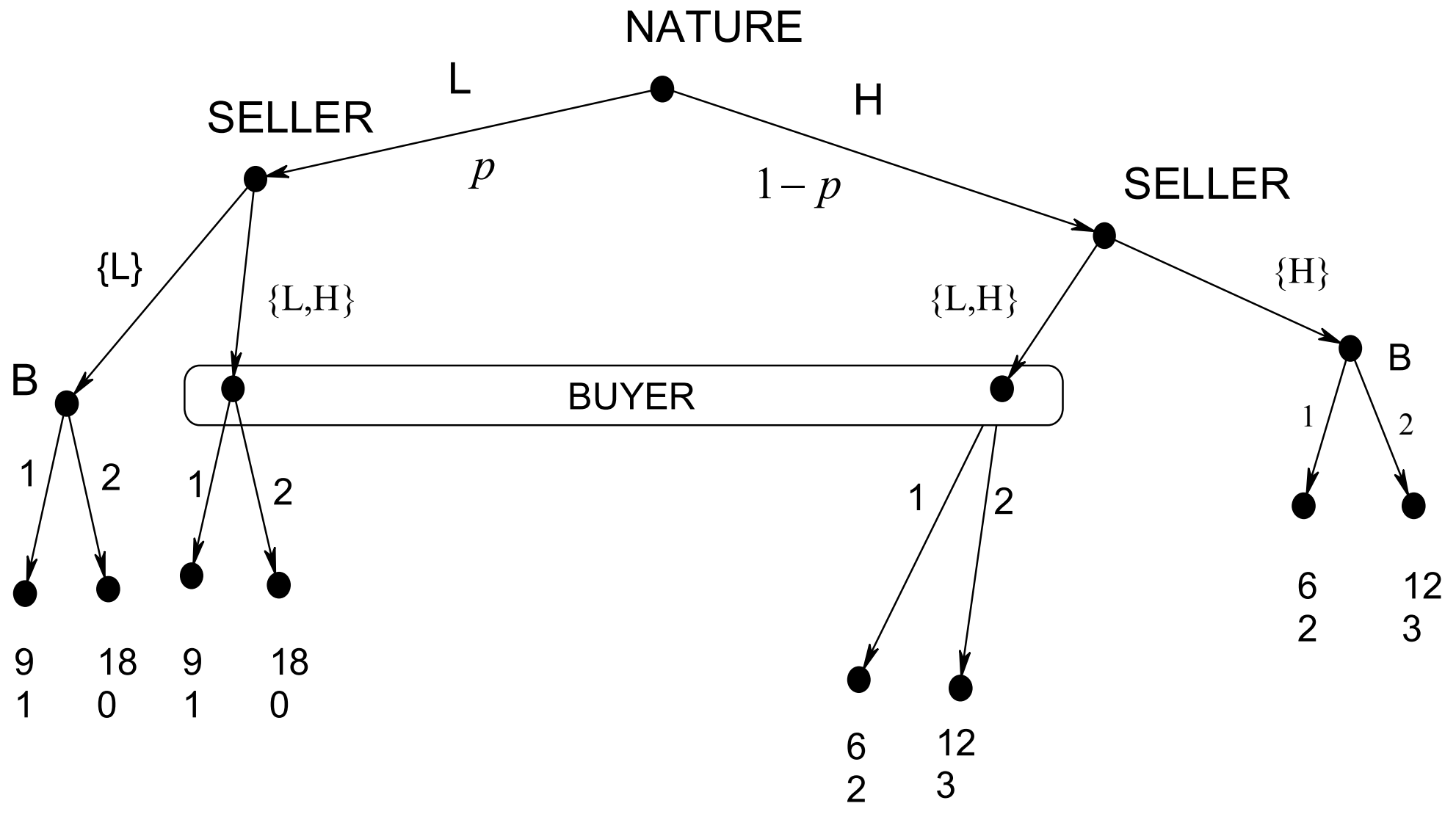


**(b)** At any subgame-perfect equilibrium, the buyer buys one unit at the singleton information set on the left and two units at the singleton information set on the right. Thus we can simplify the game by replacing the buyer's node on the left with the payoff vector (9,1) and the buyer's node on the right with the payoff vector (12,3). The strategic form corresponding to the simplified game is as follows (the buyer's strategy refers to the buyer's only information set that has remained, where he is told *LH*):

| | | **Buyer** | |
|---|---|---|---|
| | | Buy 1 unit | Buy 2 units |
| **Seller** | *L, LH* | 9p + 6(1–p), p + 2(1–p) | 18p + 12(1–p), 3(1–p) |
| | *L, H* | 9p + 12(1–p), p + 3(1–p) | 9p + 12(1–p), p + 3(1–p) |
| | *LH, LH* | 9p + 6(1–p), p + 2(1-p) | 18p + 12(1-p), 3(1–p) |
| | *LH,H* | 9p + 12(1-p), p + 3(1-p) | 18p + 12(1-p), 3(1–p) |

simplifying:

| | | Buy 1 unit | Buy 2 units |
|---|---|---|---|
| **Seller** | *L, LH* | 6 + 3p, 2 – p | 12 + 6p, 3 – 3p |
| | *L, H* | 12 – 3p, 3 – 2p | 12 – 3p, 3 – 2p |
| | *LH, LH* | 6 + 3p, 2 – p | 12 + 6p, 3 – 3p |
| | *LH,H* | 12 – 3p, 3 – 2p | 12 + 6p, 3 – 3p |

Since $p < 1$, $12{-}3p > 6 + 3p$. Thus ((*LH,H*),1) and ((*L,H*),1) are always Nash equilibria for every value or $p$. Note that $2-p > 3-3p$ if and only if $p > \frac{1}{2}$; thus if $p > \frac{1}{2}$ then there are no other pure-strategy Nash equilibria. On the other hand, if $p \leq \frac{1}{2}$ then the following are also Nash equilibria: ((*L,LH*),2) and ((*LH,LH*),2). To sum up, the pure-strategy Nash equilibria of this normal form, *which correspond to the pure-strategy subgame-perfect equilibria of the reduced extensive form*, are as follows:

- If $p > \frac{1}{2}$, ((*LH,H*),1) and ((*L,H*),1);
- If $p \leq \frac{1}{2}$, ((*LH,H*),1), ((*L,H*),1), ((*L,LH*),2) and ((*LH,LH*),2).

**Exercise 14.6.** **(a)** Let $G_1$ and $G_2$ be the following perfect-information games (whose backward induction solution has been highlighted in dark pink):

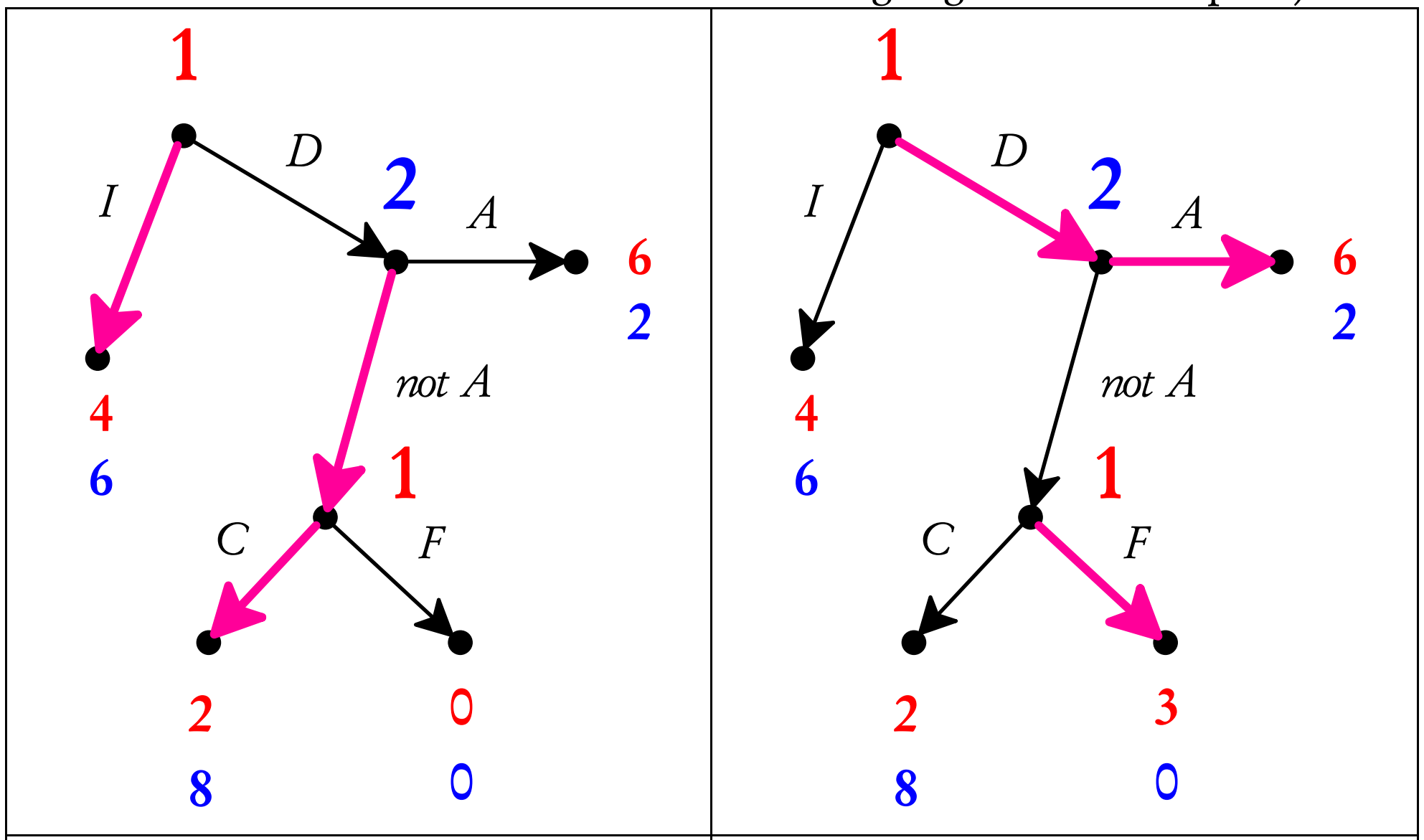


| Game $G_1$ | Game $G_2$ |
|---|---|

The interactive knowledge-belief structure is as follows:

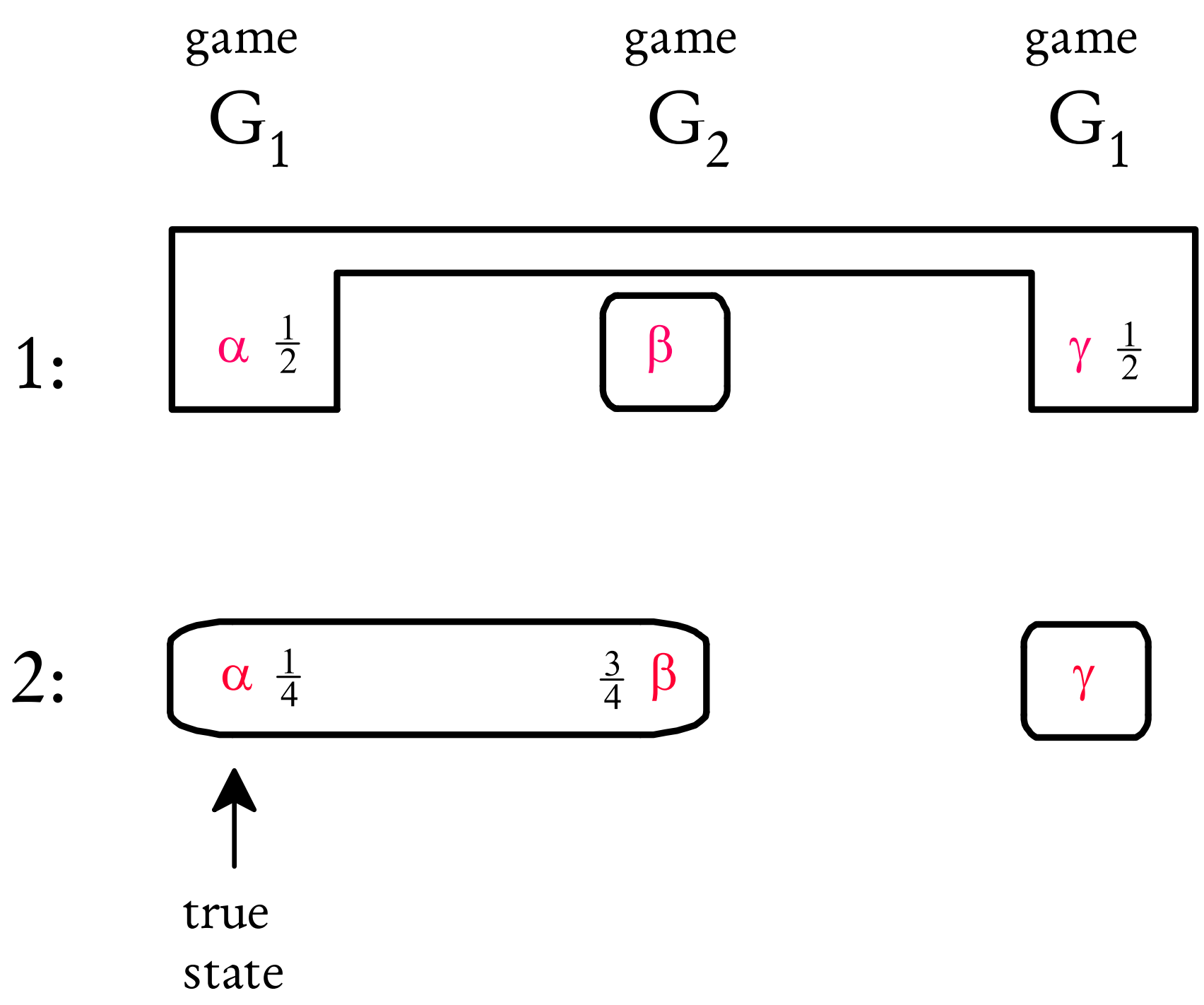

**(b)** The common prior is given by $\begin{pmatrix} \alpha & \beta & \gamma \\ \frac{1}{5} & \frac{3}{5} & \frac{1}{5} \end{pmatrix}$. The Harsanyi transformation yields the following game:

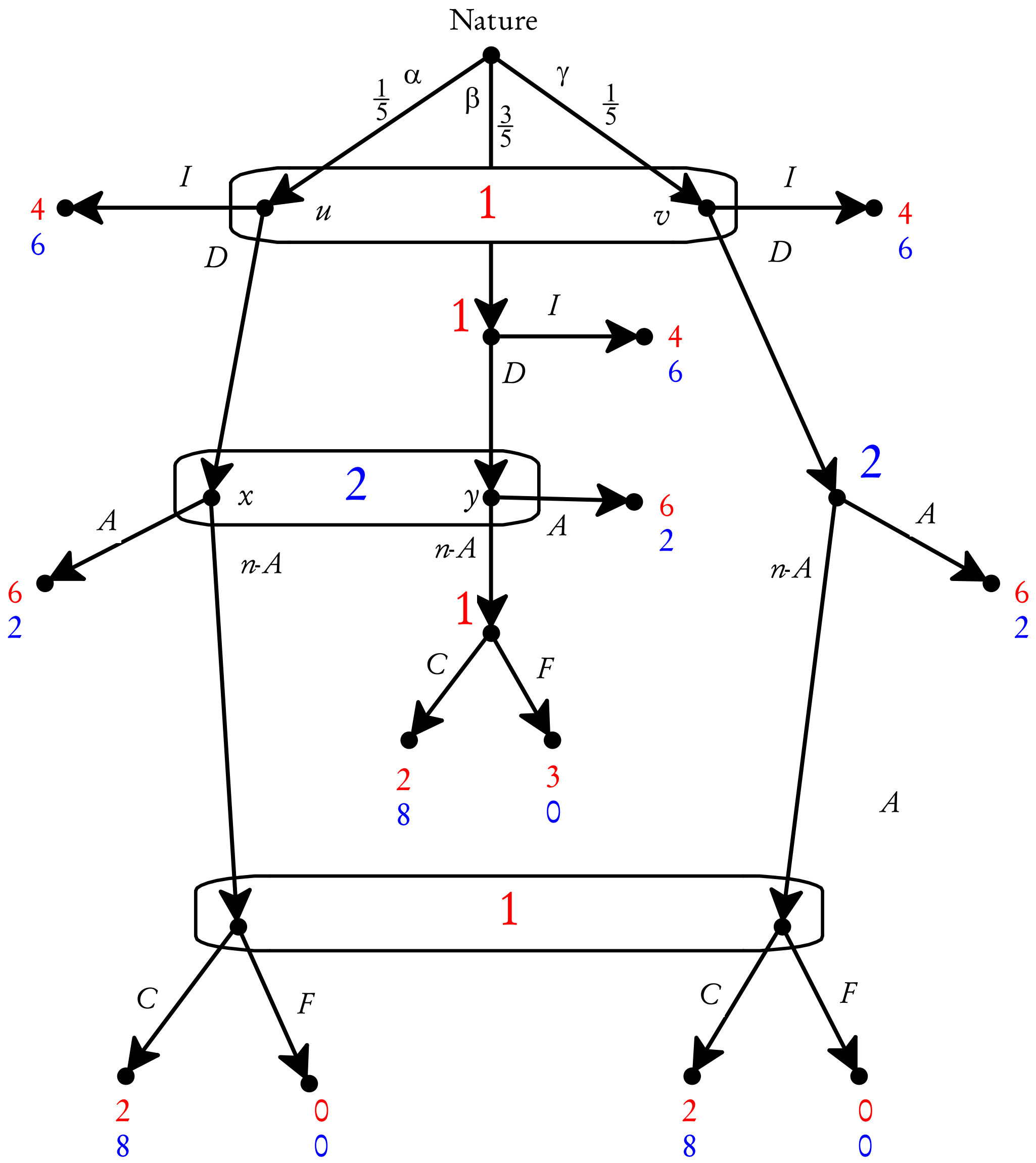


**(c)** At the bottom information set of Player 1, *C* strictly dominates *F* and thus we can replace the two nodes in that information set with the payoff vector (2,8). At the bottom singleton node of Player 1, *F* strictly dominates *C* and thus we can replace that node with the payoff vector (3,0). Thus, by appealing to sequential rationality, the game can be reduced to the following:

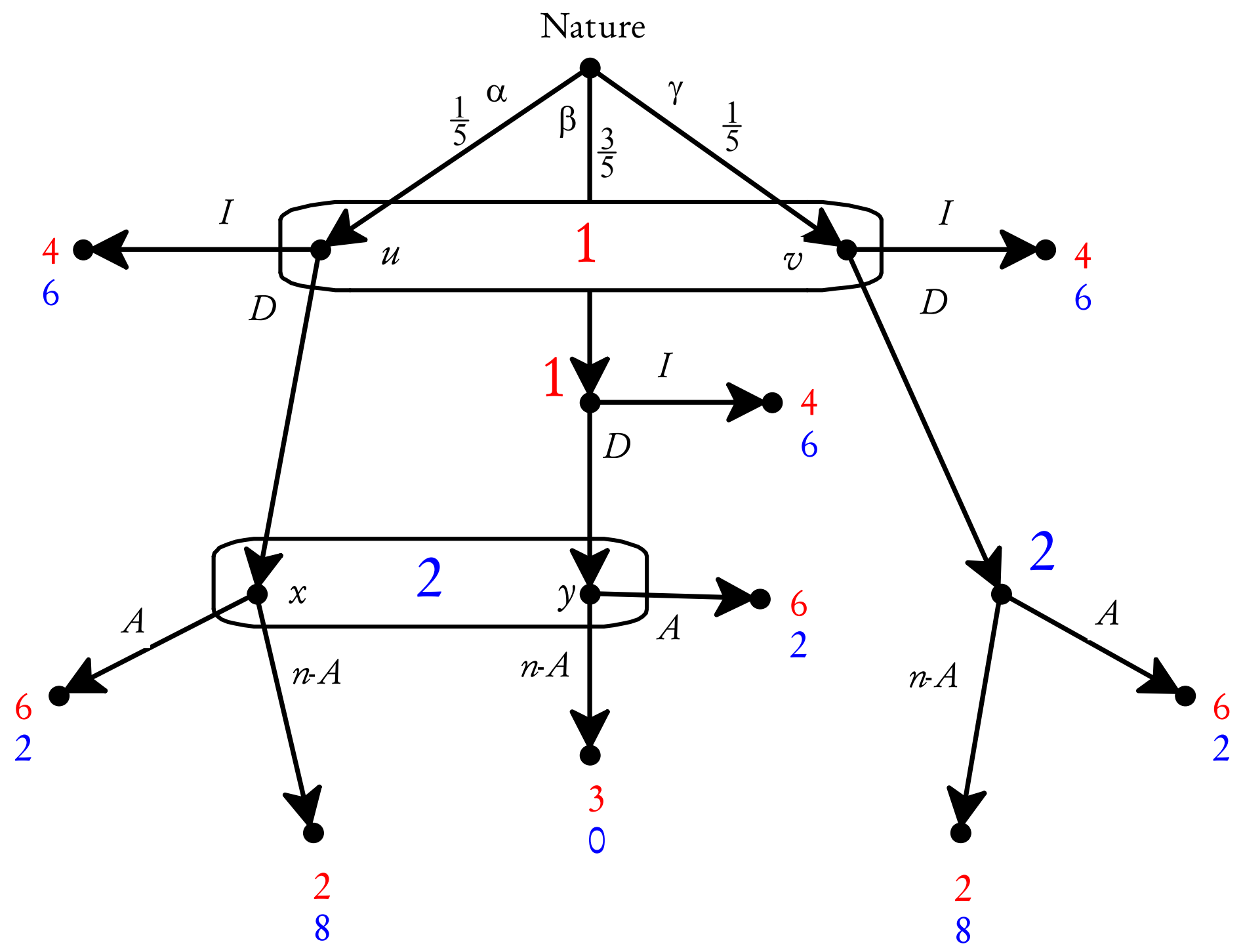


Now at the bottom-right node of Player 2, *n-A* strictly dominates *A* and thus we can replace that node with the payoff vector (2,8) and further simplify the game as follows:

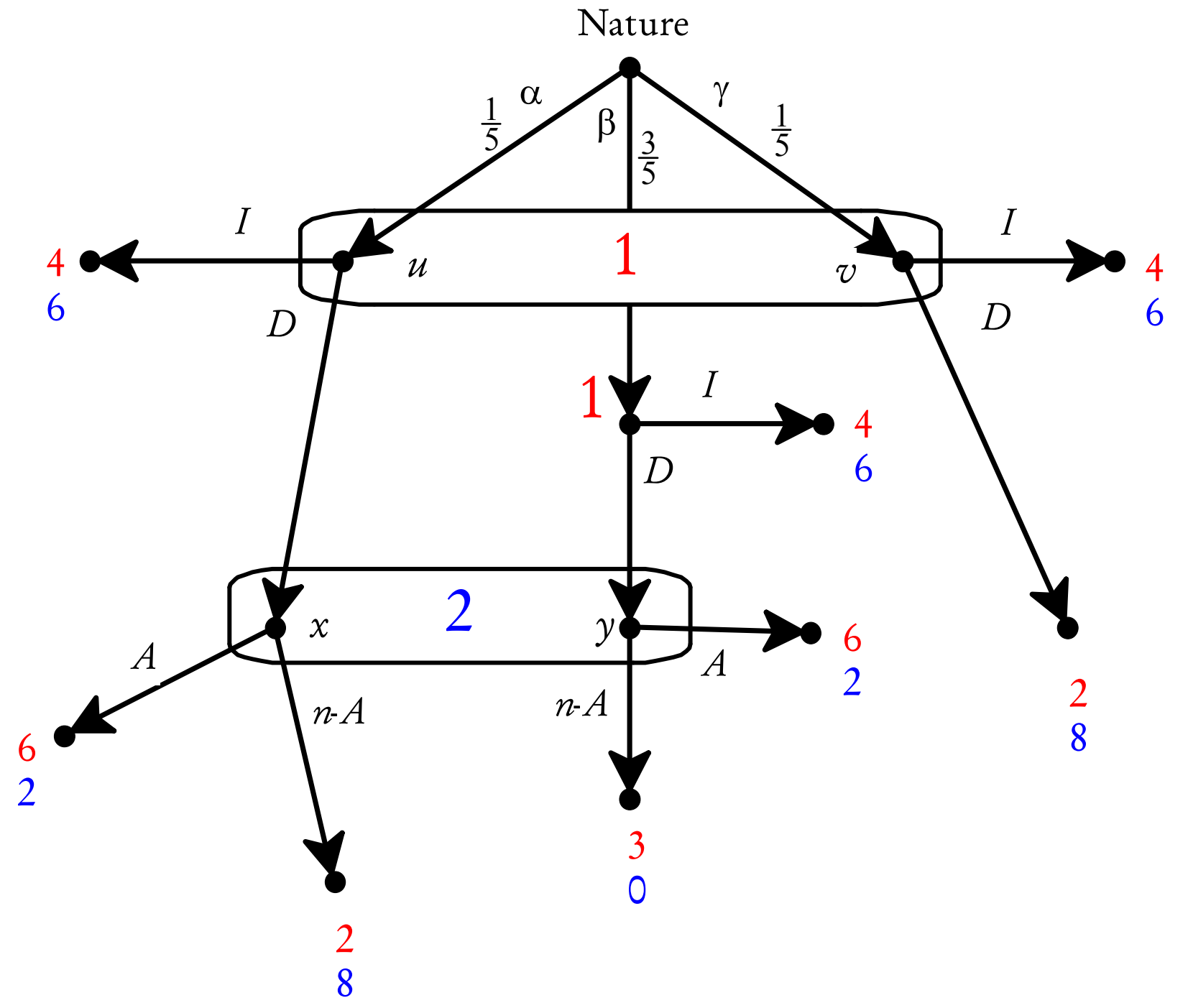

The following are the pure-strategy weak sequential equilibria of the reduced game (which can be extended into weak sequential equilibria of the original game by adding the choices that were selected during the simplification process):

- $((I,I), n\text{-}A)$ with beliefs $\mu = \left(\begin{array}{cc|cc} u & v & x & y \\ \frac{1}{2} & \frac{1}{2} & p & 1-p \end{array}\right)$ for any $p \geq \frac{1}{4}$.

- $((I,D), A)$ with beliefs $\mu = \left(\begin{array}{cc|cc} u & v & x & y \\ \frac{1}{2} & \frac{1}{2} & p & 1-p \end{array}\right)$ for any $p \leq \frac{1}{4}$.

- $((D,D), A)$ with beliefs $\mu = \left(\begin{array}{cc|cc} u & v & x & y \\ \frac{1}{2} & \frac{1}{2} & \frac{1}{4} & \frac{3}{4} \end{array}\right)$.

**Exercise 14.7 (Challenging Question).** **(a)** Let $G_1$ and $G_2$ be the following games (in $G_1$ the defendant is negligent and in $G_2$ he is not):

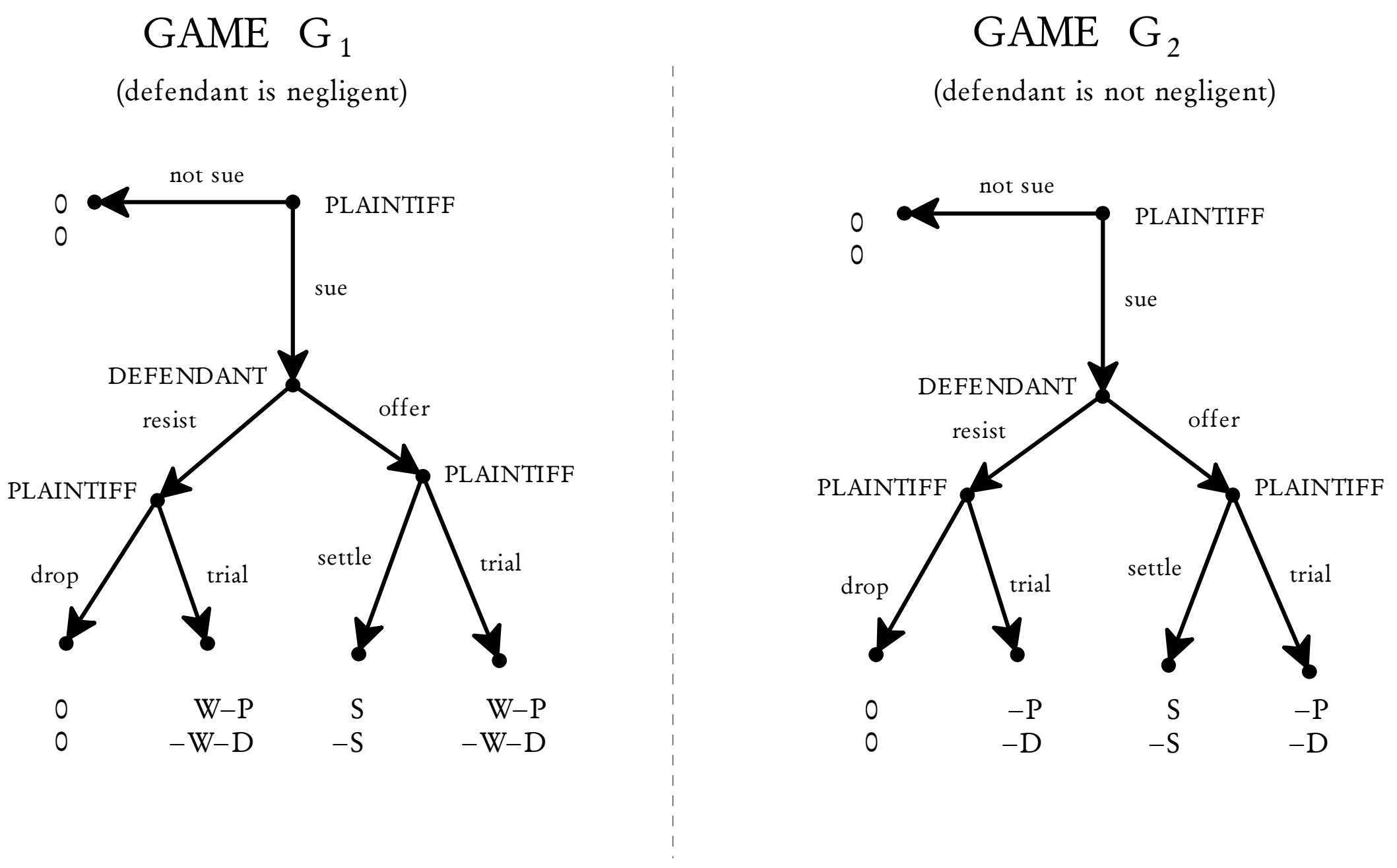


Then the situation can be represented as follows (at state $\alpha$ is the defendant is negligent and at state $\beta$ the defendant is not negligent):

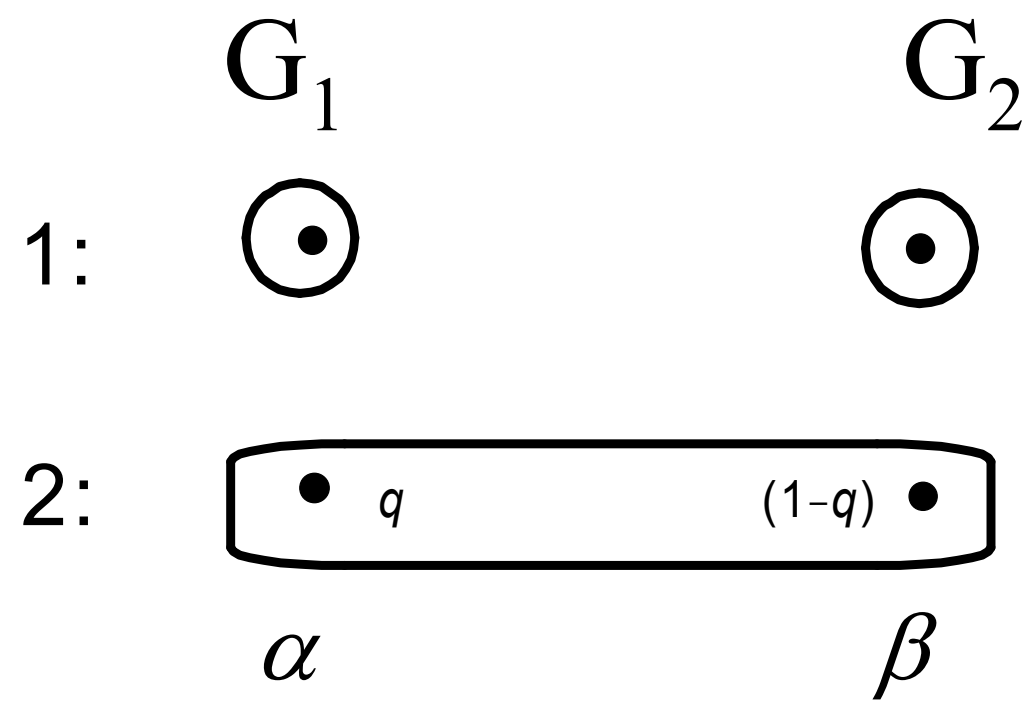


**(b)** The extensive-form game is as follows:

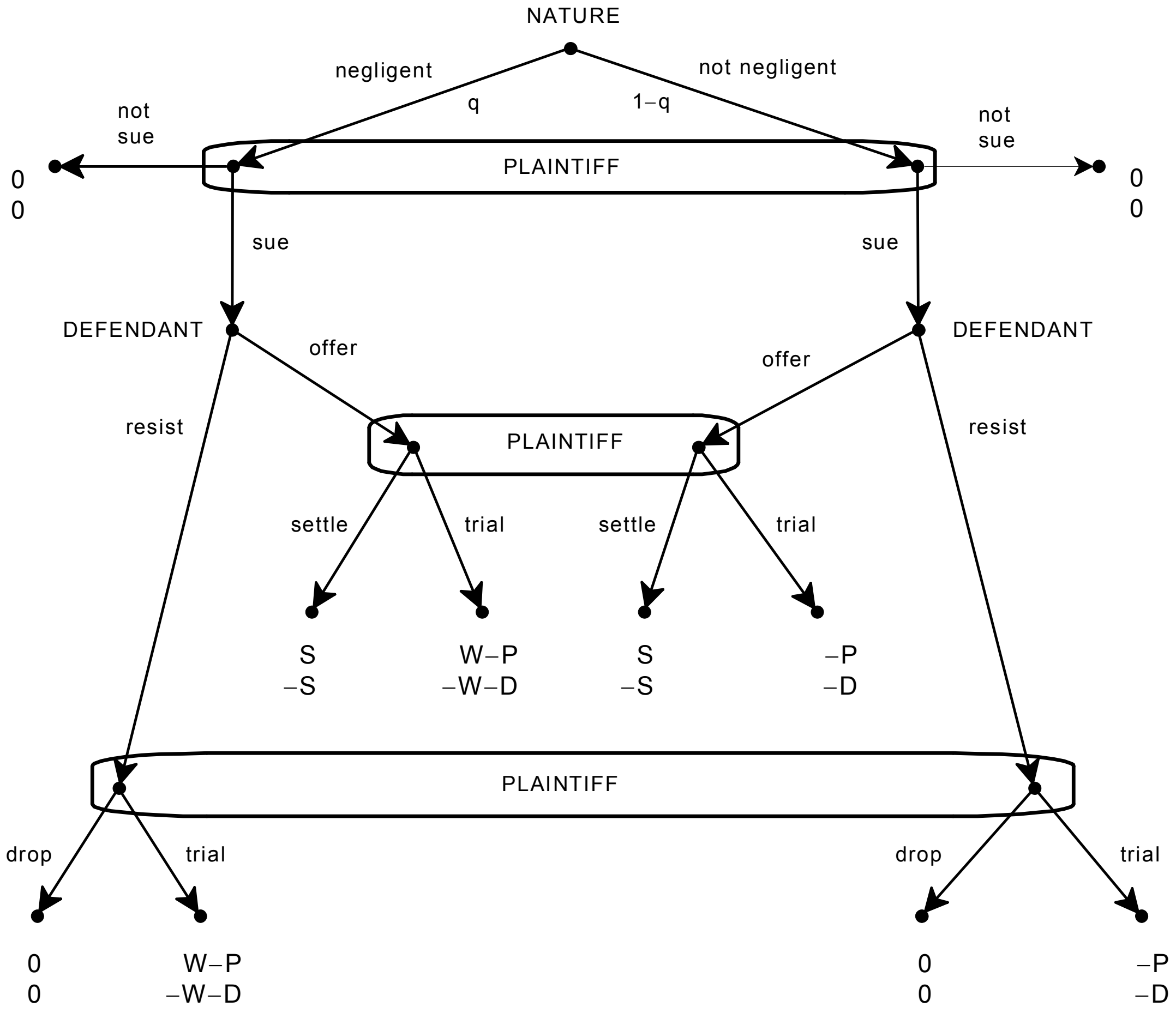

**(c)** Plaintiff's strategies: (1) (sue; if offer settle; if resist drop), (2) (sue; if offer settle; if resist go to trial), (3) (sue; if offer go to trial; if resist drop), (4) (sue; if offer go to trial; if resist go to trial), (5) (not sue; if offer settle; if resist drop), (6) (not sue; if offer settle; if resist go to trial), (7) (not sue; if offer go to trial; if resist drop), (8) (not sue; if offer go to trial; if resist go to trial).

**(d)** There are two possibilities for a separating equilibrium: (1) the defendant's strategy is "resist if negligent and offer if not negligent" and (2) the defendant's strategy is "offer if negligent and resist if not negligent". In both cases we assume that the plaintiff's strategy involves suing with probability 1.

**Consider Case 1 first.** By Bayesian updating, at the bottom information set the plaintiff must attach probability 1 to the negligent type and thus, by sequential rationality, must choose "trial" (because $W - P > 0$). Similarly, by Bayesian updating, at the middle information set the plaintiff must attach probability 1 to the non-negligent type and thus by sequential rationality must choose "settle". But then the negligent type of the defendant gets $-(W+D)$ by resisting and would get $-S$ by offering. Since, by assumption, $S < W\ (< W + D)$, the choice of resisting is not sequentially rational.

**Consider now Case 2**. By Bayesian updating, at the bottom information set the plaintiff must attach probability 1 to the non-negligent type and thus by sequential rationality must choose "drop". But then the negligent type of the defendant gets a negative payoff by offering, while he would get 0 by resisting. Hence the choice of offering is not sequentially rational.

**(e)** There are two candidates for a pure-strategy pooling equilibrium: (1) both types of the defendant choose "offer" and (2) both types of the defendant choose "resist".

**Consider Case 1 first** (both types of the defendant choose "offer"). In order for "offer" to be sequentially rational for the non-negligent type, it cannot be that the plaintiff's strategy involves "settle" at the middle information set (the non-negligent type would get either 0 or $-D$ by resisting and both payoffs are greater than $-S$) and/or "drop" at the bottom information set. That is, **it must be that the plaintiff chooses "trial" at both information sets.** By Bayesian updating, at the middle information set the plaintiff must attach probability $q$ to the negligent type and probability $(1-q)$ to the non-negligent type. Hence at the middle information set "trial" is sequentially rational if and only if $qW - P \geq S$, that is, $q \geq \frac{S+P}{W}$. In order for "trial" to be sequentially rational at the bottom information set, the plaintiff must attach

sufficiently high probability (namely $p \geq \frac{P}{W}$) to the negligent type. This is allowed by weak sequential equilibrium because the bottom information set is not reached. Finally, in order for "sue" to be sequentially rational it must be that $qW - P \geq 0$, that is, $q \geq \frac{P}{W}$, which is implied by $q \geq \frac{S+P}{W}$. Thus

| |
|---|
| there is a pooling equilibrium with ((sue,trial,trial),(offer,offer) ) if and only if $q \geq \dfrac{S+P}{W}$. |

**Now consider Case 2** (both types of the defendant choose "resist"). [Note: since in part (f) the restriction $S < W - P$ does not hold, we will carry out the analysis below at first without imposing the restriction.] If the plaintiff's strategy involves "drop" at the bottom information set, then it is indeed sequentially rational for both types of the defendant to choose "resist". Furthermore, "drop" is sequentially rational in this case if, and only if, $qW - P \leq 0$, that is, $q \leq \frac{P}{W}$. Then "sue" is also sequentially rational, since the Plaintiff's payoff is 0 no matter whether he sues or does not sue. Thus

| |
|---|
| there is a pooling equilibrium with ((sue,x,drop),(resist,resist)) if and only if $q \leq \dfrac{P}{W}$ and appropriate beliefs as follows (*p* is the probability on the left node of the unreached middle information set):<br>either $x$ = settle and any $p$ if $W \leq S+P$ or $p \leq \dfrac{S+P}{W}$ if $W > S+P$<br>or $x$ = trial and $p \geq \dfrac{S+P}{W}$, which requires $W \geq S+P$ (since $p \leq 1$). |

Since it is assumed that $W > S + P$, we can conclude that

| |
|---|
| ((sue,settle,drop),(resist,resist)) is an equilibrium if and only if $q \leq \dfrac{P}{W}$ with $p \leq \dfrac{S+P}{W}$ |
| ((sue,trial,drop),(resist,resist)) is an equilibrium if and only if $q \leq \dfrac{P}{W}$ with $p \geq \dfrac{S+P}{W}$. |

If, on he other hand, $q \geq \frac{P}{W}$, then "trial" is sequentially rational at the bottom information set. Then, in order for the non-negligent type of the defendant to choose "resist" it must be that the plaintiff's strategy involves "trial" also at the middle information set, for which we need him to assign probability $p \geq \frac{S+P}{W}$ to the negligent type (which is possible, since the middle information set is not reached); of course, this requires $W \geq S+P$. Thus,

| ((sue,trial,trial),(resist,resist)) is an equilibrium if and only if $q \geq \dfrac{P}{W}$ with $p \geq \dfrac{S+P}{W}$. |
|---|

**(f)** Note that here the restriction $W - P > S$ does **not** hold. In this case, $q < \frac{S+P}{W} = \frac{80+70}{100} = \frac{3}{2}$ and thus, by the previous analysis, there is no pooling equilibrium of type (1), that is, ((sue,trial,trial),(offer,offer)) is not an equilibrium.

Also, since $S+P > W$ there is no pooling equilibrium of type (2), that is, neither ((sue,trial,trial),(resist,resist)) nor ((sue,trial,drop),(resist,resist)) is an equilibrium.

However, there is a pooling equilibrium of type (2) with ((sue,settle,drop),(resist,resist)) with any beliefs at the middle information set, since "settle" strictly dominates "trial" there (and, of course, belief $q = \frac{1}{12}$ on the left node of the bottom information set).

Chapter

15

# Incomplete Information: The Type-Space Approach

## 15.1 Types of players

As noted in Chapter 13, the theory of "games of incomplete information"[7] was pioneered by John Harsanyi (1967-68). Harsanyi's approach was developed using a different approach from the one we employed in Chapters 13 and 14, which is based on the interactive knowledge-belief structures introduced in Chapter 8. The interactive knowledge structures of Chapter 8 are in fact a special case of the more general notion of *interactive Kripke structure,* named after the philosopher and logician Saul Kripke, whose work on this goes back to 1959 and was written while he was still an undergraduate.[8] Although well known among logicians and computer scientists, these structures were not known to game theorists. Perhaps, if Harasnyi had been aware of Kripke's work he might have developed his theory using those structures. We find the interactive Kripke structures more natural and elegant and thus prefer to explain the "theory of games of incomplete information" using those structures. In this chapter we will explain the "type-space" approach developed by Harsanyi and show that the two approaches are equivalent. We will limit ourselves to situations of incomplete information

[7] We use quotation marks because, strictly speaking, there is no such thing as a game of incomplete information. There are *situations* of incomplete information involving the playing of a game and what Harsanyi did was to suggest a way of transforming such situations into extensive-form games with imperfect information. Once the so-called "Harsanyi transformation" has been applied, the resulting game is a game of *complete* information. Thus the "theory of games of incomplete information" is a theory on how to represent a situation of incomplete information concerning the playing of a game and how to transform it into a dynamic game with imperfect, but complete, information.

[8] Kripke (1959, 1963). For further details the reader is referred to van Ditmarsch et al. (2015).

involving *static* games. We will begin in the next section with a simpler special case, which is often all that one finds in game-theory textbooks, and then explain the general case in Section 15.3.

## 15.2 Types that know their own payoffs

We take as a starting point the case of incomplete information concerning a strategic-form game-frame $\langle I,(S_1,..,S_n),O,f\rangle$.[9] In the special case considered in this section, any uncertainty Player $i$ has (if any) concerns either the von Neumann-Morgenstern utility function $U_j : O \to \mathbb{R}$ $(j \neq i)$ of another player (or several other players) or the beliefs of the other player(s) (or both). Within the approach of Chapter 13, we would represent such a situation with an interactive knowledge-belief structure,[10] by associating with every state $\omega$ a game based on the given game-frame.[11] Let $U_{i,\omega} : O \to \mathbb{R}$ be the utility function of Player $i$ at state $\omega$. The assumption that every player knows her own payoffs is captured by the following restriction:

$$\text{if } \omega,\omega' \in \Omega \text{ and } \omega' \in I_i(\omega) \text{ then } U_{i,\omega'} = U_{i,\omega}.$$

An example of such a situation of incomplete information is given in Figure 15.1, which reproduces Figure 13.3.[12]

9 Recall (Definition 1.1, Chapter 1) that the elements of $\langle I,(S_1,...,S_n),O,f\rangle$ are as follows: $I$ = {1,...,$n$} is a set of players, (2) $S_i$ is the set of *strategies* of Player $i \in I$ (and $S = S_1 \times ... \times S_n$ is the set of *strategy profiles*), $O$ is a set of *outcomes* and $f : S \to O$ is a function that associates with every strategy profile $s$ an outcome $f(s) \in O$.

10 Recall that an interactive knowledge-belief structure consists of a set of states $\Omega$, $n$ partitions $I_1,...,I_n$ of $\Omega$ and a collection of probability distributions on $\Omega$, one for each information set of each partition, whose support is a subset of that information set. Recall also that, for every state $\omega \in \Omega$, we denote by $I_i(\omega)$ the information set of Player $i$ (that is, the element of the partition $I_i$) that contains state $\omega$.

11 Recall that, given a game-frame $\langle I,(S_1,...,S_n),O,f\rangle$, a game based on it is obtained by specifying, for every Player $i$, a von Neumann-Morgenstern utility function $U_i : O \to \mathbb{R}$ on the set of outcomes.

12 In this case the game-frame is given by: $I$ = {1,2}, $S_1$ = {$T$,$B$}, $S_2$ = {$L$,$R$} (so that $S$ = {($T$,$L$), ($T$,$R$), ($B$,$L$), ($B$,$R$)}), $O = \{o_1, o_2, o_3, o_4\}$, $f((T,L)) = o_1$, $f((T,R)) = o_2$, $f((B,L)) = o_3$,

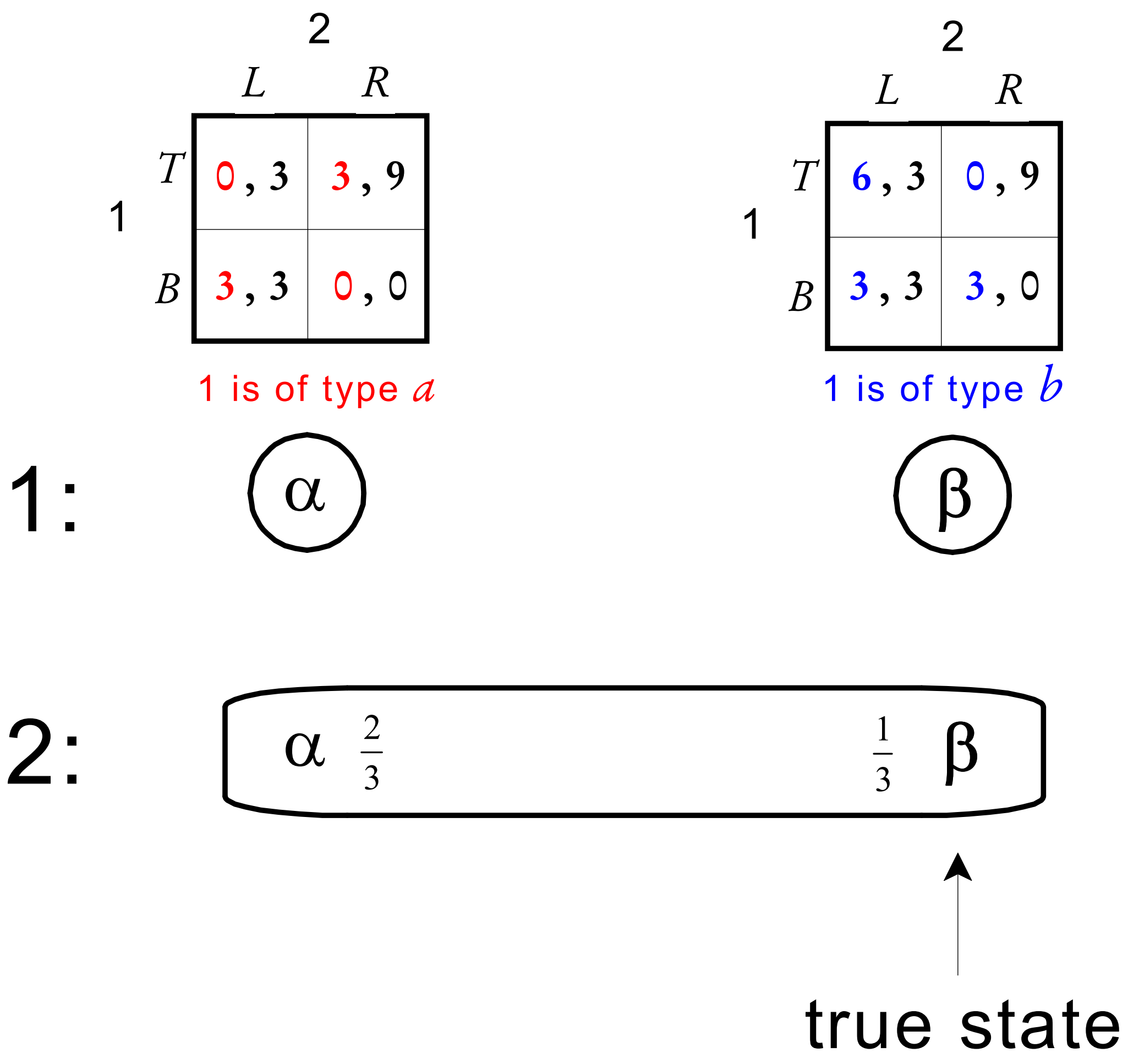


**Figure 15.1**

In the "type-space" approach the situation illustrated in Figure 15.1 would be represented using *types* of players: each type of Player $i$ represents a utility function of Player $i$ as well as Player $i$'s beliefs about the types of the other players. The formal definition is as follows.

---

$f((B,R)) = o_4$. The (state-dependent) utility functions are given by $U_{1,\alpha} = \begin{pmatrix} o_1 & o_2 & o_3 & o_4 \\ 0 & 3 & 3 & 0 \end{pmatrix}$, $U_{1,\beta} = \begin{pmatrix} o_1 & o_2 & o_3 & o_4 \\ 6 & 0 & 3 & 3 \end{pmatrix}$ and $U_{2,\alpha} = U_{2,\beta} = \begin{pmatrix} o_1 & o_2 & o_3 & o_4 \\ 3 & 9 & 3 & 0 \end{pmatrix}$.

**Definition 15.1.** A *static Bayesian game of incomplete information with knowledge of one's own payoffs* consists of the following elements:

- a set $I = \{1,...,n\}$ of players;
- for every Player $i \in I$, a set $S_i$ of strategies (as usual, we denote by $S$ the set of strategy profiles);
- for every Player $i \in I$, a set $T_i$ of possible types; we denote by $T = T_1 \times ... \times T_n$ the set of profiles of types and by $T_{-i} = T_1 \times ... \times T_{i-1} \times T_{i+1} \times ... \times T_n$ the set of profiles of types for the players other than $i$;
- for every Player $i$ and for every type $t_i \in T_i$ of Player $i$, a utility (or payoff) function $U_{i,t_i} : S \to \mathbb{R}$;
- for every Player $i$ and for every type $t_i \in T_i$ of Player $i$, a probability distribution $P_{i,t_i} : T_{-i} \to [0,1]$ representing the beliefs of type $t_i$ about the types of the other players.

The beliefs of all the types are said to be *Harsanyi consistent* if there exists a *common prior*, that is, a probability distribution $P : T \to [0,1]$ such that, for every Player $i$ and for every type $t_i \in T_i$ of Player $i$, $P_{i,t_i}$ coincides with the probability distribution obtained from $P$ by conditioning on the event $\{t_i\}$, that is (denoting a profile $t \in T$ by $(t_{-i}, t_i)$), $P_{i,t_i}(t_{-i}) = \dfrac{P\big((t_{-i}, t_i)\big)}{\sum\limits_{t'_{-i} \in T_{-i}} P\big((t'_{-i}, t_i)\big)}$.

Let us recast the situation illustrated in Figure 15.1 in the terminology of Definition 15.1. First of all, we have that $I = \{1,2\}$, $S_1 = \{T,B\}$, $S_2 = \{L,R\}$ (so that $S = \{(T,L), (T,R), (B,L), (B,R)\}$). Furthermore there are two types of Player 1 and only one type of Player 2: $T_1 = \{t_1^a, t_1^b\}$, $T_2 = \{t_2\}$, so that $T = \{(t_1^a, t_2), (t_1^b, t_2)\}$. The utility functions are given by:

$$U_{1,t_1^a} = \begin{pmatrix} (T,L) & (T,R) & (B,L) & (B,R) \\ 0 & 3 & 3 & 0 \end{pmatrix}, \quad U_{1,t_1^b} = \begin{pmatrix} (T,L) & (T,R) & (B,L) & (B,R) \\ 6 & 0 & 3 & 3 \end{pmatrix}$$

$$\text{and } U_{2,t_2} = \begin{pmatrix} (T,L) & (T,R) & (B,L) & (B,R) \\ 3 & 9 & 3 & 0 \end{pmatrix}.$$

The beliefs are given by $P_{1,t_1^a} = P_{1,t_1^b} = \begin{pmatrix} t_2 \\ 1 \end{pmatrix}$, $P_{2,t_2} = \begin{pmatrix} t_1^a & t_1^b \\ \frac{2}{3} & \frac{1}{3} \end{pmatrix}$ and the common prior is $P = \begin{pmatrix} (t_1^a, t_2) & (t_1^b, t_2) \\ \frac{2}{3} & \frac{1}{3} \end{pmatrix}$.

Before we explain in detail how to transform a "state-space" structure into a "type-space" structure and *vice versa*, we give one more example in which there is double-sided incomplete information. Consider the situation illustrated in Figure 15.2, which reproduces Figure 13.8. The true state is $\alpha$.

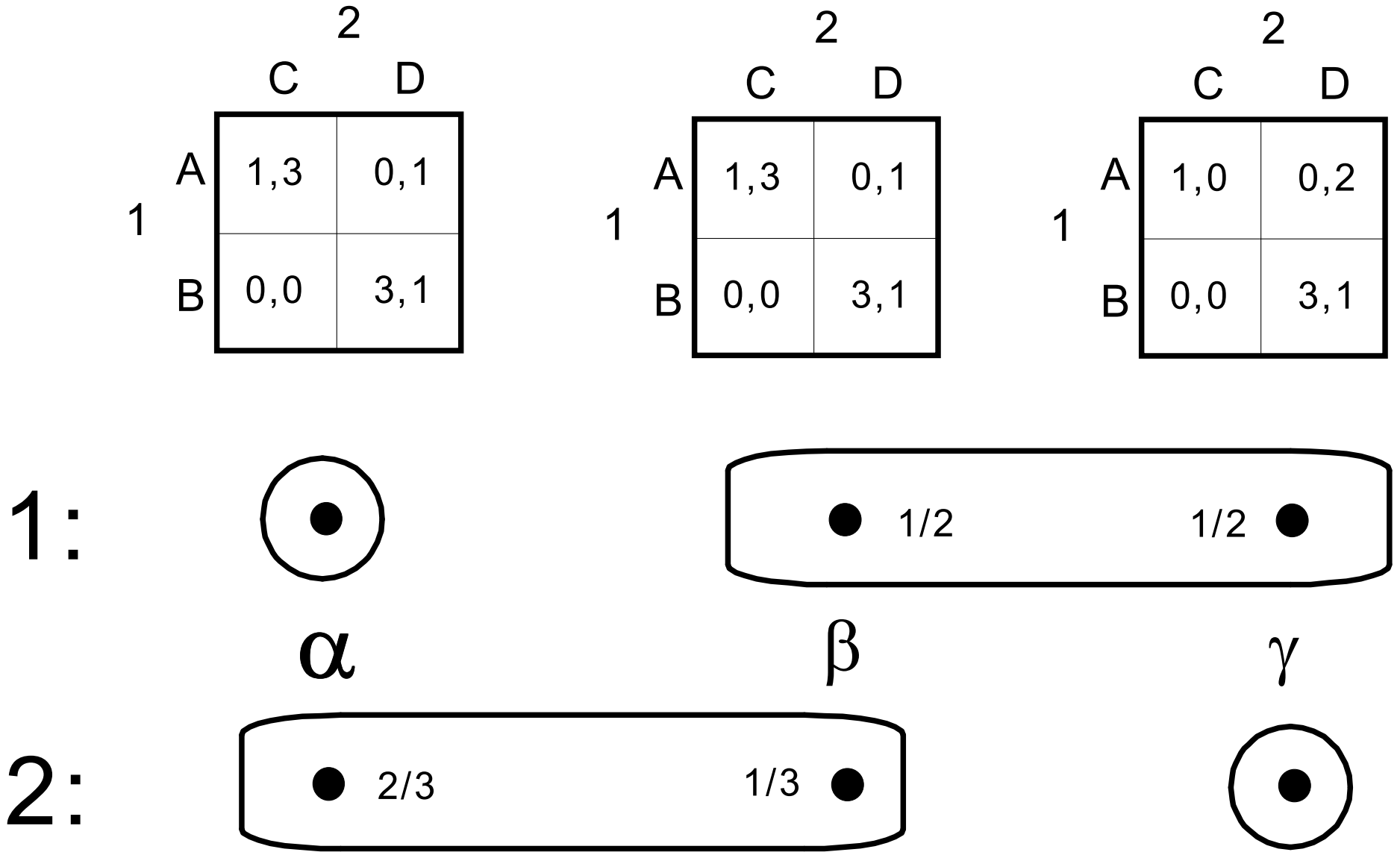


**Figure 15.2**

In this case we have that each player has two types: we can identify a type of a player with a cell of the player's information partition. Thus we have that $T_1 = \{t_1^a, t_1^b\}$, $T_2 = \{t_2^a, t_2^b\}$, $U_{1,t_1^a} = U_{1,t_1^b} = \begin{pmatrix} (A,C) & (A,D) & (B,C) & (B,D) \\ 1 & 0 & 0 & 3 \end{pmatrix}$,
$P_{1,t_1^a} = \begin{pmatrix} t_2^a \\ 1 \end{pmatrix}$, $P_{1,t_1^b} = \begin{pmatrix} t_2^a & t_2^b \\ \frac{1}{2} & \frac{1}{2} \end{pmatrix}$ $U_{2,t_2^a} = \begin{pmatrix} (A,C) & (A,D) & (B,C) & (B,D) \\ 3 & 1 & 0 & 1 \end{pmatrix}$,
$U_{2,t_2^b} = \begin{pmatrix} (A,C) & (A,D) & (B,C) & (B,D) \\ 0 & 2 & 0 & 1 \end{pmatrix}$, $P_{2,t_2^a} = \begin{pmatrix} t_1^a & t_1^b \\ \frac{2}{3} & \frac{1}{3} \end{pmatrix}$ and $P_{2,t_2^b} = \begin{pmatrix} t_1^b \\ 1 \end{pmatrix}$. The

common prior is given by $P=\begin{pmatrix}(t_1^a,t_2^a) & (t_1^a,t_2^b) & (t_1^b,t_2^a) & (t_1^b,t_2^b)\\ \frac{2}{4} & 0 & \frac{1}{4} & \frac{1}{4}\end{pmatrix}$. Thus the two types of Player 1 have the same utility function but different beliefs about the types of Player 2, while the two types of Player 2 differ both in terms of utility function and in terms of beliefs about the types of Player 1.

From the last example it should be clear how to transform a "state-space" model into a "type-space" model. First of all, for every Player *i*, create one type for every cell of Player *i*'s partition, making sure that different cells are associated with different types. In this way we have identified each state with a profile of types. Since there is a probability distribution over each information set of Player *i*, that probability distribution will give us a probability distribution for the associated type of Player *i* over some set of type-profiles for the other players. Finally, since our hypothesis (to be relaxed in the next section) is that each player knows her own utility function (that is, the utility function of a player does not vary from state to state within the same information set of that player), with each type $t_i$ of Player *i* is associated a unique utility function $U_{i,t_i}$.

Conversely, we can convert a "type-space" structure into a "state-space" structure as follows. Let the set of states be the set $T$ of profile of types. For every Player *i* and for every two states $t,t' \in T$, let $t$ and $t'$ belong to the same information set of Player $i$ (that is, to the same cell of Player $i$'s partition) if and only if Player $i$'s type is the same in $t$ and $t'$: $t' \in I_i(t)$ if and only if $t_i = t'_i$. The beliefs of each type of Player *i* then yield a probability distribution over the information set of Player *i* corresponding to that type. An example of this transformation is given in Exercise 15.2.

This is a good time to test your understanding of the concepts introduced in this section, by going through the exercises in Section 15.E.2 of Appendix 15.E at the end of this chapter.

## 15.3 The general case

As pointed out in Remark 13.1 (Chapter 13), it may very well be the case that a rational player does not know her own payoffs (for example, because she is uncertain about what outcomes might occur if she chooses a particular action). Uncertainty about a player's own payoffs is compatible with the player knowing her own preferences (that is, how she ranks the outcomes that she considers possible). Definition 15.1 is not general enough to encompass such possibilities. The following, more general, definition allows the utility function of a player to depend not only on the player's own type but also on the types of the other players. Definition 15.2 is identical to Definition 15.1 except for the fourth and fifth bulleted items, which have been highlighted in bold red characters, and the additional part in the last bulleted item, which is also highlighted in bold red characters.

**Definition 15.2.** A *static Bayesian game of incomplete information* consists of the following elements:

- a set $I = \{1,...,n\}$ of players;
- for every Player $i \in I$, a set $S_i$ of strategies (as usual we denote by $S$ the set of strategy profiles);
- for every Player $i \in I$, a set $T_i$ of possible types; we denote by $T = T_1 \times ... \times T_n$ the set of profiles of types and by $T_{-i} = T_1 \times ... \times T_{i-1} \times T_{i+1} \times ... \times T_n$ the set of profiles of types for the players other than $i$;
- **a set $Y \subseteq T$ of relevant profiles of types;**
- **for every Player $i$ and for every profile of types $t \in Y$, a utility (or payoff) function $U_{i,t} : S \to \mathbb{R}$;**
- for every Player $i$ and for every type $t_i \in T_i$ of Player $i$, a probability distribution $P_{i,t_i} : T_{-i} \to [0,1]$ representing the beliefs of type $t_i$ about the types of the other players **satisfying the restriction that if $P_{i,t_i}(t_{-i}) > 0$ then $(t_i, t_{-i}) \in Y$.**

As in the special case considered in the previous section, also in the general case one can transform a "state-space" structure into a "type-space" structure and *vice versa*.

Given a "state-space" structure, for every Player $i$ we identify the cells of Player $i$'s partition with the types of Player $i$ (one type for every information set). Since there is a probability distribution over each information set of Player $i$, that probability distribution will yield the probability distribution for the associated type of Player $i$ over some set of type-profiles for the other players. Finally, having identified each state with a profile of types (since each state belongs to one and only one information set of each player), we can assign to the corresponding type of Player $i$ the utility function of Player $i$ at that state. We shall illustrate this conversion using the "state-space" structure shown in Figure 15.3 below (taken from Exercise 13.2). The corresponding "type-space" structure (static Bayesian game of incomplete information) is the following: $I$ = {1,2} (letting Bill be Player 1 and Ann Player 2), $S_1$ = {$g$,$ng$}, $S_2$ = {$a$,$r$} (so that $S$ = {($g$,$a$), ($g$,$r$), ($ng$,$a$), ($ng$,$r$)}). Furthermore there are two types of Player 1 and only one type of Player 2: $T_1 = \{t_1^f, t_1^e\}$ ($f$ stands for 'friend' and $e$ for 'enemy'), $T_2 = \{t_2\}$, so that $T = \left\{(t_1^f, t_2), (t_1^e, t_2)\right\}$. In this case we have that $Y$ = $T$. The utility functions are given by:

$$U_{1,(t_1^f,t_2)} = U_{1,(t_1^e,t_2)} = \begin{pmatrix} (g,a) & (g,r) & (ng,a) & (ng,r) \\ 1 & -1 & 0 & 0 \end{pmatrix},$$

$$U_{2,(t_1^f,t_2)} = \begin{pmatrix} (g,a) & (g,r) & (ng,a) & (ng,r) \\ 1 & 0 & 0 & 0 \end{pmatrix}$$

$$U_{2,(t_1^e,t_2)} = \begin{pmatrix} (g,a) & (g,r) & (ng,a) & (ng,r) \\ -1 & 0 & 0 & 0 \end{pmatrix}.$$

The beliefs are given by $P_{1,t_1^f} = P_{1,t_1^e} = \begin{pmatrix} t_2 \\ 1 \end{pmatrix}$, $P_{2,t_2} = \begin{pmatrix} t_1^f & t_1^e \\ p & 1-p \end{pmatrix}$ and the common prior is $P = \begin{pmatrix} (t_1^f,t_2) & (t_1^e,t_2) \\ p & 1-p \end{pmatrix}$.

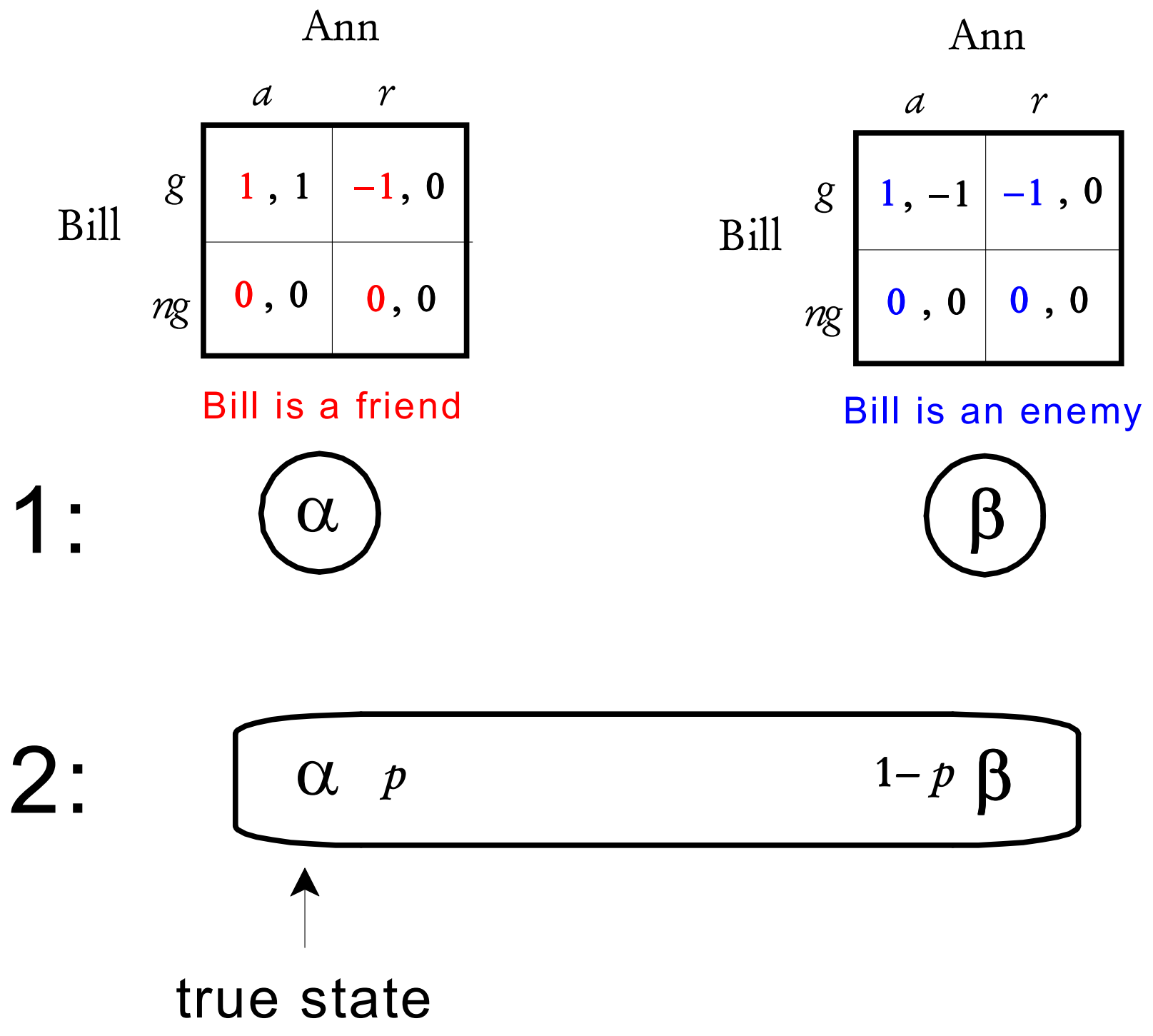


**Figure 15.3**

Conversely, we can convert a "type-space" structure into a "state-space" structure as follows. Let the set of states be the set $Y$ of relevant profiles of types. For every Player $i$ and for every two states $t, t' \in Y$, let $t$ and $t'$ belong to the same information set of Player $i$ (that is, to the same cell of Player $i$'s partition) if and only if Player $i$'s type is the same in $t$ and $t'$: $t' \in I_i(t)$ if and only if $t_i = t_i'$. The beliefs of each type of Player $i$ then yield a probability distribution over the information set of Player $i$ corresponding to that type. An example of this transformation is given in Exercise 15.3.

This is a good time to test your understanding of the concepts introduced in this section, by going through the exercises in Section 15.E.3 of Appendix 15.E at the end of this chapter.

# Appendix 15.E: Exercises

## 15.E.2. Exercises for Section 15.2: Types that know their own payoffs

The answers to the following exercises are in Appendix S.

**Exercise 15.1.** Transform the following situation of incomplete information into a "type-space" structure.

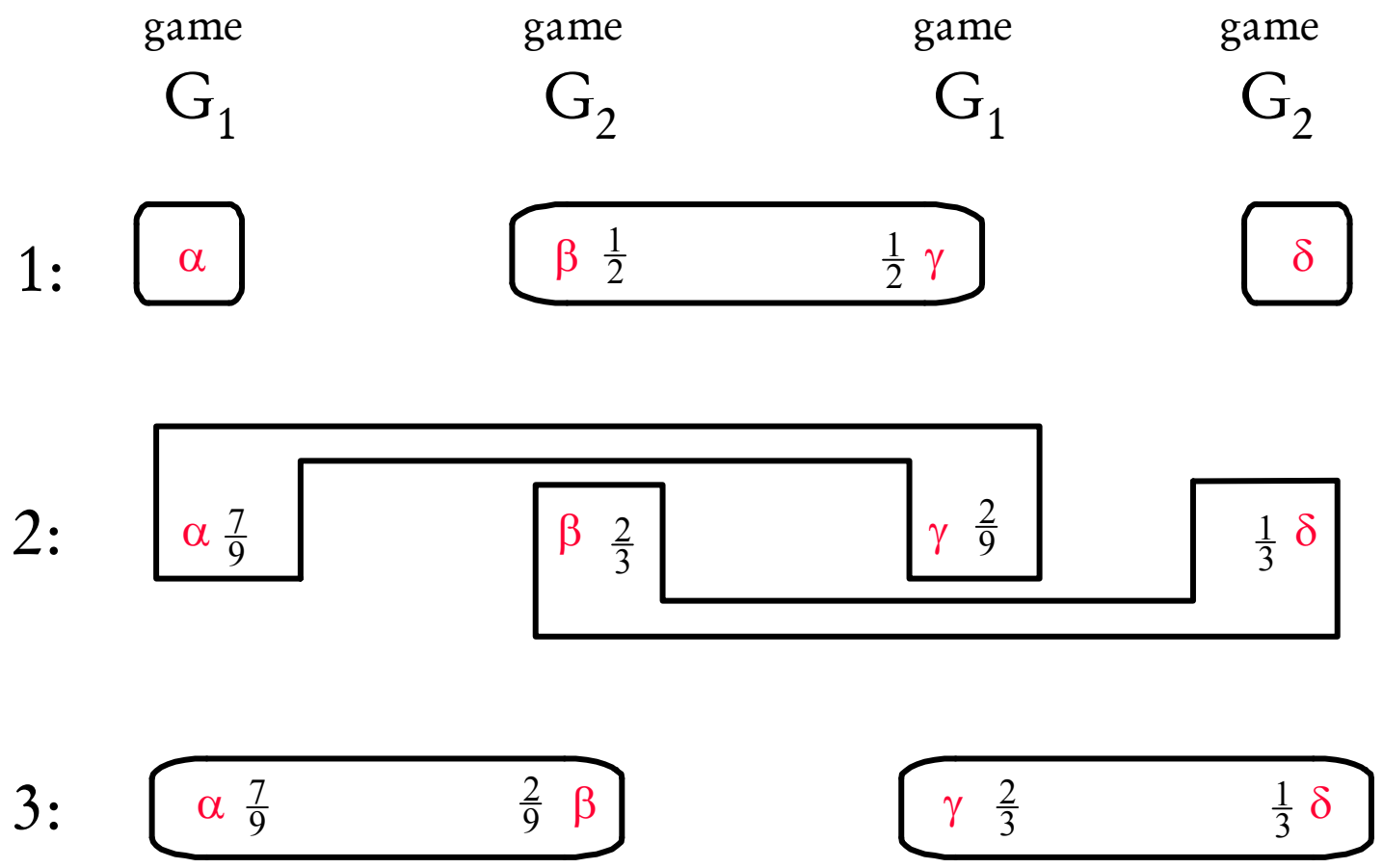


**GAME $G_1$**

Player 3 chooses E (Player 1, Player 2):

| Player 1 | C | | | D | | |
|---|---|---|---|---|---|---|
| A | 4 | 1 | 0 | 2 | 4 | 2 |
| B | 4 | 1 | 0 | 0 | 0 | 0 |

Player 3 chooses F (Player 1, Player 2):

| Player 1 | C | | | D | | |
|---|---|---|---|---|---|---|
| A | 2 | 0 | 2 | 0 | 1 | 0 |
| B | 0 | 0 | 0 | 0 | 2 | 0 |

**GAME $G_2$**

Player 3 chooses E (Player 1, Player 2):

| Player 1 | C | | | D | | |
|---|---|---|---|---|---|---|
| A | 4 | 4 | 0 | 2 | 1 | 2 |
| B | 4 | 4 | 0 | 0 | 1 | 0 |

Player 3 chooses F (Player 1, Player 2):

| Player 1 | C | | | D | | |
|---|---|---|---|---|---|---|
| A | 2 | 1 | 2 | 0 | 0 | 0 |
| B | 0 | 2 | 0 | 0 | 0 | 0 |

**Exercise 15.2.** Consider the following Bayesian game of incomplete information:

$$I = \{1,2\},\ S_1 = \{A,B\},\ S_2 = \{C,D\}\ \ T_1 = \{t_1^a, t_1^b\},\ T_2 = \{t_2^a, t_2^b, t_2^c\}$$

$$U_{1,t_1^a} = \begin{pmatrix} AC & AD & BC & BD \\ 4 & 1 & 0 & 2 \end{pmatrix},\quad U_{1,t_1^b} = \begin{pmatrix} AC & AD & BC & BD \\ 0 & 0 & 2 & 1 \end{pmatrix}$$

$$U_{2,t_2^a} = U_{2,t_2^b} = \begin{pmatrix} AC & AD & BC & BD \\ 2 & 1 & 2 & 3 \end{pmatrix},\quad U_{2,t_2^c} = \begin{pmatrix} AC & AD & BC & BD \\ 0 & 2 & 2 & 0 \end{pmatrix}$$

$$P_{1,t_1^a} = \begin{pmatrix} t_2^a & t_2^b & t_2^c \\ \frac{1}{4} & \frac{3}{4} & 0 \end{pmatrix},\quad P_{1,t_1^b} = \begin{pmatrix} t_2^a & t_2^b & t_2^c \\ \frac{4}{13} & \frac{3}{13} & \frac{6}{13} \end{pmatrix},\quad P_{2,t_2^a} = \begin{pmatrix} t_1^a & t_1^b \\ \frac{1}{3} & \frac{2}{3} \end{pmatrix},\quad P_{2,t_2^b} = \begin{pmatrix} t_1^a & t_1^b \\ \frac{2}{3} & \frac{1}{3} \end{pmatrix}$$

$$P_{2,t_2^c} = \begin{pmatrix} t_1^a & t_1^b \\ 0 & 1 \end{pmatrix}.$$

**(a)** Are the beliefs of the types consistent (that is, is there a common prior)?

**(b)** Transform this type-space structure into an interactive knowledge-belief structure.

## 15.E.3. Exercises for Section 15.3: The general case

The answers to the following exercises are in Appendix S at the end of this chapter.

**Exercise 15.3.** Consider the following Bayesian game of incomplete information.

$$I = \{1,2,3\},\quad S_1 = \{A,B\},\ S_2 = \{C,D\},\ S_3 = \{E,F\},\quad T_1 = \{t_1^a, t_1^b\},\quad T_2 = \{t_2^a, t_2^b\},$$
$$T_3 = \{t_3^a, t_3^b\},\ \ Y = \{(t_1^a, t_2^a, t_3^a), (t_1^b, t_2^b, t_3^a), (t_1^b, t_2^a, t_3^b), (t_1^b, t_2^b, t_3^b)\}$$

$$U_{1,(t_1^a,t_2^a,t_3^a)} = \begin{pmatrix} ACE & ADE & BCE & BDE & ACF & ADF & BCF & BDF \\ 2 & 2 & 3 & 0 & 0 & 0 & 1 & 1 \end{pmatrix},$$

$$U_{1,(t_1^b,t_2^b,t_3^a)} = U_{1,(t_1^b,t_2^b,t_3^b)} = \begin{pmatrix} ACE & ADE & BCE & BDE & ACF & ADF & BCF & BDF \\ 1 & 2 & 0 & 0 & 1 & 0 & 0 & 1 \end{pmatrix},$$

$$U_{1,(t_1^b,t_2^a,t_3^b)} = \begin{pmatrix} ACE & ADE & BCE & BDE & ACF & ADF & BCF & BDF \\ 0 & 0 & 2 & 1 & 2 & 0 & 0 & 1 \end{pmatrix}$$

$$U_{2,(t_1^a,t_2^a,t_3^a)} = \begin{pmatrix} ACE & ADE & BCE & BDE & ACF & ADF & BCF & BDF \\ 2 & 4 & 2 & 0 & 0 & 1 & 0 & 2 \end{pmatrix}$$

$$U_{2,(t_1^b,t_2^b,t_3^a)} = U_{2,(t_1^b,t_2^b,t_3^b)} = \begin{pmatrix} ACE & ADE & BCE & BDE & ACF & ADF & BCF & BDF \\ 4 & 1 & 4 & 1 & 1 & 0 & 2 & 0 \end{pmatrix}$$

$$U_{2,(t_1^a,t_2^a,t_3^b)} = \begin{pmatrix} ACE & ADE & BCE & BDE & ACF & ADF & BCF & BDF \\ 4 & 3 & 1 & 2 & 0 & 1 & 0 & 2 \end{pmatrix}$$

$$U_{3,(t_1^a,t_2^a,t_3^a)} = \begin{pmatrix} ACE & ADE & BCE & BDE & ACF & ADF & BCF & BDF \\ 0 & 2 & 0 & 0 & 2 & 0 & 0 & 0 \end{pmatrix}$$

$$U_{3,(t_1^b,t_2^b,t_3^a)} = U_{3,(t_1^b,t_2^b,t_3^b)} = \begin{pmatrix} ACE & ADE & BCE & BDE & ACF & ADF & BCF & BDF \\ 3 & 2 & 0 & 1 & 1 & 2 & 2 & 0 \end{pmatrix}$$

$$U_{3,(t_1^b,t_2^a,t_3^b)} = \begin{pmatrix} ACE & ADE & BCE & BDE & ACF & ADF & BCF & BDF \\ 1 & 2 & 0 & 1 & 2 & 0 & 0 & 1 \end{pmatrix}$$

$$P_{1,t_1^a} = \begin{pmatrix} (t_2^a, t_3^a) \\ 1 \end{pmatrix}, \quad P_{1,t_1^b} = \begin{pmatrix} (t_2^b, t_3^a) & (t_2^a, t_3^b) & (t_2^b, t_3^b) \\ \frac{2}{5} & \frac{2}{5} & \frac{1}{5} \end{pmatrix},$$

$$P_{2,t_2^a} = \begin{pmatrix} (t_1^a, t_3^a) & (t_1^b, t_3^b) \\ \frac{5}{7} & \frac{2}{7} \end{pmatrix}, \quad P_{2,t_2^b} = \begin{pmatrix} (t_1^b, t_3^a) & (t_1^b, t_3^b) \\ \frac{2}{3} & \frac{1}{3} \end{pmatrix}$$

$$P_{3,t_3^a} = \begin{pmatrix} (t_1^a, t_2^a) & (t_1^b, t_2^b) \\ \frac{5}{7} & \frac{2}{7} \end{pmatrix}, \quad P_{3,t_3^b} = \begin{pmatrix} (t_1^b, t_2^a) & (t_1^b, t_2^b) \\ \frac{2}{3} & \frac{1}{3} \end{pmatrix}$$

**(a)** Are the beliefs of the types consistent (that is, is there a common prior)?

**(b)** Transform this type-space structure into an interactive knowledge-belief structure.

◊◊◊◊◊◊◊◊◊◊◊◊

**Exercise 15.4: Challenging Question.**

Consider the following Bayesian game of incomplete information:

$$I=\{1,2\},\ S_1=\{T,B\},\ S_2=\{L,R\},\ T_1=\left\{t_1^A,t_1^B\right\},\ T_2=\left\{t_2^a,t_2^b\right\},$$

$$U_{1,(t_1^A,t_2^a)}=U_{1,(t_1^A,t_2^b)}=\begin{pmatrix}(T,L) & (T,R) & (B,L) & (B,R)\\ 6 & 0 & 3 & 3\end{pmatrix},$$

$$U_{1,(t_1^B,t_2^a)}=U_{1,(t_1^B,t_2^b)}=\begin{pmatrix}(T,L) & (T,R) & (B,L) & (B,R)\\ 0 & 3 & 3 & 0\end{pmatrix},$$

$$U_{2,(t_1^A,t_2^a)}=U_{2,(t_1^A,t_2^b)}=U_{2,(t_1^B,t_2^a)}=U_{2,(t_1^B,t_2^b)}=\begin{pmatrix}(T,L) & (T,R) & (B,L) & (B,R)\\ 3 & 9 & 3 & 0\end{pmatrix},$$

$$P_{1,t_1^A}=\begin{pmatrix}t_2^a & t_2^b\\ \frac{1}{2} & \frac{1}{2}\end{pmatrix},\ P_{1,t_1^B}=\begin{pmatrix}t_2^a & t_2^b\\ \frac{3}{5} & \frac{2}{5}\end{pmatrix},\ P_{2,t_2^a}=\begin{pmatrix}t_1^A & t_1^B\\ \frac{1}{4} & \frac{3}{4}\end{pmatrix}\text{ and } P_{2,t_2^b}=\begin{pmatrix}t_1^A & t_1^B\\ \frac{1}{3} & \frac{2}{3}\end{pmatrix}.$$

**(a)** Transform the situation described above into an interactive knowledge-belief structure. Assume that the true state is where Player 1 is of type $t_1^A$ and Player 2 is of type $t_2^a$.

**(b)** Apply the Harsanyi transformation to obtain an extensive-form game.

**(c)** Write the strategic-form game corresponding to the game of part (b).

**(d)** Find all the pure-strategy Bayesian Nash equilibria. Does any of these yield a Nash equilibrium in the true game being played?

**(e)** Of the pure-strategy Bayesian Nash equilibria select one where Player 1 is uncertain about Player 2's choice of action. For this equilibrium complete the structure of part (a) by turning it into a model (that is, associate with each state an action - not a strategy - for each player) and verify that at the true state there is common knowledge of rationality.

# Appendix 15.S: Solutions to exercises

**Exercise 15.1.** The structure is as follows (the elements are given as listed Definition 15.1). $I=\{1,2,3\}$, $S_1=\{A,B\}$, $S_2=\{C,D\}$, $S_3=\{E,F\}$, $T_1=\left\{t_1^a,t_1^b,t_1^c\right\}$, $T_2=\left\{t_2^a,t_2^b\right\}$, $T_3=\left\{t_3^a,t_3^b\right\}$,

$$U_{1,t_1^a}=U_{1,t_1^b}=U_{1,t_1^c}=\begin{pmatrix} ACE & ADE & BCE & BDE & ACF & ADF & BCF & BDF \\ 4 & 2 & 4 & 0 & 2 & 0 & 0 & 0 \end{pmatrix},$$

$$U_{2,t_2^a}=\begin{pmatrix} ACE & ADE & BCE & BDE & ACF & ADF & BCF & BDF \\ 1 & 4 & 1 & 0 & 0 & 1 & 0 & 2 \end{pmatrix}$$

$$U_{2,t_2^b}=\begin{pmatrix} ACE & ADE & BCE & BDE & ACF & ADF & BCF & BDF \\ 4 & 1 & 4 & 1 & 1 & 0 & 2 & 0 \end{pmatrix}$$

$$U_{3,t_3^a}=U_{3,t_3^b}=\begin{pmatrix} ACE & ADE & BCE & BDE & ACF & ADF & BCF & BDF \\ 0 & 2 & 0 & 0 & 2 & 0 & 0 & 0 \end{pmatrix}$$

$$P_{1,t_1^a}=\begin{pmatrix} (t_2^a,t_3^a) \\ 1 \end{pmatrix},\ P_{1,t_1^b}=\begin{pmatrix} (t_2^b,t_3^a) & (t_2^a,t_3^b) \\ \frac{1}{2} & \frac{1}{2} \end{pmatrix},\ P_{1,t_1^c}=\begin{pmatrix} (t_2^b,t_3^b) \\ 1 \end{pmatrix}$$

$$P_{2,t_2^a}=\begin{pmatrix} (t_1^a,t_3^a) & (t_1^b,t_3^b) \\ \frac{7}{9} & \frac{2}{9} \end{pmatrix},\quad P_{2,t_2^b}=\begin{pmatrix} (t_1^b,t_3^a) & (t_1^c,t_3^b) \\ \frac{2}{3} & \frac{1}{3} \end{pmatrix}$$

$$P_{3,t_3^a}=\begin{pmatrix} (t_1^a,t_2^a) & (t_1^b,t_2^b) \\ \frac{7}{9} & \frac{2}{9} \end{pmatrix},\quad P_{3,t_3^b}=\begin{pmatrix} (t_1^b,t_2^a) & (t_1^c,t_2^b) \\ \frac{2}{3} & \frac{1}{3} \end{pmatrix}$$

The common prior is given by

$$P=\begin{pmatrix} t_1^at_2^at_3^a & t_1^at_2^at_3^b & t_1^at_2^bt_3^a & t_1^at_2^bt_3^b & t_1^bt_2^at_3^a & t_1^bt_2^at_3^b & t_1^bt_2^bt_3^a & t_1^bt_2^bt_3^b & t_1^ct_2^at_3^a & t_1^ct_2^at_3^b & t_1^ct_2^bt_3^a & t_1^ct_2^bt_3^b \\ \frac{7}{12} & 0 & 0 & 0 & 0 & \frac{2}{12} & \frac{2}{12} & 0 & 0 & 0 & 0 & \frac{1}{12} \end{pmatrix}$$

**Exercise 15.2.** **(a)** Yes, the following is a common prior:

$$P=\begin{pmatrix} (t_1^a,t_2^a) & (t_1^a,t_2^b) & (t_1^a,t_2^c) & (t_1^b,t_2^a) & (t_1^b,t_2^b) & (t_1^b,t_2^c) \\ \frac{2}{21} & \frac{6}{21} & 0 & \frac{4}{21} & \frac{3}{21} & \frac{6}{21} \end{pmatrix}$$

**(b)** The knowledge-belief structure is as follows: the set of states is $\Omega=\{\alpha,\beta,\gamma,\delta,\varepsilon\}$, where $\alpha=(t_1^a,t_2^a)$, $\beta=(t_1^a,t_2^b)$, $\gamma=(t_1^b,t_2^a)$, $\delta=(t_1^b,t_2^b)$, $\varepsilon=(t_1^b,t_2^c)$.

| | | 2: C | 2: D |
|---|---|---|---|
| 1 | A | 4 , 2 | 1 , 1 |
| | B | 0 , 2 | 2 , 3 |

| | | 2: C | 2: D |
|---|---|---|---|
| 1 | A | 4 , 2 | 1 , 1 |
| | B | 0 , 2 | 2 , 3 |

| | | 2: C | 2: D |
|---|---|---|---|
| 1 | A | 0 , 2 | 0 , 1 |
| | B | 2 , 2 | 1 , 3 |

| | | 2: C | 2: D |
|---|---|---|---|
| 1 | A | 0 , 2 | 0 , 1 |
| | B | 2 , 2 | 1 , 3 |

| | | 2: C | 2: D |
|---|---|---|---|
| 1 | A | 0 , 0 | 0 , 2 |
| | B | 2 , 2 | 1 , 0 |

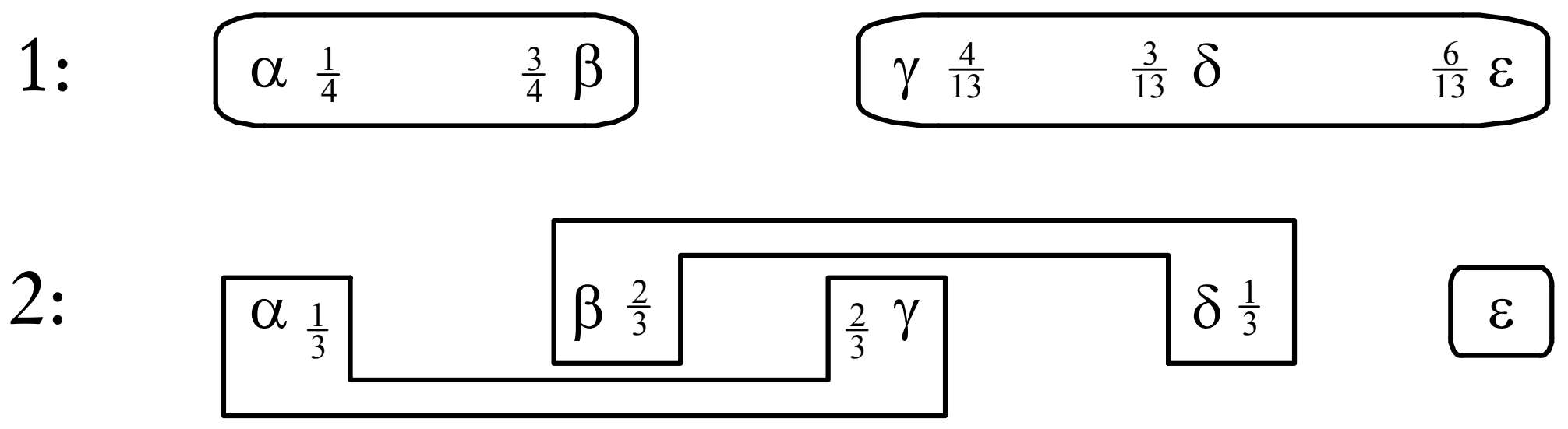


**Exercise 15.3.** (a) Yes: the following is a common prior:

$$P = \begin{pmatrix} t_1^a t_2^a t_3^a & t_1^b t_2^a t_3^b & t_1^a t_2^b t_3^a & t_1^b t_2^b t_3^b \\ \frac{1}{2} & \frac{1}{5} & \frac{1}{5} & \frac{1}{10} \end{pmatrix} \text{ and } P(t) = 0 \text{ for every other } t \in T.$$

(b) First define the following games:

**GAME $G_1$**

Player 3 chooses E

| Player 1 \ Player 2 | C | D |
|---|---|---|
| A | 2, 2, 0 | 2, 4, 2 |
| B | 3, 2, 0 | 0, 0, 0 |

Player 3 chooses F

| Player 1 \ Player 2 | C | D |
|---|---|---|
| A | 0, 0, 2 | 0, 1, 0 |
| B | 1, 0, 0 | 1, 2, 0 |

**GAME $G_2$**

Player 3 chooses E

| Player 1 \ Player 2 | C | D |
|---|---|---|
| A | 1, 4, 3 | 2, 1, 2 |
| B | 0, 4, 0 | 0, 1, 1 |

Player 3 chooses F

| Player 1 \ Player 2 | C | D |
|---|---|---|
| A | 1, 1, 1 | 0, 0, 2 |
| B | 0, 2, 2 | 1, 0, 0 |

**GAME $G_3$**

Player 3 chooses E

| Player 1 \ Player 2 | C | D |
|---|---|---|
| A | 0, 4, 1 | 0, 3, 2 |
| B | 2, 1, 0 | 1, 2, 1 |

Player 3 chooses F

| Player 1 \ Player 2 | C | D |
|---|---|---|
| A | 2, 0, 2 | 0, 1, 0 |
| B | 0, 0, 0 | 1, 2, 1 |

Then, letting $\alpha = (t_1^a, t_2^a, t_3^a)$, $\beta = (t_1^b, t_2^b, t_3^a)$, $\gamma = (t_1^b, t_2^a, t_3^b)$, $\delta = (t_1^b, t_2^b, t_3^b)$, the interactive knowledge-belief structure is as follows:

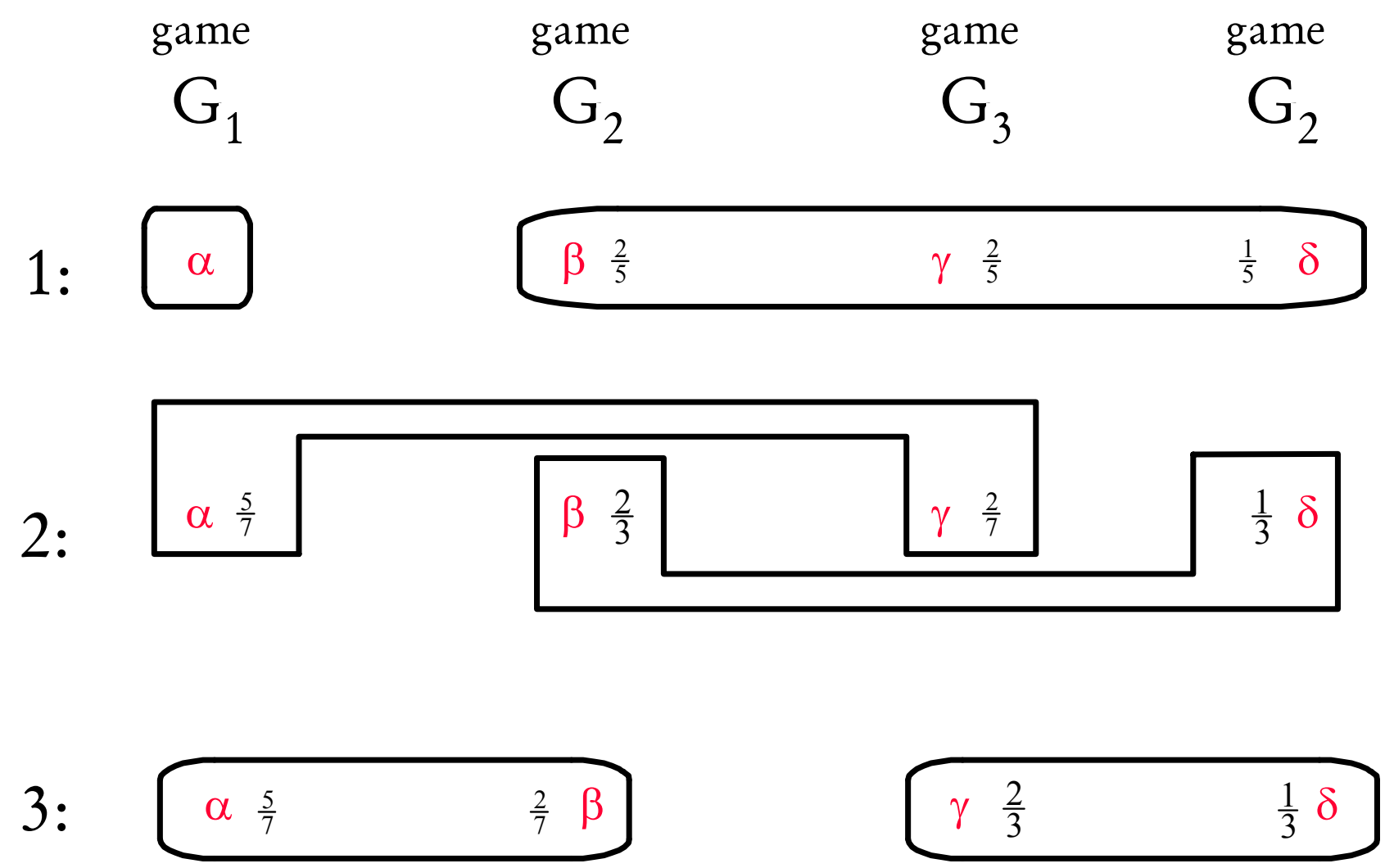

**Exercise 15.4** (Challenging Question) **(a)** The interactive knowledge-belief structure is as follows, where $\alpha = (t_1^A, t_2^a)$, $\beta = (t_1^A, t_2^b)$, $\gamma = (t_1^B, t_2^a)$, $\delta = (t_1^B, t_2^b)$:

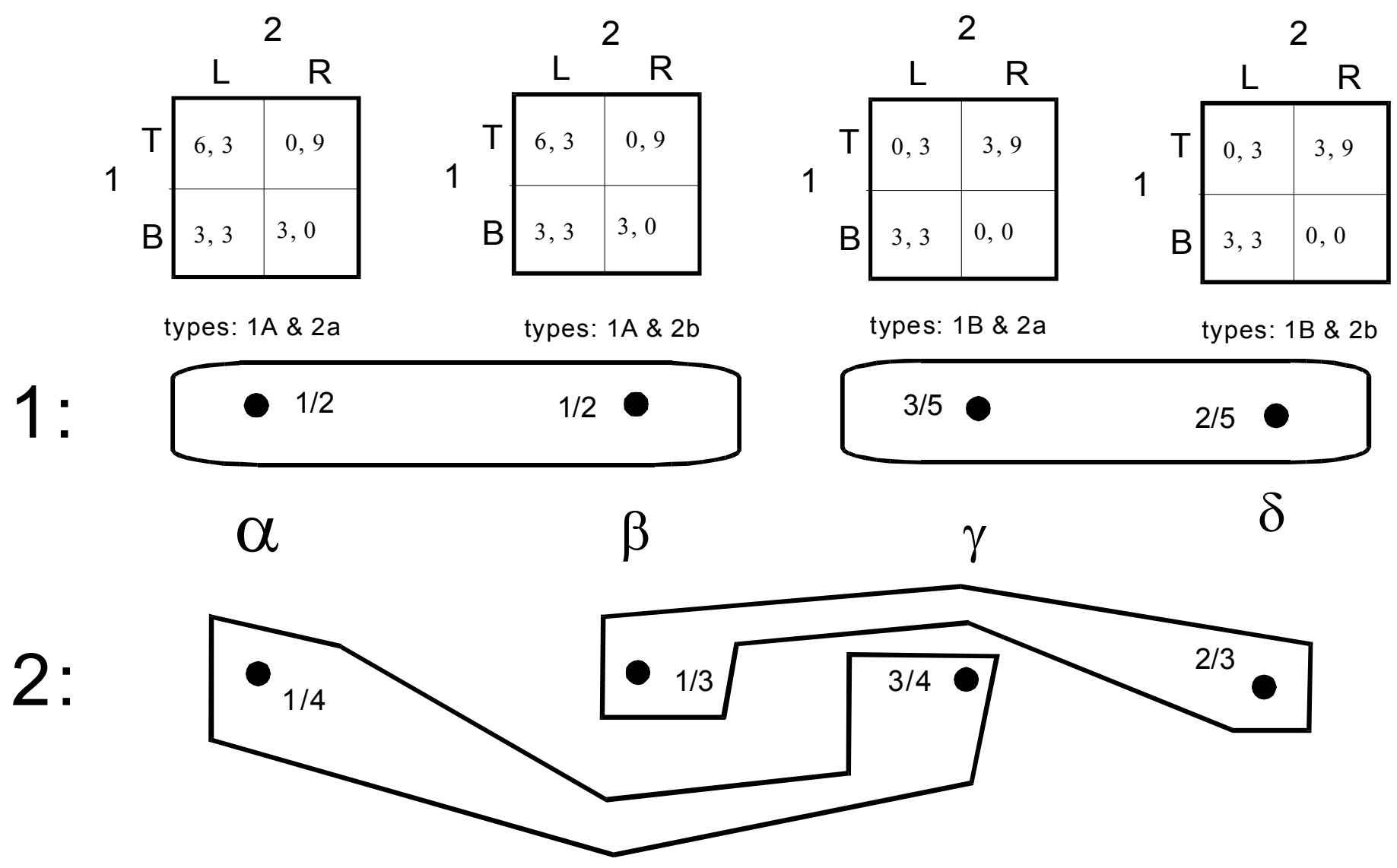

**(b)** First of all, there is a common prior: $\begin{pmatrix} \alpha & \beta & \gamma & \delta \\ \frac{1}{7} & \frac{1}{7} & \frac{3}{7} & \frac{2}{7} \end{pmatrix}$. The extensive-form game is as follows:

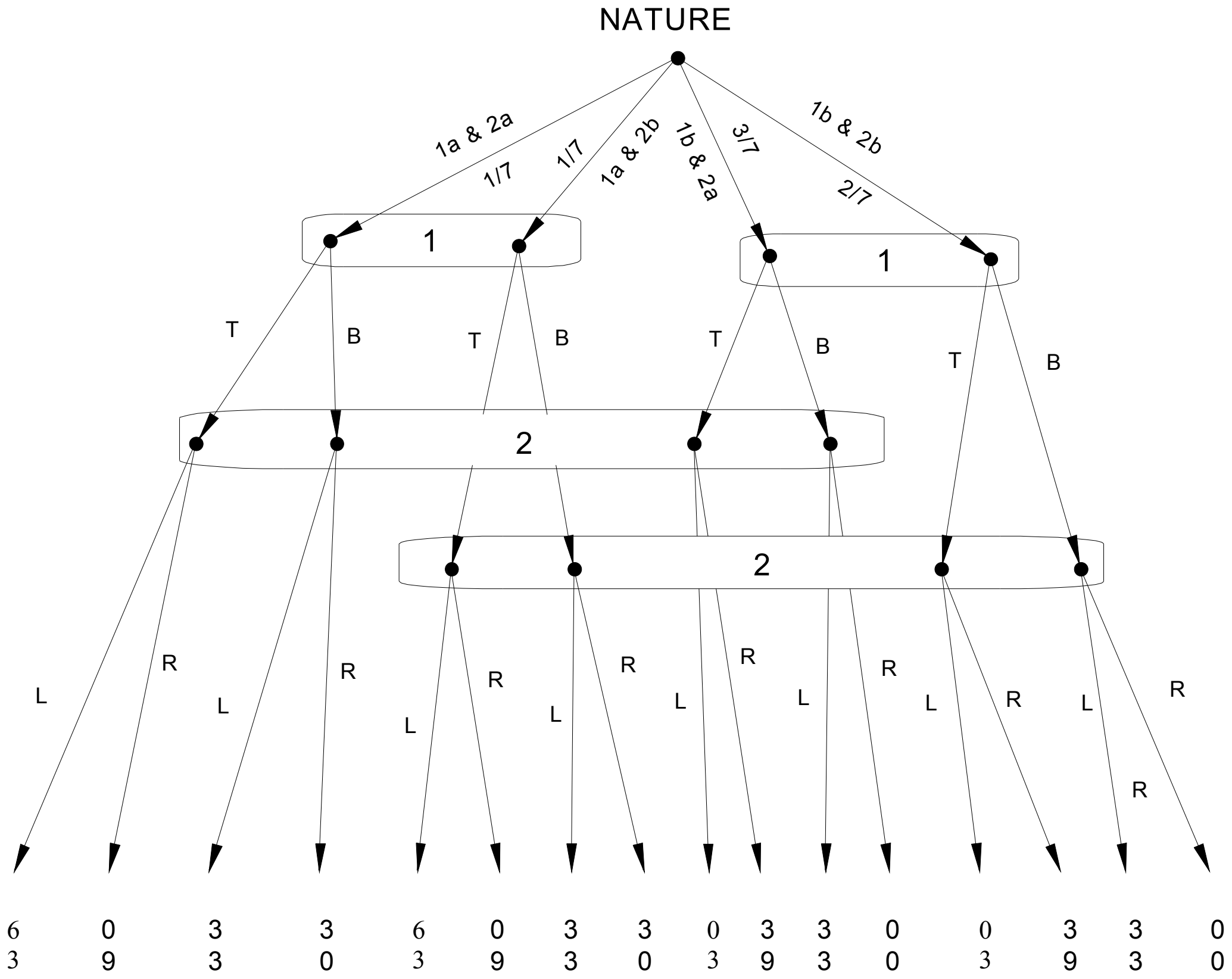

**(c)** The strategy profiles give rise to the following lotteries:

| | Player 2: L if type 2a<br>L if type 2b | L if type 2a<br>R if type 2b | R if type 2a<br>L if type 2b | R if type 2a<br>R if type 2b |
|---|---|---|---|---|
| **Player 1**: T if type 1A<br>T if type 1B | 1/7 1/7 3/7 2/7<br>*TL TL TL TL*<br>6 6 0 0<br>3 3 3 3 | 1/7 1/7 3/7 2/7<br>*TL TR TL TR*<br>6 0 0 3<br>3 9 3 9 | 1/7 1/7 3/7 2/7<br>*TR TL TR TL*<br>0 6 3 0<br>9 3 9 3 | 1/7 1/7 3/7 2/7<br>*TR TR TR TR*<br>0 0 3 3<br>9 9 9 9 |
| T if type 1A<br>B if type 1B | 1/7 1/7 3/7 2/7<br>*TL TL BL BL*<br>6 6 3 3<br>3 3 3 3 | 1/7 1/7 3/7 2/7<br>*TL TR BL BR*<br>6 0 3 0<br>3 9 3 0 | 1/7 1/7 3/7 2/7<br>*TR TL BR BL*<br>0 6 0 3<br>9 3 0 3 | 1/7 1/7 3/7 2/7<br>*TR TR BR BR*<br>0 0 0 0<br>9 9 0 0 |
| B if type 1A<br>T if type 1B | 1/7 1/7 3/7 2/7<br>*BL BL TL TL*<br>3 3 0 0<br>3 3 3 3 | 1/7 1/7 3/7 2/7<br>*BL BR TL TR*<br>3 3 0 3<br>3 0 3 9 | 1/7 1/7 3/7 2/7<br>*BR BL TR TL*<br>3 3 3 0<br>0 3 9 3 | 1/7 1/7 3/7 2/7<br>*BR BR TR TR*<br>3 3 3 3<br>0 0 9 9 |
| B if type 1A<br>B if type 1B | 1/7 1/7 3/7 2/7<br>*BL BL BL BL*<br>3 3 3 3<br>3 3 3 3 | 1/7 1/7 3/7 2/7<br>*BL BR BL BR*<br>3 3 3 0<br>3 0 3 0 | 1/7 1/7 3/7 2/7<br>*BR BL BR BL*<br>3 3 0 3<br>0 3 0 3 | 1/7 1/7 3/7 2/7<br>*BR BR BR BR*<br>3 3 0 0<br>0 0 0 0 |

Thus the strategic form is as follows:

| | Player 2: L if type 2a<br>L if type 2b | L if type 2a<br>R if type 2b | R if type 2a<br>L if type 2b | R if type 2a<br>R if type 2b |
|---|---|---|---|---|
| **Player 1**: T if type 1A<br>T if type 1B | $\frac{12}{7}, 3$ | $\frac{12}{7}, \frac{39}{7}$ | $\frac{15}{7}, \frac{45}{7}$ | $\frac{15}{7}, 9$ |
| T if type 1A<br>B if type 1B | $\boxed{\frac{27}{7}, 3}$ | $\boxed{\frac{15}{7}, 3}$ | $\frac{12}{7}, \frac{18}{7}$ | $0, \frac{18}{7}$ |
| B if type 1A<br>T if type 1B | $\frac{6}{7}, 3$ | $\frac{12}{7}, \frac{30}{7}$ | $\frac{15}{7}, \frac{36}{7}$ | $\boxed{3, \frac{45}{7}}$ |
| B if type 1A<br>B if type 1B | $3, 3$ | $\frac{15}{7}, \frac{12}{7}$ | $\frac{12}{7}, \frac{9}{7}$ | $\frac{6}{7}, 0$ |

**(d)** The pure-strategy Nash equilibria are: (*TB*,*LL*), (*TB*,*LR*) and (*BT*,*RR*). None of them yields a Nash equilibrium of the game associated with state $\alpha$, since neither (*T*,*L*) nor (*B*,*R*) are Nash equilibria of that game.

**(e)** The only pure-strategy equilibrium where Player 1 is uncertain of Player 2's choice is (*TB*,*LR*). The corresponding model is:

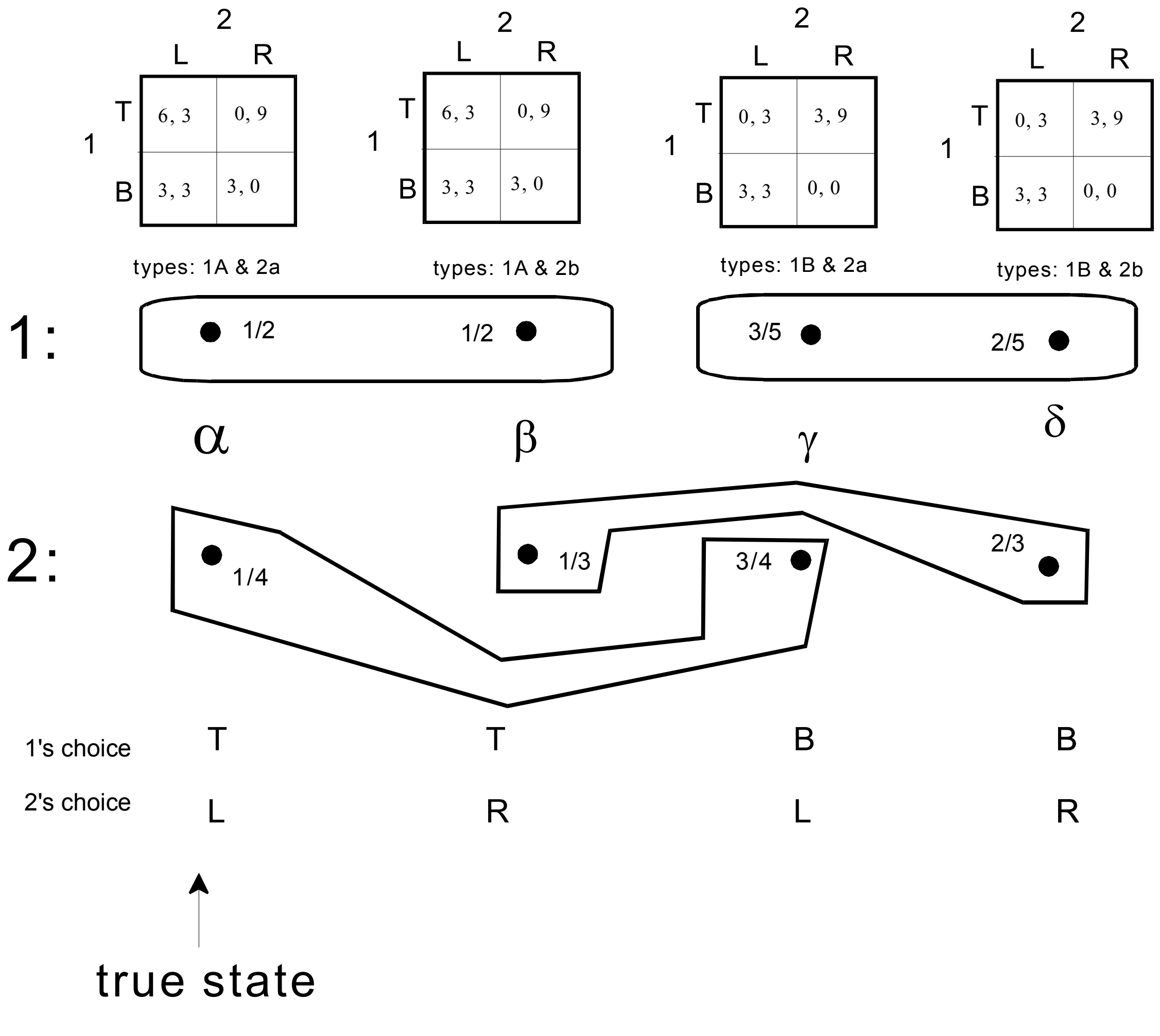

There is common knowledge of rationality because at every state both players are rational. In fact, at states $\alpha$ and $\beta$ Player 1 has an expected payoff of 3 from both *T* and *B* (thus *T* is a best reply) and at states $\gamma$ and $\delta$ Player 1 has an expected payoff of $\frac{9}{5}$ from *B* and $\frac{6}{5}$ from *T* (thus *B* is a best reply). At states $\alpha$ and $\gamma$ Player 2 has an expected payoff of 3 from *L* and $\frac{9}{4}$ from *R* (thus *L* is a best reply) and at states $\beta$ and $\delta$ Player 2 has an expected payoff of 3 from both *L* and *R* (thus *R* is a best reply).